\documentclass[a4,10pt,english]{smfbook}
\usepackage{amsmath,amscd,amssymb}
\makeindex

\newcommand{\nbiga}{\mathcal{A}}
\newcommand{\nbigb}{\mathcal{B}}
\newcommand{\nbigc}{\mathcal{C}}
\newcommand{\nbigd}{\mathcal{D}}
\newcommand{\nbige}{\mathcal{E}}
\newcommand{\nbigf}{\mathcal{F}}
\newcommand{\nbigg}{\mathcal{G}}
\newcommand{\nbigh}{\mathcal{H}}
\newcommand{\nbigi}{\mathcal{I}}
\newcommand{\nbigj}{\mathcal{J}}
\newcommand{\nbigk}{\mathcal{K}}
\newcommand{\nbigl}{\mathcal{L}}
\newcommand{\nbigm}{\mathcal{M}}
\newcommand{\nbign}{\mathcal{N}}
\newcommand{\nbigo}{\mathcal{O}}
\newcommand{\nbigp}{\mathcal{P}}
\newcommand{\nbigq}{\mathcal{Q}}
\newcommand{\nbigr}{\mathcal{R}}
\newcommand{\nbigs}{\mathcal{S}}
\newcommand{\nbigt}{\mathcal{T}}
\newcommand{\nbigu}{\mathcal{U}}
\newcommand{\nbigv}{\mathcal{V}}
\newcommand{\nbigw}{\mathcal{W}}

\newcommand{\nbigy}{\mathcal{Y}}

\newcommand{\proj}{\mathbb{P}}
\newcommand{\seisuu}{{\mathbb Z}}
\newcommand{\rnum}{{\mathbb Q}}

\newcommand{\cnum}{{\mathbb C}}
\newcommand{\real}{{\mathbb R}}
\newcommand{\hyperh}{\mathbb{H}}

\newcommand{\BB}{\mathbb{B}}

\newcommand{\Tbb}{\mathbb{T}}
\newcommand{\Abb}{\mathbb{A}}

\newcommand{\gbiga}{\mathfrak A}

\newcommand{\gbigd}{\mathfrak D}

\newcommand{\gbigf}{\mathfrak F}

\newcommand{\gbigi}{\mathfrak I}

\newcommand{\gbigk}{\mathfrak K}
\newcommand{\gbigl}{\mathfrak L}
\newcommand{\gbigm}{\mathfrak M}
\newcommand{\gbign}{\mathfrak N}
\newcommand{\gbigo}{\mathfrak O}
\newcommand{\gbigp}{\mathfrak P}
\newcommand{\gbigq}{\mathfrak Q}

\newcommand{\gbigt}{\mathfrak T}
\newcommand{\gbigu}{\mathfrak U}
\newcommand{\gbigv}{\mathfrak V}
\newcommand{\gbigw}{\mathfrak W}

\newcommand{\gminia}{\mathfrak a}
\newcommand{\gminib}{\mathfrak b}
\newcommand{\gminic}{\mathfrak c}

\newcommand{\gminim}{\mathfrak m}

\newcommand{\gminio}{\mathfrak o}

\newcommand{\gminis}{\mathfrak s}

\newcommand{\vece}{{\boldsymbol e}}

\newcommand{\vecv}{{\boldsymbol v}}

\newcommand{\vecc}{{\boldsymbol c}}

\newcommand{\vecI}{{\boldsymbol I}}

\newcommand{\vecF}{{\boldsymbol F}}

\newcommand{\vecA}{{\boldsymbol A}}
\newcommand{\vecB}{{\boldsymbol B}}
\newcommand{\vecS}{{\boldsymbol S}}

\newcommand{\llarr}{\longleftarrow}

\newcommand{\lrarr}{\longrightarrow}

\newcommand{\pf}{{\bf Proof}\hspace{.1in}}

\def\Hom{\mathop{\rm Hom}\nolimits}

\def\Cok{\mathop{\rm Cok}\nolimits}

\def\Image{\mathop{\rm Im}\nolimits}

\def\Re{\mathop{\rm Re}\nolimits}

\def\Gr{\mathop{\rm Gr}\nolimits}

\def\Tot{\mathop{\rm Tot}\nolimits}

\def\rank{\mathop{\rm rank}\nolimits}

\def\Ker{\mathop{\rm Ker}\nolimits}

\def\modulo{\mathop{\rm modulo}\nolimits}

\def\Gr{\mathop{\rm Gr}\nolimits}

\def\ord{\mathop{\rm ord}\nolimits}

\def\can{\mathop{\rm can}\nolimits}
\def\var{\mathop{\rm var}\nolimits}
\def\id{\mathop{\rm id}\nolimits}

\def\gcd{\mathop{\rm g.c.d.}\nolimits}

\def\Mer{\mathop{\rm Mer}\nolimits}
\def\RH{\mathop{\rm RH}\nolimits}

\newcommand{\del}{\partial}
\newcommand{\delbar}{\overline{\del}}

\newcommand{\nhom}{{\mathcal Hom}}

\newcommand{\vecPhi}{\boldsymbol\Phi}

\newcommand{\barz}{\overline{z}}
\newcommand{\zbar}{\barz}

\def\irr{\mathop{\rm irr}\nolimits}
\def\reg{\mathop{\rm reg}\nolimits}

\newcommand{\tildepsi}{\widetilde{\psi}}
\newcommand{\psitilde}{\tildepsi}

\newcommand{\lamda}{\lambda}

\newcommand{\distribution}{\gbigd\gminib}

\newcommand{\closedopen}[2]{[#1,#2[}
\newcommand{\openclosed}[2]{]#1,#2]}
\newcommand{\openopen}[2]{]#1,#2[}
\newcommand{\closedclosed}[2]{[#1,#2]}

\newcommand{\vhat}{\widehat{v}}

\newcommand{\rhotilde}{\widetilde{\rho}}

\newcommand{\Vhat}{\widehat{V}}
\newcommand{\nablahat}{\widehat{\nabla}}

\newcommand{\kappatilde}{\widetilde{\kappa}}

\newcommand{\nbigetilde}{\widetilde{\nbige}}

\newcommand{\ptilde}{\widetilde{p}}

\newcommand{\vtilde}{\widetilde{v}}

\newcommand{\ftilde}{\widetilde{f}}
\newcommand{\utilde}{\widetilde{u}}
\newcommand{\Ftilde}{\widetilde{F}}

\newcommand{\atilde}{\widetilde{a}}
\newcommand{\stilde}{\widetilde{s}}

\newcommand{\nbiglhat}{\widehat{\nbigl}}

\newcommand{\nbigltilde}{\widetilde{\nbigl}}

\newcommand{\Ohat}{\widehat{O}}

\newcommand{\Phihat}{\widehat{\Phi}}

\newcommand{\nbigvhat}{\widehat{\nbigv}}

\newcommand{\Utilde}{\widetilde{U}}
\newcommand{\Dtilde}{\widetilde{D}}
\newcommand{\Xtilde}{\widetilde{X}}

\newcommand{\nbigrtilde}{\widetilde{\nbigr}}

\newcommand{\Ltilde}{\widetilde{L}}

\def\ord{\mathop{\rm ord}\nolimits}

\def\Gal{\mathop{\rm Gal}\nolimits}

\def\moderate{\mathop{\rm mod}\nolimits}

\def\Hol{\mathop{\rm Hol}\nolimits}

\def\DR{\mathop{\rm DR}\nolimits}

\def\Loc{\mathop{\rm Loc}\nolimits}

\def\rd{\mathop{\rm rd}\nolimits}
\def\mg{\mathop{\rm mg}\nolimits}

\def\St{\mathop{\rm St}\nolimits}

\def\DM{\mathop{\rm DM}\nolimits}
\def\add{\mathop{\rm add}\nolimits}
\def\fin{\mathop{\rm fin}\nolimits}

\newcommand{\Dbar}{\overline{D}}

\newcommand{\Ibar}{\overline{I}}

\newcommand{\vecj}{{\boldsymbol j}}

\newcommand{\nbigutilde}{\widetilde{\nbigu}}
\newcommand{\nbigstilde}{\widetilde{\nbigs}}

\newcommand{\nbigmhat}{\widehat{\nbigm}}

\newcommand{\Ztilde}{\widetilde{Z}}
\newcommand{\nbigttilde}{\widetilde{\nbigt}}

\newcommand{\Ctilde}{\widetilde{C}}
\newcommand{\nbigntilde}{\widetilde{\nbign}}

\newcommand{\tautilde}{\widetilde{\tau}}

\newcommand{\Deltatilde}{\widetilde{\Delta}}
\newcommand{\nbigftilde}{\widetilde{\nbigf}}

\newcommand{\nbigitilde}{\widetilde{\nbigi}}

\newcommand{\vecctilde}{\widetilde{\vecc}}
\newcommand{\vecnbigi}{{\boldsymbol \nbigi}}

\newcommand{\DDD}{\boldsymbol D}
\newcommand{\jtilde}{\widetilde{j}}

\newcommand{\ctilde}{\widetilde{c}}

\newcommand{\Ktilde}{\widetilde{K}}
\newcommand{\veckappa}{{\boldsymbol \kappa}}

\newcommand{\nutilde}{\widetilde{\nu}}

\newcommand{\gminiatilde}{\widetilde{\gminia}}
\newcommand{\gminibtilde}{\widetilde{\gminib}}
\newcommand{\cnumtilde}{\widetilde{\cnum}}

\newcommand{\phitilde}{\widetilde{\phi}}

\newcommand{\Fourier}{\mathfrak{Four}}

\newcommand{\rhohat}{\widehat{\rho}}

\newcommand{\Gammatilde}{\widetilde{\Gamma}}

\newcommand{\gminiahat}{\widehat{\gminia}}

\newcommand{\nbigctilde}{\widetilde{\nbigc}}
\newcommand{\iotatilde}{\widetilde{\iota}}

\newcommand{\Jbar}{\overline{J}}
\newcommand{\vecnbigf}{\boldsymbol\nbigf}
\newcommand{\vecJ}{\boldsymbol J}
\newcommand{\Sh}{\boldsymbol{Sh}}
\newcommand{\Shcat}{\mathfrak{Sh}}

\newcommand{\vecleq}{\boldsymbol\leq}
\newcommand{\Jtilde}{\widetilde{J}}
\newcommand{\Upsilontilde}{\widetilde{\Upsilon}}
\newcommand{\Upsiloncheck}{\check{\Upsilon}}

\newcommand{\inftyhat}{\widehat{\infty}}

\newcommand{\vecnbigp}{\boldsymbol{\nbigp}}
\newcommand{\vecnbigq}{\boldsymbol{\nbigq}}
\newcommand{\vecnbigr}{\boldsymbol{\nbigr}}
\newcommand{\vecK}{\boldsymbol{K}}

\newcommand{\vecPsi}{\boldsymbol{\Psi}}

\newcommand{\projtilde}{\widetilde{\proj}}
\newcommand{\realbar}{\overline{\real}}

\newcommand{\nbiglcheck}{\check{\nbigl}}

\newcommand{\Locst}{\mathsf{Loc^{St}}}
\newcommand{\Shsf}{\mathsf{Sh}}
\newcommand{\Dsf}{\mathsf{D}}
\newcommand{\Csf}{\mathsf{C}}

\newcommand{\gbigutilde}{\widetilde{\gbigu}}
\newcommand{\vecJbar}{\overline{\vecJ}}
\newcommand{\vecL}{\boldsymbol{L}}
\newcommand{\vecnbigftilde}{\widetilde{\vecnbigf}}

\newcommand{\btilde}{\widetilde{b}}
\newcommand{\Cr}{{\rm Cr}}
\newcommand{\alphahat}{\widehat{\alpha}}
\newcommand{\nbigchat}{\widehat{\nbigc}}

\newcommand{\barast}{\underline{\ast}}
\newcommand{\barshriek}{\underline{!}}

\newcommand{\varrhobar}{\overline{\varrho}}
\newcommand{\varpibar}{\overline{\varpi}}
\newcommand{\realc}{\real{\textrm-}c}
\newcommand{\nbigibar}{\overline{\nbigi}}

\newcommand{\ttS}{{\tt S}}
\newcommand{\vecvarphi}{\boldsymbol\varphi}
\newcommand{\Sto}{\mathbb{S}\textrm{to}}
\newcommand{\varpitilde}{\widetilde{\varpi}}
\newcommand{\Jhat}{\hat{J}}
\newcommand{\Mhat}{\widehat{M}}

\newcommand{\LS}{\vecL\vecS}
\newcommand{\cantilde}{\widetilde{\can}}
\newcommand{\vartilde}{\widetilde{\var}}
\newcommand{\Kbb}{\mathbb{K}}
\newcommand{\nbigotilde}{\widetilde{\nbigo}}

\newtheorem{thm}{Theorem}[section]
\newtheorem{cor}[thm]{Corollary}

\newtheorem{rem}[thm]{Remark}
\newtheorem{lem}[thm]{Lemma}
\newtheorem{prop}[thm]{Proposition}
\newtheorem{df}[thm]{Definition}

\newtheorem{condition}[thm]{Condition}

\newtheorem{notation}[thm]{Notation}

\author{Takuro Mochizuki}
\address{Research Institute for Mathematical Sciences,
Kyoto University,
Kyoto 606-8502, Japan}
\email{takuro@kurims.kyoto-u.ac.jp}
\title{Stokes shells and Fourier transforms}

\begin{document}

\frontmatter

\begin{abstract}
Algebraic holonomic $\mathcal{D}$-modules on a complex line
are classified by the associated topological data
consisting of local systems with Stokes structure and 
the nearby and vanishing cycles at the singularities.
The Fourier transform for algebraic holonomic $\mathcal{D}$-modules 
is defined by exchanging the roles of the variable and the derivative.
It is interesting to study the induced transform for 
the associated topological data.
In particular, we closely study the local system
with Stokes structure at infinity
of the Fourier transform of a $\mathcal{D}$-module,
which also allows us to describe the remaining data. 
We introduce explicit algebraic operations
for local systems with Stokes structure,
called the local Fourier transform,
to study the case of the $\mathcal{D}$-modules
associated with basic meromorphic flat bundles.
The properties of the local Fourier transforms
are captured in terms of Stokes shells.
We also introduce the notion of extensions to study the general case.

\end{abstract}

\subjclass{14F10, 34M40, 32C38}
\keywords{meromorphic flat bundle,
$D$-module,
Fourier transform,
Stokes structure,
irregular singularity,
local Fourier transform,
Stokes shell, extension}

\maketitle

\tableofcontents

\mainmatter

\chapter{Introduction}

\section{Meromorphic flat bundles and 
local systems with Stokes structure}

\subsection{Local case}
\label{subsection;18.6.24.20}

Let $U$ be a neighbourhood of $0$ in $\cnum$
with standard coordinate $z$.
Let $(\nbigv,\nabla)$ be a meromorphic flat bundle on $(U,0)$,
i.e.,
$\nbigv$ is a locally free $\nbigo_U(\ast 0)$-module
with a flat connection $\nabla$.
Here, $\nbigo_U(\ast 0)$ denote the sheaf of
meromorphic functions allowing poles at $0$.

To understand $(\nbigv,\nabla)$,
the first goal is to know the formal structure of
$(\nbigv,\nabla)$.
For a positive integer $p$,
let $z_p$ denote a $p$-th root of $z$.
\index{formal completion $\nbigv_{|\widehat{0}}$}
Let $\nbigv_{|\widehat{0}}$
denote the formal completion of
the stalk of $\nbigv$ at $0$.
There exist
a positive integer $p>0$,
a finite subset
$\nbigi(\nbigv)\subset z_p^{-1}\cnum[z_p^{-1}]$
and a Hukuhara-Levelt-Turrittin decomposition
\begin{equation}
\label{eq;20.11.2.1}
(\nbigv,\nabla)_{|\widehat{0}}
\otimes_{\cnum(\!(z)\!)}
 \cnum(\!(z_p)\!)
=\bigoplus_{\gminia\in\nbigi(\nbigv)}
 \bigl(
 \nbigvhat_{\gminia},\nabla_{\gminia}
 \bigr)
\end{equation}
such that
$(\nbigvhat_{\gminia},\nabla_{\gminia}-d\gminia\id_{\nbigvhat_{\gminia}})$
are regular singular.
In this paper, the index set $\nbigi(\nbigv)$ is called
the set of ramified irregular values of $(\nbigv,\nabla)$.
Note that
$\nbigi(\nbigv)$ is equipped with
the $2\pi\seisuu$-action defined by
$(2\pi\ell\bullet \gminia)(z_p)
=\gminia(e^{2\pi\ell\sqrt{-1}/p}z_p)$
under the assumption that
$\nbigvhat_{\gminia}\neq 0$
for any $\gminia\in\nbigi(\nbigv)$.
\index{Hukuhara-Levelt-Turrittin decomposition}
\index{set of ramified irregular values}

Once we know the formal structure of $(\nbigv,\nabla)$,
the next goal is to study the Stokes structure of
$(\nbigv,\nabla)$.
To recall the notion of Stokes structure,
we make some preliminary.
Let $\varpi:\Utilde\lrarr U$
denote the oriented real blow up of $U$ along $0$.
\index{oriented real blow up}
Assuming $U=\{z\in\cnum\,|\,|z|<1\}$,
we have the natural identification
$\Utilde\simeq \closedopen{0}{1}\times S^1$
by the polar decomposition $z=re^{\sqrt{-1}\theta}$,
where
$\closedopen{0}{1}=\{0\leq r<1\}$.
We have the universal covering
$\varphi:\closedopen{0}{1}\times\real
 \lrarr \Utilde$
induced by $\varphi(r,\theta)=re^{\sqrt{-1}\theta}$.
We may regard $z_p=r^{1/p}\exp(\sqrt{-1}\theta/p)$,
and hence
we may naturally regard $\nbigi(\nbigv)$
as a tuple of functions on 
$\openopen{0}{1}\times\real$,
where $\openopen{0}{1}=\{0<r<1\}$.
For each $\theta$,
the partial order $\leq_{\theta}$
on $\nbigi(\nbigv)$
is determined as follows.
\begin{itemize}
\item
 $\gminia\leq_{\theta}\gminib$
 $\stackrel{\rm def}{\Longleftrightarrow}$
 there exists a neighbourhood $\nbigu_{\theta}$
of $(0,\theta)$ in 
 $\closedopen{0}{1}\times\real$
 such that
 $-\Re(\gminia)\leq -\Re(\gminib)$
 on 
 $\nbigu_{\theta}\setminus
 (\{0\}\times\real)$.
\end{itemize}
\index{order $\leq_{\theta}$}

On $U^{\ast}:=U\setminus\{0\}$,
we obtain the local system $\nbigl'$ of
the flat sections of $(\nbigv,\nabla)$.
It extends to the local system $\nbigl$ on $\Utilde$.
Let $L$ denote the $2\pi\seisuu$-equivariant
local system on $\real$
obtained as the pull back of
$\nbigl_{|\varpi^{-1}(0)}$
via the induced map $\{0\}\times\real\lrarr \varpi^{-1}(0)$.
According to the classical asymptotic analysis,
for each $\theta\in\real$,
there exists a filtration $\nbigf^{\theta}$
of $L_{|\theta}$ indexed by 
$(\nbigi(\nbigv),\leq_{\theta})$
determined as follows.
\begin{itemize}
 \item Let $(v_1,\ldots,v_m)$ be a frame of $\nbigv$ on $U$.
       An element $s\in L_{|\theta}$
       induces a flat section $\stilde$ of $\varphi^{-1}\nbigl'$ on
       $\nbigu_{\theta}\setminus \varpi^{-1}(0)$,
       where $\nbigu_{\theta}$ is a neighbourhood
       of $(0,\theta)$ in $\closedopen{0}{1}\times\real$.
       It is described as
       $\stilde=\sum_{i=1}^m s_i \varphi^{-1}(v_i)$,
       where $s_i$ are holomorphic functions on
       $\nbigu_{\theta}\setminus \varpi^{-1}(0)$.
       Then, $s$ is contained in
       $\nbigf^{\theta}_{\gminia}$
       if and only if
       there exists a neighbourhood $\nbigu_{\theta}$
       such that
       $|\exp(\gminia)s_i|=O(r^{-N})$ $(i=1,\ldots,m)$
       for a positive number $N$.
\end{itemize}
Thus, we obtain a family of filtrations $\nbigf^{\theta}$
of $L_{|\theta}$ $(\theta\in\real)$,
called the Stokes filtrations.
\index{Stokes filtration}
The family satisfies the following condition.
\begin{condition}
\label{condition;20.11.14.1}\mbox{{}}
 \begin{itemize}
 \item
 On a neighbourhood $I_{\theta}$ of $\theta$ in $\real$,
there exists a decomposition
 $L_{|I_{\theta}}=\bigoplus_{\gminia\in\nbigi(\nbigv)} 
 G_{I_{\theta},\gminia}$
 such that
the following holds for any $\theta_1\in I_{\theta}$:
\[
 \nbigf^{\theta_1}_{\gminia}(L_{|\theta_1})
=\bigoplus_{\gminib\leq_{\theta_1}\gminia}
 G_{I_{\theta},\gminib|\theta_1}.
\]
\item
     The family $\{\nbigf^{\theta}\}$ is $2\pi\seisuu$-equivariant,
     i.e.,
     $\nbigf^{\theta+2\pi}_{\gminia}(L_{|\theta+2\pi})
     =\nbigf^{\theta}_{2\pi\bullet\gminia}(L_{|\theta})$
     under the isomorphism
     $L_{|\theta}\simeq L_{|\theta+2\pi}$
     induced by the $2\pi\seisuu$-action.
     \hfill\qed
 \end{itemize}
\end{condition}
 Such a tuple $(\nbigf^{\theta}\,|\,\theta\in\real)$
is called a $2\pi\seisuu$-equivariant
local system with Stokes structure indexed by $\nbigi(\nbigv)$.
\index{Stokes structure}
We set
$\Gr^{\nbigf^{\theta}}_{\gminia}(L_{|\theta})
=\nbigf^{\theta}_{\gminia}(L_{|\theta})\big/
\sum_{\gminib<_{\theta}\gminia}
 \nbigf^{\theta}_{\gminib}(L_{|\theta})$.
The graded vector spaces
$\Gr^{\nbigf^{\theta}}(L_{|\theta})
=\bigoplus_{\gminia\in\nbigi(\nbigv)}
 \Gr_{\gminia}^{\nbigf^{\theta}}(L_{|\theta})$ 
$(\theta\in\real)$
naturally induce a $2\pi\seisuu$-equivariant
$\nbigi(\nbigv)$-graded local system
$\Gr^{\vecnbigf}(L)=
\bigoplus_{\gminia\in\nbigi(\nbigv)}\Gr^{\vecnbigf}_{\gminia}(L)$,
which is equivalent to
the right hand side
$\bigoplus_{\gminia\in\nbigi(\nbigv)}
(\nbigvhat_{\gminia},\nabla_{\gminia})$
in the Hukuhara-Levelt-Turrittin decomposition
(\ref{eq;20.11.2.1}).
\index{graded local system $\Gr^{\vecnbigf}(L)$}

The $2\pi\seisuu$-equivariant local system with Stokes structure
contains a complete information of
equivalence classes of $(\nbigv,\nabla)$.
Namely, according to the classification 
of meromorphic flat bundles due to Deligne-Malgrange-Sibuya,
the above construction induces an equivalence between 
meromorphic flat bundles on $(U,0)$
and $2\pi\seisuu$-equivariant
local systems with Stokes structure on $\real$.

\subsection{Global case}

Let $C$ be a complex curve
with a discrete subset $D$.
Let $(\nbigv,\nabla)$ be a meromorphic flat bundle
on $(C,D)$,
i.e.,
$\nbigv$ is a locally free $\nbigo_C(\ast D)$-module,
and $\nabla$ is a flat connection.
Let $\varpi:\Ctilde\lrarr C$ 
denote the oriented real blow up of $C$ along $D$.
We obtain the local system $\nbigl(\nbigv)$ on $\Ctilde$
induced by the sheaf of flat sections of $(\nbigv,\nabla)$.
For each $P\in D$,
take a holomorphic local coordinate neighbourhood
$(C_P,z_P)$ with $z_P(P)=0$.
Associated with $(\nbigv,\nabla)_{|C_P}$,
we obtain the set of the ramified irregular values
$\nbigi_P(\nbigv)$
and the $2\pi\seisuu$-equivariant
local system with Stokes structure $(L_P(\nbigv),\vecnbigf)$.
The classification of Deligne-Malgrange-Sibuya
implies that
$(\nbigv,\nabla)$ is classified by
$\nbigl$
with 
$\nbigi_P(\nbigv)$ and $(L_P(\nbigv),\vecnbigf)$ $(P\in D)$,
up to isomorphisms.
\index{set of ramified irregular values $\nbigi_P(\nbigv)$}
\index{local system with Stokes structure $(L_P(\nbigv),\vecnbigf)$}

\section{Main purpose in this paper}

\subsection{Fourier transform of $\nbigd$-modules}
Let $\nbigm$ be any algebraic $\nbigd$-module on 
$\cnum_z$.
We obtain its Fourier transform
$\Fourier_+(\nbigm)$ 
on $\cnum_w$.
\index{Fourier transform $\Fourier_+(\nbigm)$}
It is defined as an integral transform
as follows.
Let $p_z$ and $p_w$ denote the projections of
$\cnum_z\times\cnum_w$
onto $\cnum_z$ and $\cnum_w$, respectively.
Let $\nbige(zw)$ denote the algebraic meromorphic flat bundle
$(\nbigo_{\cnum_z\times\cnum_w},d+d(zw))$.
We obtain the algebraic $\nbigd$-module
on $\cnum_w$:
\begin{equation}
\label{eq;21.4.28.1}
 \Fourier_+(\nbigm):=
 p_{w+}^0\bigl(
 p_z^{\ast}(\nbigm)\otimes\nbige(zw)
 \bigr).
\end{equation}
It is also defined in terms of modules over Weyl algebras.
Let $W_z$ denote the algebra of algebraic differential operators
on $\cnum[z]$,
i.e.,
$W_z=\cnum[z]\langle\del_z\rangle$.
Let $M$ be any $W_z$-module.
We obtain a $W_w$-module $M^{\gbigf}$ 
as follows.
We set $M^{\gbigf}:=M$
as a $\cnum$-vector space.
We define the action of $W_w$ on $M^{\gbigf}$
by
$\del_w(m)=zm$ and $wm=-\del_zm$.
From any algebraic $\nbigd$-module $\nbigm$
on $\cnum_z$,
we obtain the $W_z$-module $H^0(\cnum,\nbigm)$,
and $\Fourier_+(\nbigm)$
is characterized as the algebraic $\nbigd$-module
corresponding to the $W_w$-module 
$H^0(\cnum,\nbigm)^{\gbigf}$.
The Fourier transform was studied 
by Malgrange \cite{Malgrange-book}
comprehensively.

\subsection{Main issue in this monograph}

We naturally regard algebraic holonomic $\nbigd$-module
$\nbigm$ on $\cnum$
as an analytic holonomic $\nbigd_{\proj^1}(\ast\infty)$-module,
i.e.,
an analytic holonomic $\nbigd_{\proj^1}$-module $\nbigm$
such that $\nbigm(\ast\infty)=\nbigm$.
Because $\nbigm$ is holonomic,
$\Fourier_+(\nbigm)$ is also holonomic.
There exists a neighbourhood $U_{\infty}$ of $\infty$
in $\proj^1_w$
such that 
$\Fourier_+(\nbigm)_{|U_{\infty}}$
is a meromorphic flat bundle
on $(U_{\infty},\infty)$.
On $U_{\infty}$, we use the coordinate $u=w^{-1}$.
Let $\nbigi_{\infty}(\Fourier_+(\nbigm))$ denote
the set of ramified irregular values of 
$\Fourier_+(\nbigm)_{|U_{\infty}}$.
We obtain 
the corresponding
$2\pi\seisuu$-equivariant local system
with Stokes structure 
\index{local system with Stokes structure $(\gbigl^{\gbigf}(\nbigm),\vecnbigf)$}
\[
 (\gbigl^{\gbigf}(\nbigm),\vecnbigf)
 =\Bigl(
 L_{\infty} \bigl(
\Fourier_+(\nbigm)
 \bigr),\vecnbigf\Bigr)
\]
indexed by
$\nbigi_{\infty}(\Fourier_+(\nbigm))$
on $\real$.
It is our main issue in this monograph
to study how
$(\gbigl^{\gbigf}(\nbigm),\vecnbigf)$
is described in terms of 
the topological data associated with $\nbigm$.

\subsection{Goal}
\label{subsection;21.4.28.2}

Let $D\subset\cnum$ be a finite subset
such that
$\nbigv=\nbigo_{\proj^1}(\ast D)\otimes_{\nbigo_{\proj^1}}\nbigm$
is a meromorphic flat bundle
on $(\proj^1,\Dbar)$, where $\Dbar=D\cup\{\infty\}$.
Around $\alpha\in D$,
we use the coordinate $z-\alpha$.
We obtain the following data
$\LS^{\fin}(\nbigm)$ associated with $\nbigm$.
\index{data $\LS^{\fin}(\nbigm)$}
\begin{itemize}
 \item The local system $\nbigl(\nbigm)$ on $\cnum\setminus D$
       obtained as
       the sheaf of flat sections of $\nbigm_{|\cnum\setminus D}$.
 \item The $2\pi\seisuu$-equivariant
       local systems with Stokes structure
       \[
       (L_{\alpha}(\nbigm),\vecnbigf):=
       (L_{\alpha}(\nbigv),\vecnbigf)\quad(\alpha\in D).
       \]
 \item The vector spaces
       $\psi_{z-\alpha}(\nbigm)=\Gr^{V^{\alpha}}_{-1}(\nbigm)$
       and
       $\phi_{z-\alpha}(\nbigm)=\Gr^{V^{\alpha}}_0(\nbigm)$,
       where $V^{\alpha}_{\bullet}(\nbigm)$
       denotes the $V$-filtration of $\nbigm$ along $z-\alpha$.
       They are equipped with the standard morphisms
\begin{equation}
\label{eq;21.4.27.10}
\begin{CD}
 \psi_{z-\alpha}(\nbigm)
 @>{\can}>>
 \phi_{z-\alpha}(\nbigm)
 @>{\var}>>
 \psi_{z-\alpha}(\nbigm),
\end{CD}
\end{equation}
       where $\can$ is induced by $\del_z$,
       and $\var$ is induced by $z-\alpha$.
\end{itemize}
Around $\infty$, we use the coordinate $x=z^{-1}$.
We obtain the local system with Stokes structure
\[
 (L_{\infty}(\nbigm),\vecnbigf).
\]
We obtain the tuple
$\LS(\nbigm)$
by adding $(L_{\infty},\vecnbigf)$
to $\LS^{\fin}(\nbigm)$.
\index{data $\LS(\nbigm)$}

It is our main goal to compute
$(\gbigl^{\gbigf}(\nbigm),\vecnbigf)$
directly from $\LS(\nbigm)$.

\begin{rem}
It also allows us to compute
$\LS^{\fin}(\Fourier_+(\nbigm))$
directly from $\LS(\nbigm)$.
\hfill\qed
\end{rem}

\subsection{Previous studies}

The associated $\nbigi_{\infty}(\Fourier_+(\nbigm))$-graded
local system
$\Gr^{\vecnbigf}\gbigl^{\gbigf}(\nbigm)$
has been completely well understood
by the local Fourier transforms
and their stationary phase formula.
For any $\alpha\in\proj^1$,
let $\nbigm_{|\widehat{\alpha}}$
denote the formal completion of $\nbigm$ at $\alpha$.
Let $D\subset\cnum$ be a finite subset
such that $\nbigm(\ast D)$
is a meromorphic flat bundle on
$(\proj^1,\Dbar)$, where $\Dbar=D\cup\{\infty\}$.
In the case $\nbigm=\nbigm(\ast D)$,
Bloch and Esnault introduced 
the local Fourier transforms
in \cite{Bloch-Esnault1},
and proved that 
$\Fourier_+(\nbigm)_{|\inftyhat}$
is decomposed into
the direct sum of 
the local Fourier transforms
of $\nbigm_{|\alphahat}$ $(\alpha\in \Dbar)$.
(See also \cite{Lopez}.)
It is generalized to the case of general
holonomic $\nbigd_{\proj^1}(\ast\infty)$-modules.
(See \cite{Sabbah-stationary}.)
Fang \cite{Fang},
Graham-Squire \cite{Graham-Squire},
and Sabbah \cite{Sabbah-stationary}
obtained the explicit description of the local Fourier transforms,
called the stationary phase formula.
D'Agnolo and Kashiwara \cite{D'Agnolo-Kashiwara-Fourier}
applied the theory of enhanced ind-sheaves
to the study of Fourier transforms.
The stationary phase formula implies that
$\nbigi_{\infty}\bigl(
\Fourier_+(\nbigm)\bigr)$
and $\Gr^{\vecnbigf}\gbigl^{\gbigf}(\nbigm)$
are explicitly described in terms of
$\nbigi_{\alpha}(\nbigm)$
and $\Gr^{\vecnbigf}L_{\alpha}(\nbigm)$
$(\alpha\in\Dbar)$.

It is Malgrange \cite{Malgrange-book}
who pioneered the study of this issue.
To the best of the author's understanding,
he especially obtained the following.
\begin{itemize}
 \item
      For $\alpha\in D$,
      $\Gr^{\vecnbigf^{(1)}}_{\alpha u^{-1}}
      \bigl(
      \gbigl^{\gbigf}(\nbigm),\vecnbigf
      \bigr)$
      is determined by
      the restriction of 
      $\nbigm$ to a neighbourhood of $\alpha$.
      If $\alpha=0$, the morphisms of the constructible sheaves
\[
      \Gr^{\vecnbigf^{(1)}}_0\bigl(
      \gbigl^{\gbigf}(\nbigm)
      \bigr)^{<0}
      \lrarr
      \Gr^{\vecnbigf^{(1)}}_0\bigl(
      \gbigl^{\gbigf}(\nbigm)
      \bigr)^{\leq 0}
      \lrarr
      \Gr^{\vecnbigf^{(1)}}_0\bigl(
      \gbigl^{\gbigf}(\nbigm)
      \bigr)
\]
      are described as
      a topologically defined Fourier transform
      of some morphisms of constructible sheaves induced by
      $(L_0(\nbigm),\vecnbigf)$
      and
      $\psi_{z-\alpha}(\nbigm)\to \phi_{z-\alpha}(\nbigm)
      \to\psi_{z-\alpha}$.
      (See \S\ref{subsection;21.4.24.13}
      for $\vecnbigf^{(\omega)}$,
      and \S\ref{subsection;18.4.18.1}
      for $L^{<0}$
      and $L^{\leq 0}$ induced by
      $(L,\vecnbigf)$.)
      Moreover,
      the filtered constructible sheaves
      $\bigl(
      \Gr^{\vecnbigf^{(1)}}_{0}\bigl(
      \gbigl^{\gbigf}(\nbigm)
      \bigr)^{<0},\vecnbigf\bigr)$
      and
      $\bigl(
      \Gr^{\vecnbigf^{(1)}}_0\bigl(
      \gbigl^{\gbigf}(\nbigm)\bigr)
      \big/
      \Gr^{\vecnbigf^{(1)}}_0\bigl(
      \gbigl^{\gbigf}(\nbigm)
      \bigr)^{\leq 0},\vecnbigf\bigr)$
      are isomorphic to
      the Legendre transform of
      the filtered constructible sheaves
      $(L_0(\nbigm)^{<0},\vecnbigf)$
      and
      $\bigl(
      L_0(\nbigm)\big/
      L_0(\nbigm)^{\leq 0},\vecnbigf\bigr)$,
      respectively.
      (See also Remark \ref{rem;25.4.12.20}.)

      There are similar formulas
      in the case $\alpha\neq 0$.
\item
      The morphisms of the constructible sheaves
\[
     \nbigstilde_1\bigl(
     \gbigl^{\gbigf}(\nbigm)
     \bigr)^{<0}
     \lrarr
     \nbigstilde_1\bigl(
     \gbigl^{\gbigf}(\nbigm)
     \bigr)^{\leq 0}
     \lrarr
     \nbigstilde_1\bigl(
     \gbigl^{\gbigf}(\nbigm)
     \bigr)
\]
     are described as a topologically defined
     transform of some morphisms of constructible sheaves
     induced by 
     $\bigl(
     L_{\infty}(\nbigm),\vecnbigf
     \bigr)$
     and
     $\DR(\nbigm)$.
     (See \S\ref{subsection;21.4.24.12}
     for $\nbigstilde_{\omega}$.)
     Moreover,
     the filtered constructible sheaves
     $\bigl(
     \nbigstilde_1\bigl(
     \gbigl^{\gbigf}(\nbigm)
     \bigr)^{<0},\vecnbigf\bigr)$
     and
     $\bigl(
     \nbigstilde_1\bigl(
     \gbigl^{\gbigf}(\nbigm)
     \bigr)\big/
     \nbigstilde_1\bigl(
     \gbigl^{\gbigf}(\nbigm)
     \bigr)^{\leq 0},\vecnbigf\bigr)$
     are described as
     the Legendre transform of
     $\bigl(
     L_{\infty}(\nbigm)
     \big/
     L_{\infty}(\nbigm)
     ^{\leq 0},\vecnbigf\bigr)$
     and
     $\bigl(
     L_{\infty}(\nbigm)^{<0},\vecnbigf
     \bigr)$,
     respectively.
     (See also Remark \ref{rem;25.4.12.21}.)
\item $\nbigttilde_1(\gbigl^{\gbigf}(\nbigm),\vecnbigf^{(1)})$
      is topologically described 
      in terms of $\DR(\nbigm)$.
      (See \S\ref{subsection;21.4.24.12}
      for $\nbigttilde_{\omega}$.)
 \item $\LS^{\fin}(\Fourier_+(\nbigm))$
       is also described
       in terms of $\LS(\nbigm)$.
\end{itemize}
See \cite{Malgrange-book} for
more detailed and precise explanation on his work.
In \cite{Mochizuki-Fourier-old},
the author studied how 
$(\gbigl^{\gbigf}(\nbigm),\vecnbigf)$
is described in the case where
$\nbigm=\nbigm(\ast D)$,
on the basis of 
the rapid decay homology theory
in \cite{Bloch-Esnault2, Hien},
and the saddle point method
in \cite{Beilinson-Bloch-Deligne-Esnault}.
Let $\nbigm^{\lor}$ denote the dual meromorphic flat bundle
on $(\proj^1,\Dbar)$,
i.e., it is obtained as
$\nbigm^{\lor}=\nhom_{\nbigo_{\proj^1}(\ast\Dbar)}
(\nbigm,\nbigo_{\proj^1}(\ast\Dbar))$.
The vector spaces
$\gbigl^{\gbigf}(\nbigm)_{|\theta}$
are naturally isomorphic to
the dual space of the rapid decay homology group
$H_1^{\rd}\bigl(\cnum\setminus D,
\nbigm^{\lor}
\otimes\nbige(-zu^{-1})
\bigr)$
for $u=|u|e^{\sqrt{-1}\theta}$,
where $|u|$ is sufficiently small.
To compute the Stokes filtration,
we need to find rapid decay $1$-cycles
in a way that the growth orders of
the induced flat sections are controlled.
In \cite{Mochizuki-Fourier-old},
we studied such choices of $1$-cycles
by following the idea of the saddle point method
due to Beilinson-Bloch-Deligne-Esnault.
However, the result was not completely explicit
in the general case.
In \cite{Hien-Sabbah},
Hien and Sabbah introduced 
topological Laplace transforms
for local systems with Stokes structure,
which is the counterpart of 
the Fourier transform of holonomic $\nbigd$-modules
(\ref{eq;21.4.28.1}),
which allows us to study the Fourier transform
in a purely topological way.
In particular, they closely studied the case of elementary 
meromorphic flat bundles.
The irregular Riemann-Hilbert correspondence
due to D'Agnolo and Kashiwara
\cite{D'Agnolo-Kashiwara-irregular-Riemann-Hilbert}
enables us to translate the integral transform
(\ref{eq;21.4.28.1})
into the transformation of enhanced ind-sheaves.
It also allows us to study the Fourier transform
in a topological way,
which was applied by
D'Agnolo, Hien, Morando and Sabbah
in \cite{D'Agnolo-Hien-Morando-Sabbah}
to study the Stokes structure of
the Fourier transform of
regular singular holonomic $\nbigd$-modules.
The both theories of topological Laplace transform
and enhanced ind-sheaves
provide us with topological counterparts of
the integral transform (\ref{eq;21.4.28.1}),
which are theoretically significant and useful.
However, in general,
they remain to contain non-trivial operations
which are not so easy to compute.
Sabbah studied the pure Gaussian case
\cite{Sabbah-pure-Gaussian}.
(See also \cite{Hohl-pure-Gaussian}.)

In this paper, we revisit the problem
by following the idea in \cite{Beilinson-Bloch-Deligne-Esnault} again.
This study is also regarded as an attempt to make the computations in
the previous works more explicit by using the homology theory.

\begin{rem}
More recently, Dou\c{c}ot and Hohl \cite{Doucot-Hohl}
studied a topological algorithm for this issue
for algebraic connections on $\cnum$
under some assumptions
based on the version $3$ of this monograph on arXiv.
\hfill\qed
\end{rem}

\subsection{Outline of the introduction}

We briefly review Stokes structures 
in \S\ref{subsection;20.11.14.31}.
We explain building blocks for our study
in \S\ref{subsection;24.4.15.100}
and \S\ref{subsection;24.4.19.10}.
We introduce the notion of {\em extension} of
local systems with Stokes structure
in \S\ref{section;20.11.4.32}.
We explain the process to reduce
the Stokes structure of
$(\gbigl^{\gbigf}(\nbigm),\vecnbigf)$
in \S\ref{subsection;24.4.15.101}.
Then,
we outline how to describe
$(\gbigl^{\gbigf}(\nbigm),\vecnbigf)$
in \S\ref{section;21.4.28.3}.
We mention some easy examples
in \S\ref{subsection;21.6.13.2}.

\section{A brief review of Stokes structures}
\label{subsection;20.11.14.31}

Before explaining our results,
we briefly recall the standard theory of Stokes structures.
Useful references are
\cite{Boalch-survey,Sabbah-Introduction-Stokes}.
The higher dimensional case was explained in
\cite{Mochizuki-good-Stokes, Mochizuki-wild}.

\subsection{Index sets and special directions}

Let $p$ be a positive integer.
Take a $p$-th root $z_p$ of the variable $z$.
Let $\Gal(p)$ be the group of $p$-th roots of $1$.
\index{group $\Gal(p)$}
We have the natural $\Gal(p)$-actions on $\cnum(\!(z_p)\!)$
and $z_p^{-1}\cnum[z_p^{-1}]$
by $(\tau\bullet f)(z_p)=f(\tau z_p)$
for $\tau\in\Gal(p)$.
There exists the natural map
$2\pi\seisuu\lrarr \Gal(p)$
by $2\pi \ell\longmapsto \exp(2\pi\ell\sqrt{-1}/p)$.
It induces a $2\pi\seisuu$-actions
on $\cnum(\!(z_p)\!)$
and $z_p^{-1}\cnum[z_p^{-1}]$.
For each $\omega=n/p\in \frac{1}{p}\seisuu_{>0}$,
we define the map
$\pi_{\omega}:
 z_p^{-1}\cnum[z_p^{-1}]\lrarr z_p^{-n}\cnum[z_p^{-1}]$
by
$\pi_{\omega}\bigl(
 \sum_{j\geq 1} a_jz_p^{-j}
 \bigr):=
 \sum_{j\geq n}a_jz_p^{-j}$.
\index{map $\pi_{\omega}$}
For any non-zero
$\gminia=\sum_{j=1}^m\gminia_jz_p^{-j}
\in z_p^{-1}\cnum[z_p^{-1}]$,
we set
$\ord(\gminia)=-\frac{1}{p}\max\{j\,|\,\gminia_j\neq 0\}$.
We set $\ord(0)=\infty$.
\index{order $\ord(\gminia)$}
 
For two distinct
$\gminia,\gminib\in z_p^{-1}\cnum[z_p^{-1}]$,
by using the expansion
$\gminia-\gminib
=\sum(\gminia-\gminib)_jz_p^{-j}$,
we set
\[
 S(\gminia,\gminib)
 =\bigl\{
 \theta\in\real\,\big|\,
 \Re\bigl(
 (\gminia-\gminib)_{-p\ord(\gminia-\gminib)}
 e^{\sqrt{-1}\ord(\gminia-\gminib)\theta}
 \bigr)=0
 \bigr\},
\]
\[
 A(\gminia,\gminib)
 =\bigl\{
 \theta\in\real\,\big|\,
 \Image\bigl(
 (\gminia-\gminib)_{-p\ord(\gminia-\gminib)}
 e^{\sqrt{-1}\ord(\gminia-\gminib)\theta}
 \bigr)=0
 \bigr\}.
\]
Any $\theta\in S(\gminia,\gminib)$
(resp. $\theta\in A(\gminia,\gminib)$)
is called the Stokes direction
(resp. the anti-Stokes direction)
with respect to $\gminia$ and $\gminib$.
\index{Stokes direction}
\index{anti-Stokes direction}

For a $\Gal(p)$-invariant finite subset
$\nbigi\subset z_p^{-1}\cnum[z_p^{-1}]$,
we set \index{order $\ord(\nbigi)$}
\[
\ord(\nbigi)
=\min\{\ord(\gminia)\,|\,\gminia\in\nbigi\}.
\]
We also set \index{set $S(\nbigi)$} \index{set $A(\nbigi)$}
\begin{equation}
 \label{eq;21.4.24.1}
S(\nbigi):=
 \bigcup_{
 \substack{\gminia,\gminib\in\nbigi\\
 \gminia\neq\gminib}}
  S(\gminia,\gminib),
  \quad\quad
 A(\nbigi):=
 \bigcup_{\substack{\gminia,\gminib\in\nbigi\\
  \gminia\neq\gminib
  }}
   A(\gminia,\gminib).
\end{equation}
Any $\theta\in S(\nbigi)$ (resp. $\theta\in A(\nbigi)$)
is called the Stokes direction
(resp. anti-Stokes direction)
with respect to $\nbigi$.

\subsection{Stokes structures}

Let $\nbigi$ be a $\Gal(p)$-invariant subset
of $z_p^{-1}\cnum[z_p^{-1}]$.
Let $L$ be a $2\pi\seisuu$-equivariant local system
on $\real$.
A family 
$\vecnbigf=(\nbigf^{\theta}\,|\,\theta\in\real)$
of filtrations
of $L_{|\theta}$
indexed by
$(\nbigi,\leq_{\theta})$
is called a $2\pi\seisuu$-equivariant
Stokes structure of $L$
indexed by $\nbigi$
if the condition \ref{condition;20.11.14.1} is satisfied.
\index{Stokes structure}
Note that
the filtrations $\nbigf^{\theta}$ are
constant on the complement of $S(\nbigi)$.

\begin{rem}
In general,
we allow
$\Gr^{\vecnbigf}_{\gminia}(L)=0$
for some $\gminia\in\nbigi$.
\hfill\qed
\end{rem}

Let $\Loc^{\St}(\nbigi)$ denote the category of
$2\pi\seisuu$-equivariant local systems with Stokes structure
indexed by $\nbigi$.
\index{category $\Loc^{\St}(\nbigi)$}
A morphism $f:(L_1,\vecnbigf)\lrarr(L_2,\vecnbigf)$
is defined to be
a morphism of $2\pi\seisuu$-equivariant local systems
$f:L_1\lrarr L_2$
such  that
$f(\nbigf^{\theta}_{\gminia}(L_{1|\theta}))
\subset
 \nbigf^{\theta}_{\gminia}(L_{2|\theta})$
for any $\theta\in\real$ and $\gminia\in\nbigi$.
It is well known that
$\Loc^{\St}(\nbigi)$ is an abelian category.
We can prove it by using 
the Riemann-Hilbert correspondence,
or more directly
by using the canonical splittings in \S\ref{subsection;18.5.13.10}.

Note that $\Loc^{\St}(0)$
is the category of $2\pi\seisuu$-equivariant local systems.

\subsection{Lower level Stokes structures and the associated
  graded objects}
\label{subsection;21.4.24.13}

Let $(L,\vecnbigf)\in\Loc^{\St}(\nbigi)$.
For each $\theta\in\real$
and $\gminib\in\pi_{\omega}(\nbigi)$,
by using a splitting
$L_{|\theta}=\bigoplus_{\gminia\in\nbigi} G_{\theta,\gminia}$
of $\nbigf^{\theta}$,
we set
\[
 \nbigf^{(\omega)\,\theta}_{\gminib}
 =\bigoplus_{\substack{\gminia\in\nbigi\\
  \pi_{\omega}(\gminia)\leq_{\theta}\gminib}}G_{\theta,\gminia},
\]
which is independent of the choice of a splitting.
We obtain a filtration
$\nbigf^{(\omega)\,\theta}$
of $L_{|\theta}$
indexed by
$(\pi_{\omega}(\nbigi),\leq_{\theta})$.
The family
$\vecnbigf^{(\omega)}=
\bigl(
\nbigf^{(\omega)\theta}\,\big|\,
\theta\in\real
\bigr)$
defines a $2\pi\seisuu$-equivariant
Stokes structure on $L$,
i.e.,
$(L,\vecnbigf^{(\omega)})
\in\Loc^{\St}(\pi_{\omega}(\nbigi))$.
\index{Stokes structure $\vecnbigf^{(\omega)}$}
We obtain
the $2\pi\seisuu$-equivariant
$\pi_{\omega}(\nbigi)$-graded local system
$\Gr^{\vecnbigf^{(\omega)}}(L)
=\bigoplus_{\gminib\in\pi_{\omega}(\nbigi)}
 \Gr^{\vecnbigf^{(\omega)}}_{\gminib}(L)$.
 For each $\gminib\in\pi_{\omega}(\nbigi)$,
 we set $\nbigi(\gminib):=\bigl\{
 \gminia\in\nbigi\,\big|\,\pi_{\omega}(\gminia)=\gminib\bigr\}$.
There exists an induced filtration
$\nbigf^{\theta}$ on $\Gr_{\gminib}^{\vecnbigf^{(\omega)}}(L)_{|\theta}$
indexed by
$(\nbigi(\gminib),\leq_{\theta})$.
As the direct sum,
we obtain the induced filtration
$\nbigf^{\theta}$
on $\Gr^{\vecnbigf^{(\omega)}}(L)_{|\theta}$
indexed by
$(\nbigi,\leq_{\theta})$.
Note that
$\Gr^{\vecnbigf^{(\omega)}}(L)$
with the induced family of filtrations $\vecnbigf$
is a $2\pi\seisuu$-equivariant local system with
Stokes structure,
denoted by
$\Gr^{\vecnbigf^{(\omega)}}(L,\vecnbigf)$.
\index{graded local system with Stokes structure
$\Gr^{\vecnbigf^{(\omega)}}(L,\vecnbigf)$}
This procedure defines a functor
$\Gr^{\vecnbigf^{(\omega)}}:
\Loc^{\St}(\nbigi)\lrarr \Loc^{\St}(\nbigi)$.
It is easy to see that
$\Gr^{\vecnbigf^{(\omega_1)}}
\Gr^{\vecnbigf^{(\omega_2)}}(L,\vecnbigf)$
is naturally isomorphic to
$\Gr^{\vecnbigf^{(\omega_3)}}(L,\vecnbigf)$,
where $\omega_3=\min\{\omega_1,\omega_2\}$.

Let us emphasize that
$(L,\vecnbigf)$ is easily recovered from
$(L,\vecnbigf^{(\omega)})$
and
$\Gr^{\vecnbigf^{(\omega)}}(L,\vecnbigf)$.

\subsection{Some induced local systems with Stokes structure}
\label{subsection;21.4.24.12}

We explain some notation and
some induced local systems with Stokes structure
which are useful in the study of Fourier transform.

Let $\omega\in\frac{1}{p}\seisuu_{\geq 0}$.
Let $\nbigi\subset z_p^{-1}\cnum[z_p^{-1}]$
be any $\Gal(p)$-invariant finite subset.
We put
\[
\nbigt_{\omega}(\nbigi):=
 \bigl\{
 \gminia\in\nbigi\,\big|\,
 \pi_{\omega}(\gminia)=0
 \bigr\}\cup\{0\},
\quad
 \nbigs_{\omega}(\nbigi):=
 \Bigl(
 \nbigi\setminus
\nbigt_{\omega}(\nbigi)
 \Bigr)
 \cup\{0\}.
\]
\index{set $\nbigt_{\omega}(\nbigi)$}
\index{set $\nbigs_{\omega}(\nbigi)$}
We also set
$\nbigttilde_{\omega}(\nbigi):=
 \nbigt_{\omega+p^{-1}}(\nbigi)$
and 
$\nbigstilde_{\omega}(\nbigi):=
 \nbigs_{\omega+p^{-1}}(\nbigi)$.
\index{set $\nbigttilde_{\omega}(\nbigi)$}
\index{set $\nbigstilde_{\omega}(\nbigi)$}
Clearly,
$\nbigt_{\omega}(\nbigi)$
and $\nbigs_{\omega}(\nbigi)$
are also $\Gal(p)$-invariant.

Let $(L,\vecnbigf)\in\Loc^{\St}(\nbigi)$.
Recall that we obtain
$(L,\vecnbigf^{(\omega)})
\in\Loc^{\St}(\pi_{\omega}(\nbigi))$.
We set
$\nbigt_{\omega}(L,\vecnbigf)
:=
\Gr^{\vecnbigf^{(\omega)}}_0(L,\vecnbigf)
\in\Loc^{\St}(\nbigt_{\omega}(\nbigi))$.
\index{local system with Stokes structure $\nbigt_{\omega}(L,\vecnbigf)$}
For each $\theta$ and $\gminia\in \nbigs_{\omega}(\nbigi)$,
we obtain a new filtration
$\nbigs_{\omega}(\nbigf)^{\theta}$
of $L_{|\theta}$
indexed by
$(\nbigs_{\omega}(\nbigi),\leq_{\theta})$
as follows:
\[
 \nbigs_{\omega}(\nbigf)^{\theta}_{\gminia}
 (L_{|\theta})
 =\left\{
\begin{array}{ll}
 \nbigf^{\theta}_{\gminia}(L_{|\theta})&
  (\gminia\neq 0)
  \\
 \nbigf^{(\omega)\theta}_0(L_{|\theta}).
\end{array}
 \right.
\]
They induce a $2\pi\seisuu$-equivariant Stokes structure
$\nbigs_{\omega}(\vecnbigf)$ of $L$.
Let $\nbigs_{\omega}(L,\vecnbigf)$
denote the $2\pi\seisuu$-equivariant
local system with Stokes structure
$(L,\nbigs_{\omega}(\vecnbigf))$.
\index{local system with Stokes structure $\nbigs_{\omega}(L,\vecnbigf)$}
We also set
$\nbigttilde_{\omega}(L,\vecnbigf):=
\nbigt_{\omega+p^{-1}}(L,\vecnbigf)
\in\Loc^{\St}(\nbigttilde_{\omega}(\nbigi))$,
and 
$\nbigstilde_{\omega}(L,\vecnbigf):=
\nbigs_{\omega+p^{-1}}(L,\vecnbigf)
\in\Loc^{\St}(\nbigstilde_{\omega}(\nbigi))$.
\index{local system with Stokes structure $\nbigstilde_{\omega}(L,\vecnbigf)$}
\index{local system with Stokes structure $\nbigttilde_{\omega}(L,\vecnbigf)$}

It is easy to observe that
we can recover
$(L,\vecnbigf)$
from 
$\nbigttilde_{\omega}(L,\vecnbigf)$
and 
$\nbigstilde_{\omega}(L,\vecnbigf)$.
Moreover,
we can recover
$\nbigttilde_{\omega}(L,\vecnbigf)$
from
$\nbigs_{\omega}\nbigttilde_{\omega}(L,\vecnbigf)=
\nbigttilde_{\omega}\nbigs_{\omega}(L,\vecnbigf)$
and
$\nbigt_{\omega}\nbigttilde_{\omega}(L,\vecnbigf)
=\nbigt_{\omega}(L,\vecnbigf)$.
This allows us to study $(L,\vecnbigf)$
in an inductive way.
Namely, let $\omega_1>\omega_2>\cdots>\omega_{\ell}\geq 0$
be the rational numbers such that
$\{\omega_1,\ldots,\omega_{\ell}\}=
\bigl\{-\ord(\gminia)\,\big|\,\gminia\in\nbigi\bigr\}$.
Then, for any $1\leq m\leq \ell$,
we can recover
$(L,\vecnbigf)$ from
$\nbigs_{\omega_i}\nbigttilde_{\omega_i}(L,\vecnbigf)$
$(i=1,\ldots,m)$
and $\nbigttilde_{\omega_m}(L,\vecnbigf)$.

\section{Building blocks from the irregular singularity}
\label{subsection;24.4.15.100}

\subsection{Basic meromorphic flat bundles}
\label{section;21.4.24.11}

Let us introduce the notion of 
basic meromorphic flat bundles
which are building blocks in our study.
\index{basic meromorphic flat bundle}

\begin{df}
Let $\omega\in\rnum_{>0}$.
A meromorphic flat bundle $(\nbigv,\nabla)$ on $(\proj^1,\{0,\infty\})$
is called basic of level $(0,\omega)$ if the following holds.
\index{basic of level $(0,\omega)$}
 \begin{itemize}
 \item $(\nbigv,\nabla)$ is regular singular at $\infty$.
  \item $\nbigs_{\omega}\nbigttilde_{\omega}(\nbigi_0(\nbigv))
	=\nbigi_0(\nbigv)$.
       \hfill\qed
 \end{itemize}
\end{df}

We use the coordinate $x=z^{-1}$ around $\infty$ in $\proj^1$.
\begin{df}
Let $\omega\in\rnum_{>0}$.
If $\omega\neq 1$,
we say that
a meromorphic flat bundle $(\nbigv,\nabla)$ on $(\proj^1,\{0,\infty\})$
is basic of level $(\infty,\omega)$ if the following holds.
\index{basic of level $(\infty,\omega)$}
\begin{itemize}
 \item $(\nbigv,\nabla)$ is regular singular at $0$.
 \item $\nbigs_{\omega}\nbigttilde_{\omega}(\nbigi_{\infty}(\nbigv))
       =\nbigi_{\infty}(\nbigv)$.
\end{itemize}
We say that 
a meromorphic flat bundle $(\nbigv,\nabla)$ on $(\proj^1,\{0,\infty\})$
is basic of level $(\infty,1)$ if the following holds.
\index{basic of level $(\infty,1)$}
\begin{itemize}
 \item $(\nbigv,\nabla)$ is regular singular at $0$.
 \item $\nbigi_{\infty}(\nbigv)\subset \cnum x^{-1}$.
  \hfill\qed
\end{itemize}
\end{df}

\subsection{Preliminary}

For a meromorphic flat bundle $(\nbigv,\nabla)$
on $(\proj^1,\{0,\infty\})$,
let $(\nbigv^{\lor},\nabla)$ denote the dual meromorphic flat bundle
on $(\proj^1,\{0,\infty\})$,
i.e.,
$\nbigv^{\lor}
=\nhom_{\nbigo_{\proj^1}(\ast\{0,\infty\})}
\bigl(
 \nbigv,\nbigo_{\proj^1}(\ast\{0,\infty\})
\bigr)$.
We obtain the $\nbigd_{\proj^1}$-module
$\nbigv(!0):=
\DDD_{\proj^1}(\nbigv^{\lor})
\otimes\nbigo_{\proj^1}(\ast\infty)$,
where $\DDD_{\proj^1}$ denotes the duality functor
for $\nbigd_{\proj^1}$-modules.
We also have $\nbigv(\ast 0)=\nbigv$.
\index{$\nbigd$-module \mbox{$\nbigv(! 0)$}}
\index{$\nbigd$-module $\nbigv(\ast 0)$}
To simplify the notation,
we set
$\gbigl^{\gbigf}_{\star}(\nbigv)
=\gbigl^{\gbigf}(\nbigv(\star 0))$
$(\star=!,\ast)$.
\index{local system $\gbigl^{\gbigf}_{\star}(\nbigv)$}

For a variable $y$,
and for $p,n\in\seisuu_{>0}$,
we set
\index{set $\nbigu_{y}(p,n)$}
\[
 \nbigu_{y}(p,n)=
\bigl\{\gminia\in y_p^{-1}\cnum[y_p^{-1}]\,\big|\,
-\ord(\gminia)=n/p
\bigr\}\cup\{0\}.
\]

\subsection{Fourier transform in the basic case of level $(0,\omega)$}

Let $\omega=n/p$
for some positive integers $n$ and $p$.

\subsubsection{Transformation of index sets}

As we will review in \S\ref{subsection;24.3.26.1},
by the stationary phase formula
\cite{Fang, Graham-Squire, Sabbah-stationary},
a $\Gal(n+p)$-invariant subset
$\gbigf^{(0,\infty)}_+(\nbigi)$
of $\gbigu_u(p+n,n)$ is attached
to any $\Gal(p)$-invariant subset $\nbigi$ of $\gbigu_z(p,n)$,
and the following holds.
\index{set $\gbigf^{(0,\infty)}_+(\nbigi)$}
\begin{itemize}
 \item If $(V,\nabla)$ is basic of
       level $(0,\omega)$,
       then
       $\Fourier_+(V(\star 0))(\ast 0)$
       $(\star=\ast,!)$
       are basic meromorphic flat bundles
       of level $(\infty,\frac{\omega}{1+\omega})$,
       and we have
\[
       \nbigi_{\infty}\bigl(
       \Fourier_+(V(\star 0))
       \bigr)
       =
       \gbigf^{(0,\infty)}_+\bigl(
       \nbigi_0(V)
       \bigr).
\]
\end{itemize}

\subsubsection{Local Fourier transform of local systems with Stokes structure}

Let $\nbigi\subset\gbigu_z(p,n)$ be any $\Gal(p)$-invariant subset.
In \S\ref{subsection;24.4.5.100},
for $\star=!,\ast$,
we shall introduce purely algebraically defined functors
\index{functors $\gbigf^{(0,\infty)}_{+,\star}$}
\[
 \gbigf^{(0,\infty)}_{+,\star}:
 \Loc^{\St}(\nbigi)
 \to
 \Loc^{\St}(\gbigf^{(0,\infty)}_+(\nbigi)),
 \quad
 (L,\vecnbigf)
 \mapsto
 \gbigf^{(0,\infty)}_{+,\star}(L,\vecnbigf)
 =(\gbigq^0_{\star}(L,\vecnbigf)_{\real},\vecnbigf),
\]
and morphisms of $2\pi\seisuu$-equivariant local systems
\begin{equation}
 \label{eq;24.4.19.20}
 L\lrarr
 \gbigq^0_{!}(L,\vecnbigf)_{\real}
 \lrarr
 \gbigq^0_{\ast}(L,\vecnbigf)_{\real}
 \lrarr
 L.
\end{equation}
Let us mention basic properties,
which are clear from the construction.
\begin{lem}
\mbox{{}}
\begin{itemize}
 \item 
The composition of the morphisms {\rm(\ref{eq;24.4.19.20})}
equals $\id-M^{-1}$,
where $M$ denotes the monodromy automorphism of $L$.
\item
We set $\omega^{\circ}=\frac{\omega}{1+\omega}$.
By the construction,
$\nbigt_{\omega^{\circ}}\bigl(
 \gbigf^{(0,\infty)}_{+,\star}(L,\vecnbigf)
 \bigr)$ $(\star=!,\ast)$
are identified with
$\nbigt_{\omega}(L)$,
and the following diagram is commutative.
\[
\begin{CD}
 \nbigt_{\omega}(L)
 @>{\id-M_0^{-1}}>>
 \nbigt_{\omega}(L)
 \\
 @V{=}VV @V{=}VV \\
 \nbigt_{\omega^{\circ}}\bigl(
 \gbigf^{(0,\infty)}_{+,!}(L,\vecnbigf)
 \bigr)
 @>>>
 \nbigt_{\omega^{\circ}}\bigl(
 \gbigf^{(0,\infty)}_{+,\ast}(L,\vecnbigf)
 \bigr).
\end{CD}
\]
Here, $M_0$ denotes the monodromy automorphism of
$\nbigt_{\omega}(L)$.
\end{itemize}
\end{lem}
 
We shall obtain the following theorem.
\begin{thm}[Theorem
\ref{thm;24.4.5.110}]
\label{thm;24.4.15.21}
For any basic meromorphic flat bundle
$(V,\nabla)$ of level $(0,\omega)$,
there exist the following commutative diagram
in $\Loc^{\St}(\gbigf^{(0,\infty)}_{+}(\nbigi_0(V)))$,
where the vertical arrows are isomorphisms:
 \begin{equation}
\begin{CD}
  \gbigf^{(0,\infty)}_{+,!}(L_0(V),\vecnbigf)
  @>{F_{\gbigq^0}}>>
  \gbigf^{(0,\infty)}_{+,\ast}(L_0(V),\vecnbigf)\\
  @V{\simeq}VV @V{\simeq}VV \\
  (\gbigl^{\gbigf}_!(V),\vecnbigf)
  @>>>
  (\gbigl^{\gbigf}_{\ast}(V),\vecnbigf).
\end{CD}
\end{equation}
\end{thm}

Let $\nbigl(V)$ denote the local system on
$\cnum^{\ast}$
associated with $(V,\nabla)$.
There exists the regular singular meromorphic flat bundle
$V^{\reg}$ on $(\proj^1,\{0,\infty\})$
corresponding to $\nbigl(V)$.
\index{meromorphic flat bundle $V^{\reg}$}
As explained in \S\ref{subsection;24.4.15.1},
there exist the natural morphisms of local systems
$\gbigl^{\gbigf}_!(V^{\reg})
\to
 \gbigl^{\gbigf}_!(V)$
and
$\gbigl^{\gbigf}_{\ast}(V)
\to
\gbigl^{\gbigf}_{\ast}(V^{\reg})$.

\begin{prop}
We also have the following commutative diagram of
the $2\pi\seisuu$-equivariant local systems
\[
\begin{CD}
 L_0(V) @>>> \gbigq^0_!(L_0(V),\vecnbigf)_{\real}
  @>{F_{\gbigq^0}}>> \gbigq^0_{\ast}(L_0(V),\vecnbigf)_{\real}
  @>>>
  L_0(V) \\
  @V{\simeq}VV @V{\simeq}VV @V{\simeq}VV @V{\simeq}VV \\
  \gbigl^{\gbigf}_!(V^{\reg})
  @>>>
  \gbigl^{\gbigf}_!(V)
  @>>>
  \gbigl^{\gbigf}_{\ast}(V)
  @>>>
  \gbigl^{\gbigf}_{\ast}(V^{\reg}).
\end{CD}  
\]
Here,
the lower horizontal arrows are
the natural morphisms.
\end{prop}

\subsubsection{Stokes shells}

To capture the property of
$(\gbigl^{\gbigf}_{\star}(V),\vecnbigf)$,
we shall introduce a useful invariant of local systems with Stokes structures
in \S\ref{section;18.6.3.11}
which is called Stokes shell.
For any $\Gal(p)$-invariant subset
$\nbigi$ of $\gbigu_z(p,n)$,
let $\Shcat(\nbigi)$ denote the category of Stokes shells
indexed by $\nbigi$.
There exists an equivalence
\[
 \Shsf:\Loc^{\St}(\nbigi)
 \to
 \Shcat(\nbigi).
\]
In \S\ref{subsection;18.5.7.10},
we shall introduce explicitly defined functors
\[
\gbigf^{(0,\infty)}_{+,\star}\Shsf:
\Loc^{\St}(\nbigi)
\to
\Shcat(\gbigf^{(0,\infty)}_{+}(\nbigi))
\quad(\star=!,\ast)
\]
with natural transformation
$F:\gbigf^{(0,\infty)}_{+,!}\Shsf\to
\gbigf^{(0,\infty)}_{+,\ast}\Shsf$.

\begin{prop}[Proposition
\ref{prop;24.3.25.60}]
\label{prop;24.4.15.60}
There exists the following commutative diagram:
\begin{equation}
\begin{CD}
 \gbigf^{(0,\infty)}_{+,!}\Shsf(L,\vecnbigf)
 @>{F}>>
 \gbigf^{(0,\infty)}_{+,\ast}\Shsf(L,\vecnbigf)
 \\
 @V{\simeq}VV @V{\simeq}VV \\
 \Shsf\bigl(
 \gbigf^{(0,\infty)}_{+,!}(L,\vecnbigf)
 \bigr)
 @>>>
 \Shsf\bigl(
 \gbigf^{(0,\infty)}_{+,\ast}(L,\vecnbigf)
 \bigr).
\end{CD}
\end{equation}
As a result,
$\Shsf(\gbigl^{\gbigf}_!(V),\vecnbigf)
 \to
 \Shsf(\gbigl^{\gbigf}_{\ast}(V),\vecnbigf)$
is identified with
$\gbigf^{(0,\infty)}_{+,!}\Shsf(L_0(V),\vecnbigf)
 \to
 \gbigf^{(0,\infty)}_{+,\ast}\Shsf(L_0(V),\vecnbigf)$.
\end{prop}

\subsection{Fourier transforms in the basic case of $(\infty,\omega)$}
Let $\omega=n/p$ for some positive integers $n$ and $p$.
We assume $\omega>1$, i.e., $n>p$.
Let $x=z^{-1}$.

\subsubsection{Transformation of the index sets}

As we will review in \S\ref{subsection;24.4.16.1},
by the stationary phase formula
\cite{Fang, Graham-Squire, Sabbah-stationary},
a $\Gal(n-p)$-invariant subset
$\gbigf^{(\infty,\infty)}_+(\nbigi)$
of $\gbigu_u(n-p,n)$ is attached
to any $\Gal(p)$-invariant subset $\nbigi$ of $\gbigu_x(p,n)$,
and the following holds.
\index{set $\gbigf^{(\infty,\infty)}_+(\nbigi)$}
\begin{itemize}
\item If $(V,\nabla)$ is basic of
      level $(\infty,\omega)$,
      then
      $\Fourier(V(\star 0))(\ast 0)$ $(\star=\ast,!)$
      are basic meromorphic flat bundles
      of level $(\infty,\frac{\omega}{\omega-1})$,
      and we have
\[
       \nbigi_{\infty}\bigl(
       \Fourier_+(V(\star 0))
       \bigr)
       =
       \gbigf^{(\infty,\infty)}_+\bigl(
       \nbigi_{\infty}(V)
      \bigr)
      \cup\{0\}.
\]
\end{itemize}

\subsubsection{Local Fourier transform of local systems with Stokes structure}

Let $\nbigi\subset\gbigu_x(p,n)$ be any $\Gal(p)$-invariant subset.
We set
$\nbigi^{\circ}=\gbigf^{(\infty,\infty)}_+(\nbigi)\cup\{0\}$.
In \S\ref{subsection;24.4.5.120},
for $\star=!,\ast$,
we shall introduce purely algebraically defined functors
\index{functors $\gbigf^{(\infty,\infty)}_{+,\star}$}
\[
 \gbigf^{(\infty,\infty)}_{+,\star}:
 \Loc^{\St}(\nbigi)
 \to
 \Loc^{\St}(\nbigi^{\circ}),
 \quad
 (L,\vecnbigf)
 \mapsto
 \gbigf^{(\infty,\infty)}_{+,\star}(L,\vecnbigf)
 =(\gbigq^{\infty}_{\star}(L,\vecnbigf)_{\real},\vecnbigf),
\]
and morphisms of $2\pi\seisuu$-equivariant local systems
\begin{equation}
\label{eq;24.4.19.30}
 c^{-1}(\nbigt_{\omega}(L))\lrarr
 \gbigq^{\infty}_{!}(L,\vecnbigf)_{\real}
 \lrarr
 \gbigq^{\infty}_{\ast}(L,\vecnbigf)_{\real}
 \lrarr
 c^{-1}(\nbigt_{\omega}(L)).
\end{equation}
Here $c:\real\to\real$ be the map defined by
$c(\theta)=-\theta$.
\begin{lem}
\mbox{{}}
\begin{itemize}
 \item The composition of the morphisms in {\rm(\ref{eq;24.4.19.30})}
       equal $\id-M_0$,
       where $M_0$ denotes the monodromy automorphism of
       $\nbigt_{\omega}(L)$.
 \item Set $\omega^{\circ}=(\omega-1)^{-1}\omega$.
       By the construction,
       $\nbigt_{\omega^{\circ}}
       \Bigl(
        \gbigf^{(\infty,\infty)}_{+,\star}(L,\vecnbigf)
       \Bigr)$
       are identified with
       $c^{-1}(L)$.
       The following diagram is commutative:
\begin{equation}
\begin{CD}
  c^{-1}(L)
  @>{\id-M}>>
  c^{-1}(L)\\
  @V{=}VV @V{=}VV \\
\nbigt_{\omega^{\circ}}
\Bigl(
\gbigf^{(\infty,\infty)}_{+,!}(L,\vecnbigf)
  \Bigr)
  @>>>
  \nbigt_{\omega^{\circ}}
\Bigl(
\gbigf^{(\infty,\infty)}_{+,\ast}(L,\vecnbigf)
  \Bigr).
\end{CD}
\end{equation}       
\end{itemize}
\end{lem}

We shall obtain the following theorem.
\begin{thm}[Theorem
\ref{thm;24.4.5.121}]
\label{thm;24.4.15.61}
Let $(V,\nabla)$ be a basic meromorphic flat bundle
of level $(\infty,\omega)$.
Then, there exist the following commutative diagram in
$\Loc^{\St}\bigl(
\gbigf^{(\infty,\infty)}_{+}(\nbigi_{\infty}(V))\cup\{0\}
 \bigr)$,
where the vertical arrows are isomorphisms:
 \begin{equation}
\begin{CD}
  \gbigf^{(\infty,\infty)}_{+,!}(L_{\infty}(V),\vecnbigf)
  @>{F_{\gbigq^{\infty}}}>>
  \gbigf^{(\infty,\infty)}_{+,\ast}(L_{\infty}(V),\vecnbigf)\\
  @V{\simeq}VV @V{\simeq}VV \\
  (\gbigl^{\gbigf}_!(V),\vecnbigf)
  @>>>
  (\gbigl^{\gbigf}_{\ast}(V),\vecnbigf).
\end{CD}
\end{equation}
\end{thm}

As explained in \S\ref{subsection;24.4.15.1},
there exist the natural morphisms of local systems
$\gbigl^{\gbigf}_!(\nbigt^{\infty}_{\omega}(V))
\to
 \gbigl^{\gbigf}_!(V)$
and
$\gbigl^{\gbigf}_{\ast}(V)
\to
\gbigl^{\gbigf}_{\ast}(\nbigt^{\infty}_{\omega}(V))$.

\begin{prop}
\label{prop;24.4.15.62}
We also have the following commutative diagram of
the $2\pi\seisuu$-equivariant local systems
by setting
 $L_1=\nbigt_{\omega}(L_{\infty}(V))$
 {\small
\[
\begin{CD}
 c^{-1}(L_1)
 @>>>
 \gbigq^{\infty}_!(L_{\infty}(V),\vecnbigf)_{\real}
 @>{F_{\gbigq^{\infty}}}>>
 \gbigq^{\infty}_{\ast}(L_{\infty}(V),\vecnbigf)_{\real}
  @>>>
  c^{-1}(L_1) \\
  @V{\simeq}VV @V{\simeq}VV @V{\simeq}VV @V{\simeq}VV \\
  \gbigl^{\gbigf}_!(\nbigt^{\infty}_{\omega}(V))
  @>>>
  \gbigl^{\gbigf}_!(V)
  @>>>
  \gbigl^{\gbigf}_{\ast}(V)
  @>>>
  \gbigl^{\gbigf}_{\ast}(\nbigt^{\infty}_{\omega}(V)).
\end{CD}  
\]}
Here,
the lower horizontal arrows are
the natural morphisms.
\end{prop}

\subsubsection{Stokes shells}

In \S\ref{subsection;20.11.14.20},
for any $\Gal(p)$-invariant subset
$\nbigi\subset\gbigu_x(p,n)$,
we shall introduce explicitly defined functors
\[
\gbigf^{(\infty,\infty)}_{+,\star}\Shsf:
\Loc^{\St}(\nbigi)
\to
\Shcat\bigl(\gbigf^{(\infty,\infty)}_+(\nbigi)\cup\{0\}\bigr)
\quad
(\star=!,\ast)
\]
with a natural transform
$\gbigf^{(\infty,\infty)}_{+,!}\Shsf
\to\gbigf^{(\infty,\infty)}_{+,\ast}\Shsf$.
\begin{prop}[Proposition
\ref{prop;24.3.26.31}]
For any $(L,\vecnbigf)\in\Loc^{\St}(\nbigi)$,
there exists the following commutative diagram:
\begin{equation}
\begin{CD}
 \gbigf^{(\infty,\infty)}_{+,!}\Shsf(L,\vecnbigf)
 @>{F}>>
 \gbigf^{(\infty,\infty)}_{+,\ast}\Shsf(L,\vecnbigf)
 \\
 @V{\simeq}VV @V{\simeq}VV \\
 \Shsf\bigl(
 \gbigf^{(\infty,\infty)}_{+,!}(L,\vecnbigf)
 \bigr)
 @>>>
 \Shsf\bigl(
 \gbigf^{(\infty,\infty)}_{+,\ast}(L,\vecnbigf)
 \bigr).
\end{CD}
 \end{equation}
As a result,
$\Shsf(\gbigl^{\gbigf}_!(V),\vecnbigf)
 \to
 \Shsf(\gbigl^{\gbigf}_{\ast}(V),\vecnbigf)$
is identified with
$\gbigf^{(\infty,\infty)}_{+,!}\Shsf(L_{\infty}(V),\vecnbigf)
 \to
 \gbigf^{(\infty,\infty)}_{+,\ast}\Shsf(L_{\infty}(V),\vecnbigf)$.
\end{prop}

\subsection{The reason to consider Stokes shells}
\label{subsection;24.4.3.20}

There are standard invariants
such as Stokes matrices, Stokes factors, etc.,
to capture the property of local systems with Stokes structure.
However,
the set of the Stokes directions and anti-Stokes directions
\[
 S\Bigl(
 \nbigi_{\infty}\bigl(\Fourier(\nbigv(\star 0))\bigr)
 \Bigr),
 \quad
 A\Bigl(
 \nbigi_{\infty}\bigl(\Fourier(\nbigv(\star 0))\bigr)
 \Bigr)
\]
are less directly related with
the sets
$S(\nbigi_0(\nbigv))$
and
$A(\nbigi_0(\nbigv))$,
or
the sets
$S(\nbigi_{\infty}(\nbigv))$
and
$A(\nbigi_{\infty}(\nbigv))$.
For example, let us consider the case
where $\nbigv$ is basic of level $(0,1)$
with $\nbigi_0(\nbigv)=\{\alpha_i z^{-1}\,|\,i=1,\ldots,m\}$
where $\alpha_i$ are mutually distinct non-zero complex numbers.
Then, $\nbigi_{\infty}\bigl(\Fourier(\nbigv(\star 0))\bigr)$
is
$\bigl\{
\pm 2\alpha_i^{1/2}u^{-1/2}
\bigr\}$.
The relation among the sets
\[
 S(2\alpha_i^{1/2}u^{-1/2},2\alpha_j^{1/2}u^{-1/2})=
\bigl\{
\theta\in\real\,\big|\,
\Re\bigl(
(\alpha_i^{1/2}-\alpha_j^{1/2})e^{-\sqrt{-1}\theta/2}
\bigr)=0\bigr\}
\]
for $1\leq i\neq j\leq m$
depend on the absolute values $|\alpha_k|$ $(k=1,\ldots,m)$.
If we use the classical invariants of Stokes structures,
for example Stokes matrices,
we need classifications 
depending on 
$S\Bigl(
 \nbigi_{\infty}\bigl(\Fourier(\nbigv(\star 0))\bigr)
 \Bigr)$
 or
$A\Bigl(
\nbigi_{\infty}\bigl(\Fourier(\nbigv(\star 0))\bigr)
\Bigr)$
which increases the complexity of the computation.
The results would be unnecessarily complicated.
It is one of the main reasons to consider
Stokes shells.

\section{Building blocks from the regular singularity}
\label{subsection;24.4.19.10}

The Fourier transform of regular singular holonomic $\nbigd$-modules
has been clearly understood.
There are two types of building blocks.

\subsection{Local systems}

Let $D\subset\cnum$ be a finite subset.
We set $\Dbar=D\cup\{\infty\}$.
Let $V$ be a regular singular meromorphic flat bundle
on $(\proj^1,\Dbar)$.
We obtain the dual meromorphic flat bundle
$V^{\lor}=\nhom_{\nbigo_{\proj^1}(\ast \Dbar)}(V,\nbigo_{\proj^1}(\ast \Dbar))$
on $(\proj^1,\Dbar)$,
and the $\nbigd_{\proj^1}(\ast\infty)$-module
$V(!D)=\DDD_{\proj^1}(V^{\lor})(\ast\infty)$.
\index{meromorphic flat bundle $V^{\lor}$}
\index{$\nbigd$-modules \mbox{$V(!D)$}}
For any map
$\varrho:D\to\{\ast,!\}$,
we set
$V(\varrho)=V(!D)\otimes\nbigo_{\proj^1}(\ast \varrho^{-1}(\ast))$.
\index{$\nbigd$-module $V(\varrho)$}
The $2\pi\seisuu$-equivariant local systems with Stokes structure
$(\gbigl^{\gbigf}(V(\varrho)),\vecnbigf)$
has been clearly understood
by the various previous studies.
Indeed,
$\Fourier_+(V(\varrho))(\ast 0)$
is basic of level $(\infty,1)$,
and we have
\[
 \nbigi_{\infty}\bigl(
 \Fourier_+(V(\varrho))
 \bigr)
 =\bigl\{
 \alpha u^{-1}\,|\,\alpha\in D\bigr\}=:\nbigi_D.
\]	
The Stokes structures of $(\gbigl^{\gbigf}(V(\varrho)),\vecnbigf)$
have been also described explicitly.

As a complement,
to describe the Stokes structure
in terms of the Stokes shells even in this case,
we shall introduce an explicit construction of
Stokes shells
$\gbigf_{\varrho}(\nbigl)$
from a local system $\nbigl$ on $\cnum\setminus D$
in \S\ref{subsection;18.5.15.10},
and we shall observe the following.
\begin{prop}[Proposition
\ref{prop;24.3.31.1}]
\label{prop;24.4.15.63}
Let $\nbigl(V)$ denote the local system on $\cnum\setminus D$
associated with $V$.
Then,
for any maps $\varrho:D\to\{!,\ast\}$,
there exist the natural isomorphisms
of Stokes shells
$\gbigf_{\varrho}(\nbigl(V))
 \simeq
 \Shsf(\gbigl^{\gbigf}_{\varrho}(V),\vecnbigf)$.
\end{prop}

\subsection{Regular singular monodromic $\nbigd$-modules}
\label{subsection;24.4.19.11}

Let $A$ be a finite dimensional vector space
equipped with an endomorphism $F$
such that any eigenvalue $\alpha$ of $F$
satisfies
either (i) $\alpha=0$,
or (ii) $\alpha\not\in\seisuu$,
$0\leq \Re(\alpha)<1$.
There exists the decomposition
\[
 (A,F)
 =(A^u,N)\oplus
 (A^{nu},F^{nu})
\]
such that
(i) $N$ is nilpotent,
(ii) any eigenvalue of $F^{nu}$ is not $0$.
Let $S(F^{nu})$ denote the set of
the eigenvalues of $F^{nu}$.

Let $\nbigv=A\otimes\nbigo_{\proj^1}(\ast\{0,\infty\})$
with the connection
$\nabla=d-F\frac{dz}{z}$.
Let $\nbigm$ be a regular holonomic $\nbigd_{\proj^1}$-modules
such that $\nbigm(\ast 0)=\nbigv$.
Corresponding to the decomposition
$(A,F)=(A^u,N)\oplus (A^{nu},F^{nu})$,
we have the decomposition
\[
 \nbigv=\nbigv^u\oplus \nbigv^{nu},\quad
 \nbigm=\nbigm^u\oplus\nbigm^{nu}.
\]
Moreover, $\nbigm^{nu}=\nbigv^{nu}$.

We obtain the $2\pi\seisuu$-equivariant local systems
$L_{0}(\nbigm)$
and
$L_{\infty}(\nbigm)$.
We have
$L_0(\nbigm)=c^{-1}(L_{\infty}(\nbigm))$,
where $c:\real\to\real$
is defined by
$c(\theta)=-\theta$.

\subsubsection{}

There exists the $V$-filtration of $\nbigm$
along $0$.
We set
\[
 \psitilde(\nbigm)
=\Gr^V_{-1} (\nbigm)
 \oplus
 \bigoplus_{\beta\in S(F^{nu})}
 \Gr^V_{\beta-1}(\nbigm),
\]
\[
 \phitilde(\nbigm)
=\Gr^V_{0} (\nbigm)
\oplus
\bigoplus_{\beta\in S(F^{nu})}
\Gr^V_{\beta}(\nbigm).
\]
\index{Nearby cycle functor $\psitilde$}
\index{Vanishing cycle functor $\phitilde$}
We obtain the morphism
$\cantilde_{\nbigm}:\psitilde(\nbigm)\to\phitilde(\nbigm)$
induced by
$-\del_z$.
We also obtain the morphism
$\vartilde_{\nbigm}:\phitilde(\nbigm)\to\psitilde(\nbigm)$
induced by $z$.

There exists a natural isomorphism
$\psitilde(\nbigm)
\simeq A$
under which
$-z\del_z$ is identified with $F$.
It is also standard that there exists a natural isomorphism
\[
\rhotilde_z:
\psitilde(\nbigm)
 \simeq
 H^0(\real,L_0(\nbigm))
\]
under which the monodromy automorphism of $L_0(\nbigm)$
equals
$\exp(2\pi\sqrt{-1}F)$.

\subsubsection{}

There exist the natural isomorphisms
(see \S\ref{subsection;25.3.18.1}):
\[
 \phitilde(\nbigm)\simeq
 \psitilde(\Fourier_{+}(\nbigm)),
 \quad\quad
 \psitilde(\nbigm)\simeq
 \phitilde(\Fourier_{+}(\nbigm)).
\]
We obtain the isomorphism
$\Psi_{\nbigm,\pm}:
\phitilde(\nbigm)
\simeq
H^0(\real,
\gbigl^{\gbigf}(\nbigm))$
as
the composition of the following isomorphisms:
\begin{multline}
 \phitilde(\nbigm)
 \simeq
 \psitilde(\Fourier_{+}(\nbigm))
 \simeq
 H^0\Bigl(\real,
 L_{0}\bigl(
 \Fourier_{+}(\nbigm)
 \bigr)
 \Bigr)
\\
 \simeq
 H^0\Bigl(
 \real,
 L_{\infty}\bigl(
 \Fourier_{+}(\nbigm)
 \bigr)
 \Bigr).
\end{multline}

\subsubsection{}

By using the rapid decay homology and
the moderate growth homology,
we obtain the following isomorphisms
\[
 \Abb^{\rd}_{+}:
 H^0(\real,L_0(\nbigv))
 \simeq
 H^0(\real,\gbigl^{\gbigf}_!(\nbigv)),
\]
\[
 \Abb^{\mg}_+:
 H^0(\real,L_0(\nbigv))
 \simeq
 H^0(\real,\gbigl^{\gbigf}_{\ast}(\nbigv)).
\]
(See \S\ref{subsection;25.3.10.40}.
See also
\S\ref{subsection;24.3.25.100} and \S\ref{subsection;24.4.19.2}.)

\subsubsection{}

Let $X=\realbar_{\geq 0}\times\real$
and $X^{\ast}=\real_{>0}\times\real$.
Let $\Gamma_{!}$ be a path connecting
$(\infty,-2\pi)$ and $(\infty,0)$
on $(X,X^{\ast})$.
We regard $X^{\ast}$ as a universal covering of $\cnum^{\ast}$
by the map
$(r,\theta)\longmapsto re^{\sqrt{-1}\theta}$.

Under the isomorphism
$\psitilde(\nbigm)\simeq A$,
we define the endomorphisms
$\Phi_!$ and $\Phi_{\ast}$
of $\psitilde(\nbigm)$
by 
\[
\Phi_!
=\frac{-1}{2\pi\sqrt{-1}}
\int_{\Gamma_{!}}
\exp(F\log\zeta)e^{-\zeta}
\frac{d\zeta}{\zeta},
\]
\[
\Phi_{\ast}
 =\frac{-1}{2\pi\sqrt{-1}}
 \int_{0}^{\infty}
 \exp(F\log t)e^{-t}dt.
\]
\begin{prop}[Proposition
\ref{prop;25.3.11.200}]
The endomorphisms
$\Phi_{\star}$ $(\star=!,\ast)$  are invertible.
Moreover,
the following diagrams are commutative:
\begin{equation}
\begin{CD}
 \psitilde(\nbigm)
 @>{\cantilde_{\nbigm}\circ \Phi_{!}}>>
 \phitilde(\nbigm)
 @>{(\Phi_{\ast})^{-1}\circ\vartilde_{\nbigm}}>>
 \psitilde(\nbigm)
 \\
 @V{\simeq}V{\Abb^{\rd}_{+}\circ\rhotilde_z}V
 @V{\simeq}V{\Psi_{\nbigm,+}}V
 @V{\simeq}V{\Abb^{\mg}_{+}\circ\rhotilde_z}V \\
 H^0(\real,\gbigl^{\gbigf}(\nbigv(!0)))
 @>>>
 H^0(\real,\gbigl^{\gbigf}(\nbigm))
 @>>>
 H^0(\real,\gbigl^{\gbigf}(\nbigv)).
\end{CD}
\end{equation}
Here, the lower horizontal arrows are the natural morphisms.
\end{prop}

\section{Extensions of local systems with Stokes structure}
\label{section;20.11.4.32}

Let us explain the notion of extension of
local systems with Stokes structure
which is useful in the inductive study of
$(\gbigl^{\gbigf}(\nbigm),\vecnbigf)$
for general holonomic $\nbigd$-module $\nbigm$ on $\proj^1$.
\index{extension}

\subsection{Simple case}
\label{subsection;21.4.29.1}

Let $\nbigi$ be a $\Gal(p)$-invariant finite subset of
$z_p^{-1}\cnum[z_p^{-1}]$.
Let $f:(L_1,\vecnbigf)\lrarr (L_2,\vecnbigf)$
be a morphism in $\Loc^{\St}(\nbigi)$
such  that the induced morphisms
$\Gr^{\vecnbigf}_{\gminib}(f)$ are isomorphisms
unless $\gminib=0$.
We obtain the induced morphism
$\Gr^{\vecnbigf}_0(f):
\Gr^{\vecnbigf}_0(L_1)
\lrarr
\Gr^{\vecnbigf}_0(L_2)$
in $\Loc^{\St}(0)$.

Let $\nbigc_1$ be the category of morphisms
\[
 (L_1,\vecnbigf)
\stackrel{a_1}{\lrarr} (M,\vecnbigf)
\stackrel{a_2}{\lrarr}
 (L_2,\vecnbigf)
\]
in $\Loc^{\St}(\nbigi)$
such that
(i) $a_2\circ a_1=f$,
(ii) $\Gr^{\vecnbigf}_{\gminia}(a_i)$ are isomorphisms
unless $\gminia=0$.
A morphism in $\nbigc_1$
is a commutative diagram in $\Loc^{\St}(\nbigi)$:
\[
 \begin{CD}
  (L_1,\vecnbigf)
  @>>>
  (M_1,\vecnbigf)
  @>>>
  (L_2,\vecnbigf)\\
  @V{\id}VV @VVV @V{\id}VV\\
  (L_1,\vecnbigf)
  @>>>
  (M_2,\vecnbigf)
  @>>>
  (L_2,\vecnbigf).\\
 \end{CD}
\]

Let $\nbigc_2$ be the category of
morphisms of
$2\pi\seisuu$-equivariant local systems
\[
 \Gr^{\vecnbigf}_0(L_1)
 \stackrel{b_1}{\lrarr}
 N
 \stackrel{b_2}{\lrarr}
 \Gr^{\vecnbigf}_0(L_2)
\]
 such that
 $b_2\circ b_1=\Gr^{\vecnbigf}_0(f)$.
 A morphism in $\nbigc_2$
 is a commutative diagram of
 $2\pi\seisuu$-equivariant local systems:
\[
 \begin{CD}
  \Gr^{\vecnbigf}_0(L_1)
  @>>>
  N_1
  @>>>
  \Gr^{\vecnbigf}_0(L_2)\\
  @V{\id}VV @VVV @V{\id}VV\\
  \Gr^{\vecnbigf}_0(L_1)
  @>>>
  N_2
  @>>>
  \Gr^{\vecnbigf}_0(L_2).
 \end{CD}
\]

Any object
$(L_1,\vecnbigf)
\lrarr (M,\vecnbigf)
\lrarr
(L_2,\vecnbigf)$
in $\nbigc_1$
induces an object
$\Gr^{\vecnbigf}_0(L_1)\lrarr
 \Gr^{\vecnbigf}_0(M)
\lrarr
 \Gr^{\vecnbigf}_0(L_2)$
in $\nbigc_2$.
Thus, we obtain a functor
$\nbigc_1\lrarr\nbigc_2$.
The following proposition is a special case of
Theorem \ref{thm;20.11.3.20}.
\begin{prop}
\label{prop;20.11.3.10}
 The functor 
 $\nbigc_1\lrarr\nbigc_2$
is an equivalence.
\end{prop}

By Proposition \ref{prop;20.11.3.10},
for any object
$\nbiga_0=\bigl(
\Gr^{\vecnbigf}_0(L_1)
\stackrel{b_1}{\lrarr}
N
\stackrel{b_2}{\lrarr}
\Gr^{\vecnbigf}_0(L_2)
\bigr)$
in $\nbigc_2$,
there exists
$\nbiga_1=\bigl(
(L_1,\vecnbigf)
\stackrel{a_1}{\lrarr}
(M,\vecnbigf)
\stackrel{a_2}{\lrarr}
(L_2,\vecnbigf)
\bigr)$
which induces $\nbiga_0$.
It is uniquely determined
by $f$ and $\nbiga_0$,
up to canonical isomorphisms.
Such $(M,\vecnbigf)$
is called an extension of
$f:(L_1,\vecnbigf)
\lrarr (L_2,\vecnbigf)$
by $\nbiga_0$.
\index{extension}
Let us emphasize that
we can explicitly construct
$(M,\vecnbigf)$
in an elementary way
by using canonical splittings.
(See \S\ref{subsection;21.5.1.1}.)

\subsection{Some categories}
\label{subsection;18.6.18.3}

To explain another case,
we introduce some convenient categories.
Let $\Dsf_1$ be the category given as follows.
\index{category $\Dsf_1$}
\begin{itemize}
\item
Objects of $\Dsf_1$ are $!$ and $\ast$.
\item
$\Hom_{\Dsf_1}(\star,\star)$ $(\star=!,\ast)$
consists of the identity morphism,
which we shall denote by $f_{\star,\star}$.
$\Hom_{\Dsf_1}(!,\ast)$ consists of 
a unique morphism $f_{\ast,!}$.
$\Hom_{\Dsf_1}(\ast,!)$ is empty.
\end{itemize}
Let $\Csf_1$ denote the category given as follows.
\index{category $\Csf_1$}
\begin{itemize}
\item Objects of $\Csf_1$ are $!$, $\circ$ and $\ast$.
\item
$\Hom_{\Csf_1}(\star,\star)$
$(\star=!,\circ,\ast)$
consist of the identity morphisms,
which we denote by $f_{\star,\star}$.
$\Hom_{\Csf_1}(\star_1,\star_2)$
consists of
a unique map $f_{\star_2,\star_1}$
if $(\star_1,\star_2)=(!,\circ),(\circ,\ast),(!,\ast)$.
Otherwise,
$\Hom_{\Csf_1}(\star_1,\star_2)$ is empty.
\end{itemize}
For any set $S$,
let $\Dsf(S)$ denote the category of maps
$S\lrarr \{!,\ast\}$,
where
\[
 \Hom_{\Dsf(S)}(\varrho_1,\varrho_2)
:=\prod_{\gminia\in S}
 \Hom_{\Dsf_1}(\varrho_1(\gminia),\varrho_2(\gminia)).
\]
\index{category $\Dsf(S)$}
Similarly, 
let $\Csf(S)$ denote the category of maps
$S\lrarr \{!,\circ,\ast\}$,
where
\[
 \Hom_{\Csf(S)}(\varrho_1,\varrho_2)
:=\prod_{\gminia\in S}
 \Hom_{\Csf_1}(\varrho_1(\gminia),\varrho_2(\gminia)).
\]
\index{category $\Csf(S)$}
For $\star=!,\circ,\ast$,
let $\underline{\star}\in \Csf(S)$ denote the object
such that $\underline{\star}(\gminia)=\star$
for any $\gminia\in S$.
If $\star=!,\ast$,
we may also naturally regard
$\underline{\star}$ as objects in $\Dsf(S)$.
\index{objects \mbox{$\underline{!}$},
$\underline{\circ}$, $\underline{\ast}$}

Let $\iota:\Dsf(S)\lrarr \Csf(S)$
denote the naturally defined functor.
For any functor $F$ from $\Csf(S)$ to a category,
let $\iota^{\ast}(F)$ denote the induced functor
from $\Dsf(S)$ to the category.
\index{functor $\iota:\Dsf(S)\lrarr\Csf(S)$}

\subsection{Another simple case}
\label{subsection;21.4.29.2}

Let $\nbigi$ be a $\Gal(p)$-invariant finite subset of
$z_p^{-1}\cnum[z_p^{-1}]$.
We set $\nbigi_0=\nbigi\cap \cnum z^{-1}$.
Let $\nbige$ be a functor
$\Dsf(\nbigi_0)\lrarr \Loc^{\St}(\nbigi)$
satisfying the following condition.
\begin{itemize}
 \item For any morphism $\varrho_1\to \varrho_2$
       in $\Dsf(\nbigi_0)$,
       the induced morphism
       $\Gr^{\vecnbigf}_{\gminia}(\nbige(\varrho_1))
       \lrarr
       \Gr^{\vecnbigf}_{\gminia}(\nbige(\varrho_2))$
       is an isomorphism
       unless
       $\gminia\in\nbigi_0$
       and $\varrho_1(\gminia)\neq\varrho_2(\gminia)$.
\end{itemize}

Let $\nbigc_1$ be the category of functors
$\nbigetilde:
\Csf(\nbigi_0)\lrarr \Loc^{\St}(\nbigi)$
equipped with an isomorphism
$a_{\nbigetilde}:\iota^{\ast}(\nbigetilde)\simeq\nbige$
satisfying the following condition.
\begin{itemize}
 \item For any morphism $\varrho_1\to \varrho_2$
       in $\Csf(\nbigi_0)$,
       the induced morphism
       $\Gr^{\vecnbigf}_{\gminia}(\nbigetilde(\varrho_1))
       \lrarr
       \Gr^{\vecnbigf}_{\gminia}(\nbigetilde(\varrho_2))$
       is an isomorphism
       unless
       $\gminia\in\nbigi_0$
       and $\varrho_1(\gminia)\neq\varrho_2(\gminia)$.
\end{itemize}
A morphism
$f:(\nbigetilde_1,a_{\nbigetilde_1})
\lrarr
(\nbigetilde_2,a_{\nbigetilde_2})$ in $\nbigc_1$
is defined to be
a natural transformation
$f:\nbigetilde_1\lrarr\nbigetilde_2$
such that
$a_{\nbigetilde_2}\circ \iota^{\ast}(f)
=a_{\nbigetilde_1}$.

Let $\nbigc_2$
be the category of functors
$\nbigg=\bigoplus_{\gminia\in\nbigi_0}\nbigg_{\gminia}$
from $\Csf(\nbigi_0)$
to the category of 
$2\pi\seisuu$-equivariant
$\nbigi$-graded local systems
$\nbigg$
equipped with an isomorphism
$b_{\nbigg}:
\iota^{\ast}\nbigg\simeq
\bigoplus_{\gminia\in\nbigi_0}
\Gr_{\gminia}^{\vecnbigf}(\nbige)$
satisfying the following.
\begin{itemize}
 \item For any morphism $\varrho_1\to \varrho_2$
       in $\Csf(\nbigi_0)$,
       the induced morphism
       $\nbigg_{\gminia}(\varrho_1)
       \lrarr
       \nbigg_{\gminia}(\varrho_2)$
       is an isomorphism
       unless
       $\varrho_1(\gminia)\neq\varrho_2(\gminia)$.
\end{itemize}
A morphism
$f:(\nbigg_1,b_{\nbigg_1})\lrarr(\nbigg_2,b_{\nbigg_2})$
in $\nbigc_2$ is defined to be
a natural transformation
$f:\nbigg_1\lrarr\nbigg_2$
such that $b_{\nbigg_2}\circ f=b_{\nbigg_1}$.

Any object $\nbigetilde$ of $\nbigc_1$
induces an object
$\bigoplus_{\gminia\in\nbigi_0}
\Gr_{\gminia}^{\vecnbigf}(\nbigetilde)$ in $\nbigc_2$.
Thus, we obtain a functor
$\nbigc_1\lrarr\nbigc_2$.
The following proposition is a special case of
Theorem \ref{thm;20.11.3.20}.

\begin{prop}
\label{prop;21.4.28.20}
The functor $\nbigc_1\lrarr\nbigc_2$
is an equivalence.
\end{prop}

Proposition \ref{prop;21.4.28.20}
implies that for any $\nbigg\in\nbigc_2$,
there exists $\nbigetilde$ in $\nbigc_1$
which induces $\nbigg$.
Such $\nbigetilde$ is uniquely determined
up to canonical isomorphisms,
and called
an extension of $\nbige$ by $\nbigg$.
\index{extension}
Let us again emphasize that
$\nbigetilde$ is explicitly constructed in an elementary way.
See the proof of Theorem \ref{thm;20.11.3.20}
in \S\ref{subsection;21.4.30.2}.

\section{Reductions}
\label{subsection;24.4.15.101}

For any finite subset $D\subset\cnum$,
let $\Hol(\proj^1,D,\infty)$
denote the category of holonomic $\nbigd_{\proj^1}$-modules $\nbigm$
such that $\nbigm(\ast D)$ is a meromorphic flat bundle
on $(\proj^1,D\cup\{\infty\})$.

\subsection{Reductions at $0$}
\label{subsection;24.4.15.70}

Let $\nbigm\in\Hol(\proj^1,0,\infty)$
which is regular singular at $\infty$.
We obtain the meromorphic flat bundle
$\nbigv=\nbigm(\ast 0)$.
We set $\omega=-\ord\nbigi(\nbigv)$.
We set $(L,\vecnbigftilde)=(L_0(\nbigv),\vecnbigf)$.

We obtain the basic meromorphic flat bundle $V$ of level $(0,\omega)$
corresponding to $\nbigs_{\omega}(L,\vecnbigf)$.
Let $\nbigt_{\omega}(V)$ denote the regular singular
meromorphic flat bundle on $(\proj^1,\{0,\infty\})$
corresponding to
$\nbigt_{\omega}\nbigs_{\omega}(L)\in\Loc^{\St}(0)$.

We also obtain the holonomic $\nbigd_{\proj^1}$-module
$\nbigt_{\omega}(\nbigm)\in\Hol(\proj^1,0,\infty)$
characterized by the following conditions.
(See \S\ref{subsection;24.4.15.11}.)
\begin{itemize}
 \item $(L_0(\nbigt_{\omega}(\nbigm)),\vecnbigf)
       =\nbigt_{\omega}(L,\vecnbigf)$.
 \item The standard morphisms
       $\psi_{z}(\nbigt_{\omega}(\nbigm))
       \to
       \phi_{z}(\nbigt_{\omega}(\nbigm))
       \to
       \psi_{z}(\nbigt_{\omega}(\nbigm))$
       are identified with
       $\psi_{z}(\nbigm)
       \to
       \phi_{z}(\nbigm)
       \to
       \psi_{z}(\nbigm)$.
\end{itemize}

We set $\omega^{\circ}=(1+\omega)^{-1}\omega$.
We shall prove the following theorem.

\begin{thm}[Theorem
\ref{thm;24.3.24.1},
Corollary \ref{cor;24.4.15.10}]\label{thm;24.4.15.20}
There exists the functorial isomorphism
\[
 \nbigt_{\omega^{\circ}}(\gbigl^{\gbigf}(\nbigm),\vecnbigf)
 \simeq
 (\gbigl^{\gbigf}(\nbigt_{\omega}\nbigm),\vecnbigf).
\]
\end{thm}

Theorem \ref{thm;24.4.15.20}
particularly implies that 
$\nbigt_{\omega^{\circ}}(\gbigl^{\gbigf}_{\star}(V))$ $(\star=!,\ast)$
are naturally identified with
$\gbigl^{\gbigf}_{\star}(\nbigt_{\omega}(V))$
though it also directly follows from the stationary phase formula.
As in \S\ref{subsection;24.4.15.1},
there exists the following natural morphisms of
the $2\pi\seisuu$-equivariant local systems:
\begin{equation}
\label{eq;24.4.15.55}
 \gbigl^{\gbigf}_!(\nbigt_{\omega}(V))
\lrarr
 \gbigl^{\gbigf}(\nbigt_{\omega}(\nbigm))
\lrarr
 \gbigl^{\gbigf}_{\ast}(\nbigt_{\omega}(V)).
\end{equation}
We shall prove the following theorem.
\begin{thm}[Theorem
\ref{thm;24.3.25.50},
Corollary \ref{cor;24.4.15.10}]
The $2\pi\seisuu$-equivariant local system with Stokes structure
$\nbigs_{\omega^{\circ}}\bigl(
\gbigl^{\gbigf}(\nbigm),\vecnbigf
 \bigr)$
is obtained as the extension of
$(\gbigl^{\gbigf}_!(V),\vecnbigf)
 \to
 (\gbigl^{\gbigf}_{\ast}(V),\vecnbigf)$
by {\rm(\ref{eq;24.4.15.55})}.
\end{thm}

Together with Theorem \ref{thm;24.4.15.21}
and Proposition \ref{prop;24.4.15.60},
these theorems provide us with
an inductive procedure to study
$(\gbigl^{\gbigf}(\nbigm),\vecnbigf)$.
(See \S\ref{subsection;24.4.2.110},
where the method is explained in the case $\nbigm=\nbigv[\star 0]$.)

\subsection{Reductions at finite place}
\label{subsection;24.4.15.71}

Let $\nbigm\in\Hol(\proj^1,D,\infty)$
which is regular singular at $\infty$.
Let $V$ denote the regular singular meromorphic flat bundle
on $(\proj^1,D\cup\{\infty\})$
associated with the local system $\nbigl(\nbigm)$ on $\cnum\setminus D$.

Let $U_{\alpha}$ be a neighbourhood of $\alpha$ in $\proj^1$.
Let $\nbigm_{\alpha}\in\Hol(\proj^1,\alpha,0)$
be the $\nbigd$-modules
such that $\nbigm_{\alpha|U_{\alpha}}$
is isomorphic to $\nbigm_{|U_{\alpha}}$,
and regular singular at $\infty$.
Let $V_{\alpha}$ denote the regular singular meromorphic 
flat bundles on $(\proj^1,\{\alpha,\infty\})$
obtained from $V$ in the same way.

We shall prove the following proposition.
(See Proposition \ref{prop;24.3.30.1} and Corollary \ref{cor;24.4.15.13}.)
\begin{prop}
\label{prop;24.4.15.30}
There exist the functorial isomorphisms
\[
 \Gr^{\vecnbigf^{(1)}}_{\alpha u^{-1}}
 (\gbigl^{\gbigf}(\nbigm),\vecnbigf)
 \simeq
 (\gbigl^{\gbigf}(\nbigm_{\alpha}),\vecnbigf).
\]
(See {\rm\S\ref{subsection;21.4.24.13}}
for $\vecnbigf^{(1)}$.)
\end{prop}

Proposition \ref{prop;24.4.15.30}
particularly implies that
$\Gr^{\vecnbigf^{(1)}}_{\alpha u^{-1}}
(\gbigl^{\gbigf}(V(\varrho)),\vecnbigf$)
are naturally identified with
$\Gr^{\vecnbigf^{(1)}}_{\alpha u^{-1}}
(\gbigl^{\gbigf}_{\varrho(\alpha)}(V_{\alpha}),\vecnbigf)$
though it directly follows from the stationary phase formula.
As in \S\ref{subsection;24.4.15.1},
there exists the following natural morphisms of
the $2\pi\seisuu$-equivariant local systems:
\begin{equation}
\label{eq;24.4.15.56}
 \gbigl^{\gbigf}_{!}(V_{\alpha})
\lrarr
 \gbigl^{\gbigf}(\nbigm_{\alpha})
\lrarr
 \gbigl^{\gbigf}_{\ast}(V_{\alpha}).
\end{equation}
We shall prove the following.
(See Proposition \ref{prop;24.4.14.50} and Corollary \ref{cor;24.4.15.13}.)
\begin{prop}
The $2\pi\seisuu$-equivariant local system with Stokes structure
$(\gbigl^{\gbigf}(\nbigm),\vecnbigf^{(1)})$
is the extension of 
 $(\gbigl^{\gbigf}_{\varrho}(V),\vecnbigf)$
$(\varrho\in \Dsf(D))$
by {\rm(\ref{eq;24.4.15.56})}.
\end{prop}

We can study
$(\gbigl^{\gbigf}_{\varrho}(V),\vecnbigf)$
by using this proposition,
Proposition \ref{prop;24.4.15.63}
and the results in \S\ref{subsection;24.4.15.70}.

\subsection{Reduction at infinity}

Let $\nbigm\in\Hol(\proj^1,D,\infty)$.
For any $\omega\in\rnum_{>0}$,
there exists
$\nbigstilde^{\infty}_{\omega}(\nbigm)\in\Hol(\proj^1,D,\infty)$
characterized by the following condition.
\index{$\nbigd$-module $\nbigstilde^{\infty}_{\omega}(\nbigm)$}
\begin{itemize}
 \item $\nbigstilde^{\infty}_{\omega}(\nbigm)_{|\cnum}
       =\nbigm_{|\cnum}$.
 \item $(L_{\infty}(\nbigstilde^{\infty}_{\omega}(\nbigm)),\vecnbigf)
       =\nbigstilde_{\omega}
       \bigl(L_{\infty}(\nbigm),\vecnbigf
       \bigr)$.
\end{itemize}
We obtain the following proposition.
(See Proposition \ref{prop;24.3.17.121}
and Proposition \ref{prop;24.4.15.40}.)
\begin{prop}
\label{prop;24.4.15.50}
There exists the following isomorphism
\[
 \bigl(
 \gbigl^{\gbigf}(\nbigstilde^{\infty}_1(\nbigm)),\vecnbigf
 \bigr)
 \simeq
  \bigl(
 \gbigl^{\gbigf}(\nbigm),\vecnbigf
 \bigr).
\]
\end{prop}

We set
$\omega=
\min\bigl\{
 -\ord(\gminia)\,|\,\gminia\in\nbigi_{\infty}(\nbigm)
 \bigr\}$.
By Proposition \ref{prop;24.4.15.50},
it is enough to study the case $\omega>1$.
We obtain the basic meromorphic flat bundle
$V_{\infty}=\nbigttilde^{\infty}_{\omega}(\nbigm)$
of level $(\infty,\omega)$
characterized by the following condition.
\begin{itemize}
 \item $V_{\infty}$ is regular singular at $0$,
       and
       $(L_{\infty}(V_{\infty}),\vecnbigf)$
       is isomorphic to
       $\nbigttilde_{\omega}(L_{\infty}(\nbigm),\vecnbigf)$.
\end{itemize}
We also obtain the regular singular meromorphic flat bundle
$V^{\reg}_{\infty}$ corresponding to $\nbigttilde_{\omega}(L)$.
Note that
$V^{\reg}_{\infty}=\nbigstilde^{\infty}_{\omega}(V)$.

We put $\omega^{\circ}=(\omega-1)^{-1}\omega$.
We shall prove the following theorem.
(See Theorem \ref{thm;24.3.29.40} and
Corollary \ref{cor;24.4.15.51}.)
\begin{thm}
\label{thm;24.4.15.52}
There exists the natural isomorphism
\[
 \nbigt_{\omega^{\circ}}\bigl(
 \gbigl^{\gbigf}(\nbigm),\vecnbigf
 \bigr)
\simeq
 \bigl(
 \gbigl^{\gbigf}(\nbigstilde^{\infty}_{\omega}(\nbigm)),
 \vecnbigf
 \bigr).
\]
\end{thm}

Theorem \ref{thm;24.4.15.52} particularly implies that
$\nbigt_{\omega^{\circ}}(\gbigl^{\gbigf}_{\star}(V_{\infty}))$
$(\star=!,\ast)$
are naturally identified with
$\gbigl^{\gbigf}_{\star}(V_{\infty}^{\reg})$
though it directly follows from the stationary phase formula.
As in \S\ref{subsection;24.4.15.1},
there exists the following natural morphisms of
the $2\pi\seisuu$-equivariant local systems:
\begin{equation}
\label{eq;24.4.15.57}
 \gbigl^{\gbigf}_!(V^{\reg}_{\infty})
 \to
 \gbigl^{\gbigf}(\nbigstilde^{\infty}_{\omega}(\nbigm))
 \to
 \gbigl^{\gbigf}_{\ast}(V^{\reg}_{\infty}).
\end{equation}
We shall prove the following theorem.
(See Theorem \ref{thm;24.3.29.41}
and Corollary \ref{cor;24.4.15.51}.)
\begin{thm}
\label{thm;24.4.20.10}
The $2\pi\seisuu$-equivariant local systems
with Stokes structure
$\nbigs_{\omega^{\circ}}\bigl(
\gbigl^{\gbigf}(\nbigm),\vecnbigf
 \bigr)$
is obtained as the extension of
$(\gbigl^{\gbigf}_{!}(V_{\infty}),\vecnbigf)
 \to
(\gbigl^{\gbigf}_{\ast}(V_{\infty}),\vecnbigf)$
by {\rm(\ref{eq;24.4.15.57})}.
\end{thm}

We can study
$(\gbigl^{\gbigf}(\nbigm),\vecnbigf)$
by using these theorems,
Theorem \ref{thm;24.4.15.61},
Proposition \ref{prop;24.4.15.62}
and the results in
\S\ref{subsection;24.4.15.70}
and \S\ref{subsection;24.4.15.71}.
(See also \S\ref{subsection;24.4.2.130},
where we explain the method
for $\nbigm=\nbigv(\varrho)$.)

\subsection{Remark about the connecting morphisms}

In each case,
we also obtain the connecting morphisms:
\[
 c^{-1}L_{\infty}(\nbigttilde^{\infty}_1(\nbigm))
 \stackrel{q_1}{\lrarr}
 \gbigl^{\gbigf}(\nbigm)
 \stackrel{q_1}{\lrarr}
 c^{-1}L_{\infty}(\nbigttilde^{\infty}_1(\nbigm)).
\]
We also obtain
the $\nbigd$-module
$\Fourier_+(\nbigstilde^{\infty}_1(\nbigm))$.
There exists the isomorphism
$\gbigl^{\gbigf}(\nbigm)
\simeq
\gbigl^{\gbigf}(\nbigstilde^{\infty}_1\nbigm)$.
Because $\Fourier_+(\nbigstilde^{\infty}_1\nbigm)_{|\cnum^{\ast}}$
is a flat bundle,
we have
$\gbigl^{\gbigf}(\nbigm)
=c^{-1}L_{0}(\Fourier_+(\nbigstilde^{\infty}_1\nbigm))$.
We obtain the isomorphism
\[
 a:\psitilde(\Fourier_+(\nbigstilde^{\infty}_1\nbigm))
 \simeq H^0(\real,c^{-1}\gbigl^{\gbigf}(\nbigm)).
\]

The connecting morphisms
are related as the canonical morphisms between
the vanishing cycle 
and the nearby cycle up to twists.
\begin{prop}[Proposition
\ref{prop;25.3.18.20}]
We set
$(\nbigstilde^{\infty}_1\nbigm)^{\gbigf}
=\Fourier_+(\nbigstilde^{\infty}_1\nbigm)$.
There exists the following commutative diagram:
 \[
 \begin{CD}
  \psitilde((\nbigstilde^{\infty}_1\nbigm)^{\gbigf})
  @>{\cantilde\circ \Phi'_{!,-}}>>
  \phitilde((\nbigstilde^{\infty}_1\nbigm)^{\gbigf})
  @>{(\Phi'_{\ast,-})^{-1}\circ\vartilde}>>
  \psitilde((\nbigstilde^{\infty}_1\nbigm)^{\gbigf})\\
  @V{\simeq}V{a}V @V{\simeq}VV @V{\simeq}V{a}V \\
  H^0(\real,c^{-1}\gbigl^{\gbigf}(\nbigm))
  @>{q_2}>>
  H^0(\real,L_{\infty}(\nbigttilde_1(\nbigm)))
  @>{q_1}>>
  H^0(\real,c^{-1}\gbigl^{\gbigf}(\nbigm)).
 \end{CD}
\]
(See {\rm\S\ref{subsection;25.3.12.41}}
and {\rm(\ref{eq;25.3.18.10})}
 for the automorphisms
$\Phi'_{\star,-}$.)
\end{prop}

We can compute the following morphisms
from $\LS^{\fin}(\nbigm)$.
\begin{equation}
\label{eq;25.3.18.20}
 c^{-1}L_{\infty}(\nbigm)
 \lrarr
 \gbigl^{\gbigf}(\nbigstilde^{\infty}_{\infty}\nbigm)
 \lrarr
  c^{-1}L_{\infty}(\nbigm).
\end{equation}
Conversely, we can recover $\LS^{\fin}(\nbigm)$
from (\ref{eq;25.3.18.20}).
See \S\ref{subsection;25.3.18.30}.

\section{Outline of the computation}
\label{section;21.4.28.3}

Let us explain an outline of the computation of
$(\gbigl^{\gbigf}(\nbigm),\vecnbigf)$
for $\nbigm\in\Hol(\proj^1,D,\infty)$
from the data in \S\ref{subsection;21.4.28.2}
by using the results in \S\ref{subsection;24.4.15.100},
\S\ref{subsection;24.4.19.10}
and \S\ref{subsection;24.4.15.101}.

\subsection{The associated $2\pi\seisuu$-equivariant
  local systems with Stokes structure}

We obtain the rational numbers
$1<\omega(\infty,1)<\cdots<\omega(\infty,\ell(\infty))$
determined by
\[
 \{\omega(\infty,j)\}=\rnum_{>1}\cap
 \bigl\{
 -\ord(\gminia)\,\big|\,
 \gminia\in \nbigi_{\infty}(\nbigm)\setminus\{0\}
 \bigr\}.
\]
We set $\upsilon(\infty,j):=
\omega(\infty,j)\bigl(\omega(\infty,j)-1\bigr)^{-1}$.
Note that
$\upsilon(\infty,1)>\upsilon(\infty,2)>\cdots
>\upsilon(\infty,\ell(\infty))>1$.

For $\alpha\in D$,
we obtain the rational numbers
$\omega(\alpha,1)>\omega(\alpha,2)>\cdots
>\omega(\alpha,\ell(\alpha))=0$
determined by
\[
\bigl\{\omega(\alpha,j)
\bigr\}
=\bigl\{
-\ord(\gminia)\,\big|\,
\gminia\in\nbigi_{\alpha}(\nbigv)\setminus\{0\}
\bigr\}\cup\{0\}.
\]
We set
$\upsilon(\alpha,j):=
\omega(\alpha,j)(\omega(\alpha,j)+1)^{-1}$.
Note that
$1>\upsilon(\alpha,1)>\upsilon(\alpha,2)>\cdots
>\upsilon(\alpha,\ell(\alpha))=0$.

For $\alpha\in\Dbar$ and
$\omega=\omega(\alpha,j)$,
we obtain the following
$2\pi\seisuu$-equivariant
local systems with Stokes structure
\[
 (L_{\alpha,\omega}(\nbigm),\vecnbigf):=
 \nbigs_{\omega}
 \nbigttilde_{\omega}
(L_{\alpha}(\nbigm),\vecnbigf).
\]

\subsection{Building blocks}
\label{subsection;21.4.28.4}
  
\subsubsection{}

For $\omega=\omega(\infty,j)$,
we obtain the morphism of $2\pi\seisuu$-equivariant
local systems with Stokes structure
\begin{equation}
\label{eq;21.4.27.32}
 \gbigf^{(\infty,\infty)}_{+!}\bigl(
 L_{\infty,\omega}(\nbigm),\vecnbigf
 \bigr)
 \lrarr
 \gbigf^{(\infty,\infty)}_{+\ast}\bigl(
 L_{\infty,\omega}(\nbigm),\vecnbigf
 \bigr).
\end{equation}
Set $\upsilon=\upsilon(\infty,j)$,
and then 
\begin{equation}
 \nbigs_{\upsilon}\nbigttilde_{\upsilon}
 \Bigl(
 \gbigf^{(\infty,\infty)}_{+\star}\bigl(
 L_{\infty,\omega}(\nbigm),\vecnbigf
 \bigr)
 \Bigr)
= \gbigf^{(\infty,\infty)}_{+\star}\bigl(
 L_{\infty,\omega}(\nbigm),\vecnbigf
 \bigr)
 \quad
 (\star=!,\ast).
\end{equation}

There exist the morphisms of
$2\pi\seisuu$-equivariant local systems:
{\small
\begin{equation}
 \label{eq;24.4.19.40}
 c^{-1}(\nbigt_{\omega}(L_{\infty,\omega}(\nbigm)))
\to
 \gbigq^{\infty}_!\bigl(
 L_{\infty,\omega}(\nbigm),\vecnbigf
 \bigr)_{\real}
\to
 \gbigq^{\infty}_{\ast}\bigl(
 L_{\infty,\omega}(\nbigm),\vecnbigf
 \bigr)_{\real}
\to
 c^{-1}(\nbigt_{\omega}(L_{\infty,\omega}(\nbigm))).
\end{equation}}
Here,
$\gbigq^{\infty}_{\star}\bigl(
 L_{\infty,\omega}(\nbigm),\vecnbigf
 \bigr)_{\real}$ $(\star=!,\ast)$
are the $2\pi\seisuu$-equivariant local systems
underlying 
$\gbigf^{(\infty,\infty)}_{+\star}
\bigl(
L_{\infty,\omega}(\nbigm),
\vecnbigf
\bigr)$.
There also exists the following commutative diagram:
\begin{equation}
\label{eq;24.4.19.41}
\begin{CD}
 c^{-1}\bigl(
 L_{\infty,\omega}(\nbigm)
 \bigr)
 @>{\id-M}>>
  c^{-1}\bigl(
 L_{\infty,\omega}(\nbigm)
 \bigr)
 \\
 @V{\simeq}VV @V{\simeq}VV
 \\
\nbigt_{\upsilon}
 \Bigl(
\gbigf^{(\infty,\infty)}_{+!}\bigl(
 L_{\infty,\omega}(\nbigm),\vecnbigf
 \bigr)
 \Bigr)
 @>>>
  \nbigt_{\upsilon}
 \Bigl(
\gbigf^{(\infty,\infty)}_{+\ast}\bigl(
 L_{\infty,\omega}(\nbigm),\vecnbigf
 \bigr)
 \Bigr).
\end{CD}
\end{equation}
Here, $M$ denotes the monodromy automorphism
of $L_{\infty,\omega}(\nbigm)$.

\subsubsection{}

For $\alpha\in D$ and $\omega=\omega(\alpha,j)$,
we obtain the morphism of
$2\pi\seisuu$-equivariant
local systems with Stokes structure
\begin{equation}
\label{eq;21.4.27.21}
 \gbigf^{(0,\infty)}_{+!}\bigl(
 L_{\alpha,\omega}(\nbigm),\vecnbigf
 \bigr)
 \lrarr
 \gbigf^{(0,\infty)}_{+\ast}\bigl(
 L_{\alpha,\omega}(\nbigm),\vecnbigf
 \bigr).
\end{equation}
Set $\upsilon=\upsilon(\alpha,j)$,
and then
\[
 \nbigs_{\upsilon}\nbigttilde_{\upsilon}
 \Bigl(
 \gbigf^{(0,\infty)}_{+\star}\bigl(
 L_{\alpha,\omega}(\nbigm),\vecnbigf
 \bigr)
 \Bigr)
= \gbigf^{(0,\infty)}_{+\star}\bigl(
 L_{\alpha,\omega}(\nbigm),\vecnbigf
 \bigr)
 \quad
 (\star=!,\ast).
\]
There exist the morphisms of
$2\pi\seisuu$-equivariant local systems:
\begin{equation}
\label{eq;24.4.19.60}
  L_{\alpha,\omega}(\nbigm)
\to
 \gbigq^{0}_!\bigl(
 L_{\alpha,\omega}(\nbigm),\vecnbigf
 \bigr)
\to
 \gbigq^{0}_{\ast}\bigl(
 L_{\alpha,\omega}(\nbigm),\vecnbigf
 \bigr)
\to
 L_{\alpha,\omega}(\nbigm).
\end{equation}
Here,
$\gbigq^{0}_{\star}\bigl(
 L_{\alpha,\omega}(\nbigm),\vecnbigf
 \bigr)$
$(\star=!,\ast)$
are the $2\pi\seisuu$-equivariant local systems
underlying 
$\gbigf^{(0,\infty)}_{+\star}\bigl(
 L_{\alpha,\omega}(\nbigm),\vecnbigf
 \bigr)$.
There also exists the following commutative diagram:
\begin{equation}
\label{eq;24.4.19.61}
\begin{CD}
 \nbigt_{\omega}(L_{\alpha,\omega}(\nbigm))
 @>{\id-M_0^{-1}}>>
 \nbigt_{\omega}(L_{\alpha,\omega}(\nbigm))
 \\
 @V{\simeq}VV @V{\simeq}VV
 \\ 
\nbigt_{\upsilon}\Bigl(
\gbigf^{(0,\infty)}_{+!}\bigl(
L_{\alpha,\omega}(\nbigm),\vecnbigf
\bigr)
\Bigr)
@>>>
\nbigt_{\upsilon}\Bigl(
\gbigf^{(0,\infty)}_{+\ast}\bigl(
L_{\alpha,\omega}(\nbigm),\vecnbigf
\bigr)
\Bigr).
\end{CD}
\end{equation}
Here, $M_0$ denotes the monodromy automorphism of
$\nbigt_{\omega}(L_{\alpha,\omega}(\nbigm))$.

\subsubsection{}

We set $\nbigi_D=\{\alpha u^{-1}\,|\,\alpha\in D\}$.
We obtain the functor
$\gbigf_{\varrho}(\nbigl(\nbigm))$
from $\Dsf(D)$ to $\Shcat(\nbigi_D)$.
It induces the functor
$\Locst\gbigf_{\varrho}(\nbigl(\nbigm))$
$(\varrho\in\Dsf(D))$
from $\Dsf(D)$
to $\Loc^{\St}(\nbigi_D)$.
Note that for $\varrho_1\to\varrho_2$ in $\Dsf(D)$,
the induced morphisms
\[
 \Gr^{\vecnbigf}_{\alpha u^{-1}}
 \Bigl(
 \Locst\gbigf_{\varrho_1}(\nbigl(\nbigm))
 \Bigr)
 \to
 \Gr^{\vecnbigf}_{\alpha u^{-1}}
 \Bigl(
 \Locst\gbigf_{\varrho_2}(\nbigl(\nbigm))
 \Bigr)
\]
are isomorphisms if
$\varrho_1(\alpha)=\varrho_2(\alpha)$.
Moreover,
if $\varrho_1(\alpha)=!$ and $\varrho_2(\alpha)=\ast$,
there exists the following commutative diagram:
\begin{equation}
\label{eq;24.4.19.50}
 \begin{CD}
  L_{\alpha}(\nbigm)
  @>{\id-M_{\alpha}^{-1}}>>
  L_{\alpha}(\nbigm)\\
  @V{\simeq}VV @V{\simeq}VV
  \\
 \Gr^{\vecnbigf}_{\alpha u^{-1}}
 \Bigl(
 \Locst\gbigf_{\varrho_1}(\nbigl(\nbigm))
 \Bigr)
@>>>
 \Gr^{\vecnbigf}_{\alpha u^{-1}}
 \Bigl(
 \Locst\gbigf_{\varrho_2}(\nbigl(\nbigm))
 \Bigr)
 \end{CD}
\end{equation}
Here, $M_{\alpha}$ denotes the monodromy automorphism of
$L_{\alpha}(\nbigm)$.

There also exists the following morphisms
of $2\pi\seisuu$-equivariant local systems:
\begin{equation}
\label{eq;24.4.19.51}
 c^{-1}(L_{\infty}(\nbigm))
 \lrarr
 \Locst\gbigf_{\underline{!}}(\nbigl(\nbigm))
 \lrarr
 \Locst\gbigf_{\underline{\ast}}(\nbigl(\nbigm))
 \lrarr
 c^{-1}(L_{\infty}(\nbigm)).
\end{equation}

\subsubsection{}

For $\alpha\in D$,
we obtain
the regular singular holonomic $\nbigd_{\proj^1}$-module
$\Gr^{\vecnbigf}_0(\nbigm_{\alpha})$
characterized by the following conditions.
\begin{itemize}
 \item $\Gr^{\vecnbigf}_0(\nbigm_{\alpha})(\ast 0)$
       is the regular singular meromorphic flat bundle
       on $(\proj^1,\{\alpha,\infty\})$
       corresponding to
       $\Gr^{\vecnbigf}_0(L_{\alpha}(\nbigm))$.
 \item $\psi_{z-\alpha}\bigl(
       \Gr^{\vecnbigf}_0(\nbigm_{\alpha})
       \bigr)
       \to
       \phi_{z-\alpha}\bigl(
       \Gr^{\vecnbigf}_0(\nbigm_{\alpha})
       \bigr)
       \to
       \psi_{z-\alpha}\bigl(
       \Gr^{\vecnbigf}_0(\nbigm_{\alpha})
       \bigr)$
       are identified with
       $\psi_{z-\alpha}\bigl(
       \nbigm_{\alpha}
       \bigr)
       \to
       \phi_{z-\alpha}\bigl(
       \nbigm_{\alpha}
       \bigr)
       \to
       \psi_{z-\alpha}\bigl(
       \nbigm_{\alpha}
       \bigr)$.
\end{itemize}
We obtain the morphisms of $2\pi\seisuu$-equivariant local systems:
\begin{equation}
\label{eq;24.4.19.100}
 \gbigl^{\gbigf}\bigl(
 \Gr^{\vecnbigf}_0(\nbigm_{\alpha})(!\alpha)
 \bigr)
 \to
 \gbigl^{\gbigf}\bigl(
 \Gr^{\vecnbigf}_0(\nbigm_{\alpha})
 \bigr)
 \to
 \gbigl^{\gbigf}\bigl(
 \Gr^{\vecnbigf}_0(\nbigm_{\alpha})(\ast\alpha)
 \bigr). 
\end{equation}
These morphisms can be computed from 
$\psi_{z-\alpha}\bigl(
\nbigm_{\alpha}\bigr)
\to
\phi_{z-\alpha}\bigl(
\nbigm_{\alpha}
\bigr)
\to
\psi_{z-\alpha}\bigl(
\nbigm_{\alpha}
\bigr)$,
and
the $2\pi\seisuu$-equivariant local system
$\Gr^{\vecnbigf}_0(L_{\alpha}(\nbigm))$.
(See \S\ref{subsection;24.4.19.11}.)
In particular,
there exists the following commutative diagram:
\begin{equation}
\label{eq;24.4.19.70}
 \begin{CD}
  \Gr^{\vecnbigf}_0L_{\alpha}(\nbigm)
  @>{\id-M_{\alpha,0}^{-1}}>>
  \Gr^{\vecnbigf}_0L_{\alpha}(\nbigm)\\
  @V{\simeq}VV @V{\simeq}VV \\
 \gbigl^{\gbigf}\bigl(
 \Gr^{\vecnbigf}_0(\nbigm_{\alpha})(!\alpha)
 \bigr)
  @>>>
 \gbigl^{\gbigf}\bigl(
 \Gr^{\vecnbigf}_0(\nbigm_{\alpha})(\ast\alpha)
 \bigr)
 \end{CD}
\end{equation}
Here, $M_{\alpha,0}$ denotes the monodromy automorphism of
$\Gr^{\vecnbigf}_0L_{\alpha}(\nbigm)$.

\subsection{Connecting morphisms}

\subsubsection{}

For $\omega(j)=\omega(\infty,j)$,
$\omega(j+1)=\omega(\infty,j+1)$
and $\upsilon(j)=\upsilon(\infty,j)$
$(j=1,\ldots,\ell(\infty)-1)$,
we obtain the following commutative diagram
of $2\pi\seisuu$-equivariant local systems
from (\ref{eq;24.4.19.40})
and (\ref{eq;24.4.19.41}),
which are not necessarily compatible with the Stokes structures:
\begin{equation}
\label{eq;24.4.19.101}
 \begin{CD}
  \nbigt_{\upsilon(j)}\bigl(
  \gbigf^{(\infty,\infty)}_{+,!}
  (L_{\infty,\omega(j)},\vecnbigf)
  \bigr)
  @>>>
  \nbigt_{\upsilon(j)}\bigl(
  \gbigf^{(\infty,\infty)}_{+,\ast}
  (L_{\infty,\omega(j)},\vecnbigf)
  \bigr)
  \\
  @VVV @AAA \\
  \gbigf^{(\infty,\infty)}_{+,!}\bigl(
 L_{\infty,\omega(j+1)},\vecnbigf
  \bigr)
  @>>>
  \gbigf^{(\infty,\infty)}_{+,\ast}\bigl(
 L_{\infty,\omega(j+1)},\vecnbigf
  \bigr).
 \end{CD}
\end{equation}

\subsubsection{}

For $\omega=\omega(\infty,\ell(\infty))$
and $\upsilon=\upsilon(\infty,\ell(\infty))$,
we obtain the following commutative diagram
of $2\pi\seisuu$-equivariant local systems,
from (\ref{eq;24.4.19.41})
and (\ref{eq;24.4.19.51}),
where the vertical morphisms are {\em not} necessarily compatible with
the Stokes structures:
\begin{equation}
\label{eq;21.4.27.38}
 \begin{CD}
 \nbigt_{\upsilon}\Bigl(
 \gbigf^{(\infty,\infty)}_{+!}\bigl(
 L_{\infty,\omega}(\nbigm),\vecnbigf
  \bigr)
  \Bigr)
  @>>>
 \nbigt_{\upsilon}\Bigl(
 \gbigf^{(\infty,\infty)}_{+\ast}\bigl(
 L_{\infty,\omega}(\nbigm),\vecnbigf
  \bigr)
  \Bigr)
  \\
  @VVV @AAA\\
 \Locst
  \gbigf_{\underline{!}}(\nbigl(\nbigm))
  @>>>
 \Locst\gbigf_{\underline{\ast}}(\nbigl(\nbigm)).
 \end{CD}
\end{equation}

\subsubsection{}
For $\alpha\in D$,
there exists the following commutative diagram
of {\em $2\pi\seisuu$-equivariant local systems}
from (\ref{eq;24.4.19.60})
and (\ref{eq;24.4.19.50}),
where the vertical morphisms are {\em not} necessarily compatible with
the Stokes structures:
\begin{equation}
 \label{eq;21.4.27.22}
 \begin{CD}
  \Gr^{\vecnbigf^{(1)}}_{\alpha u^{-1}}
  \Locst
  \gbigf_{\underline{!}}(\nbigl(\nbigm))
  @>>>
  \Gr^{\vecnbigf^{(1)}}_{\alpha u^{-1}}
  \Locst
  \gbigf_{\underline{\ast}}(\nbigl(\nbigm))
  \\
  @VVV @AAA \\
  \gbigf^{(0,\infty)}_{+!}\bigl(
 L_{\alpha,\omega(\alpha,1)}(\nbigm),\vecnbigf
 \bigr)
 @>>>
 \gbigf^{(0,\infty)}_{+\ast}\bigl(
 L_{\alpha,\omega(\alpha,1)}(\nbigm),\vecnbigf
 \bigr)
 \end{CD}
\end{equation}

\subsubsection{}

For $\omega(j)=\omega(\alpha,j)$,
$\omega(j+1)=\omega(\alpha,j+1)$
and $\upsilon(j)=\upsilon(\alpha,j)$
$(j=1,\ldots,\ell(\alpha)-1)$,
we obtain the following commutative diagram
of {\em $2\pi\seisuu$-equivariant local systems}
from (\ref{eq;24.4.19.60}) and (\ref{eq;24.4.19.61}),
where the vertical morphisms are {\em not} necessarily compatible with
Stokes structures:
\begin{equation}
\label{eq;21.4.27.39}
 \begin{CD}
 \nbigt_{\upsilon(j)}\Bigl(
  \gbigf^{(0,\infty)}_{+!}\bigl(
 L_{\omega(j)}(\nbigm,\alpha),\vecnbigf
 \bigr)\Bigr)
  @>>>
 \nbigt_{\upsilon(j)}\Bigl(
  \gbigf^{(0,\infty)}_{+\ast}\bigl(
 L_{\omega(j)}(\nbigm,\alpha),\vecnbigf
  \bigr)
  \Bigr)
  \\
  @VVV @AAA\\
 \gbigf^{(0,\infty)}_{+!}\bigl(
 L_{\omega(j+1)}(\nbigm,\alpha),\vecnbigf
 \bigr)
  @>>>
 \gbigf^{(0,\infty)}_{+\ast}\bigl(
 L_{\omega(j+1)}(\nbigm,\alpha),\vecnbigf
 \bigr).
 \end{CD}
\end{equation}
If $j=\ell(\alpha)-1$,
then the vertical morphisms are isomorphisms.

\subsubsection{}

For $\omega=\omega(\alpha,\ell(\alpha))$,
we obtain the following commutative diagram
of $2\pi\seisuu$-equivariant local systems
from (\ref{eq;24.4.19.70}).
\begin{equation}
\label{eq;21.4.27.14}
 \begin{CD}
   \gbigf^{(0,\infty)}_{+!}\bigl(
 L_{\alpha,\omega}(\nbigm),\vecnbigf
 \bigr)
  @>>>
 \gbigf^{(0,\infty)}_{+\ast}\bigl(
 L_{\alpha,\omega}(\nbigm),\vecnbigf
  \bigr)
  \\
  @V{\simeq}VV @A{\simeq}AA\\
  \gbigl^{\gbigf}(\Gr^{\vecnbigf}_0(\nbigm_{\alpha})(!\alpha))
  @>>>
  \gbigl^{\gbigf}(\Gr^{\vecnbigf}_0(\nbigm_{\alpha})(\ast\alpha)).
 \end{CD}
\end{equation}

\subsection{Description of
$(\gbigl^{\gbigf}(\nbigm),\vecnbigf)$}

Let us explain how to obtain an explicit description of
$(\gbigl^{\gbigf}(\nbigm),\vecnbigf)$
from the morphisms (\ref{eq;21.4.27.32}),
(\ref{eq;21.4.27.21}),
(\ref{eq;24.4.19.100})
and the functor 
$\Locst\gbigf_{\varrho}(\nbigl(\nbigm))$
$(\varrho\in\Dsf(D))$
together with
the commutative diagrams
(\ref{eq;24.4.19.101}),
(\ref{eq;21.4.27.38}),
(\ref{eq;21.4.27.22}),
(\ref{eq;21.4.27.39})
and (\ref{eq;21.4.27.14})
by using successively {\em the extension} of
local systems with Stokes structure.

\subsubsection{}

Let $\alpha\in D$.
By a descending induction,
we shall construct morphisms of
$2\pi\seisuu$-equivariant local systems
with Stokes structure 
\begin{equation}
\label{eq;21.4.27.12}
 \begin{CD}
  (\gbigl^{\gbigf}(\nbigm(!D))_{\alpha,\upsilon(j)},\vecnbigf)
  @>>>
  (\gbigl^{\gbigf}(\nbigm)_{\alpha,\upsilon(j)},\vecnbigf)
  @>>>
  (\gbigl^{\gbigf}(\nbigm(\ast D))_{\alpha,\upsilon(j)},\vecnbigf)
 \end{CD}
\end{equation}
for $\upsilon(j)=\upsilon(\alpha,j)$,
and commutative diagrams
of {\em $2\pi\seisuu$-equivariant local systems with Stokes structure}
\begin{equation}
\label{eq;21.4.27.13}
\begin{CD}
\gbigf^{(0,\infty)}_{+!}(L_{\alpha,\omega(j)}(\nbigm),\vecnbigf)
@>>>
\gbigf^{(0,\infty)}_{+\ast}(L_{\alpha,\omega(j)}(\nbigm),\vecnbigf)
\\
@VVV @AAA \\
\nbigs_{\upsilon(j)}
 (\gbigl^{\gbigf}(\nbigm(!D))_{\alpha,\upsilon(j)},\vecnbigf)
 @>>>
\nbigs_{\upsilon(j)}
(\gbigl^{\gbigf}(\nbigm(\ast D))_{\alpha,\upsilon(j)},\vecnbigf)
\end{CD}
\end{equation}
for $\upsilon(j)=\upsilon(\alpha,j)$
and $\omega(j)=\omega(\alpha,j)$,
as follows.

In the case $j=\ell(\alpha)$,
we set
$(\gbigl^{\gbigf}(\nbigm(\star D))_{\alpha,\upsilon(\ell(\alpha))},\vecnbigf)
=\gbigl^{\gbigf}\bigl(
 \Gr^{\vecnbigf}_0(\nbigm_{\alpha})(\star\alpha)
 \bigr)$
($\star=!,\ast$),
and 
$(\gbigl^{\gbigf}(\nbigm)_{\alpha,\upsilon(\ell(\alpha))},\vecnbigf)
=\gbigl^{\gbigf}\bigl(
 \Gr^{\vecnbigf}_0(\nbigm_{\alpha})
 \bigr)$.
The morphisms 
(\ref{eq;21.4.27.12})
and the diagram (\ref{eq;21.4.27.13})
are obtained from
(\ref{eq;24.4.19.100})
and (\ref{eq;21.4.27.14}),
respectively.
Suppose that we have already constructed
(\ref{eq;21.4.27.12})
and (\ref{eq;21.4.27.13})
for $\omega(\alpha,j+1)$
and $\upsilon(\alpha,j+1)$.
We obtain the following commutative diagram
of {\em $2\pi\seisuu$-equivariant local systems}
from (\ref{eq;21.4.27.39}),
and (\ref{eq;21.4.27.13})
for $\omega(\alpha,j+1)$
and $\upsilon(\alpha,j+1)$:
\begin{equation}
\label{eq;21.4.27.20}
 \begin{CD}
  \nbigt_{\upsilon(j)}\Bigl(
  \gbigf^{(0,\infty)}_{+!}\bigl(
 L_{\alpha,\omega(j)}(\nbigm),\vecnbigf
  \bigr)
  \Bigr)
  @>>>
 \nbigt_{\upsilon(j)}\Bigl(
  \gbigf^{(0,\infty)}_{+\ast}\bigl(
 L_{\alpha,\omega(j)}(\nbigm),\vecnbigf
  \bigr)
  \Bigr)
  \\
  @VVV @AAA\\
 \gbigl^{\gbigf}(\nbigm(!D))_{\alpha,\upsilon(j+1)}
  @>>>
 \gbigl^{\gbigf}(\nbigm(\ast D))_{\alpha,\upsilon(j+1)}.
 \end{CD}
\end{equation}
By using the extension of
$2\pi\seisuu$-equivariant local systems with Stokes structure
(see \S\ref{subsection;21.4.29.1}),
we obtain
(\ref{eq;21.4.27.12})
and (\ref{eq;21.4.27.13})
for $\omega=\omega(\alpha,j)$
from (\ref{eq;21.4.27.21}),
(\ref{eq;21.4.27.20})
and (\ref{eq;21.4.27.12})
for $\omega=\omega(\alpha,j+1)$.

\subsubsection{}

We obtain the $2\pi\seisuu$-equivariant local systems
$\gbigl^{\gbigf}(\nbigm)_{\alpha}:=
\gbigl^{\gbigf}(\nbigm)_{\alpha,\upsilon(\alpha,1)}$
and
$\gbigl^{\gbigf}(\nbigm(\star D))_{\alpha}:=
\gbigl^{\gbigf}(\nbigm(\star D))_{\alpha,\upsilon(\alpha,1)}$
$(\star=!,\ast)$.
We define the Stokes structure $\vecnbigf$
on $\gbigl^{\gbigf}(\nbigm)_{\alpha}$
and $\gbigl^{\gbigf}(\nbigm(\star D))_{\alpha}$
by setting
\[
 \nbigf^{\theta}_{\alpha u^{-1}+\gminia}
 \bigl(
 \gbigl^{\gbigf}(\nbigm)_{\alpha|\theta}
 \bigr)
 :=
 \nbigf^{\theta}_{\gminia}
 \bigl(
 \gbigl^{\gbigf}(\nbigm)_{\alpha,\upsilon(\alpha,1)|\theta}
 \bigr),
\]
\[
 \nbigf^{\theta}_{\alpha u^{-1}+\gminia}
 \bigl(
 \gbigl^{\gbigf}(\nbigm(\star D))_{\alpha|\theta}
 \bigr)
 :=
 \nbigf^{\theta}_{\gminia}
 \bigl(
 \gbigl^{\gbigf}(\nbigm(\star D))_{\alpha,\upsilon(\alpha,1)|\theta}
 \bigr).
\]
Thus, we obtain the following morphisms of
$2\pi\seisuu$-equivariant local systems with Stokes structure:
\begin{equation}
\label{eq;21.4.27.23}
 \begin{CD}
  (\gbigl^{\gbigf}(\nbigm(!D))_{\alpha},\vecnbigf)
  @>>>
  (\gbigl^{\gbigf}(\nbigm)_{\alpha},
  \vecnbigf)
  @>>>
  (\gbigl^{\gbigf}(\nbigm(\ast D))_{\alpha},\vecnbigf).
 \end{CD}
\end{equation}
We obtain the following commutative diagram of
{\em $2\pi\seisuu$-equivariant local systems}
from (\ref{eq;21.4.27.13}):
\begin{equation}
\label{eq;21.4.27.25}
   \begin{CD}
 \gbigf^{(0,\infty)}_{+!}\bigl(
 L_{\alpha,\omega(\alpha,1)}(\nbigm),\vecnbigf
 \bigr)
 @>>>
 \gbigf^{(0,\infty)}_{+\ast}\bigl(
 L_{\alpha,\omega(\alpha,1)}(\nbigm),\vecnbigf
 \bigr)\\
@VVV @AAA \\
  \gbigl^{\gbigf}(\nbigm(!D))_{\alpha}
  @>>>
  \gbigl^{\gbigf}(\nbigm(\ast D))_{\alpha}.
   \end{CD}
\end{equation}

\subsubsection{}

We obtain the following commutative diagram of
$2\pi\seisuu$-equivariant local systems
from (\ref{eq;21.4.27.22})
and (\ref{eq;21.4.27.25})
for $\omega=\omega(\alpha,1)$:
\begin{equation}
\label{eq;21.4.27.24}
 \begin{CD}
  \Gr^{\vecnbigf^{(1)}}_{\alpha u^{-1}}
  \Locst
  \gbigf_{\underline{!}}(\nbigl(\nbigm))
  @>>>
  \Gr^{\vecnbigf^{(1)}}_{\alpha u^{-1}}
  \Locst
  \gbigf_{\underline{\ast}}(\nbigl(\nbigm))
  \\
  @VVV @AAA \\
  \gbigl^{\gbigf}(\nbigm(!D))_{\alpha}
  @>>>
  \gbigl^{\gbigf}(\nbigm(\ast D))_{\alpha}.
 \end{CD}
\end{equation}
By using the extension of local systems with
Stokes structure (see \S\ref{subsection;21.4.29.2}),
together with
the functor
$\Locst\gbigf_{\varrho}(\nbigl(\nbigm))$
$(\varrho\in\Dsf(D))$,
(\ref{eq;21.4.27.23})
and (\ref{eq;21.4.27.24}),
we obtain a $2\pi\seisuu$-equivariant local system
with Stokes structure
$\bigl(
 \gbigl^{\gbigf}(\nbigm)^{(1)},
 \vecnbigf
\bigr)$
with morphisms of {\em $2\pi\seisuu$-equivariant local systems}
\begin{equation}
\label{eq;21.4.27.30}
\begin{CD}
 \Locst\gbigf_{\underline{!}}(\nbigl(\nbigm))
 @>>>
 \gbigl^{\gbigf}(\nbigm)^{(1)}
 @>>>
 \Locst\gbigf_{\underline{\ast}}(\nbigl(\nbigm))
\end{CD}
\end{equation}
such that the composition of (\ref{eq;21.4.27.30})
equals the natural morphism.

For $\omega=\omega(\infty,\ell(\infty))$
and $\upsilon=\upsilon(\infty,\ell(\infty))$,
we obtain the following morphisms of
$2\pi\seisuu$-equivariant local systems
from (\ref{eq;21.4.27.38}) and (\ref{eq;21.4.27.30}):
\begin{multline}
\label{eq;21.4.27.33}
 \nbigt_{\upsilon}\Bigl(
 \gbigf^{(\infty,\infty)}_{+!}\bigl(
 L_{\infty,\omega}(\nbigm),\vecnbigf
 \bigr)
 \Bigr)
\lrarr
 \gbigl^{\gbigf}(\nbigm)^{(1)} 
 \lrarr
  \nbigt_{\upsilon}\Bigl(
 \gbigf^{(\infty,\infty)}_{+\ast}\bigl(
 L_{\infty,\omega}(\nbigm),\vecnbigf
  \bigr)\Bigr).
\end{multline}

\subsubsection{}

For $\omega(j)=\omega(\infty,j)$
and $\upsilon(j)=\upsilon(\infty,j)$,
by a descending induction on $j$,
let us construct
a $2\pi\seisuu$-equivariant local systems
with Stokes structure
$(\gbigl^{\gbigf}(\nbigm)^{(\upsilon(j))},\vecnbigf)$
and morphisms of $2\pi\seisuu$-equivariant
local systems with Stokes structure
\begin{equation}
\label{eq;21.4.27.31}
 \gbigf^{(\infty,\infty)}_{+!}\bigl(
 L_{\infty,\omega(j)}(\nbigm),\vecnbigf
  \bigr)
 \lrarr
 \nbigs_{\upsilon(j)}
 (\gbigl^{\gbigf}(\nbigm)^{(\upsilon(j))},\vecnbigf)
\lrarr
 \gbigf^{(\infty,\infty)}_{+\ast}\bigl(
 L_{\infty,\omega(j)}(\nbigm),\vecnbigf
  \bigr),
\end{equation}
such that the composition of (\ref{eq;21.4.27.31})
is equal to the morphism (\ref{eq;21.4.27.32}).

In the case $\omega=\omega(\infty,\ell(\infty))$
and $\upsilon=\upsilon(\infty,\ell(\infty))$,
we obtain
$(\gbigl^{\gbigf}(\nbigm)^{(\upsilon)},\vecnbigf)$
and (\ref{eq;21.4.27.31})
from (\ref{eq;21.4.27.33})
and (\ref{eq;21.4.27.32})
by using the extension of
local systems with Stokes structure
(see \S\ref{subsection;21.4.29.1}).
Suppose that we have already constructed
$(\gbigl^{\gbigf}(\nbigm)^{(\upsilon(j+1))},\vecnbigf)$
and (\ref{eq;21.4.27.31})
for $\upsilon(\infty,j+1)$
and $\omega(\infty,j+1)$.
We obtain the following morphisms
of {\em $2\pi\seisuu$-equivariant local systems}
from (\ref{eq;24.4.19.101})
and (\ref{eq;21.4.27.31}):
\begin{multline}
\label{eq;21.4.27.35}
\nbigt_{\upsilon(j)}\Bigl(
 \gbigf^{(\infty,\infty)}_{+!}\bigl(
 L_{\infty,\omega(j)}(\nbigm),\vecnbigf
 \bigr)
 \Bigr)
 \lrarr
\gbigl^{\gbigf}(\nbigm)^{(\upsilon(j+1))} \\
 \lrarr
\nbigt_{\upsilon(j)}\Bigl(
 \gbigf^{(\infty,\infty)}_{+\ast}\bigl(
 L_{\infty,\omega(j)}(\nbigm),\vecnbigf
 \bigr)\Bigr).
\end{multline}
By using the extension of
$2\pi\seisuu$-equivariant local systems with
Stokes structure,
we obtain
$(\gbigl^{\gbigf}(\nbigm)^{(\upsilon(j))},\vecnbigf)$
and (\ref{eq;21.4.27.31})
from (\ref{eq;21.4.27.35})
and (\ref{eq;21.4.27.32}) for $\omega(j)$.
We obtain the following main theorem of this monograph
from the results in  \S\ref{subsection;24.4.15.100},
\S\ref{subsection;24.4.19.10}
and \S\ref{subsection;24.4.15.101}.

\begin{thm}
\label{thm;21.4.27.100}
 $(\gbigl^{\gbigf}(\nbigm),\vecnbigf)$
is naturally isomorphic to 
$(\gbigl^{\gbigf}(\nbigm)^{(\upsilon(\infty,1))},\vecnbigf)$.
\end{thm}

\begin{rem}
Note that we also constructed
the following morphisms
\[
c^{-1}(\nbigttilde_{1}L_{\infty}(\nbigm))
\lrarr
(\gbigl^{\gbigf}(\nbigm)^{(\upsilon(\infty,1))},\vecnbigf)
\lrarr
c^{-1}(\nbigttilde_{1}L_{\infty}(\nbigm)).
\]
\hfill\qed
\end{rem}

\subsection{Complement}

Let $h:\proj^1\to\proj^1$ denote the morphism
defined by $h(z)=-z$.
We set $\nbigm^{\gbigf}=\Fourier_+(\nbigm)$.
Because 
$\nbigm=
\Fourier_-(\nbigm^{\gbigf})$,
we obtain
\[
 \Fourier_+(\nbigm^{\gbigf})
 =h^{\ast}\nbigm.
\]

There exists the following commutative diagram of
$2\pi\seisuu$-equivariant local systems:
{\footnotesize
\begin{equation}
\label{eq;25.3.18.110}
 \begin{CD}
c^{-1}(L_{\infty}(\nbigm^{\gbigf}))
@>>>
 \gbigl^{\gbigf}(\nbigstilde^{\infty}_{\infty}(\nbigm^{\gbigf}))
  @>>>
  c^{-1}(L_{\infty}(\nbigm^{\gbigf}))
  \\
@VVV @VVV @VVV
  \\
 c^{-1}L_{\infty}
 \bigl(
 (\nbigstilde^{\infty}_1h^{\ast}\nbigm)^{\gbigf_-}
 \bigr)
 @>{(2\pi\sqrt{-1})^{-1}M\cdot b_2}>>
 L_{\infty}\bigl(
 \nbigttilde^{\infty}_1
 (h^{\ast}\nbigm)
 \bigr)
 @>{-(2\pi\sqrt{-1})b_1}>>
c^{-1}L_{\infty}
 \bigl(
 (\nbigstilde^{\infty}_1h^{\ast}\nbigm)^{\gbigf_-}
\bigr).
 \end{CD}
\end{equation}
 }
See \S\ref{subsection;25.3.18.101}
and \S\ref{subsection;25.3.18.100}
for the lower horizontal arrows.
We can recover $\LS^{\fin}(\nbigm^{\gbigf})$
from
\[
\bigl(
\gbigl^{\gbigf}(\nbigstilde^{\infty}_{\infty}(\nbigm^{\gbigf})),
 \vecnbigf
 \bigr)
\simeq 
 \nbigttilde_1\bigl(
 L_{\infty}((\nbigm^{\gbigf})^{\gbigf}),
 \vecnbigf
 \bigr)
=\nbigttilde_1\bigl(
 L_{\infty}(h^{\ast}\nbigm),
 \vecnbigf
 \bigr)
\]
and the upper horizontal arrows
in (\ref{eq;25.3.18.110}).
We can compute the lower horizontal arrows
in (\ref{eq;25.3.18.110})
from $\LS(h^{\ast}\nbigm)$.
In this way,
we can also compute
$\LS^{\fin}(\nbigm^{\gbigf})$
from $\LS(\nbigm)$.

\section{Examples}
\label{subsection;21.6.13.2}

\subsection{}

Let
$\nbigi=\bigl\{\alpha_1x^{-2},\ldots,\alpha_nx^{-2}\bigr\}$
for some mutually distinct positive numbers $\alpha_i$.
Let $\nbigm$ be a meromorphic flat bundle
on $(\proj^1,\infty)$
such that $\nbigi_{\infty}(\nbigm)=\nbigi$.
Let us study
$(\gbigl^{\gbigf}(\nbigm),\vecnbigf)$.

\subsubsection{}

We set
\[
 \nbigi^{\circ}_0
=\gbigf^{(\infty,\infty)}_+(\nbigi)
=\bigl\{
 (-1/4)\alpha_i u^{-2}
 \bigr\},
 \quad
 \nbigi^{\circ}=\nbigi^{\circ}_0\cup\{0\}.
\]
There exists the bijection
$\nbigi\simeq\nbigi_0^{\circ}$ given by
$\alpha_ix^{-2}\mapsto
-(1/4)\alpha_iu^{-2}$.
It induces the following isomorphism of the partially ordered sets
for any $\theta^u\in\real$:
\begin{equation}
 \label{eq;24.4.20.21}
  (\nbigi,\leq_{\theta^u-\pi})
  \simeq
  (\nbigi^{\circ}_0,\leq_{\theta^u})
\end{equation}

We set
$(L,\vecnbigf)=
(L_{\infty}(\nbigm),\vecnbigf)\in\Loc^{\St}(\nbigi)$. 
Let $\Ltilde^{\gbigf}$ denote the pull back of
$L$ by
$\theta^u\mapsto \theta^u-\pi$.
It is equipped with the Stokes structure
$\vecnbigf$ indexed by $\nbigi^{\circ}_0$,
induced by the Stokes structure of $L$ indexed by $\nbigi$
and the isomorphism of the partially ordered sets
(\ref{eq;24.4.20.21}).
We shall explain how to recover the following result
in \cite{Sabbah-pure-Gaussian}
(see also \cite{Hohl-pure-Gaussian}).
\begin{prop}
\label{prop;24.4.20.22}
$(\gbigl^{\gbigf}(\nbigm),\vecnbigf)
 \simeq
 (L^{\gbigf},\vecnbigf)$.
\end{prop}
 
\subsubsection{Local Fourier transform}

By the local Fourier transform of $(L,\vecnbigf)$
in \S\ref{subsection;24.4.5.120},
we obtain
\[
\gbigf^{(\infty,\infty)}_{+\star}(L,\vecnbigf)
=\bigl(
 \gbigq^{\infty}_{\star}(L,\vecnbigf)_{\real},
 \vecnbigf
\bigr)
\in
\Loc^{\St}(\nbigi^{\circ})
\quad(\star=!,\ast)
\]
which are isomorphic to
$(\gbigl^{\gbigf}(\nbigm(\star0)),\vecnbigf)$.
Let us describe
$\gbigf^{(\infty,\infty)}_{+!}(L,\vecnbigf)$
precisely.

\subsubsection{Local system $\gbigq^{\infty}_!(L,\vecnbigf)$}

For any $a\in\real$ and $r>0$,
we set $I(a,r)=\bigl\{t\in\real\,\bigr|\,|t-a|<r\bigr\}$.
For any integer $m$,
we set
$J_m=I(m\frac{\pi}{2},\frac{\pi}{4})$.
We have
$-\Re(\alpha_ie^{-2\sqrt{-1}\theta})>0$ on $J_{2\ell+1}$
and
$-\Re(\alpha_ie^{-2\sqrt{-1}\theta})<0$ on $J_{2\ell}$
for any integer $\ell$.

We consider the vector space
$\bigoplus_{\ell\in\seisuu}
H^0(J_{2\ell+1},L)$.
For $v\in H^0(\real,L)$,
the induced element in $H^0(J_{2\ell+1},L)$
is denoted by
$\langle J_{2\ell+1},v\rangle$.
We define
the vector space $\gbigq^{\infty}_!(L,\vecnbigf)$
as the quotient of 
$\bigoplus_{\ell\in\seisuu}
H^0(J_{2\ell+1},L)$
by the equivalence relation generated by
\[
 \langle J_{2\ell+1},v\rangle
 \sim
 \langle
 J_{2\ell-3},v\rangle.
\]
The $2\pi\seisuu$-equivariant local system
$\gbigq_!^{\infty}(L,\vecnbigf)_{\real}$
equals the $2\pi\seisuu$-equivariant local system
induced by 
$\gbigq_!^{\infty}(L,\vecnbigf)$
with the trivial action.

\subsubsection{Some maps}
For any integer $m$,
we set
$\vecJ_{0,m}=I(\frac{(m+1)\pi}{2},\frac{\pi}{4})$.
We have
$-\Re(-\frac{1}{4}\alpha_ie^{-2\sqrt{-1}\theta^u})>0$
on $\vecJ_{0,2\ell}$
and
$-\Re(-\frac{1}{4}\alpha_ie^{-2\sqrt{-1}\theta^u})<0$
on $\vecJ_{0,2\ell+1}$.

We naturally identify
$H^0(\real,L)=H^0(J_{m},L)$
for any interval $J$.
The maps
$\vecB_{\vecJ_{0,2\ell}}:
H^0(J_{-2\ell-1},L)
\to \gbigq^{\infty}_!(L,\vecnbigf)$
and
$\vecA_{(\vecJ_{0,2\ell+1})_{\pm}}:
H^0(J_{-2\ell},L)
\to \gbigq^{\infty}_!(L,\vecnbigf)$
are given as follows.
\[
 \vecB_{\vecJ_{0,2\ell}}(v)
=\vecA_{(\vecJ_{0,2\ell+1})_+}(v)
=\langle J_{-2\ell-1},v\rangle,
\quad
  \vecA_{(\vecJ_{0,2\ell+1})_-}(v)
 =-\langle J_{-2\ell+1},v\rangle.
\]
Let $\vecA_{\infty}:H^0(\real,L)\to \gbigq_!^{\infty}(L,\vecnbigf)$
be defined by
$\vecA_{\infty}(v)
:=\langle J_{-1},v\rangle
 +\langle J_1,v\rangle$
for any $v\in H^0(\real,L)$, 
which equals
$\langle J_{2\ell-1},v\rangle
+\langle J_{2\ell+1},v\rangle$.
We have
\[
 \vecA^{(\vecJ_{2\ell})_+}_{\infty}
=\vecA^{(\vecJ_{2\ell})_-}_{\infty}
=\vecA^{(\vecJ_{2\ell+1})_-}_{\infty}
=\vecA^{(\vecJ_{2\ell+1})_+}_{\infty}
=\vecA_{\infty}.
\]

\subsubsection{Stokes filtrations}

We have the isomorphism
$\vecJ_{0,2\ell+1}\simeq J_{-2\ell}$
and
$\vecJ_{0,2\ell}\simeq J_{-2\ell-1}$
given by
$\theta^u\mapsto
 \theta^u-(2\ell+1)\pi$.
There exists the natural bijection
$\nbigi\simeq\nbigi^{\circ}_0$
given by
$\alpha_ix^{-2}\mapsto
-(1/4)\alpha_iu^{-2}$.
It induces the isomorphism of the partially ordered sets
\begin{equation}
\label{eq;24.4.20.20}
\bigl(
 \nbigi,\leq_{\theta^u-(2\ell+1)\pi}
\bigr)
 \simeq
\bigl(
 \nbigi^{\circ}_0,\leq_{\theta^u}
 \bigr).
\end{equation}

For any $\theta\in\real$,
we identify
$L_{|\theta}$ with $H^0(\real,L)$.
We have the filtrations
$\nbigf^{\theta}$ on $H^0(\real,L)$.

The Stokes filtration
$\nbigf^{\theta^u}$
on $\gbigq^{\infty}_!(L,\vecnbigf)$
is given as follows:
\begin{itemize}
 \item For $\theta^u\in\vecJ_{0,2\ell}$,
we have 
\[
 \nbigf^{\theta^u}_{-(1/4)\alpha_iu^{-2}}
 \gbigq^{\infty}_!(L,\vecnbigf)
 =\vecB_{\vecJ_{0,2\ell}}
 \Bigl(
 \nbigf^{\theta^u-(2\ell+1)\pi}_{\alpha_ix^{-2}}H^0(\real,L)
 \Bigr),
\]
\[
\nbigf^{\theta^u}_0
\gbigq^{\infty}_!(L,\vecnbigf)
=\gbigq^{\infty}_!(L,\vecnbigf).
\]
 \item For $\theta^u\in\vecJ_{0,2\ell+1}$,
       we have
\[
 \nbigf^{\theta^u}_0
\gbigq^{\infty}_!(L,\vecnbigf)
       =\Image \vecA_{\infty},
\]
\[
 \nbigf^{\theta^u}_{-(1/4)\alpha_iu^{-2}}
 \gbigq^{\infty}_!(L,\vecnbigf)
       =
\Image\vecA_{\infty}\oplus
 \vecA_{(\vecJ_{0,2\ell+1})_+}
 \Bigl(
 \nbigf^{\theta^u-(2\ell+1)\pi}_{\alpha_ix^{-2}}H^0(\real,L)
 \Bigr).
\]
The latter equals
$\Image\vecA_{\infty}\oplus
\vecA_{(\vecJ_{0,2\ell+1})_-}
\Bigl(
\nbigf^{\theta^u-(2\ell+1)\pi}_{\alpha_ix^{-2}}H^0(\real,L)
\Bigr)$.
 \item Let $\theta^u\in \vecJbar_{0,2\ell}\setminus\vecJ_{0,2\ell}$.
       We have the decomposition
       \[
       H^0(\real,L)
       =\bigoplus \nbigf^{\theta^u-(2\ell+1)\pi}_{\alpha_ix^{-2}}
       H^0(\real,L).
       \]
       The Stokes filtration $\nbigf^{\theta^u}$ is given by the splitting
\[
       \gbigq^{\infty}_!(L,\vecnbigf)
       =\Image\vecA_{\infty}
       \oplus
       \bigoplus
       \vecB_{\vecJ_{0,2\ell}}\Bigl(
       \nbigf^{\theta^u-(2\ell+1)\pi}_{\alpha_ix^{-2}}
       H^0(\real,L)
       \Bigr).
\]       
\end{itemize}

\subsubsection{Proof of Proposition \ref{prop;24.4.20.22}}

By Theorem \ref{thm;24.4.20.10},
$(\gbigl^{\gbigf}(\nbigm),\vecnbigf)$
is the extension of
$(\gbigl^{\gbigf}(\nbigm(!0)),\vecnbigf)
\to
(\gbigl^{\gbigf}(\nbigm(\ast0)),\vecnbigf)$
by the trivial morphisms
\[
\Gr^{\vecnbigf}_0\gbigl^{\gbigf}(\nbigm(!0))
\to 0 \to
\Gr^{\vecnbigf}_0\gbigl^{\gbigf}(\nbigm(\ast 0)). 
\]
Hence, we obtain
\[
 \gbigl^{\gbigf}(\nbigm)
 \simeq
 \gbigf^{(\infty,\infty)}_{+!}(L,\vecnbigf)/\Image \vecA_{\infty}.
\]
Note that
$\vecB_{\vecJ_{0,2\ell}}
=
\vecA_{(\vecJ_{0,2\ell+1})_+}
=\vecA_{(\vecJ_{0,2\ell+1})_-}
=-\vecB_{\vecJ_{0,2\ell+2}}$
in this quotient.

We have the isomorphism
$\vecJbar_{0,2\ell}\cup\vecJbar_{0,2\ell+1}
\simeq \Jbar_{-2\ell-1}\cup\Jbar_{-2\ell}$
given by
$\theta^u\mapsto
\theta^u-(2\ell+1)\pi$.
Let $L^{\gbigf}_{\vecJbar_{0,2\ell}\cup\vecJbar_{0,2\ell+1}}$
denote the pull back of
$L_{|\Jbar_{-2\ell-1}\cup\Jbar_{-2\ell}}$.
It is equipped with the Stokes filtrations
induced by the Stokes structure of
$L_{|\Jbar_{-2\ell-1}\cup\Jbar_{-2\ell}}$
and the isomorphism of the partially ordered sets
(\ref{eq;24.4.20.20}).

Let $\theta^u(\ell,\ell+1)$ be the intersection of
$\vecJbar_{0,2\ell}\cup\vecJbar_{0,2\ell+1}$
and
$\vecJbar_{0,2\ell+2}\cup\vecJbar_{0,2\ell+3}$.
By the $2\pi\seisuu$-action on $L$,
there exists the natural isomorphism
\[
 \Psi(\ell,\ell+1):
 L^{\gbigf}_{\vecJbar_{0,2\ell}\cup\vecJbar_{0,2\ell+1}|\theta^u(\ell,\ell+1)}
 \simeq
 L^{\gbigf}_{\vecJbar_{0,2\ell+2}\cup\vecJbar_{0,2\ell+3}|\theta^u(\ell,\ell+1)}.
\]
We patch
$L^{\gbigf}_{\vecJbar_{0,2\ell}\cup\vecJbar_{0,2\ell+1}}$
$(\ell\in\seisuu)$
by using 
$-\Psi(\ell,\ell+1)$,
and we obtain a $2\pi\seisuu$-equivariant
local system with Stokes structure
$(L^{\gbigf},\vecnbigf)$ on $\real$.
We have
\[
 (L^{\gbigf},\vecnbigf)
 \simeq
 \gbigf^{(\infty,\infty)}_{+!}(L,\vecnbigf)/\Image\vecA_{\infty}
 \simeq
 (\gbigl^{\gbigf}(\nbigm),\vecnbigf).
\]

We can observe that
$(\Ltilde^{\gbigf},\vecnbigf)$
is isomorphic to
$(L^{\gbigf},\vecnbigf)$.

\subsection{}

Let $\nbigi=\{\alpha_1 x^{-3},\ldots,\alpha_nx^{-3}\}$
for some mutually distinct positive numbers $\alpha_i$.
Let $\nbigm$ be a meromorphic flat bundle on $(\proj^1,\infty)$
such that $\nbigi_{\infty}(\nbigm)=\nbigi$.
Let us study the canonical splittings and the Stokes matrices of
$(\gbigl^{\gbigf}(\nbigm),\vecnbigf)$.

\subsubsection{The index sets}

We set
$\langle 3\rangle'=2\cdot 3^{-\frac{3}{2}}$,
We set
\[
 \nbigi_0^{\circ}
 =\bigl\{
 \langle 3\rangle'\sqrt{-1}\alpha_i^{1/2}u^{-3/2}
 \bigr\},
 \quad
 \nbigi_1^{\circ}
 =\bigl\{
 -\langle 3\rangle'\sqrt{-1}\alpha_i^{1/2}u^{-3/2}
 \bigr\}.
\]
We set
$\nbigi^{\circ}=\nbigi_0^{\circ}\cup
\nbigi_1^{\circ}\cup\{0\}$.
We have
\[
(\gbigl^{\gbigf}(\nbigm),\vecnbigf)
\in
\Loc^{\St}(\nbigi_0^{\circ}\cup\nbigi_1^{\circ}).
\]

Let $\nu^0:\nbigi_0^{\circ}\simeq\nbigi$
and $\nu^1:\nbigi_1^{\circ}\simeq\nbigi$
denote the natural bijections
given by
\[
 \nu^0(\langle 3\rangle'\sqrt{-1}\alpha_i^{1/2}u^{-3/2})
=\nu^1(-\langle 3\rangle'\sqrt{-1}\alpha_i^{1/2}u^{-3/2})
 =\alpha_iu^{-3/2}.
\]

\subsubsection{Intervals}

We set $J_m=I(m\frac{\pi}{3},\frac{\pi}{6})$.
We have
$-\Re(\alpha_ix^{-3})<0$
on $J_{2\ell}$,
and
$-\Re(\alpha_ix^{-3})>0$
on $J_{2\ell+1}$
for any $\ell\in\seisuu$.

We set
$\vecJ_{0,m}=I(\frac{2}{3}m\pi+\frac{\pi}{3},\frac{\pi}{3})$
and
$\vecJ_{1,m}=I(\frac{2}{3}m\pi+\pi,\frac{\pi}{3})$
for any $m\in\seisuu$.
We have $\vecJ_{0,m}=\vecJ_{1,m-1}$.
For any $\gminia\in\nbigi_0^{\circ}$,
we have
$-\Re(\gminia)<0$ on $\vecJ_{0,2\ell}$,
and
$-\Re(\gminia)>0$ on $\vecJ_{0,2\ell+1}$.
For any $\gminib\in\nbigi_1^{\circ}$,
we have
$-\Re(\gminib)<0$ on $\vecJ_{1,2\ell}$,
and
$-\Re(\gminib)>0$ on $\vecJ_{1,2\ell+1}$.

For any interval $J=\{\theta_1<\theta<\theta_2\}$,
we set
$J_+=\{\theta_1<\theta\leq \theta_2\}$
and
$J_-=\{\theta_1\leq \theta<\theta_2\}$.

\subsubsection{Stokes matrices of $\nbigm$}

We set $(L,\vecnbigf)=
(L_{\infty}(\nbigm),\vecnbigf)$.
For any $\gminic\in\nbigi$,
let $\vece_{\gminic}$
denote a flat frame of $\Gr^{\vecnbigf}_{\gminic}(L)$.
For any $J_m$,
there exists the canonical splittings:
\[
 L_{|(J_m)_{\pm}}
=\bigoplus_{\gminic\in\nbigi}
 L_{(J_m)_{\pm},\gminic}.
\]
Let $\vece_{\gminic,(J_m)_{\pm}}$
denote the frame of
$L_{(J_m)_{\pm},\gminic}$
induced by $\vece_{\gminic}$.
Let $\vece_{(J_m)_{\pm}}$
denote the frame of
$L_{|(J_m)_{\pm}}$
induced by
$\vece_{\gminic,(J_m)_{\pm}}$
$(\gminic\in\nbigi)$.
We obtain the matrix $A_m$ by
\[
 \vece_{(J_m)_-}
 =\vece_{(J_m)_+}A_m.
\]
Note that
$H^0((J_m)_-,L_{(J_m)_-,\gminic})
=H^0((J_{m-1})_+,L_{(J_{m-1})_+,\gminic})$
in $H^0(\real,L)$
for any $m\in\seisuu$ and $\gminic\in\nbigi$.
We have
$\vece_{(J_m)_-,\gminic}
=\vece_{(J_{m-1})_+,\gminic}$
as tuples of sections of $H^0(\real,L)$.

\subsubsection{The induced frames}

Let $\gminia\in\nbigi_0^{\circ}$.
The restrictions
$\Gr^{\vecnbigf}_{\gminia}(\gbigl^{\gbigf}(\nbigm))_{\vecJbar_{0,2\ell}}$
and
$\Gr^{\vecnbigf}_{\gminia}(\gbigl^{\gbigf}(\nbigm))_{\vecJbar_{0,2\ell+1}}$
are  naturally isomorphic to the pull back of
$\Gr^{\vecnbigf}_{\nu^0(\gminia)}(L)_{|\Jbar_{-4\ell-1}}$
and
$\Gr^{\vecnbigf}_{\nu^0(\gminia)}(L)_{|\Jbar_{-4\ell}}$
by the map
$\kappa^0_{\ell}(\theta^u)=\frac{1}{2}(\theta^u-4\ell\pi-\pi)$.
The induced frames are denoted by
$\vece^0_{\gminia,\vecJbar_{0,2\ell}}$
and
$\vece^0_{\gminia,\vecJbar_{0,2\ell+1}}$,
respectively.
They induce the same frame of
$H^0(\real,\Gr^{\vecnbigf}_{\gminia}(\gbigl^{\gbigf}(\nbigm)))$.
We have
$\vece^0_{\gminia,\vecJbar_{0,2\ell}}
=-\vece^0_{\gminia,\vecJbar_{0,2\ell-1}}$
as tuples of sections of 
$H^0(\real,\Gr^{\vecnbigf}_{\gminia}(\gbigl^{\gbigf}(\nbigm)))$.

Let $\gminib\in\nbigi_0^{\circ}$.
The restrictions
$\Gr^{\vecnbigf}_{\gminib}(\gbigl^{\gbigf}(\nbigm))_{|\vecJ_{1,2\ell}}$
and 
$\Gr^{\vecnbigf}_{\gminib}(\gbigl^{\gbigf}(\nbigm))_{|\vecJ_{1,2\ell+1}}$
are naturally isomorphic to the pull back of
$\Gr^{\vecnbigf}_{\nu^1(\gminib)}(L)_{|\Jbar_{-4\ell-3}}$
and
$\Gr^{\vecnbigf}_{\nu^1(\gminib)}(L)_{|\Jbar_{-4\ell-2}}$
by $\kappa^1_{\ell}(\theta^u)=\frac{1}{2}(\theta^u-4\ell\pi-3\pi)$.
The induced frames are denoted by
$\vece^1_{\gminib,\vecJbar_{1,2\ell}}$
and
$\vece^1_{\gminib,\vecJbar_{1,2\ell+1}}$,
respectively.
They induce the same frame of
$H^0(\real,\Gr^{\vecnbigf}_{\gminib}(\gbigl^{\gbigf}(\nbigm)))$.
We have
$\vece^1_{\gminia,\vecJbar_{1,2\ell}}
=-\vece^1_{\gminia,\vecJbar_{1,2\ell-1}}$
as tuples of sections of 
$H^0(\real,\Gr^{\vecnbigf}_{\gminib}(\gbigl^{\gbigf}(\nbigm)))$.

\subsubsection{Stokes matrices of $\Fourier_+(\nbigm)$}

For any
$\vecJ_{0,m}=\vecJ_{1,m-1}$,
there exist the canonical splittings
\[
 \gbigl^{\gbigf}(\nbigm)
 _{|(\vecJ_{0,m})_{\pm}}
=\bigoplus_{\gminia\in\nbigi_0^{\circ}}
 \gbigl^{\gbigf}(\nbigm)
 _{(\vecJ_{0,m})_{\pm},\gminia}
 \oplus
 \bigoplus_{\gminib\in\nbigi_1^{\circ}}
 \gbigl^{\gbigf}(\nbigm)
 _{(\vecJ_{1,m-1})_{\pm},\gminib}.
\]
Let $\vece^0_{(\vecJ_{0,m})_{\pm}}$
denote the frame of
$\bigoplus_{\gminia}\gbigl^{\gbigf}(\nbigm)
_{(\vecJ_{0,m})_{\pm},\gminia}$
induced by $\vece^0_{\gminia,\vecJbar_{0,m}}$
$(\gminia\in\nbigi_0^{\circ})$.
Similarly,
let $\vece^1_{(\vecJ_{1,m-1})_{\pm}}$
denote the frame of
$\bigoplus_{\gminia}
\gbigl^{\gbigf}(\nbigm)
 _{(\vecJ_{1,m-1})_{\pm},\gminib}$
 induced by
 $\vece^1_{\gminib,\vecJbar_{1,m-1}}$
 $(\gminib\in\nbigi^{\circ}_1)$.

\begin{prop}
\label{prop;24.4.21.1}
\begin{equation}
\label{eq;24.4.21.2}
 \Bigl(\vece^0_{(\vecJ_{0,2\ell})_-},
 \vece^1_{(\vecJ_{1,2\ell-1})_-}\Bigr)
=\Bigl(\vece^0_{(\vecJ_{0,2\ell})_+},
 \vece^1_{(\vecJ_{1,2\ell-1})_+}
 \Bigr)
 \left(
 \begin{array}{cc}
  A_{-4\ell-1} & -A^{-1}_{4\ell}A^{-1}_{-4\ell+1}\\
  0 & A_{-4\ell+2}
 \end{array}
 \right), 
\end{equation}
\begin{equation}
\label{eq;24.4.21.3}
 \Bigl(\vece^0_{(\vecJ_{0,2\ell+1})_-},
 \vece^1_{(\vecJ_{1,2\ell})_-}
 \Bigr)
=\Bigl(\vece^0_{(\vecJ_{0,2\ell+1})_+},
 \vece^1_{(\vecJ_{1,2\ell})_+}
 \Bigr)
 \left(
 \begin{array}{cc}
  A_{-4\ell} & 0   \\
  -A^{-1}_{4\ell-2}A^{-1}_{-4\ell-1}
   & A_{-4\ell-3}
 \end{array}
 \right).
\end{equation}
\end{prop}

\subsubsection{Local Fourier transform}

By the local Fourier transform of $(L,\vecnbigf)$
in \S\ref{subsection;24.4.5.120},
we obtain
\[
 \gbigf^{(\infty,\infty)}_{+,!}(L,\vecnbigf)
 =\bigl(
 \gbigq^{\infty}_{!}(L,\vecnbigf),
 \vecnbigf
 \bigr)
 \in\Loc^{\St}(\nbigi^{\circ})
\]
which are isomorphic to
$(\gbigl^{\gbigf}(\nbigm(! 0)),\vecnbigf)$.

\subsubsection{}

We consider the vector space
$\bigoplus_{\ell\in\seisuu}
H^0(J_{2\ell+1},L)$.
For $v\in H^0(\real,L)$,
the induced element $H^0(J_{2\ell+1},L)$
is denoted by
$\langle J_{2\ell+1},v\rangle$.
We define
$\gbigq^{\infty}_!(L,\vecnbigf)$
as the quotient of 
$\bigoplus_{\ell\in\seisuu}
H^0(J_{2\ell+1},L)$
by the equivalence relation generated by
\[
 \langle J_{2\ell+1},v\rangle
 \sim
 \langle
 J_{2\ell-5},v\rangle.
\]
The $2\pi\seisuu$-equivariant local system
$\gbigq_!^{\infty}(L,\vecnbigf)_{\real}$
equals the $2\pi\seisuu$-equivariant local system
induced by 
$\gbigq_!^{\infty}(L,\vecnbigf)$
with the trivial action.

\subsubsection{Some maps}

We naturally identify
$H^0(\real,L)=H^0(J_{m},L)$
for any interval $J$.
The maps
$\vecB_{\vecJ_{0,2\ell}}:
H^0(J_{-4\ell-1},L)
\to \gbigq^{\infty}_!(L,\vecnbigf)$
and
$\vecA_{(\vecJ_{0,2\ell+1})_{\pm}}:
H^0(J_{-4\ell},L)
\to \gbigq^{\infty}_!(L,\vecnbigf)$
are given as follows.
\[
 \vecB_{\vecJ_{0,2\ell}}(v)
=\vecA_{(\vecJ_{0,2\ell+1})_+}(v)
=\langle J_{-4\ell-1},v\rangle,
\quad
  \vecA_{(\vecJ_{0,2\ell+1})_-}(v)
 =-\langle J_{-4\ell+1},v\rangle.
\]
The maps
$\vecB_{\vecJ_{1,2\ell}}:
H^0(J_{-4\ell-3},L)
\to \gbigq^{\infty}_!(L,\vecnbigf)$
and
$\vecA_{(\vecJ_{1,2\ell+1})_{\pm}}:
H^0(J_{-4\ell-2},L)
\to \gbigq^{\infty}_!(L,\vecnbigf)$
are given as follows.
\[
 \vecB_{\vecJ_{1,2\ell}}(v)
=\vecA_{(\vecJ_{1,2\ell+1})_+}(v)
=\langle J_{-4\ell-3},v\rangle,
\quad
  \vecA_{(\vecJ_{1,2\ell+1})_-}(v)
 =-\langle J_{-4\ell-},v\rangle.
\]
Let $\vecA_{\infty}:H^0(\real,L)\to \gbigq_!^{\infty}(L,\vecnbigf)$
be defined by
$\vecA_{\infty}(v)
:=\langle J_1,v\rangle
 +\langle J_3,v\rangle
 +\langle J_5,v\rangle$
for any $v\in H^0(\real,L)$, 
which equals
$\langle J_{2\ell+1},v\rangle
+\langle J_{2\ell+3},v\rangle
+\langle J_{2\ell+5},v\rangle$
for any $\ell\in\seisuu$.
We have
\[
 \vecA^{(\vecJ_{a,m})_+}_{\infty}
=\vecA^{(\vecJ_{a,m})_-}_{\infty}
=\vecA_{\infty}.
\]

\subsubsection{Stokes filtrations}

The map $\nu^a$ $(a=0,1)$ induce following
isomorphisms of partially ordered sets
for any $\theta^u\in \vecJbar_{a,2\ell}\cup\vecJbar_{a,2\ell+1}$:
\[
 (\nbigi^{\circ}_0,\leq_{\theta^u})
 \simeq
 (\nbigi,\leq_{\kappa^a_{\ell}(\theta^u)}).
\]

For any $\theta\in\real$,
we identify $L_{|\theta}$ with $H^0(\real,L)$.
We obtain the filtrations $\nbigf^{\theta}$
on $H^0(\real,L)$.

We identify $\gbigq^{\infty}_!(L,\vecnbigf)_{\real|\theta^u}$
with $\gbigq^{\infty}_!(L,\vecnbigf)$.
The Stokes filtration $\nbigf^{\theta^u}$
is given as follows.
\begin{itemize}
 \item 
For $\theta^u\in \vecJ_{0,2\ell}=\vecJ_{1,2\ell-1}$,
\[
       \nbigf^{\theta^u}_{\gminia}\gbigq^{\infty}_!(L,\vecnbigf)
       =\vecB_{\vecJ_{0,2\ell}}
       \bigl(
       \nbigf^{\kappa^0_{\ell}(\theta^u)}_{\nu^0(\gminia)}
       H^0(\real,L)
       \bigr)
       \quad
       (\gminia\in\nbigi^{\circ}_0),
\]
\[
       \nbigf^{\theta^u}_{0}\gbigq^{\infty}_!(L,\vecnbigf)
       =\Image \vecB_{\vecJ_{0,2\ell}}
       \oplus
       \Image\vecA_{\infty},
\]       
\[
       \nbigf^{\theta^u}_{\gminib}\gbigq^{\infty}_!(L,\vecnbigf)
       =\nbigf^{\theta^u}_{0}\gbigq^{\infty}_!(L,\vecnbigf)
       \oplus
       \vecA_{(\vecJ_{1,2\ell-1})_+}
       \Bigl(
       \nbigf^{\kappa_{\ell-1}^1(\theta^u)}_{\nu^1(\gminib)}
       H^0(\real,L)
       \Bigr)
       \quad
       (\gminib\in\nbigi^{\circ}_1).
\]
 \item For $\theta^u\in\vecJ_{0,2\ell+1}=\vecJ_{1,2\ell}$,
\[
       \nbigf^{\theta^u}_{\gminib}\gbigq^{\infty}_!(L,\vecnbigf)
       =\vecB_{\vecJ_{1,2\ell}}
       \bigl(
       \nbigf^{\kappa^1_{\ell}(\theta^u)}_{\nu^1(\gminib)}
       H^0(\real,L)
       \bigr)
       \quad
       (\gminib\in\nbigi^{\circ}_1),
\]
\[
       \nbigf^{\theta^u}_{0}\gbigq^{\infty}_!(L,\vecnbigf)
       =\Image \vecB_{\vecJ_{1,2\ell}}
       \oplus
       \Image\vecA_{\infty},
\]       
\[
       \nbigf^{\theta^u}_{\gminia}\gbigq^{\infty}_!(L,\vecnbigf)
       =\nbigf^{\theta^u}_{0}\gbigq^{\infty}_!(L,\vecnbigf)
       \oplus
       \vecA_{(\vecJ_{0,2\ell+1})_+}
       \Bigl(
       \nbigf^{\kappa_{\ell}^0(\theta^u)}_{\nu^0(\gminia)}
       H^0(\real,L)
       \Bigr)
       \quad
       (\gminia\in\nbigi^{\circ}_0).
\]
 \item Let $\theta^u\in
       \vecJbar_{0,2\ell}\cap\vecJbar_{0,2\ell+1}$.
       We have the decompositions
\[
       H^0(\real,L)
       =\bigoplus_{\gminia\in\nbigi}
       \nbigf^{\kappa^0_{\ell}(\theta^u)}_{\nu^0(\gminia)}
       H^0(\real,L)
       =\bigoplus_{\gminia\in\nbigi}
       \nbigf^{\kappa^1_{\ell}(\theta^u)}_{\nu^1(\gminia)}
       H^0(\real,L).
\]       
       Note that
       $\nbigf^{\kappa^1_{\ell}(\theta^u)}_{\nu^1(\gminia)}
       H^0(\real,L)
       =\nbigf^{\kappa^1_{\ell-1}(\theta^u)}_{\nu^1(\gminia)}
       H^0(\real,L)$
       because the monodromy of $L$ is trivial.
       The Stokes filtration of $\gbigq^{\infty}_!(L,\vecnbigf)$
       at $\theta^u$
       is given by the splitting
\begin{multline}
       \gbigq^{\infty}_!(L,\vecnbigf)
 =\\
 \Image\vecA_{\infty}
       \oplus
       \bigoplus_{\gminia\in\nbigi_0^{\circ}}
       \vecB_{\vecJ_{0,2\ell}}\bigl(
       \nbigf^{\kappa^0_{\ell}(\theta^u)}_{\nu^0(\gminia)}
       H^0(\real,L)
       \bigr)
       \oplus
       \bigoplus_{\gminib\in\nbigi_1^{\circ}}
       \vecB_{\vecJ_{1,2\ell}}\bigl(
       \nbigf^{\kappa^1_{\ell}(\theta^u)}_{\nu^1(\gminib)}
       H^0(\real,L)
       \bigr).
\end{multline}
       We have a similar description of
       the Stokes filtration at
       $\theta^u\in\vecJbar_{1,2\ell}\cap\vecJbar_{1,2\ell+1}$.
\end{itemize}

\subsubsection{Proof of Proposition \ref{prop;24.4.21.1}}

There exists the natural isomorphism
$(\gbigl^{\gbigf}(\nbigm),\vecnbigf)
\simeq
 \gbigf^{(\infty,\infty)}_{+!}(L,\vecnbigf)/\Image\vecA_{\infty}$.

The frames
$\vece^0_{(\vecJbar_{0,2\ell})_{\pm}}$
are given by
$\vecB_{\vecJ_{0,2\ell}}\bigl(
\vece_{(J_{-4\ell-1})_{\pm}}
\bigr)$.
The frames
$\vece^1_{(\vecJbar_{1,2\ell-1})_{\pm}}$
are given by
$\vecA_{(\vecJ_{1,2\ell-1})_{\mp}}\bigl(
\vece_{(J_{-4\ell+2})_{\pm}}
\bigr)$.
Note that
\begin{multline}
 \vecA_{(\vecJ_{1,2\ell-1})_+}(v)
-\vecA_{(\vecJ_{1,2\ell-1})_-}(v)
=\langle J_{-4\ell+1},v\rangle
+\langle J_{-4\ell+3},v\rangle
\\
 \equiv
-\langle J_{-4\ell-1},v\rangle
=-\vecB_{\vecJ_{0,2\ell}}(v).
\end{multline}
Because
$\vece_{(J_{-4\ell+2})_-}
=\vece_{(J_{-4\ell-1})_+}A_{-4\ell}^{-1}A_{-4\ell+1}^{-1}$,
we obtain (\ref{eq;24.4.21.2}).
We can obtain (\ref{eq;24.4.21.3}) similarly.
\hfill\qed

\subsection{}

Let $\nbigi=\{\alpha_1z^{-1},\ldots,\alpha_Nz^{-1}\}$
for some mutually distinct positive numbers $\alpha_i$.
Let $\nbigv$ be a basic meromorphic flat bundle
of level $(0,1)$ with
$\nbigi_0(\nbigv)=\nbigi$.
We have $\nbigv[!0]=\nbigv$.
Let us study
$(\gbigl_!^{\gbigf}(\nbigv),\vecnbigf)$.

\subsubsection{The index sets}

We set
\[
 \nbigi_0^{\circ}
 =\{2\alpha_i^{1/2}u^{-1/2}\},
 \quad
 \nbigi_1^{\circ}
 =\{-2\alpha_i^{1/2}u^{-1/2}\},
 \quad
 \nbigi^{\circ}
 =\nbigi_0^{\circ}
 \cup \nbigi_1^{\circ}.
\]
Let
$\nu^a:\nbigi_a^{\circ}\simeq\nbigi$ $(a=0,1)$
be the bijection given by
\[
 \nu^0(2\alpha_i^{1/2}u^{-1/2})
 =\alpha_iz^{-1},
 \quad
 \nu^1(-2\alpha_i^{1/2}u^{-1/2})
 =\alpha_iz^{-1}.
\]

\subsubsection{The intervals}

We set $J_m=I(m\pi,\frac{\pi}{2})$.
We have
$-\Re(\alpha_i e^{-\sqrt{-1}\theta})<0$ on $J_{2\ell}$
and 
$-\Re(\alpha_i e^{-\sqrt{-1}\theta})>0$ on $J_{2\ell+1}$.

We set
$\vecJ_{0,m}=I(2m\pi,\pi)$
and
$\vecJ_{1,m}=I(2\pi+2m\pi,\pi)$.
We have
$\vecJ_{0,m}=\vecJ_{1,m-1}$.
We have
$-\Re(2\alpha_i^{1/2}e^{-\sqrt{-1}\theta^u/2})<0$
on $\vecJ_{0,2\ell}$
and
$-\Re(2\alpha_i^{1/2}e^{-\sqrt{-1}\theta^u/2})>0$
on $\vecJ_{0,2\ell+1}$.
We have
$-\Re(-2\alpha_i^{1/2}e^{-\sqrt{-1}\theta^u/2})<0$
on $\vecJ_{1,2\ell}$
and
$-\Re(-2\alpha_i^{1/2}e^{-\sqrt{-1}\theta^u/2})>0$
on $\vecJ_{1,2\ell+1}$.

\subsubsection{Stokes matrices of $\nbigv$}

We set
$(L,\vecnbigf)=(L_0(\nbigv),\vecnbigf)$.
For any $\gminic\in\nbigi$,
let $\vece_{\gminic}$ denote a flat frame of
$\Gr^{\vecnbigf}_{\gminic}(L)$.
We obtain the matrices $G_{\gminic}$ determined by
$(\Tbb^{\ast})^{-1}\vece_{\gminic}=\vece_{\gminic}G_{\gminic}$.

For any $J_m$,
there exists the canonical splittings:
\[
 L_{|(J_{m})_{\pm}}
 =\bigoplus_{\gminic\in\nbigi}
 L_{(J_m)_{\pm},\gminic}.
\]
Let $\vece_{\gminic,(J_m)_{\pm}}$
denote the frame of
$L_{(J_m)_{\pm},\gminic}$
induced by $\vece_{\gminic}$.
Let $\vece_{(J_m)_{\pm}}$
denote the frame of $L_{|(J_m)_{\pm}}$
induced by $\vece_{\gminic,(J_m)_{\pm}}$
$(\gminic\in\nbigi)$.
We obtain the matrix $A_m$
determined by
\[
 \vece_{(J_m)-}
 =\vece_{(J_m)+}A_m.
\]
Note that
$H^0((J_m)_-,L_{(J_m)_-,\gminic})
=H^0((J_{m-1})_+,L_{(J_{m-1})_+,\gminic})$
in $H^0(\real,L)$
for any $m\in\seisuu$ and $\gminic\in\nbigi$.
We have
$\vece_{(J_m)_-,\gminic}
=\vece_{(J_{m-1})_+,\gminic}$
as tuples of sections of $H^0(\real,L)$.

We set $G=\bigoplus G_{\gminic}$.
We have
$(\Tbb^{\ast})^{-1}\vece_{(J_{m-4})_{\pm}}
=\vece_{(J_m)_{\pm}}G$.

\subsubsection{The induced frames}

Let $\gminia\in\nbigi_0^{\circ}$.
The restrictions
$\Gr^{\vecnbigf}_{\gminia}(\gbigl_!^{\gbigf}(\nbigv))_{\vecJbar_{0,2\ell}}$
and
$\Gr^{\vecnbigf}_{\gminia}(\gbigl^{\gbigf}_!(\nbigv))_{\vecJbar_{0,2\ell+1}}$
are  naturally isomorphic to the pull back of
$\Gr^{\vecnbigf}_{\nu^0(\gminia)}(L)_{|\Jbar_{4\ell}}$
and
$\Gr^{\vecnbigf}_{\nu^0(\gminia)}(L)_{|\Jbar_{4\ell+1}}$
by the map
$\kappa^0_{\ell}(\theta^u)=\frac{1}{2}(\theta^u+4\ell\pi)$.
The induced frames are denoted by
$\vece^0_{\gminia,\vecJbar_{0,2\ell}}$
and
$\vece^0_{\gminia,\vecJbar_{0,2\ell+1}}$,
respectively.
They induce the same frame of
$H^0(\real,\Gr^{\vecnbigf}_{\gminia}(\gbigl_!^{\gbigf}(\nbigv)))$.
We have
$\vece^0_{\gminia,\vecJbar_{0,2\ell}}
=-\vece^0_{\gminia,\vecJbar_{0,2\ell-1}}$
as tuples of sections of 
$H^0(\real,\Gr^{\vecnbigf}_{\gminia}(\gbigl_!^{\gbigf}(\nbigv)))$.

Let $\gminib\in\nbigi_1^{\circ}$.
The restrictions
$\Gr^{\vecnbigf}_{\gminib}(\gbigl_!^{\gbigf}(\nbigv))_{|\vecJ_{1,2\ell}}$
and 
$\Gr^{\vecnbigf}_{\gminib}(\gbigl_!^{\gbigf}(\nbigv))_{|\vecJ_{1,2\ell+1}}$
are naturally isomorphic to the pull back of
$\Gr^{\vecnbigf}_{\nu^1(\gminib)}(L)_{|\Jbar_{4\ell+2}}$
and
$\Gr^{\vecnbigf}_{\nu^1(\gminib)}(L)_{|\Jbar_{4\ell+3}}$
by $\kappa^1_{\ell}(\theta^u)=\frac{1}{2}(\theta^u+(4\ell+2)\pi)$.
The induced frames are denoted by
$\vece^1_{\gminib,\vecJbar_{1,2\ell}}$
and
$\vece^1_{\gminib,\vecJbar_{1,2\ell+1}}$,
respectively.
They induce the same frame of
$H^0(\real,\Gr^{\vecnbigf}_{\gminib}(\gbigl_!^{\gbigf}(\nbigv)))$.
We have
$\vece^1_{\gminia,\vecJbar_{1,2\ell}}
=-\vece^1_{\gminia,\vecJbar_{1,2\ell-1}}$
as tuples of sections of 
$H^0(\real,\Gr^{\vecnbigf}_{\gminib}(\gbigl_!^{\gbigf}(\nbigv)))$.

\subsubsection{Stokes matrices of $\Fourier_+(\nbigv(!0))$}

For any
$\vecJ_{0,m}=\vecJ_{1,m-1}$,
there exist the canonical splittings
\[
 \gbigl^{\gbigf}_!(\nbigv)
 _{|(\vecJ_{0,m})_{\pm}}
=\bigoplus_{\gminia\in\nbigi_0^{\circ}}
 \gbigl^{\gbigf}_!(\nbigv)
 _{(\vecJ_{0,m})_{\pm},\gminia}
 \oplus
 \bigoplus_{\gminib\in\nbigi_1^{\circ}}
 \gbigl^{\gbigf}_!(\nbigv)
 _{(\vecJ_{1,m-1})_{\pm},\gminib}.
\]
Let $\vece^0_{(\vecJ_{0,m})_{\pm}}$
denote the frame of
$\bigoplus_{\gminia}\gbigl^{\gbigf}_!(\nbigv)
_{(\vecJ_{0,m})_{\pm},\gminia}$
induced by $\vece^0_{\gminia,\vecJbar_{0,m}}$
$(\gminia\in\nbigi_0^{\circ})$.
Similarly,
let $\vece^1_{(\vecJ_{1,m-1})_{\pm}}$
denote the frame of
$\bigoplus_{\gminia}
\gbigl^{\gbigf}_!(\nbigv)
 _{(\vecJ_{1,m-1})_{\pm},\gminib}$
 induced by
 $\vece^1_{\gminib,\vecJbar_{1,m-1}}$
 $(\gminib\in\nbigi^{\circ}_1)$.

\begin{prop}
\label{prop;24.4.22.1}
{\small
 \begin{equation}
\label{eq;24.4.22.2}
 \Bigl(\vece^0_{(\vecJ_{0,2\ell})_-},
 \vece^1_{(\vecJ_{1,2\ell-1})_-}\Bigr)
=\Bigl(\vece^0_{(\vecJ_{0,2\ell})_+},
 \vece^1_{(\vecJ_{1,2\ell-1})_+}
 \Bigr)
 \left(
 \begin{array}{cc}
  A_{4\ell} & A_{4\ell}A_{4\ell-1}+GA_{4\ell-3}^{-1}A_{4\ell-2}^{-1}
   \\
  0 & A_{4\ell-2}
 \end{array}
 \right), 
 \end{equation}
\begin{equation}
\label{eq;24.4.22.3}
 \Bigl(\vece^0_{(\vecJ_{0,2\ell+1})_-},
 \vece^1_{(\vecJ_{1,2\ell})_-}
 \Bigr)
=\Bigl(\vece^0_{(\vecJ_{0,2\ell+1})_+},
 \vece^1_{(\vecJ_{1,2\ell})_+}
 \Bigr)
 \left(
 \begin{array}{cc}
  A_{4\ell+2} & 0   \\
  G A_{4\ell-1}^{-1}A_{4\ell}^{-1}
 +A_{4\ell+2}A_{4\ell+1}
   & A_{4\ell}
 \end{array}
 \right).
\end{equation}
 }
\end{prop}

\subsubsection{Local Fourier transform}

By the local Fourier transform of $(L,\vecnbigf)$
in \S\ref{subsection;24.4.5.100},
we obtain
\[
 \gbigf^{(0,\infty)}_{+,!}(L,\vecnbigf)
 =\bigl(
 \gbigq^{0}_{!}(L,\vecnbigf),
 \vecnbigf
 \bigr)
 \in\Loc^{\St}(\nbigi^{\circ})
\]
which are isomorphic to
$(\gbigl^{\gbigf}_!(\nbigv),\vecnbigf)$.

\subsubsection{}

We consider the vector space
$H^0(\real,L)\oplus
\bigoplus_{\ell\in\seisuu}
H^0(J_{2\ell},L)$.
We define the $\real$-vector space
$\gbigq^0_!(L,\vecnbigf)$
as the quotient of
$\bigoplus_{\ell\in\seisuu}
H^0(J_{2\ell},L)$
by the equivalence relation generated by
\[
 \langle J_{2\ell},v\rangle
-\langle
 J_{2\ell+2},(\Tbb^{\ast})^{-1}v\rangle
 \sim
 v.
\]
The $2\pi\seisuu$-equivariant local system
$\gbigq_!^{0}(L,\vecnbigf)_{\real}$
equals the natural $2\pi\seisuu$-equivariant local system
induced by 
$\gbigq_!^{0}(L,\vecnbigf)$.
(See \S\ref{subsection;24.4.5.100}.)

\subsubsection{Some maps}

We naturally identify
$H^0(\real,L)$ with $H^0(J_m,L)$ for any $J_m$.
The maps
$\vecA_{\vecJ_{0,2\ell}}:H^0(J_{4\ell},L)\to\gbigq^0_!(L,\vecnbigf)$
and
$\vecB_{(\vecJ_{0,2\ell+1})_{\pm}}:
H^0(J_{4\ell+1},L)\to\gbigq^0_!(L,\vecnbigf)$
are given as follows:
\[
 \vecA_{\vecJ_{0,2\ell}}(v)
 =\vecB_{(\vecJ_{0,2\ell+1})_-}(v)
 =\langle J_{4\ell},v\rangle,
 \quad
 \vecB_{(\vecJ_{0,2\ell+1})_+}(v)
 =-\langle J_{4\ell+4},(\Tbb^{\ast})^{-1}(v)\rangle.
\]
The maps
$\vecA_{\vecJ_{1,2\ell}}:
H^0(J_{4\ell+2},L)\to \gbigq^0_{!}(L,\vecnbigf)$
and
$\vecB_{(\vecJ_{1,2\ell+1})_{\pm}}:
H^0(J_{4\ell+3},L)\to \gbigq^0_{!}(L,\vecnbigf)$
are given as follows:
\[
 \vecA_{\vecJ_{1,2\ell}}(v)
 =\vecB_{(\vecJ_{1,2\ell+1})_-}(v)
=\langle J_{4\ell+2},v\rangle,
 \quad
 \vecB_{(\vecJ_{1,2\ell+1})_+}(v)
=-\langle J_{4\ell+6},(\Tbb^{\ast})^{-1}(v)\rangle.
\]
\subsubsection{Stokes filtrations}

The maps $\nu^a$ induce
the following isomorphisms of
partially ordered sets
for any $\theta^u\in \vecJbar_{a,2\ell}\cup\vecJbar_{a,2\ell+1}$:
\[
 (\nbigi^{\circ}_a,\leq_{\theta^u})
 \simeq
 (\nbigi,\leq_{\kappa^a_{\ell}(\theta^u)}).
\]
For any $\theta\in\real$,
we identify
$H^0(\real,L)$ with $L_{|\theta}$.
We obtain the filtrations
$\nbigf^{\theta}$ on $H^0(\real,L)$.

We identify
$\gbigq^0_{!}(L,\vecnbigf)_{\real|\theta^u}$
with $\gbigq^0_!(L,\vecnbigf)$.
The Stokes filtration $\nbigf^{\theta^u}$
is given as follows:
\begin{itemize}
 \item For $\theta^u\in\vecJ_{0,2\ell}=\vecJ_{1,2\ell-1}$,
       \[
       \nbigf^{\theta^u}_{\gminia}\gbigq^{0}_!(L,\vecnbigf)
       =\vecA_{\vecJ_{0,2\ell}}
       \bigl(
       \nbigf^{\kappa^0_{\ell}(\theta^u)}_{\nu^0(\gminia)}
       H^0(\real,L)
       \bigr)
       \quad
       (\gminia\in\nbigi^{\circ}_0),
\]
\[
       \nbigf^{\theta^u}_{\gminib}\gbigq^{0}_!(L,\vecnbigf)
       =\Image\vecA_{\vecJ_{0,2\ell}}
       \oplus
       \vecB_{(\vecJ_{1,2\ell-1})_+}
       \Bigl(
       \nbigf^{\kappa_{\ell-1}^1(\theta^u)}_{\nu^1(\gminib)}
       H^0(\real,L)
       \Bigr)
       \quad
       (\gminib\in\nbigi^{\circ}_1).
\]
 \item For $\theta^u\in\vecJ_{0,2\ell+1}=\vecJ_{1,2\ell}$,
\[
       \nbigf^{\theta^u}_{\gminib}\gbigq^{0}_!(L,\vecnbigf)
       =\vecA_{\vecJ_{1,2\ell}}
       \bigl(
       \nbigf^{\kappa^1_{\ell}(\theta^u)}_{\nu^1(\gminib)}
       H^0(\real,L)
       \bigr)
       \quad
       (\gminib\in\nbigi^{\circ}_1),
\]
\[
       \nbigf^{\theta^u}_{\gminia}\gbigq^{0}_!(L,\vecnbigf)
       =\Image \vecA_{\vecJ_{0,2\ell}}
       \oplus
       \vecB_{(\vecJ_{0,2\ell+1})_+}
       \Bigl(
       \nbigf^{\kappa_{\ell}^0(\theta^u)}_{\nu^0(\gminia)}
       H^0(\real,L)
       \Bigr)
       \quad
       (\gminia\in\nbigi^{\circ}_0).
\]
\item Let $\theta^u\in
       \vecJbar_{0,2\ell}\cap\vecJbar_{0,2\ell+1}$.
       We have the decompositions
\[
       H^0(\real,L)
       =\bigoplus_{\gminia\in\nbigi}
       \nbigf^{\kappa^0_{\ell}(\theta^u)}_{\nu^0(\gminia)}
       H^0(\real,L)
       =\bigoplus_{\gminia\in\nbigi}
       \nbigf^{\kappa^1_{\ell}(\theta^u)}_{\nu^1(\gminia)}
       H^0(\real,L).
\]       
       The Stokes filtration of $\gbigq^{0}_!(L,\vecnbigf)$
       at $\theta^u$
       is given by the splitting
\begin{multline}
       \gbigq^{0}_!(L,\vecnbigf)
 =
       \bigoplus_{\gminia\in\nbigi_0^{\circ}}
       \vecA_{\vecJ_{0,2\ell}}\bigl(
       \nbigf^{\kappa^0_{\ell}(\theta^u)}_{\nu^0(\gminia)}
       H^0(\real,L)
       \bigr)
       \oplus
       \bigoplus_{\gminib\in\nbigi_1^{\circ}}
       \vecA_{\vecJ_{1,2\ell}}\bigl(
       \nbigf^{\kappa^1_{\ell}(\theta^u)}_{\nu^1(\gminib)}
       H^0(\real,L)
       \bigr).
\end{multline}
       We have a similar description of
       the Stokes filtration at
       $\theta^u\in\vecJbar_{1,2\ell}\cap\vecJbar_{1,2\ell+1}$.
\end{itemize}

\subsubsection{Proof of Proposition \ref{prop;24.4.22.1}}

The frames
$\vece^0_{(\vecJbar_{0,2\ell})_{\pm}}$
and $\vece^1_{(\vecJbar_{1,2\ell-1})_{\pm}}$
are given by
$\vecA_{\vecJ_{0,2\ell}}\bigl(
\vece_{(J_{4\ell})_{\pm}}
\bigr)$
and 
$\vecB_{(\vecJ_{1,2\ell-1})_{\pm}}\bigl(
\vece_{(J_{4\ell-1})_{\pm}}
\bigr)$,
respectively.
Note that
\begin{multline}
 \vecB_{(\vecJ_{1,2\ell-1})_+}(v)
-\vecB_{(\vecJ_{1,2\ell-1})_-}(v)
=\langle J_{4\ell-2},v\rangle
+\langle J_{4\ell+2},(\Tbb^{\ast})^{-1}v\rangle
\\
=\Bigl(
 \langle J_{4\ell},(\Tbb^{\ast})^{-1}(v)
 \rangle
 +v
 \Bigr)
+\Bigl(
 \langle J_{4\ell},v
 \rangle
 -v
 \Bigr)
=\langle J_{4\ell},(\Tbb^{\ast})^{-1}(v)+v\rangle.
\end{multline}
We have
$\vece_{(J_{4\ell-1})_-}
=\vece_{(J_{4\ell})_+}A_{4\ell}A_{4\ell-1}$,
and
$\vece_{(J_{4\ell-1})_-}
=\vece_{(J_{4\ell-4})_+}A_{4\ell-3}^{-1}A_{4\ell-2}^{-1}$.
Then, we obtain (\ref{eq;24.4.21.2}).
We can obtain (\ref{eq;24.4.22.3}) similarly.
\hfill\qed

\section{Outline of this monograph}

We devote \S\ref{section;18.6.3.10}
to preliminaries for Stokes structures.
We introduce the concept of Stokes shells
in \S\ref{section;18.6.3.11}.
In \S\ref{chapter;21.6.10.30},
we make preliminaries for meromorphic flat bundles.
We study the transforms of the set of ramified irregular values
induced by the local Fourier transform
in \S\ref{section;18.6.3.12}.
We study the Fourier transforms of basic meromorphic flat bundles
and the reduction procedures
in \S\ref{section;18.6.3.20}--\S\ref{section;20.10.24.4}
by postponing the proof for the  comparison of 
the explicitly defined filtrations and the Stokes filtrations
to \S\ref{section;18.6.3.21},
where the growth orders of flat sections are studied.
In \S\ref{section;24.4.20.1},
we study the Stokes structure of
the Fourier transform for $\nbigd$-modules.

\section{Acknowledgement}

I thank 
Davide Barco,
Indranil Biswas,
Andrea D'Agnolo,
H\'{e}l\`{e}ne Esnault,
Marco Hien,
Kazuki Hiroe,
Masaki Kashiwara,
Giovanni Morando,
Toshio Oshima,
Takeshi Saito,
Carlos Simpson,
Szilard Szabo,
and Daisuke Yamakawa
for some discussions.
I thank Ting Xue and Kari Vilonen
for their kindness.
I thank Philip Boalch for asking about this study.
I appreciate some communications with
Andreas Hohl and Jean Dou\c{c}ot.
I thank Claude Sabbah for his kindness
and discussions on many opportunities.
I thank Akira Ishii and Yoshifumi Tsuchimoto
for their encouragement.
I appreciate the anonymous reviewers
for their comments on earlier versions of this manuscript.
I acknowledge the kind hospitality at 
the Tata Institute of Fundamental Research
and 
the International Center for Theoretical Sciences,
where a part of this study was done.

I am partially supported by
the Grant-in-Aid for Scientific Research (S) (No. 17H06127),
the Grant-in-Aid for Scientific Research (S) (No. 16H06335),
the Grant-in-Aid for Scientific Research (A) (No. 21H04429),
the Grant-in-Aid for Scientific Research (A) (No. 22H00094),
the Grant-in-Aid for Scientific Research (A) (No. 23H00083),
the Grant-in-Aid for Scientific Research (C) (No. 15K04843),
the Grant-in-Aid for Scientific Research (C) (No. 20K03609),
and 
the Grant-in-Aid for Scientific Research (C) (No. 25K06973),
Japan Society for the Promotion of Science.
I am also partially supported by the Research Institute for Mathematical
Sciences, an International Joint Usage/Research Center located in Kyoto
University.

\chapter{Preliminary for Stokes structures}
\label{section;18.6.3.10}

\section[Family of partially ordered sets and Stokes structure]{Family of partially ordered sets
 and Stokes structure of local systems}

\subsection{General case}
\label{subsection;18.4.3.2}

Let $G$ be a discrete group.
Let $Y$ be a manifold 
with a proper $G$-action
which may have a boundary.
Let $\pi:\vecnbigi\lrarr Y$
be a $G$-equivariant proper continuous map of manifolds
which is locally a homeomorphism.
The fibers $\pi^{-1}(y)$ $(y\in Y)$
are denoted by $\nbigi_y$.
\index{fiber $\nbigi_y$}
For any $y\in Y$,
there exists a small neighbourhood $U_y$
and a decomposition
$\pi^{-1}(U_y)=\coprod_{a\in\nbigi_y}V_a$
such that 
the restriction of $\pi$ to $V_a$ $(a\in\nbigi_y)$
induces a homeomorphism $V_a\simeq U_y$.
Hence, we obtain the bijection
$\varphi_{y',y}:
 \nbigi_y\simeq \nbigi_{y'}$ for $y'\in U_y$
by setting
$\varphi_{y',y}(a):=V_a\cap\pi^{-1}(y')$.
\index{bijection $\varphi_{y',y}$}

Let $\vecleq=(\leq_y\,|\,y\in Y)$ be a $G$-equivariant family of
partial orders on $\nbigi_y$ $(y\in Y)$
satisfying the following condition.
\begin{itemize}
 \item 
 For any $y\in Y$,
 there exists a small neighbourhood $U_y$
 such that
 $a\leq_y b$ implies
 $\varphi_{y',y}(a)\leq_{y'}\varphi_{y',y}(b)$
 for any $y'\in U_y$.
\end{itemize}
\index{family of partially ordered sets $\vecleq$}

Let $U$ be any simply connected subset of $Y$.
Let $\nbigi(U)$ denote 
the set of the connected components
of $\pi^{-1}(U)$.
For any $y\in U$,
there exists the natural identification
$\nbigi(U)\simeq \nbigi_y$.
We define the partial order
$\leq_U$ on $\nbigi(U)$
by
\[
 a\leq_Ub
\Longleftrightarrow
 a\leq_yb\,\,(\forall y\in U).
\]
\index{partial order $\leq_U$}

\subsubsection{Graded local systems}

Let $L$ be a $G$-equivariant local system of $\cnum$-vector spaces
on $Y$.
We say that $L$ is graded over $\vecnbigi$
if it is equipped with an isomorphism
$L\simeq \pi_{\ast}N$ 
for a local system $N$ on $\vecnbigi$.
\index{graded local system}
The condition is equivalent to the following.
\begin{itemize}
\item
$L_y$ $(y\in Y)$ are equipped with
the grading
$L_y=\bigoplus_{a\in\nbigi_y} L_{y,a}$
such that
$L_{y,a}=L_{y',\varphi_{y',y}(a)}$
for any $a\in\nbigi_y$
if $y'$ is sufficiently close to $y$.
\end{itemize}

\begin{rem}
We do not exclude the case $L_{y,a}=0$ for some $a\in\nbigi_y$.
\hfill\qed
\end{rem}

\subsubsection{Stokes structures on local systems}
Let $L$ be a $G$-equivariant local system of $\cnum$-vector spaces
on $Y$.
A $G$-equivariant Stokes structure of $L$ indexed by $\vecnbigi$
is defined to be a $G$-equivariant family of filtrations
$\vecnbigf=(\nbigf^y\,|\,y\in Y)$
of the stalks $L_y$ $(y\in Y)$ indexed by
$(\nbigi_y,\leq_y)$
satisfying the following condition.
\index{Stokes structure}
\begin{itemize}
 \item
 For each $y\in Y$,
 there exists a neighbourhood $U_y$ of $y$
 and a decomposition
\[
 L_{|U_y}=
 \bigoplus_{a\in\nbigi_y}
 L_{U_y,a}
\]
such that 
$\nbigf^{y'}_a(L_{y'})=
      \bigoplus_{b\leq a}L_{U_y,b|y'}$
      for any $y'\in U_y$.      
\end{itemize}

\subsubsection{The associated graded local systems}
Let $(L,\vecnbigf)$ be a $G$-equivariant
local system with Stokes structure 
indexed by $\vecnbigi$.
For any $y\in Y$ and $a\in \nbigi_y$,
we define 
\[
\Gr^{\nbigf^y}_a(L_y):=
 \frac{\nbigf^y_a(L_y)}{\sum_{b\lneq a}\nbigf^y_b(L_y)},
\quad
 \Gr^{\nbigf^y}(L_y):=
 \bigoplus_{a\in\nbigi_y}
 \Gr^{\nbigf^y}_a(L_y).
\]
By the condition,
there exists the natural isomorphism
$\Gr^{\nbigf^y}_a(L_y)\simeq
 \Gr^{\nbigf^{y'}}_{\varphi_{y',y}(a)}(L_{y'})$,
and hence
$\Gr^{\nbigf^y}(L_y)\simeq
 \Gr^{\nbigf^{y'}}(L_{y'})$
if $y'$ is sufficiently close to $y$.
Thus, we obtain a $G$-equivariant local system
$\Gr^{\vecnbigf}(L)$
induced by $\Gr^{\nbigf^y}(L_y)$
$(y\in Y)$ with the isomorphisms.
It is graded over $\vecnbigi$.
\index{graded local system $\Gr^{\vecnbigf}(L)$}

For any $G$-equivariant section $\rho:Y\lrarr\vecnbigi$,
we obtain the $G$-equivariant local subsystem
$\Gr^{\vecnbigf}_{\rho}(L)\subset \Gr^{\vecnbigf}(L)$
induced by
$\Gr^{\nbigf^y}_{\rho(y)}(L_y)$ $(y\in Y)$.
More generally,
for a $G$-invariant submanifold
$\vecnbigi_1\subset\vecnbigi$
such that the restriction of $\pi$ to $\vecnbigi_1$
is also proper and locally a homeomorphism,
we obtain the $G$-equivariant local subsystem
$\Gr^{\vecnbigf}_{\vecnbigi_1}(L)\subset \Gr^{\vecnbigf}(L)$
induced by
$\bigoplus_{a\in\nbigi_{1,y}}
 \Gr^{\nbigf^y}_a(L_y)$.

\subsubsection{Loosening of Stokes filtrations}

Let $\varpi_i:\vecnbigi_i\lrarr Y$ $(i=1,2)$
be $G$-equivariant continuous proper maps of manifolds
which are locally homeomorphisms.
Let $\vecleq_i$ $(i=1,2)$ be a family of partial orders
on $\vecnbigi_i$.
Let $\psi:\vecnbigi_1\lrarr\vecnbigi_2$
be a $G$-equivariant continuous map 
such that $\varpi_2\circ\psi=\varpi_1$.
We assume the following.
\begin{itemize}
\item
 If $a\leq_{1,y} b$,
 then $\psi(a)\leq_{2,y}\psi(b)$ holds.
 Moreover,
 if $\psi(a)<_{2,y}\psi(b)$
 then $a<_{1,y}b$ holds.
\end{itemize}
In this case, the partial order $\leq_{1,y}$
is recovered from
$\leq_{2,y}$ and the restriction of $\leq_{1,y}$
to $\psi^{-1}(c)$ $(c\in \nbigi_{2,y})$.

Let $L$ be a local system on $Y$.
Let $\vecnbigf$ be a Stokes structure of $L$
indexed by $\vecnbigi_1$.
For each $y\in Y$,
there exists a splitting
$L_y=\bigoplus_{a\in\nbigi_{1,y}}L_{y,a}$
of the filtration $\nbigf^y$.
For each $c\in\nbigi_{2,y}$,
we define
\[
 (\psi_{\ast}\nbigf)^y_c:=
 \bigoplus_{\substack{a\in\nbigi_{1,y}\\
 \psi(a)\leq_{2,y}c}}
 L_{y,a}.
\]
It is independent of the choice of a splitting.
It is easy to see that the family of filtrations
$\psi_{\ast}\vecnbigf:=
 \bigl(
 (\psi_{\ast}\nbigf)^y\,\big|\,y\in Y
 \bigr)$ is a $G$-equivariant Stokes structure of $L$
indexed by $\vecnbigi_2$.
\index{Stokes structure $\psi_{\ast}\vecnbigf$}

Note that for each $y\in Y$,
the associated graded vector space
$\Gr^{\psi_{\ast}\nbigf^y}(L_y)$
is equipped with the filtration $\nbigf^y$
indexed by $\nbigi_{1,y}$,
which is compatible with the grading,
i.e.,
each $\Gr^{\psi_{\ast}\nbigf^y}_c(L)$ 
$(c\in \nbigi_{2,y})$
is equipped with the filtration $\nbigf^y$
indexed by $\psi^{-1}(c)$,
and $\nbigf^y$ on $\Gr^{\psi_{\ast}\nbigf^y}(L_y)$
is the direct sum of $\nbigf^y$ on 
$\Gr^{\psi_{\ast}\nbigf^y}_c(L_y)$.
The family $\vecnbigf=(\nbigf^y\,|\,y\in Y)$
is a $G$-equivariant Stokes structure
on the associated graded local system
$\Gr^{\psi_{\ast}\vecnbigf}(L)$
indexed by $\vecnbigi$,
which is compatible with the grading.
The pair
$(\Gr^{\psi_{\ast}\vecnbigf}(L),\vecnbigf)$
is denoted by
$\Gr^{\psi_{\ast}\vecnbigf}(L,\vecnbigf)$.
\index{local system with Stokes structure
$\Gr^{\psi_{\ast}\vecnbigf}(L,\vecnbigf)$}
The following lemma is clear.
\begin{lem}
A $G$-equivariant
Stokes structure on $L$ indexed by $\vecnbigi_1$
is equivalent to
a $G$-equivariant
Stokes structure $\vecnbigf$ on $L$
indexed by $\vecnbigi_2$
together with 
a $G$-equivariant
Stokes structures on $\Gr^{\vecnbigf}(L)$
indexed by $\vecnbigi$ compatible with the grading.
\hfill\qed
\end{lem}

\subsubsection{Stokes graded local systems}

Let $\varpi_i:\vecnbigi_i\lrarr Y$ $(i=1,2)$
be $G$-equivariant continuous proper maps of manifolds
which are locally homeomorphisms.
Let $\psi:\vecnbigi_1\lrarr\vecnbigi_2$
be a $G$-equivariant continuous map 
such that $\varpi_2\circ\psi=\varpi_1$.
Let $\vecleq$ be a family of partial orders
on $\vecnbigi_1$.

\begin{df}
A $G$-equivariant local system $L$ equipped 
with a Stokes structure $\vecnbigf$ indexed by $\vecnbigi_1$
and a grading over $\vecnbigi_2$
is called a Stokes graded local system over
$(\vecnbigi_1,\vecnbigi_2)$
if the following holds for any $y\in Y$.
\index{Stokes graded local system over $(\vecnbigi_1,\vecnbigi_2)$}
\begin{itemize}
\item
$L_{y,a}$ $(a\in\nbigi_{2,y})$ are equipped with
a Stokes structure $\nbigf^y$ indexed by $\psi^{-1}(a)$,
and 
$(L_{y},\nbigf^y)$ is equal to the direct sum
$\bigoplus_{a\in\nbigi_{2,y}}
 (L_{y,a},\nbigf^{y})$.
\hfill\qed
 \end{itemize}
\end{df}

\subsection{Families of partially ordered sets on $S^1$}

\subsubsection{Unramified case}

Let $\varpi:\cnumtilde\lrarr\cnum$ be 
the real oriented blow up along $0$,
i.e.,
$\cnumtilde=\real_{\geq 0}\times S^1$,
and $\varphi(r,e^{\sqrt{-1}\theta})=re^{\sqrt{-1}\theta}$.
We can identify $S^1$ with the boundary $\del\cnumtilde$
by $e^{\sqrt{-1}\theta}\longleftrightarrow (0,e^{\sqrt{-1}\theta})$.
\index{oriented real blow up}

Let $z$ be the standard coordinate of $\cnum$.
Let $\nbigi\subset z^{-1}\cnum[z^{-1}]$
be a finite subset.
For $\gminia,\gminib\in\nbigi$,
set $F_{\gminia,\gminib}:=
 -|z|^{-\ord(\gminia-\gminib)}\Re(\gminia-\gminib)$
as a function on $\cnumtilde$.
For each $e^{\sqrt{-1}\theta}\in \del\cnumtilde$,
we define $\gminia\leq_{e^{\sqrt{-1}\theta}}\gminib$
if $F_{\gminia,\gminib}\leq 0$ on a neighbourhood of 
$(0,e^{\sqrt{-1}\theta})$.
We obtain the family of partial orders
$\vecleq=(\leq_{P}\,|\,P\in \del\cnumtilde)$ 
on $\vecnbigi=\nbigi\times\del\cnumtilde$.

\subsubsection{Ramified case}

Let $p$ be a positive integer.
We take a $p$-th root $z_p$ of $z$.
We have the ramified covering
$\rho_p:\cnum_{z_p}\lrarr \cnum_z$
given by $\rho_p(z_p)=z_p^p$.
We have the identification
$S^1\simeq \del\cnumtilde_{z_p}$
given by
$e^{\sqrt{-1}\theta_p}
\longleftrightarrow
 (0,e^{\sqrt{-1}\theta_p})$.
The induced map
$\del\cnumtilde_{z_p}
\lrarr
 \del\cnumtilde_z$
is identified with
$e^{\sqrt{-1}\theta_p}
\longmapsto
 e^{\sqrt{-1}p\theta_p}$.
Set $\Gal(p):=\{a\in\cnum^{\ast}\,|\,a^p=1\}$.
\index{group $\Gal(p)$}
We have the action of $\Gal(p)$ on $\cnum(\!(z_p)\!)$
by $(a^{\ast}f)(z_p)=f(a z_p)$.

Let $\nbigi$ be a $\Gal(p)$-invariant subset
of $z_p^{-1}\cnum[z_p^{-1}]$.
We consider the $\Gal(p)$-action 
on $\nbigi\times \del\cnumtilde_{z_p}$
induced by 
$a\cdot (\gminia,z_p)=
 ((a^{\ast})^{-1}\gminia,a\cdot z_p)$.
Let $\vecnbigi$ denote the quotient manifold
by the action.
The projection 
 $\nbigi\times\del\cnumtilde_{z_p}
\lrarr 
 \del\cnumtilde_{z_p}$
induces a proper map
$\vecnbigi\lrarr \del\cnumtilde_p$
which is locally a homeomorphism.

For any $P_p\in \del\cnumtilde_{z_p}$ and $a\in \Gal(p)$,
$\gminia\leq_{a(P_p)}\gminib$ holds
if and only if
$a^{\ast}\gminia\leq_{P_p}a^{\ast}\gminib$
holds.
Hence, for each $P\in \del\cnumtilde_z$,
there exists the well defined order
$\leq_{P}$ on $\nbigi_{P}$.

\subsection{Families of partially ordered sets on $\real$}
\label{subsection;18.4.3.1}

For any positive integer $p$,
let $\varphi_p:\real\lrarr \del\cnumtilde_{z_p}$
be given by
$\varphi_p(t)=\exp(\sqrt{-1}t/p)$.
We have the induced identifications
$\real/2\pi p\seisuu\simeq
 \del\cnumtilde_{z_p}$,
and  in particular
$\real/2\pi\seisuu\simeq
 \del\cnumtilde_{z}$.

Let $\nbigi$ be a $\Gal(p)$-invariant subset of
$z_p^{-1}\cnum[z_p^{-1}]$.
For each $t\in\real$,
we have the partial order $\leq_t:=\leq_{\varphi_p(t)}$
on $\nbigi$.
\index{partial order $\leq_t$}
The map $\varphi_p$ induces
a homomorphism
$\varphi_p:2\pi\seisuu\lrarr \Gal(p)$.
We denote $\varphi_p(a)^{\ast}(\gminia)$
by $a^{\ast}\gminia$
for $a\in 2\pi\seisuu$
and $\gminia\in\nbigi$.

By the construction,
for any $a\in 2\pi\seisuu$
and $\gminia,\gminib\in\nbigi$,
we have
$a^{\ast}\gminia\leq_{t}a^{\ast}\gminib$
if and only if
we have $\gminia\leq_{t+a}\gminib$.
Hence, the family of partial orders
$\vecleq=(\leq_{t}\,|\,t\in\real)$
is $2\pi\seisuu$-equivariant.

\begin{lem}
The following objects are naturally equivalent.
\begin{itemize}
\item
$2\pi\seisuu$-equivariant
local systems with Stokes structure on $\real$
indexed by $\nbigi$.
\item
Local systems with Stokes structure on 
$\del\cnumtilde_z$
indexed by $\vecnbigi$.
\hfill\qed
\end{itemize}
\end{lem}

Let $L$ be a $2\pi\seisuu$-equivariant local system on $\real$.
The following lemma is well known and clear
by definition of Stokes structure.
\begin{lem}
\label{lem;25.2.24.1}
Let $\vecnbigf$ be a $2\pi\seisuu$-equivariant local system of $L$
indexed by $\nbigi$.
For any $\theta\in\real$,
there exists $\epsilon>0$ such that the following holds.
 \begin{itemize}
 \item For any $\theta_-\in\openopen{\theta-\epsilon}{\theta}$
       and $\theta_+\in\openopen{\theta}{\theta+\epsilon}$,
       we have
       $\nbigf^{\theta}_{\gminia}
       =\nbigf^{\theta_-}_{\gminia}\cap\nbigf^{\theta_+}_{\gminia}$
       under the natural isomorphisms
       $L_{\theta}\simeq L_{\theta_{\pm}}$.
       \hfill\qed
 \end{itemize} 
\end{lem}

We obtain the following lemma as a consequence.
\begin{lem}
\label{lem;25.2.23.4}
Let $(L_i,\vecnbigf)$ be
$2\pi\seisuu$-equivariant local system
with a Stokes structure indexed by $\nbigi$.
Let $L_1\to L_2$ be a morphism of $2\pi\seisuu$-equivariant
local systems.
Then,
$\varphi$ gives a morphism
$(L_1,\vecnbigf)\to(L_2,\vecnbigf)$
of 
$2\pi\seisuu$-equivariant local system with a Stokes structure
if and only if 
there exists a discrete subset $Z\subset\real$
such that 
$\varphi(\nbigf^{\theta}_{\gminia}(L_{1,\theta}))
\subset \nbigf^{\theta}_{\gminia}(L_{2,\theta})$
for any $\theta\in\real\setminus Z$ and any $\gminia\in\nbigi$.
\hfill\qed
\end{lem}

\section{Some notation for index sets and intervals}
\label{subsection;18.5.2.10}

Let $\theta$ be the standard coordinate on $\real$.
For any open interval
\index{set $J_-$}
\index{set $J_+$}
\index{set $\Jbar$}
\[
 J=\openopen{\theta_0}{\theta_1}
:=\{\theta_0<\theta<\theta_1\}\subset\real,
\]
we put
\[
J_+:=\openclosed{\theta_0}{\theta_1}
=\{\theta_0\leq \theta<\theta_1\},
\quad
J_-:=\closedopen{\theta_0}{\theta_1}
=\{\theta_0<\theta\leq\theta_1\},
\]
\[
\Jbar:=\closedclosed{\theta_0}{\theta_1}
=\{\theta_0\leq\theta\leq\theta_1\}.
\]
The boundary points
$\theta_0$ and $\theta_1$
are also denoted by 
$\vartheta^J_{\ell}$
and $\vartheta^J_r$,
respectively.
The middle point
$(\theta_0+\theta_1)/2$
is denoted by $\vartheta^J_m$.
\index{points $\vartheta^J_r$, $\vartheta^J_{\ell}$, $\vartheta^J_m$}
For any real number $s$,
we set 
$J+s:=\{\theta+s\,|\,\theta\in J\}$.
\index{set $J+s$}

Let $p$ and $n$ be positive integers.
Set $\omega:=n/p$.
Let $\nbigi\subset \cnum\cdot z_p^{-n}\subset
 \cnum(\!(z_p)\!)$ be a $\Gal(p)$-invariant finite subset.
We set $\nbigi^{\ast}:=\nbigi\setminus\{0\}$.
\index{set $\nbigi^{\ast}$}
We have the partial orders
$\leq_{\theta}:=\leq_{\varphi_p(\theta)}$ $(\theta\in\real)$
on $\nbigi$ as in \S\ref{subsection;18.4.3.1}.
\index{partial order $\leq_{\theta}$}
For any open interval $J\subset\real$,
we define the partial order
$\leq_J$ on $\nbigi$
as in \S\ref{subsection;18.4.3.2}.
\index{partial order $\leq_J$}
For any $\gminia\in\nbigi^{\ast}$,
we set \index{set $S_0(\gminia)$}
\[
 S_0(\gminia):=
 \bigl\{
 \theta\in\real\,\big|\,
 \Re(\gminia(e^{\sqrt{-1}\theta/p}))=0
 \bigr\}.
\]
Let $T(\gminia)$
denote the set of connected components of
$\real\setminus S_0(\gminia)$.
Let 
$T_+(\gminia):=
 \bigl\{
 J\in T(\gminia)\,\big|\,
 \gminia>_J0
 \bigr\}$
and 
$T_-(\gminia):=
 \bigl\{
 J\in T(\gminia)\,\big|\,
 \gminia<_J0
 \bigr\}$.
\index{sets $T(\gminia)$, $T_+(\gminia)$, $T_-(\gminia)$}

We put
$S_0(\nbigi):=
 \bigcup_{\gminia\in\nbigi^{\ast}}
 S_0(\gminia)$,
and
$T(\nbigi):=
 \bigcup_{\gminia\in\nbigi^{\ast}}
 T(\gminia)$.
Let $T_2(\nbigi)$
denote the set of pairs 
$(J_1,J_2)$ in $T(\nbigi)$
satisfying
$J_1\cap J_2\neq\emptyset$
and $J_1\neq J_2$.
\index{set $S_0(\nbigi)$}
\index{set $T(\nbigi)$}
\index{set $T_2(\nbigi)$}

When $\nbigi^{\ast}\neq\emptyset$,
for any connected component $I=\openopen{\theta_0}{\theta_1}$ 
of $\real\setminus S_0(\nbigi)$,
let $T(\nbigi)_{I}$ denote the set of
$J\in T(\nbigi)$ such that
$I\subset J$,
\index{set $T(\nbigi)_I$}
i.e.,
$T(\nbigi)_I=
\bigl\{\openopen{\theta_1-\pi/\omega}{\theta_1}\bigr\}
\cup\bigl\{J\in T(\nbigi)\,|\,\theta_1\in J\bigr\}$.

For any $J\in T(\nbigi)$,
we set \index{sets $\nbigi_J$, $\nbigi_{J,<0}$, $\nbigi_{J,>0}$}
\[
 \nbigi_J:=\bigl\{
 \gminia\in\nbigi^{\ast}\,\big|\,
 J\in T(\gminia)
 \bigr\}
\cup\{0\},
\quad
 \nbigi_{J,<0}:=
 \bigl\{
 \gminia\in\nbigi^{\ast}\,\big|\,
 J\in T_-(\gminia)
 \bigr\},
\]
\[
 \nbigi_{J,>0}:=
 \bigl\{
 \gminia\in\nbigi^{\ast}\,\big|\,
 J\in T_+(\gminia)
 \bigr\}.
\]
We also put
$\nbigi^{\ast}_J:=\nbigi_J\setminus\{0\}
=\nbigi_{J,<0}\cup
 \nbigi_{J,>0}$. \index{set $\nbigi^{\ast}_J$}
We put 
$\nbigi_{J,\leq 0}:=
 \nbigi_{J,<0}\cup\{0\}$
and 
$\nbigi_{J,\geq 0}:=
 \nbigi_{J,>0}\cup\{0\}$.
\index{sets $\nbigi_{J,\leq 0}$, $\nbigi_{J,\geq 0}$}
Note that there exists the decomposition
$\nbigi^{\ast}
=\coprod_{J\in T(\nbigi)_I}
\nbigi_J^{\ast}$
for any connected component $I$ of
$\real\setminus S_0(\nbigi)$.

Let $J_0,J_1\in T(\nbigi)$.
We write $J_0<J_1$
if $\vartheta^{J_0}_{\ell}<\vartheta^{J_1}_{\ell}$.
\index{relation $J_0<J_1$}
We write $J_0\vdash J_1$
if the following conditions are satisfied;
(i) $\vartheta^{J_0}_{\ell}<\vartheta^{J_1}_{\ell}$,
(ii) $\openopen{\vartheta^{J_0}_{\ell}}{\vartheta^{J_1}_{\ell}}
 \cap S_0(\nbigi)=\emptyset$.
\index{relation $J_0\vdash J_1$}

Let $\Tbb:\real\lrarr\real$ be given by
$\Tbb(\theta)=\theta+2\pi$.
\index{map $\Tbb$}
Let $\Tbb^{\ast}:\nbigi\lrarr\nbigi$
be given by 
$\Tbb^{\ast}(\gminia)(z_p)=
 \gminia(e^{2\pi\sqrt{-1}/p}z_p)$.
\index{map $\Tbb^{\ast}$}
If $J\in T(\gminia)$,
then we have $\Tbb^{-1}(J)\in T(\Tbb^{\ast}\gminia)$.
In particular, 
$T(\nbigi)$ is invariant under the translation 
by $2\pi\seisuu$.

\begin{lem}
Let $\gminio$ be a $\Gal(p)$-orbit in $\nbigi$.
Then, for any $J\in T(\nbigi)$,
we have
$|\nbigi_{J,>0}\cap\gminio|\leq 1$ 
and 
$|\nbigi_{J,<0}\cap\gminio|\leq 1$.
\end{lem}
\pf
Suppose that
$\gminio\cap\nbigi_{J,>0}\neq\emptyset$.
We take $\gminia=\alpha z_p^{-n}\in\gminio\cap\nbigi_{J,>0}$.
We have 
$\alpha=-|\alpha|\exp\bigl(\sqrt{-1}(\vartheta^J_m/\omega)\bigr)$.
Then, the claim is clear.
\hfill\qed

\section{Local systems with Stokes structure on $\real$}

We prepare some notation and procedures for 
$2\pi\seisuu$-equivariant local systems
with Stokes structure on $\real$.

\subsection{Category}
Let $\nbigi$ be a $\Gal(p)$-invariant subset 
of $z_p^{-1}\cnum[z_p^{-1}]$.
Let $\Loc^{\St}(\nbigi)$ denote the category of
$2\pi\seisuu$-equivariant local systems with Stokes structure
indexed by $\nbigi$.
\index{category $\Loc^{\St}(\nbigi)$}
A morphism $f:(L_1,\vecnbigf)\lrarr(L_2,\vecnbigf)$
is defined to be
a morphism of $2\pi\seisuu$-equivariant local systems
$f:L_1\lrarr L_2$
such  that
$f(\nbigf_{\gminia}(L_{1|\theta}))
\subset
 \nbigf_{\gminia}(L_{2|\theta})$
for any $\theta\in\real$ and $\gminia\in\nbigi$.

\subsection{Loosening}
\label{subsection;18.4.18.1}

For $\omega=\ell/p\in \frac{1}{p}\seisuu_{>0}$,
let $\pi_{\omega}:
 z_p^{-1}\cnum[z_p^{-1}]
\lrarr
 z_p^{-\ell}\cnum[z_p^{-1}]$
denote the projection
given by
$\pi_{\omega}\bigl(
 \sum a_jz_p^{-j}
 \bigr):=
 \sum_{j\geq \ell}
 a_jz_p^{-j}$.
\index{projection $\pi_{\omega}$}
We set
$\nbigt_{\omega}(\nbigi):=
 (\pi_{\omega}^{-1}(0)\cap\nbigi)\cup\{0\}$ and
 $\nbigs_{\omega}(\nbigi):=
 (\nbigi\setminus \nbigt_{\omega}(\nbigi))
 \sqcup\{0\}$.
We also set
$\nbigttilde_{\omega}(\nbigi):=
 \nbigt_{\omega+1/p}(\nbigi)$
and 
$\nbigstilde_{\omega}(\nbigi):=
 \nbigs_{\omega+1/p}(\nbigi)$.
\index{set $\nbigt_{\omega}(\nbigi)$}
\index{set $\nbigs_{\omega}(\nbigi)$}
\index{set $\nbigttilde_{\omega}(\nbigi)$}
\index{set $\nbigstilde_{\omega}(\nbigi)$}

Let $(L,\vecnbigf)\in\Loc^{\St}(\nbigi)$.
Let $\psi_1:\nbigi\lrarr\nbigs_{\omega}(\nbigi)$
be the map defined by
$\psi_1(\gminia)=\gminia$ for $\gminia\not\in\nbigt_{\omega}(\nbigi)$
and $\psi_1(\gminia)=0$ for $\gminia\in\nbigt_{\omega}(\nbigi)$.
We denote 
$\psi_{1\ast}\vecnbigf$
by $\nbigs_{\omega}(\vecnbigf)$,
and we set
$ \nbigs_{\omega}(L,\vecnbigf):=
(L,\nbigs_{\omega}(\vecnbigf))$.
We also set
$\nbigstilde_{\omega}(\vecnbigf):=
 \nbigs_{\omega+1/p}(\vecnbigf)$
and 
$\nbigstilde_{\omega}(L,\vecnbigf):=
(L,\nbigstilde_{\omega}(\vecnbigf))$.
\index{local system with Stokes structure
$\nbigs_{\omega}(L,\vecnbigf)$}
\index{local system with Stokes structure
$\nbigstilde_{\omega}(L,\vecnbigf)$}
Thus, we obtain the functors
$\nbigs_{\omega}:\Loc^{\St}(\nbigi)
\lrarr\Loc^{\St}(\nbigs_{\omega}(\nbigi))$
and
$\nbigstilde_{\omega}:\Loc^{\St}(\nbigi)
\lrarr\Loc^{\St}(\nbigstilde_{\omega}(\nbigi))$.

Let
$\psi_2:\nbigi\lrarr \pi_{\omega}(\nbigi)$
be the map induced by $\pi_{\omega}$.
We set
$\vecnbigf^{(\omega)}:=
 \psi_{2\ast}\vecnbigf$.
\index{Stokes structure $\vecnbigf^{(\omega)}$}
Let 
$\nbigt_{\omega}(L,\vecnbigf)$
denote the $2\pi\seisuu$-equivariant local system
with the induced Stokes structure
$\bigl(
 \Gr^{\vecnbigf^{(\omega)}}_0(L),
 \vecnbigf
\bigr)$
indexed by $\nbigt_{\omega}(\nbigi)$.
\index{local system with Stokes structure
$\nbigt_{\omega}(L,\vecnbigf)$}
We also set
$\nbigttilde_{\omega}(L,\vecnbigf):=
 \nbigt_{\omega+1/p}(L,\vecnbigf)$.
\index{local system with Stokes structure
$\nbigttilde_{\omega}(L,\vecnbigf)$}
Thus, we obtain the functors
$\nbigt_{\omega}:\Loc^{\St}(\nbigi)
\lrarr \Loc^{\St}(\nbigt_{\omega}(\nbigi))$
and
$\nbigttilde_{\omega}:\Loc^{\St}(\nbigi)
\lrarr \Loc^{\St}(\nbigttilde_{\omega}(\nbigi))$.
By the construction,
$\nbigs_{\omega}\circ\nbigt_{\omega}(L,\vecnbigf)
\simeq
\nbigt_{\omega}\circ\nbigs_{\omega}(L,\vecnbigf)$
is just a $2\pi\seisuu$-equivariant local system
$\Gr^{\vecnbigf^{(\omega)}}_0(L)$
with the trivial Stokes structure
indexed by $\{0\}$.

For each $\theta\in \real$,
we have the subspaces
$L^{<0}_{\theta}:=\nbigf^{\theta}_{<0}$
and $L^{\leq 0}_{\theta}:=\nbigf^{\theta}_0$
of $L_{\theta}$.
They determine $2\pi\seisuu$-equivariant
constructible subsheaves
$L^{<0}$ and $L^{\leq 0}$ of $L$.
\index{constructible subsheaves $L^{<0}$, $L^{\leq 0}$}
More generally,
for any $\omega\in\rnum_{>0}$,
the subspaces
$L^{(\omega)\,<0}_{\theta}:=\nbigf^{(\omega)\,\theta}_{<0}$
and 
$L^{(\omega)\,\leq 0}_{\theta}:=
 \nbigf^{(\omega)\,\theta}_{0}$
induce $2\pi\seisuu$-equivariant constructible subsheaves
$L^{(\omega)\,<0}$
and 
$L^{(\omega)\,\leq 0}$.
\index{constructible subsheaves
 $L^{(\omega)\,<0}$, $L^{(\omega)\,\leq 0}$}
We naturally have
\[
L^{(\omega)\,<0}\subset L^{<0}
\subset L^{\leq 0}
\subset L^{(\omega)\,\leq 0}.
\]
Note that 
$\nbigt_{\omega}(L,\vecnbigf)$
is naturally isomorphic to
$L^{(\omega)\leq 0}\big/
 L^{(\omega)<0}$
with the induced Stokes structure.

Let $L_{S^1}$ denote the sheaf
on $S^1=\real/2\pi\seisuu$
obtained as the descent of $L$
with respect to the $2\pi\seisuu$-action.
\index{local system $L_{S^1}$}
For any $2\pi\seisuu$-equivariant constructible subsheaf
$K\subset L$,
let $K_{S^1}\subset L_{S^1}$ denote the subsheaf
obtained as the descent.
\index{constructible subsheaf $K_{S^1}$}
In particular,
we obtain constructible subsheaves
$L^{<0}_{S^1}$, $L^{(\omega)\leq 0}_{S^1}$,
etc.
\index{constructible subsheaves $L^{<0}_{S^1}$,
$L^{(\omega)\,\leq 0}_{S^1}$}

\subsection{Canonical splittings
and induced local systems}
\label{subsection;18.5.13.10}

Let $\nbigi$ be a $\Gal(p)$-invariant subset of 
$z_p^{-n}\cnum\subset z_p^{-1}\cnum[z_p^{-1}]$.
Let $(L,\vecnbigf)$ be a local system with Stokes structure
indexed by $\nbigi$.
Set $\omega:=n/p$.
Take any interval $J$
such that $\vartheta^J_{r}-\vartheta^J_{\ell}=\pi/\omega$.
There exist the unique decompositions
\begin{equation}
\label{eq;18.5.2.1}
 L_{|J_{\pm}}
=\bigoplus_{\gminia\in\nbigi}
 L_{J_{\pm},\gminia}
\end{equation}
such that 
$\nbigf^{\theta}_{\gminia}(L_{\theta})
=\bigoplus_{\gminib\leq_{\theta}\gminia}L_{J_{\pm},\gminib|\theta}$
for any $\theta\in J_{\pm}$.
(The uniqueness is clear.
The existence is also standard and well known.
For example,
see \cite[Proposition 3.16]{Mochizuki-good-Stokes},
where a higher dimensional analogue is proved.)
Such decompositions are called canonical splittings
in this monograph.
\index{canonical splitting}
\index{local subsystems $L_{J_{\pm},\gminia}$}

\subsubsection{Some decompositions}

Take any $\theta_0\in\real$.

\begin{lem}
\label{lem;18.5.2.1}
Suppose $\theta_0\in\real\setminus S_0(\nbigi)$.
We choose an interval
$J(0)\in T(\nbigi)$ such that $\theta_0\in J(0)$.
We also choose a function
$\mu:\nbigi\lrarr \{\pm\}$.
Then, we obtain the following decomposition:
\[
 L_{|\theta_0}
=L_{J(0)_{\mu(0)},0|\theta_0}
\oplus
\bigoplus_{\substack{J\in T(\nbigi)\\
 \theta_0\in J}}
\bigoplus_{\gminia\in\nbigi_J^{\ast}}
 L_{J_{\mu(\gminia)},\gminia|\theta_0}.
\]
\end{lem}
\pf
We have 
$L_{J(0)_{\mu(0)},0|\theta_0}
\subset
 \nbigf^{\theta_0}_0$,
and the induced map
$L_{J(0)_{\mu(0)},0|\theta_0}
\lrarr
 \Gr^{\nbigf^{\theta_0}}_0(L_{|\theta_0})$
is an isomorphism.
For each $\gminia\in\nbigi_J^{\ast}$,
we have 
$L_{J_{\mu(\gminia)},\gminia|\theta_0}
\subset
 \nbigf^{\theta_0}_{\gminia}$
and the induced map
$L_{J_{\mu(\gminia)},\gminia|\theta_0}
\lrarr
 \Gr^{\nbigf^{\theta_0}}_{\gminia}(L_{|\theta_0})$
is an isomorphism.
Then, the claim is clear.
\hfill\qed

\vspace{.1in}

Let us consider the case
$\theta_0\in S_0(\nbigi)$.
Set $J_0=\openopen{\theta_0}{\theta_0+\pi/\omega}$.
For $\gminia\in\nbigi_{J_0}^{\ast}$,
we choose $J(\gminia)\in\{J_0,J_0-\pi/\omega\}$.
We choose $J(0)\in T(\nbigi)$
such that $\theta_0\in\Jbar(0)$,
i.e.,
$J_0\leq J(0)\leq J_0+\pi/\omega$.
Note that
for any $\gminia\not\in\nbigi_{J_0}$,
there exists a unique $J\in T(\gminia)$
such that $\theta_0\in J$.
Let $\mu:\nbigi\lrarr \{\pm\}$ be
a map satisfying the following.
\begin{itemize}
\item
We assume
 $\mu(0)=+$ if $J(0)=J_0-\pi/\omega$,
 and $\mu(0)=-$ if $J(0)=J_0$.
\item
 For $\gminia\in\nbigi_{J_0}^{\ast}$,
 we have
 $\mu(\gminia)=+$  if $J(\gminia)=J_0-\pi/\omega$
and 
 $\mu(\gminia)=-$  if $J(\gminia)=J_0$.
\end{itemize}
We obtain the following lemma
by the argument in the proof of
Lemma \ref{lem;18.5.2.1}.
\begin{lem}
\label{lem;18.5.2.2}
We obtain the following decomposition:
\[
 L_{|\theta_0}
=L_{J(0)_{\mu(0)},0|\theta_0}
\oplus
 \bigoplus_{\gminia\in\nbigi_{J_0}^{\ast}}
 L_{J(\gminia)_{\mu(\gminia)},\gminia|\theta_0}
\oplus
 \bigoplus_{\substack{J\in T(\nbigi)\\ \theta_0\in J}}
 \bigoplus_{\gminib\in\nbigi_J^{\ast}}
 L_{J_{\mu(\gminib)},\gminib|\theta_0}.
\]
\hfill\qed
\end{lem}

\subsubsection{Induced local systems}

For any $J\in T(\nbigi)$,
we set
$L_{J_{\pm},<0}:=
 \bigoplus_{\gminia<_J0}
 L_{J_{\pm},\gminia}$
and 
$L_{J_{\pm},\leq 0}:=
 \bigoplus_{\gminia\leq_J0}
 L_{J_{\pm},\gminia}$
on $J_{\pm}$.
\index{local subsystems $L_{J_{\pm},<0}$}
\index{local subsystems $L_{J_{\pm},\leq 0}$}
We also set
$\gbiga_{J_{\pm}}(L):=
\bigoplus_{\gminia\in\nbigi_J}L_{J_{\pm},\gminia}$
on $J_{\pm}$.
\index{local subsystems $\gbiga_{J_{\pm}}(L)$}
The following lemma is easy to check.
\begin{lem}
We have
$L_{J_{+},<0|J}=L_{J_-,<0|J}$,
$L_{J_{+},\leq 0|J}=L_{J_-,\leq 0|J}$
 and
 $\gbiga_{J_{+}}(L)_{|J}=
 \gbiga_{J_{-}}(L)_{|J}$.
\hfill\qed
\end{lem}

Because 
$L_{J_{+},<0|J}=L_{J_-,<0|J}$,
we obtain
a local subsystem
$L_{\Jbar,<0}\subset L_{|\Jbar}$
by gluing $L_{J_{\pm},<0}$.
\index{local subsystem $L_{\Jbar,<0}$}
Because 
$L_{J_{+},\leq 0|J}=L_{J_-,\leq 0|J}$,
we obtain 
a local subsystem
$L_{\Jbar,\leq 0}\subset L_{|\Jbar}$
by gluing $L_{J_{\pm},\leq 0}$.
\index{local subsystem $L_{\Jbar,\leq 0}$}
The restrictions
$L_{\Jbar,< 0|J}$
and 
$L_{\Jbar,\leq 0|J}$
are also denoted by
$L_{J,<0}$
and $L_{J,\leq 0}$.
\index{local subsystem $L_{J,<0}$}
\index{local subsystem $L_{J,\leq 0}$}
Clearly,
$L_{\Jbar,<0}\subset (L^{<0})_{|\Jbar}$
and 
$L_{\Jbar,\leq 0}\subset (L^{\leq 0})_{|\Jbar}$.

We define
$L_{J,0}:=L_{J,\leq 0}/L_{J,<0}$ on $J$,
and 
$L_{\Jbar,0}:=L_{\Jbar,\leq 0}/L_{\Jbar,<0}$ on $\Jbar$.
We have
$L_{\Jbar,0|J}=L_{J,0}$.
\index{local system $L_{J,0}$}
\index{local system $L_{\Jbar,0}$}

Because
$\gbiga_{J_{-}}(L)_{|J}
=\gbiga_{J_{+}}(L)_{|J}$,
we obtain the local subsystem
$\gbiga_{\Jbar}(L)\subset L_{\Jbar}$.
The restriction
$\gbiga_{\Jbar}(L)_{|J}$ is denoted by
$\gbiga_J(L)$.
\index{local subsystem $\gbiga_{\Jbar}(L)$}
\index{local subsystem $\gbiga_{J}(L)$}

We set 
$L_{J_{\pm},>0}:=
 \bigoplus_{\gminia>_J0}L_{J_{\pm},\gminia}$
and 
$L_{J_{\pm},\geq 0}:=
 \bigoplus_{\gminia\geq_J0}L_{J_{\pm},\gminia}$
on $J_{\pm}$.
\index{local subsystems $L_{J_{\pm},>0}$}
\index{local subsystems $L_{J_{\pm},\geq 0}$}
Note that
$L_{J_{\pm},>0}$
and
$L_{J_{\pm},\geq 0}$
are naturally isomorphic to
$\gbiga_{J_{\pm}}(L)\big/L_{J_{\pm},\leq 0}$
and 
$\gbiga_{J_{\pm}}(L)\big/L_{J_{\pm},< 0}$,
respectively.

We define
$L_{\Jbar,>0}:=\gbiga_{\Jbar}(L)
 \big/L_{\Jbar,\leq 0}$
and 
$L_{\Jbar,\geq 0}:=\gbiga_{\Jbar}(L)\big/L_{\Jbar,<0}$
as the quotient sheaves on $\Jbar$.
\index{local system $L_{\Jbar},>0$}
\index{local system $L_{\Jbar},\geq 0$}
We also define
$L_{J,>0}:=\gbiga_{J}(L)
 \big/L_{J,\leq 0}
\simeq L_{\Jbar,>0|J}$
and 
$L_{J,\geq 0}:=\gbiga_{J}(L)\big/L_{J,<0}
\simeq L_{\Jbar,\geq 0|J}$
on $J$.
\index{local system $L_{J,>0}$}
\index{local system $L_{J,\geq 0}$}

\begin{rem}
\label{rem;20.9.7.1}
There exist the local subsystems
$L'_{J,<0}$ and $L'_{J,\leq 0}$
of $L$ on $\real$
determined by the conditions
$L'_{J,<0|J}=L_{J,<0}$
and 
$L'_{J,\leq 0|J}=L_{J,\leq 0}$,
respectively.
\index{local subsystem $L'_{J,<0}$}
\index{local subsystem $L'_{J,\leq 0}$}
Note that 
$L'_{J,<0|\theta}$ is not necessarily
contained in $\nbigf^{\theta}_{<0}$
if $\theta$ is not contained in $\Jbar$.
\hfill\qed
\end{rem}

\subsubsection{Relations}
We have the following inclusions
of subspaces in $L_{|\vartheta^J_r}$:
 \begin{equation}
\label{eq;20.9.7.10}
 L_{J_{+},>0|\vartheta^J_r}
 \subset
 L_{(J+\pi/\omega)_-,<0|\vartheta^J_r}
 \oplus
 \bigoplus_{
 \substack{J'\in T(\nbigi)\\ \vartheta^J_r\in J
 }}
 L_{J',<0|\vartheta^J_r}.
 \end{equation}
\begin{equation}
\label{eq;20.9.7.11}
 L_{J_{+},<0|\vartheta^J_r}
 \subset
 L_{(J+\pi/\omega)_-,>0|\vartheta^J_r}
 \oplus
 \bigoplus_{
 \substack{J'\in T(\nbigi)\\ \vartheta^J_r\in J
 }}
 L_{J',<0|\vartheta^J_r}.
\end{equation}
For $J\vdash J_1$,
we have the following equality
of subspaces in $L_{|\vartheta^J_r}$:
\begin{equation}
\label{eq;20.9.7.12}
 L_{J_+,0|\vartheta^J_r}
=L_{(J_1)_-,0|\vartheta^J_r}.
\end{equation}
Similarly,
we have the following inclusions
of subspaces in $L_{|\vartheta^J_{\ell}}$:
\begin{equation}
\label{eq;20.9.7.13}
 L_{J_{-},>0|\vartheta^J_{\ell}}
 \subset
 L_{(J-\pi/\omega)_+,<0|\vartheta^J_{\ell}}
 \oplus
 \bigoplus_{
 \substack{J'\in T(\nbigi)\\ \vartheta^J_{\ell}\in J
 }}
 L_{J',<0|\vartheta^J_{\ell}}.
 \end{equation}
\begin{equation}
\label{eq;20.9.7.14}
 L_{J_{-},<0|\vartheta^J_{\ell}}
 \subset
 L_{(J-\pi/\omega)_+,>0|\vartheta^J_{\ell}}
 \oplus
 \bigoplus_{
 \substack{J'\in T(\nbigi)\\ \vartheta^J_{\ell}\in J
 }}
 L_{J',<0|\vartheta^J_{\ell}}.
\end{equation}
For $J_2\vdash J$,
we have the following equality
of subspaces
in $L_{|\vartheta^J_{\ell}}$:
\begin{equation}
\label{eq;20.9.7.15}
 L_{J_-,0|\vartheta^J_{\ell}}
=L_{(J_2)_+,0|\vartheta^J_{\ell}}.
\end{equation}

\subsubsection{Some other decompositions}
\label{subsection;21.6.7.1}

Let $\gbigk(J_{0-})$
denote the set of $J\in T(\nbigi)$
such that 
$J_+\cap J_{0-}\neq\emptyset$.
Note that an interval $J\in \gbigk(J_{0-})$
does not necessarily contain $\theta_0$.
\index{set $\gbigk(J_{0-})$}

\begin{lem}
\label{lem;18.4.19.110}
There exists the following decomposition:
\[
 L_{|\theta_0}
=L_{J_{0-},0|\theta_0}
\oplus
 \bigoplus_{J\in\gbigk(J_{0-})}
 L'_{J,<0|\theta_0}.
\]
 (See Remark {\rm\ref{rem;20.9.7.1}}
 for $L'_{J,<0}$.)
\end{lem}
\pf
By Lemma \ref{lem;18.5.2.2},
we have the following decomposition:
\[
 L_{|\theta_0}
=L_{J_{0-},\leq 0|\theta_0}
\oplus
 L_{(J_0-\pi/\omega)_+,<0|\theta_0}
\oplus
\bigoplus_{\theta_0\in J'}
 \bigl(
 L_{(J')_+,>0}
\oplus
 L_{J',<0}
 \bigr)_{|\theta_0}.
\]
For $0<a<\pi/\omega$,
by using (\ref{eq;20.9.7.11}),
we obtain the following:
\begin{multline}
 L_{J_{0,-},\leq 0|\theta_0}
\oplus
 L_{(J_0-\pi/\omega)_+,<0|\theta_0}
\oplus
\bigoplus_{\theta_0\in J'}
 L_{J',<0|\theta_0} \\
\oplus
 \bigoplus_{J_0<J'< J_0+a}
 L'_{J',<0|\theta_0}
\oplus
 \bigoplus_{J_0+a\leq J'<J_0+\pi/\omega}
 L_{(J'-\pi/\omega)_+,>0|\theta_0}
=
\\
 L_{J_{0-},\leq 0|\theta_0}
\oplus
 L_{(J_0-\pi/\omega)_+,<0|\theta_0}
\oplus
\bigoplus_{\theta_0\in J'}
 L_{J',<0|\theta_0} \\
\oplus
 \bigoplus_{J_0<J'\leq J_0+a}
 L'_{J',<0|\theta_0}
\oplus
 \bigoplus_{J_0+a<J'<J_0+\pi/\omega}
 L_{(J'-\pi/\omega)_+,>0|\theta_0}.
\end{multline}
Then, we obtain the claim of the lemma.
\hfill\qed

\vspace{.1in}

Let $\gbigk(J_{0+})$
denote the set of $J\in T(\nbigi)$
such that 
$J_-\cap J_{0+}\neq\emptyset$.
\index{set $\gbigk(J_{0+})$}
As in the case of Lemma \ref{lem;18.4.19.110},
we obtain the decomposition:
\begin{equation}
\label{eq;18.5.2.3}
 L_{|\theta_0+\pi/\omega}
=L_{J_{0+}|\theta_0+\pi/\omega}
\oplus
\bigoplus_{J\in\gbigk(J_{0+})}
L'_{J,<0|\theta_0+\pi/\omega}.
\end{equation}

\subsection{Induced maps}
\label{subsection;24.2.19.1}

We continue to use the notation in \S\ref{subsection;18.5.13.10}.

\subsubsection{}
There exist the following decompositions:
\[
 H^0(J_{\pm},\gbiga_{\Jbar}(L))
=H^0(J_{\pm},L_{J_{\pm},<0})
\oplus
 H^0(J_{\pm},L_{J_{\pm},0})
 \oplus
 H^0(J_{\pm},L_{J_{\pm},>0}).
\]
We obtain the maps
$(\nbigq_{J_{\mp}},\nbigr^{J_{\mp}}_{J_{\pm}}):
H^0(J,L_{J,>0})\to
H^0(J,L_{J,0})\oplus
H^0(J,L_{J,<0})$
as the composition of
the natural isomorphisms,
the inclusion, and the projection:
\begin{multline}
 H^0(J,L_{J,>0})\simeq
 H^0(J_{\mp},L_{J_{\mp},>0})
 \lrarr
 H^0(J_{\mp},\gbiga_J(L))
 \simeq
 H^0(J_{\pm},\gbiga_J(L))
 \\
 \lrarr
 H^0(J_{\pm},L_{J_{\pm},0})
 \oplus
 H^0(J_{\pm},L_{J_{\pm},<0})
 \simeq
 H^0(J,L_{J,0})
 \oplus
 H^0(J,L_{J,<0}).
\end{multline}
We also obtain the maps
$\nbigp_{J_{\mp}}:
H^0(J,L_{J,0})\to
H^0(J,L_{J,<0})$
as the composition of the natural isomorphisms,
the inclusion and the projection:
\begin{multline}
 H^0(J,L_{J,0})\simeq
 H^0(J_{\mp},L_{J_{\mp},0})
 \lrarr
 H^0(J_{\mp},\gbiga_J(L))
 \simeq
 H^0(J_{\pm},\gbiga_J(L))
 \\
 \lrarr
 H^0(J_{\pm},L_{J_{\pm},<0})
 \simeq
 H^0(J,L_{J,<0}).
\end{multline}
The following lemma is easy to see.
\begin{lem}
 $\nbigp_{J_+}=-\nbigp_{J_-}$,
$\nbigq_{J_+}=-\nbigq_{J_-}$,
 and
 $\nbigr^{J_+}_{J_-}=-\nbigr^{J_-}_{J_+}+\nbigp_{J_-}\circ\nbigq_{J_-}$.
\hfill\qed
\end{lem}
\begin{rem}
$\nbigp_{J_-}$ and $\nbigq_{J_-}$
are often denoted by $\nbigp_J$ and $\nbigq_J$. 
\hfill\qed
\end{rem}

\subsubsection{}
We obtain the following morphisms
from the inclusions (\ref{eq;20.9.7.10}):
\[
 \nbigr^{J}_{J'}:
 H^0(J,L_{J,>0})
 \lrarr
 H^0(J',L_{J',<0})
 \quad
 (J<J'\leq J+\omega^{-1}\pi).
\]
Similarly, we obtain the following morphisms
from the inclusions (\ref{eq;20.9.7.13}):
\[
 \nbigr^{J}_{J'}:
 H^0(J,L_{J,>0})
 \lrarr
 H^0(J',L_{J',<0})
 \quad
 (J-\omega^{-1}\pi\leq J<J'<J).
\]
Note that
$\nbigr^{J}_{J+\omega^{-1}\pi}$
and
$\nbigr^{J}_{J-\omega^{-1}\pi}$
are isomorphisms.

\subsubsection{}
Let $\Phi_0^{J',J}:H^0(J,L_{J,0})\simeq H^0(J',L_{J',0})$
be the isomorphism obtained as
\[
 H^0(J,L_{J,0})
 \simeq
 H^0(\real,\Gr^{\vecnbigf}_0(L))
 \simeq
 H^0(J',L_{J',0}).
\]

\subsubsection{}
We introduce the maps
$\nbigrtilde^{J_-}_{J'}:
H^0(J,L_{J,>0})
\lrarr
H^0(J',L_{J',<0})$ for $J'\in T(\nbigi)$:
\[
 \nbigrtilde^{J_-}_{J'}:=
 \left\{
\begin{array}{ll}
 0& (J'\leq J-\omega^{-1}\pi) \\
 \nbigr^{J}_{J'} & (J-\omega^{-1}\pi\leq J'<J)\\
 \nbigr^{J_-}_{J_+} & (J'=J)\\
 \nbigr^{J}_{J'}+\nbigp_{J'_-}\circ\Phi^{J',J}_0\circ\nbigq_{J_-}
  & (J<J'\leq J+\omega^{-1}\pi)\\
 \nbigp_{J'_-}\circ\Phi^{J',J}_0\circ\nbigq_{J_-}
  &
  (J+\omega^{-1}\pi<J').
\end{array}
 \right.
\]
Similarly, we introduce the maps
$\nbigrtilde^{J_+}_{J'}:
H^0(J,L_{J,>0})
\lrarr
H^0(J',L_{J',<0})$ for $J'\in T(\nbigi)$:
\[
 \nbigrtilde^{J_+}_{J'}:=
 \left\{
\begin{array}{ll}
 \nbigp_{J'_+}\circ\Phi^{J',J}_0\circ\nbigq_{J_+}
  & (J'<J-\omega^{-1}\pi)\\
 \nbigr^{J}_{J'}+\nbigp_{J'_+}\circ\Phi^{J',J}_0\circ\nbigq_{J_+}
  & (J-\omega^{-1}\pi\leq J'<J)\\
  \nbigr^{J_+}_{J_-} & (J'=J)\\
 \nbigr^{J}_{J'} & (J<J'\leq J+\omega^{-1}\pi)\\
 0& (J+\omega^{-1}\pi\leq J') \\
\end{array}
 \right.
\]

\subsubsection{Hills}
\label{subsection;24.2.18.1}

For any $J\in T(\nbigi)$,
let $a_J:J\to \real$ denote the inclusion. 
There exists the natural isomorphism 
\[
 L/L^{\leq 0}
 \simeq
 \bigoplus_{J\in T(\nbigi)}
 a_{J\ast}(L_{J,>0}).
\]
We obtain the projection,
called a hill:
\index{hill} \index{map $R_J$}
\begin{equation}
\label{eq;21.6.7.3}
 R_J:H^0(\real,L)\lrarr
 H^0(J,L_{J,>0}).
\end{equation}

\begin{rem}
The map $R_J$ is also obtained as follows.
In the decompositions {\rm(\ref{eq;18.5.2.1})},
we have
$L_{J_+,\gminia}=L_{J_-,\gminia}$
unless $\gminia\in\nbigi_J$.
We obtain the local subsystem
$L_{\Jbar,\gminia}\subset L_{|\Jbar}$
by gluing $L_{J_{\pm},\gminia}$
for $\gminia\in\nbigi\setminus\nbigi_J$.
We obtain the decomposition
\[
 L_{|\Jbar}
=\gbiga_{\Jbar}(L)
 \oplus
 \bigoplus_{\gminia\in\nbigi\setminus\nbigi_J}
 L_{\Jbar,\gminia}.
\]
We obtain the projections
$L_{|\Jbar}\to \gbiga_{\Jbar}(L)
\lrarr L_{\Jbar,>0}$.
It induces $R_J$.
\hfill\qed
\end{rem}

\subsection{Appendix: Duality and hills}
\label{subsection;21.6.7.4}

Let $\nbigi\subset z_p^{-n}\cnum$
be a $\Gal(p)$-invariant subset.
Let $(L,\vecnbigf)$ be a $2\pi\seisuu$-equivariant local system
with Stokes structure indexed by $\nbigi$
on $\real$.

Let $L^{\lor}$  denote the $2\pi\seisuu$-equivariant
local system on $\real$
obtained as the dual of $L$.
We set $-\nbigi:=\{-\gminia\,|\,\gminia\in\nbigi\}$.
For each $\theta\in\real$ and for $\gminia\in\nbigi$,
let 
$\nbigf^{\theta}_{-\gminia}(L^{\lor}_{\theta})$
denote the subspace of $s\in L^{\lor}_{\theta}$
such that 
$s\bigl(\nbigf^{\theta}_{\gminib}(L_{\theta})\bigr)=0$
unless $\gminia\leq_{\theta}\gminib$.
It is easy to check the following lemma.
\begin{lem}
\label{lem;21.6.7.2}
 A splitting $L_{\theta}=\bigoplus_{\gminia\in\nbigi} G_{\theta,\gminia}$
of the filtration $\nbigf^{\theta}(L_{\theta})$
induces a splitting
$L^{\lor}_{\theta}=
\bigoplus_{-\gminia\in -\nbigi}
G^{\lor}_{\theta,-\gminia}$,
where
$G^{\lor}_{\theta,-\gminia}$
denotes the subspace of
$s\in L^{\lor}_{\theta}$
such that
$s(G_{\theta,\gminib})=0$
unless $\gminib=\gminia$.
\hfill\qed
\end{lem}
It is well known and easy to check
by using Lemma \ref{lem;21.6.7.2},
that the family of filtrations
$\nbigf^{\theta}(L^{\lor}_{\theta})$ $(\theta\in\real)$
is a $2\pi\seisuu$-equivariant Stokes structure of
$L^{\lor}$ indexed by $-\nbigi$.

Let $J\in T(\nbigi)$.
We obtain the local subsystems
$(L^{\lor})_{\Jbar,<0}
\subset
(L^{\lor})_{\Jbar,\leq 0}
\subset
\gbiga_{\Jbar}(L^{\lor})
\subset
(L^{\lor})_{|\Jbar}$.
We can easily check the following lemma
by using Lemma \ref{lem;21.6.7.2} on $J_{\pm}$.
\begin{lem}
\label{lem;21.6.7.3}
The natural perfect pairing between $L$ and $L^{\lor}$
induces a perfect paring of the local subsystems
$\gbiga_{\Jbar}(L)$ and $\gbiga_{\Jbar}(L^{\lor})$.
It induces perfect pairings of  
(i) $L_{\Jbar,<0}$
and $\gbiga_{\Jbar}(L^{\lor})/(L^{\lor})_{\Jbar,\leq 0}$,
 (ii) $L_{\Jbar,\leq 0}/L_{\Jbar,<0}$
 and $(L^{\lor})_{\Jbar,\leq 0}(L^{\lor})_{\Jbar,<0}$,
(iii) $\gbiga_{\Jbar}(L)/L_{\Jbar,\leq 0}$ and $(L^{\lor})_{\Jbar,<0}$.
\hfill\qed
\end{lem}

As a corollary of Lemma \ref{lem;21.6.7.3},
there exist the natural duality
\[
(L^{\lor})_{\Jbar,<0}
 \simeq
 (\gbiga_{\Jbar}(L)/L_{\Jbar,\leq 0})^{\lor}
=(L_{\Jbar,>0})^{\lor},
\]
\[
 (L^{\lor})_{\Jbar,0}
 =(L^{\lor})_{\Jbar,\leq 0}/(L^{\lor}_{\Jbar,<0})
 \simeq
  (L_{\Jbar,\leq 0}/L_{\Jbar,<0})^{\lor}
=(L_{\Jbar,0})^{\lor},
\]
\[
 (L^{\lor})_{\Jbar,>0}=
 \gbiga_{\Jbar}(L^{\lor})/(L^{\lor})_{\Jbar,\leq 0}
 \simeq
 (L_{\Jbar,<0})^{\lor}.
\]
In particular,
the natural inclusion
$(L^{\lor})_{\Jbar,<0}\lrarr (L^{\lor})_{|\Jbar}$
induces the projection
$L_{|\Jbar}\lrarr L_{\Jbar,>0}$.
In particular, we obtain the map
$R_J:H^0(\real,L)\lrarr
H^0(\Jbar,L_{\Jbar,>0})$
which equals the hill in \S\ref{subsection;24.2.18.1}.

\section{Extensions of local systems with Stokes structure}
\label{section;21.4.30.1}

Let $\nbigi$ be
a $\Gal(p)$-invariant finite subset of $z_p^{-1}\cnum[z_p^{-1}]$.
Let $\nbigj\subset\nbigi\cap z^{-1}\cnum[z^{-1}]$.
We use the notation in \S\ref{subsection;18.6.18.3}.
Let $\nbige$ be a functor from $\Dsf(\nbigj)$
to $\Loc^{\St}(\nbigi)$
such that
$\Gr^{\vecnbigf}_{\gminia}(\nbige(\varrho_1))
\lrarr
\Gr^{\vecnbigf}_{\gminia}(\nbige(\varrho_2))$
are isomorphisms
for any $\varrho_1,\varrho_2\in \Dsf(\nbigj)$
unless $\gminia\in\nbigj$
and $\varrho_1(\gminia)\neq\varrho_2(\gminia)$.
\begin{df}
\label{df;21.5.3.3}
Such a functor $\nbige$ is called a base tuple
in $\Loc^{\St}(\nbigi)$ with respect to $\nbigj$.
 \index{base tuple}
\hfill\qed
\end{df}

Let $\nbigc_1$ be the category of functors
$\nbigetilde:
\Csf(\nbigj)\lrarr \Loc^{\St}(\nbigi)$
equipped with an isomorphism
$a_{\nbigetilde}:\iota^{\ast}(\nbigetilde)\simeq\nbige$
such that
$\Gr^{\vecnbigf}_{\gminia}(\nbige(\varrho_1))
\lrarr
\Gr^{\vecnbigf}_{\gminia}(\nbige(\varrho_2))$
are isomorphisms
for any $\varrho_1,\varrho_2\in \Csf(\nbigj)$
unless $\gminia\in\nbigj$
and $\varrho_1(\gminia)\neq\varrho_2(\gminia)$.
A morphism
$f:(\nbigetilde_1,a_{\nbigetilde_1})
\lrarr
(\nbigetilde_2,a_{\nbigetilde_2})$ in $\nbigc_1$
is defined to be
a natural transformation
$f:\nbigetilde_1\lrarr\nbigetilde_2$
such that
$a_{\nbigetilde_2}\circ \iota^{\ast}(f)
=a_{\nbigetilde_1}$.

Let $\nbigc_2$
be the category of functors
$\nbigg$ from $\Csf(\nbigj)$
to the category of 
$2\pi\seisuu$-equivariant
$\nbigj$-graded local systems
$\nbigg$
equipped with an isomorphism
$b_{\nbigg}:
\iota^{\ast}\nbigg\simeq
\Gr^{\vecnbigf}_{\nbigj}(\nbige)$
such that
$\nbigg_{\gminia}(\varrho_1)
\lrarr
\nbigg_{\gminia}(\varrho_2)$
are isomorphisms for any
$\varrho_1,\varrho_2\in\Csf(\nbigj)$
unless $\varrho_1(\gminia)\neq\varrho_2(\gminia)$.

Any object $\nbigetilde$ of $\nbigc_1$
induces an object
$\Gr_{\nbigj}^{\vecnbigf}(\nbigetilde)$ in $\nbigc_2$.
Thus, we obtain a functor
$\nbigc_1\lrarr\nbigc_2$.
We shall prove the following proposition
in \S\ref{subsection;21.4.30.2}.
\begin{thm}
\label{thm;20.11.3.20}
The functor $\nbigc_1\lrarr\nbigc_2$
is an equivalence.
\end{thm}

Theorem \ref{thm;20.11.3.20}
implies that
for any $\nbigg$ in $\nbigc_2$,
there exists $\nbigetilde$ in $\nbigc_1$
which induces $\nbigg$.
Such $\nbigetilde$ is uniquely determined
up to canonical isomorphisms.
Such $\nbigetilde$ is called
an extension of $\nbige$ by $\nbigg$.
\index{extension}

\begin{df}
\label{df;21.5.3.1}
If $\nbigj=\{0\}$,
a functor $\nbige:\Dsf(\nbigj)\lrarr \Loc^{\St}(\nbigi)$ as above
is equivalent to
a morphism $F:(L_1,\vecnbigf)\lrarr(L_2,\vecnbigf)$ in
$\Loc^{\St}(\nbigi)$
such that
$\Gr^{\vecnbigf}_{\gminia}(L_1)\lrarr
\Gr^{\vecnbigf}_{\gminia}(L_2)$
is an isomorphism unless $\gminia=0$.
Such $((L_1,\vecnbigf),(L_2,\vecnbigf),F)$
is called a base tuple in $\Loc^{\St}(\nbigi)$.
\index{base tuple}
\hfill\qed
\end{df}

\begin{rem}
We obtain Proposition {\rm\ref{prop;20.11.3.10}}
by using Theorem {\rm\ref{thm;20.11.3.20}}
in the case $\nbigi\subset z_p^{-n}\cnum$
and $\nbigj=\{0\}$.
We obtain Proposition {\rm\ref{prop;21.4.28.20}}
by using
Theorem {\rm\ref{thm;20.11.3.20}}
in the case $\nbigi=\nbigj\subset z^{-1}\cnum$.
\hfill\qed
\end{rem}

\begin{rem}
In an earlier version of this monograph,
we originally proved Theorem {\rm\ref{thm;20.11.3.20}}
by using Stokes shells.
We explain another direct proof
on the basis of canonical splittings,
which would hopefully be easier to the readers.
\hfill\qed
\end{rem}

\subsection{Splittings}

We use the notation in \S\ref{subsection;20.11.14.31}.
Let $\nbigi$ be a $\Gal(p)$-invariant subset of
$z_p^{-1}\cnum[z_p^{-1}]$.
Let $(L,\vecnbigf)\in\Loc^{\St}(\nbigi)$.
For $\omega\in\frac{1}{p}\seisuu_{>0}$
and $\gminib\in \pi_{\omega}(\nbigi)$,
we set
$\nbigi(\gminib):=
\{\gminia\in\nbigi\,|\,\pi_{\omega}(\gminia)=\gminib\}$.
Assume that
$\bigl\{\omega'\in \frac{1}{p}\seisuu_{>0}
\,\big|\,
 |\pi_{\omega'}(\nbigi(\gminib))|\geq 2
 \bigr\}\neq \emptyset$.
 Let $\omega_1$ be the maximum of the set,
 which is strictly smaller than $\omega$.
If $J$ is an interval $\real$
such that 
$\bigl|J\cap S(\gminic_1,\gminic_2)\bigr|=1$
for any two distinct
$\gminic_1,\gminic_2\in
\pi_{\omega_1}(\nbigi(\gminib))$,
there exists a canonical splitting
\index{canonical splitting}
\begin{equation}
\label{eq;20.11.2.2}
 \Gr^{\vecnbigf^{(\omega)}}_{\gminib}(L)_{|J}
 =\bigoplus_{\gminic\in\pi_{\omega_1}(\nbigi(\gminib))}
 G_{J,\gminic},
\end{equation}
which induces a splitting of the filtration
$\nbigf^{(\omega_1)\,\theta}
\bigl(
 \Gr_{\gminib}^{\vecnbigf^{(\omega)}}(L)_{|\theta}
\bigr)$ 
for each $\theta\in J$.
Note that
$G_{J,\gminic}$
is naturally isomorphic to
$\Gr^{\vecnbigf^{(\omega_1)}}_{\gminic}(L)_{|J}$.
For any morphism 
$f:(L_1,\vecnbigf)\lrarr
(L_2,\vecnbigf)$ in $\Loc^{\St}(\nbigi)$,
the induced morphism
$\Gr^{\vecnbigf^{(\omega)}}_{\gminib}(f):
\Gr^{\vecnbigf^{(\omega)}}_{\gminib}(L_1)
\lrarr
\Gr^{\vecnbigf^{(\omega)}}_{\gminib}(L_2)$
preserves the canonical splittings.

In particular,
for any $\theta_0\in\real$
such that 
\[
\theta_0\not\in
\bigcup_{\substack{\gminic_1,\gminic_2\in\pi_{\omega_1}(\nbigi(\gminib))
 \\
\gminic_1\neq\gminic_2}}
 A(\gminic_1,\gminic_2),
\]
by setting 
$J=\bigl\{\theta\in\real\,\big|\,
|\theta-\theta_0|<\omega_1^{-1}\pi/2\bigr\}$,
we obtain the canonical splitting (\ref{eq;20.11.2.2})
of $\Gr^{\vecnbigf^{(\omega)}}_{\gminib}(L)_{|J}$.
As the restriction to $\theta_0$,
we obtain a splitting
\begin{equation}
\label{eq;20.11.3.1}
 \Gr^{\vecnbigf^{(\omega)}}_{\gminib}(L)_{|\theta_0}
 =\bigoplus_{\gminic\in\pi_{\omega_1}(\nbigi(\gminib))}
 \Gr^{\vecnbigf^{(\omega)}}_{\gminib}(L)_{\theta_0,\gminic}
\end{equation}
of the filtration
$\nbigf^{(\omega_1)\theta_0}
\Gr_{\gminib}^{\vecnbigf^{(\omega)}}(L)_{|\theta_0}$.
Note that
$\Gr^{\vecnbigf^{(\omega)}}_{\gminib}(L)_{\theta_0,\gminic}$
is naturally isomorphic to
\[
 \Gr^{\vecnbigf^{(\omega_1)}}_{\gminic}
\Gr^{\vecnbigf^{(\omega)}}(L)_{|\theta_0}
=\Gr^{\vecnbigf^{(\omega_1)}}_{\gminic}(L)_{|\theta_0}.
\]
For any morphism
$f:(L_1,\vecnbigf)\lrarr (L_2,\vecnbigf)$
in $\Loc^{\St}(\nbigi)$,
the induced morphism
$\Gr^{\vecnbigf^{(\omega)}}_{\gminib}(f)$
preserves the decompositions as in (\ref{eq;20.11.3.1}).

We may apply the construction of splittings successively.
Take $\theta_0\in\real$ such that
\begin{equation}
 \theta_0\not\in
 A(\nbigi)=
 \bigcup_{\substack{\gminia_1,\gminia_2\in\nbigi\\
  \gminia_1\neq\gminia_2
  }}
   A(\gminia_1,\gminia_2).
\end{equation}
Then, there uniquely exists a splitting
\begin{equation}
\label{eq;20.11.3.3}
 L_{|\theta_0}=
  \bigoplus_{\gminia\in\nbigi}
  L_{\theta_0,\gminia}
\end{equation}
of the filtration $\nbigf(L_{|\theta_0})$
such that the following holds.
\begin{itemize}
 \item 
For any $\gminib\in\pi_{\omega}(\nbigi)$,
under the natural isomorphism
\[
 \bigoplus_{\gminia\in\nbigi(\gminib)}
  L_{\theta_0,\gminia}
 \simeq
 \Gr^{\vecnbigf^{(\omega)}}_{\gminib}(L)_{|\theta_0},
\]
we obtain
\[
  \bigoplus_{\substack{\gminia\in\nbigi(\gminib)
   \\
    \pi_{\omega_1}(\gminia)=\gminic
   }}
   L_{\theta_0,\gminia}
   =
  \Gr^{\vecnbigf^{(\omega)}}_{\gminib}(L)_{\theta_0,\gminic}
\]
for any $\gminic\in\pi_{\omega_1}(\nbigi(\gminib))$.
\end{itemize}
       
Let $\theta_1\in A(\nbigi)$.
Let $\epsilon>0$ be so small that
$\{\theta_1-\epsilon\leq\theta\leq\theta_1+\epsilon\}
\cap A(\nbigi)=\{\theta_1\}$.
The natural isomorphism
$\Phi^{\theta_1+\epsilon,\theta_1-\epsilon}:
L_{|\theta_1-\epsilon}\simeq L_{|\theta_1+\epsilon}$
is contained in
\[
\bigoplus_{\gminia\in\nbigi}
\Hom\bigl(
 L_{\theta_1-\epsilon,\gminia},
 L_{\theta_1+\epsilon,\gminia}
\bigr)
\oplus
\bigoplus_{\substack{
\gminia_1,\gminia_2\in\nbigi\\
\theta_1\in A(\gminia_1,\gminia_2)\\
\gminia_1>_{\theta_1}\gminia_2}}
\Hom\Bigl(
 L_{\theta_1-\epsilon,\gminia_1},
 L_{\theta_1+\epsilon,\gminia_2}
\Bigr).
\]

\subsection{Some objects equivalent to local systems with Stokes structure}

There are several objects which are equivalent to
$2\pi\seisuu$-equivariant
local systems with Stokes structure,
as explained in \cite{Boalch-survey}.
Here, we explain some minor variants
of ``Stokes local systems''.
(Stokes local systems
are originally introduced by Boalch
as more geometric objects in a sophisticated way.)

\subsubsection{}
Let $L_{\bullet}=\bigoplus_{\gminia\in\nbigi}L_{\gminia}$
be an $\nbigi$-graded local system
equipped with a $2\pi\seisuu$-action
such that
$(2\pi\ell)^{\ast}L_{\gminia}
=L_{2\pi\ell\bullet\gminia}$.
For each $\theta_1\in A(\nbigi)$,
let
$\Sto_{\theta_1}(L_{\bullet})$
denote the group of the automorphisms
$\varphi$ of $L_{\bullet|\theta_1}$
such that
\[
 \varphi-\id
 \in
 \bigoplus_{\substack{
 \gminia_1,\gminia_2\in\nbigi\\
 \theta_1\in A(\gminia_1,\gminia_2)\\
 \gminia_1>_{\theta}\gminia_2
 }}
 \Hom\Bigl(
 L_{\gminia_1|\theta},L_{\gminia_2|\theta}
 \Bigr).
\] \index{group $\Sto_{\theta_1}(L_{\bullet})$}
Here, $\id$ denotes the identity map.
A tuple
\[
 \vecvarphi=
\bigl(
 \varphi_{\theta_1}\,\big|\,
 \theta_1\in A(\nbigi)
 \bigr)
 \in\prod_{\theta_1\in A(\nbigi)}
  \Sto_{\theta_1}(L_{\bullet})
\]
is called $2\pi\seisuu$-equivariant
if
$\varphi_{\theta_1+2\pi\ell}
=\varphi_{\theta_1}$
under the isomorphism
$L_{\bullet|\theta_1+2\pi\ell}
=(2\pi\ell)^{\ast}(L_{\bullet})_{|\theta_1}
\simeq L_{\bullet|\theta_1}$.

Let $\Loc^{\Sto}(\nbigi)$
denote the category of
$2\pi\seisuu$-equivariant
$\nbigi$-graded local systems $L$
equipped with a $2\pi\seisuu$-equivariant
tuple
$\vecvarphi\in
\prod_{\theta_1\in A(\nbigi)}\Sto_{\theta_1}(L_{\bullet})$.
\index{category $\Loc^{\Sto}(\nbigi)$}
A morphism
$f:(L_{1\bullet},\vecvarphi_1)
\lrarr
 (L_{2\bullet},\vecvarphi_2)$
is defined to be
a morphism of $2\pi\seisuu$-equivariant
$\nbigi$-graded local systems
$f:L_{1\bullet}\lrarr L_{2\bullet}$
such that
$(\varphi_2)_{\theta_1}\circ f_{|\theta_1}
=f_{|\theta_1}\circ(\varphi_1)_{\theta_1}$
for any $\theta_1\in A(\nbigi)$.

As explained in \cite{Boalch-survey},
$\Loc^{\St}(\nbigi)$
and $\Loc^{\Sto}(\nbigi)$
are equivalent.
For any $(L,\vecnbigf)\in \Loc^{\St}(\nbigi)$,
we obtain a $2\pi\seisuu$-equivariant
$\nbigi$-graded local system $\Gr^{\vecnbigf}(L)$.
For each $\theta_1\in A(\nbigi)$,
we obtain the automorphism
$\varphi_{\theta_1}$
of $\Gr^{\vecnbigf}(L)_{|\theta_1}$
as the composition of
\[
 \Gr^{\vecnbigf}(L)_{|\theta_1}
 \stackrel{b_1}{\simeq}
  \Gr^{\vecnbigf}(L)_{|\theta_1-\epsilon}
 \stackrel{a_1}{\simeq}
  L_{\theta_1-\epsilon}
 \stackrel{b_2}{\simeq}
  L_{\theta_1+\epsilon}
 \stackrel{a_2}{\simeq}
 \Gr^{\vecnbigf}(L)_{|\theta_1+\epsilon}
 \stackrel{b_3}{\simeq}
  \Gr^{\vecnbigf}(L)_{|\theta_1},
\]
where $a_i$ are induced
by the splittings as in (\ref{eq;20.11.3.3}),
and $b_i$ are the parallel transport.
It is easy to see that
$\varphi_{\theta_1}\in \Sto_{\theta_1}(L_{\bullet})$,
and
$\vecvarphi=(\varphi_{\theta_1})$ is $2\pi\seisuu$-equivariant.
This procedure induces
a functor
$\Loc^{\St}(\nbigi)\lrarr\Loc^{\Sto}(\nbigi)$.
It is the desired equivalence.
A quasi-inverse of the functor is also obtained
as follows.
Let $(L_{\bullet},\vecvarphi)\in\Loc^{\Sto}(\nbigi)$.
For each $\theta\in \Ibar$,
the $\nbigi$-grading and the order
$(\nbigi,\leq_{\theta})$
induce a filtration $\nbigf^{\theta}$
of $L_{\bullet|\theta}$.
Let $C(\nbigi)$ denote the set of the connected components of
$\real\setminus S(\nbigi)$.
For each $I\in C(\nbigi)$,
we obtain an $\nbigi$-graded local system
$L_{\bullet|\Ibar}$ on the closure $\Ibar$ of $I$ in $\real$.
For $\theta_1\in A(\nbigi)$,
there are two distinct
$I_j\in C(\nbigi)$ $(j=1,2)$
such that
$I_1=\{\theta_1'<\theta<\theta_1\}$
and $I_2=\{\theta_1<\theta<\theta_1''\}$
for some $\theta_1',\theta_1''$.
We may regard $\varphi_{\theta_1}$
as the isomorphism
$(L_{\bullet|\Ibar_1})_{|\theta_1}
\simeq
(L_{\bullet|\Ibar_2})_{|\theta_1}$.
We glue
$L_{\bullet|\Ibar}$ $(I\in C(\nbigi))$
by the isomorphisms,
and we obtain a $2\pi\seisuu$-equivariant
local system $L$ on $\real$.
Because $\varphi_{\theta_1}$ preserves the filtrations
$\nbigf^{\theta_1}$,
we obtain a family of filtrations
$\nbigf^{\theta}$ of $L_{|\theta}$ $(\theta\in\real)$,
which is a $2\pi\seisuu$-equivariant Stokes structure
on $L$.

\subsubsection{}

We may obviously consider intermediate objects.
Let $\omega\in\rnum_{>0}$.
Let $L_{\bullet}$ be a $2\pi\seisuu$-equivariant
$\pi_{\omega}(\nbigi)$-graded local system.
A $2\pi\seisuu$-equivariant Stokes structure $\vecnbigf$
of $L_{\bullet}$ indexed by $\nbigi$
is called $\pi_{\omega}(\nbigi)$-graded
if
\[
 \nbigf^{\theta}_{\gminia}(L_{\bullet|\theta})
 =\bigoplus_{\gminib\in\pi_{\omega}(\nbigi)}
  \nbigf^{\theta}_{\gminia}
  \bigl(
   L_{\gminib|\theta}
  \bigr)
\]
for any $\gminia\in\nbigi$,
and moreover
$\Gr^{\nbigf^{\theta}}_{\gminia}(L_{\gminib|\theta})=0$
unless $\pi_{\omega}(\gminia)=\gminib$.
\index{$\pi_{\omega}(\nbigi)$-graded Stokes structure
indexed by $\nbigi$}
For any $\theta_1\in A(\pi_{\omega}(\nbigi))$,
we define
$\Sto_{\theta_1}(L_{\bullet})$
as before by replacing $\nbigi$ with $\pi_{\omega}(\nbigi)$.

Let $\Loc^{\Sto}(\nbigi,\pi_{\omega}\nbigi)$
denote the category of
$2\pi\seisuu$-equivariant
$\pi_{\omega}(\nbigi)$-graded
local systems
$L_{\bullet}$
equipped with a $2\pi\seisuu$-equivariant
$\pi_{\omega}(\nbigi)$-graded Stokes structure $\vecnbigf$
indexed by $\nbigi$,
and a $2\pi\seisuu$-equivariant tuple
\[
\vecvarphi\in\prod_{\theta_1\in A(\pi_{\omega}(\nbigi))}
 \Sto_{\theta_1}(L_{\bullet}).
\]
\index{category $\Loc^{\Sto}(\nbigi,\pi_{\omega}\nbigi)$}
A morphism
$f:(L_{1\bullet},\vecnbigf,\vecvarphi_1)
\lrarr
(L_{2\bullet},\vecnbigf,\vecvarphi_2)$
is defined to be
a morphism of $2\pi\seisuu$-equivariant
$\pi_{\omega}(\nbigi)$-graded local systems
$f:L_{1\bullet}\lrarr L_{2\bullet}$
such that
$f:(L_{1\bullet},\vecnbigf)
\lrarr (L_{2\bullet},\vecnbigf)$
is a morphism in $\Loc^{\St}(\nbigi)$,
and that
$f:(L_{1\bullet},\vecvarphi_1)
\lrarr (L_{2\bullet},\vecvarphi_2)$
is a morphism in
$\Loc^{\Sto}(\pi_{\omega}(\nbigi))$.
It is easy to see that
$\Loc^{\St}(\nbigi)$
and $\Loc^{\Sto}(\nbigi,\pi_{\omega}(\nbigi))$
are equivalent.

\subsection{Proof of Theorem \ref{thm;20.11.3.20}}
\label{subsection;21.4.30.2}

\label{subsection;20.11.4.2}

Let us construct a quasi-inverse
$\nbigc_2\lrarr\nbigc_1$.
Let
$(\Gr^{\vecnbigf}(\nbige),\vecvarphi)$
denote the functor from
$\Dsf(\nbigj)$ to $\Loc^{\Sto}(\nbigi)$
corresponding to $\nbige$.
Let $(\nbigg,b_{\nbigg})\in\nbigc_2$.
We shall construct an object
$(\nbigp(\nbigg),\vecvarphi_{\nbigg})$
of $\nbigc_1$.
For any $\varrho\in \Dsf(\nbigj)$,
we set
\[
 \nbigp(\nbigg)_{\gminia}(\varrho)
 =\left\{
\begin{array}{ll}
 \Gr^{\vecnbigf}_{\gminia}(\nbige)(\underline{!})
   & (\gminia\not\in\nbigj)\\
 \nbigg_{\gminia}(\varrho)
   &(\gminia\in\nbigj).
\end{array}
 \right.
\]

For any $\theta_1\in A(\nbigi)$
and $\varrho\in\Csf(\nbigj)$,
we define
$\varphi_{\nbigg}(\varrho)_{\theta_1}$ as follows.
\begin{itemize}
 \item
      $(\varphi_{\nbigg}(\varrho)_{\theta_1})_{\gminia,\gminia}$
      are the identity map for any $\gminia\in\nbigi$.
 \item $(\varphi_{\nbigg}(\varrho)_{\theta_1})_{\gminia_1,\gminia_2}
       =(\varphi(\underline{!})_{\theta_1})_{\gminia_1,\gminia_2}$
       if $\gminia_i\not\in\nbigj$ $(i=1,2)$.
 \item If $\gminia_2\not\in \nbigj$ and $\gminia_1\in\nbigj$,
       $(\varphi_{\nbigg}(\varrho)_{\theta_1})_{\gminia_1,\gminia_2}$
       is the composite of the following morphisms:
       \[
\begin{CD}
       \nbigp(\nbigg)_{\gminia_2}(\varrho)
       =
       \Gr^{\vecnbigf}_{\gminia_2}(\nbige)(\underline{!})
       @>{\varphi(\underline{!})_{\theta_1}}>>
       \Gr^{\vecnbigf}_{\gminia_1}(\nbige)(\underline{!})
       @>>>
       \nbigg_{\gminia_1}(\varrho(\gminia_1))
       =\nbigp(\nbigg)_{\gminia_1}(\varrho).
\end{CD}
       \]
\item If $\gminia_1\not\in \nbigj$ and $\gminia_2\in\nbigj$,
       $(\varphi_{\nbigg}(\varrho)_{\theta_1})_{\gminia_1,\gminia_2}$
      is the composite of the following morphisms:
\begin{multline}
 \nbigp(\nbigg)_{\gminia_2}(\varrho)
 =\nbigg_{\gminia_2}(\varrho(\gminia_2))
 \lrarr
 \Gr^{\vecnbigf}_{\gminia_2}(\nbige)(\underline{\ast})
 \stackrel{\varphi(\underline{\ast})_{\theta_1}}{\lrarr}
 \Gr^{\vecnbigf}_{\gminia_1}(\nbige)(\underline{\ast})
 \stackrel{\simeq}{\lrarr} \\
 \Gr^{\vecnbigf}_{\gminia_1}(\nbige)(\underline{!})
 =\nbigp(\nbigg)_{\gminia_1}(\varrho).
\end{multline}
 \item Suppose that $\gminia_i\in\nbigj$ $(i=1,2)$.
       There exists
       $\varrho'\in \Dsf(\nbigj)$
       such that
       $\varrho'(\gminia_1)=\ast$
       and $\varrho'(\gminia_2)=!$.
       Then, we define
       $(\varphi_{\nbigg}(\varrho)_{\theta_1})_{\gminia_1,\gminia_2}$
       as the composite of the following morphisms:
\begin{multline}
  \nbigp(\nbigg)_{\gminia_2}(\varrho)
 =\nbigg_{\gminia_2}(\varrho(\gminia_2))
 \lrarr
 \Gr^{\vecnbigf}_{\gminia_2}(\nbige)(\underline{\ast})
 \simeq
  \Gr^{\vecnbigf}_{\gminia_2}(\nbige)(\varrho')
 \stackrel{\varphi(\varrho')_{\theta_1}}{\lrarr}
 \\
 \Gr^{\vecnbigf}_{\gminia_1}(\nbige)(\varrho')
 \simeq
 \Gr^{\vecnbigf}_{\gminia_1}(\nbige)(\underline{!})
\lrarr
 \nbigg_{\gminia_1}(\varrho(\gminia_1))
 =  \nbigp(\nbigg)_{\gminia_1}(\varrho).
\end{multline}
\end{itemize}
Let $\varphi_{\nbigg}(\varrho)_{\theta_1}$
be the automorphism of
$\nbigp(\nbigg)_{\bullet}(\varrho)_{|\theta_1}$
obtained as the sum of
$(\varphi_{\nbigg}(\varrho)_{\theta_1})_{\gminia_1,\gminia_2}$
for pairs $(\gminia_1,\gminia_2)$
such that
$\theta_1\in A(\gminia_1,\gminia_2)$
and $\gminia_1\leq_{\theta_1}\gminia_2$.
Thus, we obtain a functor
$(\nbigp(\nbigg),\vecvarphi_{\nbigg})$
from $\Csf(\nbigj)$
to $\Loc^{\Sto}(\nbigi)$.
By the construction,
there exists a natural isomorphism
$\iota^{\ast}(\nbigp(\nbigg),\vecvarphi_{\nbigg})
\simeq
 (\Gr^{\vecnbigf}(\nbige),\vecvarphi)$.
The corresponding object
$(\nbigetilde(\nbigg),a_{\nbigetilde(\nbigg)})$ of $\nbigc_1$
induces $(\nbigg,b_{\nbigg})$.
It is easy to see that this construction
induces a quasi-inverse $\nbigc_2\lrarr\nbigc_1$.
\hfill\qed

\subsection{A simple case}
\label{subsection;21.5.1.1}

We use the notation in \S\ref{subsection;21.4.29.1}.
Let 
$\nbiga_0=\bigl(
\nbigt_{\omega}(L_1)
\stackrel{b_1}{\lrarr}
N
\stackrel{b_2}{\lrarr}
\nbigt_{\omega}(L_2)
\bigr)$
be an object of $\nbigc_2$.
According to Proposition \ref{prop;20.11.3.10},
which is a special case of Theorem \ref{thm;20.11.3.20},
we have the corresponding object $\nbiga_1$ of $\nbigc_1$.
Let us compute it 
in terms of $\Loc^{\Sto}(\nbigi,\pi_{\omega}(\nbigi))$
in this particular case.
Let
\[
\Gr^{\vecnbigf^{(\omega)}}(f):
(\Gr^{\vecnbigf^{(\omega)}}(L_1,\vecnbigf),\vecvarphi_1)
\lrarr
(\Gr^{\vecnbigf^{(\omega)}}(L_2,\vecnbigf),\vecvarphi_2) 
\]
be the morphism
in $\Loc^{\Sto}(\nbigi,\pi_{\omega}(\nbigi))$
induced by $f$.
We obtain the following
$2\pi\seisuu$-equivariant local system on $\real$:
\[
 P(N)_{\bullet}:=
 N
 \oplus
 \bigoplus_{\gminib\neq 0}
  \Gr^{\vecnbigf^{(\omega)}}_{\gminib}(L_1).
\]
For $\theta_1\in A(\nbigi)$,
we define
$\varphi_{N,\theta_1}\in
\Sto_{\theta_1}(P(N)_{\bullet|\theta_1})$
as follows:
\begin{itemize}
 \item
      $(\varphi_{N,\theta_1})_{\gminib,\gminib}$
      are the identity for any $\gminib\in\pi_{\omega}(\nbigi)$.
 \item $(\varphi_{N,\theta_1})_{\gminib_1,\gminib_2}
       =(\varphi_{1,\theta_1})_{\gminib_1,\gminib_2}$
       if $\gminib_i\neq 0$ $(i=1,2)$.
 \item For $\gminib\neq 0$,
       $(\varphi_{N,\theta_1})_{0,\gminib}$
       is the composition of
       \[
\begin{CD}
       \Gr^{\vecnbigf^{(\omega)}}_{\gminib}(L_{1\bullet|\theta})
 @>{(\varphi_{1,\theta_1})_{0,\gminib}}>>
 \Gr^{\vecnbigf^{(\omega)}}_{0}(L_{1\bullet|\theta})
 @>>> N.
\end{CD}
       \]
 \item For
       $\gminib\neq 0$,
       $(\varphi_{N,\theta_1})_{\gminib,0}$
       is the composition of
       \[
\begin{CD}
 N@>>> \Gr^{\vecnbigf^{(\omega)}}_{0}(L_2)
 @>{(\varphi_2,\theta_1)_{\gminib,0}}>>
 \Gr^{\vecnbigf^{(\omega)}}_{\gminib}(L_2)
 @>{\Gr^{\vecnbigf^{(\omega)}}_{\gminib}(f)^{-1}}>{\simeq}>
  \Gr^{\vecnbigf^{(\omega)}}_{\gminib}(L_1).
\end{CD}	
       \]
\end{itemize}
Thus, we obtain morphisms
\[
 \begin{CD}
(\Gr^{\vecnbigf^{(\omega)}}(L_1),\vecvarphi_1)
  @>{c_1}>>
  (P(N)_{\bullet},\vecvarphi_{N})
   @>{c_2}>>
(\Gr^{\vecnbigf^{(\omega)}}(L_2),\vecvarphi_2)
 \end{CD}
\]
in $\Loc^{\Sto}(\nbigi,\pi_{\omega}(\nbigi))$.
Thus, we obtain the desired object $\nbiga_1$ of $\nbigc_1$.

\section{Recovery of Stokes filtrations}

For any
$\gminia=\sum_{j=1}^n\gminia_jz_p^{-j}
\in z_p^{-1}\cnum[z_p^{-1}]\setminus\{0\}$
with $\gminia_n\neq 0$,
we set $\ord(\gminia):=-\frac{n}{p}$.
\index{order $\ord(\gminia)$}

Let $\nbigitilde\subset z_p^{-1}\cnum[z_p^{-1}]$ be
a $\Gal(p)$-invariant finite subset
such that $\nbigitilde\neq\{0\}$.
We set
$\omega=\max\{-\ord(\gminia)\,|\,\gminia\in\nbigitilde\setminus\{0\}\}$.
We assume $\nbigitilde=\nbigt_{\omega}(\nbigitilde)$.
We set $\nbigi:=\pi_{\omega}(\nbigitilde)\subset z^{-\omega}\cnum$.
For any $J\in T(\nbigi)$,
let $a_J:J\to \real$
denote the inclusions of $J$.
For any $J\in T(\nbigi)$,
we set $\nbigitilde_{J,<0}=\pi_{\omega}^{-1}(\nbigi_{J,<0})$
and $\nbigitilde_{J,>0}=\pi_{\omega}^{-1}(\nbigi_{J,>0})$.

\subsection{The induced constructible subsheaves and filtrations}

Let $(L,\vecnbigftilde)\in\Loc^{\St}(\nbigitilde)$.
We set
$(L,\vecnbigf):=(L,\pi_{\omega}(\vecnbigftilde))
\in \Loc^{\St}(\nbigi)$.
We obtain the constructible subsheaves
$L^{<0}\subset L^{\leq 0}\subset L$
determined by $\vecnbigftilde$,
or equivalently determined by $\vecnbigf$.
There exist the decompositions
\begin{equation}
\label{eq;25.2.23.1}
 L^{<0}=\bigoplus_{J\in T(\nbigi)}
 a_{J!}\bigl(L_{J,<0}\bigr),
 \quad
 L/L^{\leq 0}
 =\bigoplus_{J\in T(\nbigi)}
 a_{J\ast}\bigl(
 L_{J,>0}
 \bigr).
\end{equation}

Let $J\in T(\nbigi)$.
By using the decompositions
$L_{J_{\pm},<0}=
\bigoplus_{\gminia\in\nbigi_{J,<0}}L_{J_{\pm},\gminia}$,
we obtain the filtrations
$\nbigf^{\theta}$ $(\theta\in J)$
on $H^0(J,L_{J,<0})$
indexed by
$(\nbigi_{J,<0},\leq_{\theta})$,
which is independent of the choice of $\pm$.
For $\gminia\in\nbigi_{J,<0}$, we have
\[
\Gr^{\nbigf^{\theta}}_{\gminia}H^0(J,L_{J,<0})
\simeq
H^0(J,\Gr^{\vecnbigf}_{\gminia}(L)).
\]
Let $\nbigitilde(\gminia):=
\{\gminib\in\nbigitilde\,|\,\pi_{\omega}(\gminib)=\gminia\}$.
The filtrations $\nbigftilde^{\theta}$
on $\Gr^{\vecnbigf}_{\gminia}(L)$
indexed by $(\nbigitilde(\gminib),\leq_{\theta})$,
which induces filtrations on 
$\Gr^{\nbigf^{\theta}}_{\gminia}H^0(J,L_{J,<0})$.
They induce filtrations
$\nbigftilde^{\theta}$ $(\theta\in J)$ on $H^0(J,L_{J,<0})$
indexed by $(\nbigitilde_{J,<0},\leq_{\theta})$.
Similarly, we obtain the induced filtrations
$\nbigftilde^{\theta}$ $(\theta\in J)$
on $H^0(J,L_{J,>0})$ indexed by
$(\nbigitilde_{J,>0},\leq_{\theta})$.

\subsection{Recovery of the Stokes filtrations}

Let $(L_i,\vecnbigftilde)\in\Loc^{\St}(\nbigitilde)$ $(i=1,2)$.
We obtain the constructible subsheaves
$L_i^{<0}\subset L_i^{\leq 0}\subset L_i$
determined by $\vecnbigftilde$.

Let $\varphi:L_1\to L_2$ be a morphism of
$2\pi\seisuu$-equivariant local systems.
\begin{lem}
\label{lem;25.2.24.2}
Suppose that $\varphi$ induces morphisms of constructible subsheaves
$L_1^{<0}\to L_2^{<0}$
and
$L_1^{\leq 0}\to L_2^{\leq 0}$.
The induced morphisms
$L_1^{<0}\to L_2^{<0}$
and 
$L_1/L_1^{\leq 0}\to L_2/L_2^{\leq 0}$
are compatible with the decompositions in {\rm(\ref{eq;25.2.23.1})}.
In particular, $\varphi$ induces morphisms
\begin{equation}
\label{eq;25.2.25.2}
H^0(J,(L_1)_{J,<0})\to
H^0(J,(L_2)_{J,<0}),
\quad
H^0(J,(L_1)_{J,>0})\to
H^0(J,(L_2)_{J,>0})
\end{equation}
for any $J\in T(\nbigi)$.
\end{lem}
\pf
Let $J_1,J_2\in T(\nbigi)$ with $J_1\neq J_2$.
For any local systems $M_i$ on $J_i$,
and for $\star=!,\ast$,
any morphism of constructible sheaves
$a_{J_1\star}(M_1)
\to
a_{J_2\star}(M_2)$
is $0$.
Hence, we obtain the claim of the lemma.
\hfill\qed

\begin{prop}
\label{prop;25.2.23.5}
$\varphi$ gives a morphism 
$(L_1,\vecnbigftilde)\to (L_2,\vecnbigftilde)$
in $\Loc^{\St}(\nbigitilde)$
if and only if the following conditions are satisfied.
\begin{itemize}
 \item $\varphi$ induces morphisms
       of constructible subsheaves
       $L_1^{<0}\to L_2^{<0}$
       and
       $L_1^{\leq 0}\to L_2^{\leq 0}$.
 \item The induced morphisms
       {\rm(\ref{eq;25.2.25.2})}
       are compatible with
       the induced Stokes filtrations
       $\nbigftilde^{\theta}$ $(\theta\in J)$.
\end{itemize}
\end{prop}
\pf
The ``only if'' part is clear.
Let us study the ``if'' part.
Let $\theta\in\real\setminus S_0(\nbigi)$.
We have the subspaces
\[
 (L_i^{<0})_{\theta}
 =\bigoplus_{\theta\in J}
 ((L_i)_{J,<0})_{\theta}
 \subset
 (L_i^{\leq 0})_{\theta}
 \subset (L_i)_{\theta}.
\]
Let $G_{i,\theta,0}$ be any subspace
of $(L_i^{\leq 0})_{\theta}$
such that the projection
$(L_i^{\leq 0})_{\theta}
\to (L_i^{\leq 0})_{\theta}/(L_i^{<0})_{\theta}$
induces an isomorphism 
$G_{i,\theta,0}\simeq
(L_i^{\leq 0})_{\theta}/(L_i^{<0})_{\theta}$.
For any $J\in T(\nbigi)$ such that $\theta\in J$,
let $G_{i,\theta,J,>0}\subset (L_i)_{\theta}$
be a subspace
such that the projection
$(L_i)_{\theta}\to (L_i)_{\theta}/(L_i^{\leq 0})_{\theta}$
induces
$G_{i,\theta,J,>0}\simeq (a_{\Jbar\ast}(L_i)_{\Jbar,>0})_{\theta}$.
We obtain the decomposition
\[
 (L_i)_{\theta}
 =\bigoplus_{\theta\in J}
 \Bigl(
 ((L_i)_{J,<0})_{\theta}
 \oplus G_{i,\theta,J,>0}
 \Bigr)
 \oplus
 G_{i,\theta,0}. 
\]
By Lemma \ref{lem;25.2.24.2},
we have
$\varphi\bigl(
((L_1)_{J,<0})_{\theta}
\bigr)
\subset
((L_2)_{J,<0})_{\theta}$,
and we may also assume that
$\varphi(G_{1,\theta,J,>0})
\subset G_{2,\theta,J,>0}$
and
$\varphi(G_{1,\theta,0})
\subset G_{2,\theta,0}$.

The Stokes filtrations $\nbigftilde^{\theta}$
on $(L_i)_{\theta}$
equal the filtrations induced by 
the filtrations $\nbigftilde^{\theta}$
on the spaces
$((L_i)_{J,<0})_{\theta}
\simeq
H^0(J,(L_i)_{J,<0})$
and 
$G_{i,\theta,J,>0}
\simeq
H^0(J,(L_i)_{J,>0})$,
and the trivial filtrations on
$G_{i,\theta,0}$.
Hence, 
$\varphi:(L_1)_{\theta}\to (L_2)_{\theta}$
preserves the filtrations
$\nbigftilde^{\theta}$.
As remarked in Lemma \ref{lem;25.2.23.4},
we obtain that 
$\varphi$ gives a morphism
$(L_1,\vecnbigftilde)\to (L_2,\vecnbigftilde)$
in $\Loc^{\St}(\nbigitilde)$.
\hfill\qed

\begin{prop}
\label{prop;25.3.16.11}
The $2\pi\seisuu$-equivariant Stokes structure $\vecnbigf$ on $L$
is determined by the following data.
\begin{itemize}
 \item The constructible subsheaves
       $L^{<0}\subset L^{\leq 0}\subset L$.
 \item The tuple of filtrations
       $\nbigftilde^{\theta}$ $(\theta\in J)$
       on $H^0(J,L_{J,<0})$
       indexed by
       $(\nbigi_{J,<0},\leq_{\theta})$.
 \item The tuple of filtrations
       $\nbigftilde^{\theta}$ $(\theta\in J)$
       on $H^0(J,L_{J,>0})$
       indexed by
       $(\nbigi_{J,>0},\leq_{\theta})$.
       \hfill\qed
\end{itemize}
\end{prop}

\chapter{Stokes shells}
\label{section;18.6.3.11}

\section{Preliminary}

\subsection{Notation}

We shall use the notation in \S\ref{subsection;18.5.2.10}.
Let $\nbigi\subset z_p^{-n}\cnum$
be a non-empty $\Gal(p)$-invariant subset.
We set $\nbigi^{\ast}=\nbigi\setminus\{0\}$.
\index{set $\nbigi^{\ast}$}
We define the equivalence relation $\sim$
on $\nbigi$  as follows.
\begin{itemize}
\item
$\gminia\sim\gminib$
if and only if 
there exists $a>0$ such that 
$\gminia=a\gminib$.
\end{itemize}
We set $[\nbigi]:=\nbigi/\sim$. \index{set $[\nbigi]$}
For each $\lambda\in[\nbigi]$,
let $\nbigi(\lambda)\subset\nbigi$
denote the inverse image of
$\lambda$ by the projection
$\nbigi\lrarr[\nbigi]$.
\index{set $\nbigi(\lambda)$}
A $\Gal(p)$-action
on $[\nbigi]$
is naturally induced by the $\Gal(p)$-action on $\nbigi$.
In particular,
the automorphism $\Tbb^{\ast}$ of $\nbigi$
induces an automorphism of $[\nbigi]$,
which is also denoted by $\Tbb^{\ast}$.
(See \S\ref{subsection;18.5.2.10}
for $\Tbb^{\ast}$.)
\index{map $\Tbb^{\ast}$}

For $\lambda\in[\nbigi^{\ast}]$,
we set $S_{0}(\lambda):=S_0(\nbigi(\lambda))$.
\index{set $S_0(\lambda)$}
Set $\omega:=n/p$.
There exists
 $\theta_{\lambda}\in\real$
such that 
$S_0(\lambda)=
 \bigl\{\theta_{\lambda}+\ell\pi/\omega
 \,\big|\,\ell\in\seisuu\bigr\}$.
Set $T(\lambda):=T(\nbigi(\lambda))$ for 
$\lambda\in[\nbigi^{\ast}]$,
which is identified with
the set of the connected components 
of $\real\setminus S_0(\lambda)$.
\index{set $T(\lambda)$}
Let $T(\lambda)_+\subset T(\lambda)$
be the set of $J\in T(\lambda)$
such that
$\gminia>_J0$ for any $\gminia\in\nbigi(\lambda)$.
Set $T(\lambda)_-:=T(\lambda)\setminus T(\lambda)_+$.
We set
$T(0):=T(\nbigi)$.
\index{sets $T(\lambda)_{\pm}$, $T(\lambda)_0$}
Let $T_2(\nbigi)$
denote the set of pairs 
$(J_1,J_2)$ in $T(\nbigi)$
satisfying
$J_1\cap J_2\neq\emptyset$
and $J_1\neq J_2$.
\index{set $T_2(\nbigi)$}

We set
$\nbiga(\nbigi):=
 \coprod_{\lambda\in[\nbigi]}
 \bigl\{
 (\lambda,J)\,\big|\,
 J\in T(\lambda)
 \bigr\}$.
\index{set $\nbiga(\nbigi)$}
We also set
\index{set $\nbiga(\nbigi)_{>0}$}
\index{set $\nbiga(\nbigi)_{<0}$}
\[
\nbiga_{>0}(\nbigi):=
 \coprod_{\lambda\in[\nbigi]}
 \bigl\{
 (\lambda,J)\,\big|\,
 J\in T(\lambda)_{+}
 \bigr\},
\quad
 \nbiga_{<0}(\nbigi):=
 \coprod_{\lambda\in[\nbigi]}
 \bigl\{
 (\lambda,J)\,\big|\,
 J\in T(\lambda)_{-}
 \bigr\}.
\]
We put
$\nbiga_0(\nbigi):=
 \nbiga(\nbigi)\setminus
 (\nbiga_{>0}(\nbigi)\cup\nbiga_{<0}(\nbigi))$.
It equals $\{(0,J)\,|\,J\in T(\nbigi)\}$ if $0\in\nbigi$.
\index{set $\nbiga_0(\nbigi)$}

We set $\nbigibar:=\nbigi\cup(-\nbigi)\cup\{0\}$.
We have $T(\nbigibar)=T(\nbigi)$
and $T_2(\nbigibar)=T_2(\nbigi)$.
There exist the natural embeddings
$\nbiga(\nbigi)\subset\nbiga(\nbigibar)$,
$\nbiga_{>0}(\nbigi)\subset\nbiga_{>0}(\nbigibar)$,
$\nbiga_{<0}(\nbigi)\subset\nbiga_{<0}(\nbigibar)$,
and 
$\nbiga_{0}(\nbigi)\subset\nbiga_{0}(\nbigibar)$.
For $J\in T(\nbigi)$,
the set $[\nbigibar_{J,>0}]$ consists of one element
$\lambda_+(J)$,
and the set $[\nbigibar_{J,<0}]$ consists of one element
 $\lambda_-(J)$.
\index{elements $\lambda_{\pm}(J)$}
We set \index{set $\nbigb_2(\nbigi)_J$}
\[
 \nbigb_2(\nbigi)_J:=
 \bigl\{
 (\lambda_+(J),0;J),
 (\lambda_+(J),\lambda_-(J);J),
 (0,\lambda_-(J);J)
 \bigr\}.
\]
We also set
$\nbigb_2(\nbigi):=
 \coprod_{J\in T(\nbigi)}
 \nbigb_2(\nbigi)_J$.
\index{set $\nbigb_2(\nbigi)$}

\section{Stokes graded local systems}

Let $\nbigitilde$ be a $\Gal(p)$-invariant subset
of $z_p^{-1}\cnum[z_p^{-1}]$
such that 
$\ord(\nbigitilde)=-\omega$.
Set $\nbigi:=\pi_{\omega}(\nbigitilde)$.
We set $\nbigibar=\nbigi\cup(-\nbigi)\cup\{0\}$.
For $\lambda\in[\nbigibar]$,
let $\nbigi(\lambda)$ denote the inverse image of $\lambda$
by $\nbigi\to[\nbigibar]$,
and let $\nbigitilde(\lambda)\subset\nbigitilde$
denote the inverse image of
$\nbigi(\lambda)$
by $\pi_{\omega}:\nbigitilde\to\nbigibar$.
\index{set $\nbigitilde(\lambda)$}
\index{set $\nbigi(\lambda)$}
Recall the notion of Stokes graded local systems
in this particular context.
(See \S\ref{subsection;18.4.3.2}.)
\begin{df}
 \index{$2\pi\seisuu$-equivariant Stokes graded local system
 over $(\nbigitilde,[\nbigi])$}  
A $2\pi\seisuu$-equivariant
Stokes graded local system
over $(\nbigitilde,[\nbigi])$ on $\real$
is a $2\pi\seisuu$-equivariant local system 
$\nbigk_{\bullet}
=\bigoplus_{\lambda\in[\nbigibar]}\nbigk_{\lambda}$
with a Stokes structure $\vecnbigf$
such that the following holds.
\begin{itemize}
\item $\nbigk_{\lambda}=0$ unless $\lambda\in[\nbigi]$.
 \item
The Stokes structure $\vecnbigf$ of $\nbigk_{\bullet}$
is the direct sum of 
Stokes structures on $\nbigk_{\lambda}$
indexed by $\nbigitilde(\lambda)$.
\end{itemize}
 
Similarly,
a $2\pi\seisuu$-equivariant
Stokes graded local system
over $(\nbigitilde,\nbigi)$ on $\real$
is a $2\pi\seisuu$-equivariant local system 
$\nbigk_{\bullet}
=\bigoplus_{\gminia\in\nbigibar}\nbigk_{\gminia}$
with a Stokes structure $\vecnbigf$
 such that the following holds.
  \index{$2\pi\seisuu$-equivariant Stokes graded local system
 over $(\nbigitilde,\nbigi)$}  
\begin{itemize}
\item $\nbigk_{\gminia}=0$ unless $\gminia\in\nbigi$.
 \item
The Stokes structure $\vecnbigf$ of $\nbigk_{\bullet}$
is the direct sum of 
Stokes structures on $\nbigk_{\gminia}$
indexed by $\pi_{\omega}^{-1}(\gminia)\subset\nbigitilde$.
\hfill\qed
\end{itemize}
 \end{df}

Let $\Loc^{\St}(\nbigitilde,[\nbigi])$
(resp. $\Loc^{\St}(\nbigitilde,\nbigi)$)
denote the category of 
$2\pi\seisuu$-equivariant
Stokes graded local systems over $(\nbigitilde,[\nbigi])$
(resp. $(\nbigitilde,\nbigi)$).
\index{category $\Loc^{\St}(\nbigitilde,[\nbigi])$}
\index{category $\Loc^{\St}(\nbigitilde,\nbigi)$}

\subsection{Expression as tuples of filtered vector spaces
and linear maps}
\label{subsection;20.10.9.2}

It is convenient for us
to consider another expression of
objects in $\Loc^{\St}(\nbigitilde,[\nbigi])$.

\begin{df}
 \index{$2\pi\seisuu$-equivariant Stokes tuple of vector spaces
 over $(\nbigitilde,[\nbigi])$}
A $2\pi\seisuu$-equivariant
Stokes tuple of vector spaces
over $(\nbigitilde,[\nbigi])$
is a tuple 
$(\vecK,\vecnbigf,\vecPhi,\vecPsi)$
as follows.
\begin{itemize}
\item
$\vecK=\Bigl(
 K_{\lambda,J}\,\Big|\,
 (\lambda,J)\in\nbiga(\nbigibar)
 \Bigr)$
denotes a tuple of vector spaces.
We impose $K_{\lambda,J}=0$ unless $(\lambda,J)\in\nbiga(\nbigi)$.
Each $K_{\lambda,J}$
is equipped with a Stokes structure
$\vecnbigf(K_{\lambda,J}):=
 \bigl\{ \nbigf^{\theta}(K_{\lambda,J})\,\big|\,
 \theta\in \Jbar\bigr\}$
indexed by $\nbigitilde(\lambda)$.
\item
$\vecPhi$ denotes a tuple of isomorphisms:
\[
 \Phi^{J+\pi/\omega,J}_{\lambda}:
 (K_{\lambda,J},\nbigf^{\vartheta^J_{r}})
\simeq
 (K_{\lambda,J+\pi/\omega},\nbigf^{\vartheta^J_r})
\quad\quad
 \bigl(
 (\lambda,J)\in\nbiga(\nbigibar)
 \bigr)
\]
\[
 \Phi_0^{J_2,J_1}:
 K_{0,J_1}
\simeq
 K_{0,J_2}
\quad
 ((0,J_i)\in\nbiga(\nbigibar),\,J_1\vdash J_2). 
\]
We assume that 
$\Phi_0^{J_2,J_1}$ preserves 
the filtrations $\nbigf^{\theta}$
for $\theta\in \overline{J_1}\cap\overline{J_2}$.
\item
$\vecPsi$ denotes 
a tuple of isomorphisms
$\Psi_{\lambda,J}:
 (K_{\lambda,J},\vecnbigf)
 \simeq 
 (K_{\Tbb^{\ast}(\lambda),\Tbb^{-1}(J)},
 \vecnbigf)$
for $(\lambda,J)\in \nbiga(\nbigibar)$,
where we use the bijection of the index sets
$\Tbb^{\ast}:\nbigi(\lambda)\simeq
 \nbigi(\Tbb^{\ast}(\lambda))$.
We impose the following compatibility condition:
\[
 \Phi^{\Tbb^{-1}(J+\pi/\omega),\Tbb^{-1}(J)}_{\Tbb^{\ast}(\lambda)}\circ
 \Psi_{\lambda,J}
=\Psi_{\lambda,J+\pi/\omega}\circ
  \Phi^{J+\pi/\omega,J}_{\lambda},
\]
\[
 \Phi_0^{\Tbb^{-1}(J_2),\Tbb^{-1}(J_1)}\circ\Psi_{0,J_1}
=\Psi_{0,J_2}\circ \Phi_0^{J_2,J_1}.
\]
The maps $\Psi_{\lambda,J}$ will often be denoted by
$\Psi_J$.
\hfill\qed
\end{itemize}
\end{df}

Let $(\nbigk_{\bullet},\vecnbigf)
 \in
 \Loc^{\St}(\nbigitilde,[\nbigi])$.
For each $(\lambda,J)\in \nbiga(\nbigibar)$,
we set
$K_{\lambda,J}:=
H^0(\Jbar,\nbigk_{\lambda})$.
We have $K_{\lambda,J}=0$ unless
$(\lambda,J)\in\nbiga(\nbigi)$.
The vector space $K_{\lambda,J}$
is equipped with a family of filtrations
$\nbigf^{\theta}$ $(\theta\in \Jbar)$
indexed by $(\nbigitilde(\lambda),\leq_{\theta})$.
For $\lambda\in[\nbigibar^{\ast}]$,
we obtain an isomorphism
$\Phi_{\lambda}^{J+\pi/\omega,J}$
as
$K_{\lambda,J}
\simeq
 \nbigk_{\lambda|\vartheta^J_{r}}
\simeq
 K_{\lambda,J+\pi/\omega}$.
We also obtain 
$\Phi_0^{J_2,J_1}$
as 
$K_{0,J_1}\simeq
 H^0(\Jbar_1\cap\Jbar_2,\nbigk_0)
\simeq
 K_{0,J_2}$.
Thus, we obtain a tuple of isomorphisms
$\vecPhi$.
By the $2\pi\seisuu$-action on $\nbigk_{\bullet}$,
we naturally obtain a tuple $\vecPsi$
of isomorphisms
$\Psi_{\lambda,J}:
 (K_{\lambda,J},\vecnbigf)\simeq
 (K_{\Tbb^{\ast}(\lambda),\Tbb^{-1}(J)},\vecnbigf)$.
It is easy to see that
$\gbigd(\nbigk_{\bullet},\vecnbigf):=
 (\vecK,\vecnbigf,\vecPhi,\vecPsi)$
is a $2\pi\seisuu$-equivariant
Stokes tuple of vector spaces over $(\nbigitilde,[\nbigi])$.
\index{tuple $\gbigd(\nbigk_{\bullet},\vecnbigf)$}
Conversely,
let $(\vecK,\vecnbigf,\vecPhi,\vecPsi)$ be
a $2\pi\seisuu$-equivariant
Stokes tuple of vector spaces over $(\nbigitilde,[\nbigi])$.
For each $\lambda\in[\nbigibar]$,
we obtain a local system
with Stokes structure $(\nbigk_{\lambda},\vecnbigf)$
on $\real$ indexed by $\nbigitilde(\lambda)$
by gluing $(K_{\lambda,J},\vecnbigf)$ $(J\in T(\lambda))$
via $\vecPhi$.
We have $\nbigk_{\lambda}=0$ unless $\lambda\in[\nbigi]$.
The direct sum
$\gbigl(\vecK,\vecnbigf,\vecPhi,\vecPsi)
=\bigoplus_{\lambda\in[\nbigibar]}(\nbigk_{\lambda},\vecnbigf)$
is naturally
equipped with a $2\pi\seisuu$-action,
and it induces
an object in 
$\Loc^{\St}(\nbigitilde,[\nbigi])$.
\index{graded local system with Stokes structure
$\gbigl(\vecK,\vecnbigf,\vecPhi,\vecPsi)$}
The following is clear by the constructions.
\begin{lem}
We naturally have
$\gbigd\circ\gbigl(\vecK,\vecnbigf,\vecPhi,\vecPsi)
=(\vecK,\vecnbigf,\vecPhi,\vecPsi)$
and 
$\gbigl\circ\gbigd(\nbigl_{\bullet},\vecnbigf)
\simeq
 (\nbigl_{\bullet},\vecnbigf)$.
\hfill\qed
\end{lem}

To simplify the notation,
we set
$K_{<0,J}:=K_{\lambda_-(J),J}$,
and
$K_{>0,J}:=K_{\lambda_+(J),J}$.

Let $(\nbigk_{\bullet},\vecnbigf)
 \in
 \Loc^{\St}(\nbigitilde,[\nbigi])$.
Set $(\vecK,\vecnbigf,\vecPhi,\vecPsi):=
 \gbigd(\nbigk_{\bullet},\vecnbigf)$.
For any $\lambda\in [\nbigibar]$
and for any $J_1,J_2\in T(\lambda)$,
we obtain the isomorphism
$\Phi^{J_2,J_1}_{\lambda}:
 K_{\lambda,J_1}\simeq
 K_{\lambda,J_2}$
induced by the parallel transport of
$\nbigk_{\lambda}$,
which is naturally obtained from $\vecPhi$.

\section{Stokes shells}

\subsection{Deformation data}

Let $(\nbigk_{\bullet},\vecnbigf)$ be 
a $2\pi\seisuu$-equivariant
Stokes graded local system
over $(\nbigitilde,[\nbigi])$.
We set
$(\vecK,\vecnbigf,\vecPhi,\vecPsi):=
 \gbigd(\nbigk_{\bullet},\vecnbigf)$.

\begin{df}
\index{deformation datum}
A deformation datum of $(\nbigk_{\bullet},\vecnbigf)$
is a tuple of morphisms  $\vecnbigr$:
\[
 \nbigr^{J_1}_{J_2}:
 K_{>0,J_1}
\lrarr
 K_{<0,J_2}
\quad
 \bigl(
 (J_1,J_2)\in T_2(\nbigi)
 \bigr),
\]
\[
 \nbigr^{\lambda_1,J_-}_{\lambda_2,J_+}:
 K_{\lambda_1,J}
\lrarr
 K_{\lambda_2,J}
\quad
((\lambda_1,\lambda_2,J)\in\nbigb_2(\nbigibar)).
\]
They are assumed to be equivariant with respect to $\vecPsi$
in the following sense:
\[
 \Psi_{\lambda_-(J_2),J_2}\circ
 \nbigr^{J_1}_{J_2}
=\nbigr^{\Tbb^{-1}(J_1)}_{\Tbb^{-1}(J_2)}\circ
 \Psi_{\lambda_+(J_1),J_1},
\,\,\,
 \Psi_{\lambda_2,J}\circ
 \nbigr^{\lambda_1,J_-}_{\lambda_2,J_+}
=\nbigr^{\Tbb^{\ast}(\lambda_1),\Tbb^{-1}(J)_-}
     _{\Tbb^{\ast}(\lambda_2),\Tbb^{-1}(J)_+}
\circ\Psi_{\lambda_1,J}.
\]
\hfill\qed
\end{df}

For a given deformation datum $\vecnbigr$
of $(\nbigk_{\bullet},\vecnbigf)$,
we also set
\[
\nbigr^{0,J_+}_{\lambda_-(J),J_-}:=
 -\nbigr^{0,J_-}_{\lambda_-(J),J_+},
\quad\quad
\nbigr^{\lambda_+(J),J_+}_{0,J_-}:=
-\nbigr^{\lambda_+(J),J_-}_{0,J_+},
\]
\[
 \nbigr^{\lambda_+(J),J_+}_{\lambda_-(J),J_-}:=
 -\nbigr^{\lamda_+(J),J_-}_{\lambda_-(J),J_+}
+\nbigr^{0,J_-}_{\lambda_-(J),J_+}
 \circ
 \nbigr^{\lambda_+(J),J_-}_{0,J_+}.
\]

We shall also use the following notation:
\index{maps $\nbigp_J$, $\nbigq_J$, $\nbigr^{J_-}_{J_+}$, $\nbigr^{J_+}_{J_-}$}
\[
\nbigp_J:=\nbigr^{0,J_-}_{\lambda_-(J),J_+},
\quad
\nbigq_J:=\nbigr^{\lambda_+(J),J_-}_{0,J_+},
\quad
\nbigr^{J_-}_{J_+}:=
 \nbigr^{\lambda_+(J),J_-}_{\lambda_-(J),J_+},
\quad
 \nbigr^{J_+}_{J_-}:=
 \nbigr^{\lambda_+(J),J_+}_{\lambda_-(J),J_-}.
\]

\subsection{Stokes shells}

We define the notion of Stokes shells as follows.

\begin{df}
\index{Stokes shell}
A Stokes shell 
$\Sh=(\nbigk_{\bullet},\vecnbigf,\vecnbigr)$
indexed by $\nbigitilde$
is a $2\pi\seisuu$-equivariant 
Stokes graded local system 
$(\nbigk_{\bullet},\vecnbigf)$
over $(\nbigitilde,[\nbigi])$
equipped with 
a deformation datum $\vecnbigr$.
\hfill\qed
\end{df}

For a shell $\Sh=(\nbigk_{\bullet},\vecnbigf,\vecnbigr)$,
we set 
$\gbigd(\Sh):=\gbigd(\nbigk_{\bullet},\vecnbigf)$.
\index{tuple $\gbigd(\Sh)$}

Let 
$\Sh_i=(\nbigk_{i,\bullet},\vecnbigf_i,\vecnbigr_i)$
be Stokes shells indexed by $\nbigitilde$.
A morphism 
$\Sh_1\lrarr\Sh_2$
is defined to be 
a morphism
$F:(\nbigk_{1,\bullet},\vecnbigf)
\lrarr
 (\nbigk_{2,\bullet},\vecnbigf)$
in $\Loc^{\St}(\nbigitilde,[\nbigi])$
such that
 $F$ is compatible with the deformation data
 in the following sense:
\[
 F\circ
 (\nbigr_1)^{J_1}_{J_2}
=(\nbigr_2)^{J_1}_{J_2}\circ
 F
\quad
 ((J_1,J_2)\in T_2(\nbigi)),
\]
\[
F\circ
 (\nbigr_1)^{\lambda_1,J_-}_{\lambda_2,J_+}
=(\nbigr_2)^{\lambda_1,J_-}_{\lambda_2,J_+}
\circ F
\quad
((\lambda_1,\lambda_2;J)\in \nbigb_2(\nbigibar)).
\]
\begin{notation}
Let $\Shcat(\nbigitilde)$ denote 
the category of Stoke shells indexed by $\nbigitilde$.
\index{category $\Shcat(\nbigitilde)$} 
\hfill\qed
\end{notation}

\begin{rem}
Let $(\nbigk_{\bullet},\vecnbigf,\vecnbigr)$
be a Stokes shell indexed by $\nbigitilde$.
We take a $\Gal(p)$-invariant subset
$\nbigitilde'\subset  z_p^{-1}\cnum[z_p^{-1}]$
such that
$\ord\nbigitilde'=-\omega$
and 
$\nbigitilde\subset\nbigitilde'$.
We naturally have
$[\nbigitilde]\subset[\nbigitilde']$.
By setting
$\nbigk'_{\lambda}:=\nbigk_{\lambda}$
for $\lambda\in[\nbigi]$
and 
$\nbigk'_{\lambda}:=0$
for $\lambda\in[\nbigi']\setminus[\nbigi]$,
we naturally obtain
a Stokes shell
$(\nbigk'_{\bullet},\vecnbigf,\vecnbigr)$
indexed by $\nbigitilde'$.
This procedure induces 
a fully faithful functor
$\Shcat(\nbigitilde)\lrarr
 \Shcat(\nbigitilde')$.
Therefore, we can freely enlarge
the index set $\nbigitilde$.
\hfill\qed
\end{rem}

\subsection{Induced maps}
\label{subsection;20.10.1.10}

Let $\Sh=(\nbigl_{\bullet},\vecnbigf,\vecnbigr)
\in\Shcat(\nbigitilde)$.
We set
$(\vecK,\vecnbigf,\vecPhi,\vecPsi)=\gbigd(\Sh)$.
For any  $J\in T(\nbigi)$,
we obtain the following automorphism
of $K_{0,J}\oplus K_{<0,J}\oplus K_{>0,J}$:
\begin{equation}
\label{eq;18.5.1.1}
 \Pi^{J_{+},J_{-}}:=
 \id+\sum_{(\lambda_1,\lambda_2,J)\in\nbigb_2(\nbigibar_J)}
 \nbigr^{\lambda_1,J_-}_{\lambda_2,J_+}.
\end{equation}
Let $J_1,J_2\in T(\nbigi)$ such that 
$\vartheta^{J_1}_{\ell}<\vartheta^{J_2}_{\ell}<\vartheta^{J_1}_r$.
The following map is contained in
$\vecnbigr$:
\begin{equation}
 \label{eq;18.4.30.1}
 \nbigr^{J_1}_{J_2}:
 K_{>0,J_1}
\lrarr
 K_{<0,J_2}.
\end{equation}
The following map is induced by
$\vecPhi$ and 
$-\nbigr^{J_1+\pi/\omega}_{J_2}$:
\begin{equation}
\label{eq;18.4.30.2}
 K_{<0,J_1}
\simeq
 K_{>0,J_1+\pi/\omega}
\lrarr
 K_{<0,J_2}.
\end{equation}
We obtain the following map
from (\ref{eq;18.4.30.1}) and (\ref{eq;18.4.30.2}):
\begin{equation}
 \Upsilontilde^{J_1}_{J_2}:
 \bigoplus_{\lambda\in[\nbigi^{\ast}_{J_1}]}
 K_{\lambda,J_1}
\lrarr
 K_{<0,J_2}.
\end{equation}

Similarly, the following map is contained in $\vecnbigr$:
\begin{equation}
 \label{eq;18.4.30.3}
\nbigr^{J_2}_{J_1}:
 K_{>0,J_2}
\lrarr
 K_{<0,J_1}.
\end{equation}
The following map is induced by
$\vecPhi$ and 
$-\nbigr^{J_2-\pi/\omega}_{J_1}$:
\begin{equation}
\label{eq;18.4.30.4}
 K_{<0,J_2}
\simeq
 K_{>0,J_2-\pi/\omega}
\lrarr
 K_{<0,J_1}.
\end{equation}
We obtain the following map
from (\ref{eq;18.4.30.3}) and (\ref{eq;18.4.30.4}):
\begin{equation}
 \Upsilontilde^{J_2}_{J_1}:
 \bigoplus_{\lambda\in[\nbigi^{\ast}_{J_2}]}
 K_{\lambda,J_2}
\lrarr
 K_{<0,J_1}.
\end{equation}

\section{The associated local systems with Stokes structure}

\subsection{Construction (1)}
\label{subsection;18.4.30.50}

Let $\nbigitilde$ be a $\Gal(p)$-invariant subset
of $z_p^{-1}\cnum[z_p^{-1}]$.
Set $\nbigi:=\pi_{\omega}(\nbigitilde)$.
Let $\Sh=(\nbigk_{\bullet},\vecnbigf,\vecnbigr)\in\Shcat(\nbigitilde)$.
Set 
$(\vecK,\vecnbigf,\vecPhi,\vecPsi):=\gbigd(\Sh)$.
We shall construct
a $2\pi\seisuu$-equivariant local system with Stokes structure
$\Locst(\Sh)=(\nbigh^{\Sh},\vecnbigf^{\Sh})$
indexed by $\nbigitilde$.
\index{$2\pi\seisuu$-equivariant local system with Stokes structure
$\Locst(\Sh)$}

Let $I=\openopen{\theta_0}{\theta_1}$
be any connected component of
$\real\setminus S_0(\nbigi)$.
We set $J_1:=\openopen{\theta_1-\pi/\omega}{\theta_1}\in T(\nbigi)$.
Let $\nbigh^{\Sh}_{I_{\pm}}$
denote the local system on $I_{\pm}$
induced by the vector space
\begin{equation}
\label{eq;18.4.30.20}
 K_{0,J_{1}}
\oplus
 \bigoplus_{J\in T(\nbigi)_I}
 \bigoplus_{\lambda\in [\nbigi^{\ast}_J]}
 K_{\lambda,J}.
\end{equation}
For each $\theta\in I_{\pm}$,
the stalk $\nbigh^{\Sh}_{I_{\pm}|\theta}$ at $\theta$
is identified with the graded vector space
(\ref{eq;18.4.30.20}).
Let $\nbigf^{\Sh,\theta}$
denote the filtration of
$\nbigh^{\Sh}_{I_{\pm}|\theta}$ indexed by
$(\nbigitilde,\leq_{\theta})$
obtained as the direct sum of the filtrations
$\nbigf^{\theta}$ on $K_{\lambda,J}$.
We have the automorphism
$\gbigp_{I}:=
 \id\oplus\Pi^{J_{1+},J_{1-}}$
of the following vector space
(see (\ref{eq;18.5.1.1}) for $\Pi^{J_{1+},J_{1-}}$):
\begin{equation}
\label{eq;18.4.30.21}
\Bigl(
 \bigoplus_{J\in T(\nbigi)_I\setminus\{J_1\}}
 \bigoplus_{\lambda\in [\nbigi^{\ast}_{J}]}
 K_{\lambda,J}
\Bigr)
\oplus
\bigoplus_{\lambda\in[\nbigi_{J_1}]} 
 K_{\lambda,J_1}
\end{equation}
It induces an isomorphism
of the local systems
$\nbigh^{\Sh}_{I_-|I}\simeq
 \nbigh^{\Sh}_{I_+|I}$.
It preserves the filtrations
$\nbigf^{\Sh\,\theta}$ $(\theta\in I)$.
Hence, we obtain a local system
with Stokes structure
$(\nbigh^{\Sh}_{\Ibar},\vecnbigf^{\Sh})$
on $\Ibar$
by gluing 
$(\nbigh^{\Sh}_{I_{\pm}},\vecnbigf^{\Sh})$.

Let $I$ be as above,
and let $I':=\openopen{\theta_1}{\theta_2}$
be the connected component of
$\real\setminus S_0(\nbigi)$ next to $I$.
Let us construct an isomorphism
$F_{\theta_1}:
 \nbigh^{\Sh}_{\Ibar|\theta_1}
\lrarr
 \nbigh^{\Sh}_{\Ibar'|\theta_1}$
of the stalks at $\theta_1$.
Set $J_2:=\openopen{\theta_2-\pi/\omega}{\theta_2}$.
We have the identifications:
\[
 \nbigh^{\Sh}_{\Ibar|\theta_1}
=\nbigh^{\Sh}_{I_+|\theta_1}
=\bigoplus_{\lambda\in[\nbigi_{J_1}]}
 K_{\lambda,J_1}
\oplus
 \bigoplus_{J\in T(\nbigi)_I\setminus\{J_1\}}
 \bigoplus_{\lambda\in [\nbigi^{\ast}_{J}]}
 K_{\lambda,J},
\]
\[
  \nbigh^{\Sh}_{\Ibar'|\theta_1}
=\nbigh^{\Sh}_{I'_{-}|\theta_1}
=K_{0,J_2}
\oplus
 \bigoplus_{\lambda\in[\nbigi_{J_1}^{\ast}]}
 K_{\lambda,J_1+\pi/\omega}
\oplus
 \bigoplus_{J\in T(\nbigi)_{I}\setminus\{J_1\}}
 \bigoplus_{\lambda\in [\nbigi^{\ast}_{J}]}
 K_{\lambda,J}.
\]
We have the following map
induced by
$\vecPhi$
and $\Upsilontilde^{J_1}_{J}$
for $J\ni\theta_1$:
\begin{equation}
\label{eq;18.4.30.30}
 K_{0,J_{1}}
\oplus
 \bigoplus_{\lambda\in [\nbigi^{\ast}_{J_1}]}
 K_{\lambda,J_{1}}
\lrarr
 K_{0,J_2}
 \oplus
 \bigoplus_{\lambda\in [\nbigi^{\ast}_{J_1}]}
 K_{\lambda,J_1+\pi/\omega}
 \oplus
 \bigoplus_{J\in T(\nbigi)_I\setminus\{J_1\}}
 \bigoplus_{\lambda\in[\nbigi_{J,<0}]}K_{\lambda,J}.
\end{equation}
Then, 
$F_{\theta_1}$ is defined as the map
induced by the morphism
(\ref{eq;18.4.30.30})
and the identity on
$\bigoplus_{J\in T(\nbigi)_I\setminus\{J_1\}}
 \bigoplus_{\lambda\in [\nbigi^{\ast}_{J}]}
 K_{\lambda,J}$.
The following is clear by the construction.
\begin{lem}
$F_{\theta_1}$
preserves
the filtrations $\nbigf^{\Sh\,\theta_1}$.
\hfill\qed
\end{lem}

By gluing $(\nbigh^{\Sh}_{\Ibar},\vecnbigf^{\Sh})$
for connected components $I$ 
of $\real\setminus S_0(\nbigi)$,
we obtain a $2\pi\seisuu$-equivariant
local system with Stokes structure
$\Locst(\Sh)=(\nbigh^{\Sh},\vecnbigf^{\Sh})$ on $\real$
indexed by $\nbigitilde$.

\subsection{Construction (2)}
\label{subsection;18.4.30.51}

Let $(L,\vecnbigf)$ be a $2\pi\seisuu$-equivariant
local system with Stokes structure indexed by $\nbigitilde$
on $\real$.
We obtain the Stokes structure
$\pi_{\omega\ast}(\vecnbigf)$ indexed by $\nbigi$.
There exist the canonical splittings (\ref{eq;18.5.2.1})
for any interval $J$ with $\vartheta^J_{r}-\vartheta^J_{\ell}=\pi/\omega$.
Take $J\in T(\nbigi)$.
Moreover, we have
the local subsystems
$L_{J,<0}\subset L_{J,\leq 0}\subset
 \gbiga_J(L)\subset L_{|J}$
as in \S\ref{subsection;18.5.13.10}.
We set
$K_{\lambda_-(J),J}:=
 H^0(J,L_{J,<0})$,
$K_{0,J}:=
 H^0(J,L_{J,\leq 0})\big/
 H^0(J,L_{J,<0})$,
 and
$K_{\lambda_+(J),J}:=
 H^0(J,\gbiga_J(L))\big/H^0(J,L_{J,\leq 0})$.

There exists a natural isomorphism
$K_{0,J}\simeq
H^0\bigl(J,L_{J_-,0|J}\bigr)$.
Because 
$L_{J_-,0|J}\subset
 L_{J_+,0|J}\oplus
 L_{J_+,<0|J}$,
 we obtain the map
$\nbigr^{0,J_-}_{\lambda_-(J),J_+}:
  K_{0,J}\lrarr
  K_{\lambda_-(J),J}$
  as the composite of the following natural maps:
\begin{multline}
 \nbigr^{0,J_-}_{\lambda_-(J),J_+}:
 K_{0,J}
 \simeq
 H^0\bigl(J,L_{J_-,0|J}\bigr)
 \subset
 H^0\bigl(J,
 L_{J_+,0|J}\oplus
 L_{J_+,<0|J}\bigr)
 \lrarr
  H^0\bigl(J,L_{J,<0}\bigr) \\
=K_{\lambda_-(J),J}.
\end{multline}
There exists a natural isomorphism
$K_{\lambda_+(J),J}
\simeq
 H^0(J,L_{J_-,>0|J})$.
Because
\[
 L_{J_-,>0|J}
\subset
 \gbiga_{J_+}(L)_{|J}
=L_{J_+,>0|J}
\oplus
 L_{J_+,0|J}
\oplus
 L_{J_+,<0|J},
\]
we obtain the maps
$\nbigr^{\lambda_+(J),J_-}_{0,J_+}:
 K_{\lambda_+(J),J}
\lrarr
 K_{0,J}$
and
$\nbigr^{\lambda_+(J),J_-}_{\lambda_-(J),J_+}:
 K_{\lambda_+(J),J}
\lrarr
 K_{\lamda_-(J),J}$.

By the construction,
there exist the following natural isomorphisms:
\[
 K_{\lambda_-(J),J}\simeq
 L_{J_-,<0|\vartheta^J_{\ell}}
=\bigoplus_{\gminia\in\nbigi_{J,<0}}
 L_{J_{-},\gminia|\vartheta^J_{\ell}}
\simeq
 L_{J_+,<0|\vartheta^J_{r}}
= \bigoplus_{\gminia\in\nbigi_{J,<0}}
 L_{J_{+},\gminia|\vartheta^J_{r}},
\]
\[
 K_{0,J}\simeq
 L_{J_{-},0|\vartheta^J_{\ell}}
\simeq
 L_{J_{+},0|\vartheta^J_{r}},
\]
\[
 K_{\lambda_+(J),J}\simeq
 L_{J_-,>0|\vartheta^J_{\ell}}
=\bigoplus_{\gminia\in\nbigi_{J,>0}}
 L_{J_{-},\gminia|\vartheta^J_{\ell}}
\simeq
 L_{J_+,>0|\vartheta^J_{r}}
=\bigoplus_{\gminia\in\nbigi_{J,>0}}
 L_{J_{+},\gminia|\vartheta^J_{r}}.
\]

For each $\gminia\in \nbigi_{J-\pi/\omega}^{\ast}$,
we have the following:
\begin{equation}
\label{eq;18.5.1.2}
 L_{(J-\pi/\omega)_+,\gminia|\vartheta^J_{\ell}}
\subset
 L_{J_-,\gminia|\vartheta^J_{\ell}}
\oplus
 \bigoplus_{
 \substack{J'\in T(\nbigi)\\
 \vartheta^J_{\ell}\in J' }}
 L_{J'_-,<0|\vartheta^J_{\ell}}.
\end{equation}
For $\lambda\in[\nbigi_J^{\ast}]$,
we obtain the isomorphisms
$\Phi_{\lambda}^{J,J-\pi/\omega}:
 K_{\lambda,J-\pi/\omega}
\simeq
K_{\lambda,J}$
from the isomorphisms
 $L_{(J-\pi/\omega)_+,\gminia|\vartheta^J_{\ell}}
\simeq
  L_{J_-,\gminia|\vartheta^J_{\ell}}$
induced by $(\ref{eq;18.5.1.2})$.
For $J'\in T(\nbigi)$ such that 
$\vartheta^J_{\ell}=\vartheta^{J-\pi/\omega}_r\in J'$,
we obtain the morphism
\[
 \Upsilontilde^{J-\pi/\omega}_{J'}:
 \bigoplus_{\lambda\in [\nbigi^{\ast}_{J-\pi/\omega}]}
 K_{\lambda,J-\pi/\omega}
\lrarr
 K_{<0,J'}
\]
from the morphisms
$L_{(J-\pi/\omega)_+,\gminia|\vartheta^J_{\ell}}
\lrarr
  \bigoplus_{\gminib\in\nbigi_{J',<0}}
 L_{J',\gminib|\vartheta^J_{\ell}}$
induced by (\ref{eq;18.5.1.2}).
Let $\nbigr^{J-\pi/\omega}_{J'}$
be the restriction of 
$\Upsilontilde^{J-\pi/\omega}_{J'}$
to $K_{>0,J-\pi/\omega}$.
Let $\nbigr^{J}_{J'}$ be the composite of
the following maps:
\[
 K_{>0,J}
\stackrel{a_1}{\simeq}
 K_{<0,J-\pi/\omega}
\stackrel{a_2}{\lrarr}
 K_{<0,J'}.
\]
Here,
$a_1$ is induced by
$\Phi^{J,J-\pi/\omega}_{\lambda_-(J-\pi/\omega)}$,
and $a_2$ is the restriction of
$-\Upsilontilde^{J-\pi/\omega}_{J'}$
to $K_{<0,J-\pi/\omega}$.

Let $J_1,J_2\in T(\nbigi)$ such that 
$J_1\vdash J_2$.
On $J_{1+}\cap J_{2-}$,
$L_{J_{1+},0|J_{1+}\cap J_{2-}}
=L_{J_{2-},0|J_{1+}\cap J_{2-}}$
holds.
Hence, we obtain an isomorphism
$\Phi^{J_2,J_1}_0:
 K_{0,J_1}\simeq K_{0,J_2}$.

Let $\nbigk_{\lambda,J}$
denote the local system on $\Jbar$
induced by $K_{\lambda,J}$.
It is naturally equipped with
a Stokes structure $\vecnbigf$
indexed by $\nbigitilde(\lambda)$.
By gluing $(K_{\lambda,J},\vecnbigf)$ $(J\in T(\lambda))$
by the tuple of isomorphisms $\vecPhi$,
we obtain  local systems
$(\nbigk_{\lambda},\vecnbigf)$
with Stokes structure indexed by $\nbigitilde(\lambda)$.
The direct sum
$(\nbigk_{\bullet},\vecnbigf)
=\bigoplus(\nbigk_{\lambda},\vecnbigf)$
is a $2\pi\seisuu$-equivariant
Stokes graded local system over
$(\nbigitilde,[\nbigi])$.
Together with the deformation datum $\vecnbigr$,
we obtain a Stokes shell
$\Shsf(L,\vecnbigf)$
indexed by $\nbigitilde$.
\index{Stokes shell $\Shsf(L,\vecnbigf)$}

\begin{lem}
$\Locst\Shsf(L,\vecnbigf)$
is naturally isomorphic to $(L,\vecnbigf)$.
\end{lem}
\pf
By the constructions of
$\Locst$ and $\Shsf$,
there exists a natural isomorphism
$\Locst\Shsf(L,\vecnbigf)
\lrarr
 (L,\vecnbigf)$
of $2\pi\seisuu$-equivariant
Stokes graded local systems
over $(\nbigitilde,[\nbigi])$.
\hfill\qed

\subsection{Equivalence}

Let 
$\Loc^{\St}(\nbigitilde)$ 
denote the category
of $2\pi\seisuu$-equivariant local systems
with Stokes structure indexed by $\nbigitilde$
on $\real$.
The construction in \S\ref{subsection;18.4.30.50}
induces a functor
$\Locst:\Shcat(\nbigitilde)
\lrarr
\Loc^{\St}(\nbigitilde)$.
The construction in \S\ref{subsection;18.4.30.51}
induces a functor
$\Shsf:
 \Loc^{\St}(\nbigitilde)
\lrarr
 \Shcat(\nbigitilde)$.

\begin{prop}
\label{prop;18.4.30.42}
$\Locst$ is an equivalence,
and $\Shsf$ is a quasi-inverse.
\end{prop}
\pf
Let $\Sh\in\Shcat(\nbigitilde)$.
For any
$J\in T(\nbigi)$,
we have the unique decompositions
\[
 \nbigh^{\Sh}_{|J_{\pm}}=
 \bigoplus_{\gminia\in\nbigi}
 \nbigh^{\Sh}_{J_{\pm},\gminia}
\]
which are compatible with the filtrations
$\pi_{\omega\ast}\nbigf^{\Sh,\theta}$ $(\theta\in J_{\pm})$.
Let $I$ be the connected component of $\real\setminus S_0(\nbigi)$
such that $\vartheta^I_{r}=\vartheta^J_r$.
By the construction of $\nbigh^{\Sh}$,
we have the natural isomorphism
$\nbigh^{\Sh}_{|I}\simeq\nbigh^{\Sh}_{I_-|I}$,
which gives the following natural identification
for any $\theta\in I_-$:
\[
 \nbigh^{\Sh}_{|\theta}
=\bigoplus_{\lambda\in [\nbigi_J]}K_{\lambda,J}
\oplus
 \bigoplus_{\substack{J'\in T(\nbigi)\\ 
 \vartheta^J_r\in J'}}
 \bigoplus_{\lambda\in[\nbigi^{\ast}_{J'}]}
 K_{\lambda,J'}.
\]
We have the decomposition
$K_{\lambda,J}
=\bigoplus_{\gminia\in\nbigi(\lambda)}
 K_{\gminia,J_-}$
which is a splitting of
$\pi_{\omega\ast}\nbigf^{\Sh\theta}$
$(\theta\in J_-)$
on $K_{\lambda,J}$.
For each $\gminia\in \nbigi_J$,
we obtain the local subsystem
$\nbigh_{\gminia,J_{-}}
\subset
 \nbigh^{\Sh}_{|J_{-}}$
determined by the condition
$\nbigh_{\gminia,J_{-}|\theta}
 =K_{\gminia,J_{-}}$
for $\theta\in I$.

\begin{lem}
 \label{lem;18.4.30.40}
For any $\gminia\in \nbigi_J$,
$\nbigh_{\gminia,J_{-}}
=\nbigh^{\Sh}_{J_{-},\gminia}$
holds.
Namely,
for any $\gminia\in \nbigi_J$
and for any 
$\theta\in J_{-}$,
$\nbigh_{\gminia,J_{-}|\theta}
\subset
\pi_{\omega\ast}\nbigf^{\theta}_{\gminia}$
holds.
\end{lem}
\pf
By the construction of $\nbigh^{\Sh}$
and the filtrations $\nbigf^{\Sh\,\theta}$,
$\nbigh_{\gminia,J_-|\theta}
\subset
\pi_{\omega\ast}\nbigf^{\Sh\,\theta}_{\gminia}$
clearly holds for any $\gminia\in\nbigi_{J}^{\ast}$.
Let us prove the claim for  $\nbigh_{0,J_{-}}$.
Let $\vartheta^J_{\ell}<\varphi_1<\varphi_2<\cdots<\varphi_{N-1}<
 \varphi_N:=\vartheta^J_r$
be the points of
$S_0(\nbigi)\cap\Jbar$.
We set 
$J_i:=\openopen{\varphi_i-\pi/\omega}{\varphi_i}$.
We have the local subsystem
$\nbigh_{0,J_{i,-}}\subset
 \nbigh^{\Sh}_{|J_{i,-}}$
determined by $K_{0,J_{i,-}}$.
We shall prove the claim
$\nbigh_{0,J_{i,-}|\theta}
 \subset
 \pi_{\omega\ast}\nbigf_0^{\Sh\,\theta}$
for any $\theta\in J_{-}\cap J_{i,-}$
by an induction of $i$.
If $i=1$, the claim is clear by the construction of
$\pi_{\omega\ast}\nbigf^{\Sh\,\theta}$.
Let us prove the claim for $i$ by assuming $i-1$.
Clearly,
$\nbigh_{0,J_{i,-}|\theta}
 \subset
 \pi_{\omega\ast}\nbigf^{\Sh\,\theta}$
holds
for $\theta\in\closedopen{\varphi_{i-1}}{\varphi_i}$.
For $\theta\in\closedopen{\varphi_{i-2}}{\varphi_{i-1}}$,
the construction implies
\[
 \nbigh_{0,J_{i,-}|\theta}
\subset
 \nbigh_{0,J_{i-1,-}|\theta}
\oplus
 \bigoplus_{\gminib\in \nbigi_{J_{i-1},<0}}
 \nbigh_{\gminib,J_{i-1,-}|\theta}.
\]
By the assumption of the induction,
we obtain the following for any 
$\theta\in J_{i-1,-}\cap J_-$:
\[
  \nbigh_{0,J_{i-1,-}|\theta}
\oplus
 \bigoplus_{\gminib\in \nbigi_{J_{i-1},<0}}
 \nbigh_{\gminib,J_{i-1,-}|\theta}
\subset
 \pi_{\omega\ast}\nbigf^{\Sh\,\theta}_0.
\]
Hence, we obtain the claim for $i$.
\hfill\qed

\vspace{.1in}

There exists the decomposition
$K_{\lambda,J}=
 \bigoplus_{\gminia\in\nbigi(\lambda)}
 K_{\gminia,J_+}$
which is a splitting of
$\pi_{\omega\ast}\nbigf^{\Sh\,\theta}$
$(\theta\in J_+)$.
By the isomorphism
$\nbigh^{\Sh}_{|J_+}
\simeq
 \nbigh^{\Sh}_{J_+}$,
we obtain the inclusion
$\iota_{\gminia,\vartheta^J_{r}}:
 K_{\gminia,J_+}
\lrarr
 \pi_{\omega\ast}
 \nbigf^{\Sh\,\vartheta^J_r}_{\gminia}$
for any $\gminia\in \nbigi_J$.
We obtain the local subsystem
$\nbigh_{\gminia,J_+}
\subset
 \nbigh^{\Sh}_{|J_+}$
for any $\gminia\in \nbigi_J$
determined by the condition
$\nbigh_{\gminia,J_+|\vartheta^J_r}
=\iota_{\gminia,\vartheta^J_r}(K_{\gminia,J_+})$.

\begin{lem}
\label{lem;18.4.30.41}
We have
$\nbigh_{\gminia,J_+}
=\nbigh^{\Sh}_{J_+,\gminia}$
for any $\gminia\in\nbigi_J$.
Namely,
we have
$\nbigh_{\gminia,J_+|\theta}
\subset\nbigf_{\gminia}^{\Sh\,\theta}$
for any $\theta\in J_+$
and for any $\gminia\in\nbigi_J$.
\end{lem}
\pf
The claim is clear if $\theta=\vartheta^J_r$.
Because $\gbigp_I$ 
preserves the filtrations,
we obtain the following for any 
$\gminia\in\nbigi_J$:
\[
 \nbigh_{\gminia,J_+|J}
\subset
 \bigoplus_{\substack{
 \gminib\in\nbigi_J\\
 \gminib\leq_J\gminia}}
 \nbigh_{\gminib,J_-|J}.
\]
Then, the claim of the lemma follows from
Lemma \ref{lem;18.4.30.40}.
\hfill\qed

\vspace{.1in}

By Lemma \ref{lem;18.4.30.40}
and Lemma \ref{lem;18.4.30.41},
$\Shsf(\Locst(\Sh))$
is naturally isomorphic to $\Sh$.
Thus, we obtain the claim of 
Proposition \ref{prop;18.4.30.42}.
\hfill\qed

\begin{df}
\label{df;21.5.3.2}
A functor
$\nbige:\Dsf(\nbigj)\lrarr \Shcat(\nbigitilde)$
is called a base tuple with respect to $\nbigj\subset\nbigitilde$
if the induced functor
$\Dsf(\nbigj)\lrarr \Loc^{\St}(\nbigitilde)$
is a base tuple
 in the sense of {\rm Definition \ref{df;21.5.3.3}.}
Similarly,
a morphism $F:\Sh_1\lrarr\Sh_2$ in $\Shcat(\nbigitilde)$
is called a base tuple
if the induced morphisms
 $\Gr^{\vecnbigf}_{\gminia}(F):
 \Gr^{\vecnbigf}_{\gminia}(\Sh_1)
 \lrarr
 \Gr^{\vecnbigf}_{\gminia}(\Sh_2)$
 are isomorphisms unless $\gminia=0$
 (see Definition {\rm\ref{df;21.5.3.1}}).
 \index{base tuple}
\hfill\qed
\end{df}

\section{Hills}
\label{subsection;20.10.9.20}

Let $\Sh$ be a Stokes shell.
Let
$\Locst(\Sh)=(\nbigh^{\Sh},\vecnbigf^{\Sh})$
denote the associated local system
with Stokes structure.
Let $H^0(\real,\nbigh^{\Sh})$ denote the space of
global sections of $\nbigh^{\Sh}$ on $\real$.
Set 
$(\vecK,\vecnbigf,\vecPhi,\vecPsi):=\gbigd(\Sh)$.

Recall that
$K_{<0,J}$ is identified with
the space of sections of
$\nbigh^{\Sh}_{<0,J}$,
and that
$K_{>0,J}$ is identified with
the space of the sections of
$\nbigh^{\Sh}_{>0,J}$.
(See \S\ref{subsection;18.5.13.10}
for $\nbigh^{\Sh}_{<0,J}$
and $\nbigh^{\Sh}_{>0,J}$.)
We may regard
$K_{<0,J}$ as a subspace of
$H^0(\real,\nbigh^{\Sh})$.
As we explained in \S\ref{subsection;21.6.7.4},
by the duality,
we may regard
$K_{>0,J}$
as the quotient of
$H^0(\real,\nbigh^{\Sh})$.
We shall describe the quotient map
$R_J:
H^0(\real,\nbigh^{\Sh})
\lrarr K_{>0,J}$
in (\ref{eq;21.6.7.3})
more concretely.
\index{map $R_J$}

Let $v\in H^0(\real,\nbigh^{\Sh})$.
Let $J\in T(\nbigi)$.
For any connected component $I$
of $J\setminus S_0(\gminia)$,
there exists a non-unique decomposition
\begin{equation}
\label{eq;18.4.5.50}
 v_{|I}=\sum_{J'\in T(\nbigi)_I}u_{I,J'_{\pm}|I},
\end{equation}
where 
$u_{I,J'_{\pm}}$ are sections of
$\gbiga_{J'_{\pm}}(\nbigh^{\Sh})$.
In particular,
we obtain 
sections $u_{I,J_{\pm}}
 \in\gbiga_{J_{\pm}}(\nbigh^{\Sh})$.
It is easy to observe that
the induced sections
$[u_{I,J_{\pm}}]$
of $\nbigh^{\Sh}_{>0,J}$
are independent of 
the choice of decompositions
(\ref{eq;18.4.5.50}).
Moreover, we have
$[u_{I,J_{+}}]
=[u_{I,J_-}]$.
Thus, we obtain the elements
$R_{I,J}(v):=[u_{I,J_{\pm}}]
\in K_{>0,J}$.
We can easily observe the following lemma.
\begin{lem}
$R_{I,J}(v)$
are independent of 
the choice of a connected component 
$I\subset J\setminus S_0(\nbigi)$.
\hfill\qed
\end{lem}

It is easy to see $R_J(v)=R_{I,J}(v)$
by taking a connected component $I$ of
$J\setminus S_0(\gminia)$.
By the construction, the following holds.
\begin{lem}
Let $I$ be a connected component of
$\real\setminus S_0(\nbigi)$.
For any $v\in H^0(\real,\nbigh^{\Sh})$,
we take a decomposition
\[
 v_{|I}=\sum_{\substack{
 J'\in T(\nbigi)_I}} 
 u_{J'_{\pm}|I},
\quad
\mbox{\rm where}\,\,
 u_{J'_{\pm}}\in
 H^0\bigl(
 J'_{\pm},\gbiga_{J'_{\pm}}(\nbigh^{\Sh})
 \bigr).
\]
Then, $R_{J}(v)=u_{J_{\pm}}$ in $K_{J,>0}$.
\hfill\qed
\end{lem}

\section{Appendix: Duality}
\label{subsection;18.5.1.3}

We clarify the relation between
the duality of Stokes shells and
the duality of local systems with Stokes structure
(Proposition \ref{prop;20.9.8.40}).
The reader can skip this section.

We set $\nbigitilde^{\lor}:=-\nbigitilde$
and $\nbigi^{\lor}:=-\nbigi$.
We have
$\pi_{\omega}(\nbigitilde^{\lor})
=\nbigi^\lor$.
Under the assumption
$[\nbigi]=[-\nbigi]$,
we obtain
$[\nbigi^{\lor}]=[-\nbigi^{\lor}]$.
Let $\Sh=(\nbigk_{\bullet},\vecnbigf,\vecnbigr)\in
 \Shcat(\nbigitilde)$.
We obtain the $2\pi\seisuu$-equivariant
Stokes graded local system
$(\nbigk_{\bullet},\vecnbigf)^{\lor}$
over $(\nbigitilde^{\lor},[\nbigi^{\lor}])$.
Set $(\vecK^{\lor},\vecnbigf^{\lor},\vecPhi^{\lor},\vecPsi^{\lor}):=
 \gbigd((\nbigk_{\bullet},\vecnbigf)^{\lor})$.
We may naturally
$K^{\lor}_{\lambda,J}$
as the dual space of
$K_{-\lambda,J}$.
In particular,
$(K^{\lor})_{>0,J}$,
$(K^{\lor})_{<0,J}$
and 
$(K^{\lor})_{0,J}$
are the dual spaces of 
$K_{<0,J}$,
$K_{>0,J}$
and $K_{0,J}$,
respectively.

For any $(J_1,J_2)\in T_2(\nbigi)$,
we obtain the morphism
\[
 (\nbigr^{\lor})^{J_1}_{J_2}:
 (K^{\lor})_{>0,J_1}
\lrarr
 (K^{\lor})_{<0,J_2}
\]
as the dual of 
$\nbigr^{J_2}_{J_1}$.
There exists the natural bijection
$\nbigb_2(\nbigi^{\lor})
\simeq
 \nbigb_2(\nbigi)$
induced by
$(\lambda_1,\lambda_2;J)\longmapsto
 (-\lambda_2,-\lambda_1;J)$.
Hence,
for any $(\lambda_1,\lambda_2;J)\in \nbigb_2(\nbigi^{\lor})$,
we obtain the morphism
\[
 \bigl(
 \nbigr^{\lor}
 \bigr)^{\lambda_1,J_-}_{\lambda_2,J_+}:
 (K^{\lor})_{\lambda_1,J}
\lrarr
 (K^{\lor})_{\lambda_2,J}
\]
as the dual of
$\nbigr^{-\lambda_2,J_+}_{-\lambda_1,J_-}$.
Thus,
we obtain
$\Sh^{\lor}:=\bigl(
 (\nbigk_{\bullet},\vecnbigf)^{\lor},
 \vecnbigr^{\lor}
 \bigr)
\in \Shcat(\nbigitilde)$.
We shall prove the following proposition
in \S\ref{subsection;20.10.1.3}.
\begin{prop}
\label{prop;20.9.8.40}
There exists a natural isomorphism
$\Locst(\Sh^{\lor})
\simeq
 \Locst(\Sh)^{\lor}$
in $\Loc^{\St}(\nbigitilde^{\lor})$.
\end{prop}

For the proof of Proposition \ref{prop;20.9.8.40},
we shall explain another construction of
the functor $\Locst$
in \S\ref{subsection;20.10.1.5}
after preliminary in \S\ref{subsection;20.10.1.4}.

\subsection{Preliminary}
\label{subsection;20.10.1.4}

Let us give a complement to \S\ref{subsection;20.10.1.10}.
The following map is induced by 
 $-\nbigr^{J_1}_{J_2-\pi/\omega}$
and $\vecPhi$:
\begin{equation}
\label{eq;18.4.30.10}
 K_{>0,J_1}
\lrarr
 K_{<0,J_2-\pi/\omega}
\simeq
 K_{>0,J_2}.
\end{equation}
We obtain the following map
from (\ref{eq;18.4.30.1})
and (\ref{eq;18.4.30.10}):
\begin{equation}
\Upsiloncheck^{J_1}_{J_2}:
 K_{>0,J_1}
\lrarr
  \bigoplus_{\lambda\in[\nbigi^{\ast}_{J_2}]}
 K_{\lambda,J_2}. 
\end{equation}

The following map is induced by 
$-\nbigr^{J_2}_{J_1+\pi/\omega}$
and $\vecPhi$:
\begin{equation}
\label{eq;18.4.30.11}
 K_{>0,J_2}
\lrarr
 K_{<0,J_1+\pi/\omega}
\simeq
 K_{>0,J_1}.
\end{equation}
We obtain the following map
from (\ref{eq;18.4.30.3})
and (\ref{eq;18.4.30.11}):
\begin{equation}
 \Upsiloncheck^{J_2}_{J_1}:
 K_{>0,J_2}
\lrarr
  \bigoplus_{\lambda\in[\nbigi^{\ast}_{J_1}]}
 K_{\lambda,J_1}.
\end{equation}

\subsection{Another description of the associated 
 local systems with Stokes structure}
\label{subsection;20.10.1.5}

Let $\Sh\in \Shcat(\nbigi)$.
We shall give another description of
$\Loc^{\St}(\Sh)$
for the proof of Proposition \ref{prop;20.9.8.40}.

Let $I=\openopen{\theta_0}{\theta_1}$ 
be any connected component of
$\real\setminus S_0(\nbigi)$.
Let $(\nbigh^{\Sh}_{\Ibar},\vecnbigf^{\Sh})$
denote the local system with Stokes structure on $\Ibar$
obtained as the gluing of
$(\nbigh^{\Sh}_{I_{\pm}},\vecnbigf^{\Sh})$
in \S\ref{subsection;18.4.30.50}.
For distinction,
we denote them by
$(\nbigh^{\prime\Sh}_{\Ibar},\vecnbigf^{\prime\Sh})$
in this construction.
Let $I_1:=\openopen{\theta_1}{\theta_2}$
be the connected component of $\real\setminus S_0(\nbigi)$
next to $I$.
Let us construct an isomorphism
\[
 F'_{\theta_1}:
 \nbigh^{\prime\Sh}_{\Ibar|\theta_1}
=\nbigh^{\Sh}_{I_+|\theta_1}
\lrarr
 \nbigh^{\prime\Sh}_{\Ibar_{1}|\theta_1}
=\nbigh^{\Sh}_{I_{1-}|\theta_1}
\]
which preserves the filtrations $\nbigf^{\Sh\,\theta_1}$.
Set $J_1:=\openopen{\theta_1-\pi/\omega}{\theta_1}$.
We have the following morphism
induced by
the identity of
$\bigoplus_{\substack{J'\in T(\nbigi)\\ \theta_1\in J'}}
 K_{>0,J'}$
and the morphisms 
$\Upsiloncheck^{J'}_{J_1+\pi/\omega}$:
\begin{equation}
\label{eq;18.4.30.60}
 \bigoplus_{\substack{J'\in T(\nbigi)\\ \theta_1\in J'}}
 K_{>0,J'}
\lrarr
\bigoplus_{\substack{J'\in T(\nbigi)\\ \theta_1\in J'}}
 K_{>0,J'}
\oplus
 \bigoplus_{\lambda\in[\nbigi_{J_1}^{\ast}]}
 K_{\lambda,J_1+\pi/\omega}.
\end{equation}
Set $J_2:=\openopen{\theta_2-\pi/\omega}{\theta_2}$.
We also have the following isomorphism
induced by $\vecPhi$:
\begin{equation}
\label{eq;18.4.30.61}
 \bigoplus_{\lambda\in [\nbigi_{J_1}]}
 K_{\lambda,J_1}
\simeq
 K_{0,J_2}
\oplus
 \bigoplus_{\lambda\in[\nbigi^{\ast}_{J_1}]}
 K_{\lambda,J_1+\pi/\omega}
\end{equation}
We have the following identity map:
\begin{equation}
\label{eq;18.4.30.62}
 \bigoplus_{\substack{J'\in T(\nbigi)\\ \theta_1\in J'}}
 K_{<0,J'}
\simeq
 \bigoplus_{\substack{J'\in T(\nbigi)\\ \theta_1\in J'}}
 K_{<0,J'}.
\end{equation}
We obtain the desired isomorphism $F'_{\theta_1}$
from (\ref{eq;18.4.30.60}),
(\ref{eq;18.4.30.61})
and 
(\ref{eq;18.4.30.62}).
By the construction,
it preserves the Stokes filtrations
$\nbigf^{\Sh\,\theta_1}$.
By gluing
$(\nbigh^{\prime\Sh}_{\Ibar},\vecnbigf^{\Sh})$
for connected components $I$ of
$\real\setminus S_0(\nbigi)$,
we obtain a $2\pi\seisuu$-equivariant
local system with Stokes structure
$\Locst'(\Sh):=
 (\nbigh^{\prime\,\Sh},\vecnbigf^{\prime\Sh})$
indexed by $\nbigitilde$ on $\real$.

\begin{prop}
\label{prop;18.4.30.100}
There exists an isomorphism
$\Locst(\Sh)
\simeq
 \Locst'(\Sh)$.
\end{prop}
\pf
Let $I=\openopen{\theta_0}{\theta_1}$ 
be any connected component of
$\real\setminus S_0(\nbigi)$.
For any $J\in T(\nbigi)_{I}$,
we set
\[
G_{J,I}:=
\sum_{
 \substack{J'\in T(\nbigi)_{I}
 \\  (J,J')\in T_2(\nbigi)}}
\nbigr^{J}_{J'}:
 K_{>0,J}
\lrarr
 \bigoplus_{
 \substack{J'\in T(\nbigi)_{I}
 \\  (J,J')\in T_2(\nbigi)}}
 K_{<0,J'}.
\]
Set $J_1:=\openopen{\theta_1-\pi/\omega}{\theta_1}$.
Let $G_I$ be the automorphism
of the vector space
\begin{equation}
\label{eq;20.10.1.1}
 K_{0,J_{1}}
\oplus
 \bigoplus_{J\in T(\nbigi)_I}
 \bigoplus_{\lambda\in [\nbigi^{\ast}_J]}
 K_{\lambda,J}
\end{equation}
obtained as
$G_I:=\id-\sum_{J\in T(\nbigi)_I}G_{J,I}$.
It induces an automorphism of 
$\nbigh^{\Sh}_{I_{\pm}}$.
It induces the following isomorphism
\[
 G_I:
 \nbigh^{\prime\Sh}_{\Ibar|I_{\pm}}
=\nbigh^{\Sh}_{I_{\pm}}
\lrarr
 \nbigh^{\Sh}_{I_{\pm}}
=\nbigh^{\Sh}_{\Ibar|I_{\pm}}.
\]
Let $H_I$ be the automorphism of
(\ref{eq;20.10.1.1})
induced by
$\Pi^{J_{1+},J_{1-}}$
on $K_{0,J_1}\oplus K_{<,J_{1}}\oplus K_{>0,J_1}$
and the identity map on the complement
$\bigoplus_{J\in T(\nbigi)_I\setminus\{J_1\}}
 \bigoplus_{\lambda\in [\nbigi^{\ast}_J]}
 K_{\lambda,J}$.
It is easy to check that
$G_I$ and $H_I$ are commutative.
Hence we obtain the induced isomorphism
$G_I:
 \nbigh^{\prime\Sh}_{\Ibar}
\lrarr
 \nbigh^{\Sh}_{\Ibar}$.

Let $I=\openopen{\theta_0}{\theta_1}$
and $I_1=\openopen{\theta_1}{\theta_2}$
be connected components of
$\real\setminus S_0(\nbigi)$.
The proof of Proposition \ref{prop;18.4.30.100}
is reduced to the following lemma.
\begin{lem}
\label{lem;18.4.30.101}
We have
$F_{\theta_1}\circ G_I
=G_{I_1}\circ F'_{\theta_1}$.
\end{lem}
\pf
Set $f_1:=F_{\theta_1}\circ G_I$
and $f_2:=G_{I_1}\circ F'_{\theta_1}$.
Let us prove that $f_1=f_2$.
Set $J_1:=\openopen{\theta_1-\pi/\omega}{\theta_1}$
and $J_2:=\openopen{\theta_2-\pi/\omega}{\theta_2}$.
We use the following identifications in the following argument.
\begin{equation}
\label{eq;20.10.1.2}
 \nbigh^{\Sh}_{\Ibar|\theta_1}
 =\nbigh^{\prime\Sh}_{\Ibar|\theta_1}
=K_{0,J_1}\oplus K_{<0,J_1}\oplus
 K_{>0,J_1}
\oplus
 \bigoplus_{\substack{J\in T(\nbigi)\\
 \theta_1\in J}}
 \bigoplus_{\lambda\in[\nbigi_J^{\ast}]}
 K_{\lambda,J},
\end{equation}
\begin{multline}
 \nbigh^{\Sh}_{\Ibar_1|\theta_1}
 =\nbigh^{\prime\Sh}_{\Ibar_1|\theta_1}
=K_{0,J_2}\oplus K_{<0,J_2}\oplus
 K_{>0,J_2}
\oplus
 \bigoplus_{\substack{J\in T(\nbigi)\\
 \theta_2\in J}}
 \bigoplus_{\lambda\in[\nbigi_J^{\ast}]}
 K_{\lambda,J}
 \\
 =K_{0,J_2}\oplus K_{<0,J_1+\pi/\omega}\oplus
 K_{>0,J_1+\pi/\omega}
\oplus
 \bigoplus_{\substack{J\in T(\nbigi)\\
 \theta_1\in J}}
 \bigoplus_{\lambda\in[\nbigi_J^{\ast}]}
 K_{\lambda,J}
\end{multline}
Let us study the restriction of
$f_i$ to $K_{0,J_{1}}\oplus K_{<0,J_{1}}$.
Set 
\[
 A:=
 \bigoplus_{J_1<J'<J_1+\pi/\omega}
 K_{<0,J'}.
\]
Let us look at the following commutative diagram:
\begin{equation}
 \label{eq;18.4.30.120}
 \begin{CD}
 K_{0,J_{1}}
\oplus
 K_{<0,J_{1}}
 @>{a_0}>>
 K_{0,J_{2}}
\oplus
 K_{>0,J_1+\pi/\omega}
\\
  @V{\id}VV @V{a_1}VV \\
 K_{0,J_{1}}
\oplus
 K_{<0,J_{1}}
 @>{a_2}>>
 K_{0,J_{2}}
\oplus
 K_{>0,J_1+\pi/\omega}
\oplus A
 \end{CD}
\end{equation}
Here, 
$a_0$ is induced by $\vecPhi$,
$a_1$ is induced by
the identity and 
$-\nbigr^{J_1+\pi/\omega}_{J'}$,
and $a_2$ is induced by 
$\Phi$ and $\Upsilontilde^{J_1}_{J'}$.
We can observe that 
$f_1$ is identified with
the composite of 
the left vertical arrow and the lower horizontal arrow,
and that
$f_2$ is identified with
the composite of
the upper horizontal arrows and the right vertical arrow.
Hence, we obtain $f_1=f_2$
on 
$K_{0,J_1}
\oplus
 K_{<0,J_{1}}$.

Let us study the restriction of $f_i$ to $K_{>0,J_1}$.
We also have the following commutative diagram:
\begin{equation}
\label{eq;18.4.30.131}
 \begin{CD}
 K_{>0,J_{1}} 
 @>{a_0}>>
 K_{<0,J_1+\pi/\omega}\\
 @V{a_1}VV @V{a_2}VV\\
 K_{>0,J_{1}}\oplus A
 @>{a_3}>>
  K_{<0,J_1+\pi/\omega}\oplus A
 \end{CD}
\end{equation}
Here, $a_0$ is induced by $\vecPhi$,
$a_1$ is induced by
the identity and $-\nbigr^{J_1}_{J'}$,
$a_2$ is the identity,
and $a_3$ is induced by
$\vecPhi$ and $\Upsilontilde^{J_1}_{J'}$.
We obtain 
$f_1=f_2$ on $K_{>0,J_1}$
from the commutativity 
of the diagram (\ref{eq;18.4.30.131}).

Take $J'\in T(\nbigi)$
such that $\theta_1\in J'$.
The equality $f_1=f_2$ on $K_{<0,J'}$
follows from 
the obvious commutativity 
of the following diagram:
\begin{equation}
 \label{eq;18.4.30.132}
 \begin{CD}
 K_{<0,J'}
 @>{\id}>>
  K_{<0,J'}\\
 @V{\id}VV @V{\id}VV\\
 K_{<0,J'}
 @>{\id}>>
 K_{<0,J'}.
 \end{CD}
\end{equation}
Let us prove $f_1=f_2$
on $K_{>0,J'}$.
We set
\[
 C:=
\bigoplus_{
 \substack{J''\in T(\nbigi)_I\setminus \{J_1\}
 \\
 (J',J'')\in T_2(\nbigi)} }
 K_{<0,J''}.
\]
Let us study the following diagram:
\begin{equation}
\label{eq;18.4.30.140}
\begin{CD}
 K_{>0,J'}
 @>{a_0}>>
 K_{>0,J'}
\oplus
 K_{>0,J_1+\pi/\omega}
\oplus
 K_{<0,J_1+\pi/\omega}\\
 @V{a_1}VV @V{a_2}VV \\
 K_{>0,J'}
\oplus C
\oplus K_{<0,J_{1}}
 @>{a_3}>>
 K_{>0,J'}
\oplus C
\oplus K_{>0,J_1+\pi/\omega}
\oplus K_{<0,J_1+\pi/\omega}.
\end{CD}
\end{equation}
Here, $a_0$, $a_1$, $a_2$ and $a_3$
are induced by
$F'_{\theta_1}$
$G_I$,
$G_{I_1}$
and $F_{\theta}$,
respectively.
Take $v\in K_{>0,J'}$.
By the construction,
we have
\begin{equation}
a_0(v)=
v-\Phi_{\lambda_-(J_1)}^{J_1+\pi/\omega,J_1}\circ
 \nbigr^{J'}_{J_1}(v)
+\nbigr^{J'}_{J_1+\pi/\omega}(v),
\end{equation}
\begin{equation}
 a_1(v)=
 v
-\sum_{\substack{J''\in T(\nbigi)_I\setminus J_1\\
 (J'',J')\in T_2(\nbigi)}} 
\nbigr^{J'}_{J''}(v)
-\nbigr^{J'}_{J_1}(v).
\end{equation}
We have
\begin{multline}
 a_2\Bigl(
-\Phi_{\lambda_-(J_1)}^{J_1+\pi/\omega,J_1}\circ
 \nbigr^{J'}_{J_1}(v)
 \Bigr)
= \\
-\Phi_{\lambda_-(J_1)}^{J_1+\pi/\omega,J_{1}}
 \circ\nbigr^{J'}_{J_1}(v)
+\sum_{\substack{J''\in T(\nbigi)_{I_1}\\
 (J'',J_1+\pi/\omega)\in T_2(\nbigi)}}
 \nbigr^{J_1+\pi/\omega}_{J''}\circ
 \Phi_{\lambda_-(J_1)}^{J_1+\pi/\omega,J_{1}}
 \circ\nbigr^{J'}_{J_1}(v)
\\
=
-\Phi_{\lambda_-(J_1)}^{J_1+\pi/\omega,J_{1}}
 \circ\nbigr^{J'}_{J_1}(v)
+\sum_{\substack{J''\in T(\nbigi)\\
 J_1<J''<J_1+\pi/\omega }}
 \nbigr^{J_1+\pi/\omega}_{J''}\circ
 \Phi_{\lambda_-(J_1)}^{J_1+\pi/\omega,J_{1}}
 \circ\nbigr^{J'}_{J_1}(v).
\end{multline}
We have
\begin{multline}
 a_3\Bigl(
-\nbigr^{J'}_{J_1}(v)
 \Bigr)
=
-\Phi^{J_1+\pi/\omega,J_{1}}_{\lambda_-(J_1)}
 \circ\nbigr^{J'}_{J_1}(v)
-\sum_{\substack{J''\in T(\nbigi)\\
 \theta_1\in J''}}
 \Upsilontilde^{J_1}_{J''}\circ
 \nbigr^{J'}_{J_1}(v)
\\
=
-\Phi^{J_1+\pi/\omega,J_{1}}_{\lambda_-(J_1)}
 \circ
 \nbigr^{J'}_{J_1}(v)
-\sum_{\substack{J''\in T(\nbigi)\\
 \theta_1\in J''}}
 (-\nbigr^{J_1+\pi/\omega}_{J''})
 \circ
 \Phi_{\lambda_-(J_1)}^{J_1+\pi/\omega,J_{1}}
 \circ
 \nbigr^{J'}_{J_1}(v).
\end{multline}
Hence, we have 
\[
 a_2\Bigl(
-\Phi_{\lambda_-(J_1)}^{J_1+\pi/\omega,J_{1}}\circ
 \nbigr^{J'}_{J_1}(v)
 \Bigr)
= a_3\Bigl(
 -\nbigr^{J'}_{J_1}(v)
 \Bigr).
\]
We also have
\begin{multline}
 a_2\Bigl(
 v+\nbigr^{J'}_{J_1+\pi/\omega}(v)
 \Bigr)
=v-\sum_{\substack{
 J''\in T(\nbigi)_{I_1}\setminus\{J_1+\pi/\omega\}\\
 (J',J'')\in T_2(\nbigi)
 }}
\nbigr^{J'}_{J''}(v)
=v-\sum_{\substack{
 J''\in T(\nbigi)_{I}\setminus\{J_1\}\\
 (J',J'')\in T_2(\nbigi)
 }}
\nbigr^{J'}_{J''}(v)\\
=a_3\Bigl(v
-\sum_{\substack{J''\in T(\nbigi)_I\setminus\{J_1\}\\
 (J',J'')\in T_2(\nbigi)}}
\nbigr^{J'}_{J''}(v)
 \Bigr).
\end{multline}
Thus, we obtain the commutativity of
(\ref{eq;18.4.30.140}),
which implies
$f_1=f_2$
on $K_{>0,J'}$.
The proof of Lemma \ref{lem;18.4.30.101}
and Proposition \ref{prop;18.4.30.100}
are completed.
\hfill\qed

\subsection{Proof of Proposition \ref{prop;20.9.8.40}}
\label{subsection;20.10.1.3}

By the construction,
we have
$\Locst(\Sh)^{\lor}
\simeq
 \Locst'(\Sh^{\lor})$,
which is isomorphic to
$\Locst(\Sh^{\lor})$
by Proposition \ref{prop;18.4.30.100}.
Thus, we obtain the claim of the proposition.
\hfill\qed

\chapter{Preliminary for meromorphic flat bundles}
\label{chapter;21.6.10.30}

\section{Asymptotic analysis and Riemann-Hilbert correspondence}

\subsection{Formal structure}

Let $\Delta_z:=\{z\in\cnum\,|\,|z|<\epsilon\}$
for a positive integer $\epsilon$.
\index{set $\Delta_z$}
We set $\Delta_z^{\ast}:=\Delta_z\setminus\{0\}$.
\index{set $\Delta_z^{\ast}$}
Let $(V,\nabla)$ be a meromorphic flat bundle on $(\Delta_z,0)$.
According to Hukuhara-Levelt-Turrittin theorem,
there exist a positive integer $p$,
a $\Gal(p)$-invariant
finite subset $\nbigi\subset z_p^{-1}\cnum[z_p^{-1}]$,
and a decomposition
\begin{equation}
 (V,\nabla)\otimes_{\nbigo_{\Delta_z,0}}
  \cnum[\![z_p]\!]
  =\bigoplus_{\gminia\in\nbigi}
  (\Vhat_{\gminia},\nablahat_{\gminia}),
\end{equation}
where $(\Vhat_{\gminia},\nablahat_{\gminia}-d\gminiahat\id)$
are regular singular.
We allow the case $\Vhat_{\gminia}=0$.
We set $\nbigi(V)=\{\gminia\in\nbigi\,|\,V_{\gminia}\neq 0\}$.
If $\nbigi(V)\subset z^{-1}\cnum[z^{-1}]$,
we say that $(V,\nabla)$ is unramified.
There exists the Hukuhara-Levelt-Turrittin decomposition
for $(V,\nabla)\otimes\cnum[\![z]\!]$ in the unramified case.
\index{Hukuhara-Levelt-Turrittin decomposition}
\index{set $\nbigi(V)$}

For a $\Gal(p)$-invariant
finite subset $\nbigi\subset z_p^{-1}\cnum[z_p^{-1}]$,
let $\Mer(\Delta_z,0,\nbigi)$ denote the category of
meromorphic flat bundles $(V,\nabla)$ on $(\Delta_z,0)$
such that $\nbigi(V)\subset\nbigi$.
\index{category $\Mer(\Delta_z,0,\nbigi)$}
Morphisms $(V_1,\nabla)\to (V_2,\nabla)$
are defined to be morphisms of
$\nbigo_{\Delta_z}(\ast 0)$-modules
compatible with the connections.

\subsection{The unramified case}
\subsubsection{Asymptotic analysis}
\label{subsection;24.3.28.10}

Let $\varpi:\Deltatilde_z\to \Delta_z$
denote the oriented real blow up along $0$,
i.e.,
$\Deltatilde_z=\closedopen{0}{\epsilon}\times S^1$,
and $\varpi(r,e^{\sqrt{-1}\theta})=re^{\sqrt{-1}\theta}$.
\index{oriented real blow up $\Deltatilde_z$}
A $C^{\infty}$-function $f$ on 
an open subset $\nbigu\subset \Deltatilde_z$
is called holomorphic if it is holomorphic on
$\nbigu\setminus\varpi^{-1}(0)$.
Let $\nbigo_{\Deltatilde_z}$ denote the sheaf of
holomorphic functions on $\Deltatilde_z$.
\index{sheaf $\nbigo_{\Deltatilde_z}$}
We also set
$\nbigo_{\Deltatilde_z}(\ast 0)
=\nbigo_{\Deltatilde_z}
\otimes_{\varpi^{-1}\nbigo_{\Delta_z}}
 \varpi^{-1}(\nbigo_{\Delta_z}(\ast 0))$.
\index{sheaf $\nbigo_{\Deltatilde_z}(\ast 0)$}
 For any section $f$ of $\nbigo_{\Deltatilde_z}$
on $\nbigu$,
we obtain the power series
$f_{|\widehat{\nbigu\cap\varpi^{-1}(0)}}\in \cnum[\![z]\!]$
as the Taylor series at any point of
$\nbigu\cap\varpi^{-1}(0)$.
For a section $f$ of $\nbigo_{\Deltatilde_z}(\ast 0)$,
we obtain $f_{|\widehat{\nbigu\cap\varpi^{-1}(0)}}
\in\cnum(\!(z)\!)$.

Let $(V,\nabla)$ be a meromorphic flat bundle on $(\Delta_z,0)$
which is unramified.
There exist a finite subset $\nbigi\subset z^{-1}\cnum[z^{-1}]$,
meromorphic flat bundles
$(V_{\gminia},\nabla_{\gminia})$ $(\gminia\in\nbigi)$
such that $(V_{\gminia},\nabla_{\gminia}-d\gminia\id)$
are regular singular,
and an isomorphism
\[
 \Phihat:
(V,\nabla)\otimes\cnum[\![z]\!]
 \simeq
 \bigoplus_{\gminia\in\nbigi}
 (V_{\gminia},\nabla_{\gminia})
 \otimes\cnum[\![z]\!].
\]
Note that $\Phihat$ is not necessarily convergent.

We set
$\varpi^{\ast}(V)=
\varpi^{-1}(V)\otimes_{\varpi^{-1}\nbigo_{\Delta_z}(\ast 0)}
\nbigo_{\Deltatilde_z}(\ast 0)$.
Let $P$ be any point of $\varpi^{-1}(0)$.
According to the classical asymptotic analysis,
there exist a neighbourhood $\nbigu_P$ of $P$ in $\Deltatilde_z$
and an isomorphism
\begin{equation}
\label{eq;24.3.28.1}
 \Phi_{\nbigu_P}:
 \varpi^{\ast}(V,\nabla)_{|\nbigu_P}
 \simeq
 \bigoplus_{\gminia}
 \varpi^{\ast}(V_{\gminia},\nabla_{\gminia})_{|\nbigu_P}
\end{equation}
such that
$\Phi_{\nbigu_P|\widehat{\nbigu\cap\varpi^{-1}(0)}}
=\Phihat$.
Note that, in general,
such an isomorphism is not unique.

\subsubsection{Stokes filtrations}
\label{subsection;24.3.28.2}

Let $\nbigl'$ be the local system on $\Delta_z^{\ast}$
obtained as the sheaf of flat sections of $(V,\nabla)$.
It extends to a local system on $\Deltatilde_z$,
denoted by $\nbigl$.
We set $L_{S^1}=\nbigl_{|\varpi^{-1}(0)}$.
\index{local system $L_{S^1}$}

For any $P\in \varpi^{-1}(0)$,
let $\nbigl_P$ be the stalk of $\nbigl$ at $P$.
Let $\nbigf^P_{\gminia}(\nbigl_P)\subset\nbigl_P$
denote the subspace of $s\in\nbigl_P$
satisfying the following condition.
\begin{itemize}
 \item Let $\vecv=(v_1,\ldots,v_r)$ be a frame of $V$
       over $\nbigo_{\Delta_z}(\ast 0)$.
       For the expression $s=\sum s_iv_i$,
       we obtain
       $|e^{\gminia}s_i|=O(|z|^{-N})$
       on a sector around $P$
       for some $N>0$.
\end{itemize}
In this way, we obtain the filtration
$\nbigf^P$ of $\nbigl_P$
indexed by $(\nbigi,\leq_P)$.
By the existence of an isomorphism (\ref{eq;24.3.28.1}),
there exists a splitting
$\nbigl_P=\bigoplus_{\gminia\in\nbigi}G_{P,\gminia}$
such that
$\nbigf^{P}_{\gminia}(\nbigl_P)
=\bigoplus_{\gminib\leq_P\gminia} G_{P,\gminib}$.
The following proposition is due to Deligne and Malgrange.
\begin{prop}
The family of the filtration
$\nbigf^P$ $(P\in\varpi^{-1}(0))$ is a Stokes structure of
the local system $L_{S^1}$. 
 \hfill\qed
\end{prop}

We obtain a $2\pi\seisuu$-equivariant local system
with Stokes structure
$(L,\vecnbigf)=\RH(V,\nabla)$
on $\real$ obtained as the pull back of $(L_{S^1},\vecnbigf)$
by the map
$\real\ni\theta\mapsto e^{\sqrt{-1}\theta}\in\varpi^{-1}(0)$.

\subsubsection{Riemann-Hilbert correspondence}

Let $\nbigi\subset z^{-1}\cnum[z^{-1}]$ be a finite subset.
We obtain a functor from
$\RH:\Mer(\Delta_z,0,\nbigi)
\to \Loc^{\St}(\nbigi)$.
The following is a version of the Riemann-Hilbert correspondence.
\begin{thm}[Deligne-Malgrange]
The functor 
$\RH:\Mer(\Delta_z,0,\nbigi)
\to \Loc^{\St}(\nbigi)$
is an equivalence.
\hfill\qed
\end{thm}

\subsubsection{Complement}
\label{subsection;24.3.28.3}

Let $(V,\nabla)$ and $(L_{S^1},\vecnbigf)$
be as in \S\ref{subsection;24.3.28.2}.
Let $v_1,\ldots,v_r$ be a frame of $V$ around $0$.

Let $P\in\varpi^{-1}(0)$.
Let $\nbigu_P$ be a simply connected neighbourhood of
$P$ in $\Deltatilde_z$.
There exists a decomposition
$\nbigl_{|\nbigu_P}=\bigoplus \nbigg_{\gminia}$
such that
it induces a splitting
$L_P=\bigoplus_{\gminia\in\nbigi} G_{P,\gminia}$
of the filtration $\nbigf^P$.
Let $s_1,\ldots,s_r$ be a base of $\nbigl_{|\nbigu_P}$
compatible with the decomposition,
i.e., $s_i\in\nbigg_{\gminia_i}$.
We obtain the matrix $A=(A_{i,j})$
of holomorphic functions on $\nbigu_P\setminus\varpi^{-1}(0)$
determined by
$e^{\gminia_j}s_j=\sum A_{i,j}v_i$.

\begin{lem}
There exists a neighbourhood
$\nbigu_P'$ of $P$ in $\nbigu_P$
such that
$|A_{i,j}|=O(|z|^{-N})$
and $|\det A|^{-1}=O(|z|^{-N})$
for some $N>0$
on $\nbigu_P'\setminus\varpi^{-1}(0)$.
\hfill\qed
\end{lem}

\subsubsection{The inverse construction}

Let $L_{S^1}$ be a local system
with Stokes structure indexed by $\nbigi$
on $\varpi^{-1}(0)$.
There exists the local system $\nbigl'$ on $\Delta_z^{\ast}$
obtained as the pull back via the projection
$\Delta_z^{\ast}\to\varpi^{-1}(0)$.
We set
$V'=\nbigo_{\Delta_z^{\ast}}
\otimes \nbigl'$,
which is equipped with the connection $\nabla$
such that any section of $\nbigl'$ are flat.
Let $U\subset\Delta_z$ be an open subset.
If $0\not\in U$,
we set $V(U)=V'(U)$.
If $0\in U$,
$V(U)$ be the space of
$f\in V'(U\setminus\{0\})$
satisfying the following conditions
for any $P\in\varpi^{-1}(0)$.
\begin{itemize}
 \item Let $s_1,\ldots,s_r$ be a frame of $L$ around $P$
       as in \S\ref{subsection;24.3.28.3}.
       We obtain the expression
       $f=\sum f_ie^{\gminia_i}s_i$.
       Then,
       $|f_i|=O(|z|^{-N})$
       for some $N>0$.
\end{itemize}
Then, we can prove that
$V$ is a locally free $\nbigo_{\Delta_z}(\ast 0)$-module,
and that
$(V,\nabla)\in \Mer(\Delta_z,0,\nbigi)$.
This is a quasi-inverse of $\RH$.

\subsection{The ramified case}

Let $\nbigi\subset z_p^{-1}\cnum[z_p^{-1}]$
be a $\Gal(p)$-invariant finite subset.
Let $(V,\nabla)\in\Mer(\Delta_z,0,\nbigi)$.
Let $\rho_p:\cnum\to\cnum$ be the map
defined by $\rho_p(z_p)=z_p^p$.
The meromorphic flat bundle
$\rho_p^{\ast}(V,\nabla)$ is unramified.

Let $L_{S^1}$ be the local system on $\varpi^{-1}(0)$
obtained from $(V,\nabla)$ as in \S\ref{subsection;24.3.28.2}.
Let $\varpi_p:\Deltatilde_{z_p}\to\Delta_{z_p}$
be the oriented real blow up.
There exists the map
$\rho_p:\Deltatilde_{z_p}\to\Delta_{z_p}$
induced by $\rho_p$,
which induces
$\rho_p:\varpi_p^{-1}(0)\to \varpi^{-1}(0)$.
We obtain the local system $\rho_p^{-1}(L_{S^1})$
on $\varpi_p^{-1}(0)$.
We obtain the family of Stokes filtrations
$\nbigf^P$ $(P\in\varpi_p^{-1}(0))$ of $\rho_p^{-1}(L_{S^1})_P$
as in \S\ref{subsection;24.3.28.2}.
This is a $\Gal(p)$-equivariant local system
of $\rho_p^{-1}(L_{S^1})$
indexed by $\nbigi$.

Let $\real\to \varpi_p^{-1}(0)$
be defined by
$\theta\mapsto e^{\sqrt{-1}\theta/p}$.
Let $L$ be the $2\pi\seisuu$-equivariant local system
obtained as the pull back of $\rho_p^{-1}(L_{S^1})$,
which equals the pull back of $L_{S^1}$
by the map $\theta\longmapsto e^{\sqrt{-1}\theta}$.
It is equipped with the $2\pi\seisuu$-equivariant
Stokes structure indexed by $\vecnbigi$
obtained as the pull back of
$\vecnbigf^P$ $(P\in\varpi_p^{-1}(0))$.

\begin{thm}[Deligne-Malgrange]
This procedure induces an equivalence
$\RH:\Mer(\Delta_z,0,\nbigi)
\to \Loc^{\St}(\nbigi)$.
\hfill\qed
\end{thm}

\section{Induced meromorphic flat bundles}
\label{subsection;18.6.23.3}

\subsection{Local case}
\label{subsection;18.7.10.1}

Let $p$ be a positive integer.
Let $\nbigi$ be a $\Gal(p)$-invariant finite subset
$z_p^{-1}\cnum[z_p^{-1}]$.
Let $(V,\nabla)\in\Mer(\Delta_z,0,\nbigi)$.
Let $\nbigl$ be 
the local system on $\Deltatilde_z$
associated to $(V,\nabla)_{|\Delta_z^{\ast}}$.
We set
$L_{S^1}:=\nbigl_{|\varpi^{-1}(0)}$.
We obtain the $2\pi\seisuu$-equivariant
local system with Stokes structure
$(L,\vecnbigf)=\RH(V,\nabla)$
indexed by $\nbigi$ on $\real$.
The descent of $L$ is naturally identified with $L_{S^1}$.
For $2\pi\seisuu$-equivariant constructible subsheaves
$K\subset L$,
let $K_{S^1}\subset L_{S^1}$ denote the subsheaf
obtained as the descent.
\index{constructible subsheaf $K_{S^1}$}

Take $\omega\in\rnum_{>0}$.
By the procedures in \S\ref{subsection;18.4.18.1},
we obtain 
the $2\pi\seisuu$-equivariant local systems with Stokes structure
$\nbigs_{\omega}(L,\vecnbigf)$
and $\nbigt_{\omega}(L,\vecnbigf)$.
Let $\nbigs_{\omega}(V,\nabla)$
and $\nbigt_{\omega}(V,\nabla)$ 
denote the corresponding meromorphic flat bundles
on $(\Delta_z,0)$.
\index{meromorphic flat bundles
$\nbigs_{\omega}(V,\nabla)$, $\nbigt_{\omega}(V,\nabla)$}
Similarly,
let $\nbigstilde_{\omega}(V,\nabla)$
and $\nbigttilde_{\omega}(V,\nabla)$
be the meromorphic flat bundles
corresponding to
$\nbigstilde_{\omega}(L,\vecnbigf)$
and $\nbigttilde_{\omega}(L,\vecnbigf)$,
respectively.
\index{meromorphic flat bundles
$\nbigstilde_{\omega}(V,\nabla)$, $\nbigttilde_{\omega}(V,\nabla)$}
By the construction,
we have the natural isomorphism 
of flat bundles
$(V,\nabla)_{|\Delta_z^{\ast}}
\simeq
 \nbigs_{\omega}(V,\nabla)_{|\Delta_z^{\ast}}$.
We also have the natural isomorphism
of regular singular meromorphic flat bundles 
$\nbigt_{\omega}\circ\nbigs_{\omega}(V,\nabla)
\simeq
 \nbigs_{\omega}\circ\nbigt_{\omega}(V,\nabla)$
on $(\Delta_z,0)$.

\subsection{Global case}
\label{subsection;20.10.24.3}

Let $C$ be a compact Riemann surface 
which may have smooth boundary.
Let $D\subset C\setminus \del C$
be a finite subset.
Set $C^{\circ}:=C\setminus D$.
\index{set $C^{\circ}$}

Let $(V,\nabla)$ be a meromorphic flat bundle
on $(C,D)$.
Take $Q\in D$ with 
a holomorphic coordinate 
neighbourhood $(C_Q,z)$ such that $z(Q)=0$.
Here, we assume that
the closure of $C_Q$ is isomorphic to
a closed disc.
Take $\omega\in\rnum_{>0}$.
By applying the procedure in \S\ref{subsection;18.7.10.1}
to $(V_1,\nabla):=(V,\nabla)_{|C_Q}$,
we obtain meromorphic flat bundles
$\nbigs_{\omega}(V_1,\nabla)$ 
and
$\nbigt_{\omega}(V_1,\nabla)$ 
on $(C_Q,Q)$.
We set
$\nbigt_{\omega}^Q(V,\nabla):=
 \nbigt_{\omega}(V_1,\nabla)$.
\index{meromorphic flat bundle $\nbigt^Q_{\omega}(V,\nabla)$}
By gluing 
$(V,\nabla)_{|C\setminus\{Q\}}$
and $\nbigs_{\omega}(V_1,\nabla)$
via the natural isomorphism 
$(V,\nabla)_{|C_Q\setminus \{Q\}}\simeq
\nbigs_{\omega}(V_1,\nabla)_{|C_Q\setminus\{Q\}}$,
we obtain a meromorphic flat bundle on $(C,D)$,
which we denote by $\nbigs^Q_{\omega}(V,\nabla)$.
\index{meromorphic flat bundle $\nbigs^Q_{\omega}(V,\nabla)$}
Similarly, we obtain
$\nbigstilde^Q_{\omega}(V,\nabla)$ on $(C,D)$,
and 
$\nbigttilde^Q_{\omega}(V,\nabla)$ on $(C_Q,Q)$.
\index{meromorphic flat bundle $\nbigttilde^Q_{\omega}(V,\nabla)$}
\index{meromorphic flat bundle $\nbigstilde^Q_{\omega}(V,\nabla)$}

\section{Constructible sheaves associated to
 meromorphic flat bundles}

\subsection{Associated constructible sheaves}
\label{subsection;18.5.13.1}

Let $C$ be a complex curve without boundary,
which is not necessarily compact.
Let $D\subset C$ be a discrete subset.
Set $C^{\circ}:=C\setminus D$.
\index{open subset $C^{\circ}$}
Let $\varpi:\Ctilde\lrarr C$  denote 
the oriented real blow up of $C$ along $D$.
Set $\Dtilde:=\varpi^{-1}(D)$.
We have $\del \Ctilde=\Dtilde$.
\index{set $C^{\circ}$}
\index{oriented real blow up $\Ctilde$}
\index{set $\Dtilde$}

Let $(V,\nabla)$ be a meromorphic flat bundle on $(C,D)$.
Let $(\nbigl(V),\vecnbigf)$
be the local system with Stokes structure on 
$(\Ctilde,\Dtilde)$,
i.e.,
$\nbigl(V)$ is a local system on $\Ctilde$
induced by the local system
associated to $(V,\nabla)_{|C^{\circ}}$,
and $\vecnbigf=(\nbigf^{P}\,|\,P\in \Dtilde)$
is the family of Stokes filtrations
of $\nbigl(V)_P$.
\index{local system $\nbigl(V)$}
\index{local system with Stokes structure $(\nbigl(V),\vecnbigf)$}

For each $P\in \Dtilde$,
we obtain the subspace 
$\nbigl(V)_P^{<0}:=\nbigf^P_{<0}(\nbigl(V)_P)$
of $\nbigl(V)_P$.
For each $P\in C^{\circ}$,
we set $\nbigl(V)^{<0}_P:=\nbigl(V)_P$.
They determine a constructible subsheaf of $\nbigl(V)$
on $\Ctilde$,
denoted by $\nbigl^{<0}(V)$.
\index{constructible subsheaf $\nbigl^{<0}(V,\nabla)$}

We have the meromorphic flat bundle
$(V^{\lor},\nabla)$ on $(C,D)$,
where we set $V^{\lor}:=
 \nhom_{\nbigo_{C}(\ast D)}\bigl(V,\nbigo_{C}(\ast D)\bigr)$,
and $\nabla$ denotes the naturally induced connection on $V^{\lor}$.
Let $V(!D)$ denote the $\nbigd_{C}$-module
obtained as the dual of the $\nbigd_{C}$-module $(V^{\lor},\nabla)$.
\index{$\nbigd$-module \mbox{$V(!D)$}}
We have the naturally defined morphism 
of $\nbigd_{C}$-modules $V(!D)\lrarr V$.
It is well known that there exists the natural isomorphism
\[
 V(!D)\otimes\Omega^{\bullet}_{C}
\simeq
R\varpi_{\ast}
 \nbigl^{<0}(V,\nabla)
\]
in the derived category of $\cnum_C$-modules.
(For example, see 
\cite{Mebkhout-positivity, Mochizuki-Betti, Sabbah-irregularity}.)

Similarly,
for $P\in\Dtilde$,
we obtain the subspace
$\nbigl^{\leq 0}(V)_{P}:=
 \nbigf^P_{\leq 0}(\nbigl(V)_P)$.
For $P\in C\setminus D$,
we set
 $\nbigl^{\leq 0}(V)_P:=\nbigl_P$.
They determine the constructible subsheaf
$\nbigl^{\leq 0}(V)$
of $\nbigl(V)$.
\index{constructible subsheaf $\nbigl^{\leq 0}(V)$}
Recall that there exists
the natural isomorphism
\[
 V\otimes\Omega^{\bullet}_{C}
\simeq
 R\varpi_{\ast}\nbigl^{\leq 0}(V)
\]
in the derived category of 
$\cnum_C$-modules.
(For example, see 
\cite{Mebkhout-positivity, Mochizuki-Betti, Sabbah-irregularity}.)

\vspace{.1in}
More generally,
for any $\varrho\in\Dsf(D)$,
let $\nbigl^{\varrho}(V)$
denote the constructible subsheaf of $\nbigl(V)$
on $\Ctilde(D)$
determined by the following conditions.
\index{constructible subsheaf $\nbigl^{\varrho}(V)$}
\begin{itemize}
\item
 $\nbigl^{\varrho}(V)_{|C^{\circ}}=\nbigl_{|C^{\circ}}$.
\item
 $\nbigl^{\varrho}(V)_{|\varpi^{-1}(Q)}
 =\nbigl^{\leq 0}(V)_{|\varpi^{-1}(Q)}$
for $Q\in D$ such that $\varrho(Q)=\ast$.
\item
 $\nbigl^{\varrho}(V)_{|\varpi^{-1}(Q)}
 =\nbigl^{<0}(V)_{|\varpi^{-1}(Q)}$
for $Q\in D$ such that $\varrho(Q)=!$.
\end{itemize}
We have the $\nbigd_C$-module 
\[
 V(\varrho):=
 V(!D)\otimes\nbigo_C(\ast\varrho^{-1}(\ast)),
\]
and the isomorphism
$\Omega^{\bullet}\otimes V(\varrho)
\simeq
 R\varpi_{\ast}\nbigl^{\varrho}(V)$
in the derived category of $\cnum_{C}$-modules.
\index{$\nbigd$-module $V(\varrho)$}

\subsection{Induced morphisms in the local case}
\label{subsection;18.4.17.10}

Let $p$ be a positive integer.
Let $\nbigi$ be a $\Gal(p)$-invariant finite subset
of $z_p^{-1}\cnum[z_p^{-1}]$.
Let $(V,\nabla)\in\Mer(\Delta_z,0,\nbigi)$.
We obtain the local system $\nbigl(V)$
on $\Deltatilde_z$.
We also obtain
$(L,\vecnbigf)=\RH(V,\nabla)\in\Loc^{\St}(\nbigi)$.

\subsubsection{Some morphisms for $\nbigl^{<0}$}
Recall that
for any $2\pi\seisuu$-equivariant constructible sheaf $K$ on $\real$,
let $K_{S^1}$ denote the constructible sheaf on $S^1$
obtained as the descent of $K$.
\index{constructible subsheaf $K_{S^1}$}
We have
$L_{S^1}=\nbigl_{|\varpi^{-1}(0)}$.
By the construction,
the restriction
$\nbigl^{<0}(V)_{|\varpi^{-1}(0)}$
is $L^{<0}_{S^1}$.
Similarly, 
$\nbigl^{<0}(\nbigs_{\omega}(V))_{|\varpi^{-1}(0)}$
is $L^{(\omega)\,<0}_{S^1}$.
Hence, 
there exists the following naturally defined monomorphism
\[
 \nbigl^{<0}(\nbigs_{\omega}(V))
\lrarr
 \nbigl^{<0}(V).
\]
Moreover,
there exists the following natural isomorphisms:
\begin{multline}
 \nbigl^{<0}(V)\big/
  \nbigl^{<0}(\nbigs_{\omega}(V))
\simeq
 \iota_{\ast}\bigl(
 L_{S^1}^{<0}\big/L_{S^1}^{(\omega)\,<0}
 \bigr)
\simeq \\
\nbigl^{<0}(\nbigt_{\omega}(V))\big/
  \nbigl^{<0}(\nbigs_{\omega}\circ\nbigt_{\omega}(V))
 =\iota_{\ast}\iota^{-1}\bigl(
 \nbigl^{<0}(\nbigt_{\omega}(V))
 \bigr),
\end{multline}
where $\iota:\varpi^{-1}(0)\lrarr \Deltatilde_z$
denotes the inclusion.

Let $q:\Deltatilde_z\lrarr \varpi^{-1}(0)$
be the projection $q(r,e^{\sqrt{-1}\theta})=e^{\sqrt{-1}\theta}$.
We obtain the natural monomorphism
$q^{-1}\bigl(L_{S^1}^{(\omega)\,\leq 0}\bigr)
\lrarr
 \nbigl$.
There exists the constructible subsheaf 
$\nbiglcheck^{<0}\bigl(
 \nbigt_{\omega}(V)\bigr)
\subset
 q^{-1}\bigl(L_{S^1}^{(\omega)\,\leq 0}\bigr)$
determined by the following conditions.
\index{constructible subsheaf $\nbiglcheck^{<0}(\nbigt_{\omega}(V))$}
\begin{itemize}
\item
$\nbiglcheck^{<0}\bigl(
 \nbigt_{\omega}(V)\bigr)_{|\Delta_z^{\ast}}$
is equal to 
$q^{-1}\bigl(L_{S^1}^{(\omega)\,\leq 0}\bigr)_{|\Delta_z^{\ast}}$.
\item
$\nbiglcheck^{<0}\bigl(
 \nbigt_{\omega}(V)\bigr)_{|\varpi^{-1}(0)}$
is equal to $L_{S^1}^{<0}$.
\end{itemize}
There exists the natural monomorphism
\[
 \nbiglcheck^{<0}\bigl(
 \nbigt_{\omega}(V)
 \bigr)
\lrarr
  \nbigl^{<0}(V).
\]
There also exists the following exact sequence:
\begin{equation}
\label{eq;18.5.12.10}
 0\lrarr
    q^{-1}\bigl(
 L^{(\omega)<0}_{S^1}
 \bigr)
\lrarr
 \nbiglcheck^{<0}\bigl(
 \nbigt_{\omega}(V)
 \bigr)
\lrarr
 \nbigl^{<0}\bigl(\nbigt_{\omega}(V)\bigr)
\lrarr 0.
\end{equation}

\subsubsection{Some morphisms for $\nbigl^{\leq 0}$}
Similarly, there exists the following natural monomorphism:
\[
 \nbigl^{\leq 0}(V)
\lrarr
 \nbigl^{\leq 0}\bigl(\nbigs_{\omega}(V)\bigr).
\]
We naturally obtain
\[
  \nbigl^{\leq 0}\bigl(\nbigs_{\omega}(V)\bigr)
 \big/ \nbigl^{\leq 0}(V)
\simeq
 \iota_{\ast}\bigl(
 L^{(\omega)\,\leq 0}_{S^1}\big/
 L^{\leq 0}_{S^1}
 \bigr)
\simeq
  \nbigl^{\leq 0}
 \bigl(\nbigs_{\omega}\nbigt_{\omega}(V)\bigr)
 \big/ \nbigl^{\leq 0}\bigl(
 \nbigt_{\omega}(V)
\bigr).
\]

There exists the constructible quotient sheaf
$\nbiglcheck^{\leq 0}(\nbigt_{\omega}(V))$
of $\nbigl^{\leq 0}(V)$
determined by the following conditions.
\index{constructible quotient sheaf $\nbiglcheck^{\leq 0}(\nbigt_{\omega}(V))$}
\begin{itemize}
\item
$\nbiglcheck^{\leq 0}(\nbigt_{\omega}(V))_{|\Delta_z^{\ast}}$
is equal to
$q^{-1}(L_{S^1}/L^{(\omega)\,<0}_{S^1})$.
\item
 $\nbiglcheck^{\leq 0}(\nbigt_{\omega}(V))
 _{|\varpi^{-1}(0)}$
is equal to
$L^{\leq 0}_{S^1}/L^{(\omega)\,<0}_{S^1}$.
\end{itemize}
By the construction,
there exists the natural morphism:
\[
\nbigl^{\leq 0}(V)
\lrarr
\nbiglcheck^{\leq 0}(\nbigt_{\omega}(V)).
\]
Let $k:\Delta_z^{\ast}\lrarr\Deltatilde_z$
denote the inclusion.
There also exists the following exact sequence:
\begin{equation}
\label{eq;18.5.12.12}
 0\lrarr
 \nbigl^{\leq 0}\bigl(\nbigt_{\omega}(V)\bigr)
\lrarr
 \nbiglcheck^{\leq 0}(\nbigt_{\omega}(V))
\lrarr
 k_!
 k^{-1}
 q^{-1}\bigl(
 L/L^{(\omega)\leq 0}
 \bigr)
\lrarr 0.
\end{equation}

\subsection{Induced morphisms in the global case}
\label{subsection;18.5.12.20}

Let $C$ be a compact Riemann surface
with smooth boundary $\del C$.
Set $C_1:=C\setminus \del C$.
Let $D\subset C_1$ be a finite subset.
Set $C_1^{\circ}:=C_1\setminus D$.
Let $\varpi:\Ctilde\lrarr C$ 
be the oriented real blow of $C$ along $D$.
Set $\Dtilde:=\varpi^{-1}(D)$
and $\Ctilde_1:=\varpi^{-1}(C_1)$.
Let $j_1:C_1\lrarr C$
and $\jtilde_1:\Ctilde_1\lrarr \Ctilde$
denote the inclusions.
Take $Q\in D$
with a holomorphic coordinate 
neighbourhood $(C_Q,z)$ such that $z(Q)=0$.
Here, we assume that
the closure of $C_Q$ is isomorphic to
a closed disc.

\begin{df}
We say that a constructible sheaf $G$ on $\Ctilde$
is acyclic with respect to the global cohomology
if $H^{\ast}(\Ctilde,G)=0$.
 \index{acyclic with respect to the global cohomology}
\hfill\qed
\end{df}
 
Let $(V,\nabla)$ be a meromorphic flat bundle on $(C_1,D)$.
Let $\varrho_1\lrarr\varrho_2$
be a morphism in $\Dsf(D)$.
Suppose that 
$\varrho_1(Q)=!$
and $\varrho_2(Q)=\ast$.
There exist the natural monomorphisms
$\nbigl^{\varrho_1}(\nbigs^Q_{\omega}(V))
\lrarr
 \nbigl^{\varrho_1}(V)$
and 
$\nbigl^{\varrho_2}(V)
\lrarr
 \nbigl^{\varrho_2}\bigl(
 \nbigs^Q_{\omega}(V)
 \bigr)$
on $\Ctilde_1$.
Hence, 
for any $\star_1\lrarr\star_2$ in $\Dsf_1$,
there exists the following natural commutative diagram:
\begin{equation}
 \label{eq;18.5.12.1}
 \begin{CD}
  \hyperh^i\bigl(
 C, j_{1\star_1}\bigl(
 \nbigs^Q_{\omega}(V)(\varrho_1)
 \otimes\Omega^{\bullet}
 \bigr)
 \bigr)
@>>>
\hyperh^i\bigl(C,
 j_{1\star_2}\bigl(
 \nbigs_{\omega}^Q(V)(\varrho_2)\otimes\Omega^{\bullet}
 \bigr)\bigr)
\\
 @VVV @AAA\\
  \hyperh^i\bigl(
 C,j_{1\star_1}\bigl(
 V(\varrho_1)
 \otimes\Omega^{\bullet}
 \bigr)
 \bigr)
@>>>
\hyperh^i\bigl(C,
 j_{1\star_2}\bigl(
 V(\varrho_2)\otimes\Omega^{\bullet}
 \bigr)\bigr).
 \end{CD}
\end{equation}

Set $\Ctilde_Q:=\varpi^{-1}(C_Q)$.
Let $\jtilde_Q:\Ctilde_Q\lrarr \Ctilde$ 
and $j_Q:C_Q\lrarr C$ denote the inclusions.
There exist the following natural morphisms:
\begin{equation}
 \jtilde_{Q!}
 \nbiglcheck^{<0}
 \bigl(
 \nbigt_{\omega}
 (V)_{|C_Q}
 \bigr)
\lrarr
 \jtilde_{Q!}\bigl(
 \nbigl^{<0}(V)_{|C_Q}
 \bigr)
\lrarr
 \nbigl^{\varrho_1}(V).
\end{equation}
Let $q_Q:\Ctilde_Q\lrarr \varpi^{-1}(Q)$ denote the projection.
Note that
$\jtilde_{Q!}\bigl(
 q_Q^{-1}L_{S^1}^{(\omega)\,<0}
\bigr)$
is acyclic for the global cohomology.
(See \S\ref{subsection;24.3.29.1}, for example.)
Hence, by the exact sequence (\ref{eq;18.5.12.10}),
the following morphism induces
an isomorphism of the global cohomology groups:
\[
 \jtilde_{Q!}
 \nbiglcheck^{<0}
 \bigl(
 \nbigt_{\omega}
 (V)_{|C_Q}
 \bigr)
\lrarr
 \jtilde_{Q!}\nbigl^{<0}\bigl(\nbigt_{\omega}(V)_{|C_Q}\bigr).
\]
Therefore, we obtain the following morphism:
\begin{equation}
 \label{eq;18.5.12.11}
 \hyperh^i\Bigl(
C,j_{Q!}\bigl(
 \nbigt_{\omega}(V)\otimes\Omega^{\bullet}_{C_Q}
 \bigr)
 \Bigr)
\lrarr
 \hyperh^i\bigl(
 C,j_{1!}\bigl(
 V(\varrho_1)\otimes\Omega^{\bullet}
 \bigr)
 \bigr).
\end{equation}

Similarly, there exist the following natural morphisms:
\[
 \nbigl^{\varrho_2}(V)
\lrarr
 \jtilde_{Q\ast}
 \nbigl^{\leq 0}\bigl(
  V_{|C_Q}
  \bigr)
\lrarr
 \jtilde_{Q\ast}
 \nbiglcheck^{\leq 0}\bigl(
 \nbigt_{\omega}(V)_{|C_Q}
 \bigr).
\]
Let $k_Q$ denote the inclusion
$C_Q\setminus Q\lrarr \Ctilde_Q$.
Because
$\jtilde_{Q\ast}\Bigl(
 k_{Q!}k_Q^{-1}q_Q^{-1}(L/L^{(\omega)\leq 0})
 \Bigr)$
is acyclic with respect to the global cohomology,
we obtain the following morphism by
the exact sequence (\ref{eq;18.5.12.12}):
\begin{equation}
\label{eq;18.5.12.13}
 \hyperh^i\bigl(
 C,j_{1\ast}\bigl(
 V(\varrho_2)\otimes\Omega^{\bullet}
 \bigr)
 \bigr)
\lrarr
 \hyperh^i\bigl(
 C,
 j_{Q\ast}\bigl(
 \nbigt_{\omega}(V)
 \otimes\Omega^{\bullet}_{C_Q}
 \bigr)
 \bigr)
\end{equation}
Note that the following diagram is commutative:
\begin{equation}
\label{eq;18.5.13.2}
 \begin{CD}
\hyperh^i\Bigl(
C,j_{Q!}\bigl(
 \nbigt_{\omega}(V)\otimes\Omega^{\bullet}_{C_Q}
 \bigr)
 \Bigr)
@>>>
 \hyperh^i\Bigl(
 C,
 j_{Q\ast}\bigl(
 \nbigt_{\omega}(V)
 \otimes\Omega^{\bullet}_{C_Q}
 \bigr)
 \Bigr)
\\
@VVV @AAA\\
 \hyperh^i\bigl(
 C,j_{1!}\bigl(
 V(\varrho_1)\otimes\Omega^{\bullet}
 \bigr)
 \bigr)
@>>>
\hyperh^i\bigl(
 C,j_{1\ast}\bigl(
 V(\varrho_2)\otimes\Omega^{\bullet}
 \bigr)
 \bigr).
 \end{CD}
\end{equation}

\subsection{Complement}

Let $U:=\{z\in\cnum\,|\,|z|<1\}$.
Let $(V_1,\nabla)$ be a meromorphic flat bundle
on $(U,0)$.
We extend it to a meromorphic flat bundle
on $(\proj^1,\{0,\infty\})$
with regular singularity at $\infty$,
which we denote by $(V_2,\nabla)$.

Set $D:=\{0,\infty\}$.
Let $\varpi:\projtilde^1\lrarr \proj^1$
denote the oriented real blow up
of $\proj^1$ along $D$.
Let $\Utilde:=\varpi^{-1}(U)$.
Let $\jtilde:\Utilde\lrarr\projtilde^1$
denote the inclusion.
We shall use the following natural isomorphisms:
\[
 \hyperh^i\bigl(
 \proj^1,V_2(!D)\otimes\Omega^{\bullet}
 \bigr)
\simeq
 H^i\bigl(
 \projtilde^1,
 \jtilde_{!}
 \nbigl^{<0}(V_1)
 \bigr),
\]
\[
 \hyperh^i\bigl(
 \proj^1,V_2(\ast D)\otimes\Omega^{\bullet}
 \bigr)
\simeq
 H^i\bigl(
 \projtilde^1,
 \jtilde_{\ast}
 \nbigl^{\leq 0}(V_1)
 \bigr).
\]

\section{Homology groups of meromorphic flat bundles}
\label{subsection;20.10.21.1}

\subsection{Homology groups with coefficient of constructible sheaves}
\label{subsection;20.9.8.1}

Let $Y$ be a differentiable manifold 
with boundary $\del Y$.
We assume that $Y$ is oriented.
Set $Y^{\circ}:=Y\setminus \del Y$.
\index{open subset $Y^{\circ}$}

Let $H$ be a closed subspace of $Y$.
For any open subset $U\subset Y$,
let $S_p(Y,(Y\setminus U)\cup H;\cnum)$
denote the group of piece-wise smooth $p$-chains of
$Y$ relative to $(Y\setminus U)\cup H$
with the $\cnum$-coefficient.
It induces a presheaf on $Y$.
Let $\nbigc_{Y,H}^{-p}$ denote the sheafification.
\index{sheaf $\nbigc_{Y,H}^{-p}$}
The boundary homomorphisms of the chain groups induce
$\del:\nbigc_{Y,H}^{-p}\lrarr \nbigc_{Y,H}^{-p+1}$
with which 
$\nbigc_{Y,H}^{\bullet}$ is a complex of sheaves.
\index{complex of sheaves $\nbigc_{Y,H}^{\bullet}$}
If $H=\emptyset$,
we denote it by $\nbigc_Y^{\bullet}$.
\index{complex of sheaves $\nbigc_{Y}^{\bullet}$}

Let $\nbigg$ be 
any $\real$-constructible $\cnum_Y$-module.
As mentioned in \cite{Hien},
$\nbigg\otimes\nbigc_{Y,H}^{\bullet}$
is homotopically fine
(see \cite{Bredon})
so that we may compute the hypercohomology group
$\hyperh^{\ast}\bigl(
 Y,\nbigg\otimes\nbigc_{Y,H}^{\bullet}
\bigr)$
by taking the global sections.

Because $Y$ is an oriented manifold with boundary,
there exists a natural isomorphism
$\nbigc_{Y,\del Y}^{\bullet}\simeq
 \cnum_Y[\dim Y]$
in the derived category $D^b_{\realc}(\cnum_Y)$
of cohomologically $\real$-constructible complexes.
Let $\iota_{Y^{\circ}}:Y^{\circ}\lrarr Y$
denote the inclusion.
Then, there exists a natural isomorphism
$\nbigc^{\bullet}_{Y}\simeq
\iota_{Y^{\circ}!}
\cnum_{Y^{\circ}}[\dim Y]$
in $D^b_{\realc}(\cnum_Y)$.
Hence, for any $\real$-constructible $\cnum_Y$-modules
$\nbigg$,
there exist the following natural isomorphisms
in $D^b_{\realc}(\cnum_Y)$:
\[
 \nbigg\otimes\nbigc_{Y,\del Y}^{\bullet}
\simeq \nbigg[\dim Y],
\quad\quad
 \nbigg\otimes\nbigc_{Y}^{\bullet}
\simeq 
 \nbigg\otimes
 \iota_{Y^{\circ}!}\cnum_{Y^{\circ}}[\dim Y]
=\iota_{Y^{\circ}!}\iota_{Y^{\circ}}^{-1}(\nbigg)[\dim Y].
\]

\subsection{Homology groups of meromorphic flat bundles}
\label{subsection;18.5.14.102}

The notion of rapid decay homology group
for meromorphic flat bundles was 
introduced by Bloch-Esnault \cite{Bloch-Esnault2}
in the one dimensional case,
and by Hien \cite{Hien} in the general case.
We recall the definition in the one dimensional case
by following \cite{Hien}.

Let $C$ be a compact complex curve which may have 
smooth boundary $\del C$.
We set
$C_1:=C\setminus\del C$.
Let $D\subset C_1$ be a finite subset.
Set $C^{\circ}_1:=C_1\setminus D$.
Let 
$\varpibar:
 \Ctilde\lrarr C$
and $\varpi:\Ctilde_1\lrarr C_1$ denote 
the oriented real blow up along $D$.
Set $\Dtilde:=\varpi^{-1}(D)$.
The boundary $\del\Ctilde$ of $\Ctilde$
is $\Dtilde\cup \del C$.
Let $j_1:C_1\lrarr C$ 
and $\jtilde_1:\Ctilde_1\lrarr\Ctilde$
denote the inclusions.

Let $(V,\nabla)$ be a meromorphic flat bundle on $(C_1,D)$.
Let $(\nbigl(V),\vecnbigf)$
denote the associated local system with Stokes structure on 
$(\Ctilde_1,\Dtilde)$.
As explained in \S\ref{subsection;18.5.13.1},
we obtain the associated $\real$-constructible sheaves
$\nbigl^{<0}(V)$
and 
$\nbigl^{\leq 0}(V)$
on $\Ctilde_1$.
The sheaf of rapid decay chains of $(V,\nabla)$ 
on $\Ctilde$
is defined as follows:
\index{sheaf of rapid decay chains $\nbigc^{\rd,\bullet}_{\Ctilde}(V)$}
\[
 \nbigc^{\rd,\bullet}_{\Ctilde}(V):=
 \nbigc^{\bullet}_{\Ctilde,\del\Ctilde}\otimes
 \jtilde_{1!}
 \nbigl^{<0}(V).
\]
The rapid decay $p$-th homology group of $(V,\nabla)$ 
is defined as follows:
\index{rapid decay homology group $H_p^{\rd}(C_1^{\circ},V)$}
\[
 H_p^{\rd}\bigl(
 C_1^{\circ},V
 \bigr):=
 \hyperh^{-p}\bigl(
 \Ctilde,\nbigc^{\rd,\bullet}_{\Ctilde}(V)
 \bigr).
\]
If $(V,\nabla)$ is regular singular,
$H_p^{\rd}\bigl(C_1^{\circ},V\bigr)$
equals the $p$-th homology group of $C_1^{\circ}$
with coefficient
$\nbigl(V)_{|C_1^{\circ}}$.

It is also natural to consider the homology groups
associated with $\nbigl^{\leq 0}(V)$.
The sheaf of moderate growth chains of $(V,\nabla)$
on $\Ctilde$ 
is defined as follows:
\index{sheaf of moderate growth chains
$\nbigc^{\mg,\bullet}_{\Ctilde}(V)$}
\[
 \nbigc^{\mg,\bullet}_{\Ctilde}(V)
:=\nbigc^{\bullet}_{\Ctilde,\del\Ctilde}
 \otimes
 \jtilde_{1\ast}\nbigl^{\leq 0} (V).
\]
The moderate growth homology group 
of $(V,\nabla)$ is defined as
$H_p^{\mg}(C_1^{\circ},V)
=\hyperh^{-p}\bigl(
 \Ctilde,\nbigc^{\mg,\bullet}_{\Ctilde}(V)
\bigr)$.
\index{moderate growth homology group
$H_p^{\mg}(C_1^{\circ},V)$}
 
For any $\varrho\in \Dsf(D)$
and for any $\star\in\{!,\ast\}$,
we define
the sheaf of 
$(\varrho,\star)$-type chains of $(V,\nabla)$ as 
\index{sheaf of $(\varrho,\star)$-type chains
$\nbigc^{(\varrho,\star),\bullet}_{\Ctilde}(V)$}
\[
 \nbigc^{(\varrho,\star),\bullet}_{\Ctilde}(V):=
 \nbigc^{\bullet}_{\Ctilde,\del\Ctilde}
 \otimes
 j_{1\star}
 \nbigl^{\varrho}(V).
\]
We define
the $(\varrho,\star)$-type homology group of $(V,\nabla)$ as
\index{$(\varrho,\star)$-type homology group
$H_p^{\varrho,\star}\bigl(
 C^{\circ},V\bigr)$}
\[
 H_p^{\varrho,\star}\bigl(
 C^{\circ},V\bigr):=
 \hyperh^{-p}\bigl(
 \Ctilde,\nbigc^{(\varrho,\star),\bullet}_{\Ctilde}(V)
 \bigr).
\]
If $\del C=\emptyset$,
it is denoted by
$H_p^{\varrho}\bigl(
 C^{\circ},V\bigr)$.
\index{homology group $H_p^{\varrho}\bigl(
 C^{\circ},V \bigr)$.}
By definition,
we have
\[
H_p^{\barshriek,!}\bigl(
 C^{\circ},V\bigr)
=H_p^{\rd}\bigl(
 C_1^{\circ},V\bigr),
\quad
 H_p^{\barast,\ast}\bigl(
 C^{\circ},V\bigr)
=H_p^{\mg}\bigl(
 C_1^{\circ},V\bigr).
\]

\begin{lem}
For any morphism
$\varrho_1\lrarr\varrho_2$ in $\Dsf(D)$
and $\star_1\lrarr\star_2$ in $\Dsf_1$,
there exists the following natural commutative diagram:
\[
 \begin{CD}
 H_p^{\varrho_1,\star_1}(C^{\circ},V)
 @>>>
 H_p^{\varrho_2,\star_2}(C^{\circ},V)
 \\
 @V{\simeq}VV @V{\simeq}VV \\
 \hyperh^{2-p}\bigl(C,
 j_{1\star_1}(V(\varrho_1)\otimes\Omega^{\bullet}_C)\bigr)
 @>>>
 \hyperh^{2-p}\bigl(C,
 j_{1\star_2}(
 V(\varrho_2)\otimes\Omega^{\bullet}_C)
 \bigr).
 \end{CD}
\]
\end{lem}
\pf
The vertical isomorphisms are induced by
the natural isomorphisms
$V(\varrho)\otimes\Omega^{\bullet}_{C_1}
\simeq
R\varpi_{\ast}
 \nbigl^{\varrho}(V)$
 in $D^b_{\realc}(\cnum_{C_1})$.
\hfill\qed

\subsection{Some general morphisms}
\label{subsection;18.5.15.40}

Let $C$ be a compact Riemann surface.
We assume $\del C=\emptyset$.
Let $D\subset C$
be a finite subset.
Set $C^{\circ}:=C\setminus D$.
Let $\varpi:\Ctilde\lrarr C$
be the oriented real blow of $C$ along $D$.
Let $(V,\nabla)$ be a meromorphic flat bundle on $(C,D)$.
We translate the morphisms in \S\ref{subsection;18.5.12.20}
to the context of homology groups.
Take $Q\in D$.
Let $\varrho_1\lrarr\varrho_2$ be a morphism
in $\Dsf(D)$
such that $\varrho_1(Q)=!$
and $\varrho_2(Q)=\ast$.

\subsubsection{}
There exists the natural monomorphism
$\nbigl^{\varrho_1}(\nbigs^Q_{\omega}(V))
\lrarr
 \nbigl^{\varrho_1}(V)$
on $\Ctilde$.
Let 
$\nbigq^{<0}_{Q,\omega}(V)$
denote the quotient sheaf
whose support is contained in $\varpi^{-1}(Q)$.
Note that
$\nbigc^{\bullet}_{\Ctilde}\otimes
 \nbigq^{<0}_{Q,\omega}(V)
\simeq
 \iota_{1!}\iota_{1}^{-1}
 \nbigq^{<0}_{Q,\omega}(V)=0$,
where $\iota_1$
denotes the inclusion
$\Ctilde\setminus\varpi^{-1}(Q)\lrarr \Ctilde$.
Hence, we obtain the following exact sequence
\begin{multline}
 \label{eq;18.4.17.10}
\hyperh^{-1}\bigl(
 \varpi^{-1}(Q),
 \nbigc^{\bullet}_{\varpi^{-1}(Q)}\otimes\nbigq^{<0}_{Q,\omega}(V)
 \bigr)
\lrarr
 H_1^{\varrho_1}\bigl(C^{\circ},\nbigs_{\omega}^Q(V)\bigr)
\lrarr
 H_1^{\varrho_1}\bigl(C^{\circ},V\bigr)
\lrarr 
\\
\hyperh^0\bigl(
 \varpi^{-1}(Q),
 \nbigc^{\bullet}_{\varpi^{-1}(Q)}\otimes\nbigq^{<0}_{Q,\omega}(V)
 \bigr)
\lrarr 
 H_0^{\varrho_1}\bigl(C^{\circ},\nbigs_{\omega}^Q(V)\bigr)
\lrarr
 H_0^{\varrho_1}\bigl(C^{\circ},V\bigr).
\end{multline}
We also remark that
$\nbigq^{<0}_{Q,\omega}\bigl(\nbigt_{\omega}(V)\bigr)
\simeq
 \nbigq^{<0}_{Q,\omega}\bigl(V\bigr)$
naturally.

Similarly, there exists the natural morphism
$\nbigl^{\varrho_2}(V)
\lrarr
 \nbigl^{\varrho_2}\bigl(
 \nbigs^Q_{\omega}(V)
 \bigr)$
on $\Ctilde$.
Let 
$\nbigq^{\leq 0}_{Q,\omega}(V)$
denote the quotient sheaf
whose support is 
contained in $\varpi^{-1}(Q)$.
We obtain the following exact sequence:
\begin{multline}
\hyperh^{-1}\bigl(
 \varpi^{-1}(Q),
 \nbigc^{\bullet}_{\varpi^{-1}(Q)}
 \otimes\nbigq^{\leq 0}_{Q,\omega}(V)
 \bigr)
\lrarr
 H_1^{\varrho_2}(C^{\circ},V)
\lrarr
 H_1^{\varrho_2}(C^{\circ},
 \nbigs^{Q}_{\omega}(V))
\lrarr
 \\
 \hyperh^0\bigl(
 \varpi^{-1}(Q),
 \nbigc^{\bullet}_{\varpi^{-1}(Q)}
 \otimes\nbigq^{\leq 0}_{Q,\omega}(V)
 \bigr)
\lrarr
 H_0^{\varrho_2}(C^{\circ},V)\lrarr
 H_0^{\varrho_2}(C^{\circ},
 \nbigs^{Q}_{\omega}(V)).
\end{multline}
Note that
$\nbigq^{\leq 0}_{Q,\omega}\bigl(\nbigt_{\omega}(V)\bigr)
\simeq
 \nbigq^{\leq 0}_{Q,\omega}\bigl(V\bigr)$
naturally.
The commutative diagram (\ref{eq;18.5.12.1})
is identified with the following diagram:
\begin{equation}
\label{eq;18.5.13.4}
 \begin{CD}
 H_1^{\varrho_1}(C^{\circ},
  \nbigs^{Q}_{\omega}(V))
 @>>>
 H_1^{\varrho_2}(C^{\circ},
  \nbigs^{Q}_{\omega}(V))\\
 @VVV @AAA \\
 H_1^{\varrho_1}(C^{\circ},V)
 @>>>
 H_1^{\varrho_2}(C^{\circ},V).
 \end{CD}
\end{equation}

\subsubsection{}
Set $\Ctilde_Q:=\varpi^{-1}(C_Q)$.
Set $C_Q^{\circ}:=C_Q\setminus\{Q\}$.
Recall $\varrho_1(Q)=!\in\Dsf(\{Q\})$.
As the translation of (\ref{eq;18.5.12.11}),
we obtain the following morphism:
\begin{equation}
\label{eq;18.5.13.22}
 H_1^{\varrho_{1}(Q),!}\bigl(C_Q^{\circ},\nbigt_{\omega}(V)\bigr)
\lrarr
 H_1^{\varrho_1}\bigl(C^{\circ},V\bigr).
\end{equation}
Here, 
$H_1^{\varrho_{1}(Q),!}\bigl(C_Q^{\circ},\nbigt_{\omega}(V)\bigr)$
is obtained from
the restriction of $\nbigt_{\omega}(V)$
to the closure of $C_Q$.
Recall $\varrho_2(Q)=\ast\in\Dsf(\{Q\})$.
As the translation of (\ref{eq;18.5.12.13}),
we obtain the following morphism:
\[
 H_1^{\varrho_2}\bigl(
 C^{\circ},V \bigr)
\lrarr
 H_1^{\varrho_{2}(Q),\ast}\bigl(
 C_Q^{\circ},
 \nbigt_{\omega}(V)
 \bigr).
\]
The commutative diagram (\ref{eq;18.5.13.2})
is identified with the following diagram:
\begin{equation}
\label{eq;18.5.13.3}
 \begin{CD}
 H_1^{\varrho_{1}(Q),!}\bigl(
 C_Q^{\circ},
 \nbigt_{\omega}(V)
 \bigr)
 @>>> 
 H_1^{\varrho_{2}(Q),\ast}\bigl(
 C_Q^{\circ},
 \nbigt_{\omega}(V)
 \bigr)
\\
 @VVV @AAA \\
 H_1^{\varrho_1}\bigl(
 C^{\circ},V \bigr)
 @>>>
 H_1^{\varrho_2}\bigl(
 C^{\circ},V \bigr).
 \end{CD}
\end{equation}

\subsection{Vanishing}
\label{subsection;24.3.29.1}

Let $(L,\vecnbigf)=\RH(V,\nabla)
\in\Loc^{\St}(\nbigi)$.
Let $\nbigl$  be the local system on $\Ctilde_Q$
induced by $L$.

Take $R_1<R_2$
and $0<\epsilon_1<\epsilon_2$.
We set $Y_0:=\closedopen{0}{\epsilon_2}\times \real$,
$Y_1:=\closedopen{0}{\epsilon_1}\times\openopen{R_1}{R_2}$,
$Y_2:=\openopen{0}{\epsilon_1}\times\openopen{R_1}{R_2}$
and $Y_3:=\{0\}\times \openopen{R_1}{R_2}$.
Let $q_{Y_i}:Y_i\lrarr\real$ denote the map
induced by the projection.
Let $j_{Y_i}:Y_i\lrarr Y_0$ denote the inclusion.
If $\epsilon_2$ is sufficiently small,
we obtain the map
$\varphi_Q:Y_0\lrarr \Ctilde_Q$
defined by
$\varphi_Q(r,\theta)=(r,e^{\sqrt{-1}\theta})$.

\begin{lem}
\label{lem;20.8.27.1}
 We obtain
 $H^{a}(Y_0,j_{Y_1!}q_{Y_1}^{-1}(L^{<0}))=0$
 for any $a$.
\end{lem}
\pf
For any open interval $I\subset \openopen{R_1}{R_2}$,
let $\iota_{I}:I\lrarr\real$
and $\iota_{\closedopen{0}{\epsilon_j}\times I}:
\closedopen{0}{\epsilon_j}\times I\lrarr Y_0$ $(j=1,2)$
denote the inclusions.
If $I$ satisfies $|\Ibar\cap S_0(\nbigi)|\leq 1$,
there exist open intervals $I_1,\ldots,I_N$ of $I$
such that
\[
 \iota_{I!}\iota_I^{-1}(L^{<0})
 \simeq
 \bigoplus_{k=1}^N
 \iota_{I_k!}(\cnum_{I_k}).
\]
We obtain
\[
 \iota_{\closedopen{0}{\epsilon_2}\times I!}
 \iota_{\closedopen{0}{\epsilon_2}\times I}^{-1}
 \bigl(
 j_{Y_1!}j_{Y_1}^{-1}q_{Y_0}^{-1}(L^{<0})
 \bigr)
 \simeq
 \bigoplus_{k=1}^N
 \iota_{\closedopen{0}{\epsilon_1}\times I_k!}
 \cnum_{\closedopen{0}{\epsilon_1}\times I_k}.
\]
Hence, we obtain
$H^{a}\Bigl(
 Y_0, \iota_{\closedopen{0}{\epsilon_2}\times I!}
 \iota_{\closedopen{0}{\epsilon_2}\times I}^{-1}
 \bigl(
 j_{Y_1!}q_{Y_1}^{-1}(L^{<0})
 \bigr)
\Bigr)=0$ for any $a$.
Then, we obtain the claim of the lemma easily.
\hfill\qed

\vspace{.1in}

There exists the $\real$-constructible subsheaf
$N_{Y_1}\subset j_{Y_1}^{-1}\varphi_Q^{-1}(\nbigl)$
determined by
$N_{Y_1|Y_3}=L^{<0|Y_3}$
and 
$N_{Y_1|Y_2}=
 q^{-1}_{Y_2}(L^{\leq 0})$.
We shall use the following lemma
for our computations.

\begin{lem}
 \label{lem;18.5.25.30}
We obtain 
$H^{1}\bigl(Y_0,j_{Y_1!}N_{Y_1}\bigr)=0$.
In particular,
any $1$-cycle for $j_{Y_1!}N_{Y_1}$
is $0$ in the homology level.
\end{lem}
\pf
There exists the natural exact sequence
\begin{equation}
\label{eq;20.8.27.2}
 0\lrarr
 j_{Y_1!}q_{Y_1}^{-1}(L^{<0})
 \lrarr
 j_{Y_1!}N_{Y_1}
 \lrarr
 j_{Y_2!}q_{Y_2}^{-1}(L^{\leq 0}/L^{<0})
 \lrarr 0.
\end{equation}
Because
$q_{Y_2}^{-1}(L^{\leq 0}/L^{<0})$
is a local system on $Y_2$,
we obtain
\[
 H^1\bigl(Y_0,
 j_{Y_2!}q_{Y_2}^{-1}(L^{\leq 0}/L^{<0})
 \bigr)=0.
\]
We obtain the claim of the lemma
from Lemma \ref{lem;20.8.27.1}
and the exact sequence (\ref{eq;20.8.27.2}).
\hfill\qed

\vspace{.1in}
Let us give a variant.
Take $0<\delta<R_2-R_1$,
and we set $I:=\openopen{R_1}{R_1+\delta}$.
Let $L_{I,0}\subset L_{|I}$ be a local subsystem
such that
$L_{I,0|\theta}\subset \nbigf^{\theta}_0$
for any $\theta\in I$,
and $L_{I,0}\simeq \Gr^{\vecnbigf}_0(L)_{|I}$.
We set $Y_4:=\{\epsilon_1\}\times\openopen{R_1}{R_1+\delta}$
and
$Y_5:=Y_1\cup Y_4$.
For $i=4,5$,
let
$q_{Y_i}\lrarr \real$ denote
the maps induced by the projection,
and let $j_{Y_i}:Y_i\lrarr Y_0$
denote the inclusions.
There exists the constructible subsheaf
$N_{Y_5}\subset
 j_{Y_5}^{-1}\varphi_Q^{-1}(\nbigl)$
determined by
$N_{Y_5|Y_1}=N_{Y_1}$
and
$N_{Y_5|Y_4}=q_{Y_4}^{-1}(L_{I,0})$.

\begin{lem}
We obtain 
$H^{a}\bigl(Y_0,j_{Y_5!}N_{Y_5}\bigr)=0$
for any $a$.
In particular,
any $1$-cycle for $j_{Y_5!}N_{Y_5}$
is $0$ in the homology level. 
\end{lem}
\pf
The quotient sheaf
$j_{Y_5!}N_{Y_5}\big/
j_{Y_1!}q_{Y_1}^{-1}(L^{<0})$
is also acyclic with respect to
the global cohomology.
Hence, the claim of the lemma follows.
\hfill\qed

\subsection{Some computations}

Let $J\subset\real$ be an open interval.
The inclusion $J\lrarr\real$ is denoted by $\iota_J$.
The following lemma is easy to see.
Let $L$ be a local system on $J$.
\begin{lem}
\label{lem;20.10.14.1}
By definition, we obtain 
$\hyperh^{i}(\real,\nbigc^{\bullet}_{\real}\otimes\iota_{J!}L)=0$
unless $i=0$,
and
\[
 \hyperh^{0}(\real,\nbigc^{\bullet}_{\real}\otimes\iota_{J!}L)
 =H_0(J,L).
\] 
Here, $H_0(J,L)$ denotes the $0$-th homology group
with $L$-coefficient.
We also obtain
$\hyperh^{i}(\real,\nbigc^{\bullet}_{\real}\otimes\iota_{J\ast}L)=0$
unless $i=-1$, and 
\[
 \hyperh^{-1}(\real,\nbigc^{\bullet}_{\real}\otimes\iota_{J\ast}L)
 \simeq H^{-1}(\real,\iota_{J\ast}L[1])
 =H^0(J,L),
\]
which depends on the orientation of $\real$.
\hfill\qed
\end{lem}

The natural orientation of $\real$ induces
\[
 \phi_{\real}:
 H^1\Bigl(
 \real,\iota_{J!}(L)
 \Bigr)
\simeq
 \hyperh^0(\real,\iota_{J!}(L)\otimes\nbigc^{\bullet}_{\real})
=H_0(J,L).
\]
Let $\rho:H^0(J,L)\simeq H_0(J,L)$
be the isomorphism induced by
$v\longmapsto v\otimes[x]$
for any $x\in J$,
where $[x]$ denotes the natural $0$-chain induced by $x$.
Let $\Omega_J^{\bullet}$ denote the sheaf of
$C^{\infty}$-differential forms on $J$.
There exists the natural isomorphism
$H^0(\real,\iota_{J!}L)\simeq
H^0\bigl(\real,\iota_{J!}(L\otimes\Omega^{\bullet}_J)\bigr)$.
By the integration,
we obtain
\[
\int_{\real}:
 H^1(\real,\iota_{J!}L)
 =H^1\bigl(\real,\iota_{J!}(L\otimes\Omega^{\bullet}_J)\bigr)
 \simeq
 H^0(J,L).
\]
\begin{lem}
\label{lem;24.2.14.2}
We have $\rho\circ\int_{\real}=-\phi_{\real}$.
\end{lem}
\pf
It is enough to study the case
$J=\openopen{0}{1}$ and $L=\cnum_J$.
Take $x_0<0<x_1<1<x_2$.
Let us consider the double complex of sheaves
$\iota_{J!}\Omega^{\bullet}_J\otimes
\nbigc^{\bullet}_{\real}[1]$.
Let $\del$ denote the differential of the total complex.
For $a<b$, let $[a,b]$ denote the $1$-chain of $\real$
induced by the natural inclusion.
For $a\in\real$, let $[a]$ denote the $0$-chain of $\real$
induced by $a$.
Let $\omega$ be a section of
$\iota_{J!}\Omega^{1}_J$,
i.e.,
a $1$-form whose support is contained in $J$.
We set
\[
 f_0=\int_{x_0}^x\omega,
 \quad
 f_1=-\int_{x}^{x_2}\omega.
\]
We obtain
\[
 \del\Bigl(
 f_0\otimes[x_0,x_1]
+f_1\otimes[x_1,x_2]
 \Bigr)
=\omega\otimes([x_0,x_1]+[x_1,x_2])
+\Bigl(
\int_{\real}\omega\Bigr)\otimes[x_1].
\]
It implies the claim of the lemma.
\hfill\qed

\vspace{.1in}

We set $I=\closedclosed{0}{1}$.
Let $\iotatilde_{J}:I\times J\to I\times \real$
denote the inclusion.
Let $q_J:I\times J\to J$ denote the projection.
For $a=0,1$,
let $k_a:\real\simeq\{a\}\times\real\to I\times\real$
denote the inclusions.
There exists the natural morphisms
\[
 \iotatilde_{J!}q_J^{-1}(L)
 \otimes
 \nbigc^{\bullet}_{I\times\real,\del I\times\real}
 \lrarr
 k_{a!}\Bigl(
 \iota_{J!}(L)
 \otimes
 \nbigc^{\bullet}_{\real}
 \Bigr)[1].
\]
They induce
\[
\del_a:
 \hyperh^{-1}\Bigl(
 I\times\real,
 \iotatilde_{J!}q_J^{-1}(L)
 \otimes
 \nbigc^{\bullet}_{I\times\real,\del I\times\real}
 \Bigr)
 \lrarr
 \hyperh^0(\real,\iota_{J!}(L)\otimes\nbigc^{\bullet}_{\real}).
\]
By the natural orientations of $I$ and $\real$,
we obtain
\[
 \phi_{I\times\real}:
 H^{1}\Bigl(
 I\times\real,
 \iotatilde_{J!}q_J^{-1}(L)
 \Bigr)
 \simeq
 \hyperh^{-1}\Bigl(
 I\times\real,
 \iotatilde_{J!}q_J^{-1}(L)
 \otimes
 \nbigc^{\bullet}_{I\times\real,\del I\times\real}
 \Bigr),
\]
Let $q_{\real}:I\times\real\to\real$ denote the projection.
There exists the natural isomorphism
\[
 q_{\real}^{\ast}:
 H^1\Bigl(
 \real,\iota_{J!}(L)
 \Bigr)
 \simeq
  H^{1}\Bigl(
 I\times\real,
 \iotatilde_{J!}q_J^{-1}(L)
 \Bigr).
\]

\begin{lem}
\label{lem;24.2.14.1}
$\del_a\circ \phi_{I\times\real}\circ q_{\real}^{\ast}
=-(-1)^a \phi_{\real}$. 
\end{lem}
\pf
By taking a simplicial decomposition
$I\times\Jbar=\bigcup\alpha_i$,
we construct a relative $2$-cycle
$[I\times\Jbar]=\sum\alpha_i$ of
$(I\times\Jbar,\del(I\times\Jbar))$
representing the fundamental class in
$H_2(I\times\Jbar,\del(I\times\Jbar);\seisuu)$.
It induces relative $1$-cycles
$\del_a[I\times\Jbar]$ $(a=0,1)$ of
$(\{a\}\times\Jbar,\{a\}\times\del\Jbar)$.
Note that
$-(-1)^a\del_a[I\times\Jbar]$
are the fundamental class of $H_1(\Jbar,\del\Jbar)$.

Let $\Phi_{I\times\real}:
\iota_{\Jtilde!}q_{J}^{-1}L\to
\iota_{\Jtilde!}q_{J}^{-1}L\otimes
\nbigc^{\bullet}_{I\times\real,\del_I\times\real}[-2]$
be the morphism induced by
$s\longmapsto s\otimes[I\times\Jbar]$.
It induces $\phi_{I\times\real}$.
The morphism
$\del_a\circ\phi_{I\times\real}\circ q_{\real}^{\ast}$
is induced by the composition $\Phi_{\real,a}$ of 
the following morphisms of complexes of sheaves:
\begin{multline}
 L\simeq
 q_{\real!}
 (\iotatilde_{J!}q_J^{-1}(L))
 \stackrel{q_{\real!}\Phi_{I\times\real}}{\lrarr}
 q_{\real!}
 \bigl(
 \iotatilde_{J!}q_J^{-1}L\otimes
 \nbigc^{\bullet}_{I\times\real,\del I\times\real}[-2]
 \bigr)
 \stackrel{q_{\real!}(\del_a)}{\lrarr}
 \\
 q_{\real!}
 \bigl(
 k_{a!}
 \iota_{J!}L\otimes
 \nbigc^{\bullet}_{\real}[-1]
 \bigr)
 =\iota_{J!}L\otimes
 \nbigc^{\bullet}_{\real}[-1].
\end{multline}
It equals the morphism induced by
$s\longmapsto s\otimes \del_a[I\times\Jbar]$.
It induces $-(-1)^a\phi_{\real}$
in the cohomology level.
\hfill\qed

\vspace{.1in}
Let $x_1\in J$.
Let $I\times[x_1]$ denote the relative
$1$-cycle of $(I,\del I)\times\real$
induced by the inclusion of
$I\times\{x_1\}$.
For $v\in H^0(J,L)$,
let $v\otimes (I\times[x_1])$ denote the induced section
of $L\otimes\nbigc^{-1}_{(I,\del I)\times\real}$,
which is a relative $1$-cycle.

\begin{cor}
$\int_{\real}(q_{\real}^{\ast})^{-1}
 (\phi_{\real\times I})^{-1}
 \bigl(
 v\otimes (I\times[x_1])
 \bigr)=-v$.
\end{cor}
\pf
We have
$\rho^{-1}\bigl(\del_1(v\otimes(I\times[x_1]))\bigr)
=v$.
By Lemma \ref{lem;24.2.14.1},
we obtain
\[
 \rho^{-1}\circ
 \phi_{\real}\circ
 (q_{\real}^{\ast})^{-1}\circ
 (\phi_{\real\times I})^{-1}
 \bigl(v\otimes([x_1]\otimes I)\bigr)
 =\rho_1\circ\del_1
  \bigl(v\otimes([x_1]\otimes I)\bigr)
=v.
\]
Because
$\rho^{-1}\circ\phi_{\real}=-\int_{\real}$
by Lemma \ref{lem;24.2.14.2},
we obtain the claim of the corollary.
\hfill\qed

\section{Fourier transforms and some induced maps}

\subsection{Fourier transforms}
\label{subsection;25.2.12.20}

Let $\nbigm$ be a coherent algebraic $\nbigd$-module on $\cnum_z$.
It is equivalent to a finitely generated
$\cnum[z]\langle\del_z\rangle$-module $M$.
We set $\Fourier_{\pm}(M):=M$ as $\cnum$-vector spaces.
We obtain the $\cnum[w]\langle\del_w\rangle$-modules
$\Fourier_{\pm}(M)$ by setting
$w\cdot m=\mp\del_zm$ and $\del_wm=\pm zm$.
We obtain the corresponding algebraic $\nbigd$-modules
$\Fourier_{\pm}(\nbigm)$ on $\cnum_w$.
\index{Fourier transforms $\Fourier_{\pm}(\nbigm)$}

Recall that they are also obtained as 
the integral transforms of $\nbigd$-modules.
Set
 $H:=(\proj^1_z\times\{\infty\})
 \cup(\{\infty\}\times\proj^1_w)$.
We set 
$\nbige(\pm zw):=\nbigo_{\proj^1_z\times\proj^1_w}(\ast H)$
with the connection given by
$d\pm d(zw)$.
Let $p_z:\proj^1_z\times\proj^1_w\lrarr \proj^1_z$
and $p_w:\proj^1_z\times\proj^1_w\lrarr\proj^1_w$
be the projections.
We obtain the $\nbigd$-modules
$p_z^{\ast}(\nbigm)\otimes
 \nbige(zw)$
on $\proj^1_z\times\proj^1_w$.
Then,
it is well known
and easy to check that
there exist natural isomorphisms:
\[
 \Fourier_{\pm}(\nbigm)
\simeq
 p^0_{w+}
\Bigl(
 p_z^{\ast}(\nbigm)\otimes
 \nbige(\pm zw)
 \Bigr)
 \simeq
 p_{w+}
\Bigl(
 p_z^{\ast}(\nbigm)\otimes
 \nbige(\pm zw)
 \Bigr).
\]

For any algebraic holonomic $\nbigd_{\cnum}$-module $\nbigm$,
let $\DDD(\nbigm)$ denote the dual holonomic $\nbigd_{\cnum}$-module.
Then,
there exists an isomorphism
$\Fourier_{\pm}(\DDD \nbigm)\simeq
 \DDD\Fourier_{\mp}(\nbigm)$.
In particular, 
for a meromorphic flat bundle $\nbigv$ on $(\proj^1,D\cup\{\infty\})$,
we naturally obtain
$\Fourier_+(\nbigv)\simeq
 \Fourier_-(\nbigv^{\lor}(!D))^{\lor}$
and 
$\Fourier_+(\nbigv(!D))\simeq
 \Fourier_-(\nbigv^{\lor})^{\lor}$
on a neighbourhood of $\infty$.

\subsection{Local systems with Stokes structure at $\infty$}

Let $D$ be a finite subset in $\cnum$.
Let $(\nbigv,\nabla)$ be a meromorphic flat bundle
on $(\proj^1_z,D\cup\{\infty\})$.
There exists a neighbourhood $U_{w,\infty}$ of $\infty$
in $\proj^1_w$
such that 
$\Fourier_+\bigl(\nbigv(\varrho)\bigr)_{|U_{w,\infty}\setminus \{\infty\}}$
are flat bundles
for any $\varrho\in\Dsf(D)$.
\begin{notation}
Let $(\gbigl^{\gbigf}_{\varrho}(\nbigv),\vecnbigf)$
denote the $2\pi\seisuu$-equivariant
local system with Stokes structure on $\real$
associated with $\Fourier_+(\nbigv(\varrho ))$.
\index{local system with Stokes structure
$(\gbigl^{\gbigf}_{\varrho}(\nbigv),\vecnbigf)$}
\hfill\qed
\end{notation}

Note that the polar decomposition $u=w^{-1}=|u|e^{\sqrt{-1}\theta^u}$
induces a coordinate $\theta^u$ of $\real$.

\subsection{Parallel transports}
\label{subsection;24.3.14.42}

We set $\Dtilde=D\cup\{\infty\}$.
On $U_{w,\infty}$,
we use the coordinate $u=w^{-1}$.
We consider
$u_i=|u_i|e^{\sqrt{-1}\theta^u_i}
\in U_{w,\infty}\setminus\{\infty\}$ $(i=1,2)$.
For a path connecting $u_1$ and $u_2$ in $U_{w,\infty}\setminus\{\infty\}$,
we obtain the isomorphism
\begin{equation}
\label{eq;24.3.14.20}
 H_1^{\varrho}\bigl(
 \proj^1\setminus\Dtilde,
 \nbigv\otimes\nbige(zu_1^{-1})
 \bigr)
 \simeq
 H_1^{\varrho}\bigl(
 \proj^1\setminus\Dtilde,
 \nbigv\otimes\nbige(zu_2^{-1})
 \bigr)
\end{equation}
induced by the parallel transport of the flat connection of
$\Fourier_+(\nbigv(\varrho))_{|U_{w,\infty}\setminus\{\infty\}}$.
Let us describe the isomorphism (\ref{eq;24.3.14.20})
in terms of the associated constructible sheaves
under the assumptions that
$0\leq \theta^u_1-\theta^u_2<\pi$,
that $|u_i|$ are sufficiently small,
for a path $r(t)e^{\sqrt{-1}\theta^u(t)}$ $(0\leq t\leq 1)$
satisfying $\theta^u_1\leq \theta^u(t)\leq \theta^u_2$.

We have the meromorphic flat bundle
$\nbigstilde^{\infty}_1(\nbigv)$ on $(\proj^1,\Dtilde)$.
For $\kappa=!,\ast$,
let $(\varrho,\kappa):\Dtilde\to\{!,\ast\}$
be the map determined by
$(\varrho,\kappa)(Q)=\varrho(Q)$ $(Q\in D)$
and $(\varrho,\kappa)(\infty)=\kappa$.
Let $\varpi_{\Dtilde}:\projtilde^1_{\Dtilde}\to \proj^1$
denote the oriented real blow up along $\Dtilde$.
We obtain the constructible sheaves
$\nbigl^{(\varrho,\kappa)}\bigl(
 \nbigstilde^{\infty}_1(\nbigv)
 \bigr)$
$(\kappa=!,\ast)$
on $\projtilde^1_{\Dtilde}$
and the projection
\begin{equation}
\label{eq;24.3.14.30}
 \nbigl^{(\varrho,\ast)}\bigl(
 \nbigstilde^{\infty}_1(\nbigv)
 \bigr)
 \lrarr
  \nbigl^{(\varrho,\ast)}\bigl(
 \nbigstilde^{\infty}_1(\nbigv)
 \bigr)
 \Big/
  \nbigl^{(\varrho,!)}\bigl(
 \nbigstilde^{\infty}_1(\nbigv)
 \bigr).
\end{equation}
Let $\iota_{\infty}:\varpitilde_{\Dtilde}^{-1}(\infty)
\to \projtilde^1_{\Dtilde}$ denote the inclusion.
There exists a local system on $L_{0,S^1}$ on
$\varpitilde_{\Dtilde}^{-1}(\infty)$
such that
\[
 \nbigl^{(\varrho,\ast)}\bigl(
 \nbigstilde^{\infty}_1(\nbigv)
 \bigr)
 \Big/
  \nbigl^{(\varrho,!)}\bigl(
 \nbigstilde^{\infty}_1(\nbigv)
 \bigr)
 \simeq
 \iota_{\infty\ast}(L_{0,S^1}).
\]

We identify $\varpitilde_{\Dtilde}^{-1}(\infty)$
with $\real/2\pi\seisuu$
by using the polar coordinate $z=|z|e^{\sqrt{-1}\theta}$.
For any interval with $J$ with
$\vartheta^J_r-\vartheta^J_{\ell}<2\pi$,
we obtain the natural inclusion $\iota_J:J\to S^1$.
We obtain the following subsheaf:
\[
 \iota_{\infty\ast}
 \bigl(
 \iota_{J!}
 \iota_{J}^{-1}L_{0,S^1}
 \bigr)
 \subset
 \iota_{\infty\ast}(L_{0,S^1}).
\]
Let $\nbigl^{(\varrho,\ast)}(\nbigstilde^{\infty}_1(\nbigv))_J$
denote the inverse image of
$\iota_{\infty\ast}
 \bigl(
 \iota_{J!}
 \iota_{J}^{-1}L_{0,S^1}
 \bigr)$
 via the projection (\ref{eq;24.3.14.30}).

We set $J_i=I(\theta_i^u,\pi/2)$.
Let $\alpha_1,\ldots,\alpha_m$ be the complex numbers
such that 
$\pi_1\nbigttilde_1(\nbigi_{\infty}(\nbigv))
=\bigl\{
\alpha_1z,\alpha_2z,\ldots,\alpha_mz\bigr\}$.
If $|u_i|$ are sufficiently large,
there exists a relatively compact interval
$J_0\subset J_1\cap J_2$
such that
$\Re\bigl((u_i^{-1}+\alpha_j)e^{\sqrt{-1}\theta}\bigr)>0$
for any $\theta\in J_0$,
any $j=1,\ldots,m$ and any $i=1,2$.

By using the flat sections $\exp(-zu_i^{-1})$
of $\nbige(zu_i^{-1})$,
we obtain the isomorphisms
\[
  \nbigl^{(\varrho,\ast)}(\nbigstilde^{\infty}_1(\nbigv))
  _{|\projtilde^1_{\Dtilde}\setminus \varpi_{\Dtilde}^{-1}(\infty)}
\simeq
 \nbigl^{(\varrho,\ast)}\bigl(
 \nbigv\otimes\nbige(zu_i^{-1})
 \bigr)_{|\projtilde^1_{\Dtilde}\setminus \varpi_{\Dtilde}^{-1}(\infty)}.
\]
It extends to a morphism
\[
 \nbigl^{(\varrho,\ast)}(\nbigstilde^{\infty}_1(\nbigv))_{J_0}\lrarr
 \nbigl^{(\varrho,\ast)}\bigl(
 \nbigv\otimes\nbige(zu_i^{-1})
 \bigr),
\]
and the cokernel are acyclic with respect to the global cohomology.
Therefore, we obtain the following isomorphisms:
\begin{multline}
\label{eq;24.3.14.31}
 H_1^{\varrho}\bigl(
 \proj^1\setminus\Dtilde,
 \nbigv\otimes\nbige(zu_1^{-1})
 \bigr)
 \stackrel{\simeq}{\llarr}
 H^1\Bigl(
 \projtilde^1_{\Dtilde},
 \nbigl^{(\varrho,\ast)}(\nbigstilde^{\infty}_1(\nbigv))_{J_0}
 \Bigr)
 \\
 \stackrel{\simeq}{\lrarr}
 H_1^{\varrho}\bigl(
 \proj^1\setminus\Dtilde,
 \nbigv\otimes\nbige(zu_2^{-1})
 \bigr).
\end{multline}
It is easy to see that
(\ref{eq;24.3.14.20})
equals
(\ref{eq;24.3.14.31}).

\subsection{A reduction}
\label{subsection;24.3.18.1}

\begin{lem}
\label{lem;25.3.6.10}
If $|u|$ is sufficiently large,
for any $\varrho\in \Dsf(D)$,
there exists the natural isomorphism
\begin{equation}
\label{eq;24.3.17.111}
 H_1^{\varrho}\bigl(
 \cnum\setminus D,
 \nbigv\otimes\nbige(zu^{-1})
 \bigr)
 \simeq
  H_1^{\varrho}\bigl(
 \cnum\setminus D,
 \nbigstilde^{\infty}_1(\nbigv)\otimes\nbige(zu^{-1})
 \bigr).
\end{equation}
As a result, we obtain the isomorphism
of $2\pi\seisuu$-equivariant local systems
\begin{equation}
\label{eq;24.3.17.110}
 \gbigl_{\varrho}^{\gbigf}(\nbigv)
\simeq
\gbigl_{\varrho}^{\gbigf}(\nbigstilde^{\infty}_1(\nbigv)). 
\end{equation}
\end{lem}
\pf
We use the notation in \S\ref{subsection;24.3.14.42}.
Let $u=|u|e^{\sqrt{-1}\theta^u}$.
Let $I\subset I(\theta^u,\pi/2)$
be a relatively compact interval.
By using the flat section $\exp(-zu^{-1})$ of $\nbige(zu^{-1})$,
we obtain the isomorphisms
\begin{equation}
\label{eq;24.3.17.100}
  \nbigl^{(\varrho,\ast)}(\nbigstilde^{\infty}_1(\nbigv))
  _{|\projtilde^1_{\Dtilde}\setminus \varpi_{\Dtilde}^{-1}(\infty)}
\simeq
 \nbigl^{(\varrho,\ast)}\bigl(
 \nbigv\otimes\nbige(zu^{-1})
 \bigr)_{|\projtilde^1_{\Dtilde}\setminus \varpi_{\Dtilde}^{-1}(\infty)},
\end{equation}
\begin{equation}
\label{eq;24.3.17.101}
   \nbigl^{(\varrho,\ast)}(\nbigstilde^{\infty}_1(\nbigv))
  _{|\projtilde^1_{\Dtilde}\setminus \varpi_{\Dtilde}^{-1}(\infty)}
\simeq
 \nbigl^{(\varrho,\ast)}\bigl(
 \nbigstilde^{\infty}_1(\nbigv)\otimes\nbige(zu^{-1})
 \bigr)_{|\projtilde^1_{\Dtilde}\setminus \varpi_{\Dtilde}^{-1}(\infty)}. 
\end{equation}
If $|u|$ is sufficiently large,
they extend to the following monomorphisms:
\begin{equation}
\label{eq;24.3.17.102}
  \nbigl^{(\varrho,\ast)}(\nbigstilde^{\infty}_1(\nbigv))_I
\lrarr
  \nbigl^{(\varrho,\ast)}\bigl(
 \nbigv\otimes\nbige(zu^{-1})
 \bigr),
\end{equation}
\begin{equation}
\label{eq;24.3.17.103}
   \nbigl^{(\varrho,\ast)}(\nbigstilde^{\infty}_1(\nbigv))_I
\lrarr
 \nbigl^{(\varrho,\ast)}\bigl(
 \nbigstilde^{\infty}_1(\nbigv)\otimes\nbige(zu^{-1})
 \bigr).
\end{equation}
The cokernel of (\ref{eq;24.3.17.102})
and (\ref{eq;24.3.17.103})
are acyclic with respect to the global cohomology.
Hence, we obtain 
the isomorphism (\ref{eq;24.3.17.111}).
\hfill\qed

\vspace{.1in}
We shall prove the following proposition
in \S\ref{section;24.3.17.120}.

\begin{prop}
\label{prop;24.3.17.121}
The isomorphism {\rm(\ref{eq;24.3.17.111})}
induces an isomorphism of  
$2\pi\seisuu$-equivariant local systems with Stokes structure.
\end{prop}

\subsection{Some induced maps}
\label{subsection;24.4.15.1}

 \begin{lem}
\label{lem;18.5.12.21}
Take $Q\in D$
and for $\omega\in\rnum_{>0}$.
Take a morphism
$\varrho_1\lrarr\varrho_2$ in $\Dsf(D)$
such that
$\varrho_1(Q)=!$
and $\varrho_2(Q)=\ast$.
There exists the following naturally defined commutative diagram
of $2\pi\seisuu$-equivariant local systems on $\real$:
\[
 \begin{CD}
 \gbigl^{\gbigf}_{\varrho_1}\bigl(
 \nbigs^Q_{\omega}(\nbigv)
 \bigr)
 @>>>
 \gbigl^{\gbigf}_{\varrho_2}\bigl(
 \nbigs^Q_{\omega}(\nbigv)
 \bigr)
\\
 @VVV @AAA\\
 \gbigl^{\gbigf}_{\varrho_1}\bigl(\nbigv\bigr)
 @>>>
 \gbigl^{\gbigf}_{\varrho_2}\bigl(\nbigv\bigr)
 \end{CD}
\]
 \end{lem}
\pf
Let $\nbige(zw)$
denote the meromorphic flat bundle
$\bigl(\nbigo_{\proj^1_z}(\ast\infty),d+d(zw)\bigr)$
on $(\proj^1_z,\infty)$.
The fibers of
$\gbigl^{\gbigf}_{\varrho}\bigl(\nbigv
 \bigr)$ 
over $\theta^u\in\real$
are identified with the cohomology groups
\[
 H^1\Bigl(
 \proj^1_z,
 \nbigv(\varrho)
 \otimes
  \nbige(zte^{-\sqrt{-1}\theta^u})
\otimes\Omega^{\bullet}_{\proj_z^1}
 \Bigr)
\]
for a sufficiently large $t>0$.
Similarly,
the fibers of
$\gbigl^{\gbigf}_{\varrho}\bigl(
 \nbigs^Q_{\omega}(\nbigv)
 \bigr)$ 
over $\theta^u\in\real$
are identified with the cohomology groups
\[
 H^1\Bigl(
 \proj^1_z,
 \nbigs^Q_{\omega}(\nbigv)(\varrho)
 \otimes
  \nbige(zte^{-\sqrt{-1}\theta^u})
\otimes\Omega^{\bullet}_{\proj_z^1}
 \Bigr)
\]
for a sufficiently large $t>0$.
Hence,
we obtain the desired morphisms 
by the consideration in 
\S\ref{subsection;18.5.12.20}.
\hfill\qed

\vspace{.1in}

Take a small neighbourhood 
$U_{z,\infty}$ of $\infty$ in $\proj^1_z$
such that
$(\nbigv,\nabla)_{|U_{z,\infty}\setminus\{\infty\}}$ 
is a flat bundle.
Take $\omega>1$.
We obtain the meromorphic flat bundle
$\nbigt_{\omega}\bigl((\nbigv,\nabla)_{|U_{z,\infty}}\bigr)$
on $U_{z,\infty}$.
It naturally extends to a meromorphic flat bundle
on $(\proj_z^1,\{0,\infty\})$
with regular singularity at $0$,
which we denote by 
$\nbigt^{\infty}_{\omega}(\nbigv,\nabla)$.

\begin{lem}
\label{lem;18.5.12.40}
Let $\varrho_1\lrarr\varrho_2$
be a morphism in $\Dsf(D)$.
There exists the following naturally defined 
commutative diagram:
\[
 \begin{CD}
 \gbigl^{\gbigf}_{!}(\nbigt^{\infty}_{\omega}(\nbigv))
 @>>>
 \gbigl^{\gbigf}_{\ast}(\nbigt^{\infty}_{\omega}(\nbigv))
 \\
 @VVV @AAA \\
 \gbigl^{\gbigf}_{\varrho_1}(\nbigv)
 @>>>
 \gbigl^{\gbigf}_{\varrho_2}(\nbigv).
 \end{CD}
\]
\end{lem}
\pf
We obtain the desired morphisms
from the consideration in 
\S\ref{subsection;18.5.12.20}.
\hfill\qed

\section{Variant for constructible sheaves}
\label{subsection;21.6.9.1}
Let $\varpi:\projtilde_{\infty}^1\lrarr\proj^1$
denote the oriented real blow up along $\infty$.
By the standard coordinate $z$ on $\cnum=\proj^1\setminus\{\infty\}$,
the fiber $\varpi^{-1}(\infty)$ is identified with
$S^1=\{e^{\sqrt{-1}\theta}\,|\,\theta\in\real\}$.

We set $\nbigy:=\projtilde^1_{\infty}\times\real$.
Let $Z\subset \varpi^{-1}(\infty)\times\real$
denote 
$\{(e^{\sqrt{-1}\theta},\theta^u)\,|\,
\Re(e^{\sqrt{-1}(\theta-\theta^u)})\leq 0\}$.
Let $\iota:\nbigy\setminus Z\lrarr \nbigy$
be the open embedding.
We set
$\gbigp:=\iota_!\cnum_{\nbigy\setminus Z}$.

Let $D\subset\cnum$ be a finite subset.
Let $\varpi_D:
 \projtilde^1_{D\cup\{\infty\}}\lrarr\projtilde^1_{\infty}$
denote the oriented real blow up along $D$.
We set
$\nbigy_D:=\projtilde^1_{D\cup\{\infty\}}\times \real$.
Let $q_i$ $(i=1,2)$ denote the projections of $\nbigy_D$
onto the $i$-th component.
For any constructible sheaf $\nbign$ on 
$\projtilde^1_{D\cup\infty}$,
we set
$\gbigf(\nbign):=
 Rq_{2\ast}\bigl(
 q_1^{-1}(\nbign)\otimes
 (\varpi_D\times\id_{\real})^{-1}\gbigp\bigr)[1]$
in the derived category of 
$2\pi\seisuu$-equivariant
cohomologically constructible sheaves on $\real$.
\index{sheaf $\gbigf(\nbign)$}

\vspace{.1in}

Let $(\nbigv,\nabla)$ be a meromorphic flat bundle
on $(\proj^1,D\cup\{\infty\})$
with regular singularity at $\infty$.
Then,
for any $\varrho\in \Dsf(D)$,
$\Fourier_+(\nbigv(\varrho))(\ast 0)$
is the meromorphic flat bundle
on $(\proj^1,\{0,\infty\})$
with regular singularity at $0$.
The following lemma is obvious.
\begin{lem}
\label{lem;24.2.12.5}
$\gbigl^{\gbigf}_{\varrho}(\nbigv)$
is naturally isomorphic to
$\gbigf\bigl(
 \nbigl^{\varrho}(\nbigv)
 \bigr)$.
\hfill\qed
\end{lem}

For $R\geq 0$,
let $U_R:=\{|z|>R\}\cup\{\infty\}\subset\proj^1$.
Let $\Utilde_R:=\varpi^{-1}(U_R)$.
Let $j:\Utilde_R\lrarr \projtilde^1_{\infty}$
denote the inclusion.
Let $\nbigl$ be a local system on $\Utilde_R$.
For later use,
we study
$\gbigf(j_{\star}\nbigl)$.
Let $F$ denote
the automorphism of $\nbigl$
obtained as the monodromy along the loop
$e^{2\pi\sqrt{-1}t}z$ $(0\leq t\leq 1)$
for any $z\in \Utilde(R)$.
Let $\varphi:\real\lrarr\Utilde_R$
be the map defined
by $\varphi(\theta^u)=(\infty,e^{\sqrt{-1}\theta^u})$.

\begin{lem}
There exist isomorphisms
of $2\pi\seisuu$-equivariant local systems
$\gbigf\bigl(
 j_{\star}\nbigl
 \bigr)
 \simeq
 \varphi^{-1}\nbigl$ $(\star=!,\ast)$
such that the natural morphism
$\gbigf\bigl(
j_{!}\nbigl
 \bigr)
 \lrarr
 \gbigf\bigl(
j_{\ast}\nbigl
 \bigr)$
is identified with
$\varphi^{-1}(\id-F^{-1})$.
 \end{lem}
\pf
Let $\theta^u\in\real$.
We set
$Z(\theta^u):=
\bigl\{e^{\sqrt{-1}\theta}\,\big|\,
\Re(e^{\sqrt{-1}(\theta-\theta^u)})\leq 0
\bigr\}$.
We set
$W(\theta^u):=\Utilde_R\setminus(\{\infty\}\times Z(\theta^u))$.
Let $\iota^{\theta^u}:W(\theta^u)\lrarr \Utilde_R$
denote the inclusion.
We obtain
$H^a\bigl(\projtilde^1_{\infty},
j_{\star}
\iota^{\theta^u}_!(\nbigl_{|W(\theta^u)})
\bigr)=0$
unless $a=1$,
and
the stalk of
$\gbigf(j_{\star}\nbigl)$ at $\theta^u$
is naturally isomorphic to
$H^1\bigl(\projtilde^1_{\infty},
j_{\star}
\iota^{\theta^u}_!(\nbigl_{|W(\theta^u)})
\bigr)$.

Let $\Gamma_{\ast}(\theta^u)$ denote a path in $W(\theta^u)$
connecting
$(R-\epsilon,e^{\sqrt{-1}\theta^u})$
and
$(\infty,e^{\sqrt{-1}\theta^u})$
for small $\epsilon>0$
along the ray
$\{(t,e^{\sqrt{-1}\theta^u})\,|\,R-\epsilon\leq t\leq \infty\}$.
Any $s\in \nbigl_{(\infty,e^{\sqrt{-1}\theta^u})}$
induces a section $\stilde$
of $\nbigl$ along $\Gamma_{\ast}(\theta^u)$,
which induces a global section
$\stilde\otimes\Gamma_{\ast}(\theta^u)$
of
$j_{\ast}\iota^{\theta^u}_!(\nbigl_{|W(\theta^u)})
\otimes
 \nbigc^{-1}_{\projtilde^1_{\infty},\del\projtilde^1_{\infty}}$.
(See \S\ref{subsection;20.9.8.1}
for $\nbigc^{\bullet}_{Y,\del Y}$.)
 It is a cocycle,
and induces an element
$[\stilde\otimes\Gamma_{\ast}(\theta^u)]
\in H^1(\projtilde^1_{\infty},
j_{\ast}\iota^{\theta^u}_!(\nbigl_{|W(\theta^u)}))$.
It is easy to see that
the correspondence
$s\longmapsto [\stilde\otimes\Gamma_{\ast}(\theta^u)]$
induces an isomorphism
$\nbigl_{e^{\sqrt{-1}\theta^u}}
\simeq
\gbigf(j_{\ast}\nbigl)_{|\theta^u}$.
Thus, we obtain
$\varphi^{-1}(\nbigl)\simeq
\gbigf(j_{\ast}\nbigl)$.

Let $\Gamma_{!,0}(\theta^u)$ denote a path connecting
$(\infty,e^{\sqrt{-1}\theta^u})$
and $(2R,e^{\sqrt{-1}\theta^u})$
along the ray
$\{(t,e^{\sqrt{-1}\theta^u})\,|\,2R\leq t\leq\infty\}$.
Let $\Gamma_{!,1}(\theta^u)$ denote the path
$2Re^{\sqrt{-1}\theta^u}e^{2\pi\sqrt{-1}t}$
$(0\leq t\leq 1)$.
Let $\Gamma_{!,2}(\theta^u)$ denote a path
connecting
$(2R,e^{\sqrt{-1}\theta^u})$
and 
$(\infty,e^{\sqrt{-1}\theta^u})$
along the ray
$\{(t,e^{\sqrt{-1}\theta^u})\,|\,2R\leq t\leq\infty\}$.
We obtain a $1$-chain $\Gamma_!(\theta^u)$
from $\Gamma_{!,0}(\theta^u)$,
$\Gamma_{!,1}(\theta^u)$ and $\Gamma_{!,2}(\theta^u)$.
Any $s\in \nbigl_{e^{\sqrt{-1}\theta^u}}$
induces a section $\stilde_2$ along $\Gamma_{!,2}(\theta^u)$.
Let $\stilde_1$ denote the section along $\Gamma_{!,1}(\theta^u)$
which equals $\stilde_2$ at $t=1$.
Subsequently, we obtain the section $\stilde_0$
along $\Gamma_{!,0}(\theta^u)$.
Thus, we obtain a global section
$\stilde\otimes\Gamma_!(\theta^u)
:=\sum \stilde_i\otimes\Gamma_{!i}(\theta^u)$
of 
$j_{!}\iota^{\theta^u}_!(\nbigl_{|W(\theta^u)})
\otimes
\nbigc^{-1}_{\projtilde^1_{\infty},\del\projtilde^1_{\infty}}$.
It is a cocycle,
and induces an element
$[\stilde\otimes\Gamma_{!}(\theta^u)]
\in H^1(\projtilde^1_{\infty},
j_{!}\iota^{\theta^u}_!(\nbigl_{|W(\theta^u)}))$.
It is easy to see that
the correspondence
$s\longmapsto [\stilde\otimes\Gamma_{!}(\theta^u)]$
induces an isomorphism
$\nbigl_{e^{\sqrt{-1}\theta^u}}
\simeq
\gbigf(j_{!}\nbigl)_{|\theta^u}$.
Thus, we obtain
$\varphi^{-1}(\nbigl)\simeq
\gbigf(j_{!}\nbigl)$.
By the construction,
it is easy to see that
the natural morphism
$\gbigf\bigl(
j_{!}\nbigl
 \bigr)
 \lrarr
 \gbigf\bigl(
j_{\ast}\nbigl
 \bigr)$
is identified with
$\varphi^{-1}(\id-F^{-1})$.
\hfill\qed

\begin{rem}
There are several other good
isomorphisms
$\gbigf(j_{\star}\nbigl)\simeq\varphi^{-1}\nbigl$.
For instance,
we may reverse the orientation of the paths.
In the case $\star=!$,
we may first construct
$\stilde'_{0}$ along $\Gamma_{!,0}$,
and $\stilde'_i$ along $\Gamma_{!,i}$ $(i=1,2)$
subsequently.
Under different identifications,
the morphism
$\gbigf(j_!\nbigl)\lrarr \gbigf(j_{\ast}\nbigl)$
is presented in different ways. 
\hfill\qed
\end{rem}

\section{Pairings between homology classes and Rham cohomology classes}
\label{subsection;24.3.29.110}

\subsection{Variations of the de Rham complex of $\nbigv(\varrho)$}

Let $C$ be a compact Riemann surface with $\del C=\emptyset$.
Let $D$ be a finite subset of $C$.
Let $(\nbigv,\nabla)$ be a meromorphic flat bundle on $(C,D)$.
Let $\varrho\in\Dsf(D)$.
There are convenient complexes to study
the de Rham cohomology
$\hyperh^{\ast}(C,\nbigv(\varrho)\otimes\Omega^{\bullet}_C)$.

\subsubsection{Local unramified case}

Let us consider the case
$(C,D)=(\Delta_z,0)$
and $\nbigi(\nbigv)\subset z^{-1}\cnum[z^{-1}]$.
There exists the decomposition
$(\nbigv,\nabla)\otimes\cnum[\![z]\!]
 =\bigoplus_{\gminia\in\nbigi(\nbigv)}
 (\nbigvhat_{\gminia},\nabla_{\gminia})$,
where $(\nbigvhat_{\gminia},\nabla_{\gminia}-d\gminia\id)$
are regular singular.
For any $a\in\real$,
there exist lattices
$\nbigvhat_{\gminia,-a}\subset \nbigvhat_{\gminia}$
such $\nabla_{\gminia}-d\gminia$ are logarithmic
with respect to $\nbigvhat_{\gminia,-a}$
and that the eigenvalues $\alpha$ of the residues of
$\nabla_{\gminia}-d\gminia\id$
satisfy $a<\Re(\alpha)\leq a+1$.

For each $a\geq 0$,
there exists the subcomplex
$\nbigc^{\bullet}_{!0}(\nbigv)_{-a}\subset
\nbigv\otimes\Omega^{\bullet}$
determined by the conditions
$\nbigc^{\bullet}_{!0}(\nbigv)_{-a}=
\nbigv\otimes\Omega^{\bullet}$
on $\Delta_z^{\ast}$,
and 
\[
 \nbigc^{0}_{!0}(\nbigv)_{-a}
 \otimes\cnum[\![z]\!]
=\bigoplus
 \nbigvhat_{\gminia,-a+\ord(\gminia)},
 \quad
  \nbigc^{1}_{!0}(\nbigv)_{-a}
 \otimes\cnum[\![z]\!]
=\bigoplus
 \nbigvhat_{\gminia,-a}\frac{dz}{z}.
\]       
It is well known and easy to see that
there exists a natural quasi-isomorphism
$\nbigc^{\bullet}_{!0}(\nbigv)_{-a}\to
\nbigv(!0)\otimes\Omega^{\bullet}$.

\subsubsection{Local ramified case}

For $p\in\seisuu_{>0}$,
let $\rho_p:\Delta_{z_p}\to \Delta_z$
be the map defined by $\rho_p(z_p)=z_p^p$.
There exists $p\in\seisuu_{>0}$
such that
$\rho_p^{\ast}(\nbigv,\nabla)$ is unramified.
For $a\geq 0$,
we obtain the $\Gal(p)$-invariant complex
$\nbigc^{\bullet}_{!0}(\rho_p^{\ast}(\nbigv))_{-pa}$.
As the descent,
we obtain a complex
$\nbigc^{\bullet}_{!0}(\nbigv)_{-a}$.
There exists a natural inclusion
$\nbigc^{\bullet}_{!0}(\nbigv)_{-a}
\to
\nbigv(!0)\otimes\Omega^{\bullet}$
which is a quasi-isomorphism.

\subsubsection{Global case}
\label{subsection;24.3.29.120}

We set $D(!)=\varrho^{-1}(!)$.
For $a\geq 0$,
let $\nbigc^{\bullet}_{\varrho}(\nbigv)_{-a}$
denote the subcomplex of
$\Omega^{\bullet}\otimes \nbigv(\varrho)$
determined by the following conditions.
\begin{itemize}
 \item $\nbigc^{\bullet}_{\varrho}(\nbigv)_{-a}
       =\Omega^{\bullet}\otimes\nbigv(\varrho)$
       on $C\setminus D(!)$.
 \item  $\nbigc^{\bullet}_{\varrho}(\nbigv)_{-a|C_P}=
        \nbigc^{\bullet}_{!P}(\nbigv_{|C_P})_{-a}$
	for any $P\in D(!)$.
\end{itemize}
Then, there exists the natural inclusion
$\nbigc_{\varrho}^{\bullet}(\nbigv)_{-a}\to
\Omega^{\bullet}_C\otimes\nbigv(\varrho)$,
which is a quasi-isomorphism.
Let $\nbigc_{\varrho,C^{\infty}}^{\bullet}(\nbigv)_{-a}$
denote the Dolbeault resolution of
$\nbigc_{\varrho}^{\bullet}(\nbigv)_{-a}$.

\subsubsection{Infinitely decay complex}

Let $\nbigc^{\infty}_C$ denote the sheaf of $C^{\infty}$-functions
on $C$,
and we set
$\nbigv_{C^{\infty}}=\nbigc^{\infty}\otimes_{\nbigo}\nbigv$.
Let $f$ be a section of $\nbigv_{C^{\infty}}$ around $P\in D(!)$.
Let $v_1,\ldots,v_r$ be a frame of $\nbigv$ around $P$.
There exist $N\in\seisuu_{>0}$ and
$C^{\infty}$-functions $f_i$ $(i=1,\ldots,r)$ such that
$f=\sum f_iz^{-N}v_i$ around $P$.
We say that $f$ is infinitely decay at $P$
if the Taylor series of $f_i$ at $P$ are $0$.

For any open subset $U$ of $C$,
let $\nbigv_{C^{\infty},\varrho}(U)$ be the space of
sections $f$ of $\nbigv_{C^{\infty}}$
which are infinitely decay at each point of $U\cap D(!)$.
We obtain a subsheaf
$\nbigv_{C^{\infty},\varrho}\subset\nbigv_{C^{\infty}}$.
By the connection $\nabla$,
we obtain the complex of sheaves
\[
 \nbigc^{\bullet}_{\varrho,C^{\infty}}(\nbigv)_{-\infty}=
 \Tot\bigl(\nbigv_{C^{\infty},\varrho}\otimes\Omega^{\bullet,\bullet}\bigr),
\]
where $\Tot$ denotes the total complex of the double complex.
The following is a consequence of a result of
Mebkhout \cite{Mebkhout-positivity}.
(See also \cite[Proposition 2.1.4, Proposition 3.2.1]{Mochizuki-Betti}.)
In this case, we can check it directly.

\begin{prop}
The natural inclusion
$\nbigc^{\bullet}_{\varrho,C^{\infty}}(\nbigv)_{-\infty}
 \to
 \nbigc^{\bullet}_{\varrho,C^{\infty}}(\nbigv)_{-a}$ 
is a quasi-isomorphism.
\hfill\qed
\end{prop}

\subsection{Pairings}

Let $(\nbigv^{\lor},\nabla)$ denote the dual
meromorphic flat bundle of $\nbigv$,
i.e.,
$\nbigv^{\lor}=\nhom_{\nbigo_C(\ast D)}(\nbigv,\nbigo_{C}(\ast D))$
equipped with the naturally induced connection $\nabla$.

Let $\varrhobar\in\Dsf(D)$ be defined by the condition
$\{\varrho(P),\varrhobar(P)\}=\{!,\ast\}$
for any $P\in D$.
Let $j:C\setminus D\lrarr \Ctilde$ denote the inclusion.
There exists the following naturally defined morphism
of the sheaves on $\Ctilde(D)$:
\begin{equation}
\label{eq;24.3.28.30}
  \nbigl^{\varrhobar}(\nbigv^{\lor})
  \otimes
  \nbigl^{\varrho}(\nbigv)
  \lrarr
  j_!\cnum_{C\setminus D}.
\end{equation}
It induces a perfect pairing
\[
\langle\cdot,\cdot\rangle:
 H_1^{\varrhobar}\bigl(
 C\setminus D,
 \nbigv^{\lor}
 \bigr) 
\otimes
\hyperh^1\bigl(C,\Omega^{\bullet}_C\otimes \nbigv(\varrho)\bigr)
 \lrarr
 \cnum.
\]

We shall describe it as the integration
of the $1$-forms along the $1$-chains
by following Hien \cite{Hien}.

\subsubsection{Integrability}

Let $(V,\nabla)$ be a meromorphic flat bundle on $(\Delta_z,0)$.
Let $\gamma:[0,1]\to \Deltatilde_z$ be a $C^{\infty}$-map
such that
$\gamma(0)\in \varpi^{-1}(0)$
and $\gamma(\openclosed{0}{1})\subset\Delta_z\setminus\{0\}$.
Let $c$ be a section of
$\gamma^{-1}(\nbigl^{\leq 0}(V^{\lor},\nabla))$.
Let $\tau$ be a section of
$\nbigc^1_{!,C^{\infty}}(V,\nabla)_{-a}$
for some $a\geq 0$.
We obtain the $\gamma^{-1}(V)$-valued $1$-form
$\gamma^{\ast}(\tau)$.
We obtain the $1$-form
$(c,\tau)$ on $\openclosed{0}{1}$
as the pairing of
$c$ and $\gamma^{\ast}(\tau)$.

\begin{lem}
\label{lem;24.3.28.20}
$(c,\tau)$ is integrable on $[0,1]$.
\end{lem}
\pf
We may assume that the image of $\gamma$ is contained
in a small sector.
By considering the pull back via an appropriate ramified covering,
we have only to study the case
where $\nbigi(V)\subset z^{-1}\cnum[z^{-1}]$.
By the asymptotic analysis in \S\ref{subsection;24.3.28.10},
it is enough to study the case where
there exist
$\gminia\in z^{-1}\cnum[z^{-1}]$
and a regular singular meromorphic flat bundle
$(V,\nabla_1)$ on $(\Delta_z,0)$
such that $(V,\nabla)=(V,\nabla_1+d\gminia\id)$.

There exists a lattice $V_{-a}\subset V$
such that $\nabla_1$ is logarithmic with respect to $V_{-a}$
and that the eigenvalues $\alpha$ of the residue of $\nabla_1$
satisfies $a< \Re(\alpha)\leq a+1$.
Let $v_1,\ldots,v_r$ be a frame of $V_{-a}$.
We note that $\tau$ is expressed as
$\tau=\sum
 \bigl(
 \tau_i^{1,0}dz/z
 +
 \tau_i^{0,1}d\zbar
 \bigr)v_i$,
where $\tau_i^{p,q}$ are $C^{\infty}$-functions.
Let $v_i^{\lor}$ denote the frame of $V^{\lor}$
obtained as the dual of $v_1,\ldots,v_r$.
We may regard $c$ as a section of $\nbigl(V^{\lor},\nabla)$
around $\gamma(0)$.
We have the expression $c=\sum c_iv_i^{\lor}$.

Let us study the case $\gminia=0$.
Let $A=(A_{i,j})$ be determined by
$\nabla(v_j^{\lor})=\sum A_{i,j}v_i^{\lor}\,dz/z$.
Then, $A$ is holomorphic at $z=0$,
and the eigenvalues $\beta$ of $A(0)$
satisfy $-a-1\leq \Re(\beta)<-a$.
Hence, there exists $\delta>0$ such that
$|c_i|=O\bigl(|z|^{a+\delta}\bigr)$.
We obtain
$|(c,\tau)|=O(|z|^{a+\delta-1})$
and the desired integrability.

Let us study the case $\gminia\neq 0$.
If $\gminia<_{\gamma(0)}0$ does not hold,
we obtain $c=0$
by the moderate growth condition.
If $\gminia<_{\gamma(0)}0$,
i.e., $-\Re(\gminia)<0$ around $\gamma(0)$,
then $|c_i|=O(\exp(-\delta_1|z|^{-\delta_2}))$
for some $\delta_i>0$.
Then, we obtain
$|(c,\tau)|=O(\exp(-\delta_3|z|^{-\delta_4}))$
for some $\delta_i>0$,
and hence the desired integrability.
\hfill\qed

\vspace{.1in}
Similarly, we obtain the following lemma.
\begin{lem}
Let $\tau$ be a section of
$\nbigc^0_{!,C^{\infty}}(V,\nabla)_{-a}$.
We obtain 
the function $(c,\tau)$ on $\openclosed{0}{1}$
as the pairing of $\gamma^{\ast}(c)$
and $\gamma^{\ast}(\tau)$,
and it is integrable.
\hfill\qed
\end{lem}

\subsubsection{Pairings and integrations}

For a sheaf $F$,
let $\Gamma(F)$ denote the space of global sections of $F$.
Let $\beta=\sum c_i\otimes \gamma_i\in
\Gamma\bigl(
\nbigc^{-1}_{\Ctilde,\del \Ctilde}
\otimes
\nbigl^{\varrho}(\nbigv^{\lor})\bigr)$
and
$\tau\in
\Gamma(\nbigc^1_{\varrho,C^{\infty}}(\nbigv)_{-a})$.
We may assume that 
$\gamma_i$ are $C^{\infty}$-functions
$[0,1]\to\Ctilde$
such that
$\gamma_i(\openopen{0}{1})\subset C\setminus D$.
We obtain the $1$-form
$(c_i,\tau)$ on $\openopen{0}{1}$
as the pairing of
$\gamma_i^{\ast}(\nbigv)$-valued $1$-form
$\gamma_i^{\ast}(\tau)$ 
and the section $c_i$ of $\gamma_i^{\ast}(\nbigv^{\lor})$.
By Lemma \ref{lem;24.3.28.20},
$(c_i,\tau)$ are integrable on $[0,1]$.
We obtain
\[
 \langle \beta,\tau\rangle'
=\sum_i\int_{\gamma_i} (c_i,\tau)\in\cnum.
\]
By the Stokes formula,
if $\beta$ is a cycle,
we obtain 
$\langle\beta,d\tau_0\rangle'=0$
for any 
$\tau_0\in
\Gamma\bigl(
\nbigc^0_{\varrho,C^{\infty}}(\nbigv)_{-a}
\bigr)$.
If $\tau$ is $1$-cocycle,
we obtain
$\langle\del\beta_2,\tau\rangle'=0$
for any $\beta_0\in
\Gamma\bigl(
\nbigc^{-2}_{\Ctilde,\del \Ctilde}
\otimes
\nbigl^{\varrho}(\nbigv^{\lor})\bigr)$.
The $1$-cohomology of the complex
$\Gamma\bigl(
\nbigc^{\bullet}_{\varrho,C^{\infty}}(\nbigv)_{-a}
\bigr)$
is isomorphic to 
$\hyperh^1\bigl(C,\Omega^{\bullet}_C\otimes \nbigv(\varrho)\bigr)$,
and 
the $1$-cohomology of
$\Gamma\bigl(
 \nbigc^{\bullet}_{\Ctilde,\del \Ctilde}
\otimes
\nbigl^{\varrho}(\nbigv^{\lor})[-2]
\bigr)$
is isomorphic to
$H_1^{\varrho}\bigl(
 C\setminus D,\nbigv^{\lor}
\bigr)$.
We obtain
\[
 \langle\cdot,\cdot\rangle':
 H_1^{\varrhobar}\bigl(
 C\setminus D,
 (\nbigv^{\lor},\nabla)
 \bigr) 
\otimes
\hyperh^1\bigl(C,\Omega^{\bullet}_C\otimes \nbigv(\varrho)\bigr)
 \lrarr
 \cnum.
\]
The following proposition is essentially due to
Bloch-Esnault \cite{Bloch-Esnault2}
and Hien \cite{Hien}.
\begin{prop}
$\langle\cdot,\cdot\rangle=\langle\cdot,\cdot\rangle'$.
\end{prop}
\pf
The case $D(!)=\emptyset$
is studied in \cite{Hien}.
We need only some minor modification.
We set $D(\ast)=\varrho^{-1}(\ast)$.
Let $\nbigctilde_{\varrho,C^{\infty}}^{\bullet}(\nbigv)_{-\infty}$
denote the complex of sheaves on $\Ctilde$
determined by the following conditions.
\begin{itemize}
 \item On $C\setminus D$,
       $\nbigctilde_{\varrho,C^{\infty}}^{\bullet}(\nbigv)_{-\infty}
       =\nbigc^{\bullet}_{\varrho,C^{\infty}}(\nbigv)_{-\infty}$.
 \item Let $P\in D(!)$.
       On $\varpi^{-1}(C_P)$,
\[
       \nbigctilde_{\varrho,C^{\infty}}^{j}(\nbigv)_{-\infty}
       =
       \nbigp^{<P}\otimes_{\varpi^{-1}\nbigc^{\infty}_{C_P}}
       \varpi^{-1}\bigl(
       \nbigc_{\varrho,C^{\infty}}^{\bullet}(\nbigv)_{-a}
       \bigr)
\]
       for some $0\leq a<\infty$,
       where $\nbigp^{<P}$
       denotes the sheaf of $C^{\infty}$-functions
       on $\varpi^{-1}(C_P)$
       whose Taylor series are $0$
       at any point of $\varpi^{-1}(P)$.
       It is independent of $a$.
 \item Let $P\in D(\ast)$.
       On $\varpi^{-1}(C_P)$,
\[
       \nbigctilde_{\varrho,C^{\infty}}^{\bullet}(\nbigv)
       =
       \nbigp^{\moderate}
       \otimes_{\varpi^{-1}(\nbigc^{\infty}_{C_P})}
       \varpi^{-1}
       \bigl(
       \nbigc_{\varrho,C^{\infty}}^{\bullet}(\nbigv)_{-a}
       \bigr)
\]
       for some $0\leq a<\infty$,
       where $\nbigp^{\moderate}$
       denotes the sheaf of $C^{\infty}$-functions 
       of moderate growth
       on $\varpi^{-1}(C_P)$.
       Here,
       for an open subset $\nbigu$ of
       $\varpi^{-1}(C_P)$,
       a $C^{\infty}$-function $f$ on
       $\nbigu\setminus\varpi^{-1}(P)$
       is called of moderate growth
       if any derivatives of $f$ with respect to
       coordinate systems of $\varpi^{-1}(C_P)$
       have moderate growth around any point of
       $\varpi^{-1}(P)\cap\nbigu$.      
\end{itemize}
There exists the natural inclusion
$\nbigl^{\varrho}(\nbigv)
\lrarr \nbigctilde^{\bullet}_{\varrho,C^{\infty}}(\nbigv)$
by which
$\nbigctilde^{\bullet}_{\varrho,C^{\infty}}(\nbigv)$
is a $c$-soft resolution of 
$\nbigl^{\varrho}(\nbigv)$.
By the construction,
there exists a natural quasi-isomorphism
$\nbigc^{\bullet}_{\varrho,C^{\infty}}(\nbigv)_{-\infty}
\lrarr
\varpi_{\ast}\nbigctilde^{\bullet}_{\varrho,C^{\infty}}(\nbigv)_{-\infty}$.
The pairing
between
$\Gamma(\nbigc^{-1}_{\Ctilde,\del\Ctilde}\otimes
\nbigl^{\varrho}(\nbigv^{\lor}))$
and
$\Gamma\bigl(
 \nbigc^{\bullet}_{\varrho}(\nbigv)_{-\infty}
\bigr)$
naturally extends to a pairing
between 
$\Gamma(\nbigc^{-1}_{\Ctilde,\del\Ctilde}\otimes
\nbigl^{\varrho}(\nbigv^{\lor}))$
and
$\Gamma\bigl(
\varpi_{\ast}\nbigctilde^{\bullet}_{\varrho}(\nbigv)_{-\infty}
\bigr)$,
and it induces $\langle\cdot,\cdot\rangle'$
in the cohomology level.

There exists the natural pairing
between
$\nbigc^{\bullet}_{\Ctilde,\del\Ctilde}
\otimes\nbigl^{\varrhobar}(\nbigv^{\lor})[2]$
and
$\nbigctilde^{\bullet}_{\varrho,C^{\infty}}(\nbigv)_{-\infty}$
to the sheaf of complex
$\distribution^{\rd,-\bullet}_{\Ctilde}$
of distributions with rapid decay
(see \cite{Hien})
such that
\[
  \begin{CD}
   \nbigl^{\varrhobar}(\nbigv^{\lor})
   \otimes
   \nbigl^{\varrho}(\nbigv)
   @>>>
   j_!\cnum_{C\setminus D}\\
   @VVV @VVV \\
   \bigl(
  \nbigc^{\bullet}_{\Ctilde,\del\Ctilde}
   \otimes\nbigl^{\varrhobar}(\nbigv^{\lor})
   \bigr)
  \otimes
  \nbigctilde^{\bullet}_{\varrho,C^{\infty}}(\nbigv)
  @>>>
  \distribution^{\rd,-\bullet}_{\Ctilde}
  \end{CD}
\]
is commutative.
The vertical arrows are quasi-isomorphisms.
Then, we obtain
$\langle\cdot,\cdot\rangle'=\langle\cdot,\cdot\rangle$.
\hfill\qed

\chapter[Transformations of numerical data]{Transformations of numerical associated with Fourier transform}
\label{section;18.6.3.12}

\section{Local Fourier transforms and their explicit expression}
\label{subsection;24.4.9.1}

The local Fourier transform was introduced 
in \cite{Bloch-Esnault1}.
Let $D$ be a finite subset of $\cnum$.
Let $V$ be a meromorphic flat bundle
on $(\proj^1,D\cup\{\infty\})$.
We naturally regard $V$ as a $\nbigd$-module on $\proj^1$.
It is known that
$\Fourier_{\pm}(V)_{|\inftyhat}$ depend
only on 
$V_{|\widehat{\alpha}}$ $(\alpha\in D\cup\{\infty\})$.
More precisely,
according to \cite{Bloch-Esnault1},
there exists a functor 
$\gbigf_{\pm}^{(0,\infty)}$
from the category of
$\cnum(\!(z)\!)\langle \del_z\rangle$-modules
to the category of 
$\cnum(\!(w^{-1})\!)\langle\del_{w^{-1}}\rangle$-modules,
and a functor
$\gbigf_{\pm}^{(\infty,\infty)}$
from the category of
$\cnum(\!(z^{-1})\!)\langle \del_{z^{-1}}\rangle$-modules
to the category of 
$\cnum(\!(w^{-1})\!)\langle\del_{w^{-1}}\rangle$-modules,
such that there exists a natural isomorphism
\begin{equation}
\label{eq;25.2.13.1}
 \Fourier_{\pm}(V)_{|\inftyhat}
\simeq
\gbigf_{\pm}^{(\infty,\infty)}(V_{|\inftyhat})
\oplus
 \bigoplus_{\alpha\in D}
 \gbigf^{(0,\infty)}_{\pm}(V_{|\widehat{\alpha}})\otimes
 \bigl(\cnum(\!(w^{-1})\!),d\pm\alpha\,dw\bigr).
\end{equation}
\index{local Fourier transforms $\gbigf_{\pm}^{(0,\infty)}$,
 $\gbigf^{(\infty,\infty)}_{\pm}$}
The functors $\gbigf_{\pm}^{(0,\infty)}$
and $\gbigf_{\pm}^{(\infty,\infty)}$
are called the local Fourier transforms.

The local Fourier transforms were explicitly computed 
in \cite{Fang}, \cite{Graham-Squire} 
and \cite{Sabbah-stationary}.
We recall the explicit description
of the local Fourier transforms
by following Sabbah \cite{Sabbah-stationary}.

\subsection{\mbox{$\cnum(\!(z)\!)\langle \del_z\rangle$-modules}}
\label{subsection;24.3.13.10}

Let $(V,\nabla)$
be a $\cnum(\!(z)\!)\langle\del_z\rangle$-module of finite rank.
We assume that $V$ is finite dimensional over $\cnum(\!(z)\!)$.
There exists a so called Hukuhara-Levelt-Turrittin decomposition,
i.e.,
there exist a positive integer $p$,
a $\Gal(p)$-invariant subset $\nbigi(V)\subset z_p^{-1}\cnum[z_p]$
and a decomposition
\[
 (V,\nabla)\otimes\cnum(\!(z_p)\!)
 =\bigoplus_{\gminia\in\nbigi}
 (V_{\gminia},\nabla_{\gminia})
\]
such that 
$(V_{\gminia},\nabla_{\gminia}-d\gminia\id_{V_{\gminia}})$
are regular singular.
Let $\nbigi(V)=\coprod \gbigo$
denote the decomposition
into the $\Gal(p)$-orbits.
Because
$\bigoplus_{\gminia\in\gbigo}
 (V_{\gminia},\nabla_{\gminia})$
is $\Gal(p)$-equivariant,
there exists $\cnum(\!(z)\!)\langle\del_z\rangle$-module
$(V_{\gbigo},\nabla_{\gbigo})$
such that
$(V_{\gbigo},\nabla_{\gbigo})
\otimes\cnum(\!(z_p)\!)
=\bigoplus_{\gminia\in\gbigo}
 (V_{\gminia},\nabla_{\gminia})$.
We obtain the decomposition
\begin{equation}
\label{eq;24.2.23.30}
(V,\nabla)=
\bigoplus_{\gbigo\in \nbigi(V)/\Gal(p)}
(V_{\gbigo},\nabla_{\gbigo}).
\end{equation}

For an orbit $\gbigo$ in $\nbigi(V)$,
there exists
$\gminib=\sum_{j=1}^n\gminib_jz_p^{-j}$
such that
$\gbigo=\Gal(p)\cdot\gminib$.
We set
$r=\gcd\bigl(
\{j\,|\,\gminib_j\neq 0\}
\cup\{p\}
\bigr)$,
and $p_0:=p/r$.
Any $\gminia\in\gbigo$ is contained in $z_{p_0}^{-1}\cnum[z_{p_0}^{-1}]$,
and $(V_{\gminia},\nabla_{\gminia})$
are equivariant with respect
to the natural action of the Galois group
of $\cnum(\!(z_p)\!)$ over $\cnum(\!(z_{p_0})\!)$.
We obtain the decomposition
\[
 (V_{\gbigo},\nabla_{\gbigo})
\otimes\cnum(\!(z_{p_0})\!)
=\bigoplus_{\gminia\in\gbigo}
 (V'_{\gminia},\nabla'_{\gminia}).
\]
The action of $\Gal(p_0)$ on $\gbigo$
is free and transitive.
For any $\gminia\in\gbigo$,
we can naturally regard
$(V'_{\gminia},\nabla'_{\gminia})$
as a $\cnum(\!(z)\!)\langle\del_z\rangle$-module.
Then, it is isomorphic to 
$(V_{\gbigo},\nabla_{\gbigo})$
as a $\cnum(\!(z)\!)\langle\del_z\rangle$-module.
In other words,
$(V_{\gbigo},\nabla_{\gbigo})$
is isomorphic to the push-forward of
$(V'_{\gminia},\nabla'_{\gminia})$
via the ramified covering
$z_{p_0}\longmapsto z_{p_0}^{p_0}$.

For any $\gbigo$,
we set
$\gminim(V,\gbigo)=\rank(V_{\gbigo})/|\gbigo|$,
which equals $\rank V_{\gminia}'$ for any $\gminia\in\gbigo$.
By taking $p$ such that $\gbigo\subset z_p^{-1}\cnum[z_p^{-1}]$,
we set
$\deg_{z^{-1}}(\gbigo)=
(\deg_{z_p^{-1}}\gminia)/p$.

\subsection{Expression of $\gbigf_{\pm}^{(0,\infty)}$}
\label{subsection;24.3.26.1}

For any nonzero $\rho\in \zeta\cnum[\![\zeta]\!]$
and $\gminia\in \cnum(\!(\zeta)\!)\setminus\cnum[\![z]\!]$,
we set 
\[
 \widehat{\rho}^{(0)}_{\pm}(\zeta):=
 \mp\frac{\del_{\zeta}\rho(\zeta)}{\del_{\zeta}\gminia(\zeta)},
\quad
 \widehat{\gminia}^{(0)}_{\pm}(\zeta):=
 \gminia(\zeta)-
 \frac{\rho(\zeta)}{\del_{\zeta}\rho(\zeta)}\del_{\zeta}\gminia(\zeta)
=\gminia(\zeta)\pm
 \frac{\rho(\zeta)}{\widehat{\rho}_{\pm}^{(0)}(\zeta)}.
\]
For any non-zero $g=\sum g_j\zeta^j\in \cnum(\!(\zeta)\!)$,
let $\ord(g)$ denote the minimum of
$\{j\,|\,g_j\neq 0\}$.
\index{order $\ord(g)$}
Set $p:=\ord(\rho)$
and $n:=-\ord(\gminia)$.
Then, we obtain
$\ord\widehat{\rho}^{(0)}_{\pm}=n+p$
and $\ord(\widehat{\gminia}^{(0)}_{\pm})=-n$.
If $\rho=\zeta$ and $\gminia=0$,
we set $\rhohat=\zeta$ and $\gminiahat^{(0)}_{\pm}=0$.

Let $R$ be any regular singular differential
$\cnum(\!(\zeta)\!)$-module.
Let $V$ be 
the $\cnum(\!(z)\!)\langle\del_z\rangle$-module
induced by 
$(\cnum(\!(\zeta)\!),d+d\gminia)\otimes R$
and $z=\rho(\zeta)$,
i.e.,
$V$ is obtained as the push-forward of 
$(\cnum(\!(\zeta)\!),d+d\gminia)\otimes R$
by $\rho$.
Then, 
$\gbigf_{\pm}^{(0,\infty)}(V)$
is isomorphic to
the $\cnum(\!(w^{-1})\!)\langle\del_{w^{-1}}\rangle$-module
obtained as the push-forward of
$(\cnum(\!(\zeta)\!),
d+d\widehat{\gminia}^{(0)}_{\pm}+(n/2)d\zeta/\zeta)\otimes R$
by $w^{-1}=\widehat{\rho}^{(0)}_{\pm}(\zeta)$.
By using the decomposition (\ref{eq;24.2.23.30})
we obtain the explicit expression
of $\gbigf^{(0,\infty)}_+(V)$
for general $V$.

\subsubsection{}
\label{subsection;24.2.23.40}

We may regard the above construction as follows.
We explain only the case of $\gbigf^{(0,\infty)}_+$.
For simplicity, we assume
$\rho(\zeta)=\zeta^p$,
i.e., $\zeta$ is a $p$-th root of $z$.
We consider
$\gminia=\sum_{j=1}^n\gminia_j\zeta^{-j}$
with $\gminia_n\neq 0$.
Choose an $(n+p)$-th root $\eta$ of
the variable $u=w^{-1}$,
i.e.,
$\eta^{n+p}=u$.
We set
\begin{equation}
\label{eq;20.10.4.2}
 F_{\gminia,\eta}(\zeta):=
 \gminia(\zeta)+ \eta^{-n-p}\zeta^p.
\end{equation}
We fix an $(n+p)$-th root
$\bigl(\frac{n}{p}\gminia_n\bigr)^{\frac{1}{n+p}}$
of $\frac{n}{p}\gminia_n$.
\begin{lem}
\label{lem;24.2.23.20}
There exists a convergent power series
$g(\eta)$ such that
(i) $g(0)=\bigl(\frac{n}{p}\gminia_n\bigr)^{\frac{1}{n+p}}$,
(ii) $\zeta_0(\eta)=\eta\cdot g(\eta)$ satisfies
$\del_{\zeta}F_{\gminia,\eta}(\zeta_0(\eta))=0$.
 Moreover,
 the set of solutions of $\del_{\zeta}F_{\gminia,\eta}(\zeta)=0$
equals
$\bigl\{
 \zeta_0(a\eta)\,\big|\,
 a^{n+p}=1
 \bigr\}$.
\end{lem}
\pf
Because 
$\zeta^{n+1}\del_{\zeta}F_{\gminia,\eta}(\zeta)
=p(\zeta/\eta)^{n+p}
-\sum_{j=1}^{n-1}j\gminia_j\eta^{n-j}(\zeta/\eta)^{n-j}
-n\gminia_n$,
there exists a formal power solution $g(\eta)$
satisfying the conditions (i) and (ii).
Because this is an algebraic equation,
$g(\eta)$ is convergent.
Because the equation $\del_{\zeta}F_{\gminia,\eta}(\zeta)=0$
depends only on $\eta^{n+p}$,
the set of the solutions is 
$\bigl\{
 \zeta_0(a\eta)\,\big|\,
 a^{n+p}=1
 \bigr\}$.
\hfill\qed

\vspace{.1in}
The pull back of
$(\cnum(\!(\zeta)\!),d+d\gminiahat^{(0)}_++(n/2)d\zeta/\zeta)\otimes R$
by $\zeta_0(\eta)$ is isomorphic to
\[
 \Bigl(
 \cnum(\!(\eta)\!),
 d+dF_{\gminia,\eta}(\zeta_0(\eta))+(n/2)d\eta/\eta
 \Bigr)
 \otimes R.
\]
Let $\Gal(n+p)$ denote the Galois group of the extension
$\cnum(\!(\eta)\!)$ over $\cnum(\!(u)\!)$.
The pull back of $\gbigf_+^{(0,\infty)}(V)$
by the ramified covering $u=\eta^{n+p}$
is $\Gal(n+p)$-equivariantly isomorphic to
\[
 \bigoplus_{\substack{a\in\cnum \\ a^{n+p}=1}}
 \Bigl(
 \cnum(\!(\eta)\!),
 d+dF_{\gminia,\eta}(\zeta_0(a\eta))+(n/2)d\eta/\eta
 \Bigr)
 \otimes R.
\]

Let
$\gminia^{\circ}(\eta)
=\sum_{j=1}^{n}\gminia^{\circ}_j\eta^{-j}
\in \eta^{-1}\cnum[\eta^{-1}]$
denote the polar part of
$F_{\gminia,\eta}(\zeta_0(\eta))$.
We note that $\gminia^{\circ}_{n}\neq 0$.
The following lemma is also well known.
\begin{lem}
\label{lem;24.2.23.41}  
We have
$\gcd\bigl(\{j\,|\,\gminia_j\neq 0\}\cup\{p\}\bigr)=
\gcd\bigl(
\{j\,|\,\gminia^{\circ}_j\neq 0\}\cup\{n+p\}
\bigr)$. 
\end{lem}
\pf
Let us check it by a direct computation.
We set
\[
 a=\gcd\bigl(\{j\,|\,\gminia_j\neq 0\}\cup\{p\}\bigr),
 \quad
b=\gcd\bigl(
\{j\,|\,\gminia^{\circ}_j\neq 0\}\cup\{n+p\}
\bigr).
\]
Both $a$ and $b$ are divisors of $\gcd(n,p)$.
By the construction,
it is easy to see that $a$ is a divisor of $b$.
Assume that $b/a\in\seisuu_{>1}$,
and we shall deduce a contradiction.
We set
\[
 j_0:=\max\{j\in a\seisuu\setminus b\seisuu\,|\,\gminia_j\neq 0\}<n.
\]
Let $g(\eta)=\sum_{i=0}^{\infty}g_i\eta^i$
be the power series in Lemma \ref{lem;24.2.23.20}.
We obtain
$g_i=0$ for any $i\in a\seisuu\setminus b\seisuu$ such that $0<i<n-j_0$,
and
$g_{n-j_0}=\frac{j_0}{p(n+p)}\gminia_{j_0}g_0^{-j_0-p+1}$.
We obtain
$\gminia_{j_0}^{\circ}
=\gminia_{j_0}\cdot g_0^{-j_0}(n+p+j_0)(n+p)^{-1}\neq 0$.
It contradicts the definition of $b$.
Then, we obtain the claim of the lemma.
\hfill\qed

\subsubsection{Transform of the index sets}
\label{subsection;24.2.23.50}

Let $z_p$ be a $p$-th root of the variable $z$,
and let $u_{n+p}$ be an $(n+p)$-th root of the variable $u$.
Let $\gbigo$ be a $\Gal(p)$-orbit in $z_p^{-1}\cnum[z_p^{-1}]$.
If $\gbigo\neq\{0\}$,
applying the construction in \S\ref{subsection;24.2.23.40}
to $\gminia(z_p)\in\gbigo$ with $\zeta=z_p$ and $\eta=u_{n+p}$,
we construct
$\gminia^{\circ}(u_{n+p})\in u_{n+p}^{-1}\cnum[u_{n+p}^{-1}]$.
We set
\[
 \gbigf^{(0,\infty)}_+(\gbigo)
 =\Gal(n+p)\cdot \gminia^{\circ}
 \subset u_{n+p}^{-1}\cnum[u_{n+p}^{-1}].
\]
\begin{lem}
\label{lem;24.3.13.1}
$\gbigf^{(0,\infty)}_+(\gbigo)$
is independent of the choice of $\gminia\in\gbigo$.
\end{lem}
\pf
Because it is the transformation of the index sets
induced by the local Fourier transform,
it is independent of the choice of $\gminia\in\gbigo$.
We can also check it by a direct computation.
\hfill\qed

\vspace{.1in}
We also set
$\gbigf^{(0,\infty)}_+(\{0\})
=\{0\}$.
For any $\Gal(p)$-invariant subset $\nbigitilde$
of $z_p^{-1}\cnum[z_p^{-1}]$,
by using the decomposition into orbits
$\nbigitilde=\coprod\gbigo$,
we obtain a $\Gal(n+p)$-invariant subset
$\gbigf^{(0,\infty)}_+(\nbigitilde)
=\bigsqcup \gbigf^{(0,\infty)}_+(\gbigo)$.
\index{set $\gbigf^{(0,\infty)}_+(\nbigitilde)$}

The following lemma is also well known,
which follows from Lemma \ref{lem;24.2.23.41}.
\begin{lem}
We have
$\gminim(V,\gbigo)
=\gminim(\gbigf^{(0,\infty)}_+(V),\gbigf^{(0,\infty)}_+(\gbigo))$.
(See {\rm\S\ref{subsection;24.3.13.10}}
for $\gminim(V,\gbigo)$.)
\hfill\qed
\end{lem}

\subsection{Expression of $\gbigf_{\pm}^{(\infty,\infty)}$}
\label{subsection;24.4.16.1}

Take any non-zero $\rho\in \zeta\cnum[\![\zeta]\!]$
and $\gminia\in \cnum(\!(\zeta)\!)$
such that
$p:=\ord(\rho)<-\ord(\gminia)=:n$.
We set
\[
 \widehat{\rho}_{\pm}^{(\infty)}(\zeta):=
 \pm\frac{\del_{\zeta}\rho(\zeta)}{\gminia(\zeta)\rho(\zeta)^2},
\quad
 \widehat{\gminia}_{\pm}^{(\infty)}(\zeta):=
 \gminia(\zeta)+\frac{\rho(\zeta)}{\del_{\zeta}\rho(\zeta)}
 \del_{\zeta}\gminia(\zeta)
=\gminia(\zeta)
\pm\frac{1}{\rho(\zeta)\cdot\rhohat_{\pm}^{(\infty)}(\zeta)}.
\]
Let $R$ be any regular singular differential $\cnum(\!(\zeta)\!)$-module.
Let $V$ be the $\cnum(\!(z^{-1})\!)\langle\del_{z^{-1}}\rangle$-module
obtained as the push-forward of
$\bigl(\cnum(\!(\zeta)\!),d+d\gminia\bigr)\otimes R$ by 
$z^{-1}=\rho(\zeta)$.
Then,
$\gbigf_{\pm}^{(\infty,\infty)}(V)$
is isomorphic to
the push-forward of
$\bigl(\cnum(\!(\zeta)\!),d+d\gminiahat_{\pm}^{(\infty)}+(n/2)d\zeta/\zeta\bigr)
 \otimes R$ by 
$w^{-1}=\widehat{\rho}_{\pm}^{(\infty)}(\zeta)$.

Let $V$ be a general $\cnum(\!(z)\!)\langle\del_z\rangle$-module
which is finite dimensional over $\cnum(\!(z)\!)$.
We obtain the decomposition (\ref{eq;24.2.23.30}).
If $\deg_{z}(\gbigo)\leq 1$,
we have
$\gbigf^{(\infty,\infty)}_{\pm}(V_{\gbigo})=0$.
For $\deg_{z}\gbigo>1$,
we obtain the explicit expression of
$\gbigf^{(\infty,\infty)}_{\pm}(V_{\gbigo})$
by the above procedure.
In this way,
we obtain the explicit expression
of $\gbigf^{(\infty,\infty)}_+(V)$
for general $V$.

\subsubsection{}
\label{subsection;24.2.23.43}

We may regard the construction in the following way.
We explain the case of $\gbigf^{(\infty,\infty)}_+$.
For simplicity, we assume
$\rho(\zeta)=\zeta^p$,
i.e., $\zeta$ is a $p$-th root of $z^{-1}$.
We consider
$\gminia(\zeta)=\sum_{j=1}^n \gminia_j\zeta^{-j}$
with $\gminia_n\neq 0$
and $n>p$.
Choose an $(n-p)$-th root $\eta$ of $u=w^{-1}$.
We set
\begin{equation}
\label{eq;20.10.5.11}
 G_{\gminia,\eta}(\zeta):=
 \gminia(\zeta)+\frac{1}{\eta^{n-p}\zeta^p}.
\end{equation}
We fix an $(n-p)$-th root
$(-\frac{n}{p}\gminia_n)^{\frac{1}{n-p}}$
of $-\frac{n}{p}\gminia_n$.
\begin{lem}
There exists a convergent power series
$g(\eta)$ such that
(i) $g(0)=(-\frac{n}{p}\gminia_n)^{\frac{1}{n-p}}$,
(ii) $\zeta_0(\eta)=\eta g(\eta)$ satisfies
$\del_{\zeta}G_{\gminia,\eta}(\zeta_0(\eta))=0$.
Moreover,
the set of solutions of $\del_{\zeta}G_{\gminia,\eta}(\zeta)=0$
equals
$\{\zeta_0(a\eta)\,|\,a^{n-p}=1\}$. 
\hfill\qed
\end{lem}

The pull back of
$\gbigf^{(\infty,\infty)}_+(V)$
by $\eta\longmapsto\eta^{n-p}$
is isomorphic to
\[
 \bigoplus_{\substack{a\in\cnum\\ a^{n-p}=1}}
 \Bigl(
 \cnum(\!(\eta)\!),
 d+dG_{\gminia,\eta}(\zeta_0(a\eta))+(n/2)d\eta/\eta
 \Bigr)\otimes R. 
\]

Let
$\gminia^{\circ}(\eta)
=\sum_{j=1}^{n}\gminia^{\circ}_j\eta^{-j}
\in \eta^{-1}\cnum[\eta^{-1}]$
denote the polar part of
$G_{\gminia,\eta}(\zeta_0(\eta))$.
We note that $\gminia^{\circ}_{n}\neq 0$.
\begin{lem}
\label{lem;24.3.13.11}
We have
$\gcd\bigl(\{j\,|\,\gminia_j\neq 0\}\cup\{p\}\bigr)
=\gcd\bigl(
\{j\,|\,\gminia^{\circ}_j\neq 0\}\cup\{n-p\}\bigr)$. 
\hfill\qed
\end{lem}

\subsubsection{Transform of the index sets}
\label{subsection;24.2.23.51}

Let $x_p$ be a $p$-th root of the variable $x=z^{-1}$,
and let $u_{n-p}$ be an $(n-p)$-th root of the variable $u$.
Let $\gbigo$ be a $\Gal(p)$-orbit in $x_p^{-1}\cnum[x_p^{-1}]$
If $\deg_{x^{-1}}\gbigo\leq 1$ or $\gbigo=\{0\}$,
we set
$\gbigf^{(\infty,\infty)}_+(\gbigo)=\emptyset
\subset u_{n-p}^{-1}\cnum[u_{n-p}^{-1}]$.
If $\deg_{x^{-1}}\gbigo>1$,
by applying the construction in \S\ref{subsection;24.2.23.43}
to $\gminia\in\gbigo$
with $\zeta=x_p$ and $\eta=u_{n-p}$,
we obtain $\gminia^{\circ}\in u_{n-p}^{-1}\cnum[u_{n-p}^{-1}]$,
and we set
\[
 \gbigf^{(\infty,\infty)}_+(\gbigo)
 =\Gal(n-p)\cdot \gminia^{\circ}(u_{n-p})
 \subset
 u_{n-p}^{-1}\cnum[u_{n-p}^{-1}].
\]
The following lemma is similar to Lemma \ref{lem;24.3.13.1}.
\begin{lem}
$\gbigf^{(\infty,\infty)}_+(\gbigo)$
is independent of the choice of $\gminia\in\gbigo$.
\hfill\qed
\end{lem}
For any $\Gal(p)$-invariant subset $\nbigitilde$
of $x_p^{-1}\cnum[x_p^{-1}]$,
we define $\gbigf^{(\infty,\infty)}_+(\nbigitilde)$
by using the orbit decomposition of $\nbigitilde$.
\index{set $\gbigf^{(\infty,\infty)}_+(\nbigitilde)$}

The following lemma is also well known
which follows from Lemma \ref{lem;24.3.13.11}.
\begin{lem}
If $\deg_{x^{-1}}\gbigo>1$,
we have
$\gminim(V,\gbigo)
=\gminim(\gbigf^{(\infty,\infty)}_+(V),\gbigf^{(\infty,\infty)}_+(\gbigo))$.
\hfill\qed
\end{lem}

\section{Notation}

For any $\vartheta_0\in\real$
and $L>0$,
we set $I(\vartheta_0,L):=
\bigl\{
\theta\in\real\,\big|\,
|\vartheta_0-\theta|<L
\bigr\}$.
\index{interval $I(\vartheta_0,L)$}

Let $n$ and $p$ be any positive integers.
We take a $p$-th root $z_p$ of $z$.
We set  \index{set $\gbigu_z(p,n)$}
\[
\gbigu_z(p,n):=z_p^{-n}\cnum\setminus\{0\}.
\]
For any interval $J\subset\real$,
we set \index{sets $\gbigu^{\pm}_z(p,n,J)$}
\[
\gbigu_z^-(p,n,J):=\bigl\{
 \gminia\in\gbigu_z(p,n)\,\big|\,
 \gminia<_J0
 \bigr\},
\quad
\gbigu_z^+(p,n,J):=\bigl\{
 \gminia\in\gbigu_z(p,n)\,\big|\,
 \gminia>_J0
 \bigr\}.
\]
We also set \index{set $\gbigutilde_z(p,n)$}
\[
\gbigutilde_z(p,n):=
 \left\{
 \sum_{j=1}^{n}\gminia_jz_p^{-j}\,\big|\,
 \gminia_n\neq 0
 \right\}
 \subset \cnum[z_p^{-1}]
 \simeq
 \cnum(\!(z_{p})\!)\big/
 \cnum[\![z_{p}]\!].
\]
There exists the natural map
$q_{z,p,n}:
 \gbigutilde_z(p,n)\lrarr \gbigu_z(p,n)$
defined by
$\sum\gminia_jz_p^{-j}
\longmapsto
 \gminia_nz_p^{-n}$.
For any interval $J\subset\real$,
we set \index{maps $q_{z,p,n}$} \index{sets $\gbigutilde_z^{\pm}(p,n,J)$} 
\[
\gbigutilde_z^{\pm}(p,n,J):=
 q_{z,p,n}^{-1}\bigl(
 \gbigu_z^{\pm}(p,n,J)
 \bigr).
\]

\section{From $0$ to $\infty$}

We shall refine the construction
in \S\ref{subsection;24.2.23.40}--\S\ref{subsection;24.2.23.50}.

\subsection{Preliminary computations (1)}
\label{subsection;18.4.6.10}

Let $n$ and $p$ be any positive integers.
We set $\omega:=n/p$,
and
$\langle\omega\rangle:=
\omega^{\frac{-\omega}{1+\omega}}
+\omega^{\frac{1}{1+\omega}}$.
\index{number $\langle\omega\rangle$}

Let $\alpha$ be a non-zero complex number.
We set $\gminia:=\alpha\zeta^{-n}$,
and we consider
\[
 F_{\gminia,\eta}(\zeta)=
 \gminia(\zeta)+\eta^{-n-p}\zeta^p
=\alpha\zeta^{-n}+\eta^{-n-p}\zeta^p.
\]

We obtain
$\del_{\zeta}F_{\gminia,\eta}
=-n\alpha\zeta^{-n-1}+p\zeta^{p-1}\eta^{-n-p}$.
By setting $h_{\omega}(\eta):=\omega^{\frac{1}{n+p}}\eta$,
we obtain
\[
\del_{\zeta}F_{\gminia,\eta}(\zeta_0)=0
\Longleftrightarrow
 \zeta_0\in
 \bigl\{
 h_{\omega}(\beta\eta)
 \big|\,
 \beta\in\cnum,\,
 \beta^{n+p}=\alpha
 \bigr\}.
\]

Set $J=I(\vartheta_0,\omega^{-1}\pi/2)$.
For $m\in\seisuu$,
we define the intervals $J^u(m,\pm)$
as follows:
\index{interval $J^u(m,\pm)$}
\begin{equation}
\label{eq;24.2.24.1}
\left\{
\begin{array}{l}
 J^u(m,-)
 =I\bigl(\vartheta_0-2m\pi,(1+\omega^{-1})\pi/2\bigr),
\\
  J^u(m,+)
 =I\bigl(\vartheta_0-(2m-1)\pi,(1+\omega^{-1})\pi/2\bigr).
\end{array}
\right.
\end{equation}

\subsubsection{Case 1}

Suppose that
$\Re\gminia(e^{\sqrt{-1}\theta/p})>0$ for any  $\theta \in J$.
We obtain
\begin{equation}
\label{eq;18.4.6.20}
 \alpha
 =|\alpha|
 \exp\bigl(
 \sqrt{-1}\omega\vartheta_0
 \bigr).
\end{equation}
For any $m\in\seisuu$,
we set
\begin{equation}
\label{eq;20.10.4.3}
 \beta_{J,m,-}:=|\alpha|^{\frac{1}{n+p}}
  \exp\Bigl(
 \frac{\sqrt{-1}}{n+p}
 \bigl(
 \omega\vartheta_0
 +2m\pi
 \bigr)
 \Bigr)
\in \bigl\{\beta\in\cnum\,\big|\,\beta^{n+p}=\alpha\bigr\}.
\end{equation}
The set of roots of $\del_{\zeta}F_{\gminia,\eta}$
is
$\{h(\beta_{J,m,-}\eta)\,|\,m\in\seisuu\}$.
We set \index{map $\gbigf^{(0,\infty)}_{(J,m,-)}$}
\begin{multline}
\label{eq;18.5.3.12}
 \gbigf^{(0,\infty)}_{(J,m,-)}(\gminia)(\eta):=
F_{\gminia,\eta}\bigl(h_{\omega}(\beta_{J,m,-}\eta)\bigr)
=\langle\omega\rangle
 |\alpha|^{1/(1+\omega)}
 \exp\Bigl(
 \frac{\sqrt{-1}}{1+\omega}
 \bigl(
 \omega\vartheta_0
 +2m\pi
 \bigr)
 \Bigr)\cdot \eta^{-n}
 \\
 =\langle\omega\rangle
 |\alpha|^{1/(1+\omega)}
 \exp\Bigl(
 \sqrt{-1}\Bigl(
 \frac{\omega}{1+\omega}
 (\vartheta_0-2m\pi)
 \Bigr)
 \Bigr)\eta^{-n}.
\end{multline}
\begin{lem}
\label{lem;18.5.3.10}
We set
 \[
 \arg\bigl(
 h_{\omega}(\beta_{J,m,-}e^{\sqrt{-1}\theta^u/(n+p)})
 \bigr)
= \frac{1}{p(1+\omega)}
 \bigl(
 \theta^u
+ \omega\vartheta_0
+2m\pi
 \bigr).
\] 
Then,
$p\cdot\arg\bigl(
 h_{\omega}(\beta_{J,m,-}e^{\sqrt{-1}\theta^u/(n+p)})
 \bigr)
\in J$ 
if and only if $\theta^u\in J^u(m,-)$.
Moreover,
we obtain
$\Re\gbigf^{(0,\infty)}_{(J,m,-)}(\gminia)(\eta)>0$
for $\eta=|\eta|e^{\sqrt{-1}\theta^u/(n+p)}$
with $\theta^u\in J^u(m,-)$.
\end{lem}
\pf
We obtain
\[
 p\arg\bigl(
  h_{\omega}(\beta_{J,m,-}e^{\sqrt{-1}\theta^u/(n+p)})
 \bigr)
 -\vartheta_0
 =\frac{1}{1+\omega}
 (\theta^u-\vartheta_0+2m\pi).
\]
Then, the first claim is clear.
The second claim is clear from (\ref{eq;18.5.3.12}).
\hfill\qed

\vspace{.1in}
We note the following obvious lemma.
\begin{lem}
\label{lem;18.5.24.1}
For any integer $\ell$,
we obtain
$(-1)^{\ell}\Re(\gminia(e^{\sqrt{-1}\theta/p}))>0$
on $J+\ell \omega^{-1}\pi$,
and 
$(-1)^{\ell}\Re\gbigf^{(0,\infty)}_{(J,m,-)}
(\gminia)(e^{\sqrt{-1}\theta^u/(n+p)})>0$
on $J^u(m,-)+\ell(1+\omega^{-1})\pi$.
\hfill\qed
\end{lem}

\subsubsection{Case 2}

Suppose that
$\Re\gminia(e^{\sqrt{-1}\theta/p})<0$ for any $\theta \in J$.
We obtain
\begin{equation}
 \label{eq;18.4.6.21}
  \alpha
=-|\alpha|
 \exp\bigl(
 \sqrt{-1}\omega\vartheta_0
 \bigr)
=|\alpha|  
 \exp\bigl(
 \sqrt{-1}(\omega\vartheta_0-\pi)
 \bigr).
\end{equation}
For any $m\in\seisuu$,
we set
\begin{equation}
\label{eq;20.10.4.4}
 \beta_{J,m,+}:=
 |\alpha|^{1/(n+p)}
  \exp\Bigl(
 \frac{\sqrt{-1}}{n+p}
 \bigl(
 \omega\vartheta_0
+(2m-1)\pi
\bigr)
\Bigr)
\in\{\beta\in\cnum\,|\,\beta^{n+p}=\alpha\}.
\end{equation}
The set of roots of $\del_{\zeta}F_{\gminia,\eta}$
is
$\bigl\{
 h_{\omega}(\beta_{J,m,+}\eta)\,\big|\,m\in\seisuu
\bigr\}$.
We set \index{map $\gbigf^{(0,\infty)}_{(J,m,+)}$}
\begin{multline}
\label{eq;18.5.3.13}
 \gbigf^{(0,\infty)}_{(J,m,+)}(\gminia)(\eta):=
 F_{\gminia,\eta}\bigl(
 h_{\omega}(\beta_{J,m,+}\eta)
 \bigr) 
\\
 =\langle\omega\rangle
 |\alpha|^{1/(1+\omega)}
 \exp\Bigl(
 \frac{\sqrt{-1}}{1+\omega}
 \bigl(
 \omega\vartheta_0
+(2m-1)\pi
 \bigr)
 \Bigr)\cdot \eta^{-n}
 \\
 =-\langle \omega\rangle
 |\alpha|^{1/(1+\omega)}
 \exp\Bigl(
 \sqrt{-1}\Bigl(
 \frac{\omega}{1+\omega}
 \bigl(\vartheta_0-(2m-1)\pi\bigr)
 \Bigr)
 \Bigr)\eta^{-n}.
\end{multline}

The following lemma is similar to
Lemma \ref{lem;18.5.3.10}.
\begin{lem}
\label{lem;18.5.3.11}
We set
\[
 \arg\bigl(
 h_{\omega}(\beta_{J,m,+}e^{\sqrt{-1}\theta^u/(n+p)})
 \bigr)
=\frac{1}{p(1+\omega)}\bigl(
 \theta^u+\omega\vartheta_0
+(2m-1)\pi
 \bigr).
\]
Then, 
 $p\cdot\arg\bigl(
 h_{\omega}(\beta_{J,m,+}e^{\sqrt{-1}\theta^u/(n+p)})
 \bigr)
\in J$
if and only if $\theta^u\in J^u(m,+)$.
Moreover,
we obtain
$\Re\gbigf^{(0,\infty)}_{(J,m,+)}(\gminia)(\eta)<0$
for $\eta=|\eta|e^{\sqrt{-1}\theta^u/(n+p)}$
with  $\theta^u\in J^u(m,+)$.
\hfill\qed
\end{lem}

We note the following obvious lemma.
\begin{lem}
\label{lem;18.5.24.2}
We obtain
$(-1)^{\ell}\Re(\gminia(e^{\sqrt{-1}\theta/p}))<0$
on $J+\ell\omega^{-1} \pi$,
and 
\[
 (-1)^{\ell}\Re\gbigf^{(0,\infty)}_{(J,m,+)}
(\gminia)(e^{\sqrt{-1}\theta^u/(n+p)})<0
\]
on $J^u(m,+)+\ell(1+\omega^{-1})\pi$
for any integer $\ell$.
\hfill\qed
\end{lem}

\subsection{Preliminary computations (2)}
\label{subsection;20.10.5.10}

For $\gminiatilde=
 \alpha \zeta^{-n}+\sum_{j=1}^{n-1}\gminiatilde_j\zeta^{-j}$,
we set
\[
 F_{\gminiatilde,\eta}(\zeta):=
 \gminiatilde(\zeta)+\eta^{-n-p}\zeta^p.
\]
We obtain
$\del_{\zeta}F_{\gminiatilde,\eta}(\zeta)
=-n\alpha\zeta^{-n-1}+p\zeta^{p-1}\eta^{-n-p}
-\sum_{j=1}^n j\gminiatilde_j\zeta^{-j-1}$.
The following lemma is 
a reformulation of Lemma \ref{lem;24.2.23.20}.

 \begin{lem}
\label{lem;18.4.13.10}
There exists a unique convergent power series
$a_{\gminiatilde}(\eta)=
 1+\sum_{j=1}^{\infty} a_{\gminiatilde,j}\eta^j$
such that 
the following holds for any $\beta\in\cnum$
with $\beta^{n+p}=\alpha$:
\begin{equation}
\label{eq;18.4.12.1}
 \del_{\zeta}F_{\gminiatilde,\eta}
 \Bigl(
 h_{\omega}(\beta\eta)a_{\gminiatilde}(\beta\eta)
 \Bigr)=0.
\end{equation}
Moreover, any root of $\del_{\zeta}F_{\gminiatilde,\eta}$
is described as
$h_{\omega}(\beta\eta)a_{\gminiatilde}(\beta\eta)$
for some $\beta$ with $\beta^{n+p}=\alpha$.
 \end{lem}
\pf
The condition (\ref{eq;18.4.12.1})
is equivalent to
\begin{equation}
 p\alpha\omega^{p/(n+p)}
 \Bigl(
 -a_{\gminiatilde}(\beta\eta)^{-n}+a_{\gminiatilde}(\beta\eta)^{p}
 \Bigr)
-\sum_{j=1}^{n-1}
 j\gminiatilde_j\omega^{-j/(n+p)}
  (\beta\eta)^{n-j}a_{\gminiatilde}(\beta\eta)^{-j}=0.
\end{equation}
Here, we have used $\beta^{n+p}=\alpha$.
Hence, it is enough to study the following equation:
\begin{equation}
\label{eq;18.4.12.2}
  p\alpha\omega^{p/(n+p)}
 \Bigl(
 -a_{\gminiatilde}(\eta)^{-n}+a_{\gminiatilde}(\eta)^{p}
 \Bigr)
-\sum_{j=1}^{n-1}
 j\gminiatilde_j\omega^{-j/(n+p)}
  \eta^{n-j}a_{\gminiatilde}(\eta)^{-j}=0.
\end{equation}
It is easy to check that
there exists a unique formal power series
$a_{\gminiatilde}$  of the desired form
satisfying (\ref{eq;18.4.12.2}).
It is convergent because it is a solution
of the algebraic equation
(\ref{eq;18.4.12.2}).
\hfill\qed

\vspace{.1in}
There exists the convergent power series
$1+\sum_{j=1}^{\infty} b_{\gminiatilde,j}\eta^j$
such that 
\begin{multline}
\label{eq;20.10.5.1}
 F_{\gminiatilde,\eta}\bigl(
 h_{\omega}(\beta\eta)\cdot a_{\gminiatilde}(\beta\eta)
 \bigr)
=F_{\gminia,\eta}\bigl(
 h_{\omega}(\beta\eta)
 \bigr)\cdot
 \Bigl(
 1+\sum_{j=1}^{\infty}b_{\gminiatilde,j}(\beta\eta)^j
 \Bigr) \\
 =\langle\omega\rangle
 \alpha
 (\beta\eta)^{-n}\cdot
  \Bigl(
 1+\sum_{j=1}^{\infty}b_{\gminiatilde,j}(\beta\eta)^j
 \Bigr).
\end{multline}

\begin{lem}
\label{lem;18.4.13.1}
Suppose that 
$\gminiatilde_i=
 \alpha\zeta^{-n}
+\sum_{j=1}^{n-1}\gminiatilde_{i,j}\zeta^{-j}$ $(i=1,2)$
satisfy
 $\gminiatilde_{1,j}=\gminiatilde_{2,j}$ $(j=k+1,\ldots,n-1)$
for some $0\leq k\leq n-1$.
Then, we obtain
\begin{equation}
\label{eq;18.4.12.10}
 b_{\gminiatilde_1,j}=b_{\gminiatilde_2,j}
\quad
 (j=1,\ldots,n-k-1).
\end{equation}
Moreover, we obtain
\[
 F_{\gminiatilde_1,\eta}\bigl(
 h_{\omega}(\beta\eta)\cdot a_{\gminiatilde_1}(\beta\eta)
 \bigr)
- F_{\gminiatilde_2,\eta}\bigl(
 h_{\omega}(\beta\eta)\cdot a_{\gminiatilde_2}(\beta\eta)
 \bigr)
\equiv
 \bigl(
 \gminiatilde_{1,k}
-\gminiatilde_{2,k}
 \bigr)\cdot h_{\omega}(\beta\eta)^{-k}
\]
modulo $\eta^{-k+1}\cnum[\![\eta]\!]$.
\end{lem}
\pf
For $1\leq \ell<n$,
we can easily observe that
$a_{\gminiatilde,\ell}$ depend
only on $\gminiatilde_{n-j}$ $(1\leq j \leq \ell)$.
Hence, we obtain (\ref{eq;18.4.12.10}).
Moreover, 
the dependence of $a_{\gminiatilde,\ell}$
on $\gminiatilde_{n-\ell}$ is linear,
i.e.,
$a_{\gminiatilde,\ell}=
 A_{\ell}\gminiatilde_{n-\ell}+
 Q_{\ell}(\gminiatilde_{n-1},\ldots,\gminiatilde_{n-\ell+1})$.
Hence, the following holds modulo $\eta^{-k+1}$:
\begin{multline}
 F_{\gminiatilde_1,\eta}\bigl(
 h_{\omega}(\beta\eta)\cdot a_{\gminiatilde_1}(\beta\eta)
 \bigr)
- F_{\gminiatilde_2,\eta}\bigl(
 h_{\omega}(\beta\eta)\cdot a_{\gminiatilde_2}(\beta\eta)
 \bigr)
\equiv\\
 \del_{\zeta}F_{\gminia_1,\eta}
 \bigl(h_{\omega}(\beta\eta)
 a_{\gminiatilde_2}(\beta y)
 \bigr)\cdot 
 h_{\omega}(\beta\eta)\cdot
 A_{n-k}\cdot(\gminiatilde_{1,k}-\gminiatilde_{2,k})\eta^{k}
 \\
+\gminiatilde_{1,k}
 \bigl(
 h_{\omega}(\beta\eta)a_{\gminiatilde_2}(\beta\eta)
 \bigr)^{-k}
-\gminiatilde_{2,k}
 \bigl(
 h_{\omega}(\beta\eta)a_{\gminiatilde_2}(\beta\eta)
 \bigr)^{-k}
\\
\equiv
 \bigl(
 \gminiatilde_{1,k}
-\gminiatilde_{2,k}
 \bigr)\cdot h_{\omega}(\beta\eta)^{-k}.
\end{multline}
We have used
$\del_{\zeta}F_{\gminia_1,\eta}
 \bigl(h_{\omega}(\beta\eta)
 \bigr)=0$.
Thus, we obtain the claim of the lemma.
\hfill\qed

\subsection{Direct consequences of preliminary computations}
\label{subsection;20.10.23.1}

We take a $p$-th root $z_p$ of $z$,
and an $(n+p)$-th root $u_{n+p}$ of $u=w^{-1}$.
We set $J=I(\vartheta_0,\omega^{-1}\pi/2)$.
We define the maps
$\gbigf^{(0,\infty)}_{(J,m,\pm)}:
 \gbigu^{\pm}_z(p,n,J)
\lrarr
 \gbigu_u(n+p,n)$
by the formulas (\ref{eq;18.5.3.12})
and (\ref{eq;18.5.3.13})
with $\zeta=z_{p}$ and $\eta=u_{n+p}$.
We obtain the following lemma from
Lemma \ref{lem;18.5.3.10} and Lemma \ref{lem;18.5.3.11}.
\begin{lem}
\label{lem;20.10.4.5}
They induce the following isomorphisms 
of the partially ordered sets:
\[
 \gbigf^{(0,\infty)}_{(J,m,\pm)}:
(\gbigu_z^{\pm}(p,n,J),\leq_J)
\simeq
\bigl(\gbigu_u^{\pm}(p+n,n,J^u(m,\pm)),\leq_{J^u(m,\pm)}\bigr).
\] 
\hfill\qed
\end{lem}

For
$\gminiatilde(z_p)
=\gminiatilde_n z_p^{-n}+\sum_{j=1}^{n-1}\gminiatilde_jz_p^{-j}
\in\gbigu_z^-(p,n,J)$,
we set
\[
 \beta_{\gminiatilde,J,m,-}:=
|\gminiatilde_n|^{1/(n+p)}
\exp\Bigl(
\frac{\sqrt{-1}}{n+p}
(\omega\vartheta_0+2m\pi)
\Bigr)
\]
as in (\ref{eq;20.10.4.3})
with $\alpha=\gminiatilde_n$.
Similarly,
for 
$\gminiatilde(z_p)
\in\gbigu_z^+(p,n,J)$,
we set
\[
\beta_{\gminiatilde,J,m,+}:=
|\gminiatilde_n|^{1/(n+p)}
\exp\Bigl(
\frac{\sqrt{-1}}{n+p}
\bigl(\omega\vartheta_0+(2m-1)\pi\bigr)
\Bigr)
\]
as in (\ref{eq;20.10.4.4}).
For $\gminiatilde\in \gbigutilde^{\pm}_z(p,n)$,
we define
$\gbigf^{(0,\infty)}_{(J,m,\pm)}(\gminiatilde)
\in u_{n+p}^{-1}\cnum[u_{n+p}^{-1}]$ by
\index{maps $\gbigf^{(0,\infty)}_{(J,m,\pm)}$}
\[
\gbigf^{(0,\infty)}_{(J,m,\pm)}(\gminiatilde)
:=F_{\gminiatilde,u_{p+n}}\Bigl(
 h_{\omega}(\beta_{\gminiatilde,J,m,\pm}u_{n+p})
 a_{\gminiatilde}(\beta_{\gminiatilde,J,m,\pm}u_{n+p})
 \Bigr)\,\,\,
\modulo\,\,\, \cnum[\![u_{n+p}]\!],
\]
where we set
$h_{\omega}(\beta_{\gminiatilde,J,m,\pm}\eta)
=\omega^{1/(n+p)}\beta_{\gminiatilde,J,m,\pm}\eta$,
and $a_{\gminiatilde}(\eta)$
are the convergent power series
in Lemma \ref{lem;18.4.13.10}.
Thus, we obtain  the maps
$\gbigf^{(0,\infty)}_{(J,m,\pm)}:
 \gbigutilde_z^{\pm}(p,n,J)\lrarr
 \gbigutilde_u(p+n,n)$.
By (\ref{eq;20.10.5.1}),
the following holds:
\begin{equation}
\label{eq;20.10.4.6}
 \gbigf^{(0,\infty)}_{(J,m,\pm)}
 \bigl(q_{z,p,n}(\gminiatilde)\bigr)
 =q_{u,n+p,n}(\gbigf^{(0,\infty)}_{J,m,\pm}
 \bigl(\gminiatilde)\bigr).
\end{equation}
 
\begin{lem}
\label{lem;20.10.5.23}
The maps $\gbigf^{(0,\infty)}_{(J,m,\pm)}$
induce the following bijections:
\[
 \gbigf^{(0,\infty)}_{(J,m,\pm)}:
 \gbigutilde^{\pm}_{z}(p,n,J)
\simeq
 \gbigutilde^{\pm}_{u}(p+n,n,J^u(m,\pm)).
\]
\end{lem}
\pf
By the formula (\ref{eq;20.10.4.6}),
we obtain the map
\begin{equation}
\label{eq;20.10.5.2}
 \gbigf^{(0,\infty)}_{(J,m,\pm)}:
 \gbigutilde^{\pm}_{z}(p,n,J)
\lrarr
 \gbigutilde^{\pm}_{u}(p+n,n,J^u(m,\pm)).
\end{equation}
By Lemma \ref{lem;18.4.13.1},
the map (\ref{eq;20.10.5.2}) is bijective.
\hfill\qed

\vspace{.1in}

For $\theta\in\real$,
we define $\theta^u(m,\pm)$
as follows: \index{points $\theta^u(m,\pm)$}
\[
 \theta^u(m,-):=(1+\omega)\theta-\omega\vartheta_0-2m\pi,
 \quad
 \theta^u(m,+):=(1+\omega)\theta-\omega\vartheta_0-(2m-1)\pi.
\]

\begin{prop}
\label{prop;20.10.4.10}
$\gbigf^{(0,\infty)}_{(J,m,\pm)}$ induce
the following isomorphisms of the partially
ordered sets
for any $\theta\in \real$:
\[
 \gbigf^{(0,\infty)}_{(J,m,\pm)}:
 \bigl(
 \gbigutilde^{\pm}_{z}(p,n,J),\leq_{\theta}
 \bigr)
\simeq
 \bigl(
 \gbigutilde^{\pm}_{u}(p+n,n,J^u(m,\pm)),\leq_{\theta^u(m,\pm)}
 \bigr).
\]
We also obtain the following commutative diagram:
{\small
 \begin{equation}
\label{eq;20.10.5.3}
 \begin{CD}
  \gbigu^{\pm}_{z}(p,n,J)
  @>{\iota_z}>>
  \gbigutilde^{\pm}_{z}(p,n,J)
  @>{q_{z,p,n}}>>
  \gbigu^{\pm}_{z}(p,n,J)\\
  @V{\gbigf^{(0,\infty)}_{(J,m,\pm)}}V{\simeq}V
  @V{\gbigf^{(0,\infty)}_{(J,m,\pm)}}V{\simeq}V
  @V{\gbigf^{(0,\infty)}_{(J,m,\pm)}}V{\simeq}V \\
  \gbigu^{\pm}_{u}(p+n,n,J^u(m,\pm))
  @>{\iota_u}>>
  \gbigutilde^{\pm}_{u}(p+n,n,J^u(m,\pm))
  @>{q_{u,p+n,n}}>>
  \gbigu^{\pm}_{u}(p+n,n,J^u(m,\pm)).
 \end{CD}
 \end{equation}
}
Here, $\iota_z$ and $\iota_u$
denote the natural inclusions.
\end{prop}
\pf
Let $\gminia_i=\sum_{j=1}^n\gminia_{i,j}z_p^{-j}
\in\gbigutilde^{\pm}_z(p,n,J)$.
We shall prove that
$\gminia_1<_{\theta}\gminia_2$
if and only if
$\gbigf^{(0,\infty)}_{(J,m,\pm)}(\gminia_1)
<_{\theta^u(m,\pm)}
 \gbigf^{(0,\infty)}_{(J,m,\pm)}(\gminia_2)$.
First,
let us study the case $\gminia_{1,n}\neq \gminia_{2,n}$.
Note that for each integer $\ell$,
we have
$\theta\in J+\ell \omega^{-1}\pi$
if and only if
$\theta^u(m,\pm)\in J^u(m,\pm)
 +\ell (1+\omega^{-1})\pi$.
Hence, by Lemma \ref{lem;18.5.24.1} and
Lemma \ref{lem;18.5.24.2},
we have
$\gminia_1<_{\theta}\gminia_2$
if and only if
$\gbigf^{(0,\infty)}_{(J,m,\pm)}(\gminia_1)
<_{\theta^u(m,\pm)}
 \gbigf^{(0,\infty)}_{(J,m,\pm)}(\gminia_2)$.
Second,
let us study the case
$\gminia_{1,n}=\gminia_{2,n}$.
We obtain $\beta_{J,m,\pm}$
as in (\ref{eq;20.10.4.3}) or (\ref{eq;20.10.4.4})
with $\alpha=\gminia_{1,n}=\gminia_{2,n}$.
We note that
\[
\theta=
 p\arg\bigl(
 h_{\omega}(\beta_{J,m,\pm}e^{\sqrt{-1}\theta^u(m,\pm)/(n+p)})
 \bigr).
\]
Then, by Lemma \ref{lem;18.4.13.1},
we obtain
$\gminia_1<_{\theta}\gminia_2$
if and only if
$\gbigf^{(0,\infty)}_{(J,m,\pm)}(\gminia_1)
<_{\theta^u(m,\pm)}
 \gbigf^{(0,\infty)}_{(J,m,\pm)}(\gminia_2)$.

We obtain the commutativity of (\ref{eq;20.10.5.3})
by the construction
of $\gbigf^{(0,\infty)}_{(J,m,\pm)}$
and (\ref{eq;20.10.4.6}).
\hfill\qed

\vspace{.1in}

Note that
$\gbigutilde^{\pm}_z(p,n,J+\omega^{-1}\pi)
=\gbigutilde^{\mp}_z(p,n,J)$.
We also note that
$(J+\omega^{-1}\pi)^u(m-1,-)
=J^u(m,+)+(1+\omega^{-1})\pi$
and 
$(J+\omega^{-1}\pi)^u(m,+)
=J^u(m,-)+(1+\omega^{-1})\pi$,
which imply
\[
 \gbigutilde^{-}_u\bigl(
 p+n,n,(J+\omega^{-1}\pi)^u(m-1,-)
 \bigr)
 =\gbigutilde^+_u\bigl(
 p+n,n,J^u(m,+)
 \bigr),
\]
\[
 \gbigutilde^{+}_u\bigl(
 p+n,n,(J+\omega^{-1}\pi)^u(m,+)
 \bigr)
 =\gbigutilde^-_u\bigl(
 p+n,n,J^u(m,-)
 \bigr).
\]
The following lemma can be checked by computation.
\begin{lem}
\label{lem;20.10.3.1}
 For $\gminiatilde\in \gbigutilde^{+}_z(p,n,J)$,
we obtain
$\gbigf^{(0,\infty)}_{(J,m,+)}(\gminiatilde)
\!=\!\gbigf^{(0,\infty)}_{(J+\omega^{-1}\pi,m-1,-)}(\gminiatilde)$.
For $\gminiatilde\in \gbigutilde^{-}_z(p,n,J)$,
we obtain
$\gbigf^{(0,\infty)}_{(J,m,-)}(\gminiatilde)
 =\gbigf^{(0,\infty)}_{(J+\omega^{-1}\pi,m,+)}(\gminiatilde)$.
\hfill\qed
\end{lem}

\subsection{Reformulation}

Let $\vecJ=I(\vartheta^u_0,(1+\omega^{-1})\pi/2)$.
For any integer $m$,
we set \index{intervals $\nu^{\pm}(\vecJ)$}
\begin{equation}
\label{eq;20.10.13.1}
\nu^-_m(\vecJ)=
I(\vartheta^u_0+2m\pi,\omega^{-1}\pi/2),
\quad
\nu^+_m(\vecJ)=
I(\vartheta^u_0+(2m-1)\pi,\omega^{-1}\pi/2).
\end{equation}
We define the isomorphisms of the partially ordered sets
\index{isomorphisms $\nu^{\pm}_{m,\vecJ}$}
\[
\nu^{\pm}_{m,\vecJ}:
 \bigl(
 \gbigu^{\pm}_{u}(n+p,n,\vecJ),\leq_{\vecJ}
 \bigr)
\lrarr
 \bigl(
 \gbigu^{\pm}_{z}(p,n,\nu^{\pm}_m(\vecJ)),\leq_{\nu^{\pm}_m(\vecJ)}
\bigr) 
\]
as the inverse of
$\gbigf^{(0,\infty)}_{(\nu^{\pm}_m(\vecJ),m,\pm)}$.
We define the bijections
\index{isomorphisms $\nutilde^{\pm}_{m,\vecJ}$}
\[
  \nutilde^{\pm}_{m,\vecJ}:
 \gbigutilde^{\pm}_{u}\bigl(
 p+n,n,\vecJ
 \bigr)
 \simeq
 \gbigutilde^{\pm}_z\bigl(
 p,n,\nu^{\pm}_{m}(\vecJ)
 \bigr)
\]
as the inverse of 
$\gbigf^{(0,\infty)}_{(\nu^{\pm}_m(\vecJ),m,\pm)}$.
We define the maps
$\kappa_{m,\vecJ}^{\pm}:\real\lrarr\real$ by the following formulas:
\index{maps $\kappa_{m,\vecJ}^{\pm}$}
\begin{equation}
\label{eq;20.10.9.1}
 \kappa_{m,\vecJ}^-(\theta^u)
 =\frac{1}{1+\omega}(\theta^u+\omega\vartheta^u_0)+2m\pi,
\quad
 \kappa_{m,\vecJ}^+(\theta^u)
 =\frac{1}{1+\omega}(\theta^u+\omega\vartheta^u_0)+(2m-1)\pi.
\end{equation}
Note that
$\kappa_{m,\vecJ}^{\pm}$
induces bijections
$\vecJ+\ell(1+\omega^{-1})\pi\simeq
\nu_m^{\pm}(\vecJ)+\ell \omega^{-1}\pi$
for each integer $\ell$.

\begin{prop}
\label{prop;20.10.4.11}
The maps
$\nutilde^{\pm}_{m,\vecJ}$ induce isomorphisms of
the following partially ordered sets
for any $\theta^u\in\real$:
\[
\nutilde^{\pm}_{m,\vecJ}:
(\gbigutilde^{\pm}_u(n+p,p,\vecJ),\leq_{\theta^u})
\simeq
 \bigl(
 \gbigutilde^{\pm}_z(n,p,\nu^{\pm}_m(\vecJ)),
 \leq_{\kappa^{\pm}_{m,\vecJ}(\theta^u)}\bigr).
\]
We also obtain the following commutative diagram:
\[
  \begin{CD}
  \gbigu^{\pm}_u(n+p,p,\vecJ)
  @>{\iota_u}>>   
  \gbigutilde^{\pm}_u(n+p,p,\vecJ)
  @>{q_{u,p+n,n}}>>
  \gbigu^{\pm}_u(n+p,p,\vecJ)
  \\
   @V{\nu^{\pm}_{m,\vecJ}}V{\simeq}V
   @V{\nutilde^{\pm}_{m,\vecJ}}V{\simeq}V
   @V{\nu^{\pm}_{m,\vecJ}}V{\simeq}V \\
 \gbigu^{\pm}_z(n,p,\nu^{\pm}_m(\vecJ))
@>{\iota_z}>>
 \gbigutilde^{\pm}_z(n,p,\nu^{\pm}_m(\vecJ))
@>{q_{z,p,n}}>>
 \gbigu^{\pm}_z(n,p,\nu^{\pm}_m(\vecJ)).
  \end{CD}
\]
 Here, $\iota_u$  and $\iota_z$
 denote the natural inclusions.
\end{prop}
\pf
It follows from Proposition \ref{prop;20.10.4.10}.
\hfill\qed

\vspace{.1in}

We obtain the following lemma
from Lemma \ref{lem;20.10.3.1}.
\begin{lem}
We obtain
$\nu_m^+(\vecJ+(1+\omega^{-1})\pi)
 =\nu_m^-(\vecJ)+\omega^{-1}\pi$
 and
$\nu_m^-(\vecJ+(1+\omega^{-1})\pi)
=\nu_{m+1}^+(\vecJ)+\omega^{-1}\pi$.
Moreover, we obtain
$\nutilde^{+}_{m,\vecJ+(1+\omega^{-1})\pi}(\gminibtilde)
=\nutilde^{-}_{m,\vecJ}(\gminibtilde)$
for
 $\gminibtilde\in
 \gbigutilde_u^{+}\bigl(n+p,p,\vecJ+(1+\omega^{-1})\pi\bigr)
 =\gbigutilde_u^{-}(n+p,p,\vecJ)$,
and 
$\nutilde^{-}_{m,\vecJ+(1+\omega^{-1})\pi}(\gminibtilde)
=\nutilde^{+}_{m+1,\vecJ}(\gminibtilde)$
for
 $\gminibtilde\in
 \gbigutilde_u^{-}\bigl(n+p,p,\vecJ+(1+\omega^{-1})\pi\bigr)
 =\gbigutilde_u^{+}(n+p,p,\vecJ)$.
\hfill\qed
 \end{lem}

For $\gminib\in\gbigu^{\pm}_u(n+p,n,\vecJ)$,
we can describe
$\nu^{\pm}_{m,\vecJ}(\gminib)$ explicitly.
Indeed, for any 
\[
 \gminib(u_{n+p})=\mp a\exp\Bigl(
 \sqrt{-1}\Bigl(
\frac{\omega}{1+\omega}\vartheta^u_0
 \Bigr)
 \Bigr)\cdot u_{n+p}^{-n}
\in  \gbigu^{\pm}_{u}(n+p,n,\vecJ)
\quad (a>0),
\]
we obtain
\[
 \nu_{m,\vecJ}^-(\gminib)(z_{p})=
 \bigl(
 \langle\omega\rangle^{-1}a
 \bigr)^{1+\omega}
  \exp\Bigl(
 \sqrt{-1}
 \omega\bigl(\vartheta^u_0+2m\pi\bigr)
 \Bigr) z_p^{-n}
\in \gbigu_z^{-}(p,n,\nu^{-}_m(\vecJ)),
\]
\[
 \nu_{m,\vecJ}^+(\gminib)(z_{p})=
-\bigl(
 \langle\omega\rangle^{-1}a
 \bigr)^{1+\omega}
  \exp\Bigl(
 \sqrt{-1}
 \omega\bigl(\vartheta^u_0+(2m-1)\pi\bigr)
 \Bigr) z_p^{-n}
\in \gbigu_z^{+}(p,n,\nu^{+}_m(\vecJ)).
\]

\subsection{Transformation of index sets
induced by the local Fourier transform}
\label{subsection;18.5.5.1}

As explained in \S\ref{subsection;24.2.23.40}--\ref{subsection;24.2.23.50},
the local Fourier transform induces a transformation
of any $\Gal(p)$-invariant subset
$\nbigitilde\subset z_p^{-1}\cnum[z_p^{-1}]$
to
a $\Gal(n+p)$-invariant subset
$\gbigf^{(0,\infty)}_+(\nbigitilde)
\subset u_{n+p}^{-1}\cnum[u_{n+p}^{-1}]$.
By the construction,
we have
$\gbigf_+^{(0,\infty)}(q_{z,p,n}(\nbigitilde))
=q_{u,n+p,n}(\gbigf_+^{(0,\infty)}(\nbigitilde))$.

Let $\nbigitilde$ be a $\Gal(p)$-invariant subset of 
$\gbigutilde_z(n,p)$.
Set $\nbigitilde^{\circ}:=\gbigf_+^{(0,\infty)}(\nbigitilde)$.
We also set
$\nbigi:=q_{z,p,n}(\nbigitilde)$
and $\nbigi^{\circ}:=q_{u,n+p,n}(\nbigitilde^{\circ})$.
For $J\in T(\nbigi)$,
we set
$\nbigitilde_{J,>0}:=q_{z,p,n}^{-1}(\nbigi_{J,>0})$
and 
$\nbigitilde_{J,<0}:=q_{z,p,n}^{-1}(\nbigi_{J,<0})$.
\index{sets $\nbigitilde_{J,<0}$, $\nbigitilde_{J,>0}$}
Similarly,
for $\vecJ\in T(\nbigi^{\circ})$,
we set
$\nbigitilde^{\circ}_{\vecJ,>0}:=
 q_{u,n+p,n}^{-1}(\nbigi^{\circ}_{\vecJ,>0})$
and 
$\nbigitilde^{\circ}_{\vecJ,<0}:=
 q_{u,n+p,n}^{-1}(\nbigi^{\circ}_{\vecJ,<0})$.

\begin{prop}
\label{prop;18.5.5.40}
For any $\vecJ\in T(\nbigi^{\circ})$,
the maps $\nu^{\pm}_{m,\vecJ}$ induce
isomorphisms of
the following partially ordered sets
for any $\theta^u\in\real$:
\[
\nu^{-}_{m,\vecJ}:
(\nbigi^{\circ}_{\vecJ,<0},\leq_{\theta^u})
\simeq
 \bigl(
 \nbigi_{\nu^{-}_m(\vecJ),<0},
 \leq_{\kappa^{-}_{m,\vecJ}(\theta^u)}\bigr),
\]
\[
 \nu^{+}_{m,\vecJ}:
(\nbigi^{\circ}_{\vecJ,>0},\leq_{\theta^u})
\simeq
 \bigl(
 \nbigi_{\nu^{+}_m(\vecJ),>0},
 \leq_{\kappa^{+}_{m,\vecJ}(\theta^u)}\bigr).
\]
The maps $\nutilde^{\pm}_{m,\vecJ}$ induce
isomorphisms of the following partially ordered sets
for any $\theta^u\in\real$:
\[
 \nutilde_{m,\vecJ}^{-}:
 \bigl(
 \nbigitilde^{\circ}_{\vecJ,<0},
 \leq_{\theta^u}
 \bigr)
\lrarr
 \bigl(
 \nbigitilde_{\nu_m^-(\vecJ),<0},
 \leq_{\kappa^-_{m,\vecJ}(\theta^u)}
 \bigr),
\]
\[
\nutilde_{m,\vecJ}^{+}:
 \bigl(
 \nbigitilde^{\circ}_{\vecJ,>0},\leq_{\theta^u}
\bigr)
\lrarr
 \bigl(
 \nbigitilde_{\nu_m^+(\vecJ),>0},
 \leq_{\kappa^+_{m,\vecJ}(\theta^u)}
\bigr).
\]
We also have the commutativity
$q_{u,p+n,n}\circ\nutilde^{\pm}_{m,\vecJ}
=\nu^{\pm}_{m,\vecJ}\circ q_{z,p,n}$.
\end{prop}
\pf
It follows from Proposition \ref{prop;20.10.4.11}.
\hfill\qed

\section{From $\infty$ to $\infty$}
\label{subsection;20.10.9.30}

We shall refine the construction in
\S\ref{subsection;24.2.23.43}--\S\ref{subsection;24.2.23.51}.

\subsection{Preliminary computations (1)}
\label{subsection;18.4.7.3}

Let $n>p$ be two positive integers.
We set $\omega=n/p$
and 
$\langle\omega\rangle':=
\omega^{\frac{-1}{\omega-1}}
-\omega^{\frac{-\omega}{\omega-1}}>0$.
\index{number $\langle\omega\rangle'$}
Let $\alpha$ be a non-zero complex number.
We set $\gminia:=\alpha\zeta^{-n}$,
and consider
\[
 G_{\gminia,\eta}(\zeta):=
 \gminia(\zeta)+\eta^{-n+p}\zeta^{-p}
=\alpha\zeta^{-n}+\eta^{-n+p}\zeta^{-p}.
\]

We obtain
$\del_{\zeta}G_{\gminia,\eta}(\zeta)
=-n\alpha\zeta^{-n-1}-p\eta^{-n+p}\zeta^{-p-1}$.
By setting
$h_{\omega}(\eta)=\omega^{1/(n-p)}\eta$,
we obtain
\[
 \del_{\zeta}G^+_{\gminia,\eta}(\zeta_0)=0
\Longleftrightarrow
 \zeta_0\in
 \bigl\{
 h_{\omega}(\beta\eta)\,\big|\,
 \beta\in\cnum,\,\,
 \beta^{n-p}=-\alpha
 \bigr\}.
\]

Take an interval
$J=I(\vartheta_0,\omega^{-1}\pi/2)$.
We set \index{interval $J^u(m,\pm)$}
\[
\left\{
\begin{array}{l}
J^u(m,-):=
I\bigl(-\vartheta_0-(2m-1)\pi,(1-\omega^{-1})\pi/2\bigr), \\
J^u(m,+):=
 I\bigl(
 -\vartheta_0-2m\pi,
 (1-\omega^{-1})\pi/2
 \bigr).
\end{array}
\right.
\]

\subsubsection{Case 1}
If $\Re\gminia(e^{\sqrt{-1}\theta/p})>0$ for any $\theta\in J$,
we obtain
\[
 -\alpha=|\alpha|
\exp\bigl(
\sqrt{-1}\bigl(
\omega\vartheta_0
-\pi
\bigr)
\bigr).
\]
For any integer $m$,
we set
\begin{equation}
 \beta_{J,m,-}:=
   |\alpha|^{1/(n-p)}
 \exp\Bigl(
 \frac{\sqrt{-1}}{n-p}
 \bigl(
 \omega\vartheta_0+(2m-1)\pi
 \bigr)
 \Bigr)
 \in\bigl\{
 \beta\,\big|\,
 \beta^{n-p}=-\alpha
\bigr\}.
\end{equation}
The set of roots of $\del_{\zeta}G_{\gminia,\eta}$
is
$\{h_{\omega}(\beta_{J,m,-}\eta)\,|\,m\in\seisuu\}$.
For $m\in\seisuu$,
we set \index{map $\gbigf^{(\infty,\infty)}_{(J,m,-)}$}
\begin{multline}
\label{eq;18.5.5.12}
 \gbigf^{(\infty,\infty)}_{(J,m,-)}(\gminia)(\eta):=
 G_{\gminia,\eta}(h_{\omega}(\beta_{J,m,-}\eta))
\\
 =
 \langle\omega\rangle'
 |\alpha|^{\frac{-1}{\omega-1}}
 \exp\Bigl(
 \frac{-\sqrt{-1}}{\omega-1}
 \bigl(\omega\vartheta_0+(2m-1)\pi\bigr)
 \Bigr)\cdot \eta^{-n} \\
=-\langle\omega\rangle'
 |\alpha|^{\frac{-1}{\omega-1}}
 \exp\Bigl(
 \frac{-\sqrt{-1}\omega}{\omega-1}
 \bigl(\vartheta_0+(2m-1)\pi\bigr)
 \Bigr)\cdot \eta^{-n}.
\end{multline}

\begin{lem}
\label{lem;18.5.5.10}
We set
\[
 \arg
 \bigl(
 h_{\omega}(\beta_{J,m,-}e^{\sqrt{-1}\theta^u/(n-p)})
 \bigr)
 =\frac{1}{p(\omega-1)}
 \bigl(
 \theta^u+\omega\vartheta_0+(2m-1)\pi
 \bigr).
\]
Then,
$p\cdot\arg
 \bigl(
 h_{\omega}(\beta_{J,m,-}e^{\sqrt{-1}\theta^u/(n-p)})
 \bigr)
\in J$
if and only if
$\theta^u\in J^u(m,-)$.
Moreover,
we obtain
$\Re\gbigf^{(\infty,\infty)}_{(J,m,-)}(\gminia)(\eta)<0$
for $\eta=|\eta|\exp(\sqrt{-1}\theta^u/(n-p))$
with $\theta^u\in J^u(m,-)$.
\end{lem}
\pf
We obtain
\[
  p\arg
 \bigl(
 h_{\omega}(\beta_{m,-}e^{\sqrt{-1}\theta^u/(n-p)})
 \bigr)
-\vartheta_0
=
\frac{1}{\omega-1}
\bigl(
\theta^u+
\vartheta_0+(2m-1)\pi
\bigr).
\]
Then, the first claim is clear.
The second claim is clear
by the formula (\ref{eq;18.5.5.12}).
\hfill\qed

\vspace{.1in}
We remark the following obvious lemma.
\begin{lem}
\label{lem;18.5.24.20}
For any integer $\ell$,
we obtain
$(-1)^{\ell}\Re\gminia(e^{\sqrt{-1}\theta/p})>0$ 
on $J+\ell \omega^{-1}\pi$,
and
$(-1)^{\ell}\Re\gbigf^{(\infty,\infty)}_{(J,m,-)}(\gminia)
 (e^{\sqrt{-1}\theta^u/(n-p)})<0$
on $J^u(m,-)+\ell (1-\omega^{-1})\pi$.
\hfill\qed
\end{lem}

\subsubsection{Case 2}

If $\Re\gminia(e^{\sqrt{-1}\theta/p})<0$ for any $\theta\in J$,
we obtain
\[
 -\alpha=|\alpha|
 \exp\bigl(
 \sqrt{-1}
 \omega\vartheta_0
 \bigr).
\]
For any integer $m$, we set
\[
 \beta_{J,m,+}:=
 |\alpha|^{1/(n-p)}
 \exp\Bigl(
 \frac{\sqrt{-1}}{n-p}
 \bigl(
 \omega\vartheta_0+2m\pi
 \bigr)
 \Bigr).
\]
The set of roots of $\del_{\zeta}G_{\gminia,\eta}$ is
$\{h_{\omega}(\beta_{J,m,+}\eta)\,|\,m\in\seisuu\}$.
We set \index{map $\gbigf^{(\infty,\infty)}_{(J,m,+)}$}
\begin{multline}
\label{eq;18.5.5.13}
 \gbigf^{(\infty,\infty)}_{(J,m,+)}(\gminia)(\eta):=
 G_{\gminia,\eta}\bigl(h_{\omega}(\beta_{J,m,+}\eta)\bigr)=
 \\
 \langle\omega\rangle'
 |\alpha|^{-p/(n-p)}
 \exp\Bigl(
 \frac{-\sqrt{-1}}{\omega-1}
 \bigl(\omega\vartheta_0+2m\pi\bigr)
 \Bigr) \eta^{-n}=
 \\
 \langle\omega\rangle'
 |\alpha|^{-p/(n-p)}
 \exp\Bigl(
 \frac{-\sqrt{-1}\omega}{\omega-1}
 \bigl(\vartheta_0+2m\pi\bigr)
 \Bigr)\eta^{-n}.
\end{multline}
The following lemma is similar to Lemma \ref{lem;18.5.5.10}.
\begin{lem}
\label{lem;18.5.5.11}
We set
\[
 \arg\bigl(
  h_{\omega}(\beta_{m,+}e^{\sqrt{-1}\theta^u/(n-p)})
 \bigr)
 =\frac{1}{p(\omega-1)}
 \bigl(
 \theta^u+\omega\vartheta_0+2m\pi
 \bigr).
\]
Then,
 $ p\arg\bigl(
  h_{\omega}(\beta_{m,+}e^{\sqrt{-1}\theta^u/(n-p)})
 \bigr)
\in J$ 
if and only if
$\theta^u\in J^u(m,+)$.
Moreover,
we obtain
$\Re\gbigf^{(\infty,\infty)}_{(m,J,+)}(\gminia)(\eta)>0$
for $\eta=|\eta|\exp(\sqrt{-1}\theta^u/(n-p))$
with $\theta^u\in J^u(m,+)$.
\hfill\qed
\end{lem}

We note the following obvious lemma.
\begin{lem}
\label{lem;18.5.24.21}
For any integer $\ell$,
we obtain
$(-1)^{\ell}\Re\gminia(e^{\sqrt{-1}\theta/p})<0$ 
on  $J+\ell \omega^{-1}\pi$,
and
$(-1)^{\ell}\Re\gbigf^{(\infty,\infty)}_{(J,m,+)}
 (\gminia)(e^{\sqrt{-1}\theta^u/(n-p)})>0$
on $J^u(m,+)+\ell (1-\omega^{-1})\pi$.
\hfill\qed
\end{lem}

\subsection{Preliminary computations (2)}

For 
$\gminiatilde=
 \alpha \zeta^{-n}+\sum_{j=1}^{n-1}\gminiatilde_j\zeta^{-j}$,
we set
\[
 G_{\gminiatilde,\eta}(\zeta):=
 \gminiatilde(\zeta)+\eta^{-n+p}\zeta^{-p}.
\]
We obtain
$\del_{\zeta}G_{\gminiatilde,\eta}
=-n\alpha\zeta^{-n-1}-p\zeta^{-p-1}\eta^{-n+p}
-\sum j\gminiatilde_j\zeta^{-j-1}$.
The following lemma is similar to
Lemma \ref{lem;18.4.13.10}.
\begin{lem}
\label{lem;18.5.5.22}
There exists a unique convergent power series
$a_{\gminiatilde}(\eta)=
 1+\sum_{j=1}^{\infty} a_{\gminiatilde,j}\eta^j$
such that 
the following holds for any $\beta\in\cnum$
with $\beta^{n-p}=-\alpha$:
\begin{equation}
\label{eq;18.4.13.3}
 \del_{\zeta}G_{\gminiatilde,\eta}
 \Bigl(
 h_{\omega}(\beta\eta)a_{\gminiatilde}(\beta\eta)
 \Bigr)=0.
\end{equation}
Conversely,
any root of $\del_{\zeta}G_{\gminiatilde,\eta}$
is $h_{\omega}(\beta\eta)a_{\gminiatilde}(\beta\eta)$
for some $\beta$ with $\beta^{n-p}=-\alpha$.
\hfill\qed
\end{lem}

There exists the convergent formal power series
$1+\sum_{j=1}^{\infty} b_{\gminiatilde}\eta^j$
such that 
\begin{multline}
\label{eq;20.10.5.20}
 G_{\gminiatilde,\eta}\bigl(
 h_{\omega}(\beta\eta)\cdot a_{\gminiatilde}(\beta\eta)
 \bigr)
=G_{\gminia,\eta}\bigl(
 h_{\omega}(\beta\eta)
 \bigr)\cdot
 \Bigl(
 1+\sum_{j=1}^{\infty}b_{\gminiatilde,j}(\beta\eta)^j
 \Bigr)
  \\
=-\alpha\langle\omega\rangle'
 (\beta\eta)^{-n}\cdot
  \Bigl(
 1+\sum_{j=1}^{\infty}b_{\gminiatilde,j}(\beta\eta)^j
 \Bigr).
\end{multline}
We obtain the following lemma
as in the case of Lemma \ref{lem;18.4.13.1}.
\begin{lem}
\label{lem;18.4.13.31}
If
$\gminiatilde_i=
 \alpha\zeta^{-n}
 +\sum_{j=1}^{n-1}\gminiatilde_{i,j}\zeta^{-j}$ $(i=1,2)$
satisfies
 $\gminiatilde_{1,j}=\gminiatilde_{2,j}$ $(j=k+1,\ldots,n-1)$
for some $1\leq k\leq n-1$,
then we obtain
$b_{\gminiatilde_1,j}=b_{\gminiatilde_2,j}$
$(j=1,\ldots,n-k-1)$.
Moreover, we obtain
\[
 G_{\gminiatilde_1,\eta}\bigl(
 h_{\omega}(\beta\eta)\cdot a_{\gminiatilde_1}(\beta\eta)
 \bigr)
-G_{\gminiatilde_2,\eta}\bigl(
 h_{\omega}(\beta\eta)\cdot a_{\gminiatilde_2}(\beta\eta)
 \bigr)
\equiv
 \bigl(
 \gminiatilde_{1,k}
-\gminiatilde_{2,k}
 \bigr)\cdot h_{\omega}(\beta\eta)^{-k}
\]
modulo $\eta^{-k+1}\cnum[\![\eta]\!]$.
\hfill\qed
\end{lem}

\subsection{Direct consequences of preliminary computations}

We take a $p$-th root $x_p$ of $x=z^{-1}$
and $(n-p)$-th root $u_{n-p}$ of $u=w^{-1}$.
\index{variable $x_p$}
Let $J=I(\vartheta_0,\omega^{-1}\pi/2)$.
We define the maps
$\gbigf^{(\infty,\infty)}_{(J,m,\pm)}:
 \gbigu^{\pm}_x(p,n,J)
\lrarr
 \gbigu_u(n-p,n)$
by the formulas (\ref{eq;18.5.5.12}) and 
(\ref{eq;18.5.5.13})
with $\zeta=z_p$ and $\eta=u_{n-p}$.
\begin{lem}
\label{lem;20.10.5.22}
They induce the following isomorphisms of
the partially ordered sets:
\[
 \gbigf^{(\infty,\infty)}_{(J,m,\pm)}:
 (\gbigu^{\pm}_{x}(p,n,J),\leq_J)
\simeq
 (\gbigu^{\mp}_{u}(n-p,n,J^u(m,\pm)),
 \leq_{J^u(m,\pm)}).
\]
\end{lem}
\pf
It follows from Lemma \ref{lem;18.5.5.10}
and Lemma \ref{lem;18.5.5.11}.
\hfill\qed

\vspace{.1in}

For
$\gminiatilde(x_p)
=\sum_{j=1}^{n}\gminiatilde_jx_p^{-j}
\in\gbigutilde_x^-(p,n,J)$,
we set
\[
 \beta_{\gminiatilde,J,m,-}:=
|\gminiatilde_n|^{1/(n-p)}
\exp\Bigl(
\frac{\sqrt{-1}}{n-p}
(\omega\vartheta_0+(2m-1)\pi)
\Bigr).
\]
Similarly,
for 
$\gminiatilde(x_p)
\in\gbigutilde_x^+(p,n,J)$,
we set
\[
\beta_{\gminiatilde,J,m,+}:=
|\gminiatilde_n|^{1/(n-p)}
\exp\Bigl(
\frac{\sqrt{-1}}{n-p}
(\omega\vartheta_0+2m\pi)
\Bigr).
\]
For $\gminiatilde\in\gbigutilde^{\pm}_x(n,p,J)$,
we define \index{maps $ \gbigf^{(\infty,\infty)}_{(J,m,\pm)}$}
\[
 \gbigf^{(\infty,\infty)}_{(J,m,\pm)}(\gminiatilde)
=G_{\gminiatilde,u_{n-p}}\bigl(
 h_{\omega}(\beta_{\gminiatilde,J,m,\pm}\eta)\cdot
 a_{\gminiatilde}(\beta_{\gminiatilde,J,m,\pm}\eta)
 \bigr)
\,\,\,\modulo\,\,\,
\cnum[\![u_{n-p}]\!],
\]
where
we set
$h_{\omega}(\beta_{\gminiatilde,J,m,\pm}\eta)
=\omega^{1/(n-p)}\beta_{\gminiatilde,J,m,\pm}\eta$,
and $a_{\gminiatilde}$ is
the convergent power series
as in Lemma \ref{lem;18.5.5.22}.
Thus, we obtain the map
$\gbigf^{(\infty,\infty)}_{(J,m,\pm)}:
 \gbigutilde^{\pm}_x(p,n,J)\lrarr
 \gbigutilde_u(n-p,n)$.
By (\ref{eq;20.10.5.20}),
the following holds:
\begin{equation}
 \label{eq;20.10.5.21}
 \gbigf^{(\infty,\infty)}_{(J,m,\pm)}
 \bigl(q_{x,p,n}(\gminiatilde)\bigr)
 =q_{u,n-p,n}(\gbigf^{(\infty,\infty)}_{(J,m,\pm)}
 \bigl(\gminiatilde)\bigr).
\end{equation}

We obtain the following lemma
from (\ref{eq;20.10.5.21}),
Lemma \ref{lem;20.10.5.22}
and Lemma \ref{lem;18.4.13.31}
as in the case of Lemma \ref{lem;20.10.5.23}.
\begin{lem}
The maps $\gbigf_{(J,m,\pm)}^{(\infty,\infty)}$ induce
bijections
\[
 \gbigf^{(\infty,\infty)}_{(J,m,\pm)}:
  \gbigutilde^{\pm}_x(p,n,J)
\simeq
 \gbigutilde^{\mp}_u\bigl(n-p,n,J^u(m,\pm)\bigr).
\]
\hfill\qed
\end{lem}

For $\theta\in\real$,
we define
$\theta^u(m,\pm)$
by the following formulas:
\index{points $\theta^u(m,\pm)$}
\[
\theta^u(m,-)
=(\omega-1)\theta+\omega\vartheta_0-(2m-1)\pi,
\quad\quad
\theta^u(m,+)
=(\omega-1)\theta+\omega\vartheta_0-2m\pi.
\]
If $\theta\in J+\ell\omega^{-1}\pi$,
then 
$\theta^u(m,\pm)
\in J^u(m,\pm)+\ell(1-\omega^{-1})\pi$.
The following proposition is similar to
Proposition \ref{prop;20.10.4.10}.
\begin{prop}
\label{prop;20.10.5.30}
There exist the following isomorphisms of the partially
ordered sets
for any $\theta\in \real$:
\[
 \gbigf^{(\infty,\infty)}_{(J,m,\pm)}:
 \bigl(
 \gbigutilde^{\pm}_{x}(p,n,J),\leq_{\theta}
 \bigr)
\simeq
 \bigl(
 \gbigutilde^{\mp}_{u}(n-p,n,J^u(m,\pm)),\leq_{\theta^u(m,\pm)}
 \bigr).
\]
We also obtain the following commutative diagram:
{\small
\begin{equation}
\label{eq;20.10.5.24}
 \begin{CD}
  \gbigu^{\pm}_{x}(p,n,J)
  @>{\iota_x}>>
  \gbigutilde^{\pm}_{x}(p,n,J)
  @>{q_{x,p,n}}>>
  \gbigu^{\pm}_{x}(p,n,J)\\
  @V{\gbigf^{(\infty,\infty)}_{(J,m,\pm)}}V{\simeq}V
  @V{\gbigf^{(\infty,\infty)}_{(J,m,\pm)}}V{\simeq}V
  @V{\gbigf^{(\infty,\infty)}_{(J,m,\pm)}}V{\simeq}V \\
  \gbigu^{\mp}_{u}(n-p,n,J^u(m,\pm))
  @>{\iota_u}>>
  \gbigutilde^{\mp}_{u}(n-p,n,J^u(m,\pm))
  @>{q_{u,n-p,n}}>>
  \gbigu^{\mp}_{u}(n-p,n,J^u(m,\pm)).
 \end{CD}
 \end{equation}
}
Here, $\iota_x$ and $\iota_u$
denote the natural inclusions.
\hfill\qed
\end{prop}
 
Note that
$(J+\omega^{-1}\pi)^u(m,-)
=J^u(m,+)+(1-\omega^{-1})\pi$
and 
$(J+\omega^{-1}\pi)^u(m-1,+)
=J^u(m,-)+(1-\omega^{-1})\pi$.
We have
\[
 \gbigutilde^{-}_u\bigl(
 n-p,n,(J+\omega^{-1}\pi)^u(m,-)
 \bigr)
 =\gbigutilde^+_u\bigl(
 n-p,n,J^u(m,+)
 \bigr),
\]
\[
 \gbigutilde^{+}_u\bigl(
 n-p,n,(J+\omega^{-1}\pi)^u(m-1,+)
 \bigr)
 =\gbigutilde^-_u\bigl(
 n-p,n,J^u(m,-)
 \bigr).
\]
The following lemma can be checked by computation.
\begin{lem}
\label{lem;20.10.4.1}
For $\gminiatilde\in \gbigutilde^{+}_x(p,n,J)$,
we obtain
$\gbigf^{(\infty,\infty)}_{(J,m,+)}(\gminiatilde)
 =\gbigf^{(\infty,\infty)}_{(J+\omega^{-1}\pi,m,-)}(\gminiatilde)$.
For $\gminiatilde\in \gbigutilde^{-}_x(p,n,J)$,
we obtain
$\gbigf^{(\infty,\infty)}_{(J,m,-)}(\gminiatilde)
 =\gbigf^{(\infty,\infty)}_{(J+\omega^{-1}\pi,m-1,+)}(\gminiatilde)$.
\hfill\qed
\end{lem}

\subsection{Reformulation}
\label{subsection;25.2.6.1}

Let $\vecJ=I\bigl(\vartheta^u_0,(1-\omega^{-1})\pi/2\bigr)$.
For any integer $m$,
we set
\index{intervals $\nu^{\pm}_m(\vecJ)$}
\[
 \nu^-_m(\vecJ)
 =I\bigl(
 -\vartheta^u_0-(2m-1)\pi,
 \omega^{-1}\pi/2
 \bigr),
 \quad
  \nu^+_m(\vecJ)
 =I\bigl(
 -\vartheta^u_0-2m\pi,
 \omega^{-1}\pi/2
 \bigr).
\]
We define the isomorphisms of the partially ordered sets
\index{isomorphisms $\nu^{\pm}_{m,\vecJ}$}
\[
 \nu^{\pm}_{m,\vecJ}:
 \bigl(
 \gbigu^{\mp}_{u}(n-p,n,\vecJ),\leq_{\vecJ}
 \bigr)
\simeq
 \bigl(
 \gbigu^{\pm}_{x}(p,n,\nu^{\pm}_m(\vecJ)),\leq_{\nu^{\pm}_m(\vecJ)}
\bigr)
\]
as the inverse of
$\gbigf^{(\infty,\infty)}_{(\nu^{\pm}_m(\vecJ),m,\pm)}$.
We define the bijection
\index{isomorphisms $\nutilde^{\pm}_{m,\vecJ}$}
\[
 \nutilde^{\pm}_{m,\vecJ}:
 \gbigutilde^{\mp}_{u}(n-p,n,\vecJ)
\lrarr
 \gbigutilde^{\pm}_{x}(p,n,\nu_m^{\pm}(\vecJ))
\]
as the inverse of
$\gbigf^{(\infty,\infty)}_{(\nu^{\pm}_m(\vecJ),m,\pm)}$.

We define the maps
$\kappa_{m,\vecJ}^{\pm}:
\real\lrarr\real$
by the following formulas:
\index{maps $\kappa_{m,\vecJ}^{\pm}$}
\[
 \kappa^-_{m,\vecJ}(\theta^u)=
 \frac{1}{\omega-1}(\theta^u-\omega\vartheta^u_0)
 -(2m-1)\pi,
 \quad
 \kappa^+_{m,\vecJ}(\theta^u)=
 \frac{1}{\omega-1}(\theta^u-\omega\vartheta^u_0)
 -2m\pi.
\]
Note that
$\kappa^{\pm}_{m,\vecJ}$
induces bijections
$\vecJ+\ell(1-\omega^{-1})\pi
\simeq
\nu_{m}^{\pm}(\vecJ)
+\ell \omega^{-1}\pi$
for any integer $\ell$.

\begin{prop}
\label{prop;20.10.5.31}
The maps
$\nutilde^{\pm}_{m,\vecJ}$ induce isomorphisms of
the following partially ordered sets
for any $\theta^u\in\real$:
\[
\nutilde^{\pm}_{m,\vecJ}:
(\gbigutilde^{\mp}_u(n-p,p,\vecJ),\leq_{\theta^u})
\simeq
 \bigl(
 \gbigutilde^{\pm}_x(n,p,\nu^{\pm}_m(\vecJ)),
 \leq_{\kappa^{\pm}_{m,\vecJ}(\theta^u)}\bigr).
\]
We also obtain the following commutative diagram:
\[
  \begin{CD}
  \gbigu^{\mp}_u(n-p,p,\vecJ)
  @>{\iota_u}>>   
  \gbigutilde^{\mp}_u(n-p,p,\vecJ)
  @>{q_{u,n-p,n}}>>
  \gbigu^{\mp}_u(n-p,p,\vecJ)
  \\
   @V{\nu^{\pm}_{m,\vecJ}}V{\simeq}V
   @V{\nutilde^{\pm}_{m,\vecJ}}V{\simeq}V
   @V{\nu^{\pm}_{m,\vecJ}}V{\simeq}V \\
 \gbigu^{\pm}_x(n,p,\nu^{\pm}_m(\vecJ))
@>{\iota_z}>>
 \gbigutilde^{\pm}_x(n,p,\nu^{\pm}_m(\vecJ))
@>{q_{z,p,n}}>>
 \gbigu^{\pm}_x(n,p,\nu^{\pm}_m(\vecJ)).
  \end{CD}
\]
 Here, $\iota_u$  and $\iota_x$
 denote the natural inclusions.
\end{prop}
\pf
It follows from Proposition \ref{prop;20.10.5.30}.
\hfill\qed

\vspace{.1in}

We obtain the following lemma
from Lemma \ref{lem;20.10.4.1}.
\begin{lem}
We obtain
$\nu_m^+(\vecJ+(1-\omega^{-1})\pi)
 =\nu_{m+1}^-(\vecJ)+\omega^{-1}\pi$
and
$\nu_m^-(\vecJ+(1-\omega^{-1})\pi)
=\nu_{m}^+(\vecJ)+\omega^{-1}\pi$.
Moreover,
we obtain
$\nutilde^{+}_{m,\vecJ+(1-\omega^{-1})\pi}(\gminibtilde)
=\nutilde^{-}_{m+1,\vecJ}(\gminibtilde)$
for
 $\gminibtilde\in
 \gbigutilde_u^{-}\bigl(n-p,p,\vecJ+(1-\omega^{-1})\pi\bigr)
 =\gbigutilde_u^{+}(n-p,p,\vecJ)$,
and 
$\nutilde^{-}_{m,\vecJ+(1-\omega^{-1})\pi}(\gminibtilde)
=\nutilde^{+}_{m,\vecJ}(\gminibtilde)$
for
 $\gminibtilde\in
 \gbigutilde_u^{+}\bigl(n-p,p,\vecJ+(1-\omega^{-1})\pi\bigr)
 =\gbigutilde_u^{-}(n-p,p,\vecJ)$.
\hfill\qed
 \end{lem}

For $\gminib\in\gbigutilde^{\mp}_u(n-p,n,\vecJ)$,
we can describe
$\nu^{\pm}_{m,\vecJ}(\gminib)$
explicitly.
Indeed, for 
\[
 \gminib=
 \pm a\exp\Bigl(
 \frac{\sqrt{-1}\omega}{\omega-1}
 \vartheta_0
 \Bigr)\cdot u_{n-p}^{-n}
 \in\gbigu^{\mp}_u(\vecJ)
\quad (a>0),
\]
we obtain
\[
  \nu_m^{-}(\gminib)=
-\Bigl(
 \frac{a}{\langle\omega\rangle'}
 \Bigr)^{-(\omega-1)}
 \exp\Bigl(
 \sqrt{-1}
 \omega\bigl(-\vartheta_0-(2m-1)\pi\bigr)
 \Bigr)
\cdot x_p^{-n}
\in \gbigu^{+}_x(p,n,\nu^{-}_{m}(\vecJ)),
\]
\[
  \nu_m^{+}(\gminib)=
 \Bigl(
 \frac{a}{\langle\omega\rangle'}
 \Bigr)^{-(\omega-1)}
 \exp\Bigl(
 \sqrt{-1}
 \omega\bigl(-\vartheta_0-2m\pi\bigr)
 \Bigr)
\cdot x_p^{-n}
\in \gbigu^{-}_x(p,n,\nu^{+}_{m}(\vecJ)).
\]

\subsection{Transformation of index sets induced by
the local Fourier transform}

As explained in
\S\ref{subsection;24.2.23.43}--\S\ref{subsection;24.2.23.51},
the local Fourier transform induces
a transformation of any
$\Gal(p)$-invariant subset
$\nbigitilde\subset x_p^{-1}\cnum[x_p^{-1}]$
to
a $\Gal(n+p)$-invariant subset
$\gbigf^{(\infty,\infty)}_+(\nbigitilde)
\subset
 u_{n-p}^{-1}\cnum[u_{n-p}^{-1}]$.
By the construction,
we have
$\gbigf_+^{(\infty,\infty)}(q_{x,p,n}(\nbigitilde))
=q_{u,n-p,n}(\gbigf_+^{(\infty,\infty)}(\nbigitilde))$.

Let $\nbigitilde$ be a $\Gal(p)$-invariant subset of
$\gbigutilde_x(n,p)$.
Set $\nbigitilde^{\circ}:=
\gbigf_+^{(\infty,\infty)}(\nbigitilde)$.
We also set
$\nbigi:=q_{x,p,n}(\nbigitilde)$
and $\nbigi^{\circ}:=q_{u,n-p,n}(\nbigitilde^{\circ})$.
For $J\in T(\nbigi)$,
we set
$\nbigitilde_{J,>0}:=q_{x,p,n}^{-1}(\nbigi_{J,>0})$
and 
$\nbigitilde_{J,<0}:=q_{x,p,n}^{-1}(\nbigi_{J,<0})$.
Similarly,
for $\vecJ\in T(\nbigi^{\circ})$,
we set
$\nbigitilde^{\circ}_{\vecJ,>0}:=
 q_{u,n-p,n}^{-1}(\nbigi^{\circ}_{\vecJ,>0})$
and 
$\nbigitilde^{\circ}_{\vecJ,<0}:=
 q_{u,n-p,n}^{-1}(\nbigi^{\circ}_{\vecJ,<0})$.

\begin{prop}
\label{prop;18.5.6.1}
For any $\vecJ\in T(\nbigi^{\circ})$,
the maps $\nu^{\pm}_{m,\vecJ}$ induce
the following isomorphisms
of the partially ordered sets
for any $\theta^u\in\real$:
\[
\nu^{-}_{m,\vecJ}:
 \bigl(
 \nbigi^{\circ}_{\vecJ,>0},\leq_{\theta^u}
 \bigr)
\simeq
 \bigl(
 \nbigi_{\nu^-_m(\vecJ),<0},\leq_{\kappa^-_{m,\vecJ}(\theta^u)}
 \bigr),
\]
\[
 \nu^{+}_{m,\vecJ}:
 \bigl(
 \nbigi^{\circ}_{\vecJ,<0},\leq_{\theta^u}
 \bigr)
\simeq
 \bigl(
 \nbigi_{\nu^+_m(\vecJ),>0},\leq_{\kappa^+_{m,\vecJ}(\theta^u)}
 \bigr).
\]
The maps $\nutilde^{\pm}_{m,\vecJ}$ induce
the following isomorphisms 
of the partially ordered sets
for any $\theta^u\in\real$:
\[
 \nutilde^-_{m,\vecJ}:
 \bigl(
 \nbigitilde^{\circ}_{\vecJ,>0},
 \leq_{\theta^u}
 \bigr)
\simeq
 \bigl(
 \nbigitilde_{\nu^-_m(\vecJ),<0},
 \leq_{\kappa^-_{m,\vecJ}(\theta^u)}
 \bigr),
\]
\[
 \nutilde^+_{m,\vecJ}:
 \bigl(
 \nbigitilde^{\circ}_{\vecJ,<0},
 \leq_{\theta^u}
 \bigr)
\simeq
 \bigl(
 \nbigitilde_{\nu^+_m(\vecJ),>0},
 \leq_{\kappa^+_{m,\vecJ}(\theta^u)}
 \bigr).
\]
\end{prop}
\pf
It follows from Proposition \ref{prop;20.10.5.31}.
\hfill\qed

\chapter[Local Fourier transform and reductions at $0$]{Local Fourier transform and reductions  at $0$}
\label{section;18.6.3.20}

\section{Introduction to \S\ref{section;18.6.3.20}}

Let $(\nbigv,\nabla)$ be a meromorphic flat bundle
on $(\proj^1,\{0,\infty\})$
with regular singularity at $\infty$.
Let $\nbigi(\nbigv)$ be the set of ramified irregular values at $0$
of $(\nbigv,\nabla)$.
Suppose that $\nbigi(\nbigv)\neq\{0\}$,
and we set $\omega=-\ord(\nbigi(\nbigv))>0$.
We obtain the meromorphic flat bundle
$(V,\nabla):=\nbigs_{\omega}(\nbigv,\nabla)$ on $(\proj^1,\{0,\infty\})$.
Note that
$\nbigi(V)=\nbigs_{\omega}(\nbigi(\nbigv))$.
We also obtain
the meromorphic flat bundle
$\nbigt_{\omega}(\nbigv,\nabla)$ on $(\Delta,0)$,
which extends to the meromorphic flat bundle on
$(\proj^1,\{0,\infty\})$ with regular singularity at $\infty$,
denoted by the same notation $\nbigt_{\omega}(\nbigv,\nabla)$.
We have
$\nbigi(\nbigt_{\omega}(\nbigv))=\nbigt_{\omega}(\nbigi(\nbigv))$.

\subsection{Reduction of $\gbigl^{\gbigf}_{\star}(\nbigv)$}
\label{subsection;24.3.25.30}

We set $\omega^{\circ}=(1+\omega)^{-1}\omega$.
\begin{thm}
\label{thm;24.3.24.1}
There exists the following commutative diagram
of local systems with Stokes structure:
 \[
\begin{CD}
 \nbigt_{\omega^{\circ}}\bigl(
 \gbigl^{\gbigf}_{!}(\nbigv),\vecnbigf
 \bigr)
@>>>
 \nbigt_{\omega^{\circ}}\bigl(
 \gbigl^{\gbigf}_{\ast}(\nbigv),\vecnbigf
 \bigr)
 \\
 @V{\simeq}VV @V{\simeq}VV \\ 
 \bigl(
 \gbigl^{\gbigf}_{!}(\nbigt_{\omega}\nbigv),
 \vecnbigf
 \bigr)
 @>>>
 \bigl(
 \gbigl^{\gbigf}_{\ast}(\nbigt_{\omega}\nbigv),
 \vecnbigf
 \bigr).
\end{CD}
\]
\end{thm}
When $\nbigv=V$,
the theorem says that the morphism of
the $2\pi\seisuu$-equivariant local systems
$\nbigt_{\omega^{\circ}}\bigl(
 \gbigl^{\gbigf}_{!}(V)
 \bigr)
\lrarr
\nbigt_{\omega^{\circ}}\bigl(
 \gbigl^{\gbigf}_{\ast}(V)
 \bigr)$
is identified with
$ \gbigl^{\gbigf}_{!}(\nbigt_{\omega}(V))
 \lrarr
 \gbigl^{\gbigf}_{\ast}(\nbigt_{\omega}(V))$,
which also directly follows from  the stationary phase formula
in \S\ref{subsection;24.3.26.1}.
Note that
$\nbigt_{\omega}(V)=\nbigs_{\omega}\nbigt_{\omega}(\nbigv)$
is regular singular.

\begin{thm}
\label{thm;24.3.25.50}
The $2\pi\seisuu$-equivariant local systems with Stokes structure
$\nbigs_{\omega^{\circ}}\bigl(
\gbigl^{\gbigf}_{\star}(\nbigv),\vecnbigf
 \bigr)$ ($\star=!,\ast$)
are obtained as the extension of
the base tuple
$(\gbigl^{\gbigf}_!(V),\vecnbigf)
 \to
 (\gbigl^{\gbigf}_{\ast}(V),\vecnbigf)$
by the morphisms of
the $2\pi\seisuu$-equivariant local systems
\begin{equation}
\label{eq;24.3.25.21}
 \gbigl^{\gbigf}_!(\nbigt_{\omega}(V))
 \to
 \gbigl^{\gbigf}_{\star}(\nbigt_{\omega}(\nbigv))
 \to
 \gbigl^{\gbigf}_{\ast}(\nbigt_{\omega}(V)).
\end{equation}
\end{thm}

\subsection{Stokes structures of $(\gbigl^{\gbigf}_{\star}(V),\vecnbigf)$}
\label{subsection;24.4.5.111}

It is fundamental for us
to study
$(\gbigl^{\gbigf}_{\star}(V),\vecnbigf)$.
Let $(L,\vecnbigftilde)$ be the $2\pi\seisuu$-equivariant
local system with Stokes structure on $\real$
indexed by $\nbigi(V)$,
corresponding to $(V,\nabla)$.
We shall give two types of explicit descriptions
of $(\gbigl^{\gbigf}_{\star}(V),\vecnbigf)$.

\subsubsection{Local systems with Stokes structure}
In \S\ref{subsection;24.4.5.100},
from $(L,\vecnbigftilde)$,
we shall explicitly construct
$2\pi\seisuu$-equivariant
local systems with Stokes structure
$\gbigf^{(0,\infty)}_{+,\star}(L,\vecnbigftilde)=
(\gbigq^0_{\star}(L,\vecnbigftilde)_{\real},\vecnbigf)$ $(\star=!,\ast)$
and morphisms of $2\pi\seisuu$-equivariant local systems
\[
L\lrarr \gbigq^0_!(L,\vecnbigftilde)_{\real}
\lrarr \gbigq^0_{\ast}(L,\vecnbigftilde)_{\real}
\lrarr L.
\]

\begin{thm}
\label{thm;24.4.5.110}
There exists the following commutative diagram of
$2\pi\seisuu$-equivariant local systems with Stokes structure:
\[
 \begin{CD}
   \gbigf^{(0,\infty)}_{+,!}(L,\vecnbigftilde)
  @>{F_{\gbigq^0}}>>
  \gbigf^{(0,\infty)}_{+,\ast}(L,\vecnbigftilde)\\
  @V{\simeq}VV @V{\simeq}VV \\
  (\gbigl^{\gbigf}_!(V),\vecnbigftilde)
  @>>>
  (\gbigl^{\gbigf}_{\ast}(V),\vecnbigftilde).
 \end{CD}
\]
We also have the following commutative diagram of
the local systems
\[
\begin{CD}
 L @>>> \gbigq^0_!(L,\vecnbigftilde)_{\real}
  @>{F_{\gbigq^0}}>> \gbigq^0_{\ast}(L,\vecnbigftilde)_{\real}
  @>>>
  L \\
  @V{\simeq}VV @V{\simeq}VV @V{\simeq}VV @V{\simeq}VV \\
  \gbigl^{\gbigf}_!(V^{\reg})
  @>>>
  \gbigl^{\gbigf}_!(V)
  @>>>
  \gbigl^{\gbigf}_{\ast}(V)
  @>>>
  \gbigl^{\gbigf}_{\ast}(V^{\reg}).
\end{CD}  
\]
In the diagrams, the lower horizontal arrows are
the natural morphisms.
\end{thm}

\begin{rem}
In {\rm\S\ref{subsection;25.2.25.20}},
we describe the constructible subsheaves
$\gbigq^0_{\star}(L,\vecnbigftilde)^{<0}_{\real}
\subset
\gbigq^0_{\star}(L,\vecnbigftilde)^{\leq 0}_{\real}
\subset
\gbigq^0_{\star}(L,\vecnbigftilde)_{\real}$,
and the induced filtrations
on $H^0(\vecJ,\gbigq^0_{\star}(L,\vecnbigftilde)_{\vecJ,>0})$
and 
$H^0(\vecJ,\gbigq^0_{\star}(L,\vecnbigftilde)_{\vecJ,<0})$.
\hfill\qed
\end{rem}

\subsubsection{Stokes shells of
$\gbigf^{(0,\infty)}_{+,\star}(L,\vecnbigftilde)$}
In \S\ref{subsection;18.5.7.10},
we shall introduce an explicit construction of
a base tuple of Stokes shells
$\bigl(
\gbigf^{(0,\infty)}_{!}(\Sh),
\gbigf^{(0,\infty)}_{\ast}(\Sh),F
\bigr)$
in $\Shcat(\gbigf^{(0,\infty)}_+(\nbigi(V)))$
from a Stokes shell $\Sh$
in 
$\Shcat(\nbigi(V))$.

\begin{prop}
\label{prop;24.3.25.60}
There exists the following commutative diagram:
\begin{equation}
\begin{CD}
 \gbigf^{(0,\infty)}_{+,!}(\Shsf(L,\vecnbigftilde))
 @>{F}>>
 \gbigf^{(0,\infty)}_{+,\ast}(\Shsf(L,\vecnbigftilde))
 \\
 @V{\simeq}VV @V{\simeq}VV \\
 \Shsf\bigl(
 \gbigf^{(0,\infty)}_{+,!}(L,\vecnbigftilde)
 \bigr)
 @>>>
 \Shsf\bigl(
 \gbigf^{(0,\infty)}_{+,\ast}(L,\vecnbigftilde)
 \bigr).
\end{CD}
\end{equation}
As a result,
we can identify
$\Shsf(\gbigl^{\gbigf}_!(V))
 \to
 \Shsf(\gbigl^{\gbigf}_{\ast}(V))$
 with
$\gbigf^{(0,\infty)}_{+,!}(\Shsf(L,\vecnbigftilde))
 \to
 \gbigf^{(0,\infty)}_{+,\ast}(\Shsf(L,\vecnbigftilde))$.
\end{prop}

\subsection{Inductive procedure}
\label{subsection;24.4.2.110}

These theorems provide us
with the following procedure to study
$(\gbigl^{\gbigf}_{\star}(\nbigv),\vecnbigf)$ inductively.
\begin{itemize}
 \item
$(\gbigl^{\gbigf}_{\star}(\nbigv),\vecnbigf)$
are obtained from 
$\nbigs_{\omega^{\circ}}\bigl(
\gbigl^{\gbigf}_{\star}(\nbigv),\vecnbigf
 \bigr)$
and 
\[
      \nbigt_{\omega^{\circ}}(\gbigl^{\gbigf}_{\star}(\nbigv),\vecnbigf)
\simeq(\gbigl^{\gbigf}_{\star}(\nbigt_{\omega}(\nbigv)),\vecnbigf). 
\]
      Note that
      $-\ord\nbigi(\nbigt_{\omega}(\nbigv))<\omega$.
 \item $\nbigs_{\omega^{\circ}}\bigl(
\gbigl^{\gbigf}_{\star}(\nbigv),\vecnbigf
 \bigr)$ are explicitly described as
as the extensions of the base tuple
\[
\begin{CD}
 \gbigf^{(0,\infty)}_{+,!}(L,\vecnbigftilde)
 @>>>
\gbigf^{(0,\infty)}_{+,\ast}(L,\vecnbigftilde)
\end{CD}
\]
by (\ref{eq;24.3.25.21}).
\end{itemize}

As the complement to this procedure,
we note that the morphisms of the local systems
\[
 \gbigl^{\gbigf}_!(V^{\reg})
 \lrarr
 \gbigl^{\gbigf}_!(V)
 \lrarr
  \gbigl^{\gbigf}_{\ast}(V)
 \lrarr
  \gbigl^{\gbigf}_{\ast}(V^{\reg})
\]
are explicitly described
by Theorem \ref{thm;24.4.5.110}.
It allows us to describe explicitly
the morphisms
\[
 \gbigl^{\gbigf}_!(V^{\reg})
 \lrarr
 \gbigl^{\gbigf}_{\star}(\nbigv)
 \lrarr 
 \gbigl^{\gbigf}_{\star}(V^{\reg}).
\]
We remark
$V^{\reg}=\nbigv^{\reg}$
and 
$\nbigt_{\omega}(V)=\nbigt_{\omega}(\nbigv)^{\reg}$.

\subsection{Extensions and the recovery of Stokes structure $\vecnbigftilde$}
\label{subsection;25.3.17.2}

Let $M_0$ denote the monodromy automorphism of
$\nbigt_{\omega}(L)$.
There exists the following commutative diagram:
\[
\begin{CD}
 \nbigt_{\omega^{\circ}}\bigl(
 \gbigq^0_{!}(L,\vecnbigf)_{\real}
 \bigr)
@>>>
 \nbigt_{\omega^{\circ}}\bigl(
 \gbigq^0_{\ast}(L,\vecnbigf)_{\real}
 \bigr)
 \\
 @V{\simeq}VV @V{\simeq}VV \\
 \nbigt_{\omega}(L)
 @>{\id-M_{0}^{-1}}>>
 \nbigt_{\omega}(L).
\end{CD}
\]
Let
$\nbigt_{\omega}(L)\stackrel{a}{\lrarr}
L_1
\stackrel{b}{\lrarr}
\nbigt_{\omega}(L)$
be morphisms such that $b\circ a=\id-M_{0}^{-1}$.
We obtain the extension
of $2\pi\seisuu$-equivariant local systems with Stokes structure:
\[
\gbigf^{(0,\infty)}_{+,!}(L,\vecnbigftilde)
\stackrel{u_1}{\lrarr}
(\Ltilde_1,\vecnbigf)
\stackrel{u_2}{\lrarr}
\gbigf^{(0,\infty)}_{+,\ast}(L,\vecnbigftilde).
\]
We also obtain morphisms of $2\pi\seisuu$-equivariant local systems
$L
\stackrel{\atilde}{\lrarr}
\Ltilde_1
\stackrel{\btilde}{\lrarr}
L$.
\begin{prop}[Proposition
\ref{prop;25.2.23.10}]
Let $M_{L_1}$ and $M_{\Ltilde_1}$
denote the monodromy automorphism of $L_1$ and $\Ltilde_1$,
respectively. 
If $b\circ a=\id-M_{L_1}^{-1}$ holds,
then $\btilde\circ\atilde=\id-M_{\Ltilde_1}^{-1}$ holds. 
\end{prop}

In \S\ref{subsection;25.2.23.20},
we explain how to recover
the constructible subsheaves
$L^{<0}\subset L^{\leq 0}\subset L$
and the filtrations on
$H^0(J,L_{J,<0})$
and
$H^0(J,L_{J,>0})$
from $(\Ltilde_1,\vecnbigf)$
and the morphisms
$L\stackrel{\atilde}{\lrarr}\Ltilde_1
\stackrel{\btilde}{\lrarr} L$.

\subsection{Homology groups}

For $u\in\cnum^{\ast}$,
let $\nbige(zu^{-1})$ be the meromorphic flat bundle
$\bigl(
 \nbigo_{\proj^1}(\ast\{\infty\}),d+d(zu^{-1})
 \bigr)$.
We take $\theta^u\in\real$
such that $\theta^u=\arg(u)$,
i.e.,
$\exp(\sqrt{-1}\theta^u)=|u|^{-1}u$.
There exist the natural isomorphisms
\[
 \gbigl^{\gbigf}_{!}(\nbigv)_{|\theta^u}
 \simeq
 H_1^{\rd}\bigl(\cnum^{\ast},
 \nbigv\otimes\nbige(zu^{-1})\bigr),
 \quad
 \gbigl^{\gbigf}_{\ast}(\nbigv)_{|\theta^u}
 \simeq
 H_1^{\mg}\bigl(\cnum^{\ast},
 \nbigv\otimes\nbige(zu^{-1})\bigr).
\]
To obtain the theorems in
\S\ref{subsection;24.3.25.30}--\ref{subsection;24.4.5.111},
we study the rapid decay homology groups
and the moderate growth homology groups
of $(V,\nabla)\otimes\nbige(zu^{-1})$
and $(\nbigv,\nabla)\otimes\nbige(zu^{-1})$.
We set
$\nbigi:=\pi_{\omega}(\nbigi(V))$
and $\nbigi^{\circ}=\gbigf^{(0,\infty)}_+(\nbigi)$.

\subsubsection{Rapid decay homology group
$H_1^{\rd}(\cnum^{\ast},V\otimes\nbige(zu^{-1}))$}

We shall introduce the following maps in
\S\ref{subsection;24.3.25.100} and \S\ref{subsection;24.2.18.12}:
\[
 \Abb^{\rd}_{\infty,\theta^u}:
 H^0(\real,L)
 \lrarr
 H_1^{\rd}\bigl(\cnum^{\ast},
 V\otimes\nbige(zu^{-1})\bigr),
\]
\[
 \Abb^{\rd}_{J,\theta^u}:
 H^0(J,L_{J,<0})
 \lrarr
 H_1^{\rd}\bigl(\cnum^{\ast},
 V\otimes\nbige(zu^{-1})\bigr)
 \quad (J\in T(\nbigi)).
\]
These induce the isomorphism
\[
 H_1^{\rd}(\cnum^{\ast},V\otimes\nbige(zu^{-1}))
 \simeq
\Bigl(
 H^0(\real,L)
 \oplus
 \bigoplus_{J\in T(\nbigi)}
 H^0(J,L_{J,<0})
 \Bigr)
 \Big/\!\!\sim.
\]
(See \S\ref{subsection;24.3.21.3}
for the equivalence relation.)
The $2\pi\seisuu$-action is also defined naturally
on the right hand side.
This is the isomorphism of the local systems
$\gbigl^{\gbigf}_!(V)\simeq \gbigq^0_!(V)_{\real}$
in Theorem \ref{thm;24.4.5.110}.

To study the Stokes structure,
in \S\ref{subsection;24.3.25.110},
we shall introduce the following maps for
any $J\in T(\nbigi)$:
\[
 B^{J_{\pm}}_{\infty,\theta^u}:
 H^0(\real,\nbigt_{\omega}(L))
 \lrarr
 H_1^{\rd}\bigl(
 \cnum^{\ast},V\otimes\nbige(zu^{-1})
 \bigr),
\]
\[
 B_{J_{\pm},\theta^u}:
 H^0(J,L_{J,>0})
 \lrarr
  H_1^{\rd}\bigl(
 \cnum^{\ast},V\otimes\nbige(zu^{-1})
 \bigr).
\]
We obtain the isomorphisms of vector spaces
(\ref{eq;18.4.19.40})
and (\ref{eq;18.5.13.20})
(Proposition \ref{prop;18.4.19.100}).
The both hand sides of
(\ref{eq;18.4.19.40})
and (\ref{eq;18.5.13.20})
are equipped with the filtrations
indexed by the partially ordered set
$\bigl(
\nbigi(\Fourier(V)),\leq_{\theta^u}
\bigr)$.
As in Theorem \ref{thm;24.3.15.10},
(\ref{eq;18.4.19.40})
and (\ref{eq;18.5.13.20})
are isomorphisms of filtered vector spaces.
(The proof of Theorem \ref{thm;24.3.15.10}
will be given in \S\ref{section;20.11.21.1}.)
This gives the isomorphism
$\gbigf^{(0,\infty)}_{+,\star}(L,\vecnbigftilde)
\simeq
(\gbigl^{\gbigf}_!(V),\vecnbigftilde)$
in Theorem \ref{thm;24.4.5.110}.
It also provides us with
the following isomorphisms of the filtered vector spaces
\[
 H^0(\nu_0^-(\vecJ),L_{\nu_0^-(\vecJ),<0})
 \simeq
 H^0(\vecJ,\gbigl^{\gbigf}_{!}(V)_{\vecJ,<0}),
\] 
\[
 H^0(\nu_0^+(\vecJ)_{\pm},L_{\nu_0^+(\vecJ)_{\pm},>0})
 \simeq
 H^0(\vecJ_{\pm},\gbigl^{\gbigf}_!(V)_{\vecJ_{\pm},>0}),
\] 
\[
 H^0(\nu_0^-(\vecJ)_{\pm},L_{\nu_0^-(\vecJ)_{\pm},0})
 \simeq
 H^0(\vecJ_{\pm},\gbigl^{\gbigf}_!(V)_{\vecJ_{\pm},0}).
\]
By the relations among
$\Abb^{\rd}_{J,\theta^u}$,
$B_{J,\theta^u}$
and $B^{J_{\pm}}_{\infty,\theta^u}$
$(J\in T(\nbigi))$,
we obtain that
the Stokes shell of
$(\gbigl^{\gbigf}_!(V),\vecnbigf)
\simeq
\gbigf^{(0,\infty)}_{+,!}(L,\vecnbigftilde)$
is isomorphic to
$\gbigf^{(0,\infty)}_{+,!}(\Shsf(L,\vecnbigftilde))$
as in Proposition \ref{prop;24.3.25.60}.

\subsubsection{Moderate growth homology group
$H_1^{\mg}(\cnum^{\ast},V\otimes\nbige(zu^{-1}))$}

We shall introduce the following maps for $J\in T(\nbigi)$
in \S\ref{subsection;24.3.25.111} and \S\ref{subsection;24.3.25.112}:
\[
 \BB^{\mg}_{J,u}:
 H^0(J,L_{J,>0})
 \lrarr
 H_1^{\mg}\bigl(
 \cnum^{\ast},V\otimes\nbige(zu^{-1})
 \bigr),
\]
\[
 \Abb^{\mg,J_{\pm}}_{\infty,\theta^u}:
 H^0(\real,L)
 \lrarr
  H_1^{\mg}\bigl(
 \cnum^{\ast},V\otimes\nbige(zu^{-1})
 \bigr).
\]
They induce the isomorphism:
\[
 H_1^{\mg}\bigl(
 \cnum^{\ast},
 V\otimes\nbige(zu^{-1})
 \bigr)
 \simeq
 \Bigl(
 \bigoplus_{\pm}
 \bigoplus_{J\in T(\nbigi)}
 H^0(\real,L)
 \oplus
 \bigoplus_{J\in T(\nbigi)}
 H^0(J,L_{J,>0})
 \Bigr)
 \Big/
 \!\!\sim.
\]
(See \S\ref{subsection;24.3.21.4}
for the equivalence relation.)
The $2\pi\seisuu$-action is defined naturally.
This gives the isomorphism of
the $2\pi\seisuu$-equivariant local systems
$\gbigl^{\gbigf}_{\ast}(V)\simeq
 \gbigq^0_{\ast}(V)_{\real}$
in Theorem \ref{thm;24.4.5.110}.

To study the Stokes structure,
in \S\ref{subsection;24.3.21.13},
we shall introduce
\[
 B^{\mg,J_{\pm}}_{\infty,\theta^u}:
 H^0(\real,\nbigt_{\omega}(L))
 \lrarr
 H_1^{\mg}(\cnum^{\ast},V\otimes\nbige(zu^{-1})).
\]
We obtain the isomorphisms of the vector spaces
{\rm(\ref{eq;24.2.20.51})}
and {\rm(\ref{eq;24.2.20.52})}.
The both hand sides of
{\rm(\ref{eq;24.2.20.51})}
and {\rm(\ref{eq;24.2.20.52})}
are equipped with the filtrations
indexed by the partially ordered set
$\bigl(
 \nbigi(\Fourier(V)),\leq_{\theta^u}
 \bigr)$.
We shall prove that  
{\rm(\ref{eq;24.2.20.51})}
and {\rm(\ref{eq;24.2.20.52})}
are isomorphisms of filtered vector spaces
(Theorem \ref{thm;24.3.15.10}).
It gives an isomorphism
$\gbigf^{(0,\infty)}_{+,\ast}(L,\vecnbigftilde)
\simeq
(\gbigl^{\gbigf}_{\ast}(V),\vecnbigftilde)$
in Theorem \ref{thm;24.4.5.110}.
It also provides us with
the following isomorphisms of the filtered vector spaces
\[
 H^0(\nu_0^-(\vecJ),L_{\nu_0^-(\vecJ),<0})
 \simeq
 H^0(\vecJ,\gbigl^{\gbigf}_{\ast}(V)_{\vecJ,<0}),
\] 
\[
 H^0(\nu_0^+(\vecJ)_{\pm},L_{\nu_0^+(\vecJ)_{\pm},>0})
 \simeq
 H^0(\vecJ_{\pm},\gbigl^{\gbigf}_{\ast}(V)_{\vecJ_{\pm},>0}),
\] 
\[
 H^0(\nu_0^-(\vecJ)_{\pm},L_{\nu_0^-(\vecJ)_{\pm},0})
 \simeq
 H^0(\vecJ_{\pm},\gbigl^{\gbigf}_{\ast}(V)_{\vecJ_{\pm},0}).
\]
By the relations among
$\Abb^{\rd}_{J,\theta^u}$,
$B_{J,\theta^u}$
and $B^{\mg,J_{\pm}}_{\infty,\theta^u}$
$(J\in T(\nbigi))$,
we obtain that
the Stokes shell of
$(\gbigl^{\gbigf}_{\ast}(V),\vecnbigf)$
is isomorphic to
$\gbigf^{(0,\infty)}_{\ast}(\Sh(V))$
as in Proposition \ref{prop;24.3.25.60}.
 
\subsubsection{Homology groups
$H_1^{\rd}(\cnum^{\ast},\nbigv\otimes\nbige(zu^{-1}))$
and $H_1^{\mg}(\cnum^{\ast},\nbigv\otimes\nbige(zu^{-1}))$}

We shall introduce the following maps in
\S\ref{subsection;24.2.20.30} and \S\ref{subsection;24.3.21.13}:
\[
 C^{J_{\pm}}_{\infty,\theta^u}:
 H_1^{\mg}(\cnum^{\ast},\nbigt_{\omega}(\nbigv)\otimes\nbige(zu^{-1}))
 \lrarr
 H_1^{\rd}(\cnum^{\ast},\nbigv\otimes\nbige(zu^{-1})),
\]
\[
 C^{\mg,J_{\pm}}_{\infty,\theta^u}:
 H_1^{\mg}(\cnum^{\ast},\nbigt_{\omega}(\nbigv)\otimes\nbige(zu^{-1}))
 \lrarr
 H_1^{\mg}(\cnum^{\ast},\nbigv\otimes\nbige(zu^{-1})).
\]
We obtain the isomorphisms of the vector spaces
(\ref{eq;24.2.20.21}), (\ref{eq;24.2.20.25}),
(\ref{eq;24.2.20.51})
and (\ref{eq;24.2.20.52}).
The both hand sides of
(\ref{eq;24.2.20.21}), (\ref{eq;24.2.20.25}),
{\rm(\ref{eq;24.2.20.51})}
and {\rm(\ref{eq;24.2.20.52})}
are equipped with the filtrations
indexed by the partially ordered set
$\bigl(
 \nbigi(\Fourier(\nbigv)),\leq_{\theta^u}
 \bigr)$.
We shall prove that
(\ref{eq;24.2.20.21}), (\ref{eq;24.2.20.25}),
{\rm(\ref{eq;24.2.20.51})}
and {\rm(\ref{eq;24.2.20.52})}
are isomorphisms of filtered vector spaces
(Theorem \ref{thm;24.3.15.10}).
It implies Theorem \ref{thm;24.3.24.1}
and Theorem \ref{thm;24.3.25.50}.
We also obtain the isomorphisms
(\ref{eq;24.3.25.52}) and (\ref{eq;24.3.25.51}).

\subsection{Notation}
\label{subsection;18.5.13.50}

Let $\projtilde^1$
be the oriented real blow up
of $\proj^1$ along $\{0,\infty\}$.
\index{oriented real blow up $\projtilde^1$}
Let $\realbar_{\geq 0}:=\closedclosed{0}{\infty}$.
\index{set $\realbar_{\geq 0}$}
We identify 
$\projtilde^1$
with $\realbar_{\geq 0}\times S^1$,
which preserves the natural orientations.
We set 
$X:=\realbar_{\geq 0}\times\real$
and $X^{\ast}:=\real_{>0}\times\real$.
\index{spaces $X$ and $X^{\ast}$}
For any subset $Z\subset X$,
let $\iota_Z$ denote the inclusion $Z\lrarr X$,
and $q_Z$ denote the projection $Z\lrarr \real$.
\index{maps $\iota_Z$, $q_Z$}
Let $\varphi:X
\lrarr\realbar_{\geq 0}\times S^1$
be the map given by
$\varphi(r,\theta)=(r,e^{\sqrt{-1}\theta})$.
Similarly, let $\varphi_1:\real\lrarr S^1$
be given by $\varphi_1(\theta)=e^{\sqrt{-1}\theta}$.
\index{maps $\varphi$, $\varphi_1$}
For any subset $A\subset \real$,
let $a_A$ denote the inclusion $A\lrarr\real$.
\index{map $a_A$}

In the following,
a path on $(X^{\ast},X)$ means 
a continuous map
$\gamma:[0,1]\lrarr X$
such that 
$\gamma(\openopen{0}{1})\subset
 X^{\ast}$.
\index{path on $(X^{\ast},X)$}
We say that a path $\gamma$ on $(X^{\ast},X)$
connects
$P$ to $Q$
if $\gamma(0)=P$
and $\gamma(1)=Q$.

Let $(L,\vecnbigftilde)$ be the $2\pi\seisuu$-equivariant
local system with Stokes structure on $\real$
indexed by $\nbigitilde$,
corresponding to $(V,\nabla)$.
We set
$\nbigi:=\pi_{\omega}(\nbigitilde)$
and 
$\vecnbigf:=\pi_{\omega\ast}(\vecnbigftilde)$.
Let $\Tbb:\real\to\real$ denote the map
defined by $\Tbb(\theta)=\theta+2\pi$.
We have the isomorphism
$\Tbb^{\ast}(L)\simeq L$
\index{map $\Tbb$}
because $L$ is $2\pi\seisuu$-equivariant.
For any $\theta\in\real$,
let $M:L_{|\theta}\to L_{|\theta}$
denote the monodromy automorphism induced by
the parallel transport $L_{|\theta}\simeq L_{|\theta+2\pi}$
and the isomorphism
$L_{|\theta+2\pi}=\Tbb^{\ast}(L)_{|\theta}\simeq L_{|\theta}$.
It also induces the automorphism $M$ on $H^0(\real,L)$.
Similarly,
let $M_0$ denote the automorphism of
$H^0(\real,\nbigt_{\omega}(L))$
obtained as the monodromy.

Let $\nbigl$ denote the local system
on $\projtilde^1$
corresponding to 
$(V,\nabla)_{|\cnum^{\ast}}
=(\nbigv,\nabla)_{|\cnum^{\ast}}$.
We set
$\nbigitilde:=\nbigi(V)$.
We have the natural identifications
$\varphi^{\ast}(\nbigl)_{|\{0\}\times\real}
=L$
and 
$\varphi^{\ast}(\nbigl)
=q_X^{-1}(L)$.

We take $\theta^u\in\real$
such that $\theta^u=\arg(u)$,
i.e.,
$\exp(\sqrt{-1}\theta^u)=|u|^{-1}u$.
We set
\index{interval $\vecI(\theta^u)$}
\[
\vecI(\theta^u)=I(\theta^u-\pi,\pi/2)=
\openopen{\theta^u-3\pi/2}{\theta^u-\pi/2}.
\]
We have
$\Re(zu^{-1})>0$
if and only if
$\arg(z)\in\vecI(\theta^u)+(2m+1)\pi$ for some $m\in\seisuu$,
and $\Re(zu^{-1})<0$
if and only if
$\arg(z)\in\vecI(\theta^u)+2m\pi$ for some $m\in\seisuu$.

\section{Rapid decay homology group of $V\otimes\nbige(zu^{-1})$}

\subsection{Description of
$H_1^{\rd}\bigl(\cnum^{\ast},V^{\reg}\otimes\nbige(zu^{-1})\bigr)$}
\label{subsection;24.3.25.100}

Let $\Gamma_{\infty,\theta^u}$ be any path
on $(X^{\ast},X)$
connecting $(\infty,\theta^u-2\pi)$
and $(\infty,\theta^u)$.
\index{path $\Gamma_{\infty,\theta^u}$}
Any $v\in H^0(\real,L)$
induces
a flat section of $\varphi^{\ast}\nbigl$
along $\Gamma_{\infty,\theta^u}$,
which is also denoted by $v$.
We can naturally regard 
$\varphi_{\ast}\bigl(
 \Gamma_{\infty,\theta^u}\otimes v
 \bigr)$
as a rapid decay $1$-cycle for 
$V^{\reg}\otimes \nbige(zu^{-1})$.
It is easy to see the following lemma.
\begin{lem}
\label{lem;18.4.17.20}
The above procedure induces
an isomorphism of the vector spaces:
\begin{equation}
 \label{eq;18.4.19.2}
 H^0(\real,L)\simeq
  H_1^{\rd}\bigl(\cnum^{\ast},
 V^{\reg}\otimes \nbige(zu^{-1})\bigr).
\end{equation}
We can also regard it as an isomorphism
\[
 \varphi^{\ast}(\nbigl)_{|(\infty,\theta^u)}
=L_{|\theta^u}
\simeq
 H_1^{\rd}\bigl(\cnum^{\ast},
 V^{\reg}\otimes \nbige(zu^{-1})\bigr)
\]
by using the natural isomorphism
$H^0(\real,L)\simeq L_{|\theta^u}$.
\hfill\qed
\end{lem}

Let $\Abb^{\rd}_{\infty,\theta^u}$  denote the composition
of the following maps: \index{map $\Abb^{\rd}_{\infty,\theta^u}$}
\[
 H^0(\real,L)
 \to
H_1^{\rd}(\cnum^{\ast},V^{\reg}\otimes\nbige(zu^{-1}))
\to
 H_1^{\rd}(\cnum^{\ast},V\otimes\nbige(zu^{-1})).
\]
Because
$\varphi_{\ast}(\Gamma_{\infty,\theta^u+2\pi}\otimes v)
=\varphi_{\ast}(\Gamma_{\infty,\theta^u}\otimes M(v))$
for any $v\in H^0(\real,L)$,
we obtain the following lemma.

\begin{lem}
\label{lem;24.4.5.2}
We have
$\Abb^{\rd}_{\infty,\theta^u+2\pi}
=\Abb^{\rd}_{\infty,\theta^u}\circ M$.
\hfill\qed
\end{lem}

\subsection{Rapid decay $1$-homology classes
$\Abb^{\rd}_{J,\theta^u}(v)$}
\label{subsection;24.2.18.12}

For any $J\in T(\nbigi)$,
we take  a path $\Gamma_{J,\theta^u}$
on $(X^{\ast},X)$
connecting a point in $\{0\}\times J$
and $(\infty,\theta^u)$.
\index{path $\Gamma_{J,\theta^u}$}
Any $v\in H^0(J,L_{J,<0})$
induces a flat section 
$\varphi^{\ast}\nbigl$
along $\Gamma_{J,\theta^u}$,
which is also denoted by $v$.
We obtain the rapid decay $1$-cycle
$\varphi_{\ast}\bigl(
 v\otimes \Gamma_{J,\theta^u}
 \bigr)$
for $(V,\nabla)$.
Let $\Abb^{\rd}_{J,\theta^u}(v)$ denote the homology class.
We obtain the following maps: \index{map $\Abb^{\rd}_{J,\theta^u}$}
\[
 \Abb^{\rd}_{J,\theta^u}:
H^0(J,L_{J,<0})
\lrarr
 H_1^{\rd}\bigl(
 \cnum^{\ast},V\otimes \nbige(zu^{-1})\bigr).
\]
The following lemma is clear by the construction.
\begin{lem}
\label{lem;24.4.5.1}
By the isomorphism
$\Tbb^{\ast}:H^0(J+2\pi,L_{J+2\pi,<0})\simeq H^0(J,L_{J,<0})$,
we obtain the following equality on $H^0(J+2\pi,L_{J+2\pi,<0})$:
\[
 \Abb^{\rd}_{J+2\pi,\theta^u+2\pi}
=\Abb^{\rd}_{J,\theta^u}\circ\Tbb^{\ast}.
\]
By the natural inclusion
$\rho_J:H^0(J,L_{J,<0})\lrarr H^0(\real,L)$,
we obtain the following equality on $H^0(J,L_{J,<0})$:
\index{map $\rho_J$}
\[
\Abb^{\rd}_{J,\theta^u}-\Abb^{\rd}_{J,\theta^u-2\pi}
=\Abb^{\rd}_{\infty,\theta^u}\circ\rho_J.
\] 
As a result,
we obtain
$\Abb^{\rd}_{J,\theta^u}
-\Abb^{\rd}_{J+2\pi,\theta^u}\circ(\Tbb^{\ast})^{-1}
=\Abb^{\rd}_{\infty,\theta^u}\circ\rho_J$
on $H^0(J,L_{J,<0})$.
\hfill\qed
\end{lem}

\subsection{Exact sequence and splittings}
\label{subsection;18.5.13.20}

For each $J\in T(\nbigi)$,
there exists the natural morphism
$\varphi_{1!}
 a_{J!}(L_{J,<0})
\lrarr
 L_{S^1}$.
Let $\gbigt(\nbigi,\theta^u)$ be the set of
the intervals $J\in T(\nbigi)$ such that 
$\vartheta^{\vecI(\theta^u)}_{\ell}<\vartheta^J_{\ell}\leq
\vartheta^{\vecI(\theta^u)}_{\ell}+2\pi$.
\index{set $\gbigt(\nbigi,\theta^u)$}
We obtain the following isomorphism:
\[
 \bigoplus_{J\in \gbigt(\nbigi,\theta^u)}
 \varphi_{1!}
 a_{J!}(L_{J,<0})
\simeq
 L^{<0}_{S^1}.
\]

\begin{lem}
We obtain the following exact sequence:
\begin{multline}
\label{eq;18.4.19.1}
 0\lrarr
 H_1^{\rd}\bigl(
 \cnum^{\ast},V^{\reg}\otimes \nbige(zu^{-1})\bigr)
\stackrel{c_{1,u}}{\lrarr}
 H_1^{\rd}\bigl(\cnum^{\ast},V\otimes \nbige(zu^{-1})\bigr)
\stackrel{c_{2,u}}{\lrarr} \\
\bigoplus_{J\in \gbigt(\nbigi,\theta^u)}
 H^0(J,L_{J,<0})
\lrarr 0.
\end{multline}
\end{lem}
\pf
Take $\omega'>\omega$.
Let $b:\varpi^{-1}(0)\lrarr \projtilde^1$
denote the inclusion.
We have
\[
 \nbigq_{\omega',0}^{<0}(V,\nabla)=
 b_{\ast}\Bigl(
 \bigoplus_{J\in \gbigt(\nbigi,\theta^u)}
 \varphi_{1!}a_{J!}L_{J,<0}
 \Bigr).
\]
(See \S\ref{subsection;18.5.15.40} for
$\nbigq_{\omega',0}^{<0}(V,\nabla)$.)
By Lemma \ref{lem;20.10.14.1},
we obtain
\[
 \hyperh^{-i}\bigl(
 \varpi^{-1}(0),
 \nbigc_{\varpi^{-1}(0)}^{\bullet}
 \otimes\nbigq_{\omega',0}^{<0}(V,\nabla)
\bigr)
\simeq
 \bigoplus_{J\in\gbigt(\nbigi,\theta^u)}
 H_i\bigl(J,L_{J,<0}\bigr).
\]
Clearly, we obtain
$H_i\bigl(J,L_{J,<0}\bigr)=0$ unless $i=0$.
It is easy to see
$H_0^{\rd}\bigl(\cnum^{\ast},
V^{\reg}\otimes \nbige(zu^{-1})\bigr)=0$.
We also have
$H_0(J,L_{J,<0})=H^0(J,L_{J,<0})$.
Hence, we obtain the desired exact sequence
from (\ref{eq;18.4.17.10}).
\hfill\qed

\vspace{.1in}
The following lemma gives 
a splitting of the exact sequence (\ref{eq;18.4.19.1}).

\begin{lem}
\label{lem;18.4.22.121}
The maps $\Abb^{\rd}_{\infty,\theta^u}$
and $\Abb^{\rd}_{J,\theta^u}$ $(J\in\gbigt(\nbigi,\theta^u))$
induce an isomorphism
\begin{equation}
\label{eq;18.4.19.3}
 H^0(\real,L)
\oplus
\bigoplus_{J\in \gbigt(\nbigi,\theta^u)}
H^0(J,L_{J,<0})
\lrarr
  H_1^{\rd}\bigl(
 \cnum^{\ast},V\otimes \nbige(zu^{-1})\bigr).
\end{equation}
\end{lem}
\pf
Let us look at the following morphisms:
\[
H^0(J,L_{J,<0})
\stackrel{\Abb^{\rd}_{J,\theta^u}}{\lrarr}
   H_1^{\rd}\bigl(
 \cnum^{\ast},V\otimes \nbige(zu^{-1})\bigr)
\stackrel{c_{2,u}}{\lrarr}
 \bigoplus_{J'\in \gbigt(\nbigi,\theta^u)}
 H_0(J',L_{J',<0}).
\]
By the construction,
the induced map
$H^0(J,L_{J,<0})
\lrarr
  H_0(J,L_{J,<0})$
is an isomorphism.
If $J'\neq J$,
the induced map
$H^0(J,L_{J,<0})
\lrarr
  H_0(J',L_{J',<0})$ is $0$.
Then, we obtain the claim of the lemma.
\hfill\qed

\subsection{Some useful classes}
\label{subsection;24.3.25.110}

We introduce some rapid decay homology classes
which are useful in our study of
$(\gbigl^{\gbigf}_{\star}(V),\vecnbigf)$ $(\star=!,\ast)$.
(See \S\ref{subsection;24.3.14.41})

\subsubsection{Rapid decay homology classes
$B^{J_{\pm}}_{\infty,\theta^u}(v)$}
\label{subsection;25.2.11.30}

For any $J\in T(\nbigi)$,
there exist the natural isomorphisms:
\[
 H^0(\real,\nbigt_{\omega}(L))
\simeq H^0(J,L_{J,0})
\simeq H^0(J_{\pm},L_{J_{\pm},0}).
\]
For $v\in H^0(\real,\nbigt_{\omega}(L))$,
let $v_J$ denote the image in $H^0(J,L_{J,0})$,
and let $v_{J_{\pm}}$ denote the images in
$H^0(J_{\pm},L_{J_{\pm},0})$.
Let $\rho_{J_{\pm}}:H^0(J_{\pm},L_{J_{\pm},0})\lrarr H^0(\real,L)$
denote the natural inclusions.
\index{maps $\rho_{J_{\pm}}$}
We set
\begin{equation}
\label{eq;25.2.22.1}
 B^{J_{+}}_{\infty,\theta^u}(v)
=\Abb^{\rd}_{\infty,\theta^u}(\rho_{J_{+}}(v_{J_{+}}))
+\sum_{J-2\pi<J'\leq J}
 \Abb^{\rd}_{J',\theta^u-2\pi}(\nbigp_{J'_+}(v_{J'})),
\end{equation}
\begin{equation}
\label{eq;25.2.22.2}
 B^{J_{-}}_{\infty,\theta^u}(v)
=\Abb^{\rd}_{\infty,\theta^u}(\rho_{J_{-}}(v_{J_{-}}))
+\sum_{J-2\pi\leq J'<J}
 \Abb^{\rd}_{J',\theta^u-2\pi}(\nbigp_{J'_+}(v_{J'})).
\end{equation}
We obtain the following linear maps:
\index{maps $B^{J_{\pm}}_{\infty,\theta^u}$}
\[
  B^{J_{\pm}}_{\infty,\theta^u}:
 H^0(\real,\nbigt_{\omega}(L))
 \lrarr
 H_1^{\rd}(\cnum^{\ast},V\otimes\nbige(zu^{-1})).
\]
Because $\nbigt_{\omega}(V,\nabla)$ is regular singular at $\{0,\infty\}$,
there exists the natural isomorphism
as in Lemma \ref{lem;18.4.17.20}:
\begin{equation}
\label{eq;24.2.21.1}
 H_1^{\rd}(\cnum^{\ast},\nbigt_{\omega}(V)\otimes\nbige(zu^{-1}))
 \simeq
 H^0(\real,\nbigt_{\omega}(L)).
\end{equation}
We can regard $B^{J_{\pm}}_{\infty,\theta^u}$
as maps
\[
 B^{J_{\pm}}_{\infty,\theta^u}:
 H_1^{\rd}(\cnum^{\ast},\nbigt_{\omega}(V)\otimes\nbige(zu^{-1}))
 \lrarr
 H_1^{\rd}(\cnum^{\ast},V\otimes\nbige(zu^{-1})).
\]
The following lemma is clear by the construction.
\begin{lem}
If $J_2\vdash J_1$,
$B^{J_{2+}}_{\infty,\theta^u}=B^{J_{1-}}_{\infty,\theta^u}$.
\hfill\qed
\end{lem}

Let $M_0$ denote the automorphism of
$H^0(\real,\nbigt_{\omega}(L))$
obtained as the monodromy.
\begin{lem}
\label{lem;25.2.11.10}
For any $v\in H^0(\real,\nbigt_{\omega}(L))$,
we obtain
\begin{equation}
\label{eq;24.2.19.10}
 B^{J_-}_{\infty,\theta^u}(v)
 -B^{J_+}_{\infty,\theta^u}(v)
 =\Abb^{\rd}_{J,\theta^u}\bigl(
 \nbigp_{J}(v-M_0^{-1}(v))
 \bigr).
\end{equation}
\end{lem}
\pf
We recall
$\nbigp_{J}=\nbigp_{J_-}=-\nbigp_{J_+}$.
We omit to denote $\rho_{J_{\pm}}$.
We have
\begin{multline}
 B^{J_-}_{\infty,\theta^u}(v)
-B^{J_+}_{\infty,\theta^u}(v)
 =\\
 \Abb^{\rd} _{\infty,\theta^u}(v_{J_-})
+\Abb^{\rd}_{J-2\pi,\theta^u-2\pi}(\nbigp_{(J-2\pi)_+}(v))
-\Abb^{\rd}_{\infty,\theta^u}(v_{J_+})
 -\Abb^{\rd}_{J,\theta^u-2\pi}(\nbigp_{J_+}(v))\\
= \Abb^{\rd} _{\infty,\theta^u}(v_{J_-})
-\Abb^{\rd}_{J-2\pi,\theta^u-2\pi}(\nbigp_{J-2\pi}(v))
-\Abb^{\rd}_{\infty,\theta^u}(v_{J_+})
+\Abb^{\rd}_{J,\theta^u-2\pi}(\nbigp_{J}(v)).
\end{multline}
We have
\[
 \Abb^{\rd}_{\infty,\theta^u}(v_{J_-})
-\Abb^{\rd}_{\infty,\theta^u}(v_{J_+})
=\Abb^{\rd}_{\infty,\theta^u}(\nbigp_J(v))
=\Abb^{\rd}_{J,\theta^u}(\nbigp_{J}(v))
-\Abb^{\rd}_{J,\theta^u-2\pi}(\nbigp_{J}(v)).
\]
We obtain
\[
 B^{J_-}_{\infty,\theta^u}(v)
-B^{J_+}_{\infty,\theta^u}(v)
=\Abb^{\rd}_{J,\theta^u}(\nbigp_{J}(v))
-\Abb^{\rd}_{J-2\pi,\theta^u-2\pi}(\nbigp_{J-2\pi}(v)).
\]
Then, we obtain (\ref{eq;24.2.19.10}).
\hfill\qed

\begin{cor}
\label{cor;25.2.10.11}
For any $J_1<J_2$,
we obtain
\[
  B^{J_{1-}}_{\infty,\theta^u}
-B^{J_{2-}}_{\infty,\theta^u}
=\sum_{J_1\leq J'<J_2}
 \Abb^{\rd}_{J',\theta^u}
 \circ\nbigp_{J'}\circ
 (\id-M_0^{-1}).
\]
\hfill\qed
 \end{cor}

The following lemma is clear by the construction.
Note that the isomorphism
$H^0(\real,\nbigt_{\omega}(L))
\simeq
H_1^{\rd}(\cnum^{\ast},\nbigt_{\omega}(V)\otimes\nbige(zu^{-1}))$
depends on the choice of $\theta^u$.
\begin{lem}
 $B^{(J+2\pi)_{\pm}}_{\infty,\theta^u+2\pi}
=B^{J_{\pm}}_{\infty,\theta^u}\circ M_0$
as maps on $H^0(\real,\nbigt_{\omega}(L))$.
We have  $B^{(J+2\pi)_{\pm}}_{\infty,\theta^u+2\pi}
=B^{J_{\pm}}_{\infty,\theta^u}$
as maps on 
$H_1^{\rd}(\cnum^{\ast},\nbigt_{\omega}(V)\otimes\nbige(zu^{-1}))$.
\hfill\qed
\end{lem}

\begin{lem}
\label{lem;25.2.10.10}
For any $J-2\pi\leq J_1\leq J$
and any $v\in H^0(\real,\nbigt_{\omega}(L))$,
{\small
\begin{multline}
\label{eq;24.2.19.20}
 B^{J_+}_{\infty,\theta^u}(v)
=\Abb^{\rd}_{\infty,\theta^u}(\rho_{J_{1+}}(v_{J_{1+}}))
+\!\!\!\sum_{J-2\pi<J'\leq J_1}\!\!\!\!\!
 \Abb^{\rd}_{J',\theta^u-2\pi}(\nbigp_{J'_+}(v))
 -\sum_{J_1< J'\leq J}\!\!\!
  \Abb^{\rd}_{J',\theta^u}(\nbigp_{J'_-}(v))
\\
=\Abb^{\rd}_{\infty,\theta^u}(\rho_{J_{1-}}(v_{J_{1-}}))
+\sum_{J-2\pi<J'< J_1}\!\!\!
 \Abb^{\rd}_{J',\theta^u-2\pi}(\nbigp_{J'_+}(v))
 -\sum_{J_1\leq J'\leq J}\!\!\!
  \Abb^{\rd}_{J',\theta^u}(\nbigp_{J'_-}(v)).
 \end{multline}}
\end{lem}
\pf
If $J_1=J$,
we obtain the first equality by definition.
Because
\begin{multline}
 \Abb^{\rd}_{\infty,\theta^u}(\rho_{J_+}(v_{J_+}))
=\Abb^{\rd}_{\infty,\theta^u}(\rho_{J_-}(v_{J_-}))
+\Abb^{\rd}_{J,\theta^u}(\nbigp_{J_+}(v))
 -\Abb^{\rd}_{J,\theta^u-2\pi}(\nbigp_{J_+}(v))
 \\
=\Abb^{\rd}_{\infty,\theta^u}(\rho_{J_-}(v_{J_-}))
-\Abb^{\rd}_{J,\theta^u}(\nbigp_{J_-}(v))
-\Abb^{\rd}_{J,\theta^u-2\pi}(\nbigp_{J_+}(v)),
\end{multline}
we obtain the second equality in the case $J_1=J$.
Suppose we have already proved the claim for $J_1$.
If $J_2\vdash J_1$,
we have $v_{J_{2+}}=v_{J_{1-}}$,
and the first equality in the case $J_2$
is the same as the second equality in the case $J_1$.
We obtain the second equality in the case $J_2$
from the first equality as in the case of $J$.
\hfill\qed

\subsubsection{Rapid decay homology classes
$B_{J_{\pm},\theta^u}(v)$}
\label{subsection;25.2.11.31}

For $J\in T(\nbigi)$ and $v\in H^0(J,L_{J,>0})$,
we set
\index{maps $B_{J_{\pm},\theta^u}$}
\begin{multline}
\label{eq;24.3.16.1}
 B_{J_-,\theta^u}(v)
 :=\sum_{J-\omega^{-1}\pi\leq J'<J}
 \Abb^{\rd}_{J',\theta^u}(\nbigrtilde^{J_-}_{J'}(v))
 -\sum_{J\leq J'<J+\pi}
 \Abb^{\rd}_{J',\theta^u}(\nbigrtilde^{J_-}_{J'}(v))
 \\
 -\sum_{J-\omega^{-1}\pi\leq J'<J-\pi}
 \Abb^{\rd}_{J',\theta^u-2\pi}(\nbigrtilde^{J_-}_{J'}(v)),
\end{multline}
\begin{multline}
\label{eq;24.3.16.2}
 B_{J_+,\theta^u}(v):=
 -\!\!\!\sum_{J<J'\leq J+\omega^{-1}\pi}\!\!\!
 \Abb^{\rd}_{J',\theta^u-2\pi}(\nbigrtilde^{J_+}_{J'}(v)) 
+\!\!\!\sum_{J-\pi<J'\leq J}\!\!\!
 \Abb^{\rd}_{J',\theta^u-2\pi}(\nbigrtilde^{J_+}_{J'}(v)) 
\\
 +\sum_{J+\pi<J'\leq J+\omega^{-1}\pi}
 \Abb^{\rd}_{J',\theta^u}(\nbigrtilde^{J_+}_{J'}(v)).
\end{multline}
(See \S\ref{subsection;24.2.19.1} for
the maps $\nbigrtilde^{J_{\pm}}_{J'}$.)
The following lemma is obvious
by the construction.
\begin{lem}
\label{lem;25.2.17.1}
$B_{(J+2\pi)_{\pm},\theta^u+2\pi}
=B_{J_{\pm},\theta^u}\circ\Tbb^{\ast}$
on $H^0(J+2\pi,L_{J+2\pi,>0})$.
\hfill\qed
\end{lem}

We shall prove the following proposition in
\S\ref{subsection;24.3.25.1}--\ref{subsection;24.3.25.2}.
\begin{prop}
\label{prop;24.2.19.20}
For $v\in H^0(J,L_{J,>0})$,
we obtain 
\begin{multline}
\label{eq;24.2.19.30}
B_{J_-,\theta^u}(v)
-B_{J_+,\theta^u}(v)
=B^{(J+\pi)_+}_{\infty,\theta^u}(\nbigq_J(v))\\
 +\Abb^{\rd}_{J+\pi,\theta^u}
 \bigl(
 \nbigrtilde^{J_-}_{J+\pi}(v)
 \bigr)
 +\Abb^{\rd}_{J-\pi,\theta^u-2\pi}
 \bigl(\nbigrtilde^{J_-}_{J-\pi}(v)\bigr).
\end{multline}
\end{prop}

\subsubsection{Proof of
Proposition \ref{prop;24.2.19.20} in the case $\omega>1$}
\label{subsection;24.3.25.1}

We formally set $\nbigr^{J}_{J'}=0$
if $J'< J-\omega^{-1}\pi$
or $J'> J+\omega^{-1}\pi$.

We obtain
\begin{multline}
\label{eq;25.2.11.3}
B_{J_-,\theta^u}(v)
 =\sum_{J-\omega^{-1}\pi\leq J'<J}
 \Abb^{\rd}_{J',\theta^u}(\nbigr^{J}_{J'}(v))
-\Abb^{\rd}_{J,\theta^u}(\nbigr^{J_-}_{J_+}(v))
\\
 -\sum_{J< J'<J+\pi}
 \Bigl(
 \Abb^{\rd}_{J',\theta^u}(\nbigr^{J}_{J'}(v))
+\Abb^{\rd}_{J',\theta^u}(\nbigp_{J'_-}\nbigq_{J_-}(v)
 \Bigr).
\end{multline}
\begin{multline}
 B_{J_+,\theta^u}(v)
 =\sum_{J-\pi<J'< J}
 \Bigl(
 \Abb^{\rd}_{J',\theta^u-2\pi}(\nbigr^{J}_{J'}(v))
+\Abb^{\rd}_{J',\theta^u-2\pi}(\nbigp_{J'_+}\nbigq_{J_+}(v)
 \Bigr)
\\
+\Abb^{\rd}_{J,\theta^u-2\pi}(\nbigr^{J_+}_{J_-}(v))
-\sum_{J<J'\leq J+\omega^{-1}\pi}
 \Abb^{\rd}_{J',\theta^u-2\pi}(\nbigr^{J}_{J'}(v)). 
\end{multline}
We note $\nbigq_{J_+}=-\nbigq_J$,
and
$\nbigr^{J_+}_{J_-}(v)
=-\nbigr^{J_-}_{J_+}(v)+\nbigp_{J_-}\nbigq_{J_-}(v)$.
We obtain
\begin{multline}
 B_{J_-,\theta^u}(v)
-B_{J_+,\theta^u}(v)
=\sum_{J-\omega^{-1}\pi\leq J'<J}
\Abb^{\rd}_{\infty,J'}(\nbigr^{J}_{J'}(v))
 -\sum_{J<J'\leq J+\omega^{-1}\pi}
 \Abb^{\rd}_{\infty,J'}(\nbigr^{J}_{J'}(v))
 \\
-\sum_{J<J'<J+\pi}
 \Abb^{\rd}_{J',\theta^u}(\nbigp_{J'}\nbigq_{J}(v))
+\sum_{J-\pi<J'<J}
 \Abb^{\rd}_{J',\theta^u-2\pi}(\nbigp_{J'+}\nbigq_{J}(v))
 \\
-\Abb^{\rd}_{\infty,\theta^u}(\nbigr^{J_-}_{J_+}(v))
+\Abb^{\rd}_{\infty,\theta^u-2\pi}(\nbigp_{J_+}\nbigq_J(v)).
\end{multline}
We have
\begin{multline}
 \sum_{J-\omega^{-1}\pi\leq J'<J}
\Abb^{\rd}_{\infty,J'}(\nbigr^{J}_{J'}(v))
 -\sum_{J<J'\leq J+\omega^{-1}\pi}
 \Abb^{\rd}_{\infty,J'}(\nbigr^{J}_{J'}(v))
-\Abb^{\rd}_{\infty,\theta^u}(\nbigr^{J_-}_{J_+}(v))
\\
=\Abb^{\rd}_{\infty,\theta^u}(v_{J_-})
 -\Abb^{\rd}_{\infty,\theta^u}(v_{J_+})
 -\Abb^{\rd}_{\infty,\theta^u}(\nbigr^{J_-}_{J_+}(v))
=\Abb^{\rd}_{\infty,\theta^u}(\nbigq_J(v)_{J_+}).
\end{multline}
Then, we obtain (\ref{eq;24.2.19.30})
in this case.

\subsubsection{Proof of 
Proposition \ref{prop;24.2.19.20} in the case $\omega<1$}

We have
\begin{multline}
B_{J_-,\theta^u}(v)
=\sum_{J-\omega^{-1}\pi\leq J'<J-\pi}
 \Abb^{\rd}_{\infty,\theta^u}(\nbigr^{J}_{J'}(v))
+\sum_{J-\pi\leq J'<J}
 \Abb^{\rd}_{J',\theta^u}(\nbigr^{J}_{J'}(v))
 \\
-\Abb^{\rd}_{J,\theta^u}(\nbigr^{J_-}_{J_+}(v))
-\sum_{J<J'<J+\pi}
 \Bigl(
 \Abb^{\rd}_{J',\theta^u}
 \bigl(\nbigr^{J}_{J'}(v)\bigr)
 +\Abb^{\rd}_{J',\theta^u}\bigl(
 \nbigp_{J'}\nbigq_{J}(v)
 \bigr)
 \Bigr),
\end{multline}
\begin{multline}
B_{J_+,\theta^u}(v)=
 \sum_{J-\pi<J'<J}\Bigl(
 \Abb^{\rd}_{J',\theta^u-2\pi}
 \bigl(\nbigr^J_{J'}(v)\bigr)
+\Abb^{\rd}_{J',\theta^u-2\pi}
 \bigl(\nbigp_{J'_+}\nbigq_{J_+}(v)\bigr)
 \Bigr)
 \\
 +\Abb^{\rd}_{J,\theta-2\pi}(\nbigr^{J_+}_{J_-}(v))
 -\sum_{J<J'\leq J+\pi}
 \Abb^{\rd}_{J',\theta^u-2\pi}
 (\nbigr^{J}_{J'}(v))
 +\sum_{J+\pi<J'\leq J+\omega^{-1}\pi}
 \Abb^{\rd}_{\infty,\theta^u}(\nbigr^{J}_{J'}(v)).
\end{multline}
We obtain
\begin{multline}
B_{J_-,\theta^u}(v)
-B_{J_+,\theta^u}(v)
 =\\
 \sum_{J-\omega^{-1}\pi\leq J'<J-\pi}
 \Abb^{\rd}_{\infty,\theta^u}(\nbigr^J_{J'}(v))
- \sum_{J+\pi< J'\leq J+\omega^{-1}\pi}
 \Abb^{\rd}_{\infty,\theta^u}(\nbigr^J_{J'}(v))\\
+\sum_{J-\pi< J'<J}
 \Abb^{\rd}_{\infty,\theta^u}(\nbigr^J_{J'}(v))
- \sum_{J<J'<J+\pi}
 \Abb^{\rd}_{\infty,\theta^u}(\nbigr^J_{J'}(v))\\
 \\
-\sum_{J<J'<J+\pi}
 \Abb^{\rd}_{J',\theta^u}
 (\nbigp_{J'}\nbigq_{J}(v))
+\sum_{J-\pi<J'<J}
 \Abb^{\rd}_{J',\theta^u-2\pi}
 (\nbigp_{J'_+}\nbigq_{J}(v))
 \\
+\Abb^{\rd}_{J-\pi,\theta^u}(\nbigr^{J}_{J-\pi}(v))
+\Abb^{\rd}_{J+\pi,\theta^u-2\pi}(\nbigr^{J}_{J+\pi}(v))
-\Abb^{\rd}_{J,\theta^u}(\nbigr^{J_-}_{J_+}(v))
-\Abb^{\rd} _{J,\theta^u-2\pi}(\nbigr^{J_+}_{J_-}(v)).
\end{multline}
We have
\begin{multline}
\label{eq;24.2.20.1}
 \Abb^{\rd}_{J-\pi,\theta^u}(\nbigr^{J}_{J-\pi}(v))
+\Abb^{\rd}_{J+\pi,\theta^u-2\pi}(\nbigr^{J}_{J+\pi}(v))
 =\\
 \Abb^{\rd}_{\infty,\theta^u}(\nbigr^J_{J-\pi}(v))
-\Abb^{\rd}_{\infty,\theta^u}(\nbigr^J_{J+\pi}(v))
+\Abb^{\rd}_{J-\pi,\theta^u-2\pi}(\nbigr^J_{J-\pi}(v))
+\Abb^{\rd}_{J+\pi,\theta^u}(\nbigr^J_{J+\pi}(v)),
\end{multline}
\begin{equation}
\label{eq;24.2.20.2}
 -\Abb^{\rd}_{J,\theta^u}(\nbigr^{J_-}_{J_+}(v))
 -\Abb^{\rd} _{J,\theta^u-2\pi}(\nbigr^{J_+}_{J_-}(v))
=-\Abb^{\rd}_{\infty,\theta^u}(\nbigr^{J_-}_{J_+}(v))
+\Abb^{\rd}_{J,\theta^u-2\pi}(\nbigp_{J_+}\nbigq_J(v)).
\end{equation}
Then, we obtain (\ref{eq;24.2.19.30})
as in \S\ref{subsection;24.3.25.1}.

\subsubsection{Proof of 
Proposition \ref{prop;24.2.19.20} in the case $\omega=1$}
\label{subsection;24.3.25.2}

We have
\begin{multline}
 B_{J_-,\theta^u}(v)
 =\sum_{J-\pi\leq J'<J}
 \Abb^{\rd}_{J',\theta^u}(\nbigr^J_{J'}(v))
 -\Abb^{\rd}_{J,\theta^u}(\nbigr^{J_-}_{J_+}(v))
\\
 -\sum_{J<J'<J+\pi}
 \Abb^{\rd}_{J',\theta^u}
 \bigl(
 \nbigr^{J}_{J'}(v)
+\nbigp_{J'} \nbigq_{J}(v)
 \bigr).
\end{multline}
\begin{multline}
 B_{J_+,\theta^u}(v)
 =\sum_{J-\pi<J'<J}
 \bigl(
 \Abb^{\rd}_{J',\theta^u-2\pi}
 \bigl(\nbigr^J_{J'}(v)\bigr)
+\Abb^{\rd} _{J',\theta^u-2\pi}(\nbigp_{J'_+}\nbigq_{J_+}(v))
 \bigr)
 \\
 +\Abb^{\rd}_{J,\theta^u-2\pi}(\nbigr^{J_+}_{J_-}(v))
 -\sum_{J<J'\leq J+\pi}
 \Abb^{\rd}_{J',\theta^u-2\pi}\bigl(
 \nbigr^{J}_{J'}(v)
 \bigr).
\end{multline}
We obtain
\begin{multline}
 B_{J_-,\theta^u}(v)
-B_{J_+,\theta^u}(v)
=\sum_{J-\pi<J'<J}
 \Abb^{\rd}_{\infty,\theta^u}(\nbigr^{J}_{J'}(v))
-\sum_{J<J'<J+\pi}
 \Abb^{\rd}_{\infty,\theta^u}(\nbigr^{J}_{J'}(v))\\
 -\sum_{J<J'<J+\pi}
 \Abb^{\rd}_{J',\theta^u}(\nbigp_{J'_-}\nbigq_{J_-}(v))
 +\sum_{J-\pi<J'<J}
 \Abb^{\rd}_{J',\theta^u-2\pi}(\nbigp_{J'_+}\nbigq_{J}(v))
\\
 +\Abb^{\rd}_{J-\pi,\theta^u}(\nbigr^{J}_{J-\pi}(v))
+\Abb^{\rd}_{J+\pi,\theta^u-2\pi}(\nbigr^{J}_{J+\pi}(v))
-\Abb^{\rd}_{J,\theta^u}(\nbigr^{J_-}_{J_+}(v))
-\Abb^{\rd}_{J,\theta^u-2\pi}(\nbigr^{J_+}_{J_-}(v)).
\end{multline}
By using (\ref{eq;24.2.20.1}) and (\ref{eq;24.2.20.2}),
we obtain (\ref{eq;24.2.19.30})
in this case.
Thus,
the proof of Proposition \ref{prop;24.2.19.20} is completed.
\hfill\qed

\subsection{Decompositions of
$H_1^{\rd}(\cnum^{\ast},V\otimes\nbige(zu^{-1}))$}
\label{subsection;24.2.20.13}

Let $\gbigw_{1}(\nbigi,\vecI(\theta^u)_{\pm})$
be the set of 
$J\in T(\nbigi)$
such that 
$J_{\mp}\cap\vecI(\theta^u)_{\pm}\neq\emptyset$.
Let $\gbigw_{2}(\nbigi,\vecI(\theta^u)_{\pm})$
be the set of
$J \in T(\nbigi)$
such that
$J_{\mp}\cap(\vecI(\theta^u)+\pi)_{\pm}\neq\emptyset$.
\index{sets $\gbigw_{1}(\nbigi,\vecI(\theta^u)_{\pm})$}
\index{sets $\gbigw_{2}(\nbigi,\vecI(\theta^u)_{\pm})$}

Take $J_1\in\gbigw_2(\vecI(\theta^u)_+)$.
We obtain the following map
induced by
$B_{J_{-},\theta^u}$ $(J\in \gbigw_{1}(\nbigi,\vecI(\theta^u)_{+}))$,
$\Abb^{\rd}_{J,\theta^u}$ $(J\in\gbigw_{2}(\nbigi,\vecI(\theta^u)_{+}))$,
and $B^{J_{1,-}}_{\infty,\theta^u}$:
{\small
\begin{multline}
\label{eq;18.4.19.40}
 \bigoplus_{J\in\gbigw_{1}(\nbigi,\vecI(\theta^u)_+)}
 \!\!\!\!\!\!H^0(J_-,L_{J_-,>0})
\oplus\!\!\!\!
 \bigoplus_{J\in\gbigw_{2}(\nbigi,\vecI(\theta^u)_+)}
 \!\!\!\!\!\!H^0(J,L_{J<0})
\oplus
 H_1^{\rd}\bigl(
 \cnum^{\ast},
 \nbigt_{\omega}(V)\otimes
 \nbige(zu^{-1})
 \bigr) \\
\lrarr
 H_1^{\rd}\bigl(
 \cnum^{\ast},V\otimes\nbige(zu^{-1})\bigr).
\end{multline}
}
Similarly,
we take $J_2\in\gbigw_2(\nbigi,\vecI(\theta^u)_-)$.
We obtain the following map
induced by
$B_{J_{+},\theta^u}$ $(J\in \gbigw_{1}(\nbigi,\vecI(\theta^u)_{-}))$,
$\Abb^{\rd}_{J,\theta^u}$ $(J\in\gbigw_{2}(\nbigi,\vecI(\theta^u)_{-}))$,
and $B^{J_{2,+}}_{\infty,\theta^u}$:
{\small
\begin{multline}
\label{eq;18.5.13.20}
 \bigoplus_{J\in\gbigw_{1}(\nbigi,\vecI(\theta^u)_-)}
 \!\!\!\!\!\!\!H^0(J_+,L_{J_+,>0})
\oplus\!\!\!
 \bigoplus_{J\in\gbigw_{2}(\nbigi,\vecI(\theta^u)_-)}
 \!\!\!\!\!\!\!H^0(J,L_{J,<0})
\oplus
 H_1^{\rd}\bigl(
 \cnum^{\ast},
 \nbigt_{\omega}(V)\otimes
 \nbige(zu^{-1})
 \bigr) \\
\lrarr
 H_1^{\rd}\bigl(
 \cnum^{\ast},V\otimes\nbige(zu^{-1})\bigr).
\end{multline}
}
\begin{prop}
\label{prop;18.4.19.100}
The morphisms {\rm(\ref{eq;18.4.19.40})}
and {\rm(\ref{eq;18.5.13.20})}
are isomorphisms.
\end{prop}

We shall prove the claim for (\ref{eq;18.4.19.40})
in \S\ref{subsection;25.4.6.1}.
The claim for (\ref{eq;18.5.13.20})
can be proved similarly.
We set
$\theta_0=\theta^u-3\pi/2=\vartheta^{\vecI(\theta^u)}_{\ell}$.
There exists $J_0\in T(\nbigi)$
such that $\theta_0\in (J_{0})_-$,
and that
$\openclosed{\vartheta^{J_0}_{\ell}}{\theta_0}
\cap S_0(\nbigi)=\emptyset$.
For $J\in T(\nbigi)$,
we obtain $J_0<J$ if and only if
$\theta_0<\vartheta^{J}_{\ell}$.
To simplify the description,
we use the notation
$\gbigw_{i}(\nbigi)$
instead of
$\gbigw_{i}(\nbigi,\vecI(\theta^u)_+)$.

\subsubsection{Preliminary}

Because
\begin{multline}
 L_{|\vartheta^{J_0}_{\ell}}=
 L'_{J_{0-},0|\vartheta^{J_0}_{\ell}}
\oplus
 L'_{J_0,<0|\vartheta^{J_0}_{\ell}}
\oplus
 L'_{(J_0-\omega^{-1}\pi),<0|\vartheta^{J_0}_{\ell}}
\\
 \oplus
\bigoplus_{\vartheta^{J_0}_{\ell}\in J}
 \bigl(
 L'_{J,<0|\vartheta^{J_0}_{\ell}}
\oplus
 L'_{J_-,>0|\vartheta^{J_0}_{\ell}}
 \bigr),
\end{multline}
we obtain
\begin{equation}
 L_{|\vartheta^{J_0}_{\ell}}=
  L'_{J_{0-},0|\vartheta^{J_0}_{\ell}}
\oplus
 \bigoplus_{
 J_0-2\omega^{-1}\pi<J\leq J_0}
 L'_{J,<0|\vartheta^{J_0}_{\ell}}.
\end{equation}
We set
$H=\bigoplus_{J\in \gbigt(\nbigi)}
\Image \Abb^{\rd}_{J,\theta^u}$.
(See \S\ref{subsection;18.5.13.20}
for $\gbigt(\nbigi)$.)
\begin{lem}
For any $0\leq a<2\omega^{-1}$,
we obtain
\begin{equation}
\label{eq;25.2.11.1}
 \Abb^{\rd}_{\infty,\theta^u}
\Bigl(
\bigoplus_{J_0-a\pi\leq J\leq J_0}
 L'_{J,<0|\vartheta_{\ell}^{J_0}}
\Bigr)
 \oplus H
=
 \bigoplus_{J_0-a\pi\leq J\leq J_0}
 \Image \Abb^{\rd}_{J,\theta^u}
 \oplus H
\end{equation}
As a result, we obtain
\begin{equation}
\label{eq;25.2.11.2}
 \Abb^{\rd}_{\infty,\theta^u}
\Bigl(
\bigoplus_{J_0-2\omega^{-1}\pi<J\leq J_0}
 L'_{J,<0|\vartheta_{\ell}^{J_0}}
\Bigr)
 \oplus H
=
 \bigoplus_{J_0-2\omega^{-1}\pi<J\leq J_0}
 \Image \Abb^{\rd}_{J,\theta^u}
 \oplus H
\end{equation}
\end{lem}
\pf
For $0\leq a\leq 2\omega^{-1}$,
let $H_a$ denote the right hand side of (\ref{eq;25.2.11.1}).
We set $H_{<a}=\sum_{b<a} H_b$.
Note that
$\Image \Abb^{\rd}_{J_0-a\pi,\theta^u-2\pi}
\subset H_{<a}$.
Hence, we obtain
\[
\Image \Abb^{\rd}_{J_0-a\pi,\theta^u}\oplus H_{<a}
\oplus
\Abb_{\infty,\theta^u}
 \bigl(
 L'_{J_0-a\pi|\vartheta^{J_0}_{\ell}}
 \bigr)
\oplus H_{<a}
\]
by Lemma \ref{lem;24.4.5.1}.
Then, we obtain (\ref{eq;25.2.11.1}) by an easy induction.
\hfill\qed

\begin{cor}
We obtain
\[
 H_1^{\rd}\bigl(\cnum^{\ast},
 V\otimes\nbige(zu^{-1})\bigr)
 =\Abb^{\rd}_{\infty,\theta^u}
 \Bigl(
 L'_{J_{0-},0|\vartheta^{J_0}_{\ell}}
 \Bigr)
 \oplus
 \bigoplus_{J_0-2\omega^{-1}\pi<J\leq J_0}
 \Image \Abb^{\rd}_{J,\theta^u}
 \oplus H.
\] 
\hfill\qed
\end{cor}

\subsubsection{Proof of Proposition \ref{prop;18.4.19.100}}
\label{subsection;25.4.6.1}

Let $\gbigw_{2}'(\nbigi)$ be the set of
$J\in T(\nbigi)$
such that
$J_0-2\omega^{-1}\pi<J\leq J_0+(1-\omega^{-1})\pi$.
Note that
\[
 \gbigw_2(\nbigi)\sqcup
 \gbigw_2'(\nbigi)
 =\gbigt(\nbigi)
 \sqcup
 \{J\in \nbigt(\nbigi)\,|\,J_0-2\omega^{-1}\pi<J\leq J_0\}.
\]
Let $J'\in\gbigw_1(\nbigi)$.
Note that
$\nbigrtilde^{J'_-}_{J'-\omega^{-1}\pi}$ is an isomorphism.
By the construction in (\ref{eq;24.3.16.1}),
we have
\[
 \Image B_{J'_-,\theta^u}
 \oplus
 \bigoplus_{J'-\omega^{-1}\pi<J\leq J_0}
 \Image \Abb^{\rd}_{J,\theta^u}
 \oplus H
=\bigoplus_{J'-\omega^{-1}\pi\leq J\leq J_0}
 \Image \Abb^{\rd}_{J,\theta^u}
 \oplus H.
\]

Therefore, we can obtain
\begin{multline}
\label{eq;18.4.21.50}
 \bigoplus_{J\in \gbigw_1(\nbigi)}
\Image B_{J_-,\theta^u}
 \oplus
 \bigoplus_{J\in\gbigw_2(\nbigi)}
\Image \Abb^{\rd}_{J,\theta^u}
= \\
\Abb^{\rd}_{\infty,\theta^u}
\Bigl(
\bigoplus_{J_0-2\omega^{-1}\pi<J'\leq J_0}
 L_{J',<0|\vartheta_{\ell}^{J_0}}
\Bigr)
\oplus
 \bigoplus_{J\in\gbigt(\nbigi)}
\Image \Abb^{\rd}_{J,\theta^u}.
\end{multline}

Let $K$ denote the vector space (\ref{eq;18.4.21.50}).
For any $v\in L_{J_{0-},0}$,
$\Abb^{\rd}_{\infty,\theta^u}(v)$
equals
$B^{(J_0+2\pi)_{-}}_{\infty,\theta^u}(v')$
modulo $K$,
where $v'$ is the section of $L_{(J_0+2\pi)_-,0}$
induced by $v$ and the parallel transport of
$\nbigt_{\omega}(L)=\Gr^{\vecnbigf}_0(L)$.
(See Lemma \ref{lem;25.2.10.10}.)
Hence, we can conclude that
the morphism {\rm(\ref{eq;18.4.19.40})}
is an isomorphism
in the case $J_0+2\pi$.

Take any $J_1\in T(\nbigi)$
such that $(J_1)_-\cap(\vecI+\pi)_+\neq\emptyset$.
Take any $v\in L_{(J_1)-,0}$.
We obtain the section
$v'\in L_{(J_0+2\pi)_-,0}$
induced by $v$
and the parallel transport of
$\Gr^{\vecnbigf}_0(L)$.
Because
$B^{(J_1)_-}_{\infty,\theta^u}(v)-B^{(J_0+2\pi)_-}_{\infty,\theta^u}(v')$
is contained in $K$
(see Corollary \ref{cor;25.2.10.11}),
we obtain that
(\ref{eq;18.4.19.40})
is an isomorphism for any $J_1$ as above.
Thus the proof of Proposition \ref{prop;18.4.19.100}
is completed.
\hfill\qed

\vspace{.1in}
From the proof and (\ref{eq;25.2.11.2}),
we obtain the following corollary.
\begin{cor}
\label{cor;25.2.15.100}
\[
  \bigoplus_{J\in \gbigw_1(\nbigi)}
\Image B_{J_-,\theta^u}
 \oplus
 \bigoplus_{J\in\gbigw_2(\nbigi)}
\Image \Abb^{\rd}_{J,\theta^u}
= 
 \!\!\!\!\!
 \bigoplus_{J_0-2\omega^{-1}\pi<J'\leq J_0}
 \!\!\!\!\!
\Image \Abb^{\rd}_{J,\theta^u}
\oplus
 \bigoplus_{J\in\gbigt(\nbigi)}
\Image \Abb^{\rd}_{J,\theta^u}.
\]
\hfill\qed
\end{cor}

\subsection{Appendix: Another description of
$\Abb^{\rd}_{J,\theta^u}$}

Let us give another but equivalent description of 
the map $\Abb^{\rd}_{J,\theta^u}$.
Take $\epsilon>0$.
Set $I_1^{\circ}:=\openopen{-\epsilon}{\epsilon}$
and $I_2:=\closedclosed{0}{1}$.
We take an embedding
$F_J:
 I_1^{\circ}\times I_2
\lrarr X$
such that
(i) $F_{J|I_1^{\circ}\times\{0\}}
 \subset \{0\}\times J$,
(ii) $F_{J|I_1^{\circ}\times\{1\}}
\subset
 \{\infty\}\times(\vecI(\theta^u)+\pi)$,
(iii)
$F_{J|I_1^{\circ}\times\openopen{0}{1}}
\subset X^{\ast}$.
We have the local subsystem 
$\nbigl_{J,<0}\subset
 F_J^{-1}\varphi^{\ast}(\nbigl)$
on $I_1^{\circ}\times I_2$ 
induced by $L_{J,<0}$.
We obtain the constructible subsheaf
\[
 F_{J!}\nbigl_{J,<0}
\subset
 \varphi^{\ast}\bigl(
 \nbigl^{<0}\bigl(
 V\otimes\nbige(zu^{-1})
 \bigr)\bigr).
\]
Let $j_{I_1^{\circ}}$ denote the embedding of
$I_1^{\circ}$ to $\real$
obtained as the restriction
$F_{J|I_1^{\circ}\times\{0\}}$.
There exist the following isomorphisms:
\[
 H^1\bigl(X,
 F_{J!}\nbigl_{J,<0}\bigr)
\simeq
 H^1\bigl(
 \real,
 j_{I_1^{\circ}!}
 j_{I_1^{\circ}}^{-1}L_{J,<0}
 \bigr)
\simeq
 H_0\bigl(I_1^{\circ},
 j_{I_1^{\circ}}^{-1}L_{J,<0}
 \bigr)
\simeq
 H_0(J,L_{J,<0}).
\]
Hence,
we obtain
$H_0(J,L_{J,<0})
\lrarr
 H_1^{\rd}\bigl(
 \cnum^{\ast},V\otimes\nbige(zu^{-1})
 \bigr)$
 induced by
$\varphi_!
 F_{J!}\nbigl_{J,<0}
\subset
\nbigl^{<0}\bigl(V\otimes\nbige(zu^{-1})
 \bigr)$.
It equals $\Abb^{\rd}_{J,\theta^u}$
up to the signature.

\section{Rapid decay homology group of
$\nbigv\otimes\nbige(zu^{-1})$}

\subsection{Lifting maps}
\label{subsection;24.2.20.30}

From $(\nbigv,\nabla)$,
we obtain the local system with Stokes structure
$(L_{S^1},\vecnbigf^{\nbigv})$ on $\varpi^{-1}(0)$.
We obtain the constructible subsheaf
$L_{S^1}^{\nbigv,<0}$ of $L_{S^1}$.

For $J\in T(\nbigi)$,
let us construct maps
\index{maps $C^{J_{\pm}}_{\infty,\theta^u}$}
\[
 C^{J_{\pm}}_{\infty,\theta^u}:
 H_1^{\rd}\bigl(
 \cnum^{\ast},
 \nbigt_{\omega}(\nbigv)
 \otimes\nbige(zu^{-1})
 \bigr)
 \lrarr
 H_1^{\rd}\bigl(
 \cnum^{\ast},
 \nbigv\otimes\nbige(zu^{-1})
 \bigr).
\]

Let $I_1=\openopen{0}{1}$,
$I_2=\closedclosed{0}{1}$
and $I_2^{\circ}=I_2\setminus\{0,1\}$.
Let $F:I_1\times I_2\to X$ be an embedding
such that
(i) $F(I_1\times\{0\})\subset \{2\}\times J$,
(ii) $F(I_1\times I_2^{\circ})\subset
\{2<r<\infty\}\times \real$,
(iii) $F(I_1\times \{1\})\subset \{\infty\}\times(\vecI(\theta^u)+\pi)$.
We consider the following subsets of
$\projtilde^1=\real_{\geq 0}\times S^1$:
\[
Z_0:=\openopen{0}{2}\times S^1,
\quad
Z_1:=\closedopen{0}{2}\times S^1,
\quad
Z=Z_1\cup\Image (\varphi\circ F).
\]
Let $q_i:Z_i\to \varpi^{-1}(0)$ denote the projection.
Let $\nbign_{J_{\kappa},!}(\nbigv)$ $(\kappa=\pm)$
be the constructible subsheaves of
$\nbigl^{<0}\bigl(\nbigv\otimes\nbige(zu^{-1})\bigr)_{|Z}$
determined by the following conditions.
\begin{itemize}
 \item $\nbign_{J_{\kappa},!}(\nbigv)_{|\varpi^{-1}(0)}=L_{S^1}^{\nbigv,<0}$.
 \item $\nbign_{J_{\kappa},!}(\nbigv)_{|Z_0}
       =q_0^{-1}L_{S^1}^{\leq 0}$.
 \item $\nbign_{J_{\kappa},!}(\nbigv)_{|\Image (\varphi\circ F)}
       =\varphi_{\ast}(q_{\Image(F)}^{-1}(L'_{J_{\kappa},0}))$.
       Here,
       $L'_{J_{\kappa},0}$ denote the local subsystems
       of $L$
       determined by
       $L'_{J_{\kappa},0|J_{\kappa}}=L_{J_{\kappa},0}$.
\end{itemize}
Let $j_Z:Z\to \projtilde^1$ and $j_{Z_i}:Z_i\to \projtilde^1$
denote the inclusions.
We obtain the following exact sequence:
\[
 0\lrarr
  j_{Z_1!}(q_{Z_1}^{-1}L_{S^1}^{<0})
  \lrarr
  j_{Z!}\nbign_{J_{\pm},!}(\nbigv)
  \lrarr
  j_{Z!}\Bigl(
  \nbigl^{<0}\bigl(
  \nbigt_{\omega}(\nbigv)
  \otimes\nbige(zu^{-1})
  \bigr)_{|Z}
  \Bigr)
  \lrarr 0.
\]
The constructible sheaf
$j_{Z_1!}(q_{Z_1}^{-1}L_{S^1}^{<0})$
is acyclic with respect to the global cohomology.
The quotient of the natural monomorphism
\[
   j_{Z!}\Bigl(
  \nbigl^{<0}\bigl(
  \nbigt_{\omega}(\nbigv)
  \otimes\nbige(zu^{-1})
  \bigr)_{|Z}
 \Bigr)
 \lrarr
   \nbigl^{<0}\bigl(
  \nbigt_{\omega}(\nbigv)
  \otimes\nbige(zu^{-1})
  \bigr)
\]
is acyclic with respect to the global cohomology.
As a result,
there exists the natural isomorphism
\begin{equation}
 H^1(\projtilde^1,j_{Z!}\nbign_{J_{\pm},!}(\nbigv))
 \simeq
 H_1^{\rd}\bigl(
 \cnum^{\ast},
 \nbigt_{\omega}(\nbigv)\otimes\nbige(zu^{-1})
\bigr).
\end{equation}

We obtain the maps
$C^{J_{\pm}}_{\infty,\theta^u}$
from the natural morphisms
\[
 j_{Z!}\nbign_{J_{\pm},!}(\nbigv)\to
 \nbigl^{<0}\bigl(\nbigv\otimes\nbige(zu^{-1})\bigr).
\]
\begin{lem}
We have
$C^{(J+2\pi)_{\pm}}_{\infty,\theta^u+2\pi}
=C^{J_{\pm}}_{\infty,\theta^u}$.
\hfill\qed
\end{lem}

See Proposition \ref{prop;24.3.21.20}
for the difference
$C^{J_+}_{\infty,\theta^u}-C^{J_-}_{\infty,\theta^u}$.

\subsection{Basic properties}

\begin{lem}
\label{lem;25.2.15.10}
If $\nbigv=V$,
the maps 
$C^{J_{\pm}}_{\infty,\theta^u}$
equal
$B^{J_{\pm}}_{\infty,\theta^u}$.
\end{lem}
\pf
Let us study the case of $J_-$.
The other case can be argued similarly.
Let us give a description of
$B^{J_{-}}_{\infty,\theta^u}$
in terms of $1$-cycles.
Let $\gamma_1$ be a path on $(X^{\ast},X)$
connecting
$(1,\vartheta^J_{\ell})$
and $(\infty,\theta^u)$.
Let $\gamma_2$ be a path on $(X^{\ast},X)$
connecting
$(\infty,\theta_{\infty}-2\pi)$
and 
$(1,\vartheta^J_{\ell}-2\pi)$.

Let $b_0=\vartheta^J_{\ell}>b_1>\cdots>b_N=\vartheta^J_{\ell}-2\pi$
be the intersection of
$S_0(\nbigi)$ and 
$\closedclosed{\vartheta^J_{\ell}-2\pi}{\vartheta^J_{\ell}}$.
Set $J_i:=\openopen{b_i}{b_i+\omega^{-1}\pi}$ $(i=0,\ldots,N)$.
Take paths $I_i$ $(i=1,\ldots,N)$ on $X^{\ast}$
connecting
$(1,b_{i})$ and $(1,b_{i-1})$.
Let $K_i$ $(i=0,\ldots,N)$ be paths on $(X^{\ast},X)$
connecting
$(1,b_i)$ and 
a point in $\{0\}\times J_{i+1}$.

Because $\nbigt_{\omega}(V)$ is regular singular at $\{0,\infty\}$,
there exists the natural isomorphism (\ref{eq;24.2.21.1}).
For $v\in H^0(\real,\nbigt_{\omega}(L))$,
we obtain the corresponding elements
$v_i\in H^0(J_{i-},L_{J_{i-},0})\simeq H^0(\real,\nbigt_{\omega}(L))$
$(i=0,\ldots,N)$.
Note that
$v_{i+1}-v_i\in H^0(J_{i+1},L_{J_{i+1},<0})$.
We obtain the following $1$-cycle
for $V\otimes\nbige(zu^{-1})$:
\begin{equation}
\label{eq;24.2.20.3}
\varphi_{\ast}\Bigl(
 v_0\otimes\gamma_1
+\sum_{i=1}^{N}
 v_i\otimes I_i
+\sum_{i=0}^{N-1}
 (v_{i+1}-v_{i})\otimes K_{i}
 +v_N\otimes \gamma_2
 \Bigr).
\end{equation}
The homology class equals
$B^{J_{-}}_{\infty,\theta^u}(v)$.
Because (\ref{eq;24.2.20.3})
is a $1$-cocycle of
$j_{Z!}\nbign_{J_-,!}(V)\otimes
\nbigc^{\bullet}_{\projtilde^1,\del\projtilde^1}[-2]$,
we obtain the claim of the lemma.
\hfill\qed

\vspace{.1in}

There exist the following commutative diagrams:
\begin{equation}
\label{eq;24.2.20.10}
 \begin{CD}
  j_{Z!}\nbign_{J_{\pm},!}(V)
  @>{a_0}>>
  \nbigl^{<0}(V\otimes\nbige(zu^{-1}))
  \\
  @V{a_1}VV @V{a_2}VV \\
  j_{Z!}\nbign_{J_{\pm},!}(\nbigv)
  @>{a_3}>>
  \nbigl^{<0}(\nbigv\otimes\nbige(zu^{-1})).  
 \end{CD}
\end{equation}
From (\ref{eq;24.2.20.10}),
we obtain the following commutative diagrams:
\begin{equation}
 \begin{CD}
  H_1^{\rd}(\cnum^{\ast},\nbigt_{\omega}(V)\otimes\nbige(zu^{-1}))
  @>{B^{J_{\pm}}_{\infty,\theta^u}}>>
  H_1^{\rd}(\cnum^{\ast},V\otimes\nbige(zu^{-1}))
  \\
  @V{b_1}VV @V{b_2}VV \\
  H_1^{\rd}(\cnum^{\ast},\nbigt_{\omega}(\nbigv)\otimes\nbige(zu^{-1}))
  @>{C^{J_{\pm}}_{\infty,\theta^u}}>>
  H_1^{\rd}(\cnum^{\ast},\nbigv\otimes\nbige(zu^{-1})).
 \end{CD}
\end{equation}

\begin{prop}
\label{prop;24.2.20.22}
We obtain the following exact sequences:
\begin{multline}
 \label{eq;24.2.20.12}
 0\lrarr H_1^{\rd}(\cnum^{\ast},\nbigt_{\omega}(V)\otimes\nbige(zu^{-1}))
 \stackrel{B^{J_{\pm}}_{\infty,\theta^u}+b_1}{\lrarr}
 \\
 H_1^{\rd}(\cnum^{\ast},V\otimes\nbige(zu^{-1}))
\oplus 
  H_1^{\rd}(\cnum^{\ast},\nbigt_{\omega}(\nbigv)\otimes\nbige(zu^{-1}))
\\
 \stackrel{C^{J_{\pm}}_{\infty,\theta^u}-b_2}{\lrarr}
  H_1^{\rd}(\cnum^{\ast},\nbigv\otimes\nbige(zu^{-1}))\lrarr 0.
\end{multline}
\end{prop}
\pf
We obtain the following exact sequences from (\ref{eq;24.2.20.10}):
\begin{multline}
\label{eq;24.2.20.11}
 0\lrarr j_{Z!}\nbign_{J_{\pm},!}(V)
  \stackrel{a_1+a_0}{\lrarr}
 j_{Z!}\nbign_{J_{\pm},!}(\nbigv)
\oplus\nbigl^{<0}(V\otimes\nbige(zu^{-1})) \\
 \stackrel{a_3-a_2}{\lrarr}
 \nbigl^{<0}(\nbigv\otimes\nbige(zu^{-1}))
 \lrarr 0.
\end{multline}
For the constructible sheaves $\nbigf$ in  (\ref{eq;24.2.20.11}),
we have
$H^j(\projtilde^1,\nbigf)=0$ unless $j=1$.
Hence, we obtain the exact sequence (\ref{eq;24.2.20.12}).
\hfill\qed

\subsection{Decompositions of
$H_1^{\rd}(\cnum^{\ast},\nbigv\otimes\nbige(zu^{-1}))$}

Because $(V,\nabla)=\nbigs_{\omega}(\nbigv,\nabla)$,
there exists the natural morphism
\begin{equation}
\label{eq;24.2.20.20}
 H_1^{\rd}\bigl(\cnum^{\ast},V\otimes
  \nbige(zu^{-1})\bigr)
  \lrarr
  H_1^{\rd}\bigl(\cnum^{\ast},\nbigv\otimes
 \nbige(zu^{-1})\bigr).
\end{equation}

For $J\in T(\nbigi)$,
$\Abb^{\rd}_{J,\theta^u}$ and (\ref{eq;24.2.20.20})
induce
the following morphisms,
which are also denoted by $\Abb^{\rd}_{J,\theta^u}$:
\index{maps $\Abb^{\rd}_{J,\theta^u}$}
\[
 \Abb^{\rd}_{J,\theta^u}:
 H^0(J,L_{J,<0})
 \lrarr
 H_1^{\rd}\bigl(
 \cnum^{\ast},
 \nbigv\otimes\nbige(zu^{-1})
 \bigr).
\]
We also obtain the following maps
from 
$B_{J_{\pm},\theta^u}$ and (\ref{eq;24.2.20.20})
which are also denoted by $B_{J_{\pm},\theta^u}$:
\index{maps $B_{J_{\pm},\theta^u}$}
\[
 B_{J_{\pm},\theta^u}:
 H^0(J,L_{J,>0})
 \lrarr
 H_1^{\rd}\bigl(
 \cnum^{\ast},
 \nbigv\otimes\nbige(zu^{-1})
 \bigr).
\]

Let $\gbigw_j(\nbigi,\vecI(\theta^u)_{\pm})$ $(j=1,2)$
be as in \S\ref{subsection;24.2.20.13}.
Take $J_1\in\gbigw_2(\nbigi,\vecI(\theta^u)_+)$.
We obtain the following map
induced by
$B_{J_{-},\theta^u}$ $(J\in \gbigw_{1}(\nbigi,\vecI(\theta^u)_{+}))$,
$\Abb^{\rd}_{J,\theta^u}$ $(J\in\gbigw_{2}(\nbigi,\vecI(\theta^u)_{+}))$,
and $C^{J_{1,-}}_{\infty,\theta^u}$:
{\small
\begin{multline}
\label{eq;24.2.20.21}
 \bigoplus_{J\in\gbigw_{1}(\nbigi,\vecI(\theta^u)_+)}
 \!\!\!\!\!\!\!H^0(J,L_{J,>0})
\oplus\!\!\!\!
 \bigoplus_{J\in\gbigw_{2}(\nbigi,\vecI(\theta^u)_+)}
 \!\!\!\!\!\!\!H^0(J,L_{J,<0})
\oplus
 H_1^{\rd}\bigl(
 \cnum^{\ast},
 \nbigt_{\omega}(\nbigv)\otimes
 \nbige(zu^{-1})
 \bigr) \\
\lrarr
 H_1^{\rd}\bigl(
 \cnum^{\ast},\nbigv\otimes\nbige(zu^{-1})\bigr).
\end{multline}
}
Similarly, 
taking $J_2\in \gbigw_2(\nbigi,\vecI(\theta^u)_-)$,
we obtain the following map
induced by
$B_{J_{+},\theta^u}$ $(J\in \gbigw_{1}(\nbigi,\vecI(\theta^u)_{-}))$,
$A_{J,\theta^u}$ $(J\in\gbigw_{2}(\nbigi,\vecI(\theta^u)_{-}))$,
and $C^{J_{2,+}}_{\infty,\theta^u}$:
{\small
 \begin{multline}
\label{eq;24.2.20.25}
 \bigoplus_{J\in\gbigw_{1}(\nbigi,\vecI(\theta^u)_-)}\!\!\!\!
 H^0(J,L_{J,>0})
\oplus
 \!\!\!\bigoplus_{J\in\gbigw_{2}(\nbigi,\vecI(\theta^u)_-)}\!\!\!\!
 H^0(J,L_{J,<0})
\oplus
 H_1^{\rd}\bigl(
 \cnum^{\ast},
 \nbigt_{\omega}(\nbigv)\otimes
 \nbige(zu^{-1})
 \bigr) \\
\lrarr
 H_1^{\rd}\bigl(
 \cnum^{\ast},\nbigv\otimes\nbige(zu^{-1})\bigr).
 \end{multline}
}
We obtain the following corollary from
Proposition \ref{prop;18.4.19.100} and
Proposition \ref{prop;24.2.20.22}.
\begin{cor}
\label{cor;24.3.14.20}
The morphisms {\rm(\ref{eq;24.2.20.21})}
and {\rm(\ref{eq;24.2.20.25})}
are isomorphisms.
\hfill\qed
\end{cor}

\section{Moderate growth homology groups}

\subsection{Exact sequence}
Let $\gbigt(\nbigi,\theta^u)$ be as in \S\ref{subsection;18.5.13.20}.
There exists the following isomorphism
obtained as the projection:
\[
 L_{S^1}/L_{S^1}^{\leq 0}
\simeq 
 \bigoplus_{J\in\gbigt(\nbigi,\theta^u)}
 \varphi_{1\ast}a_{J\ast}L_{J,>0}.
\]
Let $b:\varpi^{-1}(0)\lrarr\projtilde^1$ denote the inclusion.
We obtain the following exact sequence:
\begin{multline}
0\lrarr
 \nbigl^{\leq 0}\bigl(V\otimes\nbige(zu^{-1})\bigr)
\lrarr
 \nbigl^{\leq 0}\bigl(V^{\reg}\otimes\nbige(zu^{-1})\bigr)
 \lrarr
 \\
 b_{\ast}\Bigl(
 \bigoplus_{J\in\gbigt(\nbigi,\theta^u)}
 \varphi_{1\ast}a_{J\ast}L_{J,>0}
 \Bigr)
\lrarr 0.
\end{multline}
By Lemma \ref{lem;20.10.14.1},
we obtain the following isomorphism
\begin{multline}
 \hyperh^{-2}\bigl(
 \projtilde^1,
 \nbigc^{\bullet}_{\projtilde^1,\del\projtilde^1}
 \otimes
 \varphi_{1\ast}a_{J\ast}L_{J,>0}
 \bigr)
 \simeq
 \hyperh^{-1}\bigl(
 \varpi^{-1}(0),
 \nbigc^{\bullet}_{\varpi^{-1}(0)}\otimes
 \varphi_{1\ast}a_{J\ast}L_{J,>0}
 \bigr)
 \\
 \simeq
 H^0(J,L_{J,>0}).
\end{multline}
Here,
the orientation of $\varpi^{-1}(0)$
is the opposite to the natural orientation of
$\varpi^{-1}(0)\simeq S^1$,
i.e.,
the orientation obtained as the component of
the boundary of $\projtilde^1$.
We obtain the following exact sequence:
\begin{multline}
\label{eq;24.2.15.1}
 0\lrarr
 \bigoplus_{J\in\gbigt(\nbigi,\theta^u)}
 H^0(J,L_{J,>0})
\stackrel{c_{1,u}}{\lrarr}
 H_1^{\mg}(\cnum^{\ast},V\otimes\nbige(zu^{-1}))
\stackrel{c_{2,u}}{\lrarr} \\
 H_1^{\mg}(\cnum^{\ast},V^{\reg}\otimes\nbige(zu^{-1}))
\lrarr 0.
\end{multline}

\subsection{Moderate growth $1$-homology classes
$\BB^{\mg}_{J,u}(v)$}
\label{subsection;24.3.25.111}

To represent $c_{1,u}$ in terms of $1$-cycles,
for any $J\in T(\nbigi)$,
let us construct a map
depending only on $u$:
\index{maps $\BB^{\mg}_{J,u}$}
\[
 \BB^{\mg}_{J,u}:
 H^0(J,L_{J,>0})
 \lrarr
 H^{\mg}_1(\cnum^{\ast},V\otimes\nbige(zu^{-1})).
\]
Let $\delta>0$ be any sufficiently small number.
We take a path $\gamma_{1,J}$ connecting
$(0,\vartheta^J_{\ell}-\delta)$ to 
$(1,\vartheta^J_r)$ on $(X,X^{\ast})$.
We also take paths $\gamma_{2,J,\pm}$
connecting $(1,\vartheta^J_r)$
and $(0,\vartheta^J_r\pm\delta)$ on $(X,X^{\ast})$.
By using
$H^0(J,L_{J,>0})
\simeq
H^0(J_-,L_{J_-,>0})
\subset
H^0(\real,L)$,
any
$v\in 
H^0\bigl(J,L_{J,>0}\bigr)$
induces
a section 
 of $\varphi^{\ast}\nbigl$
 along $\gamma_{1,J}$,
 which is also denoted by $v_{J_-}$.
Note that there exists the natural isomorphism:
\begin{equation}
 \label{eq;20.10.14.2}
\varphi^{\ast}(\nbigl)_{|(1,\vartheta^J_r)}
=(L_{J_+,0})_{|\vartheta^J_r}
\oplus
(L_{J_+,<0})_{|\vartheta^J_r}
\oplus
(L_{(J+\omega^{-1}\pi),<0})_{|\vartheta^J_r}
\oplus
\bigoplus_{\vartheta^J_r\in J'}
\bigl(
L_{J',<0}
\oplus L_{J',>0}
\bigr)_{|\vartheta^J_r}.
\end{equation}
According to (\ref{eq;20.10.14.2}),
we obtain the decomposition
\[
 v_{J_-|\vartheta^J_r}=u_{J,0}+\sum_{J\leq J'\leq J+\omega^{-1}\pi} u_{J'},
\]
where 
$u_{J'}\in L_{J',<0|\vartheta^J_r}$ 
and $u_{J,0}\in L_{J_+,0|\vartheta^J_r}$.
They naturally induce sections of $\varphi^{\ast}(\nbigl)$.
We obtain the following moderate growth $1$-cycle
of $(V,\nabla)\otimes\nbige(zu^{-1})$:
\[
\varphi_{\ast}\Bigl(
 v_{J_-}\otimes\gamma_{1,J}
+(u_J+u_{J,0})\otimes\gamma_{2,J,-}
+u_{J+\omega^{-1}\pi}\otimes\gamma_{2,J,+}
+\!\!\!\!\!
\sum_{J<J'<J+\omega^{-1}\pi}
\!\!\!\!\!u_{J'}\otimes\gamma_{2,J,+}
\Bigr).
\]
Let $\BB^{\mg}_{J,u}(v)\in
H^{\mg}_1(\cnum^{\ast},V\otimes\nbige(zu^{-1}))$
denote the homology class.
\begin{rem}
\label{rem;25.2.11.21}
$\BB^{\mg}_{J,u}(v)$ depends on $u$,
but is independent of the choice of $\theta^u$.
\hfill\qed
\end{rem}

The following lemma is clear by the construction.
\begin{lem}
\label{lem;25.2.11.20}
For any $v\in H^0(J+2\pi,L_{J+2\pi,>0})$,
we have
$\BB^{\mg}_{J+2\pi,u}(v)
=\BB^{\mg}_{J,u}(\Tbb^{\ast}(v))$.
\hfill\qed 
\end{lem}

\begin{lem}
Let $J\in\gbigt(\nbigi,\theta^u)$.
For any
\[
 v\in H^0(J,L_{J,>0})\simeq
 H^0(\varpi^{-1}(0),\varphi_{1\ast}a_{J\ast}L_{J,>0}),
\]
we have $c_{1,u}(v)=\BB^{\mg}_{J,u}(v)$.
\end{lem}
\pf
We set $K_{\delta}(J):=(J-\delta)\cup (J+\delta)$
and $K_{1,\delta}(J):=K_{\delta}(J)\setminus J$.
Let $[K_{\delta}(J)]$ be
a relative $1$-cycle in $H_1(K_{\delta}(J),K_{1,\delta}(J))$
induced by the inclusion of $K_{\delta}(J)$.
There exists the natural isomorphism
\begin{multline}
 H^0\bigl(\varpi^{-1}(0),
 \varphi_{1\ast}a_{J\ast}L_{J,>0}
 \bigr)
\simeq
 \hyperh^{-1}\bigl(
 \varpi^{-1}(0),
 \varphi_{1\ast}a_{J\ast}L_{J,>0}
 \otimes\nbigc^{\bullet}_{\varpi^{-1}(0)}
 \bigr)
 \\
\simeq
 H_1\bigl(
 (K_{\delta}(J),K_{1,\delta}(J));
 L'_{J,>0|K_{\delta}(J)}
 \bigr).
\end{multline}
Here, the correspondence is given by
the multiplication of $-[K_{\delta}(J)]$.
We set
$\Ktilde_{\delta}(J)=\closedclosed{0}{\delta}\times K_{\delta}(J)$.
By taking a simplicial decomposition of $\Ktilde_{\delta}(J)$,
we obtain a $2$-chain $[\Ktilde_{\delta}(J)]$.
We obtain a section
$v\otimes [\Ktilde_{\delta}(J)]$
of
$\varphi_{1\ast}a_{J\ast}L_{J,>0}
\otimes
\nbigc^{\bullet}_{\projtilde^1,\del\projtilde^1}[-2]$
which induces
$-v\otimes[K_{\delta}(J)]$
in 
$H_1\bigl(
 (K_{\delta}(J),K_{1,\delta}(J));
 L'_{J,>0|K_{\delta}(J)}
 \bigr)$.
By a direct computation,
it is mapped to
$\BB_{J,u}^{\mg}(v)\in
H_1^{\mg}(\cnum^{\ast},V\otimes\nbige(zu^{-1}))$
via $c_{1,u}$.
Hence, we obtain $\BB_{J,u}^{\mg}(v)=c_{1,u}(v)$.
\hfill\qed

\vspace{.1in}

Let us explain another description
of $\BB^{\mg}_{J,u}(v)$.
Let $\gamma'_{1,J}$ be a path connecting
$(1,\vartheta^J_{\ell})$ to $(0,\vartheta^J_r+\delta)$.
We also take paths $\gamma'_{2,J\pm}$
connecting
$(0,\vartheta^J_{\ell}\pm\delta)$ and $(1,\vartheta^J_{\ell})$.
By using
$H^0(J,L_{>0})\simeq
H^0(J_+,L_{J_+,>0})$,
any $v\in H^0(J,L_{J,>0})$
induces a section
of $\varphi^{\ast}\nbigl$ along $\gamma'_{1,J}$,
which is denoted by $v_{J_+}$.
There exists the decomposition
\[
 v_{J_+|\vartheta^J_{\ell}}
 =u'_{J,0}+\sum_{J-\omega^{-1}\pi\leq J'\leq J}u'_{J'},
\]
such that
$u'_{J'}\in L_{J',<0|\vartheta^J_{\ell}}$
and
$u'_{J,0}\in L_{J_-,0|\vartheta^J_{\ell}}$.
They naturally induce sections of $\varphi^{\ast}(\nbigl)$.
We obtain the following moderate growth $1$-cycle
of $V\otimes\nbige(zu^{-1})$:
\[
\varphi_{\ast}\Bigl(
 v\otimes\gamma'_{1,J}
+(u'_J+u'_{J,0})\otimes\gamma_{2,J,+}
+u_{J-\omega^{-1}\pi}\otimes\gamma_{2,J,-}
+\sum_{J'}u_{J'}\otimes\gamma_{2,J,-}
\Bigr).
\]
It also represents $\BB^{\mg}_{J,u}(v)$.

\subsection{Description of
$H_1^{\mg}(\cnum^{\ast},V^{\reg}\otimes\nbige(zu^{-1}))$}
\label{subsection;24.4.19.2}

Let $\Gamma_{\infty,\theta^u}$ be a path on $(X,X^{\ast})$ connecting
a point in $\{0\}\times\real$
and $(\infty,\theta^u)$.
\index{path $\Gamma_{\infty,\theta^u}$}
Any $v\in H^0(\real,L)$
induces a section of $\varphi^{\ast}(\nbigl)$
along $\Gamma_{\infty,\theta^u}$,
which is also denoted by $v$.
We obtain the moderate growth $1$-cycle
$v\otimes\Gamma_{\infty,\theta^u}$.
It induces a homology class
in
$H_1^{\mg}\bigl(
 \cnum^{\ast},
 V^{\reg}\otimes\nbige(zu^{-1})
 \bigr)$.
\begin{lem}
By this correspondence,
we obtain the isomorphism
\[
 H^0(\real,L)\simeq
 H_1^{\mg}\bigl(
\cnum^{\ast},
V^{\reg}\otimes\nbige(zu^{-1})
\bigr). 
\]
We also obtain
$L_{|\theta^u}
 \simeq
H_1^{\mg}\bigl(
\cnum^{\ast},
V^{\reg}
\otimes\nbige(zu^{-1})
 \bigr)$
by
$H^0(\real,L)\simeq L_{|\theta^u}$.
\hfill\qed
\end{lem}

\subsection{Splittings of $c_{2,u}$ in (\ref{eq;24.2.15.1})}
\label{subsection;24.3.25.112}

For each $J\in T(\nbigi)$,
let us construct morphisms
\index{maps $\Abb^{\mg, J_{\pm}}_{\infty,\theta^u}$}
\[
 \Abb^{\mg, J_{\pm}}_{\infty,\theta^u}:
 H_1^{\mg}\bigl(
 \cnum^{\ast},
 V^{\reg}
 \otimes\nbige(zu^{-1})
 \bigr)
=L_{|\theta^u} 
 \lrarr
 H_1^{\mg}\bigl(
 \cnum^{\ast},
 V\otimes\nbige(zu^{-1})
 \bigr)
\]
such that
$c_{2,u}\circ \Abb^{\mg,J_{\pm}}_{\infty,\theta^u}$
are isomorphisms.

For any $J'\in T(\nbigi)$,
let $\Gamma_{J'}$ be a path on $(X,X^{\ast})$
connecting a point in $\{0\}\times J'$
and $(\infty,\theta^u)$.
For $\kappa=\pm$,
let $L_{J_{\kappa},0}'\subset L$ be the local subsystem
determined by the condition
$L'_{J_{\kappa},0|J_{\kappa}}=L_{J_{\kappa},0}$.
Similarly,
let $L'_{J,<0}\subset L$ denote the local subsystem
determined by
$L'_{J,<0|J}=L_{J,<0}$.
Recall that there exist the following decompositions:
\[
 L=
 L'_{J_+,0}\oplus
 \bigoplus_{J-\pi^{-1}\omega< J'\leq J+\omega^{-1}\pi}
 L'_{J',<0}
 =L'_{J_-,0}
 \oplus
  \bigoplus_{J-\pi^{-1}\omega\leq J'<J+\omega^{-1}\pi}
 L'_{J',<0}.
\]
For $v\in H^0(\real,L)$,
there exists the decomposition
\begin{equation}
\label{eq;24.3.25.10}
 v=u_{J,0}+\sum_{J-\pi^{-1}\omega<J'\leq J+\omega^{-1}\pi}u_{J'},
\end{equation}
where $u_{J,0}\in H^0(\real,L'_{J_+,0})$
and $u_{J'}\in H^0(\real,L'_{J',<0})$.
We obtain the following moderate growth cycle of
$V\otimes\nbige(zu^{-1})$:
\[
 \varphi_{\ast}\Bigl(
 u_{J,0}\otimes\Gamma_J
 +\sum_{J-\pi^{-1}\omega<J'\leq J+\omega^{-1}\pi}
 u_{J'}\otimes\Gamma_{J'}
 \Bigr).
\]
It induces the homology class
$\Abb^{\mg,J_+}_{\infty,\theta^u}(v)
\in H_1^{\mg}(\cnum^{\ast},V\otimes\nbige(zu^{-1}))$.
Similarly,
there exists the decomposition
\begin{equation}
 v=w_{J,0}+\sum_{J-\pi^{-1}\omega\leq J'<J+\omega^{-1}\pi}w_{J'},
\end{equation}
where $w_{J,0}\in H^0(\real,L'_{J_-,0})$
and $w_{J'}\in H^0(\real,L'_{J',<0})$.
We obtain the following moderate growth cycle of
$V\otimes\nbige(zu^{-1})$:
\[
 \varphi_{\ast}\Bigl(
 w_{J,0}\otimes\Gamma_J
 +\sum_{J-\pi^{-1}\omega\leq J'< J+\omega^{-1}\pi}
 w_{J'}\otimes\Gamma_{J'}
 \Bigr).
\]
It induces a homology class
$\Abb^{\mg,J_-}_{\infty,\theta^u}(v)
\in H_1^{\mg}(\cnum^{\ast},V\otimes\nbige(zu^{-1}))$.
The following lemma is clear by construction.
\begin{lem}
$c_{2,u}\circ \Abb^{\mg,J_{\pm}}_{\infty,\theta^u}$
are the identity
on $H_1^{\mg}(\cnum^{\ast},V^{\reg}\otimes\nbige(zu^{-1}))$.
\hfill\qed
\end{lem}

We obtain the following lemma by the construction.
\begin{lem}
\label{lem;24.4.5.3}
We have
$\Abb^{\mg,(J+2\pi)_{\pm}}_{\infty,\theta^u+2\pi}
=\Abb^{\mg,J_{\pm}}_{\infty,\theta^u}\circ M$.
If $J_1\vdash J_2$,
then we have $\Abb^{\mg,J_{1+}}_{\infty,\theta^u}
=\Abb^{\mg,J_{2-}}_{\infty,\theta^u}$. 
\hfill\qed
\end{lem}

\begin{lem}
\label{lem;24.3.25.10}
We have $\Abb^{\mg,J_-}_{\infty,\theta^u}
-\Abb^{\mg,J_+}_{\infty,\theta^u}
=\BB^{\mg}_{J,u}\circ R_J$
on $H^0(\real,L)$.
(See {\rm\S\ref{subsection;24.2.18.1}} for the maps $R_J$.)
\end{lem}
\pf
Let $v\in H^0(\real,L)$.
It is enough to consider the cases
(i) $v$ is contained in the image
$H^0(J_-,L_{J_-,>0})\to H^0(\real,L)$,
(ii) $R_J(v)=0$.
In the case (i),
we obtain
$\BB^{\mg}_{J,u}(R_J(v))
=\Abb^{\mg,J_-}_{\infty,\theta^u}(v)
-\Abb^{\mg,J_+}_{\infty,\theta^u}(v)$
by the construction.
In the case (ii),
we obtain $u_{J+\omega^{-1}\pi}=0$,
$w_{J-\omega^{-1}\pi}=0$,
and 
$u_{J'}=w_{J'}$ for
$J-\omega^{-1}\pi<J'<J+\omega^{-1}\pi$ with $J'\neq J$.
We also have
$u_{J}+u_{J,0}=w_{J}+w_{J,0}$.
Hence, we obtain
$\Abb^{\mg,J_-}_{\infty,\theta^u}(v)
=\Abb^{\mg,J_+}_{\infty,\theta^u}$
and 
$\BB^{\mg}_{J,u}(R_J(v))=0$.

Let us give another proof.
We identify $\projtilde^1$ with
$\realbar_{\geq 0}\times S^1$.
We set $Z=\closedopen{0}{1}\times S^1$.
Let $q_{Z}$ denote the projection of $Z$ to $S^1$.
Let $\iota_Z:Z\to \projtilde^1$ denote the inclusion.
We obtain the constructible subsheaf
$\iota_{Z!}\bigl(
q_Z^{-1}(L^{\leq 0})_{S^1}
\bigr)$
which is acyclic with respect to the global cohomology.
The homology class
$\Abb^{\mg,J_-}_{\infty,\theta^u}(v)
-\Abb^{\mg,J_+}_{\infty,\theta^u}
-\BB^{\mg}_{J,u}(R_J(v))$
is induced by a $1$-cocycle of
$\iota_{Z!}\bigl(
q_Z^{-1}(L^{\leq 0})_{S^1}
\bigr)
\otimes\nbigc^{\bullet}_{\projtilde^1,\del \projtilde^1}[2]$.
Hence, it is $0$.
\hfill\qed

\begin{cor}
\label{cor;24.4.5.4}
If $J_1\leq J_2$, the following holds on $H^0(\real,L)$:
\begin{equation}
\label{eq;24.2.18.11}
\Abb^{\mg,J_{1-}}_{\infty,\theta^u}
-\Abb^{\mg,J_{2-}}_{\infty,\theta^u}
=\sum_{J_1\leq J<J_2}\BB^{\mg}_{J,u}\circ R_J.
\end{equation}
\hfill\qed
\end{cor}

\subsection{Relations with rapid decay homology classes}
\label{subsection;24.3.26.3}

There exists the natural morphism
\[
 H_1^{\rd}(\cnum^{\ast},V\otimes\nbige(zu^{-1}))
 \lrarr
 H_1^{\mg}(\cnum^{\ast},V\otimes\nbige(zu^{-1})).
\]
The image of any element in
$H_1^{\rd}(\cnum^{\ast},V\otimes\nbige(zu^{-1}))$
is denoted by the same notation.
\begin{lem}
\label{lem;24.4.5.10}
For any $J_1\in T(\nbigi)$ and $v\in H^0(\real,L)$,
we obtain the following relation
in $H_1^{\mg}(\cnum^{\ast},V\otimes\nbige(zu^{-1}))$:
\begin{equation}
\label{eq;24.3.25.11}
 \Abb^{\rd}_{\infty,\theta^u}(v)
=\Abb^{\mg,J_{1+}}_{\infty,\theta^u}(v)
-\Abb^{\mg,J_{1+}}_{\infty,\theta^u}(M^{-1}(v))
 +\sum_{J_1-2\pi<J\leq J_1}
 \BB^{\mg}_{J,u}(R_{J}(v)).
\end{equation}
For any $J\in T(\nbigi)$
and $v\in H^0(J,L_{J,<0})$, we obtain
\begin{equation}
\label{eq;24.3.25.12}
 \Abb^{\rd}_{J,\theta^u}(v)
=\Abb^{\mg, J_+}_{\infty,\theta^u}(\rho_{J}(v))
=\Abb^{\mg,J_{-}}_{\infty,\theta^u}(\rho_{J}(v)),
\end{equation}
where
$\rho_{J}:H^0(J,L_{J,<0})\to H^0(\real,L)$
denotes the natural inclusion.
\end{lem}
\pf
The equalities (\ref{eq;24.3.25.12}) clear by the constructions.
We can obtain (\ref{eq;24.3.25.11})
by the second proof of Lemma \ref{lem;24.3.25.10}.
\hfill\qed

\begin{lem}
\label{lem;25.2.15.20}
For $J\in T(\nbigi)$ and $v\in H^0(\real,\nbigt_{\omega}(L))$,
we have
\[
 B^{J_{\pm}}_{\infty,\theta^u}(v)
=\Abb^{J_{\pm}}_{\infty,\theta^u}(v)
-\Abb^{J_{\pm}}_{\infty,\theta^u}(M^{-1}(v)).
\]
\end{lem}
\pf
By the description of $B^{J_-}_{\infty,\theta^u}(v)$
in the proof of Lemma \ref{lem;25.2.15.10},
we obtain
\[
  B^{J_-}_{\infty,\theta^u}(v)
=\Abb^{J_-}_{\infty,\theta^u}(v)
-\Abb^{J_--2\pi}_{\infty,\theta^u-2\pi}(v)
=\Abb^{J_-}_{\infty,\theta^u}(v)
-\Abb^{J_-}_{\infty,\theta^u}(M^{-1}(v)).
\]
We can obtain the claim for $B^{J_+}_{\infty,\theta^u}(v)$
similarly.
\hfill\qed

\subsection{Lifting maps for the moderate homology groups}
\label{subsection;24.3.21.13}

For $J\in T(\nbigi)$,
we define
\index{maps $B^{\mg,J_{\pm}}_{\infty,\theta^u}$}
\[
 B^{\mg,J_{\pm}}_{\infty,\theta^u}:
H^0(\real,\nbigt_{\omega}(L))
 \lrarr
H_1^{\mg}(\cnum^{\ast},V\otimes\nbige(zu^{-1}))
\]
by setting
\[
 B^{\mg,J_{\pm}}_{\infty,\theta^u}(v)
=\Abb^{\mg,J_{\pm}}_{\infty,\theta^u}(v_{J_{\pm}})
\]
for any $v\in H^0(\real,\nbigt_{\omega}(L))$,
where $v_{J_{\pm}}$ are obtained by
\[
H^0(\real,\nbigt_{\omega}(L))
 \simeq
 H^0(J_{\pm},L_{J_{\pm},0})
 \subset
 H^0(\real,L). 
\]
We may also regard $B^{\mg,J_{\pm}}_{\infty,\theta^u}$
as maps
\[
  B^{\mg,J_{\pm}}_{\infty,\theta^u}:
H_1^{\mg}(\cnum^{\ast},\nbigt_{\omega}(V)
 \otimes\nbige(zu^{-1}))
 \lrarr
H_1^{\mg}(\cnum^{\ast},V\otimes\nbige(zu^{-1})).
\]
We obtain the following lemmas by the construction.
Note that the isomorphism
$H^0(\real,\nbigt_{\omega}(L))
\simeq
 H_1^{\mg}\bigl(\cnum^{\ast},
 \nbigt_{\omega}(V)\otimes\nbige(zu^{-1})\bigr)$
depends on the choice of $\theta^u$.
\begin{lem}
We have
$B^{\mg,(J+2\pi)_{\pm}}_{\infty,\theta^u+2\pi}
=B^{\mg,J_{\pm}}_{\infty,\theta^u}\circ M_0$
 as maps on
$H^0(\real,\nbigt_{\omega}(L))$. 
We have
$B^{\mg,(J+2\pi)_{\pm}}_{\infty,\theta^u+2\pi}
=B^{\mg,J_{\pm}}_{\infty,\theta^u}$
 as maps on
$H_1^{\mg}(\cnum^{\ast},\nbigt_{\omega}(V)
 \otimes\nbige(zu^{-1}))$.
\hfill\qed
\end{lem}

The following lemma is clear by the construction.
\begin{lem}
\label{lem;24.4.2.100}
$B^{\mg,J_+}_{\infty,\theta^u}
-B^{\mg,J_-}_{\infty,\theta^u}
=\Abb^{\rd}_{J,\theta^u}\circ\nbigp_J$
as maps
$H^0(\real,\nbigt_{\omega}(L))
 \lrarr 
 H_1^{\mg}\bigl(\cnum^{\ast},V\otimes\nbige(zu^{-1})\bigr)$.
\hfill\qed
\end{lem}

We obtain the following lemma from
Lemma \ref{lem;25.2.15.20}
\begin{lem}
\label{lem;25.2.15.21}
For $v\in H^0(\real,\nbigt_{\omega}(L))$
and for $J\in T(\nbigi)$,
we have
\[
 B^{J_{\pm}}_{\infty,\theta^u}(v)
=B^{\mg,J_{\pm}}_{\infty,\theta^u}(v)
-B^{\mg,J_{\pm}}_{\infty,\theta^u}(M^{-1}(v)).
\]
\hfill\qed
\end{lem}

For $J\in T(\nbigi)$,
we shall construct maps
\index{maps $C^{\mg,J_{\pm}}_{\infty,\theta^u}$}
\[
 C^{\mg,J_{\pm}}_{\infty,\theta^u}:
 H_1^{\mg}\bigl(
 \cnum^{\ast},
 \nbigt_{\omega}(\nbigv)
 \otimes\nbige(zu^{-1})
 \bigr)
 \lrarr
 H_1^{\mg}\bigl(
 \cnum^{\ast},\nbigv \otimes\nbige(zu^{-1})
 \bigr).
\]
We use the notation in \S\ref{subsection;24.2.20.30}.
Let $\nbign_{J_{\pm},\ast}$ denote the constructible subsheaf
of $\nbigl^{\leq 0}(\nbigv\otimes\nbige(zu^{-1}))_{|Z}$
determined by the following conditions:
\begin{itemize}
 \item $\nbign_{J_{\pm},\ast}(\nbigv)_{|\varpi^{-1}(0)}
       =L_{S^1}^{\nbigv,\leq 0}$.
 \item $\nbign_{J_{\pm},\ast}(\nbigv)_{|Z_0}
       =q_0^{-1}L_{S^1}^{\leq 0}$.
 \item $\nbign_{J_{\pm},\ast}(\nbigv)_{|\Image (\varphi\circ F)}
       =\varphi_{\ast}(q_{\Image(F)}^{-1}(L'_{J_{\pm},0}))$.
\end{itemize}
We obtain the following exact sequence:
\[
 0\lrarr
  j_{Z_1!}(q_{Z_1}^{-1}L_{S^1}^{<0})
  \lrarr
  j_{Z!}\nbign_{J_{\pm},\ast}(\nbigv)
  \lrarr
  j_{Z!}\Bigl(
  \nbigl^{\leq 0}\bigl(
  \nbigt_{\omega}(\nbigv)
  \otimes\nbige(zu^{-1})
  \bigr)_{|Z}
  \Bigr)
  \lrarr 0.
\]
The constructible sheaf
$j_{Z_1!}(q_{Z_1}^{-1}L_{S^1}^{<0})$
is acyclic with respect to the global cohomology.
The cokernel of the natural monomorphism
\[
  j_{Z!}\Bigl(
  \nbigl^{\leq 0}\bigl(
  \nbigt_{\omega}(\nbigv)
  \otimes\nbige(zu^{-1})
  \bigr)_{|Z}
 \Bigr)
 \lrarr
   \nbigl^{\leq 0}\bigl(
  \nbigt_{\omega}(\nbigv)
  \otimes\nbige(zu^{-1})
  \bigr)
\]
is acyclic with respect to the global cohomology.
As a result,
there exists the natural isomorphism
\begin{equation}
 H^1(\projtilde^1,j_{Z!}\nbign_{J_{\pm},\ast}(\nbigv))
 \simeq
 H_1^{\mg}\bigl(
 \cnum^{\ast},
 \nbigt_{\omega}(\nbigv)\otimes\nbige(zu^{-1})
\bigr).
\end{equation}
The maps are induced by
$C^{\mg,J_{\pm}}_{\infty,\theta^u}$
the natural morphisms
\[
j_{Z!}\nbign_{J_{\pm},\ast}(\nbigv)\to
 \nbigl^{\leq 0}\bigl(\nbigv\otimes\nbige(zu^{-1})\bigr).
\]

\subsubsection{Basic properties}

The following lemma is clear by the construction.
\begin{lem}
If $\nbigv=V$,
then
$C^{\mg,J_{\pm}}_{\infty,\theta^u}=B^{\mg,J_{\pm}}_{\infty,\theta^u}$. 
\hfill\qed
\end{lem}

We obtain the following commutative diagrams:
\begin{equation}
\label{eq;24.2.20.40}
 \begin{CD}
  j_{Z!}\nbign_{J_{\pm},!}(V)
  @>{a_0}>>
  \nbigl^{<0}(V\otimes\nbige(zu^{-1}))
  \\
  @V{a_1}VV @V{a_2}VV \\
  j_{Z!}\nbign_{J_{\pm},\ast}(\nbigv)
  @>{a_3}>>
  \nbigl^{\leq 0}(\nbigv\otimes\nbige(zu^{-1})).  
 \end{CD}
\end{equation}

From (\ref{eq;24.2.20.40}),
we obtain the following commutative diagrams:
\begin{equation}
\label{eq;24.2.20.41}
 \begin{CD}
  H_1^{\rd}(\cnum^{\ast},\nbigt_{\omega}(V)\otimes\nbige(zu^{-1}))
  @>{B^{J_{\pm}}_{\infty,\theta^u}}>>
  H_1^{\rd}(\cnum^{\ast},V\otimes\nbige(zu^{-1}))
  \\
  @V{b_1}VV @V{b_2}VV \\
  H_1^{\mg}(\cnum^{\ast},\nbigt_{\omega}(\nbigv)\otimes\nbige(zu^{-1}))
  @>{C^{\mg,J_{\pm}}_{\infty,\theta^u}}>>
  H_1^{\mg}(\cnum^{\ast},\nbigv\otimes\nbige(zu^{-1})).
 \end{CD}
\end{equation}
The following proposition is similar to
Proposition \ref{prop;24.2.20.22}.
\begin{prop}
\label{prop;24.2.20.42}
We obtain the following exact sequence
\begin{multline}
 \label{eq;24.2.20.43}
 0\lrarr H_1^{\rd}(\cnum^{\ast},\nbigt_{\omega}(V)\otimes\nbige(zu^{-1}))
 \stackrel{B^{J_{\pm}}_{\infty,\theta^u}+b_1}{\lrarr}
 \\
 H_1^{\rd}(\cnum^{\ast},V\otimes\nbige(zu^{-1}))
\oplus 
  H_1^{\mg}(\cnum^{\ast},\nbigt_{\omega}(\nbigv)\otimes\nbige(zu^{-1}))
\\
 \stackrel{C^{\mg,J_{\pm}}_{\infty,\theta^u}-b_2}{\lrarr}
  H_1^{\mg}(\cnum^{\ast},\nbigv\otimes\nbige(zu^{-1}))\lrarr 0.
\end{multline}
\hfill\qed
\end{prop}

The commutative diagram (\ref{eq;24.2.20.41})
is a part of the following.
\begin{equation}
 \begin{CD}
  H_1^{\rd}(\cnum^{\ast},\nbigt_{\omega}(V)\otimes\nbige(zu^{-1}))
  @>{B^{J_{\pm}}_{\infty,\theta^u}}>>
  H_1^{\rd}(\cnum^{\ast},V\otimes\nbige(zu^{-1}))
  \\
  @VVV @VVV \\
  H_1^{\rd}(\cnum^{\ast},\nbigt_{\omega}(\nbigv)\otimes\nbige(zu^{-1}))
  @>{C^{J_{\pm}}_{\infty,\theta^u}}>>
  H_1^{\rd}(\cnum^{\ast},\nbigv\otimes\nbige(zu^{-1}))\\
  @VVV @VVV \\
  H_1^{\mg}(\cnum^{\ast},\nbigt_{\omega}(\nbigv)\otimes\nbige(zu^{-1}))
  @>{C^{\mg,J_{\pm}}_{\infty,\theta^u}}>>
  H_1^{\mg}(\cnum^{\ast},\nbigv\otimes\nbige(zu^{-1}))
  \\
  @VVV @VVV \\
  H_1^{\mg}(\cnum^{\ast},\nbigt_{\omega}(V)\otimes\nbige(zu^{-1}))
  @>{B^{\mg,J_{\pm}}_{\infty,\theta^u}}>>
  H_1^{\mg}(\cnum^{\ast},V\otimes\nbige(zu^{-1})).
 \end{CD}
\end{equation}

\subsubsection{Decompositions}

There exists the natural morphism
\begin{equation}
\label{eq;24.2.20.50}
 H_1^{\rd}\bigl(\cnum^{\ast},V\otimes
  \nbige(zu^{-1})\bigr)
  \lrarr
  H_1^{\mg}\bigl(\cnum^{\ast},\nbigv\otimes
 \nbige(zu^{-1})\bigr).
\end{equation}
For $J\in T(\nbigi)$,
$\Abb^{\rd}_{J,\theta^u}$ and (\ref{eq;24.2.20.50})
induce the following morphisms,
which are also denoted by $\Abb^{\rd}_{J,\theta^u}$:
\[
 \Abb^{\rd}_{J,\theta^u}:
 H^0(J,L_{J,<0})
 \lrarr
 H_1^{\mg}\bigl(
 \cnum^{\ast},
 \nbigv\otimes\nbige(zu^{-1})
 \bigr).
\]
We also obtain the following maps from 
$B_{J_{\pm},\theta^u}$ and (\ref{eq;24.2.20.20})
which are also denoted by $B_{J_{\pm},\theta^u}$:
\[
 B_{J_{\pm},\theta^u}:
 H^0(J,L_{J,>0})
 \lrarr
 H_1^{\mg}\bigl(
 \cnum^{\ast},
 \nbigv\otimes\nbige(zu^{-1})
 \bigr).
\]

Let $\gbigw_j(\nbigi,\vecI(\theta^u)_{\pm})$ $(j=1,2)$
be as in \S\ref{subsection;24.2.20.13}.
Take $J_1\in\gbigw_2(\nbigi,\vecI(\theta^u)_+)$.
We obtain the following map
induced by
$B_{J_{-},\theta^u}$ $(J\in \gbigw_{1}(\nbigi,\vecI(\theta^u)_{+}))$,
$\Abb^{\rd}_{J,\theta^u}$ $(J\in\gbigw_{2}(\nbigi,\vecI(\theta^u)_{+}))$,
and $C^{\mg,J_{1,-}}_{\infty,\theta^u}$:
{\small
\begin{multline}
\label{eq;24.2.20.51}
 \bigoplus_{J\in\gbigw_{1}(\nbigi,\vecI(\theta^u)_+)}
 \!\!\!\!H^0(J,L_{J,>0})
\oplus
 \bigoplus_{J\in\gbigw_{2}(\nbigi,\vecI(\theta^u)_+)}
 \!\!\!\!H^0(J,L_{J<0})
\oplus
 H_1^{\mg}\bigl(
 \cnum^{\ast},
 \nbigt_{\omega}(\nbigv)\otimes
 \nbige(zu^{-1})
 \bigr) \\
\lrarr
 H_1^{\mg}\bigl(
 \cnum^{\ast},\nbigv\otimes\nbige(zu^{-1})\bigr).
\end{multline}
}
Similarly, 
we take $J_2\in \gbigw_2(\nbigi,\vecI(\theta^u)_-)$.
We obtain the following map
induced by
$B_{J_{+},\theta^u}$
$(J\in \gbigw_{1}(\nbigi,\vecI(\theta^u)_{-}))$,
$\Abb^{\rd}_{J,\theta^u}$
$(J\in\gbigw_{2}(\nbigi,\vecI(\theta^u)_{-}))$,
and $C^{\mg,J_{2,+}}_{\infty,\theta^u}$:
{\small
 \begin{multline}
\label{eq;24.2.20.52}
 \bigoplus_{J\in\gbigw_{1}(\nbigi,\vecI(\theta^u)_-)}\!\!\!\!
 H^0(J,L_{J,>0})
\oplus
 \!\!\!\bigoplus_{J\in\gbigw_{2}(\nbigi,\vecI(\theta^u)_-)}\!\!\!\!
 H^0(J,L_{J,<0})
\oplus
 H_1^{\mg}\bigl(
 \cnum^{\ast},
 \nbigt_{\omega}(\nbigv)\otimes
 \nbige(zu^{-1})
 \bigr) \\
\lrarr
 H_1^{\mg}\bigl(
 \cnum^{\ast},\nbigv\otimes\nbige(zu^{-1})\bigr).
 \end{multline}}
We obtain the following corollary from
Proposition \ref{prop;18.4.19.100} and
Proposition \ref{prop;24.2.20.42}.
\begin{cor}
\label{cor;24.3.14.21}
The morphisms {\rm(\ref{eq;24.2.20.51})}
and {\rm(\ref{eq;24.2.20.52})}
are isomorphisms.
\hfill\qed
\end{cor}

\subsubsection{Difference of lifting maps}

There exist the natural morphisms
\begin{multline}
 H_1^{\rd}\bigl(\cnum^{\ast},
  \nbigt_{\omega}(\nbigv)
  \otimes\nbige(zu^{-1})
  \bigr)
  \stackrel{a_1}{\lrarr}
  H_1^{\mg}\bigl(\cnum^{\ast},
  \nbigt_{\omega}(\nbigv)
  \otimes\nbige(zu^{-1})
 \bigr)
 \\
 \stackrel{a_2}{\lrarr}
  H_1^{\mg}\bigl(\cnum^{\ast},
  \nbigt_{\omega}(V)
  \otimes\nbige(zu^{-1})
 \bigr)
 \simeq H^0(\real,\nbigt_{\omega}(L)).
\end{multline}
\begin{prop}
\label{prop;24.3.21.20}
The following equality holds
as maps to $H_1^{\rd}(\cnum^{\ast},\nbigv\otimes\nbige(zu^{-1}))$: 
\begin{equation}
\label{eq;24.3.21.10}
 C^{J_+}_{\infty,\theta^u}
-C^{J_-}_{\infty,\theta^u}
=\Abb^{\rd}_{J,\theta^u}\circ\nbigp_J
\circ(a_2\circ a_1),
\end{equation}
The following equality holds
as maps to $H_1^{\mg}(\cnum^{\ast},\nbigv\otimes\nbige(zu^{-1}))$: 
 \begin{equation}
\label{eq;24.3.21.11}
 C^{\mg,J_+}_{\infty,\theta^u}
-C^{\mg,J_-}_{\infty,\theta^u}
=\Abb^{\rd}_{J,\theta^u}\circ\nbigp_J\circ a_2.
 \end{equation}
\end{prop}
\pf
This is reduced to Lemma \ref{lem;24.4.2.100}
by Theorem \ref{thm;24.3.25.50}.
We explain a direct sheaf theoretic argument to the issue
to Lemma \ref{lem;24.4.2.100}.
We use the notation in \S\ref{subsection;24.2.20.30}
and \S\ref{subsection;24.3.21.13}.
For $\varrho=!,\ast$ and $\kappa=\pm$,
there exist the following epimorphisms:
\[
 b_{\kappa,\varrho}^{\nbigv}:
 \nbign_{J_{\kappa},\varrho}(\nbigv)
\lrarr
  \nbigl^{\varrho}\bigl(
  \nbigt_{\omega}(\nbigv)
  \otimes\nbige(zu^{-1})
  \bigr)_{|Z}.
\]
Let $\nbigk_{J,\varrho}(\nbigv)$ denote the kernel of
the following morphism:
\[
\begin{CD}
 \nbign_{J_{+},\varrho}(\nbigv)
 \oplus
 \nbign_{J_{-},\varrho}(\nbigv)
 @>{b^{\nbigv}_{+,\varrho}-b^{\nbigv}_{-,\varrho}}>>
   \nbigl^{\varrho}\bigl(
  \nbigt_{\omega}(\nbigv)
  \otimes\nbige(zu^{-1})
  \bigr)_{|Z}.
\end{CD}
\]
The composition of the morphisms
\[
 j_{Z!}\nbigk_{J,\varrho}(\nbigv)
 \lrarr
 j_{Z!}\nbign_{J_+,\varrho}(\nbigv)
 \lrarr
    \nbigl^{\varrho}\bigl(
  \nbigt_{\omega}(\nbigv)
  \otimes\nbige(zu^{-1})
  \bigr)
\]
induce the isomorphisms
\[
 H^1(\projtilde^1,j_{Z!}\nbigk_{J,\varrho}(\nbigv))
 \simeq
 H_1^{\varrho}
 \bigl(
 \cnum^{\ast},
 \nbigt_{\omega}(\nbigv)\otimes\nbige(zu^{-1})
 \bigr).
\]

There exist the natural monomorphisms
\[
 c_{\kappa,\varrho}^{\nbigv}:
 j_{Z!}\nbign_{J_{\kappa},\varrho}(\nbigv)
 \lrarr
 \nbigl^{\varrho}\bigl(
 \nbigv\otimes\nbige(zu^{-1})
 \bigr).
\]
The map
$C^{J_+}_{\infty,\theta^u}-C^{J_-}_{\infty,\theta^u}$
(resp. $C^{\mg,J_+}_{\infty,\theta^u}-C^{\mg,J_-}_{\infty,\theta^u}$)
is induced by the following morphism
in the case
$\varrho=!$ (resp. $\varrho=\ast$):
\[
 c_{+,\varrho}^{\nbigv}-c_{-,\varrho}^{\nbigv}:
 j_{Z!}\nbigk_{J,\varrho}(\nbigv)
 \lrarr
 \nbigl^{\varrho}\bigl(
 \nbigv\otimes\nbige(zu^{-1})
 \bigr).
\]
Let $\nbigntilde_{J,\varrho}(\nbigv)$ be
the constructible subsheaf of
$\nbigl^{\varrho}\bigl(
 \nbigv\otimes\nbige(zu^{-1})
 \bigr)_{|Z}$
determined by the following conditions.
\begin{itemize}
 \item $\nbigntilde_{J,\varrho}(\nbigv)_{|Z_1}
       =\nbign_{J_+,\varrho}(\nbigv)_{|Z_1}
       =\nbign_{J_-,\varrho}(\nbigv)_{|Z_1}$.
 \item $\nbigntilde_{J,\varrho}(\nbigv)_{|\Image(\varphi\circ F)}
       =\varphi_{\ast}\bigl(
        q_{\Image F}^{-1}(L'_{J,\leq 0})
       \bigr)$.
       Here, $L'_{J,\leq 0}$ denotes the local subsystem of $L$
       determined by
       $L'_{J,\leq 0|J=L_{J,\leq 0}}$.
\end{itemize}
Let $\nbigntilde'_{J}(\nbigv)$ be the constructible subsheaf
of $\nbigntilde_{J,\varrho}(\nbigv)$
determined by the following conditions.
\begin{itemize}
 \item $\nbigntilde'_{J}(\nbigv)_{|Z_1}
       =q_1^{-1}(L_{S^1}^{<0})$.
 \item $\nbigntilde'_{J}(\nbigv)_{|\Image(\varphi\circ F)}
       =\varphi_{\ast}\bigl(
        q_{\Image F}^{-1}(L'_{J,<0})
       \bigr)$.
       Here, $L'_{J,<0}$ denotes the local subsystem of $L$
       determined by
       $L'_{J,<0|J=L_{J,< 0}}$.
\end{itemize}
We have
$\nbigntilde'_J(\nbigv)
\subset
\nbigntilde_{J,!}(\nbigv)
\subset
\nbigntilde_{J,\ast}(\nbigv)$.
There exists the following exact sequences:
\[
 0\lrarr
 j_{Z!}\nbigntilde'_{J}(\nbigv)
 \lrarr
 j_{Z!}\nbigntilde_{J,\varrho}(\nbigv)
 \lrarr
   j_{Z!}\Bigl(
  \nbigl^{\varrho}\bigl(
  \nbigt_{\omega}(\nbigv)
  \otimes\nbige(zu^{-1})
  \bigr)_{|Z}
  \Bigr)
  \lrarr 0,
\]
There exist the natural monomorphisms
$d^{\nbigv}_{\pm,\varrho}:
\nbign_{J_{\pm},\varrho}(\nbigv)\lrarr  \nbigntilde_{J,\varrho}(\nbigv)$.
The morphism
$c^{\nbigv}_{+,\varrho}-c^{\nbigv}_{-,\varrho}$
is the composition of the following morphisms:
\[
 \begin{CD}
  j_{Z!}\nbigk_{J,\varrho}(\nbigv)
  @>{d^{\nbigv}_{+,\varrho}-d^{\nbigv}_{-,\varrho}}>>
  j_{Z!}\nbigntilde_{J,\varrho}(\nbigv)
  @>>>
  \nbigl^{\varrho}
  \bigl(\nbigv\otimes\nbige(zu^{-1})
  \bigr).
 \end{CD}
\]
The morphism 
$d^{\nbigv}_{+,\varrho}-d^{\nbigv}_{-,\varrho}$
is the composition of the morphism
$j_{Z!}\nbigk_{J,\varrho}(\nbigv)
\lrarr
j_{Z!}\nbigntilde'_{J}(\nbigv)$,
and the inclusion
$j_{Z!}\nbigntilde'_{J}(\nbigv)
\lrarr
j_{Z!}\nbigntilde_{J,\varrho}(\nbigv)$,
and there exists the following commutative diagram:
\[
 \begin{CD}
j_{Z!}\nbigk_{J,!}(V)
@>>>
j_{Z!}\nbigntilde'_{J}(V)
@>>>
j_{Z!}\nbigntilde_{J,!}(V).
  \\
@VVV @V{=}VV @VVV \\ 
j_{Z!}\nbigk_{J,\varrho}(\nbigv)
@>>>
j_{Z!}\nbigntilde'_{J}(\nbigv)
@>>>
j_{Z!}\nbigntilde_{J,\varrho}(\nbigv)
  \\
@VVV @V{=}VV @VVV \\ 
j_{Z!}\nbigk_{J,\ast}(V)
@>>>
j_{Z!}\nbigntilde'_{J}(V)
@>>>
j_{Z!}\nbigntilde_{J,\ast}(V).
 \end{CD}
\]  
Then, we obtain (\ref{eq;24.3.21.10})
and (\ref{eq;24.3.21.11}).
\hfill\qed

\section{Stokes filtrations}
\label{subsection;24.3.14.41}

\subsection{}

Let $u=|u|e^{\sqrt{-1}\theta^u}$.
There exist the isomorphisms:
\[
 \gbigl^{\gbigf}_!(\nbigv)_{|\theta^u}
 \simeq
 H_1^{\rd}\bigl(
 \cnum^{\ast},
 \nbigv\otimes\nbige(zu^{-1})
 \bigr),
 \quad
 \gbigl^{\gbigf}_{\ast}(\nbigv)_{|\theta^u}
 \simeq
 H_1^{\mg}\bigl(
 \cnum^{\ast},
 \nbigv\otimes\nbige(zu^{-1})
 \bigr).
\]
The Stokes filtrations of
$\gbigl^{\gbigf}_{\varrho}(\nbigv)_{|\theta^u}$
$(\varrho=!,\ast)$
induce
the filtrations
$\nbigf^{\circ\,\theta^u}$
on the spaces
$H_1^{\varrho}\bigl(
 \cnum^{\ast},
 \nbigv\otimes\nbige(zu^{-1})
 \bigr)$ $(\varrho=\rd,\mg)$
indexed by the partially ordered set
$\Bigl(
 \gbigf^{(0,\infty)}_+\bigl(
 \nbigi(\nbigv)
 \bigr),
 \leq_{\theta^u}
 \Bigr)$.
Similarly,
we obtain the filtrations
$\nbigf^{\circ\,\theta^u}$
on the spaces
$H_1^{\varrho}\bigl(
 \cnum^{\ast},
 \nbigt_{\omega}(\nbigv)\otimes\nbige(zu^{-1})
 \bigr)$
($\varrho=\rd,\mg$)
indexed by the partially ordered set
$\Bigl(
\gbigf^{(0,\infty)}_+(\nbigt_{\omega}(\nbigi(\nbigv))),
\leq_{\theta^u}
\Bigr)$.

The following lemma is obvious by the constructions.
(See \S\ref{subsection;24.3.14.42}
for the isomorphism
$\gbigl^{\gbigf}_{\varrho}(\nbigv)_{|\theta^u_1}
 \simeq
 \gbigl^{\gbigf}_{\varrho}(\nbigv)_{|\theta^u_2}$.)
\begin{lem}
Let $\vecJ\in T(\nbigi^{\circ})$.
For any $\theta^u_1,\theta^u_2\in \vecJbar_{\pm}$,
we have 
\[
 \Abb^{\rd}_{\nu_0^-(\vecJ),\theta^u_1}
 =\Abb^{\rd}_{\nu_0^-(\vecJ),\theta^u_2},
 \quad
 B_{\nu_0^+(\vecJ)_{\pm},\theta^u_1}
 =B_{\nu_0^+(\vecJ)_{\pm},\theta^u_2},
\]
\[
 C^{\nu_0^-(\vecJ)_{\pm}}_{\infty,\theta^u_1}
 =C^{\nu_0^-(\vecJ)_{\pm}}_{\infty,\theta^u_2},
\quad
  C^{\mg,\nu_0^-(\vecJ)_{\pm}}_{\infty,\theta^u_1}
 =C^{\mg,\nu_0^-(\vecJ)_{\pm}}_{\infty,\theta^u_2}
\]
under the natural isomorphisms
$\gbigl^{\gbigf}_{\varrho}(\nbigv)_{|\theta^u_1}
 \simeq
 \gbigl^{\gbigf}_{\varrho}(\nbigv)_{|\theta^u_2}$.
\hfill\qed
\end{lem}

\subsection{}
\label{subsection;24.4.4.1}

Recall
$V=\nbigs_{\omega}(\nbigv)$,
$\nbigitilde=\nbigi(V)$
and $\nbigi=\pi_{\omega}(\nbigitilde)$.
We set $\nbigi^{\circ}=\gbigf^{(0,\infty)}_+(\nbigi)$
and $\nbigitilde^{\circ}=\gbigf^{(0,\infty)}_+(\nbigitilde)$.
We set
\index{sets $\gbigm_{\pm}(\nbigi^{\circ},\theta^u)$}
\[
 \gbigm_-(\nbigi^{\circ},\theta^u)
=\bigl\{
 \vecJ\in T(\nbigi^{\circ})\,\big|\,
 \theta^u\in \vecJ_-
 \bigr\},
\quad
  \gbigm_+(\nbigi^{\circ},\theta^u)
=\bigl\{
 \vecJ\in T(\nbigi^{\circ})\,\big|\,
 \theta^u\in \vecJ_+
 \bigr\}.
\]

\begin{lem}
For any $\vecJ\in T(\nbigi^{\circ})$,
the following conditions are equivalent. 
\begin{itemize}
 \item $\vecJ\in\gbigm_-(\nbigi^{\circ},\theta^u)$
 \item $\nu_0^-(\vecJ)_-\cap (\vecI(\theta^u)+\pi)_+\neq\emptyset$.
 \item $\nu_0^+(\vecJ)_-\cap \vecI(\theta^u)_+\neq\emptyset$.
\end{itemize}
In the case,
$\kappa^-_{0,\vecJ}(\theta^u)
 \in
 \nu^-_{0}(\vecJ)_-
 \cap
 (\vecI(\theta^u)+\pi)_+$
and
$\kappa^+_{0,\vecJ}(\theta^u)
 \in
 \nu^+_{0}(\vecJ)_-
 \cap
 \vecI(\theta^u)_+$
hold.

Similarly, the following conditions
are equivalent.
\begin{itemize}
 \item $\vecJ\in\gbigm_+(\nbigi^{\circ},\theta^u)$.
 \item $\nu_0^-(\vecJ)_+\cap (\vecI(\theta^u)+\pi)_-\neq\emptyset$.
 \item $\nu_0^+(\vecJ)_+\cap \vecI(\theta^u)_-\neq\emptyset$.
\end{itemize} 
In the case,
$\kappa^-_{0,\vecJ}(\theta^u)
 \in
 \nu^-_{0}(\vecJ)_+
 \cap
 (\vecI(\theta^u)+\pi)_-$
and
$\kappa^+_{0,\vecJ}(\theta^u)
 \in
 \nu^+_{0}(\vecJ)_+
 \cap
 \vecI(\theta^u)_-$
hold.
\hfill\qed
\end{lem}

Recall that there exist the isomorphisms
of the partially ordered sets
in Proposition \ref{prop;18.5.5.40}:
\begin{equation}
\label{eq;24.3.14.40}
(\nbigitilde^{\circ}_{\vecJ,<0},\leq_{\theta^u})
\simeq
(\nbigitilde_{\nu^-_0(\vecJ),<0},\leq_{\kappa^-_{0,\vecJ}(\theta^u)}),
\quad
(\nbigitilde^{\circ}_{\vecJ,>0},\leq_{\theta^u})
\simeq
(\nbigitilde_{\nu^+_0(\vecJ),>0},\leq_{\kappa^+_{0,\vecJ}(\theta^u)}).
\end{equation}
When $\theta_u\in \vecJbar$,
we obtain the filtration $\nbigf^{\prime \theta^u}$
of
\[
 H^0(\nu^-_0(\vecJ),L_{\nu_0^-(\vecJ),<0})
 \simeq
 H^0(\overline{\nu^-_0(\vecJ)},L_{\nu_0^-(\vecJ),<0})
\]
indexed by the partially ordered set
$(\nbigitilde^{\circ}_{\vecJ,<0},\leq_{\theta^u})$
from the filtration
$\nbigftilde^{\kappa^-_{0,\vecJ}(\theta^u)}$
indexed by
$(\nbigitilde_{\nu^-_0(\vecJ),<0},\leq_{\kappa^-_{0,\vecJ}(\theta^u)})$.
We also obtain the filtration $\nbigf^{\prime \theta^u}$
of
\[
 H^0(\nu^+_0(\vecJ),L_{\nu_0^+(\vecJ),>0})
 \simeq
 H^0(\overline{\nu^+_0(\vecJ)},L_{\nu_0^+(\vecJ),>0})
\]
indexed by the partially ordered set
$(\nbigitilde^{\circ}_{\vecJ,>0},\leq_{\theta^u})$
from the filtration
$\nbigftilde^{\kappa^+_{0,\vecJ}(\theta^u)}$
indexed by
$(\nbigitilde_{\nu^+_0(\vecJ),>0},\leq_{\kappa^+_{0,\vecJ}(\theta^u)})$.

\subsection{Isomorphisms of filtered vector spaces}

Let $\vecJ_1\in \gbigm_-(\nbigi^{\circ},\theta^u)$.
According to Corollary \ref{cor;24.3.14.20},
we obtain the following isomorphism
induced by
$\Abb^{\rd}_{\nu_0^-(\vecJ),\theta^u}$,
$B_{\nu_0^+(\vecJ)_-,\theta^u}$
$(\vecJ\in \gbigm_-(\nbigi^{\circ},\theta^u))$
and
$C^{\nu_0^-(\vecJ_1)_-}_{\infty,\theta^u}$:
{\small
\begin{multline}
\label{eq;24.3.14.10}
\!\!\!\!\!
 \bigoplus_{\vecJ\in\gbigm_-(\nbigi^{\circ},\theta^u)}\!\!
 \Bigl(
 H^0(\nu^-_0(\vecJ),L_{\nu_0^-(\vecJ),<0})
\oplus
 H^0(\nu^+_0(\vecJ),L_{\nu_0^+(\vecJ),>0})
  \Bigr)
 \oplus
 H_1^{\rd}\bigl(
 \cnum^{\ast},\nbigt_{\omega}(\nbigv)\otimes
 \nbige(zu^{-1})
 \bigr)
 \\
\stackrel{\simeq}{\lrarr}
  H_1^{\rd}\bigl(
 \cnum^{\ast},\nbigv\otimes
 \nbige(zu^{-1})
 \bigr).
\end{multline}
}
According to Corollary \ref{cor;24.3.14.21},
we also obtain the following isomorphism
induced by
$\Abb^{\rd}_{\nu_0^-(\vecJ),\theta^u}$,
$B_{\nu_0^+(\vecJ)_-,\theta^u}$
$(\vecJ\in \gbigm_-(\nbigi^{\circ},\theta^u))$
and
$C^{\mg,\nu_0^-(\vecJ_1)_-}_{\infty,\theta^u}$:
{\small
\begin{multline}
\label{eq;24.3.14.11}
 \!\!\!\!\!
 \bigoplus_{\vecJ\in\gbigm_-(\nbigi^{\circ},\theta^u)}\!\!
 \Bigl(
 H^0(\nu^-_0(\vecJ),L_{\nu_0^-(\vecJ),<0})
\oplus
 H^0(\nu^+_0(\vecJ),L_{\nu_0^+(\vecJ),>0})
  \Bigr)
 \oplus
 H_1^{\mg}\bigl(
 \cnum^{\ast},\nbigt_{\omega}(\nbigv)\otimes
 \nbige(zu^{-1})
 \bigr)
 \\
\stackrel{\simeq}{\lrarr}
  H_1^{\mg}\bigl(
 \cnum^{\ast},\nbigv\otimes\nbige(zu^{-1})
 \bigr).
\end{multline}
}
Similarly, 
for $\vecJ_2\in\gbigm_+(\nbigi^{\circ},\theta^u)$,
we obtain the following isomorphism
induced by
$\Abb^{\rd}_{\nu_0^-(\vecJ),\theta^u}$,
$B_{\nu_0^+(\vecJ)_+,\theta^u}$
$(\vecJ\in \gbigm_+(\nbigi^{\circ},\theta^u))$
and
$C^{\nu_0^-(\vecJ_2)_+}_{\infty,\theta^u}$:
{\small
\begin{multline}
\label{eq;24.3.14.12}
 \!\!\!\!\!
 \bigoplus_{\vecJ\in\gbigm_+(\nbigi^{\circ},\theta^u)}\!\!
 \Bigl(
 H^0(\nu^-_0(\vecJ),L_{\nu_0^-(\vecJ),<0})
\oplus
 H^0(\nu^+_0(\vecJ),L_{\nu_0^+(\vecJ),>0})
  \Bigr)
 \oplus
 H_1^{\rd}\bigl(
 \cnum^{\ast},\nbigt_{\omega}(\nbigv)\otimes
 \nbige(zu^{-1})
 \bigr)
 \\
\stackrel{\simeq}{\lrarr}
  H_1^{\rd}\bigl(
 \cnum^{\ast},\nbigv\otimes
 \nbige(zu^{-1})
 \bigr).
\end{multline}
}
We also obtain the following
induced by
$\Abb^{\rd}_{\nu_0^-(\vecJ),\theta^u}$,
$B_{\nu_0^+(\vecJ)_+,\theta^u}$
$(\vecJ\in \gbigm_+(\nbigi^{\circ},\theta^u))$
and
$C^{\mg,\nu_0^-(\vecJ_2)_+}_{\infty,\theta^u}$:
{\small
\begin{multline}
\label{eq;24.3.14.13}
 \!\!\!\!\!
 \bigoplus_{\vecJ\in\gbigm_+(\nbigi^{\circ},\theta^u)}\!\!
 \Bigl(
 H^0(\nu^-_0(\vecJ),L_{\nu_0^-(\vecJ),<0})
\oplus
 H^0(\nu^+_0(\vecJ),L_{\nu_0^+(\vecJ),>0})
  \Bigr)
 \oplus
 H_1^{\mg}\bigl(
 \cnum^{\ast},\nbigt_{\omega}(\nbigv)\otimes
 \nbige(zu^{-1})
 \bigr)
 \\
\stackrel{\simeq}{\lrarr}
  H_1^{\mg}\bigl(
 \cnum^{\ast},\nbigv\otimes\nbige(zu^{-1})
 \bigr).
\end{multline}
}
Note that
\[
 \gbigf^{(0,\infty)}_+(\nbigi(\nbigv))
 =\bigl(
 \nbigitilde^{\circ}\setminus\{0\}
 \bigr)
 \sqcup
 \gbigf^{(0,\infty)}_+\bigl(
 \nbigt_{\omega}(\nbigi(\nbigv))
 \bigr).
\]
The left hand sides of
(\ref{eq;24.3.14.10}),
(\ref{eq;24.3.14.11}),
(\ref{eq;24.3.14.12})
and (\ref{eq;24.3.14.13})
are equipped with the filtration
$\nbigf^{\prime\theta^u}$
obtained from the filtrations
$\nbigf^{\prime\theta^u}$
on $H^0(\nu_{0}^{-}(\vecJ),L_{\nu^-_0(\vecj),<0})$
and $H^0(\nu_{0}^{+}(\vecJ),L_{\nu^+_0(\vecj),<0})$
$(\vecJ\in \gbigm_{\pm}(\nbigi^{\circ},\theta^u))$,
and $\nbigf^{\circ\,\theta^u}$
on 
$H_1^{\varrho}\bigl(
 \cnum^{\ast},\nbigt_{\omega}(\nbigv)\otimes
 \nbige(zu^{-1})
 \bigr)$ $(\varrho=\rd,\mg)$.
The right hand sides are also equipped with the filtration
$\nbigf^{\circ\theta^u}$ 
induced by the Stokes filtrations of
$\gbigl^{\gbigf}_{\varrho}(\nbigv)$
$(\varrho=!,\ast)$.
The following is one of the main theorem,
which we shall prove in \S\ref{section;20.11.21.1}.
\begin{thm}
\label{thm;24.3.15.10}
 The isomorphisms {\rm(\ref{eq;24.3.14.10}), (\ref{eq;24.3.14.11})
(\ref{eq;24.3.14.12}), (\ref{eq;24.3.14.13})}
are isomorphisms of  
filtered vector spaces. 
\end{thm}

\subsection{Some canonically defined subspaces}
\label{subsection;24.3.24.10}

By Theorem \ref{thm;24.3.15.10},
we obtain the following corollary.
\begin{cor}
\label{cor;25.2.22.11}
For any $\theta^u\in\vecJbar$,
$\Abb^{\rd}_{\nu_0^-(\vecJ),\theta^u}$ induce isomorphisms
of filtered vector spaces
\[
 H^0(\nu_0^-(\vecJ),L_{\nu_0^-(\vecJ),<0})
 \simeq
 H^0(\vecJ,\gbigl^{\gbigf}_{\varrho}(\nbigv)_{\vecJ,<0}).
\] 
 For any $\theta^u\in\vecJ_{\pm}$,
$B_{\nu_0^+(\vecJ)_{\pm},\theta^u}$ induce isomorphisms
of filtered vector spaces
 \[
 H^0(\nu_0^+(\vecJ)_{\pm},L_{\nu_0^+(\vecJ)_{\pm},>0})
 \simeq
 H^0(\vecJ_{\pm},\gbigl^{\gbigf}_{\varrho}(\nbigv)_{\vecJ_{\pm},>0}).
\] 
Here, we use the isomorphism of
the partially ordered sets in {\rm(\ref{eq;24.3.14.40})}
to identify the index sets of the filtrations.
\hfill\qed
\end{cor}

\begin{cor}
\label{cor;25.2.15.30}
For any $\theta^u\in\vecJ_{\pm}$,
$C^{\nu_0^-(\vecJ)_{\pm}}_{\infty,\theta^u}$ induce
isomorphisms of filtered vector spaces
\begin{equation}
\label{eq;24.3.25.52}
 H^0(\vecJ_{\pm},\gbigl^{\gbigf}_{!}(\nbigt_{\omega}\nbigv))
 \simeq
 H^0(\vecJ_{\pm},\gbigl^{\gbigf}_{!}(\nbigv)_{\vecJ_{\pm},0}).
\end{equation}
and
$C^{\mg,\nu_0^-(\vecJ)_{\pm}}_{\infty,\theta^u}$ induce isomorphisms
of filtered vector spaces
\begin{equation}
\label{eq;24.3.25.51}
 H^0(\vecJ_{\pm},\gbigl^{\gbigf}_{\ast}(\nbigt_{\omega}\nbigv))
 \simeq
 H^0(\vecJ_{\pm},\gbigl^{\gbigf}_{\ast}(\nbigv)_{\vecJ_{\pm},0}).
\end{equation}
As a result,
by setting $\omega^{\circ}=(1+\omega)^{-1}\omega$,
we obtain
$\gbigl^{\gbigf}_{\star}(\nbigt_{\omega}(\nbigv))
\simeq
 \nbigt_{\omega^{\circ}}\bigl(
 \gbigl^{\gbigf}_{\star}(\nbigv)
 \bigr)$. 
\hfill\qed
\end{cor}

\subsection{Transformations of cycles adapted to Stokes filtrations}
\label{subsection;24.3.24.11}

Let $\vecJ\in T(\nbigi^{\circ})$.

\subsubsection{}

We set
$\theta^u=\vartheta^{\vecJ}_{\ell}
=\vartheta^{\vecJ-(1+\omega^{-1})\pi}_{r}$.
We have
$\nu_0^-(\vecJ-(1+\omega^{-1})\pi)
=\nu_0^+(\vecJ)-\omega^{-1}\pi$,
and
\[
 \kappa^+_{0,\vecJ}(\theta^u)
=\kappa^-_{0,\vecJ-(1+\omega^{-1})\pi}(\theta^u)
=\theta^u-\pi/2
=\vartheta^{\nu_0^+(\vecJ)}_{\ell}.
\]
For $v\in H^0(J,L_{J,>0})$,
by the construction in (\ref{eq;24.3.16.1}),
we obtain
\begin{multline}
 B_{\nu_0^+(\vecJ)_-,\theta^u}(v)
 =\Abb^{\rd}_{\nu_0^-(\vecJ-(1+\omega^{-1})\pi),\theta^u}
 \bigl(
 \nbigrtilde^{\nu_0^+(\vecJ)_-}_{\nu_0^-(\vecJ-(1+\omega^{-1})\pi)}(v)
 \bigr)
\\
 +\sum_{\vecJ-(1+\omega^{-1})\pi<\vecJ'<\vecJ-\pi}
 \Abb^{\rd}_{\nu_0^-(\vecJ'),\theta^u}
 \bigl(
 \nbigrtilde^{\nu_0^+(\vecJ)_-}_{\nu_0^-(\vecJ')}(v)
 \bigr)
-\sum_{\vecJ-\pi\leq \vecJ'<\vecJ}
 \Abb^{\rd}_{\nu_0^-(\vecJ'),\theta^u}
 \bigl(
 \nbigrtilde^{\nu_0^+(\vecJ)_-}_{\nu_0^-(\vecJ')}(v)
 \bigr)
 \\
-\sum_{\vecJ-(\omega^{-1}-1)\pi\leq \vecJ'<\vecJ}
 \Abb^{\rd}_{\nu_{-1}^-(\vecJ'),\theta^u-2\pi}
 \bigl(
 \nbigrtilde^{\nu_0^+(\vecJ)_-}_{\nu_{-1}^-(\vecJ')}(v)
 \bigr).
\end{multline}
Note that
$\nbigrtilde^{\nu_0^+(\vecJ)_-}_{\nu_0^-(\vecJ-(1+\omega^{-1})\pi)}$
is an isomorphism
\[
H^0\bigl(\nu_0^+(\vecJ)_-,L_{\nu_0^+(\vecJ)_-,>0}\bigr)
\simeq
H^0\Bigl(
\nu_0^+(\vecJ)-\omega^{-1}\pi,
L_{(\nu_0^+(\vecJ)-\omega^{-1}\pi)_+,<0}
\Bigr)
\]
which preserves the Stokes filtrations
$\nbigf^{\prime\theta^u}$.

\subsubsection{}

We set
$\theta^u
=\vartheta^{\vecJ}_{r}
=\vartheta^{\vecJ+(1+\omega^{-1}\pi)}_{\ell}$.
We have
$\nu_{-1}^-(\vecJ+(1+\omega^{-1})\pi)
=\nu_0^+(\vecJ)+\omega^{-1}\pi$,
and
\[
 \kappa^+_{0,\vecJ}(\theta^u)
=\kappa^-_{-1,\vecJ+(1+\omega^{-1})\pi}(\theta^u)
=\theta^u-3\pi/2
=\vartheta^{\nu_0^+(\vecJ)}_r.
\]
For $v\in H^0(J,L_{J,>0})$,
by the construction in (\ref{eq;24.3.16.2}),
we obtain
\begin{multline}
 B_{\nu_0^+(\vecJ)_+,\theta^u}(v)
=-\Abb^{\rd}_{\nu_{-1}^-(\vecJ+(1+\omega^{-1})\pi),\theta^u-2\pi}
 \bigl(
 \nbigrtilde^{\nu_{0}^+(\vecJ)_+}_{\nu_{-1}^-(\vecJ+(1+\omega^{-1})\pi)}(v)
 \bigr)
\\
-\!\!\!\!\!\sum_{\vecJ+\pi<\vecJ'<\vecJ+(1+\omega^{-1})\pi}\!\!\!\!\!\!\!\!\!
 \Abb^{\rd}_{\nu_{-1}^-(\vecJ'),\theta^u-2\pi}
 \bigl(
 \nbigrtilde^{\nu_0^+(\vecJ)_+}_{\nu_{-1}^-(\vecJ')}(v)
 \bigr)
+\!\!\!\sum_{\vecJ<\vecJ'\leq\vecJ+\pi}\!\!\!\!
 \Abb^{\rd}_{\nu_{-1}^-(\vecJ'),\theta^u-2\pi}
 \bigl(
 \nbigrtilde^{\nu_0^+(\vecJ)_+}_{\nu_{-1}^-(\vecJ')}(v)
 \bigr)
 \\
+\sum_{\vecJ<\vecJ'\leq \vecJ+(\omega^{-1}-1)\pi}
 \Abb^{\rd}_{\nu_{0}^-(\vecJ'),\theta^u}
 \bigl(
 \nbigrtilde^{\nu_0^+(\vecJ)_+}_{\nu_{0}^-(\vecJ')}(v)
 \bigr).
\end{multline}
Note that
$\nbigrtilde^{\nu_0^+(\vecJ)_+}_{\nu_0^-(\vecJ+(1+\omega^{-1})\pi)}$
is an isomorphism
\[
H^0\bigl(\nu_0^+(\vecJ)_+,L_{\nu_0^+(\vecJ)_+,>0}\bigr)
\simeq
H^0\Bigl(
\nu_0^+(\vecJ)+\omega^{-1}\pi,
L_{(\nu_0^+(\vecJ)+\omega^{-1}\pi),<0}
\Bigr)
\]
which preserves the Stokes filtrations
$\nbigf^{\prime\theta^u}$.

\subsubsection{}

For $v\in H^0(\nu_0^+(\vecJ),L_{\nu_0^+(\vecJ),>0})$,
by Proposition \ref{prop;24.2.19.20},
we obtain
\begin{multline}
B_{\nu_0^+(\vecJ)_-,\theta^u}(v)
 - B_{\nu_0^+(\vecJ)_+,\theta^u}(v)
 =B^{\nu_0^-(\vecJ)_+} _{\infty,\theta^u}
 \bigl(
 \nbigq_{\nu_0^+(\vecJ)}(v)
 \bigr)
 \\
 +\Abb^{\rd}_{\nu_0^-(\vecJ),\theta^u}
 \bigl(
 \nbigrtilde^{\nu_0^+(\vecJ)_-}_{\nu_0^-(\vecJ)}(v)
 \bigr)
+\Abb^{\rd}_{\nu_{-1}^-(\vecJ),\theta^u-2\pi}
 \bigl(
 \nbigrtilde^{\nu_0^+(\vecJ)_-}_{\nu_{-1}^-(\vecJ)}(v)
 \bigr).
\end{multline}

\subsubsection{}

For $v\in H^0(\real,\nbigt_{\omega}(L))$,
by Lemma \ref{lem;25.2.11.10} and Lemma \ref{lem;24.4.2.100},
we obtain
\begin{equation}
 B^{\nu_0^-(\vecJ)_-}_{\infty,\theta^u}(v)
 -B^{\nu_0^-(\vecJ)_+}_{\infty,\theta^u}(v)
 =\Abb^{\rd}_{\nu_0^-(\vecJ),\theta^u}\bigl(
 \nbigp_{\nu_0^-(\vecJ)}(v-M_0^{-1}(v))
\bigr),
\end{equation}
\begin{equation}
B^{\mg,\nu_0^-(\vecJ)_-}_{\infty,\theta^u}(v)
-B^{\mg,\nu_0^-(\vecJ)_+}_{\infty,\theta^u}(v)
=\Abb^{\rd}_{\nu_0^-(\vecJ),\theta^u}\bigl(
\nbigp_{\nu_0^-(\vecJ)}(v)
\bigr).
\end{equation}

\subsection{The induced constructible subsheaves and filtrations}

\subsubsection{Constructible subsheaves}

Let $\theta^u\in\real$.
There exist the following isomorphisms for $\star=!,\ast$
induced by $\Abb^{\rd}_{\nu_0^-(\vecJ),\theta^u}$:
\begin{equation}
\label{eq;25.2.24.10}
 \gbigl^{\gbigf}_{\star}(V)^{<0}_{\theta^u}
=\bigoplus_{\theta^u\in\vecJ}
 H^0(\vecJ,\gbigl^{\gbigf}_{\star}(V)_{\vecJ,<0})
 \simeq
\bigoplus_{\theta^u\in\vecJ}
 H^0(\nu_0^-(\vecJ),L_{\nu_0^-(\vecJ),<0}).
\end{equation}
Take $\vecJ_1\in T(\nbigi^{\circ})$
such that $\theta^u\in\vecJ_{1+}$.
The isomorphisms (\ref{eq;25.2.24.10})
extend to the following isomorphisms
by
$B^{\nu_0^-(\vecJ_1)+}_{\infty,\theta^u}$
or
$B^{\mg,\nu_0^-(\vecJ_1)+}_{\infty,\theta^u}$:
\begin{multline}
 \gbigl^{\gbigf}_{\star}(V)^{\leq 0}_{\theta^u}
=\bigoplus_{\theta^u\in\vecJ}
 H^0(\vecJ,\gbigl^{\gbigf}_{\star}(V)_{\vecJ,<0})
 \oplus
 H^0(\vecJ_{1+},\gbigl^{\gbigf}_{\star}(V)_{\vecJ_{1+},0})
\\
 \simeq
 \bigoplus_{\theta^u\in\vecJ}
 H^0(\nu_0^-(\vecJ),L_{\nu_0^-(\vecJ),<0})
 \oplus
 H^0(\nu_0^-(\vecJ_1),L_{\nu_0^-(\vecJ_1)_+,0}).
\end{multline}

\subsubsection{The filtrations on
$H^0(\vecJ,\gbigl^{\gbigf}_{\star}(V)_{\vecJ,<0})$}

For $\vecJ\in T(\nbigi)$,
we have the isomorphism
induced by $\Abb^{\rd}_{\nu_0^-(\vecJ),\theta^u}$
for any $\theta^u\in\vecJ$:
\[
 H^0(\nu_0^-(\vecJ),L_{\nu_0^-(\vecJ),<0})
 \simeq
 H^0(\vecJ,\gbigl^{\gbigf}_{\star}(V)_{<0}).
\]
It is an isomorphism of filtered vector spaces
where we use the isomorphism of
the partially ordered sets in {\rm(\ref{eq;24.3.14.40})}
to identify the index sets of the filtrations.

\subsubsection{The filtrations on
$H^0(\vecJ,\gbigl^{\gbigf}_{\star}(V)_{\vecJ,>0})$}

For $\vecJ\in T(\nbigi)$,
we have the isomorphism
induced by $B_{\nu_0^+(\vecJ)_{\pm},\theta^u}$
for any $\theta^u\in\vecJ$:
\begin{equation}
\label{eq;25.2.24.11}
 H^0(\nu_0^+(\vecJ),L_{\nu_0^+(\vecJ),>0})
 \simeq
 H^0(\vecJ,\gbigl^{\gbigf}_{\star}(V)_{>0}).
\end{equation}
It is independent of the choice of $\pm$.
It is an isomorphism of filtered vector spaces
where we use the isomorphism of
the partially ordered sets in {\rm(\ref{eq;24.3.14.40})}
to identify the index sets of the filtrations.

Because
$\gbigl^{\gbigf}_{\star}(V)/\gbigl^{\gbigf}_{\star}(V)^{\leq 0}
=\bigoplus_{\vecJ\in T(\nbigi^{\circ})}
a_{\vecJ\ast}
\gbigl_{\star}^{\gbigf}(V)_{\vecJ,>0}$,
there exists the following morphisms
(see \S\ref{subsection;24.2.18.1} for $R_{\vecJ}$):
\begin{multline}
\label{eq;25.2.24.12}
H^0(\nu_0^+(\vecJ),L_{\nu_0^+(\vecJ),>0})
 \simeq
 H^0(\nu_0^+(\vecJ)_-,L_{\nu_0^+(\vecJ)_-,>0})
 \lrarr
 H^0(\real,L)
 \\
 \stackrel{f}{\lrarr}
 H^0(\real,\gbigl^{\gbigf}_{\star}(V))
 \stackrel{R_{\vecJ}}{\lrarr}
 H^0\bigl(
 \vecJ,\gbigl^{\gbigf}_{\star}(V)_{\vecJ,>0}
 \bigr).
\end{multline}
\begin{prop}
\label{prop;25.2.25.2}
The composition of {\rm(\ref{eq;25.2.24.12})}
equals {\rm(\ref{eq;25.2.24.11})}.
\end{prop}
\pf
We set
$\theta^u=\vartheta^{\vecJ}_{\ell}$.
We set $\theta_1=\theta^u-\pi/2$.
We set
$J_0=\nu_0^+(\vecJ)$
and $J_1=J_0-\omega^{-1}\pi$.
Note that
$\theta_1=
\vartheta^{J_0}_{\ell}
=\vartheta^{\vecI(\theta^u)}_r$.
We have
\[
 \gbigl^{\gbigf}_{!}(V)^{<0}_{\theta^u}=
 \bigoplus_{J_1<J'\leq J_0+\pi}
 \Image \Abb^{\rd}_{J',\theta^u}.
\]

Let $v\in H^0(J_{0-},L_{J_{0-},>0})$.
We set
$v'=\nbigr^{J_{0-}}_{J_1}(v)\in H^0(J_{1},L_{J_1,<0})$.
(See \S\ref{subsection;24.2.19.1} for the map $\nbigr^{J_{0-}}_{J_1}$.)
By (\ref{eq;24.3.16.1}),
we have
$B_{J_{0-},\theta^u}(v)
 \equiv
 \Abb^{\rd}_{J_1,\theta^u}(v')$
in
$\gbigl^{\gbigf}_{\star}(V)_{\theta^u}\big/
 \gbigl^{\gbigf}_{\star}(V)^{\leq 0}_{\theta^u}$.
We also have
\[
\Abb^{\rd}_{\infty,\theta^u}(v)
=
\sum_{J_1\leq J<J_0}
\bigl(
 \Abb^{\rd}_{J,\theta^u}(\nbigr^{J_0}_J(v))
-\Abb^{\rd}_{J,\theta^u-2\pi}(\nbigr^{J_0}_J(v))
\bigr).
\]
We have
\[
 \sum_{J_1<J<J_0}
 \Abb^{\rd}_{J,\theta^u}(\nbigr^{J_0}_J(v))
 \in\gbigl^{\gbigf}_!(V)^{<0}_{\theta^u}.
\]
\begin{lem}
For $J_0-\pi<J\leq J_0+\omega^{-1}\pi$,
we have 
\begin{equation}
\label{eq;25.2.25.1}
 \Image \Abb^{\rd}_{J,\theta^u-2\pi}
 \subset
 \gbigl^{\gbigf}_{!}(V)^{<0}_{\theta^u}
 \oplus
 \bigoplus_{J_1-\pi<J'\leq J-\omega^{-1}\pi}
 \Image B_{J'_+,\theta^u}.
\end{equation}
\end{lem}
\pf
By using (\ref{eq;24.3.16.2}),
we can prove (\ref{eq;25.2.25.1})
for $J_0-\pi<J\leq J_0-\pi+a\pi$
$(0\leq a\leq 1+\omega^{-1})$
by an induction on $a$.
\hfill\qed

\vspace{.1in}

Let $v\in H^0(J_+,L_{J_+,>0})$.
We set $\vecJ_3=\vecJ-(1+\omega^{-1})\pi$.
We have
\[
 B_{J_{0+},\theta^u}(v)-
 \Abb^{\rd}_{\theta^u,v}
 \in
 \gbigl^{\gbigf}_!(V)^{<0}_{\theta^u}
 \oplus
 \bigoplus_{\vecJ_3<\vecJ'<\vecJ-\omega^{-1}\pi}
 H^0(\vecJ'_+,\gbigl^{\gbigf}_{!}(V)_{\vecJ'_+,>0}).
\]
Therefore, we obtain the claim of
Proposition \ref{prop;25.2.25.2}.
\hfill\qed

\begin{rem}
\label{rem;25.4.12.20}
To the best of the author's understanding,
we may also obtain the above descriptions of
$\gbigl_{\star}^{\gbigf}(V)^{<0}
\subset
\gbigl_{\star}^{\gbigf}(V)^{\leq 0}
\subset 
\gbigl_{\star}^{\gbigf}(V)$
and the Stokes filtrations 
on $\gbigl_{\star}^{\gbigf}(V)^{<0}$
and $\gbigl_{\star}^{\gbigf}(V)/\gbigl_{\star}^{\gbigf}(V)^{\leq 0}$
by applying the results
in {\rm\cite[VII, VIII]{Malgrange-book}}
to the cases $V(\star 0)$ $(\star=!,\ast)$.
\hfill\qed 
\end{rem}

\section{Extensions and the recovery of the Stokes structure}
\subsection{Preliminary}

There exists the following natural morphisms
of $2\pi\seisuu$-equivariant local systems.
\begin{equation}
\label{eq;25.2.15.2} 
\gbigl^{\gbigf}_!(V^{\reg})
\stackrel{a_1}{\lrarr}
\gbigl^{\gbigf}_!(V)
\stackrel{a_2}{\lrarr}
\gbigl^{\gbigf}_{\ast}(V)
\stackrel{a_3}{\lrarr}
\gbigl^{\gbigf}_{\ast}(V^{\reg}).
\end{equation}
As explained
in \S\ref{subsection;24.3.25.100} and
\S\ref{subsection;24.4.19.2},
there exist the natural isomorphisms
\begin{equation}
\label{eq;25.2.15.1}
 \gbigl^{\gbigf}_{!}(V^{\reg})\simeq
 L\simeq
 \gbigl^{\gbigf}_{\ast}(V^{\reg}).
\end{equation}
The following lemma is obvious.
\begin{lem}
Under the isomorphisms {\rm(\ref{eq;25.2.15.1})},
the induced endomorphism $a_3\circ a_2\circ a_1$ of $L$
is $\id-M^{-1}$.
\hfill\qed
\end{lem}

Let $M^{\gbigf}_!$ and $M^{\gbigf}_{\ast}$
denote the monodromy automorphisms
$\gbigl^{\gbigf}_!(V)$
and 
$\gbigl^{\gbigf}_{\ast}(V)$,
respectively.
\index{automorphisms $M^{\gbigf}_{\star}$}
We set $f_!=a_1\circ a_3\circ a_2$
and 
$f_{\ast}=a_2\circ a_1\circ a_3$.

\begin{lem}
Under the identification {\rm(\ref{eq;25.2.15.1})},
we have 
$f_!=\id-(M^{\gbigf}_!)^{-1}$
and 
$f_{\ast}=\id-(M^{\gbigf}_{\ast})^{-1}$.
\end{lem}
\pf
Let $\theta^u\in\real$.
By Lemma \ref{lem;24.4.5.10} and Lemma \ref{lem;24.4.5.2},
we obtain
\[
 f_!(\Abb^{\rd}_{\infty,\theta^u}(v))
 =\Abb^{\rd}_{\infty,\theta^u}(v-M^{-1}(v))
=(\id-(M^{\gbigf}_!)^{-1})(\Abb^{\rd}_{\infty,\theta^u}(v)).
\]
By using Lemma \ref{lem;24.4.5.1}
and Lemma \ref{lem;24.4.5.10}, we obtain
\[
 f_!(\Abb^{\rd}_{J,\theta^u}(v))
 =\Abb^{\rd}_{\infty,\theta^u}(v)
 =\Abb^{\rd}_{J,\theta^u}(v)-\Abb^{\rd}_{J,\theta^u-2\pi}(v)
=(\id-(M^{\gbigf}_!)^{-1})(\Abb^{\rd}_{J,\theta^u}(v)).
\]
Thus, we obtain
$f_!=\id-(M^{\gbigf}_!)^{-1}$.
By Remark \ref{rem;25.2.11.21},
we have
\[
 f_{\ast}\bigl(
 \BB^{\mg}_{J,u}(v)
 \bigr)
 =0
 =(\id-(M^{\gbigf}_{\ast})^{-1})(\BB^{\mg}_{J,u}(v)).
\]
Let $J_1\in T(\nbigi)$.
By using Lemma \ref{lem;24.4.5.10},
we obtain the following equality:
\begin{multline}
 f_{\ast}\bigl(
 \Abb^{\mg,J_{1+}}_{\infty,\theta^u}(v)
 \bigr)
=a_2(\Abb^{\rd}_{\infty,\theta^u}(v))
\\
 = \Abb^{\mg,J_{1+}}_{\infty,\theta^u}(v)
-\Abb^{\mg,J_{1+}}_{\infty,\theta^u}(M^{-1}(v))
 +\sum_{J_1-2\pi<J'\leq J_1}
 \BB^{\mg}_{J',u}(R_{J'}(v)).
\end{multline}
By using Lemma \ref{lem;24.4.5.3} and Lemma \ref{lem;24.3.25.10},
we obtain
\begin{multline}
(\id-(M^{\gbigf}_{\ast})^{-1})(\Abb^{\mg,J_{1+}}_{\infty,\theta^u}(v))
=\Abb^{\mg,J_{1+}}_{\infty,\theta^u}(v)
 -\Abb^{\mg,J_{1+}}_{\infty,\theta^u-2\pi}(v)\\
 =\Abb^{\mg,J_{1+}}_{\infty,\theta^u}(v)
 -\Abb^{\mg,J_{1+}}_{\infty,\theta^u}(M^{-1}(v))
+\sum_{J_1-2\pi<J'\leq J_1}
 \BB^{\mg}_{J',u}(R_{J'}(v)).
\end{multline}
Hence, $f_{\ast}=\id-(M^{\gbigf}_{\ast})^{-1}$.
\hfill\qed

\begin{lem}
Let $\omega^{\circ}=\omega(1+\omega)^{-1}$.
The following diagram is commutative:
\[
\begin{CD}
 \nbigt_{\omega^{\circ}}\bigl(
 \gbigl^{\gbigf}_{!}(V)
 \bigr)
 @>>>
 \nbigt_{\omega^{\circ}}\bigl(
 \gbigl^{\gbigf}_{\ast}(V)
 \bigr)\\
 @V{\simeq}VV @V{\simeq}VV \\
 \nbigt_{\omega}(L)
 @>{\id-M_0^{-1}}>>
 \nbigt_{\omega}(L).
\end{CD} 
\] 
Here, the vertical arrows are the isomorphisms
in Corollary {\rm\ref{cor;25.2.15.30}}.
\end{lem}
\pf
It follows from Lemma \ref{lem;25.2.15.21}.
\hfill\qed

\subsection{Extension}

Let $L_1$ be a $2\pi\seisuu$-equivariant local system.
Let $M_{L_1}$ be the monodromy automorphism.
We consider morphisms of $2\pi\seisuu$-equivariant local systems
\begin{equation}
\label{eq;25.2.15.40}
\nbigt_{\omega}(L)
\stackrel{a}{\lrarr} L_1
\stackrel{b}{\lrarr} \nbigt_{\omega}(L)
\end{equation}
such that $b\circ a=\id-M_0^{-1}$.
As the extension of
$(\gbigl^{\gbigf}_!(V),\vecnbigf)
\to
(\gbigl^{\gbigf}_{\ast}(V),\vecnbigf)$
by (\ref{eq;25.2.15.40}),
we obtain 
\begin{equation}
 (\gbigl^{\gbigf}_!(V),\vecnbigf)
\stackrel{u_1}{\lrarr}
 (\Ltilde_1,\vecnbigf)
\stackrel{u_2}{\lrarr}
(\gbigl^{\gbigf}_{\ast}(V),\vecnbigf)
\end{equation}
in $\Loc^{\St}(\nbigi^{\circ})$.
(See Theorem \ref{thm;20.11.3.20}.)
Together with
(\ref{eq;25.2.15.2}) and
(\ref{eq;25.2.15.1}),
we obtain the induced morphisms
of $2\pi\seisuu$-equivariant local systems:
\begin{equation}
\label{eq;25.3.16.10}
 L
\stackrel{\atilde}{\lrarr}  
 \Ltilde_1
\stackrel{\btilde}{\lrarr}
 L.
\end{equation}

\subsection{The induced endomorphisms}
By the construction,
$\btilde\circ\atilde=\id-M^{-1}$.
Let $M_{\Ltilde_1}$ denote the monodromy automorphism of
$\Ltilde_1$.
\begin{prop}
\label{prop;25.2.22.10}
If $a\circ b=\id-M_{L_1}^{-1}$ holds,
then
$\atilde\circ\btilde
=\id-M_{\Ltilde_1}^{-1}$ holds.
\end{prop}
\pf
Let $\theta^u\in T(\nbigi^{\circ})$.
Let $\vecJ_1\in\gbigm_+(\nbigi^{\circ},\theta^u)$
such that
$\theta^u=\vartheta^{\vecJ_1}_r$.
It is equivalent to the condition that
$J< \nu_0^-(\vecJ_1)
\Longleftrightarrow
 \vartheta^{J}_r<\vartheta^{\vecI(\theta^u)}_{r}$.
We set
\[
 W=
 \bigoplus_{\vecJ\in\gbigm_+(\nbigi^{\circ},\theta^u)}\!\!
 \Bigl(
 H^0(\nu^-_0(\vecJ),L_{\nu_0^-(\vecJ),<0})
\oplus
 H^0(\nu^+_0(\vecJ),L_{\nu_0^+(\vecJ),>0})
  \Bigr).
\]
By 
$\Abb^{\rd}_{\nu_0^-(\vecJ),\theta^u}$,
$B_{\nu_0^+(\vecJ)_+,\theta^u}$
$(\vecJ\in \gbigm_+(\nbigi^{\circ},\theta^u))$,
$B^{\nu_0^-(\vecJ_1)_+}_{\infty,\theta^u}$
and
$B^{\mg,\nu_0^-(\vecJ_1)_+}_{\infty,\theta^u}$,
we identify the morphisms
$\gbigl^{\gbigf}_!(\nbigv)_{|\theta^u}
\stackrel{u_1}{\lrarr}
\Ltilde_{1|\theta^u}
\stackrel{u_2}{\lrarr}
\gbigl^{\gbigf}_{\ast}(\nbigv)_{|\theta^u}$
with the morphisms
\begin{equation}
\label{eq;25.2.15.41}
W\oplus  H^0(\real,\nbigt_{\omega}(L))
\lrarr
W\oplus  H^0(\real,L_1)
\lrarr
W\oplus  H^0(\real,\nbigt_{\omega}(L)),
\end{equation}
which are the direct sum of the identity map of $W$
and the morphisms
\[
H^0(\real,\nbigt_{\omega}(L))
\to
H^0(\real,L_1)
\to
H^0(\real,\nbigt_{\omega}(L))
\]
induced by (\ref{eq;25.2.15.40}).
We describe elements
of $\gbigl^{\gbigf}_{\star}(\nbigv)_{|\theta^u}$ $(\star=!,\ast)$
and $\Ltilde_{1|\theta^u}$
as $(s_1,s_2)$
according to the decompositions in (\ref{eq;25.2.15.41}).

Let $s=(s_1,0)\in \Ltilde_{1|\theta^u}$.
We have
$s'=(s_1,0)\in \gbigl^{\gbigf}_!(V)_{|\theta^u}$
which satisfies $u_1(s')=s$.
Because
$f_{!}(s')
=\bigl(
\id-(M^{\gbigf}_{!})^{-1}
\bigr)(s')$,
we obtain
\begin{multline}
 \atilde\circ\btilde(s)
 =
 (u_1\circ a_1)
 \circ
 (a_3\circ u_2)(u_1(s'))
 =u_1\circ f_!(s')=u_1\bigl(
 (\id-(M^{\gbigf}_!)^{-1}(s'))
 \bigr)
\\ 
=(\id-M_{\Ltilde_1}^{-1})(s).
\end{multline}

Let $s=(0,s_2)\in\Ltilde_{1|\theta^u}$.
It induces
$u_2(s)=(0,b(s_2))
\in \gbigl^{\gbigf}_{\ast}(\nbigv)_{|\theta^u}$.
The element
$a_3\circ u_2(s)\in H^0(\real,L)$
equals the element 
induced from $b(s_2)\in H^0(\real,\nbigt_{\omega}(L))$
by 
\[
H^0(\real,\nbigt_{\omega}(L))
\simeq
H^0(\nu_0^-(\vecJ_1)_+,L_{0,\nu^0_-(\vecJ_1)_+})
\subset
H^0(\real,L).
\]
We set
$v:=a_3\circ u_2(s)
\in H^0(\nu_0^-(\vecJ_1)_+,L_{0,\nu_0^-(\vecJ_1)_+})$.
Under the isomorphism
$\gbigl^{\gbigf}_!(V)_{|\theta^u}
\simeq H_1^{\rd}(\cnum^{\ast},V\otimes\nbige(zu^{-1}))$,
we obtain the following by (\ref{eq;25.2.22.1}):
\begin{multline}
\label{eq;25.2.16.1}
 a_1\circ a_3\circ u_2(s)
=\Abb_{\infty,\theta^u}\bigl(
 v
 \bigr)
 \\
=B^{\nu_0^-(\vecJ_1)_+}_{\infty,\theta^u}(v)
 -\!\!\!\!\!\sum_{\nu_0^-(\vecJ_1)-2\pi<J\leq \nu_0^-(\vecJ_1)}
 \!\!\!\!\!\!
\Abb^{\rd}_{J,\theta^u-2\pi}\bigl(
 \nbigp_{J}(v_{J})
\bigr).
\end{multline}

\begin{lem}
\begin{equation}
\label{eq;25.2.17.2}
 \bigoplus_{\nu_0^-(\vecJ_1)-2\pi<J\leq \nu_0^-(\vecJ_1)}
 \Image \Abb^{\rd}_{J,\theta^u-2\pi}
 \subset W.
\end{equation}
\end{lem}
\pf
By using the notation in \S\ref{subsection;24.2.20.13},
we rewrite $W$ as
\[
W=
 \bigoplus_{J\in\gbigw_1(\nbigi,\vecI(\theta^u)_-)}
 \Image B_{J_+,\theta^u}
 \oplus
 \bigoplus_{J\in\gbigw_2(\nbigi,\vecI(\theta^u)_-)}
 \Image \Abb^{\rd}_{J,\theta^u}.
\]
We set $J_1=\nu_0^-(\vecJ_1)$.
By Lemma \ref{lem;24.4.5.1},
we have
\[
 \bigoplus_{J_1-2\pi<J<J_1-\pi+\omega^{-1}\pi}
 \Image \Abb^{\rd}_{J,\theta^u-2\pi}
\subset W.
\]
If $\omega<1$, then it implies (\ref{eq;25.2.17.2}).
Let us consider the case where $\omega\geq 1$.
For $0\leq a\leq 1-\omega^{-1}$,
we set
\[
 K_a=
 \bigoplus_{J_1-\pi+\omega^{-1}\pi
 \leq J\leq J_1-\pi+\omega^{-1}\pi+a\pi}
 \Image \Abb^{\rd}_{J,\theta^u-2\pi}.
\]
By using (\ref{eq;24.3.16.2}),
we can prove $K_a\subset W$ by using an induction on $a$.
\hfill\qed

\vspace{.1in}

There exists $t_1\in W$
such that
$a_1\circ a_3\circ u_2(s)
 =
 (t_1,b(s_2))
 \in\gbigl^{\gbigf}_!(\nbigv)_{|\theta^u}$.
We obtain
\[
\atilde\circ\btilde(s)
= u_1\circ a_2\circ a_1\circ a_3(0,b(s_2))
=(t_1,a\circ b(s_2)).
\]
As remarked in Lemma \ref{lem;25.2.17.3} below,
because $\omega^{\circ}<1$,
there exists $t_2\in W$ such that
\[
 M_{\Ltilde_1}^{-1}(s)
 =(t_2,M_{L_1}^{-1}(s)).
\]
We have
\[
 u_2\circ \atilde\circ\btilde(s)
=f_{\ast}(u_2(s))
=(\id-(M^{\gbigf}_{\ast})^{-1})(u_2(s))
=u_2\bigl((\id-M_{\Ltilde_1}^{-1})(s)\bigr).
\]
We obtain $t_1=-t_2$,
and $\atilde\circ\btilde(s)=(\id-M_{\Ltilde_1}^{-1})(s)$.
\hfill\qed

\subsubsection{Appendix}
We consider any object $(L_2,\vecnbigf)\in\Loc^{\St}(\nbigi^{\circ})$.
We note that $\omega^{\circ}=-\ord(\nbigi^{\circ})<1$.
Let $M_{L_2}$ denote the monodromy automorphism of $L_2$.
Let $M_{L_2,0}$ denote the monodromy automorphism of
$\nbigt_{\omega^{\circ}}(L_2)$.

Let $\theta^u\in\real$.
Take any $\vecJ_1\in T(\nbigi^{\circ})$
such that $\theta^u\in (\vecJ_1)_+$.
We set
\[
 W_{\theta^u}=\bigoplus_{\vecJ\in \gbigm_+(\nbigi^{\circ},\theta^u)}
\bigl(
 (L_2)'_{\vecJ,<0|\theta^u}
 \oplus
 (L_2)'_{\vecJ,>0|\theta^u}
 \bigr).
\]
We have the decomposition
\[
L_{2|\theta^u}
=
 W_{\theta^u}
 \oplus
 (L_2)'_{(\vecJ_1)_+,0|\theta^u}.
\]
An element of $L_{2|\theta^u}$ is denoted by $(s_1,s_2)$
according to the decomposition.
It is easy to observe the following.
\begin{lem}
\label{lem;25.2.17.3}
For any $s_2\in (L_2)'_{|(\vecJ_1)_+,0|\theta^u}$,
there exists $t_2\in W_{\theta^u}$ such that
$M^{-1}_{L_2}(0,s_2)=(t_2,M^{-1}_{L_2,0}(s_2))$. 
\hfill\qed
\end{lem}

\subsection{The recovery of the Stokes structure of $L$}

Let us observe that
we can recover the Stokes structure $\vecnbigf$ on $L$
from the induced Stokes structure $\vecnbigf$ on $\Ltilde_1$.

\subsubsection{Recovery of $(L^{<0},\vecnbigftilde)$}

Let $\theta^u\in S_0(\nbigi^{\circ})$.
Let $\vecJ_1\in T(\nbigi^{\circ})$
such that $\vartheta^{\vecJ_1}_r=\theta^u$.
We have
\begin{multline}
 H^0(\real,\Ltilde_1)=
 \Ltilde_{1|\theta^u}= 
\\
 \bigoplus_{\theta^u\in\vecJ_+}
 \Bigl(
 H^0(\vecJ,\gbigl^{\gbigf}_!(V)_{\vecJ,<0})
 \oplus
 H^0(\vecJ_+,\gbigl^{\gbigf}_!(V)_{\vecJ_+,>0})
 \Bigr)
 \oplus
 H^0(\vecJ_{1+},(L_1)_{\vecJ_{1+},0}).
\end{multline}
The morphism
$a_3\circ u_2:H^0(\real,\Ltilde_1)\to H^0(\real,L)$
induces an isomorphism
\[
H^0(\vecJ,\gbigl^{\gbigf}_!(V)_{\vecJ,<0})
\simeq
H^0(\nu_0^-(\vecJ),L_{\nu_0^-(\vecJ),<0}).
\]

Note that
$L^{<0}
 =\bigoplus_{J\in T(\nbigi)}
 \iota_{J!}(L_{J,<0})$.
We can recover
$\iota_{J!}(L_{J,<0})\subset L^{<0}$
from $H^0(J,L_{J,<0})\subset H^0(\real,L)$.
According to Corollary \ref{cor;25.2.22.11},
the Stokes filtrations $\vecnbigftilde$ on
$H^0(\nu_0^-(\vecJ),L_{\nu_0^-(\vecJ),<0})$
is recovered from
$\bigl(
H^0(\vecJ,\gbigl^{\gbigf}_!(V)_{\vecJ,<0}),
\vecnbigf
\bigr)$,
where we use the identification
of the index sets in (\ref{eq;24.3.14.40}).
In this way,
$(L^{<0},\vecnbigftilde)$ is recovered
from $(\Ltilde_1,\vecnbigf)$.

\subsubsection{Recovery of $L^{\leq 0}$ and the positive parts}

We set $\theta_1=\vartheta^{\vecI(\theta^u)}_r=\theta^u-\pi/2$.
We set $J_1=\nu_0^-(\vecJ_1)$
and $J_2=J_1+\omega^{-1}\pi$.
We have
\[
 H^0(\real,L)=L_{|\theta_1}
 =\bigoplus_{J_1\leq J<J_2}
 \Bigl(
 H^0(J,L_{J,<0})
 \oplus
 H^0(J_+,L_{J_+,>0})
 \Bigr)
 \oplus
 H^0(J_{1+},L_{J_{1+},0}).
\]
We have
\[
 H^0(\real,\Ltilde_1)
 =\Ltilde_{1|\theta^u}
 =
 \!\!\!
 \bigoplus_{J\in \gbigm_1(\nbigi,\vecI(\theta^u)_-)}
 \!\!\!
 \Image B_{J_+,\theta^u}
 \oplus
 \!\!\!
 \bigoplus_{J\in \gbigm_2(\nbigi,\vecI(\theta^u)_-)}
 \!\!\!
 \Image \Abb^{\rd}_{J,\theta^u}
 \oplus
 H^0(\real,L_1).
\]
We set
\[
 K_{\theta^u}:=
 \bigoplus_{J_1-\pi\leq J<J_1}
 \Image B_{J_+,\theta^u}
 \oplus
 \!\!\!
\bigoplus_{J_1<J<J_2+\pi}
 \!\!\!
 \Image \Abb^{\rd}_{J,\theta^u}
 \oplus
 H^0(\real,L_1)
\subset H^0(\real,\Ltilde_1).
\]
Note that
\[
 H^0(\real,\Ltilde_1)
 =K_{\theta^u}\oplus
 \bigoplus_{J_1\leq J<J_2}
 \Image B_{J_+,\theta^u}
 \oplus
 \Image \Abb^{\rd}_{J_1,\theta^u}.
\]
\begin{lem}
\label{lem;25.2.22.20}
If $J_1-2\pi< J<J_2$,
then 
$\Image \Abb^{\rd}_{J,\theta^u-2\pi}
 \subset K_{\theta^u}$.
\end{lem}
\pf
If $J_1-2\pi\leq J<J_2-\pi$,
we have $J+2\pi\in \gbigm_2(\nbigi,\vecI(\theta^u)_-)$.
Hence, the claim follows from Lemma \ref{lem;24.4.5.1}.
For $J_2-\pi\leq J<J_2-\pi+a\pi$ $(0\leq a<1)$,
we obtain 
$\Image \Abb^{\rd}_{J,\theta^u-2\pi}
 \subset K_{\theta^u}$.
by using (\ref{eq;24.3.16.2})
and an easy induction on $a$.
\hfill\qed

\vspace{.1in}
Let $h:H^0(\real,L)\to H^0(\real,\Ltilde_1)$ denote
the morphism induced by $a_1$ and $u_1$.
\begin{lem}
\label{lem;25.2.22.21}
\[
 h\Bigl(
 \bigoplus_{J_1<J<J_2}
 H^0(J,L_{J,<0})
 \oplus
 H^0(J_{1+},L_{J_{1+},0})
 \Bigr)
 \subset K_{\theta^u}.
\]
\end{lem}
\pf
It is enough to study the case
$\Ltilde_1=\gbigl^{\gbigf}_!(V)$.
Let $J_1<J<J_2$.
Because
$\Abb_{\infty,\theta^u}(v)
=\Abb^{\rd}_{J,\theta^u}(v)
-\Abb^{\rd}_{J,\theta^u-2\pi}(v)$
for $v\in H^0(J,L_{J,<0})$,
we obtain
$\Abb_{\infty,\theta^u}(v)\in K_{\theta^u}$
by Lemma \ref{lem;25.2.22.20}.
By using (\ref{eq;25.2.22.1}) and Lemma \ref{lem;25.2.22.20},
we obtain 
$A_{\infty,\theta^u}(v)\in K_{\theta^u}$ for
$v\in H^0(J_{1+},L_{J_{1+},0})$.
\hfill\qed

\begin{lem}
\label{lem;25.2.22.23}
We have
$h\bigl(
H^0(J_{1},L_{J_1,<0})
\bigr)
\subset
K_{\theta^u}\oplus
\Image \Abb^{\rd}_{J_1,\theta^u}$.
The induced map
\[
 H^0(J_{1},L_{J_1,<0})
 \lrarr
\bigl(
 K_{\theta^u}\oplus
 \Image \Abb^{\rd}_{J_1,\theta^u}
 \bigr)\Big/K_{\theta^u}
 \simeq
 \Image \Abb^{\rd}_{J_1,\theta^u}
\]
equals $\Abb^{\rd}_{J_1,\theta^u}$.
\end{lem}
\pf
The first claim is similar to Lemma \ref{lem;25.2.22.21}.
The second claim follows from the construction.
\hfill\qed

\begin{lem}
\label{lem;25.2.22.22}
For any $J_1\leq J<J_2$,
we obtain 
\[
 h\bigl(
 H^0(J_{+},L_{J_{+},>0})
 \bigr)
 \subset
 K_{\theta^u}
 \oplus
 \bigoplus_{J_1\leq J'\leq J}
 \Image B_{J'_+,\theta^u}.
\]
Moreover, the induced map
$H^0(J_+,L_{J_+,>0})\to
\Image B_{J_+,\theta^u}$
equals $B_{J_+,\theta^u}$.
\end{lem}
\pf
It follows from (\ref{eq;25.2.22.1})
and Lemma \ref{lem;25.2.22.20}.
\hfill\qed

\vspace{.1in}

By Lemma \ref{lem;25.2.22.21},
Lemma \ref{lem;25.2.22.23}
and Lemma \ref{lem;25.2.22.22},
under the isomorphism
$L_{\theta_1}=H^0(\real,L)$,
we have
\[
 L^{\leq 0}_{\theta_1}
 =h^{-1}\bigl(
 K_{\theta^u}
 \bigr).
\]
Let $\theta_1'\in S_0(\nbigi)$
determined by 
$\openopen{\theta_1'}{\theta_1}\cap S_0(\nbigi)=\emptyset$.
For any $\theta\in \openopen{\theta_1'}{\theta_1}$,
under the isomorphism
$L_{\theta}=H^0(\real,L)$,
by Lemma \ref{lem;25.2.22.21},
Lemma \ref{lem;25.2.22.23}
and Lemma \ref{lem;25.2.22.22},
we have
\[
 (L^{\leq 0})_{\theta}
 =h^{-1}\Bigl(
 K_{\theta^u}
 \oplus
 \Image\Abb^{\rd}_{J_1,\theta^u}
 \Bigr).
\]
Let $\theta_1''\in S_0(\nbigi)$
determined by 
$\openopen{\theta_1}{\theta_1''}\cap S_0(\nbigi)=\emptyset$.
For any $\theta\in\openopen{\theta_1}{\theta_1''}$,
under the isomorphism
$L_{\theta}=H^0(\real,L)$,
by Lemma \ref{lem;25.2.22.21},
Lemma \ref{lem;25.2.22.23}
and Lemma \ref{lem;25.2.22.22},
we have
\[
 (L^{\leq 0})_{\theta}
 =h^{-1}\Bigl(
 K_{\theta^u}
 \oplus
 \Image B_{J_1,\theta^u}
 \Bigr).
\]
Thus,
the constructible subsheaf
$L^{\leq 0}\subset L$
is recovered.

Note that
\[
L/L^{\leq 0}\simeq
\bigoplus_{J\in T(\nbigi)}
\iota_{\Jbar\ast}(L_{\Jbar,L_{\Jbar,>0}}).
\]
The Stokes filtrations
$\nbigftilde$
on $H^0(\nu_0^+(\Jbar),L_{\nu_0^+(\Jbar),>0})$
are recovered from the Stokes filtrations
$\vecnbigf$
on $H^0(\vecJbar,\gbigl^{\gbigf}_!(V)_{\Jbar,>0})$
and the isomorphism
in Lemma \ref{lem;25.2.22.22}.
By Proposition \ref{prop;25.3.16.11},
we can recover $(L,\vecnbigftilde)$
from $(\Ltilde_1,\vecnbigf)$
with morphisms (\ref{eq;25.3.16.10}).

\section{Local Fourier transforms of Stokes structure
from $0$ to $\infty$}
\label{subsection;24.4.5.100}

To describe 
$(\gbigl^{\gbigf}_{\star}(V),\vecnbigf)$
it is convenient to introduce the local Fourier transform
of a Stokes structure.

\subsection{$2\pi\seisuu$-equivariant
local system $\gbigq^0_!(L,\vecnbigftilde)_{\real}$}
\label{subsection;24.3.21.3}

We consider the vector space
\begin{equation}
\label{eq;24.3.21.1}
 H^0(\real,L)
 \oplus
 \bigoplus_{J\in T(\nbigi)}
 H^0(J,L_{J,<0}).
\end{equation}
An element of $v\in H^0(J,L_{J,<0})$ is denoted as
a pair $\langle J,v\rangle$.

Let $\gbigq^0_!(L,\vecnbigftilde)$ denote
the quotient space of
(\ref{eq;24.3.21.1}) by the equivalence relation
generated by the following
(see Lemma \ref{lem;24.4.5.1}):
\index{vector space \mbox{$\gbigq^0_!(L,\vecnbigftilde)$}}
\[
 \langle J,v\rangle-\langle J+2\pi,(\Tbb^{\ast})^{-1}(v)\rangle
 \sim
 \rho_J(v)
 \in H^0(\real,L).
\]
Here, $\rho_{J}:H^0(J,L_{J,<0})\to H^0(\real,L)$
denote the natural inclusions.

Let $\Tbb^{\ast}_{\gbigq^0,!}$ denote the automorphism on
$\gbigq^0_!(L,\vecnbigftilde)$
induced by 
$M$ on $H^0(\real,L)$,
and the maps
(see Lemma \ref{lem;24.4.5.2} and Lemma \ref{lem;24.4.5.1}):
\[
\Tbb^{\ast}:H^0(J+2\pi,L_{J+2\pi,<0})\simeq H^0(J,L_{J,<0}),
\quad
\langle J+2\pi,v\rangle
 \longmapsto
 \langle J,\Tbb^{\ast}(v)\rangle.
\]

Let $\gbigq^0_!(L,\vecnbigftilde)_{\real}$ denote
the local system on $\real$ induced by $\gbigq^0_!(L,\vecnbigftilde)$.
\index{local system \mbox{$\gbigq^0_!(L,\vecnbigftilde)_{\real}$}}
We naturally identify
$H^0(\real,\gbigq^0_!(L,\vecnbigftilde)_{\real})$
with $\gbigq^0_!(L,\vecnbigftilde)$.
There exists the $2\pi\seisuu$-action
on $\gbigq^0_!(L,\vecnbigftilde)_{\real}$
such that
the pull back
$\Tbb^{\ast}:H^0(\real,\gbigq^0_{!}(L,\vecnbigftilde)_{\real})
\simeq
 H^0(\real,\gbigq^0_{!}(L,\vecnbigftilde)_{\real})$
equals $\Tbb^{\ast}_{\gbigq^0,!}$.

\begin{prop}
There exists the isomorphism
of $2\pi\seisuu$-equivariant local systems
$\gbigq^0_{!}(L,\vecnbigftilde)_{\real}
\simeq
\gbigl^{\gbigf}_!(V)$ 
induced by
$\Abb^{\rd}_{\infty,\theta^u}$
and
$\Abb^{\rd}_{J,\theta^{u}}$
$(\theta^u\in\real,J\in T(\nbigi))$. 
\end{prop}
\pf
It follows from Lemma \ref{lem;24.4.5.2},
Lemma \ref{lem;24.4.5.1} and
Lemma \ref{lem;18.4.22.121}.
\hfill\qed

\subsubsection{Another expression and the monodromy}
\label{subsection;24.2.21.2}

Fix $u(0)\in\cnum^{\ast}$ and $\theta^{u}_0\in\real$
such that
$\theta^u_0=\arg(u(0))$.
For $J\in \gbigt(\nbigi,\theta^u_0)$,
we obtain the constant $2\pi\seisuu$-equivariant local system
$H^0(J,L_{J,<0})_{\real}$ on $\real$
induced by $H^0(J,L_{J,<0})$.
We obtain the following exact sequence
of $2\pi\seisuu$-equivariant local systems:
\[
 0\lrarr L
 \lrarr
 \gbigq^0_{!}(L,\vecnbigftilde)_{\real}
  \lrarr
 \bigoplus_{J\in \gbigt(\nbigi,\theta^u_0)}
 H^0(J,L_{J,<0})_{\real}
 \lrarr 0.
\]
There exists the natural isomorphism
\[
 \gbigq^0_{!}(L,\vecnbigftilde)_{\real|\theta^{u}_0}
 \simeq
 H^0(\real,L)\oplus
 \bigoplus_{J\in\gbigt(\nbigi,\theta^u_0)}
 H^0(J,L_{J,<0})
\]
under which the monodromy of
$\gbigq^0_{!}(L,\vecnbigftilde)_{\real}$
is described as
\[
\Bigl(
w,\sum_{J}v_J
\Bigr)
\longmapsto
\Bigl(
M(w)+\sum_JM\circ\rho_J(v_J),
\sum_Jv_J
\Bigr).
\]

\subsection{$2\pi\seisuu$-equivariant
local system $\gbigq^0_{\ast}(L,\vecnbigftilde)_{\real}$}
\label{subsection;24.3.21.4}

We consider the vector space
\begin{equation}
 \label{eq;24.3.21.2}
  \bigoplus_{\pm}\bigoplus_{J\in T(\nbigi)}
  H^0(\real,L)
  \oplus
  \bigoplus_{J\in T(\nbigi)}
  H^0(J,L_{J,>0}).
\end{equation}
An element of $w\in H^0(\real,L)$ corresponding
to the $(\kappa,J)$-component
$((\kappa,J)\in \{\pm\}\times T(\nbigi))$
is denoted as $\langle J_{\kappa},w\rangle^{\mg}$.
An element of $v\in H^0(J,L_{J,>0})$
is denoted as a pair $\langle J,v\rangle^{\mg}$.

Let $\gbigq^0_{\ast}(L,\vecnbigftilde)$
denote the quotient space of
(\ref{eq;24.3.21.2}) by the equivalence relation
generated by the following
(see Lemma \ref{lem;25.2.11.20},
Lemma \ref{lem;24.3.25.10} and Corollary \ref{cor;24.4.5.4}).
\index{vector space $\gbigq^0_{\ast}(L,\vecnbigftilde)$}
\begin{itemize}
 \item  $\langle J+2\pi,v\rangle^{\mg}
	\sim\langle J,\Tbb^{\ast}(v)\rangle^{\mg}$
	for any $J\in T(\nbigi)$
	and $v\in H^0(J+2\pi,L_{J+2\pi,>0})$.
 \item $\langle J_-,w\rangle^{\mg}-\langle J_+,w\rangle^{\mg}
       \sim\langle J,R_J(w)\rangle^{\mg}$
       for any $J\in T(\nbigi)$
       and $w\in H^0(\real,L)$.
 \item $\langle J_{1-},w\rangle^{\mg}
       -\langle J_{2-},w\rangle^{\mg}
       \sim
       \sum_{J_1\leq J'<J_2}
       \langle J',R_{J'}(w)\rangle^{\mg}$
       for any
       $J_1\leq J_2$ in $T(\nbigi)$
       and $w\in H^0(\real,L)$.
\end{itemize}

Let $\Tbb^{\ast}_{\gbigq^0,\ast}$ denote the automorphism of
$\gbigq^0_{\ast}(L,\vecnbigftilde)$
induced by 
\[
 \langle (J+2\pi)_{\pm},w\rangle^{\mg}
 \longmapsto
 \langle J_{\pm},M(w)
 \rangle^{\mg},
 \quad\quad
 \langle J,v\rangle^{\mg}
 \longmapsto
 \langle J,v\rangle^{\mg}.
\]
(See Remark \ref{rem;25.2.11.21} and Lemma \ref{lem;24.4.5.3}.)
Let $\gbigq^0_{\ast}(L,\vecnbigftilde)_{\real}$
denote the local system on $\real$
induced by $\gbigq^0_{\ast}(L,\vecnbigftilde)$.
\index{local system $\gbigq^0_{\ast}(L,\vecnbigftilde)_{\real}$}
We naturally identify
$H^0(\real,\gbigq^0_{\ast}(L,\vecnbigftilde)_{\real})$
with $\gbigq^0_{\ast}(L,\vecnbigftilde)$.
There exists the $2\pi\seisuu$-action
on $\gbigq^0_{\ast}(L,\vecnbigftilde)_{\real}$
such that the pull back
$\Tbb^{\ast}:H^0(\real,\gbigq^0_{\ast}(L,\vecnbigftilde)_{\real})
\simeq
 H^0(\real,\gbigq^0_{+,\ast}(L,\vecnbigftilde)_{\real})$
equals $\Tbb^{\ast}_{\gbigq^0,\ast}$.

\begin{prop}
There exists the isomorphism
of $2\pi\seisuu$-equivariant local systems
$\gbigq^0_{\ast}(L,\vecnbigftilde)_{\real}
 \simeq
\gbigl^{\gbigf}_{\ast}(V)$ 
induced by
$\Abb^{\mg,J_{\pm}}_{\infty,\theta^u}$
and
$\BB^{\mg}_{J,u}$
$(\theta^u\in\real,J\in T(\nbigi))$.
\end{prop}
\pf
It follows from Lemma \ref{lem;24.4.5.3},
Lemma \ref{lem;24.3.25.10}
and Corollary \ref{cor;24.4.5.4}.
\hfill\qed

\subsubsection{Another expression and the monodromy}
\label{subsection;24.3.26.2}
Fix $u(0)\in\cnum^{\ast}$ and
$\theta^u_0\in\real$ such that
$\theta^u_0=\arg(u(0))$.
For $J\in \gbigt(\nbigi,\theta^u_0)$,
we obtain the constant $2\pi\seisuu$-equivariant local system
$H^0(J,L_{J,>0})_{\real}$ on $\real$
induced by $H^0(J,L_{J,>0})$.
We obtain the following exact sequence
of $2\pi\seisuu$-equivariant local systems on $\real$:
\[
 0\lrarr
\bigoplus_{J\in \gbigt(\nbigi,\theta^u_0)}
 H^0(J,L_{J,>0})_{\real}
 \lrarr
\gbigq^0_{\ast}(L,\vecnbigftilde)_{\real}
 \lrarr
 L
 \lrarr 0.
\]
Choosing $J_0\in T(\nbigi)$,
and considering
$\langle (J_0)_+,v\rangle^{\mg}$
for $v\in H^0(\real,L)$,
we obtain the isomorphism
\begin{equation}
\label{eq;24.2.15.3}
\gbigq^0_{\ast}(L,\vecnbigftilde)_{\real|\theta^{u}_0}
 \simeq
 H^0(\real,L)
 \oplus
 \bigoplus_{J\in\gbigt(\nbigi,\theta^u_0)}
 H^0(J,L_{J,>0}),
\end{equation}
under which the monodromy of
$\gbigq^0_{\ast}(L,\vecnbigftilde)_{\real}$
is described as
\[
\Bigl(
w,\sum_{J}v_J
\Bigr)
\longmapsto
\Bigl(
M(w),
\sum_J\bigl(v_J+R_J(w)\bigr)
\Bigr).
\]
(See \S\ref{subsection;24.2.18.1}
for the maps $R_J:H^0(\real,L)\to H^0(J,L_{J,>0})$.)

\subsection{Morphisms}
\label{subsection;24.4.4.2}
Let
$F_{\gbigq^0}\colon
\gbigq^0_!(L,\vecnbigftilde)
\to\gbigq^0_{\ast}(L,\vecnbigftilde)$
be the morphism obtained as follows
(see Lemma \ref{lem;24.4.5.10}):
\index{morphism $F_{\gbigq^0}$}
\begin{itemize}
 \item  For any $J\in T(\nbigi)$ and $v\in H^0(J,L_{J,<0})$,
\[
 \langle J,v\rangle
 \longmapsto
 \langle J_+,\rho_J(v)\rangle^{\mg}
=\langle J_-,\rho_J(v)\rangle^{\mg}.
\]
\item For any $w\in H^0(\real,L)$,
\[
 w\longmapsto \langle J_{1+},w-M^{-1}(w)\rangle^{\mg}
 +\sum_{J_1-2\pi<J'\leq J_1}
 \langle J',R_{J'}(w)\rangle^{\mg}.
\]
      The right hand side is independent of $J_1\in T(\nbigi)$
      in $\gbigq^0_{\ast}(L,\vecnbigftilde)$.
\end{itemize}
It induces the morphism of $2\pi\seisuu$-equivariant local systems
$F_{\gbigq^0}:\gbigq^0_{!}(L,\vecnbigftilde)_{\real}
\to \gbigq^0_{\ast}(L,\vecnbigftilde)_{\real}$.

Let $H^0(\real,L)\to \gbigq^0_!(L,\vecnbigftilde)$
denote the morphism induced by
the inclusion of $H^0(\real,L)$ into the space (\ref{eq;24.3.21.1}).
Let $\gbigq^0_{\ast}(L,\vecnbigftilde)\to H^0(\real,L)$
denote the morphism induced by
the projection of the space (\ref{eq;24.3.21.2})
onto $H^0(\real,L)$.
They induce the morphisms
of $2\pi\seisuu$-equivariant local systems
$d_1:L\to\gbigq^0_{!}(L,\vecnbigftilde)_{\real}$
and
$d_2:\gbigq^0_{\ast}(L,\vecnbigftilde)_{\real}\to L$.
\begin{prop}
We obtain the following commutative diagram:
\begin{equation}
\label{eq;25.2.22.30}
 \begin{CD}
  L @>{d_1}>>
  \gbigq^0_{!}(L,\vecnbigftilde)_{\real}
  @>{F_{\gbigq^0}}>>
  \gbigq^0_{\ast}(L,\vecnbigftilde)_{\real}
  @>{d_2}>>
  L \\
  @V{\simeq}VV @V{\simeq}VV @V{\simeq}VV @V{\simeq}VV \\
  \gbigl^{\gbigf}_!(V^{\reg})
  @>>>
  \gbigl^{\gbigf}_!(V)
  @>>>
  \gbigl^{\gbigf}_{\ast}(V)
  @>>>
  \gbigl^{\gbigf}_{\ast}(V^{\reg})
 \end{CD}
\end{equation}
\end{prop}
\pf
We obtain of the commutativity of
the middle square from Lemma \ref{lem;24.4.5.10}.
The commutativity of the left and right squares are clear
by the construction.
\hfill\qed

\vspace{.1in}
Let $M_{\gbigq^0_!}$
and $M_{\gbigq^0_{\ast}}$
denote the monodromy automorphisms
of $\gbigq^0_!(L,\vecnbigftilde)$
and $\gbigq^0_{\ast}(L,\vecnbigftilde)$,
respectively.
We have
\[
d_2\circ F_{\gbigq^0}\circ d_1
=\id-M^{-1},
\quad
d_1\circ d_2\circ F_{\gbigq^0}
=\id-M_{\gbigq^0_!}^{-1},
\quad
F_{\gbigq^0}\circ d_1\circ d_2
=\id-M_{\gbigq^0_{\ast}}^{-1}.
\]

\subsection{Stokes structure of $\gbigl^{\gbigf}_{!}(L,\vecnbigftilde)$}
\label{subsection;24.4.4.3}

Let $\vecJ\in T(\nbigi^{\circ})$.
We define the map
$\vecA_{\vecJ}:H^0(\nu_0^-(\vecJ),L_{\nu_0^-(\vecJ),<0})
\lrarr
\gbigq^0_!(L,\vecnbigftilde)$
by \index{map $\vecA_{\vecJ}$}
\[
 \vecA_{\vecJ}(v)
:=\langle \nu_0^-(\vecJ),v\rangle.
\]
For any $v\in H^0(\real,\nbigt_{\omega}(L))$
and any $J\in T(\nbigi)$,
we obtain
$v_{J_{\pm}}\in H^0(J_{\pm},L_{J_{\pm},0})\subset H^0(\real,L)$.
Then, we define
$\vecB^{\vecJ_{\pm}}_{\infty}:
H^0(\real,\nbigt_{\omega}(L))
\lrarr
\gbigq^0_!(L,\vecnbigftilde)$ by
\index{maps $\vecB^{\vecJ_{\pm}}_{\infty}$}
\[
 \vecB^{\vecJ_+}_{\infty}(v):=
 v_{\nu_0^-(\vecJ)_+}
 -\sum_{\nu_0^-(\vecJ)-2\pi<J'\leq \nu_0^-(\vecJ)}
 \langle J',\nbigp_{J'}(v)\rangle,
\]
\[
  \vecB^{\vecJ_-}_{\infty}(v):=
 v_{\nu_0^-(\vecJ)_-}
 -\sum_{\nu_0^-(\vecJ)-2\pi\leq J'< \nu_0^-(\vecJ)}
 \langle J',\nbigp_{J'}(v)\rangle.
\]
(See \S\ref{subsection;25.2.11.30}.)
We define
$\vecB_{\vecJ_{\pm}}:
H^0(\nu_0^+(\vecJ),L_{\nu_0^+(\vecJ),>0})
\lrarr
\gbigq^0_!(L,\vecnbigftilde)$
by
\index{maps $\vecB_{\vecJ_{\pm}}$}
\begin{multline}
 \vecB_{\vecJ_-}(v)
 :=
\!\!\!
 \sum_{\nu_0^+(\vecJ)-\omega^{-1}\pi\leq J'<\nu_0^+(\vecJ)}
\!\!\!\!
\bigl\langle J',\nbigrtilde^{\nu_0^+(\vecJ)_-}_{J'}(v)\bigr\rangle
\\
 -
\!\!\!\!\!\!
 \sum_{\nu_0^+(\vecJ)\leq J'<\nu_0^+(\vecJ)+\pi}
\!\!\!\!\!\!
\bigl\langle J',\nbigrtilde^{\nu_0^+(\vecJ)_-}_{J'}(v)\bigr\rangle
 -
 \!\!\!\!\!\!
 \sum_{\nu_0^+(\vecJ)-\omega^{-1}\pi\leq J'<\nu_0^+(\vecJ)-\pi}
 \!\!\!\!\!
\bigl\langle J'+2\pi,
 (\Tbb^{\ast})^{-1}
 \nbigrtilde^{\nu_0^+(\vecJ)_-}_{J'}(v)\bigr\rangle,
\end{multline}
\begin{multline}
 \vecB_{\vecJ_+}(v)
 :=-
\!\!\!\!\!\!
 \sum_{\nu_0^+(\vecJ)<J'\leq \nu_0^+(\vecJ)+\omega^{-1}\pi}
\!\!\!\!\!\!
 \bigl\langle J'+2\pi,(\Tbb^{\ast})^{-1}
 \nbigrtilde^{\nu_0^+(\vecJ)_+}_{J'}(v)\bigr\rangle
\\
 +
\!\!\!\!\!\!
 \sum_{\nu_0^+(\vecJ)-\pi< J'\leq \nu_0^+(\vecJ)}
\!\!\!\!\!\!
 \bigl\langle J'+2\pi,
 (\Tbb^{\ast})^{-1}
 \nbigrtilde^{\nu_0^+(\vecJ)_+}_{J'}(v)\bigr\rangle
 +
 \!\!\!\!\!\!\!\!
 \sum_{\nu_0^+(\vecJ)+\pi<J'\leq \nu_0^+(\vecJ)+\omega^{-1}\pi}
 \!\!\!\!\!\!\!\!
 \bigl\langle J',
  \nbigrtilde^{\nu_0^+(\vecJ)_+}_{J'}(v)\bigr\rangle.
\end{multline}
(See \S\ref{subsection;25.2.11.31}.)
By Theorem \ref{thm;24.3.15.10},
we obtain the following.
\begin{prop}
\label{prop;24.4.6.2}
Let $\theta^u\in\real$.
Choose $\vecJ_1\in\gbigm_{-}(\nbigi^{\circ},\theta^u)$.
Then,
$\vecA_{\vecJ}$, $\vecB_{\vecJ_{-}}$
$(\vecJ\in \gbigm_-(\nbigi^{\circ},\theta^u))$
and $\vecB^{\vecJ_{1-}}_{\infty}$
induce the isomorphism of the vector spaces:
\begin{multline}
\label{eq;24.4.6.1}
 \!\!\!\!\!\!\!\!
 \bigoplus_{\vecJ\in\gbigm_-(\nbigi^{\circ},\theta^u)}
 \!\!
 \Bigl(
 H^0(\nu_0^-(\vecJ),L_{\nu_0^-(\vecJ),<0})
\oplus 
 H^0(\nu_0^+(\vecJ),L_{\nu_0^+(\vecJ),>0})
 \Bigr)
 \oplus
 H^0(\real,\nbigt_{\omega}(L))
 \\
\simeq\gbigq^0_!(L,\vecnbigftilde)\simeq
 \gbigl^{\gbigf}_!(V)_{|\theta^u}.
 \end{multline}
Moreover, 
if we consider the filtrations $\nbigf^{\prime\theta^u}$
on the spaces
$H^0(\nu_0^-(\vecJ),L_{\nu_0^-(\vecJ),<0})$
and
$H^0(\nu_0^+(\vecJ),L_{\nu_0^+(\vecJ),>0})$
defined in {\rm\S\ref{subsection;24.4.4.1}},
the trivial filtration on $H^0(\real,\nbigt_{\omega}(L))$
indexed by $0$,
and the Stokes filtration $\nbigf^{\theta^u}$
on $\gbigl^{\gbigf}_!(V)$,
then {\rm(\ref{eq;24.4.6.1})}
induces the isomorphism of filtered vector spaces.

We also obtain a similar isomorphism
by choosing
$\vecJ_1\in\gbigm_+(\nbigi^{\circ},\theta^u)$
and using
$\vecA_{\vecJ}$,
$\vecB_{\vecJ_+}$
$(\vecJ\in \gbigm_+(\nbigi^{\circ},\theta^u))$
and $\vecB^{\vecJ_{1+}}_{\infty}$.
\hfill\qed
\end{prop}

By Theorem \ref{thm;24.3.15.10},
we also obtain the following.
\begin{prop}
\label{prop;24.4.4.21}
Under the isomorphism
$\gbigq^0_!(L,\vecnbigftilde)\simeq H^0(\real,\gbigl^{\gbigf}_!(V))$,
we have
\[
\Image\vecA_{\vecJ}
=H^0(\vecJ,\gbigl^{\gbigf}_{!}(V)_{\vecJ,<0}),
\quad
\Image \vecB^{\vecJ_{\pm}}_{\infty}
=H^0(\vecJ_{\pm},\gbigl^{\gbigf}_{!}(V)_{\vecJ,0}),
\]
\[
\Image \vecB_{\vecJ_{\pm}}
=H^0(\vecJ_{\pm},\gbigl^{\gbigf}_{!}(V)_{\vecJ,>0}).
\]
\hfill\qed 
\end{prop}

\subsection{Stokes structure of $\gbigl^{\gbigf}_{\ast}(L,\vecnbigftilde)$}
\label{subsection;24.4.4.20}

For $\vecJ\in T(\nbigi^{\circ})$,
we obtain
\[
\vecA^{\mg}_{\vecJ}:
H^0(\nu_0^-(\vecJ),L_{\nu_0^-(\vecJ),<0})
\lrarr
\gbigq^0_{\ast}(L,\vecnbigftilde),
\]
\[
\vecB^{\mg}_{\vecJ_{\pm}}:
H^0(\nu_0^+(\vecJ),L_{\nu_0^+(\vecJ),>0})
\lrarr
\gbigq^0_{\ast}(L,\vecnbigftilde)
\]
as the composition of $\vecA_{\vecJ}$ and $\vecB_{\vecJ}$
in \S\ref{subsection;24.4.4.3}
and the morphism
$\gbigq^0_!(L,\vecnbigftilde)\to\gbigq^0_{\ast}(L,\vecnbigftilde)$.
We define
$\vecB^{\mg,\vecJ_{\pm}}_{\infty}:
H^0(\real,\nbigt_{\omega}(L))
\to\gbigq^0_{\ast}(L,\vecnbigftilde)$
by
\index{maps $\vecB^{\mg,\vecJ_{\pm}}_{\infty}$}
\[
 \vecB^{\mg,\vecJ_{\pm}}_{\infty}(v)
 =
 \bigl\langle
 \nu_0^-(\vecJ)_{\pm},
 v_{\nu_0^-(\vecJ)_{\pm}}
 \bigr\rangle^{\mg}.
\]
(See \!\ref{subsection;24.4.4.3}
for $v_{J_{\pm}}$.)
By Theorem \ref{thm;24.3.15.10},
we obtain the following.
\begin{prop}
\label{prop;24.4.6.20}
Let $\theta^u\in\real$.
Choose $\vecJ_1\in\gbigm_{-}(\nbigi^{\circ},\theta^u)$.
Then,
$\vecA^{\mg}_{\vecJ}$, $\vecB^{\mg}_{\vecJ_{-}}$
$(\vecJ\in \gbigm_-(\nbigi,\theta^u))$
and $\vecA^{\mg,\vecJ_{1-}}_{\infty}$
induce the isomorphism of the vector spaces:
 \begin{multline}
\label{eq;24.4.6.10}
\!\!\! \bigoplus_{\vecJ\in\gbigm_-(\nbigi^{\circ},\theta^u)}\!\!
\Bigl(
 H^0(\nu_0^-(\vecJ),L_{\nu_0^-(\vecJ),<0})
\oplus 
 H^0(\nu_0^+(\vecJ),L_{\nu_0^+(\vecJ),>0})
 \Bigr)
 \oplus
 H^0(\real,\nbigt_{\omega}(L))
 \\
\simeq\gbigq^0_{\ast}(L,\vecnbigftilde)
 \simeq
 \gbigl^{\gbigf}_{\ast}(V)_{|\theta^u}.
 \end{multline}
Moreover, if we consider the filtrations $\nbigf^{\prime\theta^u}$
on the spaces
$H^0(\nu_0^-(\vecJ),L_{\nu_0^-(\vecJ),<0})$
and
$H^0(\nu_0^+(\vecJ),L_{\nu_0^+(\vecJ),>0})$
defined in {\rm\S\ref{subsection;24.4.4.1}},
the trivial filtration on $H^0(\real,\nbigt_{\omega}(L))$
indexed by $0$,
and the Stokes filtration $\nbigf^{\theta^u}$
on $\gbigl^{\gbigf}_{\ast}(V)_{|\theta^u}$,
then {\rm(\ref{eq;24.4.6.10})}
induces an isomorphism of filtered vector spaces.

We also obtain a similar isomorphism
by choosing $\vecJ_1\in\gbigm_+(\nbigi^{\circ},\theta^u)$
and by using
$\vecA^{\mg}_{\vecJ}$,
$\vecB^{\mg}_{\vecJ_{-}}$ 
$(\vecJ\in \gbigm_+(\nbigi^{\circ},\theta^u))$
and $\vecB^{\mg,\vecJ_{1+}}_{\infty}$.
\hfill\qed
\end{prop}
By Theorem \ref{thm;24.3.15.10},
we also obtain the following.
\begin{prop}
\label{prop;24.4.4.22}
Under the isomorphism
$\gbigq^0_{\ast}(L,\vecnbigftilde)\simeq H^0(\real,\gbigl^{\gbigf}_{\ast}(V))$,
we have
\[
\Image\vecA^{\mg}_{\vecJ}
=H^0(\vecJ,\gbigl^{\gbigf}_{\ast}(V)_{\vecJ,<0}),
\quad
\Image \vecB^{\mg,\vecJ_{\pm}}_{\infty}
=H^0(\vecJ_{\pm},\gbigl^{\gbigf}_{\ast}(V)_{\vecJ,0}),
\]
\[
\Image \vecB^{\mg}_{\vecJ_{\pm}}
=H^0(\vecJ_{\pm},\gbigl^{\gbigf}_{\ast}(V)_{\vecJ,>0}).
\]

\hfill\qed 
\end{prop}

\subsection{Isomorphisms}

For any $\theta^u\in\real$,
we define the filtrations $\nbigf^{\theta^u}$
on $\gbigq^0_{\star}(L,\vecnbigftilde)
=\gbigq^0_{\star}(L,\vecnbigftilde)_{\real|\theta^u}$
$(\star=!,\ast)$
indexed by
$(\nbigitilde^{\circ},\leq_{\theta^u})$
by using the isomorphisms
(\ref{eq;24.4.6.1})
and (\ref{eq;24.4.6.10})
and the filtrations
$\nbigf^{\prime\theta^u}$
on 
$H^0(\nu_0^-(\vecJ),L_{\nu_0^-(\vecJ),<0})$
and
$H^0(\nu_0^+(\vecJ),L_{\nu_0^+(\vecJ),>0})$
defined in {\rm\S\ref{subsection;24.4.4.1}},
and the trivial filtration on $H^0(\real,\nbigt_{\omega}(L))$
indexed by $0$.
It is independent of the choice of $\vecJ_1$.
We obtain the $2\pi\seisuu$-equivariant family of
filtrations $\vecnbigf=(\nbigf^{\theta^u}\,|\,\theta^u\in\real)$
of $\gbigq^0_{\star}(L,\vecnbigftilde)_{\real}$.
By Proposition \ref{prop;24.4.6.2} and Proposition \ref{prop;24.4.6.20}
we obtain the following.
\begin{thm}
$(\gbigq^0_{\star}(L,\vecnbigftilde),\vecnbigf)$
are local systems with Stokes structure indexed by 
$\nbigitilde^{\circ}$.
Moreover, there exists the following commutative diagram 
in $\Loc^{\St}(\nbigitilde^{\circ})$:
\[
\begin{CD}
(\gbigq^0_!(L,\vecnbigftilde)_{\real},\vecnbigf)
 @>{F}>>
 (\gbigq^0_{\ast}(L,\vecnbigftilde)_{\real},\vecnbigf)
 \\
 @V{\simeq}VV @V{\simeq}VV \\ 
 (\gbigl^{\gbigf}_!(V),\vecnbigftilde)
 @>>>
 (\gbigl^{\gbigf}_{\ast}(V),\vecnbigftilde),
\end{CD}
\]
where the lower horizontal arrow is induced by
$V(!0)\to V$.
\hfill\qed
\end{thm}

\begin{df}
We set
$\gbigf^{(0,\infty)}_{+,\star}\bigl(
L,\vecnbigftilde
\bigr)
 :=(\gbigq^0_{\star}(L,\vecnbigftilde),\vecnbigf)$,
called the local Fourier transform of
$(L,\vecnbigftilde)$.
\hfill\qed
\end{df}

\subsection{The induced constructible sheaves and filtrations}
\label{subsection;25.2.25.20}

For $\star=!,\ast$,
we have the constructible subsheaves
$\gbigq^0_{\star}(L,\vecnbigftilde)^{<0}_{\real}
\subset
\gbigq^0_{\star}(L,\vecnbigftilde)^{\leq 0}_{\real}
\subset
\gbigq^0_{\star}(L,\vecnbigftilde)_{\real}$.
For $\theta^u\in\real$,
we have
\[
 \gbigq^0_{!}(L,\vecnbigftilde)^{<0}_{\theta^u}
 =\bigoplus_{\theta^u\in \vecJ}
 \Image \vecA_{\vecJ},
 \quad\quad
 \gbigq^0_{\ast}(L,\vecnbigftilde)^{<0}_{\theta^u}
 =\bigoplus_{\theta^u\in \vecJ}
 \Image \vecA^{\mg}_{\vecJ}.
\]
By choosing $\vecJ_1\in T(\nbigi^{\circ})$
such that $\theta^u\in \vecJ_{1+}$,
we have
\[
 \gbigq^0_{!}(L,\vecnbigftilde)^{\leq 0}_{\theta^u}
 =\bigoplus_{\theta^u\in \vecJ}
 \Image \vecA_{\vecJ}
 \oplus
 \Image \vecB^{\vecJ_{1-}}_{\infty},
 \quad
 \gbigq^0_{\ast}(L,\vecnbigftilde)^{\leq 0}_{\theta^u}
 =\bigoplus_{\theta^u\in \vecJ}
 \Image \vecA^{\mg}_{\vecJ}
 \oplus
 \Image \vecB^{\mg\vecJ_{1-}}_{\infty}.
\]

For any $\vecJ\in T(\nbigi^{\circ})$,
we have the isomorphism induced by $\vecA_{\vecJ}$
or $\vecA^{\mg}_{\vecJ}$:
\begin{equation}
\label{eq;25.2.25.10}
 H^0(\nu_0^-(\vecJ),L_{\nu_0^-(\vecJ),<0})
 \simeq
 H^0(\vecJ,\gbigq^0_{\star}(L,\vecnbigftilde)_{\vecJ,<0}).
\end{equation}
We also have the isomorphism induced by $\vecB_{\vecJ_-}$
or $\vecB^{\mg}_{\vecJ_-}$:
\begin{equation}
\label{eq;25.2.25.11}
 H^0\bigl(\nu_0^+(\vecJ),L_{\nu_0^+(\vecJ),>0}\bigr)
 \simeq
 H^0\bigl(\vecJ,\gbigq^0_{\star}(L,\vecnbigftilde)_{\vecJ,>0}\bigr).
\end{equation}
It is also obtained as the composition of the following maps:
\begin{multline}
 H^0(\nu_0^+(\vecJ),L_{\nu_0^+(\vecJ),>0})
 \simeq
 H^0(\nu_0^+(\vecJ)_-,L_{\nu_0^+(\vecJ)_-,>0})
 \subset
 H^0(\real,L)
 \stackrel{a}{\lrarr}
 \\
 \gbigq_{\star}^0(L,\vecnbigftilde)
 \stackrel{R_{\vecJ}}{\lrarr}
 H^0(\vecJ,\gbigq^0_{\star}(L,\vecnbigftilde)_{\vecJ,>0}). 
\end{multline}
Here, $a$ is induced by the natural morphism
$L\to \gbigq^0_{\star}(L,\vecnbigftilde)_{\real}$.
The Stokes filtrations $\vecnbigf$
on $H^0(\vecJ,\gbigq^0_{\star}(L,\vecnbigftilde)_{<0})$
and $H^0(\vecJ,\gbigq^0_{\star}(L,\vecnbigftilde)_{>0})$
equal to the filtrations $\nbigftilde$
on $H^0(\nu_0^-(\vecJ),L_{\nu_0^-(\vecJ),<0})$
and $H^0(\nu_0^+(\vecJ),L_{\nu_0^+(\vecJ),>0})$
by the isomorphisms (\ref{eq;25.2.25.10})
and (\ref{eq;25.2.25.11}), respectively.
Here, we use the isomorphisms of
the partially ordered sets in {\rm(\ref{eq;24.3.14.40})}
to identify the index sets of the filtrations.

\subsection{Extensions}

Let $M$ and $M_0$ denote the monodromy automorphisms
of $L$ and $\nbigt_{\omega}(L)$, respectively.
Let $L_1$ be a $2\pi\seisuu$-equivariant local system
with morphisms
\[
\nbigt_{\omega}(L)\stackrel{a}{\lrarr}
L_1\stackrel{b}{\lrarr} \nbigt_{\omega}(L).
\]
Together with $\gbigq_!(L,\vecnbigftilde)
\to
\gbigq_{\ast}(L,\vecnbigftilde)$,
we obtain the extension
$\Ltilde_1$.
We have the induced Stokes structure $\vecnbigf$
of $\Ltilde_1$.
We also have the induced morphisms
$\gbigq_!(L,\vecnbigftilde)
\stackrel{u_1}{\lrarr}
\Ltilde_1
\stackrel{u_2}{\lrarr}
\gbigq_{\ast}(L,\vecnbigftilde)$,
and
$L
\stackrel{\atilde}{\lrarr}
\Ltilde_1
\stackrel{\btilde}{\lrarr}
L$.
Let $M_{L_1}$ and $M_{\Ltilde_1}$
denote the monodromy automorphisms of $L_1$
and $\Ltilde_1$, respectively.
We obtain the following proposition from
Proposition \ref{prop;25.2.22.10}.
\begin{prop}
\label{prop;25.2.23.10}
If $b\circ a=\id-M_{L_1}^{-1}$,
then we have
$\btilde\circ\atilde=\id-M_{\Ltilde_1}^{-1}$.
\hfill\qed
\end{prop}

\subsection{The recovery of the Stokes filtrations}
\label{subsection;25.2.23.20}

Let us recover the Stokes structure $\vecnbigftilde$ of $L$
from $(\Ltilde_1,\vecnbigf)$.

\subsubsection{The recovery of $L^{<0}$}

For any $\vecJ\in T(\nbigi^{\circ})$,
the morphism $H^0(\real,\Ltilde_1)\to H^0(\real,L)$
induces an isomorphism
\begin{equation}
\label{eq;25.2.22.40}
 H^0\bigl(\vecJ,(\Ltilde_1)_{\vecJ,<0}\bigr)
 \lrarr
 H^0(\nu_0^-(\vecJ),L_{\nu_0^-(\vecJ),<0}). 
\end{equation}
Because
$L^{<0}=
\bigoplus_{J\in T(\nbigi)}
a_{J!}(L_{J,<0})$,
we can recover $L^{<0}$
by (\ref{eq;25.2.22.40}).
The Stokes filtrations $\vecnbigftilde$ on
$H^0(\nu_0^-(\vecJ),L_{\nu_0^-(\vecJ),<0})$
are recovered from the Stokes filtrations on
$H^0\bigl(\vecJ,(\Ltilde_1)_{\vecJ,<0}\bigr)$,
where we use the isomorphism of
the partially ordered sets in {\rm(\ref{eq;24.3.14.40})}
to identify the index sets of the filtrations.

\subsubsection{The recovery of $L^{\leq 0}$ and the positive parts}

Let $\theta_1\in T(\nbigi)$.
We set $\theta^u=\theta_1+\pi/2$.
Let $\vecJ_1\in T(\nbigi^{\circ})$
such that $\vartheta^{\vecJ_1}_r=\theta^u$.
We set $\vecJ_2=\vecJ_1+(1+\omega^{-1})\pi$
and $\vecJ_3=\vecJ_1+\pi$.
We have the isomorphism
\begin{multline}
 H^0(\real,\Ltilde_1)
 =\Ltilde_{1|\theta^u}
 \simeq
\\
 \bigoplus_{\vecJ_1\leq \vecJ<\vecJ_2}
 \Bigl(
 H^0(\vecJ,(\Ltilde_1)_{\vecJ,<0})
 \oplus
 H^0(\vecJ_+,(\Ltilde_1)_{\vecJ_+,>0})
 \Bigr)
 \oplus
 H^0(\vecJ_{1+},(\Ltilde_1)_{\vecJ_{1+},0}).
\end{multline}
We set
\[
 K_{\theta^u}
 =\!\!\!\!\!\!\bigoplus_{\vecJ_1<\vecJ<\vecJ_2}
 \!\!\!\!\!
 H^0(\vecJ,(\Ltilde_1)_{\vecJ,<0})
 \oplus
 \!\!\!\!\!
 \bigoplus_{\vecJ_1\leq \vecJ<\vecJ_3}
 \!\!\!\!\!
 H^0(\vecJ_+,(\Ltilde_1)_{\vecJ_+,>0})
 \oplus
 H^0(\vecJ_{1+},(\Ltilde_1)_{\vecJ_{1+},0}).
\]
We have the natural map
$h:
 H^0(\real,L)\lrarr
 H^0(\real,\Ltilde_1)$.
Under the natural isomorphism $H^0(\real,L)\simeq L_{\theta_1}$,
we have
\[
 (L^{\leq 0})_{\theta_1}
 =h^{-1}\bigl(
 K_{\theta^u}
 \bigr).
\]
Let $\theta_1'\in S_0(\nbigi)$
determined by 
$\openopen{\theta_1'}{\theta_1}\cap S_0(\nbigi)=\emptyset$.
For any $\theta\in \openopen{\theta_1'}{\theta_1}$,
under the isomorphism
$L_{\theta}=H^0(\real,L)$,
we have
\[
 (L^{\leq 0})_{\theta}
 =h^{-1}\Bigl(
 K_{\theta^u}
 \oplus
 H^0(\vecJ_1,(\Ltilde_1)_{\vecJ_1,<0})
 \Bigr).
\]
Let $\theta_1''\in S_0(\nbigi)$
determined by 
$\openopen{\theta_1}{\theta_1''}\cap S_0(\nbigi)=\emptyset$.
For any $\theta\in\openopen{\theta_1}{\theta_1''}$,
under the isomorphism
$L_{\theta}=H^0(\real,L)$,
we have
\[
 (L^{\leq 0})_{\theta}
 =h^{-1}\Bigl(
 K_{\theta^u}
 \oplus
 H^0(\vecJ_{3+},(\Ltilde_1)_{\vecJ_{3+,>0}})
 \Bigr).
\]
Thus,
the constructible subsheaf
$L^{\leq 0}\subset L$
is recovered.

Let $J_1\in T(\nbigi)$ such that
$\vartheta^{J_1}_r=\theta_1$.
We set $J_2=J_1+\omega^{-1}\pi$.
Note that
$\nu_0^{+}(\vecJ_3)=J_1$
and
$\nu_0^{+}(\vecJ_2)=J_2$.
We have the decomposition
\[
 L_{\theta_1}
 =\bigoplus_{J_1\leq J<J_3}
 \Bigl(
 H^0(J,L_{J,<0})
 \oplus
 H^0(J_{+},L_{J_{+},>0})
 \Bigr)
 \oplus
 H^0(J_{1+},L_{J_{1+},0})
\]
We have
$h\bigl(
 H^0(J_{1+},L_{J_{1+},>0})
 \bigr)
 \subset
 K_{\theta^u}
 \oplus
 H^0(\vecJ_{3+},(\Ltilde_1)_{\vecJ_{3+},>0})$,
and it induces an isomorphism
\[
 H^0(J_{1},L_{J_{1},>0})
 \simeq
 H^0(\vecJ_{3},(\Ltilde_1)_{\vecJ_{3},>0}).
\]
The Stokes filtrations $\vecnbigftilde$ on
$H^0(J_1,L_{J_1,>0})$
are recovered from the Stokes filtrations on
$H^0\bigl(\vecJ_3,(\Ltilde_1)_{\vecJ_3,>0}\bigr)$,
where we use the isomorphism of
the partially ordered sets in {\rm(\ref{eq;24.3.14.40})}
to identify the index sets of the filtrations.

\section{Stokes shells}
\label{subsection;18.5.7.10}

To explain the formula for
$\Shsf\bigl(
\gbigf^{(0,\infty)}_{+,\star}(L,\vecnbigftilde)
\bigr)$,
we introduce transformations for
$\Sh=(\nbigk_{\bullet},\vecnbigf,\vecnbigr)
\in \Shcat(\nbigitilde)$.
We set
$(\vecK,\vecnbigf,\vecPhi,\vecPsi):=
\gbigd(\nbigk_{\bullet},\vecnbigf)$.
We use the notation
$\nbigp_J=\nbigr^{0,J_-}_{\lambda_-(J),J_+}$, 
$\nbigq_J=\nbigr^{\lambda_+(J),J_-}_{0,J_+}$,
$\nbigr^{J_-}_{J_+}=\nbigr^{\lambda_+(J),J_-}_{\lambda_-(J),J_+}$
and 
$\nbigr^{J_+}_{J_-}=\nbigr^{\lambda_+(J),J_+}_{\lambda_-(J),J_-}$
for $\Sh$ and $J\in T(\nbigi)$.

\subsection{Stokes graded local systems}

Take $\vecJ=I(\vartheta^u_0,(1+\omega^{-1})\pi/2)\in T(\nbigi^{\circ})$.
We obtain the intervals
$\nu^{\pm}_{m}(\vecJ)\in T(\nbigi)$ $(m\in\seisuu)$
as in (\ref{eq;20.10.13.1}).
There exist the isomorphisms
$\kappa^{\pm}_{m,\vecJ}:
 \vecJbar\simeq
 \nu^{\pm}_m(\vecJbar)$
as in (\ref{eq;20.10.9.1}).
By Proposition \ref{prop;18.5.5.40},
we obtain the following local systems with Stokes structure
indexed by $\nbigitilde^{\circ}_{\vecJ,<0}$
on $\vecJbar$:
\index{local system with Stokes structure
 $(\nbigk^{\circ}_{\lambda_-(\vecJ),\vecJ},\vecnbigf^{\circ})$}
\[
 (\nbigk^{\circ}_{\lambda_-(\vecJ),\vecJ},\vecnbigf^{\circ}):=
 \bigl(
 \kappa_{0,\vecJ}^{-}
 \bigr)^{-1}
 \bigl(\nbigk_{\lambda_-(\nu^-_0(\vecJ))},
 \vecnbigf\bigr)_{|\nu^{-}_{0}(\vecJbar)}.
\]
We also obtain the following local systems with Stokes structure
indexed by $\nbigitilde^{\circ}_{\vecJ,>0}$
on $\vecJbar$:
\[
 (\nbigk^{\circ}_{\lambda_+(\vecJ),\vecJ},\vecnbigf^{\circ}):=
 \bigl(
 \kappa_{0,\vecJ}^{+}
 \bigr)^{-1}\bigl(
 \nbigk_{\lambda_+(\nu^+_0(\vecJ))},
 \vecnbigf\bigr)_{|\nu^{+}_{0}(\vecJbar)}.
\]
We obtain the following local system on $\vecJbar$:
\[
 \nbigk^{\circ}_{0,\vecJ}:=
 (\kappa_{0,\vecJ}^-)^{-1}
 \bigl(
 \nbigk_0
 \bigr)_{|\nu^{-}_{0}(\vecJbar)}.
\]
The spaces of the global sections 
of $\nbigk^{\circ}_{\lambda,\vecJ}$
are denoted by
$K^{\circ}_{\lambda,\vecJ}$.
\index{vector spaces $K^{\circ}_{\lambda,\vecJ}$}
There exist the natural identifications:
\[
 K^{\circ}_{\lambda_-(\vecJ),\vecJ}=
 K_{<0,\nu^-_0(\vecJ)},
\quad\quad
 K^{\circ}_{\lambda_+(\vecJ),\vecJ}=
 K_{>0,\nu^+_0(\vecJ)},
\quad\quad
 K^{\circ}_{0,\vecJ}=
 K_{0,\nu^-_0(\vecJ)}.
\]
By the construction
and the relation
$\kappa^{\pm}_{0,\vecJ}\circ\Tbb
=\Tbb\circ \kappa^{\pm}_{0,\Tbb^{-1}(\vecJ)}$,
there exist the natural isomorphisms
$\Tbb^{-1}\nbigk^{\circ}_{\lambda,\vecJ}
\simeq
 \nbigk^{\circ}_{\Tbb^{\ast}(\lambda),\Tbb^{-1}(\vecJ)}$,
which induce
$\Psi^{\circ}_{\lambda,\vecJ}:
 K^{\circ}_{\lambda,\vecJ}
\simeq
 K^{\circ}_{\Tbb^{\ast}(\lambda),\Tbb^{-1}(\vecJ)}$.

Because
$\nu_0^{+}(\vecJ+(1+\omega^{-1})\pi)
=\nu_0^-(\vecJ)+\omega^{-1}\pi$,
we obtain the following isomorphisms:
\[
 (\Phi^{\circ})^{\vecJ+(1+\omega^{-1})\pi,\vecJ}_{\lambda_-(\vecJ)}:
=\Phi^{\nu_0^-(\vecJ)+\omega^{-1}\pi,\nu_0^-(\vecJ)}
 _{\lambda_-(\nu_0^-(\vecJ))}:
 K^{\circ}_{\lambda_-(\vecJ),\vecJ}
\simeq
 K^{\circ}_{\lambda_-(\vecJ),\vecJ+(1+\omega^{-1})\pi}.
\]
Because 
$\nu_{-1}^{-}(\vecJ+(1+\omega^{-1})\pi)
=\nu_0^+(\vecJ)+\omega^{-1}\pi$,
we obtain the following isomorphisms:
\begin{multline}
  (\Phi^{\circ})^{\vecJ+(1+\omega^{-1})\pi,\vecJ}_{\lambda_+(\vecJ)}:
=-\Psi^{-1}_{\nu_0^-(\vecJ+(1+\omega^{-1})\pi)}\circ
 \Phi^{\nu_0^+(\vecJ)+\omega^{-1}\pi,\nu_0^+(\vecJ)}
 _{\lambda_+\bigl(\nu_0^+(\vecJ)\bigr)}: \\
 K^{\circ}_{\lambda_+(\vecJ),\vecJ}
\simeq
 K^{\circ}_{\lambda_+(\vecJ),\vecJ+(1+\omega^{-1})\pi}.
\end{multline}
For $\vecJ_1\vdash \vecJ_2$ in $T(\nbigi^{\circ})$,
because $\nu^-_0(\vecJ_1)\vdash\nu^-_0(\vecJ_2)$,
we obtain the following isomorphisms:
\index{maps $(\Phi^{\circ})^{\vecJ_1,\vecJ_2}_{\lambda}$}
\[
 \bigl(
 \Phi^{\circ}
 \bigr)^{\vecJ_2,\vecJ_1}_0:=
\Phi_0^{\nu^-_0(\vecJ_2),\nu^-_0(\vecJ_1)}:
 K^{\circ}_{0,\vecJ_1}\simeq
 K^{\circ}_{0,\vecJ_2}.
\]
By gluing
$(\nbigk^{\circ}_{\lambda,\vecJ},\vecnbigf^{\circ})$ 
via the tuple of the isomorphisms
$\vecPhi^{\circ}$,
we obtain a Stokes graded local system
$(\nbigk^{\circ}_{\bullet},\vecnbigf^{\circ})$
over $(\nbigitilde^{\circ},[\nbigi^{\circ}])$.
By the construction,
it is naturally $2\pi\seisuu$-equivariant.
\index{local systems with Stokes structure
$(\nbigk^{\circ}_{\bullet},\vecnbigf^{\circ})$}

For both $\star=!,\ast$,
we set
$(\nbigk_{\star\,\bullet},\vecnbigf^{\circ}):=
 (\nbigk_{\bullet},\vecnbigf^{\circ})$.
We naturally have
$\gbigd(\nbigk_{\star\bullet},\vecnbigf^{\circ})
=(\vecK^{\circ},\vecnbigf^{\circ},\vecPhi^{\circ},\vecPsi^{\circ})$.
\index{local systems with Stokes structure
$(\nbigk^{\circ}_{\bullet\ast},\vecnbigf^{\circ})$}
\index{local systems with Stokes structure
$(\nbigk^{\circ}_{\bullet!},\vecnbigf^{\circ})$}

\subsection{Morphisms $\vecnbigp^{\circ}_{\star}$ 
and $\vecnbigq^{\circ}_{\star}$ $(\star=!,\ast)$}

Let $M_0$ denote the automorphism of $\nbigk_0$
obtained as the composition of the isomorphisms
$\nbigk_0\stackrel{a}{\simeq}
 \Tbb^{-1}\nbigk_0
 \stackrel{b}{\simeq}\nbigk$,
where $a$ is induced by the parallel transport,
and $b$ is induced by the $2\pi\seisuu$-equivariance.
\index{automorphism $M_0$}
For $\vecJ\in T(\nbigi^{\circ})$,
we set
\[
 (\nbigp^{\circ}_!)_{\vecJ}:=
 \nbigp_{\nu_0^-(\vecJ)}\circ
 (\id-M_0^{-1}),
\quad\quad
 (\nbigp^{\circ}_{\ast})_{\vecJ}:=
 \nbigp_{\nu_0^-(\vecJ)},
\]
\[
 (\nbigq^{\circ}_!)_{\vecJ}:=
 \Phi_0^{\nu_0^-(\vecJ),\nu_0^+(\vecJ)}\circ
 \nbigq_{\nu_0^+(\vecJ)},
\quad\quad
  (\nbigq^{\circ}_{\ast})_{\vecJ}:=
 (\id-M_0^{-1}) \circ
 \Phi_0^{\nu_0^-(\vecJ),\nu_0^+(\vecJ)}\circ
 \nbigq_{\nu_0^+(\vecJ)}.
\]
\index{maps \mbox{$(\nbigp^{\circ}_!)_{\vecJ}$},
$(\nbigp^{\circ}_{\ast})_{\vecJ}$,
\mbox{$(\nbigq^{\circ}_!)_{\vecJ}$},
$(\nbigq^{\circ}_{\ast})_{\vecJ}$}
(Recall that there exist the isomorphisms
$\Phi_{\lambda}^{J_2,J_1}:K_{\lambda,J_1}\simeq K_{\lambda,J_2}$
for any $J_1,J_2\in T(\lambda)$
induced by the parallel transport of $\nbigk_{\lambda}$,
as in \S\ref{subsection;20.10.9.2}.))

\subsection{Morphisms $\vecnbigr^{\circ}$}

We set
\[
 (\nbigr^{\circ})^{\vecJ_-}_{\vecJ_+}
 :=
 \nbigrtilde^{\nu_0^+(\vecJ)_-}_{\nu_0^-(\vecJ)}
+\nbigrtilde^{\nu_0^+(\vecJ)_-}_{\nu_0^-(\vecJ)-2\pi}.
\]
For $\vecJ'<\vecJ$,
we set
\[
 (\nbigr^{\circ}_1)^{\vecJ}_{\vecJ'}:=
 \left\{
\begin{array}{ll}
 \nbigrtilde^{\nu_0^+(\vecJ)_-}_{\nu_0^-(\vecJ)}
  & (\vecJ-(1+\omega^{-1}\pi)<\vecJ'<\vecJ-\pi)
  \\
  -\nbigrtilde^{\nu_0^+(\vecJ)_-}_{\nu_0^-(\vecJ)}
  & (\vecJ-\pi\leq\vecJ'<\vecJ),
\end{array}
 \right.
\]
\[
 (\nbigr^{\circ}_2)^{\vecJ}_{\vecJ'}:=
 \left\{
\begin{array}{ll}
 0
  & (\vecJ'<\vecJ-(\omega^{-1}-1)\pi)
  \\
 -\Psi^{-1}\circ
  \nbigrtilde^{\nu_0^+(\vecJ)_-}_{\nu_0^-(\vecJ)-2\pi}
  & (\vecJ-(\omega^{-1}-1)\pi\leq\vecJ'<\vecJ).
\end{array}
 \right.
\]
For $\vecJ<\vecJ'$,
we set
\[
 (\nbigr^{\circ}_1)^{\vecJ}_{\vecJ'}
 :=\left\{
\begin{array}{ll}
 -\Psi^{-1}\circ \nbigrtilde^{\nu_0^+(\vecJ)}_{\nu_0^-(\vecJ')-2\pi}&
 (\vecJ+\pi'<\vecJ'<\vecJ+(1+\omega^{-1})\pi)
 \\
 \Psi^{-1}\circ \nbigrtilde^{\nu_0^+(\vecJ)}_{\nu_0^-(\vecJ')-2\pi}&
 (\vecJ<\vecJ'\leq \vecJ+\pi)
\end{array}
 \right.
\]
\[
 (\nbigr^{\circ}_2)^{\vecJ}_{\vecJ'}
 :=\left\{
\begin{array}{ll}
 0 &(\vecJ+(\omega^{-1}-1)\pi<\vecJ')
 \\
 \nbigrtilde^{\nu_0^+(\vecJ)_+}_{\nu_0^-(\vecJ')}
 & (\vecJ<\vecJ'\leq \vecJ+(\omega^{-1}-1)\pi).
\end{array}
 \right.
\]
Then, we set
\[
 (\nbigr^{\circ})^{\vecJ}_{\vecJ'}
 =(\nbigr^{\circ}_1)^{\vecJ}_{\vecJ'}
 +(\nbigr^{\circ}_2)^{\vecJ}_{\vecJ'}.
\]

\subsection{Isomorphisms}

For $\star=!,\ast$,
let $\gbigf^{(0,\infty)}_{+\star}(\Sh)$
be the Stokes shell obtained as
$(\nbigk^{\circ}_{\star \bullet},\vecnbigf^{\circ})$
with the tuple of the morphisms
$\bigl(
 \vecnbigp^{\circ}_{\star},
 \vecnbigq^{\circ}_{\star},
 \vecnbigr \bigr)$.
They are objects in $\Shcat(\nbigitilde^{\circ})$.
\index{Stokes shells
$\gbigf^{(0,\infty)}_{+\ast}(\Sh)$,
\mbox{$\gbigf^{(0,\infty)}_{+!}(\Sh)$}}
We have the morphism
$F:
 \gbigf^{(0,\infty)}_{+!}(\Sh)
\lrarr
 \gbigf^{(0,\infty)}_{+\ast}(\Sh)$
induced by the identity maps
$\nbigk^{\circ}_{!\lambda}
=\nbigk^{\circ}_{\ast\lambda}$ 
$(\lambda\neq 0)$
and 
$\id-M_0^{-1}:
 \nbigk^{\circ}_{!0}\lrarr
 \nbigk^{\circ}_{\ast 0}$.
Then, by the construction,
$(\gbigf^{(0,\infty)}_{+!}(\Sh),
 \gbigf^{(0,\infty)}_{+\ast}(\Sh),
 F)$
is a base tuple 
in $\Shcat(\nbigitilde^{\circ})$.
\index{base tuple
\mbox{$(\gbigf^{(0,\infty)}_{+!}(\Sh),
 \gbigf^{(0,\infty)}_{+\ast}(\Sh),F)$}}
As the translation of the results in
{\rm\S\ref{subsection;24.3.24.10}--\ref{subsection;24.3.24.11}},
we obtain the following.
\begin{prop}
There exists the commutative diagram:
\[
\begin{CD}
 \gbigf^{(0,\infty)}_{+,!}\bigl(\Shsf(L,\vecnbigftilde)\bigr)
 @>>>
 \gbigf^{(0,\infty)}_{+,\ast}\bigl(\Shsf(L,\vecnbigftilde)\bigr)
 \\
 @V{\simeq}VV @V{\simeq}VV \\
 \Shsf\bigl(
 \gbigf^{(0,\infty)}_{+,!}(L,\vecnbigftilde)
 \bigr)
 @>>>
 \Shsf\bigl(
 \gbigf^{(0,\infty)}_{+,\ast}(L,\vecnbigftilde)
 \bigr).
\end{CD}
\]
\hfill\qed
\end{prop}

\subsection{Another description of the Stokes graded local systems}
\label{subsubsection;18.5.8.1}

For each $\lambda\in[(\nbigi^{\circ})^{\ast}]$,
we take
$\vecJ_{\lambda}=
I(\vartheta^u_{0,\lambda},(1+\omega^{-1})\pi/2)
\in T(\lambda)_{<0}$.
We define the map
$\kappa_{\vecJ_{\lambda}}:\real\lrarr\real$
as
\[
 \kappa_{\vecJ_{\lambda}}(\theta^u)
 =\frac{1}{1+\omega}(\theta^u+\omega\vartheta^u_{0,\lambda}). 
\]
We obtain the local system with Stokes structure
$(\nbigk^{\circ\circ}_{\lambda},\vecnbigf^{\circ\circ}):=
 \kappa_{\vecJ_{\lambda}}^{-1}
 \bigl(\nbigk_{\lambda_-(\nu_0^-(\vecJ_{\lambda}))},\vecnbigf\bigr)$.
 \index{local system with Stokes structure
 $(\nbigk^{\circ\circ}_{\lambda},\vecnbigf^{\circ\circ})$}
By the construction,
there exists the natural isomorphism
\[
 (\nbigk^{\circ\circ}_{\lambda},\vecnbigf^{\circ\circ})
 _{|\vecJ_{\lambda}\cup (\vecJ_{\lambda}+(1+\omega^{-1})\pi)}
\simeq
 (\nbigk^{\circ}_{\lambda},\vecnbigf^{\circ})
 _{|\vecJ_{\lambda}\cup(\vecJ_{\lambda}+(1+\omega^{-1})\pi)}.
\]
 Because
 \[
\kappa_{\vecJ_{\lambda}} =
\Tbb^{-m}\circ\kappa^{-}_{0,\vecJ_{\lambda}+2m(1+\omega^{-1})\pi}
=\Tbb^{-m}\circ\kappa^{+}_{0,\vecJ_{\lambda}+(2m+1)(1+\omega^{-1})\pi}
\quad (m\in\seisuu),
 \]
it uniquely extends to an isomorphism
$b_{\lambda}:(\nbigk^{\circ\circ}_{\lambda},\vecnbigf^{\circ\circ})
\simeq
 (\nbigk^{\circ}_{\lambda},\vecnbigf^{\circ})$,
where the restriction of $b_{\lambda}$ to
$(\vecJ_{\lambda}+2m(1+\omega^{-1})\pi)
\cup
 (\vecJ_{\lambda}+(2m+1)(1+\omega^{-1})\pi)$
is induced by $(-1)^m\Psi^{-m}$.

We also set
$\nbigk^{\circ\circ}_{0}:=\nbigk_0$.
Let us observe that 
there exists a natural isomorphism
$b_0:\nbigk^{\circ\circ}_0\simeq
 \nbigk^{\circ}_0$.
Take any $\vecJ\in T(\nbigi^{\circ})$.
If $\theta^u\in\vecJ$,
we obtain
$\nu_{0}^-(\vecJ)\cap
 \openopen{\theta^u-\pi/2}{\theta^u+\pi/2}
\neq\emptyset$.
Hence, there exists an isomorphism
\[
\nbigk^{\circ\circ}_{0|\theta^u}
:=\nbigk_{0|\theta^u}
\simeq
 K_{0,\nu_0^-(\vecJ)},
\]
which induces the desired isomorphism
$b_0:\nbigk^{\circ\circ}_0\simeq
 \nbigk^{\circ}_0$.

 We set
$(\nbigk^{\circ\circ}_{\bullet},\vecnbigf^{\circ\circ}):=
 \bigoplus_{\lambda\in[\nbigi^{\circ}]}
 (\nbigk^{\circ\circ}_{\lambda},\vecnbigf^{\circ\circ})$.
There exists the isomorphism
$b:(\nbigk^{\circ\circ}_{\bullet},\vecnbigf^{\circ\circ})
\simeq
 (\nbigk^{\circ}_{\bullet},\vecnbigf)$
 induced by $b_{\lambda}$ $(\lambda\in[\nbigi^{\circ}])$.
An action of $2\pi\seisuu$ 
on 
$(\nbigk^{\circ\circ}_{\bullet},\vecnbigf^{\circ\circ})$
is induced
by the isomorphism $b$
and the $2\pi\seisuu$-action
on $(\nbigk^{\circ}_{\bullet},\vecnbigf^{\circ})$.

\begin{rem}
There exist positive integers $n_1,p_1$
such that $n_1/p_1=\omega$ with $\gcd(n_1,p_1)=1$.
For $\lambda\in[(\nbigi^{\circ})^{\ast}]$,
we obtain the following isomorphism:
\[
 a_0:
 (\Tbb^{n_1+p_1})^{\ast}(\nbigk_{\lambda}^{\circ\circ},
 \vecnbigf^{\circ\circ})
\simeq
   (\Tbb^{n_1+p_1})^{\ast}(\nbigk_{\lambda}^{\circ},\vecnbigf^{\circ})
\simeq
 (\nbigk_{\lambda}^{\circ},\vecnbigf^{\circ})
\simeq 
 (\nbigk_{\lambda}^{\circ\circ},\vecnbigf^{\circ\circ}).
\]
We also have the following natural isomorphism:
\begin{multline}
a_1:
 (\Tbb^{n_1+p_1})^{\ast}(\nbigk_{\lambda}^{\circ\circ},
 \vecnbigf^{\circ\circ})
=\kappa_{\vecJ_{\lambda}}^{-1}
 \Bigl(
 (\Tbb^{p_1})^{\ast}(\nbigk_{\lambda_-(\nu_{0}^-(\vecJ_{\lambda}))})
 \Bigr)
\simeq
 \kappa_{\vecJ_{\lambda}}^{-1}
 (\nbigk_{\lambda_-(\nu_{0}^-(\vecJ_{\lambda}))})
 \\
=(\nbigk_{\lambda}^{\circ\circ},\vecnbigf^{\circ\circ}).
\end{multline}
Note that $a_0=(-1)^{n_1}a_1$.
\hfill\qed
\end{rem}

\subsection{Example}

Let $\omega\in\seisuu_{>1}$.
Let
$\nbigi=\{\alpha_iz^{-\omega}\mid i=1,\ldots,N\}
\subset \real_{>0}z^{-\omega}$
be a finite subset.
We set $J_m:=I(m\omega^{-1}\pi,\omega^{-1}\pi/2)$.
We have
$T(\nbigi)=
\bigl\{J_m\,\big|\,
m\in\seisuu
\bigr\}$,
and
$\nbigi=\nbigi_{J_{2\ell},<0}
=\nbigi_{J_{2\ell+1},>0}$.

Let $(V,\nabla)$ be a basic meromorphic flat bundle
of level $(0,\omega)$ 
such that
$\nbigi_0(V)\subset\nbigi$.
Let $(L,\vecnbigf)\in\Loc^{\St}(\nbigi)$
be the corresponding local system with Stokes structure.
In this case,
the associated Stokes shell consists of
$(\nbigk_{\bullet},\vecnbigf)=(L,\vecnbigf)$
and $\vecnbigr=\emptyset$.
We note that $V(\ast 0)=V(!0)$ in this case.
Let $(\gbigl^{\gbigf}(V),\vecnbigf)$
denote the local system with Stokes structure
corresponding to $\Fourier(V)$ at $\infty$.
Let us describe 
the associated Stokes shell
$(\nbigk^{\gbigf}_{\bullet},\vecnbigf,\vecnbigr^{\gbigf})$.

For $k\in\seisuu$,
we set
$\beta_k=\exp\bigl(2\pi\sqrt{-1}k/(1+\omega)\bigr)$.
We set
$\nbigi^{\circ}_k:=\bigl\{
\langle\omega\rangle\alpha_i^{\frac{1}{1+\omega}}
\beta_k u^{-\frac{\omega}{1+\omega}}
\bigr\}$,
and
$\nbigi^{\circ}=
\bigcup_{k=0}^{\omega}\nbigi^{\circ}_k$.
Note that $\nbigi_{\infty}(\Fourier(V))\subset\nbigi^{\circ}$.
 
We set
$\vecJ_{k,m}:=
I\bigl(2\pi\omega^{-1}k+m(1+\omega^{-1})\pi,
(1+\omega^{-1})\pi/2\bigr)$.
For $k\in\seisuu$, we have
$T(\nbigi_k^{\circ})
=\{\vecJ_{k,m}\,|\,m\in\seisuu\}$.
Hence,
$T(\nbigi^{\circ})=
 \bigcup_{k=0}^{\omega}\bigl\{
  \vecJ_{k,m}\,\big|\,m\in\seisuu
  \bigr\}$.
We have
$\nbigi^{\circ}_{\vecJ_{k,2\ell},<0}
=\nbigi^{\circ}_{\vecJ_{k,2\ell+1},>0}
=\nbigi^{\circ}_k$.
We have
$\nu_0^+(\vecJ_{k,2\ell+1})
 =J_{2(k+\ell\omega+\ell)+1}$
and
$\nu_0^-(\vecJ_{k,2\ell})
=J_{2(k+\ell\omega+\ell)}$.

We obtain local systems with filtrations
$(\kappa_{0,\vecJ_{k,2\ell+1}}^+)^{-1}
(L_{|\overline{\nu_0^+(\vecJ_{k,2\ell+1})}},\vecnbigf)$
on $\overline{\vecJ_{k,2\ell+1}}$
and 
$(\kappa_{0,\vecJ_{k,2\ell}}^-)^{-1}
(L_{|\overline{\nu_0^-(\vecJ_{k,2\ell})}},\vecnbigf)$
on $\overline{\vecJ_{k,2\ell}}$.
The index sets are $\nbigi^{\circ}_k$.
Because
$\nu_0^+(\vecJ_{k,2\ell+1})=\nu_0^-(\vecJ_{k,2\ell})+\omega^{-1}\pi$,
we have the natural isomorphism
at $\vartheta^u_1=\overline{\vecJ_{k,2\ell}}\cap\overline{\vecJ_{k,2\ell+1}}$:
\[
(\kappa_{0,\vecJ_{k,2\ell}}^-)^{-1}
(L_{|\overline{\nu_0^-(\vecJ_{k,2\ell})}},\vecnbigf)
 _{|\vartheta^u_1}
 \simeq
 (\kappa_{0,\vecJ_{k,2\ell+1}}^+)^{-1}
 (L_{|\overline{\nu_0^+(\vecJ_{k,2\ell+1})}},\vecnbigf)
 _{|\vartheta^u_1}.
\]
Because 
$\nu_0^-(\vecJ_{k,2\ell+2})=
\nu_0^+(\vecJ_{k,2\ell+1})+\omega^{-1}\pi+2\pi$,
we obtain the isomorphism
\[
(\kappa_{0,\vecJ_{k,2\ell+1}}^+)^{-1}
(L_{|\overline{\nu_0^+(\vecJ_{k,2\ell+1})}},\vecnbigf)
 _{|\vartheta^u_2}
 \simeq
 (\kappa_{0,\vecJ_{k,2\ell+2}}^-)^{-1}
 (L_{|\overline{\nu_0^-(\vecJ_{k,2\ell+2})}},\vecnbigf)
 _{|\vartheta^u_2}
\]
at $\vartheta^u_2\in
\overline{\vecJ_{k,2\ell+1}}\cap
\overline{\vecJ_{k,2\ell+2}}$,
as the $-1$ times the natural isomorphism.
By patching them, we obtain
a local system with filtrations
$(\nbigk^{\circ}_k,\vecnbigf)$ on $\real$.
There exist the natural isomorphisms
$\Tbb^{-1}(\nbigk^{\circ}_k,\vecnbigf)
\simeq
(\nbigk^{\circ}_{k-\omega},\vecnbigf)$,
where $k-\omega$ is considered in $\seisuu/(\omega+1)\seisuu$.
They induce
the natural $2\pi\seisuu$-action
on $(\nbigk^{\circ}_{\bullet},\vecnbigf)=
\bigoplus_{k=0}^{\omega}(\nbigk^{\circ}_k,\vecnbigf)$.
By Proposition \ref{prop;24.3.25.60},
we have the isomorphism
$(\nbigk^{\circ}_{\bullet},\vecnbigf)
\simeq(\nbigk^{\gbigf}_{\bullet},\vecnbigf)$.

Let us compute $\vecnbigr^{\gbigf}$.
By Proposition \ref{prop;24.3.25.60},
the non-trivial terms are
$(\nbigr^{\gbigf})^{\vecJ}_{\vecJ'}$
in the cases
$\vecJ'=\vecJ-(1-\omega^{-1})\pi$
or
$\vecJ'=\vecJ+(1-\omega^{-1})\pi$.

We have
$\vecJ_{k,2\ell+1}-(1-\omega^{-1})\pi=\vecJ_{k+1,2\ell}$.
If $k=\omega$, we regard $\vecJ_{\omega+1,2\ell}=\vecJ_{0,2(\ell+1)}$.
We have
$\nu_0^-(\vecJ_{k+1,2\ell})
=\nu_0^+(\vecJ_{k,2\ell+1})+\omega^{-1}\pi$.
By Proposition \ref{prop;24.3.25.60},
$(\nbigr^{\gbigf})^{\vecJ_{k,2\ell+1}}_{\vecJ_{k+1,2\ell}}$
equals the $-1$ times 
the natural isomorphism
$H^0(\nu_0^+(\vecJ_{k,2\ell+1}),L)
\simeq
H^0(\nu_0^-(\vecJ_{k+1,2\ell}),L)$.

We have
$\vecJ_{k,2\ell+1}+(1-\omega^{-1})\pi
=\vecJ_{k-1,2\ell+2}$.
We regard
$\vecJ_{-1,2\ell+2}=\vecJ_{\omega,2\ell}$.
We have
$\nu_0^-(\vecJ_{k-1,2\ell+2})
=\nu_0^+(\vecJ_{k,2\ell+1})-\omega^{-1}\pi+2\pi$.
By Proposition \ref{prop;24.3.25.60},
$(\nbigr^{\gbigf})^{\vecJ_{k,2\ell+1}}_{\vecJ_{k-1,2\ell+2}}$
equals the natural isomorphism
$H^0(\nu_0^+(\vecJ_{k,2\ell+1}),L)
\simeq
H^0(\nu_0^-(\vecJ_{k-1,2\ell+2}),L)$.

\chapter{Reduction at finite place}
\label{section;20.10.30.2}

\section{Introduction to \S\ref{section;20.10.30.2}}

Let $D\subset\cnum$ be a finite subset.
Let $(\nbigv,\nabla)$ be a meromorphic flat bundle
on $(\proj^1,D\cup\{\infty\})$
with regular singularity at $\infty$.
Let $(V,\nabla)$ be the regular singular
meromorphic flat bundle on $(\proj^1,D\cup\{\infty\})$
associated with the local system
corresponding to $(\nbigv,\nabla)$.

For each $\alpha\in D$,
set $U_{\alpha}:=\bigl\{z\in\cnum\,\big|\,|z-\alpha|<\epsilon\bigr\}$
for a small positive number $\epsilon$
such that
$U_{\alpha}\cap D=\{\alpha\}$.
\index{open set $U_{\alpha}$}
Each restriction $(\nbigv,\nabla)_{|U_{\alpha}}$
induces a meromorphic flat bundle
$(\nbigv_{\alpha},\nabla)$ on $(\proj^1,\{\alpha,\infty\})$
with regular singularity at $\infty$.
Similarly, each restriction $(V,\nabla)_{|U_{\alpha}}$
induces a regular meromorphic flat bundle
$(V_{\alpha},\nabla)$ on $(\proj^1,\{\alpha,\infty\})$.

Note that $-\ord(\Fourier_+(\nbigv))\leq 1$, and
\[
\pi_1\bigl(\nbigi(\Fourier_+(\nbigv))\bigr)
=\nbigi(\Fourier_+(V))=\{\alpha u^{-1}\,\big|\,\alpha\in D\}
=:\nbigi^{\circ}.
\]
(Here, $\pi_1$ denotes the projection as in \S\ref{subsection;18.4.18.1}.)

\subsection{Reduction of $\gbigl^{\gbigf}_{\varrho}(\nbigv)$}
\label{subsection;24.4.1.100}
We set $\vecnbigf^{(1)}=\pi_{1\ast}(\vecnbigf)$
on $\gbigl^{\gbigf}_{\varrho}(\nbigv)$.

\begin{prop}
\label{prop;24.3.30.1}
For each $\alpha\in D$,
there exists an isomorphism 
\[
  \Gr^{\vecnbigf^{(1)}}_{\alpha u^{-1}}
(\gbigl^{\gbigf}_{\varrho}(\nbigv),\vecnbigf)
\simeq
(\gbigl^{\gbigf}_{\varrho(\alpha)}(\nbigv_{\alpha}),\vecnbigf).
\]
They induce an isomorphism of functors
from $\Dsf(D)$ to the category 
of local systems with Stokes structure,
i.e.,
for any $\varrho_1\to\varrho_2$ in $\Dsf(D)$,
the following diagram is commutative:
\[
 \begin{CD}
 \Gr^{\vecnbigf^{(1)}}_{\alpha u^{-1}}
  (\gbigl^{\gbigf}_{\varrho_1}(\nbigv),\vecnbigf)
@>>>
   \Gr^{\vecnbigf^{(1)}}_{\alpha u^{-1}}
(\gbigl^{\gbigf}_{\varrho_2}(\nbigv),\vecnbigf)
  \\
  @V{\simeq}VV @V{\simeq}VV \\  
  (\gbigl^{\gbigf}_{\varrho_1(\alpha)}(\nbigv_{\alpha}),\vecnbigf)
 @>>>
  (\gbigl^{\gbigf}_{\varrho_2(\alpha)}(\nbigv_{\alpha}),\vecnbigf).
 \end{CD}
\]
\end{prop}

Let $\rho_{\alpha}:\proj^1\lrarr\proj^1$
be the map determined by
$\rho_{\alpha}(z)=z+\alpha$.
The following lemma is easy to see.
Note that we can apply the results in \S\ref{section;18.6.3.20} to
each
$\bigl(
\gbigl^{\gbigf}_{\star}\bigl(\rho_{\alpha}^{\ast}(\nbigv_{\alpha})\bigr),
\vecnbigf\bigr)$.
\begin{lem}
There exists a natural isomorphism of local systems
$\gbigl^{\gbigf}_{\varrho_1(\alpha)}(\nbigv_{\alpha})
\simeq
 \gbigl^{\gbigf}_{\varrho_1(\alpha)}(\rho_{\alpha}^{\ast}\nbigv_{\alpha})$,
which preserves the Stokes filtrations
under the bijection of the index sets 
$\nbigi_{\infty}\bigl(
\Fourier_+(\nbigv_{\alpha}(\varrho_1(\alpha)))
\bigr)
\simeq
\nbigi_{\infty}\bigl(
\Fourier_+(\rho_{\alpha}^{\ast}\nbigv_{\alpha}(\varrho_1(\alpha)))
\bigr)$
defined by $\gminia\longmapsto \gminia-\alpha u^{-1}$.
\hfill\qed
\end{lem}

When $\nbigv=V$,
Proposition \ref{prop;24.3.30.1} implies that
$\Gr^{\vecnbigf^{(1)}}_{\alpha u^{-1}}\gbigl^{\gbigf}_{\varrho}(V)$
are functorially identified with
$\gbigl^{\gbigf^{(1)}}_{\varrho}(V_{\alpha})$,
which also follows from the stationary phase formula.

\begin{prop}
\label{prop;24.4.14.50}
The $2\pi\seisuu$-equivariant local system with Stokes structure
$\bigl(
 \gbigl^{\gbigf}_{\varrho}(\nbigv),
 \vecnbigf^{(1)}
\bigr)$
are obtained as the extension of
the base tuple
$(\gbigl^{\gbigf}_{\varrho}(V),\vecnbigf)$
$(\varrho\in\Dsf(D))$
by the natural morphisms of
the $2\pi\seisuu$-equivariant local systems:
\begin{equation}
\label{eq;24.3.30.2}
 \gbigl^{\gbigf}_{!}(V_{\alpha})
\lrarr
 \gbigl^{\gbigf}_{\varrho(\alpha)}(\nbigv_{\alpha})
\lrarr
 \gbigl^{\gbigf}_{\ast}(V_{\alpha}).
\end{equation}
\end{prop}

In \S\ref{subsection;18.5.15.10},
we shall introduce an explicit construction of
a base tuple of Stokes shells
$\gbigf_{\varrho}(\nbigl)$ $(\varrho\in\Dsf(D))$
from a local system $\nbigl$ on $\cnum\setminus D$.
\index{shells $\gbigf_{\varrho}(\nbigl)$}

\begin{prop}
\label{prop;24.3.31.1}
Let $\nbigl(V)$ denote the local system on $\cnum\setminus D$
associated with $(V,\nabla)$.
Then, there exists the isomorphism
of base tuples
in the category of local systems with Stokes structure. 
\[
 \Locst\bigl(
\gbigf_{\varrho}(\nbigl(V))
\bigr)\simeq
 (\gbigl^{\gbigf}_{\varrho}(V),\vecnbigf)
 \quad (\varrho\in\Dsf(D)).
\]
\end{prop}

\subsubsection{}
\label{subsection;24.4.2.120}
These propositions provide us with the following procedure
to study $(\gbigl^{\gbigf}_{\varrho}(\nbigv),\vecnbigf)$.
\begin{itemize}
 \item $\bigl(\gbigl^{\gbigf}_{\varrho}(\nbigv),\vecnbigf\bigr)$
       are recovered from
       $\bigl(\gbigl^{\gbigf}_{\varrho}(\nbigv),\pi_{1\ast}\vecnbigf\bigr)$
       and the Stokes filtrations of
        $\Gr^{\vecnbigf}_{\alpha u^{-1}}(\gbigl^{\gbigf}_{\varrho}(\nbigv))
\simeq
       \gbigl^{\gbigf}_{\varrho(\alpha)}(\nbigv_{\alpha})$.
       We can apply the results in \S\ref{section;18.6.3.20} to
       $\bigl(
       \gbigl^{\gbigf}_{\star}\bigl(\rho_{\alpha}^{\ast}(\nbigv_{\alpha})\bigr),
       \vecnbigf\bigr)$.
 \item $\bigl(\gbigl^{\gbigf}_{\varrho}(\nbigv),
       \pi_{1\ast}\vecnbigf\bigr)$
       are explicitly described as
       the extension of the base tuple
       $\gbigf_{\varrho}(\nbigl)$ $(\varrho\in\Dsf(D))$
       by the morphisms of the local systems
       (\ref{eq;24.3.30.2}).
\end{itemize}

\subsubsection{Complement}
\label{subsection;24.4.3.1}

Let $U$ be a small neighbourhood of $\infty$ in $\proj^1$.
We obtain the regular singular meromorphic flat bundle
$(V,\nabla)_{|U}$ on $(U,\infty)$,
which extends to a regular singular meromorphic flat bundle
$(V_{\infty},\nabla)=\nbigttilde^{\infty}_0(V,\nabla)$
on $(\proj^1,\{0,\infty\})$.
By Lemma \ref{lem;18.5.12.40},
there exist the following natural morphisms:
\begin{equation}
\label{eq;24.4.1.41}
 \gbigl^{\gbigf}_!(V_{\infty})
 \lrarr
 \gbigl^{\gbigf}_{\varrho}(\nbigv)
 \lrarr
 \gbigl^{\gbigf}_{\ast}(V_{\infty}).
\end{equation}
In \S\ref{section;24.4.1.40},
we shall explicitly describe
(\ref{eq;24.4.1.41}).

\subsection{Extensions}

For each $\alpha\in D$,
let $U_{\alpha}$ denote a small neighbourhood.
We set $U_{\alpha}^{\ast}=U_{\alpha}\setminus\{\alpha\}$.
Let $L_{\alpha}$ be the local systems
on $U_{\alpha}^{\ast}$ obtained as the restriction of
$\nbigl(V)$.
Let $M_{\alpha}$ denote the monodromy automorphism.
We consider morphisms of local systems
\begin{equation}
\label{eq;25.3.1.30}
 L_{\alpha}\stackrel{a_{\alpha}}{\lrarr}
 L_{1,\alpha}
 \stackrel{b_{\alpha}}{\lrarr}
 L_{\alpha}
\end{equation}
such that $b_{\alpha}\circ a_{\alpha}=\id-M_{\alpha}^{-1}$.
We obtain the extension $(\Ltilde_1,\vecnbigf)$
of the base tuple
$(\gbigl^{\gbigf}_{\varrho}(V),\vecnbigf)$
by (\ref{eq;25.3.1.30}).
There exist the natural morphisms
\begin{equation}
\label{eq;25.3.1.50}
 \gbigl^{\gbigf}_!(V_{\infty})
 \stackrel{\atilde}{\lrarr}
 \Ltilde_1
 \stackrel{\btilde}{\lrarr}
 \gbigl_{\ast}^{\gbigf}(V_{\infty}).
\end{equation}
Let $M_{\Ltilde_1}$ denote the monodromy automorphism of
$\Ltilde_1$.
Let $M_{1,\alpha}$ denote the monodromy automorphisms of
$L_{1,\alpha}$.
\begin{prop}[Proposition
\ref{prop;25.3.1.20}]
If
$a_{\alpha}\circ b_{\alpha}=\id-M_{1,\alpha}^{-1}$
for any $\alpha$,
we obtain 
$a_{\Ltilde_1}\circ b_{\Ltilde_1}=
\id-M^{-1}_{\Ltilde_1}$.
\end{prop}

We shall also observe that
$\nbigl(V)$ is recovered from (\ref{eq;25.3.1.50}).

\subsection{Homology groups}

To prove the propositions in \S\ref{subsection;24.4.1.100},
we shall study homology groups of
$(\nbigv,\nabla)$ and $(V,\nabla)$.
In \S\ref{subsection;20.10.24.2},
when $\theta^u\in\vecJ_{\pm}$,
we shall introduce
the following maps
\index{maps $C^{\varrho}_{\vecJ_{\pm},\alpha}$}
\[
 C^{\varrho}_{\vecJ_{\pm},\alpha}:
 H^{\varrho(\alpha)}_1\bigl(
 \cnum\setminus\{\alpha\},
 \nbigv_{\alpha}
 \otimes\nbige(zu^{-1})
 \bigr)
 \lrarr
 H^{\varrho}_1\bigl(
 \cnum\setminus D,
 \nbigv
 \otimes\nbige(zu^{-1})
 \bigr).
\]
We obtain the commutative diagram
(\ref{eq;18.5.15.50}),
where the horizontal arrows are isomorphisms.
The left hand side and the right hand side of
(\ref{eq;18.5.15.50})
are equipped with the filtrations.
As stated in Proposition \ref{prop;24.3.16.10},
they are isomorphisms of filtered vector spaces,
which we shall prove in \S\ref{section;20.11.21.3}.
It implies the propositions in \S\ref{subsection;24.4.1.100}.

\section{Decompositions of homology groups}

\subsection{Construction of maps}
\label{subsection;20.10.24.2}

For any $\vecJ=I(\vartheta^{\vecJ}_0,\pi/2)\in T(\nbigi^{\circ})$,
let $D_{\vecJ}$ denote the set of $\alpha\in D$
such that 
$\alpha u^{-1}\in\nbigi^{\circ}_{\vecJ}$.
\index{set $D_{\vecJ}$}
Any element $\alpha\in D_{\vecJ}$
has the expression
$\alpha=-a\cdot \exp\bigl(\sqrt{-1}\vartheta^{\vecJ}_0\bigr)$
for some $a\in\real$.

Take $\vecJ=I(\vartheta^{\vecJ}_0,\pi/2)\in T(\nbigi^{\circ})$
such that $\arg(u)=\theta^u\in\vecJbar$.
We shall construct the following morphisms
for any $\alpha\in D_{\vecJ}$
and $\varrho\in \Dsf(D)$:
\index{map $C^{\varrho}_{\vecJ_{\pm},\alpha}$}
\begin{equation}
\label{eq;24.3.15.1}
 C^{\varrho}_{\vecJ_{\pm},\alpha}:
 H_1^{\varrho(\alpha)}\bigl(
 \cnum\setminus\{\alpha\},
 \nbigv_{\alpha}\otimes\nbige(zu^{-1})
 \bigr)
\lrarr
 H_1^{\varrho}\bigl(
 \cnum\setminus D,
 \nbigv\otimes\nbige(zu^{-1})
 \bigr).
\end{equation}
We mean that
we construct 
$C^{\varrho}_{\vecJ_-,\alpha}$ 
if $\theta^u\in\vecJ_-$,
and
$C^{\varrho}_{\vecJ_+,\alpha}$ 
if $\theta^u\in\vecJ_+$.
Similarly, we shall construct
the following maps:
\begin{equation}
 \label{eq;24.3.15.2}
 C^{\varrho}_{\vecJ_{\pm},\alpha}:
 H_1^{\varrho(\alpha)}\bigl(
 \cnum\setminus\{\alpha\},
 V_{\alpha}\otimes\nbige(zu^{-1})
 \bigr)
\lrarr
 H_1^{\varrho}\bigl(
 \cnum\setminus D,
 V\otimes\nbige(zu^{-1})
 \bigr).
\end{equation}

\subsection{The case of ``$-$''}

Suppose that $\theta^u\in \vecJ_-$.
Let $\varpi:\projtilde^1_{\infty}\lrarr\proj^1$ denote
the oriented real blow up along $\infty$.
Let $\varpi_D:
 \projtilde^1_{\infty\cup D}\lrarr\projtilde_{\infty}^1$
denote the oriented real blow up along $D$.
For $\alpha\in D$,
let $\varpi_{\alpha}:\projtilde^1_{\infty\cup\alpha}
\lrarr \projtilde^1_{\infty}$
denote the oriented real blow up along $\alpha$.
We set $\vartheta^{\vecJ}_{\ell}=\vartheta^{\vecJ}_0-\pi/2$
and $\vartheta^{\vecJ}_{r}=\vartheta^{\vecJ}_0+\pi/2$.

Take $0<\delta<\!<\vartheta^{\vecJ}_{r}-\theta^u$.
Take a small $\epsilon>0$.
We have the following subset of $\projtilde^1_{\infty}$:
\begin{multline}
 \nbigu_{\vecJ_-,\theta^u}:=
 \bigl\{
 s_1e^{\sqrt{-1}\vartheta^{\vecJ}_{0}}
+s_2e^{\sqrt{-1}\vartheta^{\vecJ}_{\ell}}\,\big|\,
 s_1\in\real,\,\,
 0<s_2<2\epsilon
 \bigr\}
\cup \\
 \bigl\{
 re^{\sqrt{-1}\theta}\,\big|\,
 0\leq r\leq\infty,\,\,
 \vartheta^{\vecJ}_0-\delta<\theta<\vartheta^{\vecJ}_0
 \bigr\}.
\end{multline}
Put $I_1:=\openopen{0}{1}$
and $I_2:=\closedclosed{0}{1}$.
For each $\alpha\in D_{\vecJ}$,
we take an embedding
$F_{\alpha}:I_1\times I_2\lrarr
 \nbigu_{\vecJ_-,\theta^u}$
such that
(i) $F_{\alpha}(I_1\times\{0\})\subset
 \del U_{\alpha}(\epsilon)\cap \nbigu_{\vecJ_-,\theta^u}$,
(ii) $F_{\alpha}(I_1\times\{1\})
 \subset
 \nbigu_{\vecJ_-,\theta^u}\cap
 \varpi^{-1}(\infty)$,
(iii)
 $F_{\alpha}(I_1\times(I_2\setminus\del I_2))
\subset
 \nbigu_{\vecJ_-,\theta^u}\setminus
\bigl(
 \varpi^{-1}(\infty)
 \cup
 \bigcup_{\beta\in D_{\vecJ}}
 \overline{U_{\beta}}
 \bigr)$.
We may naturally regard 
$F_{\alpha}$ as a map to
$\projtilde^1_{\infty\cup D}$.

Let $Y_{\alpha,\vecJ_-}$
denote the union of $\varpi_D^{-1}(U_{\alpha})$
and $F_{\alpha}(I_1\times I_2)$
in $\projtilde^1_{\infty\cup D}$.
It is an open subset in $\projtilde^1_{\infty\cup D}$.
Let $j_{Y_{\alpha,\vecJ_-}}$ denote the inclusion
$Y_{\alpha,\vecJ_-}\lrarr 
\projtilde^1_{\infty\cup D}$.
We also have the natural inclusion
$j'_{Y_{\alpha,\vecJ_-}}:Y_{\alpha,\vecJ_-}
 \lrarr
 \projtilde^1_{\infty\cup\alpha}$.
We set
\[
 N_{\alpha}^{\varrho}\bigl(
 \nbigv\otimes
 \nbige(zu^{-1})
 \bigr)
:=
 j_{Y_{\alpha,\vecJ_-}}^{-1}
 \nbigl^{\varrho}\bigl(
 \nbigv\otimes
 \nbige(zu^{-1})
 \bigr).
\]
There exists the natural monomorphism:
\[
 j_{Y_{\alpha,\vecJ_-}!}
   N_{\alpha}^{\varrho}\bigl(
 \nbigv\otimes
 \nbige(zu^{-1})
 \bigr)
\lrarr 
 \nbigl^{\varrho}\bigl(
 \nbigv\otimes
 \nbige(zu^{-1})
 \bigr).
\]
There also exists the natural monomorphism:
\begin{equation}
\label{eq;18.5.15.20}
 j'_{Y_{\alpha,\vecJ_-}!}
  N_{\alpha}^{\varrho}\bigl(
 \nbigv\otimes
 \nbige(zu^{-1})
 \bigr)
\lrarr 
 \nbigl^{\varrho(\alpha)}\bigl(
 \nbigv_{\alpha}\otimes
 \nbige(zu^{-1})
 \bigr).
\end{equation}
The cokernel of the morphisms (\ref{eq;18.5.15.20})
is acyclic with respect to the global cohomology.
Hence, we obtain the desired morphism
$C^{\varrho}_{\vecJ_-,\alpha}$
in (\ref{eq;24.3.15.1}).
Applying the same constructions to $(V,\nabla)$,
we obtain the map
$C^{\varrho}_{\vecJ_-,\alpha}$
in (\ref{eq;24.3.15.2}).

\subsubsection{Explicit $1$-cycles in the case of $(V,\nabla)$}

Let us describe $C^{\varrho}_{\vecJ_{-},\alpha}$
for $(V,\nabla)$
in terms of explicit $1$-cycles.
For each $\alpha\in D_{\vecJ}$,
we set \index{points $\alpha(\vecJ_{\pm})$}
\[
 \alpha(\vecJ_-):=
 \alpha+\epsilon\exp(\sqrt{-1}\vartheta^{\vecJ}_{\ell}).
\]
Here,
$\epsilon$ denotes a small positive number.
Let $\gamma_{\vecJ_-,\alpha,1}$ be a path
from $\alpha(\vecJ_-)$
to $(\infty,\vartheta^{\vecJ}_0-\delta/2)$
in $\nbigu_{\vecJ_-,\theta^u}$.
Let $\gamma_{\vecJ_-,\alpha,2}$
be the path given by
$\alpha+\epsilon e^{\sqrt{-1}(\vartheta^{\vecJ}_{\ell}+t)}$
$(-2\pi\leq t\leq 0)$.
Take $v\in \nbigl_{|\alpha(\vecJ_-)}$.
We have the section $\vtilde$ of 
$\gamma_{\vecJ_-,\alpha,2}^{\ast}\nbigl$
such that $\vtilde(0)=v$.
Let $v'$ be the element of
$\nbigl_{|\alpha(\vecJ_-)}$
obtained as $\vtilde(-2\pi)$.
We have the sections $\check{v}$
and $\check{v}'$ of $\nbigl$
along $\gamma_{\vecJ_-,\alpha,1}$
induced by $v$ and $v'$, respectively.
If $\varrho(\alpha)=!$,
we obtain the following cycle
for $\nbigl^{\varrho}\bigl(V\otimes\nbige(zu^{-1})\bigr)$:
\[
 \nbigc^{\varrho}_{\vecJ_-,\alpha}(v):=
 \vtilde\otimes\gamma_{\vecJ_-,\alpha,2}
+(\check{v}-\check{v}')\otimes
 \gamma_{\vecJ_-,\alpha,1}.
\]
Let $\gamma_{\vecJ_-,\alpha,3}$ be
a path
from 
a point of $\varpi_D^{-1}(\alpha)$
to $\alpha(\vecJ_-)$
in $\varpi_D^{-1}(U_{\alpha})$.
We have the section $\hat{v}$
along $\gamma_{\vecJ_-,\alpha,3}$
induced by $v$.
If $\varrho(\alpha)=\ast$,
we obtain the following cycle
for  $\nbigl^{\varrho}\bigl(V\otimes\nbige(zu^{-1})\bigr)$:
\[
 \nbigc^{\varrho}_{\vecJ_-,\alpha}(v):=
\check{v}\otimes\gamma_{\vecJ_-,\alpha,1}
+\hat{v}\otimes\gamma_{\vecJ_-,\alpha,3}.
\]

In the case of $(V,\nabla)=(V_{\alpha},\nabla)$,
these constructions induce isomorphisms
\begin{equation}
\label{eq;24.4.2.1}
 \nbigl_{\alpha(\vecJ_-)}\simeq
H_1^{\varrho(\alpha)}\bigl(\cnum\setminus\{\alpha\},
V_{\alpha}\otimes\nbige(zu^{-1})\bigr).
\end{equation}
Under these identifications (\ref{eq;24.4.2.1}),
the cycles
$\nbigc^{\varrho}_{\vecJ_-,\alpha}(v)$
for $\nbigl^{\varrho}(V\otimes\nbige(zu^{-1}))$
represent
$C^{\varrho}_{\vecJ_-,\alpha}(v)
\in H_1^{\varrho}(\cnum\setminus D,V\otimes\nbige(zu^{-1}))$.

\subsection{The case of ``$+$''}

Suppose that $\theta^u\in \vecJ_+$.
Take $0<\delta<\!<\theta^u-\vartheta^{\vecJ}_{\ell}$.
Take a small $\epsilon>0$.
We consider the following subset
of $\projtilde^1_{\infty}$:
\begin{multline}
 \nbigu_{\vecJ_+,\theta^u}:=
 \bigl\{
 s_1e^{\sqrt{-1}\vartheta^{\vecJ}_0}
+s_2e^{\sqrt{-1}\vartheta^{\vecJ}_r}\,\big|\,
 s_1\in\real,\,\,
 0<s_2<2\epsilon
 \bigr\}
\cup \\
 \bigl\{
 re^{\sqrt{-1}\theta}\,\big|\,
 0\leq r\leq\infty,\,\,
 \vartheta^{\vecJ}_0<\theta<\vartheta_0^{\vecJ}+\delta
 \bigr\}.
\end{multline}
For each $\alpha\in D_{\vecJ}$,
we take an embedding
$F_{\alpha}:I_1\times I_2\lrarr
 \nbigu_{\vecJ_+,\theta^u}$
such that
(i) $F_{\alpha}(I_1\times\{0\})\subset
 \del U_{\alpha}\cap \nbigu_{\vecJ_+,\theta^u}$,
(ii) $F_{\alpha}(I_1\times\{1\})
 \subset
 \nbigu_{\vecJ_+,\theta^u}\cap
 \varpi^{-1}(\infty)$,
(iii)
 $F_{\alpha}(I_1\times(I_2\setminus\del I_2))
\subset
 \nbigu_{\vecJ_+,\theta^u}\setminus
\bigl(
 \varpi^{-1}(\infty)
 \cup
 \bigcup_{\beta\in D_{\vecJ}}
 \overline{U_{\beta}}
 \bigr)$.
Let $Y_{\alpha,\vecJ_+}$
denote the union of $U_{\alpha}$
and $F_{\alpha}(I_1\times I_2)$.
Let $j_{Y_{\alpha,\vecJ_+}}$
denote the inclusion
$Y_{\alpha,\vecJ_+}\lrarr 
\projtilde^1_{\infty\cup D}$.
Let $j'_{Y_{\alpha,\vecJ_+}}$
denote the inclusion
$Y_{\alpha,\vecJ_+}\lrarr 
\projtilde^1_{\infty\cup\alpha}$.
By using 
$Y_{\alpha,\vecJ_+}$ 
with embeddings
$j_{Y_{\alpha,\vecJ_+}}$
and 
$j'_{Y_{\alpha,\vecJ_+}}$
instead of
$Y_{\alpha,\vecJ_-}$
with
$j_{Y_{\alpha,\vecJ_-}}$
and 
$j'_{Y_{\alpha,\vecJ_-}}$,
we construct the morphisms
$C^{\varrho}_{\vecJ_+,\alpha}$
for $(\nbigv,\nabla)$
and $(V,\nabla)$.

Let us describe $C^{\varrho}_{\vecJ_{+},\alpha}$
for $(V,\nabla)$
in terms of explicit $1$-cycles.
We set
\[
\alpha(\vecJ_+):=
  \alpha+\epsilon\exp(\sqrt{-1}\vartheta^{\vecJ}_r).
\]
We take a path $\gamma_{\vecJ_+,\alpha,1}$
from $\alpha(\vecJ_+)$ to 
$(\infty,\vartheta_0^{\vecJ}+\delta/2)$
in $\nbigu_{\vecJ_+,\theta^u}$.
Let $\gamma_{\vecJ_+,\alpha,2}$
be the path given as
$\alpha+\epsilon e^{\sqrt{-1}(\vartheta^{\vecJ}_{r}+t)}$
$(-2\pi\leq t\leq 0)$.
We take a path $\gamma_{\vecJ_+,\alpha,3}$
from a point of $\varpi_D^{-1}(\alpha)$
to $\alpha(\vecJ_+)$
in $\varpi_D^{-1}(U_{\alpha})$.
Then, 
for $v\in \nbigl_{|\alpha(\vecJ_+)}$,
we construct cycles
$\nbigc^{\varrho}_{\vecJ_+,\alpha}(v)$
of $(V,\nabla)$
by using $\gamma_{\vecJ_+,\alpha,i}$
as in the case of ``$-$''.
These constructions induce isomorphisms
$\nbigl_{\alpha(\vecJ_+)}\simeq
H_1^{\varrho}\bigl(\cnum\setminus\{\alpha\},
 V_{\alpha}\otimes\nbige(zu^{-1})\bigr)$
in the case of $(V,\nabla)=(V_{\alpha},\nabla)$.
By these identifications,
the cycles
$\nbigc^{\varrho}_{\vecJ_+,\alpha}(v)$
represent
$C^{\varrho}_{\vecJ_+,\alpha}(v)$.

\subsection{Commutativity of the morphisms}

We obtain the following diagrams
for $\alpha\in D_{\vecJ}$
and for $\varrho\in \Dsf(D)$:
\begin{equation}
\label{eq;18.5.15.30}
\begin{CD}
  H_1^{\rd}(\cnum\setminus\{\alpha\},
 V_{\alpha}\otimes\nbige(zu^{-1})\bigr)
 @>{C^{\barshriek}_{\vecJ_{\pm},\alpha}}>>
  H_1^{\rd}\bigl(\cnum\setminus D,
 V\otimes\nbige(zu^{-1})\bigr)
 \\
 @V{d_1}VV @V{d_2}VV\\
  H_1^{\varrho(\alpha)}(\cnum\setminus\{\alpha\},
 \nbigv_{\alpha}\otimes\nbige(zu^{-1}))
 @>{C^{\varrho}_{\vecJ_{\pm},\alpha}}>>
  H_1^{\varrho}(\cnum\setminus D,
 \nbigv\otimes\nbige(zu^{-1}))
\\
 @V{d_3}VV @V{d_4}VV \\
  H_1^{\mg}(\cnum\setminus\{\alpha\},
 V_{\alpha}\otimes\nbige(zu^{-1})\bigr)
 @>{C^{\barast}_{\vecJ_{\pm},\alpha}}>>
  H_1^{\mg}\bigl(\cnum\setminus D,
 V\otimes\nbige(zu^{-1})\bigr).
\end{CD}
\end{equation}
Here, $d_i$ are the natural morphisms
in \S\ref{subsection;18.5.15.40}.

\begin{lem}
The diagram {\rm(\ref{eq;18.5.15.30})}
is commutative.
\end{lem}
\pf
By the construction,
there exist the following commutative diagrams:
\[
 \begin{CD}
j_{Y_{\alpha,\vecJ_-}!}
N_{\alpha}^{\barshriek}\bigl(
V\otimes\nbige(zu^{-1})
 \bigr)
@>>>
\nbigl^{\barshriek}\bigl(
V\otimes\nbige(zu^{-1})
 \bigr)\\
@VVV @VVV \\
j_{Y_{\alpha,\vecJ_-}!}
N_{\alpha}^{\varrho}\bigl(
\nbigv\otimes\nbige(zu^{-1})
 \bigr)
@>>>
\nbigl^{\varrho}\bigl(
\nbigv\otimes\nbige(zu^{-1})
 \bigr)\\
@VVV @VVV\\
j_{Y_{\alpha,\vecJ_-}!}
N_{\alpha}^{\barast}\bigl(
V\otimes\nbige(zu^{-1})
 \bigr)
@>>>
\nbigl^{\barast}\bigl(
V\otimes\nbige(zu^{-1})
 \bigr).
 \end{CD}
\]
We obtain the claims in the case of $-$.
The case of $+$ can be argued similarly.
\hfill\qed

\subsection{Decompositions}

For any $\alpha\in D$,
let us choose
$(\vecJ_{\alpha},\nu(\alpha))
\in T(\nbigi^{\circ})\times\{\pm\}$
satisfying the following condition.
\begin{itemize}
 \item $\theta^u\in \vecJbar_{\alpha}$
       and $\alpha\in D_{\vecJ_{\alpha}}$.
 \item If $\theta^u=\vartheta^{\vecJ_{\alpha}}_{\ell}$,
       then $\nu(\alpha)=-$.
       If $\theta^u=\vartheta^{\vecJ_{\alpha}}_{r}$,
       then $\nu(\alpha)=+$.
\end{itemize}
We obtain the following commutative diagram:
\begin{equation}
\label{eq;18.5.15.50}
 \begin{CD}
\bigoplus_{\alpha\in D}
H_1^{\rd}(\cnum\setminus\{\alpha\},
 V_{\alpha}\otimes\nbige(zu^{-1})
 )
 @>{a_1}>>
  H_1^{\rd}(\cnum\setminus D,
 V\otimes\nbige(zu^{-1}))
 \\
 @V{d_1}VV @V{d_2}VV\\
\bigoplus_{\alpha\in D}
  H_1^{\varrho}(\cnum\setminus\{\alpha\},
 \nbigv_{\alpha}\otimes\nbige(zu^{-1}))
 @>{a_2}>>
  H_1^{\varrho}(\cnum\setminus D,
 \nbigv\otimes\nbige(zu^{-1}))
\\
 @V{d_3}VV @V{d_4}VV \\
\bigoplus_{\alpha\in D}
H_1^{\mg}(\cnum\setminus\{\alpha\},
 V_{\alpha}
 \otimes\nbige(zu^{-1})
 )
 @>{a_3}>>
  H_1^{\mg}(\cnum\setminus D,
 V\otimes\nbige(zu^{-1})).
 \end{CD}
\end{equation}
Here,
$a_1$ is induced by
$C^{\barshriek}_{(\vecJ_{\alpha})_{\nu(\alpha)},\alpha}$,
$a_2$ is induced by
$C^{\varrho}_{(\vecJ_{\alpha})_{\nu(\alpha)},\alpha}$
and 
$a_3$ is induced by
$C^{\barast}_{(\vecJ_{\alpha})_{\nu(\alpha)},\alpha}$.

\begin{lem}
We obtain the following exact sequence:
\begin{multline}
\label{eq;18.5.15.51}
0\lrarr
 \bigoplus_{\alpha\in D}
   H_1^{\rd}(\cnum\setminus\{\alpha\},
 V_{\alpha}\otimes\nbige(zu^{-1}))
\stackrel{a_1+d_1}\lrarr
 \\
  H_1^{\rd}(\cnum\setminus D,
 V\otimes\nbige(zu^{-1}))
\oplus
 \bigoplus_{\alpha\in D}
   H_1^{\varrho}(\cnum\setminus\{\alpha\},
 \nbigv_{\alpha}\otimes\nbige(zu^{-1}))
\stackrel{a_2-d_2}\lrarr 
 \\
  H_1^{\varrho}(\cnum\setminus D,
 \nbigv\otimes\nbige(zu^{-1}))
\lrarr 0.
\end{multline}
\end{lem}
\pf
We obtain the following naturally defined exact sequence
from the diagrams (\ref{eq;18.5.15.50}):
\begin{multline}
0\lrarr
\bigoplus_{\alpha\in D}
j_{Y_{\alpha,\vecJ_{\pm}}}
N_{\alpha}^{\barshriek}\bigl(
V\otimes\nbige(zu^{-1})
 \bigr)
\lrarr \\
\nbigl^{\barshriek}\bigl(
V\otimes\nbige(zu^{-1})
 \bigr)
\oplus
\bigoplus_{\alpha\in D}
j_{Y_{\alpha,\vecJ_{\pm}}}
N_{\alpha}^{\varrho}\bigl(
\nbigv\otimes\nbige(zu^{-1})
 \bigr)
\lrarr \\
\nbigl^{\varrho}\bigl(
\nbigv\otimes\nbige(zu^{-1})
 \bigr)
\lrarr 0.
\end{multline}
Thus, we obtain the exactness of
(\ref{eq;18.5.15.51}).
\hfill\qed

\begin{prop}
\label{prop;18.6.21.20}
The morphisms $a_i$ $(i=1,2,3)$
in {\rm(\ref{eq;18.5.15.50})}
are isomorphisms.
\end{prop}
\pf
The claims for $a_1$ and $a_3$ are easy.
We obtain the claim for $a_2$
from the claim for $a_1$
and the exact sequence (\ref{eq;18.5.15.51}).
\hfill\qed

\section{Stokes filtrations}

\subsection{}

There exist the isomorphisms:
\[
 \gbigl^{\gbigf}_{\varrho}(\nbigv)_{|\theta^u}
 \simeq
 H_1^{\varrho}(\cnum\setminus D,\nbigv\otimes\nbige(zu^{-1})).
\]
The Stokes filtration of
$\gbigl^{\gbigf}_{\varrho}(\nbigv)_{|\theta^u}$
induces the filtration $\nbigf^{\circ\,\theta^u}$
on 
$H_1^{\varrho}(\cnum\setminus D,\nbigv\otimes\nbige(zu^{-1}))$
indexed by
the partially ordered set
$\bigl(
 \nbigi(\Fourier_+(\nbigv)),\leq_{\theta^u}
 \bigr)$.

\subsection{}

Let $\rho_{\alpha}:\proj^1\to\proj^1$
be given by $\rho_{\alpha}(z)=z+\alpha$.
There exists the isomorphism
\[
 \gbigl^{\gbigf}_{\varrho(\alpha)}(\rho_{\alpha}^{\ast}(\nbigv))_{|\theta^u}
 \simeq
 H_1^{\varrho(\alpha)} \bigl(
 \cnum^{\ast},
 \rho_{\alpha}^{\ast}(\nbigv_{\alpha})
 \otimes\nbige(zu^{-1})
 \bigr).
\]
The Stokes filtration of
$\gbigl^{\gbigf}_{\varrho(\alpha)}(\rho_{\alpha}^{\ast}(\nbigv))_{|\theta^u}$
induces a filtration
$\nbigf^{\circ \theta^u}$ of the space
$H_1^{\varrho(\alpha)} \bigl(
 \cnum^{\ast},
 \rho_{\alpha}^{\ast}\nbigv_{\alpha}
 \otimes\nbige(zu^{-1})
 \bigr)$
 indexed by
$\bigl(
\gbigf^{(0,\infty)}_+(\nbigi(\rho_{\alpha}^{\ast}\nbigv_{\alpha})),
\leq_{\theta^u}
\bigr)$.

We have the isomorphism
$\rho_{\alpha}^{\ast}\nbige(zu^{-1})
\simeq
 \nbige(zu^{-1})$
given by
$\rho_{\alpha}^{\ast}(\exp(-zu^{-1}))
=\exp(-\alpha u^{-1}) \exp(-zu^{-1})$.
It induces the following isomorphisms:
\[
 H_1^{\varrho(\alpha)}\bigl(\cnum\setminus\{\alpha\},
 \nbigv_{\alpha}
 \otimes\nbige(zu^{-1})
 \bigr)
\simeq
  H_1^{\varrho(\alpha)}\bigl(\cnum^{\ast},
 \rho_{\alpha}^{\ast}(\nbigv_{\alpha})
 \otimes\nbige(zu^{-1})
 \bigr).
\]
There exists the natural isomorphism
of the partially ordered sets
\[
\gbigf^{(0,\infty)}_+
(\nbigi(\rho_{\alpha}^{\ast}\nbigv_{\alpha}))
\simeq
\gbigf^{(\alpha,\infty)}_+(\nbigi(\nbigv_{\alpha}))
:=
 \bigl\{
 \alpha u^{-1}+\gminib\,\big|\,
 \gminib\in \gbigf^{(0,\infty)}_+
  (\nbigi(\rho_{\alpha}^{\ast}\nbigv_{\alpha}))
 \bigr\}
\]
equipped with $\leq_{\theta^u}$.
We obtain the filtrations
$\nbigf^{\circ\theta^u}$
on 
$H_1^{\varrho(\alpha)}\bigl(\cnum\setminus\{\alpha\},
 \nbigv_{\alpha}
 \otimes\nbige(zu^{-1})
 \bigr)$
 indexed by the partially ordered set
 $\bigl(
 \gbigf^{(\alpha,\infty)}_+(\nbigi(\nbigv_{\alpha})),
 \leq_{\theta^u}
 \bigr)$. 
We note that
\[
 \nbigi(\Fourier_+(\nbigv))
 =\bigsqcup_{\alpha\in D}
 \bigl(
 \gbigf^{(\alpha,\infty)}_+(\nbigi(\nbigv_{\alpha})).
\]
We obtain the filtration $\nbigf^{\prime\theta^u}$ of the space
\[
 \bigoplus_{\alpha\in D}H_1^{\varrho(\alpha)}
 \bigl(
 \cnum\setminus\{\alpha\},
 \nbigv_{\alpha}\otimes\nbige(zu^{-1})
 \bigr)
\]
indexed by
$\bigl(
\nbigi(\Fourier_+(\nbigv)),\leq_{\theta^u}
\bigr)$
from the filtrations $\nbigf^{\circ\theta^u}$
on the direct summands
$H_1^{\varrho(\alpha)}\bigl(\cnum\setminus\{\alpha\},
 \nbigv_{\alpha}\otimes\nbige(zu^{-1})
 \bigr)$.

\subsection{Isomorphisms of filtered vector spaces}
We shall prove the following proposition
in \S\ref{section;20.11.21.3}.
Note that $a_1$ and $a_3$ are the special cases of $a_2$.
\begin{prop}
\label{prop;24.3.16.10}
The isomorphism $a_2$ in {\rm(\ref{eq;18.5.15.50})}
is an isomorphism of filtered vector space.
\end{prop}

\subsection{Some canonically defined spaces}

By Proposition \ref{prop;24.3.16.10},
we obtain the following corollary.
\begin{cor}
For any $\theta^u\in\vecJ_{\pm}$,
 the maps $C^{\varrho}_{\vecJ_{\pm},\alpha}$
$(\alpha\in D_{\vecJ})$
induce the following isomorphisms of filtered vector spaces:
\[
 \bigoplus_{\alpha<_{\vecJ}0}
 \nbigl_{\alpha(\vecJ_{\pm})}
 \simeq
 H^0(\vecJ_{\pm},\gbigl^{\gbigf}_{\varrho}(V)_{\vecJ_{\pm,<0}}),
\quad
 \bigoplus_{\alpha>_{\vecJ}0}
 \nbigl_{\alpha(\vecJ_{\pm})}
 \simeq
 H^0(\vecJ_{\pm},\gbigl^{\gbigf}_{\varrho}(V)_{\vecJ_{\pm,>0}}).
\] 
We also obtain
$\nbigl_{0(\vecJ_{\pm})}
 \simeq
 H^0(\vecJ_{\pm},\gbigl^{\gbigf}_{\varrho}(V)_{\vecJ_{\pm,0}})$.
\hfill\qed
\end{cor}

\section{Fourier transform of local systems in terms of Stokes shells}
\label{subsection;18.5.15.10}

Let $D\subset\cnum$ be a finite subset.
For simplicity, we assume $0\in D$.
Let $\nbigl$ be a local system on $\cnum\setminus D$.
Set $\nbigi^{\circ}:=\{\alpha u^{-1}\,|\,\alpha\in D\}$.
Let us construct a functor
$\gbigf_{\varrho}(\nbigl)$ $(\varrho\in\Dsf(\nbigi^{\circ}))$
from $\Dsf(\nbigi^{\circ})$
to $\Shcat(\nbigi^{\circ})$.

\subsection{Preliminary}

Take $\vecJ=I(\vartheta^{\vecJ}_0,\pi/2)$.
Let $D_{\vecJ}$ be the set of $\alpha\in D$
such that 
$\alpha u^{-1}\in\nbigi^{\circ}_{\vecJ}$.
\index{set $D_{\vecJ}$}
Any element $\alpha\in D_{\vecJ}$
has the expression
$\alpha=-a\cdot \exp\bigl(
 \sqrt{-1}\vartheta^{\vecJ}_0
\bigr)$
for $a\in\real$.
Let $D_{\vecJ,>0}$ (resp. $D_{\vecJ,<0}$)
be the set of $\alpha$ such that $a>0$
(resp. $a<0$).
\index{sets $D_{\vecJ,>0}$, $D_{\vecJ,<0}$}
We have 
$\alpha\in D_{\vecJ,>0}$
(resp. $\alpha\in D_{\vecJ,<0}$)
if and only if
$\alpha u^{-1}>_{\vecJ}0$
(resp. $\alpha u^{-1}<_{\vecJ}0$).

We define the order $\leq_{\vecJ}$
on $D_{\vecJ}$
by 
$\alpha\leq_{\vecJ}\beta
\Longleftrightarrow
 \alpha u^{-1}\leq_{\vecJ}\beta u^{-1}$.
\index{order $\leq_{\vecJ}$}
 
Take $\vecJ\in T(\nbigi^{\circ})$.
For each $\alpha\in D_{\vecJ}$,
we set \index{points $\alpha(\vecJ_{\pm})$}
\[
 \alpha(\vecJ_-):=
 \alpha+\epsilon\exp(\sqrt{-1}\vartheta^{\vecJ}_{\ell}),
\quad
 \alpha(\vecJ_+):=
  \alpha+\epsilon\exp(\sqrt{-1}\vartheta^{\vecJ}_{r}).
\]
Here, $\epsilon$ denotes a small positive number.
Let $\gamma^{\alpha(\vecJ_-)}_{\alpha(\vecJ_+)}$
be the path given by
$\alpha+\epsilon e^{\sqrt{-1}(\vartheta^{\vecJ}_{\ell}+t)}$
$(0\leq t\leq \pi)$.
\index{path $\gamma^{\alpha(\vecJ_-)}_{\alpha(\vecJ_+)}$}
We have the isomorphism
$G^{\alpha(\vecJ_-)}_{\alpha(\vecJ_+)}:
 \nbigl_{|\alpha(\vecJ_-)}
\lrarr
 \nbigl_{|\alpha(\vecJ_+)}$
obtained as the parallel transport
along the path $\gamma^{\alpha(\vecJ_-)}_{\alpha(\vecJ_+)}$.
\index{isomorphism $G^{\alpha(\vecJ_-)}_{\alpha(\vecJ_+)}$}

Let $\alpha_1,\alpha_2\in D_{\vecJ}$.
We have the expressions
$\alpha_i=-a_i\exp(\sqrt{-1}\vartheta^{\vecJ}_0)$.
Suppose $\alpha_1u^{-1}>_{\vecJ}\alpha_2u^{-1}$.
We have the path
$\gamma^{\alpha_1(\vecJ_-)}_{\alpha_2(\vecJ_+)}$
from
$\alpha_1(\vecJ_-)$
to $\alpha_2(\vecJ_+)$
obtained as the union of the following.
\index{path $\gamma^{\alpha_1(\vecJ_-)}_{\alpha_2(\vecJ_+)}$}
\begin{itemize}
\item
$\gamma_1$ is the segment
connecting
$\alpha_1(\vecJ_-)$
and $\alpha_2(\vecJ_-)$.
\item
$\gamma_2$ is the path connecting 
$\alpha_2(\vecJ_-)$
and 
$\alpha_2(\vecJ_+)$
given by
$\gamma_2(t)=
 \alpha_2+e^{\sqrt{-1}(\vartheta^{\vecJ}_{\ell}-t)}$
$(0\leq t\leq \pi)$.
\end{itemize}
Let $G^{\alpha_1(\vecJ_-)}_{\alpha_2(\vecJ_+)}$
denote the isomorphism
$\nbigl_{|\alpha_1(\vecJ_-)}
\simeq
\nbigl_{|\alpha_2(\vecJ_+)}$
along the path 
$\gamma^{\alpha_1(\vecJ_-)}_{\alpha_2(\vecJ_+)}$.
\index{isomorphism $G^{\alpha_1(\vecJ_-)}_{\alpha_2(\vecJ_+)}$}

Suppose
$\vecJ-\pi<\vecJ'=I(\vartheta^{\vecJ'}_0,\pi/2)<\vecJ$.
Let $\beta=-b\exp(\sqrt{-1}\vartheta^{\vecJ'}_0)$
with $b<0$.
We consider the path
$\gamma^{\alpha((\vecJ-\pi)_+)}_{\beta(\vecJ'_-)}$
connecting
$\alpha((\vecJ-\pi)_+)$
and 
$\beta(\vecJ'_-)$
given as the union of the following.
We note $\vartheta^{\vecJ}_{\ell}=\vartheta^{\vecJ-\pi}_r$
and $\alpha(\vecJ_-)=\alpha((\vecJ-\pi)_+)$.
\index{path $\gamma^{\alpha((\vecJ-\pi)_+)}_{\beta(\vecJ'_-)}$}
\begin{itemize}
\item a path connecting
      $\alpha((\vecJ-\pi)_+)$
      and
      $\beta(\vecJ'_+)$
      on the union of the lines
     $\real e^{\sqrt{-1}\vartheta^{\vecJ}_0}
     +\epsilon e^{\sqrt{-1}\vartheta^{\vecJ}_{\ell}}$
     and
     $\real e^{\sqrt{-1}\vartheta^{\vecJ'}_0}
     +\epsilon e^{\sqrt{-1}\vartheta^{\vecJ'}_{r}}$.
 \item the path
       $\beta+\epsilon e^{\sqrt{-1}(\vartheta^{\vecJ'}_{r}+\pi s)}$
       $(0\leq s\leq 1)$.
\end{itemize}
We obtain the isomorphism
$G^{\alpha((\vecJ-\pi)_+)}_{\beta(\vecJ'_-)}:
 \nbigl_{\alpha((\vecJ-\pi)_+)}
\lrarr
 \nbigl_{\beta(\vecJ'_-)}$
as the parallel transport along the path
$\gamma^{\alpha((\vecJ-\pi)_+)}_{\beta(\vecJ'_-)}$.
\index{isomorphism $G^{\alpha((\vecJ-\pi)_+)}_{\beta(\vecJ'_-)}$}

Let $\vecJ_1\vdash\vecJ_2$ in $T(\nbigi^{\circ})$.
Let 
$\gamma^{0(\vecJ_{1+})}_{0(\vecJ_{2-})}$
be the path connecting 
$0(\vecJ_{1+})$ and $0(\vecJ_{2-})$
given by
$\epsilon \exp\bigl(
 s\vartheta^{\vecJ_2}_{\ell}+(1-s)\vartheta^{\vecJ_1}_r
 \bigr)$ $(0\leq s\leq 1)$.
\index{path $\gamma^{0(\vecJ_{1+})}_{0(\vecJ_{2-})}$}
We obtain the isomorphism
$G^{0(\vecJ_{1+})}_{0(\vecJ_{2-})}:
 \nbigl_{|0(\vecJ_{1+})}
\lrarr
 \nbigl_{|0(\vecJ_{2-})}$
obtained as the parallel transport
along the path
$\gamma^{0(\vecJ_{1+})}_{0(\vecJ_{2-})}$.
\index{isomorphism $G^{0(\vecJ_{1+})}_{0(\vecJ_{2-})}$}

For each $\alpha\in D$,
let $M_{\alpha}$ denote the automorphisms 
of $\nbigl_{|\alpha(\vecJ_{\pm})}$
obtained as the monodromy along the loop
$\alpha\mp\epsilon\exp(\sqrt{-1}(\vartheta^{\vecJ}_{\ell}+s))$
$(0\leq s\leq 2\pi)$.
\index{automorphism $M_{\alpha}$}

\subsection{Stokes-graded local systems}

For $\vecJ\in T(\nbigi^{\circ})$,
let $\nbigk^{\circ}_{\lambda_+(\vecJ),\vecJ_{\pm}}$
be the local systems on $\vecJ_{\pm}$
induced by the graded vector spaces
$\bigoplus_{\alpha\in D_{\vecJ,>0}}
 \nbigl_{|\alpha(\vecJ_{\pm})}$.
The grading and the orders
$(\nbigi_{\vecJ,>0},\leq_{\theta})$ $(\theta\in\vecJ_{\pm})$
induce a Stokes structure
$\bigl(\nbigf^{\theta^u}\,|\,\theta^u\in \vecJ_{\pm}\bigr)$
indexed by $\nbigi^{\circ}_{\vecJ,>0}$.
For $\varrho\in\Dsf(\nbigi^{\circ})$,
we have the isomorphisms
$\mu_{\varrho}:
 \bigl(
 \nbigk^{\circ}_{\lambda_+(\vecJ),\vecJ_{-}},
 \vecnbigf
 \bigr)_{|\vecJ}
\simeq
 \bigl(
 \nbigk^{\circ}_{\lambda_+(\vecJ),\vecJ_{+}},
 \vecnbigf
 \bigr)_{|\vecJ}$
induced by
\[
 \mu_{\varrho}:=
 \sum_{\alpha\in D_{\vecJ,>0}}
 G^{\alpha(\vecJ_-)}_{\alpha(\vecJ_+)}
-\sum_{
 \substack{\alpha_i\in D_{\vecJ,>0}\\
 \alpha_1>_{\vecJ}\alpha_2
 }}
 (\id-M_{\alpha_2}^{-1})^{\delta(\ast,\varrho,\alpha_2)}\circ
 G^{\alpha_1(\vecJ_-)}_{\alpha_2(\vecJ_+)}
 \circ(\id-M_{\alpha_1}^{-1})^{\delta(!,\varrho,\alpha_1)}.
\]
Here,
$\delta(\star,\varrho,\alpha)=1$ if $\varrho(\alpha u^{-1})=\star$
and $\delta(\star,\varrho,\alpha)=0$ if
$\varrho(\alpha u^{-1})\neq \star$.
By gluing
$(\nbigk^{\circ}_{\lambda_+(\vecJ),\vecJ_{\pm}},\vecnbigf)$
via the isomorphism,
we obtain the local systems with Stokes structure
$(\nbigk^{\circ}_{\varrho,\lambda_+(\vecJ),\vecJbar},\vecnbigf)$
on $\vecJbar$.

Similarly,
let $\nbigk^{\circ}_{\lambda_-(\vecJ),\vecJ_{\pm}}$
be the local systems on $\vecJ_{\pm}$
induced by
$\bigoplus_{\alpha\in D_{\vecJ,<0}}
 \nbigl_{|\alpha(\vecJ_{\pm})}$,
which are equipped with the Stokes structure
indexed by $\nbigi^{\circ}_{\vecJ,<0}$.
For $\varrho\in\Dsf(\nbigi^{\circ})$,
we have the isomorphism
$\mu_{\varrho}:
 \bigl(
 \nbigk^{\circ}_{\lambda_-(\vecJ),\vecJ_{-}},\vecnbigf
 \bigr)_{|\vecJ}
\simeq
\bigl(
 \nbigk^{\circ}_{\lambda_-(\vecJ),\vecJ_{+}},\vecnbigf
 \bigr)_{|\vecJ}$
induced by
 \[
 \mu_{\varrho}:=
 \sum_{\alpha\in D_{\vecJ,<0}}
 G^{\alpha(\vecJ_-)}_{\alpha(\vecJ_+)}
-\sum_{
 \substack{\alpha_i\in D_{\vecJ,<0}\\
 \alpha_1>_{\vecJ}\alpha_2
 }}
 (\id-M_{\alpha_2}^{-1})^{\delta(\ast,\varrho,\alpha_2)}\circ
 G^{\alpha_1(\vecJ_-)}_{\alpha_2(\vecJ_+)}
 \circ(\id-M_{\alpha_1}^{-1})^{\delta(!,\varrho,\alpha_1)}.
\]
By gluing
$(\nbigk^{\circ}_{\lambda_-(\vecJ),\vecJ_{\pm}},\vecnbigf)$
via the isomorphism,
we obtain local system with Stokes structure
$(\nbigk^{\circ}_{\varrho,\lambda_-(\vecJ),\vecJbar},\vecnbigf)$
on $\vecJbar$.

We have
$\lambda_{\pm}(\vecJ+\pi)
=\lambda_{\mp}(\vecJ)$
and
$\alpha((\vecJ+\pi)_{\pm})
=\alpha(\vecJ_{\mp})$.
At $\vartheta^{\vecJ}_r=\vartheta^{\vecJ+\pi}_{\ell}$,
we have the isomorphism
$(\nbigk^{\circ}_{\varrho,\lambda_{\mp}(\vecJ),\vecJbar},
 \vecnbigf)_{|\vartheta^{\vecJ}_r}
\simeq
(\nbigk^{\circ}_{\varrho,\lambda_{\pm}(\vecJ+\pi),\vecJbar+\pi},
 \vecnbigf)_{|\vartheta^{\vecJ+\pi}_{\ell}}$
induced by the identity on
$\nbigl_{|\alpha(\vecJ_{\mp})}$.
By gluing them,
we obtain local systems with Stokes structure
$(\nbigk^{\circ}_{\varrho,\lambda},\vecnbigf)$
$(\lambda\in[\nbigi^{\circ\ast}],\star=!,\ast)$
on $\real$.
\index{local systems with Stokes structure $(\nbigk^{\circ}_{\varrho,\lambda},\vecnbigf)$}

Let $\nbigk^{\circ}_{0,\vecJ_{\pm}}$
be the local system on
$\vecJ_{\pm}$
induced by
$\nbigl_{|0(\vecJ_{\pm})}$.
We have the isomorphism
$\nbigk^{\circ}_{0,\vecJ_{-}|\vecJ}
\simeq
\nbigk^{\circ}_{0,\vecJ_{+}|\vecJ}$
induced by 
$G^{0(\vecJ_-)}_{0(\vecJ_+)}$.
For $\vecJ_1\vdash\vecJ_2$,
we have the isomorphism
$\nbigk^{\circ}_{0,\vecJ_{1+}|\vecJ_1\cap\vecJ_2}
\simeq
\nbigk^{\circ}_{0,\vecJ_{2-}|\vecJ_1\cap\vecJ_2}$
induced by
$G^{0(\vecJ_{1+})}_{0(\vecJ_{2-})}$.
By gluing them,
we obtain a local system
$\nbigk^{\circ}_{0}$ on $\real$.
We set
$\nbigk^{\circ}_{\varrho,0}:=\nbigk^{\circ}_0$
for any $\varrho\in\Dsf(\nbigi^{\circ})$.
\index{local system $\nbigk^{\circ}_{0}$}

Thus, we obtain $2\pi\seisuu$-equivariant
Stokes graded local systems
$(\nbigk^{\circ}_{\varrho,\bullet},\vecnbigf)$
$(\varrho\in\Dsf(\nbigi^{\circ}))$
over $(\nbigi^{\circ},[\nbigi^{\circ}])$.
\index{local system with Stokes structure
$(\nbigk^{\circ}_{\varrho,\bullet},\vecnbigf)$}
For any morphism
$f_{\varrho_2,\varrho_1}:\varrho_1\lrarr\varrho_2$
in $\Dsf(\nbigi^{\circ})$,
we have the morphism
$\gbigf(f_{\varrho_2,\varrho_1}):
(\nbigk^{\circ}_{\varrho_1,\bullet},\vecnbigf)
\lrarr
(\nbigk^{\circ}_{\varrho_2,\bullet},\vecnbigf)$
induced by
$\bigoplus_{\alpha}
(\id-M_{\alpha}^{-1})^{\epsilon(\varrho_1,\varrho_2,\alpha)}$,
where
$\epsilon(\varrho_1,\varrho_2,\alpha)=1$
if
$\varrho_1(\alpha u^{-1})=!$
and $\varrho_2(\alpha u^{-1})=\ast$,
or 
$\epsilon(\varrho_1,\varrho_2,\alpha)=0$
otherwise.

We have the isomorphisms by the constructions:
\[
 K^{\circ}_{\varrho,\lambda_{+}(\vecJ),\vecJ}
\simeq
\bigoplus_{\alpha\in D_{\vecJ,>0}}
 \nbigl_{|\alpha(\vecJ_{-})}
\simeq
 \bigoplus_{\alpha\in D_{\vecJ,>0}}
 \nbigl_{|\alpha(\vecJ_+)},
\]
\[
 K^{\circ}_{\varrho,\lambda_{-}(\vecJ),\vecJ}
\simeq
 \bigoplus_{\alpha\in D_{\vecJ,<0}}
 \nbigl_{|\alpha(\vecJ_{-})}
\simeq
 \bigoplus_{\alpha\in D_{\vecJ,<0}}
 \nbigl_{|\alpha(\vecJ_{+})},
\]
\[
  K^{\circ}_{\varrho,0,\vecJ}
\simeq
 \nbigl_{|0(\vecJ_{-})}
\simeq
  \nbigl_{|0(\vecJ_{+})}.
\]
Let 
$N_{\star,\lambda,\vecJ_{\pm }}(\varrho)$ $(\star=\ast,!)$
denote the endomorphisms
of the vector spaces $K^{\circ}_{\varrho,\lambda,\vecJ}$
induced by
$\bigoplus_{\alpha}
(\id-M_{\alpha}^{-1})^{\delta(\star,\varrho,\alpha)}$
of $\bigoplus \nbigl_{|\alpha(\vecJ_{\pm})}$.

\subsection{Deformation data}
\index{tuple of morphisms $\vecnbigr^{\circ}$}

For $\vecJ\in T(\nbigi^{\circ})$,
we set 
\[
 G^{\lambda_+(\vecJ),\vecJ_-}_{\lambda_-(\vecJ),\vecJ_+}:
=\sum_{\alpha_1\in D_{\vecJ,>0}}
 \sum_{\alpha_2\in D_{\vecJ,<0}}
 G^{\alpha_1(\vecJ_-)}_{\alpha_2(\vecJ_+)}:
 \bigoplus_{\alpha_1\in D_{\vecJ,>0}}
 \nbigl_{|\alpha_1(\vecJ_-)}
\lrarr
 \bigoplus_{\alpha_2\in D_{\vecJ,<0}}
 \nbigl_{|\alpha_2(\vecJ_+)},
\]
\[
 G^{\lambda_+(\vecJ),\vecJ_-}_{0,\vecJ_+}:=
 \sum_{\alpha_1\in D_{\vecJ,>0}}
 G^{\alpha_1(\vecJ_-)}_{0(\vecJ_+)}:
 \bigoplus_{\alpha_1\in D_{\vecJ,>0}}
 \nbigl_{|\alpha_1(\vecJ_-)}
\lrarr
 \nbigl_{|0(\vecJ_+)},
\]
\[
 G^{0,\vecJ_-}_{\lambda_-(\vecJ),\vecJ_+}:=
 \sum_{\alpha_2\in D_{\vecJ,<0}}
 G^{0(\vecJ_-)}_{\alpha_2(\vecJ_+)}:
 \nbigl_{|0(\vecJ_-)}
\lrarr
 \bigoplus_{\alpha_2\in D_{\vecJ,<0}}
 \nbigl_{|\alpha_2(\vecJ_+)}.
\]
For $(\lambda_1,\lambda_2,\vecJ)\in\nbigb_2(\nbigi^{\circ})$,
we set
\[
  (\nbigr^{\circ}_{\varrho})^{\lambda_1,\vecJ_-}_{\lambda_2,\vecJ_+}:=
-N_{\ast,\lambda_2,\vecJ_+}(\varrho)\circ
  G^{\lambda_1,\vecJ_-}_{\lambda_2,\vecJ_+}\circ
N_{!,\lambda_1,\vecJ_-}(\varrho).
\]

For $(\vecJ_1,\vecJ_2)\in T_2(\nbigi^{\circ})$
with $\vecJ_1<\vecJ_2$,
we set
\begin{equation}
 G^{\vecJ_1}_{\vecJ_2}:=
 \sum_{\alpha_1\in D_{\vecJ_1,>0}}
 \sum_{\alpha_2\in D_{\vecJ_2,<0}}
 G^{\alpha_1(\vecJ_{1+})}_{\alpha_2(\vecJ_{2-})}:
 \bigoplus_{\alpha_1\in D_{\vecJ_1,>0}}
 \nbigl_{|\alpha_1(\vecJ_{1+})}
 \lrarr
 \bigoplus_{\alpha_2\in D_{\vecJ_2,<0}}
 \nbigl_{|\alpha_2(\vecJ_{2-})},
\end{equation}
{\small
\begin{equation}
\label{eq;18.5.7.1}
 G^{\vecJ_2}_{\vecJ_1}:=
\!\!\!\!\sum_{\alpha_1\in D_{\vecJ_1,<0}}
 \sum_{\alpha_2\in D_{\vecJ_2,>0}}
 G^{\alpha_2((\vecJ_2-\pi)_{+})}_{\alpha_1(\vecJ_{1-})}:
 \!\!\!
 \bigoplus_{\alpha_2\in D_{\vecJ_2,>0}}
 \!\!\!
 \nbigl_{|\alpha_2((\vecJ_2-\pi)_+)}
\lrarr \!\!\!\!
 \bigoplus_{\alpha_1\in D_{\vecJ_1,<0}}
 \nbigl_{|\alpha_1(\vecJ_{1-})}.
\end{equation}
}

Moreover, we set 
\[
 (\nbigr^{\circ}_{\varrho})^{\vecJ_1}_{\vecJ_2}:=
 N_{\ast,\lambda_-(\vecJ_2),\vecJ_{2-}}(\varrho)
 \circ
 G^{\vecJ_1}_{\vecJ_2}
 \circ
  N_{!,\lambda_+(\vecJ_1),\vecJ_{1+}}(\varrho),
\]
\[
 (\nbigr^{\circ}_{\varrho})^{\vecJ_2}_{\vecJ_1}:=
 -N_{\ast,\lambda_-(\vecJ_1),\vecJ_{1-}}(\varrho)
 \circ
 G^{\vecJ_2}_{\vecJ_1}
 \circ
  N_{!,\lambda_+(\vecJ_2-\pi),(\vecJ_{2}-\pi)_+}(\varrho).
\]

We obtain Stokes shells
$\gbigf_{\varrho}(\nbigl):=
 \bigl(
 \nbigk_{\varrho,\bullet},\vecnbigf,
 \vecnbigr^{\circ}_{\varrho}
 \bigr)$ for $(\varrho\in\Dsf(\nbigi^{\circ}))$.
\index{Stokes shells $\nbigf_{\varrho}(\nbigl)$}
For $\varrho_1\lrarr\varrho_2$ in $\Dsf(\nbigi^{\circ})$,
$\gbigf(f_{\varrho_2,\varrho_1})$
induces a morphism
$\gbigf_{\varrho_1}(\nbigl)
\lrarr
\gbigf_{\varrho_2}(\nbigl)$.
Thus, we obtain a functor
$\gbigf(\nbigl):\Dsf(\nbigi^{\circ})
\lrarr\Shcat(\nbigi^{\circ})$.
\index{Stokes shells $\gbigf_{\varrho}(\nbigl)$}

To simplify the description,
we denote
$\gbigf_{\underline{\star}}(\nbigl)
=(\nbigk_{\underline{\star},\bullet},\vecnbigf,
\vecnbigr^{\circ}_{\underline{\star}})$
by
$\gbigf_{\star}(\nbigl)
=(\nbigk_{\star,\bullet},\vecnbigf,
\vecnbigr^{\circ}_{\star})$.

\subsection{Proof of Proposition \ref{prop;24.3.31.1}}

By the construction,
we obtain the desired isomorphism
$\Locst\bigl(\gbigf_{\varrho}(\nbigl)\bigr)
\simeq 
(\gbigl^{\gbigf}_{\varrho}(V),\vecnbigf)$
by using $C^{\varrho}_{\vecJ_{\pm},\alpha}$.
\hfill\qed

\subsection{The associated graded local systems}

Let $\varpi:\cnumtilde(D)\lrarr \cnum$
be the oriented real blow up of $\cnum$ along $D$.
Let $\iota:\cnum\setminus D\lrarr \cnumtilde(D)$
be the inclusion.
We obtain the local system
$\iota_{\ast}(\nbigl)$ on $\cnumtilde(D)$.

For each $\alpha\in D$,
$\varpi^{-1}(\alpha)$
is identified with $\{e^{\sqrt{-1}\theta}\,|\,\theta\in\real\}$
by the coordinate $z-\alpha$.
We have
$\varphi_{\alpha}:
 \real\lrarr \varpi^{-1}(\alpha)$
given by $\theta\longmapsto \exp(\sqrt{-1}\theta)$.
We obtain the $2\pi\seisuu$-equivariant
local systems
$\varphi_{\alpha}^{-1}(\nbigl)$
on $\real$.

For $\lambda\in [\nbigi^{\circ}]$,
we set $D_{\lambda}:=
 \{\alpha\in D\,|\,\alpha u^{-1}\in 
 \nbigi^{\circ}_{\lambda}\}$.
We have the associated graded local systems
\[
\Gr^{\vecnbigf}(\nbigk^{\circ}_{\star,\lambda})
=\bigoplus_{\alpha \in D_{\lambda}}
 \Gr^{\vecnbigf}_{\alpha u^{-1}}
 (\nbigk^{\circ}_{\star,\lambda}).
\]
By the construction,
we have the natural $2\pi\seisuu$-equivariant
isomorphism
$b_{\alpha}:
 \varphi_{\alpha}^{-1}(\nbigl)
\simeq
  \Gr^{\vecnbigf}_{\alpha u^{-1}}
 (\nbigk^{\circ}_{\star,\lambda})$.
The induced endomorphism
\[
\varphi_{\alpha}^{-1}(\nbigl)
\simeq
\Gr^{\vecnbigf}_{\alpha u^{-1}}
(\nbigk^{\circ}_{\underline{!},\lambda})
\lrarr
\Gr^{\vecnbigf}_{\alpha u^{-1}}
(\nbigk^{\circ}_{\underline{\ast},\lambda})
\simeq
\varphi_{\alpha}^{-1}(\nbigl)
\]
is identified with $\id-M_{\alpha}^{-1}$.
 
\section{The local systems}
\label{section;24.4.1.40}

Let us give descriptions of the local systems
$\gbigl^{\gbigf}_{\varrho}(V)$.
We continue to assume $0\in D$.
Let $\nbigl$ be the local system on $\cnum\setminus D$
associated with $(V,\nabla)$.
Let $L_{\infty}=\varphi_{\infty}^{-1}(\nbigl_{|\varpi^{-1}(\infty)})$
be the $2\pi\seisuu$-equivariant local system on $\real$.
Let $M$ denote the automorphism of $L_{\infty}$.

\subsection{Basic homology classes (1)}

Let $\varpi:\projtilde^1\to\proj^1$ denote the oriented real blow up
along $\{0,\infty\}$.
We identify $\projtilde^1$ with $\realbar_{\geq 0}\times S^1$.
The points $(0,e^{\sqrt{-1}\psi})$,
$(r,e^{\sqrt{-1}\psi})$ $(0<r<\infty)$
and $(\infty,e^{\sqrt{-1}\psi})$
are denoted by
$0e^{\sqrt{-1}\psi}$,
$re^{\sqrt{-1}\psi}$ and
$\infty e^{\sqrt{-1}\psi}$, respectively.
Let $\varphi:\realbar_{\geq 0}\times \real\to \projtilde^1$
be the map given by
$\varphi(r,\theta^u)=re^{\sqrt{-1}\theta^u}$.
Let $\varphi_{\infty}:\real\to \varpi^{-1}(\infty)$
be the map given by 
$\theta^u\longmapsto \infty e^{\sqrt{-1}\theta^u}$.

Let $R$ be a sufficiently large number such that
$R>|\alpha|$ for any $\alpha\in D$.
Let $\epsilon$ denote a sufficiently small positive number.

Let $\arg(D)\subset\real$ be the set of
$\psi$ such that $e^{\sqrt{-1}\psi}=|\alpha|^{-1}\alpha$
for some $\alpha\in D\setminus\{0\}$.
\index{set $\arg(D)$}

\subsubsection{}

Let $\Gamma_{\infty,\theta^u}$ be a path
on $\realbar_{\geq R}\times\real$
connecting $(\infty,\theta^u-2\pi)$ and $(\infty,\theta^u)$.
\index{path $\Gamma_{\infty,\theta^u}$}

Let $\psi\in\arg(D)$.
Let $\delta>0$ be sufficiently small.
Let $\Gamma_{\infty,\psi,\pm,\theta^u}$
be a path connecting
$(0,\psi\pm\delta)$
and
$(\infty,\theta^u)$ on $\realbar_{\geq 0}\times\real$
obtained as the union of the following:
\index{paths $\Gamma_{\infty,\psi,\pm,\theta^u}$}
\begin{itemize}
 \item the segment connecting 
       $(0,\psi\pm\delta)$ and $(R,\psi\pm\delta)$,
 \item a path connecting
       $(R,\psi\pm\delta)$
       and $(\infty,\theta^u)$
       on $\realbar_{\geq R}\times\real$.
\end{itemize}

\subsubsection{}

For $v\in H^0(\real,L_{\infty})$,
we obtain flat sections along $\Gamma_{\infty,\theta^u}$
and $\Gamma_{\infty,\psi,\pm,\theta^u}$
which are also denoted by $v$.

Let $\Abb_{\infty,\theta^u}(v)$
denote the homology class of
$\varphi_{\ast}\bigl(
v\otimes\Gamma_{\infty,\theta^u}
\bigr)$
in $H_1^{\varrho}(\cnum\setminus D,V\otimes\nbige(zu^{-1}))$
for any $\varrho\in\Dsf(\nbigi^{\circ})$.
\index{map $\Abb_{\infty,\theta^u}$}

Let $\Abb^{\mg,(\psi,\pm)}_{\infty,\theta^u}(v)$
denote the homology class of
$\varphi_{\ast}\bigl(
v\otimes\Gamma_{\infty,\psi,\pm,\theta^u}
\bigr)$
in $H^{\varrho}_{1}(\cnum\setminus D,V\otimes\nbige(zu^{-1}))$
for any $\varrho\in\Dsf(\nbigi^{\circ})$
such that $\varrho(0)=\ast$.
\index{map $\Abb^{\mg,(\psi,\pm)}_{\infty,\theta^u}$}

\begin{lem}
\label{lem;25.2.14.4}
\label{lem;25.2.14.15}
We have
$\Abb_{\infty,\theta^u+2\pi}
=\Abb_{\infty,\theta^u}\circ M$
 and
$\Abb^{\mg,(\psi+2\pi,\pm)}_{\infty,\theta^u+2\pi}
=\Abb^{\mg,(\psi,\pm)}_{\infty,\theta^u}\circ M$.
\hfill\qed
\end{lem}

\subsection{Basic homology classes (2)}

Let $\alpha\in D\setminus\{0\}$
and $\psi\in\real$ such that
$\alpha=|\alpha|e^{\sqrt{-1}\psi}$.
Let
$(\alpha,\psi,\pm)$
denote 
$\alpha\pm\epsilon\sqrt{-1} e^{\sqrt{-1}\psi}$.
Let $\varpi_D:\projtilde^1_{D\cup\{\infty\}}\to \proj^1$
denote the oriented real blow up along $D\cup\{\infty\}$.

\subsubsection{}

Let $\gamma_{1,(\alpha,\psi,\pm),\theta^u}$
be a path
from $(\alpha,\psi,\pm)$
to $\infty e^{\sqrt{-1}\theta^u}$
in $\projtilde^1$ 
obtained as the union of the following paths.
\begin{itemize}
\item the segment
      connecting
      $\alpha\pm\epsilon\sqrt{-1}e^{\sqrt{-1}\psi}$
      and
      $Re^{\sqrt{-1}\psi}\pm \epsilon \sqrt{-1}e^{\sqrt{-1}\psi}$.
 \item the segment
       connecting
       $Re^{\sqrt{-1}\psi}\pm \epsilon \sqrt{-1}e^{\sqrt{-1}\psi}$
       and
       $Re^{\sqrt{-1}\psi}$.
 \item the path
       $R\exp\bigl(\sqrt{-1}((1-s)\psi+s\theta^u)
       \bigr)$ $(0\leq s\leq 1)$
       connecting
       $Re^{\sqrt{-1}\psi}$
       and
       $Re^{\sqrt{-1}\theta^u}$.
 \item the path
       $(1-s)^{-1}Re^{\sqrt{-1}\theta^u}$ $(0\leq s\leq 1)$
       connecting
       $Re^{\sqrt{-1}\theta^u}$
       and
       $\infty e^{\sqrt{-1}\theta^u}$.
\end{itemize}

Let $\gamma_{2,(\alpha,\psi,\pm)}$
be the loop
$\alpha\pm \epsilon \sqrt{-1}e^{\sqrt{-1}(\psi+2\pi t}$
$(-2\pi\leq t\leq 0)$.

Let $\gamma_{3,(\alpha,\psi,\pm)}$
be the segment connecting
$(\alpha,\psi,\pm)$
and a point in $\varpi^{-1}(0)$.

Let $\gamma_{4,(\alpha,\psi,\pm)}$
be the segment connecting
a point in $\varpi_{D}^{-1}(\alpha)$
and $(\alpha,\psi,\pm)$.

\subsubsection{}

Let $v\in\nbigl_{(\alpha,\psi,\pm)}$.
It induces sections along
$\gamma_{i,(\alpha,\psi,\pm)}$ $(i=1,3,4)$,
which are also denoted by $v$.
Let $\vtilde$ denote the section of $\nbigltilde$
along $\gamma_{2,(\alpha,\psi,\pm)}$
such that $\vtilde(0)=v$.
We obtain $M_{\alpha}(v)=\vtilde(2\pi)\in \nbigl_{(\alpha,\psi,\pm)}$.

\subsubsection{}
\label{subsection;24.4.1.2}

Let $\Abb_{(\alpha,\psi,\pm),\theta^u}(v)$
denote the homology class of the cycle
\[
 \vtilde\otimes \gamma_{2,(\alpha,\psi,\pm)}
+(v-M_{\alpha}^{-1}(v)) \otimes\gamma_{1,(\alpha,\psi,\pm),\theta^u}
\]
in $H_1^{\varrho}(\cnum\setminus D,V\otimes\nbige(zu^{-1}))$
for any $\varrho\in \Dsf(\nbigi^{\circ})$.
\index{maps $\Abb_{(\alpha,\psi,\pm),\theta^u}$}

Let $\Abb^{\mg}_{(\alpha,\psi,\pm),\theta^u}(v)$
denote the homology class of the cycle
\[
 v\otimes\gamma_{1,(\alpha,\psi,\pm),\theta^u}
+v\otimes\gamma_{4,(\alpha,\psi,\pm)}
\]
in $H_1^{\varrho}(\cnum\setminus D,V\otimes\nbige(zu^{-1}))$
for any $\varrho\in \Dsf(\nbigi^{\circ})$
such that $\varrho(\alpha)=\ast$.
\index{maps $\Abb^{\mg}_{(\alpha,\psi,\pm),\theta^u}$}

Let $\BB_{(\alpha,\psi,\pm),\theta^u}(v)$
denote the homology class of the cycle
\[
 \vtilde\otimes \gamma_{2,(\alpha,\psi,\pm)}
+v\otimes \gamma_{3,(\alpha,\psi,\pm)}
\]
in $H_1^{\varrho}(\cnum\setminus D,V\otimes\nbige(zu^{-1}))$
for any $\varrho\in \Dsf(\nbigi^{\circ})$
such that $\varrho(0)=\ast$.
\index{maps $\BB_{(\alpha,\psi,\pm),\theta^u}$}

Let $\BB^{\mg}_{(\alpha,\psi,\pm),\theta^u}(v)$
denote the homology class of the cycle
\[
 v\otimes \gamma_{3,(\alpha,\psi,\pm)}
+v\otimes \gamma_{4,(\alpha,\psi,\pm)}
\]
in $H_1^{\varrho}(\cnum\setminus D,V\otimes\nbige(zu^{-1}))$
for any $\varrho\in \Dsf(\nbigi^{\circ})$
such that $\varrho(\alpha)=\ast$ and $\varrho(0)=\ast$.
\index{maps $\BB^{\mg}_{(\alpha,\psi,\pm),\theta^u}$}

We have the natural identification
$\nbigl_{(\alpha,\psi+2\pi,\pm)}=\nbigl_{(\alpha,\psi,\pm)}$.
The following lemma is clear by the construction.
\begin{lem}
\label{lem;25.2.14.3}
\label{lem;25.2.14.10}
We have the following equalities:
\[
\Abb_{(\alpha,\psi+2\pi,\pm),\theta^u+2\pi}
 =\Abb_{(\alpha,\psi,\pm),\theta^u},
 \quad
\Abb^{\mg}_{(\alpha,\psi+2\pi,\pm),\theta^u+2\pi}
=\Abb^{\mg}_{(\alpha,\psi,\pm),\theta^u},
\]
\[
 \BB_{(\alpha,\psi,\pm),\theta^u}
=\BB_{(\alpha,\psi+2\pi,\pm),\theta^u}
=\BB_{(\alpha,\psi,\pm),\theta^u+2\pi},
\]
\[
\BB^{\mg}_{(\alpha,\psi,\pm),\theta^u}
=\BB^{\mg}_{(\alpha,\psi+2\pi,\pm),\theta^u}
=\BB^{\mg}_{(\alpha,\psi,\pm),\theta^u+2\pi}.
\]
\hfill\qed
\end{lem}

\subsection{Description of Homology groups
$H_1^{\varrho}(\cnum\setminus D,V\otimes\nbige(zu^{-1}))$}

\label{subsection;24.4.2.11}

For $\alpha\in D\setminus\{0\}$,
we choose $\arg(\alpha)$
such that $\alpha=|\alpha|\exp(2\pi\sqrt{-1}\arg(\alpha))$.
The following lemma is easy to see.
\begin{lem}
\label{lem;25.3.16.1}
If $\varrho(0)=!$,
the maps
$\Abb_{\infty,\theta^u}$,
$\Abb_{(\alpha,\arg(\alpha),-),\theta^u}$
$(\alpha\in D\setminus\{0\},\varrho(\alpha)=!)$
and
$\Abb^{\mg}_{(\alpha,\arg(\alpha),-),\theta^u}$
$(\alpha\in D\setminus\{0\},\varrho(\alpha)=\ast)$ 
induce the following isomorphism:
\[
 H^0(\real,L_{\infty})
 \oplus
 \bigoplus_{\alpha\in D\setminus\{0\}}
 \nbigl_{(\alpha,\arg(\alpha),-)}
 \lrarr
 H_1^{\varrho}(\cnum\setminus D,V\otimes\nbige(zu^{-1})).
\]
We also obtain such an isomorphism
by using 
$\Abb_{\infty,\theta^u}$,
$\Abb_{(\alpha,\arg(\alpha),+),\theta^u}$
$(\alpha\in D\setminus\{0\},\varrho(\alpha)=!)$
and
$\Abb^{\mg}_{(\alpha,\arg(\alpha),+),\theta^u}$
$(\alpha\in D\setminus\{0\},\varrho(\alpha)=\ast)$.
\hfill\qed
\end{lem}

\begin{lem}
If $\varrho(0)=\ast$,
the maps
$\Abb^{\mg}_{\infty,\theta^u}$,
$\BB_{(\alpha,\arg(\alpha),-),\theta^u}$
$(\alpha\in D\setminus\{0\},\varrho(\alpha)=!)$
and
$\BB^{\mg}_{(\alpha,\arg(\alpha),-),\theta^u}$
$(\alpha\in D\setminus\{0\},\varrho(\alpha)=\ast)$ 
induce the following isomorphism:
\[
 H^0(\real,L_{\infty})
 \oplus
 \bigoplus_{\alpha\in D\setminus\{0\}}
 \nbigl_{(\alpha,\arg(\alpha),-)}
 \lrarr
 H_1^{\varrho}(\cnum\setminus D,V\otimes\nbige(zu^{-1})).
\]
We also obtain such an isomorphism
by using 
$\Abb^{\mg}_{\infty,\theta^u}$,
$\BB_{(\alpha,\arg(\alpha),+),\theta^u}$
$(\alpha\in D\setminus\{0\},\varrho(\alpha)=!)$
and
$\BB^{\mg}_{(\alpha,\arg(\alpha),+),\theta^u}$
$(\alpha\in D\setminus\{0\},\varrho(\alpha)=\ast)$.
\hfill\qed
\end{lem}

\subsection{Relations among the basic homology classes}

\subsubsection{}

Let $G_{\infty}^{(\alpha,\psi,\pm)}:
\nbigl_{(\alpha,\psi,\pm)}\to H^0(\real,L_{\infty})$
denote the isomorphism
induced by the parallel transport
along the path
$|\alpha|(1-s)^{-1}e^{\sqrt{-1}\psi}\pm\epsilon\sqrt{-1}e^{\sqrt{-1}\psi}$
$(0\leq s\leq 1)$
connecting
$(\alpha,\psi,\pm)$
and $\infty e^{\sqrt{-1}\psi}\in \varpi^{-1}(\infty)$,
and the identification
$H^0(\real,L_{\infty})\simeq
(L_{\infty})_{\psi}\simeq \nbigl_{\infty e^{\sqrt{-1}\psi}}$.
\index{isomorphisms $G_{\infty}^{(\alpha,\psi,\pm)}$}

Let $G_{(\alpha,\psi,{\pm})}^{(\alpha,\psi,\mp)}:
\nbigl_{(\alpha,\psi,\mp)}
\simeq
\nbigl_{(\alpha,\psi,\pm)}$
be the isomorphism
obtained as the parallel transport along
the path
$\alpha\mp\epsilon\sqrt{-1} e^{\sqrt{-1}\psi+\sqrt{-1}\pi s}$
$(0\leq s\leq 1)$.
\index{isomorphisms $G_{(\alpha,\psi,{\pm})}^{(\alpha,\psi,\mp)}$}

Let $G_{(\beta,\psi,\pm)}^{(\alpha,\psi,\pm)}:
 \nbigl_{(\alpha,\psi,\pm)}
 \simeq
 \nbigl_{(\beta,\psi,\pm)}$
 denote the isomorphism induced by the segment
connecting 
$(\alpha,\psi,\pm)$ and $(\beta,\psi,\pm)$.
\index{isomorphisms $G_{(\beta,\psi,\pm)}^{(\alpha,\psi,\pm)}$}

\subsubsection{Rapid decay and moderate growth conditions}

\begin{lem}
\label{lem;24.4.2.10}
We obtain the following relations by the construction
if $\varrho(0)=\ast$:
\begin{equation}
\label{eq;25.2.27.1}
 \Abb_{\infty,\theta^u}
=\Abb^{\mg,(\psi,-)}_{\infty,\theta^u}
 -\Abb^{\mg,(\psi,-)}_{\infty,\theta^u-2\pi}
= \Abb^{\mg,(\psi,-)}_{\infty,\theta^u}
 -\Abb^{\mg,(\psi+2\pi,-)}_{\infty,\theta^u}\circ M^{-1},
\end{equation}
\[
 \Abb_{(\alpha,\psi,-),\theta^u}
=\BB_{(\alpha,\psi,-),\theta^u}
 +\Abb^{\mg,(\psi,-)}_{\infty,\theta^u}
 \circ
 G_{\infty}^{(\alpha,\psi,-)}\circ
 (\id-M_{\alpha}^{-1}),
\]
\[
 \Abb^{\mg}_{(\alpha,\psi,-),\theta^u}
=\BB^{\mg}_{(\alpha,\psi,-),\theta^u}
 +\Abb^{\mg,(\psi,-)}_{\infty,\theta^u}
 \circ
 G_{\infty}^{(\alpha,\psi,-)}.
\]
We also have the following relations if $\varrho(\alpha)=\ast$:
\[
 \Abb_{(\alpha,\psi,-),\theta^u}
=\Abb^{\mg}_{(\alpha,\psi,-),\theta^u}\circ(\id-M_{\alpha}^{-1}),
\]
\[
 \BB_{(\alpha,\psi,-),\theta^u}
=\BB^{\mg}_{(\alpha,\psi,-),\theta^u}\circ(\id-M_{\alpha}^{-1}).
\]
\hfill\qed
\end{lem}
 
\subsubsection{Change from $\theta^u$ to $\theta^u-2\pi$}

\begin{lem}
\label{lem;25.2.14.1}
We have
\begin{equation}
\label{eq;25.2.27.2}
  \Abb_{(\alpha,\psi,\pm),\theta^u}
 -\Abb_{(\alpha,\psi,\pm),\theta^u-2\pi}
=\Abb_{\infty,\theta^u}\circ
 G_{\infty}^{(\alpha,\psi,\pm)}\circ(\id-M_{\alpha}^{-1}),
\end{equation}
\begin{equation}
\label{eq;25.2.27.3}
  \Abb^{\mg}_{(\alpha,\psi,\pm),\theta^u}
 -\Abb^{\mg}_{(\alpha,\psi,\pm),\theta^u-2\pi}
=\Abb_{\infty,\theta^u}\circ
 G_{\infty}^{(\alpha,\psi,\pm)}.
\end{equation}
As a result, we obtain
\[
 \Abb_{(\alpha,\psi+2\pi,\pm),\theta^u}
=\Abb_{(\alpha,\psi,\pm),\theta^u}
 -\Abb_{\infty,\theta^u}\circ
  G_{\infty}^{(\alpha,\psi,\pm)}\circ(\id-M_{\alpha}^{-1}),
\]
\[
 \Abb^{\mg}_{(\alpha,\psi+2\pi,\pm),\theta^u}
=\Abb^{\mg}_{(\alpha,\psi,\pm),\theta^u}
 -\Abb_{\infty,\theta^u}\circ
  G_{\infty}^{(\alpha,\psi,\pm)}.
\] 
\hfill\qed
\end{lem}

\subsubsection{Change of $\pm$} 

\begin{lem}
\label{lem;25.2.14.2}
\label{lem;25.2.14.11}
 We obtain
\begin{multline}
\Abb_{(\alpha,\psi,-),\theta^u}
 =\Abb_{(\alpha,\psi,+),\theta^u}\circ
 G_{(\alpha,\psi,+)}^{(\alpha,\psi,-)}
\\
 +\sum_{|\beta|>|\alpha|}
 \Abb_{(\beta,\psi,+),\theta^u}\circ
 G_{(\beta,\psi,+)}^{(\beta,\psi,-)}
 \circ
 G_{(\beta,\psi,-)}^{(\alpha,\psi,-)}
 \circ (\id-M_{\alpha}^{-1}),
\end{multline}
\begin{multline}
 \Abb^{\mg}_{(\alpha,\psi,-),\theta^u}
 =\Abb^{\mg}_{(\alpha,\psi,+),\theta^u}\circ
 G_{(\alpha,\psi,+)}^{(\alpha,\psi,-)}
 +\!\!\sum_{|\beta|>|\alpha|}\!\!
 \Abb_{(\beta,\psi,+),\theta^u}\circ
 G_{(\beta,\psi,+)}^{(\beta,\psi,-)}\circ
 G_{(\beta,\psi,-)}^{(\alpha,\psi,-)}.
\end{multline}
We have the following relations:
\begin{multline}
 \BB_{(\alpha,\psi,-),\theta^u}
-\BB_{(\alpha,\psi,+),\theta^u}\circ
 \bigl(G_{(\alpha,\psi,-)}^{(\alpha,\psi,+)}\bigr)^{-1}
 =\\
-\sum_{|\beta|<|\alpha|}
 \BB_{(\beta,\psi,+),\theta^u}
 \circ
 G_{(\beta,\psi,+)}^{(\beta,\psi,-)}
 \circ
 G_{(\beta,\psi,-)}^{(\alpha,\psi,-)}
 \circ(\id-M_{\alpha}^{-1}),
 \end{multline}
\begin{multline}
 \BB^{\mg}_{(\alpha,\psi,-),\theta^u}
-\BB^{\mg}_{(\alpha,\psi,+),\theta^u}\circ
 \bigl(G_{(\alpha,\psi,-)}^{(\alpha,\psi,+)}\bigr)^{-1}
 =\\
-\sum_{|\beta|<|\alpha|}
 \BB_{(\beta,\psi,+),\theta^u}
 \circ
 G_{(\beta,\psi,+)}^{(\beta,\psi,-)}
 \circ
 G_{(\beta,\psi,-)}^{(\alpha,\psi,-)}.
 \end{multline}
We also have the following relations:
\begin{equation}
 \Abb^{\mg,(\psi,-)}_{\infty,\theta^u}
=\Abb^{\mg,(\psi,+)}_{\infty,\theta^u}
-\sum_{\arg\alpha=\psi}
\BB_{(\alpha,\psi,+),\theta^u}\circ
 G_{(\alpha,\psi,+)}^{(\alpha,\psi,-)}
 \!\circ\!
 \bigl(G_{\infty}^{(\alpha,\psi,-)}\bigr)^{-1}\!\!\!.
\end{equation}
 
For $\psi_1<\psi_2$ in $\arg(D)$,
\begin{multline}
 \Abb^{\mg,(\psi_1,-)}_{\infty,\theta^u}
 =\Abb^{\mg,(\psi_2,-)}_{\infty,\theta^u}
 \\
 \!-\!\sum_{\psi_1\leq \psi<\psi_2}
 \sum_{\arg\alpha=\psi}
\BB_{(\alpha,\psi,+),\theta^u}\circ
 G_{(\alpha,\psi,+)}^{(\alpha,\psi,-)}
 \circ
 \bigl(G_{\infty}^{(\alpha,\psi,-)}\bigr)^{-1}.
 \end{multline}
\hfill\qed
\end{lem}

\subsection{The homology classes adapted to Stokes structures}
\label{subsection;24.4.2.12}
\label{subsection;24.4.2.20}

Let $\vecJ\in T(\nbigi^{\circ})$.

\subsubsection{}
For $\alpha\in D_{\vecJ}$,
let
$G_{(\beta,\psi,+)}^{\alpha((\vecJ-\pi)_+)}:
\nbigl_{\alpha((\vecJ-\pi)_+)}
\simeq
\nbigl_{(\beta,\psi,+)}$
denote the isomorphisms
induced by the path along
the union of the lines
$\real e^{\sqrt{-1}\vartheta^{\vecJ}_0}
+\epsilon e^{\sqrt{-1}\vartheta^{\vecJ-\pi}_{r}}$
and 
$\real e^{\sqrt{-1}\psi}
+\epsilon\sqrt{-1}e^{\sqrt{-1}\psi}$. 
\index{isomorphisms $G_{(\beta,\psi,+)}^{\alpha((\vecJ-\pi)_+)}$}

Let
$G_{\infty}^{0(\vecJ_-)}:
 \nbigl_{0(\vecJ_-)}
 \simeq
 \nbigl_{\infty e^{\sqrt{-1}\vartheta^{\vecJ}_0}}
 \simeq
 (L_{\infty})_{\vartheta^{\vecJ}_0}
\simeq H^0(\real,L_{\infty})$
denote the isomorphism
induced by the path along the line
$\real e^{\sqrt{-1}\vartheta_0^{\vecJ}}
 +\epsilon e^{\sqrt{-1}\vartheta_{\ell}^{\vecJ}}$.
\index{isomorphisms $G_{\infty}^{0(\vecJ_-)}$}
 
Let
$G_{(\beta,\psi,\pm)}^{0(\vecJ_-)}$
denote the isomorphism
induced by the path obtained as the union of
the arc
\[
 \epsilon
 \exp\Bigl(
 \sqrt{-1}\bigl(
 (1-s)\vartheta^{\vecJ}_{\ell}
 +s(\psi\pm\pi/2)
 \bigr)
 \Bigr)\quad (0\leq s\leq 1)
\]
and the segment connecting
$\pm\epsilon \sqrt{-1}e^{\sqrt{-1}\psi}$
and $\beta\pm\epsilon \sqrt{-1}e^{\sqrt{-1}\psi}$. 
\index{isomorphisms $G_{(\beta,\psi,\pm)}^{0(\vecJ_-)}$}

\subsubsection{}

Let $\theta^u\in\vecJ_-$.

\begin{lem}
\label{lem;25.3.1.1}
If $\alpha<_{\vecJ}0$, we obtain
\[
 C^{\varrho}_{\vecJ_{-},\alpha}
 =
\left\{
\begin{array}{ll}
 \Abb_{(\alpha,\vartheta^{\vecJ}_0,-),\theta^u}
  & (\varrho(\alpha)=!)
  \\
 \Abb^{\mg}_{(\alpha,\vartheta^{\vecJ}_0,-),\theta^u}
  & (\varrho(\alpha)=\ast)
\end{array}
\right.
\]
\end{lem}

\begin{lem}
\label{lem;25.3.1.2}
Suppose $\alpha>_{\vecJ}0$.
If $\varrho(\alpha)=!$,
we obtain
\begin{multline}
 C^{\varrho}_{\vecJ_-,\alpha}
=\Abb_{(\alpha,\vartheta^{\vecJ}_0-\pi,+),\theta^u}
\\
 -\sum_{\vartheta^{\vecJ}_0-\pi<\psi<\vartheta^{\vecJ}_0}
 \sum_{\arg(\beta)=\psi}
 \Abb_{(\beta,\psi,-),\theta^u}
 \circ
  G_{(\beta,\psi,-)}^{(\beta,\psi,+)}
\circ
 G_{(\beta,\psi,+)}^{\alpha((\vecJ-\pi)_+)}
 \circ(\id-M_{\alpha}^{-1}).
\end{multline}
If $\varrho(\alpha)=\ast$,
we obtain
\begin{multline}
 C^{\varrho}_{\vecJ_-,\alpha}
=\Abb^{\mg}_{(\alpha,\vartheta^{\vecJ}_0-\pi,+),\theta^u}
\\
 -\sum_{\vartheta^{\vecJ}_0-\pi<\psi<\vartheta^{\vecJ}_0}
 \sum_{\arg(\beta)=\psi}
 \Abb_{(\beta,\psi,-),\theta^u}
 \circ
  G_{(\beta,\psi,-)}^{(\beta,\psi,+)}
\circ
 G_{(\beta,\psi,+)}^{\alpha((\vecJ-\pi)_+)}.
\end{multline}
If moreover $\varrho(0)=\ast$,
we obtain
\[
 C^{\varrho}_{\vecJ_-,\alpha}
 =
\left\{
\begin{array}{ll}
 \BB_{(\alpha,\vartheta^{\vecJ}_0-\pi,+),\theta^u}
 +\Abb^{(\vartheta^{\vecJ}_0,-)}_{\infty,\theta^u}
 \circ
 G_{\infty}^{\alpha((\vecJ-\pi)_+)}
 \circ(\id-M_{\alpha}^{-1})&
 (\varrho(\alpha)=!)
\\
 \BB^{\mg}_{(\alpha,\vartheta^{\vecJ}_0-\pi,+),\theta^u}
 +\Abb^{\mg,(\vartheta^{\vecJ}_0,-)}_{\infty,\theta^u}
 \circ
 G_{\infty}^{\alpha((\vecJ-\pi)_+)}  &
 (\varrho(\alpha)=\ast).
\end{array}
\right.
\]
Here,
 $G_{\infty}^{\alpha((\vecJ-\pi)_+)}:
 \nbigl_{\alpha((\vecJ-\pi)_+)}
 \simeq
 \nbigl_{\infty e^{\sqrt{-1}\vartheta_0^{\vecJ}}}
 \simeq
 (L_{\infty})_{\vartheta_0^{\vecJ}}
 \simeq H^0(\real,L_{\infty})$
denotes the isomorphism induced by the path along the line
$\real e^{\sqrt{-1}\vartheta_0^{\vecJ}}
+\epsilon e^{\sqrt{-1}\vartheta_r^{\vecJ-\pi}}$.
\hfill\qed
\end{lem}

\begin{lem}
If $\varrho(0)=!$, we obtain
\begin{equation}
 C^{\varrho}_{\vecJ_-,0}=
 \Abb_{\infty,\theta^u}\circ
  G_{\infty}^{0(\vecJ_-)}
 -\!\!\!\!\!\!\!\sum_{
 \substack{\vartheta^{\vecJ}_0-2\pi\leq \psi<\vartheta^{\vecJ}_0
 \\ \arg(\beta)=\psi
 }}\!\!\!\!\!\!\!
 \Abb_{(\beta,\psi,+),\theta^u}
 \circ
 G_{(\beta,\psi,+)}^{(\beta,\psi,-)}
 \circ
 G_{(\beta,\psi,-)}^{0(\vecJ_-)}.
 \end{equation}
We also have
\begin{equation}
 C^{\varrho}_{\vecJ_-,0}=
 \Abb_{\infty,\theta^u}\circ
 M\circ G_{\infty}^{0(\vecJ_-)}\circ M_0^{-1}
 +\!\!\!\!\!\!\!\sum_{
 \substack{
 \vartheta^{\vecJ}_0-2\pi\leq \psi<\vartheta^{\vecJ}_0
 \\
  \arg(\beta)=\psi}
 }\!\!\!\!\!\!\!
 \Abb_{(\beta,\psi,-),\theta^u}
 \circ
 G_{(\beta,\psi,-)}^{(\beta,\psi,+)}
 \circ
 G_{(\beta,\psi,+)}^{0(\vecJ_-)}.
 \end{equation}
\hfill\qed
\end{lem}

\begin{lem}
If $\varrho(0)=\ast$, we obtain
$C^{\varrho}_{\vecJ_-,0}
=\Abb^{\mg,(\vartheta^{\vecJ}_0,-)}_{\infty,\theta^u}
 \circ G_{\infty}^{0(\vecJ_-)}$.
\hfill\qed
\end{lem}

\subsubsection{}

Let $\vecJ\in T(\nbigi^{\circ})$.
Let $\theta^u\in\real$ such that
$\vartheta^{\vecJ}_0
-(\theta^u-\pi/2)$
is a sufficiently small positive number.

\begin{lem}
\label{lem;25.3.1.10}
Let $\alpha\in D_{\vecJ,<0}$
and $v\in\nbigl_{\alpha(\vecJ_-)}$.
For 
$\beta\in D_{\vecJ}\setminus\{\alpha\}$
or 
$\beta \in D_{\vecJ'}\setminus\{0\}$
$(\vecJ-2\pi<\vecJ'<\vecJ)$,
there uniquely exist $s_{\beta}(v)\in\nbigl_{\beta(\vecJ')_+}$
such that 
\begin{multline}
 \label{eq;25.3.1.4}
\Abb_{\infty,\theta^u}\bigl(
G^{(\alpha,\vartheta^{\vecJ}_0,-)}_{\infty}(v)
\bigr)
-C^{\underline{!}}_{\vecJ_-,\alpha}(v)
= \\
 \sum_{\beta\in D_{\vecJ}\setminus\{\alpha\}}
 C^{\underline{!}}_{\vecJ_+,\beta}(s_{\beta}(v))
 +\sum_{\vecJ-2\pi<\vecJ'<\vecJ}
 \sum_{\beta\in D_{\vecJ'}\setminus\{0\}}
 C^{\underline{!}}_{\vecJ'_+,\beta}(s_{\beta}(v)).
\end{multline} 
\end{lem}
\pf
There exist $s'_{\beta}$
such that
 \begin{multline}
\label{eq;25.3.1.3}
  C^{\underline{!}}_{\vecJ_-,\alpha}(v)
-\Abb_{\infty,\theta^u}\bigl(
G^{(\alpha,\vartheta^{\vecJ}_0,-)}_{\infty}(v)
\bigr)
 =\\
 \sum_{\substack{\beta\in D_{\vecJ,<0}\\ \beta\neq\alpha}}
 \Abb_{(\beta,\vartheta^{\vecJ}_0,+),\theta^u}(s'_{\beta})
 +C^{\underline{!}}_{\vecJ_+,0}(s'_0)
 +\sum_{\substack{\vartheta^{\vecJ}_0-2\pi<\psi<\vartheta^{\vecJ}_0\\
 \arg(\beta)=\psi}}
 \Abb_{(\beta,\psi,+),\theta^u-2\pi}(s'_{\beta}).
 \end{multline}
By using Lemma \ref{lem;25.3.1.1} and
Lemma \ref{lem;25.3.1.2},
we rewrite (\ref{eq;25.3.1.3}) to (\ref{eq;25.3.1.4}).
The uniqueness follows from Lemma \ref{lem;25.3.16.1}.
\hfill\qed

\subsection{The monodromy automorphisms and the induced morphisms}

There exist the following natural morphisms
\[
L_{\infty}
\stackrel{a_{\varrho}}{\lrarr}
\gbigl^{\gbigf}_{\varrho}(V)
\stackrel{b_{\varrho}}{\lrarr}
L_{\infty}.
\]
The morphism $a_{\varrho}$ equals
the morphism induced by
$\Abb_{\infty,\theta^u}$.

\begin{lem}
\label{lem;25.2.27.4}
We have the following for $v\in H^0(\real,L_{\infty})$:
\[
 b_{\varrho}\bigl(
 \Abb_{\infty,\theta^u}(v)
 \bigr)
=(\id-M^{-1})(v),
\quad\quad
 b_{\varrho}\bigl(
 \Abb^{\mg}_{\infty,\theta^u}(v)
 \bigr)
=v.
\]
For $\alpha\in D\setminus\{0\}$
and $v\in \nbigl_{(\alpha,\psi,\pm)}$, we obtain
\[
 b_{\varrho}
 \bigl(
 \Abb_{(\alpha,\psi,\pm),\theta^u}(v)
 \bigr)
=G^{(\alpha,\psi,\pm)}_{\infty}\circ
 (\id-M_{\alpha}^{-1})(v),
\quad
 b_{\varrho}
 \bigl(
 \Abb^{\mg}_{(\alpha,\psi,\pm),\theta^u}(v)
 \bigr)
=G^{(\alpha,\psi,\pm)}_{\infty}(v),
\]
\[
 b_{\varrho}
 \bigl(
 \BB_{(\alpha,\psi,\pm),\theta^u}(v)
 \bigr)
=0,
\quad\quad
 b_{\varrho}
 \bigl(
 \BB^{\mg}_{(\alpha,\psi,\pm),\theta^u}(v)
 \bigr)
=0.
\]
\hfill\qed
\end{lem}

\begin{lem}
$b_{\varrho}\circ a_{\varrho}
=\id-M^{-1}$.
\hfill\qed 
\end{lem}

Let $M^{\gbigf}_{\varrho}$ denote
the monodromy automorphism of
$\gbigl^{\gbigf}_{\varrho}(V)$.
\index{automorphism $M^{\gbigf}_{\varrho}$}
\begin{lem}
\label{lem;25.2.28.1}
$a_{\varrho}\circ b_{\varrho}=\id-(M^{\gbigf}_{\varrho})^{-1}$.
\end{lem}
\pf
In the case $\varrho(0)=!$,
the claim follows from (\ref{eq;25.2.25.2}), (\ref{eq;25.2.27.3})
and Lemma \ref{lem;25.2.27.4}.
In the case $\varrho(0)=\ast$,
the claim follows from Lemma \ref{lem;25.2.14.3}
and Lemma \ref{lem;25.2.27.4}.
\hfill\qed

\subsection{Extension}

Let
$\varphi_{\alpha}^{-1}(\nbigl)
\stackrel{a_{\alpha}}\lrarr
L^{\alpha}_{1}
\stackrel{b_{\alpha}}{\lrarr}
\varphi_{\alpha}^{-1}(\nbigl)$
be morphisms of local systems
such that
$b_{\alpha}\circ a_{\alpha}=\id-M_{\alpha}^{-1}$.
Together with the functor 
$\gbigl_{\varrho}^{\gbigf}(V)$,
we obtain the local system
$\Ltilde_1$ with the morphisms
\[
\gbigl^{\gbigf}_{\underline{!}}(V)
\lrarr
\Ltilde_1
\lrarr
\gbigl^{\gbigf}_{\underline{\ast}}(V).
\]
By the natural identifications
$\gbigl^{\gbigf}_{\star}(V_{\infty})
=L_{\infty}$,
we also obtain the following morphisms
\[
 L_{\infty}
 \stackrel{a_{\Ltilde_1}}{\lrarr}
 \Ltilde_1
 \stackrel{b_{\Ltilde_1}}{\lrarr}
 L_{\infty}.
\]
Let $M_{\Ltilde_1}$
and $M_{L^{\alpha}_1}$
are the monodromy automorphisms of $\Ltilde_1$
and $L^{\alpha}_1$, respectively.

\subsubsection{The induced endomorphisms}

\begin{prop}
\label{prop;25.3.1.20}
If
$a_{\alpha}\circ b_{\alpha}=\id-M_{L^{\alpha}_1}^{-1}$
for any $\alpha$,
we obtain 
$a_{\Ltilde_1}\circ b_{\Ltilde_1}=
\id-M^{-1}_{\Ltilde_1}$.
\end{prop}
\pf
Let $D'\subset D$ be the subset of
$\alpha\in D$
such that 
one of $a_{\alpha}$ or $b_{\alpha}$ is
an isomorphism.
If $D'=D$, the claim follows from
Lemma \ref{lem;25.2.28.1}.
We shall use an induction on $m=|D\setminus D'|$.

Let $\alpha\in D\setminus D'$.
Let $\Ltilde_{1\star}$ $(\star=!,\ast)$
denote the local systems
from
$\varphi_{\beta}^{-1}(\nbigl)
\to L_{\beta}\to\varphi_{\beta}^{-1}(\nbigl)$
$(\beta\in D\setminus\{\alpha\})$
and
\[
\varphi_{\alpha}^{-1}(\nbigl)
\stackrel{=}{\lrarr} \varphi_{\alpha}^{-1}(\nbigl)
\stackrel{b_{\alpha}\circ a_{\alpha}}
{\lrarr}
\varphi_{\alpha}^{-1}(\nbigl)
\quad (\star=!),
\]
\[
\varphi_{\alpha}^{-1}(\nbigl)
\stackrel{b_{\alpha}\circ a_{\alpha}}{\lrarr}
\varphi_{\alpha}^{-1}(\nbigl)
\stackrel{=}{\lrarr}
\varphi_{\alpha}^{-1}(\nbigl)
\quad(\star=\ast).
\]
There exist the natural morphisms
$\Ltilde_{1!}\stackrel{u_1}{\lrarr}
\Ltilde_1
\stackrel{u_2}{\lrarr}
\Ltilde_{1\ast}$.
To simplify the notation,
we set
$\atilde=a_{\Ltilde_1}$,
$\btilde=b_{\Ltilde_1}$,
$\atilde_{\star}=a_{\Ltilde_{1\star}}$,
and $\btilde_{\star}=b_{\Ltilde_{1\star}}$.
We shall identify
$\nbigi^{\circ}$ and $D$
by $\alpha u^{-1}\leftrightarrow \alpha$.

Let $\vecJ(\alpha)\in T(\nbigi^{\circ})$
such that $\alpha\in D_{\vecJ(\alpha),<0}$.
Let $\theta^u\in\real$ such that
$\vartheta^{\vecJ(\alpha)}_{0}-(\theta^u-\pi/2)$
is a sufficiently small positive number.
Let $S(\theta^u)$ denote the set of $\vecJ\in T(\nbigi)$
such that $\theta^u\in\vecJ$.
We obtain the decomposition
\begin{equation}
\label{eq;25.3.1.11}
 \Ltilde_{1|\theta^u}
 =
 \bigoplus_{\vecJ\in S(\theta^u)}
 \Bigl(
 (\Ltilde_{1})_{\vecJ_+,>0|\theta^u}
 \oplus
 (\Ltilde_{1})_{\vecJbar,<0|\theta^u}
 \Bigr)
 \oplus
 (\Ltilde_{1})_{\vecJ(\alpha)_+,0|\theta^u}.
\end{equation}
Let $q_{\vecJ_+,>0}$, $q_{\vecJbar,<0}$
and $q_{\vecJ(\alpha)_+,0}$
denote the projections onto the corresponding components.
We have the decomposition
\[
 (\Ltilde_1)_{\vecJbar(\alpha),<0|\theta^u}
 =\bigoplus_{\beta\in D_{\vecJ(\alpha),<0}\setminus\{\alpha\}}
 (\Ltilde_1)_{\vecJ(\alpha)_+,\beta|\theta^u}
 \oplus
 (\Ltilde_1)_{\vecJ(\alpha)_-,\alpha|\theta^u}.
\]
Let $q_{\beta}$ denote the projection onto the component.
If $q_{\alpha}\circ q_{\vecJ(\alpha)_+,<0}(s)=0$,
then there exists
$s'\in \Ltilde_{1!}$
such that $s=u_1(s')$.
By using the assumption of the induction,
we obtain
$(\atilde_!\circ\btilde_!)(s')
=(\id-M^{-1}_{\Ltilde_{1!}})(u_1(s'))$.
We can easily check that
$(\atilde\circ\btilde)(u_1(s'))
=(\id-M^{-1}_{\Ltilde_{1}})(u_1(s'))$.

We consider
$t_{\alpha}\in
(\Ltilde_1)_{\vecJ(\alpha)_-,\alpha|\theta^u}
\subset
\Ltilde_{1|\theta^u}$.
Note that
\[
 \atilde_{\ast}\circ\btilde_{\ast}\bigl(
 u_2(t'_{\alpha})
 \bigr)
=(\id-M^{-1}_{\Ltilde_{1\ast}})(u_2(t'_{\alpha})).
\]
It implies that
$\atilde\circ\btilde\bigl(
t'_{\alpha}
\bigr)
=(\id-M^{-1}_{\Ltilde_{1}})(t'_{\alpha})$
except for the $q_{\alpha}\circ q_{\vecJ(\alpha)_+,<0}$-component.
We obtain
\[
q_{\alpha}\circ
q_{\vecJbar(\alpha),<0}\bigl(
M_{\Ltilde_1}^{-1}(t_{\alpha})
\bigr)
=M_{L^{\alpha}_1}^{-1}\bigl(
 q_{\alpha}(t_{\alpha})
\bigr).
\]
We identify
$(\Ltilde_{1\star})_{\vecJ(\alpha)_-,\alpha|\theta^u}
=\nbigl_{\alpha(\vecJ_-)}$.
By Lemma \ref{lem;25.3.1.10},
we have
\[
 q_{\alpha}\circ
 q_{\vecJbar(\alpha),<0}\bigl(
 \atilde_!\circ\btilde(t_{\alpha})
 \bigr)
 =b_{\alpha}(t_{\alpha}).
\]
Hence, we obtain
\[
 q_{\alpha}\circ
 q_{\vecJbar(\alpha),<0}\bigl(
 \atilde\circ\btilde(t_{\alpha})
 \bigr)
 =a_{\alpha}\circ b_{\alpha}(t_{\alpha}).
\]
Then, we obtain the desired equality for
the $q_{\alpha}\circ q_{\vecJ(\alpha)_+,<0}$-components.
\hfill\qed

\subsubsection{The recovery of $\nbigl$}
\label{subsection;25.3.17.3}

Let us observe that
the local system $\nbigl$ is recovered from
$L_{\infty}\to \Ltilde_1\to L_{\infty}$.

Let $\alpha\in D\setminus\{0\}$.
We consider $\theta^j$
and the decomposition (\ref{eq;25.3.1.11})
in the proof of Proposition \ref{prop;25.3.1.20}.
We obtain the following morphisms
\[
 L_{\infty|\theta^u}
 \stackrel{\atilde}{\lrarr}
 \Ltilde_{1|\theta^u}
 \stackrel{c_1(\alpha)}{\lrarr}
 (\Ltilde_1)_{\vecJ(\alpha)_-,\alpha|\theta^u}
 \stackrel{c_2(\alpha)}{\lrarr}
 \Ltilde_{1|\theta^u}
 \stackrel{\btilde}{\lrarr}
 L_{\infty|\theta^u},
\]
where $c_1(\alpha)$ denotes the projection,
and $c_2(\alpha)$ denotes the inclusion.
Under the isomorphism
$G^{\alpha(\vecJ(\alpha)_-)}_{\infty}:
\nbigl_{\alpha(\vecJ(\alpha)_-)}
\simeq
L_{\infty,\theta^u}$,
we have
\[
M_{\alpha}^{-1}
=
\id-
\btilde\circ c_2(\alpha)\circ c_1(\alpha)\circ\atilde.
\]

\subsection{Local systems}

We translate the results
in \S\ref{subsection;24.4.2.11}--\S\ref{subsection;24.4.2.12}.
\subsubsection{Local system $\gbigl^{\gbigf}_{\underline{!}}(V)$}

We consider the vector space
\begin{equation}
\label{eq;24.4.1.1}
 H^0(\real,L_{\infty})
 \oplus
 \bigoplus_{\pm}
 \bigoplus_{\alpha\in D\setminus\{0\}}
 \bigoplus_{e^{\sqrt{-1}\psi}=\alpha/|\alpha|}
 \nbigl_{(\alpha,\psi,\pm)}.
\end{equation}
An element $v\in\nbigl_{(\alpha,\psi,\pm)}$
is denoted by
$\langle(\alpha,\psi,\pm),v\rangle$.
For any $\theta^u\in\real$,
$\gbigl^{\gbigf}_{\underline{!}}(V)_{|\theta^u}$
is identified with the quotient space of 
(\ref{eq;24.4.1.1})
by the equivalence relation generated by
(\ref{eq;24.4.1.10}) and (\ref{eq;24.4.1.11}) below
(see Lemma \ref{lem;25.2.14.1} and Lemma \ref{lem;25.2.14.2}):
\begin{equation}
\label{eq;24.4.1.10}
 \langle(\alpha,\psi+2\pi,\pm),v\rangle
=\langle
 (\alpha,\psi,\pm),v
 \rangle
+G_{\infty}^{(\alpha,\psi,\pm)}\bigl(
v-M_{\alpha}(v)\bigr),
\end{equation}
\begin{multline}
\label{eq;24.4.1.11}
 \langle(\alpha,\psi,-),v\rangle
=\langle (\alpha,\psi,+),G_{(\alpha,\psi,+)}^{(\alpha,\psi,-)}(v)
 \rangle
\\
+\sum_{|\beta|>|\alpha|}
 \langle (\beta,\psi,+),
 \bigl(
 G_{(\beta,\psi,-)}^{(\beta,\psi,+)}
 \bigr)^{-1}
 \circ G_{(\beta,\psi,-)}^{(\alpha,\psi,-)}
 (v-M_{\alpha}(v))
 \rangle.
\end{multline}
The $2\pi\seisuu$-action
is induced by the monodromy automorphism $M$
on $H^0(\real,L_{\infty})$
and the natural shift
$\langle (\alpha,\psi+2\pi,\pm),v
\rangle
\longmapsto
\langle (\alpha,\psi,\pm),v
\rangle$
(see Lemma \ref{lem;25.2.14.4} and Lemma \ref{lem;25.2.14.3}).

\subsubsection{Local system $\gbigl^{\gbigf}_{\underline{\ast}}(V)$}

We consider the following vector space
\begin{equation}
\label{eq;24.4.1.3}
\bigoplus_{\pm}
\bigoplus_{\arg(D)}
H^0(\real,L_{\infty})
\oplus
\bigoplus_{\pm}  
\bigoplus_{\alpha\in D\setminus\{0\}}
\bigoplus_{e^{\sqrt{-1}\psi}=|\alpha|^{-1}\alpha}
\nbigl_{(\alpha,\psi,\pm)}.
\end{equation}
An element of $H^0(\real,L_{\infty})$
corresponding to the $(\pm,\psi)$-component
is denoted by
$\langle \pm,\psi,w\rangle^{\mg}$.
An element of
$\nbigl_{(\alpha,\psi,\pm)}$
is denoted by
$\langle (\alpha,\psi,\pm),v\rangle^{\mg}$.
For any $\theta^u\in\real$,
the space $\gbigl^{\gbigf}_{\underline{\ast}}(V)_{|\theta^u}$
is identified with
the quotient of (\ref{eq;24.4.1.3})
by the equivalence relation generated by
(\ref{eq;24.4.1.23}),
(\ref{eq;24.4.1.22}),
(\ref{eq;24.4.1.20})
and
(\ref{eq;24.4.1.21})
(see Lemma \ref{lem;25.2.14.10},
Lemma \ref{lem;25.2.14.11}):
\begin{equation}
\label{eq;24.4.1.23}
 \langle
  (\alpha,\psi+2\pi,\pm),
  v\rangle
  =\langle
  (\alpha,\psi,\pm),v
  \rangle.
\end{equation}
\begin{multline}
\label{eq;24.4.1.22}
 \langle (\alpha,\psi,-),v\rangle^{\mg}
 -\langle (\alpha,\psi,+),(G_{(\alpha,\psi,-)}^{(\alpha,\psi,+)})^{-1}(v)
 \rangle^{\mg}
 =\\
+\sum_{|\beta|<|\alpha|}
 \langle
 (\beta,\psi,+),
 (\id-M_{\beta})\circ
 \bigl(
 G_{(\beta,\psi,-)}^{(\beta,\psi,+)}
 \bigr)^{-1}
 \circ
 G_{(\beta,\psi,-)}^{(\alpha,\psi,-)}(v)
 \rangle^{\mg}.
\end{multline}
\begin{multline}
\label{eq;24.4.1.20}
 \langle -,\psi,w\rangle^{\mg}
 -\langle +,\psi,w\rangle^{\mg}
 =\\
 +\sum \langle (\alpha,\psi,+),
 (\id-M_{\alpha})\circ
 \bigl(
 G_{(\alpha,\psi,-)}^{(\alpha,\psi,+)}
 \bigr)^{-1}
 \circ
 (G_{\infty}^{(\alpha,\psi,-)})^{-1}(v)
 \rangle^{\mg}.
\end{multline}
For $\psi_1<\psi_2$,
\begin{multline}
\label{eq;24.4.1.21}
 \langle -,\psi_1,w\rangle^{\mg}
 -\langle -,\psi_2,w\rangle^{\mg}=
 \\
 +\sum_{\psi_1\leq \psi<\psi_2} \langle (\alpha,\psi,+),
 (\id-M_{\alpha})\circ
 \bigl(
 G_{(\alpha,\psi,-)}^{(\alpha,\psi,+)}
 \bigr)^{-1}
 \circ
 (G_{\infty}^{(\alpha,\psi,-)})^{-1}(v)
 \rangle^{\mg}.
\end{multline}
The $2\pi\seisuu$-action is induced by
the automorphism of (\ref{eq;24.4.1.3})
obtained from
the monodromy $M$ on $H^0(\real,L_{\infty})$
and the shift
$\langle (\alpha,\psi+2\pi,\pm),v\rangle^{\mg}
\longmapsto
\langle (\alpha,\psi,\pm),v\rangle^{\mg}$.
(see Lemma \ref{lem;25.2.14.15} and Lemma \ref{lem;25.2.14.10}).

\subsubsection{Morphism}

The morphism
$\gbigl^{\gbigf}_{\underline{!}}(V)\to
\gbigl^{\gbigf}_{\underline{\ast}}(V)$
is described as
\begin{equation}
 \langle
 (\alpha,\psi,-),v
 \rangle
 \longmapsto
 \langle
 \psi,-,G_{\infty}^{(\alpha,\psi,-)}\circ(\id-M_{\alpha})(v)
 \rangle^{\mg}
 +\langle
  (\alpha,\psi,-),
  (\id-M_{\alpha})(v)
 \rangle^{\mg}, 
\end{equation}
and 
\begin{equation}
\label{eq;24.4.1.30}
 w\longmapsto
 \langle
 \psi,-,w
 \rangle^{\mg}
- \langle
  \psi+2\pi,-,
  M^{-1}(w)
 \rangle^{\mg}.
\end{equation}
(See Lemma \ref{lem;24.4.2.10}.)
The image of the right hand side (\ref{eq;24.4.1.30})
in the quotient of (\ref{eq;24.4.1.3})
is independent of $\psi$.

\begin{rem}
We can explicitly describe the isomorphisms
$\gbigl^{\gbigf}_{\underline{\star}}(V)
 \simeq
 \Loc^{\St}(\gbigf_{\underline{\star}}(\nbigl(V)))$
by the relations in {\rm\S\ref{subsection;24.4.2.20}}.
\hfill\qed
\end{rem}

\chapter{Local Fourier transform and reductions at $\infty$}
\label{section;20.10.24.4}

\section{Introduction to \S\ref{section;20.10.24.4}}
\label{section;25.2.5.1}

Let $D$ be a finite subset in $\cnum$.
We set $\Dtilde=D\cup\{\infty\}$.
Let $(\nbigv,\nabla)$ be a meromorphic flat bundle on
$(\proj_z^1,\Dtilde)$.
Let $\nbigi_{\infty}(\nbigv)$
denote the set of ramified irregular values of $\nbigv$ at $\infty$.
When $\nbigi_{\infty}(\nbigv)\neq \{0\}$,
we set
\[
\omega:=
\min\{-\ord(\gminia)\,|\,\gminia\in\nbigi_{\infty}(\nbigv)\setminus\{0\}\}
=\min\bigl\{
 \omega'\,\big|\,
     \nbigstilde_{\omega'}(\nbigi_{\infty}(\nbigv))
     \neq\nbigi_{\infty}(\nbigv)
 \bigr\}.
\]
We study the case $1<\omega$.
(See Proposition \ref{prop;24.3.17.121} for the case $\omega\leq 1$.)

Let $U$ be a small neighbourhood of $\infty$ in $\proj_z^1$
such that $D\cap U=\emptyset$.
We obtain the meromorphic flat bundle
$(V_{\infty},\nabla):=
\nbigttilde^{\infty}_{\omega}(\nbigv,\nabla)$ on $(U,0)$.
(See \S\ref{subsection;18.6.23.3} for
$\nbigttilde^{\infty}_{\omega}(\nbigv,\nabla)$.)
It extends to a meromorphic flat bundle
on $(\proj^1,\{0,\infty\})$
with regular singularity at $0$.
The extended bundle is also denoted by 
$(V_{\infty},\nabla)$.
We set $\nbigitilde:=\nbigi_{\infty}(V_{\infty})=
\nbigttilde_{\omega}(\nbigi_{\infty}(\nbigv))$.
Note $\nbigs^{\infty}_{\omega}(\nbigitilde)=\nbigitilde$.

For any $\varrho\in\Dsf(D)$,
let $\bigl(
 \gbigl^{\gbigf}_{\varrho}(\nbigv),\vecnbigf
 \bigr)$
denote the local system with Stokes structure
corresponding to $\Fourier_+(\nbigv(\varrho))$ at $\infty$. 
We obtain the functor
$\Dsf(D)\to \Loc^{\St}(\nbigi(\Fourier_+(\nbigv)))$
given by
$\varrho\longmapsto
(\gbigl^{\gbigf}_{\varrho}(\nbigv),\vecnbigf)$.
We also obtain the local systems with Stokes structure
$(\gbigl^{\gbigf}_{\star}(V_{\infty}),\vecnbigf)$ $(\star=!,\ast)$.

\subsection{Reduction of $\gbigl^{\gbigf}_{\varrho}(\nbigv)$}
\label{subsection;24.3.26.100}

We set $\omega^{\circ}=(\omega-1)^{-1}\omega$.
\begin{thm}
\label{thm;24.3.29.40}
For any $\varrho\in\Dsf(D)$,
there exists 
the isomorphism of local systems with Stokes structure
$\nbigt_{\omega^{\circ}}\bigl(
 \gbigl^{\gbigf}_{\varrho}(\nbigv),\vecnbigf
 \bigr)
\simeq
 \bigl(
 \gbigl^{\gbigf}_{\varrho}(\nbigstilde^{\infty}_{\omega}\nbigv),
 \vecnbigf
 \bigr)$.
They induce an isomorphism of functors from $\Dsf(D)$
to the category of local systems with Stokes structure,
i.e., for any $\varrho_1\to\varrho_2$ in $\Dsf(D)$,
the following diagram is commutative:
\[
\begin{CD}
  \nbigt_{\omega^{\circ}}\bigl(
 \gbigl^{\gbigf}_{\varrho_1}(\nbigv),\vecnbigf
 \bigr)
@>>>
  \nbigt_{\omega^{\circ}}\bigl(
 \gbigl^{\gbigf}_{\varrho_2}(\nbigv),\vecnbigf
 \bigr)
\\
 @V{\simeq}VV @V{\simeq}VV \\ 
 \bigl(
 \gbigl^{\gbigf}_{\varrho_1}(\nbigstilde^{\infty}_{\omega}\nbigv),
 \vecnbigf
 \bigr)
 @>>>
 \bigl(
 \gbigl^{\gbigf}_{\varrho_2}(\nbigstilde^{\infty}_{\omega}\nbigv),
 \vecnbigf
 \bigr).
\end{CD}
\]
\end{thm}

When $\nbigv=V_{\infty}$,
the theorem says that
the morphism of the local systems
$\nbigt_{\omega^{\circ}}\gbigl^{\gbigf}_!(V_{\infty})
\to \nbigt_{\omega}^{\circ}\gbigl^{\gbigf}_{\ast}(V_{\infty})$
is identified with
$\gbigl^{\gbigf}_!(V_{\infty}^{\reg})
\to
\gbigl^{\gbigf}_{\ast}(V_{\infty}^{\reg})$,
where
$V_{\infty}^{\reg}=\nbigstilde^{\infty}_{\omega}(V_{\infty})
=\nbigttilde^{\infty}_{\omega}\nbigstilde^{\infty}_{\omega}(\nbigv)$
is regular singular at $\{0,\infty\}$.
It also directly follows from the stationary phase formula
in \S\ref{subsection;24.3.26.1}.

\begin{thm}
\label{thm;24.3.29.41}
The functor
from $\Dsf(D)$ to the category of
$2\pi\seisuu$-equivariant local systems with Stokes structure
$\nbigs_{\omega^{\circ}}\bigl(
\gbigl^{\gbigf}_{\varrho}(\nbigv),\vecnbigf
 \bigr)$
is obtained as the extension of
$(\gbigl^{\gbigf}_{!}(V_{\infty}),\vecnbigf)
 \to
(\gbigl^{\gbigf}_{\ast}(V_{\infty}),\vecnbigf)$
by the following natural morphisms of
$2\pi\seisuu$-equivariant local systems:
\begin{equation}
\label{eq;24.3.26.30}
 \gbigl^{\gbigf}_!(V^{\reg}_{\infty})
 \to
 \gbigl^{\gbigf}_{\varrho}(\nbigstilde^{\infty}_{\omega}(\nbigv))
 \to
 \gbigl^{\gbigf}_{\ast}(V^{\reg}_{\infty}).
\end{equation}
\end{thm}

\subsection{Stokes structure of
$(\gbigl^{\gbigf}_{\star}(V_{\infty}),\vecnbigf)$}
\label{subsection;24.4.5.122}

It is fundamental for us to study
$(\gbigl^{\gbigf}_{\star}(V_{\infty}),\vecnbigf)$.
Let $(L,\vecnbigftilde)$ denote the $2\pi\seisuu$-equivariant
local system with Stokes structure indexed by $\nbigitilde$
on $\real$
corresponding to $(V_{\infty},\nabla)$ at $\infty$.
We shall give two types of explicit descriptions
of $(\gbigl^{\gbigf}_{\star}(V_{\infty}),\vecnbigf)$.

\subsubsection{Local system with Stokes structure}

In \S\ref{subsection;24.4.5.120},
from $(L,\vecnbigftilde)$,
we shall explicitly construct
$2\pi\seisuu$-equivariant
local systems with Stokes structure
$\gbigf^{(\infty,\infty)}_{+,\star}(L,\vecnbigftilde)=
(\gbigq^{\infty}_{\star}(V_{\infty})_{\real},\vecnbigf)$ $(\star=!,\ast)$
and morphisms of local systems
\index{local systems with Stokes structure
$\gbigf^{(\infty,\infty)}_{+,\star}(L,\vecnbigftilde)$}
\[
c^{-1}(\nbigt_{\omega}(L))\to \gbigq^{\infty}_!(V_{\infty})_{\real}
\to \gbigq^{\infty}_{\ast}(V_{\infty})_{\real}
\to
c^{-1}(\nbigt_{\omega}(L)).
\]
Here,
$c:\real\to\real$ be the map
defined by $c(\theta^u)=-\theta^u$.
\index{map $c$}
\index{local systems $\gbigq^{\infty}_{\star}(V_{\infty})_{\real}$}

\begin{thm}
\label{thm;24.4.5.121}
There exists the following commutative diagram of
$2\pi\seisuu$-equivariant local systems with Stokes structure:
\[
 \begin{CD}
  \gbigf^{(\infty,\infty)}_{+,!}(L,\vecnbigftilde)
  @>{F_{\gbigq^{\infty}}}>>
  \gbigf^{(\infty,\infty)}_{+,\ast}(L,\vecnbigftilde)
  \\
  @V{\simeq}VV @V{\simeq}VV \\
  (\gbigl^{\gbigf}_!(V_{\infty}),\vecnbigftilde)
  @>>>
  (\gbigl^{\gbigf}_{\ast}(V_{\infty}),\vecnbigftilde).
 \end{CD}
\]
We also have the following commutative diagram of
the local systems
\[
\begin{CD}
 c^{-1}(\nbigt_{\omega}(L)) @>>> \gbigq^{\infty}_!(L,\vecnbigftilde)_{\real}
 @>{F_{\gbigq^{\infty}}}>>
 \gbigq^{\infty}_{\ast}(L,\vecnbigftilde)_{\real}
  @>>>
  c^{-1}(\nbigt_{\omega}(L)) \\
  @V{\simeq}VV @V{\simeq}VV @V{\simeq}VV @V{\simeq}VV \\
  \gbigl^{\gbigf}_!(\nbigt^{\infty}_{\omega}(V_{\infty}))
  @>>>
  \gbigl^{\gbigf}_!(V_{\infty})
  @>>>
  \gbigl^{\gbigf}_{\ast}(V_{\infty})
  @>>>
  \gbigl^{\gbigf}_{\ast}(\nbigt^{\infty}_{\omega}(V_{\infty})).
\end{CD}  
\]
In the diagrams, the lower horizontal arrows are
the natural morphisms.
\end{thm}

\subsubsection{Stokes shells of
$(\gbigl^{\gbigf}_{\star}(V_{\infty}),\vecnbigf)$}

In \S\ref{subsection;20.11.14.20},
we introduce an explicit construction of
a base tuple of Stokes shells
$\bigl(
 \gbigf^{(\infty,\infty)}_{+,!}(\Shsf(L,\vecnbigftilde)),
 \gbigf^{(\infty,\infty)}_{+,\ast}(\Shsf(L,\vecnbigftilde)),
 F
\bigr)$
in $\Shcat\bigl(\gbigf^{(\infty,\infty)}_+(\nbigi_{\infty}(V_{\infty}))
\cup\{0\}\bigr)$
from any Stokes shell $(L,\vecnbigftilde)$.
\begin{prop}
\label{prop;24.3.26.31}
There exists the following commutative diagram
of Stokes shells: 
\begin{equation}
\label{eq;24.3.26.23}
\begin{CD}
  \gbigf_{+,!}^{(\infty,\infty)}(\Shsf(L,\vecnbigftilde))
  @>{F}>>
 \gbigf_{+,\ast}^{(\infty,\infty)}(\Shsf(L,\vecnbigftilde))
 \\
  @V{\simeq}VV @V{\simeq}VV \\
  \Shsf\bigl(
  \gbigf^{(\infty,\infty)}_{+,!}(L,\vecnbigftilde)
  \bigr)
  @>>>
  \Shsf\bigl(
  \gbigf^{(\infty,\infty)}_{+,\ast}(L,\vecnbigftilde)
  \bigr).
 \end{CD}
\end{equation}
As a result, the base tuple
$\Shsf\bigl(
  \gbigl^{\gbigf}_!(V_{\infty}),\vecnbigf
 \bigr)
 \to
  \Shsf\bigl(
  \gbigl^{\gbigf}_{\ast}(V_{\infty}),\vecnbigf
  \bigr)$
can be identified with 
$\gbigf_{+,!}^{(\infty,\infty)}(\Shsf(L,\vecnbigftilde))
 \to
 \gbigf_{+,\ast}^{(\infty,\infty)}(\Shsf(L,\vecnbigftilde))$.
\end{prop}

\subsection{Inductive procedure}
\label{subsection;24.4.2.130}

These theorems provide us with the following procedure
to study $(\gbigl^{\gbigf}_{\varrho}(\nbigv),\vecnbigf)$ inductively.
\begin{itemize}
 \item $(\gbigl^{\gbigf}_{\varrho}(\nbigv),\vecnbigf)$
are recovered from
$\nbigs_{\omega^{\circ}}\bigl(
\gbigl^{\gbigf}_{\varrho}(\nbigv),\vecnbigf
 \bigr)$
and
\[
      \nbigt_{\omega^{\circ}}(\gbigl^{\gbigf}_{\varrho}(\nbigv),\vecnbigf)
       \simeq
       (\gbigl^{\gbigf}_{\varrho}(\nbigstilde^{\infty}_{\omega}(\nbigv)),
       \vecnbigf). 
\]
       Note that either
       $\min\bigl\{
       -\ord(\gminia)\big|\,
       \gminia\in
       \nbigi_{\infty}(\nbigstilde^{\infty}_{\omega}(\nbigv))
       \bigr\}>\omega$
       or
       $\nbigi_{\infty}(\nbigstilde^{\infty}_{\omega}(\nbigv))=\{0\}$
       holds.
       If $\nbigi_{\infty}(\nbigstilde^{\infty}_{\omega}(\nbigv))=\{0\}$,
       we may apply the results in \S\ref{section;20.10.30.2}
       to study
       $\gbigl_{\varrho}^{\gbigf}(\nbigstilde^{\infty}_{\omega}(\nbigv))$.
 \item
      $\nbigs_{\omega^{\circ}}\bigl(
\gbigl^{\gbigf}_{\varrho}(\nbigv),\vecnbigf
      \bigr)$
      is explicitly described as
      the extension of the base tuple
\[
 \begin{CD}
  (\gbigq^{\infty}_!(L,\vecnbigftilde)_{\real},\vecnbigf)
  @>>>
  (\gbigq^{\infty}_{\ast}(L,\vecnbigftilde)_{\real},\vecnbigf)
 \end{CD}
\]
by (\ref{eq;24.3.26.30}).
\end{itemize}

As the complement to this procedure,
we note that
the morphisms of the local systems
\[
 \gbigl^{\gbigf}_!(\nbigt^{\infty}_{\omega}(V_{\infty}))
 \lrarr
 \gbigl^{\gbigf}_!(V_{\infty})
 \lrarr
 \gbigl^{\gbigf}_{\ast}(V_{\infty})
 \lrarr
 \gbigl^{\gbigf}_{\ast}(\nbigt_{\omega}^{\infty}(V_{\infty}))
\]
are explicitly described
by Theorem \ref{thm;24.4.5.121}.
It allows us to describe explicitly
the morphisms of local systems
\[
 \gbigl^{\gbigf}_!(\nbigt^{\infty}_{\omega}(V_{\infty}))
 \lrarr
 \gbigl^{\gbigf}_{\varrho}(\nbigv)
 \lrarr
 \gbigl^{\gbigf}_{\ast}(\nbigt^{\infty}(V_{\infty})).
\]
We remark that
$\nbigt^{\infty}_{\omega}(V_{\infty})
=\nbigt^{\infty}_{\omega}(\nbigv)$,
and it is regular singular at $\{0,\infty\}$.

\subsection{Homology groups}

Let $u=|u|\exp(\sqrt{-1}\theta^u)\in\cnum^{\ast}$.
When $|u|$ is sufficiently small,
there exist the natural isomorphisms
\[
 \gbigl^{\gbigf}_{\varrho}(\nbigv)_{|\theta^u}
 \simeq
 H_1^{\varrho}\bigl(
 \cnum\setminus D,
 \nbigv\otimes\nbige(zu^{-1})
 \bigr).
\]
There also exist the following isomorphisms:
\[
 \gbigl^{\gbigf}_{!}(V_{\infty})_{|\theta^u}
 \simeq
 H_1^{\rd}\bigl(
 \cnum^{\ast},
 V_{\infty}\otimes\nbige(zu^{-1})
 \bigr),
 \quad
  \gbigl^{\gbigf}_{\ast}(V_{\infty})_{|\theta^u}
 \simeq
 H_1^{\mg}\bigl(
 \cnum^{\ast},
 V_{\infty}\otimes\nbige(zu^{-1})
 \bigr).
\]
To obtain the theorems in
\S\ref{subsection;24.3.26.100}--\S\ref{subsection;24.4.5.122},
we shall study these homology groups.

Set $\nbigi:=\pi_{\omega\ast}(\nbigitilde)$.
We set
$\vecI_x(\theta^u)=\openopen{-\theta^u+\pi/2}{-\theta^u+3\pi/2}$.

\subsubsection{Homology group
$H_1^{\rd}(\cnum^{\ast},V_{\infty}\otimes\nbige(zu^{-1}))$}

In \S\ref{subsection;24.4.2.31}--\S\ref{subsection;24.2.23.10},
we shall introduce the following maps
\[
 \Abb^{\rd}_{\infty,\theta^u}:
 H^0(\real,\nbigt_{\omega}(L))
 \lrarr
 H_1^{\rd}\bigl(\cnum^{\ast},
 V_{\infty}\otimes\nbige(zu^{-1})\bigr),
\]
\[
 \BB^{\rd}_{J,\theta^u}:
 H^0(J,L_{J,>0})
 \lrarr
 H_1^{\rd}\bigl(
 \cnum^{\ast},
 V_{\infty}\otimes\nbige(zu^{-1})
 \bigr)
 \quad
 (J\in T(\nbigi)).
\]
They induce the isomorphism
\[
 H_1^{\rd}\bigl(
 \cnum^{\ast},
 V_{\infty}\otimes\nbige(zu^{-1})
 \bigr)
 \simeq
 \Bigl(
 H^0(\real,\nbigt_{\omega}(L))
 \oplus
 \bigoplus_{J\in T(\nbigi)}
 H^0(J,L_{J,>0})
 \Bigr)\Big/\!\!\sim.
\]
(See \S\ref{subsection;24.3.26.20}
for the equivalence relation.)
The $2\pi\seisuu$-action is also defined naturally
on the right hand side.
This gives the isomorphism of $2\pi\seisuu$-equivariant
local systems
$\gbigl^{\gbigf}_{!}(V_{\infty})
\simeq
 \gbigq^{\infty}_{!}(V_{\infty})_{\real}$
in Theorem \ref{thm;24.4.5.121}.
 
To study the Stokes structure,
in \S\ref{subsection;24.2.22.20},
we shall introduce maps 
\[
 A_{J_+,\theta^u}:
 H^0(J,L_{J,<0})
 \lrarr
 H_1^{\rd}\bigl(
 \cnum^{\ast},V_{\infty}\otimes\nbige(zu^{-1})
 \bigr)
\]
for $J\in T(\nbigi)$ such that $J_+\subset \vecI_x(\theta^u)$,
and 
\[
 A_{J_-,\theta^u}:
 H^0(J,L_{J,<0})
 \lrarr
 H_1^{\rd}\bigl(
 \cnum^{\ast},V_{\infty}\otimes\nbige(zu^{-1})
 \bigr)
\]
for $J\in T(\nbigi)$ such that $J_-\subset \vecI_x(\theta^u)$.
We shall also construct
\begin{equation}
 A^{J_{1+}}_{\infty,\theta^u}:
 H^0(\real,L)
\lrarr
H_1^{\rd}\bigl(\cnum^{\ast},
V_{\infty}\otimes
 \nbige(zu^{-1})\bigr)
\end{equation}
for $J_1\in T(\nbigi)$
such that $J_{1+}\subset\vecI_x(\theta^u)-\pi$,
and
\begin{equation}
 A^{J_{2-}}_{\infty,\theta^u}:
 H^0(\real,L)
\lrarr
 H_1^{\rd}\bigl(\cnum^{\ast},
 V_{\infty}\otimes\nbige(zu^{-1})\bigr)
\end{equation}
for $J_2\in T(\nbigi)$
such that $J_{2-}\subset\vecI_x(\theta^u)-\pi$.
Then, we obtain the isomorphism of the vector spaces
(\ref{eq;18.4.22.100})
(Proposition \ref{prop;18.4.20.31}).
The both sides of (\ref{eq;18.4.22.100})
are equipped with the filtrations
indexed by
$\bigl(\nbigi_{\infty}(\Fourier_+(V_{\infty})),\leq_{\theta^u}\bigr)$.
As in Theorem \ref{thm;24.3.16.20},
they are isomorphisms of filtered vector spaces.
(The proof of Theorem \ref{thm;24.3.16.20} will be given in
\S\ref{subsection;18.5.23.1}.)
This gives the isomorphism
$(\gbigq^{\infty}_!(V_{\infty})_{\real},\vecnbigf)
\simeq
 (\gbigl^{\gbigf}_!(V_{\infty}),\vecnbigf)$
in Theorem \ref{thm;24.4.5.121}.
It also provides us with the following isomorphisms of
the filtered vector spaces:
\[
 H^0(\nu_0^-(\vecJ),L_{\nu_0^-(\vecJ),<0})
 \simeq
 H^0(\vecJ_{\mp},\gbigl^{\gbigf}_{!}(V_{\infty})_{\vecJ_{\mp},>0}),
\]
\[
 H^0(\nu_0^+(\vecJ),L_{\nu_0^+(\vecJ),>0})
 \simeq
 H^0(\vecJ,\gbigl^{\gbigf}_{!}(V_{\infty})_{\vecJ,<0}),
\] 
\[
 H^0(\real,L)
 \simeq
 H^0(\vecJ_{\pm},\gbigl^{\gbigf}_{!}(V_{\infty})_{\vecJ_{\pm,0}}).
\]
By the relation among
$\BB^{\rd}_{J,\theta^u}$,
$A_{J_{\pm},\theta^u}$ $(J\in T(\nbigi))$
and
$A^{J_{1\pm}}_{\infty,\theta^u}$,
we obtain that
the Stokes shell of
$(\gbigl^{\gbigf}_!(V_{\infty}),\vecnbigf)$
is isomorphic to
$\gbigf^{(\infty,\infty)}_{!}(\Shsf(L,\vecnbigftilde))$
as in Proposition \ref{prop;24.3.26.31}.

\subsubsection{Homology group
$H_1^{\mg}(\cnum^{\ast},V_{\infty}\otimes \nbige(zu^{-1}))$}

We shall introduce 
\[
 \Abb^{\mg}_{J,\theta^u}:
 H^0(J,L_{J,<0})
 \lrarr
   H_1^{\mg}\bigl(
  \cnum^{\ast},
  V\otimes\nbige(x^{-1}u^{-1})
  \bigr),
  \quad (J\in T(\nbigi)),
\]
\[
 \Abb^{\mg,J_{\pm}}_{\infty,\theta^u}:
 H^0(\real,\nbigt_{\omega}(L))
 \lrarr
    H_1^{\mg}\bigl(
  \cnum^{\ast},
  V\otimes\nbige(x^{-1}u^{-1})
  \bigr),
  \quad
  (J\in T(\nbigi)).
\]
They induce the isomorphism
\[
 H_1^{\mg}(\cnum^{\ast},V_{\infty}\otimes\nbige(zu^{-1}))
 \simeq
 \Bigl(
 \bigoplus_{\pm}
 \bigoplus_{J\in T(\nbigi)}
 H^0(\real,\nbigt_{\omega}(L))
 \oplus
 \bigoplus_{J\in T(\nbigi)}
 H^0(J,L_{J,<0})
 \Bigr)
 \Big/\!\!\sim.
\]
(See \S\ref{subsection;24.3.26.21}
for the equivalence class.)
The $2\pi\seisuu$-action is naturally defined
on the right hand side.
This gives the isomorphism
of the $2\pi\seisuu$-equivariant local systems
$\gbigl^{\gbigf}_{\ast}(V_{\infty})
\simeq
\gbigq^{\infty}_{\ast}(V_{\infty})_{\real}$
in Theorem \ref{thm;24.4.5.121}.

To study the Stokes structure,
we shall introduce
\begin{equation}
 A^{\mg,J_{1+}}_{\infty,\theta^u}:
 H^0(\real,L)
\lrarr
H_1^{\rd}\bigl(\cnum^{\ast},
V_{\infty}\otimes
 \nbige(zu^{-1})\bigr)
\end{equation}
for $J_1\in T(\nbigi)$
such that $J_{1+}\subset\vecI_x(\theta^u)-\pi$,
and
\begin{equation}
 A^{\mg,J_{2-}}_{\infty,\theta^u}:
 H^0(\real,L)
\lrarr
 H_1^{\rd}\bigl(\cnum^{\ast},
 V_{\infty}\otimes\nbige(zu^{-1})\bigr)
\end{equation}
for $J_2\in T(\nbigi)$
such that $J_{2-}\subset\vecI_x(\theta^u)-\pi$.
Then, we obtain the isomorphism of vector spaces
(\ref{eq;18.5.26.100})
(Proposition \ref{prop;24.3.29.10}).
The both sides of (\ref{eq;18.5.26.100})
are equipped with the filtrations
indexed by
$\bigl(\nbigi_{\infty}(\Fourier_+(V_{\infty})),\leq_{\theta^u}\bigr)$.
As in Theorem \ref{thm;24.3.16.20},
they are isomorphisms of filtered vector spaces,
which will be proved in \S\ref{section;20.11.21.3}.
It gives the isomorphism
$(\gbigq^{\infty}_{\ast}(V_{\infty})_{\real},\vecnbigf)
\simeq
 (\gbigl^{\gbigf}_{\ast}(V_{\infty}),\vecnbigf)$.
It also provides us with the following isomorphisms of
the filtered vector spaces:
\[
 H^0(\nu_0^-(\vecJ),L_{\nu_0^-(\vecJ),<0})
 \simeq
 H^0(\vecJ_{\mp},\gbigl^{\gbigf}_{\ast}(V_{\infty})_{\vecJ_{\mp},>0}),
\]
\[
 H^0(\nu_0^+(\vecJ),L_{\nu_0^+(\vecJ),>0})
 \simeq
 H^0(\vecJ,\gbigl^{\gbigf}_{\ast}(V_{\infty})_{\vecJ,<0}),
\] 
\[
 H^0(\real,L)
 \simeq
 H^0(\vecJ_{\pm},\gbigl^{\gbigf}_{\ast}(V_{\infty})_{\vecJ_{\pm,0}}).
\]
By the relation among
$\BB_{J,\theta^u}$,
$A_{J_{\pm},\theta^u}$ $(J\in T(\nbigi))$
and
$A^{\mg,J_{1\pm}}_{\infty,\theta^u}$,
we obtain that
the Stokes shell of
$(\gbigl^{\gbigf}_{\ast}(V_{\infty}),\vecnbigf)$
is isomorphic to
$\gbigf^{(\infty,\infty)}_{\ast}(\Shsf(L,\vecnbigftilde))$
as in Proposition \ref{prop;24.3.26.31}.

\subsubsection{Homology groups
$H^{\varrho}_1(\cnum\setminus D,\nbigv\otimes\nbige(zu^{-1}))$}

In \S\ref{subsection;24.3.29.20},
we shall construct
\begin{equation}
C^{J_{1+}}_{\infty,\theta^u}:
H_1^{\varrho}\bigl(
 \cnum\setminus D,
 \nbigstilde^{\infty}_{\omega}(\nbigv)\otimes\nbige(u^{-1}z)
 \bigr)
\lrarr
H_1^{\varrho}\bigl(
 \cnum\setminus D,
 \nbigv\otimes\nbige(u^{-1}z)
 \bigr),
\end{equation}
for $J_1\in T(\nbigi)$
such that $J_{1+}\subset\vecI_x(\theta^u)-\pi$,
and
\begin{equation}
C^{J_{1-}}_{\infty,\theta^u}:
H_1^{\varrho}\bigl(
 \cnum\setminus D,
 \nbigstilde^{\infty}_{\omega}(\nbigv)\otimes\nbige(u^{-1}z)
 \bigr)
\lrarr
H_1^{\varrho}\bigl(
 \cnum\setminus D,
 \nbigv\otimes\nbige(u^{-1}z)
 \bigr),
\end{equation}
for $J_1\in T(\nbigi)$
such that $J_{1-}\subset \vecI_x(\theta^u)-\pi$.
We obtain the isomorphism of vector spaces
(\ref{eq;24.2.22.23}).
The both hand sides of
(\ref{eq;24.2.22.23})
are equipped with the filtrations
indexed by
$\bigl(
 \nbigi(\Fourier_+(\nbigv)),\leq_{\theta^u}
 \bigr)$.
We shall prove that
(\ref{eq;24.2.22.23})
are isomorphisms of filtered vector spaces
(Theorem \ref{thm;24.3.16.20}).
It implies Theorem \ref{thm;24.3.29.40}
and Theorem \ref{thm;24.3.29.41}.

\subsection{Some notation}
\label{subsection;25.2.5.2}

\subsubsection{}
Let $\psi:\proj^1_x\to\proj^1$ be defined by $\psi(x)=x^{-1}$.
We set $D'=\psi^{-1}(D)$
and $\Dtilde'=D'\cup\{0\}$.
\index{sets $D'$ and $\Dtilde'$}
We obtain the meromorphic flat bundle
$(\nbigv',\nabla)=\psi^{\ast}(\nbigv,\nabla)$
on $(\proj^1_x,\Dtilde')$.
\index{map $\varrho'$}
We also obtain
$(V,\nabla)=\psi^{\ast}(V_{\infty},\nabla)$
on $(\proj^1_x,\{0,\infty\})$.

Let $\nbige(u^{-1}x^{-1})$
denote the meromorphic flat bundle
$\bigl(
 \nbigo_{\proj^1}(\ast 0),d+d(u^{-1}x^{-1})
 \bigr)$.
Let $\varrho'\in\Dsf(\Dtilde')$ be defined by
$\varrho'(P)=\varrho(\psi_x(P))$
and $\varrho'(0)=\ast$.
We study
\[
 H_1^{\varrho'}(\proj^1\setminus D',\nbigv'\otimes\nbige(x^{-1}u^{-1})),
 \quad
 H_1^{\kappa}(\cnum^{\ast},V\otimes\nbige(x^{-1}u^{-1}))
 \quad(\kappa=\rd,\mg).
\]

For  $\theta=\arg(x)$,
we have $\Re(x^{-1}u^{-1})<0$
if and only if 
$\theta\in \bigcup_{m\in\seisuu}(\vecI_x(\theta^u)+2m\pi)$.
We have $\Re(x^{-1}u^{-1})>0$
if and only if 
$\theta\in
 \bigcup_{m\in\seisuu}
 \bigl(\vecI_x(\theta^u)+(2m+1)\pi\bigr)$.

\subsubsection{}
We set $X:=\realbar_{\geq 0}\times\real$
and $X^{\ast}:=\real_{>0}\times\real$.
\index{spaces $X$ and $X^{\ast}$}
For any subset $Z\subset X$,
let $q_Z:Z\lrarr\real$ denote the projection,
and let $\iota_Z$ denote the inclusion
$Z\lrarr X$.

Let $\varpi:\projtilde^1\lrarr\proj^1$
denote the oriented real blow up
along $\{0,\infty\}$.
\index{oriented real blow up $\projtilde^1$}
We identify $\projtilde^1$ with
$\realbar_{\geq 0}\times S^1$
by using the coordinate $x$.
Let $\varphi:X\lrarr \projtilde^1$
be given by
$\varphi(r,\theta)=(r,e^{\sqrt{-1}\theta})$.
Let $\varphi_1:\real\lrarr S^1$ be 
given by $\varphi_1(\theta)=e^{\sqrt{-1}\theta}$,
which is identified with the restriction of
$\varphi$ to $\real\times\{0\}$.
\index{maps $\varphi$, $\varphi_1$}
For any subset $A\subset\real$,
let $a_A:A\to \real$ denote the inclusion.
\index{map $a_A$}

\subsubsection{}
\label{subsection;25.2.6.4}
Let $L^{<0}\subset L$ and $L^{\leq 0}$
be the $2\pi\seisuu$-equivariant
constructible subsheaves determined by
$(L^{<0})_{\theta}=\nbigf^{\theta}_{<0}$
and 
$(L^{\leq 0})_{\theta}=\nbigf^{\theta}_{\leq 0}$.
We obtain the constructible subsheaves
$L^{<0}_{S^1}\subset L^{\leq 0}_{S^1}
\subset L_{S^1}$
on $\varpi^{-1}(0)$
as the descent.

We have the meromorphic flat bundle
$\pi_{\omega\ast}(V,\nabla)$
corresponding to $(L,\vecnbigf)$.
Let $\vecnbigf^F$ be the Stokes structure of $L$ 
corresponding to 
$\pi_{\omega\ast}(V,\nabla)\otimes \nbige(x^{-1}u^{-1})$.
\index{Stokes structure $\vecnbigf^F$}
Let $L^{F\,<0}\subset L$ denote the constructible subsheaf
determined by
$(L^{F\,<0})_{\theta}
=\nbigf^{F,\theta}_{<0}(L_{\theta})$.
\index{constructible sheaf $L^{F\,<0}$}
Note that 
$\nbigf^{F,\theta}_{<0}=\nbigf^{F,\theta}_{\leq 0}$.
We obtain the constructible subsheaf
$L^{F\,<0}_{S^1}\subset L_{S^1}$.
Note that
$L_{S^1}^{<0}\subset L_{S^1}^{F\,<0}
\subset
 L_{S^1}^{\leq 0}$.
The cokernel
$L_{S^1}^{F\,< 0}\big/L_{S^1}^{<0}$
is isomorphic to
$\varphi_{1!}
 \bigl(
 a_{(\vecI_x(\theta^u)-\pi)!}
 \nbigt_{\omega}(L)_{|\vecI_x(\theta^u)-\pi}
\bigr)$.

Let $\nbigl$ be the local system on $\projtilde^1$
corresponding to $(V,\nabla)$.
The restriction $\nbigl_{|\varpi^{-1}(0)}$
is identified with $L_{S^1}$.
We have the natural
$2\pi\seisuu$-equivariant isomorphism
$\varphi^{-1}(\nbigl)\simeq q_X^{-1}(L)$.

\section{Rapid decay homology group of
$(V,\nabla)\otimes\nbige(x^{-1}u^{-1})$}
\label{subsection;24.4.2.30}

\subsection{Exact sequence}
\label{subsection;24.2.21.6}

We set
$\cnumtilde=\projtilde^1\setminus\varpi^{-1}(\infty)$,
which is identified with $\real_{\geq 0}\times \varpi^{-1}(0)$.
Let $q_0:\cnumtilde\lrarr \varpi^{-1}(0)$
denote the projection,
and let $j_0:\cnumtilde\lrarr \projtilde^1$ denote the inclusion.
Let $q_1:\cnum^{\ast}\to \varpi^{-1}(0)$ denote the projection,
and let $j_1:\cnum^{\ast}\to\projtilde^1$ denote the inclusion.

\subsubsection{}

Let $\nbign_0$ denote the constructible subsheaf of
$\nbigl^{<0}(V\otimes\nbige(x^{-1}u^{-1}))_{|\cnumtilde}$
determined by the following conditions:
\[
 \nbign_{0|\varpi^{-1}(0)}=L_{S^1}^{F\,<0},
 \quad\quad
 \nbign_{0|\cnum^{\ast}}=q_1^{-1}(L_{S^1}^{\leq 0}).
\]
There exists the following exact sequence:
\begin{equation}
 0\lrarr
  j_{0!}q_0^{-1}(L_{S^1}^{<0})
  \lrarr
  j_{0!}\nbign_0
  \lrarr
  \nbigl^{<0}(\nbigt_{\omega}(V)\otimes\nbige(x^{-1}u^{-1}))
  \lrarr 0.
\end{equation}
Because $j_{0!}q_0^{-1}(L_{S^1}^{<0})$
is acyclic with respect to the global cohomology,
we obtain
\[
 H^1\bigl(\projtilde^1,
 j_{0!}\nbign_0
 \bigr)
 \simeq
 H_1^{\rd}\bigl(
 \cnum^{\ast},
 \nbigt_{\omega}(V)\otimes
 \nbige(x^{-1}u^{-1})
 \bigr).
\]

\subsubsection{}

Let $\delta>0$.
For any $J\in T(\nbigi)$,
let $\gamma_J$ be a path connecting
$(1,\vartheta^J_{r}+\delta)$ and $(1,\vartheta^J_{\ell}-\delta)$.
For any $v\in H^0(J,L_{J,>0})$,
we obtain the section $v\otimes\gamma_J$
of
$\nbigc^{-1}_{\proj^1,\del\proj^1}
\otimes j_{1!}q_1^{-1}(a_{J\ast}L_{J,>0})$.
It induces an isomorphism
\begin{equation}
\label{eq;24.2.21.3}
 H^0(J,L_{J,>0})
 \simeq
 \hyperh^{-1}\bigl(
 \projtilde^1,
 \nbigc^{\bullet}_{\projtilde^1,\del\projtilde^1}
 \otimes
 j_{1!}q_1^{-1}\bigl(
 a_{J\ast}L_{J,>0}
 \bigr)
 \bigr).
\end{equation}
We shall identify them by this isomorphism.

\subsubsection{}
\label{subsection;25.2.9.1}
There exists the following exact sequence:
\begin{equation}
 0\lrarr
 j_{0!}\nbign_0
 \lrarr
 \nbigl^{<0}(V\otimes\nbige(x^{-1}u^{-1}))
 \lrarr
 j_{1!}
 q_1^{-1}
 \bigl(L_{S^1}/L_{S^1}^{\leq 0}\bigr)
 \lrarr 0.
\end{equation}
Let $\gbigt(\nbigi,\theta^u)$
denote the set of $J\in T(\nbigi)$
such that
$-\theta^u-\pi/2< \vartheta^J_{\ell}
\leq -\theta^u+3\pi/2$.
There exists the natural isomorphism
\[
  j_{1!}
 q_1^{-1}
 \bigl(L_{S^1}/L_{S^1}^{\leq 0}\bigr)
 =\bigoplus_{J\in \gbigt(\nbigi,\theta^u)}
 j_{1!}q_1^{-1}\bigl(
 a_{J\ast}L_{J,>0}
 \bigr).
\]
We obtain the following exact sequence:
\begin{multline}
\label{eq;24.2.21.4}
 0\lrarr
 H_1^{\rd}\bigl(
 \cnum^{\ast},
 \nbigt_{\omega}(V)\otimes
 \nbige(x^{-1}u^{-1})
 \bigr)
 \stackrel{c_{1,u}}{\lrarr}
 H_1^{\rd}\bigl(
 \cnum^{\ast},
 V\otimes
 \nbige(x^{-1}u^{-1})
 \bigr)\\
 \stackrel{c_{2,u}}\lrarr
 \bigoplus_{J\in\gbigt(\nbigi,\theta^u)}
 H^0(J,L_{J,>0})\lrarr 0.
\end{multline}

\subsection{Description of
 $H_1^{\rd}\bigl(\cnum^{\ast},
 \nbigt_{\omega}(V)
 \otimes\nbige(x^{-1}u^{-1})
 \bigr)$}

Let $\Gamma_{\theta^u}$ be a path on $(X,X^{\ast})$
connecting $(0,-\theta^u+2\pi)$ and $(0,-\theta^u)$.
\index{path $\Gamma_{\theta^u}$}
For $v\in H^0(\real,\nbigt_{\omega}(L))$,
we obtain the rapid decay $1$-cycle
$\varphi_{\ast}(v\otimes\Gamma_{\theta^u})$ of
$\nbigt_{\omega}(V)\otimes\nbige(x^{-1}u^{-1})$.
This procedure induces an isomorphism,
depending on the choice of $\theta^u$:
\begin{equation}
\label{eq;24.2.21.40}
 H^0(\real,\nbigt_{\omega}(L))\simeq
 H_1^{\rd}\bigl(\cnum^{\ast},
 \nbigt_{\omega}(V)
 \otimes\nbige(x^{-1}u^{-1})
 \bigr).
\end{equation}
Let $M_0$ denote
the automorphism of $H^0(\real,\nbigt_{\omega}(L))$
obtained as the monodromy of $\nbigt_{\omega}(L)$.
\begin{lem}
$\varphi_{\ast}(v\otimes \Gamma_{\theta^u-2\pi})
=\varphi_{\ast}\bigl(
 M_0(v)\otimes\Gamma_{\theta^u}
 \bigr)$
in  
$H_1^{\rd}\bigl(\cnum^{\ast},
 \nbigt_{\omega}(V)
 \otimes\nbige(x^{-1}u^{-1})
 \bigr)$.
\hfill\qed
\end{lem}

\subsection{Rapid decay classes $\Abb^{\rd}_{\infty,\theta^u}(v)$}
\label{subsection;24.4.2.31}

Let us describe $c_{1,u}$ in terms of
$1$-cycles.
Take $a_1\in S_0(\nbigi)\cap (\vecI_x(\theta^u)-\pi)$
and $a_2\in S_0(\nbigi)\cap(\vecI_x(\theta^u)+\pi)$.
Let $a_1=b_0<b_1<\cdots<b_N=a_2$
be the set
$S_0(\nbigi)\cap \closedclosed{a_1}{a_2}$.
We set $J_i=\openopen{b_i-\omega^{-1}\pi}{b_i}$ $(i=0,\ldots,N-1)$
and $I_i=\{1\}\times\openopen{b_{i}}{b_{i+1}}$.
Let $\gamma_{i}$ be a path connecting
$(1,b_i)$ and a point in $\{0\}\times J_i$.
For $v\in H^0(\real,\nbigt_{\omega}(L))$,
we obtain
$v_{i}\in H^0(J_{i+},L_{J_{i+},0})\simeq H^0(\real,\nbigt_{\omega}(L))$
$(i=0,\ldots,N)$.
The induced sections of $\varphi^{\ast}(\nbigl)$
are also denoted by $v_{i}$.
Note that
$v_i-v_{i-1}\in H^0(J_{i},L_{J_i,<0})$
for $i=1,\ldots,N$.
We obtain the following rapid decay $1$-cycle
of $V\otimes\nbige(x^{-1}u^{-1})$:
\[
\varphi_{\ast}\Bigl(
 v_0\otimes \gamma_0
-\sum_{i=0}^{N-1} v_i\otimes I_i
-\sum_{i=1}^{N-1}(v_{i}-v_{i-1})\otimes \gamma_i
-v_{N-1}\otimes\gamma_{N}
\Bigr).
\]
The homology class
is denoted by
$\Abb^{\rd}_{\infty,\theta^u}(v)$.
\index{map $\Abb^{\rd}_{\infty,\theta^u}$}
It equals $c_{1,u}(v)$
under the identification (\ref{eq;24.2.21.40}).
\begin{lem}
\label{lem;24.2.23.4}
 $\Abb^{\rd}_{\infty,\theta^u-2\pi}
 =\Abb^{\rd}_{\infty,\theta^u}\circ M_0$
on $H^0(\real,\nbigt_{\omega}(L))$. 
\hfill\qed
\end{lem}

\subsection{Rapid decay classes $\BB^{\rd}_{J,\theta^u}(v)$}
\label{subsection;24.2.23.10}

For any $J\in T(\nbigi)$,
let us construct
a map
\index{map $\BB^{\rd}_{J,\theta^u}$}
\[
 \BB^{\rd}_{J,\theta^u}:
 H^0(J,L_{J,>0})
 \lrarr
 H_1^{\rd}\bigl(
 \cnum^{\ast},
 V\otimes\nbige(x^{-1}u^{-1})
 \bigr).
\]
Let us consider the case
$-\theta^u-\pi/2<\vartheta^J_{\ell}$.
We take any $a_1\in S_0(\nbigi)\cap(\vecI_x(\theta^u)-\pi)$
such that $a_1\leq \vartheta^J_{\ell}$.
Let $a_1=b_N<b_{N-1}<\cdots<b_0=\vartheta^J_{\ell}$
be $S_0(\nbigi)\cap\closedclosed{a_1}{\vartheta^J_{\ell}}$.
We set $J_i=\openopen{b_i}{b_i+\omega^{-1}\pi}$.
We have $J=J_0$.
Let $I_i$ $(i=0,\ldots,N-1)$ be paths connecting
$(1,b_{i})$ to $(1,b_{i+1})$.
Let $\gamma_i$ be a path connecting
$(1,b_i)$ and a point in $\{0\}\times J_i$.
For $J'$ such that $J-\omega^{-1}\pi\leq J'\leq J$,
let $\delta_{J'}$ be a path connecting
$(1,\vartheta^J_{\ell})$ and a point in $\{0\}\times J'$.
Let $\Gamma_J$ be a path connecting
$(0,\vartheta^J_{r}+\delta)$ and $(1,\vartheta^J_{\ell})$,
where $\delta$ denotes any sufficiently small positive number.
For $v\in H^0(J,L_{J,>0})$,
we obtain $v_{J_+}\in H^0(J_+,L_{J_+,>0})\subset H^0(\real,L)$.
There exists the decomposition
\[
 v_{J_+}
 =u_{J,0}+\sum_{J-\omega^{-1}\pi\leq J'\leq J}
 u_{J'},
\]
where $u_{J,0}$ is a section of $L'_{J_-,0}$,
and $u_{J'}$ are sections of $L'_{J',<0}$.
We obtain
$(u_{J,0})_i\in
H^0(J_{i-},L_{J_{i-},0})\simeq
H^0(\real,\nbigt_{\omega}(L))$
induced by $u_{J,0}$.
Note that $(u_{J,0})_{i-1}-(u_{J,0})_{i}\in H^0(J_{i},L_{J_i,<0})$.
We obtain the following rapid decay $1$-cycle
of $V\otimes\nbige(x^{-1}u^{-1})$:
\begin{multline}
\label{eq;24.2.23.1}
 \varphi_{\ast}\Bigl(
 v_{J_+}\otimes\Gamma_J
 +\!\!\!\!\sum_{J-\omega^{-1}\pi\leq J'\leq J}\!\!\!\!
 u_{J'}\otimes\delta_{J'}
 +\sum_{i=0}^{N-1}
 (u_{J,0})_{i}
 \otimes I_i
\\
 +\sum_{i=1}^{N-1}
 \bigl(
(u_{J,0})_{i-1}
 -(u_{J,0})_{i}
 \bigr)\otimes\gamma_i
 +(u_{J,0})_{N-1}\otimes\gamma_N
 \Bigr).
\end{multline}
Let $\BB^{\rd}_{J,\theta^u}(v)\in
H^{\rd}_1(\cnum^{\ast},V\otimes\nbige(x^{-1}u^{-1}))$
denote the homology class.

Let us consider the case
$\vartheta^J_{r}<-\theta^u+\pi/2$.
We take any $a_1\in S_0(\nbigi)\cap(\vecI_x(\theta^u)-\pi)$
such that $a_1\geq \vartheta^J_{r}$.
Let $b_0=\vartheta^J_{\ell}<b_1<\cdots<b_{N-1}<b_N=a_1$
be the set $S_0(\nbigi)\cap\closedclosed{\vartheta^J_{r}}{a_1}$.
We set $J_i=\openopen{b_i-\omega^{-1}\pi}{b_i}$.
We have $J=J_0$.
Let $I_i$ $(i=0,\ldots,N-1)$ be paths connecting
$(1,b_{i})$ to $(1,b_{i+1})$.
Let $\gamma_i$ be a path connecting
$(1,b_i)$ and a point in $\{0\}\times J_i$.
For $J'$ such that $\vartheta^J_r\in J'$,
let $\delta_{J'}$ be a path connecting
$(1,\vartheta^J_{r})$ and a point in $\{0\}\times J'$.
Let $\Gamma_J$ be a path connecting
$(0,\vartheta^J_{\ell}-\delta)$ and $(1,\vartheta^J_r)$,
where $\delta$ denotes any sufficiently small positive number.
For $v\in H^0(J,L_{J,>0})$,
we obtain $v_{J_-}\in H^0(J_-,L_{J_-,>0})\subset H^0(\real,L)$.
There exists the decomposition
\[
 v_{J_-}
 =u_{J,0}+\sum_{J\leq J'\leq J+\omega^{-1}\pi}
 u_{J'},
\]
where $u_{J,0}$ is a section of $L'_{J_+,0}$,
and $u_{J'}$ are sections of $L'_{J',<0}$.
We obtain
$(u_{J,0})_i\in
H^0(J_{i+},L_{J_{i+},0})\simeq
H^0(\real,\nbigt_{\omega}(L))$
induced by $u_{J,0}$.
We obtain the following rapid decay $1$-cycle of
$V\otimes\nbige(x^{-1}u^{-1})$:
\begin{multline}
\label{eq;24.2.23.2}
 \varphi_{\ast}\Bigl(
 -v_{J_-}\otimes\Gamma_J
 -\!\!\!\!\sum_{J\leq J'\leq J+\omega^{-1}\pi}\!\!\!\!
 u_{J'}\otimes\delta_{J'}
 -\sum_{i=0}^{N-1}
 (u_{J,0})_{i}
 \otimes I_i
 \\
 -\sum_{i=1}^{N-1}
 \bigl(
(u_{J,0})_{i-1}
 -(u_{J,0})_{i}
 \bigr)\otimes\gamma_i
 -(u_{J,0})_{N-1}\otimes\gamma_N
 \Bigr).
\end{multline}
Let $\BB^{\rd}_{J,\theta^u}(v)\in
H^{\rd}_1(\cnum^{\ast},V\otimes\nbige(x^{-1}u^{-1}))$
denote the homology class.
The following lemma is easy to see.
\begin{lem}
In the constructions, the homology classes are independent of
the choice of $a_1$.
If $\Jbar\subset \vecI_x(\theta^u)-\pi$,
the two constructions give the same homology class. 
\end{lem}
\pf
Let us explain another construction of
$\BB^{\rd}_{J,\theta^u}(v)$
in the case $-\theta^u-\pi/2<\vartheta^J_{\ell}$.
Let $\delta>0$ be sufficiently small.
We set $W=\openopen{-\theta^u_0-\pi/2}{\vartheta^J_{r}+\delta}$,
and
\[
Z_0=\closedopen{0}{\epsilon}\times W,
\quad
Z_1=\openopen{0}{\epsilon}\times W,
\quad
Z_2=\{0\}\times W.
\]
Let $M_{J,1}$ be the constructible subsheaf of $L$
determined by
$M_{J,1}=L^{\leq 0}$ on $\real\setminus \Jbar$
and
$M_{J,1}=L^{\leq 0}+\gbiga_J(L)$ on $\Jbar$.
Let $M_{J,2}$ be the constructible subsheaf of $M_{J,1}$
determined by
$M_{J,2}=L^{<0}$ on $\real\setminus(\vecI_x(\theta^u)-\pi)$
and
$M_{J,2}=L^{\leq 0}$ on $\vecI_x(\theta^u)-\pi$.
Let $K$ be the constructible subsheaf
of $q_{Z_0}^{-1}(L)$ on $Z_0$ determined by the following conditions:
\[
 K_{|Z_2}=M_{J,2|W},
 \quad
 K_{|Z_1}=q_{Z_1}^{-1}(M_{J,1}).
\]
We have the constructible subsheaves
$K_0$ and $K_1$ of $K$ on $Z_0$
determined as follows:
\[
 K_0=q_{Z_0}^{-1}(L^{<0}),
 \quad
 K_{1|Z_2}=M_{J,2|Z_2},
 \quad
 K_{1|Z_1}=q_{Z_1}^{-1}(L^{\leq 0}). 
\]
We obtain the following constructible subsheaves of
$\varphi^{-1}(\nbigl^{<0}((V,\nabla)\otimes\nbige(x^{-1}u^{-1})))$:
\[
 \iota_{Z_0!}(K_0)
 \subset
 \iota_{Z_0!}(K_1)
 \subset
 \iota_{Z_0!}(K).
\]
The constructible sheaves
$\iota_{Z_0!}(K_0)$
and $\iota_{Z_0!}(K_1/K_0)$
are acyclic with respect to the global cohomology.
We have
\[
 \iota_{Z_0!}(K/K_1)
 =\iota_{Z_0!}\bigl(
 q_{Z_1}^{-1}(a_{J\ast}(L_{J,>0}))
 \bigr).
\]
Hence, we obtain
\[
 H^0(J,L_{J,>0})
 \simeq
 H^0(X,\iota_{Z_0!}(K))
 \lrarr
 H_1^{\rd}\bigl(
 \cnum^{\ast},
 (V,\nabla)\otimes\nbige(x^{-1}u^{-1})
 \bigr).
\]
It equals $\BB^{\rd}_{J,\theta^u}$.
In particular, we obtain that
$\BB^{\rd}_{J,\theta^u}$
is independent of the choice of $a_1$
in the case $-\theta^u-\pi/2<\vartheta^J_{\ell}$.
The other case can be shown similarly.

Suppose that $\Jbar\subset \vecI_x(\theta^u)-\pi$.
Let $\BB^{\rd}_{J,\theta^u,1}(v)$
and $\BB^{\rd}_{J,\theta^u,2}(v)$
denote the homology classes obtained from
(\ref{eq;24.2.23.1}) and (\ref{eq;24.2.23.2}),
respectively.
We set
$Z'=\closedopen{0}{\epsilon}\times(\vecI_x(\theta^u)-\pi)$.
The difference
$\BB^{\rd}_{J,\theta^u,1}(v)-
\BB^{\rd}_{J,\theta^u,2}(v)$
can be represented by a rapid decay $1$-cycle
obtained from a $1$-cocycle of
$\nbigc^{\bullet}_{X,\del X}\otimes
 \iota_{Z'!}q_{Z'}^{-1}\bigl(L^{\leq 0}\bigr)[-2]$.
Because
$\iota_{Z'!}q_{Z'}^{-1}\bigl(L^{\leq 0}\bigr)$ is acyclic
with respect to the global cohomology,
we obtain 
$\BB^{\rd}_{J,\theta^u,1}(v)-
\BB^{\rd}_{J,\theta^u,2}(v)=0$.
\hfill\qed

\vspace{.1in}

Let $\Tbb:\real\to\real$ be defined by $\Tbb(\theta)=\theta+2\pi$.
\index{map $\Tbb$}
We have the isomorphism
$\Tbb^{\ast}:H^0(J+2\pi,L_{J+2\pi,>0})
\simeq H^0(J,L_{J,>0})$.
The following lemma is clear by the construction.
\begin{lem}
\label{lem;24.4.5.30}
$\BB^{\rd}_{J+2\pi,\theta^u-2\pi}
=\BB^{\rd}_{J,\theta^u}\circ\Tbb^{\ast}$
on $H^0(J+2\pi,L_{J+2\pi,>0})$. 
\hfill\qed
\end{lem}

\begin{lem}
\label{lem;24.4.5.31}
For any $v\in H^0(J,L_{J,>0})$, we obtain
 \begin{equation}
\label{eq;24.2.23.3}
 \BB^{\rd}_{J,\theta^u+2\pi}(v)
=\BB^{\rd}_{J,\theta^u}(v)
 +\Abb^{\rd}_{\infty,\theta^u}
 (M_0^{-1}\circ\nbigq_{J_+}(v)).
 \end{equation}
As a result, we obtain
$\BB^{\rd}_{J+2\pi,\theta^u}
=\BB^{\rd}_{J,\theta^u}\circ (\Tbb)^{\ast}
+\Abb^{\rd}_{\infty,\theta^u}\circ
\nbigq_{(J+2\pi)_+}$.
\end{lem}
\pf
If $-\theta^u-\pi/2<\vartheta^J_{\ell}$,
we obtain
\[
 \BB^{\rd}_{J,\theta^u+2\pi}(v)
=\BB^{\rd}_{J,\theta^u}(v)
+\Abb^{\rd}_{\infty,\theta^u+2\pi}(\nbigq_{J_+}(v))
\]
by the construction.
If $\vartheta^J_{\ell}\leq -\theta^u-\pi/2$,
we obtain
\[
\BB^{\rd}_{J,\theta^u+2\pi}(v)=
\BB^{\rd}_{J,\theta^u}(v)
-\Abb^{\rd}_{\infty,\theta^u+2\pi}(\nbigq_{J_-}(v))
=\BB^{\rd}_{J,\theta^u}(v)
+\Abb^{\rd}_{\infty,\theta^u+2\pi}(\nbigq_{J_+}(v))
\]
by the construction.
Then, we obtain (\ref{eq;24.2.23.3})
from Lemma \ref{lem;24.2.23.4}.
\hfill\qed

\subsection{Splitting}

\begin{prop}
The maps
$\Abb^{\rd}_{\infty,\theta^u}$ and
$\BB^{\rd}_{J,\theta^u}$ $(J\in \gbigt(\nbigi,\theta^u))$
induce an isomorphism
\[
H^0(\real,\nbigt_{\omega}(L))
 \oplus
 \bigoplus_{J\in\gbigt(\nbigi,\theta^u)}
 H^0(J,L_{J,>0})
 \lrarr
 H_1^{\rd}\bigl(
 \cnum^{\ast},
 V\otimes\nbige(x^{-1}u^{-1})
 \bigr).
\]
\end{prop}
\pf
It is easy to check that
the tuple $\BB^{\rd}_{J,u}$ $(J\in \gbigt(\nbigi,\theta^u))$
induces a splitting of the exact sequence (\ref{eq;24.2.21.4}).
\hfill\qed

\subsection{Rapid decay homology classes $A_{J_{\pm},\theta^u}(v)$}
\label{subsection;24.2.22.20}

For $J\in T(\nbigi)$ such that $J_+\subset \vecI_x(\theta^u)$,
we shall construct a map
\index{maps $A_{J_{\pm},\theta^u}$}
\[
 A_{J_+,\theta^u}:
 H^0(J,L_{J,<0})
 \lrarr
 H_1^{\rd}\bigl(
 \cnum^{\ast},V\otimes\nbige(x^{-1}u^{-1})
 \bigr).
\]
For $J\in T(\nbigi)$ such that $J_-\subset \vecI_x(\theta^u)$,
we shall construct a map
\[
 A_{J_-,\theta^u}:
 H^0(J,L_{J,<0})
 \lrarr
 H_1^{\rd}\bigl(
 \cnum^{\ast},V\otimes\nbige(x^{-1}u^{-1})
 \bigr).
\]
They will be useful in our study of the Stokes structure of
$(\gbigl^{\gbigf}_{\star}(V_{\infty}),\vecnbigf)$ $(\star=!,\ast)$.

\subsubsection{}
\label{subsection;25.2.9.10}

For any $J_1\in T(\nbigi)$,
let $\gbigk((J_1)_+)$
denote the set of $J\in T(\nbigi)$
such that $J_-\cap (J_1)_+\neq\emptyset$,
i.e.,
$J_1-\omega^{-1}\pi<J\leq J_1+\omega^{-1}\pi$.
Similarly,
let $\gbigk((J_1)_-)$
denote the set of $J\in T(\nbigi)$
such that $J_+\cap (J_1)_-\neq\emptyset$,
i.e.,
$J_1-\omega^{-1}\pi\leq J<J_1+\omega^{-1}\pi$.
\index{sets $\gbigk(J_{\pm})$}
There exist the decompositions of the local system
\begin{equation}
\label{eq;24.2.22.2}
 L
=
 L'_{(J_1)_{+},0}
\oplus
\bigoplus_{J\in\gbigk((J_1)_{+})}
 L'_{J,<0}
 =L'_{(J_1)_{-},0}
 \oplus
\bigoplus_{J\in\gbigk((J_1)_{-})}
 L'_{J,<0}
\end{equation}
(See Remark \ref{rem;20.9.7.1}
for the local systems
$L'_{(J_1)_{\pm},0}$ and $L'_{J,<0}$.)

\subsubsection{Construction of $A_{J_+,\theta^u}$}

Let $J\in T(\nbigi)$ such that $J_+\subset \vecI_x(\theta^u)$.
Let $\gamma_1$ be a path
connecting a point in $\{0\}\times J$
and $(1,\vartheta^{J}_{r}-\pi)$.
For each $J'\in \gbigk((J-\pi)_+)$,
the intersection
$J'\cap(\vecI_x(\theta^u)-\pi)$ is not empty.
Let $\gamma_{J'}$ be a path 
connecting $(1,\vartheta^J_r-\pi)$
and a point in $\{0\}\times (J'\cap (\vecI_x(\theta^u)-\pi))$.

Recall the decomposition (\ref{eq;24.2.22.2})
with $J_1=J-\pi$.
For $v\in H^0(J,L_{J,<0})$,
there exists the decomposition 
\[
 v=
 u_{J-\pi,0}
+\sum_{J'\in\gbigk((J-\pi)_+)}
 u_{J'},
\]
where $u_{J-\pi,0}$
is a section of $L'_{(J-\pi)_{+},0}$,
and $u_{J'}$ are sections of $L'_{J',<0}$.
We obtain the following rapid decay $1$-cycle
of $(V,\nabla)\otimes\nbige(x^{-1}u^{-1})$:
\[
\varphi_{\ast}\Bigl(
 v\otimes \gamma_1
+u_{J-\pi,,0}\otimes\gamma_{J-\pi}
+\sum_{J'}
 u_{J'}\otimes\gamma_{J'}
 \Bigr).
\]
Let $A_{J_+,\theta^u}(v)$ denote the homology class.

\subsubsection{Construction of $A_{J_-,\theta^u}$}

Let $J\in T(\nbigi)$ such that $J_-\subset \vecI_x(\theta^u)$.
Let $\gamma_1$ be a path
connecting a point in $\{0\}\times J$
and $(1,\vartheta^{J}_{\ell}+\pi)$.
For each $J'\in \gbigk((J+\pi)_-)$,
the intersection
$J'\cap(\vecI_x(\theta^u)+\pi)$ is not empty.
Let $\gamma_{J'}$ be a path
connecting a point in $\{0\}\times (J'\cap (\vecI_x(\theta^u)+\pi))$
and $(1,\vartheta^J_{\ell}+\pi)$.

For $v\in H^0(J,L_{J,<0})$,
there exists the decomposition
\[
 v=
 u_{J+\pi,0}
 +\sum_{J'\in\gbigk((J+\pi)_-)}
  u_{J'},
\]
where $u_{J+\pi,0}$
is a section of $L'_{(J+\pi)_{-},0}$,
and $u_{J'}$ are sections of $L'_{J',<0}$.
We obtain the following rapid decay $1$-cycle
of $(V,\nabla)\otimes\nbige(x^{-1}u^{-1})$:
\[
\varphi_{\ast}\Bigl(
 v\otimes \gamma_1
+u_{J+\pi,0}\otimes\gamma_{J+\pi}
+\sum_{J'}
 u_{J'}\otimes\gamma_{J'}
 \Bigr).
\]
Let $A_{J_-,\theta^u}(v)$ denote the homology class.

\subsubsection{}

The following lemma is clear by the construction.
\begin{lem}
$A_{(J+2\pi)_{\pm},\theta^u-2\pi}(v)
=A_{J_{\pm},\theta^u}(\Tbb^{\ast}(v))$. 
\hfill\qed
\end{lem}

We express
$A_{J_-,\theta^u}(v)$
in terms of standard homology classes
in \S\ref{subsection;24.2.23.10}.

\begin{prop}
\label{prop;24.3.17.1}
For $J\in T(\nbigi)$ such that $J_+\subset\vecI_x(\theta^u)$
and for $v\in H^0(J,L_{J,<0})$,
we obtain
\begin{equation}
\label{eq;24.2.21.30}
 A_{J_+,\theta^u}(v)
 =\sum_{J-\pi<J'\leq J-\omega^{-1}\pi}
 \BB^{\rd}_{J',\theta^u}\bigl(
 R_{J'}(v)
 \bigr).
\end{equation}
(See {\rm\S\ref{subsection;24.2.18.1}} for the maps $R_{J'}$.)
 For $J\in T(\nbigi)$ such that $J_-\subset\vecI_x(\theta^u)$
and for $v\in H^0(J,L_{J,<0})$,
we obtain
\begin{equation}
\label{eq;24.2.21.31}
 A_{J_-,\theta^u}(v)
=\sum_{J+\omega^{-1}\pi\leq J'<J+\pi}
-\BB^{\rd}_{J',\theta^u-2\pi}\bigl(
 R_{J'}(v)
 \bigr).
\end{equation}
\end{prop}
\pf
We explain a proof for (\ref{eq;24.2.21.30}).
The other case can be argued similarly.
Let $\delta>0$ be sufficiently small.
We set
$W=\openopen{\vartheta^{\vecI_x(\theta^u)}_{\ell}-\pi}{\vartheta^J_r+\delta}$,
and
\[
 Z_0=\closedopen{0}{\epsilon}\times W,
 \quad
 Z_1=\openopen{0}{\epsilon}\times W,
 \quad
 Z_2=\{0\}\times W.
\]
Let $M$ be the constructible subsheaf of $L$
determined by
$M=L^{<0}$ on $\real\setminus (J-\pi)$,
and $M=L^{\leq 0}$ on $J-\pi$.
We consider the constructible subsheaf
$K$ of $q_{Z_0}^{-1}(L)$
determined by
\[
 K_{|Z_2}=M,
 \quad
 K_{|Z_1}=q_{Z_1}^{-1}(L^{\leq 0}).
\]
It is easy to see that
$\iota_{Z_0!}K$ is acyclic with respect
to the global cohomology.
The homology class
\[
\alpha=A_{J_+,\theta^u}(v)
-\sum_{J-\pi<J'\leq J-\omega^{-1}\pi}
\BB^{\rd}_{J',\theta^u}\bigl(
R_{J'}(v)
\bigr)
\]
is induced by a $1$-cocycle of
$\nbigc^{\bullet}_{X,\del X}\otimes
 \iota_{Z_0!}K[-2]$.
Hence, we obtain $\alpha=0$.
\hfill\qed

\begin{rem}
See Proposition {\rm\ref{prop;24.3.22.31}}
for the difference
$A_{J_+,\theta^u}-A_{J_-,\theta^u}$.
\hfill\qed
\end{rem}

\section{Moderate growth homology of $(V,\nabla)\otimes\nbige(x^{-1}u^{-1})$}

\subsection{Exact sequence}

\subsubsection{Moderate growth homology classes
$\Abb^{\mg}_{J,\theta^u}(v)$}
We use the notation in \S\ref{subsection;24.2.21.6}.
Let $J\in T(\nbigi)$.
Let $\Gamma_J$ be a path on $(X,X^{\ast})$
connecting a point in $\{0\}\times J$
and a point in $\{\infty\}\times J$.
The image of $\Gamma_J$ is assumed to be contained in
$\real_{\geq 0}\times J$.
For $v\in H^0(J,L_{J,<0})$,
we obtain a $1$-cocycle
$\varphi_{\ast}(v\otimes\Gamma_J)$
of 
$\nbigc^{\bullet}_{\projtilde^1,\del\projtilde^1}
\otimes
j_{0\ast}q_0^{-1}a_{J!}L_{J,<0}[-2]$.
This procedure induces an isomorphism
\begin{equation}
\label{eq;24.2.23.11}
H^0(J,L_{J,<0})
 \simeq
 \hyperh^{-1}\bigl(
 \projtilde^1,
 \nbigc^{\bullet}_{\projtilde^1,\del\projtilde^1}
 \otimes
  j_{0\ast}q_0^{-1}a_{J!}L_{J,<0}
 \bigr).
\end{equation}
We shall identify them by this isomorphism.
There exists the natural morphism
\begin{equation}
\label{eq;24.2.23.12}
  \hyperh^{-1}\bigl(
 \projtilde^1,
 \nbigc^{\bullet}_{\projtilde^1,\del\projtilde^1}
 \otimes
  j_{0\ast}q_0^{-1}a_{J!}L_{J,<0}
  \bigr)
  \lrarr
  H_1^{\mg}\bigl(
  \cnum^{\ast},
  V\otimes\nbige(x^{-1}u^{-1})
  \bigr).
\end{equation}
The image of $v$ via (\ref{eq;24.2.23.11})
and (\ref{eq;24.2.23.12})
is denoted by $\Abb^{\mg}_{J,\theta^u}(v)$.
\index{maps $\Abb^{\mg}_{J,\theta^u}$}
Thus, we obtain
\[
 \Abb^{\mg}_{J,\theta^u}:
 H^0(J,L_{J,<0})
 \lrarr
   H_1^{\mg}\bigl(
  \cnum^{\ast},
  V\otimes\nbige(x^{-1}u^{-1})
  \bigr).
\]
\begin{lem}
\label{lem;24.3.22.1}
$\Abb^{\mg}_{J,\theta^u+2\pi}
=\Abb^{\mg}_{J,\theta^u}$.
\hfill\qed
\end{lem}

\begin{lem}
\label{lem;24.4.5.210}
We have
$\Abb^{\mg}_{J+2\pi,\theta^u-2\pi}
=\Abb^{\mg}_{J,\theta^u}\circ \Tbb^{\ast}$
on $H^0(J+2\pi,L_{J+2\pi,<0})$,
and hence
$\Abb^{\mg}_{J+2\pi,\theta^u}
=\Abb^{\mg}_{J,\theta^u}\circ\Tbb^{\ast}$.
\hfill\qed
\end{lem}

\subsubsection{Expression of
$H_1^{\mg}(\cnum^{\ast},\nbigt_{\omega}(V)\otimes
\nbige(x^{-1}u^{-1}))$}

Let $\Gamma_{\theta^u}$ be a path
connecting a point in $\{\infty\}\times\real$
and $(0,-\theta^u)$.
For any $v\in H^0(\real,\nbigt_{\omega}(L))$,
we obtain
the moderate growth cycle
$\varphi_{\ast}(v\otimes\Gamma_{\theta^u})$
of $\nbigt_{\omega}(V)\otimes\nbige(x^{-1}u^{-1})$.
This procedure induces an isomorphism
depending on $\theta^u$.
\begin{equation}
\label{eq;24.2.22.10}
 H^0(\real,\nbigt_{\omega}(L))
 \simeq
 H_1^{\mg}\bigl(
 \cnum^{\ast},
 \nbigt_{\omega}(V)
 \otimes\nbige(x^{-1}u^{-1})
 \bigr).
\end{equation}
We shall identify them by this isomorphism.
\begin{lem}
\label{lem;25.2.9.11}
The homology classes of
$\varphi_{\ast}(v\otimes\Gamma_{\theta^u-2\pi})$
and
$\varphi_{\ast}(M_0(v)\otimes\Gamma_{\theta^u})$
are the same.
\hfill\qed
\end{lem}

\subsubsection{Moderate growth homology classes
$\Abb^{\mg,J_{\pm}}_{\infty,\theta^u}(v)$}

We shall construct the following maps
for any $J\in T(\nbigi)$:
\index{maps $\Abb^{\mg,J_{\pm}}_{\infty,\theta^u}$}
\begin{equation}
\label{eq;24.3.26.10}
 \Abb^{\mg,J_{\pm}}_{\infty,\theta^u}:
 H^0(\real,\nbigt_{\omega}(L))
 \lrarr
    H_1^{\mg}\bigl(
  \cnum^{\ast},
  V\otimes\nbige(x^{-1}u^{-1})
  \bigr).
\end{equation}

Let $J\in T(\nbigi)$
such that
$\vartheta^J_r\in (\vecI_x(\theta^u)-\pi)_-$.
For $v\in H^0(\real,\nbigt_{\omega}(L))$,
we have
$v_{J_{+}}\in H^0(J_{+},L_{J_{+},0})
\subset H^0(\real,L)$.
Let $\Gamma_{J}$ be
a path connecting a point in $\{\infty\}\times\real$
and $(0,\vartheta^J_r+\delta)\in \{0\}\times(\vecI_x(\theta^u)-\pi)$,
where $\delta>0$ denotes a sufficiently small number.
We obtain the moderate growth $1$-cycle
$\varphi_{\ast}\bigl(
v_{J_+}\otimes\Gamma_J
\bigr)$
of $(V,\nabla)\otimes\nbige(x^{-1}u^{-1})$.
The homology class is denoted by
$\Abb^{\mg,J_{+}}_{\infty,\theta^u}(v)$.

Similarly, 
let $J\in T(\nbigi)$
such that
$\vartheta^J_{\ell}\in (\vecI_x(\theta^u)-\pi)_+$.
For $v\in H^0(\real,\nbigt_{\omega}(L))$,
we have
$v_{J_-}\in H^0(J_-,L_{J_-,0})\subset H^0(\real,L)$.
Let $\Gamma_J$ be
a path connecting a point in $\{\infty\}\times\real$
and $(0,\vartheta^J_{\ell}-\delta)\in \{0\}\times(\vecI_x(\theta^u)-\pi)$
for any sufficiently small $\delta>0$.
We obtain the moderate growth $1$-cycle
$\varphi_{\ast}\bigl(
 v_{J_-}\otimes\Gamma_J
 \bigr)$
of $(V,\nabla)\otimes\nbige(x^{-1}u^{-1})$.
The homology class is denoted by
$\Abb^{\mg,J_-}_{\infty,\theta^u}(v)$.

Let $J\in T(\nbigi)$.
Choosing $J_1\in T(\nbigi)$ such that
$\vartheta^{J_1}_{r}\in (\vecI_x(\theta^u)-\pi)_-$,
we set
\[
 \Abb^{\mg,J_{+}}_{\infty,\theta^u}(v)
 =
\left\{
\begin{array}{ll}
 \Abb^{\mg,J_{1+}}_{\infty,\theta^u}(v)
+\sum_{J_1<J'\leq J}\Abb^{\mg}_{J',\theta^u}\circ\nbigp_{J'}(v)
& (J_1\leq J)
  \\
\Abb^{\mg,J_{1+}}_{\infty,\theta^u}(v)
-\sum_{J< J'\leq J_1}\Abb^{\mg}_{J',\theta^u}\circ\nbigp_{J'}(v)
& (J\leq J_1).
\end{array}
 \right.
\]
Choosing $J_2\in T(\nbigi)$ such that
$\vartheta^{J_2}_{\ell}\in (\vecI_x(\theta^u)-\pi)_+$,
we set
\[
 \Abb^{\mg,J_-}_{\infty,\theta^u}(v)
 =
\left\{
\begin{array}{ll}
 \Abb^{\mg,J_{2-}}_{\infty,\theta^u}(v)
+\sum_{J_2\leq J'<J}\Abb^{\mg}_{J',\theta^u}\circ\nbigp_{J'}(v)
& (J_2\leq J)
  \\
\Abb^{\mg,J_{2-}}_{\infty,\theta^u}(v)
-\sum_{J\leq J'<J_2}\Abb^{\mg}_{J',\theta^u}\circ\nbigp_{J'}(v)
& (J\leq J_2).
\end{array}
 \right.
\]
They are independent of the choices of $J_1$ and $J_2$.
Therefore, we obtain (\ref{eq;24.3.26.10}).

The following lemma is clear by the construction.
\begin{lem}
\mbox{{}}\label{lem;24.3.22.2}
\begin{itemize}
 \item
       $\Abb^{\mg,(J+2\pi)_{\pm}}_{\infty,\theta^u-2\pi}
=\Abb^{\mg,J_{\pm}}_{\infty,\theta^u}\circ M_0$.
 \item
 $\Abb^{\mg,J_{-}}_{\infty,\theta^u}
-\Abb^{\mg,J_{+}}_{\infty,\theta^u}
=-\Abb^{\mg}_{J,\theta^u}\circ\nbigp_{J}$.      
 \item
For $J_1<J_2$ in $T(\nbigi)$,
we obtain 
$\Abb^{\mg,J_{1-}}_{\infty,\theta^u}
-\Abb^{\mg,J_{2-}}_{\infty,\theta^u}
=-\sum_{J_1\leq J<J_2} \Abb^{\mg}_{J,\theta^u}\circ\nbigp_{J}$
and
$\Abb^{\mg,J_{1+}}_{\infty,\theta^u}
-\Abb^{\mg,J_{2+}}_{\infty,\theta^u}
=-\sum_{J_1< J\leq J_2} \Abb^{\mg}_{J,\theta^u}\circ\nbigp_{J}$.
 \hfill\qed      
\end{itemize}
\end{lem}

\subsubsection{Exact sequence}

We use the notation in \S\ref{subsection;24.2.21.6}.
Let $k:\cnum^{\ast}\to\cnumtilde$ denote the inclusion.
We obtain the following exact sequence:
\[
 0\lrarr j_{0\ast}\nbign_0
 \lrarr
 \nbigl^{\leq 0}(V\otimes\nbige(x^{-1}u^{-1}))
 \lrarr
 j_{0\ast}
 \bigl(
 k_!(L_{S^1}/L_{S^1}^{\leq 0})
 \bigr)
 \lrarr 0.
\]
Note that
$j_{0\ast}
 \bigl(
 k_!(L_{S^1}/L_{S^1}^{\leq 0})
 \bigr)$
 is acyclic with respect to the global cohomology.
We obtain the natural isomorphism
\[
 H^1\bigl(\projtilde^1,
j_{0\ast}\nbign_0
  \bigr)
  \simeq
  H_1^{\mg}\bigl(
  \cnum^{\ast},
  V\otimes\nbige(x^{-1}u^{-1})
  \bigr).
\]
There exists the following natural exact sequence:
\[
 0\lrarr
 j_{0\ast}
 q_0^{-1}(L_{S_1}^{<0})
 \lrarr
 j_{0\ast}\nbign_0
 \lrarr
 \nbigl^{\leq 0}\bigl(
 \nbigt_{\omega}(V)
 \otimes\nbige(x^{-1}u^{-1})
 \bigr)
 \lrarr 0.
\]
There exists the following isomorphism:
\[
  j_{0\ast}
  q_0^{-1}(L_{S_1}^{<0})
=\bigoplus_{J\in \gbigt(\nbigi,\theta^u)}
 j_{0\ast}q_0^{-1}a_{J!}L_{J,<0}.
\]

We obtain the following exact sequence:
\begin{multline}
\label{eq;24.2.21.10}
 0\lrarr
 \bigoplus_{J\in \gbigt(\nbigi,\theta^u)}
 H^0(J,L_{J,<0})
 \stackrel{c_{1,u}}\lrarr
 H_1^{\mg}\bigl(
 \cnum^{\ast},
 V\otimes\nbige(x^{-1}u^{-1})
 \bigr)
 \\
 \stackrel{c_{2,u}}
 \lrarr
 H_1^{\mg}\bigl(
 \cnum^{\ast},
 \nbigt_{\omega}(V)
 \otimes
 \nbige(x^{-1}u^{-1})
 \bigr)
\lrarr 0.
\end{multline}
For $J_1\in T(\nbigi)$ such that
$\vartheta^{J_1}_{r}\in (\vecI_x(\theta^u)-\pi)_-$,
the map
$\Abb^{\mg,J_{1+}}_{\infty,\theta^u}$
is a splitting of (\ref{eq;24.2.21.10}).
The maps
$\Abb^{\mg,J_{1+}}_{\infty,\theta^u}$
and
$\Abb^{\mg}_{J,\theta^u}$
$(J\in \gbigt(\nbigi,\theta^u))$
induce an isomorphism:
\[
 \bigoplus_{J\in \gbigt(\nbigi,\theta^u)}
 H^0(J,L_{J,<0})
 \oplus
 H^0(\real,\nbigt_{\omega}(L))
 \simeq
 H_1^{\mg}(\cnum^{\ast},
 V\otimes\nbige(x^{-1}u^{-1})).
\]
Similarly,
for $J_2\in T(\nbigi)$ such that
$\vartheta^{J_2}_{\ell}\in (\vecI_x(\theta^u)-\pi)_+$,
the map
$\Abb^{\mg,J_{2-}}_{\infty,\theta^u}$
is a splitting of (\ref{eq;24.2.21.10}).
The maps
$\Abb^{\mg,J_{2-}}_{\infty,\theta^u}$
and
$\Abb^{\mg}_{J,\theta^u}$
$(J\in \gbigt(\nbigi,\theta^u))$
induce an isomorphism:
\[
 \bigoplus_{J\in \gbigt(\nbigi,\theta^u)}
 H^0(J,L_{J,<0})
 \oplus
 H^0(\real,\nbigt_{\omega}(L))
 \simeq
 H_1^{\mg}(\cnum^{\ast},
 V\otimes\nbige(x^{-1}u^{-1})).
\]

\subsection{Relations with rapid decay cycles}
\label{subsection;24.3.26.22}

There exists the natural morphism
$H_1^{\rd}(\cnum^{\ast},(V,\nabla)\otimes\nbige(x^{-1}u^{-1}))
\to
H_1^{\mg}(\cnum^{\ast},(V,\nabla)\otimes\nbige(x^{-1}u^{-1}))$.
The image of an element of
$H_1^{\rd}(\cnum^{\ast},(V,\nabla)\otimes\nbige(x^{-1}u^{-1}))$
is denoted by the same notation.
We obtain the following lemmas by the construction.
\begin{lem}
\label{lem;24.4.5.40}
For any $J\in T(\nbigi)$, and for $v\in H^0(\real,\nbigt_{\omega}(L))$,
\begin{multline}
 \Abb^{\rd}_{\infty,\theta^u}(v)
=\Abb^{\mg,(J+2\pi)_+}_{\infty,\theta^u}(v)
-\Abb^{\mg,(J+2\pi)_+}_{\infty,\theta^u-2\pi}(v)
=\Abb^{\mg,(J+2\pi)_-}_{\infty,\theta^u}(v)
-\Abb^{\mg,(J+2\pi)_-}_{\infty,\theta^u-2\pi}(v)
 \\
=\Abb^{\mg,J_+}_{\infty,\theta^u}(v-M_0(v))
+\sum_{J<J'\leq J+2\pi}
\Abb^{\mg}_{J',\theta^u}(\nbigp_{J'}(v))
\\
=\Abb^{\mg,J_-}_{\infty,\theta^u}(v-M_0(v))
+\sum_{J\leq J'< J+2\pi}
\Abb^{\mg}_{J',\theta^u}(\nbigp_{J'}(v)).
\end{multline}
\end{lem}
\pf
For $J\in T(\nbigi)$ such that 
$\vartheta^J_{r}\in(\vecI_x(\theta^u)-\pi)$,
we obtain
$\Abb^{\rd}_{\infty,\theta^u}(v)
=\Abb^{\mg,(J+2\pi)_+}_{\infty,\theta^u}(v)
-\Abb^{\mg,(J+2\pi)_+}_{\infty,\theta^u-2\pi}(v)$
by the construction.
By using Lemma \ref{lem;24.3.22.1}
and Lemma \ref{lem;24.3.22.2},
we obtain 
$\Abb^{\rd}_{\infty,\theta^u}(v)
=\Abb^{\mg,(J+2\pi)_+}_{\infty,\theta^u}(v)
-\Abb^{\mg,(J+2\pi)_+}_{\infty,\theta^u-2\pi}(v)
=\Abb^{\mg,(J+2\pi)_-}_{\infty,\theta^u}(v)
-\Abb^{\mg,(J+2\pi)_-}_{\infty,\theta^u-2\pi}(v)$
for any $J\in T(\nbigi)$.
We obtain the other equalities from Lemma \ref{lem;24.3.22.2}.
\hfill\qed

\begin{lem}
\label{lem;24.4.5.41}
For $J\in T(\nbigi)$
and $v\in H^0(J,L_{J,>0})$,
we obtain
\begin{multline}
\label{eq;24.3.22.10}
 \BB^{\rd}_{J,\theta^u}(v)
 =\sum_{J<J'\leq J+\omega^{-1}\pi}
 \Abb^{\mg}_{J',\theta^u}(\nbigr^{J}_{J'}(v))
 -\sum_{J-\omega^{-1}\pi\leq J'<J}
 \Abb^{\mg}_{J',\theta^u}(\nbigr^{J}_{J'}(v))
 \\
-\Abb^{\mg}_{J,\theta^u}(\nbigr^{J_+}_{J_-}(v))
+\Abb^{\mg,J_{-}}_{\infty,\theta^u} (\nbigq_{J_+}(v)).
\end{multline}
\end{lem}
\pf
It $-\theta^u-\pi/2<\vartheta^J_{\ell}$,
we obtain (\ref{eq;24.3.22.10})
by the construction.
If $\vartheta^J_r<-\theta^u+\pi/2$,
we obtain the following:
\begin{multline}
 \BB^{\rd}_{J,\theta^u}(v)
 =\sum_{J<J'\leq J+\omega^{-1}\pi}
 \Abb^{\mg}_{J',\theta^u}(\nbigr^{J}_{J'}(v))
 -\sum_{J-\omega^{-1}\pi\leq J'<J}
 \Abb^{\mg}_{J',\theta^u}(\nbigr^{J}_{J'}(v))
 \\
+\Abb^{\mg}_{J,\theta^u}(\nbigr^{J_-}_{J_+}(v))
-\Abb^{\mg,J_+}_{\infty,\theta^u} (\nbigq_{J_-}(v)).
\end{multline}
By using
$\nbigr^{J_-}_{J_+}
=-\nbigr^{J_+}_{J_-}
+\nbigp_{J_+}\circ\nbigq_{J_+}$,
$\nbigq_{J_+}=-\nbigq_{J_-}$
and Lemma \ref{lem;24.3.22.2},
we obtain (\ref{eq;24.3.22.10}).
\hfill\qed

\section{Lifting maps for $(V,\nabla)\otimes\nbige(x^{-1}u^{-1})$}

\subsection{Some constructible sheaves on $S^1$}
\label{subsection;18.4.23.1}

We take 
$J_1\in T(\nbigi)$
such that 
$J_{1+}\subset\vecI_x(\theta^u)-\pi$.
For each $J\in \gbigk(J_{1+})$,
we take an interval
$I_{10}(J)\subset
 J\cap (\vecI_x(\theta^u)-\pi)\neq\emptyset$.
We obtain the following constructible subsheaves of $L$:
\[
K_0^{J_{1+}}:=
 \bigoplus_{J\in \gbigk(J_{1+})}
 a_{I_{10}(J)!}
 L_{J,<0|I_{10}(J)},
\quad\quad
K_1^{J_{1+}}:=
 a_{J_1!}L_{J_{1+},0|J_1}
\oplus
 K_0^{J_{1+}}.
\]
We obtain the constructible subsheaves
$\varphi_{1!}K_0^{J_{1+}}\subset
\varphi_{1!}K_1^{J_{1+}}\subset L_{S^1}$.
There exists the following commutative diagram:
{\small
\[
 \begin{CD}
 0@>>>
 \varphi_{1!}K_0^{J_{1+}} @>>> 
 \varphi_{1!}K_1^{J_{1+}} @>>>
 \varphi_{1!}\bigl(
 a_{J_1!}(L_{J_{1+},0|J_1})
 \bigr)
 @>>> 0
 \\
 @.
 @V{c_1}VV @V{c_2}VV @V{c_3}VV 
 @. \\
 0@>>>
 L^{<0}_{S^1}
@>>>
 L^{F,<0}_{S^1}
@>>>
 \varphi_{1!}
 a_{(\vecI_x(\theta^u)-\pi)!}\bigl(
 \nbigt_{\omega}(L)_{|\vecI_x(\theta^u)-\pi}
 \bigr)
 @>>> 0.
 \end{CD}
\]}
The rows are exact.
The morphisms $c_i$ are monomorphisms.
Note that $\Cok c_3$ is acyclic with respect to
the global cohomology.
Let $\gbign(J_{1+})$ denote the set of
$J\in T(\nbigi)$
satisfying
$-\theta^u-\pi/2
\leq
\vartheta^J_r
\leq \vartheta^{J_1}_{\ell}$
or 
$\vartheta^{J_1}_{r}+\omega^{-1}\pi<\vartheta^J_r
< -\theta^u+3\pi/2$.
We have 
\[
 \Cok(c_1)=
 \bigoplus_{J\in \gbign(J_{1+})}
 \varphi_{1!}(L_{J,<0|J})
\oplus
 \bigoplus_{J\in\gbigk(J_{1+})}
 \varphi_{1!}
 \Bigl( 
 \bigl(
 a_{J!}L_{J,<0|J}
 \bigr)\big/
 a_{I_{10}(J)!}L_{J,<0|I_{10}(J)}
\Bigr).
\]
The second term in the right hand side
is acyclic with respect to the global cohomology.

Similarly, let $J_1\in T(\nbigi)$.
For each $J\in \gbigk(J_{1-})$,
we take an interval
$I_{11}(J)\subset J\cap (\vecI_x(\theta^u)-\pi)$.
We obtain the following constructible sheaves
\[
K_0^{J_{1-}}:=
 \bigoplus_{J\in \gbigk(J_{1-})}
 a_{I_{11}(J)!}
 L_{J,<0|I_{10}(J)},
\quad\quad
K_1^{J_{1-}}:=
 a_{J_1!}L_{J_{1-},0|J_1}
\oplus
 K_0^{J_{1-}}.
\]

\subsection{Rapid decay case}

Let $(V^{\reg},\nabla)=\nbigstilde_{\omega}(V,\nabla)$
be the regular singular meromorphic flat bundle
on $(\proj^1,\{0,\infty\})$
corresponding to $\nbigl$.
For an interval 
$J_1\in T(\nbigi)$
such that $J_{1+}\subset\vecI_x(\theta^u)-\pi$,
we shall construct the following map:
\index{maps $A^{J_{1+}}_{\infty,\theta^u}$}
\begin{equation}
 A^{J_{1+}}_{\infty,\theta^u}:
 H_1^{\rd}\bigl(\cnum^{\ast},V^{\reg}\otimes
 \nbige(u^{-1}x^{-1})\bigr)
\lrarr
 H_1^{\rd}\bigl(\cnum^{\ast},V\otimes
 \nbige(u^{-1}x^{-1})\bigr).
\end{equation}
Similarly,
for an interval $J_2\in T(\nbigi)$
such that $J_{2-}\subset\vecI_x(\theta^u)-\pi$,
we shall construct the following map:
\index{map $ A^{J_{2-}}_{\infty,\theta^u}$}
\begin{equation}
\label{eq;24.2.22.3}
 A^{J_{2-}}_{\infty,\theta^u}:
 H_1^{\rd}\bigl(\cnum^{\ast},V^{\reg}\otimes
 \nbige(u^{-1}x^{-1})\bigr)
\lrarr
 H_1^{\rd}\bigl(\cnum^{\ast},
 V\otimes\nbige(u^{-1}x^{-1})\bigr).
\end{equation}
Because $(V^{\reg},\nabla)$ is regular singular at $\{0,\infty\}$,
there exists the isomorphism
as in (\ref{eq;24.2.21.40}):
\begin{equation}
 H^0(\real,L)\simeq
 H_1^{\rd}(\cnum^{\ast},V^{\reg}\otimes\nbige(x^{-1}u^{-1}))
\end{equation}
We shall also obtain the maps
\[
 A^{J_{1+}}_{\infty,\theta^u}:
 H^0(\real,L)
 \lrarr
 H_1^{\rd}\bigl(\cnum^{\ast},
 V\otimes\nbige(x^{-1}u^{-1})
 \bigr),
\]
\[
 A^{J_{2-}}_{\infty,\theta^u}:
 H^0(\real,L)
 \lrarr
 H_1^{\rd}\bigl(\cnum^{\ast},
 V\otimes\nbige(x^{-1}u^{-1})
 \bigr).
\]

\subsubsection{Construction in the case of ``$+$''}
\label{subsection;25.2.5.3}

There exists the constructible subsheaf
\begin{equation}
\label{eq;24.2.12.1}
 \nbigm^{J_{1+}}\bigl(
 V\otimes\nbige(u^{-1}x^{-1})
 \bigr)
 \subset 
 \nbigl^{<0}\bigl(
 V\otimes\nbige(u^{-1}x^{-1})
 \bigr)
\end{equation}
determined by the following conditions.
\begin{itemize}
\item
$\nbigm^{J_{1+}}\bigl(
 V\otimes\nbige(u^{-1}x^{-1})
 \bigr)_{|\projtilde^1\setminus\varpi^{-1}(0)}
=
 \nbigl^{<0}\bigl(
 V\otimes\nbige(u^{-1}x^{-1})
 \bigr)_{|\projtilde^1\setminus\varpi^{-1}(0)}$.
\item
$\nbigm^{J_{1+}}\bigl(
 V\otimes\nbige(u^{-1}x^{-1})
 \bigr)_{|\varpi^{-1}(0)}
=\varphi_{1!}K_1^{J_{1+}}$.
\end{itemize}
The inclusion (\ref{eq;24.2.12.1}) induces a morphism
\begin{equation}
 \label{eq;24.2.12.2}
   H^1\Bigl(\projtilde^1,
  \nbigm^{J_{1+}}\bigl(
 V\otimes\nbige(u^{-1}x^{-1})
 \bigr)\Bigr)
\lrarr
 H_1^{\rd}\bigl(\cnum^{\ast},V\otimes
 \nbige(u^{-1}x^{-1})\bigr).
\end{equation}
By the construction,
there exists a natural monomorphism
\begin{equation}
\label{eq;18.4.23.11}
 g:
 \nbigm^{J_{1+}}\bigl(
 V\otimes\nbige(u^{-1}x^{-1})
 \bigr)
\lrarr
 \nbigl^{<0}\bigl(
 V^{\reg}\otimes\nbige(u^{-1}x^{-1})
 \bigr).
\end{equation}

\begin{lem}
\label{lem;18.4.23.41}
$\Cok(g)$ is acyclic with respect to
the global cohomology.
Therefore,
the morphism {\rm(\ref{eq;18.4.23.11})}
induces an isomorphism
\begin{equation}
 \label{eq;24.2.12.3}
 H^1\Bigl(\projtilde^1,
  \nbigm^{J_{1+}}\bigl(
 V\otimes\nbige(u^{-1}x^{-1})
 \bigr)\Bigr)
 \simeq
 H_1^{\rd}\bigl(\cnum^{\ast},V^{\reg}\otimes
 \nbige(u^{-1}x^{-1})\bigr).
\end{equation}
\end{lem}
\pf
There exists the decomposition (\ref{eq;24.2.22.2}).
Let $i_0:\varpi^{-1}(0)\lrarr\projtilde^1$
denote the inclusion.
Then,
$\Cok(g)$ is isomorphic to
\begin{multline}
 i_{0!}
 \varphi_{1!}
\left(
 a_{(\vecI_x(\theta^u)-\pi)!}
 \bigl(
  L'_{J_{1+,0}|\vecI_x(\theta^u)-\pi}
 \bigr)
 \Big/
 a_{J_1!}\bigl(
 L'_{J_{1+},0|J_1}
 \bigr)
\right)
\\
 \oplus
 \bigoplus_{J\in\gbigk(J_{1+})}
 i_{0!}
 \varphi_{1!}
\left(
 a_{(\vecI_x(\theta^u)-\pi)!}\bigl(
 L'_{J,<0|\vecI_x(\theta^u)-\pi}
 \bigr)
\Big/
 a_{I_{10}(J)!}\bigl(
 L'_{J,<0|I_{10}(J)}
 \bigr)
 \right).
\end{multline}
Hence, we obtain the claim of the lemma.
\hfill\qed

\vspace{.1in}
We obtain the desired map
$A^{J_{1+}}_{\infty,\theta^u}$
from (\ref{eq;24.2.12.2}) and (\ref{eq;24.2.12.3}).

\subsubsection{Explicit $1$-cycles}

Let us describe $A^{J_{1+}}_{\infty,\theta^u}$
in terms of explicit $1$-cycles.

Let $\gamma_1$ be a path
connecting  $(1,\vartheta^{J_1}_r+2\pi)$ 
and $(1,\vartheta^{J_1}_r)$ on $(X^{\ast},X)$.
We take $\theta_J\in J\cap(\vecI_x(\theta^u)-\pi)$
for each $J\in \gbigk(J_{1+})$,
and path $\gamma_{J,2}$
connecting $(1,\vartheta^{J_1}_r)$
to $(0,\theta_J)$ on $(X^{\ast},X)$.
By shifting $\gamma_{J,2}$ by $2\pi$,
we obtain paths
$\gamma_{J,3}$ connecting
$(1,\vartheta^{J_1}_r+2\pi)$
to $(0,\theta_J+2\pi)$ on $(X^{\ast},X)$.
Let $v\in H^0(\real,L)$.
We have the decompositions
\begin{equation}
\label{eq;18.6.22.1}
 v=u_{J_1+2\pi,0}
+\sum_{J\in \gbigk((J_1+2\pi)_+)} u_{J},
\quad\quad
 v=u_{J_1,0}
+\sum_{J\in\gbigk(J_{1+})}u_J,
\end{equation}
where $u_J$ are sections of $L'_{J,<0}$,
$u_{J_1+2\pi,0}$ is a section of
$L_{(J_1+2\pi)_+,0}$,
and
$u_{J_1,0}$ is a section of
$L_{J_{1+},0}$.
We obtain the $1$-cycle
\begin{multline}
\nbiga^{J_{1+}}_{\infty,\theta^u}(v)
=
 v\otimes\gamma_1
\\
 -\sum_{J\in\gbigk((J_1+2\pi)_+)}
 u_{J}\otimes\gamma_{3,J}
-u_{J_1+2\pi,0}\otimes\gamma_{3,J_1+2\pi}
+\sum_{J\in\gbigk(J_{1+})}
 u_{J}\otimes\gamma_{2,J}
+u_{J_1,0}\otimes\gamma_{2,J_1}.
\end{multline}
Then,
$\varphi_{\ast}(\nbiga^{J_{1+}}_{\infty,\theta^u}(v))$
represents $A^{J_{1+}}_{\infty,\theta^u}(v)$.

\subsubsection{The case of ``$-$''}

Let $J_2\in T(\nbigi)$
such that $J_{2-}\subset\vecI_x(\theta^u)-\pi$.
By using $K^{J_{2-}}_1$
instead of $K^{J_{1+}}_1$,
we obtain the constructible subsheaf
$\nbigm^{J_{2-}}\bigl(
 V\otimes\nbige(u^{-1}x^{-1})
 \bigr)
\subset
 \nbigl^{<0}\bigl(
 V\otimes\nbige(u^{-1}x^{-1})
 \bigr)$.
It induces the desired morphism
$A^{J_{2-}}_{\infty,\theta^u}$ in (\ref{eq;24.2.22.3}).

Let us describe it in terms of $1$-cycles.
Let $\gamma_1$ be a path connecting
$(1,\vartheta^{J_2}_{\ell}+2\pi)$
and $(1,\vartheta^{J_2}_{\ell})$.
We take $\theta_J\in J\cap(\vecI_x(\theta^u)-\pi)$
for each $J\in \gbigk(J_{2-})$,
and path $\gamma_{J,2}$
connecting $(1,\vartheta^{J_2}_{\ell})$
to $(0,\theta_J)$ on $(X^{\ast},X)$.
By shifting $\gamma_{J,2}$ by $2\pi$,
we obtain paths
$\gamma_{J,3}$ connecting
$(1,\vartheta^{J_2}_{\ell}+2\pi)$
to $(0,\theta_J+2\pi)$ on $(X^{\ast},X)$.

Let $v\in H^0(\real,L)$.
We have the decompositions
\begin{equation}
\label{eq;24.3.22.30}
 v'=u_{J_2+2\pi,0}
+\sum_{J\in \gbigk((J_2+2\pi)_-)} u_{J},
\quad\quad
 v=u_{J_2,0}
+\sum_{J\in\gbigk(J_{2-})}u_J,
\end{equation}
where $u_J$ are sections of $L'_{J,<0}$,
$u_{J_2+2\pi,0}$ is a section of
$L_{(J_2+2\pi)_-,0}$,
and
$u_{J_2,0}$ is a section of
$L_{J_{2-},0}$.
We obtain the $1$-cycle
\begin{multline}
\nbiga^{J_{2-}}_{\infty,\theta^u}(v)=
 \vtilde\otimes\gamma_1
-\sum_{J\in\gbigk((J_2+2\pi)_-)}
 u_{J}\otimes\gamma_{3,J}
 -u_{J_2+2\pi,0}\otimes\gamma_{3,J_2+2\pi}
 \\
+\sum_{J\in\gbigk(J_{2-})}
 u_{J}\otimes\gamma_{2,J}
+u_{J_2,0}\otimes\gamma_{2,J_1}.
\end{multline}
Then,
$\varphi_{\ast}(\nbiga^{J_{2-}}_{\infty,\theta^u}(v))$
represents $A^{J_{2-}}_{\infty,\theta^u}(v)$.

\subsubsection{}

\begin{lem}
For $J_1\in T(\nbigi)$ such that
$J_{1+}\subset\vecI_x(\theta^u)-\pi$
and for $v\in H^0(\real,L)$,
\[
 A^{J_{1+}}_{\infty,\theta^u}(v)
 =\sum_{J_1<J'\leq J_1+2\pi}
 \BB_{J',\theta^u}(R_{J'}(v))
 +\Abb^{\rd}_{\infty,\theta^u}(u_{J_1+2\pi,0}).
\]
Here, $u_{J_1+2\pi,0}$ is the section in {\rm(\ref{eq;18.6.22.1})}.
(See {\rm\S\ref{subsection;24.2.18.1}} for the maps $R_{J'}$.)
There exists a similar expression for
$A^{J_{1-}}_{\infty,\theta^u}(v)$.
\hfill\qed
\end{lem}

\subsection{Moderate growth case}

We obtain a constructible subsheaf
\[
\nbigm^{\mg,J_{1\pm}}\bigl(
 V\otimes\nbige(u^{-1}x^{-1})
 \bigr)
\subset
 \nbigl^{\leq 0}
 \bigl(
  V\otimes\nbige(u^{-1}x^{-1})
 \bigr)
\]
by replacing
$\nbigl^{<0}\bigl(
 V\otimes\nbige(u^{-1}x^{-1})
  \bigr)$
with 
$\nbigl^{\leq 0}\bigl(
 V\otimes\nbige(u^{-1}x^{-1})
  \bigr)$
in the construction of 
$\nbigm^{J_{1\pm}}\bigl(
 V\otimes\nbige(u^{-1}x^{-1})
 \bigr)$.
Note that 
$\nbigm^{\mg,J_{1\pm}}\bigl(
 V\otimes\nbige(u^{-1}x^{-1})
 \bigr)$
and 
$\nbigm^{J_{1\pm}}\bigl(
 V\otimes\nbige(u^{-1}x^{-1})
 \bigr)$
are the same outside of
$\varpi^{-1}(\infty)$.
By using $\nbigm^{\mg,J_{1\pm}}\bigl(
 V\otimes\nbige(u^{-1}x^{-1})
 \bigr)$,
we obtain the maps
\index{maps $ A^{\mg,J_{1\pm}}_{\infty,\theta^u}$}
\begin{equation}
 A^{\mg,J_{1\pm}}_{\infty,\theta^u}:
 H_1^{\mg}\bigl(
 \cnum^{\ast},V^{\reg}\otimes
 \nbige(u^{-1}x^{-1})\bigr)
\lrarr
 H_1^{\mg}\bigl(\cnum^{\ast},
 V\otimes\nbige(u^{-1}x^{-1})\bigr).
\end{equation}
As in (\ref{eq;24.2.22.10}),
there exists the isomorphism:
\[
 H^0(\real,L)
 \simeq
 H_1^{\mg}\bigl(
 \cnum^{\ast},
 V\otimes\nbige(x^{-1}u^{-1})
 \bigr).
\]
Hence, we also obtain the maps
\begin{equation}
 A^{\mg,J_{1\pm}}_{\infty,\theta^u}:
H^0(\real,L)
\lrarr
 H_1^{\mg}\bigl(\cnum^{\ast},
 V\otimes\nbige(u^{-1}x^{-1})\bigr).
\end{equation}

Let us describe
$A^{\mg,J_{1+}}_{\infty,\theta^u}$
in terms of $1$-cycles.
Let $\gamma_0$ be a path 
connecting $(\infty,\vartheta^{J_1}_r)$
and $(1,\vartheta^{J_1}_r)$.
Let $v\in H^0(\real,L)$.
There exits a decomposition of $v$
as in (\ref{eq;18.6.22.1}).
Then, 
$A^{\mg,J_{1+}}_{\infty,\theta^u}(v)$
is represented by
\[
 v\otimes\gamma_0
+\sum_{J\in\gbigk(J_{1+})}
 u_{J}\otimes \gamma_{2,J}
+u_{J_1,0}\otimes\gamma_{2,J_1}.
\]
There exists a similar expression
for $A^{\mg,J_{1-}}_{\infty,\theta^u}(v)$.

\begin{lem}
\label{lem;25.2.9.20}
Let $J_1\in T(\nbigi)$ such that $J_{1+}\subset\vecI_x(\theta^u)-\pi$.
For $J\in \gbigk(J_{1+})$,
let $\Gamma_{J}$  be a path connecting
a point in $\{\infty\}\times\real$
and $\{0\}\times J$.
Then,
\[
A^{\mg,J_{1+}}_{\infty,\theta^u}(v)
=-\sum_{J\in \gbigk(J_{1+})}
\Abb^{\mg}_{J,\theta^u}(u_{J})
+\Abb^{\mg,J_{1+}}_{\infty,\theta^u}(u_{J_1,0}).
\] 
There exists a similar expression for
$A^{\mg,J_{1-}}_{\infty,\theta^u}(v)$.
\hfill\qed
\end{lem}

\subsection{Relations}

Let $M$ be the automorphism of $H^0(\real,L)$
obtained as the monodromy.
The following lemmas are clear by the construction.

\begin{lem}
$A^{J_{1\pm}}_{\infty,\theta^u}\circ M
=A^{(J_1+2\pi)_{\pm}}_{\theta^u-2\pi}$
and
$A^{\mg,J_{1\pm}}_{\infty,\theta^u}\circ M
=A^{\mg,(J_1+2\pi)_{\pm}}_{\infty,\theta^u-2\pi}$.
\hfill\qed
\end{lem}

\begin{lem}
\label{lem;25.2.8.1}
In $H_1^{\mg}(\cnum^{\ast},V\otimes\nbige(x^{-1}u^{-1}))$,
for any $v\in H_0(\real,L)$,
we have
\[
 A^{J_{1\pm}}_{\infty,\theta^u}(v)
=A^{\mg,J_{1\pm}}_{\infty,\theta^u}(v)
-A^{\mg,J_{1\pm}}_{\infty,\theta^u}(M(v)).
\]
\hfill\qed
\end{lem}

\begin{prop}
\label{prop;24.3.23.4}
Let $J\in T(\nbigi)$ such that
$\Jbar\subset \vecI_x(\theta^u)-\pi$.
For any $v\in H^0(\real,L)$,
we have
\begin{equation}
\label{eq;24.2.22.11}
 A^{J_{+}}_{\infty,\theta^u}(v)
-A^{J_{-}}_{\infty,\theta^u}(v)
=-\BB_{J,\theta^u}(R_J(v))
+\BB_{J,\theta^u}(R_J(M(v))),
\end{equation}
\begin{equation}
\label{eq;24.2.22.12}
 A^{\mg,J_{+}}_{\infty,\theta^u}(v)
-A^{\mg,J_{-}}_{\infty,\theta^u}(v)
=-\BB_{J,\theta^u}(R_J(v)).
\end{equation}
\end{prop}
\pf
We set 
$W=
\openopen{\vartheta^{J}_{\ell}-\omega^{-1}\pi}{\vartheta^J_r+\omega^{-1}\pi}$.
We consider
\[
 Z_0=\closedopen{0}{\epsilon}\times W,
 \quad
 Z_1=\openopen{0}{\epsilon}\times W,
 \quad
 Z_2=\{0\}\times W.
\]
Let $M$ be the constructible subsheaf of $L$
determined by
$M=L^{<0}$ on $\real\setminus J$
and
$M=L^{\leq 0}$ on $J$.
Let $K$ be the constructible subsheaf of $q_{Z_0}^{-1}(L)$
determined by
\[
 K_{|Z_2}=M_{|Z_2},
 \quad
 K_{Z_1}=q_{Z_1}^{-1}(L^{\leq 0}).
\]
We set $W'=W+2\pi$, and consider
\[
 Z'_0=\closedopen{0}{\epsilon}\times W',
 \quad
 Z'_1=\openopen{0}{\epsilon}\times W',
 \quad
 Z'_2=\{0\}\times W'.
\]
Let $M'$ be the constructible subsheaf of $L$
determined by
$M'=L^{<0}$ on $\real\setminus (J+2\pi)$
and
$M'=L^{\leq 0}$ on $J+2\pi$.
Let $K'$ be the constructible subsheaf of $q_{Z'_0}^{-1}(L)$
determined by
\[
 K'_{|Z'_2}=M'_{|Z'_2},
 \quad
 K'_{Z'_1}=q_{Z'_1}^{-1}(L^{\leq 0}).
\]
We obtain the constructible sheaf
$\iota_{Z_0!}(K)\oplus
\iota_{Z'_0!}(K')$.
It is acyclic with respect to the global cohomology.

The homology class
\[
\alpha=A^{J_{1+}}_{\infty,\theta^u}(v)
-A^{J_{1-}}_{\infty,\theta^u}(v)
+\BB_{J,\theta^u}(R_J(v))
-\BB_{J+2\pi,\theta^u-2\pi}(R_{J+2\pi}(v))
\]
is induced by a $1$-cocycle of
$\nbigc^{\bullet}_{X,\del X}\otimes
\bigl(
\iota_{Z_0!}(K)\oplus
\iota_{Z'_0!}(K')
\bigr)[-2]$.
Hence, we obtain $\alpha=0$
which is (\ref{eq;24.2.22.11}).
We obtain (\ref{eq;24.2.22.12}) similarly.
\hfill\qed

\vspace{.1in}
We obtain the following proposition similarly.
\begin{prop}
\label{prop;24.3.22.31}
Let $J\in T(\nbigi)$ such that $\Jbar\subset \vecI_x(\theta^u)$.
For $v\in H^0(J,L_{J,<0})\subset H^0(\real,L)$,
we have
\[
 A_{J_+,\theta^u}(v)
-A_{J_-,\theta^u}(v) 
=A^{(J-\pi)_-}_{\infty,\theta^u}(v)
-\BB^{\rd}_{J-\pi,\theta^u}(R_{J-\pi}(v)).
\] 
\hfill\qed
\end{prop}

\subsection{Auxiliary isomorphism}

Take $J_1\in T(\nbigi)$
 such that $J_{1+}\subset \vecI_x(\theta^u)-\pi$.
Let $i_0:\varpi^{-1}(0)\to\projtilde^1$ denote the inclusion.
 We have
\[
 \nbigl^{<0}\bigl(
 V\otimes\nbige(u^{-1}x^{-1})
 \bigr)
 \Big/
 \nbigm^{J_{1+}}\bigl(
 V\otimes\nbige(u^{-1}x^{-1})
 \bigr)
\simeq
i_{0!}\bigl(
 \Cok(c_2)
 \bigr).
\]
By the consideration in \S\ref{subsection;18.4.23.1},
we have
$H^i(\varpi^{-1}(0),\Cok(c_2))=0$ for $i\neq 1$,
and 
\[
 H^1\bigl(\varpi^{-1}(0),\Cok(c_2)\bigr)
=\bigoplus_{J\in\gbign(J_{1+})}
 H_0(J,L_{J,<0}).
\]
We obtain the following exact sequence:
\begin{multline}
\label{eq;18.4.22.110}
 0\lrarr
 H_1^{\rd}\bigl(\cnum^{\ast},
 V^{\reg}\otimes\nbige(u^{-1}x^{-1})
 \bigr)
\stackrel{A^{J_{1+}}_{\infty,\theta^u}}{\lrarr}
 H_1^{\rd}\bigl(\cnum^{\ast},
 V\otimes\nbige(u^{-1}x^{-1})
 \bigr)
 \\
\lrarr
 \bigoplus_{J\in\gbign(J_{1+})}
 H_0(J,L_{J,<0})
\lrarr 0.
\end{multline}
Let us construct
a splitting of the exact sequence
(\ref{eq;18.4.22.110}).

For each $J\in \gbign(J_{1+})$,
we take $\theta_J\in J$.
We take a path $\gamma_{10,J}$
connecting $(0,\theta_J)$
and $(1,\vartheta^{J_1}_r)$ on $(X^{\ast},X)$.
There exists the decomposition
\[
 v=u_{J_1,0}+\sum_{J'\in\gbigk(J_{1+})} u_{J'},
\]
where $u_{J_1,0}$ is a section of
$L_{J_{1+},0}$,
and $u_{J'}$ are sections of
$L_{J',<0}$.
Then, we obtain the following
rapid decay cycle of $(V,\nabla)\otimes\nbige(x^{-1}u^{-1})$:
\[
 \varphi_{\ast}\Bigl(
 v\otimes\gamma_{10,J}
+u_{J_{1},0}\otimes\gamma_{2,J_1}
+\sum_{J'\in \gbigk(J_{1+})}
 u_{J'}\otimes \gamma_{2,J'}
 \Bigr).
\]
The homology classes 
are denoted by
$A^{J_{1+}}_{J,\theta^u}(v)$.
Thus, we obtain
\[
 A^{J_{1+}}_{J,\theta^u}:
 H_0(J,L_{J,<0})
\lrarr
 H_1^{\rd}\bigl(\cnum^{\ast},
 V\otimes\nbige(u^{-1}x^{-1})\bigr).
\]
We obtain the following morphism induced by 
$A^{J_{1+}}_{\infty,\theta^u}$
and $A^{J_{1+}}_{J,\theta^u}$ $(J\in\gbign(J_{1+}))$:
\begin{multline}
\label{eq;18.4.22.120}
H_1^{\rd}\bigl(\cnum^{\ast},
 V^{\reg}\otimes\nbige(u^{-1}x^{-1})\bigr)
\oplus
 \bigoplus_{J\in\gbign(J_{1+})}
 H_0(J,L_{J,<0}) \\
\lrarr
 H_1^{\rd}\bigl(\cnum^{\ast},
 V\otimes\nbige(u^{-1}x^{-1})
 \bigr).
\end{multline}
We can show the following lemma
by an argument as in the proof of
Lemma \ref{lem;18.4.22.121}.
\begin{lem}
\label{lem;18.4.22.130}
The morphism {\rm(\ref{eq;18.4.22.120})}
is an isomorphism.
\hfill\qed
\end{lem}

\subsubsection{}

There exists the following exact sequence:
\begin{multline}
\label{eq;18.6.22.10}
0\lrarr
 H_1^{\mg}\bigl(\cnum^{\ast},
 V^{\reg}\otimes\nbige(u^{-1}x^{-1})
 \bigr)
\stackrel{A^{\mg,J_{1+}}_{\infty,\theta^u}}{\lrarr}
 H_1^{\mg}\bigl(\cnum^{\ast},
 V\otimes\nbige(u^{-1}x^{-1})
 \bigr)
 \\
\lrarr
 \bigoplus_{J\in\gbign(J_{1+})}
 H_0(J,L_{J,<0})
\lrarr 0.
\end{multline}
A splitting of the exact sequence
is given by the composition of 
$A^{J_{1+}}_{J,\theta^u}$
and the natural morphism
$H_1^{\rd}\bigl(\cnum^{\ast},
 V\otimes\nbige(u^{-1}x^{-1})
 \bigr)
\lrarr
 H_1^{\mg}\bigl(\cnum^{\ast},
 V\otimes\nbige(u^{-1}x^{-1})
 \bigr)$.

\section{Decomposition of the homology groups of
$(V,\nabla)\otimes\nbige(x^{-1}u^{-1})$}

\subsection{Rapid decay homology group}
\label{subsection;24.2.22.22}

Let $\gbigw_{1}(\nbigi,\theta^u,\pm)$
be the sets of $J\in T(\nbigi)$
such that $J_{\pm}\subset \vecI_x(\theta^u)-\pi$.
\index{sets $\gbigw_1(\nbigi,\theta^u,\pm)$}
Let $\gbigw_{2}(\nbigi,\theta^u,\pm)$
be the sets of $J\in T(\nbigi)$
such that $J_{\pm}\subset \vecI_x(\theta^u)$.
\index{sets $\gbigw_2(\nbigi,\theta^u,\pm)$}

Take $J_1\in T(\nbigi)$
such that $J_{1+}\subset \vecI_x(\theta^u)-\pi$
or $J_{1-}\subset\vecI_x(\theta^u)-\pi$.
We consider the following morphism induced by
$\BB^{\rd}_{J_{\pm},\theta^u}$ $(J\in \gbigw_1(\nbigi,\theta^u,\pm))$,
$A_{J_{\pm},\theta^u}$ $(J\in\gbigw_2(\nbigi,\theta^u,\pm))$
and $A^{J_{1\pm}}_{\infty,\theta^u}$:
\begin{multline}
\label{eq;18.4.22.100}
F^{J_{1\pm}}_{\theta^u}:
 \!\!\!\!\!\bigoplus_{J\in \gbigw_1(\nbigi,\theta^u,\pm)}\!\!\!\!\!\!
 H^0(J,L_{J,>0})
\oplus\!\!\!\!\!\!
 \bigoplus_{J\in\gbigw_2(\nbigi,\theta^u,\pm)}\!\!\!\!\!\!
 H^0(J,L_{J,<0})
\oplus
 H^0(\real,L)
 \\
\lrarr
 H_1^{\rd}\bigl(\cnum^{\ast},
 V\otimes
 \nbige(u^{-1}x^{-1})\bigr).
\end{multline}
We mean that
we consider $F_{\theta^u}^{J_{1+}}$
if $J_{1+}\subset \vecI_x(\theta^u)-\pi$,
and that 
we consider $F_{\theta^u}^{J_{1-}}$
if $J_{1-}\subset \vecI_x(\theta^u)-\pi$.
We may consider both $F_{\theta^u}^{J_{1+}}$ and
$F_{\theta^u}^{J_{1-}}$
if $\Jbar\subset\vecI_x(\theta^u)-\pi$.

\begin{prop}
\label{prop;18.4.20.31}
The morphisms {\rm(\ref{eq;18.4.22.100})}
are isomorphisms.
\end{prop}
\pf
We prove the case of ``$+$''.
The other case is similar.
By Lemma \ref{lem;18.4.22.130},
we obtain the following decomposition:
\[
 H_1^{\rd}\bigl(
 \cnum^{\ast},
 V\otimes\nbige(u^{-1}x^{-1})
 \bigr)
=\Image A^{J_{1+}}_{\infty,\theta^u}
\oplus
 \bigoplus_{J\in\gbign(J_{1+})}
 \Image A^{J_{1+}}_{J,\theta^u}.
\]

We set $\theta_0=-\theta^u+\pi/2$.
Let $\gbign'(J_{1+})$ be the set of
$J\in T(\nbigi)$
such that
$\theta_0-\pi
\leq\vartheta^J_{r}
\leq\vartheta^{J_1}_{\ell}$.
Let $\gbign''(J_{1+})$ be the set of
$J\in T(\nbigi)$
such that
$\vartheta^{J_1}_{r}+\omega^{-1}\pi
<\vartheta^J_r<\theta_0+\omega^{-1}\pi$.
Let $\gbign'''(J_{1+})$ be the set of
$J\in T(\nbigi)$
such that
$\theta_0+\omega^{-1}\pi
\leq\vartheta^J_r<\theta_0+\pi$.
We have 
$\gbign(J_{1+})
=\gbign'(J_{1+})
\sqcup
 \gbign''(J_{1+})
\sqcup
 \gbign'''(J_{1+})$.

Let $\gbigw_1'(\nbigi,+)$ be the set of
$J\in T(\nbigi)$
such that
$\theta_0-\pi\leq \vartheta^J_{\ell}
\leq\vartheta^{J_1}_{\ell}$.
There exits the bijection
$\gbigw_1'(\nbigi,+)\simeq
 \gbign'(J_{1+})$
defined by $J\longmapsto J-\omega^{-1}\pi$.
We can easily observe that
\[
 \bigoplus_{J\in \gbigw_1'(\nbigi,+)}
 \Image B_{J_+,u}
=\bigoplus_{J\in\gbign'(J_{1+})}
 \Image A^{J_{1+}}_{J,u}.
\]
Let $\gbigw_1''(\nbigi,+)$
be the set of
$J\in T(\nbigi)$
such that
$\vartheta^{J_1}_{\ell}<\vartheta^{J}_{\ell}
< \theta_0-\omega^{-1}\pi$.
We have
$\gbigw_1(\nbigi,\theta^u,+)
=\gbigw_1'(\nbigi,+)
\sqcup
 \gbigw_1''(\nbigi,+)$.
We have the bijection
$\gbigw_1''(\nbigi,+)
\simeq
 \gbign''(J_{1+})$
given by
$J\longmapsto J+\omega^{-1}\pi$.
We can easily observe 
\[
 \bigoplus_{J\in\gbigw_1''(\nbigi,+)}
\Image B_{J_+,u}
=\bigoplus_{J\in\gbign''(J_{1+})}
\Image A^{J_{1+}}_{J,u}.
\]
Hence, we obtain the following:
\[
 H_1^{\rd}\bigl(
 \cnum^{\ast},V\otimes\nbige(u^{-1}x^{-1})
 \bigr)
=\Image A^{J_{1+}}_{\infty,u}
\oplus\!\!\!
 \bigoplus_{J\in \gbign'''(J_{1+})}\!\!\!
 \Image A^{J_{1+}}_{J,u}
\oplus\!\!\!
 \bigoplus_{J\in\gbigw_1(\nbigi,\theta^u,+)}\!\!\!
 \Image B_{J_+,u}.
\]
We have
$\Image A^{J_{1+}}_{J,u}
\equiv
 \Image A_{J_+,u}$ modulo
$\bigoplus_{J\in\gbigw_1(\nbigi,\theta^u,+)}
 \Image B_{J_+,u}$.
Thus, the proof of Proposition \ref{prop;18.4.20.31}
is completed.
\hfill\qed

\subsection{Moderate growth homology group}

Take $J_1\in T(\nbigi)$
such that $J_{1+}\subset \vecI_x(\theta^u)-\pi$
or $J_{1-}\subset\vecI_x(\theta^u)-\pi$.
We obtain the following morphism
induced by
$\BB_{J_{\pm},\theta^u}$ $(J\in \gbigw_1(\nbigi,\theta^u,\pm))$,
$A^{\mg}_{J_{\pm},\theta^u}$ $(J\in\gbigw_2(\nbigi,\theta^u,\pm))$
and $A^{\mg,J_{1\pm}}_{\infty,\theta^u}$:
\begin{multline}
\label{eq;18.5.26.100}
F^{\mg,J_{1\pm}}_{\theta^u}:\!\!\!
 \bigoplus_{J\in \gbigw_1(\nbigi,\theta^u,\pm)}
 \!\!\!\!\!H^0(J,L_{J,>0})
\oplus\!\!\!\!
 \bigoplus_{J\in\gbigw_2(\nbigi,\theta^u,\pm)}
 \!\!\!\!\!H^0(J,L_{J,<0})
\oplus
H^0(\real,L)
 \\
\lrarr
 H_1^{\mg}\bigl(\cnum^{\ast},
 V\otimes
 \nbige(u^{-1}x^{-1})\bigr).
\end{multline}
We obtain the following proposition
as a corollary of 
Proposition \ref{prop;18.4.20.31}
and the exact sequences 
(\ref{eq;18.4.22.110})
and (\ref{eq;18.6.22.10}).
\begin{prop}
\label{prop;24.3.29.10}
The morphisms {\rm(\ref{eq;18.5.26.100})}
are isomorphisms.
\hfill\qed
\end{prop}

\section{Homology groups of $(\nbigv,\nabla)\otimes\nbige(zu^{-1})$}
\label{section;24.3.17.2}

We use the notation in \S\ref{section;25.2.5.1},
in particular \S\ref{subsection;25.2.5.2}.
There exists the natural isomorphism
\[
 H_1^{\varrho}(\cnum\setminus D,
 \nbigv\otimes\nbige(zu^{-1}))
 \simeq
 H_1^{\varrho'}(\proj^1\setminus \Dtilde',
 \nbigv'\otimes\nbige(x^{-1}u^{-1})).
\]

\subsection{Lifting maps}
\label{subsection;24.3.29.20}

For an interval $J_1\in T(\nbigi)$
such that $J_{1+}\subset\vecI_x(\theta^u)-\pi$,
we shall construct the following maps:
\index{maps $C^{J_{1\pm}}_{\infty,\theta^u}$}
\begin{equation}
\label{eq;24.3.16.30}
C^{J_{1+}}_{\infty,\theta^u}:
H_1^{\varrho}\bigl(
 \cnum\setminus D,
 \nbigstilde^{\infty}_{\omega}(\nbigv)\otimes\nbige(u^{-1}z)
 \bigr)
\lrarr
H_1^{\varrho}\bigl(
 \cnum\setminus D,
 \nbigv\otimes\nbige(u^{-1}z)
 \bigr),
\end{equation}
or equivalently
\begin{equation}
\label{eq;18.4.21.11}
C^{J_{1+}}_{\infty,\theta^u}:
H_1^{\varrho'}\bigl(
 \proj_x^1\setminus \Dtilde',
 \nbigstilde^0_{\omega}(\nbigv')\otimes\nbige(u^{-1}x^{-1})
 \bigr)
\lrarr
H_1^{\varrho'}\bigl(
 \proj_x^1\setminus \Dtilde',
 \nbigv'\otimes\nbige(u^{-1}x^{-1})
 \bigr).
\end{equation}
For an interval $J_1\in T(\nbigi)$
such that $J_{1-}\subset \vecI_x(\theta^u)-\pi$,
we shall construct
\begin{equation}
C^{J_{1-}}_{\infty,\theta^u}:
H_1^{\varrho}\bigl(
 \cnum\setminus D,
 \nbigstilde^{\infty}_{\omega}(\nbigv)\otimes\nbige(u^{-1}z)
 \bigr)
\lrarr
H_1^{\varrho}\bigl(
 \cnum\setminus D,
 \nbigv\otimes\nbige(u^{-1}z)
 \bigr),
\end{equation}
or equivalently,
\begin{equation}
C^{J_{1-}}_{\infty,\theta^u}:
H_1^{\varrho'}\bigl(
  \proj^1\setminus \Dtilde',
 \nbigstilde^0_{\omega}(\nbigv')\otimes\nbige(u^{-1}x^{-1})
 \bigr)
\lrarr
H_1^{\varrho'}\bigl(
 \proj^1\setminus \Dtilde',
 \nbigv'\otimes\nbige(u^{-1}x^{-1})
 \bigr).
\end{equation}

\subsubsection{Construction of $C^{J_{1+}}_{\infty,\theta^u}$}

Take $J_1\in T(\nbigi)$
such that $J_{1+}\subset\vecI_x(\theta^u)-\pi$.
Let $(L(\nbigv'),\vecnbigf)$
be the local system with Stokes structure 
indexed by $\nbigi_0(\nbigv')$
on $\real$ corresponding to $(\nbigv',\nabla)$ at $x=0$.
Take $\omega_1>\omega$ such that
$\omega_1-\omega$ is sufficiently small
so that
$\nbigstilde_{\omega}(\nbigi_0(\nbigv'))
=\nbigs_{\omega_1}(\nbigi_0(\nbigv'))$.
It implies that
$\nbigs^0_{\omega_1}(\nbigv',\nabla)
=\nbigstilde^0_{\omega}(\nbigv',\nabla)$.

Let $\varpi_{\Dtilde'}:\projtilde_x^1(\Dtilde')\lrarr\proj_x^1$
be the oriented real blow up
of $\proj_x^1$ along $\Dtilde'$.
There exist the constructible subsheaves
$L(\nbigv')_{S^1}^{(\omega_1)\,<0}\subset
 L(\nbigv')_{S^1}^{(\omega_1)\,\leq 0}\subset
 L(\nbigv')_{S^1}$ on $\varpi_{\Dtilde'}^{-1}(0)$
with respect to
$\pi_{\omega_1\ast}(\vecnbigf)$,
and the following holds:
\begin{equation}
\label{eq;18.4.23.30}
 L_{S^1}\simeq
 L(\nbigv')_{S^1}^{(\omega_1)\,\leq 0}\big/
 L(\nbigv')_{S^1}^{(\omega_1)\,<0}.
\end{equation}
(See \S\ref{subsection;18.4.18.1}
for the notation.)
There exists the constructible subsheaf
$\varphi_{1!}K^{J_{1+}}_1\subset L_{S^1}$
as in \S\ref{subsection;18.4.23.1}.
By using (\ref{eq;18.4.23.30}),
we obtain the constructible subsheaf
$K^{J_{1+},\nbigv'}_1\subset L(\nbigv')_{S^1}$
with the exact sequence:
\[
 0\lrarr
 L(\nbigv')_{S^1}^{(\omega_1)\,<0}\lrarr
 K^{J_{1+},\nbigv'}_1
 \lrarr
 \varphi_{1!}K_1^{J_{1+}}
 \lrarr 0.
\]

We have the constructible sheaf
$\nbigl^{\varrho'}\bigl(
 \nbigv'\otimes\nbige(u^{-1}x^{-1})
 \bigr)$
on $\projtilde^1_x(\Dtilde')$.
We also have the constructible subsheaf
\[
 \nbigm^{\varrho',J_{1+}}\bigl(
 \nbigv'\otimes\nbige(u^{-1}x^{-1})
\bigr)
\subset
 \nbigl^{\varrho'}\bigl(
 \nbigv'\otimes\nbige(u^{-1}x^{-1})
 \bigr)
\]
on $\projtilde^1_x(\Dtilde')$
determined by the following conditions.
\begin{equation}
\nbigm^{\varrho',J_{1+}}\bigl(
 \nbigv'\otimes\nbige(u^{-1}x^{-1})
\bigr)_{|\projtilde_x^1(\Dtilde')\setminus\varpi_{\Dtilde'}^{-1}(0)}
=
 \nbigl^{\varrho'}\bigl(
 \nbigv'\otimes\nbige(u^{-1}x^{-1})
 \bigr)_{|\projtilde_x^1(\Dtilde')\setminus\varpi_{\Dtilde'}^{-1}(0)}
\end{equation}
\begin{equation}
 \nbigm^{\varrho',J_{1+}}\bigl(
 \nbigv'\otimes\nbige(u^{-1}x^{-1})
\bigr)_{|\varpi_{\Dtilde'}^{-1}(0)}
=K_1^{J_{1+},\nbigv'}. 
\end{equation}
By the construction,
there exists the following natural monomorphism:
\begin{equation}
\nbigm^{\varrho',J_{1+}}\bigl(
 \nbigv'\otimes\nbige(u^{-1}x^{-1})
 \bigr)
\lrarr
 \nbigl^{\varrho'}\bigl(
\nbigv'\otimes\nbige(u^{-1}x^{-1})
 \bigr).
\end{equation}
There also exists the following natural monomorphism:
\begin{equation}
\label{eq;18.4.21.10}
 \nbigm^{\varrho',J_{1+}}\bigl(
 \nbigv'\otimes\nbige(u^{-1}x^{-1})
 \bigr)
\lrarr
 \nbigl^{\varrho'}\bigl(
 \nbigs^0_{\omega_1}(\nbigv')
 \otimes\nbige(u^{-1}x^{-1})
 \bigr).
\end{equation}
The cokernel of (\ref{eq;18.4.21.10})
is acyclic with respect to the global cohomology,
which we can show by an argument
in the proof of Lemma \ref{lem;18.4.23.41}.
Hence, we obtain the desired map $C_{\infty,\theta^u}^{J_{1+}}$
as the composite of the following maps:
\begin{multline}
 H_1^{\varrho'}\bigl(
 \proj^1\setminus \Dtilde',
 \nbigs^0_{\omega_1}(\nbigv')\otimes\nbige(u^{-1}x^{-1})
 \bigr)
\simeq
 H^1\Bigl(\projtilde_x^1(\Dtilde'),
 \nbigm^{\varrho',J_{1+}}
 \bigl(\nbigv'\otimes\nbige(u^{-1}x^{-1})\bigr)
 \Bigr)
\\
\lrarr
 H_1^{\varrho}\bigl(
 \proj_x^1\setminus \Dtilde',
 \nbigv'\otimes\nbige(u^{-1}x^{-1})
 \bigr).
\end{multline}

\subsubsection{Construction of $C^{J_{1-}}_{\infty,\theta^u}$}

For $J_1\in T(\nbigi)$
such that $J_{1-}\subset\vecI_x(\theta^u)-\pi$,
we construct the constructible subsheaf
$\nbigm^{\varrho',J_{1-}}\bigl(
 \nbigv'\otimes\nbige(u^{-1}x^{-1})
 \bigr)
\subset
 \nbigl^{\varrho'}\bigl(
 \nbigv'\otimes
 \nbige(u^{-1}x^{-1})\bigr)$
by using $K^{J_{1-}}_1$ instead of $K^{J_{1+}}_1$.
Then, we obtain the desired map
$C^{J_{1-}}_{\infty,\theta^u}$
by a similar argument.

\subsection{Some commutative diagrams}

The following lemma is clear by the construction.
\begin{lem}
For any morphism $\varrho_1\lrarr\varrho_2$
in $\Dsf(D)$,
the following diagrams are commutative:
\begin{equation}
\begin{CD}
H_1^{\varrho_1}\bigl(
  \cnum\setminus D,
 \nbigstilde^{\infty}_{\omega}(\nbigv)\otimes\nbige(u^{-1}z)
 \bigr)
@>{C^{J_{1\pm}}_{\infty,\theta^u}}>>
H_1^{\varrho_1}\bigl(
 \cnum\setminus D,
 \nbigv\otimes\nbige(u^{-1}z)
 \bigr)
 \\
 @VVV @VVV \\
H_1^{\varrho_2}\bigl(
 \cnum\setminus D,
 \nbigstilde^{\infty}_{\omega}(\nbigv)\otimes\nbige(u^{-1}z)
 \bigr)
@>{C^{J_{1\pm}}_{\infty,\theta^u}}>>
H_1^{\varrho_2}\bigl(
 \cnum\setminus D,
 \nbigv\otimes\nbige(u^{-1}z)
 \bigr).
\end{CD}
\end{equation}
\hfill\qed
\end{lem}

Recall $\psi(x)=x^{-1}$.
Because
$\psi^{\ast}\Bigl(
\nbigttilde^{\infty}_{\omega}\bigl(
\nbigv
\otimes\nbige(u^{-1}z)
\bigr)
 \Bigr)
=V\otimes\nbige(u^{-1}x^{-1})$,
there exist the following natural morphisms
as explained in \S\ref{subsection;18.5.14.102}:
\begin{multline}
 H_1^{\rd}\bigl(\cnum^{\ast},
 V\otimes\nbige(u^{-1}x^{-1})
 \bigr)
\lrarr
 H_1^{\varrho}\bigl(
 \cnum\setminus D, 
 \nbigv\otimes\nbige(u^{-1}z)
 \bigr)
 \\
\lrarr
  H_1^{\mg}\bigl(\cnum^{\ast},
 V\otimes\nbige(u^{-1}x^{-1})
 \bigr).
\end{multline}
Similarly,
we obtain
\begin{multline}
H_1^{\rd}\bigl(\cnum^{\ast},
 \nbigstilde^0_{\omega}(V)\otimes\nbige(u^{-1}x^{-1})
 \bigr)
\lrarr
 H_1^{\varrho}\bigl(
 \cnum\setminus D, 
\nbigstilde^{\infty}_{\omega}(\nbigv)\otimes\nbige(u^{-1}z)
 \bigr)
 \\
\lrarr
 H_1^{\mg}\bigl(\cnum^{\ast},
 \nbigstilde^0_{\omega}(V)\otimes\nbige(u^{-1}x^{-1})
 \bigr).
\end{multline}
Note that
$(V^{\reg},\nabla)
=\nbigstilde^0_{\omega}(V,\nabla)$.

\begin{prop}
\label{prop;18.4.22.200}
The following diagram is commutative.
\begin{equation}
\label{eq;18.4.23.60}
 \begin{CD}
H_1^{\rd}\bigl(
 \cnum^{\ast},
 V^{\reg}\otimes\nbige(u^{-1}x^{-1})
 \bigr)
@>{A^{J_{1\pm}}_{\infty,\theta^u}}>>
H_1^{\rd}\bigl(
 \cnum^{\ast},
 V\otimes\nbige(u^{-1}x^{-1})
 \bigr)\\
@V{d_1}VV @V{d_2}VV \\
H_1^{\varrho}\bigl(
  \cnum\setminus D,
 \nbigstilde^{\infty}_{\omega}(\nbigv)\otimes\nbige(u^{-1}z)
 \bigr)
@>{C^{J_{1\pm}}_{\infty,\theta^u}}>>
H_1^{\varrho}\bigl(
 \cnum\setminus D,
 \nbigv\otimes\nbige(u^{-1}z)
 \bigr)
\\
 @VVV @VVV \\
H_1^{\mg}\bigl(
 \cnum^{\ast},
V^{\reg}\otimes\nbige(u^{-1}x^{-1})
 \bigr)
@>{A^{\mg,J_{1\pm}}_{\infty,\theta^u}}>>
H_1^{\mg}\bigl(
 \cnum^{\ast},
 V\otimes\nbige(u^{-1}x^{-1})
 \bigr).
 \end{CD}
\end{equation}
The morphisms $C^{J_{1\pm}}_{\infty,\theta^u}$ are injective,
and the following are exact:
\begin{multline}
 \label{eq;18.4.23.61}
0\lrarr
 H_1^{\rd}\bigl(
 \cnum^{\ast},
 V^{\reg}\otimes\nbige(u^{-1}x^{-1})
 \bigr)
\stackrel{f_1^{\pm}}\lrarr \\
H_1^{\varrho}\bigl(
 \cnum\setminus D,
 \nbigstilde^{\infty}_{\omega}(\nbigv)\otimes\nbige(u^{-1}z)
 \bigr)
\oplus
 H_1^{\rd}\bigl(
 \cnum^{\ast},
 V\otimes\nbige(u^{-1}x^{-1})
 \bigr)
\stackrel{f_2^{\pm}}{\lrarr}
 \\
H_1^{\varrho}\bigl(
 \cnum\setminus D,
 \nbigv\otimes\nbige(u^{-1}z)
 \bigr)
\lrarr 0.
\end{multline}
Here, $f_1^{\pm}=d_1+A^{J_{1\pm}}_{\infty,\theta^u}$
and $f_2^{\pm}=C^{J_{1\pm}}_{\infty,\theta^u}-d_2$.
\end{prop}
\pf
We explain the proof in the case of ``$+$''.
The other case can be proved similarly.
To simplify the description,
we omit to denote the superscript ``$J_{1+}$''.
We can replace (\ref{eq;24.3.16.30})
with (\ref{eq;18.4.21.11}).

We take a small positive number $\epsilon>0$.
We consider the subspaces
$Y_{0,\epsilon}:=\closedopen{0}{\epsilon}\times S^1$
and
$Y_{1,\epsilon}:=\closedopen{\epsilon}{\infty}\times S^1$
of $\projtilde^1$.
Let $j_{Y_{i,\epsilon}}\lrarr \projtilde^1$
denote the inclusions
$Y_{i,\epsilon}\lrarr \projtilde^1$.
Let $q_{Y_{i,\epsilon}}$
denote the projection
$Y_{i,\epsilon}\lrarr S^1$.

\subsubsection{}
There exists the following exact sequence on $Y_{0,\epsilon}$:
\begin{equation}
\label{eq;24.3.23.10}
 0\lrarr
 q_{Y_{0,\epsilon}}^{-1}
 \bigl(
 L(\nbigv')_{S^1}^{(\omega_1)\,< 0}
 \bigr)
\lrarr
 q_{Y_{0,\epsilon}}^{-1}
 \bigl(
 L(\nbigv')_{S^1}^{(\omega_1)\,\leq 0}
 \bigr)
\stackrel{h}{\lrarr}
 q_{Y_{0,\epsilon}}^{-1}
 \bigl(
 L_{S^1}
 \bigr)
\lrarr 0.
\end{equation}
We have the constructible subsheaf
$j_{Y_{0,\epsilon}}^{-1}
 \nbigm\bigl(
 V\otimes\nbige(u^{-1}x^{-1})
 \bigr)
\subset
 q_{Y_{0,\epsilon}}^{-1}(L_{S^1})$
(see \S\ref{subsection;25.2.5.3}),
and we obtain the constructible subsheaf
$\check{\nbigm}
 \bigl(V\otimes\nbige(u^{-1}x^{-1})\bigr)$
of 
$q_{Y_{0,\epsilon}}^{-1}
 \bigl(
 L(\nbigv')_{S^1}^{(\omega_1)\,\leq 0}
 \bigr)$
as the pull back of
$j_{Y_{0,\epsilon}}^{-1}
 \nbigm\bigl(
 V\otimes\nbige(u^{-1}x^{-1})
 \bigr)$ 
by $h$.
Similarly,
we have the constructible subsheaf
$j_{Y_{0,\epsilon}}^{-1}
 \nbigl^{<0}\bigl(
 V\otimes\nbige(u^{-1}x^{-1})
 \bigr)
\subset
 q_{Y_{0,\epsilon}}^{-1}(L_{S^1})$,
and we obtain 
$\check{\nbigl}^{<0}
 \bigl(V\otimes\nbige(u^{-1}x^{-1})\bigr)$
as the pull back of 
$j_{Y_{0,\epsilon}}^{-1}
 \nbigl^{<0}\bigl(
 V\otimes\nbige(u^{-1}x^{-1})
 \bigr)$
by $h$.
There exists the following commutative diagram:
\[
 \begin{CD}
 j_{Y_{0,\epsilon}!}
 \check{\nbigm}
 \bigl(V\otimes\nbige(u^{-1}x^{-1})\bigr)
@>{a_1}>>
 \nbigm
 \bigl(V\otimes\nbige(u^{-1}x^{-1})\bigr)
 \\
@VVV @VVV 
 \\ 
 j_{Y_{0,\epsilon}!}
 \check{\nbigl}^{<0}
 \bigl(V\otimes\nbige(u^{-1}x^{-1})\bigr)
 @>{a_2}>>
 \nbigl^{<0}
 \bigl(V\otimes\nbige(u^{-1}x^{-1})\bigr).
 \end{CD}
\]
We have the following:
\[
 \Ker(a_1)=\Ker(a_2)
=j_{Y_{0,\epsilon}!}
 q_{Y_{0,\epsilon}}^{-1}
 \bigl(
 L(\nbigv')_{S^1}^{(\omega_1)\,< 0}
 \bigr),
\quad
 \Cok(a_1)=\Cok(a_2)
=j_{Y_{1,\epsilon}!}
 q_{Y_{1,\epsilon}}^{-1}
 L_{S^1}.
\]
Hence, $\Ker(a_i)$ and $\Cok(a_i)$
are acyclic with respect to the global cohomology.
We obtain the following commutative diagram:
\[
 \begin{CD}
H^1\bigl(
\projtilde^1,
j_{Y_{0,\epsilon}!}
 \check{\nbigm}
 \bigl(V\otimes\nbige(u^{-1}x^{-1})\bigr)
\bigr)
@>{\simeq}>>
H_1^{\rd}\bigl(\cnum^{\ast},
V^{\reg}\otimes\nbige(u^{-1}x^{-1})\bigr)
 \\
@VVV @VVV 
 \\ 
H^1\bigl(
\projtilde^1,
 j_{Y_{0,\epsilon}!}
 \check{\nbigl}^{<0}
 \bigl(V\otimes\nbige(u^{-1}x^{-1})\bigr)
 \bigr)
 @>{\simeq}>>
 H_1^{\rd}\bigl(\cnum^{\ast},
 V\otimes\nbige(u^{-1}x^{-1})\bigr).
 \end{CD}
\]

We may naturally regard $Y_{0,\epsilon}$
as subspaces of $\projtilde^1(\Dtilde')$.
Let $j'_{Y_{0,\epsilon}}$
denote the inclusions
$Y_{0,\epsilon}\lrarr \projtilde^1(\Dtilde')$.
We have the following natural commutative diagram:
{\small
\[
 \begin{CD}
 H^1\bigl(
\projtilde^1\!\!\!,
j_{Y_{0,\epsilon}!}
 \check{\nbigm}
 \bigl(V\otimes\nbige(u^{-1}x^{-1})\bigr)
\bigr)
@>{\simeq}>>
H^1\bigl(\projtilde^1(\Dtilde'),
j'_{Y_{0,\epsilon}!}
 \check{\nbigm}
 \bigl(V\otimes\nbige(u^{-1}x^{-1})\bigr)
\bigr)
\bigr)
 \\
@VVV @VVV 
 \\ 
H^1\bigl(
\projtilde^1\!\!\!,
 j_{Y_{0,\epsilon}!}
 \check{\nbigl}^{<0}
 \bigl(V\otimes\nbige(u^{-1}x^{-1})\bigr)
 \bigr)
 @>{\simeq}>>
H^1\bigl(
\projtilde^1(\Dtilde'),
 j'_{Y_{0,\epsilon}!}
 \check{\nbigl}^{<0}
 \bigl(V\otimes\nbige(u^{-1}x^{-1})\bigr)
 \bigr)
 \end{CD}
\]
}

There exists the following natural commutative diagram:
\[
 \begin{CD}
 j'_{Y_{0,\epsilon}!}
 \check{\nbigm}\bigl(
 V\otimes\nbige(u^{-1}x^{-1})
 \bigr)
 @>{c_1}>>
 \nbigm^{\varrho'}\bigl(
 \nbigv'\otimes\nbige(u^{-1}x^{-1})
 \bigr)
 \\
 @V{c_2}VV @V{c_3}VV \\
 j'_{Y_{0,\epsilon}!}
 \check{\nbigl}^{<0}\bigl(
 V\otimes\nbige(u^{-1}x^{-1})
 \bigr)
 @>{c_4}>>
 \nbigl^{\varrho'}\bigl(
 \nbigv'\otimes\nbige(u^{-1}x^{-1})
 \bigr).
 \end{CD}
\]
The morphisms $c_i$ are monomorphisms,
and the following is exact:
\begin{multline}
0\lrarr
  j'_{Y_{0,\epsilon}!}
  \check{\nbigm}\bigl(
 V\otimes\nbige(u^{-1}x^{-1})
 \bigr)
 \stackrel{c_1+c_2}{\lrarr} \\
  \nbigm^{\varrho'}\bigl(
 \nbigv'\otimes\nbige(u^{-1}x^{-1})
 \bigr)
\oplus
  j'_{Y_{0,\epsilon}!}
 \check{\nbigl}^{<0}\bigl(
 V\otimes\nbige(u^{-1}x^{-1})
 \bigr)
\stackrel{c_3-c_4}{\lrarr}
 \\
 \nbigl^{\varrho'}\bigl(
 \nbigv'\otimes\nbige(u^{-1}x^{-1})
 \bigr)
\lrarr 0.
\end{multline}
We have
$H_i^{\rd}\bigl(
 \cnum^{\ast},
 V\otimes\nbige(u^{-1}x^{-1})
 \bigr)
=H_i^{\rd}\bigl(
 \cnum^{\ast},
 V^{\reg}\otimes\nbige(u^{-1}x^{-1})
 \bigr)=0$ unless $i=1$.
We also have
$H_i^{\varrho'}\bigl(
 \proj^1\setminus \Dtilde',
 \nbigv'\otimes\nbige(u^{-1}x^{-1})
 \bigr)
=H_i^{\varrho'}\bigl(
 \proj^1\setminus \Dtilde',
\nbigstilde^0_{\omega}(\nbigv')\otimes\nbige(u^{-1}x^{-1})
 \bigr)=0$ unless $i=1$.
Hence, we obtain the commutativity 
of the upper square of (\ref{eq;18.4.23.60}),
and the exact sequence (\ref{eq;18.4.23.61}).
Because $A^{J_{1+}}_{\infty,\theta^u}$ is injective,
we obtain the injectivity of $C^{J_{1+}}_{\infty,\theta^u}$.

\subsubsection{}
Let us study the commutativity of the lower square
of the diagram (\ref{eq;18.4.23.60}).
There exists the constructible subsheaf
\[
 j_{Y_{0,\epsilon}}^{-1}\nbigm\bigl(
 V\otimes\nbige(u^{-1}x^{-1})
 \bigr)
\subset
 q_{Y_{0,\epsilon}}^{-1}(L_{S^1})
\subset
 q_{Y_{0,\epsilon}}^{-1}\bigl(
 L(\nbigv')_{S^1}\big/
 L(\nbigv')^{(\omega_1)<0}_{S^1}
 \bigr).
\]
We set
$Y_{2,\epsilon}:=
 \openopen{0}{\epsilon}\times S^1$.
Let $j_{Y_{2,\epsilon}}:Y_{2,\epsilon}\lrarr \projtilde^1$
denote the inclusion.
Let 
$\nbigmhat\bigl(
 V\otimes\nbige(u^{-1}x^{-1})
 \bigr)$
 be the constructible subsheaf
of $q_{Y_{0,\epsilon}}^{-1}
 \bigl(
  L(\nbigv')_{S^1}\big/
 L(\nbigv')^{(\omega_1)\,<0}_{S^1}
 \bigr)$
determined by the following conditions.
\begin{itemize}
\item
$\nbigmhat\bigl(
 V\otimes\nbige(u^{-1}x^{-1})
     \bigr)_{|Y_{2,\epsilon}}
= q_{Y_{0,\epsilon}}^{-1}
 \bigl(
  L(\nbigv')_{S^1}\big/
 L(\nbigv')^{(\omega_1)\,<0}_{S^1}
 \bigr)_{|Y_{2,\epsilon}}$.
\item
$\nbigmhat\bigl(
 V\otimes\nbige(u^{-1}x^{-1})
 \bigr)_{|\{0\}\times S^1}
=\nbigm\bigl(
 V\otimes\nbige(u^{-1}x^{-1})
 \bigr)_{|\{0\}\times S^1}$.
\end{itemize}
Let $j:Y_{2,\epsilon}\lrarr Y_{0,\epsilon}$
denote the inclusion.
We have the following exact sequence:
\begin{multline}
0\lrarr
 j'_{Y_{0,\epsilon}\ast}
 j_{Y_{0,\epsilon}}^{-1}
 \nbigm\bigl(V\otimes\nbige(u^{-1}x^{-1})\bigr)
\lrarr
 j'_{Y_{0,\epsilon}\ast}
\nbigmhat\bigl(
 V\otimes\nbige(u^{-1}x^{-1})
 \bigr) 
\lrarr \\
  j'_{Y_{0,\epsilon}\ast}
 j_!\bigl(
 q_{Y_{2,\epsilon}}^{-1}\bigl(
 L(\nbigv')_{S^1}/L(\nbigv')^{(\omega_1)\,\leq 0}_{S^1}
 \bigr)
 \bigr)
\lrarr 0.
\end{multline}
Note that 
$j'_{Y_{0,\epsilon}\ast}
 j_!\bigl(
 q_{Y_{2,\epsilon}}^{-1}\bigl(
 L(\nbigv')_{S^1}/L(\nbigv')^{\leq 0}_{S^1}
 \bigr)
 \bigr)$
is acyclic with respect to the global cohomology.

Similarly,
we have the constructible subsheaf
\[
\nbiglhat^{\leq 0}(V\otimes\nbige(u^{-1}x^{-1}))
\subset
 q_{Y_{0,\epsilon}}^{-1}\bigl(
 L(\nbigv')_{S^1}/L(\nbigv')^{(\omega_1)\,<0}_{S^1}
 \bigr)
\]
determined by the following conditions.
\begin{itemize}
\item
$\nbiglhat^{\leq 0}(V\otimes\nbige(u^{-1}x^{-1}))
     _{|Y_{2,\epsilon}}
= q_{Y_{0,\epsilon}}^{-1}
 \bigl(
  L(\nbigv')_{S^1}\big/
 L(\nbigv')^{(\omega_1)\,<0}_{S^1}
 \bigr)_{|Y_{2,\epsilon}}$.
\item
$\nbiglhat^{\leq 0}(V\otimes\nbige(u^{-1}x^{-1}))
 _{|\{0\}\times S^1}
= 
 \nbigl^{\leq 0}(V\otimes\nbige(u^{-1}x^{-1}))
 _{|\{0\}\times S^1}$.
\end{itemize}
Then, the cokernel of the natural morphism
\[
 j'_{Y_{0,\epsilon}\ast}
 j_{Y_{0,\epsilon}}^{-1}
  \nbigl^{\leq 0}(V\otimes\nbige(u^{-1}x^{-1}))
\lrarr
  j'_{Y_{0,\epsilon}\ast}
\nbiglhat^{\leq 0}(V\otimes\nbige(u^{-1}x^{-1}))
\]
is acyclic with respect to the global cohomology.
Hence,
the morphism
\begin{multline}
 H^1\bigl(
 \projtilde^1(\Dtilde'),
  j'_{Y_{0,\epsilon}\ast}
\nbigmhat\bigl(
 V\otimes\nbige(u^{-1}x^{-1})
 \bigr) 
 \bigr)
\lrarr \\
 H^1\bigl(
 \projtilde^1(\Dtilde'),
   j'_{Y_{0,\epsilon}\ast}
\nbiglhat^{\leq 0}(V\otimes\nbige(u^{-1}x^{-1}))
 \bigr)
\end{multline}
is identified with 
$A^{\mg,J_{1+}}_{\infty,\theta^u}$.
Then,
we obtain the 
the commutativity of the lower square
of the diagram (\ref{eq;18.4.23.60})
from the following commutative diagram:
\[
 \begin{CD}
 \nbigm^{\varrho'}\bigl(
 \nbigv'\otimes\nbige(u^{-1}x^{-1})
 \bigr)
@>>>
 \nbigl^{\varrho'}\bigl(
 \nbigv'\otimes\nbige(u^{-1}x^{-1})
 \bigr)
 \\ 
 @VVV @VVV \\
   j'_{Y_{0,\epsilon}\ast}
\nbigmhat\bigl(
 V\otimes\nbige(u^{-1}x^{-1})
 \bigr) 
@>>>
  j'_{Y_{0,\epsilon}\ast}
\nbiglhat^{\leq 0}(V\otimes\nbige(u^{-1}x^{-1})).
 \end{CD}
\]
Thus, the proof of Proposition \ref{prop;18.4.22.200}
is completed.
\hfill\qed

\subsection{Decompositions of the homology groups
of $(\nbigv,\nabla)\otimes\nbige(x^{-1}u^{-1})$}

Let $\gbigw_{1}(\nbigi,\theta^u,\pm)$
and $\gbigw_{2}(\nbigi,\theta^u,\pm)$
be as in \S\ref{subsection;24.2.22.22}.
By the composition with
$H_1^{\rd}(\cnum^{\ast},V\otimes
\nbige(u^{-1}x^{-1}))
\to
H_1^{\varrho}(\cnum\setminus D,
\nbigv\otimes
\nbige(u^{-1}z))$,
we obtain the following morphisms:
\[
\BB^{\rd}_{J_{\pm},\theta^u}:
 H^0(J,L_{J,>0})
 \lrarr
 H_1^{\varrho}(\cnum\setminus D,
\nbigv\otimes
\nbige(u^{-1}z))
\quad
(J\in\gbigw_1(\nbigi,\theta^u,\pm)),
\]
\[
A_{J_{\pm},\theta^u}:
 H^0(J,L_{J,<0})
 \lrarr
 H_1^{\varrho}(\cnum\setminus D,
\nbigv\otimes
\nbige(u^{-1}z))
\quad
(J\in\gbigw_2(\nbigi,\theta^u,\pm)).
\]

Take $J_1\in T(\nbigi)$
such that $J_{1+}\subset \vecI_x(\theta^u)-\pi$
or $J_{1-}\subset\vecI_x(\theta^u)-\pi$.
We obtain the following morphism induced by
$\BB^{\rd}_{J_{\pm},\theta^u}$ $(J\in \gbigw_1(\nbigi,\theta^u,\pm))$,
$A_{J_{\pm},\theta^u}$ $(J\in\gbigw_2(\nbigi,\theta^u,\pm))$
and $C^{J_{1\pm}}_{\infty,\theta^u}$:
\begin{multline}
\label{eq;24.2.22.23}
F^{J_{1\pm}}_{\theta^u}:
 \!\!\!\!\!\!\bigoplus_{J\in \gbigw_1(\nbigi,\theta^u,\pm)}\!\!\!\!\!\!\!\!\!
 H^0(J,L_{J,>0})
\oplus\!\!\!\!\!\!\!\!\!
 \bigoplus_{J\in\gbigw_2(\nbigi,\theta^u,\pm)}\!\!\!\!\!\!\!\!
 H^0(J,L_{J,<0})
\oplus
 H_1^{\varrho}\bigl(\cnum\setminus D,
 \nbigstilde^{\infty}_{\omega}(\nbigv)\otimes
 \nbige(u^{-1}z)\bigr)
 \\
\lrarr
 H_1^{\varrho}\bigl(\cnum\setminus D,
 \nbigv\otimes
 \nbige(u^{-1}z)\bigr).
\end{multline}

We obtain the following proposition
from Proposition \ref{prop;18.4.20.31}
and Proposition \ref{prop;18.4.22.200}.

\begin{prop}
The morphisms {\rm(\ref{eq;24.2.22.23})}
are isomorphisms.
\hfill\qed
\end{prop}

\subsection{Difference of lifting maps}

Suppose that $\Jbar_1\subset \vecI_x(\theta^u)-\pi$.
We study the map:
\[
 C^{J_{1+}}_{\infty,\theta^u}
-C^{J_{1-}}_{\infty,\theta^u}:
 H_1^{\varrho}\bigl(
 \cnum\setminus D,
 \nbigstilde^{\infty}_{\omega}(\nbigv)
 \otimes
 \nbige(zu^{-1})
 \bigr)
 \lrarr
 H_1^{\varrho}\bigl(
 \cnum\setminus D,
 \nbigv
 \otimes
 \nbige(zu^{-1})
 \bigr).
\]
\begin{prop}
\label{prop;25.4.6.10}
The map
$C^{J_{1+}}_{\infty,\theta^u}-C^{J_{1-}}_{\infty,\theta^u}$
is identified with the composition of the following morphisms:
\begin{multline}
 H_1^{\varrho}\bigl(
 \cnum\setminus D,
 \nbigstilde^{\infty}_{\omega}(\nbigv)
 \otimes
 \nbige(zu^{-1})
 \bigr)
\stackrel{a_1}{\lrarr}
 H_1^{\mg}\bigl(
 \cnum^{\ast},
 V^{\reg}\otimes\nbige(x^{-1}u^{-1})
 \bigr)
 \simeq
 H^0(\real,L)
\stackrel{R_{J_1}}{\lrarr}
 \\
 H^0(J_1,L_{J_1,>0})
\stackrel{\Abb^{\rd}_{J_1,\theta^u}}{\lrarr}
 H_1^{\rd}\bigl(
 \cnum^{\ast},
 V\otimes\nbige(x^{-1}u^{-1})
 \bigr)
\stackrel{a_2}{\lrarr}
  H_1^{\varrho}\bigl(
 \cnum\setminus D,
 \nbigv
 \otimes
 \nbige(zu^{-1})
 \bigr). 
\end{multline}
Here, $a_1$ and $a_2$ are the natural morphisms. 
\end{prop}
\pf
We use the notation in the proof of
Proposition \ref{prop;18.4.22.200}.
The sheaves
$\check{\nbigm}(V\otimes\nbige(x^{-1}u^{-1}))$
and $\nbigmhat(V\otimes\nbige(x^{-1}u^{-1}))$
are denoted by
$\check{\nbigm}^{J_{1\pm}}$
and $\nbigmhat^{J_{1\pm}}$
to denote the dependence on $J_{1\pm}$.
The sheaves
$\nbigm^{\varrho',J_{1\pm}}(\nbigv'\otimes\nbige(u^{-1}x^{-1}))$
are denoted by
$\nbigm^{\varrho',J_{1\pm}}$.

\subsubsection{}
We consider
$\Jtilde_1=\openopen{\vartheta^{J_1}_{\ell}-\delta}{\vartheta^{J_1}_r+\delta}
\subset \vecI_x(\theta^u)-\pi$
for a sufficiently small $\delta>0$
such that
$\openopen{\vartheta^{J_1}_{\ell}-\delta}{\vartheta^{J_1}_{\ell}}
\cap S_0(\nbigi)=\emptyset$
and
$\openopen{\vartheta^{J_1}_{r}}{\vartheta^{J_1}_{r}+\delta}
\cap S_0(\nbigi)=\emptyset$.
We may assume that
the intervals $I_{10}(J)$ and $I_{11}(J)$
for $K_0^{J_{1\pm}}$ and $K_1^{J_{1\pm}}$
in \S\ref{subsection;18.4.23.1}
are contained in $\Jtilde_1$.

Let $L'_{J_{1\pm},>0}$ be the local subsystems of $L$
determined by
$L'_{J_{1\pm},>0|J_{1\pm}}=L_{J_{1\pm},>0}$.
We obtain the constructible subsheaf
$\gbiga_{J_1}(L)
=L^{\leq 0}+L'_{J_{1+},>0}
=L^{\leq 0}+L'_{J_{1-},>0}$
of $L$.
We note that
$L'_{J_{1\pm},>0|\Jtilde_1\setminus\Jbar_1}\subset
L^{\leq 0}_{\Jtilde_1\setminus \Jbar_1}$.
We also note that
$K_1^{J_{1\pm}}\subset
a_{\Jtilde_1!}a_{\Jtilde_1}^{-1}(\gbiga_{J_1}(L))$.

\begin{lem}
\label{lem;25.2.6.2}
The cokernel of the natural inclusion
$a_{\Jtilde_1!}a_{\Jtilde_1}^{-1}(\gbiga_{J_1}(L))
\to
(a_{\vecI_x(\theta^u)-\pi})_!a_{\vecI_x(\theta^u)-\pi}^{-1}(L)$
is acyclic with respect to the global cohomology. 
\end{lem}
\pf
The cokernel of the natural morphism
$a_{\Jtilde_1!}a_{\Jtilde_1}^{-1}L
\to
(a_{\vecI_x(\theta^u)-\pi})_!a_{\vecI_x(\theta^u)-\pi}^{-1}L$
is acyclic with respect to the global cohomology.
Let $S$ denote the set of $J\in T(\nbigi)$
such that $J\cap J_1\neq\emptyset$ and $J\neq J_1$.
For any $J\in S$,
let $a_{\Jtilde_1,J\cap \Jtilde_1}:J\cap\Jtilde_1\to \Jtilde_1$
and $a_{J,J\cap \Jtilde_1}:J\cap\Jtilde_1\to J$
denote the inclusions.
The cokernel of the natural inclusion
$a_{\Jtilde!}a_{\Jtilde}^{-1}(\gbiga_J(L))
\to
a_{\Jtilde!}a_{\Jtilde}^{-1}(L)$
is isomorphic to
\[
\bigoplus_{J\in S}
 a_{\Jtilde_1!}\bigl(
 (a_{\Jtilde_1,J\cap \Jtilde_1})_{\ast}
 \circ
 (a_{J,J\cap\Jtilde_1})^{-1}(L_{J,>0})
 \bigr).
\]
Thus, we obtain the claim of the lemma.
\hfill\qed

\subsubsection{}

We obtain the constructible subsheaf
$\varphi_{1!}
a_{\Jtilde_1!}a_{\Jtilde_1}^{-1}(\gbiga_{J_1}(L))
\subset L_{S^1}$.
By using (\ref{eq;18.4.23.30}),
we obtain the constructible subsheaf
$\varphi_{1!}
a_{\Jtilde_1!}a_{\Jtilde_1}^{-1}(\gbiga_{J_1}(L))^{\nbigv'}
\subset
L(\nbigv')_{S^1}$
with the following exact sequence
\[
 0\lrarr
 L(\nbigv')_{S^1}^{(\omega_1)\,<0}
 \lrarr
 \varphi_{1!}
 a_{\Jtilde_1!}a_{\Jtilde_1}^{-1}(\gbiga_{J_1}(L))^{\nbigv'}
 \lrarr
 \varphi_{1!}
 a_{\Jtilde_1!}a_{\Jtilde_1}^{-1}(\gbiga_{J_1}(L))
 \lrarr 0.
\]
Let $\nbigm_1^{\varrho',\Jtilde_1}
\subset
\nbigl^{\varrho'}\bigl(\nbigv'\otimes\nbige(u^{-1}x^{-1})\bigr)$
be the constructible subsheaf on
$\projtilde^1_x(\Dtilde')$
which equals
$\nbigl^{\varrho'}\bigl(\nbigv'\otimes\nbige(u^{-1}x^{-1})\bigr)$
outside $\varpi_{\Dtilde'}^{-1}(0)$,
and equals
$\varphi_{1!}
a_{\Jtilde_1!}a_{\Jtilde_1}^{-1}(\gbiga_{J_1}(L))^{\nbigv'}$
on $\varpi_{\Dtilde'}^{-1}(0)$.
There exists the natural monomorphism
\[
\nbigm_1^{\varrho',\Jtilde_1}
\to
\nbigl^{\varrho'}
\bigl(\nbigstilde_{\omega}^0\nbigv'\otimes\nbige(u^{-1}x^{-1})\bigr),
\]
and the cokernel is acyclic with respect to the global cohomology
by Lemma \ref{lem;25.2.6.2}.

\subsubsection{}
Let $\nbigm_1^{\Jtilde_1}\subset
\nbigl^{<0}(V\otimes\nbige(u^{-1}x^{-1}))$
determined by the following conditions.
\begin{itemize}
 \item $\nbigm^{\Jtilde_1}_{1|\projtilde^1\setminus\varpi^{-1}(0)}
       =\nbigl^{<0}(V\otimes\nbige(u^{-1}x^{-1}))
       _{|\projtilde^1\setminus\varpi^{-1}(0)}$.
 \item $\nbigm^{\Jtilde_1}_{1|\varpi^{-1}(0)}
       =\varphi_{1!}a_{\Jtilde_1!}a_{\Jtilde_1}^{-1}(\gbiga_{J_1}(L))$.
\end{itemize}
We obtain the constructible subsheaf
$\check{\nbigm}_1^{\Jtilde_1}$
of
$q_{Y_{0,\epsilon}}^{-1}(L(\nbigv')_{S^1}^{(\omega_1)\,\leq 0})$
as the pull back of
$j_{Y_{0,\epsilon}}^{-1}\nbigm_1^{\Jtilde_1}
\subset q_{Y_{0,\epsilon}}^{-1}(L_{S^1})$
by $h$ in (\ref{eq;24.3.23.10}).

\subsubsection{}
Let $\nbigmhat_1^{\Jtilde_1}$
be the constructible subsheaf 
of $q_{Y_{0,\epsilon}}^{-1}\bigl(L(\nbigv')_{S^1}
\big/L(\nbigv')_{S^1}^{(\omega_1)\,<0}\bigr)$
on $Y_{0,\epsilon}$
such that 
it equals
$q_{Y_{0,\epsilon}}^{-1}\bigl(L(\nbigv')_{S^1}
\big/L(\nbigv')_{S^1}^{(\omega_1)\,<0}\bigr)$
outside of $\varpi^{-1}(0)$,
and that it equals 
$\varphi_{1!}a_{\Jtilde_1!}a_{\Jtilde_1}^{-1}(\gbiga_{J_1}(L))$
on $\varpi^{-1}(0)$.

\subsubsection{}

There exist the natural monomorphisms
$f_{J_{1\pm}}:
 \nbigm^{\varrho',J_{1\pm}}
 \to
 \nbigm^{\varrho',\Jtilde_1}_1$,
and the cokernel are acyclic with respect to the global cohomology.
 We obtain the complex
$\nbigc^{\bullet}(\nbigstilde^0_{\omega}(\nbigv'),\varrho')$:
\[
\begin{CD}
 \nbigm^{\varrho',J_{1+}}
\oplus 
 \nbigm^{\varrho',J_{1-}}
 @>{f_{J_{1+}}-f_{J_{1-}}}>>
 \nbigm^{\varrho',\Jtilde_1}_1.
\end{CD}
\]
Here, the first term sits in the degree $0$.
The projections onto $\nbigm^{\varrho',J_{1\pm}}$
and the inclusions 
$\nbigm^{\varrho',J_{1\pm}}
\to
\nbigl^{\varrho'}\bigl(
 \nbigstilde^0_{\omega}(\nbigv')
 \otimes\nbige(x^{-1}u^{-1})
 \bigr)$
induce the following morphism of complexes:
\[
 m_{\pm}:
\nbigc^{\bullet}(\nbigstilde^0_{\omega}\nbigv',\varrho')
\lrarr
\nbigl^{\varrho'}\bigl(
 \nbigstilde^0_{\omega}(\nbigv')
 \otimes\nbige(x^{-1}u^{-1})
 \bigr).
\]
They induce the isomorphisms of the global cohomology groups
\[
 H(m_{\pm}):H^1\bigl(
 \projtilde^1(\Dtilde'),
\nbigc^{\bullet}(\nbigstilde^0_{\omega}\nbigv',\varrho')
 \bigr)
 \simeq
 H_1^{\varrho'}(\proj^1\setminus \Dtilde',
 \nbigstilde^{0}_{\omega}(\nbigv')
 \otimes\nbige(x^{-1}u^{-1})).
\]
\begin{lem}
$H(m_{+})=H(m_-)$.
\end{lem}
\pf
We consider the complex
$\nbigc^{\prime\bullet}(\nbigstilde_{\omega}\nbigv',\varrho')$
given by
\begin{multline}
 \nbigl^{\varrho'}\bigl(
 \nbigstilde^0_{\omega}(\nbigv')
 \otimes\nbige(x^{-1}u^{-1})
 \bigr)
 \oplus
 \nbigl^{\varrho'}\bigl(
 \nbigstilde^0_{\omega}(\nbigv')
 \otimes\nbige(x^{-1}u^{-1})
 \bigr)
 \\
 \stackrel{\id-\id}{\lrarr}
 \nbigl^{\varrho'}\bigl(
 \nbigstilde^0_{\omega}(\nbigv')
 \otimes\nbige(x^{-1}u^{-1})
 \bigr).
\end{multline}
There exists the natural morphism
$\nbigc^{\bullet}(\nbigstilde_{\omega}\nbigv',\varrho')\to
\nbigc^{\prime\bullet}(\nbigstilde_{\omega}\nbigv',\varrho')$,
which induces the isomorphism of the global cohomology groups.
The projections onto the $j$-th 
$\nbigl^{\varrho'}\bigl(
\nbigstilde^0_{\omega}(\nbigv')
\otimes\nbige(x^{-1}u^{-1})
\bigr)$
in the degree $0$
induce quasi-isomorphisms of complexes
$m'_j:
\nbigc^{\prime\bullet}(\nbigstilde_{\omega}\nbigv',\varrho')
\to
 \nbigl^{\varrho'}\bigl(
 \nbigstilde^0_{\omega}(\nbigv')
 \otimes\nbige(x^{-1}u^{-1})
 \bigr)$,
 which induce the isomorphisms of
the global cohomology groups
\[
 H(m'_j):
 H^1\Bigl(
 \projtilde^1(\Dtilde'),
  \nbigc^{\prime\bullet}(\nbigstilde_{\omega}\nbigv',\varrho')
 \Bigr)
 \simeq
 H^1\Bigl(
 \projtilde^1(\Dtilde'),
 \nbigl^{\varrho'}\bigl(
 \nbigstilde^0_{\omega}(\nbigv')
 \otimes\nbige(x^{-1}u^{-1})
 \bigr)
 \Bigr).
\]
The diagonal embedding of
$\nbigl^{\varrho'}\bigl(
 \nbigstilde^0_{\omega}(\nbigv')
 \otimes\nbige(x^{-1}u^{-1})
 \bigr)$
to the degree $0$ part of 
$\nbigc^{\prime\bullet}(\nbigstilde_{\omega}\nbigv',\varrho')$
induces
a quasi-isomorphism of complexes
$k:
\nbigl^{\varrho'}\bigl(
 \nbigstilde^0_{\omega}(\nbigv')
 \otimes\nbige(x^{-1}u^{-1})
 \bigr)
\to 
\nbigc^{\prime\bullet}(\nbigstilde_{\omega}\nbigv',\varrho')$,
which induces
\[
 H(k):
 H^1\Bigl(
 \projtilde^1(\Dtilde'),
 \nbigl^{\varrho'}\bigl(
 \nbigstilde^0_{\omega}(\nbigv')
 \otimes\nbige(x^{-1}u^{-1})
 \bigr)
 \Bigr)
 \simeq
H^1\Bigl(
 \projtilde^1(\Dtilde'),
  \nbigc^{\prime\bullet}(\nbigstilde_{\omega}\nbigv',\varrho')
 \Bigr).
\]
Both $H(m'_{j})\circ H(k)$ $(j=1,2)$ equal the identity,
and hence we obtain $H(m'_1)=H(m_2')$
from which we obtain $H(m_+)=H(m_-)$.
\hfill\qed

\vspace{.1in}
There exist the natural morphisms
$g_{J_{1\pm}}:
 \nbigm^{\varrho',J_{1\pm}}
 \lrarr
 \nbigl^{\varrho'}\bigl(
 \nbigv'
 \otimes\nbige(x^{-1}u^{-1})
 \bigr)$.
We obtain the morphism of complexes
\begin{equation}
\label{eq;24.3.23.1}
 \nbigc^{\bullet}(\nbigstilde^0_{\omega}\nbigv',\varrho')
 \lrarr
  \nbigl^{\varrho'}\bigl(
 \nbigv'
 \otimes\nbige(x^{-1}u^{-1})
 \bigr)
\end{equation}
induced by $g_{J_{1+}}-g_{J_{1-}}$ at the degree $0$.
It induces the morphism
$C^{J_{1+}}_{\infty,\theta^u}-C^{J_{1-}}_{\infty,\theta^u}$
in the level of the global cohomology.

\subsubsection{}
There exist the natural monomorphisms
\[
 f^{\rd}_{J_{1\pm}}:
j'_{Y_{0,\epsilon}!}\check{\nbigm}^{J_{1\pm}}
\lrarr
j'_{Y_{0,\epsilon}!}\check{\nbigm}^{\Jtilde_1}_1
\]
whose cokernel is acyclic with respect to the global cohomology,
and their cohomology groups are isomorphic to
$H_1^{\rd}\bigl(\cnum^{\ast},
V^{\reg}\otimes\nbige(x^{-1}u^{-1})
\bigr)$.
We obtain the complex
$\nbigc^{\bullet}(V^{\reg},\rd)$
\[
j'_{Y_{0,\epsilon}!}\check{\nbigm}^{J_{1+}}
\oplus
j'_{Y_{0,\epsilon}!}\check{\nbigm}^{J_{1-}}
\stackrel{f^{\rd}_{J_{1+}}-f^{\rd}_{J_{1-}}}{\lrarr}
j'_{Y_{0,\epsilon}!}\check{\nbigm}^{\Jtilde_1}_1.
\]
There exist the following morphisms
\[
  g^{\rd}_{J_{1\pm}}:
j'_{Y_{0,\epsilon}!}\check{\nbigm}^{J_{1\pm}}
\lrarr
j'_{Y_{0,\epsilon}!}\check{\nbigl}^{<0}
\bigl(
 V\otimes\nbige(u^{-1}x^{-1})
\bigr).
\]
We obtain the following morphism of complexes by
$g^{\rd}_{J_{1+}}-g^{\rd}_{J_{1-}}$:
\begin{equation}
\label{eq;24.3.23.2}
 \nbigc^{\bullet}(V^{\reg},\rd)
 \lrarr
 j'_{Y_{0,\epsilon}!}\check{\nbigl}^{<0}
\bigl(
 V\otimes\nbige(u^{-1}x^{-1})
\bigr).
\end{equation}
It induces
$A^{J_{1+}}_{\infty,\theta^u}-A^{J_{1-}}_{\infty,\theta^u}$.

\subsubsection{}

There exist the natural monomorphisms
\[
 f^{\mg}_{J_{1\pm}}:
j'_{Y_{0,\epsilon}\ast}\nbigmhat^{J_{1\pm}}
\lrarr
j'_{Y_{0,\epsilon}\ast}\nbigmhat^{\Jtilde_1}_1
\]
whose cokernel is acyclic with respect to the global cohomology,
and their global cohomology groups are isomorphic to
$H_1^{\mg}\bigl(\cnum^{\ast},
V^{\reg}\otimes\nbige(x^{-1}u^{-1})\bigr)$.
We obtain the complex
$\nbigc^{\bullet}(V^{\reg},\mg)$:
\[
j'_{Y_{0,\epsilon}\ast}\nbigmhat^{J_{1+}}
\oplus
j'_{Y_{0,\epsilon}\ast}\nbigmhat^{J_{1-}}
\stackrel{f^{\mg}_{J_{1+}}-f^{\mg}_{J_{1-}}}{\lrarr}
j'_{Y_{0,\epsilon}\ast}\nbigmhat^{\Jtilde_1}_1.
\]
There exist the following morphisms
\[
  g^{\mg}_{J_{1\pm}}:
j'_{Y_{0,\epsilon}\ast}\nbigmhat^{J_{1\pm}}
\lrarr
j'_{Y_{0,\epsilon}\ast}\nbiglhat^{\leq 0}
\bigl(
 V\otimes\nbige(u^{-1}x^{-1})
\bigr).
\]
We obtain the following morphism of complexes by
$g^{\mg}_{J_{1+}}-g^{\mg}_{J_{1-}}$:
\begin{equation}
\label{eq;24.3.23.3}
 \nbigc^{\bullet}(V^{\reg},\mg)
 \lrarr
 j'_{Y_{0,\epsilon}\ast}\nbiglhat^{\leq 0}
\bigl(
 V\otimes\nbige(u^{-1}x^{-1})
\bigr).
\end{equation}
It induces
$A^{\mg,J_{1+}}_{\infty,\theta^u}-A^{\mg,J_{1-}}_{\infty,\theta^u}$.

\subsubsection{}
We set
$W_0=\closedopen{0}{\epsilon/2}\times \Jtilde_1$
and $W_1=\openopen{0}{\epsilon/2}\times \Jtilde_1$.
Let $M_{\Jtilde_1}$ denote the constructible subsheaf of
$q_{W_0}^{-1}(L)$
determined by the following conditions:
\[
 M_{J_1|W_1}
 =q_{W_1}^{-1}(\gbiga_{J_1}(L)),
 \quad
 M_{J_1|\{0\}\times \Jtilde_1}
 =L^{\leq 0}.
\]
It induces the constructible subsheaf
$\Mhat_{J_1}=
 j_{Y_{0,\epsilon}}^{-1}
\varphi_{\ast}\bigl(
 \iota_{W_0!}(M_{J_1})
\bigr)$
of
$q_{Y_{0,\epsilon}}^{-1}(L_{S^1})$.
(See the proof of Proposition \ref{prop;18.4.22.200}).
Let $\check{M}_{J_1}$ denote the inverse image of $\Mhat_{J_1}$
by the projection $h$ in (\ref{eq;24.3.23.10}).
There exist the following commutative diagram:
\[
\begin{CD}
 j'_{Y_{0,\epsilon}!}
 \check{M}_{J_1}
 @>{d_1}>>
 j'_{Y_{0,\epsilon}!}
 \check{\nbigl}^{<0}
 \bigl(
 V\otimes\nbige(u^{-1}x^{-1})
 \bigr)
 \\
 @V{=}VV @VVV \\
 j'_{Y_{0,\epsilon}!}
 \check{M}_{J_1}
 @>{d_2}>>
 \check{\nbigl}^{\varrho'}
 \bigl(
 \nbigv\otimes\nbige(u^{-1}x^{-1})
 \bigr)
 \\
 @V{b_1}VV @VVV \\
 j'_{Y_{0,\epsilon}!}
 \Mhat_{J_1}
 @>{d_3}>>
 j_{Y_{0,\epsilon},\ast}\nbiglhat^{\leq 0}
 \bigl(
 V\otimes\nbige(u^{-1}x^{-1})
 \bigr).
\end{CD}
\]
We obtain the following complexes
\[
 \nbigc_1^{\bullet}:
  j'_{Y_{0,\epsilon}!}
 \check{\nbigl}^{<0}
 \bigl(
 V\otimes\nbige(u^{-1}x^{-1})
 \bigr)
 \lrarr
 \Cok(d_1),
\]
\[
 \nbigc_2^{\bullet}:
  \check{\nbigl}^{\varrho'}
 \bigl(
 \nbigv\otimes\nbige(u^{-1}x^{-1})
 \bigr)
 \lrarr
 \Cok(d_2),
\]
\[
 \nbigc_3^{\bullet}:
  j_{Y_{0,\epsilon},\ast}\nbiglhat^{\leq 0}
 \bigl(
 V\otimes\nbige(u^{-1}x^{-1})
 \bigr)
 \lrarr
 \Cok(d_3).
\]
There exist the natural quasi-isomorphisms
$j'_{Y_{0,\epsilon}!}
 \check{M}_{J_1}
 \lrarr \nbigc_1^{\bullet}$,
$j'_{Y_{0,\epsilon}!}
 \check{M}_{J_1}
 \lrarr\nbigc_2^{\bullet}$,
and 
$j'_{Y_{0,\epsilon}!}
 \hat{M}_{J_1}
\lrarr \nbigc_3^{\bullet}$.

\subsubsection{}

Let $K_2\subset L_{S^1}$ be the constructible subsheaf
determined by the following conditions.
\begin{itemize}
 \item $K_2=L_{S^1}^{F\,<0}$
       on $S^1\setminus\overline{\varphi_1(\Jtilde_1)}$.
       (See \S\ref{subsection;25.2.6.4} for $L_{S^1}^{F\,<0}$.)
 \item $K_2=\gbiga_{J_1}(L)$
       on $\overline{\varphi_1(\Jtilde_1)}
       \simeq \overline{\Jtilde_1}$.
\end{itemize}
Let $\nbigl_1$ be the constructible subsheaf
of $q_{Y_{0,\epsilon}}^{-1}(L_{S^1})$ on $Y_{0,\epsilon}$
such that
\[
\nbigl_{1|\openopen{0}{\epsilon}\times S^1}
=q_{Y_{0,\epsilon}}^{-1}(L_{S^1})_{1|\openopen{0}{\epsilon}\times S^1},
\quad\quad
\nbigl_{1|\{0\}\times S^1}
=K_2.
\]

We have the constructible subsheaf
$M'_{J_1}=q_{W_0}^{-1}(\gbiga_{J_1}(L))\subset
q_{W_0}^{-1}(L)$.
It induces a constructible subsheaf
$\Mhat'_{J_1}=j_{Y_{\epsilon,0}}^{-1}\varphi_{\ast}\bigl(
\iota_{W_0!}M'_{J_1}
\bigr)$
of $q_{Y_{0,\epsilon}}^{-1}(L_{S^1})$.
There exists the following natural commutative diagram:
\[
\begin{CD}
 \Mhat_{J_1} @>>>
 j_{Y_{0,\epsilon}}^{-1}\nbigl^{<0}(V\otimes\nbige(u^{-1}x^{-1}))
 \\
@VVV @VVV \\
\Mhat'_{J_1}
 @>>>
 \nbigl_1.
 \end{CD}  
\]
The horizontal arrows are monomorphisms.
It induces an isomorphism
\[
j_{Y_{0,\epsilon}}^{-1}\nbigl^{<0}(V\otimes\nbige(u^{-1}x^{-1}))
\big/
\Mhat_{J_1}
\simeq
\nbigl_1\big/
\Mhat'_{J_1}.
\]
There exists the natural morphism
$\nbigm_1^{\Jtilde_1}
\to \nbigl_1$
which induces
\[
 \nbigm_1^{\Jtilde_1}
 \lrarr
 j_{Y_{0,\epsilon}}^{-1}\nbigl^{<0}(V\otimes\nbige(u^{-1}x^{-1}))
\big/
\Mhat_{J_1}.
\]
By using the above consideration,
we obtain the following lemma.

\begin{lem}
There exist natural morphisms
$\nbigm^{\varrho',\Jtilde_1}_1\to
\Cok(d_2)$,
$j'_{Y_{0,\epsilon}!}\check{\nbigm}^{\Jtilde_1}_1\to
\Cok(d_1)$,
and 
$j'_{Y_{0,\epsilon}\ast}
 \nbigmhat^{\Jtilde_1}_1\to
\Cok(d_3)$.
\hfill\qed
\end{lem}

\subsubsection{}
We obtain the following commutative diagram:
\[
 \begin{CD}
  \nbigc^{\bullet}(V^{\reg},\rd)
  @>>>
  \nbigc_1^{\bullet}
@>>>
  j'_{Y_{0,\epsilon}!}\check{\nbigl}^{<0}
\bigl(
 V\otimes\nbige(u^{-1}x^{-1})
  \bigr)
  \\
  @VVV @VVV @VVV \\
 \nbigc^{\bullet}(\nbigstilde^0_{\omega}\nbigv',\varrho)
  @>>>
\nbigc_2^{\bullet}
  @>>>
  \nbigl^{\varrho'}\bigl(
 \nbigv'
 \otimes\nbige(x^{-1}u^{-1})
 \bigr)
  \\
  @VVV @VVV @VVV \\
 \nbigc^{\bullet}(V^{\reg},\mg)
  @>>>
  \nbigc_3^{\bullet}
 @>>>
  j'_{Y_{0,\epsilon}\ast}\nbiglhat^{\leq 0}
\bigl(
 V\otimes\nbige(u^{-1}x^{-1})
\bigr).
 \end{CD}
\]
Here, the composite of the morphisms in the rows
are equal to the morphisms
in (\ref{eq;24.3.23.1}),
(\ref{eq;24.3.23.2}) and (\ref{eq;24.3.23.3}).
We obtain the following commutative diagram:
{
\footnotesize
\[
 \begin{CD}
 H_1^{\rd}(\cnum^{\ast},V^{\reg}\otimes\nbige(x^{-1}u^{-1}))
  @>>>
  H^0(J_1,L_{J_1,>0})
  @>>>
  H_1^{\rd}(\cnum^{\ast},V^{\reg}\otimes\nbige(x^{-1}u^{-1}))
  \\
  @VVV @V{=}VV @VVV \\
  H_1^{\varrho'}(\proj^1\setminus \Dtilde',
  \nbigstilde^0_{\omega}(\nbigv)\otimes\nbige(x^{-1}u^{-1}))
  @>>>
  H^0(J_1,L_{J_1,>0})
  @>>>
  H_1^{\varrho'}(\proj^1\setminus \Dtilde',
  \nbigv\otimes\nbige(x^{-1}u^{-1}))
  \\
  @VVV @V{=}VV @VVV \\
  H_1^{\mg}(\cnum^{\ast},V^{\reg}\otimes\nbige(x^{-1}u^{-1}))
  @>>>
  H^0(J_1,L_{J_1,>0})
  @>>>
  H_1^{\mg}(\cnum^{\ast},V^{\reg}\otimes\nbige(x^{-1}u^{-1})).
 \end{CD}
\]
}
The composition of the morphisms in the rows are
$A^{J_{1+}}_{\infty,\theta^u}-A^{J_{1-}}_{\infty,\theta^u}$,
$C^{J_{1+}}_{\infty,\theta^u}-C^{J_{1-}}_{\infty,\theta^u}$
and
$A^{\mg,J_{1+}}_{\infty,\theta^u}-A^{\mg,J_{1-}}_{\infty,\theta^u}$.
Then, we obtain the claim of Proposition \ref{prop;25.4.6.10}
from Proposition \ref{prop;24.3.23.4}.
\hfill\qed

\section{Stokes filtrations}

\subsection{}

Let $u=|u|e^{\sqrt{-1}\theta^u}$.
If $|u|$ is sufficiently small,
there exist the natural isomorphisms:
\[
 \gbigl^{\gbigf}_{\varrho}(\nbigv)_{|\theta^u}
 \simeq
 H_1^{\varrho}\bigl(
 \cnum\setminus D,
 \nbigv\otimes\nbige(u^{-1}z)
 \bigr).
\]
The Stokes filtrations of
$\gbigl^{\gbigf}_{\varrho}(\nbigv)_{|\theta^u}$
induce
the filtrations
$\nbigf^{\circ\,\theta^u}$
on the spaces
\[
H_1^{\varrho}\bigl(
 \cnum\setminus D,
 \nbigv\otimes\nbige(u^{-1}z)
 \bigr)
\]
indexed by the partially ordered set
$\Bigl(
 \nbigi(\Fourier_+(\nbigv)),
 \leq_{\theta^u}
 \Bigr)$.
Similarly,
we obtain the filtrations
$\nbigf^{\circ\,\theta^u}$
on the space
\[
H_1^{\varrho}\bigl(
 \cnum\setminus D,
 \nbigstilde_{\omega}^{\infty}(\nbigv)\otimes\nbige(u^{-1}z)
 \bigr)
\]
indexed by the partially ordered set
$\Bigl(
\nbigi\bigl(
\Fourier(\nbigstilde^{\infty}_{\omega}(\nbigv))
\bigr),
\leq_{\theta^u}
\Bigr)$.

The following lemma is obvious by the constructions.
(See \S\ref{subsection;24.3.14.42}
for the isomorphism
$\gbigl^{\gbigf}_{\varrho}(\nbigv)_{|\theta^u_1}
 \simeq
 \gbigl^{\gbigf}_{\varrho}(\nbigv)_{|\theta^u_2}$.)
\begin{lem}
Let $\vecJ\in T(\nbigi^{\circ})$.
For any $\theta^u_1,\theta^u_2\in \vecJ_{\mp}$,
we have 
\[
 A_{\nu_0^-(\vecJ)_{\pm},\theta^u_1}
=A_{\nu_0^-(\vecJ)_{\pm},\theta^u_2},
 \quad
 \BB^{\rd}_{\nu_0^+(\vecJ)_{\pm},\theta^u_1}
=\BB^{\rd}_{\nu_0^+(\vecJ)_{\pm},\theta^u_2},
\quad
 C^{\nu_0^+(\vecJ)_{\pm}}_{\infty,\theta^u_1}
 =C^{\nu_0^+(\vecJ)_{\pm}}_{\infty,\theta^u_2}
\]
under the natural isomorphisms
$\gbigl^{\gbigf}_{\varrho}(\nbigv)_{|\theta^u_1}
 \simeq
 \gbigl^{\gbigf}_{\varrho}(\nbigv)_{|\theta^u_2}$.
(See {\rm\S\ref{subsection;25.2.6.1}} for $\nu_0^{\pm}$.)
\hfill\qed
\end{lem}

\subsection{}
\label{subsection;24.4.4.10}

Recall $V=\nbigttilde^{\infty}_{\omega}(\nbigv)$,
$\nbigitilde=\nbigi(V)$
and $\nbigi=\pi_{\omega}(\nbigi)$.
We set
$\nbigi^{\circ}=\gbigf^{(\infty,\infty)}_+(\nbigi)\cup\{0\}$
and
$\nbigitilde^{\circ}=\gbigf^{(\infty,\infty)}_+(\nbigitilde)\cup\{0\}$.
We set
\index{sets $\gbigm_{\pm}(\nbigi^{\circ},\theta^u)$}
\[
 \gbigm_-(\nbigi^{\circ},\theta^u)
 =\bigl\{
 \vecJ\in T(\nbigi^{\circ})\,\big|\,
 \theta^u\in\vecJ_-
 \bigr\},
 \quad
  \gbigm_+(\nbigi^{\circ},\theta^u)
 =\bigl\{
 \vecJ\in T(\nbigi^{\circ})\,\big|\,
 \theta^u\in\vecJ_+
 \bigr\}.
\]

\begin{lem}
For any $\vecJ\in T(\nbigi^{\circ})$,
the following conditions are equivalent. 
\begin{itemize}
 \item $\vecJ\in\gbigm_-(\nbigi^{\circ},\theta^u)$
 \item $\nu_0^-(\vecJ)_+\subset \vecI_x(\theta^u)$.
 \item $\nu_0^+(\vecJ)_+\subset \vecI_x(\theta^u)-\pi$.
\end{itemize}
In the case,
$\kappa^-_{0,\vecJ}(\theta^u)
 \in
 \nu^-_{0}(\vecJ)_-$
and
$\kappa^+_{0,\vecJ}(\theta^u)
 \in
 \nu^+_{0}(\vecJ)_-$.

Similarly, the following conditions
are equivalent.
\begin{itemize}
 \item $\vecJ\in\gbigm_+(\nbigi^{\circ},\theta^u)$.
 \item $\nu_0^-(\vecJ)_-\subset \vecI_x(\theta^u)$.
 \item $\nu_0^+(\vecJ)_-\subset \vecI_x(\theta^u)-\pi$.
\end{itemize} 
In the case,
$\kappa^-_{0,\vecJ}(\theta^u)
 \in
 \nu^-_{0}(\vecJ)_+$
and
$\kappa^+_{0,\vecJ}(\theta^u)
 \in
 \nu^+_{0}(\vecJ)_+$.
\hfill\qed
\end{lem}

Recall that there exist the isomorphisms
of the partially ordered sets
in Proposition \ref{prop;18.5.6.1}:
\begin{equation}
\label{eq;24.3.14.140}
(\nbigitilde^{\circ}_{\vecJ,>0},\leq_{\theta^u})
\simeq
(\nbigitilde_{\nu^-_0(\vecJ),<0},\leq_{\kappa^-_{0,\vecJ}(\theta^u)}),
\quad
(\nbigitilde^{\circ}_{\vecJ,<0},\leq_{\theta^u})
\simeq
(\nbigitilde_{\nu^+_0(\vecJ),>0},\leq_{\kappa^+_{0,\vecJ}(\theta^u)}).
\end{equation}
When $\theta_u\in \vecJbar$,
we obtain the filtration $\nbigf^{\prime \theta^u}$
of
\[
 H^0(\nu^-_0(\vecJ),L_{\nu_0^-(\vecJ),<0})
 \simeq
 H^0(\overline{\nu^-_0(\vecJ)},L_{\nu_0^-(\vecJ),<0})
\]
indexed by the partially ordered set
$(\nbigitilde^{\circ}_{\vecJ,>0},\leq_{\theta^u})$
from the filtration
$\nbigf^{\kappa^-_{0,\vecJ}(\theta^u)}$
indexed by
$(\nbigitilde_{\nu^-_0(\vecJ),<0},\leq_{\kappa^-_{0,\vecJ}(\theta^u)})$.
We also obtain the filtration $\nbigf^{\prime \theta^u}$
of
\[
 H^0(\nu^+_0(\vecJ),L_{\nu_0^+(\vecJ),>0})
 \simeq
 H^0(\overline{\nu^+_0(\vecJ)},L_{\nu_0^+(\vecJ),>0})
\]
indexed by the partially ordered set
$(\nbigitilde^{\circ}_{\vecJ,<0},\leq_{\theta^u})$
from the filtration
$\nbigf^{\kappa^+_{0,\vecJ}(\theta^u)}$
indexed by
$(\nbigitilde_{\nu^+_0(\vecJ),>0},\leq_{\kappa^+_{0,\vecJ}(\theta^u)})$.

\subsection{Isomorphisms of the filtered vector spaces}

For $\vecJ_1\in \gbigm_-(\nbigi^{\circ},\theta^u)$,
we obtain the following isomorphism
induced by
$\BB^{\rd}_{\nu_0^+(\vecJ)_+,\theta^u}$,
$A_{\nu_0^-(\vecJ)_+,\theta^u}$
$(\vecJ\in \gbigm_-(\nbigi^{\circ},\theta^u))$
and
$C^{\nu_0^+(\vecJ_1)_+}_{\infty,\theta^u}$:
{\footnotesize
\begin{multline}
\label{eq;24.3.14.110}
\!\!\!\!\!
 \bigoplus_{\vecJ\in\gbigm_-(\nbigi^{\circ},\theta^u)}\!\!
 \Bigl(
 H^0(\nu^+_0(\vecJ),L_{\nu_0^-(\vecJ),>0})
\oplus
 H^0(\nu^-_0(\vecJ),L_{\nu_0^+(\vecJ),<0})
  \Bigr)
 \oplus
 H_1^{\varrho}\bigl(
 \cnum\setminus D,
 \nbigstilde^{\infty}_{\omega}(\nbigv)\otimes
 \nbige(u^{-1}z)
 \bigr)
 \\
\stackrel{\simeq}{\lrarr}
 H_1^{\varrho}\bigl(
 \cnum\setminus D,
 \nbigv\otimes
 \nbige(u^{-1}z)
 \bigr).
\end{multline}
}
Similarly,
for $\vecJ_1\in \gbigm_+(\nbigi^{\circ},\theta^u)$,
we obtain the following isomorphism
induced by
$\BB^{\rd}_{\nu_0^+(\vecJ)_-,\theta^u}$,
$A_{\nu_0^-(\vecJ)_-,\theta^u}$
$(\vecJ\in \gbigm_+(\nbigi^{\circ},\theta^u))$
and
$C^{\nu_0^+(\vecJ_1)_-}_{\infty,\theta^u}$:
{\footnotesize
\begin{multline}
\label{eq;24.3.14.111}
\!\!\!\!\!
 \bigoplus_{\vecJ\in\gbigm_+(\nbigi^{\circ},\theta^u)}\!\!
 \Bigl(
 H^0(\nu^+_0(\vecJ),L_{\nu_0^-(\vecJ),>0})
\oplus
 H^0(\nu^-_0(\vecJ),L_{\nu_0^+(\vecJ),<0})
  \Bigr)
 \oplus
 H_1^{\varrho}\bigl(
 \cnum\setminus D,
 \nbigstilde^{\infty}_{\omega}(\nbigv)\otimes
 \nbige(u^{-1}z)
 \bigr)
 \\
\stackrel{\simeq}{\lrarr}
 H_1^{\varrho}\bigl(
 \cnum\setminus D,
 \nbigv\otimes
 \nbige(u^{-1}z)
 \bigr).
\end{multline}
}
Note that
$\nbigi\bigl(
 \Fourier(\nbigv(\varrho))
 \bigr)
 \subset
 \nbigi\Bigl(
 \Fourier\bigl(\nbigstilde^{\infty}_{\omega}(\nbigv(\varrho))\bigr)
 \Bigr)
 \sqcup
 \nbigitilde^{\circ}$.
The left hand side of (\ref{eq;24.3.14.110})
and (\ref{eq;24.3.14.111})
are equipped with the filtrations
$\nbigf^{\prime\theta^u}$
obtained from the filtrations
$\nbigf^{\prime\theta^u}$
on  $H^0(\nu_{0}^{-}(\vecJ),L_{\nu^-_0(\vecj),<0})$
and $H^0(\nu_{0}^{+}(\vecJ),L_{\nu^+_0(\vecj),<0})$
$(\vecJ\in \gbigm_{\pm}(\nbigi^{\circ},\theta^u))$,
and $\nbigf^{\circ\,\theta^u}$
on 
$H_1^{\varrho}\bigl(
 \cnum\setminus D,
 \nbigstilde^{\infty}_{\omega}(\nbigv)\otimes
 \nbige(u^{-1}z)
 \bigr)$.
The right hand side of (\ref{eq;24.3.14.110})
and (\ref{eq;24.3.14.111})
are equipped with the filtrations
$\nbigf^{\circ\theta^u}$
induced by the Stokes filtrations of
$\gbigl^{\gbigf}_{\varrho}(\nbigv)$.
The following is one of the main theorems,
which we shall prove in \S\ref{section;18.6.3.21}.
\begin{thm}
\label{thm;24.3.16.20}
The isomorphisms {\rm(\ref{eq;24.3.14.110})}
and {\rm(\ref{eq;24.3.14.111})}
are isomorphisms of filtered vector spaces. 
\end{thm}

\subsection{Some canonically defined spaces}
\label{subsection;24.3.29.100}

\begin{cor}
\label{cor;25.4.12.10}
For any $\theta^u\in \vecJ_{\mp}$,
$A_{\nu_0^-(\vecJ)_{\pm},\theta^u}$ induces an isomorphism
of filtered vector spaces
\[
 H^0(\nu_0^-(\vecJ),L_{\nu_0^-(\vecJ),<0})
 \simeq
 H^0(\vecJ_{\mp},\gbigl^{\gbigf}_{\varrho}(\nbigv)_{\vecJ_{\mp},>0}),
\] 
and 
$\BB^{\rd}_{\nu_0^+(\vecJ)_{\pm},\theta^u}$
induces an isomorphism of filtered vector spaces
\[
 H^0(\nu_0^+(\vecJ),L_{\nu_0^+(\vecJ),>0})
 \simeq
 H^0(\vecJ,\gbigl^{\gbigf}_{\varrho}(\nbigv)_{\vecJ,<0}).
\] 
Here, we use the isomorphisms of
the partially ordered sets in {\rm(\ref{eq;24.3.14.140})}
to identify the index sets of the filtrations.
\hfill\qed
\end{cor}

\begin{cor}
For any $\theta^u\in\vecJ_{\mp}$,
$C^{\nu_0^-(\vecJ)_{\pm}}_{\infty,\theta^u}$ induces
an isomorphism of filtered vector spaces
\[
 H^0(\vecJ_{\mp},\gbigl^{\gbigf}_{\varrho}(\nbigstilde_{\omega}(\nbigv)))
 \simeq
 H^0(\vecJ_{\mp},\gbigl^{\gbigf}_{\varrho}(\nbigv)_{\vecJ_{\mp},0}).
\] 
\hfill\qed
\end{cor}

\subsection{Transformations of cycles adapted to the Stokes filtrations}
\label{subsection;24.3.29.101}

Let $\vecJ\in T(\nbigi^{\circ})$.

\subsubsection{}

Let $\theta^u\in\vecJ$.
By Lemma \ref{lem;25.2.8.1} and Proposition \ref{prop;24.3.23.4},
for any $v\in H^0(\real,L)$,
we obtain
\begin{equation}
 A^{\nu_0^+(\vecJ)_{\pm}}_{\infty,\theta^u}(v)
 =
 A^{\mg,\nu_0^+(\vecJ)_{\pm}}_{\infty,\theta^u}(v)
- A^{\mg,\nu_0^+(\vecJ)_{\pm}}_{\infty,\theta^u}(M(v)),
\end{equation}
\begin{equation}
 A^{\nu_0^+(\vecJ)_+}_{\infty,\theta^u}(v)
-A^{\nu_0^+(\vecJ)_-}_{\infty,\theta^u}(v)
=-\BB^{\rd}_{\nu_0^+(\vecJ),\theta^u}(R_{\nu_0^+(\vecJ)}(v))
+\BB^{\rd}_{\nu_0^+(\vecJ),\theta^u}(R_{\nu_0^+(\vecJ)}(M(v))),
\end{equation}
\begin{equation}
 A^{\nu_0^+(\vecJ)_+,\mg}_{\infty,\theta^u}(v)
-A^{\nu_0^+(\vecJ)_-,\mg}_{\infty,\theta^u}(v)
=-\BB^{\rd}_{\nu_0^+(\vecJ),\theta^u}(R_{\nu_0^+(\vecJ)}(v)).
\end{equation}

\subsubsection{}

Let $\theta^u\in\vecJ$.
For $v\in H^0(\nu_0^-(\vecJ),L_{\nu_0^-(\vecJ),<0})$,
Proposition \ref{prop;24.3.22.31} implies
\begin{equation}
 A_{\nu_0^-(\vecJ)_+,\theta^u}(v)
-A_{\nu_0^-(\vecJ)_-,\theta^u}(v)
=A^{\nu_0^+(\vecJ)_-}_{\infty,\theta^u}(v)
-\BB^{\rd}_{\nu_0^+(\vecJ),\theta^u}(R_{\nu_0^+(\vecJ)}(v)).
\end{equation}

\subsubsection{}

Let $\theta^u=\vartheta^{\vecJ}_{\ell}$.
For any $v\in H^0(\nu_0^-(\vecJ),L_{\nu_0^-(\vecJ),<0})$,
Proposition \ref{prop;24.3.17.1} implies
\begin{multline}
 A_{\nu_0^-(\vecJ)_+,\theta^u}(v)
 =\BB^{\rd}_{\nu_0^+(\vecJ-(1-\omega^{-1})\pi),\theta^u}
 \bigl(
 R_{\nu_0^+(\vecJ-(1-\omega^{-1})\pi)}(v)
 \bigr)
 \\
 +\sum_{\vecJ-(1-\omega^{-1})\pi<\vecJ'<\vecJ}
 \BB^{\rd}_{\nu_0^+(\vecJ'),\theta^u}
 \bigl(
 R_{\nu_0^+(\vecJ')}(v)
 \bigr).
\end{multline}
Note that
$\nu_0^+(\vecJ-(1-\omega^{-1})\pi)
=\nu_0^-(\vecJ)-\omega^{-1}\pi$,
and that
$R_{\nu_0^+(\vecJ-(1-\omega^{-1})\pi)}$ induces an isomorphism
\begin{multline}
R_{\nu_0^+(\vecJ-(1-\omega^{-1})\pi)}:
H^0(\nu_0^-(\vecJ),L_{\nu_0^-(\vecJ),<0})
 \simeq
 \\
H^0(\nu_0^+(\vecJ-(1-\omega^{-1})\pi),
 L_{\nu_0^+(\vecJ-(1-\omega^{-1})\pi),>0})
\end{multline}
which preserves the filtrations $\nbigf^{\prime\theta^u}$.

\subsubsection{}

Let $\theta^u=\vartheta^{\vecJ}_{r}$.
For any $v\in H^0(\nu_0^-(\vecJ),L_{\nu_0^-(\vecJ),<0})$,
Proposition \ref{prop;24.3.17.1} also implies
\begin{multline}
 A_{\nu_0^-(\vecJ)_-,\theta^u}(v)
=-\BB^{\rd}_{\nu_{-1}^+(\vecJ+(1-\omega^{-1})\pi),\theta^u-2\pi}
 \bigl(
 R_{\nu_{-1}^+(\vecJ+(1-\omega^{-1})\pi)}(v)
 \bigr)
 \\
 -\sum_{\vecJ<\vecJ'<\vecJ+(1-\omega^{-1})\pi}
 \BB^{\rd}_{\nu_{-1}^+(\vecJ'),\theta^u-2\pi}
 \bigl(
 R_{\nu_{-1}^+(\vecJ')}(v)
 \bigr).
\end{multline}
Note that
$\nu_{-1}^+(\vecJ+(1-\omega^{-1})\pi)
=\nu_0^-(\vecJ)+\omega^{-1}\pi$,
and that
$R_{\nu_0^+(\vecJ+(1-\omega^{-1})\pi)}$ induces an isomorphism
\begin{multline}
R_{\nu_0^+(\vecJ+(1-\omega^{-1})\pi)}:
H^0(\nu_0^-(\vecJ),L_{\nu_0^-(\vecJ),<0})
\simeq \\ 
 H^0(\nu_{-1}^+(\vecJ+(1-\omega^{-1})\pi),
L_{\nu_{-1}^+(\vecJ+(1-\omega^{-1})\pi),>0})
\end{multline}
which preserves the filtrations $\nbigf^{\prime\theta^u}$.

\subsection{The induced constructible sheaves and filtrations}

Let $\theta^u\in\real$.
There exist the following isomorphisms for $\star=!,\ast$
induced by $\BB^{\rd}_{\nu_0^+(\vecJ)_+,\theta^u}$:
\[
 \gbigl^{\gbigf}_{\star}(V)^{<0}_{\theta^u}
=\bigoplus_{\theta^u\in\vecJ}
 H^0(\vecJ,\gbigl^{\gbigf}_{\star}(V)_{\vecJ,<0})
 \simeq
 \bigoplus_{\theta^u\in\vecJ}
 H^0(\nu_0^+(\vecJ),L_{\nu_0^+(\vecJ),>0}).
\]
If $\theta^u\in\real\setminus S_0(\nbigi^{\circ})$,
there exist the following isomorphisms
induced by $A_{\nu_0^-(\vecJ)_{\pm},\theta^u}$:
\[
\Bigl(
\gbigl^{\gbigf}_{\star}(V)\big/
\gbigl^{\gbigf}_{\star}(V)^{\leq 0}
\Bigr)_{\theta^u}
\simeq
\bigoplus_{\theta^u\in\vecJ}
 H^0(\vecJ,\gbigl^{\gbigf}_{\star}(V)_{\vecJ,>0})
\simeq
 \bigoplus_{\theta^u\in\vecJ}
 H^0(\nu_0^-(\vecJ),L_{\nu_0^-(\vecJ),<0}).
\]
If $\theta^u\in S_0(\nbigi^{\circ})$,
we set
$\vecJ^1(\theta^u)=
\openopen{\theta^u-\frac{\omega\pi}{(\omega-1)}}{\theta^u}$
and
$\vecJ^2(\theta^u)=
\openopen{\theta^u}{\theta^u+\frac{\omega\pi}{\omega-1}}$.
There exist the following isomorphisms induced by
$A_{\nu_0^-(\vecJ)_{\pm},\theta^u}$,
$A_{\nu_0^-(\vecJ^1(\theta^u))_{-},\theta^u}$
and
$A_{\nu_0^-(\vecJ^2(\theta^u))_{+},\theta^u}$:
\begin{multline}
\Bigl(
\gbigl^{\gbigf}_{\star}(V)\big/
\gbigl^{\gbigf}_{\star}(V)^{\leq 0}
\Bigr)_{\theta^u}
 \simeq
 \bigoplus_{\theta^u\in\vecJ}
 H^0\Bigl(\nu_0^-(\vecJ),L_{\nu_0^-(\vecJ),<0}\Bigr)
\\
 \oplus 
 H^0\Bigl(
 \nu_0^-(\vecJ^1(\theta^u))_-,
 L_{\nu_0^-(\vecJ^1(\theta^u))_-,<0}\Bigr)
 \oplus
 H^0\Bigl(
 \nu_0^-(\vecJ^2(\theta^u))_+,
 L_{\nu_0^-(\vecJ^2(\theta^u))_+,<0}\Bigr).
\end{multline}
The isomorphisms also induces isomorphisms of
the filtered vector spaces
as in Corollary {\rm\ref{cor;25.4.12.10}}.

\begin{rem}
\label{rem;25.4.12.21}
To the best of the author's understanding,
we may deduce the above descriptions of
the constructible sheaves 
$\gbigl^{\gbigf}_{\star}(V)^{<0}$
and
$\gbigl^{\gbigf}_{\star}(V)\big/
\gbigl^{\gbigf}_{\star}(V)^{\leq 0}$
and the Stokes filtrations 
by applying the results in
{\rm\cite[\S X]{Malgrange-book}}
to the cases $V(\star 0)$ $(\star=!,\ast)$.
\hfill\qed
\end{rem}

\section{Local Fourier transforms of Stokes structure
from $\infty$ to $\infty$}
\label{subsection;24.4.5.120}

To describe $(\gbigl^{\gbigf}_{\star}(V_{\infty}),\vecnbigf)$,
it is convenient to introduce the local Fourier transform
of a Stokes structure.

\subsection{$2\pi\seisuu$-equivariant
local system $\gbigq^{\infty}_!(L,\vecnbigftilde)_{\real}$}
\label{subsection;24.3.26.20}

We consider the vector space
\begin{equation}
\label{eq;24.3.21.30}
 H^0(\real,\nbigt_{\omega}(L))
 \oplus
 \bigoplus_{J\in T(\nbigi)}
 H^0(J,L_{J,>0}).
\end{equation}
An element $v\in H^0(J,L_{J,>0})$
is denoted as a pair
$\langle J,v\rangle$.

Let $\gbigq^{\infty}_!(L,\vecnbigftilde)$
denote the quotient space of (\ref{eq;24.3.21.30})
by the equivalence relation generated by
\index{vector space \mbox{$\gbigq^{\infty}_!(L,\vecnbigftilde)$}}
\[
\langle J,v\rangle
\sim
\langle J-2\pi,\Tbb^{\ast}(v)\rangle
+\nbigq_{J_+}(v)
\quad
(J\in T(\nbigi),\,\,v\in H^0(J,L_{J,>0})).
\]
(See Lemma \ref{lem;24.4.5.31}.)
Let $\Tbb^{\ast}_{\gbigq^{\infty},!}$ be the automorphism of
$\gbigq^{\infty}_!(L,\vecnbigftilde)$
induced by
the automorphism of
(\ref{eq;24.3.21.30})
given by
$M_0^{-1}$
on $H^0(\real,\nbigt_{\omega}(L))$,
and $(\Tbb^{\ast})^{-1}:
H^0(J-2\pi,L_{J-2\pi,>0})
\simeq
H^0(J,L_{J,>0})$.
(See Lemma \ref{lem;24.2.23.4} and Lemma \ref{lem;24.4.5.30}.)

Let $\gbigq^{\infty}_!(L,\vecnbigftilde)_{\real}$
denote the local system on $\real$
induced by $\gbigq^{\infty}_!(L,\vecnbigftilde)$.
\index{local system \mbox{$\gbigq^{\infty}_!(L,\vecnbigftilde)_{\real}$}}
We naturally identify
$H^0(\real,\gbigq^{\infty}_!(L,\vecnbigftilde)_{\real})$
with $\gbigq^{\infty}_!(L,\vecnbigftilde)$.
There exists the $2\pi\seisuu$-action on
$\gbigq^{\infty}_!(L,\vecnbigftilde)$
such that
the pull back
$\Tbb^{\ast}:H^0(\real,\gbigq^{\infty}_!(L,\vecnbigftilde)_{\real})
\simeq
H^0(\real,\gbigq^{\infty}_!(L,\vecnbigftilde)_{\real})$
equals $\Tbb^{\ast}_{\gbigq^{\infty},!}$.

\begin{prop}
There exists the isomorphism of
$2\pi\seisuu$-equivariant local systems
 $\gbigq^{\infty}_!(L,\vecnbigftilde)_{\real}
 \simeq
\gbigl^{\gbigf}_!(V_{\infty})$
induced by 
$\Abb^{\rd}_{\infty,\theta^u}$
and $\BB^{\rd}_{J,\theta^u}$.
\end{prop}
\pf
It follows from
Lemma \ref{lem;24.2.23.4},
Lemma \ref{lem;24.4.5.30}
and Lemma \ref{lem;24.4.5.31}.
\hfill\qed

\subsubsection{Another expression and the monodromy}
\label{subsection;24.2.22.1}
Let $c:\real\to\real$ be the map
defined by $c(\theta^u)=-\theta^u$.
Let $H^0(J,L_{J,>0})_{\real}$ denote the constant local system
on $\real$ induced by $H^0(J,L_{J,>0})$.
We fix $u(0)\in\cnum^{\ast}$ and $\theta^u_0\in\real$
such that $\arg(u(0))=\theta^u_0$.
Let $\gbigt(\nbigi,\theta^u)$ be as in \S\ref{subsection;25.2.9.1}.
We obtain the following exact sequence
(compare it with (\ref{eq;24.2.21.4})):
\begin{equation}
 0\lrarr
  c^{-1}(\nbigt_{\omega}(L))
  \lrarr
  \gbigq^{\infty}_!(L,\vecnbigftilde)_{\real}
  \lrarr
  \bigoplus_{J\in \gbigt(\nbigi,\theta_0^u)}
  H^0(J,L_{J,>0})_{\real}
  \lrarr 0.
\end{equation}
There exists the natural isomorphism
\begin{equation}
\label{eq;24.2.21.5}
  \gbigq^{\infty}_!(L,\vecnbigftilde)_{\real|\theta_0^u}
 \simeq
H^0(\real,\nbigt_{\omega}(L))
 \oplus
 \bigoplus_{J\in\gbigt(\nbigi,\theta_0^u)}
 H^0(J,L_{J,>0})
\end{equation}
under which 
the monodromy of 
$\gbigq^{\infty}_!(L,\vecnbigftilde)_{\real}$
is represented by the automorphism of
(\ref{eq;24.2.21.5})
given as follows:
\[
\Bigl(
v,\sum_{J}v_J
\Bigr)
\longmapsto
\Bigl(
 M_0^{-1}(v)
+\sum_JM_0^{-1}\circ\nbigq_{J_+}(v_J),\,
\sum_Jv_J
\Bigr).
\]
Here,
for any $v_J\in H^0(J,L_{J,>0})$,
we regard
$\nbigq_{J_+}(v)\in H^0(J,L_{J,0})\simeq H^0(\real,\nbigt_{\omega}(L))$.

\subsection{$2\pi\seisuu$-equivariant local system
$\gbigq^{\infty}_{\ast}(L,\vecnbigftilde)$}
\label{subsection;24.3.26.21}

We consider the vector space
\begin{equation}
\label{eq;24.3.21.40}
 \bigoplus_{\pm}
 \bigoplus_{J\in T(\nbigi)}
 H^0(\real,\nbigt_{\omega}(L))
 \oplus
 \bigoplus_{J\in T(\nbigi)}
 H^0(J,L_{J,<0}).
\end{equation}
An element $w\in H^0(\real,\nbigt_{\omega}(L))$
corresponding to the $(\kappa,J)$-component
($(\kappa,J)\in \{\pm\}\times T(\nbigi)$)
is denoted as
$\langle J_{\kappa},w\rangle^{\mg}$.
An element $v\in H^0(J,L_{J,<0})$
is denoted as $\langle J,v\rangle^{\mg}$.

Let $\gbigq^{\infty}_{\ast}(L,\vecnbigftilde)$ denote
the quotient space of (\ref{eq;24.3.21.40})
by the equivalence relation generated by
the following
(see Lemma \ref{lem;24.4.5.210} and Lemma \ref{lem;24.3.22.2}):
\index{vector space $\gbigq^{\infty}_{\ast}(L,\vecnbigftilde)$}
\begin{itemize}
 \item $\langle J+2\pi,v\rangle^{\mg}\sim\langle J,\Tbb^{\ast}(v)\rangle^{\mg}$
       for any $J\in T(\nbigi)$ and $v\in H^0(J,L_{J,<0})$
 \item $\langle J_+,w\rangle^{\mg}
       -\langle J_-,w\rangle^{\mg}
       \sim \langle J,\nbigp_J(w)\rangle^{\mg}$
       for any $J\in T(\nbigi)$
       and $w\in H^0(\real,\nbigt_{\omega}(L))$.
 \item $\langle J_{1-},w\rangle^{\mg}
       -\langle J_{2-},w\rangle^{\mg}
       \sim \sum_{J_2\leq J<J_1}
       \langle J,\nbigp_J(w)\rangle^{\mg}$
       for any $J_2<J_1$ in $T(\nbigi)$
       and $w\in H^0(\real,\nbigt_{\omega}(L))$.       
\end{itemize}
Let $\Tbb^{\ast}_{\gbigq^{\infty},\ast}$
denote the automorphism of
$\gbigq^{\infty}_{\ast}(L,\vecnbigftilde)$
induced by
\[
 \langle
 J_{\pm},w
 \rangle^{\mg}
 \longmapsto
 \langle
 (J+2\pi)_{\pm},M_0^{-1}(w)
 \rangle^{\mg},
 \quad\quad
 \langle J,v\rangle^{\mg}
 \longmapsto
 \langle J+2\pi,(\Tbb^{\ast})^{-1}(v)\rangle^{\mg}.
\]
(See Lemma \ref{lem;24.4.5.210} and Lemma \ref{lem;24.3.22.2}.)

Let $\gbigq^{\infty}_{\ast}(L,\vecnbigftilde)_{\real}$
denote the local system on $\real$
induced by $\gbigq^{\infty}_{\ast}(L,\vecnbigftilde)$.
\index{local system $\gbigq^{\infty}_{\ast}(L,\vecnbigftilde)_{\real}$}
We naturally identify
$H^0(\real,\gbigq^{\infty}_{\ast}(L,\vecnbigftilde)_{\real})$
with $\gbigq^{\infty}_{\ast}(L,\vecnbigftilde)$.
There exists the $2\pi\seisuu$-action on
$\gbigq^{\infty}_{\ast}(L,\vecnbigftilde)$
such that
the pull back
$\Tbb^{\ast}:H^0(\real,\gbigq^{\infty}_{\ast}(L,\vecnbigftilde)_{\real})
\simeq
H^0(\real,\gbigq^{\infty}_{\ast}(L,\vecnbigftilde)_{\real})$
equals $\Tbb^{\ast}_{\gbigq^{\infty},\ast}$.

\begin{prop}
There exists the isomorphism of
$2\pi\seisuu$-equivariant local systems
 $\gbigq^{\infty}_{\ast}(L,\vecnbigftilde)_{\real}
 \simeq
\gbigl^{\gbigf}_{\ast}(V_{\infty})$
induced by 
$\Abb^{\mg,J_{\pm}}_{\infty,\theta^u}$
and 
$\Abb^{\mg}_{J,\theta^u}$.
\end{prop}
\pf
It follows from Lemma \ref{lem;24.4.5.210}
and Lemma \ref{lem;24.3.22.2}.
\hfill\qed

\subsubsection{Another expression and the monodromy}

Let $H^0(J,L_{J,<0})_{\real}$
denote the constant $2\pi\seisuu$-equivariant local system
on $\real$ induced by $H^0(J,L_{J,<0})$.
We fix $u(0)\in\cnum^{\ast}$
and $\theta^u_0\in\real$
such that
$\arg(u(0))=\theta^u_0$.
Let $\gbigt(\nbigi,\theta^u_0)$ be as in \S\ref{subsection;25.2.9.1}.
We obtain the following exact sequence
(compare it with (\ref{eq;24.2.21.10})):
\[
 0\lrarr
 \bigoplus_{J\in \gbigt(\nbigi,\theta^u_0)}
 H^0(J,L_{J,<0})_{\real}
 \lrarr
 \gbigq^{\infty}_{\ast}(L,\vecnbigftilde)
  \lrarr
 c^{-1}(\nbigt_{\omega}(L))
 \lrarr 0.
\]
Let $J_1\in T(\nbigi)$
be determined by
$\theta^u-\pi/2<\vartheta^J_{\ell}$
and
$\openopen{\theta^u-\pi/2}{\vartheta^J_{\ell}}
\cap S_0(\nbigi)=\emptyset$.
By considering
$\langle J_{1-},w\rangle^{\mg}$
for $w\in H^0(\real,\nbigt_{\omega}(L))$,
we obtain
the isomorphism
\begin{equation}
\label{eq;24.2.21.11}
  \gbigq^{\infty}_{\ast}(L,\vecnbigftilde)_{\real|\theta_0^u}
  \simeq
  H^0(\real,\nbigt_{\omega}(L))
  \oplus
  \bigoplus_{J\in \gbigt(\nbigi,\theta^u_0)}
  H^0(J,L_{J,<0}),
\end{equation}
under which 
the monodromy of $\gbigq^{\infty}_{\ast}(L,\vecnbigftilde)_{\real}$
is represented by
\[
\Bigl(
w,\sum_Jv_J
\Bigr)
\longmapsto
\Bigl(
M_0^{-1}(w),
 \sum_J(v_J+\nbigp_{J}(M_0^{-1}w))
 \Bigr).
\]

\subsection{Morphisms}
\label{subsection;24.4.4.11}

Let $F_{\gbigq^{\infty}}:\gbigq^{\infty}_!(L,\vecnbigftilde)
\to\gbigq^{\infty}_{\ast}(L,\vecnbigftilde)$
be the morphism obtained as follows
(see Lemma \ref{lem;24.4.5.40} and Lemma \ref{lem;24.4.5.41}):
\index{morphism $F_{\gbigq^{\infty}}$}
\begin{itemize}
 \item For any $J\in T(\nbigi)$ and $v\in H^0(J,L_{J,>0})$,
\begin{multline}
 \langle J,v\rangle
 \longmapsto
 \langle J_-,\nbigq_{J_+}(v)\rangle^{\mg}
-\langle J,\nbigr^{J_+}_{J_-}(v)\rangle^{\mg}
\\
 +\!\!\!\!\sum_{J<J'\leq J+\omega^{-1}\pi}
 \!\!\!\!
 \langle J',\nbigr^J_{J'}(v)\rangle^{\mg}
-\!\!\!\!\sum_{J-\omega^{-1}\pi\leq J'<J}
 \!\!\!\!
 \langle J',\nbigr^{J}_{J'}(v)\rangle^{\mg}.
\end{multline}
 \item For any $w\in H^0(\real,\nbigt_{\omega}(L))$,
\[
 w\longmapsto
 \langle J_{1-},w-M_0(w)\rangle^{\mg}
 +\sum_{J_1\leq J'<J_1+2\pi}
 \langle J',\nbigp_{J'}(w)\rangle^{\mg}.
\]
The right hand side is independent of $J_1$.
\end{itemize}
It induces the morphism of $2\pi\seisuu$-equivariant
local systems
$F_{\gbigq^{\infty}}:\gbigq^{\infty}_{!}(L,\vecnbigftilde)_{\real}
\to \gbigq^{\infty}_{\ast}(L,\vecnbigftilde)_{\real}$.

The morphism $c^{-1}(\nbigt_{\omega}(L))
\to\gbigq^{\infty}_!(L,\vecnbigftilde)_{\real}$
is induced by the inclusion of
$H^0(\real,\nbigt_{\omega}(L))$
into the space (\ref{eq;24.3.21.30}).
The morphism 
$\gbigq^{\infty}_{\ast}(L,\vecnbigftilde)_{\real}
\to
c^{-1}(\nbigt_{\omega}(L))$
is induced by the projection of (\ref{eq;24.3.21.40})
onto $H^0(\real,\nbigt_{\omega}(L))$.
\begin{prop}
We obtain the following commutative diagram:
\[
 \begin{CD}
  c^{-1}(\nbigt_{\omega}(L))@>>>
  \gbigq^{\infty}_!(L,\vecnbigftilde)_{\real}
  @>>>
  \gbigq^{\infty}_{\ast}(L,\vecnbigftilde)_{\real}
  @>>>
  c^{-1}(\nbigt_{\omega}(L))
  \\
  @V{\simeq}VV @V{\simeq}VV @V{\simeq}VV @V{\simeq}VV \\
  \gbigl^{\gbigf}_!(\nbigt_{\omega}(V_{\infty}))
  @>>>
  \gbigl^{\gbigf}_!(V_{\infty})
  @>>>
  \gbigl^{\gbigf}_{\ast}(V_{\infty})
  @>>>
  \gbigl^{\gbigf}_{\ast}(\nbigt_{\omega}(V_{\infty})).
 \end{CD}
\]
\end{prop}
\pf
We obtain the commutativity of the middle square
from Lemma \ref{lem;24.4.5.40}
and Lemma \ref{lem;24.4.5.41}.
\hfill\qed

\subsection{Stokes filtrations of $\gbigq^{\infty}_!(L,\vecnbigftilde)_{\real}$}
\label{subsection;24.4.4.40}

Let $\vecJ\in T(\nbigi^{\circ})$.
We define the map
$\vecB_{\vecJ}:
H^0(\nu_0^+(\vecJ),L_{\nu_0^+(\vecJ),>0})
\to \gbigq^{\infty}_!(L,\vecnbigftilde)$
by \index{map $\vecB_{\vecJ}$}
\[
 \vecB_{\vecJ}(v)
 =\langle \nu_0^+(\vecJ),v\rangle.
\]
We define the map
$\vecA_{\vecJ_{\pm}}:
H^0(\nu_0^-(\vecJ),L_{\nu_0^-(\vecJ),<0})
\to \gbigq^{\infty}_!(L,\vecnbigftilde)$
by \index{maps $\vecA_{\vecJ_{\pm}}$}
\[
 \vecA_{\vecJ_+}(v)
 =\sum_{\nu_0^-(\vecJ)-\pi<J'\leq \nu_0^-(\vecJ)-\omega^{-1}\pi}
 \langle J',R_{J'}(v)\rangle,
\]
\[
 \vecA_{\vecJ_-}(v)
 =-\sum_{\nu_0^-(\vecJ)+\omega^{-1}\pi\leq J'<\nu_0^-(\vecJ)+\pi}
 \langle
 J'-2\pi,
 \Tbb^{\ast}(R_{J'}(v))
 \rangle.
\]
(See Lemma \ref{lem;24.4.5.30} and Proposition \ref{prop;24.3.17.1}.)
We introduce the maps
$\vecA^{\vecJ_{\pm}}_{\infty}:
H^0(\real,L)\to\gbigq^{\infty}_!(L,\vecnbigftilde)$
as follows.
\index{maps $\vecA^{\vecJ_{\pm}}_{\infty}$}
For any $v\in H^0(\real,L)$,
we have the decomposition
\[
 v=u_{(\nu_0^+(\vecJ)+2\pi)_+,0}
 +\sum_{J'\in \gbigk(\nu_0^+(\vecJ)_+)}
 u_{J'+2\pi},
\]
where $u_{(\nu_0^+(\vecJ)+2\pi)_+,0}$ is
a section of $L'_{(\nu_0^+(\vecJ)+2\pi)_+,0}$,
and $u_{J'+2\pi}$ are sections of $L'_{J'+2\pi,<0}$.
(See \S\ref{subsection;25.2.9.10}.)
We obtain a section of $\nbigt_{\omega}(L)$ on $\real$
from $u_{(\nu_0^+(\vecJ)+2\pi)_+,0}$,
which is also denoted by $u_{(\nu_0^+(\vecJ)+2\pi)_+,0}$.
We set
\[
 \vecA^{\vecJ_+}_{\infty}(v)
=u_{(\nu_0^+(\vecJ)+2\pi)_+,0}
+\sum_{\nu_0^+(\vecJ)<J'\leq \nu_0^+(\vecJ)+2\pi}
\langle
J',R_{J'}(v)
\rangle.
\]
(See Lemma \ref{lem;25.2.9.11}.)
Similarly,
we have the decomposition
\[
 v=w_{(\nu_0^+(\vecJ)+2\pi)_-,0}
 +\sum_{J'\in \gbigk(\nu_0^+(\vecJ)_-)}
 w_{J'+2\pi},
\]
where $w_{(\nu_0^+(\vecJ)+2\pi)_-,0}$ is
a section of $L'_{(\nu_0^+(\vecJ)+2\pi)_-,0}$,
and $w_{J'+2\pi}$ are sections of $L'_{J'+2\pi,<0}$.
We set
\[
 \vecA^{\vecJ_-}_{\infty}(v)
=w_{(\nu_0^+(\vecJ)+2\pi)_-,0}
+\sum_{\nu_0^+(\vecJ)\leq J'< \nu_0^+(\vecJ)+2\pi}
\langle
J',R_{J'}(v)
\rangle.
\]
Thus, we obtain the maps
$\vecA^{\vecJ_{\pm}}_{\infty}:
H^0(\real,L)\to\gbigq^{\infty}_!(L,\vecnbigftilde)$.
By Theorem \ref{thm;24.3.16.20},
we obtain the following.
\begin{prop}
\label{prop;24.4.6.210}
Let $\theta^u\in\real$.
Choose $\vecJ_1\in\gbigm_{-}(\nbigi^{\circ},\theta^u)$.
Then,
$\vecA_{\vecJ_{+}}$, $\vecB_{\vecJ}$
$(\vecJ\in \gbigm_-(\nbigi,\theta^u))$
and $\vecA^{\vecJ_{1+}}_{\infty}$
induce the isomorphism of the vector spaces:
 \begin{multline}
\label{eq;24.4.6.200}
 \bigoplus_{\vecJ\in\gbigm_-(\nbigi^{\circ},\theta^u)}
\Bigl(
 H^0(\nu_0^-(\vecJ),L_{\nu_0^-(\vecJ),<0})
\oplus 
 H^0(\nu_0^+(\vecJ),L_{\nu_0^+(\vecJ),>0})
 \Bigr)
 \oplus
 H^0(\real,L)
 \\
\simeq
\gbigq^{\infty}_!(L,\vecnbigftilde)
\simeq
\gbigl^{\gbigf}_!(V_{\infty})_{|\theta^u}.
\end{multline}
Moreover,
if we consider the filtrations $\nbigf^{\prime\theta^u}$
on the spaces
$H^0(\nu_0^-(\vecJ),L_{\nu_0^-(\vecJ),<0})$
and
$H^0(\nu_0^+(\vecJ),L_{\nu_0^+(\vecJ),>0})$
defined in {\rm\S\ref{subsection;24.4.4.10}},
the trivial filtration on $H^0(\real,L)$
indexed by $0$,
and the Stokes filtration $\nbigf^{\theta^u}$
on $\gbigl^{\gbigf}_!(V_{\infty})_{|\theta^u}$,
then {\rm(\ref{eq;24.4.6.200})} induces
an isomorphism of the filtered spaces.

We also obtain a similar isomorphism
by choosing 
$\vecJ_1\in\gbigm_+(\nbigi^{\circ},\theta^u)$
and by using 
$\vecA_{\vecJ_{-}}$, $\vecB_{\vecJ}$
$(\vecJ\in \gbigm_+(\nbigi,\theta^u))$
and $\vecA^{\vecJ_{1-}}_{\infty}$.
\hfill\qed
\end{prop}

We also obtain the following
by Theorem \ref{thm;24.3.16.20}.

\begin{prop}
\label{prop;24.4.4.50}
Under the isomorphism
$\gbigq^{\infty}_!(L,\vecnbigftilde)\simeq
H^0(\real,\gbigl^{\gbigf}_!(V_{\infty}))$,
we have
\[
 \Image \vecB_{\vecJ}
 =H^0(\vecJ,\gbigl^{\gbigf}_{!}(V_{\infty})_{\vecJ,<0}),
 \quad
 \Image \vecA_{\vecJ_{\pm}}
 =H^0(\vecJ_{\mp},\gbigl^{\gbigf}_{!}(V_{\infty})_{\vecJ_{\mp},>0}),
\]
\[
 \Image \vecA^{\vecJ_{\pm}}_{\infty}
=H^0(\vecJ_{\mp},\gbigl^{\gbigf}_!(V_{\infty})_{\vecJ_{\mp},0}).
\]

\hfill\qed  
\end{prop}

\subsection{Stokes filtrations of $\gbigq^{\infty}_{\ast}(L,\vecnbigftilde)$}
\label{subsection;24.4.4.41}

For $\vecJ\in T(\nbigi^{\circ})$,
we obtain the maps
\index{maps $\vecB^{\mg}_{\vecJ}$}
\index{maps $\vecA^{\mg}_{\vecJ_{\pm}}$}
\[
\vecB^{\mg}_{\vecJ}:
H^0(\nu_0^+(\vecJ),L_{\nu_0^+(\vecJ),>0})
\to
\gbigq^{\infty}_{\ast}(L,\vecnbigftilde),
\]
\[
\vecA^{\mg}_{\vecJ_{\pm}}:
H^0(\nu_0^+(\vecJ),L_{\nu_0^+(\vecJ),>0})
\to
\gbigq^{\infty}_{\ast}(L,\vecnbigftilde)
\]
as the composition of
$\vecB_{\vecJ}$ and
$\vecA_{\vecJ_{\pm}}$
with the map
$F_{\gbigq^{\infty}}:\gbigq^{\infty}_!(L,\vecnbigftilde)
\to \gbigq^{\infty}_{\ast}(L,\vecnbigftilde)$.
We introduce the maps
$\vecA^{\mg,\vecJ_{\pm}}_{\infty}:
H^0(\real,L)\to\gbigq^{\infty}_{\ast}(L,\vecnbigftilde)$
as follows.
\index{maps $\vecA^{\mg,\vecJ_{\pm}}_{\infty}$}
For $v\in H^0(\real,L)$,
as in \S\ref{subsection;25.2.9.10},
we have the decomposition
\[
 v=u_{\nu_0^+(\vecJ)_+,0}
 +\sum_{J'\in\gbigk(\nu_0^+(\vecJ)_+)}
 u_{J'},
\]
where $u_{\nu_0^+(\vecJ)_+,0}$
is a section of $L'_{\nu_0^+(\vecJ)_+,0}$
and $u_{J'}$ are sections of $L'_{J',<0}$.
We set
\[
 \vecA^{\mg,\vecJ_+}_{\infty}(v)
=\langle
  \nu_0^+(\vecJ)_+,
  u_{\nu_0^+(\vecJ)_+,0}
  \rangle^{\mg}
  -\sum_{J'\in\gbigk(\nu_0^+(\vecJ)_+)}
  \langle
  J',u_{J'}
  \rangle^{\mg}.
\]
(See Lemma \ref{lem;25.2.9.20}.)
Similarly,
there exists the decomposition
\[
 v=w_{\nu_0^+(\vecJ)_-,0}
 +\sum_{J'\in\gbigk(\nu_0^+(\vecJ)_-)}
 w_{J'},
\]
where $w_{\nu_0^+(\vecJ)_-,0}$
is a section of $L'_{\nu_0^+(\vecJ)_-,0}$
and $w_{J'}$ are sections of $L'_{J',<0}$.
We set
\[
 \vecA^{\mg,\vecJ_-}_{\infty}(v)
=\langle
  \nu_0^+(\vecJ)_-,
  w_{\nu_0^+(\vecJ)_-,0}
  \rangle^{\mg}
  -\sum_{J'\in\gbigk(\nu_0^+(\vecJ)_-)}
  \langle
  J',w_{J'}
  \rangle^{\mg}.
\]
Thus, we obtain the maps
$\vecA^{\mg,\vecJ_{\pm}}_{\infty}:
H^0(\real,L)\to
\gbigq^{\infty}_{\ast}(L,\vecnbigftilde)$.
We obtain the following
by Theorem \ref{thm;24.3.16.20}.

\begin{prop}
\label{prop;24.4.6.211}
Let $\theta^u\in\real$.
Choose $\vecJ_1\in\gbigm_{-}(\nbigi^{\circ},\theta^u)$.
Then,
$\vecA^{\mg}_{\vecJ_{+}}$, $\vecB^{\mg}_{\vecJ}$
$(\vecJ\in \gbigm_-(\nbigi,\theta^u))$
and $\vecA^{\mg,\vecJ_{1+}}_{\infty}$
induce the isomorphism of the vector spaces:
 \begin{multline}
\label{eq;24.4.6.201}
 \bigoplus_{\vecJ\in\gbigm_-(\nbigi^{\circ},\theta^u)}
\Bigl(
 H^0(\nu_0^-(\vecJ),L_{\nu_0^-(\vecJ),<0})
\oplus 
 H^0(\nu_0^+(\vecJ),L_{\nu_0^+(\vecJ),>0})
 \Bigr)
 \oplus
 H^0(\real,L)
 \\
\simeq
\gbigq^{\infty}_{\ast}(L,\vecnbigftilde)\simeq
  \gbigl^{\gbigf}_{\ast}(V_{\infty})_{|\theta^u}.
 \end{multline}
If we consider the filtrations $\nbigf^{\prime\theta^u}$
on the spaces
$H^0(\nu_0^-(\vecJ),L_{\nu_0^-(\vecJ),<0})$
and
$H^0(\nu_0^+(\vecJ),L_{\nu_0^+(\vecJ),>0})$
defined in {\rm\S\ref{subsection;24.4.4.10}},
the trivial filtration on $H^0(\real,L)$
indexed by $0$,
and the Stokes filtration $\nbigf^{\theta^u}$
on $\gbigl^{\gbigf}_{\ast}(V_{\infty})_{|\theta^u}$,
{\rm(\ref{eq;24.4.6.201})} induces an isomorphism of
filtered vector spaces.
 
We obtain a similar isomorphisms
by choosing $\vecJ_1\in\gbigm_+(\nbigi^{\circ},\theta^u)$
and by using
$\vecA^{\mg}_{\vecJ_{-}}$, $\vecB^{\mg}_{\vecJ}$
$(\vecJ\in \gbigm_+(\nbigi,\theta^u))$
and $\vecA^{\mg,\vecJ_{1-}}_{\infty}$.
\hfill\qed 
\end{prop}

We also obtain the following
by Theorem \ref{thm;24.3.16.20}.
\begin{prop}
\label{prop;24.4.4.51}
Under the isomorphism
$\gbigq^{\infty}_{\ast}(L,\vecnbigftilde)\simeq
H^0(\real,\gbigl^{\gbigf}_{\ast}(V_{\infty}))$,
we have
\[
 \Image \vecB^{\mg}_{\vecJ}
 =H^0(\vecJ,\gbigl^{\gbigf}_{\ast}(V_{\infty})_{\vecJ,<0}),
 \quad
 \Image \vecA^{\mg}_{\vecJ_{\pm}}
 =H^0(\vecJ_{\mp},\gbigl^{\gbigf}_{\ast}(V_{\infty})_{\vecJ_{\mp},>0}),
\]
\[
 \Image \vecA^{\mg,\vecJ_{\pm}}_{\infty}
=H^0(\vecJ_{\mp},\gbigl^{\gbigf}_{\ast}(V_{\infty})_{\vecJ_{\mp},0}).
\]

\hfill\qed  
\end{prop}

\subsection{Isomorphisms}

For any $\theta^u\in\real$,
we define the filtrations $\nbigf^{\theta^u}$
on $\gbigq^{\infty}_{\star}(L,\vecnbigftilde)
=\gbigq^{\infty}_{\star}(L,\vecnbigftilde)_{\real|\theta^u}$
$(\star=!,\ast)$
indexed by
$(\nbigitilde^{\circ},\leq_{\theta^u})$
by using the isomorphisms
(\ref{eq;24.4.6.200})
and (\ref{eq;24.4.6.201}),
and the filtrations
$\nbigf^{\prime\theta^u}$
on 
$H^0(\nu_0^-(\vecJ),L_{\nu_0^-(\vecJ),<0})$
and
$H^0(\nu_0^+(\vecJ),L_{\nu_0^+(\vecJ),>0})$
defined in {\rm\S\ref{subsection;24.4.4.10}},
and the trivial filtration on $H^0(\real,\nbigt_{\omega}(L))$
indexed by $0$.
It is independent of the choice of $\vecJ_1$.
We obtain the $2\pi\seisuu$-equivariant family of
filtrations $\vecnbigf=(\nbigf^{\theta^u}\,|\,\theta^u\in\real)$
of $\gbigq^{\infty}_{\star}(L,\vecnbigftilde)_{\real}$.
By Proposition \ref{prop;24.4.6.210} and Proposition \ref{prop;24.4.6.211},
we obtain the following.
\begin{thm}
$(\gbigq^{\infty}_{\star}(L,\vecnbigftilde),\vecnbigf)$
are local systems with Stokes structure indexed by 
$\nbigitilde^{\circ}$.
Moreover, there exists the following commutative diagram 
in $\Loc^{\St}(\nbigitilde^{\circ})$:
\[
\begin{CD}
(\gbigq^{\infty}_!(L,\vecnbigftilde)_{\real},\vecnbigf)
 @>{F}>>
 (\gbigq^{\infty}_{\ast}(L,\vecnbigftilde)_{\real},\vecnbigf)
 \\
 @V{\simeq}VV @V{\simeq}VV \\ 
 (\gbigl^{\gbigf}_!(V_{\infty}),\vecnbigftilde)
 @>>>
 (\gbigl^{\gbigf}_{\ast}(V_{\infty}),\vecnbigftilde),
\end{CD}
\]
where the lower horizontal arrow is induced by
$V_{\infty}(!0)\to V_{\infty}$.
\hfill\qed
\end{thm}

\begin{df}
We set
$\gbigf^{(\infty,\infty)}_{+,\star}\bigl(
L,\vecnbigftilde
\bigr)
 :=(\gbigq^{\infty}_{\star}(L,\vecnbigftilde),\vecnbigf)$,
called the local Fourier transform of
$(L,\vecnbigftilde)$.
\index{local system with Stokes structure
$\gbigf^{(\infty,\infty)}_{+,\star}\bigl(
L,\vecnbigftilde
\bigr)$}
\hfill\qed
\end{df}

\section{Stokes shells}
\label{subsection;20.11.14.20}

Let us describe
$\Shsf(\gbigf^{(\infty,\infty)}_{+,\star}(L,\vecnbigftilde))$.
We set
$(\nbigk_{\bullet},\vecnbigf,\vecnbigr)=
\Shsf(L,\vecnbigftilde)
\in\Shcat(\nbigitilde)$.
We set
$(\vecK,\vecnbigf,\vecPhi,\vecPsi):=
\gbigd(\nbigk_{\bullet},\vecnbigf)$.
We use the notation
$\nbigp_J=\nbigr^{0,J_-}_{\lambda_-(J),J_+}$, 
$\nbigq_J=\nbigr^{\lambda_+(J),J_-}_{0,J_+}$,
$\nbigr^{J_-}_{J_+}=\nbigr^{\lambda_+(J),J_-}_{\lambda_-(J),J_+}$
and 
$\nbigr^{J_+}_{J_-}=\nbigr^{\lambda_+(J),J_+}_{\lambda_-(J),J_-}$
for $\Shsf(L,\vecnbigftilde)$.
For any $J_1,J_2\in T(\nbigi)$,
let $\Phi^{J_1,J_2}$ denote the isomorphism
$H^0(J_2,L)\simeq H^0(J_1,L)$
induced by the parallel transport of $L$.
The restriction of $\Phi^{J_1,J_2}$
to a subspace of
$H^0(J_2,L)$ is also denoted by
$\Phi^{J_1,J_2}$.
Let $M$ denote the automorphism of $L$
obtained as 
$L\stackrel{c_1}{\simeq}
 \Tbb^{-1}(L)
 \stackrel{c_2}{\simeq}L$,
where $c_1$ is induced by the parallel transport,
and $c_2$ is induced by the $2\pi\seisuu$-equivariance.
\index{automorphism $M$}

\subsection{Stokes graded local systems}

We use the notation in \S\ref{subsection;20.10.9.30}.
For $\vecJ\in T(\nbigi^{\circ})$,
we obtain the intervals
$\nu^{\pm}_m(\vecJ)\in T(\nbigi)$
and the isomorphisms
$\kappa_{m,\vecJ}^{\pm}:
 \vecJbar\simeq
 \nu_m^{\pm}(\vecJbar)$.
By Proposition \ref{prop;18.5.6.1},
we obtain the following local system
with Stokes structure indexed by
$\nbigitilde^{\circ}_{\vecJ,<0}$ on $\vecJbar$:
\[
 \bigl(
 \nbigk^{\circ}_{\lambda_-(\vecJ),\vecJ},
 \vecnbigf^{\circ}
 \bigr)
:=
 (\kappa_{0,\vecJ}^+)^{-1}\bigl(
 \nbigk_{\lambda_+(\nu_0^+(\vecJ))},
 \vecnbigf
 \bigr)_{|\nu_0^+(\vecJbar)}.
\]
We also obtain the following local system
with Stokes structure indexed by
$\nbigitilde^{\circ}_{\vecJ,>0}$ on $\vecJbar$:
\[
 \bigl(
 \nbigk^{\circ}_{\lambda_+(\vecJ),\vecJ},
 \vecnbigf^{\circ}
 \bigr)
:=
 (\kappa_{0,\vecJ}^-)^{-1}
 \bigl(
 \nbigk_{\lambda_-(\nu_0^-(\vecJ))},
 \vecnbigf
 \bigr)_{|\nu_0^-(\vecJbar)}.
\]
We obtain the following local system on $\vecJbar$:
\[
 \nbigk^{\circ}_{0,\vecJ}:=
 (\kappa_{0,\vecJ}^+)^{-1}\bigl(
 L_{|\nu_0^+(\vecJbar)}
 \bigr).
\]
The spaces of the global sections of
$\nbigk^{\circ}_{\lambda,\vecJ}$
are denoted by 
$K^{\circ}_{\lambda,\vecJ}$.
\index{vector space $K^{\circ}_{\lambda,\vecJ}$}
By the construction and
the relation
$\kappa^{\pm}_{0,\vecJ}\circ\Tbb
=\Tbb^{-1}\circ\kappa^{\pm}_{0,\Tbb^{-1}(\vecJ)}$,
there exist natural isomorphisms
$\Tbb^{-1}(\nbigk^{\circ}_{\lambda,\vecJ},\vecnbigf)
\simeq
 (\nbigk^{\circ}_{\Tbb^{\ast}(\lambda),\Tbb^{-1}(\vecJ)},\vecnbigf)$,
 which induces
 $\Psi^{\circ}_{\lambda,\vecJ}:
 K^{\circ}_{\lambda,\vecJ}
 \simeq
 K^{\circ}_{\Tbb^{\ast}(\lambda),\Tbb^{-1}(\vecJ)}$.
 We have the natural identifications:
\[
 K^{\circ}_{\lambda_-(\vecJ),\vecJ}
=K_{\lambda_+(\nu_0^+(\vecJ)),\nu_0^+(\vecJ)},
\quad\quad
 K^{\circ}_{\lambda_+(\vecJ),\vecJ}
=K_{\lambda_-(\nu_0^-(\vecJ)),\nu_0^-(\vecJ)}.
\]
By the construction, we have
$K^{\circ}_{0,\vecJ}
 =H^0\bigl(
 \nu_0^{+}(\vecJbar),
 L
 \bigr)$.

Because
$\nu_0^-(\vecJ+(1-\omega^{-1})\pi)=\nu_0^+(\vecJ)+\omega^{-1}\pi$,
we obtain the following isomorphism:
\[
 \bigl(
 \Phi^{\circ}
 \bigr)^{\vecJ+(1-\omega^{-1})\pi,\vecJ}
 _{\lambda_-(\vecJ)}
:=\Phi_{\lambda_+(\nu_0^+(\vecJ))}
 ^{\nu_0^+(\vecJ)+\omega^{-1}\pi,\nu_0^+(\vecJ)}:
 K^{\circ}_{\lambda_-(\vecJ),\vecJ}
\simeq
 K^{\circ}_{\lambda_-(\vecJ),\vecJ+(1-\omega^{-1})\pi}.
\]
Because
$\nu_{-1}^+(\vecJ+(1-\omega^{-1})\pi)=\nu_0^-(\vecJ)+\omega^{-1}\pi$,
we obtain the following isomorphism:
\[
 \bigl(
 \Phi^{\circ}
 \bigr)^{\vecJ+(1-\omega^{-1})\pi,\vecJ}
 _{\lambda_+(\vecJ)}
:=
 -\Psi\circ
 \Phi_{\lambda_-(\nu_0^-(\vecJ))}
 ^{\nu_0^-(\vecJ)+\omega^{-1}\pi,\nu_0^-(\vecJ)}:
 K^{\circ}_{\lambda_+(\vecJ),\vecJ}
\simeq
 K^{\circ}_{\lambda_+(\vecJ),\vecJ+(1-\omega^{-1})\pi}.
\]
For $\vecJ_1\vdash\vecJ_2$ in $T(\nbigi^{\circ})$,
we obtain the following natural isomorphism
induced by the parallel transport of $L$:
\index{maps $(\Phi^{\circ})^{\vecJ_1,\vecJ_2}_{\lambda}$}
\[
 \bigl(
 \Phi^{\circ}
 \bigr)^{\vecJ_2,\vecJ_1}_0:
 K^{\circ}_{0,\vecJ_1}
\simeq
 K^{\circ}_{0,\vecJ_2}.
\]
By gluing
$(\nbigk^{\circ}_{\lambda,\vecJ},\vecnbigf^{\circ})$ 
via the tuple of the isomorphisms
$\vecPhi^{\circ}$,
we obtain a Stokes graded local system
$(\nbigk^{\circ}_{\bullet},\vecnbigf^{\circ})$
over $(\nbigitilde^{\circ},[\nbigi^{\circ}])$.
By the construction,
it is naturally $2\pi\seisuu$-equivariant.
\index{local system with Stokes structure
$(\nbigk^{\circ}_{\bullet},\vecnbigf^{\circ})$}

For $\star=!,\ast$,
we set
$(\nbigk_{\star\,\bullet},\vecnbigf^{\circ}):=
 (\nbigk_{\bullet},\vecnbigf^{\circ})$.
We naturally have
$\gbigd(\nbigk_{\star\bullet},\vecnbigf^{\circ})
=(\vecK^{\circ},\vecnbigf^{\circ},\vecPhi^{\circ},\vecPsi^{\circ})$.
\index{local system with Stokes structure
$(\nbigk^{\circ}_{\bullet\ast},\vecnbigf^{\circ})$}
\index{local system with Stokes structure
\mbox{$(\nbigk^{\circ}_{\bullet!},\vecnbigf^{\circ})$}}

\subsection{Morphisms $\vecnbigp^{\circ}_{\star}$
and $\vecnbigq^{\circ}_{\star}$ 
($\star=!,\ast$)}

For any $J\in T(\nbigi)$,
there exists the morphism
$R_J:H^0(\real,L)
\lrarr
H^0(J,L_{J,>0})=K_{J,>0}$.
(See \S\ref{subsection;20.10.9.20}.)
For $\vecJ\in T(\nbigi^{\circ})$,
we set
\index{maps \mbox{$(\nbigp^{\circ}_!)_{\vecJ}$},
\mbox{$(\nbigp^{\circ}_{\ast})_{\vecJ}$},
\mbox{$(\nbigq^{\circ}_!)_{\vecJ}$},
\mbox{$(\nbigq^{\circ}_{\ast})_{\vecJ}$}}
\[
 (\nbigp^{\circ}_!)_{\vecJ}:=
 -R_{\nu_0^+(\vecJ)}\circ(\id-M),
\quad
 (\nbigp^{\circ}_{\ast})_{\vecJ}:=
 -R_{\nu_0^+(\vecJ)},
\]
\[
 (\nbigq^{\circ}_!)_{\vecJ}:=
 \Phi^{\nu_0^+(\vecJ),\nu_0^-(\vecJ)},
\quad
 (\nbigq^{\circ}_{\ast})_{\vecJ}:=
 (\id-M)\circ
 \Phi^{\nu_0^+(\vecJ),\nu_0^-(\vecJ)}.
\]

\subsection{Morphisms $\vecnbigr^{\circ}$}
\index{tuple of morphisms $\vecnbigr^{\circ}$}
\label{subsection;25.2.10.1}

For $\vecJ\in T(\nbigi^{\circ})$,
we set
\[
 (\nbigr^{\circ})^{\vecJ_-}_{\vecJ_+}
=-R_{\nu_0^+(\vecJ)}\circ
 \Phi^{\nu_0^+(\vecJ),\nu_0^-(\vecJ)}.
\]
We remark the following
\begin{itemize}
\item
 $\vecJ-(1-\omega^{-1})\pi<\vecJ'<\vecJ$
if and only if
 $\nu_0^{+}(\vecJ)
<\nu_0^{+}(\vecJ')
<\nu_0^{-}(\vecJ)-\omega^{-1}\pi$.
\item
$\vecJ<\vecJ'<\vecJ+(1-\omega^{-1})\pi$
if and only if
$\nu_0^-(\vecJ)+\omega^{-1}\pi
<\nu_{-1}^+(\vecJ')
<\nu_{-1}^+(\vecJ)$.
\end{itemize}
For $\vecJ-(1-\omega^{-1})\pi<\vecJ'<\vecJ+(1-\omega^{-1})\pi$
with $\vecJ'\neq \vecJ$,
we define
\[
 (\nbigr^{\circ})^{\vecJ}_{\vecJ'}:=
 \left\{
 \begin{array}{ll}
 R_{\nu_0^+(\vecJ')}\circ
 \Phi^{\nu_0^+(\vecJ'),\nu_0^-(\vecJ)}
 & (\vecJ-(1-\omega^{-1})\pi<\vecJ'<\vecJ)
 \\
 -\Psi\circ R_{\nu_{-1}^+(\vecJ')}\circ
 \Phi^{\nu_{-1}^+(\vecJ'),\nu_0^-(\vecJ)}
 & (\vecJ<\vecJ'<\vecJ+(1-\omega^{-1})\pi).
 \end{array}
 \right.
\]

\subsection{Description}

Let $\gbigf^{(\infty,\infty)}_{+,\star}(\Shsf(L,\vecnbigftilde))$
$(\star=!,\ast)$
be the shells consisting of
$(\nbigk^{\circ}_{\star \bullet},\vecnbigf^{\circ})$
and the tuple of the morphisms
$\bigl(
\vecnbigp^{\circ}_{\star},
\vecnbigq^{\circ}_{\star},
\vecnbigr \bigr)$.
\index{shells $\gbigf^{(\infty,\infty)}_{+,\star}(\Shsf(L,\vecnbigftilde))$}
Let 
$\gbigf^{(\infty,\infty)}_{+,!}(\Shsf(L,\vecnbigftilde))
\to
\gbigf^{(\infty,\infty)}_{+,\ast}(\Shsf(L,\vecnbigftilde))$
be the morphism induced by
the identity maps
$\nbigk^{\circ}_{!\lambda}
=\nbigk^{\circ}_{\ast\lambda}$ 
$(\lambda\neq 0)$
and
$\id-M:
 \nbigk^{\circ}_{!0}\lrarr
 \nbigk^{\circ}_{\ast 0}$.
We obtain the following
as the translation of the results in
\S\ref{subsection;24.3.29.100}--\S\ref{subsection;24.3.29.101}.
\begin{prop}
There exists the following commutative diagram
in $\Shcat(\nbigitilde^{\circ})$:
\[
 \begin{CD}
  \gbigf^{(\infty,\infty)}_{+,!}(\Shsf(L,\vecnbigftilde))
  @>{F}>>
  \gbigf^{(\infty,\infty)}_{+,\ast}(\Shsf(L,\vecnbigftilde))
  \\
  @V{\simeq}VV @V{\simeq}VV \\  
  \Shsf(\gbigf^{(\infty,\infty)}_{+,!}(L,\vecnbigftilde))
  @>>>
  \Shsf(\gbigf^{(\infty,\infty)}_{+,\ast}(L,\vecnbigftilde)).
 \end{CD}
\]
\hfill\qed
\end{prop}
 
\subsection{Another description of the Stokes graded local systems}
\label{subsection;18.6.24.10}

For $\lambda\in [(\nbigi^{\circ})^{\ast}]$,
we take 
$\vecJ_{\lambda}
=I(\vartheta^u_{0,\lambda},(1-\omega^{-1})\pi/2)
\in T(\lambda)_{<0}$.
We define the map
$\kappa_{\vecJ_{\lambda}}:\real\lrarr\real$
by
\[
 \kappa_{\vecJ_{\lambda}}(\theta^u)
 =\frac{1}{\omega-1}
 (\theta^u-\omega\vartheta^u_{0,\lambda}).
\]
We obtain the local system with Stokes structure
$\bigl(\nbigk_{\lambda}^{\circ\circ},\vecnbigf^{\circ\circ}\bigr)
:=
 \kappa_{\vecJ_{\lambda}}^{-1}\bigl(
 \nbigk_{\lambda_+(\nu_0^+(\vecJ_{\lambda}))},\vecnbigf
 \bigr)$.
Because
\[
 \kappa_{\vecJ_{\lambda}}
 =\Tbb^m\circ
 \kappa^+_{0,\vecJ_{\lambda}+2m(1+\omega^{-1})\pi}
 =\Tbb^m\circ
 \kappa^-_{0,\vecJ_{\lambda}+(2m+1)(1+\omega^{-1})\pi}.
\]
Hence, there exists an isomorphism
$b_{\lambda}:
 (\nbigk_{\lambda}^{\circ\circ},\vecnbigf^{\circ\circ})
\simeq
 (\nbigk_{\lambda}^{\circ},\vecnbigf^{\circ})$
whose restriction to 
$(\vecJ+2m\pi)\cup(\vecJ+(2m+1)\pi)$
are induced by
$(-1)^m\Psi^{m}$.

Let $c:\real\lrarr\real$ be given by
$c(\theta^u)=-\theta^u$.
We set
$\nbigk^{\circ\circ}_0:=
 c^{-1}(\nbigh^{\Sh})$.
Let us observe that
there exists a natural isomorphism
$b_0:\nbigk^{\circ\circ}_0\simeq
 \nbigk^{\circ}_0$.
Take $\vecJ\in T(\nbigi^{\circ})$.
If $\theta^u\in\vecJ$,
we have 
$\nu_0^+(\vecJ)
\subset
 \openopen{-\theta^u-\pi/2}{-\theta^u+\pi/2}$.
Hence, we obtain the natural isomorphisms
\[
 (\nbigk^{\circ\circ}_0)_{|\theta^u}
:=(\nbigk_0)_{|-\theta^u}
\simeq
 K_{0,\nu_0^+(\vecJ)},
\]
which induce the desired isomorphism $b_0$.

We set 
$(\nbigk^{\circ\circ}_{\bullet},\vecnbigf^{\circ\circ})
:=
 \bigoplus_{\lambda\in[\nbigi]}
 (\nbigk_{\lambda}^{\circ\circ},\vecnbigf^{\circ\circ})$.
We constructed the natural isomorphism
$(\nbigk^{\circ\circ}_{\bullet},\vecnbigf^{\circ\circ})
\simeq
(\nbigk^{\circ}_{\bullet},\vecnbigf^{\circ})$.
An action of $2\pi\seisuu$
on 
$(\nbigk^{\circ\circ}_{\bullet},\vecnbigf^{\circ\circ})$
is induced
by the isomorphism $b$
and the $2\pi\seisuu$-action
on $(\nbigk^{\circ}_{\bullet},\vecnbigf^{\circ})$.
\index{local system with Stokes structure
$(\nbigk^{\circ\circ}_{\bullet},\vecnbigf^{\circ\circ})$}

\begin{rem}
There exist positive integers $n_1,p_1$
such that $n_1/p_1=\omega$ with $\gcd(n_1,p_1)=1$.
For $\lambda\in[(\nbigi^{\circ})^{\ast}]$,
we obtain the following isomorphism:
\[
 a_0:
 (\Tbb^{n_1-p_1})^{\ast}
  (\nbigk_{\lambda}^{\circ\circ},\vecnbigf^{\circ\circ})
\simeq
   (\Tbb^{n_1-p_1})^{\ast}(\nbigk_{\lambda}^{\circ},\vecnbigf^{\circ})
\simeq
 (\nbigk_{\lambda}^{\circ},\vecnbigf^{\circ})
\simeq 
 (\nbigk_{\lambda}^{\circ\circ},\vecnbigf^{\circ\circ}).
\]
We also have the following natural isomorphism:
\begin{multline}
a_1:
 (\Tbb^{n_1-p_1})^{\ast}
 (\nbigk_{\lambda}^{\circ\circ},\vecnbigf^{\circ\circ})
=\kappa_{\vecJ_{\lambda}}^{-1}
 \Bigl(
 (\Tbb^{p_1})^{\ast}(\nbigk_{\lambda_+(\nu_{0}^+(\vecJ_{\lambda}))})
 \Bigr)
\simeq
 \kappa_{\vecJ_{\lambda}}^{-1}
(\nbigk_{\lambda_+(\nu_{0}^+(\vecJ_{\lambda})})) \\
=(\nbigk_{\lambda}^{\circ\circ},\vecnbigf^{\circ\circ}).
\end{multline}
We have $a_0=(-1)^{n_1}a_1$.  
\hfill\qed
\end{rem}

\subsection{Example}

\label{subsection;21.6.13.1}

Let $\omega\in\seisuu_{>1}$.
Let $\nbigi=
\{\alpha_ix^{-\omega}\,|\,i=1,\ldots,N\}
\subset \real_{>0}x^{-\omega}$
be a finite subset.
We have
$T(\nbigi)=
\bigl\{J_m\,\big|\,
m\in\seisuu
\bigr\}$,
where
$J_m:=I(m\omega^{-1}\pi,\omega^{-1}\pi/2)$.
On $J_{2\ell}$,
we have $-\Re(\alpha_ix^{-\omega})<0$,
and hence we have $\nbigi=\nbigi_{J_{2\ell},<0}$.
Similarly,
we have $\nbigi=\nbigi_{J_{2\ell+1},>0}$.

Let $(\nbigv,\nabla)$ be a meromorphic flat bundle
on $(\proj^1,\infty)$
such that $\nbigi_{\infty}(\nbigv)\subset\nbigi$.
Let $(L,\vecnbigf)\in\Loc^{\St}(\nbigi)$
be the corresponding local system with Stokes structure.
In this case,
the associated Stokes shell consists of
$(\nbigk_{\bullet},\vecnbigf)=(L,\vecnbigf)$
and $\vecnbigr=\emptyset$.
Let $(\gbigl^{\gbigf}(\nbigv),\vecnbigf)$
denote the local system with Stokes structure
corresponding to $\Fourier(\nbigv)$ at $\infty$.
Let $(\nbigk^{\gbigf}_{\bullet}(\nbigv),\vecnbigf,\vecnbigr^{\gbigf})$
denote the associated Stokes shell.
For $k\in\seisuu$,
we set
$\beta_k=\exp((2k+1)\pi\sqrt{-1}/(\omega-1))$.
We set
$\nbigi^{\circ}_k:=
\bigl\{
 \langle\omega\rangle'
 \alpha_i^{\frac{1}{\omega-1}}\beta_k
 u^{-\frac{\omega}{\omega-1}}
 \bigr\}$
and $\nbigi^{\circ}=\bigcup_{k=0}^{\omega-2}\nbigi_k^{\circ}$.
We have
$\nbigi_{\infty}\bigl(
\Fourier(\nbigv)
\bigr)\subset\nbigi^{\circ}$.

We set $V=\nbigv(\ast 0)$,
which is a basic meromorphic flat bundle
of level $(\infty,\omega)$
with $\nbigi_{\infty}(V)=\nbigi_{\infty}(\nbigv)$.
Let $(\gbigl^{\gbigf}_{\star}(V),\vecnbigf)$ $(\star=!,\ast)$
denote the local systems with Stokes structure
corresponding to $\Fourier(V(\star 0))$ at $\infty$.
We obtain the associated Stokes shells
$(\nbigk^{\gbigf}_{\bullet}(V(\star 0)),\vecnbigf,\vecnbigr)$.

There exist the natural morphisms
\[
 (\gbigl^{\gbigf}_{!}(V),\vecnbigf)
 \lrarr
 (\gbigl^{\gbigf}(\nbigv),\vecnbigf)
 \lrarr
 (\gbigl^{\gbigf}_{\ast}(V),\vecnbigf).
\]
Note that
$\Gr^{\vecnbigf}_0\gbigl^{\gbigf}(\nbigv)=0$,
and that the morphisms
\[
 \Gr^{\vecnbigf}_{\gminia}(\gbigl^{\gbigf}_{!}(V),\vecnbigf)
 \lrarr
 \Gr^{\vecnbigf}_{\gminia}(\gbigl^{\gbigf}(\nbigv),\vecnbigf)
 \lrarr
 \Gr^{\vecnbigf}_{\gminia}(\gbigl^{\gbigf}_{\ast}(V),\vecnbigf)
\]
are isomorphisms if $\gminia\neq 0$.
Hence, $(\gbigl^{\gbigf}(\nbigv),\vecnbigf)$
is the extension of the base tuple
$(\gbigl^{\gbigf}_{!}(V),\vecnbigf)
\lrarr
(\gbigl^{\gbigf}_{\ast}(V),\vecnbigf)$
by the trivial maps
$\Gr^{\vecnbigf}_{0}\gbigl^{\gbigf}_{!}(V)
 \lrarr
0
 \lrarr
 \Gr^{\vecnbigf}_{0}\gbigl^{\gbigf}_{\ast}(V)$.
Equivalently, in terms of Stokes shells,
there exist natural morphisms
\[
(\nbigk^{\gbigf}_{\bullet}(V(!0)),\vecnbigf,\vecnbigr)
\lrarr
(\nbigk^{\gbigf}_{\bullet}(\nbigv),\vecnbigf,\vecnbigr)
\lrarr
(\nbigk^{\gbigf}_{\bullet}(V(\ast 0)),\vecnbigf,\vecnbigr),
\]
and we have
$\nbigk^{\gbigf}_0(\nbigv)=0$,
and
$\nbigk^{\gbigf}_{\lambda}(V(!0))
\simeq
\nbigk^{\gbigf}_{\lambda}(\nbigv)
\simeq
\nbigk^{\gbigf}_{\lambda}(V(\ast 0))$
for $\lambda\in [\nbigi^{\circ}]$.
In this way,
we can compute
$(\nbigk^{\gbigf}_{\bullet}(\nbigv),\vecnbigf,\vecnbigr^{\gbigf})$
from 
$(\nbigk^{\gbigf}_{\bullet}(V(\star 0)),\vecnbigf,\vecnbigr^{\gbigf})$.

We set
$\vecJ_{k,m}=I\bigl((1-\omega^{-1})m\pi+(2k+1)\omega^{-1}\pi,
(1-\omega^{-1})\pi/2
\bigr)$ for $k=0,\ldots,\omega-2$
and $m\in\seisuu$.
Then,
$T(\nbigi_k^{\circ})=\{\vecJ_{k,m}\,|\,m\in\seisuu\}$
and
$T(\nbigi^{\circ})=\bigcup_{k=0}^{\omega-2}T(\nbigi_k^{\circ})$.
We have
$\nbigi^{\circ}_{\vecJ_{k,2\ell},<0}
=\nbigi_k^{\circ}
=\nbigi^{\circ}_{\vecJ_{k,2\ell+1},>0}$.
We have
$\nu_0^+(\vecJ_{k,2\ell})=J_{-2(k+\ell(\omega-1))-1}$
and
$\nu_0^-(\vecJ_{k,2\ell+1})=J_{-2(k+\ell(\omega-1))}$.

We obtain local systems with filtrations
$(\kappa_{0,\vecJ_{k,2\ell+1}}^-)^{-1}
\bigl(
L_{|\overline{\nu_0^-(\vecJ_{k,2\ell+1})}},\vecnbigf
\bigr)$
on $\vecJ_{k,2\ell+1}$
and 
$(\kappa_{0,\vecJ_{k,2\ell}}^+)^{-1}
\bigl(
L_{|\overline{\nu_0^+(\vecJ_{k,2\ell})}},\vecnbigf
\bigr)$
on $\vecJ_{k,2\ell}$.
The index sets are $\nbigi_k^{\circ}$.
Because
$\nu_0^-(\vecJ_{k,2\ell+1})
=\nu_0^+(\vecJ_{k,2\ell})+\omega^{-1}\pi$,
we have the natural isomorphism
at $\vartheta^u_1\in
\overline{\vecJ_{k,2\ell}}
\cap\overline{\vecJ_{k,2\ell+1}}$:
\[
(\kappa_{0,\vecJ_{k,2\ell}}^+)^{-1}
\bigl(
L_{|\overline{\nu_0^+(\vecJ_{k,2\ell})}},\vecnbigf
\bigr)_{|\vartheta^u_1}
\simeq
 (\kappa_{0,\vecJ_{k,2\ell+1}}^-)^{-1}
\bigl(
L_{|\overline{\nu_0^-(\vecJ_{k,2\ell+1})}},\vecnbigf
\bigr)_{|\vartheta^u_1}.
\]
Because
$\nu_0^+(\vecJ_{k,2(\ell+1)})
=\nu_0^-(\vecJ_{k,2\ell+1})+\omega^{-1}\pi-2\pi$,
we obtain the isomorphism
\[
 (\kappa_{0,\vecJ_{k,2\ell+1}}^-)^{-1}
\bigl(
L_{|\overline{\nu_0^-(\vecJ_{k,2\ell+1})}},\vecnbigf
\bigr)_{|\vartheta^u_2}
\simeq
 (\kappa_{0,\vecJ_{k,2\ell}}^+)^{-1}
\bigl(
L_{|\overline{\nu_0^+(\vecJ_{k,2\ell})}},\vecnbigf
\bigr)_{|\vartheta^u_2}
\]
at $\vartheta^u_2\in
\overline{\vecJ_{k,2\ell+1}}
\cap
\overline{\vecJ_{k,2(\ell+1)}}$,
as the $(-1)$ times the natural isomorphism.
By patching them,
we obtain a local system with a family of filtrations
$(\nbigk^{\circ}_k,\vecnbigf)$ on $\real$.
There exist natural isomorphism
$\Tbb^{-1}(\nbigk^{\circ}_k,\vecnbigf)
\simeq
(\nbigk^{\circ}_{k-\omega},\vecnbigf)$,
where we consider $k-\omega$ in $\seisuu/(\omega-1)\seisuu$,
we have the natural $2\pi\seisuu$-action
on
$(\nbigk^{\circ}_{\bullet},\vecnbigf)
=\bigoplus_{k=0}^{\omega-2}
(\nbigk^{\circ}_k,\vecnbigf)$.
According to Proposition \ref{prop;24.3.26.31},
we have
$(\nbigk^{\gbigf}_{\bullet}(\nbigv),\vecnbigf)
\simeq
(\nbigk^{\circ}_{\bullet},\vecnbigf)$.

Let us compute $\vecnbigr^{\gbigf}$.
The non-trivial terms are
$(\nbigr^{\gbigf})^{\vecJ}_{\vecJ'}$
with
\begin{equation}
\label{eq;21.6.12.3}
 \vecJ-(1-\omega^{-1})\pi<\vecJ'<\vecJ+(1-\omega^{-1})\pi
\end{equation}
for
$\vecJ=\vecJ_{k_1,2\ell_1+1}$
and $\vecJ'=\vecJ_{k_2,2\ell_2}$
for some $0\leq k_1,k_2\leq \omega-2$
and $\ell_1,\ell_2\in\seisuu$.
Here, if $\vecJ=\vecJ'$, we regard it as
$(\nbigr^{\gbigf})^{\vecJ_-}_{\vecJ_+}$.
The condition (\ref{eq;21.6.12.3})
is equivalent to
\begin{equation}
\label{eq;25.2.10.1}
 (\omega-1)\ell_1+k_1<(\omega-1)\ell_2+k_2
<(\omega-1)\ell_1+k_1+\omega-1.
\end{equation}
Note that for any integer $m$
satisfying
$(\omega-1)\ell_1+k_1<m
<(\omega-1)\ell_1+k_1+\omega-1$,
there uniquely exist $(k_2,\ell_2)\in\seisuu^2$
such that
$m=k_2+(\omega-1)\ell_2$
and $0\leq k_2\leq \omega-2$.

Let $(k_1,\ell_1)$ and $(k_2,\ell_2)$
be pairs satisfying (\ref{eq;25.2.10.1}).
If 
$\ell_1(\omega-1)+k_1<(\omega-1)\ell_2+k_2
<\ell_1(\omega-1)+k_1+(\omega-1)/2$,
i.e.,
$\vecJ-(1-\omega^{-1})\pi<\vecJ'<\vecJ$,
we have
$\nu_0^-(\vecJ_{k_1,\ell_1})-\pi
<\nu_0^+(\vecJ_{k_2,\ell_2})
<\nu_0^-(\vecJ_{k_1,\ell_1})-\omega^{-1}\pi$.
As explained in \S\ref{subsection;25.2.10.1},
$(\nbigr^{\gbigf})^{\vecJ_{k_1,\ell_1}}_{\vecJ_{k_2,\ell_2}}$
is equal to the following isomorphism induced by the parallel transport:
\[
 H^0(\nu_0^-(\vecJ_{k_1,\ell_1}),L)
 \simeq
 H^0(\nu_0^+(\vecJ_{k_2,\ell_2}),L).
\]
If 
$\ell_1(\omega-1)+k_1+(\omega-1)/2<(\omega-1)\ell_2+k_2
<\ell_1(\omega-1)+k_1+\omega-1$,
i.e.,
$\vecJ<\vecJ'<\vecJ+(1-\omega^{-1})\pi$,
we have
$\nu_0^-(\vecJ_{k_1,\ell_1})+\omega^{-1}\pi
<\nu_{-1}^+(\vecJ_{k_2,\ell_2})
=\nu_0^+(\vecJ_{k_2,\ell_2})+2\pi
<\nu_0^-(\vecJ_{k_1,\ell_1})+\pi$.
As explained in \S\ref{subsection;25.2.10.1},
$(\nbigr^{\gbigf})^{\vecJ_{k_1,\ell_1}}_{\vecJ_{k_2,\ell_2}}$
is equal to the isomorphism
\[
 H^0(\nu_0^-(\vecJ_{k_1,\ell_1}),L)
 \stackrel{a}{\simeq}
 H^0(\nu_{-1}^+(\vecJ_{k_2,\ell_2}),L)
 \stackrel{-\Psi}\simeq
 H^0(\nu_0^+(\vecJ_{k_2,\ell_2}),L),
\]
where $a$ is the parallel transport,
and $\Psi$ is induced by the $2\pi\seisuu$-action on $L$.
If 
$\ell_1(\omega-1)+k_1+(\omega-1)/2=(\omega-1)\ell_2+k_2$,
we have
$\vecJ_{k_2,\ell_2}=\vecJ_{k_1,\ell_1}=:\vecJ$,
and
$\nu_0^+(\vecJ)=\nu_0^-(\vecJ)-\pi$.
By Proposition \ref{prop;24.3.26.31},
$(\nbigr^{\gbigf})^{\vecJ_-}_{\vecJ_+}$
is equal to the isomorphism
\[
 H^0(\nu_0^-(\vecJ),L)
\simeq
 H^0(\nu_0^+(\vecJ),L)
\]
obtained as the $-1$ times the parallel transport.

\begin{rem}
If $\omega=2$,
there is no integer $m$
satisfying
$(\omega-1)\ell_1+k_1<m<(\omega-1)\ell_1+k_1+\omega-1$.
Hence, $\vecnbigr^{\gbigf}=\emptyset$.
Moreover,
$(\nbigk^{\circ}_{\bullet},\vecnbigf)
=(\nbigk^{\circ}_1,\vecnbigf)$
is isomorphic to
the pull back of
$(L,\vecnbigf)$.
This recovers a result in {\rm\cite{Sabbah-pure-Gaussian}}.
(See Proposition {\rm\ref{prop;24.4.20.22}}.)
\hfill\qed
\end{rem}

\chapter{Estimate of growth orders}
\label{section;18.6.3.21}

\section{Preliminaries}

\subsection{Horizontal and vertical paths}

Let $\realbar_{\geq 0}:=\realbar_{\geq 0}\cup\{\infty\}$.
Set $X:=\realbar_{\geq 0}\times\real$.
For $\theta_1,\theta_2\in\real$ and $0<r<\infty$,
let $\gamma_{h}(r;\theta_1,\theta_2):[0,1]\lrarr X$
denote the path given as 
$\gamma_{h}(r;\theta_1,\theta_2)(s)=
 (r,s\theta_2+(1-s)\theta_1)$.
\index{path $\gamma_h(r;\theta_1,theta_2)$}
For $\theta\in\real$ and $0\leq r_2<r_1<\infty$,
let $\gamma_{v}(r_1,r_2;\theta):[0,1]\lrarr X$
denote the path given as
$\gamma_{v}(r_1,r_2;\theta)(s)=((1-s)r_1+sr_2,\theta)$.
\index{path $\gamma_v(r_1,r_2;\theta)$}
For $\theta\in\real$ and $0\leq r<\infty$,
let $\gamma_v(\infty,r;\theta)$ denote the path
given as 
$\gamma_v(\infty,r;\theta)(s)
=(r+s(1-s)^{-1},\theta)$.
\index{path $\gamma_v(\infty,r;\theta)$}

We identify $\projtilde^1$
with $\realbar_{\geq 0}\times S^1$
by the polar coordinate.
We have the morphism $\varphi:X\lrarr\projtilde^1$
induced by
$\theta\longmapsto e^{\sqrt{-1}\theta}$.
We use the same notation
to denote the induced paths on $\projtilde^1$.

\subsection{Metrics}

Let $C$ be a complex curve with a discrete subset $D$.
Let $\nbigv$ be a locally free $\nbigo_C(\ast D)$-module.
A Hermitian metric $h_{\nbigv}$
of $\nbigv_{|C\setminus D}$
is called adapted to the meromorphic structure of $\nbigv$
if the following holds:
\index{adapted metric}
\begin{itemize}
\item
Take any point $P$ of $D$.
Take a frame 
$\vecv=(v_1,\ldots,v_r)$ of $\nbigv$ 
on a neighbourhood $C_P$ of $P$ in $C$
with a holomorphic coordinate $z_P$
such that $z_P(P)=0$.
Let $H$ denote the Hermitian-matrix valued 
function on $C_P\setminus\{P\}$
determined by
$H_{i,j}=h(v_i,v_j)$.
Then, 
$A^{-1}|z|^{N}I_r\leq H\leq A|z|^{-N}I_r$
for some positive constants $A$ and $N$,
where $I_r$ denote the $r$-th identity matrix.
\end{itemize}
Take $P\in D$ with a neighbourhood as above.
If an adapted metric $h_{\nbigv}$ is given,
a section $f$ of $\nbigv$ on $C_P\setminus P$
is a section of $\nbigv$ on $C_P$
if and only if
$|f|_{h^{\nbigv}}=O(|z_P|^{-N})$ for some $N$.

\subsection{Stokes filtrations and adapted metrics}

Let $(\nbigv,\nabla)$ be a meromorphic flat bundle
on $(\Delta,0)$.
Let $h_{\nbigv}$ be an adapted metric for $\nbigv$.
Let $\nbigi(\nbigv)$ denote the set of ramified irregular values.
Let $(L,\vecnbigf)$ be the associated
$2\pi\seisuu$-equivariant local system
with Stokes structure.
Let $\varphi_1:\real\lrarr\varpi^{-1}(0)$
be given by $\varphi_1(\theta)=e^{\sqrt{-1}\theta}$.

Let $\varpi:\Deltatilde\lrarr \Delta$
denote the oriented real blow up
along $0$.
We have the local system $\nbigl$ on $\Deltatilde$
corresponding to the flat bundle 
$(\nbigv,\nabla)_{|\Delta\setminus\{0\}}$.

Take $\theta\in\real$.
Let $s\in L_{\theta}$.
We have the induced flat section
$\stilde$ of $\nbigl$
on a neighbourhood 
$U$ of $\varphi_1(\theta)$ in $\Deltatilde$.
Note that
we may naturally regard
elements of $\nbigi(\nbigv)$ as functions on
$U\setminus\varpi^{-1}(0)$
by the choice of $\theta$.
The following lemma is obvious.

\begin{lem}
$s$ is contained in $\nbigf^{\theta}_{\gminia}$
for $\gminia\in\nbigi(\nbigv)$
if and only if
$|s|_{h^{\nbigv}}=
 O\Bigl(
 \exp\bigl(
 -\Re(\gminia)
 \bigr)\cdot|z|^{-N}
 \Bigr)$
for some $N$.
\hfill\qed
\end{lem}

\subsection{A property of Stokes filtrations}

We shall also use the following well known and easy lemma
for Stokes structures.
\begin{lem}
\label{lem;20.10.23.1}
For any $\theta\in\real$,
there exists $\epsilon_0>0$
such that
 $\nbigf^{\theta}_{\gminia}
 =\nbigf^{\theta+\epsilon}_{\gminia}
 \cap\nbigf^{\theta-\epsilon}_{\gminia}$
for any $\gminia\in\nbigi(\nbigv)$
and for any $0<\epsilon<\epsilon_0$.
 \end{lem}

\subsection{Estimate of the norms of the induced sections}

Let $D$ be a finite subset of $\cnum$.
Let $(\nbigv,\nabla)$ be a meromorphic flat bundle
on $(\proj^1,D\cup\{\infty\})$.
We take a metric $h_{\nbigv}$
of $\nbigv_{|\cnum\setminus D}$
which is adapted to the meromorphic structure of $\nbigv$.
For $\varrho\in\Dsf(D)$,
we set
$\nbigv^{\gbigf}_{\varrho}:=
 \Fourier_+(\nbigv(\varrho))$.
\index{$\nbigd$-module $\nbigv^{\gbigf}_{\varrho}$}
We take a neighbourhood $U_{\infty}$ of $\infty$
such that
$\nbigv^{\gbigf}_{\varrho|U_{\infty}}$
are meromorphic flat bundles
on $(U_{\infty},\infty)$.
We take metrics $h_{\varrho}$ of
$\nbigv^{\gbigf}_{\varrho|U_{\infty}\setminus\{\infty\}}$
which are adapted to the meromorphic structure.
\index{metric $h_{\varrho}$}

For any $u\in\cnum^{\ast}$,
we have the natural metric $h_u$ of 
$\nbige(zu^{-1})=\nbigo_{\proj^1}(\ast\infty)$
given by $|1|_{h_u}=1$.
\index{metric $h_u$}
Let $h^{\nbigv}_u$ denote the induced metric
on $\nbigv\otimes\nbige(zu^{-1})$.
\index{metric $h^{\nbigv}_u$}

Let $\projtilde^1_{D\cup\infty}\lrarr\proj^1$
denote the oriented real blow up along $D\cup\{\infty\}$.
\index{oriented real blow up $\projtilde^1_{D\cup\infty}$}
Take $\theta^u\in\real$ and $C>0$,
and we put $u:=Ce^{\sqrt{-1}\theta^u}$.

Let $d_{\proj^1}$ be the distance induced by
a Riemannian metric of $\proj^1$.
For any $\ell\in\real$ and $\star\in\{\ast,!\}$,
we set
$W_{\varrho,\star,\ell}(\beta):=
 \prod_{\alpha\in \varrho^{-1}(\star)}
 d_{\proj^1}(\alpha,\beta)^{-\ell}$.

From $\varrho\in \Dsf(D)$,
we obtain $\varrho_1\in\Dsf(D\cup\{\infty\})$
by setting
$\varrho_1(P)=\varrho(P)$ $(P\in D)$
and
$\varrho_1(\infty)=!$.
If $u\in U_{\infty}\setminus\{\infty\}$,
a section of
$\nbigc^{-1}_{\projtilde^1_{D\cup\infty},\del\projtilde^1_{D\cup\infty}}
\otimes
\nbigl^{\varrho_1}\bigl(
\nbigv\otimes\nbige(zu^{-1})
\bigr)$
is called a $\varrho$-type $1$-chain of
$\nbigv\otimes\nbige(zu^{-1})$.
(See \S\ref{subsection;20.10.21.1}
for the notation.)
It is called a $\varrho$-type $1$-cycle
of $\nbigv\otimes\nbige(zu^{-1})$
if it is a cycle.

\begin{rem}
If $u\in U_{\infty}\setminus\{\infty\}$,
the regular part of
$\bigl(\nbigv\otimes\nbige(zu^{-1})\bigr)\otimes\cnum[\![z^{-1}]\!]$
is $0$,
and hence
$(\nbigv\otimes\nbige(zu^{-1}))(!\infty)
\simeq
(\nbigv\otimes\nbige(zu^{-1}))(\ast\infty)$.
\hfill\qed
\end{rem}

Let $\vecc(t)=\sum_{i=1}^N c_{i,t}\otimes\gamma_{i,t}$
$(0<t<t_0)$ be a family of $\varrho$-type $1$-cycles for
$(\nbigv,\nabla)\otimes\nbige(zu^{-1}t^{-1})$,
i.e.,
$\gamma_{i,t}$ are families of paths on
$\projtilde^1_{D\cup\infty}$,
and $c_{i,t}$ are flat sections of
$\nbigv\otimes\nbige(zu^{-1}t^{-1})$
along $\gamma_{i,t}$
satisfying the condition associated with $\varrho$.
Suppose that there exist
$Q_i(t)\in \real[t^{-1/e}]$ $(i=1,\ldots,N)$
and $m\in\real$
such that the following holds.
\begin{itemize}
\item
For any $\ell\in\real$,
there exist $M(\ell)>0$ and $C(\ell)>0$
such that 
\[
 \int_{\gamma_{i,t}}
 \bigl|
 c_{i,t}
 \bigr|_{h_{ut}^{\nbigv}}
 W_{\varrho,\ast,m}
 W_{\varrho,!,\ell}
 \cdot(1+|z|^2)^{\ell}|dz|
\leq
 C(\ell)\exp\bigl(Q_i(t)\bigr)t^{-M(\ell)}.
\]
\end{itemize}

Let $[\vecc(t)]$
denote the element of $\nbigv^{\gbigf}_{\varrho|ut}$
induced by $\vecc(t)$.

\begin{lem}
\label{lem;18.5.27.10}
There exist $C>0$ and $M>0$
such that 
$\bigl|
 [\vecc(t)]
 \bigr|_{h_{\varrho}}
\leq
 C\exp\bigl(\max_iQ_i(t)\bigr) t^{-M}$.
\end{lem}
\pf
Let $\varrhobar_1\in\Dsf(D\cup\{\infty\})$
be determined by
$\{
 \varrhobar_1(\alpha),\varrho_1(\alpha)\}
=\{!,\ast\}$ for any $\alpha\in D\cup\{\infty\}$.
We set
$\nbigv_1:=
 \Fourier_-(\nbigv^{\lor}(\varrhobar_1))$.
There exists a natural isomorphism
$\nbigv_{1|U_{\infty}}
 \simeq
 (\nbigv_{\varrho|U_{\infty}}^{\gbigf})^{\lor}$
as a meromorphic flat bundle.
Set
$U_{\infty}^{\ast}:=U_{\infty}\setminus\{\infty\}$.
The natural pairing
$\langle\!\langle\cdot,\cdot \rangle\!\rangle:
 \nbigv^{\gbigf}_{\varrho|U_{\infty}^{\ast}}
\otimes
 \nbigv_{1|U_{\infty}^{\ast}}
\lrarr
\nbigo_{U_{\infty}^{\ast}}$ at $tu$
is induced by the following pairing
(see \S\ref{subsection;24.3.29.110}):
\[
 H_1^{\varrho_1}\bigl(
 \cnum\setminus D,
 \nbigv\otimes\nbige(zu^{-1}t^{-1})
 \bigr)
\otimes
 \hyperh^1\bigl(
 \proj^1,
 \nbigv^{\lor}(\varrhobar_1)\otimes\nbige(-zu^{-1}t^{-1})
 \otimes\Omega^{\bullet}
 \bigr)
\lrarr
 \cnum.
\]

We have the natural section
$1$ of $\nbige(-zu^{-1})$,
which we denote by $e_{-u^{-1}}$.
Set $D_1:=\varrhobar_1^{-1}(!)$
and $D_2=\varrhobar_1^{-1}(\ast)\ni\infty$.
For $a\geq 0$,
let $\nbigc_{\varrhobar_1}^{\bullet}(\nbigv^{\lor})_{-a}$
be the complex as in \S\ref{subsection;24.3.29.120}.
If $a$ is sufficiently large compared with $|m|$,
the following holds:
\begin{itemize}
\item
     Let $f\,dz$ be a section of
     $\nbigc^{1}_{\varrhobar_1}(\nbigv^{\lor})_{-a}$.
     For any $u\in U_{\infty}^{\ast}$,
     we obtain the section
     $f\otimes e_{-u^{-1}}\,dz$
     of
     $\nbigc^{1}_{\varrhobar_1}(\nbigv^{\lor}\otimes\nbige(-u^{-1}z))_{-a}$,
     which
     induces an element of 
$\nbigv_{1|u}=
 \hyperh^1\bigl(
 \proj^1,
 \nbigv^{\lor}(\varrho_1)\otimes\nbige(-zu^{-1})
 \otimes\Omega^{\bullet}
\bigr)$.
We obtain 
the induced section of
$\nbigv_{1|U_{\infty}^{\ast}}$
denoted by $F_f$.
\item
     An appropriate tuple of such sections
     $f_1,\ldots,f_L$
     induces a frame
     $F_{f_1},\ldots,F_{f_L}$ of
     $\nbigv_{1|U_{\infty}}$.
\end{itemize}
Because
$\bigl|
 \bigl(
 c_{i,t},
 f\cdot e_{-u^{-1}t^{-1}}\,dz
 \bigr)
\bigr|
\leq
 C_1
 |c_{i,t}|_{h^{\nbigv}_{tu}}
\cdot
 |f|_{h^{\nbigv^{\lor}}}
 |dz|$,
we obtain
\begin{equation}
\label{eq;18.5.27.1}
\Bigl|
 \langle\!\langle
 [\vecc(t)],
 F_f
 \rangle\!\rangle_{|tu}
\Bigr|
\leq
 \sum_{i=1}^N
 C_i\exp(Q_i(t))t^{-N_i}.
\end{equation}
We obtain the claim of the lemma
from the estimates (\ref{eq;18.5.27.1}).
\hfill\qed

\section{Some estimates}

\subsection{Critical points of some functions on $\real$}

The calculations in this subsection
are essentially contained in 
\cite[\S3, \S4]{Mochizuki-Fourier-old}.
Take $\veckappa:=(\kappa_1,\kappa_2)\in\real^2$.
Take $\omega\in\rnum_{>0}$.
We consider the following function:
\index{function $H_{\veckappa}$}
\[
 H_{\veckappa}(\theta)=
-\frac{1}{\omega}\cos(\omega\theta-\kappa_1)
-\cos(\theta-\kappa_2).
\]
We have
$\del_{\theta}H_{\veckappa}(\theta)
=\sin(\omega\theta-\kappa_1)+\sin(\theta-\kappa_2)$.
We have
$\del_{\theta}H_{\veckappa}(\theta)=0$
if and only if
one of the following holds:
(i) $\omega\theta-\kappa_1=-(\theta-\kappa_2)+2m\pi$
 for an integer $m$,
(ii) $\omega\theta-\kappa_1=\theta-\kappa_2+(2q+1)\pi$
 for an integer $q$.
We set
\index{point $[\omega,m;\veckappa]$}
\[
 [\omega,m;\veckappa]:=
 \frac{1}{\omega+1}(\kappa_1+\kappa_2+2m\pi).
\]
If $\omega\neq 1$,
we also set
\index{point $(\omega,m;\veckappa)$}
\[
 (\omega,m;\veckappa):=
 \frac{1}{\omega-1}\bigl(
 \kappa_1-\kappa_2+(2m+1)\pi
 \bigr).
\]
Then, the condition (i) is equivalent to
$\theta=[\omega,m;\veckappa]$ for an integer $m$,
and the condition (ii) is equivalent to
$\theta=(\omega,m;\veckappa)$ for an integer $m$.
We have the following:
\begin{equation}
\label{eq;24.2.16.1}
 \cos\bigl(
 \omega[\omega,m;\veckappa]-\kappa_1
 \bigr)
=
 \cos\bigl(
 [\omega,m;\veckappa]-\kappa_2
 \bigr),
\end{equation}
\begin{equation}
\label{eq;24.2.16.2}
 \cos\bigl(
 \omega(\omega,m;\veckappa)-\kappa_1
 \bigr)
=
-\cos\bigl(
 (\omega,m;\veckappa)-\kappa_2
 \bigr).  
\end{equation}
Let $\Cr_1(\omega,\veckappa)$ 
denote the set of $[\omega,m;\veckappa]$ $(m\in\seisuu)$.
\index{set $\Cr_1(\omega,\veckappa)$}
Let $\Cr_2(\omega,\veckappa)$
denote the set of $(\omega,m;\veckappa)$ $(m\in\seisuu)$.
\index{set $\Cr_2(\omega,\veckappa)$}
When we consider $\Cr_2(\omega,\veckappa)$,
we implicitly assume that $\omega\neq 1$.
If $\theta_0\in\Cr_1(\omega,\veckappa)$,
we have
\[
 \del_{\theta}^2H_{\veckappa}(\theta_0)
=(\omega+1)\cos(\omega\theta_0-\kappa_1)
=-\omega H_{\veckappa}(\theta_0).
\]
If $\theta_0\in\Cr_2(\omega,\veckappa)$,
we have
\[
 \del_{\theta}^2H_{\veckappa}(\theta_0)
=-(-\omega+1)\cos(\omega\theta_0-\kappa_1)
=\omega H_{\veckappa}(\theta_0).
\]
Hence, for 
$\theta_0\in \Cr_1(\omega,\veckappa)$,
$H_{\veckappa}(\theta_0)$ is maximal
(resp. minimal)
if $H_{\veckappa}(\theta_0)>0$
(resp. $H_{\veckappa}(\theta_0)<0$).
Similarly, for 
$\theta_0\in \Cr_2(\omega,\veckappa)$,
$H_{\veckappa}(\theta_0)$ is maximal
(resp. minimal)
if $H_{\veckappa}(\theta_0)<0$
(resp. $H_{\veckappa}(\theta_0)>0$).

For any $\kappa\in\real$ and $\ell>0$,
we set
\index{sets $\nbigt_{\ell}(\kappa)_{\pm}$}
\[
 \nbigt_{\ell}(\kappa)_+:=\Bigl\{
 \openopen{\ell^{-1}(\kappa-\pi/2+2m\pi)}
 {\ell^{-1}(\kappa+\pi/2+2m\pi)}\,\Big|\,
 m\in\seisuu
 \Bigr\},
\]
\[
  \nbigt_{\ell}(\kappa)_-:=\Bigl\{
 \openopen{\ell^{-1}(\kappa+\pi/2+2m\pi)}
 {\ell^{-1}(\kappa+3\pi/2+2m\pi)}\,\Big|\,
 m\in\seisuu
 \Bigr\}.
\]
We set 
$\nbigt_{\ell}(\kappa):=
\nbigt_{\ell}(\kappa)_+
\sqcup
\nbigt_{\ell}(\kappa)_-$.
\index{set $\nbigt_{\ell}(\kappa)$}
We have
$\pm\cos(\ell\theta-\kappa)>0$
on $\nbigt_{\ell}(\kappa)_{\pm}$.

\begin{lem}
Take $\theta_0\in
 \Cr_1(\omega,\veckappa)\cup
 \Cr_2(\omega,\veckappa)$.
Then, $\theta_0$
is an end point of an interval $J\in \nbigt_{\omega}(\kappa_1)$
if and only if 
$\theta_0$ is an end point of an interval 
$J'\in \nbigt_1(\kappa_2)$.
Moreover, we have
\[
 \cos(\omega(\theta_0+a)-\kappa_1)
 \cos(\theta_0+a-\kappa_2)<0
\]
if $0<|a|$ is sufficiently small.
\hfill\qed
\end{lem}

\begin{lem}
\label{lem;18.5.18.1}
Take $J_1\in \nbigt_{\omega}(\kappa_1)_{\pm}$
and $J_2\in\nbigt_1(\kappa_2)_{\pm}$.
\begin{itemize}
\item
We have $J_1\cap J_2\neq\emptyset$
if and only if
$J_1\cap J_2\cap\Cr_1(\omega,\veckappa)\neq\emptyset$.
Moreover, 
$J_1\cap J_2\cap\Cr_1(\omega,\veckappa)$
consists of one element.
\item
If $\Jbar_1\cap\Jbar_2=\{\theta_0\}$,
then $\theta_0\in \Cr_1(\omega,\veckappa)$.
If $\omega\neq 1$, it is also an element of
$\Cr_2(\omega,\veckappa)$.
\end{itemize}
\end{lem}
\pf
The second claim can be checked by a direct computation.
We can prove the first claim
by using (\ref{eq;24.2.16.1}) and the continuity argument 
with varying $\kappa_2$.
\hfill\qed

\begin{lem}
Suppose $\omega>1$.
Take $J_1\in \nbigt_{\omega}(\kappa_1)_{\pm}$
and $J_2\in\nbigt_1(\kappa_2)_{\mp}$.
\begin{itemize}
\item
We have $\Jbar_1\subset J_2$
if and only if
$J_1\cap J_2\cap\Cr_2(\omega,\veckappa)\neq\emptyset$.
Moreover, 
$J_1\cap J_2\cap\Cr_2(\omega,\veckappa)$
consists of one element.
\item
If $J_1\subset J_2$ 
and $\Jbar_1\setminus J_2=\{\theta_0\}$,
then $\theta_0\in \Cr_2(\omega,\veckappa)$.
It is also an element of $\Cr_1(\omega,\veckappa)$.
\end{itemize}
In particular,
for $J\in \nbigt_{\omega}(\kappa_1)_{\pm}$,
we have
$J\cap\Cr_2(\omega,\veckappa)\neq\emptyset$
if and only if
$J\subset J'$
for some $J'\in\nbigt_1(\kappa_2)_{\mp}$.
\end{lem}
\pf
The second claim is implied by Lemma \ref{lem;18.5.18.1}.
We can prove the first claim
by using the continuity
with varying $\kappa_2$.
\hfill\qed

\vspace{.1in}
Similarly,
we obtain the following.

\begin{lem}
Suppose $\omega<1$.
Take $J_1\in \nbigt_{\omega}(\kappa_1)_{\pm}$
and $J_2\in\nbigt_1(\kappa_2)_{\mp}$.
\begin{itemize}
\item
We have $\Jbar_2\subset J_1$
if and only if
$J_1\cap J_2\cap\Cr_2(\omega,\veckappa)\neq\emptyset$.
Moreover, 
$J_1\cap J_2\cap\Cr_2(\omega,\veckappa)$
consists of one element.
\item
If $J_2\subset J_1$ 
and $\Jbar_2\setminus J_1=\{\theta_0\}$,
then $\theta_0\in \Cr_2(\omega,\veckappa)$.
It is also an element of $\Cr_1(\omega,\veckappa)$.
\end{itemize}
In particular,
for $J\in \nbigt_1(\kappa_2)_{\pm}$,
we have
$J\cap\Cr_2(\omega,\veckappa)\neq\emptyset$
if and only if
$J\subset J'$
for some $J'\in\nbigt_{\omega}(\kappa_1)_{\mp}$.
\hfill\qed
\end{lem}

\begin{cor}
\label{cor;18.5.19.11}
Take $J_1\in \nbigt_{\omega}(\kappa_1)_{\pm}$
and $J_2\in\nbigt_1(\kappa_2)_{\mp}$.
If $J_1\setminus J_2\neq\emptyset$
 and $J_2\setminus J_1\neq\emptyset$,
then we have
$J_1\cap J_2\cap\Cr_2(\omega,\veckappa)=\emptyset$.
In particular,
$H_{\veckappa|J_1\cap J_2}$ is monotonic
on $J_1\cap J_2$.
\hfill\qed
\end{cor}

\begin{lem}
Take $J\in \nbigt_{\omega}(\kappa_1)_{-}$
such that
$J\cap \Cr_1(\omega;\veckappa)\neq\emptyset$.
For  $\theta_0\in J\cap\Cr_1(\omega;\veckappa)$,
there exists $J'\in \nbigt_1(\kappa_2)_{-}$
such that $\theta_0\in J'$.
Moreover, the following holds.
\begin{itemize}
\item
If $\omega>1$,
$\theta_0$ is the unique maximum point of 
$H_{\veckappa|\Jbar}$.
\item
If $\omega<1$,
$\theta_0$ is the unique maximum point of
$H_{\veckappa|\Jbar'}$.
\item
If $\omega=1$,
$\theta_0$ is the unique maximum point of
$H_{\veckappa|\Jbar'\cup\Jbar}$.
\end{itemize}
 \end{lem}
\pf
We can prove the first claim by
using the second claim of Lemma \ref{lem;18.5.18.1}
and the continuity varying $\kappa_2$.
If $\omega>1$,
by the previous lemmas,
we obtain that
$\theta_0$ is the unique critical point of $H_{\veckappa|J}$.
Similarly, if $\omega<1$,
we obtain that 
$\theta_0$ is the unique critical point of $H_{\veckappa|J'}$.
If $\omega=1$,
because $\Cr_2(\omega,\veckappa)=\emptyset$,
$\theta_0$
is the unique critical point of
$H_{\veckappa|J\cup J'}$.
Then, the claim of the lemma follows.
\hfill\qed

\subsection{Behaviour of some functions along paths (1)}
\label{subsection;18.5.20.1}

Take $\veckappa=(\kappa_1,\kappa_2)\in\real^2$
and $\omega\in\rnum_{>0}$.
We consider the following function on
$\real_{>0}\times\real$:
\index{function $F_{\veckappa}$}
\[
 F_{\veckappa}(r,\theta):=
 \frac{r^{-\omega}}{\omega}
 e^{\sqrt{-1}(-\omega\theta+\kappa_1)}
+re^{\sqrt{-1}(\theta-\kappa_2)}.
\]
We obtain the following function:
\index{function $f_{\veckappa}$}
\[
 f_{\veckappa}(r,\theta):=
-\Re F_{\veckappa}(r,\theta)=
-\frac{r^{-\omega}}{\omega}\cos(\omega\theta-\kappa_1)
-r\cos(\theta-\kappa_2).
\]
We have 
$d_{(r,\theta)}f_{\veckappa}=0$
if and only if
$r=1$ and $\theta\in \Cr_1(\omega,\veckappa)$.

\subsubsection{Paths which contain a critical point}
\index{path $\Gamma_{\theta_0}$}

Set $\theta_0:=[\omega,m;\veckappa]$.
Suppose that
$\cos(\omega\theta_0-\kappa_1)>0$.
Let $\Gamma_{\theta_0}$ be the path
$(t,\theta_0)$ $(0<t<\infty)$.
We have
$f_{\veckappa|\Gamma_{\theta_0}}(r)=
-(r^{-\omega}/\omega+r)\cos(\omega\theta_0-\kappa_1)$.
Then, it is easy to see that
$r=1$ is the unique maximum point of
$f_{\veckappa|\Gamma_{\theta_0}}(r)$.
We have
$f_{\veckappa|\Gamma_{\theta_0}}(r)
 \sim
 -\cos(\omega\theta_0-\kappa_1) r$ as $r\to\infty$,
and 
$f_{\veckappa|\Gamma_{\theta_0}}(r)\sim
 -\omega^{-1}\cos(\omega\theta_0-\kappa_1) r^{-\omega}$
as $r\to 0$.

Suppose that
$\cos(\omega\theta_0-\kappa_1)<0$.
There exist the intervals 
$J_1=\openopen{\theta_1}{\theta_1+\omega^{-1}\pi}
 \in\nbigt_{\omega}(\kappa_1)$
and
$J_2=\openopen{\theta_2}{\theta_2+\pi}
 \in\nbigt_1(\kappa_2)$
such that $\theta_0\in J_i$.
Take a small $\delta>0$.
If $\omega>1$,
we set
$\Gamma_{\theta_0}:=
 \gamma_h(1;\theta_1-\delta,\theta_1+\omega^{-1}\pi+\delta)$.
If $\omega<1$,
we set
$\Gamma_{\theta_0}:=
 \gamma_h(1;\theta_2-\delta,\theta_2+\pi+\delta)$.
If $\omega=1$,
we set
$\Gamma_{\theta_0,-}:=
\gamma_h(1;\theta_1-\delta,\theta_2+\pi+\delta)$
and 
$\Gamma_{\theta_0,+}:=
\gamma_h(1;\theta_2-\delta,\theta_1+\pi+\delta)$.
Then, $(1,\theta_0)$ is the unique maximum point
of the restriction of $f_{\veckappa}$ to the paths.

Suppose that 
$\cos(\omega\theta_0-\kappa_1)=0$.
Take a small neighborhood $\nbigu$
of $(1,\theta_0)$.
Take a small positive number $\epsilon>0$.
Note that
$\cos\bigl(\omega(\theta_0+\epsilon)-\kappa_1\bigr)
 \cos\bigl((\theta_0+\epsilon)-\kappa_2\bigr)<0$.
Set $v:=(-1,1)$ 
if $\cos((\theta_0+\epsilon)-\kappa_2)<0$,
or $v:=(1,1)$ 
if $\cos((\theta_0+\epsilon)-\kappa_2)>0$.
Let $\Gamma_{\theta_0}$ be the paths
given as
$(1,\theta_0)+tv$ $(-\epsilon\leq t\leq \epsilon)$.
Then, $(1,\theta_0)$ is the unique maximum point of
$f_{\veckappa|\Gamma_{\theta_0}}$.

\subsubsection{Vertical paths}

Take $\theta_1$ such that
$\cos(\omega\theta_1-\kappa_1)
 \cos(\theta_1-\kappa_2)\neq 0$.
Let us consider the restriction of
$f_{\veckappa}$
to $L_{\theta_1}=\{(r,\theta_1)\,|\,0<r<\infty\}$.
We have the following obvious lemma.
\begin{lem}
\label{lem;20.10.23.10}
If 
$\cos(\omega\theta_1-\kappa_1)
 \cos(\theta_1-\kappa_2)<0$,
$f_{\veckappa|L_{\theta_1}}$ is monotonic.
If $\cos(\omega\theta_1-\kappa_1)>0$
and $\cos(\theta_1-\kappa_2)>0$,
then $f_{\veckappa|L_{\theta_1}}<0$.
\hfill\qed
\end{lem}

\subsubsection{Perturbation of functions}

Let $A$ be a finite subset in
$\{a\in\rnum\,|\,0<a<\omega\}$.
For $\gminic=(\gminic_j)_{j\in A}\in \cnum^A$,
let us consider the following function on 
$\real_{>0}\times\real$:
\index{function $ F_{\veckappa,\gminic}$}
\[
 F_{\veckappa,\gminic}(r,\theta)
= \frac{r^{-\omega}}{\omega}
 e^{\sqrt{-1}(-\omega\theta+\kappa_1)}
+r e^{\sqrt{-1}(\theta-\kappa_2)}
+\sum_{j\in A}
 \gminic_jr^{-j}e^{-\sqrt{-1}j\theta}.
\]
We obtain the following function:
\index{function $f_{\veckappa,\gminic}$}
\[
 f_{\veckappa,\gminic}(r,\theta):=
-\Re F_{\veckappa,\gminic}
=
-\frac{r^{-\omega}}{\omega}
 \cos(\omega\theta-\kappa_1)
-r\cos(\theta-\kappa_2)
-\sum_{j\in A}
 \Re\bigl(
 \gminic_jr^{-j} e^{-\sqrt{-1}j\theta}
 \bigr).
\]
We may regard $F_{\veckappa,\gminic}$
as a holomorphic function of
$\eta=\log r+\sqrt{-1}\theta$.

Take 
$\eta_0=\sqrt{-1}\theta_0$
with $\theta_0=[\omega,m,\veckappa]$.
Let $\Gamma_{\theta_0}$ be the paths as above.
The following lemma is easy to see.

\begin{lem}
For any $\epsilon_1>0$,
there exists $\delta>0$
such that the following holds
if $|\gminic|<\delta$.
\begin{itemize}
\item
There exists a unique root $\eta_{\gminic}$ of 
the function
 $\del_{\eta}F_{\veckappa,\gminic}$
in $\{|\eta-\eta_0|<\epsilon_1\}$.
\item
There exists a path $\Gamma_{\theta_0,\gminic}$
 such that 
(i)  $\Gamma_{\theta_0,\gminic}$ contains $\eta_{\gminic}$,
(ii)  $\Gamma_{\theta_0,\gminic}$ and $\Gamma_{\theta_0}$
are equal on the outside of 
$\Gamma_{\theta_0}^{-1}\bigl(
 \{|\eta-\eta_0|<\epsilon_1\}\bigr)$,
(iii) $\eta_{\gminic}$ is the unique maximum point of
 the restriction of $\Re F_{\veckappa,\gminic}$
     to $\Gamma_{\theta_0,\gminic}$.
     \index{path $\Gamma_{\theta_0,\gminic}$}
\hfill\qed
\end{itemize}
\end{lem}

\subsection{Behaviour of some functions along paths (2)}

Take $\omega\in\rnum_{>1}$.
Take $\veckappa=(\kappa_1,\kappa_2)\in\real^2$.
We consider the following function on 
$\real_{>0}\times\real$:
\index{function $G_{\veckappa}$}
\[
 G_{\veckappa}(r,\theta)=
 \frac{r^{-\omega}}{\omega}
 e^{\sqrt{-1}(-\omega\theta+\kappa_1)}
+r^{-1}e^{\sqrt{-1}(-\theta+\kappa_2)}
\]
We obtain the following function:
\index{function $g_{\veckappa}$}
\[
 g_{\veckappa}(r,\theta):=
-\Re G_{\veckappa}(r,\theta)
=-\frac{r^{-\omega}}{\omega}\cos(\omega\theta-\kappa_1)
-r^{-1}\cos(\theta-\kappa_2).
\]

\subsubsection{Paths which contain a critical point}

\index{path $\Gamma_{\theta_0}$}

Let $\theta_0=(\omega,m;\veckappa)$.
Suppose that
$\cos(\omega\theta_0-\kappa_1)>0$.
Let $\Gamma_{\theta_0}$ be the path
$(t,\theta_0)$ $(0<t<C)$ for some $C>1$.
We have
$g_{\veckappa|\Gamma_{\theta_0}}(t)=
 -(t^{-\omega}/\omega-t^{-1})\cos(\omega\theta_0-\kappa_1)$.
It is easy to see that
$t=1$ is the unique maximum point of
$g_{\veckappa|\Gamma_{\theta_0}}(t)$.
Note that we have
$g_{\veckappa}(1,\theta_0)=
 (1-\omega^{-1})\cos(\omega\theta_0-\kappa_1)>0$.
We have
$g_{\veckappa|\Gamma_{\theta_0}}(r)\sim
 -\omega^{-1}\cos(\omega\theta_0-\kappa_1) r^{-\omega}$
as $r\to 0$.
If $C$ is sufficiently large,
we have
$0<g_{\veckappa}(C,\theta_0)
 <\!<g_{\veckappa}(1,\theta_0)$.

Suppose that
$\cos(\omega\theta_0-\kappa_1)<0$.
There exist the intervals
$J_1=\openopen{\theta_1}{\theta_1+\omega^{-1}\pi}
 \in\nbigt_{\omega}(\kappa_1)_-$
and 
$J_2=\openopen{\theta_2}{\theta_2+\pi}
 \in\nbigt_1(\kappa_2)_+$
such that $\theta_0\in J_1\subset J_2$.
Take a small $\delta>0$.
We set
$\Gamma_{\theta_0}:=
 \gamma_h(1;\theta_1-\delta,\theta_1+\omega^{-1}\pi+\delta)$.
Then, $(1,\theta_0)$ is the unique maximum point
of the restriction of $g_{\veckappa}$ to the path.

Suppose that 
$\cos(\omega\theta_0-\kappa_1)=0$.
Take a small neighborhood $\nbigu$
of $(1,\theta_0)$.
Take a small positive number $\epsilon>0$,
then 
$\cos(\omega(\theta_0+\epsilon)-\kappa_1)
 \cos((\theta_0+\epsilon)-\kappa_2)<0$.
Set $v:=(-1,1)$ 
if $\cos((\theta_0+\epsilon)-\kappa_2)<0$,
or $v:=(1,1)$ 
if $\cos((\theta_0+\epsilon)-\kappa_2)>0$.
Let $\Gamma_{\theta_0}$ be the paths
given as
$(1,\theta_0)+tv$ $(-\epsilon\leq t\leq \epsilon)$.
Then, $(1,\theta_0)$ is the unique maximum point of
$g_{\veckappa|\Gamma_{\theta_0}}$.

\subsubsection{Vertical paths}

Take $\theta_1$ such that
$\cos(\omega\theta_1-\kappa_1)>0$
and 
$\cos(\theta_1-\kappa_2)>0$.
It is easy to check the following.
\begin{lem}
Let $L_{\theta_1}=\{(r,\theta_1)\,|\, 0<r<\infty\}$.
Then, $g_{|L_{\theta_1}}$ is negative and monotonically
increasing with respect to $r$.
\hfill\qed
\end{lem}

\subsubsection{Perturbation of functions}

Let $A$ be a finite subset in
$\{a\in\rnum\,|\,0<a<\omega\}$.
For $\gminic=(\gminic_j)_{j\in A}\in \cnum^A$,
let us consider the following function on $S$:
\index{function $G_{\veckappa,\gminic}$}
\[
 G_{\veckappa,\gminic}(r,\theta):=
 G_{\veckappa}(r,\theta)
+\sum_{j\in A}
 \gminic_jr^{-j}e^{\sqrt{-1}j\theta}.
\]
We may naturally regard
$G_{\veckappa,\gminic}$
as a holomorphic function of
$\eta=\log r+\sqrt{-1}\theta$.
We set
\index{function $g_{\veckappa,\gminic}$}
\[
 g_{\veckappa,\gminic}:=
-\Re G_{\veckappa,\gminic}.
\]
Take $\eta_0=\sqrt{-1}\theta_0$
with $\theta_0=(\omega,m;\veckappa)$.

\begin{lem}
For any $\epsilon_1>0$,
there exists $\delta>0$
such that the following holds
if $|\gminic|<\delta$.
\begin{itemize}
\item
There exists a unique root $\eta_{\gminic}$ of 
 $\del_{\eta}G_{\veckappa,\gminic}$
in $\{|\eta-\eta_0|<\epsilon_1\}$.
\item
There exists a path $\Gamma_{\theta_0,\gminic}$
 such that 
(i)  $\Gamma_{\theta_0,\gminic}$ contains $\zeta_{\gminic}$,
(ii)  $\Gamma_{\theta_0,\gminic}$ and $\Gamma_{\theta_0}$
 are the same 
on the outside of 
$\Gamma_{\theta_0}^{-1}\bigl(
 \{|\eta-\eta_0|<\epsilon_1\}\bigr)$,
(iii) $\eta_{\gminic}$ is the unique maximum point of
 the restriction of $\Re G_{\veckappa,\gminic}$ 
to $\Gamma_{\theta_0,\gminic}$.
\index{path $\Gamma_{\theta_0,\gminic}$}
\hfill\qed
\end{itemize}
\end{lem}

\subsection{Scaling}

Let $\alpha\neq 0$ and $u\neq 0$.
Let us consider the following function
\[
F(s,\theta)=
 \alpha s^{-\omega}e^{-\sqrt{-1}\omega\theta}
+t^{-1}u^{-1}se^{\sqrt{-1}\theta}.
\]
By a scaling, we obtain
\begin{multline}
 F\bigl(
 (|\alpha||u|t\omega )^{1/(1+\omega)}r,\theta\bigr)=
\\
 (\omega|\alpha|)^{1/(1+\omega)}
 (t|u|)^{-\omega/(1+\omega)}\cdot
 \Bigl(
 \omega^{-1}r^{-\omega}
 e^{-\sqrt{-1}(\omega\theta-\arg(\alpha))}
+re^{\sqrt{-1}(\theta-\arg(u))}
 \Bigr).
\end{multline}

Suppose $\omega>1$.
Let us consider the following function
\[
G(s,\theta)=
 \alpha s^{-\omega}e^{-\sqrt{-1}\omega\theta}
+t^{-1}u^{-1}
 s^{-1}e^{-\sqrt{-1}\theta}.
\]
By a scaling 
we obtain
\begin{multline}
 G\bigl(
 (|\alpha||u|t\omega)^{1/(\omega-1)}r,\theta
 \bigr)
= \\
\bigl(
 \omega|\alpha|
\bigr)^{-1/(\omega-1)}
 \bigl(
 t|u|
 \bigr)^{-\omega/(\omega-1)}
 \Bigl(
 \omega^{-1}
 r^{-\omega}
 e^{\sqrt{-1}(-\omega\theta-\arg(\alpha))}
+r^{-1}
 e^{-\sqrt{-1}(\theta+\arg(u))}
 \Bigr).
\end{multline}

\section{Proof of Theorem \ref{thm;24.3.15.10}}
\label{section;20.11.21.1}

\subsection{Families of cycles}
\label{subsection;18.5.21.10}

Let $(\nbigv,\nabla)$ be a meromorphic flat bundle
on $(\proj^1,\{0,\infty\})$
with regular singularity at $\infty$.
We take a Hermitian metric $h^{\nbigv}$
of $\nbigv_{|\cnum^{\ast}}$
adapted to the meromorphic structure of $\nbigv$.
Set $\omega:=-\ord(\nbigi(\nbigv))$.
Let $\varrho\in \Dsf(\{0\})$.

Take $u\in\cnum^{\ast}$.
Set $\theta^u:=\arg(u)$.
Let $d\in\rnum$ such that $0<d\leq 1$.
Let 
$\vecc(t)$ $(0<t\leq t_0)$
be a family of $\varrho$-type $1$-chains
of $(\nbigv,\nabla)\otimes\nbige(zu^{-1})$
of the following form:
{\small
\begin{multline}
\label{eq;18.5.17.2}
 \vecc(t)= \\
 \Bigl(
 \sum_{i=1}^{N_0} a_{0,i}
 \otimes\gamma_{0,i,t}
+\sum_{j=1}^{N_1}a_{1,j}
 \otimes \gamma_{1,j,t}
+\sum_{j=1}^{N_2}a_{2,j}
 \otimes \gamma_{2,j,t}
+\sum_{k=1}^{N_3}
 b_k\otimes\eta_{k}
+\sum_{i=1}^{N_4} c_i\otimes\Gamma_i
\Bigr)\exp(-zu^{-1}).
\end{multline}
}
Here, we impose the following condition
by using the polar coordinate $z=re^{\sqrt{-1}\theta}$.
\begin{itemize}
\item
$\gamma_{0,i,t}=t^d\gamma_{i,t}$
for a continuous family of paths
$\gamma_{i,t}$ $(0\leq t\leq t_0)$
in $U_0\setminus\{0\}$
whose end points are independent of $t$.     
Note that the family $\gamma_{i,t}$
is assumed to extend at $t=0$.
\item
 $\gamma_{1,j,t}$ are paths of the form
 $\gamma_{h}(t^{d_j}r_{1,j};\phi_{1,j,1},\phi_{1,j,2})$
where $d_j\geq d$.
\item
 $\gamma_{2,j,t}$ are paths of the form
 $\gamma_{v}(t^{d_{2,j,1}}r_{2,j,1},t^{d_{2,j,2}}r_{2,j,2};\phi_{2,j})$
where we impose that
$d_{2,j,1}=0$ or $d_{2,j,1}\geq d$
 and that $d_{2,j,2}\geq d$.
We admit $r_{2,j,2}=0$.
If $d_{2,j,1}=0$,
then
$\gamma_{2,j,t}$
are contained in
$\{z\,|\,\Re(zu^{-1})>0\}$.

\item
Each $\gamma_{p,i,t}$ is contained in a small sector,
and $a_{p,i}$ is a flat section of $\nbigv$ on the sector.
\item
 $\eta_k$ are of the form
 $\gamma_{h}(\epsilon;\psi_{k,1},\psi_{k,2})$,
and they are contained in 
 $\bigl\{
 z\,\big|\,\Re(zu^{-1})>0
 \bigr\}$,
and $b_k$ are flat sections of $\nbigv$
along $\eta_k$.
\item
 $\Gamma_i$ are paths connecting
 $\epsilon e^{\sqrt{-1}\varphi_{i,1}}$
and
 $\infty e^{\sqrt{-1}\varphi_{i,2}}$
in $\{z\in\cnum\,|\, \Re(zu^{-1})>0\}$,
and $c_i$ are flat sections of $\nbigv$ along $\Gamma_i$.
Here, $\varphi_{i,b}$ $(b=1,2)$ satisfy
$\theta^u-\pi/2<\varphi_{i,b}<\theta^u+\pi/2$,
which implies that
$\cos(\varphi_{i,b}-\theta^u)>0$.
\end{itemize}

Note that 
for any $N>0$
there exists $\delta>0$ such that 
\[
 \int_{\Gamma_j}
 |c_j|_{h^{\nbigv}}\cdot
 \exp\bigl(
 -t^{-1}\Re(zu^{-1})
 \bigr)\cdot |z|^N 
=O\Bigl(
 \exp(-\delta t^{-1})
 \Bigr).
\]
We also have the following estimate for some $\delta>0$:
\[
 \int_{\eta_j}
 |b_j|_{h^{\nbigv}}\cdot
 \exp\bigl(
 -t^{-1}\Re(zu^{-1})
 \bigr)
=O\Bigl(
 \exp(-\delta t^{-1})
 \Bigr).
\]
Let $Q\in \real[t^{-1/e}]$
such that
$|t^{1-d}Q(t)|$ are bounded as $t\to 0$.
We say that 
the growth order of $\vecc(t)$ 
is less than $Q$ 
if the following holds.
\index{growth order}
\begin{itemize}
 \item
For any $N>0$,
there exist $C>0$ and $M>0$
such that 
\[
 \int_{\gamma_{0,i,t}}
 \bigl|a_{0,i}\bigr|_{h^{\nbigv}}
 \exp\bigl(-t^{-1}\Re(zu^{-1})\bigr)
 |z|^{-N}
 \leq 
 C\exp\bigl(Q(t)\bigr)t^{-M}.
\]
\item
For any $N>0$, there exist 
$\delta>0$ and $C>0$
such that 
\[
\int_{\gamma_{1,i,t}}
 \bigl|a_{1,i}\bigr|_{h^{\nbigv}}
 \exp\bigl(-t^{-1}\Re(zu^{-1})\bigr)
 |z|^{-N}
\leq
C\exp\bigl(Q(t)-\delta t^{-(1-d)}\bigr).
\]
\item
In the case $\varrho=!$,
for any $N>0$ there exist 
$\delta>0$ and $C>0$
such that
\[
 \int_{\gamma_{2,i,t}}
 \bigl|a_{2,i}\bigr|_{h^{\nbigv}}
 \exp\bigl(-t^{-1}\Re(zu^{-1})\bigr)
 |z|^{-N}
\leq
C\exp\bigl(Q(t)-\delta t^{-(1-d)}\bigr).
\]
In the case $\varrho=\ast$,
there exist $N>0$, $\delta>0$ and $C>0$ such that
\[
 \int_{\gamma_{2,i,t}}
 \bigl|a_{2,i}\bigr|_{h^{\nbigv}}
 \exp\bigl(-t^{-1}\Re(zu^{-1})\bigr)
 |z|^{N}
\leq
C\exp\bigl(Q(t)-\delta t^{-(1-d)}\bigr).
\]
\end{itemize}
Let $\nbigc_0^{\varrho}((\nbigv,\nabla),u,d,Q)$ 
be the set of
families of $\varrho$-type $1$-chains
for $(\nbigv,\nabla)\otimes\nbige(zu^{-1})$
whose growth order is less than $Q$.
\index{set $\nbigc_0^{\varrho}((\nbigv,\nabla),u,d,Q)$}

Let $\vecc(t)\in \nbigc_0^{\varrho}((\nbigv,\nabla),u,d,Q)$
such that $\vecc(t)$ are cycles for any $t$.
\begin{lem}
\label{lem;20.10.24.1}
 The homology classes of
$\vecc(t)$ in
$H^{\varrho}_{1}\bigl(
\cnum^{\ast},
(\nbigv,\nabla)\otimes\nbige(zu^{-1})
\bigr)$
are constant. 
\end{lem}
\pf
Let $0<t_1<t_2\leq t_0$.
For $q=0,1,2$,
we obtain the paths
$\kappa_{0,q,i}(s)=\gamma_{q,i,s}(0)$
and
$\kappa_{1,q,i}(s)=\gamma_{q,i,s}(1)$
for $t_1\leq s\leq t_2$.
It is easy to see that
\[
 \vecc(t_2)-\vecc(t_1)+
   \left(
\sum_{k=0,1,2}
  \sum_{i=1}^{N_k}
  a_{k,i} \otimes (\kappa_{0,k,i}-\kappa_{1,k,i})
  \right)
 \exp(-zu^{-1})
\]
is homologue to $0$.
Because $\vecc(t)$ are cycles for any $t$,
we obtain
\[
  \left(
\sum_{k=0,1,2}
  \sum_{i=1}^{N_k}
  a_{k,i} \otimes (\kappa_{0,k,i}-\kappa_{1,k,i})
  \right)
 \exp(-zu^{-1})
 =0.
\]
Thus, the claims of the lemma follows.
\hfill\qed

\vspace{.1in}

We obtain the family of $\varrho$-type $1$-cycles 
\[
 \vecctilde(t):=\vecc(t)\cdot
 \exp\bigl(-(t^{-1}-1)zu^{-1}\bigr)
\]
of $(\nbigv,\nabla)\otimes\nbige(t^{-1}zu^{-1})$.
They induce a flat section $[\vecctilde(t)]$ of
$\nbigv^{\gbigf}_{\varrho}=\Fourier_+(\nbigv(\varrho))$
along the path $tu$ $(0<t\leq 1)$.
We obtain the following as a special case of
Lemma \ref{lem;18.5.27.10}.

\begin{lem}
\label{lem;20.10.27.20}
 $\bigl|
 [\vecc(t)]
 \bigr|_{h_{\varrho}}=O\bigl(
 \exp(Q) t^{-N}\bigr)$ for some $N>0$
as $t\to 0$.
\hfill\qed
\end{lem}

\subsection{Statements}

We use the notation in \S\ref{section;18.6.3.20}.
Take $u=|u|e^{\sqrt{-1}\theta^u}\in\cnum^{\ast}$,
and set
$\vecI(\theta^u)=\openopen{\theta^u-3\pi/2}{\theta^u-\pi/2}$.
Take $J_{\pm}\in T(\nbigi)$
such that $J_{\pm}\cap(\vecI(\theta^u)+\pi)_{\mp}\neq\emptyset$.
There exist splittings
\[
 L_{J_{\pm},<0}
=\bigoplus_{\gminia\in\nbigitilde_{J,<0}}
 L_{J_{\pm},\gminia}
\]
of the Stokes filtrations
$\nbigftilde^{\theta}$ $(\theta\in J_{\pm})$.
For any $\gminia\in\nbigitilde_{J,<0}$,
we obtain the following map
as the composition of
the restriction of $\Abb^{\rd}_{J,\theta^u}$
and the natural morphism
$H^{\rd}_1\bigl(
\cnum^{\ast},(V,\nabla)\otimes\nbige(zu^{-1})
\bigr)
\lrarr
H^{\rd}_1\bigl(
\cnum^{\ast},(\nbigv,\nabla)\otimes\nbige(zu^{-1})
\bigr)$:
\[
 \Abb^{\rd}_{J_{\pm},\theta^u,\gminia}:
 H^0(J_{\pm},L_{J_{\pm},\gminia})
\lrarr
 H_1^{\rd}\bigl(\cnum^{\ast},
 (\nbigv,\nabla)\otimes\nbige(zu^{-1})\bigr).
\]
For $u_1=|u|e^{\sqrt{-1}\theta^u_1}$
with $|\theta^u-\theta^u_1|<\pi/2$,
there exists the isomorphism
\begin{equation}
\label{eq;18.5.28.1}
  H_1^{\rd}\bigl(\cnum^{\ast},
 (\nbigv,\nabla)\otimes\nbige(zu^{-1})\bigr)
\simeq
   H_1^{\rd}\bigl(\cnum^{\ast},
 (\nbigv,\nabla)\otimes\nbige(zu_1^{-1})\bigr)
\end{equation}
induced by the parallel transport
as in \S\ref{subsection;24.3.14.42}.
We obtain the following morphism:
\[
 \Abb^{\rd}_{J_{\pm},(\theta^u,\theta_1^{u}),\gminia}:
 H^0(J_{\pm},L_{J_{\pm},\gminia})
\lrarr
 H_1^{\rd}\bigl(\cnum^{\ast},
 (\nbigv,\nabla)\otimes\nbige(zu_1^{-1})\bigr).
\]

For any $\gminia\in\nbigitilde_{J,<0}$,
we set
$\gminia^{\circ}:=
 \gbigf^{(0,\infty)}_{(J,0,-)}(\gminia)
 \in\nbigitilde^{\circ}
 =\gbigf^{(0,\infty)}_+(\nbigitilde)$,
 and $d(\omega):=(1+\omega)^{-1}$.
 (See \S\ref{subsection;20.10.23.1}
 for $\gbigf^{(0,\infty)}_{(J,0,-)}(\gminia)$.)
Note that by the choice of $\theta^u=\arg(u)$
we may naturally regard
$\gminia^{\circ}$
as a function on a sector which contains $u$.

\begin{prop}
 \label{prop;18.5.20.30}
\mbox{{}}
\begin{itemize}
 \item If $J\cap (\vecI(\theta^u)+\pi)\neq\emptyset$,
       then
       for any
       $v\in
       H^0(J_{\pm},L_{J_{\pm},\gminia})$,
       $\Abb^{\rd}_{J_{\pm},\theta^u,\gminia}(v)$
       is represented by
       a family of cycles contained in
       $\nbigc^{!}_0
       \Bigl((\nbigv,\nabla),u,d(\omega),
       -\Re\bigl(\gminia^{\circ}(ut)\bigr)
       \Bigr)$.
 \item
      Suppose that
      $J_{\pm}\cap (\vecI(\theta^u)+\pi)_{\mp}$ consists of one point.
      If
      $|\theta^u-\theta^u_1|\neq 0$ is sufficiently small,
      for any
      $v\in
      H^0(J_{\pm},L_{J_{\pm},\gminia})$,
      $\Abb^{\rd}_{J_{\pm},(\theta^u,\theta^u_1),\gminia}(v)$
      is represented by
      a family of cycles in
      $\nbigc^{!}_0
      \Bigl((\nbigv,\nabla),u_1,d(\omega),
      -\Re\bigl(\gminia^{\circ}(u_1t)\bigr)
      \Bigr)$.
\end{itemize}
\end{prop}

Take $J\in T(\nbigi)$ such that 
$J_{\pm}\cap \vecI(\theta^u)_{\mp}\neq\emptyset$.
There exist splittings
\[
 L_{J_{\pm},>0}
=\bigoplus_{\gminia\in \nbigitilde_{J,>0}}
 L_{J_{\pm},\gminia}
\]
of $\nbigftilde^{\theta}$ $(\theta\in J_{\pm})$.
For any $\gminia\in\nbigitilde_{J,>0}$,
and for any $u_1=|u|e^{\sqrt{-1}\theta^u_1}$
with $|\theta^u-\theta^u_1|<\pi/2$,
we have the following map
induced by 
$B_{J_{\pm},\theta^u}$,
the natural morphism
$H^{\rd}_1\bigl(\cnum^{\ast},(V,\nabla)\otimes\nbige(zu^{-1})\bigr)
\lrarr
H^{\rd}_1\bigl(\cnum^{\ast},(\nbigv,\nabla)\otimes\nbige(zu^{-1})\bigr)$
and (\ref{eq;18.5.28.1}):
\[
 B_{J_{\pm},(\theta^u,\theta^u_1),\gminia}:
 H^0(J_{\pm},L_{J_{\pm},\gminia})
\lrarr
 H_1^{\rd}\bigl(\cnum^{\ast},
 (\nbigv,\nabla)\otimes\nbige(zu_1^{-1})\bigr).
\]

For any $\gminia\in\nbigitilde_{J,>0}$,
we set
$\gminia^{\circ}:=
 \gbigf^{(0,\infty)}_{(J,0,+)}(\gminia)
 \in\nbigitilde^{\circ}$.
Note that by the choice of $\theta^u=\arg(u)$
we may naturally regard
$\gminia^{\circ}$
as a function on a sector which contains $u$.
 
\begin{prop}
\label{prop;18.5.20.31}
\mbox{{}}
\begin{itemize}
 \item If $J\cap\vecI(\theta^u)\neq\emptyset$,
 for any 
$v\in H^0(J_{\pm},L_{J_{\pm},\gminia})$,
$B_{J_{\pm},\theta^u,\gminia}(v)$ is represented by a family of cycles
contained in 
 $\nbigc^{!}_0\Bigl(
 (\nbigv,\nabla),u,d(\omega),
  -\Re\bigl(\gminia^{\circ}(ut)\bigr)
  \Bigr)$.
 \item
      Suppose that $J_{\pm}\cap\vecI(\theta^u)_{\mp}$
      consists of a point.
 If
$|\theta^u-\theta^u_1|\neq 0$ is sufficiently small,
for any 
$v\in 
 H^0(J_{\pm},L_{J_{\pm},\gminia})$,
      $B_{J_{\pm},(\theta^u,\theta^u_1),\gminia}(v)$
      is represented by a family of cycles
contained in 
 $\nbigc^{!}_0\Bigl(
 (\nbigv,\nabla),u_1,d(\omega),
  -\Re\bigl(\gminia^{\circ}(u_1t)\bigr)
  \Bigr)$.
\end{itemize}
\end{prop}

\begin{rem}
By modifying the constructions appropriately,
we may also construct 
$1$-cycles representing
$A_{J_{\pm},\theta^u,\gminia}(v)$
and
$B_{J_{\pm},\theta^u,\gminia}(v)$
in the critical cases,
i.e.,
$J_{\pm}\cap (\vecI(\theta^u)+\pi)_{\mp}$
or
$J_{\pm}\cap \vecI(\theta^u)_{\mp}$
consists of one point.
We omit it to simplify the explanations.
\hfill\qed
\end{rem}

Let $J_1\in T(\nbigi)$
such that
$J_{1\pm}\cap (\vecI(\theta^u)+\pi)_{\mp}\neq\emptyset$.
Let $y\in H_1^{\varrho}\bigl(\cnum^{\ast},
 \nbigt_{\omega}(\nbigv,\nabla)\otimes\nbige(zu^{-1})\bigr)$.
We have 
$C^{(J_{1\pm})}_{\infty,\theta^u}(y)
 \in
 H^{\varrho}_1\bigl(
 \cnum^{\ast},(\nbigv,\nabla)\otimes\nbige(zu^{-1})
 \bigr)$.

\begin{prop}
\label{prop;18.5.20.20}
Suppose that $y$ is represented 
by a family of $1$-cycles
contained in $\nbigc_0^{\varrho}(\nbigt_{\omega}(\nbigv,\nabla),u,d,Q)$.
Assume $\omega>(1-d)/d$.
Then, $C^{(J_{1\pm})}_{\infty,\theta^u}(y)$ is represented 
by a family of $1$-cycles
in $\nbigc_0^{\varrho}((\nbigv,\nabla),u,d,Q)$.
\end{prop}

\subsection{Proof of Theorem \ref{thm;24.3.15.10}}

Let us prove Theorem \ref{thm;24.3.15.10}
together with the following proposition.
\begin{prop}
\label{prop;24.3.15.20}
There exists a finite subset
 $\ttS\subset
 \bigl\{\alpha\in\cnum\,\big|\,|\alpha|=1\bigr\}$
 such that the following holds
 unless $|u|^{-1}u\in\ttS$.
\begin{itemize}
\item
For any $\gminia^{\circ}\in \gbigf^{(0,\infty)}_+(\nbigi(\nbigv))$,
any element of
$\nbigf^{\circ\theta^u}_{\gminia^{\circ}}
 H_1^{\varrho}\bigl(
     \cnum^{\ast},(\nbigv,\nabla)\otimes\nbige(zu^{-1})
 \bigr)$
is represented as a sum $\sum c_i$,
where $c_i$ are families of cycles
contained in 
$\nbigc_0^{\varrho}\bigl(
 (\nbigv,\nabla),u,d_i, Q_i \bigr)$
such that
$Q_i(t)\leq -\Re(\gminia^{\circ}(ut))$
for any sufficiently small $t>0$.
\end{itemize}
\end{prop}

We shall prove the claims
of Theorem  \ref{thm;24.3.15.10}
and Proposition \ref{prop;24.3.15.20}
by an induction on
$-\ord(\nbigi)$.
If $\ord(\nbigi)=0$,
Theorem \ref{thm;24.3.15.10} is trivial
because both the filtrations
$\nbigf^{\prime\theta^u}$ and $\nbigf^{\circ\theta^u}$
are indexed
by the trivial partially ordered set $(\{0\})$.
Proposition \ref{prop;24.3.15.20} is restated as follows,
which is easy to see.
\begin{lem}
If $(\nbigv,\nabla)$ is regular singular at $0$,
then 
any $y\in H^{\varrho}_1\bigl(
 \cnum^{\ast},
 (\nbigv,\nabla)\otimes\nbige(zu^{-1})
 \bigr)$
is represented by a cycle
in $\nbigc^{\varrho}_0\bigl(
 (\nbigv,\nabla),u,1,0
 \bigr)$.
\hfill\qed
\end{lem}

We assume that both
Theorem  \ref{thm;24.3.15.10}
and Proposition \ref{prop;24.3.15.20}
are proved in the case $-\ord(\nbigi(\nbigv))<\omega$,
and let us prove them in the case $-\ord(\nbigi(\nbigv))=\omega$.

Let us study the isomorphism (\ref{eq;24.3.14.10}).
The other isomorphisms can be studied similarly.
By Proposition \ref{prop;18.5.20.30},
Proposition \ref{prop;18.5.20.31}
together with
Lemma \ref{lem;20.10.23.1}
and Lemma \ref{lem;20.10.27.20},
for any
$\vecJ\in\gbigm_-(\nbigi^{\circ},\theta^u)$,
we obtain
\[
 \Abb^{\rd}_{\nu_0^-(\vecJ),\theta^u}
 \bigl(
 \nbigf^{\prime\,\theta^u}_{\gminib}
 H^0(\nu_0^-(\vecJ),L_{\nu_0^-(\vecJ),<0})
 \bigr)
 \subset
 \nbigf^{\circ\theta^u}_{\gminib}
 H_1^{\rd}\bigl(
 \cnum^{\ast},
 (\nbigv,\nabla)\otimes\nbige(zu^{-1})
 \bigr),
\]
\[
 B_{\nu_0^+(\vecJ),\theta^u}
 \bigl(
 \nbigf^{\prime\,\theta^u}_{\gminib}
 H^0(\nu_0^+(\vecJ),L_{\nu_0^-(\vecJ),>0})
 \bigr)
 \subset
 \nbigf^{\circ\theta^u}_{\gminib}
 H_1^{\rd}\bigl(
 \cnum^{\ast},
 (\nbigv,\nabla)\otimes\nbige(zu^{-1})
 \bigr).
\]
By the hypothesis of the induction,
Proposition \ref{prop;24.3.15.20} holds
$\nbigt_{\omega}(\nbigv,\nabla)$.
Then,
by Proposition \ref{prop;18.5.20.20}
together with
Lemma \ref{lem;20.10.23.1} and Lemma \ref{lem;20.10.27.20},
we obtain
\[
 C^{(\nu_0^-(\vecJ_{1})_-)}_{\infty,\theta^u}
 \Bigl(
 \nbigf^{\circ\theta^u}_{\gminib}
 H_1^{\rd}\bigl(\cnum^{\ast},
 \nbigt_{\omega}(\nbigv,\nabla)\otimes\nbige(zu^{-1})
 \bigr)
 \Bigr)
 \subset
 \nbigf^{\circ\theta^u}_{\gminib}
 H_1^{\rd}\bigl(\cnum^{\ast},
 (\nbigv,\nabla)\otimes\nbige(zu^{-1})\bigr).
\]
Hence, 
we obtain
$\nbigftilde^{\prime\theta^u}_{\gminib}
\subset
 \nbigftilde^{\circ\theta^u}_{\gminib}$
for any 
$\gminib\in\gbigf^{(0,\infty)}_+(\nbigi(\nbigv))$
under the isomorphism (\ref{eq;24.3.14.10}).
We obtain that
$\nbigftilde^{\prime\theta^u}_{\gminib}
=\nbigftilde^{\circ\theta^u}_{\gminib}$
for any $\gminib\in\gbigf^{(0,\infty)}_+(\nbigi(\nbigv))$,
because the dimension of the associated 
graded spaces of the filtrations are the same.
The claim of Proposition \ref{prop;24.3.15.20} also follows.
\hfill\qed

\subsection{Preliminary}
To simplify the notation,
we denote $\vecI(\theta^u)$ by $\vecI$
in the rest of \S\ref{section;20.11.21.1}.

For any 
$\gminia
=\sum_{0<j\leq\omega}
 \gminia_j r^{-j}e^{-\sqrt{-1}j \theta}
\in\nbigitilde$,
we set
$\veckappa(\gminia,u):=
 (\arg(\gminia_{\omega}),\theta_0^u)$
and
$\gminis(\gminia,u):=
 \bigl|\omega \gminia_{\omega} u\bigr|^{1/(1+\omega)}$.
We also set
\[
\gminic(\gminia,u)=
\bigl(
 \gminia_j\cdot
 \bigl|\omega\gminia_{\omega}\bigr|^{-(j+1)/(\omega+1)}
\cdot
 \bigl|u\bigr|^{(\omega-j)/(\omega+1)}
 \bigr)_{0<j<\omega}.
\]

Set
$F_{\gminia,u}(r,\theta)
:=\gminia(re^{\sqrt{-1}\theta})
+u^{-1}re^{\sqrt{-1}\theta}$.
We shall use the following rescaling:
\[
 F_{\gminia,u}\Bigl(
 \gminis(\gminia,u) r,
 \theta
 \Bigr)
=\bigl(\omega|\gminia_{\omega}|\bigr)^{1/(1+\omega)}
 |u|^{-\omega/(1+\omega)}
 F_{\veckappa(\gminia,u),\gminic(\gminia,u)}(r,\theta).
\]

We also remark the following,
which allows us to avoid the study of the critical cases.
\begin{lem}
\label{lem;18.5.30.10}
 The first claims of
 Proposition {\rm\ref{prop;18.5.20.30}}
 and Proposition {\rm\ref{prop;18.5.20.31}}
 imply the second claims of 
 Proposition {\rm\ref{prop;18.5.20.30}}
 and Proposition {\rm\ref{prop;18.5.20.31}}.
\end{lem}
\pf
Suppose that we have already proved
the first claims of
Proposition {\rm\ref{prop;18.5.20.30}}
and Proposition {\rm\ref{prop;18.5.20.31}}.

Let us prove the second claim of
Proposition {\rm\ref{prop;18.5.20.30}}
in the case
$J_+\cap(\vecI+\pi)_-=\{\vartheta^J_{r}\}$.
Take $u_1=|u|e^{\sqrt{-1}\theta^u_1}$
such that $|\theta^u_1-\theta^u|$ is sufficiently small.
Take $\gminia\in\nbigitilde_{J,<0}$
and $v\in H^0(J_+,L_{J_+,\gminia})$.
If $\theta^u_1-\theta^u<0$,
by the first claim of Proposition {\rm\ref{prop;18.5.20.30}},
$\Abb^{\rd}_{J_-,(\theta^u,\theta^u_1),\gminia}(v)$
is represented by a family of cycles
contained in
$\nbigc_0^{!}\bigl(
 (\nbigv,\nabla),u_1,d,-\Re(\gminia^{\circ}(u_1t))
 \bigr)$.
Let us consider the case $\theta^u_1-\theta^u>0$.
We set $\hat{J}=J+\omega^{-1}\pi$.
There exists
$\hat{v}\in H^0(\hat{J}_-,L_{\hat{J}_-,>0})$
such that
$v=\nbigrtilde^{\hat{J}_-}_{J}(\vhat)$.
By the formula (\ref{eq;24.3.16.1}), we obtain
\begin{multline}
 B_{\Jhat_-,\theta^u}(\vhat)
-\Abb^{\rd}_{J,\theta^u}(v)
=\sum_{\Jhat-\omega^{-1}\pi< J'<\Jhat}
 \Abb^{\rd}_{J',\theta^u}(\nbigrtilde^{\Jhat_-}_{J'}(\vhat))
 -\sum_{\Jhat\leq J'<\Jhat+\pi}
 \Abb^{\rd}_{J',\theta^u}(\nbigrtilde^{\Jhat_-}_{J'}(\vhat))
 \\
 -\sum_{\Jhat-\omega^{-1}\pi\leq J'<\Jhat-\pi}
 \Abb^{\rd}_{J'+2\pi,\theta^u}
 \bigl((\Tbb^{\ast})^{-1}\nbigrtilde^{\Jhat_-}_{J'}(\vhat)\bigr).
\end{multline}
Note that
$J'\cap (\vecI+\pi)\neq\emptyset$
for any $\Jhat-\omega^{-1}\pi<J'<\Jhat+\pi$.
We also note that 
$(J'+2\pi)\cap (\vecI+\pi)\neq\emptyset$
for any $\Jhat-\omega^{-1}\pi\leq J'<\Jhat-\pi$
in the case $\omega<1$.
There exist the expressions
\[
 \nbigr^{\Jhat_-}_{J'}(\vhat) 
 =\sum_{\gminib\in\nbigitilde_{J',<0}}
 \nbigr^{\Jhat_-}_{J'}(\vhat)_{\gminib},
 \quad
 (\Tbb^{\ast})^{-1}\nbigrtilde^{\Jhat_-}_{J'}(\vhat)
 =\sum_{\gminib\in \nbigitilde_{J'+2\pi,<0}}
 (\Tbb^{\ast})^{-1}\nbigrtilde^{\Jhat_-}_{J'}(\vhat)_{\gminib},
\]
where
$\nbigr^{\Jhat_-}_{J'}(\vhat)_{\gminib}
\in H^0\bigl(
 J'_-,L_{J'_-,\gminib}
 \bigr)$
and 
$(\Tbb^{\ast})^{-1}\nbigrtilde^{\Jhat_-}_{J'}(\vhat)_{\gminib}
\in H^0\bigl(
 (J'+2\pi)_-,L_{(J'+2\pi)_-,\gminib}
 \bigr)$.
Note that 
$-\Re\gminib^{\circ}_-(u_1t)t^{(1+\omega)^{-1}\omega}$
is convergent to a negative number as $t\to 0$
for any $\gminib\in\nbigitilde_{J',<0}$
or $\gminib\in\nbigitilde_{J'+2\pi,<0}$.
There exists the expression
\[
 \vhat
 =\sum_{\substack{\gminib\in\nbigitilde_{\Jhat,>0}}}
 \vhat_{\gminib},
\]
where $\vhat_{\gminib}\in H^0(\Jhat_-,L_{\Jhat_-,\gminib})$.
We have
$\vhat_{\gminib}=0$ unless $\gminib\leq_{\vartheta^J_r}\gminia$.
For any $\gminib\in\nbigitilde_{\Jhat,>0}$,
\[
-\Re \gbigf^{(0,\infty)}_{J,0,+}(\gminib)(u_1t)t^{(1+\omega)^{-1}\omega}
\]
is convergent to a positive number $C(\gminib)$ as $t\to 0$.
If $\gminib\leq_{\vartheta^J_r}\gminia$ and $\gminib\neq\gminia$,
we have
$C(\gminib)<C(\gminia)$.
We also note that
$\gbigf^{(0,\infty)}_{\Jhat,0,+}(\gminia)
=\gbigf^{(0,\infty)}_{J,0,-}(\gminia)$.
Therefore, by the first claims of
Proposition {\rm\ref{prop;18.5.20.30}}
and Proposition {\rm\ref{prop;18.5.20.31}},
we obtain the first claim of
Proposition {\rm\ref{prop;18.5.20.30}}
in the case where $J_+\cap(\vecI+\pi)$
consists of one point.
By the same argument,
we can prove the second claim of Proposition {\rm\ref{prop;18.5.20.31}}
in the case where
$J_-\cap\vecI_+$ consists of one point.
We can prove the other cases of
Proposition {\rm\ref{prop;18.5.20.30}}
and Proposition {\rm\ref{prop;18.5.20.31}}
by using the formula (\ref{eq;24.3.16.2}).
\hfill\qed

\subsection{Proof of
the first claim of Proposition \ref{prop;18.5.20.30}}

We take $\gminia\in\nbigitilde_{J,<0}$.
Let us study the claim for
$\Abb^{\rd}_{J_-,\theta^u,\gminia}$
in the case $J\cap(\vecI+\pi)\neq\emptyset$.
The claim for $\Abb^{\rd}_{J_+,\theta^u,\gminia}$
can be argued similarly.

There exists $\theta_1\in J\cap(\vecI+\pi)$
such that
$\theta_1\in \Cr_1(\omega,\veckappa(\gminia,u))$.
Let $\Gamma_{\theta_1}$ be the path
on $\projtilde^1$
defined by
$s\longmapsto (s(1-s)^{-1},\theta_1)$ $(0\leq s\leq 1)$.
For a sufficiently small $t_0>0$,
we construct a continuous family of paths
$\Gamma_{\theta_1,\gminic(\gminia,tu)}$
$(0\leq t\leq t_0)$
for $F_{\veckappa(\gminia,u),\gminic(\gminia,tu)}$
and $\theta_1$
by modifying $\Gamma_{\theta_1}$
as in \S\ref{subsection;18.5.20.1}.
Any element 
$v\in H^0(J_-,L_{J_-,\gminia})$
naturally induces a flat section
$\vtilde$ of $\nbigv$
on a sector which contains
$\Gamma_{\theta_1,\gminic(\gminia,tu)}$.
We obtain the following family of rapid decay $1$-cycles
for $(\nbigv,\nabla)\otimes\nbige(zu^{-1})$,
which represents $\Abb^{\rd}_{J_-,u}(v)$:
\begin{equation}
 \label{eq;18.5.28.10}
  \Bigl(
\vtilde
\cdot\exp(-zu^{-1})
\Bigr)
\otimes
\Bigl(
 \gminis(\gminia,ut)
 \cdot
 \Gamma_{\theta_1,\gminic(\gminia,tu)}
 \Bigr).
\end{equation}
Here,
$\gminis(\gminia,ut)\Gamma_{\theta_1,\gminic(\gminia,tu)}$
is a family of paths in $\cnum^{\ast}$
obtained as
the multiplication of $\gminis(\gminia,ut)$
to $\Gamma_{\theta_1,\gminic(\gminia,tu)}$.

\begin{lem}
\label{lem;18.5.28.11}
We can divide $\Gamma_{\theta_1,\gminic(\gminia,tu)}$
into a sum of $1$-chains such that
the family of cycles {\rm(\ref{eq;18.5.28.10})}
is contained in
\[
 \nbigc_0^{!}\bigl(
 (\nbigv,\nabla),
 u,d(\omega),-\Re(\gminia^{\circ}(ut))
 \bigr).
\]
\end{lem}
\pf
We may naturally regard $\Gamma_{\theta_1,\gminic(\gminia,tu)}$
as a path on $X$.
(See \S\ref{section;18.6.3.20} for $X$.)
By the estimates in \S\ref{subsection;18.5.20.1},
there exist $C,N>0$
such that 
$|\vtilde|_{h^{\nbigv}}
 \exp\bigl(-\Re(zu^{-1}t^{-1})\bigr)
\leq
 C\exp\bigl( 
-\Re\gminia^{\circ}(ut)
\bigr)t^{-N}$
along $\gminis(\gminia,ut)\Gamma_{\theta_1,\gminic(\gminia,tu)}$.
Moreover, for any $\epsilon>0$,
there exist a neighbourhood 
$U_{\epsilon}$ of $(1,\theta_1)$
in $X$
such that the following holds
on $\gminis(\gminia,ut)\cdot 
\bigl(
 \Gamma_{\theta_1,\gminic(\gminia,ut)}
 \setminus
 U_{\epsilon}
 \bigr)$
for some $C_{i,\epsilon}>0$ $(i=1,2)$:
\begin{multline}
 |\vtilde|_{h^{\nbigv}}
 \exp\bigl(-\Re(zu^{-1}t^{-1})\bigr)
 =\\
 O\left(
 \exp\Bigl[
 -\Re\gminia^{\circ}(ut)
 -t^{-\omega/(1+\omega)}
\bigl(
 \epsilon+C_{1,\epsilon}r+C_{2,\epsilon}r^{-\omega}
 \bigr)
 \Bigr]
 \right).
\end{multline}
Then, we obtain the claim of the lemma.
\hfill\qed

\vspace{.1in}
We immediately obtain the first claim of the proposition
from the lemma.

\subsection{Proof of
the first claim of Proposition \ref{prop;18.5.20.31}
in the case $\omega>1$}
\label{subsection;18.5.28.100}

Take $\gminia\in \nbigitilde_{J,>0}$.
Let us study the claim for $B_{J_-,\theta^u,\gminia}$
in the case $J\cap \vecI\neq\emptyset$.
The claim for $B_{J_+,\theta^u,\gminia}$ can be argued similarly.

There exists
$\theta_1\in
 \Cr_1(\omega;\veckappa(\gminia,u))$
such that
$\theta_1\in J_-\cap\vecI_+$.
Take a small $\delta>0$.
Let $\Gamma_{\theta_1}$ be the path
obtained as 
$\gamma_h(1;\vartheta^J_{\ell}-\delta,\vartheta^J_r)$
for $F_{\veckappa(\gminia,u)}$.
For a small $t_0>0$,
we construct a continuous family of paths
$\Gamma_{\theta_1.\gminic(\gminia,ut)}$
$(0\leq t\leq t_0)$
for $F_{\veckappa(\gminia,u),\gminic(\gminia,ut)}$
by modifying $\Gamma_{\theta_1}$
as in \S\ref{subsection;18.5.20.1}.
By adding the segment
$\gamma_v(1,0;\vartheta^J_{\ell}-\delta)$
to $\Gamma_{\theta_1,\gminic(\gminia,ut)}$,
we obtain a continuous family of paths
$\Gammatilde_{\theta_1,\gminic(\gminia,ut)}$
$(0\leq t\leq t_0)$
connecting
$(0,\vartheta^J_{\ell}-\delta)$
and $(1,\vartheta^J_r)$.

Any element
$v\in H^0(J_-,L_{J_-,\gminia})$ 
induces a flat section $\vtilde$ of $\nbigl$
on $\{(r,\theta)\,|\,0\leq r\leq\infty,\,\,
 \theta\in \openclosed{\vartheta^J_{\ell}-\delta}{\vartheta^J_r} \}$.
At $\vartheta^J_{r}$,
we have the decomposition
$v=u_{J_+,0}+\sum_{J\leq J'\leq J+\omega^{-1}\pi}
 u_{J'}$,
where $u_{J_+,0}$ is a section of
$L_{J_+,0}$
and $u_{J'}$ are sections of $L_{J',<0}$.

We take a sufficiently small $\delta>0$.
Let $\Gamma_{2,\pm}$ be the path
connecting
$(1,\vartheta^J_r)$
and $(0,\vartheta^J_r\pm\delta)$,
obtained as the union of 
$\gamma_h(1;\vartheta^J_{r},\vartheta^J_r\pm\delta)$
and 
$\gamma_v(1,0;\vartheta^J_r\pm\delta)$.
We have the flat sections 
$\utilde_{J'}$ induced by $u_{J'}$
along $\Gamma_{2,+}$
if $J<J'\leq J+\omega^{-1}\pi$.
We have the section
$\utilde_{J}$ induced by $u_J$
on the sector which contains
$\Gamma_{2,-}$ if $J'=J$.

Let $a_0<a_1<\cdots<a_N$
be the intersection of 
$S_0(\nbigi)\cap
 \closedopen{\vartheta^J_{r}}{\vartheta^J_r+\pi}$.
We take 
$a_{N+1}\in\openopen{a_N}{\vartheta^J_r+\pi}$.
We set $J_i:=\openopen{a_i-\omega^{-1}\pi}{a_i}$.
We obtain the sections 
$u_{J_{i+},0}\in H^0(J_{i+},L_{J_{i+},0})$
induced by $u_{J_+,0}$ and the parallel transport of
$\Gr^{\vecnbigf}_0(L)$.

We set
$\omega'=
\max\{-\ord(\gminia)\,|\,\gminia\in\nbigt_{\omega}(\nbigi(\nbigv))\}
<\omega$.
Let $\beta>0$ be sufficiently small that
\begin{equation}
\label{eq;24.3.14.1}
 \omega'
 \left(
  \frac{1}{1+\omega}+\beta
  \right)
 <\frac{\omega}{1+\omega}.
\end{equation}
Let $\Gamma_{3}$ be the path
$\gamma_v(1,t^{\beta};\vartheta^J_r)$.
Let $\Gamma_{4,i}$ $(i=0,\ldots,N)$ be the paths
$\gamma_h(t^{\beta};a_i,a_{i+1})$.
Let $\Gamma_{5}$ be the paths
$\gamma_v(\infty,t^{\beta};a_{N+1})$.
Let $\Gamma_{6,i}$ be the paths
$\gamma_v(t^{\beta},0;a_i)$.

We obtain the following continuous family of cycles
which represents $B_{J_-,\theta^u,\gminia}(v)$:
\begin{multline}
\label{eq;18.5.20.40}
 \Bigl[
 \vtilde
 \otimes
 \bigl(
 \gminis(\gminia,ut)
 \Gammatilde_{\theta_1,\gminic(\gminia,ut)}
 \bigr)
+\utilde_J
 \otimes
 \bigl(
 \gminis(\gminia,ut)\Gamma_{2,-}
 \bigr)
+\sum_{J<J'\leq J+\pi/\omega}
 \utilde_{J'}
 \otimes
 \bigl(
  \gminis(\gminia,ut)
 \Gamma_{2,+}
 \bigr)
 \Bigr.
 \\
 \Bigl.
+u_{J_+,0}
 \otimes
 \bigl(
 \gminis(\gminia,ut)\Gamma_{3}
 \bigr)
+\sum_{i=0}^N
 u_{J_{i+},0}
 \otimes
 \bigl(
 \gminis(\gminia,ut)\Gamma_{4,i}
 \bigr)
+u_{J_{N_+},0}
 \otimes  
 \gminis(\gminia,ut) 
 \Gamma_5
 \Bigr.\\
 \Bigl.
+\sum_{i=1}^N
 \bigl(u_{J_{i+},0}-u_{(J_{i-1})_+,0}\bigr)
\otimes\bigl(
 \gminis(\gminia,ut)\Gamma_{6,i}
 \bigr)
 \Bigr]
 \exp(-zu^{-1}).
\end{multline}

\begin{lem}
\label{lem;18.5.28.30}
We can divide paths into sums of $1$-chains
such that
the family {\rm(\ref{eq;18.5.20.40})}
is contained in 
\[
 \nbigc_0^{!}\bigl(
 (\nbigv,\nabla),u,d(\omega),
 -\Re(\gminia^{\circ}(ut))
 \bigr).
\]
\end{lem}
\pf
By the estimates in \S\ref{subsection;18.5.20.1},
there exist $C,N>0$
such that
\[
\bigl|
 \vtilde
 \bigr|_{h^{\nbigv}}
 \exp\bigl(-\Re(zu^{-1}t^{-1})\bigr)
\leq
 C\exp(-\Re\gminia^{\circ}(tu)) t^{-N}
\]
on 
$\gminis(\gminia,ut)\cdot
 \Gamma_{\theta_1,\gminic(\gminia,ut)}$.
 Moreover, for any $\epsilon>0$,
there exists
a neighbourhood $U_{\epsilon}$
of $(1,\theta_1)$ such that
the following holds
on $\gminis(\gminia,ut)\cdot
\bigl(
 \Gamma_{\gminia,\gminic(\gminia,ut)}
\setminus
 U_{\epsilon}
\bigr)$:
\begin{equation}
 \bigl|
 \vtilde
 \bigr|_{h^{\nbigv}}
 \exp\bigl(-\Re(zu^{-1}t^{-1})\bigr)
=O\Bigl(
 \exp\Bigl(-\Re(\gminia^{\circ}(tu))
 -\epsilon t^{-\omega/(1+\omega)}
 \Bigr)
 \Bigr).
\end{equation}

In the following, $C_1$ and $\epsilon_1$
denote positive constants,
which can vary.
Note that 
$-\Re(\gminia^{\circ}(ut))t^{\omega/(1+\omega)}$
is convergent to a positive number as $t\to 0$.
We shall use the estimates in \S\ref{subsection;18.5.20.1}.
On $\gminis(\gminia,ut)\cdot
\gamma_v(1,0;\vartheta^J_{\ell}-\delta)$,
by using Lemma \ref{lem;20.10.23.10},
we obtain
\[
 \bigl|
 \vtilde
 \bigr|_{h^{\nbigv}}
 \exp\bigl(-\Re(zu^{-1}t^{-1})\bigr)
=O\Bigl(
 \exp\Bigl(
 -\Re(\gminia^{\circ}(tu))
 -\epsilon_1 t^{-\omega/(1+\omega)}
 -C_1\gminis(\gminia,ut)^{\omega}r^{-\omega}
 \Bigr)
 \Bigr).
\]
Similarly,
on $\gminis(\gminia,ut)\cdot \Gamma_{2,-}$,
by using Lemma \ref{lem;20.10.23.10},
we obtain
\[
 \bigl|
 \utilde_J
 \bigr|_{h^{\nbigv}}\cdot
 \exp(-\Re(zu^{-1}t^{-1}))
 \!=\!
 O\Bigl(
  \exp\Bigl(-\Re(\gminia^{\circ}_+(tu))
 -\epsilon_1 t^{-\omega/(1+\omega)}
 -C_1\gminis(\gminia,ut)^{\omega}r^{-\omega}
 \Bigr)
 \Bigr).
\]
We obtain similar estimates for
$\bigl|
 \utilde_{J'}
 \bigr|_{h^{\nbigv}}
  \exp(-\Re(zu^{-1}t^{-1}))$
on $\gminis(\gminia,ut)\cdot \Gamma_{2,+}$.
\begin{lem}
We have the following estimates
for some $\epsilon_1>0$
on $\gminis(\gminia,ut)\Gamma_3$:
\begin{equation}
\label{eq;18.5.28.20}
 \bigl|
 u_{J_+,0}
 \bigr|_{h^{\nbigv}}
 \exp(-\Re(zu^{-1}t^{-1}))
=O\Bigl(
 \exp\Bigl(
 -\Re\gminia^{\circ}(ut)
-\epsilon_1t^{-\omega/(1+\omega)}
 \Bigr)
\Bigr).
\end{equation}
\end{lem}
\pf
If $-\Re(zu^{-1})<0$ on along $\arg(z)=\vartheta^J_r$,
the claim is clear.
If $-\Re(zu^{-1})>0$,
then $\exp\bigl(-\Re(zu^{-1}t^{-1})\bigr)$ is 
monotonously increasing with respect to $|z|$.
We also have the following for
$t^{\beta}\gminis(\gminia,ut)\leq |z|\leq\gminis(\gminia,ut)$:
\[
\log\left(
 \frac
{\bigl(
 |u_{J_+,0}|_{h^{\nbigv}}\bigr)_{|(|z|,\vartheta^J_r)}
 }
 {\bigl(
 |u_{J_+,0}|_{h^{\nbigv}}\bigr)_{|(\gminis(\gminia,ut),\vartheta^J_r)}
  }
\right)
=O\bigl(
t^{-\omega'((1+\omega)^{-1}+\beta)}
\bigr).
\]
The following holds on $\gminis(\gminia,ut)\Gamma_3$:
\begin{multline}
\log\Bigl(
 \bigl|
 u_{J_+,0}
 \bigr|_{h^{\nbigv}}
 \exp(-\Re(zu^{-1}t^{-1}))
 \Bigr)
 =\\
\log\Bigl(
 \bigl|
 u_{J_+,0}
 \bigr|_{h^{\nbigv}}
 \exp(-\Re(zu^{-1}t^{-1}))
 _{|(\gminis(\gminia,ut),\vartheta^J_r)}
 \Bigr)
+
 O\bigl(
 t^{-\omega'((1+\omega)^{-1}+\beta))}
 \bigr).
\end{multline}
Because the estimate (\ref{eq;18.5.28.20}) holds at 
$(\gminis(\gminia,ut),\vartheta^J_{r})$,
we obtain (\ref{eq;18.5.28.20}) on 
$\gminis(\gminia,ut)\Gamma_3$.
\hfill\qed

\vspace{.1in}

On $\gminis(\gminia,ut)\cdot\Gamma_{4,i}$,
we have
$\Re(zu^{-1}t^{-1})
=
O\bigl(
t^{-\omega/(1+\omega)+\beta}
\bigr)$
and
\[
 \log|u_{J_{i+},0}|_{h^{\nbigv}}
=O\bigl(
 t^{-\omega'((1+\omega)^{-1}+\beta)}
\bigr).
\]
Hence, we obtain the following
because $z(ut)^{-1}$ is bounded on
$\gminis(\gminia,ut)\cdot\Gamma_{4,i}$:
\begin{equation}
\label{eq;18.5.28.21}
 \bigl|
 u_{J_{i+},0}
 \bigr|_{h^{\nbigv}}
  \exp(-\Re(zu^{-1}t^{-1}))
=O\Bigl(
 \exp\Bigl(
 -\Re\gminia^{\circ}(ut)
-\epsilon_1 t^{-\omega/(1+\omega)}
 \Bigr)
 \Bigr).
\end{equation}
On $\gminis(\gminia,ut)\cdot\Gamma_5$,
we obtain the following estimate
because $\Re(zu^{-1})<0$ around $\Gamma_5$:
\[
 \bigl|
 u_{J_{N+},0}
 \bigr|_{h^{\nbigv}}
 \exp(-\Re(zu^{-1}t^{-1}))
=O\Bigl(
  \exp\Bigl(
 -\Re\gminia^{\circ}(ut)
-\epsilon_1 t^{-\omega/(1+\omega)}
-C_1|z|
 \Bigr)
 \Bigr).
\]
On $\gminis(\gminia,ut)\cdot\Gamma_{6,i}$,
we obtain the following estimate
because
$\Re(z(ut)^{-1})=O(t^{-\omega'(1+\omega)^{-1}+\beta})$
on $\gminis(\gminia,ut)\cdot\Gamma_{6,i}$:
\begin{multline}
 \bigl|
 u_{J_{i+},0}-u_{(J_{i-1})_+,0}
 \bigr|_{h^{\nbigv}}
 \exp\bigl(-\Re(zu^{-1}t^{-1})\bigr)
 =\\
 O\Bigl(
   \exp\Bigl(
 -\Re\gminia^{\circ}(ut)
-\epsilon_1 t^{-\omega/(1+\omega)}
-C_1
 |z|^{-\omega}
 \Bigr)
 \Bigr).
\end{multline}
Then, we obtain the claim of 
Lemma \ref{lem;18.5.28.30}.
\hfill\qed

\vspace{.1in}
Thus, we obtain the claim of Proposition \ref{prop;18.5.20.31}
in the case $\omega>1$.

\subsection{Proof of the first claim of
 Proposition  \ref{prop;18.5.20.31}
 in the case $\omega<1$}

Take $\gminia\in\nbigitilde_{J,>0}$.
We shall explain the proof for $B_{J_-,\theta^u,\gminia}$
in the case $J\cap\vecI\neq\emptyset$.
The proof for $B_{J_+,\theta^u,\gminia}$ is similar.
There exists
$\theta_1\in J\cap \vecI$
such that
$\theta_1\in \Cr_1(\omega;\veckappa(\gminia,u))$.

For $v\in H^0(J_-,L_{J_-,>0})$,
we have the expression
\[
 v=\sum_{J-\omega^{-1}\pi\leq J'<J} u_{J'},
\]
where $u_{J'}$ are sections of
$L_{J',<0}$.
We also have the expression
\[
 v=u_{J_+,0}
+\sum_{J\leq J'\leq J+\omega^{-1}\pi}
 u_{J'},
\]
where 
$u_{J_+,0}$ is a section of $L_{J_+,0}$,
and $u_{J'}$ are sections of $L_{J',<0}$.
For $J<J'<J+\pi$,
we obtain the sections
$u_{J'_{\pm},0}$ of $L_{J'_{\pm},0}$
induced by $u_{J_{+},0}$
and the parallel transport of
$\Gr^{\vecnbigf}_0(L)$.
We obtain the sections 
$u_{J'_+,0}-u_{J'_-,0}$
of $L_{J'_+,<0}$.

\vspace{.1in}

\subsubsection{The case $\vecI\subset J$}
Let us consider the case $\vecI\subset J$.
Take a small $\delta>0$.
Let $\Gamma_{\theta_1}$ be the path
for $F_{\veckappa(\gminia)}$
obtained as
$\gamma_h(1;\vartheta^{\vecI}_{\ell}-\delta,\vartheta^{\vecI}_r+\delta)$.
We modify it to
$\Gamma_{\theta_1,\gminic(\gminia,ut)}$
as in \S\ref{subsection;18.5.20.1}.
By adding 
$\gamma_v(\infty,1;\vartheta^{\vecI}_{\ell}-\delta)$
and 
$\gamma_v(\infty,1;\vartheta^{\vecI}_{r}+\delta)$,
we obtain a family of paths
$\Gammatilde_{\theta_1,\gminic(\gminia,ut)}$
connecting
$(\infty,\vartheta^{\vecI}_{\ell}-\delta)$
and 
$(\infty,\vartheta^{\vecI}_{r}+\delta)$.

For any $J'\in T(\nbigi)$
such that $J-\pi\leq J'<J$,
we have $J'\cap(\vecI-\pi)\neq\emptyset$.
We take $\theta_{J'}\in J'\cap(\vecI-\pi)$.
Let $\Gamma_{J'}$ be the path 
$\gamma_v(\infty,0;\theta_{J'})$.

For any $J'\in T(\nbigi)$
such that $J<J'<J+\pi$,
we have $J'\cap(\vecI+\pi)\neq\emptyset$.
We take $\theta_{J'}\in J'\cap(\vecI+\pi)$.
Let $\Gamma_{J'}$ be the path 
$\gamma_v(\infty,0;\theta_{J'})$.

If $J\cap(\vecI+\pi)\neq\emptyset$,
we take $\theta_J\in J\cap(\vecI+\pi)$,
and let $\Gamma_J$ be the path
$\gamma_v(\infty,0;\theta_{J})$.
If $J\cap(\vecI+\pi)=\emptyset$,
we have
$J_+\cap(\vecI+\pi)_-=\{\vartheta^J_r\}$.
We take a sufficiently small $\delta>0$,
and we consider the path
$\Gamma_{J}$ connecting
$(0,\vartheta^J_r-\delta)$
and $(\infty,\vartheta^J_{r}+\delta)$
obtained as the union of
$\gamma_v(1,0;\vartheta^J_r-\delta)$,
$\gamma_h(1;\vartheta^J_r-\delta,\vartheta^J_r+\delta)$
and 
$\gamma_v(\infty,1;\vartheta^J_r+\delta)$.

Let $\vtilde$ denote the flat section induced by $v$.
We obtain the following family of cycles
for $(\nbigv,\nabla)\otimes\nbige(zu^{-1})$
which represents $B_{J_-,\theta^u,\gminia}(v)$:
\begin{multline}
\label{eq;18.5.29.1}
 \Bigl[
\vtilde\otimes
\gminis(\gminia,ut)
\Gammatilde_{\theta_1,\gminic(\gminia,ut)}
+\sum_{J-\pi\leq J'<J}
 u_{J'}\otimes
\gminis(\gminia,ut)\Gamma_{J'}
 \\
-u_{J}\otimes
\gminis(\gminia,ut)\Gamma_{J}
-\sum_{J<J'<J+\pi}
 \bigl(
  u_{J'}+(u_{J_+',0}-u_{J'_-,0})
 \bigr)
\otimes
\gminis(\gminia,ut)\Gamma_{J'}
\Bigr]\exp(-zu^{-1}).
\end{multline}

\begin{lem}
\label{lem;18.5.29.4}
We can divide the paths into sums of $1$-chains
such that
the family {\rm(\ref{eq;18.5.29.1})}
is contained in
\[
 \nbigc_0^{!}\bigl(
 (\nbigv,\nabla),u,d(\omega),
-\Re(\gminia^{\circ}(ut))
 \bigr).
\]
\end{lem}
\pf
On $\gminis(\gminia,ut)\Gamma_{\theta_1,\gminic(\gminia,ut)}$,
we obtain the following for some $N>0$
by using the estimates in \S\ref{subsection;18.5.20.1}:
\[
 \bigl|
 \vtilde
 \bigr|_{h^{\nbigv}}
 \exp\bigl(
 -\Re(zu^{-1}t^{-1})
 \bigr)
=O\Bigl(
 \exp\Bigl(
 -\Re(\gminia^{\circ}(ut))
 \Bigr) t^{-N}
 \Bigr).
\]
 Moreover, for any $\epsilon>0$,
there exists
a neighbourhood $U_{\epsilon}$
of $(1,\theta_1)$ such that
the following holds
on $\gminis(\gminia,ut)\cdot
\bigl(
 \Gamma_{\gminia,\gminic(\gminia,ut)}
\setminus
 U_{\epsilon}
\bigr)$:
\begin{equation}
 \bigl|
 \vtilde
 \bigr|_{h^{\nbigv}}
 \exp\bigl(-\Re(zu^{-1}t^{-1})\bigr)
=O\Bigl(
 \exp\Bigl(-\Re(\gminia^{\circ}(tu))
 -\epsilon t^{-\omega/(1+\omega)}
 \Bigr)
 \Bigr).
\end{equation}

In the following,
$C_1$ and $\epsilon_1$ denote positive constants,
which can vary.
Note that 
$-\Re(\gminia^{\circ}(ut))t^{\omega/(1+\omega)}$
is convergent to a positive number as $t\to 0$.

On $\gminis(\gminia,ut)
\gamma_v(\infty,1;\vartheta^J_{\ell}-\delta)$
and
$\gminis(\gminia,ut)
\gamma_v(\infty,1;\vartheta^J_{r}+\delta)$,
we obtain
\[
  \bigl|
 \vtilde
 \bigr|_{h^{\nbigv}}
 \exp\bigl(
 -\Re(zu^{-1}t^{-1})
 \bigr)
=O\Bigl(
 \exp\Bigl(
 -\Re(\gminia^{\circ}(ut))
-\epsilon_1 t^{-\omega/(1+\omega)}
-C_1 \gminis(\gminia,ut)^{-1}|z|
 \Bigr)
 \Bigr).
\]
Hence, we obtain the following estimates on $\Gamma_{J'}$
(we use Lemma \ref{lem;20.10.23.10}
in the case $\Gamma'=\Gamma$.):
\begin{multline}
 \bigl|
 u_{J'}
 \bigr|_{h^{\nbigv}}
 \exp\bigl(
 -\Re(zu^{-1}t^{-1})
 \bigr)
=\\
 O\Bigl(
  \exp\Bigl(
 -\Re(\gminia^{\circ}(ut))
-\epsilon_1 t^{-\omega/(1+\omega)}
-C_1
 \bigl(
 \gminis(\gminia,ut)^{-1}|z|
+\gminis(\gminia,ut)^{\omega}|z|^{-\omega}
\bigr)
 \Bigr)
 \Bigr).
\end{multline}
\begin{multline}
 \bigl|
 u_{J'_+,0}-u_{J'_-,0}
 \bigr|_{h^{\nbigv}}
 \exp\bigl(
 -\Re(zu^{-1}t^{-1})
 \bigr)
=\\
 O\Bigl(
  \exp\Bigl(
 -\Re(\gminia^{\circ}(ut))
-\epsilon_1 t^{-\omega/(1+\omega)}
-C_1 
\bigl(
 \gminis(\gminia,ut)^{-1}|z|
+\gminis(\gminia,ut)^{\omega}|z|^{-\omega}
\bigr)
 \Bigr)
 \Bigr).
\end{multline}
Then, we obtain the claim of the lemma.
\hfill\qed

\subsubsection{The case 
$\vartheta^J_r<\vartheta^{\vecI}_{r}<\vartheta^J_r+\pi$}

Let us study the case
$\vartheta^J_r<\vartheta^{\vecI}_{r}<\vartheta^J_r+\pi$.
Let $\Gamma_{\theta_1}$ be the path
$\gamma_h(1;\vartheta^{\vecI}_{\ell}-\delta,\vartheta^J_r)$.
We modify it to $\Gamma_{\theta_1,\gminic(\gminia,ut)}$
as in \S\ref{subsection;18.5.20.1}.
By adding
$\gamma_v(\infty,1;\vartheta^{\vecI}_{\ell}-\delta)$,
we obtain a family of paths
$\Gammatilde_{\theta_1,\gminic(\gminia,ut)}$
connecting
$(\infty,\vartheta^{\vecI}_{\ell}-\delta)$
and $(1,\vartheta^J_r)$.

Let $\Gamma_{2,\pm}$ be the paths
connecting $(1,\vartheta^J_r)$
and $(0,\vartheta^J_r\pm\delta)$
obtained as the union of
$\gamma_h(1;\vartheta^J_r,\vartheta^J_r\pm\delta)$
and $\gamma_v(1,0;\vartheta^J_r\pm\delta)$.

Let $a_0<a_1<\cdots<a_N$ be the intersection of
$S_0(\nbigi)\cap\closedopen{\vartheta^J_r}{\vartheta^J_r+\pi}$.
We take
$a_{N+1}\in \openopen{a_N}{\vartheta^{J}_r+\pi}$.
We set $J_i:=\openopen{a_i-\omega^{-1}\pi}{a_i}$.
We obtain the sections 
$u_{J_{i+},0}\in H^0(J_{i+},L_{J_{i+},0})$
induced by $u_{J_+,0}$ and the parallel transport of
$\Gr^{\vecnbigf}_0(L)$.

Let $\beta>0$ be sufficiently
satisfying (\ref{eq;24.3.14.1}).
Let $\Gamma_{3}$ be the path
$\gamma_v(1,t^{\beta};\vartheta^J_r)$.
Let $\Gamma_{4,i}$ $(i=0,\ldots,N)$ be the paths
$\gamma_h(t^{\beta};a_i,a_{i+1})$.
Let $\Gamma_{5}$ be the paths
$\gamma_v(\infty,t^{1-d(\omega)};a_{N+1})$.
Let $\Gamma_{6,i}$ be the paths
$\gamma_v(t^{\beta},0;a_i)$.

We obtain the following continuous family of cycles:
\begin{multline}
\label{eq;18.5.29.2}
 \Bigl[
\vtilde\otimes
 \gminis(\gminia,ut)\Gammatilde_{\theta_1,\gminic(\gminia,ut)}
+u_J\otimes
 \gminis(\gminia,ut)\Gamma_{2,-}
+\sum_{J<J'\leq J+\omega^{-1}\pi}
 u_{J'}\otimes
 \gminis(\gminia,ut)\Gamma_{2,+}
 \\
+u_{J_+,0}\otimes
 \gminis(\gminia,ut)\Gamma_{3}
+\sum_{i=0}^Nu_{J_{i+},0}\otimes
 \gminis(\gminia,ut)\Gamma_{4,i}
+u_{J_N,0}\otimes
 \gminis(\gminia,ut)\Gamma_5
 \\
+\sum_{i=1}^N
 (u_{J_{i+},0}-u_{(J_{i-1})_+,0})\otimes
 \gminis(\gminia,ut)
 \Gamma_{6,i}
\Bigr]
 \exp(-zu^{-1}).
\end{multline}

For $J'$ such that $J-\pi\leq J'<J$,
we have $J'\cap(\vecI-\pi)\neq\emptyset$.
We take $\theta_{J'}\in J'\cap(\vecI-\pi)$,
and we set
$\Gamma_{J'}:=\gamma_v(\infty,0;\theta_{J'})$.
For $J'$ such that
$J+\pi\leq J'<J+\omega^{-1}\pi$,
we have 
$J'\cap(\vecI+\pi)\neq\emptyset$.
We take $\theta_{J'}\in J'\cap(\vecI+\pi)$,
and we set $\Gamma_{J'}:=\gamma_v(\infty,0;\theta_{J'})$.
We obtain the following family of cycles:
\begin{equation}
\label{eq;18.5.29.3}
\Bigl(
\sum_{J-\pi\leq J'<J}
 u_{J'}\otimes
 \gminis(\gminia,ut)
 \Gamma_{J'}
+\sum_{J+\pi\leq J'<J+\omega^{-1}\pi}
 u_{J'}\otimes
 \gminis(\gminia,ut)
 \Gamma_{J'}
 \Bigr)\exp(-zu^{-1}).
\end{equation}
The sum of (\ref{eq;18.5.29.2})
and (\ref{eq;18.5.29.3}) represents $B_{J_-,\theta^u,\gminia}(v)$.

\begin{lem}
We can divide the paths into sums of $1$-chains such that
the families {\rm(\ref{eq;18.5.29.2})}
and {\rm(\ref{eq;18.5.29.3})}
are contained in
\[
 \nbigc_0^{!}\bigl(
(\nbigv,\nabla),u,d(\omega),
-\Re(\gminia^{\circ}(ut))
 \bigr).
\]
\end{lem}
\pf
The estimate for 
$\vtilde\exp(-zu^{-1}t^{-1})
 \otimes
 \gminis(\gminia,ut)
 \Gammatilde_{\theta_1,\gminic(\gminia,ut)}$
is similar to that in Lemma \ref{lem;18.5.29.4}.
The contributions of
the other terms in (\ref{eq;18.5.29.2})
are dominated as in the case of
Lemma \ref{lem;18.5.28.30},
and the terms in (\ref{eq;18.5.29.3})
are dominated as in the case of
Lemma \ref{lem;18.5.29.4}.
\hfill\qed

\subsubsection{The case
$\vartheta^J_{\ell}-\pi<
 \vartheta^{\vecI}_{\ell}<
 \vartheta^J_{\ell}$}
Let us study the case
$\vartheta^J_{\ell}-\pi<
 \vartheta^{\vecI}_{\ell}<
 \vartheta^J_{\ell}$.
Let $\Gamma_{\theta_1}$ be the path
$\gamma_h(1;\vartheta^J_{\ell}-\delta,\vartheta^{\vecI}_r+\delta)$.
We modify it to
$\Gamma_{\theta_1,\gminic(\gminia,ut)}$
as in \S\ref{subsection;18.5.20.1}.
By adding 
$\gamma_v(1,0;\vartheta^J_{\ell}-\delta)$
and 
$\gamma_v(\infty,1;\vartheta^{\vecI}_r+\delta)$,
we obtain the path
$\Gammatilde_{\theta_1,\gminic(\gminia,ut)}$
connecting
$(0,\vartheta^J_{\ell}-\delta)$
and $(\infty,\vartheta^{\vecI}_{r}+\delta)$.

For $J'$ such that
$J-\omega^{-1}\pi\leq J'<J-\pi$,
we have 
$J'\cap(\vecI-\pi)\neq\emptyset$.
We take $\theta_{J'}\in J'\cap(\vecI-\pi)$,
and we put
$\Gamma_{J'}:=\gamma_v(\infty,0;\theta_{J'})$.

For $J'$ such that
$J\leq J'<J+\pi$,
we have
$J'\cap (\vecI+\pi)\neq\emptyset$.
We take $\theta_{J'}\in J'\cap(\vecI+\pi)$,
and we put
$\Gamma_{J'}:=\gamma_v(\infty,0;\theta_{J'})$.

Then, the following family of cycles
represents $B_{J_-,\theta^u,\gminia}(v)$:
\begin{multline}
\label{eq;18.5.29.5}
 \Bigl[
 \vtilde\otimes
 \gminis(\gminia,ut)\Gammatilde_{\theta_1,\gminic(\gminia,ut)}
+\sum_{J-\omega^{-1}\pi\leq J'<J}
 u_{J'}\otimes
 \gminis(\gminia,ut)\Gamma_{J'}
 \\
-u_{J}\otimes
 \gminis(\gminia,ut)\Gamma_{J}
-\sum_{J< J'<J+\pi}
 \bigl(
 u_{J'}+(u_{J'_+,0}-u_{J'_-,0})
 \bigr)
 \otimes
 \gminis(\gminia,ut)\Gamma_{J'}
 \Bigr]
 \exp(-zu^{-1}).
\end{multline}
By similar arguments,
it is proved that
we can divided the paths into sums of $1$-chains
such that
the family {\rm(\ref{eq;18.5.29.5})}
is contained in
$\nbigc_0^{!}\bigl(
 (\nbigv,\nabla),u,d(\omega),
 -\Re\gminia^{\circ}(ut)
\bigr)$.
Thus, the first claim of Proposition \ref{prop;18.5.20.31}
is proved in the case $\omega<1$.

\subsection{Proof of the first claim of
Proposition \ref{prop;18.5.20.31}
 in the case $\omega=1$}

Take $\gminia\in\nbigitilde_{J,>0}$.
We shall explain the proof for
$B_{J_-,\theta^u,\gminia}$
in the case $J\cap\vecI\neq\emptyset$.
The proof for $B_{J_+,\theta^u,\gminia}$ is similar.
There exists
$\theta_1\in J\cap \vecI$
such that
$\theta_1\in \Cr_1(\omega;\veckappa(\gminia,u))$.

For $v\in H^0(J_-,L_{J_-,>0})$,
we have the expression
\[
 v=u_{J_+,0}
+\sum_{J\leq J'\leq J+\pi}
 u_{J'},
\]
where 
$u_{J_+,0}$ is a section of $L_{J_+,0}$,
and $u_{J'}$ are sections of $L_{J',<0}$.

\vspace{.1in}

We take a small $\delta>0$.
Let $\Gamma_{\theta_1}$ be the path
obtained as
$\gamma_h(1,\vartheta^J_1-\delta,\vartheta^J_2)$.
We obtain a continuous family
$\Gamma_{\theta_1,\gminic(\gminia,ut)}$
by modifying $\Gamma_{\theta_1}$
as in \S\ref{subsection;18.5.20.1}.
By adding the paths
$\gamma_v(1,0;\vartheta^J_{\ell}-\delta)$,
we obtain a family of paths
$\Gammatilde_{\theta_1,\gminic(\gminia,ut)}$
connecting
$(0,\vartheta^J_{\ell}-\delta)$
and $(1,\vartheta^J_r)$.

Let $\Gamma_{2_{\pm}}$ be the paths obtained as
the union of 
$\gamma_{h}(1;\vartheta^J_{r},\vartheta^J_r\pm\delta)$
and 
$\gamma_v(1,0;\vartheta^J_r\pm\delta)$.

Let $a_0<a_1<\cdots<a_N$
be the intersection of 
$S_0(\nbigi)\cap
 \closedopen{\vartheta^J_{r}}{\vartheta^J_r+\pi}$.
We take
$a_{N+1}\in \openopen{a_N}{\vartheta^J_r+\pi}$.
We set $J_i:=\openopen{a_i-\pi/\omega}{a_i}$.
We have the sections 
$u_{J_{i+},0}\in H^0(J_{i+},L_{J_{i+},0})$
induced by $u_{J,0}$ and the parallel transport of
$\Gr^{\vecnbigf}_0(L)$.

Let $\beta>0$ be sufficiently small satisfying
(\ref{eq;24.3.14.1}).
Let $\Gamma_{3}$ be the path
$\gamma_v(1,t^{\beta};\vartheta^J_r)$.
Let $\Gamma_{4,i}$ $(i=0,\ldots,N)$ be the paths
$\gamma_h(t^{\beta};a_i,a_{i+1})$.
Let $\Gamma_{5}$ be the path
$\gamma_v(\infty,t^{\beta};a_{N+1})$.
Let $\Gamma_{6,i}$ be the paths
$\gamma_v(t^{\beta},0;a_i)$.
We take $\theta_{J+\pi}\in (J+\pi)\cap (\vecI+\pi)\neq\emptyset$.
Let $\Gamma_{7}$ be the path
$\gamma_v(\infty,0;\theta_{J+\pi})$.

We obtain the following family of cycles
which represents
$B_{J_-,\theta^u,\gminia}(v)$:
\begin{multline}
\label{eq;18.5.29.6}
 \Bigl[
 v
 \otimes 
 \gminis(\gminia,ut)
\Gammatilde_{\theta_1,\gminic(\gminia,ut)}
+u_J
 \otimes
  \gminis(\gminia,ut)\Gamma_{2,-}
+\sum_{J<J'\leq J+\pi}
 u_{J'}
 \otimes
  \gminis(\gminia,ut)
 \Gamma_{2,+}
 \\
+u_{J,0}
 \otimes
  \gminis(\gminia,ut)\Gamma_{3}
+\sum_{i=0}^N
 u_{J_i,0}
 \otimes  
 \gminis(\gminia,ut)\Gamma_{4,i}
+u_{J_N,0}
 \otimes 
 \gminis(\gminia,ut) 
 \Gamma_5
 \\
+\sum_{i=1}^N
 (u_{J_{i},0}-u_{J_{i-1},0})
\otimes
  \gminis(\gminia,ut)\Gamma_{6,i}
+u_{J+\pi}\otimes
 \gminis(\gminia,ut)\Gamma_{7}
\Bigr]
 \exp(-zu^{-1}).
\end{multline}

\begin{lem}
The family {\rm(\ref{eq;18.5.29.6})}
is contained in
\[
 \nbigc^{!}_0\bigl(
 (\nbigv,\nabla),u,d(\omega),
 -\Re(\gminia^{\circ}(ut))
 \bigr).
\]
\end{lem}
\pf
Note that 
$-\Re(\gminia^{\circ}(ut))t^{\omega/(1+\omega)}$
converges to a positive numbers
as $t\to 0$.
Hence, the contribution of
$u_{J+\pi}\exp(zu^{-1}t^{-1})
 \otimes
 \gminis(\gminia,ut)\Gamma_7$
can be ignored as in the case of 
Lemma \ref{lem;18.5.29.4}.
We obtain the estimate for
the contributions of the other terms 
by the argument
in the proof of Lemma \ref{lem;18.5.28.30}.
\hfill\qed

\vspace{.1in}

In all,
by Lemma \ref{lem;18.5.30.10},
Proposition \ref{prop;18.5.20.30}
and Proposition \ref{prop;18.5.20.31}
are proved.
\hfill\qed

\subsection{Proof of Proposition \ref{prop;18.5.20.20}}

Let us explain the proof for $J_{1-}$.
The proof for $J_{1+}$ is similar.
We take a small $\delta>0$
such that
$\openopen{\vartheta^{J_1}_{\ell}-\delta}{\vartheta^{J_1}_r}
 \cap(\vecI+\pi)\neq\emptyset$,
and $J_{10}\subset
\openopen{\vartheta^{J_1}_{\ell}-\delta}{\vartheta^{J_1}_r}
 \cap(\vecI+\pi)\neq\emptyset$.

Let $\nbigi(\nbigv)$
denote the set of ramified irregular values of
$(\nbigv,\nabla)$.
Let $(L,\vecnbigftilde)$ be 
the $2\pi\seisuu$-equivariant local system 
with Stokes structure over $\nbigi(\nbigv)$ on $\real$
corresponding to $(\nbigv,\nabla)$.

Let $\projtilde^1\lrarr\proj^1$
denote the oriented real blow up of $\proj^1$ along 
$\{0,\infty\}$.
Let $\nbigl$ denote the local system 
on $\projtilde^1$ associated to $(\nbigv,\nabla)$.

We take a representative
$\vecc(t)$ of $y$ contained in
$\nbigc^{\varrho}_0(\nbigt_{\omega}(\nbigv,\nabla),u,d,Q)$
expressed as in (\ref{eq;18.5.17.2}).
We may assume the following.
\begin{itemize}
\item
There exist intervals $J_{p,j,t}\subset\real$
such that
(i)
$\vartheta^{J_{p,j,t}}_{\ell}-\vartheta^{J_{p,j,t}}_r
<\omega^{-1}\pi$,
(ii) $\gamma_{p,j,t}$
are contained in $\real_{\geq 0}\times J_{p,j,t}$.
\item
There exist intervals $J_{k}\subset\real$
such that
(i)
$\vartheta^{J_k}_{\ell}-\vartheta^{J_k}_r
<\omega^{-1}\pi$,
(ii) $\eta_k$
are contained in $\real_{\geq 0}\times J_k$.
\item
$N_4=1$.
Moreover, $\Gamma_1$ is contained in 
$\realbar_{\geq 0}\times J_{10}$.
Let $P_0$ denote the end point of $\Gamma_1$
contained in $\cnum^{\ast}$.
\end{itemize}

There exist splittings
$L_{|J_{p,j,t}}=
 \bigoplus_{\gminia\in\nbigi}
 L_{J_{p,j,t},\gminia}$
of the Stokes filtrations
$\pi_{\omega\ast}(\nbigftilde)^{\theta}$
$(\theta\in J_{p,j,t})$.
By using the isomorphism
$L_{J_{p,j,t},0}\simeq
 \nbigt_{\omega}(L)_{|J_{p,j,t}}$,
we construct flat sections
$\atilde_{p,j}$
of $\nbigl$
on $\real_{\geq 0}\times J_{p,j,t}$.
Similarly,
by using a splitting
$L_{|J_k}=\bigoplus_{\gminia\in\nbigi}
 L_{k,\gminia}$
 of $\pi_{\omega\ast}(\vecnbigf)$ on $J_k$,
we construct sections
$\btilde_k$ of 
$\nbigl$ on $\real_{\geq 0}\times J_k$
from $b_k$.

By using the canonical splitting
$L_{|J_{1-}}=\bigoplus_{\gminia\in\nbigi}
 L_{J_{1-},\gminia}$
of $\pi_{\omega\ast}(\vecnbigftilde)$,
we construct a section $\ctilde_1$
of $\nbigl$ on $\real_{\geq 0}\times J_{10}$
from $c_1$.

We obtain the following family of chains:
\[
 \vecctilde(t):=
 \Bigl(
 \sum_{\ell=0,1,2}
 \sum_{i=1}^{N_{\ell}}
 \atilde_{\ell,i}\otimes\gamma_{\ell,i,t}
+\sum_{k=1}^{N_3}
 \btilde_{k}\otimes\eta_k
+\ctilde_1\otimes\Gamma_1
 \Bigr)
 \exp(-zu^{-1}).
\]
We obtain
$\del \vecctilde(t)
=\sum e_j\exp(-zu^{-1})\otimes P_{j,t}$,
where the following holds.
\begin{itemize}
\item
 $P_{j,t}=(t^{d_j}s_j,\theta_j)$  for some  $d_j\geq 0$.
We have $d_j=0$ or $d_j\geq d$.
\item
If $d_j=0$,
then $(s_j,\theta_j)\in \{z\,|\,\Re(zu^{-1})>0\}$.
\item
$e_j$ are sections of $q^{-1}(L_{S^1}^{<0})$
around the segments 
     $Z_{j,t}:=\gamma_v(t^{d_j}s_j,0;\theta_j)$,
     where $q:\cnum^{\ast}\lrarr\varpi^{-1}(0)$ is the projection.
\end{itemize}
We obtain the following family of cycles
for $(\nbigv,\nabla)$,
which represents $C^{(J_{1-})}_u(y)$:
\begin{equation}
\label{eq;18.5.29.20}
 \vecctilde(t)
-\sum_{j}e_j\exp(-zu^{-1})
 \otimes Z_{j,t}.
\end{equation}
We obtain the desired estimate for
the components of $\vecctilde(t)$
from the estimate for the components of
$\vecc(t)$.
If $d_{j}\geq d$,
we obtain
$d_j\omega> 1-d\geq 1-d_j$,
and hence the following holds for some $\epsilon>0$:
\[
 \int_{Z_{j,t}}
 \bigl|e_j\bigr|_{h^{\nbigv}}
 \exp(-\Re(zu^{-1})t^{-1})
 =O\Bigl(
 \exp\Bigl(
 -\epsilon t^{-d_j\omega}
 \Bigr)
 \Bigr).
\]
If $d_j=0$,
we have $P_j\in\{z\,|\,\Re(zu^{-1})>0\}$.
Hence, we can easily obtain
\[
  \int_{Z_j}
 \bigl|e_j\bigr|_{h^{\nbigv}}
 \exp(-\Re(zu^{-1})t^{-1})
 =O\Bigl(
 \exp\Bigl(
 -\epsilon t^{-\omega/(1+\omega)}
 \Bigr)
 \Bigr).
\]
Note that $\omega/(1+\omega)>1-d$.
Hence, we can conclude that
the family (\ref{eq;18.5.29.20})
is contained in
$\nbigc_0^{!}\bigl(
 (\nbigv,\nabla),u,d,Q
 \bigr)$.

\hfill\qed

\section{Proof of Proposition \ref{prop;24.3.16.10}}

\label{section;20.11.21.3}

Let $D\subset\cnum$ be a finite subset.
Let $(\nbigv,\nabla)$ be a meromorphic flat bundle
on $(\proj^1,D\cup\{\infty\})$
with regular singularity at $\infty$.
Let $\rho_{\alpha}:\cnum\lrarr \cnum$
be given by $\rho_{\alpha}(z)=z+\alpha$.
We set $\nbigi^{\circ}=\{\alpha u^{-1}\,|\,\alpha\in D\}$.
Let $U_{\alpha}$ be a neighbourhood of $\alpha\in D$.
Let $U_0$ be a neighbourhood of $0$.

\subsection{Families of cycles}
\label{subsection;18.5.23.10}

Let $u\in\cnum^{\ast}$.
Let $\varrho\in\Dsf(D)$.
Let $\alpha\in D$.
Let $\epsilon>0$ be sufficiently small.
Let $0<d\leq 1$.
Let $\vecc^{(\alpha)}(t)$
be a family of $\varrho$-type chains
for $(\nbigv,\nabla)\otimes\nbige(zu^{-1})$
of the following form:
\begin{multline}
\label{eq;18.5.23.2}
 \vecc^{(\alpha)}(t)=
 \vecc_1^{(\alpha)}(t)
+c_{\infty}\exp(-zu^{-1})\otimes\gamma_{\infty}
= \\
\Bigl(
\sum_{\ell=0,1,2}
 \sum_{i=1}^{N_{\ell}}
  a_{\ell,i}\otimes
 \gamma_{\ell,i,t}
+\sum_{i=1}^{N_3}
 b_{i}\otimes
 \eta_i
+c\otimes\Gamma
+c_{\infty}\otimes\gamma_{\infty}
\Bigr)\exp(-zu^{-1}).
\end{multline}
We impose the following conditions:
\begin{itemize}
\item
$\gamma_{0,i,t}$ are paths of the form
$\rho_{\alpha\ast}(t^d\gamma_{i,t})$
for a continuous family of paths
 $\gamma_{i,t}$ $(0\leq t\leq t_0)$
in $U_0\setminus\{0\}$
whose end points are independent of $t$.
\item
$\gamma_{1,i,t}$
are paths of the form
$\rho_{\alpha\ast}\bigl(
 \gamma_h(t^{d_i}r_{1,i};
 \phi_{1,i,1},
 \phi_{1,i,2})
 \bigr)$,
where $d_i\geq d$.
\item
$\gamma_{2,i,t}$ 
are paths of the form
$\rho_{\alpha\ast}\bigl(
 \gamma_v\bigl(
 \epsilon,
 t^{d_{2,i}}r_{2,i};
 \phi_{2,i}
 \bigr)
\bigr)$
where $r_{2,i}\geq 0$,
     $d_{2,i}\geq d$,
and $\gamma_v\bigl(
 \epsilon,
 t^{d_{2,i}}r_{2,i};
 \phi_{2,i}
     \bigr)$ is contained in
     $\{\Re(zu^{-1})>0\}$,
or of the form
$\rho_{\alpha\ast}\bigl(
 \gamma_v\bigl(
 t^{d_{2,i,1}}r_{2,i,1},
 t^{d_{2,i,2}}r_{2,i,2};
 \phi_{2,i}
 \bigr)
 \bigr)$
where $r_{2,i,1}>0$, $r_{2,i,2}\geq 0$
and 
$d\leq d_{2,i,1}\leq d_{2,i,2}$.
     
\item
$\eta_k$ are of the form
$\rho_{\alpha\ast}\bigl(
 \gamma_h(\epsilon;\psi_{k,1},\psi_{k,2})
 \bigr)$
such that $\gamma_h(\epsilon;\psi_{k,1},\psi_{k,2})$
are contained in
$\{\Re(zu^{-1})>0\}$.
\item
Each $\gamma_{p,i,t}$ is contained in a small sector,
and $\rho_{\alpha}^{\ast}(a_{p,i})$ is
a flat section of
$\rho_{\alpha}^{\ast}\bigl(
 \nbigv
\bigr)$ on the sector.
\item
$\rho_{\alpha}^{\ast}b_k$ 
are flat sections of $\rho_{\alpha}^{\ast}(\nbigv)$
on sectors which contain $\eta_k$.
\item
$\Gamma$ is a path connecting
$\alpha+\epsilon e^{\sqrt{-1}\varphi_2}$
and 
$Re^{\sqrt{-1}\varphi_1}$
in $\bigl\{
 \Re((z-\alpha)u^{-1})>0
 \bigr\}\setminus
 \bigcup_{\alpha\in D} U_{\alpha}$,
where
$R$ is a large number,
and 
$\varphi_i$ are chosen such that
$\Re(e^{\sqrt{-1}\varphi_i}u^{-1})>0$.
\item
$\gamma_{\infty}$
is a path connecting
$Re^{\sqrt{-1}\varphi_1}$
and 
$\infty e^{\sqrt{-1}\varphi_1}$
in $\bigl\{
 \Re((z-\alpha)u^{-1})>0
 \bigr\}\setminus
 \bigcup_{\alpha\in D} U_{\alpha}$.
\item
$c$ and $c_{\infty}$
are flat sections of $\nbigv$
along $\Gamma$ and $\Gamma_{\infty}$,
respectively.
\end{itemize}
\begin{rem}
$\vecc_1^{(\alpha)}(t)$
and $c_{\infty}\otimes\gamma_{\infty}$
are divided for the use in 
{\rm\S\ref{subsection;18.5.23.1}}.
\hfill\qed
\end{rem}

Let $Q\in \real[t^{-1/e}]$
such that 
$|t^{1-d}Q(t)|\leq C$ for $C>0$
 as $t\to 0$.
We also consider the following condition.
\begin{itemize}
\item
For any $N$,
there exist $M>0$ and $C>0$ such that
\[
 \int_{\rho_{\alpha}^{-1}\circ\gamma_{0,i,t}}
 \bigl|\rho_{\alpha}^{\ast}a_{0,i}\bigr|_{\rho_{\alpha}^{\ast}(h^{\nbigv})}
 \exp\bigl(-t^{-1}\Re(zu^{-1})\bigr)
 |z|^{-N}
 \leq 
 C\exp\bigl(Q(t)\bigr)t^{-M}.
\]
\item
For any $N>0$,
there exist $\delta>0$ and $C>0$ such that
\[
 \int_{\rho_{\alpha}^{-1}\circ\gamma_{1,i,t}}
 \bigl|\rho_{\alpha}^{\ast}a_{1,i}\bigr|_{\rho_{\alpha}^{\ast}(h^{\nbigv})}
 \exp\bigl(-t^{-1}\Re(zu^{-1})\bigr)
 |z|^{-N}
\leq
C\exp\bigl(Q-\delta t^{-(1-d)}\bigr). 
\]
\item
If $\varrho(\alpha)=!$,
for any $N>0$, there exist $\delta>0$ and $C>0$
such that
\[
 \int_{\rho_{\alpha}^{-1}\circ\gamma_{2,i,t}}
 \bigl|\rho_{\alpha}^{\ast}a_{2,i}\bigr|_{\rho_{\alpha}^{\ast}(h^{\nbigv})}
 \exp\bigl(-t^{-1}\Re(zu^{-1})\bigr)
 |z|^{-N}
\leq
C\exp(Q-\delta t^{-(1-d)}).
\]
If $\varrho(\alpha)=\ast$,
there exist $N>0$, $\delta>0$ and $C>0$
such that
\[
 \int_{\rho_{\alpha}^{-1}\circ\gamma_{2,i,t}}
 \bigl|\rho_{\alpha}^{\ast}a_{2,j}\bigr|_{\rho_{\alpha}^{\ast}(h^{\nbigv})}
 \exp\bigl(-t^{-1}\Re(zu^{-1})\bigr)
 |z|^{N}
\leq
C\exp(Q-\delta t^{-(1-d)}).
\]
\end{itemize}

Let $\nbigc_{\alpha}^{\varrho}((\nbigv,\nabla),u,d,Q)$ 
be the set of such families of $\varrho$-type $1$-chains
for $(\nbigv,\nabla)\otimes\nbige(zu^{-1})$.
\index{set $\nbigc_{\alpha}^{\varrho}((\nbigv,\nabla),u,d,Q)$}

Let $\vecc^{(\alpha)}(t)\in 
 \nbigc_{\alpha}^{\varrho}((\nbigv,\nabla),u,d,Q)$.
If each $\vecc^{(\alpha)}(t)$ is a cycle,
the homology classes of
$\vecc^{(\alpha)}(t)$ are independent of $t$.
(See Lemma \ref{lem;20.10.24.1}.)
We obtain the family of $1$-cycles 
$\vecc^{(\alpha)}(t)\cdot
 \exp\bigl(-(t^{-1}-1)zu^{-1}\bigr)$
of $(\nbigv,\nabla)\otimes\nbige(t^{-1}zu^{-1})$.
They induce a flat section $[\vecc^{(\alpha)}(t)]$ of
$\nbigv^{\circ}_{\varrho}=\Fourier_+(\nbigv(\varrho))$
along $\{tu\,|\,0<t\leq t_0\}$.
We obtain the following 
as a special case of Lemma \ref{lem;18.5.27.10}.

\begin{lem}
$\bigl|
 [\vecc^{(\alpha)}(t)]
 \bigr|_{h_{\varrho}}
=O\Bigl(
 \exp\bigl(-t^{-1}\Re(\alpha u^{-1})+Q\bigr)
 \cdot t^{-N}
 \Bigr)$ for some $N>0$
as $t\to 0$.
\hfill\qed
\end{lem}

\subsection{Proof of Proposition \ref{prop;24.3.16.10}}

\begin{prop}
\label{prop;24.3.16.11}
There exists a finite subset
$\ttS\subset\bigl\{\alpha\in\cnum\,\big|\,|\alpha|=1\bigr\}$
such that the following holds
unless $|u|^{-1}u\in\ttS$.
\begin{itemize}
\item
Any $y\in \nbigf^{\circ\theta^u}_{\gminia^{\circ}}
H_1^{\varrho}(\cnum\setminus D,
(\nbigv,\nabla)\otimes\nbige(zu^{-1}))$
is expressed as a sum
$\sum \vecc^{(\alpha_i)}(t)$
where
$\vecc^{(\alpha_i)}(t)$
are represented by families of cycles
in $\nbigc^{\varrho}_{\alpha_i}((\nbigv,\nabla),u,d_i,Q_i)$
such that
$-\Re(\alpha_i u^{-1}t^{-1})+Q_i(t)
\leq
 -\Re(\gminia^{\circ}(ut))$
for any sufficiently small $t>0$.
\end{itemize}
\end{prop}
\pf
Let
$\theta^u\in \real\setminus S_0(\nbigi^{\circ})$.
Suppose that 
$\rho_{\alpha}^{\ast}y\in 
 H_1^{\varrho(\alpha)}\bigl(\cnum\setminus\{0\},
 \rho_{\alpha}^{\ast}(\nbigv_{\alpha},\nabla)
 \otimes\nbige(zu^{-1})
 \bigr)$
is represented by a family of cycles
$\vecc(t)$
contained in 
$\nbigc_0^{\varrho(\alpha)}
 (\rho_{\alpha}^{\ast}(\nbigv_{\alpha},\nabla),u,d,Q)$
as in \S\ref{subsection;18.5.21.10}.
We obtain the family of 
$\varrho(\alpha)$-type cycles
$\rho_{\alpha\ast}(\vecc(t))\exp(-\alpha u^{-1})$
for $(\nbigv_{\alpha},\nabla)\otimes\nbige(zu^{-1})$.
By modifying $\Gamma$,
we may assume that the underlying chains of
$\rho_{\alpha\ast}(\vecc(t))
\exp(-\alpha u^{-1})$ are contained in
the union of $U_{\alpha}$
and $\nbigu_{\vecJ_+,u}$ or $\nbigu_{\vecJ_-,u}$.
We have an isomorphism
$(\nbigv,\nabla)\simeq(\nbigv_{\alpha},\nabla)$
on the union of $U_{\alpha}$
and a neighbourhood of $\Gamma$.
We naturally regard 
$\rho_{\alpha\ast}(\vecc(t))
\exp(-\alpha u^{-1})$ as a family of
$\varrho$-type $1$-cycles
for $(\nbigv,\nabla)\otimes\nbige(zu^{-1})$,
which we denoted by
$\vecc_{\alpha}(t)$.
Then, we can easily see that
$C^{\varrho}_{\vecJ_{\pm,\alpha}}(y)$
is represented by $\vecc_{\alpha}(t)$,
and 
it is contained in 
$\nbigc^{\varrho}_{\alpha}\bigl(
 (\nbigv,\nabla),u,d,Q
 \bigr)$.

The claim of Proposition \ref{prop;24.3.16.11}
follows from Proposition \ref{prop;24.3.15.20}
and the above consideration.
\hfill\qed

\vspace{.1in}

For $u$ such that $|u|^{-1}u\not\in\ttS$,
we obtain
$\nbigf^{\prime\theta^u}_{\gminia}
\subset
 \nbigf^{\circ\theta^u}_{\gminia}$
for any $\gminia\in\nbigi^{\circ}$
from Proposition \ref{prop;24.3.16.11}.
We obtain $\nbigf^{\prime\theta^u}=\nbigf^{\circ\theta^u}$
because the dimensions of the associated graded pieces
are the same.
Then, we obtain
$\nbigf^{\prime\theta^u}=\nbigf^{\circ\theta^u}$
for any $u$ by Lemma \ref{lem;20.10.23.1}.
Thus, we obtain Proposition \ref{prop;24.3.16.10}.
\hfill\qed

\section{Proof of Theorem \ref{thm;24.3.16.20}}
\label{subsection;18.5.23.1}

Let $(\nbigv,\nabla)$ be a meromorphic flat bundle
on $(\proj^1,D\cup\{\infty\})$.
Let $\nbigi_{\infty}(\nbigv)$
be the set of ramified irregular values of $(\nbigv,\nabla)$
at $\infty$.
Let $h_{\nbigv}$ be a metric of $\nbigv_{|\cnum\setminus D}$ 
adapted to the meromorphic structure.
Take $u\in\cnum^{\ast}$.
Set $\theta^u:=\arg(u)$.
We set
\[
 \omega(\nbigv):=
 \min\{\omega\,|\,\nbigstilde^{\infty}_{\omega}(\nbigv)\neq\nbigv\}
 =\min\{\omega\,|\,\nbigstilde_{\omega}(\nbigi_{\infty}(\nbigv))
 \neq\nbigi_{\infty}(\nbigv)\}.
\]
 
Let $U_{\infty}$ be a neighbourhood of $\infty$
with the coordinate $x=z^{-1}$.
Let $\Utilde_{\infty}\lrarr U_{\infty}$
denote the oriented real blow up.
We use the polar coordinate induced by $x$.

\subsection{Families of cycles}

Take $0<d$ such that
$\omega(\nbigv)\leq (d+1)/d$.
We consider families of 
$\varrho$-type $1$-chains $\vecc^{(\infty)}(t)$ $(0<t\leq 1)$ of
$(\nbigv,\nabla)\otimes\nbige(zu^{-1})$
of the following form:
\begin{equation}
 \label{eq;18.5.17.20}
 \vecc^{(\infty)}(t)=
\Bigl(
 \sum_{i=1}^{N_{-1}}c_i\otimes \nu_{i,t}
+\sum_{i=1}^{N_0}
 a_{i}\otimes\gamma_{i,t}
+\sum_{j=1}^{N_1}
 b_j\otimes\eta_{j,t}
\Bigr)
 \exp(-zu^{-1}).
\end{equation}
We impose the following conditions
by using the polar coordinate $x=re^{\sqrt{-1}\theta}=z^{-1}$.
\begin{itemize}
\item
$\nu_{i,t}$ are continuous families of paths
of the form
$t^d\nu_{i,t}'$ on $U_{\infty}\setminus\{\infty\}$,
where $\nu'_{i,t}$ $(0\leq t\leq t_0)$
are continuous families of paths
     on $U_{\infty}\setminus\{\infty\}$
     whose end points are independent of $t$.
We assume that 
each $\nu_{i,t}$ is contained in a small sector $S$,
and $c_i$ is a flat section on $S$.
Moreover, we assume that
$c_i$ is a section of
a direct summand 
$\nbigv_{S,\gminia_i}$
for a splitting of the Stokes filtration of $\nbigv$
on $S$,
where
$0<C_1<|x|^{(1+d)/d}|\gminia_i|<C_2$ 
for some constants $C_b$.
\item
     For each $\gamma_{i,t}$,
     one of the following holds:
     $\gamma_{i,t}$ is of the form
     $\gamma_{v}(r_{i,1},t^{d_{i,2}}r_{i,2};\phi_i)$,
     where $d< d_{i,2}$, $r_{i,1}>0$, $r_{i,2}\geq 0$,
     and contained in
     $\{\Re(x^{-1}u^{-1})>0\}$;
     or $\gamma_{i,t}$ is of the form
     $\gamma_{v}(t^{d_{i,1}}r_{i,1},t^{d_{i,2}}r_{i,2};\phi_i)$,
     where 
     $d\leq d_{i,1}<d_{i,2}$,
     $r_{i,1}>0$ and $r_{i,2}\geq 0$.
     We assume
     $\omega(\nbigv)\leq (d_{i,1}+1)/d_{i,1}$ and
     $\omega(\nbigv)\leq (d_{i,2}+1)/d_{i,2}$
     in any case.
\item
 $\eta_{i,t}$
are of the form
$\gamma_h(t^{d_i}r_i;\psi_{i,1},\psi_{i,2})$
     where we assume
     $d_i=0$, or
     $d_i\geq d$
     and $\omega(\nbigv)\leq (d_i+1)/d_i$.
\item
$\gamma_{i,t}$ and $\eta_{j,t}$
are contained in a small sector
in $(U_{\infty},\infty)$,
and $a_i$ and $b_j$ are flat sections of
$(\nbigv,\nabla)$ on the sectors.
\end{itemize}
Let $Q\in \real[t^{-1/e}]$ such that 
$|t^{1+d}Q(t)|\leq C$ for some $C>0$ as $t\to 0$.
We say that the growth order of $\vecc^{(\infty)}(t)$ is
less than $Q$ if the following holds:
\index{growth order}
\begin{itemize}
\item
For any $N>0$,
there exists $M>0$ and $C>0$
such that
\[
 \int_{\nu_{i,t}}|c_i|_{h^{\nbigv}}
 \exp\Bigl(
 -t^{-1}\Re(zu^{-1})
 \Bigr)|z|^N
\leq C \exp(Q(t))t^{-M}.
\]
\item
For any $N$,
there exist $C>0$ and $\delta>0$
such that
\[
 \int_{\gamma_{i,t}}
 |a_{i}|_{h^{\nbigv}}
 \exp\Bigl(
 -t^{-1}\Re(zu^{-1})
 \Bigr)|z|^N
\leq C
 \exp\Bigl(Q(t)-\delta t^{-(1+d)}\Bigr).
\]
Moreover, we also impose
\[
  |a_{i}|_{h^{\nbigv}}
 \exp\bigl(
 -t^{-1}\Re(zu^{-1})
 \bigr)
\leq 
  C\exp\Bigl(Q(t)-\delta_1t^{-(1+d)}-\delta_2 t^{-1}|x|^{-1}\Bigr)
\]
on $\gamma_{i,t}$
for some $C>0$ and $\delta_i>0$.
\item
For any $N>0$,
there exist $C>0$ and $\delta>0$
such that
\[
 \int_{\eta_{i,t}}
 |b_{i}|_{h^{\nbigv}}
 \exp\Bigl(
 -t^{-1}\Re(zu^{-1})
 \Bigr)|z|^N
\leq
 C\exp\Bigl(Q(t)-\delta t^{-(1+d_i)}\Bigr).
\]
Note that if $d_i>d$,
it implies the following for some 
$C'>0$ and $\delta'>0$:
\[
 \int_{\eta_{i,t}}
 |b_{i}|_{h^{\nbigv}}
 \exp\Bigl(
 -t^{-1}\Re(zu^{-1})
 \Bigr)|z|^N
\leq C'
 \exp\bigl(-\delta' t^{-(1+d_i)}\bigr).
\]
\end{itemize}
Let $\nbigc^{(\infty)\,\varrho}_{\infty}\bigl(
 (\nbigv,\nabla),u,d,Q
 \bigr)$
be the set of such families of $1$-cycles
for $(\nbigv,\nabla)\otimes\nbige(zu^{-1})$.
\index{set $\nbigc^{(\infty)\,\varrho}_{\infty}\bigl(
 (\nbigv,\nabla),u,d,Q
 \bigr)$}

\vspace{.1in}

Let $\alpha\in D$.
We also consider families of $1$-cycles 
for $(\nbigv,\nabla)\otimes\nbige(zu^{-1})$
of the form
\[
 \vecc^{(\alpha)}(t)=
\Bigl(
 \sum_{i=1}^{N_0}
 a_{i}\otimes\gamma_{i,t}
+\sum_{j=1}^{N_1}
 b_j\otimes\eta_{j,t}
\Bigr)\exp(-zu^{-1})
+\vecc^{(\alpha)}_1(t).
\]
Here, 
$\sum_{i=1}^{N_0}
 a_{i}\otimes\gamma_{i,t}
+\sum_{j=1}^{N_1}
 b_j\otimes\eta_{j,t}$
is as in (\ref{eq;18.5.17.20}),
and $\vecc^{(\alpha)}_1(t)$
is as in \S\ref{subsection;18.5.23.10}
for some $0<d(\alpha)<1$.
Let $Q\in \real[t^{-1/e}]$
such that $\deg Q\leq 1-d(\alpha)$.
We say that
the growth order of $\vecc^{(\alpha)}(t)$
is less than
$-\Re(\alpha u^{-1}t^{-1})+Q(t)$
if the following holds.
\begin{itemize}
\item
$\vecc^{(\alpha)}_1(t)+c_{\infty}\exp(-zu^{-1})\otimes\gamma_{\infty}$
is contained in
     $\nbigc^{\varrho}_{\alpha}
     \bigl(\nbigs^{\infty}_0(\nbigv,\nabla),u,d(\alpha),Q\bigr)$.
Here, see \S\ref{subsection;18.5.23.10}
for $c_{\infty}$ and $\gamma_{\infty}$,
and \S\ref{subsection;20.10.24.3} for
$\nbigs^{\infty}_0(\nbigv,\nabla)$.
 \item
For any $N>0$, there exist $C>0$ and $\delta>0$
such that
\[
 \int_{\gamma_{i,t}}
 |a_{i}|_{h^{\nbigv}}
 \exp\Bigl(
 -t^{-1}\Re(zu^{-1})
 \Bigr)|z|^N
\leq 
 C\exp\Bigl(-\Re(\alpha u^{-1}t^{-1})+Q(t)-\delta t^{-1}\Bigr).
\]
Moreover, we also impose
\[
 |a_{i}|_{h^{\nbigv}}
 \exp\bigl(
 -t^{-1}\Re(zu^{-1})
 \bigr)
\leq 
  C\exp\Bigl(-\Re(\alpha u^{-1}t^{-1})
+Q(t)-\delta t^{-1}|x|^{-1}\Bigr)
\]
on $\gamma_{i,t}$
for some $C>0$ and $\delta>0$.
\item
There exist $C>0$ and $\delta>0$ such that
\[
 \int_{\eta_{i,t}}
 |b_{i}|_{h^{\nbigv}}
 \exp\bigl(
 -t^{-1}\Re(zu^{-1})
 \bigr)|z|^N
\leq C
 \exp\bigl(-\Re(\alpha u^{-1}t^{-1})+Q(t)-\delta t^{-(1+d_i)}\bigr)
 \bigr)
\]
Note that we have
\[
 \int_{\eta_{i,t}}
 |b_{i}|_{h^{\nbigv}}
 \exp\bigl(
 -t^{-1}\Re(zu^{-1})
 \bigr)|z|^N
\leq C'
 \exp\bigl(-\delta' t^{-(1+d_i)}\bigr)
\]
for some $C'>0$ and $\delta'>0$.
\end{itemize}
Let $\nbigc^{(\alpha)\,\varrho}_{\infty}((\nbigv,\nabla),
 u,d,Q)$
denote the set of such families of cycles.
\index{set $\nbigc^{(\alpha)\,\varrho}_{\infty}((\nbigv,\nabla),
 u,d,Q)$}

\subsubsection{}
Let $\vecc^{(\alpha)}(t)\in
\nbigc^{(\alpha)\,\varrho}_{\infty}
\bigl(
(\nbigv,\nabla),u,d,Q
\bigr)$,
where $\alpha\in D\cup\{\infty\}$.
If each $\vecc(t)$ is a $1$-cycle,
the homology classes of $\vecc(t)$
are constant
as in the case of Lemma \ref{lem;20.10.24.1}.
We obtain the family of  $\varrho$-type $1$-cycles
$\vecc(t)\cdot 
 \exp\Bigl(
 -(t^{-1}-1)zu^{-1}\Bigr)$
of $(\nbigv,\nabla)\otimes\nbige(t^{-1}zu^{-1})$.
They induce a flat section
$[\vecc(t)]$ of 
along $\{tu\,|\,0<t<t_0\}$.
\begin{lem}
We obtain the following the estimate for some $N$ as $t\to\ 0$.
\begin{itemize}
 \item 
$\bigl|
[\vecc(t)]
\bigr|_{h_{\varrho}}
 =O\Bigl(
 \exp(Q)\cdot t^{-N}
       \Bigr)$
        in the case $\alpha=\infty$.
 \item
      $\bigl|
      [\vecc(t)]
      \bigr|_{h_{\varrho}}
      =O\Bigl(
      \exp(-t^{-1}\Re(\alpha u^{-1})+Q)\cdot t^{-N}
      \Bigr)$
      in the case $\alpha\in D$.
\end{itemize}
\end{lem}
\pf
It follows from Lemma \ref{lem;18.5.27.10}.
\hfill\qed

\subsection{Statements}

Suppose that
$\omega:=\omega(\nbigv)>1$.
Let $(V,\nabla):=\nbigttilde^{\infty}_{\omega}(\nbigv,\nabla)$.
Set $\nbigitilde:=\nbigi(V)$.
We have $\ord(\nbigitilde)=-\omega$.
We set $\nbigi:=\pi_{\omega}(\nbigitilde)$.

Let $(L,\vecnbigftilde)$ be the $2\pi\seisuu$-equivariant
local system with Stokes structure on $\real$
indexed by $\nbigitilde$ associated to $(V,\nabla)$.
For $J\in \gbigw_2(\nbigi,\theta^u,\pm)$,
there exist splittings
\[
 L_{J_{\mp},<0}=
 \bigoplus_{\gminia\in\nbigitilde_{J,<0}}
 L_{J_{\mp},\gminia}
\]
of the Stokes filtrations
$\nbigftilde^{\theta}$ $(\theta\in J_{\mp})$.
For any $\gminia\in\nbigitilde_{J,<0}$,
we obtain the following map
induced by $A_{J_{\pm},\theta^u}$
and the natural morphism
$H_1^{\rd}\bigl(
 \cnum^{\ast},
 (V,\nabla)\otimes
 \nbige(x^{-1}u^{-1})
 \bigr)
\lrarr
 H_1^{\varrho}\bigl(
 \cnum\setminus D,
 (\nbigv,\nabla)\otimes
 \nbige(zu^{-1})
 \bigr)$
in \S\ref{subsection;18.5.15.40}:
\[
 A_{J_{\pm},\theta^u,\gminia}:
 H^0(J_{\mp},L_{J_{\mp},\gminia})
\lrarr
 H_1^{\varrho}\bigl(
 \cnum\setminus D,
 (\nbigv,\nabla)\otimes
 \nbige(zu^{-1})
 \bigr).
\]
For any $u_1=|u|e^{\sqrt{-1}\theta_1^u}$
such that $|\theta_1^u-\theta^u|<\pi/2$,
there exists the natural isomorphism
\[
  H_1^{\varrho}\bigl(
 \cnum\setminus D,
 (\nbigv,\nabla)\otimes
 \nbige(zu^{-1})
 \bigr)
\simeq
   H_1^{\varrho}\bigl(
 \cnum\setminus D,
 (\nbigv,\nabla)\otimes
 \nbige(zu_1^{-1})
 \bigr).
\]
Let $A_{J_{\pm},(\theta^u,\theta^u_1),\gminia}$
denote the following induced map
\[
 H^0(J_{\mp},L_{J_{\mp},\gminia})
\lrarr
 H_1^{\varrho}\bigl(
 \cnum^{\ast},
 (V,\nabla)\otimes
 \nbige(zu^{-1})
 \bigr)
\lrarr 
 H_1^{\varrho}\bigl(
 \cnum\setminus D,
 (\nbigv,\nabla)\otimes
 \nbige(zu_1^{-1})
 \bigr).
\]
For any $\gminia\in\nbigitilde_{J,<0}$,
we set
$\gminia^{\circ}:=
 \gbigf^{(\infty,\infty)}_{(J,0,-)}(\gminia)
 \in\nbigitilde^{\circ}:=
 \gbigf^{(\infty,\infty)}_+(\nbigitilde)$.
We also put $d(\omega):=(\omega-1)^{-1}$.

\begin{prop}
\label{prop;18.5.21.30}
 \mbox{{}}
\begin{itemize}
 \item
      If $\Jbar\subset\vecI_x(\theta^u)$, then
      for any $v\in H^0(J_{\mp},L_{J_{\mp},\gminia})$,
      $A_{J_{\pm},\theta^u,\gminia}(v)$
is represented by a family of cycles contained in
$\nbigc^{(\infty)\,\varrho}_{\infty}\bigl(
 (\nbigv,\nabla),
 u,d(\omega),-\Re(\gminia^{\circ}(ut))
 \bigr)$.
 \item Suppose $J\subset\vecI_x(\theta^u)$
       but $\Jbar\not\subset\vecI_x(\theta^u)$.
       If $|\theta^u-\theta^u_1|\neq 0$ is sufficiently small,
       then for any $v\in H^0(J_{\mp},L_{J_{\mp},\gminia})$,
$A_{J_{\pm},(\theta^u,\theta^u_1),\gminia}(v)$
is represented by a family of cycles in
$\nbigc^{(\infty)\,\varrho}_{\infty}\bigl(
 (\nbigv,\nabla),
 u_1,d(\omega),-\Re(\gminia^{\circ}(u_1t))
 \bigr)$.
\end{itemize}
\end{prop}

Take $J\in\gbigw_1(\nbigi,\theta^u,\pm)$.
There exist splittings
\[
 L_{J_{\mp},>0}=
 \bigoplus_{\gminia\in\nbigitilde_{J,>0}}
 L_{J_{\mp},\gminia}
\]
of the Stokes filtrations
$\nbigftilde^{\theta}$ $(\theta\in J_{\mp})$.
By using $\BB^{\rd}_{J_{\pm},\theta^u}$
and the natural morphism in \S\ref{subsection;18.5.15.40},
we obtain the following morphism
$\BB^{\rd}_{J_{\pm},\theta^u,\gminia}$
for any $\gminia\in\nbigitilde_{J,>0}$:
\[
   H^0(J_{\mp},L_{J_{\mp},\gminia})
 \lrarr
 H_1^{\rd}\bigl(\cnum^{\ast},(V,\nabla)\otimes\nbige(x^{-1}u^{-1})
 \bigr)
\lrarr
  H_1^{\varrho}\bigl(\cnum\setminus D,(\nbigv,\nabla)\otimes\nbige(zu^{-1})
 \bigr)
\]
Moreover, for
any $u_1=|u|e^{\sqrt{-1}\theta_1^u}$
such that $|\theta^u_1-\theta^u|<\pi/2$,
we obtain the following morphism
$\BB^{\rd}_{J_{\pm},(\theta^u,\theta^u_1),\gminia}$:
\[
 H^0(J_{\mp},L_{J_{\mp},\gminia})
\lrarr
   H_1^{\varrho}\bigl(\cnum\setminus D,
   (\nbigv,\nabla)\otimes\nbige(zu_1^{-1})
 \bigr).
\]
For any $\gminia\in\nbigitilde_{J,>0}$,
we set
$\gminia^{\circ}:=
 \gbigf^{(\infty,\infty)}_{(J,0,+)}(\gminia)
 \in\nbigitilde^{\circ}$.

\begin{prop}
\label{prop;18.5.21.31}
\mbox{}
\begin{itemize}
 \item If $\Jbar\subset\vecI_x(\theta^u)-\pi$,
 for any $v\in H_0(J_{\mp},L_{J_{\mp},\gminia})$,
       $\BB^{\rd}_{J_{\pm},\theta^u,\gminia}(v)$
       is represented by a family of cycles
       contained in 
$\nbigc^{(\infty)\,\varrho}_{\infty}
 \bigl((\nbigv,\nabla),u,d(\omega),
 -\Re(\gminia^{\circ}(ut))
 \bigr)$.
 \item
      Suppose that $J\subset\vecI_x(\theta^u)-\pi$
      but $\Jbar\not\subset\vecI_x(\theta^u)-\pi$.
      If $|\theta^u-\theta^u_1|\neq 0$,
 for any $v\in H_0(J_{\mp},L_{J_{\mp},\gminia})$,
$\BB^{\rd}_{J_{\pm},(\theta^u,\theta^u_1),\gminia}(v)$
is represented by a family of cycles contained in 
$\nbigc^{(\infty)\,\varrho}_{\infty}
 \bigl((\nbigv,\nabla),u_1,d(\omega),
 -\Re(\gminia^{\circ}(u_1t))
 \bigr)$.       
\end{itemize}
\end{prop}

We remark the following.

\begin{lem}
\label{lem;20.10.25.22}
The first claims of 
Proposition {\rm\ref{prop;18.5.21.30}}
and Proposition {\rm\ref{prop;18.5.21.31}}
imply
the second claims of  
Proposition {\rm\ref{prop;18.5.21.30}}
and Proposition {\rm\ref{prop;18.5.21.31}}.
\end{lem}
\pf
Suppose $\vartheta^J_{\ell}=\vartheta^{\vecI_x(\theta^u)}_{\ell}$.
We set $\Jhat=J-\omega^{-1}\pi$.
By Proposition \ref{prop;24.3.17.1},
for any $v\in H^0(J_{-},L_{J,\gminia})$,
we obtain
\[
A_{J,\theta^u}(v)
=\BB^{\rd}_{\Jhat,\theta^u}
\bigl(
\nbigr^{J_-}_{\Jhat}(v)
\bigr)
+\sum_{J-\pi<J'<\Jhat}
\BB^{\rd}_{J',\theta^u}\bigl(
\nbigr^{J_-}_{J'}(v)
\bigr).
\]
For any $J-\pi<J'<\Jhat$,
we have
$J'\cap (\vecI_x(\theta^u)-\pi)\neq\emptyset$,
and $\nbigr^{J_-}_{J'}(v)\in H^0(J',L_{J',>0})$.
We also have
$\gminia^{\circ}=\gbigf^{(\infty,\infty)}_{J,0,-}(\gminia)
=\gbigf^{(\infty,\infty)}_{\Jhat,0,+}(\gminia)$
and
$\nbigr^{J_-}_{\Jhat}(v)
\in \nbigf^{\prime\theta^u}_{\gminia^{\circ}}
H^0(\Jhat,L_{\Jhat_+,>0})$.
Then, we obtain the claim of Lemma \ref{lem;20.10.25.22}
by using the arguments for Lemma \ref{lem;18.5.30.10}.
\hfill\qed

\vspace{.1in}

Take $J_1\in T(\nbigi)$
such that 
$J_{1\pm}\subset \vecI_x(\theta^u)-\pi$.
Let
$y\in
 H_1^{\varrho}\bigl(
 \cnum\setminus D,
 \nbigstilde^{\infty}_{\omega}(\nbigv,\nabla)
 \otimes\nbige(zu^{-1})
 \bigr)$.
We obtain
$C^{J_{1\pm}}_{\infty,\theta^u}(y)
\in
 H_1^{\varrho}\bigl(
 \cnum\setminus D,
 (\nbigv,\nabla)\otimes\nbige(zu^{-1})
 \bigr)$
as in \S\ref{section;24.3.17.2}.

\begin{prop}
\label{prop;18.5.21.20}
Suppose that $y$ is represented by 
a family of $1$-cycles contained in
$\nbigc^{(\alpha)\,\varrho}_{\infty}\bigl(
 \nbigstilde_{\omega}^{\infty}(\nbigv,\nabla),
 u,d,Q
 \bigr)$ for $\alpha\in D\cup\{\infty\}$.
Suppose $\omega< d^{-1}(1+d)$.
Then, $C^{J_{1\pm}}_{\infty,\theta^u}(y)$
are represented by a family of $1$-cycles
contained in 
$\nbigc^{(\alpha)\,\varrho}_{\infty}\bigl(
 (\nbigv,\nabla),
 u,d,Q
 \bigr)$.
\end{prop}

We shall prove
Propositions
\ref{prop;18.5.21.30}, \ref{prop;18.5.21.31}
and \ref{prop;18.5.21.20}
in \S\ref{subsection;20.10.25.1}--\S\ref{subsection;20.10.25.3}
after preliminaries in
\S\ref{subsection;20.10.25.4}--\S\ref{subsection;18.5.21.100}.

\subsection{Proof of Theorem \ref{thm;24.3.16.20}}

By using
Propositions \ref{prop;18.5.21.30}--\ref{prop;18.5.21.20}
together with an argument in the proof of 
Proposition \ref{prop;24.3.15.20},
we can prove Theorem \ref{thm;24.3.16.20}
and the following proposition.

\begin{prop}
\label{prop;24.3.18.4}
Any element of
$\nbigf^{\circ\theta^u}_{\gminia^{\circ}}
 H_1^{\varrho}\bigl(
 \cnum\setminus D,\nbigv\otimes\nbige(zu^{-1})
 \bigr)$
is represented as a sum $\sum c_i$,
where $c_i$ are families of cycles
in 
$\nbigc_{\infty}^{(\alpha_i)\,\varrho}\bigl(
 (\nbigv,\nabla),u,d_i, Q_i \bigr)$
satisfying the following condition.
\begin{itemize}
 \item If $\alpha_i=\infty$
       we have $Q_i(t)\leq -\Re(\gminia^{\circ}(ut))$
       for any sufficiently small $t>0$
 \item If $\alpha_i\in D$,
       we have
       $-\Re(\alpha_i u^{-1}t^{-1})+Q_i(t)\leq -\Re(\gminia^{\circ}(ut))$
       for any sufficiently small $t>0$.
       \hfill\qed
\end{itemize}
\end{prop}

\subsection{Scaling}
\label{subsection;20.10.25.4}

For any 
$\gminia
=\sum_{0<j\leq\omega}
 \gminia_j r^{-j}e^{\sqrt{-1}j \theta}
\in\nbigitilde$,
we set
$\veckappa(\gminia,u):=
 (\arg(\gminia_{\omega}),-\theta^u)$
and
$\gminis(\gminia,u):=
 \bigl|\omega \gminia_{\omega} u\bigr|^{1/(\omega-1)}$.
We also set
\[
\gminic(\gminia,u)=
\bigl(
 \gminia_j\cdot
 \bigl|\omega\gminia_{\omega}\bigr|^{(-j+1)/(\omega-1)}
\cdot
 \bigl|u\bigr|^{(\omega-j)/(\omega-1)}
 \bigr)_{0<j<\omega}.
\]

Set
$G_{\gminia,u}(r,\theta)
:=\gminia(re^{\sqrt{-1}\theta})
+u^{-1}r^{-1}e^{-\sqrt{-1}\theta}$.
We remark
\[
 G_{\gminia,u}\Bigl(
 \gminis(\gminia,u) r,
 \theta
 \Bigr)
=\bigl(\omega|\gminia_{\omega}|\bigr)^{-1/(\omega-1)}
 |u|^{-\omega/(\omega-1)}
 G_{\veckappa(\gminia,u),\gminic(\gminia,u)}(r,\theta).
\]

\subsection{Lift of $1$-chains}
\label{subsection;18.5.21.100}

We set $\Delta_x=\{|x|<1\}$.
Let $(\nbigv_0,\nabla)$ be
a meromorphic flat bundle on $(\Delta,0)$.
Take $\omega_0>0$.
We obtain the induced meromorphic flat bundle
$\nbigt_{\omega_0}(\nbigv_0,\nabla)$ on $(\Delta,0)$.
Let $\varpi:\Deltatilde\lrarr\Delta$
be the oriented real blow up along $0$.
Let $\nbigl(\nbigv_0)$ and
$\nbigt_{\omega_0}(\nbigl(\nbigv_0))$
be the local systems on $\Deltatilde$
associated to $\nbigv_0$
and $\nbigt_{\omega_0}(\nbigv_0)$.

Let $\nbigi(\nbigv_0)$ be the set of ramified irregular values
of $(\nbigv_0,\nabla)$.
For simplicity,
we assume $0\not\in\nbigi(\nbigv_0)$.
Let $(L(\nbigv_0),\vecnbigftilde)$ be
the $2\pi\seisuu$-equivariant
local system with Stokes structure indexed by
$\nbigi(\nbigv_0)$
on $\real$ corresponding to $(\nbigv_0,\nabla)$.
Set $\nbigi_0=\pi_{\omega_0}(\nbigi(\nbigv_0))$
and $\vecnbigf:=\pi_{\omega_0\ast}(\vecnbigftilde)$.

Let $\gamma:\closedclosed{0}{1}\lrarr \Deltatilde$ 
be a path
such  that $\gamma(\openopen{0}{1})\subset\Delta\setminus\{0\}$.
There exists a sequence
$t_0=0<t_1<\ldots<t_{N-1}<t_N=1$
such that 
each $\gamma(\closedclosed{t_i}{t_{i+1}})$
is contained in 
$\{re^{\sqrt{-1}}\,|\,0\leq r<1,\,\theta\in I_i\}$
for sectors $I_i$
with $\vartheta^{I_i}_{\ell}-\vartheta^{I_i}_r<\omega_0^{-1}\pi$.
There exist splittings
$L_{|I_i}=\bigoplus_{\gminia\in \nbigi} L_{I_i,\gminia}$
of Stokes filtrations $\nbigf^{\theta}$ $(\theta\in \nbigi)$.
Let $\gamma_i$ denote the restriction of $\gamma$
to $\closedclosed{t_i}{t_{i+1}}$.
Let $Z_i$ be the segments
connecting $\gamma(t_i)$ and $0$.
They naturally induce paths on $\Deltatilde$,
which are also denoted by $Z_i$.

Let $v$ be any element of
$\nbigt_{\omega_0}(\nbigl(\nbigv_0))_{|\gamma(0)}$.
If $\gamma(0)\in\varpi^{-1}(0)$,
we assume that
$v\in \nbigt_{\omega_0}(L)^{<0}_{|\gamma(0)}$.
It induces a section $\vtilde$ of 
$\nbigt_{\omega_0}(\nbigl(\nbigv_0))$ along $\gamma$.
If $\gamma(1)\in\varpi^{-1}(0)$,
we assume that
$\vtilde_{|\gamma(1)}\in
\nbigt_{\omega_0}(L)^{<0}_{|\gamma(1)}$.
Then, 
$\vtilde\otimes\gamma$ is 
a rapid decay $1$-chain $\vtilde\otimes\gamma$
of $\nbigt_{\omega_0}(\nbigl(\nbigv_0))$.
By using the splitting on $I_i$,
we obtain flat sections $\vtilde_i$  of $\nbigl(\nbigv_0)$
along $\gamma_i$ from the restriction of $\vtilde$
to $\gamma_i$.
Thus, we obtain the following rapid decay $1$-chain
of $(\nbigv_0,\nabla)$:
\begin{equation}
\label{eq;20.10.25.10}
\mu(\vtilde\otimes\gamma)=
 \sum_{i=0}^{N-1} \vtilde_i\otimes\gamma_i
 +\sum_{i=1}^{N-1}
 (\vtilde_{i-1|\gamma(t_i)}-\vtilde_{i|\gamma(t_i)})\otimes Z_i.
\end{equation}
It is called a lift of $\vtilde\otimes\gamma$
to a rapid decay $1$-chain of $(\nbigv_0,\nabla)$.
We set
$\mu(\vtilde\otimes\gamma)_{|\gamma(0)}
=\vtilde_{0|\gamma(0)}
\in\nbigl(\nbigv_0)_{|\gamma(0)}$
and 
$\mu(\vtilde\otimes\gamma)_{|\gamma(1)}
=\vtilde_{N-1|\gamma(1)}
\in\nbigl(\nbigv_0)_{|\gamma(1)}$.
The term
$\mu^{\add}(\vtilde\otimes\gamma)
=\sum_{i=1}^{N-1}
 (\vtilde_{i-1|\gamma(t_i)}-\vtilde_{i|\gamma(t_i)})\otimes Z_i$
is called the additional term.

Let $\sum_{j=1}^M \vtilde^{(j)}\otimes\gamma^{(j)}$
be a $1$-chain for $(\nbigt_{\omega_0}(\nbigv_0),\nabla)$.
By applying the above procedure to each
$\vtilde^{(j)}\otimes\gamma^{(j)}$,
we construct $1$-chains
$\mu(\vtilde^{(j)}\otimes\gamma^{(j)})$.
As the boundary of
$\sum_{j=1}^M\mu(\vtilde^{(j)}\otimes\gamma^{(j)})$,
we obtain a $0$-chain
$\sum u_k\otimes P_k$,
where $u_k\in\nbigl_{P_k}$.
If $P_k=(r_ke^{\sqrt{-1}\theta_k})\in\Delta\setminus\{0\}$,
let $Z(P_k)$ be the segment connecting $P_k$ and $0$.
It induces a path in $\Deltatilde$,
which is also denoted by $Z(P_k)$.
We obtain the following rapid decay $1$-chain of
$(\nbigv,\nabla)$:
\[
 \mu\Bigl(
  \sum_{j=1}^M \vtilde^{(j)}\otimes\gamma^{(j)}
  \Bigr)
 =\sum_{j=1}^M\mu(\vtilde^{(j)}\otimes\gamma^{(j)})
+\sum u_k\otimes Z(P_k).
\]
The term
$\sum_{j=1}^M
\mu^{\add}(\vtilde^{(j)}\otimes\gamma^{(j)})
+\sum u_k\otimes Z(P_k)$
is called the additional term.
If $\sum_{j=1}^M\vtilde^{(j)}\otimes\gamma^{(j)}$
is a rapid decay $1$-cycle,
then 
$\mu\Bigl(
  \sum_{j=1}^M \vtilde^{(j)}\otimes\gamma^{(j)}
  \Bigr)$
is also a rapid decay $1$-cycle.
 
\subsection{Proof of the first claim of
Proposition \ref{prop;18.5.21.30}}
\label{subsection;20.10.25.1}

We take $\gminia\in \nbigitilde_{J,<0}$.
We explain the proof for $A_{J_+,\theta^u,\gminia}$
in the case $\Jbar\subset\vecI_x(\theta^u)$.
The proof for $A_{J_-,\theta^u,\gminia}$ is similar.
There exists
$\theta_1\in J\cap
 \Cr_2(\omega,\veckappa(\gminia,u))$.

There exists a large $C>0$ such that
$C^{-1}<\!<
 g_{\veckappa(\gminia,u)}(1,\theta_1)$.
By using the coordinate $x=z^{-1}$,
we set $\Gamma_{\theta_1}:=\gamma_v(C,0;\theta_1)$.
By modifying $\Gamma_{\theta_1}$
as in \S\ref{subsection;18.5.20.1},
we obtain
a family of paths $\Gamma_{\theta_1,\gminic(\gminia,ut)}$
for $g_{\veckappa(\gminia,u)}$.
By adding $\gamma_h(C;\theta_1,\vartheta^J_r-\pi)$,
we obtain a path
$\Gammatilde_{\theta_1,\gminic(\gminia,ut)}$
connecting
$(0,\theta_1)$
and $(C,\vartheta^J_r-\pi)$.
Any $v\in H^0(J_{+},L_{J_+,\gminia})$
induces a section $\vtilde$ 
along $\Gammatilde_{\theta_1,\gminic(\gminia,ut)}$.

There exists the decomposition
\[
 v=u_{(J-\pi)_+,0}
+\sum_{J-(1+\omega^{-1})\pi<J'\leq J-(1-\omega^{-1})\pi}
 u_{J'},
\]
where $u_{(J-\pi)_+,0}$
is a section of $L_{(J-\pi)_+,0}$,
and $u_{J'}$ are sections of $L_{J',<0}$.

If $J-(1+\omega^{-1})\pi<J'\leq J+(-1+\omega^{-1})\pi$,
we have $J'\cap (\vecI_x(\theta^u)-\pi)\neq\emptyset$.
We take $\theta_{J'}\in J'\cap (\vecI_x(\theta^u)-\pi)$.
Let $\Gamma_{J'}$ be the paths 
obtained as the union of
$\gamma_{h}(C;\vartheta^J_{r}+\pi,\theta_{J'})$
and 
$\gamma_v(C,0;\theta_{J'})$.
Then, we obtain the following continuous family of
rapid decay $1$-cycles
for $(V,\nabla)\otimes\nbige(zu^{-1})$
which represents $A_{J_+,\theta^u,\gminia}(v)$
in the case $\nbigv=V$:
\begin{multline}
\label{eq;18.5.21.40}
 \langle v\rangle_t:=
 \Bigl(
 \vtilde\otimes
 \gminis(\gminia,ut)
 \Gammatilde_{\theta_1,\gminic(\gminia,ut)}
+u_{(J-\pi)_+,0}\otimes 
 \gminis(\gminia,ut)
 \Gamma_{J-\pi}
 \\
+\sum_{J-(1+\omega^{-1})\pi<J'\leq J+(-1+\omega^{-1})\pi}
 u_{J'}\otimes
 \gminis(\gminia,ut)\Gamma_{J'}
\Bigr)
 \exp(-zu^{-1}).
\end{multline}

\begin{lem}
\label{lem;20.10.25.20}
We can divide the paths such that
the family {\rm (\ref{eq;18.5.21.40})}
is contained in
\[
 \nbigc^{(\infty)\,\varrho}_{\infty}\bigl(
 (V,\nabla),u,d(\omega),
-\Re \gminia^{\circ}(ut)
 \bigr).
 \]
\end{lem}
\pf
Let $h^V$ be a Hermitian metric of
$V_{|\cnum^{\ast}}$ adapted to the meromorphic structure
of $V$.
In the following, $C_1$ and $N_1$
denotes positive constants.
The following holds on 
$\gminis(\gminia,ut)
 \Gammatilde_{\theta_1,\gminic(\gminia,ut)}$:
\[
 \bigl|
 \vtilde
 \bigr|_{h^{V}}
 \exp\Bigl(
 -\Re(zu^{-1}t^{-1})
 \Bigr)
\leq
 C_1\exp\Bigl(
 -\Re\gminia^{\circ}(ut)
 \Bigr)t^{-N_1}.
\]
Moreover,
for any $\epsilon>0$,
there exist a small neighbourhood
$U_{\epsilon}$ of $(1,\theta_1)$
such that the following holds on 
$\gminis(\gminia,ut)\bigl(
 \Gammatilde_{\theta_1,\gminic(\gminia,ut)}
 \setminus U_{\epsilon}
 \bigr)$:
\[
 \bigl|
 \vtilde
 \bigr|_{h^{V}}
 \exp\Bigl(
 -\Re(zu^{-1}t^{-1})
 \Bigr)
\leq
 C_1\exp\Bigl(
 -\Re\gminia^{\circ}(ut)
-\epsilon t^{-\omega/(\omega-1)}
-\epsilon
 t^{-\omega/(\omega-1)}
 |x|^{-\omega}
 \Bigr).
\]
Note that 
$-t^{\omega/(\omega-1)}\Re\gminia^{\circ}(ut)$
converges to a positive number as $t\to 0$.

On $\gminis(\gminia,ut)\gamma_h(C;\theta_1,\vartheta^J_r-\pi)$,
we have the following:
\[
 \bigl|
 \vtilde
 \bigr|_{h^{V}}
 \exp\Bigl(
 -\Re(zu^{-1}t^{-1})
 \Bigr)
\leq
 C_1\exp\Bigl(
 -\Re\gminia^{\circ}(ut)
-\epsilon t^{-\omega/(\omega-1)}
 \Bigr).
\]
We have the following estimate
on $\gminis(\gminia,ut)\Gamma_{J-\pi}$
for some $\delta>0$ and $C_1>0$:
\[
 \bigl|
 u_{(J-\pi)_+,0}
 \bigr|_{h^{V}}
 \exp\Bigl(
 -\Re(zu^{-1}t^{-1})
 \Bigr)
 \leq
  C_1\exp\Bigl(
 -\Re\gminia^{\circ}(ut)
-\epsilon t^{-\omega/(\omega-1)}
-\delta
 t^{-1}|x|^{-1}
 \Bigr).
\]
We have the following estimate
on $\Gamma_{J'}$ for some $C_1>0$ and $\delta>0$:
\[
 \bigl|
  u_{J'}
 \bigr|_{h^{V}}
 \exp\Bigl(
 -\Re(zu^{-1}t^{-1})
 \Bigr)
\leq
 C_1\exp\Bigl(
 -\Re\gminia^{\circ}(ut)
-\epsilon t^{-\omega/(\omega-1)}
-\delta
 t^{-\omega/(\omega-1)}
 |x|^{-\omega}
 \Bigr).
\]
Thus, we obtain the claim of the lemma.
\hfill\qed

\vspace{.1in}

By applying the lifting procedure
in \S\ref{subsection;18.5.21.100}
to
the family of cycles $\langle v\rangle_t$
to obtain a family
$\mu(\langle v\rangle_t)$ of rapid decay $1$-cycles
for $(\nbigv,\nabla)\otimes\nbige(zu^{-1})$.
It represents $A_{J_+,\theta^u,\gminia}(v)$ for $\nbigv$.

We set $d(\omega)=(\omega-1)^{-1}$.
Note that the additional term
of $\mu(\langle v\rangle_t)$
is the sum of the $1$-chains of the form
$c\otimes \gamma_v(rt^{d(\omega)},0;\phi)
\cdot \exp(-zu^{-1})$,
where $c$ is a flat section of $\nbigv$
along $\gamma_v(rt^{d(\omega)},0;\phi)$
such that
$\bigl|
 c
 \bigr|_{h^{\nbigv}}
 \leq
 C\exp(-\delta |x|^{-\omega'})$
for some $C,\delta>0$ and $\omega'>\omega$.
Note that
$|x|^{-1}t^{-1}=O(|x|^{-\omega})$
and
$t^{-\omega/(\omega-1)}=O(|x|^{-\omega})$
on $\gamma_v(rt^{d(\omega)},0;\phi)$.
Hence, from Lemma \ref{lem;20.10.25.20},
we obtain that
$\mu(\langle v\rangle_t)$
is contained in
$\nbigc^{(\infty)\,\varrho}_{\infty}\bigl(
 (\nbigv,\nabla),u,d(\omega),
-\Re \gminia^{\circ}(ut)
 \bigr)$.
Thus, we obtain the first claim of
Proposition \ref{prop;18.5.21.30}.

\subsection{Proof of
the first claim of Proposition \ref{prop;18.5.21.31}}
\label{subsection;20.10.25.2}

We take $\gminia\in \nbigitilde_{J,<0}$.
We explain the proof for $\BB^{\rd}_{J_+,\theta^u,\gminia}$
in the case $\Jbar\subset \vecI_x(\theta^u)-\pi$.
The proof for $\BB^{\rd}_{J_-,\theta^u,\gminia}$ is similar.

There exists $\theta_1\in J$
such that
$\theta_1\in \Cr_2(\omega,\veckappa(\gminia,u))$.
We use the coordinate $x=re^{\sqrt{-1}\theta}=z^{-1}$.
Let $\Gamma_{\theta_1}$ be the path
$\gamma_h(1;\vartheta^J_{\ell}-\delta,\vartheta^J_{r})$.
By modifying it as in \S\ref{subsection;18.5.20.1},
we obtain a continuous family of paths
$\Gamma_{\theta_1,\gminic(\gminia,ut)}$.
By adding
$\gamma_v(1,0;\vartheta^J_{\ell}-\delta)$,
we obtain a family of paths
$\Gammatilde_{\theta_1,\gminic(\gminia,ut)}$
connecting
$(1,\vartheta^J_{r})$ and $(0,\vartheta^J_{\ell}-\delta)$.

Any $v\in H^0(J_-,L_{J_-,\gminia})$
induces a section $\vtilde$ 
along $\Gammatilde_{\theta_1,\gminic(\gminia,ut)}$.
There exists the decomposition
\[
 v=u_{J,0}+\sum_{J\leq J'\leq J+\omega^{-1}\pi}
 u_{J'},
\]
where $u_{J,0}$ is a section of $L_{J_+,0}$,
and $u_{J'}$ are sections of $L_{J',<0}$.
Let $\Gamma_{2,\pm}$ denote the paths
connecting
$(1,\vartheta^J_{r})$
and $(0,\vartheta^J_{r}\pm\delta)$,
obtained as the union of
$\gamma_h(1;\vartheta^J_r,\vartheta^J_{r}\pm\delta)$
and 
$\gamma_v(1,0;\vartheta^J_{r}\pm\delta)$.

We obtain the following family of cycles,
which represents $\BB^{\rd}_{J_+,\theta^u}(v)$
in the case $\nbigv=V$:
\begin{multline}
\label{eq;18.5.29.30}
\langle v\rangle_t=
\Bigl(
\vtilde\otimes
 \gminis(\gminia,ut)\Gammatilde_{\theta_1,\gminic(\gminia,ut)}
-u_{J,0}\otimes
 \gminis(\gminia,ut)\Gamma_{2,-}
-u_J\otimes
 \gminis(\gminia,ut)\Gamma_{2,-}
 \\
-\sum_{J<J'\leq J+\omega^{-1}\pi}
 u_{J'}\otimes
 \gminis(\gminia,ut)\Gamma_{2,+}
\Bigr)
 \exp(-zu^{-1}).
\end{multline}

\begin{lem}
\label{lem;20.10.25.21}
We can divide the paths such that
the family {\rm(\ref{eq;18.5.29.30})}
is contained in
\[
 \nbigc_{\infty}^{(\infty)\,\varrho}\bigl(
 (V,\nabla),u,d(\omega),
 -\Re\gminia^{\circ}(ut)
 \bigr)
\]
\end{lem}
\pf
Let $h^V$ be a Hermitian metric of $V_{|\cnum^{\ast}}$
which is adapted to the meromorphic structure of $V$.
On $\gminis(\gminia,ut)\Gamma_{\theta_1,\gminic(\gminia,ut)}$,
we have the following estimate
for some $C_1>0$ and $N_1>0$:
\[
 \bigl|
 \vtilde
 \bigr|_{h^{V}}
 \exp\bigl(-\Re(zu^{-1}t^{-1})\bigr)
\leq
 C_1\exp\Bigl(
 -\Re\gminia^{\circ}(ut)
 \Bigr)t^{-N_1}.
\]
Moreover,
for any $\epsilon>0$,
there exists a neighbourhood $U_{\epsilon}$
of $(1,\theta_1)$
such that the following estimate holds
on 
$\gminis(\gminia,ut)\bigl(
 \Gamma_{\theta_1,\gminic(\gminia,ut)}
 \setminus U_{\epsilon}
 \bigr)$
for some $C_1>0$:
\[
 \bigl|
 \vtilde
 \bigr|_{h^{V}}
 \exp\bigl(-\Re(zu^{-1}t^{-1})\bigr)
\leq
 C_1\exp\Bigl(
 -\Re\gminia^{\circ}(ut)
-\epsilon t^{-\omega/(\omega-1)}
 \Bigr).
\]
On $\gminis(\gminia,ut)\gamma_v(1,0;\vartheta^J_{\ell}-\delta)$,
we have the following
for some $C_1>0$ and $\epsilon>0$:
\[
 \bigl|
 \vtilde
 \bigr|_{h^{V}}
 \exp\bigl(-\Re(zu^{-1}t^{-1})\bigr)
\leq
 C_1\exp\Bigl(
 -\Re\gminia^{\circ}(ut)
-\epsilon t^{-\omega/(\omega-1)}
-\epsilon
 t^{-\omega/(\omega-1)}
 |x|^{-\omega}
 \Bigr).
\]
On $\gminis(\gminia,ut)\Gamma_{2,-}$,
we have the following estimates
for some $C_1>0$ and $\epsilon>0$:
\[
 \bigl|
 u_{J,0}
\bigr|_{h^{V}}
  \exp\bigl(-\Re(zu^{-1}t^{-1})\bigr)
\leq
  C_1\exp\Bigl(
 -\Re\gminia^{\circ}(ut)
-\epsilon t^{-\omega/(\omega-1)}
-\epsilon
 t^{-1}|x|^{-1}
 \Bigr),
\]
\[
 \bigl|
 u_{J}
 \bigr|_{h^{V}}
 \exp\bigl(-\Re(zu^{-1}t^{-1})\bigr)
\leq
 C_1\exp\Bigl(
 -\Re\gminia^{\circ}(ut)
-\epsilon t^{-\omega/(\omega-1)}
-\epsilon
 t^{-\omega/(\omega-1)}
 |x|^{-\omega}
 \Bigr).
\]
We have similar estimates for
$\bigl|
 u_{J'}
 \bigr|_{h^{V}}
 \exp\bigl(-\Re(zu^{-1}t^{-1})\bigr)$
on $\gminis(\gminia,ut)\Gamma_{2,+}$.
Hence, we obtain the claim of the lemma.
\hfill\qed

\vspace{.1in}

By applying the lifting procedure
in \S\ref{subsection;18.5.21.100}
to $\langle v\rangle_t$,
we obtain a continuous family of
rapid decay $1$-cycles
$\mu(\langle v\rangle_t)$
for $\nbigv\otimes\nbige(zu^{-1})$.
As in the last part of \S\ref{subsection;20.10.25.1},
by Lemma \ref{lem;20.10.25.21},
$\mu(\langle v\rangle_t)$
is contained in
$\nbigc_{\infty}^{(\infty)\,\varrho}\bigl(
 (V,\nabla),u,d(\omega),
 -\Re\gminia^{\circ}(ut)
 \bigr)$.
Thus, we obtain the first claim of
Proposition \ref{prop;18.5.21.31}.
 
\vspace{.1in}
 
By Lemma \ref{lem;20.10.25.22},
Proposition {\rm\ref{prop;18.5.21.30}}
and Proposition {\rm\ref{prop;18.5.21.31}}
are proved.
\hfill\qed

\subsection{Proof of Proposition \ref{prop;18.5.21.20}}
\label{subsection;20.10.25.3}

Let us explain the case $\alpha=\infty$
and for $C^{J_{1+}}_{\infty,\theta^u}(y)$.
The other cases can be argued similarly.
We set $\omega=\omega(\nbigv)$.

We describe $y$ by a cycle
$\vecc(t)\in \nbigc^{(\infty)\,\varrho}_{\infty}\bigl(
 \nbigstilde_{\omega}(\nbigv,\nabla),u,d,Q
 \bigr)$
as in (\ref{eq;18.5.17.20}).
We naturally regard 
$c_i\otimes\nu_{i,t}$ and
$b_i\otimes\eta_{i,t}$
are $1$-chains for 
$(\nbigv,\nabla)$.
Because $c_i$ are assumed to be sections of
a direct summand $\nbigv_{S,\gminia_i}$ 
for a splitting of Stokes filtrations,
the growth condition of $c_i\otimes\nu_{i,t}$
is satisfied also for $(\nbigv,\nabla)$.
Let us consider $b_i\otimes\eta_{i,t}$.
If $d_i=0$,
the growth condition of $b_i\otimes\eta_{i,t}$
is clearly satisfied for $(\nbigv,\nabla)$.
If $d_i\geq d$,
because
$\omega<\omega(\nbigstilde_{\omega}(\nbigv))\leq(d_i+1)/d_i$
is also assumed,
the growth condition for 
$b_i\otimes \eta_{i,t}$ is satisfied
for $(\nbigv,\nabla)$.

Let us consider
$a_i\otimes\gamma_{i,t}$.
If $\gamma_{i,t}=
 \gamma_v(t^{d_{i,1}}r_{i,1},t^{d_{i,2}}r_{i,2};\phi)$
with $r_{i,2}\neq 0$,
then 
$a_i\otimes\gamma_{i,t}$ is naturally regarded
as a $1$-cycle for $(\nbigv,\nabla)$.
We have
$|x|^{-\omega}=O\bigl(|x|^{-1}t^{-d_{i,2}(\omega-1)}\bigr)$
on $\gamma_{i,t}$.
Because $d_{i,2}(\omega-1)<1$,
the growth condition is also satisfied.
Let us consider the case
$\gamma_{i,t}=
 \gamma_v(t^{d_{i,1}}r_{i,1},0;\phi)$.
We take a small sector $S$ which contains $\gamma_{i,t}$.
Note that one of the following holds.
\begin{description}
\item[(A1)]
 $a_i$ is a section of
 $\nbigf^S_{<0}\nbigstilde_{\omega}^{\infty}(\nbigv)$.
\item[(A2)]
 $a_i$ is a section of
 $\nbigf^{S}_{\leq 0}\nbigstilde_{\omega}^{\infty}(\nbigv)$
 and $u^{-1}z<_{S}0$.
\end{description}
In the case {\bf (A1)},
$a_i\otimes\gamma_{i,t}$
is naturally regarded as
a family of cycles for $(\nbigv,\nabla)$,
and the growth condition is satisfied.

Let us consider the case {\bf (A2)}.
We shall replace $t_0$ with a smaller number if it is necessary.
We set $d(\omega)=(\omega-1)^{-1}>d_{i,1}$.
We take $0<r_0$
such that
$t_0^{d(\omega)}r_0<t_0^{d_{i,1}}r_{i,1}$.
In the following,
we shall replace $r_0$ with a larger number
if it is necessary,
which is possible by replacing $t_0$ with a small number.
We divide $\gamma_{i,t}$
into the union of
$\gamma_v(t^{d(\omega)}r_0,0;\phi)$
and 
$\gamma_v(t^{d_i,1}r_{i,1},t^{d(\omega)}r_0;\phi)$.
Here, we may assume that 
$\phi$ is contained in $\vecI_x(\theta^u)-\pi$.
We have the following
on $\gamma_v(t^{d_i,1}r_{i,1},t^{d(\omega)}r_0;\phi)$
for some $C_j>0$ and $\delta_j>0$:
\[
 |a_i|_{h^{\nbigv}}
 \exp\bigl(-\Re(x^{-1}u^{-1}t^{-1})\bigr)
\leq 
 C_1\exp\Bigl(
 Q(t)-\delta_1t^{-(1+d)}
 -\delta_2t^{-1}|x|^{-1}
 +C_2|x|^{-\omega}
 \Bigr).
\]
Because $d(\omega)\cdot\omega=1+d(\omega)$,
if $r_0$ is sufficiently large,
we have the following on 
$\gamma_v(t^{d_i,1}r_{i,1},t^{d(\omega)}r_0;\phi)$
for some $\delta_3>0$:
\[
 -\delta_2t^{-1}|x|^{-1}
 +C_2|x|^{-\omega}
\leq
 -\delta_3t^{-1}|x|^{-1}.
\]
We obtain the following
on $\gamma_v(t^{d_i,1}r_{i,1},t^{d(\omega)}r_0;\phi)$:
\[
 |a_i|_{h^{\nbigv}}
 \exp\bigl(-\Re(x^{-1}u^{-1}t^{-1})\bigr)
\leq 
 C_1\exp\Bigl(
-Q(t)-\delta_1t^{-(1+d)}
-\delta_3t^{-1}|x|^{-1}
 \Bigr).
\]
Hence, the chain 
$a_i\otimes \gamma_v(t^{d_i,1}r_{i,1},t^{d(\omega)}r_0;\phi)$
satisfies the desired growth condition.

We shall replace $a_i\otimes\gamma_{v}(t^{d(\omega)}r_0,0;\phi)$
with another chain in the following.
Let $J_1$ be as in the statement of Proposition \ref{prop;18.5.21.20}.
We study the case  $(J_1)_+\subset \vecI_x(\theta^u)-\pi$.

Recall
$(V,\nabla):=\nbigttilde^{\infty}_{\omega}(\nbigv,\nabla)$.
Let $(L,\vecnbigftilde)$ be the local system
corresponding to $(V,\nabla)$.
We obtain the section $[a_i]$ of $L$
induced by $a_i$.
We have the expression
\[
 [a_i]=u_{J_{1,+},0}
+\sum_{J_1-\omega^{-1}\pi<J'\leq J_1+\omega^{-1}\pi}
 u_{J'},
\]
where 
$u_{J_{1,+},0}$
is a section of $L_{J_{1+},0}$,
and 
$u_{J'}$ are sections of $L_{J',<0}$.

Let $\Gamma_1$ be the path 
$\gamma_h(t^{d(\omega)}r_0;\phi,\vartheta^{J_1}_{r})$.
For each $J'$ such that
$J_1-\omega^{-1}\pi<J'\leq J_1+\omega^{-1}\pi$,
we take $\theta_{J'}\in (\vecI_x(\theta^u)-\pi)\cap J'$.
Let $\Gamma_{J'}$ be the path
obtained as the union of
$\gamma_h(t^{d(\omega)}r_0;\vartheta^{J_1}_{r},
 \theta_{J'})$
and $\gamma_v(t^{d(\omega)}r_0,0;\theta_{J'})$.
We obtain the following family of $1$-chains
for $(V,\nabla) \otimes\nbige(x^{-1}u^{-1})$:
\begin{multline}
\label{eq;18.5.21.120}
 \langle
 a_i\exp(-zu^{-1})
  \otimes\gamma_v(t^{d(\omega)}r_0,0;\phi)
 \rangle
= \\
 \Bigl(
 [a_i]\otimes \Gamma_1
+u_{J_{1+},0}
 \otimes \Gamma_{J_1}
+\sum_{J_1-\pi/\omega<J'\leq J_1+\pi/\omega}
 u_{J'}\otimes \Gamma_{J'}
 \Bigr)
 \exp(-zu^{-1}).
\end{multline}

\begin{lem}
\label{lem;18.5.30.20}
The family of $1$-chains
$\langle
 a_i\exp(-zu^{-1})
  \otimes\gamma_v(t^{d(\omega)}r_0,0;\phi)
 \rangle
 $
is contained in 
$\nbigc_{\infty}^{(\infty)\,\varrho}
 \bigl(
 (V,\nabla),
 u,d,Q-\delta t^{-(1+d)}
 \bigr)$
for some $\delta>0$.
\end{lem}
\pf
Let $h^{V}$ be a Hermitian metric of $V_{|\cnum^{\ast}}$
adapted to the meromorphic structure of $V$.
In the following, $C_j$ and $\delta_j$
denote positive constants.
On $\Gamma_1$,
we have the following:
\[
\bigl|
 [a_{i}]
 \bigr|_{h^{V}}
 \exp\bigl(-\Re(x^{-1}u^{-1}t^{-1})\bigr)
\leq
 C_1\exp\Bigl(
 -\delta_5 r_0^{-1}t^{-1-d(\omega)}
+C_2 r_0^{-1-d(\omega)}t^{-1-d(\omega)}
 \Bigr).
\]
Hence, if $r_0$ is sufficiently large,
we have the following on $\Gamma_1$:
\[
 \bigl|
 [a_{i}]
 \bigr|_{h^{V}}
 \exp\bigl(-\Re(x^{-1}u^{-1}t^{-1})\bigr)
\leq
 C_1\exp\Bigl(
 -\delta_6 t^{-1-d(\omega)}
 \Bigr).
\]
Similarly,
if $r_0$ is sufficiently large,
we obtain
\[
\bigl|
 u_{J'}
 \bigr|_{h^{V}}
 \exp\bigl(-\Re(x^{-1}u^{-1}t^{-1})\bigr) 
\leq
 C_1\exp\bigl(
 -\delta_6 t^{-1-d(\omega)}
 \bigr)
\]
on $\gamma_h(t^{d(\omega)}r_0;\vartheta^J_{\ell},\theta_{J'})$,
and 
\[
\bigl|
 u_{J_{1,+},0}
 \bigr|_{h^{V}}
 \exp\bigl(-\Re(x^{-1}u^{-1}t^{-1})\bigr) 
\leq
 C_1\exp\bigl(
 -\delta_6 t^{-1-d(\omega)}
 \bigr)
\]
on $\gamma_h(t^{d(\omega)}r_0;\vartheta^J_{\ell},\theta_{J_1})$.

We have the following on 
$\gamma_v(t^{d(\omega)}r_0,0;\theta_{J'})$:
\[
 \bigl|
 u_{J'}
 \bigr|_{h^V}
 \exp\bigl(-\Re(x^{-1}u^{-1}t^{-1})\bigr)
\leq
  C_1\exp\bigl(
 -\delta_7|x|^{-1}|t^{-1}|
 -\delta_8|x|^{-\omega}
 \bigr).
\]
We have the following on 
$\gamma_v(t^{d(\omega)}r_0,0;\theta_{J_1})$:
\[
 \bigl|
 u_{J_{1+},0}
 \bigr|_{h^V}
 \exp\bigl(-\Re(x^{-1}u^{-1}t^{-1})\bigr)
\leq
  C_1\exp\bigl(
 -\delta_7|x|^{-1}|t^{-1}|
 \bigr).
\]
Then, we obtain the claim of 
Lemma \ref{lem;18.5.30.20}.
\hfill\qed

\vspace{.1in}

By applying the lifting procedure in \S\ref{subsection;18.5.21.100}
to the chain
$\langle
a_i\exp(-zu^{-1})
\otimes \gamma_v(t^{d(\omega)}r_0,0;\phi)
\rangle$,
we obtain a family of $1$-chains
$\mu\bigl(
\langle
a_i\exp(-zu^{-1})
\otimes \gamma_v(t^{d(\omega)}r_0,0;\phi)
\rangle
\bigr)$
for $(\nbigv,\nabla)\otimes\nbige(zu^{-1})$.
We may assume that
$\mu([a_i]\otimes\Gamma_1)_{|\Gamma_1(0)}=a_i$.
Then,
by replacing
$a_i\exp(-zu^{-1})
\otimes \gamma_v(t^{d(\omega)}r_0,0;\phi)$
with
$\mu\bigl(
\langle
a_i\exp(-zu^{-1})
\otimes \gamma_v(t^{d(\omega)}r_0,0;\phi)
\rangle
\bigr)$ for each $i$,
we obtain a family of $1$-cycles,
which represents
$C^{J_{1+}}_{\infty,\theta^u}(y)$.

The additional term of each
$\mu\bigl(
\langle
a_i\exp(-zu^{-1})
\otimes \gamma_v(t^{d(\omega)}r_0,0;\phi)
\rangle
\bigr)$
is of the form
$(c\otimes \gamma_v(rt^{d(\omega)},0;\phi))\exp(-zu^{-1})$,
where
 $c$ is a flat section of $\nbigv$
along $\gamma_v(rt^{d(\omega)},0;\phi)$
such that
$\bigl|
 c
 \bigr|_{h^{\nbigv}}
 \leq
 C\exp(-\delta |x|^{-\omega'})$
for some $C,\delta>0$ and $\omega'>\omega$.
Note that
$|x|^{-1}t^{-1}=O(|x|^{-\omega})$
and
$t^{-(1+d)}=O(|x|^{-\omega})$
on $\gamma_v(rt^{d(\omega)},0;\phi)$.
Hence, by Lemma \ref{lem;18.5.30.20},
we obtain that 
the family of $1$-cycles satisfies the estimate
as desired in Proposition \ref{prop;18.5.21.20}.
\hfill\qed

\section{Proof of Proposition \ref{prop;24.3.17.121}}
\label{section;24.3.17.120}

If $|u|$ is sufficiently small,
there exist the following commutative diagram of
isomorphisms for any $0<t\leq 1$,
as explained in
\S\ref{subsection;24.3.14.42}--\S\ref{subsection;24.3.18.1}.
\begin{equation}
 \begin{CD}
 H_1^{\varrho}\bigl(
 \cnum\setminus D,
 \nbigstilde_{1}(\nbigv,\nabla)\otimes\nbige(zu^{-1})
 \bigr)
 @>{a_1}>{\simeq}>  
 H_1^{\varrho}\bigl(
 \cnum\setminus D,
 (\nbigv,\nabla)\otimes\nbige(zu^{-1})
 \bigr)
  \\
  @V{b_{1,t}}V{\simeq}V @V{b_{2,t}}V{\simeq}V \\
 H_1^{\varrho}\bigl(
 \cnum\setminus D,
 \nbigstilde_{1}(\nbigv,\nabla)\otimes\nbige(z(tu)^{-1})
 \bigr)
 @>{a_t}>{\simeq}>  
 H_1^{\varrho}\bigl(
 \cnum\setminus D,
 (\nbigv,\nabla)\otimes\nbige(z(tu)^{-1})
 \bigr)
 \end{CD}
\end{equation}

Let $y\in 
 H_1^{\varrho}\bigl(
 \cnum\setminus D,
 \nbigstilde_1(\nbigv,\nabla)\otimes\nbige(zu^{-1})
 \bigr)$.

\begin{lem}
\label{lem;24.3.18.3}
Suppose that $y$ is represented by
a family of $1$-cycles
contained in 
$\nbigc^{(\alpha)\,\varrho}_{\infty}
 \bigl(
 \nbigstilde_1(\nbigv,\nabla)\otimes\nbige(zu^{-1}),
 u,d,Q
 \bigr)$ 
for some $\alpha\in D\cup\{\infty\}$.
Then, there exists $0<t_0\leq 1$ such that
$b_{2,t_0}\circ a_1$ is also represented by
a family of $1$-cycles
contained in 
$\nbigc^{(\alpha)\,\varrho}_{\infty}
 \bigl(
 \nbigv,\nabla\otimes\nbige(z(t_0u)^{-1}),
 u,d,Q
 \bigr)$.
\end{lem}
\pf
Let $\vecc(t)$ $(0\leq t\leq 1)$
be a family of $1$-cycles in
$\nbigc^{(\infty)\,\varrho}_{\infty}
 \bigl(
 \nbigstilde_1(\nbigv,\nabla)\otimes\nbige(zu^{-1}),
 u,d,Q
 \bigr)$ 
 which represents $y$,
as in (\ref{eq;18.5.17.20}).
We identify $\varpi^{-1}_{\Dtilde}(\infty)\simeq \real/2\pi\seisuu$
by the polar coordinate $x=re^{\sqrt{-1}\theta}=z^{-1}$. 
There exists a relatively compact interval
$I\subset
\openopen{-\theta^u-\pi/2}{-\theta^u+\pi/2}
\subset\varpi^{-1}_{\Dtilde}(\infty)$
such that 
any $\gamma_{i,t}=\gamma_v(r_{i,1},t^{d_{i,2}}r_{i,2};\phi_i)$
or $\gamma_{i,t}=\gamma_v(t^{d_{i,1}}r_{i,1},t^{d_{i,2}}r_{i,2};\phi_i)$
are contained in the sector corresponding to $I$,
i.e., $\phi_i\in I$ modulo $2\pi\seisuu$.
Let $\beta_1,\ldots,\beta_m$
be complex numbers such that
such that 
$\pi_1\nbigttilde_1(\nbigi_{\infty}(\nbigv))
=\bigl\{
\beta_1x^{-1},\beta_2x^{-1},\ldots,\beta_mx^{-1}\bigr\}$.
There exists $t_0$ such that
$\Bigl(
\bigl(
 (|u|t)^{-1}e^{-\sqrt{-1}\theta^u}-\beta_i
\bigr)e^{-\sqrt{-1}\theta}
\Bigr)>0$
for any $\theta\in I$.
For any $0<t<t_0$,
$c(t)\exp\bigl(-(t^{-1}-1)u^{-1}z\bigr)$
are $\varrho$-type $1$-cycles of
$(\nbigv,\nabla)\otimes\nbige(z(tu)^{-1})$.
Then, we obtain the claim of Lemma \ref{lem;24.3.18.3}
in the case $\alpha=\infty$.
The case $\alpha\in D$ can be argued similarly.
\hfill\qed

\vspace{.1in}

Let $f:\gbigl^{\gbigf}_{\varrho}(\nbigstilde_1\nbigv)
\to \gbigl^{\gbigf}_{\varrho}(\nbigv)$
denote the isomorphism of $2\pi\seisuu$-equivariant local systems
in (\ref{eq;24.3.17.110}).
By Proposition \ref{prop;24.3.18.4},
there exists a finite subset
$\ttS\subset
\bigl\{a\in\cnum\,\big|\,
|a|=1
 \bigr\}$
such that
$f_{\theta^u}(\nbigf^{\theta^u}_{\gminib})
\subset
 \nbigf^{\theta^u}_{\gminib}$
 for any $\gminib\in \nbigi(\Fourier_+(\nbigv))$
 unless $e^{\sqrt{-1}\theta}\in \ttS$.
By the comparison of the dimensions
of the associated graded spaces,
we obtain that 
$f_{\theta^u}(\nbigf^{\theta^u}_{\gminib})
=\nbigf^{\theta^u}_{\gminib}$ for any
$\gminib\in\nbigi(\Fourier_+(\nbigv))$
unless $e^{\sqrt{-1}\theta^u}\in \ttS$.
We obtain Proposition \ref{prop;24.3.17.121}
by Lemma \ref{lem;20.10.23.1}.
\hfill\qed

\chapter{Fourier transform of $D$-modules and Stokes structures}
\label{section;24.4.20.1}

\section{Holonomic $\nbigd$-modules on a punctured disc}

\subsection{Local description of holonomic $\nbigd$-modules}
\label{subsection;24.4.11.1}

Set $C:=\{z\in\cnum\,|\,|z|<1\}$.
Let $O$ denote the origin.
Set $V_0\nbigd_C:=\nbigo_C\langle z\del_z\rangle
\subset\nbigd_C$.
Let $\nbigm$ be any holonomic $\nbigd_C$-module
on $C$.

Let $\leq_{\cnum}$ be the total order on $\cnum$
induced by the lexicographic order
and the identification 
$\cnum\simeq\real\times\real$
obtained as $a+\sqrt{-1}b\longleftrightarrow(a,b)$.
According to Kashiwara and Malgrange,
there exists an increasing filtration $V_{\bullet}(\nbigm)$
indexed by $(\cnum,\leq_{\cnum})$
characterized by the following conditions.
\index{$V$-filtration}
\begin{itemize}
\item
Each $V_{\alpha}(\nbigm)$ is a coherent 
$V_0\nbigd_C$-submodule of $\nbigm$
such that
$V_{\alpha}(\nbigm)(\ast O)=\nbigm(\ast O)$.
\item
We have
$\nbigm=\bigcup_{\alpha\in\cnum}V_{\alpha}(\nbigm)$,
and
$V_{\alpha}(\nbigm)=
 \bigcap_{\beta>_{\cnum}\alpha}V_{\beta}(\nbigm)$
for any $\alpha$.  
\item
There exists a finite subset
$S\subset \{\alpha\in\cnum\,|\,0\leq \Re(\alpha)<1\}$
such that
\[
 \Gr^V_{\alpha}(\nbigm):=
 V_{\alpha}(\nbigm)\bigl/
\bigcup_{\beta<_{\cnum}\alpha}
 V_{\beta}(\nbigm)
\]
is $0$ unless $\alpha\in S+\seisuu$.
\item
We have
$zV_{\alpha}(\nbigm)\subset V_{\alpha-1}\nbigm$
for any $\alpha$.
Moreover,
$zV_{\alpha}(\nbigm)=V_{\alpha-1}(\nbigm)$     
for $\alpha<_{\cnum}0$     .
\item
We have
$\del_zV_{\alpha}(\nbigm)\subset V_{\alpha+1}(\nbigm)$
for any $\alpha\in\cnum$.
\item
The induced actions of
$\del_zz+\alpha$ on $\Gr^{V}_{\alpha}(\nbigm)$
are nilpotent.
\end{itemize}
It is easy to see that
$z:\Gr^V_{\alpha}(\nbigm)\simeq 
 \Gr^V_{\alpha-1}(\nbigm)$
is an isomorphism 
unless $\alpha=0$,
and 
$\del_z:\Gr^V_{\alpha}(\nbigm)\simeq 
 \Gr^V_{\alpha+1}(\nbigm)$
is an isomorphism unless $\alpha=-1$.
Let $\var$ and $\can$ denote the maps
$z:\Gr^V_0(\nbigm)\lrarr \Gr^V_{-1}(\nbigm)$
and
$-\del_z:\Gr^{V}_{-1}(\nbigm)\lrarr\Gr^V_0(\nbigm)$,
respectively.
By the construction,
$\var\circ\can$ is 
the nilpotent map $N$ on $\Gr^V_{-1}(\nbigm)$
induced by $-z\del_z=-\del_zz+1$.
It is easy to see that
$\nbigm\lrarr\nbigm(\ast O)$
induces
$\Gr^V_{\alpha}(\nbigm)
\simeq
\Gr^V_{\alpha}(\nbigm(\ast O))$
for any $\alpha\in\cnum\setminus\seisuu_{\geq 0}$.
The $V$-filtrations are functorial
in the sense that
a morphism of holonomic $\nbigd_C$-modules
$f:\nbigm_1\lrarr\nbigm_2$
preserves the $V$-filtrations,
i.e.,
$f(V_{\alpha}\nbigm_1)\subset V_{\alpha}\nbigm_2$.

Let $\nbigc$ denote the category of the tuples
$(\nbigv,\nbigq;a,b)$ as follows.
\begin{itemize}
\item
$\nbigv$ is a meromorphic flat bundle on $(C,O)$.
\item
$\nbigq$ is a $\cnum$-vector space.
\item
$a:\Gr^V_{-1}(\nbigv)\lrarr \nbigq$
and 
$b:\nbigq\lrarr\Gr^V_{-1}(\nbigv)$ 
are $\cnum$-linear maps
such that $b\circ a=\var\circ\can$.
\item
A morphism
$(\nbigv_1,\nbigq_1;a_1,b_1)
\lrarr
 (\nbigv_2,\nbigq_2;a_2,b_2)$
is defined as a tuple of morphisms
$f:\nbigv_1\lrarr\nbigv_2$
and $g:\nbigq_1\lrarr \nbigq_2$
such that
the following induced diagram is commutative:
\[
 \begin{CD}
 \Gr^{V}_{-1}(\nbigv_1)
 @>{a_1}>>
 \nbigq_1
 @>{b_1}>>
 \Gr^{V}_{-1}(\nbigv_1)
\\
 @V{\Gr^V_{-1}f}VV @V{g}VV @V{\Gr^V_{-1}f}VV \\
 \Gr^{V}_{-1}(\nbigv_2)
 @>{a_2}>>
 \nbigq_2
 @>{b_2}>>
 \Gr^{V}_{-1}(\nbigv_2).
 \end{CD}
\]
\end{itemize}

Let $\Hol(C,O)$ denote the category of holonomic
$\nbigd_C$-modules $\nbigm$
such that $\nbigm(\ast O)$ is a meromorphic flat bundle
on $(C,O)$.
\index{category $\Hol(C,O)$}
We obtain the functor
$\Psi:\Hol(C,O)\lrarr \nbigc$
defined as 
$\Psi(\nbigm)=
(\nbigm(\ast O),\Gr^V_0(\nbigm);\can,\var)$.
The following is well known due to 
Beilinson, Kashiwara and Malgrange.
\begin{prop}
\label{prop;25.3.12.70}
The functor $\Psi$ is an equivalence.
\hfill\qed
\end{prop}

For a meromorphic flat bundle $\nbigv$ on $(C,O)$,
the holonomic $\nbigd_C$-module
$\nbigv\in \Hol(C,O)$ corresponds to
$(\nbigv,\Gr^V_{-1}(\nbigv);\var\circ\can,\id)$,
and
$\nbigv(!O)\in \Hol(C,O)$ corresponds to
$(\nbigv,\Gr^V_{-1}(\nbigv);\id,\var\circ\can)$.
Hence,
$\Gr^V_{-1}(\nbigm)\stackrel{\can}{\lrarr}
 \Gr^V_0(\nbigm)\stackrel{\var}{\lrarr}
 \Gr^V_{-1}(\nbigm)$
are identified with
the morphisms:
\[
\Gr^V_{0}(\nbigm(! O))\lrarr
 \Gr^V_0(\nbigm)\lrarr
 \Gr^V_0(\nbigm(\ast O)).
\]

\subsection{Formal completion}
\label{subsection;24.4.14.3}

For any $\nbigo_C$-module $\nbigm$,
let $\nbigm_O$  denote the stalk of $\nbigm$ at $O$,
and let $\nbigm_{|\Ohat}$ denote the formal completion of $\nbigm$
at $O$,
i.e.,
$\nbigm_{|\Ohat}=\nbigm_O\otimes_{\nbigo_{C,O}}\cnum[\![z]\!]$.
\index{formal completion $\nbigm_{|\Ohat}$}

Let $(\nbigv,\nabla)$ be a meromorphic flat bundle on $(C,O)$.
There exist meromorphic flat bundles
$(\nbigv_i,\nabla)$ $(i=1,2)$
and an isomorphism
\begin{equation}
\label{eq;18.6.1.1}
 (\nbigv_{|\Ohat},\nabla)
\simeq
 (\nbigv_{1|\Ohat},\nabla)
\oplus
 (\nbigv_{2|\Ohat},\nabla)
\end{equation}
such that
(i) $(\nbigv_1,\nabla)$ is regular singular,
(ii) the set of ramified irregular values of
 $(\nbigv_2,\nabla)$ does not contain $0$.
The isomorphism (\ref{eq;18.6.1.1})
preserves the $V$-filtrations,
i.e.,
$V_a(\nbigv)_{|\Ohat}\simeq
V_a(\nbigv_1)_{|\Ohat}\oplus
V_a(\nbigv_2)_{|\Ohat}$.
We also have
$V_{\alpha}(\nbigv_2)=\nbigv_2$
for any $\alpha$.
Hence, 
we obtain the natural isomorphisms
\begin{equation}
\label{eq;18.6.1.2}
 \Gr^{V}_{\alpha}(\nbigv)
\simeq
 \Gr^V_{\alpha}(\nbigv_1).
\end{equation}

Let $(\nbigv,\nbigq;a,b)$ be an object in $\nbigc$.
By the isomorphisms (\ref{eq;18.6.1.2}),
we obtain the morphisms
\[
 \begin{CD}
 \Gr^V_{-1}(\nbigv_1)
@>{a_1}>>
 \nbigq
@>{b_1}>>
 \Gr^V_{-1}(\nbigv_1).
 \end{CD}
\]
There exists a regular holonomic $\nbigd_C$-module
$\nbigm_1$
corresponding to
$(\nbigv_1,\nbigq;a_1,b_1)$.

\begin{lem}
Let $\nbigm$ be a holonomic $\nbigd_C$-module
corresponding to $(\nbigv,\nbigq;a,b)$.
Then, there exists a natural isomorphism
$\nbigm_{|\Ohat}
\simeq
 \nbigm_{1|\Ohat}
\oplus
 \nbigv_{2|\Ohat}$.
\end{lem}
\pf
There exist the natural morphisms
$f:
 \nbigv_{O}\lrarr
 \nbigv_{|\Ohat}$
and $g:
 \nbigm_{1|\Ohat}\oplus
 \nbigv_{2|\Ohat}
 \lrarr
 \nbigv_{1|\Ohat}\oplus
 \nbigv_{2|\Ohat}=
 \nbigv_{|\Ohat}$.
We obtain the following morphism
induced by $f$ and $-g$:
\begin{equation}
\label{eq;18.6.1.3}
  \nbigv_{O}\oplus
 \Bigl(
  \nbigm_{1|\Ohat}\oplus
 \nbigv_{2|\Ohat}
 \Bigr)
\lrarr
 \nbigv_{|\Ohat}.
\end{equation}
We also obtain the induced morphisms
\begin{equation}
\label{eq;18.6.1.4}
 V_{\alpha}
 (\nbigv)_{O}
\oplus
 \Bigl(
  V_{\alpha}(\nbigm_{1|\Ohat})\oplus
 \nbigv_{2|\Ohat}
 \Bigr)
 \lrarr
 V_{\alpha}(\nbigv)_{|\Ohat}.
\end{equation}
If $\alpha<_{\cnum}0$,
then 
$V_{\alpha}(\nbigm_{1|\Ohat})\oplus
 \nbigv_{2|\Ohat}
\simeq
 V_{\alpha}(\nbigv)_{|\Ohat}$ holds.
It is easy to see that the induced morphisms
\[ 
 V_{\beta}(\nbigv)_O\big/
 V_{\alpha}(\nbigv)_O
\lrarr
 V_{\beta}(\nbigv)_{|\Ohat}\big/
 V_{\alpha}(\nbigv)_{|\Ohat}
\]
are isomorphism.
Hence, the morphisms (\ref{eq;18.6.1.3})
and (\ref{eq;18.6.1.4})
are surjective.
Let $\nbigk_O$ denote the kernel of (\ref{eq;18.6.1.3}),
which is a finitely generated $\nbigd_{C,O}$-module.
It induces a $\nbigd_{C'}$-module $\nbigk'$
on a neighbourhood $C'$ of $O$ in $C$
such that $\nbigk'_O=\nbigk_O$.
Note that
$\nbigk_O(\ast O)=\nbigv_O$,
and hence $\nbigk'(\ast O)\simeq\nbigv(\ast O)_{|C'}$.
Therefore, we obtain $\nbigk$ in $\Hol(C,O)$
whose stalk at $O$ is $\nbigk_O$.
It is equipped with a $V$-filtration
$V_{\alpha}(\nbigk)$ $(\alpha\in\cnum)$
such that each
 $V_{\alpha}(\nbigk)$ is isomorphic to
the kernel of (\ref{eq;18.6.1.4}).
By the construction,
there exists a natural isomorphism
\begin{equation}
\label{eq;18.6.1.5}
\nbigk_{|\Ohat}
\simeq
 \nbigm_{1|\Ohat}
\oplus
 \nbigv_{2|\Ohat}.
\end{equation}
By the isomorphism (\ref{eq;18.6.1.5}),
$\nbigk$ is an object in
$\Hol(C,O)$
corresponding to $(\nbigv,\nbigq;a,b)$.
Hence, there exists an isomorphism
$\nbigk\simeq\nbigm$.
It implies the claim of the lemma.
\hfill\qed

\subsection{Reduction with respect to the Stokes structure}
\label{subsection;24.4.15.11}

Let $\nbigm\in\Hol(C,O)$.
We obtain
$(\nbigv,\nbigq;a,b):=\Psi(\nbigm)$.
For $\omega\in\rnum_{>0}$,
we obtain the meromorphic flat bundle
$\nbigt_{\omega}(\nbigv)$ on $(C,O)$.
There exists the natural isomorphism
\[
 \Gr^{V}_{-1}(\nbigt_{\omega}(\nbigv))
\simeq
 \Gr^V_{-1}(\nbigv).
\]
Hence, we obtain the induced morphisms:
\[
\begin{CD}
 \Gr^V_{-1}(\nbigt_{\omega}(\nbigv))
@>{a_1}>>
 \nbigq
@>{b_1}>>
  \Gr^V_{-1}(\nbigt_{\omega}(\nbigv)).
\end{CD}
\]
We obtain a holonomic $\nbigd$-module
$\nbigt_{\omega}(\nbigm)$
corresponding to
$(\nbigt_{\omega}(\nbigv),\nbigq;a_1,b_1)$.
Note that
$\nbigt_{\omega}(\nbigv)(!O)
\simeq
 \nbigt_{\omega}\bigl(\nbigv(!O)\bigr)$.
 
\subsection{Beilinson functors and the gluing}
\label{subsection;24.4.14.1}

We set
$A^a=s^a\cnum[\![s]\!]$ $(a\in\seisuu)$
and
$A^{a,b}=A^a/A^b$ for $a\leq b$.
The multiplication of $s$
induces a nilpotent endomorphism $N_A$ of $A^{a,b}$.
For $\alpha\in\cnum$,
we set
\[
 \gbigi^{a,b}
 =\nbigo_{C}(\ast O)\otimes_{\cnum}A^{a,b}.
\]
It is equipped with the connection
defined by
$\nabla(g)=N_A(\alpha)dz/z$
for any $g\in A^{a,b}$.
We have natural morphisms
$\gbigi^{a,b}\to \gbigi^{c,d}$
for any $a\geq c$ and $b\geq d$
which are compatible with the connections.
We have the natural isomorphism
$\gbigi^{a,a+1}\simeq\nbigo_{C}(\ast O)$
by $s^a\longleftrightarrow 1$.

\subsubsection{Nearby cycle functor, maximal functor and gluing}

Let $\nbigm\in\Hol(C,O)$.
We set
$\Pi^{a,b}(\nbigm)
 =\nbigm\otimes\gbigi^{a,b}$.
We obtain the $\nbigd_C$-modules
$\Pi^{a,b}_{\star}(\nbigm)
 :=\Pi^{a,b}(\nbigm)(\star O)$
$(\star=!,\ast)$.
\index{$\nbigd$-modules $\Pi^{a,b}_{\star}(\nbigm)$}
We define
\[
 \Pi^{a,b}_{\ast,!}(\nbigm)
 :=\varprojlim_{N\to\infty}
 \Cok
 \Bigl(
 \Pi^{b,N}_{!}(\nbigm)
 \lrarr
 \Pi^{a,N}_{\ast}(\nbigm)
 \Bigr).
\]
There exists the following natural isomorphism:
\[
 \Pi^{a,b}_{\ast,!}(\nbigm)
 \simeq
 \varinjlim_{N\to\infty}
 \Ker\Bigl(
  \Pi^{-N,b}_{!}(\nbigm)
 \lrarr
 \Pi^{-N,a}_{\ast}(\nbigm)
 \Bigr).
\]
The following lemma is easy to see.
\begin{lem}
If $N$ is sufficiently large,
the natural morphisms
\[
  \Cok
 \Bigl(
 \Pi^{b,N+1}_{!}(\nbigm)
 \lrarr
 \Pi^{a,N+1}_{\ast}(\nbigm)
 \Bigr)
\lrarr
 \Cok
 \Bigl(
 \Pi^{b,N}_{!}(\nbigm)
 \lrarr
 \Pi^{a,N}_{\ast}(\nbigm)
 \Bigr),
\]
\[
  \Ker\Bigl(
  \Pi^{-N,b}_{!}(\nbigm)
 \lrarr
 \Pi^{-N,a}_{\ast}(\nbigm)
 \Bigr)
 \to
 \Ker\Bigl(
  \Pi^{-N-1,b}_{!}(\nbigm)
 \lrarr
 \Pi^{-N-1,a}_{\ast}(\nbigm)
 \Bigr)
\]
are isomorphisms.
\hfill\qed 
\end{lem}

The nearby cycle functor and the maximal functors are defined as
\index{nearby cycle functor $\psi^{(a)}$}
\index{maximal functor $\Xi$}
\[
 \psi^{(a)}(\nbigm)
 =\Pi^{a,a}_{\ast !} (\nbigm),
 \quad
 \Xi(\nbigm)
 =\Pi^{0,1}_{\ast !} (\nbigm).
\]
The multiplication of $s$ induces 
$\psi^{(a)}(\nbigm)\simeq\psi^{(a+1)}(\nbigm)$.
There exist the exact sequences:
\[
 0\lrarr
 \nbigm(!O)
 \stackrel{c_1}{\lrarr}
 \Xi(\nbigm)
 \stackrel{c_2}{\lrarr}
 \psi^{(0)}(\nbigm)
 \lrarr 0,
\]
\[
0\lrarr
 \psi^{(1)}(\nbigm)
 \stackrel{d_1}{\lrarr}
 \Xi(\nbigm)
 \stackrel{d_2}{\lrarr}
 \nbigm(\ast O)
 \lrarr 0. 
\]
The multiplication of $s$
and $c_2\circ d_1$
induce a nilpotent endomorphism of $\psi^{(1)}(\nbigm)$.

\subsubsection{Vanishing cycle functor}

There exist the natural morphisms
\[
 \nbigm(!O)
 \stackrel{e_1}{\lrarr}
 \nbigm
 \stackrel{e_2}{\lrarr}
 \nbigm(\ast O).
\]
Note that
$e_2\circ e_1
=d_2\circ c_1$.
We obtain the following complex:
\begin{equation}
\label{eq;24.4.13.1}
\begin{CD}
 \nbigm(!O)
 @>{c_1+e_1}>>
 \Xi(\nbigm)
 \oplus
 \nbigm
 @>{d_2-e_2}>>
 \nbigm(\ast O).
\end{CD}
\end{equation}
Beilinson defined the vanishing cycle functor
$\phi(\nbigm)$
as the cohomology of the complex (\ref{eq;24.4.13.1}).
\index{vanishing cycle functor $\phi$}
The morphisms
$d_1$
and $c_2$
induce $\can$ and $\var$:
\[
\begin{CD}
 \psi^{(1)}(\nbigm) 
 @>{\can}>>
 \phi(\nbigm)
 @>{\var}>>
 \psi^{(0)}(\nbigm).
\end{CD}
\]
We recall that $\nbigm$
is reconstructed as the cohomology of
the complex:
\begin{equation}
\label{eq;24.4.13.30}
\begin{CD}
 \psi^{(1)}(\nbigm)
 @>{d_1+\can}>>
 \Xi(\nbigm)
 \oplus
 \phi(\nbigm)
 @>{c_2-\var}>>
 \psi^{(0)}(\nbigm).
\end{CD}
\end{equation}

\subsubsection{Cohomology of the nearby cycle sheaf}

Let $\nbigv$ be a meromorphic flat bundle on $(C,O)$.
Let $\varpi:\Ctilde\to C$
denote the oriented real blow up along $O$.
There exist the sheaves
$\nbigl^{<0}(\Pi^{a,b}\nbigv)$
and
$\nbigl^{\leq 0}(\Pi^{a,b}\nbigv)$
on $\Ctilde$.
Let $\nbigc^{a,b}(\nbigv)$
denote the cokernel of the natural monomorphism
$\nbigl^{<0}(\Pi^{a,b}\nbigv)
\to
\nbigl^{\leq 0}(\Pi^{a,b}\nbigv)$.
\begin{lem}
For any sufficiently large $N$,
there exists the following natural isomorphism
\[
 H^1\bigl(
 C,
 \Omega^{\bullet}\otimes
 \psi^{(a)}(\nbigv)
 \bigr)
 \simeq
 H^1\bigl(
 \Ctilde,
 \nbigc^{a,N}(\nbigv)
 \bigr).
\] 
\end{lem}
\pf
There exists the natural isomorphism
\begin{multline}
 H^1\bigl(
 C,\Omega^{\bullet}\otimes
 \psi^{(a)}(\nbigv)
 \bigr)
 \simeq \\
 \Cok\Bigl(
 H^1\bigl(
 C,
 \Omega^{\bullet}\otimes
 \Pi^{a,N}_{!}(\nbigv)
 \bigr)
\to
 H^1\bigl(
 C,
 \Omega^{\bullet}\otimes
 \Pi^{a,N}_{\ast}(\nbigv)
 \bigr)
 \Bigr).
\end{multline}
We also have the natural isomorphisms
$H^1(C,\Omega^{\bullet}\otimes
 \Pi^{a,N}_{!}(\nbigm))
 \simeq
 H^1\bigl(
 \Ctilde,\nbigl^{<0}(\Pi^{a,N}\nbigm)
 \bigr)$
and 
$H^1(C,\Omega^{\bullet}\otimes
 \Pi^{a,N}_{\ast}(\nbigm))
 \simeq
 H^1\bigl(
 \Ctilde,\nbigl^{\leq 0}(\Pi^{a,N}\nbigm)
 \bigr)$.
Because
$H^1\bigl(
\Ctilde,
\nbigc^{a,N}(\nbigv)
\bigr)$
is isomorphic to the cokernel of
$H^1\bigl(
 \Ctilde,\nbigl^{<0}(\Pi^{a,N}\nbigm)
 \bigr)
 \lrarr
  H^1\bigl(
 \Ctilde,\nbigl^{\leq 0}(\Pi^{a,N}\nbigm)
 \bigr)$,
we obtain the desired isomorphism.
\hfill\qed

\subsubsection{Comparison of isomorphisms}
\label{subsection;24.4.13.20}

Let $\iota:\{O\}\to C$
denote the inclusion.
We recall that there exist the natural isomorphisms
\[
 \iota_{\dagger}
 \Gr^V_{-1}(\nbigv)
 \simeq
 \psi^{(a)}(\nbigv).
\]
Equivalently, 
$\Gr^V_{-1}(\nbigv)
\simeq
 \Gr^V_0(\psi^{(a)}(\nbigv))$.
Indeed, 
\begin{multline}
 \Gr^V_0\psi^{(a)}(\nbigv)
 =
\Cok\Bigl(
 \Gr^V_0\bigl(
 \Pi^{a,N}_{!}(\nbigv)
 \bigr)
 \to
 \Gr^V_0\bigl(
 \Pi^{a,N}_{\ast}(\nbigv)
 \bigr)
 \Bigr)
 \\
 =\Cok\Bigl(
 \Gr^V_{-1}(\nbigv)\otimes A^{a,N}
 \stackrel{s+N_A}{\lrarr}
 \Gr^V_{-1}(\nbigv)\otimes A^{a,N}
 \Bigr)
 \simeq
 \Gr^V_{-1}(\nbigv).
\end{multline}
We recall the decomposition (\ref{eq;18.6.1.1}).
Because there exists the natural isomorphism
$\Gr^V_{-1}(\nbigv)
\simeq
\Gr^V_{-1}(\nbigv_1)$,
we obtain
\begin{equation}
\label{eq;24.4.13.3}
\psi^{(a)}(\nbigv_1)\simeq
\psi^{(a)}(\nbigv).
\end{equation}
From (\ref{eq;24.4.13.3}),
we obtain
\begin{equation}
 f_1:
 H^1\bigl(C,
 \psi^{(a)}(\nbigv_1)
 \otimes\Omega^{\bullet}
 \bigr)
\simeq
 H^1\bigl(C,
  \psi^{(a)}(\nbigv)
  \otimes
  \Omega^{\bullet}
  \bigr).
\end{equation}

We also obtain the isomorphism
\begin{equation}
 f_2:
 H^1\bigl(C,
 \psi^{(a)}(\nbigv_1)
 \otimes\Omega^{\bullet}
 \bigr)
 \simeq
H^1\bigl(C,
  \psi^{(a)}(\nbigv)
  \otimes
  \Omega^{\bullet}
  \bigr)
\end{equation}
from the natural isomorphism
$\nbigc^{a,b}(\nbigv_1)\simeq
\nbigc^{a,b}(\nbigv)$.

\begin{lem}
\label{lem;24.4.14.21}
We have $f_1=f_2$.
\end{lem}
\pf
Because
$H^1(C,\Omega^{\bullet}\otimes\Pi^{a,N}_!\nbigv_1)=0$,
we have
$H^1(C,\psi^{(a)}(\nbigv_1)\otimes\Omega^{\bullet})
\simeq
H^1(C,\Pi^{a,N}\nbigv_1\otimes\Omega^{\bullet})$.

Let $\nbigc^{\infty}$ denote the sheaf of $C^{\infty}$-functions on $C$.
Let $\nbiga^{\bullet}(\Pi^{a,N}\nbigv)$
denote the Dolbeault resolution of
$\Pi^{a,N}\nbigv\otimes\Omega^{\bullet}$.
Let $\nbigp^{<O}$ denote the sheaf of $C^{\infty}$-functions on $C$
which are infinitely decay at $0$.
We set
$\nbiga^{\rd,\bullet}(\Pi^{a,N}\nbigv)\bigr)
=\nbigp^{<O}\otimes_{\nbigc^{\infty}}
 \nbiga^{\bullet}(\Pi^{a,N}\nbigv)$.
There exist the natural isomorphisms:
\[
 H^1\bigl(
 C,\Pi^{a,N}\nbigv\otimes\Omega^{\bullet}
 \bigr)
\simeq
 H^1\bigl(
 C,\nbiga^{\bullet}(\Pi^{a,N}\nbigv)\bigr),
\]
\[
 H^1\bigl(
 C,\Pi^{a,N}\nbigv[!O]\otimes\Omega^{\bullet}
 \bigr)
 \simeq
 H^1\bigl(
 C,\nbiga^{\rd,\bullet}(\Pi^{a,N}\nbigv)\bigr)
 \bigr).
\]

Let $\nbigl$ denote the local system on $\Ctilde$
associated with $\Pi^{a,N}\nbigv$.
We set $L_{S^1}=\nbigl_{|\varpi^{-1}(O)}$.
We obtain the constructible subsheaves
$L_{S^1}^{<0}\subset L_{S^1}^{\leq 0}\subset L_{S^1}$.
Let $q:\Ctilde\to\varpi^{-1}(0)$ be the projection
induced by the polar coordinate.
We obtain the constructible subsheaves
$q^{-1}(L_{S^1}^{<0})\subset
q^{-1}(L_{S^1}^{\leq 0})
\subset\nbigl$.

Let $\nbigl_1$ denote the local system on $\Ctilde$
associated with $\Pi^{a,N}\nbigv_1$.
We set $L_{1,S^1}:=\nbigl_{1|\varpi^{-1}(O)}$.
We have
$L_{1,S^1}=L_{S^1}^{\leq 0}/L_{S^1}^{<0}$
and
$\nbigl_1=q^{-1}(L_{S^1}^{\leq 0})/q^{-1}(L_{S^1}^{<0})$.

Let $I\subset S^1$ be any small interval
on which there exists a splitting
$L_{1,S^1}\to L_{S^1}^{\leq 0}$.
It induces a splitting
$\nbigl_{1}\to q^{-1}(L_{S^1}^{\leq 0})$
on the sector $q^{-1}(I)$.

Let $\tau$ be a holomorphic section of
$\Pi^{a,N}\nbigv_1\otimes\Omega^1$
For any small interval $I\subset S^1$,
by using a splitting as above,
we construct a holomorphic section $\tau_I$
of 
$q^{-1}(L^{\leq 0}_{S^1})\otimes \Omega^1$
on $q^{-1}(I)$
which induces $\tau_{|q^{-1}(I)}$.
Let $S^1=\bigcup I_i$ be a covering by sectors
such that $I_i\cap I_j$ $(i\neq j)$
do not include the Stokes directions of $\nbigv$.
Let $\chi_{i}$ be a partition of unity of $S^1$
subordinating the covering.
We obtain the $C^{\infty}$-section
$\tautilde_1
=\sum \chi_i \tau_{S_i}$
of $\Pi^{a,N}\nbigv\otimes\Omega^{1,0}$ on $C$.
We obtain
the $C^{\infty}$-section
$\delbar\tautilde_1$
of $\Pi^{a,N}\nbigv\otimes\Omega^{1,1}$ on $C$.
Note that on $C\setminus\{O\}$,
it is a section of
$q^{-1}(L_{S^1}^{<0})\otimes\Omega^{1,1}$.
By the integration in the radial direction,
we obtain
the $C^{\infty}$-section $\tautilde_2$
of $L^{<0}\otimes
\bigl(
 \Omega^{1,0}\oplus \Omega^{0,1}
 \bigr)$
such that
$\nabla\tautilde_2
+\delbar\tautilde_1=0$.
Note that $\tautilde_2$ induces a $C^{\infty}$-section of
$\Pi^{a,N}\nbigv\otimes(\Omega^{1,0}\oplus \Omega^{0,1})$
which is infinitely decay at $O$.
We obtain a $C^{\infty}$-section
$\tautilde=\tautilde_1+\tautilde_2$
of $\nbigv\otimes(\Omega^{1,0}\oplus\Omega^{0,1})$
such that the restriction to $C\setminus O$
is a $C^{\infty}$-section of
$q^{-1}(L^{\leq 0}_{S^1})\otimes(\Omega^{1,0}\oplus \Omega^{0,1})$
and that the induced section of
$\nbigl_1\otimes(\Omega^{1,0}\oplus\Omega^{0,1})$ equals $\tau$.
For another such $\tautilde'$,
$\kappatilde=\tautilde-\tautilde'$
is a $C^{\infty}$-section of
$\Pi^{a,N}\nbigv\otimes(\Omega^{1,0}\oplus\Omega^{0,1})$
satisfying
$\nabla(\kappatilde)=0$
which is infinitely decay at $O$.
Hence, we obtain the map
\begin{multline}
H^1(C,\Pi^{a,N}\nbigv_1\otimes\Omega^{\bullet})
\to
 \Cok\Bigl(
 H^1\bigl(C,
 \nbiga^{\rd\bullet}(\Pi^{a,N}\nbigv)
 \bigr)
\lrarr
 H^1\bigl(
 C,\nbiga^{\bullet}(\Pi^{a,N}\nbigv)
 \bigr)
 \Bigr)
 \\
\simeq\Cok\Bigl(
 H^1\bigl(
 C,\Omega^{\bullet}\otimes\Pi^{a,N}_{!}\nbigv
 \bigr)
\lrarr
 H^1\bigl(
 C,\Omega^{\bullet}\otimes\Pi^{a,N}\nbigv
 \bigr)
 \Bigr)
 \simeq
 H^1\bigl(
 C,
 \psi^{(a)}(\nbigv)\otimes\Omega^{\bullet}
 \bigr).
\end{multline}
It equals both $f_1$ and $f_2$.
\hfill\qed

\subsection{Regular holonomic $\nbigd$-modules and local systems}

We consider
the full subcategory 
$\Hol^{\reg}(C,O)\subset\Hol(C,O)$
of regular holonomic $\nbigd$-modules.
We fix a total order $\leq_{\cnum}$ on $\cnum$
as in \S\ref{subsection;24.4.11.1}.
We fix a subset
$T\subset\cnum\setminus\seisuu$
such that the projection
$\cnum\to \cnum/\seisuu$
induces a bijection
$T\simeq (\cnum\setminus\seisuu)/\seisuu$.
We assume that $\alpha<_{\cnum}0$ for any $\alpha\in T$
for simplicity.

\subsubsection{The nearby cycle functor and the vanishing cycle functor}
\label{subsection;25.4.5.1}

Let $\nbigm\in\Hol^{\reg}(C,O)$.
We set
\index{nearby cycle functor $\psitilde$}
\index{vanishing cycle functor $\phitilde$}
\[
 \psitilde(\nbigm):=
 \Gr^V_{-1}(\nbigm)
 \oplus
 \bigoplus_{\alpha\in T}
 \Gr^V_{\alpha-1}(\nbigm),
 \quad
 \quad
 \phitilde(\nbigm):=
 \Gr^V_{0}(\nbigm)
 \oplus
 \bigoplus_{\alpha\in T}
 \Gr^V_{\alpha}(\nbigm).
\]
Let $\cantilde:\psitilde(\nbigm)\to\phitilde(\nbigm)$
be the map defined by
$-\del_z:\Gr^V_{\beta}(\nbigm)\to\Gr^V_{\beta+1}(\nbigm)$.
\index{morphism $\cantilde$}
Let $\vartilde:\phitilde(\nbigm)\to\psitilde(\nbigm)$
be the map defined by
$z:\Gr^V_{\beta+1}(\nbigm)\to\Gr^V_{\beta}(\nbigm)$.
\index{morphism $\vartilde$}

We have the endomorphisms
$f_{\beta}$ on $\Gr^V_{\beta}(\nbigm)$
induced by $-z\del_z$.
Let $M_{\psitilde(\nbigm)}$
and $M_{\phitilde(\nbigm)}$
denote the automorphisms of
$\psitilde(\nbigm)$
and $\phitilde(\nbigm)$
obtained as the direct sum of
$\exp\bigl(2\pi\sqrt{-1}f_{\beta}\bigr)$.
They are called the monodromy automorphisms.

\subsubsection{The associated regular meromorphic flat bundles}
For any $\nbigm\in\Hol^{\reg}(C,O)$,
we set
\[
 \nbigv=
 \psitilde(\nbigm)\otimes\nbigo_C(\ast O).
\]
Let $f$ be the endomorphism
of $\psitilde(\nbigv)$
obtained as
\[
f=
 f_{-1}
 \oplus
 \bigoplus_{\alpha\in T}
 f_{\alpha-1}.
\]
We consider the connection
$\nabla=d-f\frac{dz}{z}$.
We recall the following lemma.
\begin{lem}
There exists an isomorphism
$\nbigm(\ast O)\simeq\nbigv$ 
which induces
$\psitilde(\nbigm)\simeq \psitilde(\nbigv)$.
\hfill\qed
\end{lem}

\subsubsection{The associated local systems and the nearby cycle functor}

We set $C^{\ast}=C\setminus\{O\}$.
By using the polar decomposition
$z=|z|e^{\sqrt{-1}\theta}$,
we obtain
$C^{\ast}=\openopen{0}{1}\times S^1$.
We define the map
$\varphi_{z,0}:\real\to C^{\ast}$
by $\varphi_{z,0}=\epsilon e^{\sqrt{-1}\theta}$
for any small positive number $\epsilon$.

Let $\nbigl(\nbigm)$ denote the local system on
$C^{\ast}$ obtained as
the sheaf of flat sections of
$\nbigm_{|C^{\ast}}$.
We obtain the $2\pi\seisuu$-equivariant
local system
$L_0(\nbigm)=\varphi_{z,0}^{-1}(\nbigl(\nbigm))$
on $\real$.
Let $M_{L_0(\nbigm)}$ denote
the monodromy automorphism of $L_0(\nbigm)$.
\index{local system $L_0(\nbigm)$}

Let $v_1,\ldots,v_r$ be a frame of $\nbigm(\ast O)$
such that the following holds.
\begin{itemize}
 \item There exists $\beta_i\in T\cup\{-1\}$
       such that
       $v_i\in V_{\beta_i}(\nbigm)$.
       Moreover,
       $\{v_i\,|\,\beta_i=\beta\}$
       induces a frame of
       $\Gr^V_{\beta}(\nbigm)$.
\end{itemize}
For $s\in H^0(\real,L_0(\nbigm))$,
there exist holomorphic functions
$g_{i,k}$ $(1\leq i\leq r,\,\,k\in\seisuu_{\geq 0})$ such that
\[
 s=\sum g_{i,k}(z)z^{\beta_i+1}(\log z)^kv_i.
\]
Here, $g_{i,k}=0$  for any sufficiently large $k$.
We obtain
\[
 \upsilon(s)\in \sum g_{i,0}(0)v_i
 \in\psitilde(\nbigm).
\]
The following lemma is well known.
\begin{lem}
The above procedure induces a well defined isomorphism
\[
 \upsilon:H^0(\real,L_0(\nbigm))
 \simeq
 \psitilde(\nbigm).
\]
Under the isomorphism,
we have
$M_{L_0(\nbigm)}=M_{\psitilde(\nbigm)}$.
\hfill\qed
\end{lem}

\subsubsection{Appendix: Topological vanishing cycle functor}
\index{nearby cycle functor $\psitilde^t$}
\index{vanishing cycle functor $\phitilde^t$}

We mention the topological vanishing cycle functor
for $\DR(\nbigm)$ though we do not use it.
See \cite{Malgrange-book} for more detail and precise.
Let $\pi:\Ctilde^{\ast}\to C^{\ast}$ denote a universal covering.
Let $j:C^{\ast}\to C$ denote the inclusion.
We obtain the sheaf
$j_{\ast}(\pi_{\ast}(\nbigo_{\Ctilde^{\ast}}))$.
Let $\iota:\{O\}\to C$ denote the inclusion.
We set
$\nbigotilde_O=\iota^{-1}\bigl(
j_{\ast}(\pi_{\ast}(\nbigo_{\Ctilde^{\ast}}))
\bigr)$.
Let $\nbigo_O$ denote the stalk of $\nbigo_C$ at $O$.
We also set
$\nbigctilde_O=\nbigotilde_O/\nbigo_O$.

Let $(\nbigm\otimes\Omega^{\bullet})_O$
denote the stalk of $\nbigm\otimes\Omega^{\bullet}$ at $O$.
We obtain the following:
\begin{equation}
\label{eq;25.4.5.1}
 \nbigm_O\otimes\nbigotilde_O
 \lrarr
 (\nbigm\otimes
 \Omega^1)_O\otimes\nbigotilde_O.
\end{equation}
\begin{equation}
\label{eq;25.4.5.2}
 \nbigm_O\otimes\nbigctilde_O
 \lrarr
 (\nbigm\otimes\Omega^1)_O\otimes\nbigctilde_O.
\end{equation}
The morphisms (\ref{eq;25.4.5.1})
and (\ref{eq;25.4.5.2}) are epimorphisms.
Let $\psitilde^t(\nbigm)$
and $\phitilde^t(\nbigm)$
denote the kernels of (\ref{eq;25.4.5.1})
and (\ref{eq;25.4.5.2}),
respectively.
It is easy to see that
$\psitilde^t(\nbigm)\simeq
H^0(\real,L_0(\nbigm))$.
The projection
$\nbigotilde_O\to\nbigctilde_O$
induces
$\can^t:\psitilde^t(\nbigm)
\to\phitilde^t(\nbigm)$.
There exists the automorphism
of $\varphi:\Ctilde^{\ast}$
given by $\log x\to\log x+2\pi$.
The pull back by $\varphi$
induces an automorphism $T:\nbigotilde_O\to\nbigotilde_O$.
The endomorphism $T-1:\nbigotilde_O\to\nbigotilde_O$ induces
a morphism $\nbigctilde_O\to \nbigotilde_O$.
It induces
$\var^t:\phitilde^t(\nbigm)\to\psitilde^t(\nbigm)$.

Let $G(t)=t^{-1}(e^{-2\pi\sqrt{-1}t}-1)$.
Let $F$ denote the endomorphism of
$\psitilde(\nbigm)$
induced by $z\del_z$,
which equals $-\vartilde\circ\cantilde$.
According to \cite{Malgrange-book},
there exist natural isomorphisms
\[
\psitilde(\nbigm)\simeq\psitilde^t(\nbigm),
\quad
\phitilde(\nbigm)\simeq\phitilde^t(\nbigm)
\]
for which the following diagram is commutative:
\[
 \begin{CD}
  \psitilde(\nbigm)
  @>{-\cantilde}>>
  \phitilde(\nbigm)
  @>{G(F)\circ\vartilde}>>
  \psitilde(\nbigm)
  \\
  @V{\simeq}VV @V{\simeq}VV @V{\simeq}VV \\
  \psitilde^t(\nbigm)
  @>{\can^t}>>
  \phitilde^t(\nbigm)
  @>{\var^t}>>
  \psitilde^t(\nbigm).
 \end{CD}
\]

\section{Monodromic regular holonomic $\nbigd$-modules}
\label{subsection;25.3.12.42}

Let $A$ be a finite dimensional vector space
equipped with an endomorphism $F$.
Let $\nbigv=A\otimes\nbigo_{\proj^1}(\ast\{0,\infty\})$
with the connection
$\nabla=d-F\frac{dz}{z}$.
Let $\nbigm$ be a regular holonomic $\nbigd_{\proj^1}$-modules
such that $\nbigm(\ast 0)=\nbigv$.

There exists the decomposition
\[
 (A,F)=(A^u,F^u)\oplus (A^{nu},F^{nu}),
\]
where any eigenvalue of $F^u$ is an integer,
and any eigenvalue of $F^{nu}$ is not an integer.
We have the corresponding decompositions
\[
 \nbigv=\nbigv^u\oplus \nbigv^{nu},
 \quad
 \nbigm=\nbigm^u\oplus\nbigm^{nu}.
\]
Moreover, we have $\nbigm^{nu}=\nbigv^{nu}$.

Let $S(F^{nu})$ denote the set of
the eigenvalues of $F^{nu}$.
We may assume the following conditions for $F$
without loss of generality.
\begin{itemize}
 \item For any two distinct elements
       $\alpha,\beta\in S(F^{nu})$,
       $\alpha-\beta$ is not an integer.
 \item $F^u$ is nilpotent.
\end{itemize}

\subsection{The generalized eigen decompositions and the $V$-filtrations}

\subsubsection{The generalized eigen decompositions}

There exists the generalized eigen decomposition
\[
 H^0(\proj^1,\nbigm)
 =\bigoplus_{\beta\in\cnum\setminus\seisuu}
 H^0(\proj^1,\nbigm)_{\beta}
\]
with respect to the action of $-z\del_z$,
i.e.,
$H^0(\proj^1,\nbigm)_{\beta}$
is the kernel of
$(-z\del_z-\beta)^m$ for a sufficiently large $m$.
Let $F_{\beta}$ denote the endomorphism of
$H^0(\proj^1,\nbigm)_{\beta}$
induced by $-z\del_z$.
There exist the generalized eigen decompositions
\[
 H^0(\proj^1,\nbigm^{nu})
 =\bigoplus_{\beta\in\cnum\setminus\seisuu}
 H^0(\proj^1,\nbigm^{nu})_{\beta},
 \quad
 H^0(\proj^1,\nbigm^{u})
 =\bigoplus_{n\in\seisuu}
 H^0(\proj^1,\nbigm^{u})_{n}.
\]

Let $N_{\nbigm}$ denote the nilpotent endomorphism
on $H^0(\proj^1,\nbigm^u)_{1}$
induced by $-\del_zz$,
and
let $N$ denote the nilpotent endomorphism
on $H^0(\proj^1,\nbigm^u)_{0}$
induced by $-z\del_z$.
Under the identification
$H^0(\proj^1,\nbigm^u)_{0}=A^u$,
we have $N=F^u$.

\subsubsection{$V$-filtrations}

There exists the $V$-filtration of $\nbigm$ along $0$.
There exists the decomposition
\[
 V_{\beta}(\nbigm)
 =
 V_{\beta}(\nbigm^u)
 \oplus
 V_{\beta}(\nbigm^{nu}).
\]

Because $-\del_zz=-z\del_z-1$,
there exists the natural isomorphism
\[
 H^0(\proj^1,V_{\beta}(\nbigm))
=\bigoplus_{\alpha\leq_{\cnum} \beta+1}
 H^0(\proj^1,\nbigm)_{\alpha}.
\]
(See \S\ref{subsection;24.4.11.1} for $\leq_{\cnum}$.)
In particular, we have
\[
\Gr^V_{\beta}(\nbigm)
\simeq
H^0(\proj^1,\nbigm)_{\beta+1}.
\]

\subsubsection{The nearby cycle functor and the vanishing cycle functor}
Recall
\index{nearby cycle functor $\psitilde$}
\index{vanishing cycle functor $\phitilde$}
\[
 \psitilde(\nbigm)
 :=\Gr^V_{-1}(\nbigm)
 \oplus
 \bigoplus_{\beta\in S(F^{nu})}
 \Gr^V_{\beta-1}(\nbigm),
\]
\[
  \phitilde(\nbigm)
 :=\Gr^V_{0}(\nbigm)
 \oplus
 \bigoplus_{\beta\in S(F^{nu})}
 \Gr^V_{\beta}(\nbigm).
\]

There exist the natural isomorphisms
\begin{equation}
\label{eq;25.3.12.10}
 \psitilde(\nbigm)
 \simeq
 H^0(\proj^1,\nbigm^u)_{0}
 \oplus
 \bigoplus_{\beta\in S(F^{nu})}
 H^0(\proj^1,\nbigm^{nu})_{\beta},
\end{equation}
\begin{equation}
\label{eq;25.3.12.11}
 \phitilde(\nbigm)
 \simeq
 H^0(\proj^1,\nbigm^u)_{1}
 \oplus
 \bigoplus_{\beta\in S(F^{nu})}
 H^0(\proj^1,\nbigm^{nu})_{\beta+1}. 
\end{equation}
We obtain the morphism
$\cantilde_{\nbigm}:\psitilde(\nbigm)\to\phitilde(\nbigm)$
induced by
$-\del_z$.
We also obtain the morphism
$\vartilde_{\nbigm}:\phitilde(\nbigm)\to\psitilde(\nbigm)$
induced by $z$.

\subsection{The associated local systems}
\label{subsection;25.3.12.40}

Let $\nbigl(\nbigm)$ denote the local system on
$\cnum^{\ast}$ obtained as
the sheaf of flat sections of
$\nbigm_{|\cnum^{\ast}}=\nbigv_{|\cnum^{\ast}}$.
By using the polar decomposition
$z=|z|e^{\sqrt{-1}\theta}$,
we obtain
$\cnum^{\ast}=\real_{>0}\times \real$.
We define the map
$\varphi_{z,0}:\real\to\cnum^{\ast}$
by $\varphi_{z,0}=e^{\sqrt{-1}\theta}$.
We obtain the $2\pi\seisuu$-equivariant
local system
$L_0(\nbigm)=\varphi_{z,0}^{-1}(\nbigl(\nbigm))$
on $\real$.
\index{local system $L_0(\nbigm)$}
Let $M_{L_0(\nbigm)}$ denote
the monodromy automorphism of $L_0(\nbigm)$.
\index{automorphism $M_{L_0(\nbigm)}$}

Let $c:\real\to\real$ be defined by $c(\theta)=-\theta$.
\index{map $c$}
We set
$L_{\infty}(\nbigm)=c^{-1}L_0(\nbigm)$.
\index{local system $L_{\infty}(\nbigm)$}
Let $M_{L_{\infty}(\nbigm)}$ denote
the monodromy automorphism of $L_{\infty}(\nbigm)$.
\index{automorphism $M_{L_{\infty}(\nbigm)}$}
We have
$M_{L_{\infty}(\nbigm)}=c^{-1}(M^{-1}_{L_0(\nbigm)})$.

Let $\varphi_{z,\infty}:\real\to \cnum^{\ast}$
be the map defined by $\varphi_{z,\infty}(\theta')=e^{\sqrt{-1}\theta'}$
with respect to the polar decomposition
$z^{-1}=|z|^{-1}e^{\sqrt{-1}\theta'}$.
We have
$L_{\infty}(\nbigv)=\varphi_{z,\infty}^{-1}(\nbigl)$.

For $\kappa=0,\infty$
there exist the generalized eigen decompositions
with respect to $M_{L_{\kappa}(\nbigm)}$:
\[
 H^0(\real,L_{\kappa}(\nbigm))
 =\bigoplus_{b\in\cnum^{\ast}}
 H^0(\real,L_{\kappa}(\nbigm))_b.
\]
There exists the isomorphisms
\[
 H^0(\real,L_{\kappa}(\nbigm^u))
 \simeq
 H^0(\real,L_{\kappa}(\nbigm))_1,
 \quad
  H^0(\real,L_{\kappa}(\nbigm^{nu}))
 =\bigoplus_{b\neq 1}
 H^0(\real,L_{\kappa}(\nbigm))_b.
\]
There exist the natural isomorphisms
\[
 H^0(\real,L_{\infty}(\nbigm))_b
 \simeq
  H^0(\real,L_{0}(\nbigm))_{b^{-1}}.
\]

\subsubsection{Isomorphisms}

For any $\beta\in\cnum\setminus\seisuu_{>0}$
and any $v\in H^0(\proj^1,\nbigm)_{\beta}$,
we obtain
\index{map $\rho_{z,\beta}$}
\[
 \rho_{z,\beta}(v)=\exp(F_{\beta}\log z)v
 \in H^0(\real,L_0(\nbigm)).
\]
It induces an isomorphism
$H^0(\proj^1,\nbigm)_{\beta}
 \simeq
 L_0(\nbigm)_{\exp(2\pi\sqrt{-1}\beta)}$.
The monodromy automorphism on
$L_0(\nbigm)_{\exp(2\pi\sqrt{-1}\beta)}$
equals
$\exp(2\pi\sqrt{-1}F_{\beta})$
under the isomorphism.

We obtain the isomorphism
\index{map $\rho_z$}
\[
 \rho_z:
 H^0(\proj^1,\nbigm^u)_{0}
 \oplus
 \bigoplus_{\beta\in S(F^{nu})}
 H^0(\proj^1,\nbigm^{nu})_{\beta}
\simeq
 H^0(\real,L_0(\nbigm)).
\]
As the composition of
(\ref{eq;25.3.12.10})
and $\rho_z$,
we obtain the isomorphism
\index{map $\rhotilde_z$}
\[
 \rhotilde_z:
 \psitilde(\nbigm)
 \simeq
 H^0(\real,L_0(\nbigm)).
\]

\subsection{Fourier transforms}

We consider the Fourier transforms
$\Fourier_{\pm}(\nbigm)$ of $\nbigm$,
which are regular holonomic $\nbigd$-modules
on $\proj^1$
such that
$\Fourier_{\pm}(\nbigm)(\ast 0)$
are meromorphic flat bundles on $(\proj^1,\{0,\infty\})$.
We have the decomposition
\[
 \Fourier_{\pm}(\nbigm)
 =
 \Fourier_{\pm}(\nbigm)^u
 \oplus
 \Fourier_{\pm}(\nbigm)^{nu}
=\Fourier_{\pm}(\nbigm^u)
 \oplus
 \Fourier_{\pm}(\nbigm^{nu}).
\]
There exist the natural isomorphisms
\index{maps $s_{\nbigm,\pm}$}
\[
 s_{\nbigm,\pm}:
 H^0(\proj^1,\nbigm)
\simeq
 H^0(\proj^1,\Fourier_{\pm}\nbigm).
\]
For $v\in H^0(\proj^1,\nbigm)$,
we have
$-\del_w(ws_{\nbigm,\pm}(v))
=s_{\nbigm,\pm}(z\del_zv)
=s_{\nbigm,\pm}(\del_z(zv)-v)$
and
$-w\del_w(s_{\nbigm,\pm}(v))
=s_{\nbigm,\pm}(\del_z (zv))
=s_{\nbigm,\pm}\bigl(
 z\del_zv+v
 \bigr)$.
Hence, we obtain the following isomorphisms
\[
 s_{\nbigm,\pm}:
 H^0(\proj^1,\nbigm^{nu})_{\beta}
 \simeq
 H^0(\proj^1,\Fourier_{\pm}\nbigm^{nu})_{-\beta+1}.
\]

\subsubsection{Comparison of the nearby cycle functors and the vanishing cycle functors}
\label{subsection;25.3.18.1}
There exist the natural isomorphisms:
\[
\psitilde(\Fourier_{\pm}(\nbigm)^u)
\simeq H^0(\proj^1,\Fourier_{\pm}(\nbigm^u))_{0}
 \simeq H^0(\proj^1,\nbigm^u)_{1}
 \simeq
 \phitilde(\nbigm^u),
\]
\[
 \phitilde(\Fourier_{\pm}(\nbigm)^u)
\simeq H^0(\proj^1,\Fourier_{\pm}(\nbigm^u))_{1}
 \simeq H^0(\proj^1,\nbigm^u)_{0}
 \simeq
 \psitilde(\nbigm^u).
\]
Because
$-w\del_ws_{\nbigm^{nu},\pm}(z^{-1}v)=
-s_{\nbigm^{nu},\pm}(-\del_zv)
=-s_{\nbigm^{nu},\pm}(z^{-1}F^{nu}(v))$,
there exist the natural isomorphisms
\[
 \Fourier_{\pm}(\nbigm^{nu})
 \simeq
 (z^{-1}A^{nu})\otimes\nbigo_{\proj^1}(\ast\{0,\infty\})
\]
under which
the connection equals
$d+F^{nu}\frac{dw}{w}$.
We set
$F^{nu,\gbigf}=-F^{nu}$.
We have
$S(F^{nu,\gbigf})=
\bigl\{
 -\beta\,\big|\,\beta\in S(F^{nu})
\bigr\}$.
There also exist the following isomorphisms:
\begin{multline}
 \psitilde(\Fourier_{\pm}(\nbigm)^{nu})
:=
 \bigoplus_{\gamma\in S(F^{nu,\gbigf})}
 \Gr^V_{\gamma-1}(\Fourier_{\pm}(\nbigm)^{nu})
 \\
 \simeq
 \bigoplus_{\gamma\in S(F^{nu,\gbigf})}
 H^0\bigl(
 \proj^1,\Fourier_{\pm}(\nbigm)^{nu}
 \bigr)_{\gamma}
 \simeq
 \bigoplus_{\beta\in S(F^{nu})}
 H^0\bigl(
 \proj^1,\nbigm^{nu}
 \bigr)_{\beta+1}
 \simeq
 \phitilde(\nbigm^{nu}).
\end{multline}
Similarly, we obtain the following isomorphisms:
\begin{multline}
 \phitilde(\Fourier_{\pm}(\nbigm)^{nu})
:=
 \bigoplus_{\gamma\in S(F^{nu,\gbigf})}
 \Gr^V_{\gamma}(\Fourier_{\pm}(\nbigm)^{nu})
 \\
 \simeq
 \bigoplus_{\gamma\in S(F^{nu,\gbigf})}
 H^0\bigl(
 \proj^1,\Fourier_{\pm}(\nbigm)^{nu}
 \bigr)_{\gamma+1}
 \simeq
 \bigoplus_{\beta\in S(F^{nu})}
 H^0\bigl(
 \proj^1,\nbigm^{nu}
 \bigr)_{\beta}
 \simeq
 \psitilde(\nbigm^{nu}).
\end{multline}
Therefore,
we obtain the following isomorphisms:
\[
 \phitilde(\nbigm)\simeq
 \psitilde(\Fourier_{\pm}(\nbigm)),
 \quad\quad
 \psitilde(\nbigm)\simeq
 \phitilde(\Fourier_{\pm}(\nbigm)).
\]
We have
\[
\cantilde_{\Fourier_{\pm}(\nbigm)}
=\mp \vartilde_{\nbigm},
\quad\quad
\vartilde_{\Fourier_{\pm}(\nbigm)}
=\pm\cantilde_{\nbigm}.
\]

\subsubsection{The induced isomorphisms}
\label{subsection;25.3.12.30}

We obtain the isomorphism
\index{isomorphisms $\Psi_{\nbigm,\pm}$}
\[
 \Psi_{\nbigm,\pm}:
\phitilde(\nbigm)
\simeq
H^0(\real,
L_{\infty}(\Fourier_{\pm}(\nbigm)))
\]
as
the composition of the following isomorphisms:
\begin{multline}
 \phitilde(\nbigm)
 \simeq
 \psitilde(\Fourier_{\pm}(\nbigm))
 \simeq
 H^0\Bigl(\real,
 L_{0}\bigl(
 \Fourier_{\pm}(\nbigm)
 \bigr)
 \Bigr)
\\
 \simeq
 H^0\Bigl(
 \real,
 L_{\infty}\bigl(
 \Fourier_{\pm}(\nbigm)
 \bigr)
 \Bigr).
\end{multline}

\subsection{Rapid decay and moderate growth homology}
\label{subsection;25.3.10.40}

Let $X=\realbar_{\geq 0}\times\real$
and $X^{\ast}=\real_{>0}\times\real$.
For $\theta^u\in\real$,
we consider paths
$\Gamma_{\star,\pm,\theta^u}$ $(\star=!,\ast)$
on $(X,X^{\ast})$.
\index{paths $\Gamma_{\star,\pm,\theta^u}$}
\begin{itemize}
 \item 
       $\Gamma_{!,+,\theta^u}$ is a path connecting
       $(\infty,\theta^u-2\pi)$ and $(\infty,\theta^u)$.
 \item $\Gamma_{\ast,+,\theta^u}$
       is a path connecting
       $(0,\theta^u)$ and $(\infty,\theta^u)$.
 \item $\Gamma_{!,-,\theta^u}$
       is a path connecting
       $(\infty,\theta^u-\pi)$ and $(\infty,\theta^u+\pi)$.
 \item $\Gamma_{\ast,-,\theta^u}$
       is a path connecting
       $(0,\theta^u+\pi)$ and $(\infty,\theta^u+\pi)$.
\end{itemize}

Let $\varpi:\projtilde^1\to\proj^1$
denote the oriented real blow up of $\proj^1$ at $\{0,\infty\}$.
Let $\varphi:\Xtilde\to \projtilde^1$
denote the map given by
$\varphi(r,\theta)=re^{\sqrt{-1}\theta}$.

We use the polar decomposition
$u=w^{-1}=|u|^{-1}\exp(\sqrt{-1}\theta^u)$.
Let $t\in H^0(\real,L_0(\nbigv))$.
We obtain the following
rapid decay $1$-cycles for
$\nbigv\otimes\nbige(\pm zw)$:
\[
 \varphi_{\ast}\bigl(
 t\cdot 
 \exp(\mp z u^{-1})
 \otimes
 \Gamma_{!,\pm,\theta^u}
 \bigr).
\]
By the isomorphisms
$\Fourier_{\pm}(\nbigv(!0))_{|u}
=H_1^{\rd}(\cnum^{\ast},\nbigv\otimes\nbige(zu^{-1}))$,
they induce sections of
$H^0(\real,L_{\infty}(\Fourier_{\pm}(\nbigv(!0))))$
denoted by
$\Abb^{\rd}_{\nbigv,\pm}(t)$.
We obtain the isomorphism
\index{isomorphisms $\Abb^{\rd}_{\nbigv,\pm}$}
\[
 \Abb^{\rd}_{\nbigv,\pm}:
 H^0(\real,L_0(\nbigv))
 \simeq
 H^0(\real,L_{\infty}(\Fourier_{\pm}\nbigv(\star 0))).
\]
We also obtain the following
moderate growth $1$-cycles for
$\nbigv\otimes\nbige(\pm zw)$:
\[
 \varphi_{\ast}\bigl(
 t\cdot 
 \exp(\mp z u^{-1})
 \otimes
 \Gamma_{\ast,\pm,\theta}
 \bigr).
\]
They induce the sections of
$H^0(\real,L_{\infty}(\Fourier_{\pm}(\nbigv)))$
denoted by
$\Abb^{\mg}_{\nbigv,\pm}(t)$.
We obtain the isomorphisms
\index{isomorphisms $\Abb^{\mg}_{\nbigv,\pm}$}
\[
 \Abb^{\mg}_{\nbigv,\pm}:
 H^0(\real,L_0(\nbigv))
 \simeq
 H^0(\real,L_{\infty}(\Fourier_{\pm}\nbigv)).
\]
The following lemma is clear
by the construction.
\begin{lem}
\label{lem;25.3.12.50}
Let $a:L_{\infty}(\Fourier_{\pm}(\nbigv[!0]))
\to L_{\infty}(\Fourier_{\pm}\nbigv)$
denote the natural morphism.
Then, $a\circ\Abb^{\rd}_{\nbigv,\pm}
=\Abb^{\mg}_{\nbigv,\pm}
 \circ(\id-M_{L_0(\nbigv)}^{-1})$.
\hfill\qed
\end{lem}

\subsection{Some endomorphisms}
\label{subsection;25.3.12.41}

We set
$\Gamma_{\star,\pm}:=
\Gamma_{\star,\pm,0}$.
\index{paths $\Gamma_{\star,\pm}$}

For any $\beta\in\cnum$,
let $F_{\beta}$
denote the endomorphism of
$\Gr^V_{\beta-1}(\nbigm)
=H^0(\proj^1,\nbigm)_{\beta}$
induced by $-z\del_z$.
We define the endomorphisms
$\Phi_{\beta,!,\pm}$
on $H^0(\proj^1,\nbigm)_{\beta}$
by
\begin{equation}
 \Phi_{\beta,!,\pm}
 =
 \frac{-1}{2\pi\sqrt{-1}}
 \int_{\Gamma_{!,\pm}}
 \exp\bigl(F_{\beta}\log \zeta\bigr)
 e^{\mp\zeta}\frac{d\zeta}{\zeta}.
\end{equation}
For $\beta\in\cnum\setminus\seisuu_{<0}$,
we define the endomorphisms
$\Phi_{\beta,\ast,\pm}$
on $H^0(\proj^1,\nbigm)_{\beta}$
by
\begin{equation}
 \Phi_{\beta,\ast,\pm}
=\frac{-(\pm 1)^{n-1}}{2\pi\sqrt{-1}}\int_{\Gamma_{\ast,\pm}}
 \exp\bigl((F_{\beta}+n\id)\log \zeta\bigr)
\prod_{j=1}^{n}
(F_{\beta}+j\id)^{-1}
 e^{\mp\zeta}\frac{d\zeta}{\zeta}.
\end{equation}
Here, $n$ denotes any non-negative integer such that
$\Re(\beta)+n>-1$.
The endomorphisms $\Phi_{\beta,\ast,\pm}$
are independent of the choice of $n$.
If $\Re(\beta)>-1$,
then
\[
 \Phi_{\beta,\ast,\pm}
 =\frac{\mp 1}{2\pi\sqrt{-1}}
 \int_{\Gamma_{\ast,\pm}}
 \exp\bigl(F_{\beta}\log \zeta\bigr)
 e^{\mp\zeta}d\zeta.
\]

We obtain the following endomorphisms
of $\psitilde(\nbigm)$:
\index{endomorphisms \mbox{$\Phi_{\star,\pm}$}}
\[
 \Phi_{\star,\pm}
 =
 \Phi_{0,\star,\pm}
 \oplus
 \bigoplus_{\beta\in S(F^{nu})}
 \Phi_{\beta,\star,\pm}.
\]

\subsection{Statements}
We explain some results which will be proved in 
\S\ref{subsection;25.3.12.23}--\S\ref{subsection;25.3.12.24}.

\subsubsection{Commutative diagrams}
\label{subsection;25.3.11.210}

We set
$\gbigl^{\gbigf}_{\pm}(\nbigm)
:=L_{\infty}(\Fourier_{\pm}(\nbigm))$
to simplify the notation.
We obtain the following proposition
from
Proposition \ref{prop;25.3.10.2},
Proposition \ref{prop;25.3.11.201},
Proposition \ref{prop;25.3.10.110},
and 
Proposition \ref{prop;25.3.11.202}
below.

\begin{prop}
\label{prop;25.3.11.200}
The endomorphisms
$\Phi_{\star,\pm}$  are invertible.
Moreover,
the following diagrams are commutative:
\begin{equation}
\begin{CD}
 \psitilde(\nbigm)
 @>{\cantilde\circ \Phi_{!,\pm}}>>
 \phitilde(\nbigm)
 @>{(\Phi_{\ast,\pm})^{-1}\circ\vartilde_{\nbigm}}>>
 \psitilde(\nbigm)
 \\
 @V{\simeq}V{\Abb^{\rd}_{\nbigv,\pm}\circ\rhotilde_z}V
 @V{\simeq}V{\Psi_{\nbigm,\pm}}V
 @V{\simeq}V{\Abb^{\mg}_{\nbigv,\pm}\circ\rhotilde_z}V \\
 H^0(\real,\gbigl^{\gbigf}_{\pm}(\nbigv(!0)))
 @>>>
 H^0(\real,\gbigl^{\gbigf}_{\pm}(\nbigm))
 @>>>
 H^0(\real,\gbigl^{\gbigf}_{\pm}(\nbigv)).
\end{CD}
\end{equation}
Here, the lower horizontal arrows are the natural morphisms.
The monodromy automorphisms of
$H^0(\real,\gbigl^{\gbigf}_{\pm}(\nbigv(!0)))$,
$H^0(\real,\gbigl^{\gbigf}_{\pm}(\nbigm))$
and
$H^0(\real,\gbigl^{\gbigf}_{\pm}(\nbigv))$
are equal to
$M_{\psitilde(\nbigm)}$,
$M_{\phitilde(\nbigm)}$,
and 
$M_{\psitilde(\nbigm)}$,
respectively.
(See {\rm\S\ref{subsection;25.4.5.1}} for the notation.)
\end{prop}

We obtain the following proposition from 
Corollary \ref{cor;25.3.12.20} below.
\begin{prop}
\label{prop;25.3.11.150}
The composition of the natural morphisms
\begin{multline}
 H^0(\real,\gbigl^{\gbigf}_{\pm}(\nbigm))
 \to
 H^0(\real,\gbigl^{\gbigf}_{\pm}(\nbigv))
 \simeq
 \psitilde(\nbigm)
 \simeq
 H^0(\real,\gbigl^{\gbigf}_{\pm}(\nbigv(!0)))
 \\
 \to
 H^0(\real,\gbigl^{\gbigf}_{\pm}(\nbigm))
\end{multline}
equal $\id-M^{-1}$,
where $M$ denote the monodromy automorphisms of
$\gbigl^{\gbigf}_{\pm}(\nbigm)$.
\end{prop}

\subsubsection{Some isomorphisms}

Let $a_{\ast,\pm}$ denote the automorphisms of
$\psitilde(\nbigv)$
obtained as the composition of the following:
\[
\begin{CD}
 \psitilde(\nbigv)
 @>{\vartilde^{-1}}>{\simeq}>
 \phitilde\bigl(\nbigv\bigr)
 @>>{\simeq}>
 H^0(\real,\gbigl^{\gbigf}_{\pm}(\nbigv))
 @>{(\Abb^{\mg}_{\nbigv,\pm}\circ\rhotilde_z)^{-1}}>{\simeq}>
 \psitilde(\nbigv).
\end{CD}
\]
Let $a_{!,\pm}$ denote the automorphisms of
$\psitilde(\nbigv(!0))$
obtained as the composition of the following:
\[
\begin{CD}
 \psitilde(\nbigv(!0))
 @>{\cantilde}>{\simeq}>
 \phitilde\bigl(\nbigv(!0)\bigr)
 @>>{\simeq}>
 H^0(\real,\gbigl^{\gbigf}_{\pm}(\nbigv(!0)))
 @>{(\Abb^{\rd}_{\nbigv,\pm}\circ\rhotilde_z)^{-1}}>{\simeq}>
 \psitilde(\nbigv(!0)).
\end{CD}
\]
Because $\psitilde(\nbigv)=\psitilde(\nbigv(!0))$,
we obtain the automorphisms
$a_{\ast,\pm}\circ a_{!,\pm}^{-1}$
of $\psitilde(\nbigm)$.
We set $G(t)=t^{-1}(1-e^{-2\pi\sqrt{-1}t})$.
\begin{cor}
$a_{\ast,\pm}\circ a_{!,\pm}^{-1}$
equals the direct sum of $G(F_{\beta})$.
\end{cor}
\pf
We have
$a_{\star,\pm}=\Phi_{\star,\pm}^{-1}$
for $\star=!,\ast$.
Hence, $a_{\ast,\pm}\circ a_{!,\pm}^{-1}$
is the direct sum of
$\Phi_{\beta,\ast,\pm}\circ
\Phi_{\beta,!,\pm}^{-1}
=G(F_{\beta})$.
\hfill\qed

\subsubsection{The inversion}

There exist the natural isomorphisms
\[
\Fourier_{\pm}\circ\Fourier_{\mp}(\nbigm)\simeq \nbigm.
\]
We have
$s_{\Fourier_{\pm}(\nbigm),\mp}\circ s_{\nbigm,\pm}
=\id$.
We set $\vecF_{\pm}(\nbigv)=\Fourier(\nbigv)(\ast 0)$.
We obtain the following proposition from 
Proposition \ref{prop;25.3.12.21}
and Proposition \ref{prop;25.3.12.22} below.
\begin{prop}
\label{prop;25.3.13.2}
On $H^0(\real,L_0(\nbigv(!0)))$,
we have
\begin{equation}
(c^{-1}\circ\Abb^{\mg}_{\vecF_-\nbigv,+})
\circ
(c^{-1}\circ\Abb^{\rd}_{\nbigv,-})
=
-(2\pi\sqrt{-1})^{-1}\id,
\end{equation}
\begin{equation}
 (c^{-1}\circ\Abb^{\mg}_{\vecF_+\nbigv,-})
  \circ
  (c^{-1}\circ\Abb^{\rd}_{\nbigv,+})\circ\rho_z
  =(2\pi\sqrt{-1})^{-1}
  M_{L_0(\nbigv(!0))}^{-1}.
\end{equation}
On $H^0(\real,L_0(\nbigv))$,
we have 
\begin{equation}
(c^{-1}\circ\Abb^{\rd}_{\vecF_-\nbigv,+})
\circ
(c^{-1}\circ\Abb^{\mg}_{\nbigv,-})
=(2\pi\sqrt{-1})^{-1}\cdot M_{L_0(\nbigv)},
\end{equation}
\begin{equation}
(c^{-1}\circ\Abb^{\rd}_{\vecF_+\nbigv,-})
\circ
(c^{-1}\circ\Abb^{\mg}_{\nbigv,+})
=-(2\pi\sqrt{-1})^{-1}\id.
\end{equation}
\end{prop}

\begin{cor}
\label{cor;25.3.15.1}
On $H^0(\real,L_{\infty}(\nbigv(!0)))$,
we have
\begin{equation}
(\Abb^{\mg}_{\vecF_-\nbigv,+}\circ c^{-1})
\circ
(\Abb^{\rd}_{\nbigv,-}\circ c^{-1})
=
-(2\pi\sqrt{-1})^{-1}\id,
\end{equation}
\begin{equation}
 (\Abb^{\mg}_{\vecF_+\nbigv,-}\circ c^{-1})
  \circ
  (\Abb^{\rd}_{\nbigv,+}\circ c^{-1})
  =(2\pi\sqrt{-1})^{-1}
  M_{L_{\infty}(\nbigv(!0))}.
\end{equation}
On $H^0(\real,L_0(\nbigv))$,
we have 
\begin{equation}
(\Abb^{\rd}_{\vecF_-\nbigv,+}\circ c^{-1})
\circ
(\Abb^{\mg}_{\nbigv,-}\circ c^{-1})
=(2\pi\sqrt{-1})^{-1}\cdot M_{L_{\infty}(\nbigv)}^{-1},
\end{equation}
\begin{equation}
(\Abb^{\rd}_{\vecF_+\nbigv,-}\circ c^{-1})
\circ
(\Abb^{\mg}_{\nbigv,+}\circ c^{-1})
=-(2\pi\sqrt{-1})^{-1}\id.
\end{equation}
\end{cor}

\subsubsection{Complement to Proposition \ref{prop;25.3.11.200}}
\label{subsection;25.3.13.1}

We consider morphisms of regular holonomic $\nbigd$-modules
$\nbigm_1
\to
\nbigm
\to
\nbigm_2$
satisfying the following conditions.
\begin{description}
 \item[(a)] $\nbigm(\ast 0)$ and
	    $\nbigm_i(\ast 0)$ are meromorphic flat bundles
       on $(\proj^1,\{0,\infty\})$.
 \item[(b)] The kernel and the cokernel of the morphisms
       are flat bundles.
 \item[(c)] $\nbigm_1(\ast 0)=\nbigm_1$
	    and $\nbigm_2(!0)=\nbigm_2$.
\end{description}
The condition (b) is equivalent to the following condition.
\begin{itemize}
 \item The induced morphisms
       $\phitilde(\nbigm_1)\to\phitilde(\nbigm)\to\phitilde(\nbigm_2)$
       are isomorphisms.
\end{itemize}

We obtain the induced morphisms
of $2\pi\seisuu$-equivariant local systems
\[
 L_{\infty}(\nbigm_1)
 \stackrel{a}{\lrarr} L_{\infty}(\nbigm)
 \stackrel{b}{\lrarr}
 L_{\infty}(\nbigm_2).
\]

By applying $\Fourier_-$
to $\nbigm_1\to\nbigm\to\nbigm_2$,
we obtain morphisms 
$\nbign_1\to \nbign\to\nbign_2$
of regular holonomic $\nbigd$-modules.
By the inversion,
we recover
$\nbigm_1\to\nbigm\to\nbigm_2$
by applying 
$\Fourier_+$ to $\nbign_1\to\nbign\to\nbign_2$.
In particular, we have $\Fourier_+(\nbign)=\nbigm$
and $\Fourier_+(\nbign_i)=\nbigm_i$.

\begin{lem}
The following holds.
\begin{itemize}
 \item $\nbign(\ast 0)$ and $\nbign_i(\ast 0)$
       are meromorphic flat bundles
       on $(\proj^1,\{0,\infty\})$.
 \item $\nbign_1(!0)=\nbign_1$
       and $\nbign_2(\ast 0)=\nbign_2$.
 \item The induced morphisms
       $\nbign_1(\ast 0)\to \nbign(\ast 0)\to\nbign_2(\ast 0)$
       are isomorphisms.
\end{itemize}
We shall identify
$\nbign_1=\nbign(!0)$
and $\nbign_2=\nbign(\ast 0)$.
\hfill\qed
\end{lem}

We have the isomorphisms
$\Abb^{\rd}_{+}:
H^0(\real,L_0(\nbign(!0)))\simeq
H^0(\real,L_{\infty}(\nbigm_1))$
and
$\Abb^{\mg}_+:
H^0(\real,L_0(\nbign(\ast 0)))
\simeq
H^0(\real,L_{\infty}(\nbigm_2))$.
We have the natural isomorphism
$L_0(\nbign(! 0))\simeq L_0(\nbign(\ast 0))$.
Hence, we obtain
\begin{equation}
\label{eq;25.3.5.40}
 L_{\infty}(\nbigm_1)\simeq L_{\infty}(\nbigm_2).
\end{equation}

\begin{prop}
\label{prop;25.3.5.51}
We obtain $b\circ a=\id-M_{L_{\infty}(\nbigm_1)}^{-1}$
and $a\circ b=\id-M_{L_{\infty}(\nbigm)}^{-1}$
under the isomorphism {\rm(\ref{eq;25.3.5.40})}.
\end{prop}
\pf
We obtain $b\circ a=\id-M_{L_{\infty}(\nbigm_1)}^{-1}$
from the relation between
$\Abb^{\rd}_{\pm}$ and $\Abb^{\mg}_{\pm}$.
We obtain 
$a\circ b=\id-M_{L_{\infty}(\nbigm)}^{-1}$
from Proposition \ref{prop;25.3.11.150}.
\hfill\qed

\subsubsection{The recovery of the nearby cycle functor
and the vanishing cycle functor}
\label{subsection;25.3.13.12}

We continue to use the notation in \S\ref{subsection;25.3.13.1}.
We consider the maps
\[
p_1,q_1:H^0(\real,L_0(\nbign))\to
H^0(\real,L_0(\nbigm)).
\]
Here, $p_1$ is the composition of
\begin{multline}
 H^0(\real,L_0(\nbign))=
 H^0(\real,L_0(\Fourier_-(\nbigm)))
 \lrarr
 H^0(\real,L_0(\Fourier_-(\nbigm(\ast 0))))
 \\
 \stackrel{(c^{-1}\circ \Abb^{\mg}_-)^{-1}}{\lrarr}
 H^0(\real,L_0(\nbigm)),
\end{multline}
and $q_1$ is the composition of
\begin{multline}
 H^0(\real,L_0(\nbign))
 \stackrel{c^{-1}\circ \Abb^{\rd}_+}{\simeq}
 H^0(\real,L_0(\Fourier_+(\nbign(!0))))
\\
 \lrarr
 H^0(\real,L_0(\Fourier_+(\nbign)))
 =H^0(\real,L_0(\nbigm)).
\end{multline}
\begin{prop}
$p_1=q_1\circ (2\pi\sqrt{-1})M_{L_0(\nbign)}$.
\end{prop}
\pf
We have
$q_1=p_1\circ (c^{-1}\circ\Abb^{\mg}_-)\circ (c^{-1}\circ \Abb^{\rd}_{+})$.
By Proposition \ref{prop;25.3.13.2},
we obtain
$q_1=p_1\circ (2\pi\sqrt{-1})^{-1}M^{-1}_{L_0(\nbign)}$.
\hfill\qed

\vspace{.1in}
We consider the maps
\[
 p_2,q_2:
 H^0(\real,L_0(\nbigm))
 \lrarr
 H^0(\real,L_0(\nbign)).
\]
Here, $p_2$ is the composition of
\begin{multline}
 H^0(\real,L_0(\nbigm))
 \stackrel{c^{-1}\circ\Abb^{\rd}_-}{\simeq}
 H^0(\real,L_0(\Fourier_-(\nbigm(!0))))
 \lrarr
 \\
 H^0\Bigl(\real,
 L_0\bigl(
 \Fourier_-(\Fourier_+(\nbign(\ast 0))(!0))
 \bigr)
 \Bigr)
=H^0\Bigl(\real,
 L_0\bigl(
 \Fourier_-(\Fourier_+(\nbign(\ast 0)))
 \bigr)
 \Bigr)
\\
=H^0\Bigl(\real,
 L_0\bigl(
 \nbign(\ast 0)
 \bigr)
 \Bigr)
=H^0\bigl(\real,
 L_0(\nbign)
 \bigr),
\end{multline}
and $q_2$ is the composition of the following maps:
\begin{multline}
 H^0(\real,L_0(\nbigm))
 \lrarr
 H^0\Bigl(\real,
 L_0\bigl(
 \Fourier_+(\nbign(\ast 0))\bigr)
 \Bigr)
 \stackrel{(c^{-1}\circ \Abb^{\mg}_{+})^{-1}}{\simeq}
 H^0(\real,L_0(\nbign)).
\end{multline}
\begin{prop}
$p_2=(-2\pi\sqrt{-1})^{-1}q_2$.
\end{prop}
\pf
We have
$p_2=(c^{-1}\circ \Abb^{\rd}_-)\circ
(c^{-1}\Abb^{\mg}_+)\circ q_2$.
By Proposition \ref{prop;25.3.13.2},
we obtain
$p_2=(-2\pi\sqrt{-1})^{-1}q_2$.
\hfill\qed

\vspace{.1in}
We set
$\phitilde'(\nbigm)=H^0(\real,L_0(\nbign))\simeq
H^0(\real,L_0(\nbign(\star 0)))$ $(\star=!,\ast)$.
There exists the following commutative diagram
in Proposition \ref{prop;25.3.11.200}:
{\small
\begin{equation}
\label{eq;25.3.13.10}
 \begin{CD}
  \psitilde(\nbigm)
  @>{\cantilde\circ\Phi_{!,-}}>>
  \phitilde(\nbigm)
  @>{(\Phi_{\ast,-})^{-1}\circ\vartilde}>>
  \psitilde(\nbigm)
  \\
  @V{\simeq}VV @V{\simeq}VV @V{\simeq}VV \\
  H^0(\real,L_0(\nbigm))
  @>{p_2}>>
  \phitilde'(\nbigm)
  @>{p_1}>>
  H^0(\real,L_0(\nbigm)).
 \end{CD}
\end{equation}}
There also exists the following commutative diagram:
{\small
\begin{equation}
\label{eq;25.3.13.11}
\begin{CD}
H^0(\real,L_0(\Fourier_+(\nbign(!0))))
@>>> 
H^0(\real,L_0(\nbigm))
@>>> 
H^0(\real,L_0(\Fourier_+(\nbign(\ast 0))))
\\
@V{\simeq}V{(c^{-1}\circ \Abb_+^{\rd})^{-1}}V
@V{=}VV
@V{\simeq}V{(c^{-1}\circ \Abb_+^{\mg})^{-1}}V
\\
\phitilde'(\nbigm)
@>{q_1}>>
H^0(\real,L_0(\nbigm))
@>{q_2}>>
\phitilde'(\nbigm).
\end{CD}
\end{equation}}
Hence, we can recover
$\psitilde(\nbigm)\to\phitilde(\nbigm)\to\psitilde(\nbigm)$
in (\ref{eq;25.3.13.10}) from (\ref{eq;25.3.13.11}).
Namely,
by setting
$M_{\phitilde'(\nbigm)}:=M_{L_0(\nbign)}$,
we define
\begin{equation}
\label{eq;25.3.18.10}
 \Phi'_{!,-}
 =(2\pi\sqrt{-1})^{-1}\Phi_{!,-},
 \quad
 \Phi'_{\ast,-}
 =(2\pi\sqrt{-1})^{-1}
 M_{\phitilde'(\nbigm)}^{-1}
 \cdot
 \Phi_{\ast,-}.
\end{equation}

\begin{prop}
\label{prop;25.3.18.20}
We obtain the following commutative diagram:
\[
 \begin{CD}
  \psitilde(\nbigm)
  @>{\cantilde\circ \Phi'_{!,-}}>>
  \phitilde(\nbigm)
  @>{(\Phi'_{\ast,-})^{-1}\circ\vartilde}>>
  \psitilde(\nbigm)\\
  @VVV @VVV @VVV \\
  H^0(\real,L_0(\nbigm))
  @>{q_2}>>
  \phitilde'(\nbigm)
  @>{q_1}>>
  H^0(\real,L_0(\nbigm)).
 \end{CD}
\]
 \end{prop}

\section{Non-unipotent monodromic regular meromorphic flat bundles}
\label{subsection;25.3.12.23}

Let $\alpha\in\cnum\setminus\seisuu$.
Let $A$ be a finite dimensional complex vector space
equipped with an endomorphism $F$
which has a unique eigenvalue $\alpha$.
Let $\nbigv=A\otimes\nbigo_{\proj^1}(\ast\{0,\infty\})$
with the connection
$\nabla=d-F\frac{dz}{z}$.

\subsection{Some notation}
We recall some notation in \S\ref{subsection;25.3.12.42}.

\subsubsection{The generalized eigen decompositions}
\label{subsection;25.3.11.1}

For any $\beta=\alpha+n$ with $n\in\seisuu$,
we obtain the subspace
$H^0(\proj^1,\nbigv)_{\beta}
 =z^{-n}A
 \subset H^0(\proj^1,\nbigv)$.
We obtain the generalized eigen decomposition
\[
 H^0(\proj^1,\nbigv)
 =\bigoplus_{\beta\in\alpha+\seisuu}
 H^0(\proj^1,\nbigv)_{\beta}
\]
of the endomorphism $-z\del_z$,
i.e.,
$H^0(\proj^1,\nbigv)_{\beta}$
is the kernel of
$(-z\del_z-\beta)^m$ for a sufficiently large $m$.
Let $F_{\beta}$ denote the endomorphism of
$H^0(\proj^1,\nbigv)_{\beta}$
induced by $-z\del_z$.
Under the isomorphism
$H^0(\proj^1,\nbigv)_{\alpha+n}\simeq z^{-n}A$,
we have 
$F_{\alpha+n}=F+n\id_A$.

\subsubsection{The associated local systems}

We obtain the $2\pi\seisuu$-equivariant local systems
$L_0(\nbigm)$ and $L_{\infty}(\nbigm)$ on $\real$
as in \S\ref{subsection;25.3.12.40}.
It is well known and easy to check that
the monodromy automorphism $M_{L_0(\nbigv)}$ has a unique eigenvalue
$\exp(2\pi\sqrt{-1}\alpha)$.

For any $v\in H^0(\proj^1,\nbigv)_{\beta}$,
we obtain 
\[
 \rho_{z,\beta}(v)
=\exp\bigl(
 F_{\beta}\log z
 \bigr)(v)
\in H^0(\real,L_0(\nbigv)).
\]
It induces an isomorphism
\[
 H^0(\proj^1,\nbigv)_{\beta}
 \simeq
 H^0(\real,L_0(\nbigv)).
\]
Under the isomorphism,
we have
$\exp(2\pi\sqrt{-1}F_{\beta})=M_{L_0(\nbigv)}$.
It is easy to check
$\rho_{z,\beta}(v)
=\rho_{z,\beta+n}(z^{-n}v)$
for any integer $n$.

\subsubsection{Fourier transforms}

We consider
the Fourier transforms
$\Fourier_{\pm}(\nbigv)$
of $\nbigv$,
which are regular singular meromorphic flat bundles
on $(\proj^1,\{0,\infty\})$.
There exist the natural isomorphisms
\[
 s_{\nbigv,\pm}:
 H^0(\proj^1,\nbigv)
\simeq
 H^0(\proj^1,\Fourier_{\pm}\nbigv).
\]
For $v\in H^0(\proj^1,\nbigv)$,
we have
$-w\del_w(s_{\nbigv,\pm}(v))
=s_{\nbigv,\pm}(\del_z zv)
=s_{\nbigv,\pm}\bigl(
 z\del_zv+v
\bigr)$.
Hence,
we obtain the following isomorphisms
for any $\beta\in\alpha+\seisuu$:
\[
 H^0(\proj^1,\nbigv)_{\beta}
 \simeq
 H^0(\proj^1,\Fourier_{\pm}\nbigv)_{-\beta+1}.
\]
There exist the natural isomorphisms
\[
 \Fourier_{\pm}(\nbigv)
 \simeq
 H^0(\proj^1,\nbigv)_{\alpha+1}
 \otimes
 \nbigo_{\proj^1}(\ast\{0,\infty\})
\]
under which
the connection of
$\Fourier_{\pm}(\nbigv)$
is identified with
$d+F\,dw/w$.

\begin{lem}
Let $F^{\gbigf}_{-\beta}$
denote the endomorphism
on $H^0(\proj^1,\Fourier_{\pm}\nbigv)_{-\beta}$
induced by $-w\del_w$.
Under the identification
\[
 H^0(\proj^1,\Fourier_{\pm}\nbigv)_{-\beta+1}
 \simeq
 H^0(\proj^1,\nbigv)_{\beta},
\] 
we have
$F^{\gbigf}_{-\beta+1}=\id-F_{\beta}$. 
\hfill\qed
\end{lem}

\subsubsection{The induced isomorphisms}
As a special case of the isomorphism in \S\ref{subsection;25.3.12.30},
we obtain the isomorphisms
\[
\Psi_{\nbigv,\beta,\pm}:H^0(\proj^1,\nbigv)_{\beta}
\simeq
H^0(\real,L_{\infty}(\Fourier_{\pm}\nbigv))
\]
as the composition of the following isomorphisms:
\begin{multline}
 H^0(\proj^1_z,\nbigv)_{\beta}
 \stackrel{s_{\nbigv,\pm}}\simeq
 H^0(\proj^1_w,\Fourier_{\pm}\nbigv)_{-\beta+1}
 \stackrel{\rho_{w,-\beta+1}}{\simeq}
 H^0(\real,L_0(\Fourier_{\pm}\nbigv))
\\
 \stackrel{c^{-1}}{\simeq}
 H^0(\real,L_{\infty}(\Fourier_{\pm}\nbigv)).
\end{multline}

\subsection{Formulas}

Let $\Gamma_{\star,\pm}$ be paths as in
\S\ref{subsection;25.3.12.41}.
\index{endomorphisms $\Ftilde_{\beta,\star,\pm}$}
We define the endomorphisms
$\Ftilde_{\beta,!,\pm}$
on $H^0(\proj^1,\nbigv)_{\beta}$
by
\index{endomorphisms $\Ftilde_{\beta,\star,\pm}$}
\begin{equation}
\label{eq;25.3.10.20}
 \Ftilde_{\beta,!,\pm}
 =\int_{\Gamma_{!,\pm}}
 \exp\bigl(F_{\beta}\log \zeta\bigr)
 e^{\mp\zeta}\frac{d\zeta}{\zeta}.
\end{equation}
We also define
the endomorphisms
$\Ftilde_{\beta,\ast,\pm}$
on $H^0(\proj^1,\nbigv)_{\beta}$
by
\begin{equation}
 \label{eq;25.3.10.21}
 \Ftilde_{\beta,\ast,\pm}
=(\pm 1)^n\int_{\Gamma_{\ast,\pm}}
 \exp\bigl((F_{\beta}+n\id)\log \zeta\bigr)
\prod_{j=0}^{n-1}
(F_{\beta}+j\id)^{-1}
 e^{\mp\zeta}\frac{d\zeta}{\zeta}.
\end{equation}
Here, $n$ denotes any non-negative integer such that
$\Re(\beta)+n>0$.
The endomorphisms $\Ftilde_{\beta,\ast,\pm}$
are independent of the choice of $n$.
If $\Re(\beta)>0$,
then
\[
 \Ftilde_{\beta,\ast,\pm}
=\int_{\Gamma_{\ast,\pm}}
 \exp\bigl(F_{\beta}\log \zeta\bigr)
 e^{\mp\zeta}\frac{d\zeta}{\zeta}.
\]

We have the isomorphism
$c^{-1}:H^0(\real,L_{\infty}(\Fourier_{\pm}\nbigv))
\simeq H^0(\real,L_{0}(\Fourier_{\pm}\nbigv))$.
We shall prove the following proposition
in \S\ref{subsection;25.3.10.1}.
\begin{prop}
\label{prop;25.3.10.2}
For $\beta\in\alpha+\seisuu$,
$\Ftilde_{\beta,\star,\pm}$
are isomorphisms.
Moreover, 
we have the following equalities
for maps
$H^0(\proj^1,\nbigv)_{\beta}
\to H^0(\real,L_{0}(\Fourier_{\pm}(\nbigv)))$:
\begin{equation}
\label{eq;25.3.10.12}
c^{-1}\circ\Abb^{\rd}_{\nbigv,\pm}\circ \rho_{z,\beta}
=\mp(2\pi\sqrt{-1})^{-1}
 \rho_{w,-\beta+1}\circ
 s_{\nbigv,\pm}\circ
 \Ftilde_{\beta,!,\pm},
\end{equation}
\begin{equation}
\label{eq;25.3.10.13}
c^{-1}\circ\Abb^{\mg}_{\nbigv,\pm}\circ
 \rho_{z,\beta}
=\mp(2\pi\sqrt{-1})^{-1}
 \rho_{w,-\beta+1}\circ
 s_{\nbigv,\pm}
 \circ
 \Ftilde_{\beta,\ast,\pm}.
\end{equation}
In other words, the following equalities hold:
\begin{equation}
\Abb^{\rd}_{\nbigv,\pm}\circ \rho_{z,\beta}
=\mp(2\pi\sqrt{-1})^{-1}
 \Psi_{\nbigv,\beta,\pm}\circ
 \Ftilde_{\beta,!,\pm},
\end{equation}
\begin{equation}
\Abb^{\mg}_{\nbigv,\pm}\circ
 \rho_{z,\beta}
=\mp(2\pi\sqrt{-1})^{-1}
\Psi_{\nbigv,\beta,\pm}
 \circ
 \Ftilde_{\beta,\ast,\pm}.
\end{equation}
\end{prop}

\begin{rem}
We have
$\Ftilde_{\beta,!,\pm}
=(1-e^{-2\pi\sqrt{-1}F_{\beta}})
\Ftilde_{\beta,\ast,\pm}$.
It is consistent with 
the relation between
$\Abb^{\rd}_{\pm}$ and $\Abb^{\rd}_{\pm}$
in Lemma {\rm\ref{lem;25.3.12.50}}.
\hfill\qed
\end{rem}

\subsubsection{Reformulation}
\label{subsection;25.3.11.130}

\index{automorphisms $\Phi_{\beta,\star,\pm}$}
We define the automorphism
$\Phi_{\beta,!,\pm}$ 
on $H^0(\proj^1,\nbigv)_{\beta}$ by
\begin{equation}
\Phi_{\beta,!,\pm}
=\frac{-1}{2\pi\sqrt{-1}}
\int_{\Gamma_{!,\pm}}
\exp(F_{\beta}\log\zeta)e^{\mp\zeta}\frac{d\zeta}{\zeta}
=\frac{-1}{2\pi\sqrt{-1}}
 \Ftilde_{\beta,!,\pm}.
\end{equation}
We also define the automorphism
$\Phi_{\beta,\ast,\pm}$ 
on $H^0(\proj^1,\nbigv)_{\beta}$ by
\begin{equation}
 \Phi_{\beta,\ast,\pm}
=\frac{-(\pm 1)^{n-1}}{2\pi\sqrt{-1}}\int_{\Gamma_{\ast,\pm}}
 \exp\bigl((F_{\beta}+n\id)\log \zeta\bigr)
\prod_{j=1}^{n}
(F_{\beta}+j\id)^{-1}
 e^{\mp\zeta}d\zeta,
\end{equation}
where $n$ is chosen as $\Re(\beta)+n>-1$.
If $\Re(\beta)>-1$, we have
\begin{equation}
 \Phi_{\beta,\ast,\pm}
  =\frac{\mp 1}  {2\pi\sqrt{-1}}
  \int_{\Gamma_{\ast,\pm}}
  \exp(F_{\beta}\log\zeta)
  e^{\mp \zeta}d\zeta.
\end{equation}
\begin{prop}
\label{prop;25.3.11.201}
We set
$\gbigl^{\gbigf}_{\pm}(\nbigv)
:=L_{\infty}(\Fourier_{\pm}(\nbigv))$
to simplify the notation.
The following diagram is commutative:
\begin{equation}
\label{eq;25.3.11.112}
\begin{CD}
H^0(\proj^1,\nbigv)_{\beta}
 @>{(-\del)\circ\Phi_{\beta,!,\pm}}>>
H^0(\proj^1,\nbigv)_{\beta+1}
 @>{\Phi_{\beta,\ast,\pm}^{-1}\circ z}>>
H^0(\proj^1,\nbigv)_{\beta}\\
 @V{\simeq}V{\Abb^{\rd}_{\nbigv,\pm}\circ\rho_{z,\beta}}V
 @V{\simeq}V{\Psi_{\nbigv,\beta+1,\pm}}V
 @V{\simeq}V{\Abb^{\mg}_{\nbigv,\pm}\circ\rho_{z,\beta}}V \\
 H^0(\real,\gbigl^{\gbigf}_{\pm}(\nbigv))
 @>{=}>>
 H^0(\real,\gbigl^{\gbigf}_{\pm}(\nbigv))
 @>{=}>>
 H^0(\real,\gbigl^{\gbigf}_{\pm}(\nbigv)).
\end{CD}
\end{equation}
The monodromy automorphisms of
$H^0(\real,\gbigl^{\gbigf}_{\pm}(\nbigv))$
are equal to
$\exp(2\pi\sqrt{-1}F_{\beta})$
on $H^0(\proj^1,\nbigv)_{\beta}$,
and
$\exp(2\pi\sqrt{-1}F_{\beta+1})$
on $H^0(\proj^1,\nbigv)_{\beta+1}$.
\end{prop}
\pf
Because
\begin{multline}
 \frac{\mp 1}{2\pi\sqrt{-1}}
 \Psi_{\nbigv,\beta+1,\pm}
 \circ\Ftilde_{\beta+1,!,\pm}(-\del v)
=\Abb^{\rd}_{\nbigv,\pm}\circ
 \rho_{z,\beta+1}(-\del v)
\\
 =\Abb^{\rd}_{\nbigv,\pm}\circ
 \rho_{z,\beta}(-z\del v)
 =\Abb^{\rd}_{\nbigv,\pm}\circ
 \rho_{z,\beta}(F_{\beta} v),
\end{multline}
we obtain
\begin{multline}
 \Abb^{\rd}_{\pm}\circ\rho_{z,\beta}(v)
= \frac{\mp 1}{2\pi\sqrt{-1}}
 \Psi_{\nbigv,\beta+1,\pm}
 \circ\Ftilde_{\beta+1,!,\pm}(-\del F_{\beta}^{-1}v)
\\
 = \frac{\mp 1}{2\pi\sqrt{-1}}
  \Psi_{\nbigv,\beta+1,\pm}
  \circ\Ftilde_{\beta+1,!,\pm}
  \circ(F_{\beta+1}-\id)^{-1}(-\del v).
\end{multline}
Note that
\[
 \int_{\Gamma_{!,\pm}}
 \exp(F_{\beta+1}\log\zeta)
 e^{\mp \zeta}
 (F_{\beta+1}-\id)^{-1}\frac{d\zeta}{\zeta}
 =\pm
 \int_{\Gamma_{!,\pm}}
 \exp((F_{\beta+1}-\id)\log\zeta)e^{\mp \zeta}
 \frac{d\zeta}{\zeta}.
\]
Hence, we obtain
\[
 \Ftilde_{\beta+1,!,\pm}\circ(F_{\beta+1}-\id)^{-1}(-\del v)
=-\del
 \Bigl(
 \pm\Ftilde_{\beta,!,\pm}(v)
 \Bigr).
\]
We obtain
\[
  \Abb^{\rd}_{\pm}\circ\rho_{z,\beta}(v)
  =\frac{-1} {2\pi\sqrt{-1}}
  \Psi_{\nbigv,\beta+1,\pm}\circ(-\del)\circ
  \Ftilde_{\beta,!,\pm}(v).
\]
This is the commutativity of the left square.

For the commutativity of the right square,
it is enough to study the case $\Re\beta>-1$.
Because
\[
 \Abb^{\mg}_{\nbigv,\pm}\circ\rho_{z,\beta}
 (zv)
 =\mp
 (2\pi\sqrt{-1})^{-1}
 \Psi_{\nbigv,\beta+1,\pm}
 \circ
 \Ftilde_{\beta+1,\ast,\pm}(v),
\]
we obtain
\[
 \Psi_{\nbigv,\beta+1,\pm}(v)
 =\mp (2\pi\sqrt{-1})
 \Abb^{\mg} _{\nbigv,\pm}\circ
 \rho_{z,\beta}
 \Bigl(
 z\cdot (\Ftilde_{\beta+1,\ast,\pm})^{-1}(v)
 \Bigr).
\]
Note that
\begin{multline}
 \int_{\Gamma_{\ast,\pm}}
 \exp(F_{\beta}\log\zeta)e^{\mp\zeta}
  \Bigl(
 z\cdot (\Ftilde_{\beta+1,\ast,\pm})^{-1}(v)
 \Bigr)\,d\zeta
= \\
 \int_{\Gamma_{\ast,\pm}}
 \exp((F_{\beta}+\id)\log\zeta)e^{\mp\zeta}
  \Bigl(
 z\cdot (\Ftilde_{\beta+1,\ast,\pm})^{-1}(v)
 \Bigr)\,\frac{d\zeta}{\zeta}
 =
\\
 z
 \Ftilde_{\beta+1,\ast,\pm}
 (\Ftilde_{\beta+1,\ast,\pm})^{-1}(v)
=zv. 
\end{multline}
Hence, we obtain
$ \Psi_{\nbigv,\beta+1,\pm}(v)
=
\Abb^{\mg}_{\nbigv,\pm}\circ\rho_{z,\beta}
\bigl(
(\Phi_{\beta,\ast,\pm})^{-1}(zv)
\bigr)$.
This is the commutativity of the right square.
\hfill\qed

\subsection{The inversion}

There exist the natural isomorphisms
\[
 \Fourier_{\pm}\circ\Fourier_{\mp}(\nbigv)\simeq \nbigv.
\]
We have
$s_{\Fourier_{\pm}(\nbigv),\mp}\circ s_{\nbigv,\pm}
=\id$.

\begin{prop}
\label{prop;25.3.12.21}
On $H^0(\real,L_0(\nbigv))$,
we have
\begin{equation}
(c^{-1}\circ \Abb^{\mg}_{\Fourier_+\nbigv,-})
\circ  
(c^{-1}\circ \Abb^{\rd}_{\nbigv,+})
=(2\pi\sqrt{-1})^{-1}
M_{L_0}^{-1},
\end{equation}
\begin{equation}
 (c^{-1}\circ \Abb^{\rd}_{\Fourier_+\nbigv,-})
\circ  
(c^{-1}\circ \Abb^{\mg}_{\nbigv,+})
=-(2\pi\sqrt{-1})^{-1}\id,
\end{equation}
\begin{equation}
(c^{-1}\circ \Abb^{\mg}_{\Fourier_-\nbigv,+})
\circ
(c^{-1}\circ \Abb^{\rd}_{\nbigv,-})
=-(2\pi\sqrt{-1})^{-1}\id,
\end{equation}
\begin{equation}
(c^{-1}\circ \Abb^{\rd}_{\Fourier_-\nbigv,+})
\circ
(c^{-1}\circ \Abb^{\mg}_{\nbigv,-})
=(2\pi\sqrt{-1})^{-1}M_{L_0}.
\end{equation}
\end{prop}
\pf
We have only to prove the following equalities
on $H^0(\proj^1,\nbigv)_{\beta}$:
\begin{equation}
(c^{-1}\circ \Abb^{\mg}_{\Fourier_+\nbigv,-})
\circ  
(c^{-1}\circ \Abb^{\rd}_{\nbigv,+})\circ \rho_{z,\beta}
=(2\pi\sqrt{-1})^{-1}
\rho_{z,\beta}\circ
\exp\bigl(-2\pi\sqrt{-1}F_{\beta}\bigr),
\end{equation}
\begin{equation}
 (c^{-1}\circ \Abb^{\rd}_{\Fourier_+\nbigv,-})
\circ  
(c^{-1}\circ \Abb^{\mg}_{\nbigv,+})\circ \rho_{z,\beta}
=-(2\pi\sqrt{-1})^{-1}
\rho_{z,\beta},
\end{equation}
\begin{equation}
(c^{-1}\circ \Abb^{\mg}_{\Fourier_-\nbigv,+})
\circ
(c^{-1}\circ \Abb^{\rd}_{\nbigv,-})\circ\rho_{z,\beta}
=-(2\pi\sqrt{-1})^{-1}\rho_{z,\beta},
\end{equation}
\begin{equation}
(c^{-1}\circ \Abb^{\rd}_{\Fourier_-\nbigv,+})
\circ
(c^{-1}\circ \Abb^{\mg}_{\nbigv,-})\circ\rho_{z,\beta}
=(2\pi\sqrt{-1})^{-1}\rho_{z,\beta}\circ
\exp(2\pi\sqrt{-1}F_{\beta}).
\end{equation}
We have
\begin{multline}
 (c^{-1}\circ \Abb^{\mg}_{\Fourier_+\nbigv,-})
 \circ
 (c^{-1}\circ\Abb^{\rd}_{\nbigv,+})\circ\rho_{z,\beta}(v)
 =\\
 -(2\pi\sqrt{-1})^{-2}
 \rho_{z,\beta}\circ
 s_{\Fourier_+\nbigv,-}
 \Ftilde^{\gbigf}_{-\beta+1,\ast,-}
 s_{\nbigv,+}(\Ftilde_{\beta,!,+}(v)).
\end{multline}
We also have
\begin{multline}
 s_{\Fourier_+\nbigv,-}
 \Ftilde^{\gbigf}_{-\beta+1,\ast,-}
 s_{\nbigv,+}(\Ftilde_{\beta,!,+}(v))
 =
\\
 \int_{\Gamma_{\ast,-}}
 \exp\bigl(
 (\id-F_{\beta})\log\eta
 \bigr)
 e^{\eta}\frac{d\eta}{\eta}
 \cdot
 \int_{\Gamma_{!,+}}
 \exp\bigl(
 F_{\beta}\log\zeta
 \bigr)(v)\cdot
 e^{-\zeta}\frac{d\zeta}{\zeta}
\\
 =-2\pi\sqrt{-1}
 \exp\bigl(
  -2\pi\sqrt{-1}F_{\beta}
 \bigr).
\end{multline}
We obtain other formulas similarly.
\hfill\qed

\subsubsection{Appendix}

Let $f$ be any endomorphism of $A$.
We set
\[
 \ftilde_{!,\pm}
 =\int_{\Gamma_{!,\pm}}
 \exp(f\log\zeta)
 e^{\mp \zeta}\frac{d\zeta}{\zeta}.
\]
If any eigenvalue of $f$ is not a non-positive integer,
we also set
\[
 \ftilde_{\ast,\pm}
 =\int_{\Gamma_{\ast,\pm}}
 \exp((f+n\id)\log\zeta)
 \prod_{j=0}^{n-1}(f+j\id)^{-1}
 e^{\mp\zeta}
 \frac{d\zeta}{\zeta},
\]
where $n$ denotes any non negative integer
such that
$\Re\alpha+n>0$ holds
for any eigenvalue $\alpha$ of $f$.
When $n=0$,
$\prod_{j=0}^{n-1}(f+j\id)^{-1}$ means the identity.

\begin{lem}
\label{lem;25.3.11.40}
We set $f^{\gbigf}=\id-f$.
Suppose that 
any eigenvalue of $f$ is not a non-positive integer,
i.e., $\ftilde_{\ast,\pm}$ are defined.
Then, we obtain the following equalities:
\begin{equation}
\label{eq;25.3.11.20}
 \ftilde^{\gbigf}_{!,-}\circ
 \ftilde_{\ast,+}
=2\pi\sqrt{-1}\id,
\end{equation}
\begin{equation}
\label{eq;25.3.11.21}
 \ftilde^{\gbigf}_{!,+}\circ
 \ftilde_{\ast,-}
=-(2\pi\sqrt{-1})\exp(2\pi\sqrt{-1}f).
\end{equation}
\end{lem}
\pf
Because the both sides of the equalities are
complex analytic with respect to $f$,
it is enough to prove the equalities
for $f$ satisfying the following conditions.
\begin{itemize}
 \item Any eigenvalue of $f$ satisfies $0<\Re\alpha<1$.
 \item $f$ is semisimple.
       It implies that it is enough to consider the case
       $\dim A=1$.
\end{itemize}
We obtain the equalities
from the standard reflection formula
for the Gamma functions
\[
\Gamma(1-z)\Gamma(z)=\frac{\pi}{\sin \pi z}
=\frac{2\pi\sqrt{-1}}{e^{\pi z}-e^{-\pi z}}.
\]
Indeed,
we obtain (\ref{eq;25.3.11.20})
from the following.
\begin{multline}
 \int_{\Gamma_{!,-}}
 \exp((1-\alpha)\log\eta)e^{\eta}\frac{d\eta}{\eta}
 \int_{\Gamma_{\ast,+}}
 \exp(\alpha \log\zeta)e^{-\zeta}\frac{d\zeta}{\zeta}
 =\\
 (e^{(1-\alpha)\pi\sqrt{-1}}-e^{-(1-\alpha)\pi\sqrt{-1}})
 \Gamma(1-\alpha)\Gamma(\alpha)
=2\pi\sqrt{-1}.
\end{multline}
We obtain (\ref{eq;25.3.11.21})
from the following.
\begin{multline}
 \int_{\Gamma_{!,+}}
 \exp((1-\alpha)\log\eta)e^{\eta}\frac{d\eta}{\eta}
 \int_{\Gamma_{\ast,-}}
 \exp(\alpha \log\zeta)e^{-\zeta}\frac{d\zeta}{\zeta}
 =\\
 (1-e^{-2(1-\alpha)\pi\sqrt{-1}})
 e^{\alpha\pi\sqrt{-1}}
 \Gamma(1-\alpha)\Gamma(\alpha)
=-2\pi\sqrt{-1}e^{2\pi\sqrt{-1}\alpha}.
\end{multline}
\hfill\qed

\begin{cor}
\label{cor;25.3.11.41}
Suppose that
any eigenvalue of $f$ is not a positive integer,
i.e., 
$\ftilde^{\gbigf}_{\ast,\pm}$ are defined.
Then, we obtain the following equalities:
\begin{equation}
\label{eq;25.3.10.30}
 \ftilde^{\gbigf}_{\ast,-}\circ
 \ftilde_{!,+}
=-2\pi\sqrt{-1}\exp(-2\pi\sqrt{-1}f),
\end{equation}
\begin{equation}
 \ftilde^{\gbigf}_{\ast,+}\circ
 \ftilde_{!,-}
=2\pi\sqrt{-1}\id,
 \end{equation}
\end{cor}
\pf
Note that $(f^{\gbigf})^{\gbigf}=f$.
We also note that
$f$, $f^{\gbigf}$,
$\ftilde_{\star,\pm}$
and
$\ftilde^{\gbigf}_{\star,\pm}$
are mutually commuting.
Hence, we obtain the claims from Lemma \ref{lem;25.3.11.40}.
\hfill\qed

\subsection{Proof of Proposition \ref{prop;25.3.10.2}}
\label{subsection;25.3.10.1}

We explain the case $\beta=\alpha$.
We obtain the claim in the other cases
by replacing $(A,F)$ with $(z^{-n}A,F+n\id)$.

\subsubsection{Complexes}
\label{subsection;25.3.12.1}

For $m\in\seisuu$,
we set
\[
 \nbigc^{0}_{\pm,m}(\nbigv)
=A\otimes\nbigo_{\proj^1}\bigl(
 (m-1)\{0\}+(m-2)\{\infty\}
 \bigr),
\quad
 \nbigc^1_{\pm,m}(\nbigv)
=A\otimes \Omega^1_{\proj^1}(m\{0\}+m\{\infty\}).
\]
Let $\pi_z:\proj^1_z\times\cnum_w\to \proj^1_z$
and $\pi_w:\proj^1_z\times\cnum_w\to \cnum_w$
denote the projections.
We obtain the following complexes
$\nbigctilde_{\pm,m}^{\bullet}(\nbigv)$
on $\proj^1_z\times \cnum_w$:
\[
 \pi_z^{\ast}\nbigc^0_{\pm,m}(\nbigv)
 \stackrel{d\pm w\,dz}{\lrarr}
 \pi_z^{\ast}\nbigc^1_{\pm,m}(\nbigv).
\]
There exist the natural isomorphisms for any $m$:
\begin{equation}
 \Fourier_{\pm}(\nbigv)(\ast 0)
 \simeq
 R^1\pi_{w\ast}
 \bigl(
 \nbigctilde^{\bullet}_{\pm,m}(\nbigv)
 \bigr)(\ast 0).
\end{equation}

\subsubsection{Representatives}
\label{subsection;25.3.12.2}

For $w\in\cnum^{\ast}$,
we obtain the complex
$\nbigc^{\bullet}_{\pm,m}(\nbigv)_w$:
\[
 \nbigc^{0}_{\pm,m}(\nbigv)
 \stackrel{d\pm w\,dz}{\lrarr}
 \nbigc^{1}_{\pm,m}(\nbigv).
\]
There exist the isomorphisms
\begin{equation}
\label{eq;25.3.10.3}
 \Fourier_{\pm}(\nbigv)(\ast 0)_{|w}
 \simeq
 H^1\bigl(
 \proj^1,
 \nbigc^{\bullet}_{\pm,m}(\nbigv)_w
\bigr).
\end{equation}

\begin{lem}
For any $v\in A= H^0(\proj^1,\nbigv)_{\alpha}$,
$s_{\nbigv,\pm}(z^{-1}v)_{|w}$
are represented by
$v\otimes dz/z$ of
$\nbigc^1_{\pm,1}(\nbigv)$
under the isomorphisms {\rm(\ref{eq;25.3.10.3})}
with $m=1$.
\hfill\qed
\end{lem}

Let $\nbigc^{\bullet}_{\pm,0,C^{\infty}}(\nbigv)_w$
denote the Dolbeault resolution of
$\nbigc^{\bullet}_{\pm,0}(\nbigv)_w$.
Let $\chi:\proj^1\to \{0\leq a\leq 1\}$ be
a $C^{\infty}$-function
such that
$\chi(z)=0$ $(|z|<1/2)$
and
$\chi(z)=1$ $(|z|>1)$.
For any $v\in A$,
we set
\[
 B_{\pm}(v)
=v\otimes dz
\mp(d\pm w\,dz)
\Bigl(
 w^{-1}\chi v
\pm w^{-2}\chi F(v) z^{-1}
\Bigr)
\in \nbigc^1_{\pm,0,C^{\infty}}(\nbigv)_w.
\]
\begin{lem}
$s_{\nbigv,\pm}(v)_{|w}$
are represented by
$B_{\pm}(v)$
under the isomorphisms {\rm(\ref{eq;25.3.10.3})}
with $m=0$.
\end{lem}

\subsubsection{Duality}

Let $A^{\lor}$ denote the dual space,
which is equipped with the dual endomorphism $F^{\lor}$.
Let $\langle\cdot,\cdot\rangle$ denote the natural pairing of
$A$ and $A^{\lor}$.
We obtain the dual bundle
$\nbigv^{\lor}=
A^{\lor}\otimes
\nbigo_{\proj^1}(\ast\{0,\infty\})$
with the induced connection $d+F^{\lor}\frac{dz}{z}$.

We obtain the complexes
$\nbigctilde^{\bullet}_{\pm,m}(\nbigv^{\lor})$
on $\proj^1\times\cnum_w$.
There exist the natural pairings
\[
 \langle\cdot,\cdot\rangle_{\pm}:
\nbigctilde^{\bullet}_{\pm,0}(\nbigv)
 \otimes
 \nbigctilde^{\bullet}_{\mp,1}(\nbigv^{\lor})
 \lrarr
 \pi_z^{\ast}\Omega^{\bullet}_{\proj^1}.
\]
They induce the pairings:
\[
 \langle\cdot,\cdot\rangle_{\pm}:
 R^1\pi_{w\ast}\bigl(
 \nbigctilde^{\bullet}_{\pm,0}(\nbigv)
 \bigr)(\ast 0)
 \otimes
 R^1\pi_{w\ast}
 \bigl(
 \nbigctilde^{\bullet}_{\mp,1}(\nbigv^{\lor})
 \bigr)(\ast 0)
 \lrarr
 \nbigo_{\cnum_w}(\ast 0).
\]

Let $v\in A=H^0(\proj^1,\nbigv)_{\alpha}$
and $v^{\lor}\in A^{\lor}=H^0(\proj^1,\nbigv^{\lor})_{-\alpha}$.
Note that
$s_{\nbigv,\pm}(v)\in
H^0(\proj^1,\Fourier_{\pm}(\nbigv))_{-\alpha+1}$
and 
$s_{\nbigv^{\lor},\pm}(z^{-1}v)\in
H^0(\proj^1,\Fourier_{\pm}(\nbigv))_{\alpha}$.
\begin{lem}
\begin{equation}
\label{eq;25.3.10.4}
 \langle s_{\nbigv,\pm}(v),
  s_{\mp}(z^{-1}v^{\lor})
 \rangle_{\pm}
=\mp(2\pi\sqrt{-1})w^{-1}
 \langle v,v^{\lor}\rangle.
\end{equation}
\end{lem}
\pf
We have
\begin{multline}
\langle s_{\nbigv,\pm}(v),
s_{\mp}(z^{-1}v^{\lor})
\rangle_{\pm}
=\langle
B_{\pm}(v),v^{\lor}(dz/z)
\rangle
\\
=
\mp w^{-1}\langle v,v^{\lor}\rangle
 \int_{\proj^1}
\delbar\chi
 \frac{dz}{z}
-w^{-2}
 \langle
F(v),v^{\lor}
\rangle
\int_{\proj^1}
\delbar\chi
 d(z^{-1}).
\end{multline}
We obtain (\ref{eq;25.3.10.4})
from
$\int_{\proj^1}
\delbar\chi
 \frac{dz}{z}
=2\pi\sqrt{-1}$
and
$\int_{\proj^1}
\delbar\chi
 d(z^{-1})=0$.
\hfill\qed

\subsubsection{Proof of Proposition \ref{prop;25.3.10.2}}

For $v\in A$ and $v\in A^{\lor}$,
we have
\begin{multline}
\label{eq;25.3.10.10}
\bigl\langle
c^{-1}\circ\Abb^{\rd}_{\nbigv,\pm}\circ\rho_{z,\alpha}(v),
s_{\nbigv^{\lor},\mp}(z^{-1}v^{\lor})
\bigr\rangle_{\pm}
 =
\\\int_{\Gamma_{!,\theta^u,\pm}}
 \bigl\langle
 \exp(F_{\alpha}\log z)v,
 v^{\lor}
 \bigr\rangle e^{\mp zw}\,\frac{dz}{z}
=\int_{\Gamma_{!,\pm}}
 \bigl\langle
 \exp(F_{\alpha}\log(\zeta w^{-1})v,v^{\lor})
 \bigr\rangle
 e^{\mp \zeta}\frac{d\zeta}{\zeta}.
\end{multline}
By using
$\exp(F_{\alpha}\log(\zeta w^{-1}))
=w^{-1}\exp\bigl((-F_{\alpha}+\id)\log w\bigr)\cdot
\exp(F_{\alpha}\log \zeta)$,
we rewrite (\ref{eq;25.3.10.10}) as
\begin{multline}
 w^{-1}\Bigl\langle
 \exp\bigl((-F_{\alpha}+\id)\log w\bigr)
 \Ftilde_{\alpha,!,\pm}(v),
 v^{\lor}
 \Bigr\rangle 
 =\\
 \mp
 (2\pi\sqrt{-1})^{-1}\Bigl\langle
 \rho_{w,-\alpha+1}\circ
 s_{\nbigv,\pm}
 \bigl(\Ftilde_{\alpha,!,\pm}(v)
 \bigr),
 s_{\nbigv^{\lor},\mp}(z^{-1}v^{\lor})
 \Bigr\rangle_{\pm}.
\end{multline}
Hence, we obtain (\ref{eq;25.3.10.12}).
Concerning (\ref{eq;25.3.10.13}),
it is enough to consider the case $\Re(\alpha)>0$.
Then, we can prove it by the same argument.
\hfill\qed

\section{Unipotent monodromic holonomic $\nbigd$-modules}
\label{subsection;25.3.12.24}

Let $A$ be a finite dimensional complex vector space
equipped with a nilpotent endomorphism $N$.
Let $\nbigv=A\otimes\nbigo_{\proj^1}(\ast\{0,\infty\})$
with the connection
$\nabla=d-N\frac{dz}{z}$.
Let $\nbigm$ be a regular holonomic $\nbigd_{\proj^1}$-modules
such that $\nbigm(\ast 0)=\nbigv$.

\subsection{Some notation}
We recall some notation in \S\ref{subsection;25.3.12.42}.

\subsubsection{The generalized eigen decompositions and the $V$-filtrations}

We obtain the $\cnum$-linear endomorphism
$-z\del_z$ on $H^0(\proj^1,\nbigm)$,
and we obtain the generalized eigen decomposition
\[
 H^0(\proj^1,\nbigm)
 =\bigoplus_{m\in\seisuu}
 H^0(\proj^1,\nbigm)_m,
\]
where $H^0(\proj^1,\nbigm)_m$
denotes the kernel of
$(-z\del_z-m)^n$ for any sufficiently large $n$.
Let $N_{\nbigm}$ denote the nilpotent endomorphism
of $H^0(\proj^1,\nbigm)_1$ induced by $-\del_zz$.

There exists the $V$-filtration of $\nbigm$.
It is indexed by $\seisuu$ in this case.
We have
$H^0(\proj^1,V_m(\nbigm))
=\bigoplus_{p\leq m+1}H^0(\proj^1,\nbigm)_p$.

\begin{lem}
We have $H^0(\proj^1,\nbigm)_n=H^0(\proj^1,\nbigv)_n=z^{n}A$
for any $n\leq 0$.
It is naturally isomorphic to
$\Gr^V_{n-1}(\nbigm)$.
The nilpotent part of $-z\del_z$ is identified with $N$.
 \hfill\qed
\end{lem}

There exist the isomorphisms
\[
-\del_z:H^0(\proj^1,\nbigv(!0))_{0}\simeq H^0(\proj^1,\nbigv(!0))_1,
\quad\quad
z:H^0(\proj^1,\nbigv)_{1}\simeq H^0(\proj^1,\nbigv)_{0}.
\]
By using $A=H^0(\proj^1,\nbigv[\star 0])_{0}$,
we denote elements of $H^0(\proj^1,\nbigv)_1$
by $z^{-1}v$ $(v\in A)$,
and
elements of $H^0(\proj^1,\nbigv[!0])_1$
by $-\del_z\otimes v$ $(v\in A)$.
The endomorphisms
$N_{\nbigv(\star 0)}$ $(\star=!,\ast)$
are naturally identified with $N$.

\subsubsection{The associated local systems}

We obtain
the $2\pi\seisuu$-equivariant local systems
$L_0(\nbigm)$ and $L_{\infty}(\nbigm)$ on $\real$
as in \S\ref{subsection;25.3.12.40}.
It is well known and easy to check that
the monodromy automorphism $M_{L_0(\nbigm)}$
has a unique eigenvalue $1$.

For any $v\in H^0(\proj^1,\nbigm)_{0}$,
we obtain 
\[
 \rho_{z}(v)
=\exp\bigl(
 N\log z
 \bigr)(v)
\in H^0(\real,L_0(\nbigm)).
\]
It induces an isomorphism
$H^0(\proj^1,\nbigm)_{0}
 \simeq
 H^0(\real,L_0(\nbigm))$
under which
\[
\exp(2\pi\sqrt{-1}N)=M_{L_0(\nbigm)}.
\]

\subsubsection{Fourier transforms}

We consider the Fourier transforms
$\Fourier_{\pm}(\nbigm)$ of $\nbigm$,
which are regular holonomic $\nbigd$-modules
on $\proj^1$
such that
$\Fourier_{\pm}(\nbigm)(\ast 0)$
are meromorphic flat bundles on $(\proj^1,\{0,\infty\})$.
There exist the natural isomorphisms
\[
 s_{\nbigm,\pm}:
 H^0(\proj^1,\nbigm)
\simeq
 H^0(\proj^1,\Fourier_{\pm}\nbigm).
\]
It induces the isomorphisms
$s_{\nbigm,\pm}:
 H^0(\proj^1,\nbigm)_{m}
 \simeq
 H^0(\proj^1,\Fourier_{\pm}\nbigm)_{-m+1}$.
In particular,
it induces
\[
 H^0(\proj^1,\nbigm)_{1}
 \simeq
 H^0(\proj^1,\Fourier_{\pm}\nbigm)_{0}.
\]
There exist the natural isomorphisms
\[
 \Fourier_{\pm}(\nbigm)(\ast 0)
 \simeq
 H^0(\proj^1,\nbigm)_1
 \otimes\nbigo_{\proj^1}(\ast\{0,\infty\})
\]
under which
the connections of
$\Fourier_{\pm}(\nbigm)(\ast 0)$
are identified with
$d+N_{\nbigm}\cdot dw/w$.

\subsubsection{The induced isomorphisms}
As a special case of the isomorphism in \S\ref{subsection;25.3.12.30},
we obtain the isomorphisms
$\Psi_{\nbigm,\pm}:\Gr^V_0(\nbigm)\simeq
H^0(\real,L_{\infty}(\Fourier_{\pm}(\nbigv)))$
as the composition of
the following isomorphisms:
\begin{multline}
 \Gr^{V}_0(\nbigm)\simeq
 H^0(\proj^1_z,\nbigm)_1
 \stackrel{s_{\nbigm,\pm}}{\simeq}
 H^0(\proj^1_w,\Fourier_{\pm}\nbigm)_{0}
 \stackrel{\rho_{w}}{\simeq}
 H^0(\real,L_0(\Fourier_{\pm}(\nbigm)))
\\
 \stackrel{c^{-1}}{\simeq}
 H^0(\real,L_{\infty}(\Fourier_{\pm}(\nbigm))).
\end{multline}
Let $N_{\nbigm}$ denote the nilpotent endomorphisms
$\Gr^V_0(\nbigm)$
induced by $-\del_zz$.
The monodromy automorphism
on $H^0(\real,L_{\infty}(\Fourier_{\pm}(\nbigm)))$
equals $\exp(2\pi\sqrt{-1}N_{\nbigm})$.

\subsection{Formulas}

\index{endomorphisms $\Ftilde_{\star,\pm}$}

We regard $N$ as an endomorphism of
the $\nbigd$-modules $\nbigv[\star 0]$.
We define the endomorphisms
$\Ftilde_{!,\pm}$ of $\nbigv[\star 0]$
by
\begin{equation} 
\label{eq;25.3.11.140}
 \Ftilde_{!,\pm}
=\int_{\Gamma_{!,\pm}}
\exp\bigl(N\log \zeta\bigr)
e^{\mp\zeta}\frac{d\zeta}{\zeta}.
\end{equation}
We also define the endomorphisms
$\Ftilde_{\ast,\pm}$ of $\nbigv[\star 0]$
by
\begin{equation}
\label{eq;25.3.11.141}
 \Ftilde_{\ast ,\pm}
=\pm\int_{\Gamma_{\ast,\pm}}
\exp\bigl(N\log \zeta\bigr)
e^{\mp\zeta}d\zeta.
\end{equation}
We shall prove the following proposition
in \S\ref{subsection;25.3.16.2}.
\begin{prop}
\label{prop;25.3.10.110}
$\Ftilde_{\star,\pm}$
are isomorphisms.
Moreover,
for any
$v\in A=H^0(\proj^1,\nbigv(\star 0))_{0}$,
we have
\begin{equation}
\label{eq;25.3.10.130}
c^{-1}\circ\Abb^{\rd}_{\nbigv,\pm}\circ\rho_z(v)=
\frac{-1}{2\pi\sqrt{-1}}
\rho_w\circ
 s_{\nbigv(!0),\pm}(-\del_z\otimes \Ftilde_{!,\pm}(v)),
\end{equation}
\begin{equation}
\label{eq;25.3.10.131}
 c^{-1}\circ\Abb^{\mg}_{\nbigv,\pm}\circ\rho_z(v)
 =\frac{-1}{2\pi\sqrt{-1}}
  \rho_w\circ s_{\nbigv,\pm}(z^{-1}\Ftilde_{\ast,\pm}(v)).
\end{equation}
In other words,
\begin{equation}
\Abb^{\rd}_{\nbigv,\pm}\circ\rho_z(v)=
\frac{-1}{2\pi\sqrt{-1}}
\Psi_{\nbigv(!0)}
(-\del_z\otimes \Ftilde_{!,\pm}(v)),
\end{equation}
\begin{equation}
 \Abb^{\mg}_{\nbigv,\pm}\circ\rho_z(v)
 =\frac{-1}{2\pi\sqrt{-1}}
 \Psi_{\nbigv,\pm}
 (z^{-1}\Ftilde_{\ast,\pm}(v)).
\end{equation}
\end{prop}

\begin{rem}
Let $G(t)=t^{-1}(1-e^{-2\pi\sqrt{-1}t})$.
It is easy to see 
$\Ftilde_{!,\pm}
=G(N)\Ftilde_{\ast,\pm}$.
In particular,
$\Ftilde_{!,\pm}\circ N
=(\id-e^{-2\pi\sqrt{-1}}N)\Ftilde_{\ast,\pm}$.
It is consistent with
the relation between
$\Abb^{\rd}_{\nbigv,\pm}$  and $\Abb^{\mg}_{\nbigv,\pm}$
in Lemma {\rm\ref{lem;25.3.12.50}}
\hfill\qed
\end{rem}

\subsubsection{Reformulation}
\label{subsection;25.3.11.131}

\index{automorphisms $\Phi_{0,\star,\pm}$}

We define the automorphisms
$\Phi_{0,\star,\pm}$ $(\star=!,\ast)$
on $\Gr^V_{-1}(\nbigv)=H(\proj^1,\nbigv)_0$ by
\[
 \Phi_{0,\star,\pm}
=\frac{-1}{2\pi\sqrt{-1}}
 \Ftilde_{\star,\pm}.
\]
We have the following equalities on
$\Gr^V_{-1}(\nbigv(!0))$:
\begin{equation}
\label{eq;25.3.11.100}
 \Abb^{\rd}_{\nbigv,\pm}\circ\rho_z
 =\Psi_{\nbigv(!0),\pm}\circ\can_{\nbigv(!0)}\circ\Phi_{0,!,\pm}.
\end{equation}
We also have the following equalities on
$\Gr^V_0(\nbigv)$:
\begin{equation}
\label{eq;25.3.11.101}
 \Abb^{\mg}_{\nbigv,\pm}\circ\rho_z
 \circ
 (\Phi_{0,\ast,\pm})^{-1}\circ
 \var_{\nbigv}
 =\Psi_{\nbigv,\pm}.
\end{equation}

\begin{prop}
\label{prop;25.3.11.202}
\label{prop;24.4.18.4}
We set
$\gbigl^{\gbigf}_{\pm}(\nbigm)
:=L_{\infty}(\Fourier_{\pm}(\nbigm))$
to simplify the notation.
The following diagrams are commutative:
\begin{equation}
\label{eq;25.3.11.102}
\begin{CD}
 \Gr^V_{-1}(\nbigv)
 @>{\can_{\nbigm}\circ \Phi_{0,!,\pm}}>>
 \Gr^V_0(\nbigm)
 @>{(\Phi_{0,\ast,\pm})^{-1}\circ\var_{\nbigm}}>>
 \Gr^V_{-1}(\nbigv)\\
 @V{\simeq}V{\Abb^{\rd}_{\nbigv,\pm}\circ\rho_z}V
 @V{\simeq}V{\Psi_{\nbigm,\pm}}V
 @V{\simeq}V{\Abb^{\mg}_{\nbigv,\pm}\circ\rho_z}V \\
 H^0(\real,\gbigl^{\gbigf}_{\pm}(\nbigv(!0)))
 @>>>
 H^0(\real,\gbigl^{\gbigf}_{\pm}(\nbigm))
 @>>>
 H^0(\real,\gbigl^{\gbigf}_{\pm}(\nbigv)).
\end{CD}
\end{equation}
Here, the lower horizontal arrows are the natural morphisms.
The monodromy automorphisms of
$H^0(\real,\gbigl^{\gbigf}_{\pm}(\nbigv(!0)))$,
$H^0(\real,\gbigl^{\gbigf}_{\pm}(\nbigm))$
and
$H^0(\real,\gbigl^{\gbigf}_{\pm}(\nbigv))$
are equal to
$\exp(2\pi\sqrt{-1}N)$,
$\exp(2\pi\sqrt{-1}N_{\nbigm})$
and
$\exp(2\pi\sqrt{-1}N)$
on
$\Gr^V_{-1}(\nbigv)$,
$\Gr^V_0(\nbigm)$
and
$\Gr^V_{-1}(\nbigv)$. 
\end{prop}
\pf
There exists the following commutative diagram:
\[
 \begin{CD}
  \Gr^V_0(\nbigv[!0])
  @>>>
  \Gr^V_0(\nbigm)
  @>>>
  \Gr^V_0(\nbigv)
  \\
  @V{\simeq}V{\Psi_{\nbigv[!0],\pm}}V
  @V{\simeq}V{\Psi_{\nbigm,\pm}}V
  @V{\simeq}V{\Psi_{\nbigv,\pm}}V \\
  H^0(\real,\gbigl^{\gbigf}_{\pm}(\nbigv(!0)))
  @>>>
  H^0(\real,\gbigl^{\gbigf}_{\pm}(\nbigm))
  @>>>
  H^0(\real,\gbigl^{\gbigf}_{\pm}(\nbigv)).
 \end{CD}
\]
By using (\ref{eq;25.3.11.100})
and (\ref{eq;25.3.11.101}),
we obtain the commutativity of (\ref{eq;25.3.11.102}).
\hfill\qed

\begin{cor}
\label{cor;25.3.12.20}
The composition of the natural morphisms
\begin{multline}
 H^0(\real,\gbigl^{\gbigf}_{\pm}(\nbigm))
 \to
 H^0(\real,\gbigl^{\gbigf}_{\pm}(\nbigv))
 \simeq
 \\
 \Gr^V_{-1}(\nbigv)
 \simeq
 H^0(\real,\gbigl^{\gbigf}_{\pm}(\nbigv(!0)))
 \to
 H^0(\real,\gbigl^{\gbigf}_{\pm}(\nbigm))
\end{multline}
equal $1-M_{\pm}^{-1}$,
where $M_{\pm}$ denote the monodromy automorphisms of
$\gbigl^{\gbigf}_{\pm}(\nbigm)$.
\end{cor}
\pf
From $N$ on $H^0(\proj^1,\nbigm)_n$ $(n\leq 0)$
and $N_{\nbigm}$ on $H^0(\proj^1,\nbigm)_n$ $(n>0)$,
we obtain an endomorphism $N$ on $\nbigm$.
We define the endomorphism
$\Ftilde^{\nbigm}_{\star,\pm}$ $(\star=!,\ast)$
of $\nbigm$
as in the case of $\nbigv[\star 0]$
by using the formulas (\ref{eq;25.3.11.140})
and (\ref{eq;25.3.11.141}).
We also define
$\Phi^{\nbigm}_{0,\star,\pm}
=\frac{-1}{2\pi\sqrt{-1}}
\Ftilde^{\nbigm}_{\star,\pm}$.
We have
\begin{multline}
 \can_{\nbigm}
 \circ
 \Phi_{0,!,\pm}
 \circ
 (\Phi_{0,\ast,\pm})^{-1}\circ
 \var_{\nbigm}
 =
  \can_{\nbigm}
 \circ
 \var_{\nbigm}\circ
 \Phi^{\nbigm}_{0,!,\pm}
 \circ
 (\Phi^{\nbigm}_{0,\ast,\pm})^{-1}
\\
=N_{\nbigm}\circ
  \Phi^{\nbigm}_{0,!,\pm}
 \circ
 (\Phi^{\nbigm}_{0,\ast,\pm})^{-1}.
\end{multline}
Because 
$N_{\nbigm}\circ
  \Phi^{\nbigm}_{0,!,\pm}
  =(1-e^{-2\pi\sqrt{-1}N_{\nbigm}})
  \Phi^{\nbigm}_{0,\ast,\pm}$,
we obtain the claim of the corollary.
\hfill\qed

\subsection{The inversion}

There exist the natural isomorphisms
\[
\Fourier_{\pm}\circ\Fourier_{\mp}(\nbigm)\simeq \nbigm.
\]
We have
$s_{\Fourier_{\pm}(\nbigm),\mp}\circ s_{\nbigm,\pm}
=\id$.
We set $\vecF_{\star,\pm}(\nbigv)=\Fourier(\nbigv(\star 0))(\ast 0)$
$(\star=!,\ast)$.

\begin{prop}
\label{prop;25.3.12.22}
On $H^0(\real,L_0(\nbigv[!0]))$,
we have
\begin{equation}
\label{eq;25.3.11.10}
(c^{-1}\circ\Abb^{\mg}_{\vecF_{!,-}(\nbigv),+})
\circ
(c^{-1}\circ\Abb^{\rd}_{\nbigv,-})
=
-(2\pi\sqrt{-1})^{-1}\id,
\end{equation}
\begin{equation}
\label{eq;25.3.11.11}
 (c^{-1}\circ\Abb^{\mg}_{\vecF_{!,+}(\nbigv),-})
  \circ
  (c^{-1}\circ\Abb^{\rd}_{\nbigv,+})\circ\rho_z
  =(2\pi\sqrt{-1})^{-1}
  M_{L_0(\nbigv[!0])}^{-1}.
\end{equation}
On $H^0(\real,L_0(\nbigv))$,
we have 
\begin{equation}
\label{eq;25.3.11.12}
(c^{-1}\circ\Abb^{\rd}_{\vecF_{\ast,-}(\nbigv),+})
\circ
(c^{-1}\circ\Abb^{\mg}_{\nbigv,-})
=(2\pi\sqrt{-1})^{-1}\cdot M_{L_0(\nbigv)},
\end{equation}
\begin{equation}
\label{eq;25.3.11.13}
(c^{-1}\circ\Abb^{\rd}_{\vecF_{\ast,+}(\nbigv),-})
\circ
(c^{-1}\circ\Abb^{\mg}_{\nbigv,+})
=-(2\pi\sqrt{-1})^{-1}\id.
\end{equation}
\end{prop}
\pf
The equalities 
(\ref{eq;25.3.11.10})
and 
(\ref{eq;25.3.11.11})
are the translation of
the following equalities
on $H^0(\proj^1,\nbigv[!0])_{-1}$:
\begin{equation}
\label{eq;25.3.11.50}
(c^{-1}\circ\Abb^{\mg}_{\vecF_{!,-}(\nbigv),+})
\circ
(c^{-1}\circ\Abb^{\rd}_{\nbigv,-})\circ\rho_z
=
-(2\pi\sqrt{-1})^{-1}\rho_z,
\end{equation}
\begin{equation}
\label{eq;25.3.11.51}
 (c^{-1}\circ\Abb^{\mg}_{\vecF_{!,+}(\nbigv),-})
  \circ
  (c^{-1}\circ\Abb^{\rd}_{\nbigv,+})\circ\rho_z
  =(2\pi\sqrt{-1})^{-1}\rho_z
  \circ e^{-2\pi\sqrt{-1} N}.
\end{equation}
The equalities 
(\ref{eq;25.3.11.12})
and 
(\ref{eq;25.3.11.13})
are the translation of
the following equalities
on $H^0(\proj^1,\nbigv)_{-1}$:
\begin{equation}
\label{eq;25.3.11.52}
(c^{-1}\circ\Abb^{\rd}_{\vecF_{\ast,-}(\nbigv),+})
\circ
(c^{-1}\circ\Abb^{\mg}_{\nbigv,-})\circ\rho_z
=(2\pi\sqrt{-1})^{-1}\circ\rho_z\circ
e^{2\pi\sqrt{-1}N},
\end{equation}
\begin{equation}
\label{eq;25.3.11.53}
(c^{-1}\circ\Abb^{\rd}_{\vecF_{\ast,+}(\nbigv),-})
\circ
(c^{-1}\circ\Abb^{\mg}_{\nbigv,+})\circ\rho_z
=-(2\pi\sqrt{-1})^{-1}\rho_z.
\end{equation}
Let $F^{\gbigf}$ be the endomorphism of
$\vecF_{\star,\pm}(\nbigv)[\star'0]$
induced by $-N$ on
\[
 H^0(\proj^1,\vecF_{\star,\pm}(\nbigv)[\star' 0])_{-1}
\simeq H^0(\proj^1,\nbigv[\star 0])_0.
\]
We have
\begin{multline}
 (c^{-1}\circ\Abb^{\mg}_{\vecF_{!,\mp}(\nbigv),\pm})
 \circ
 (c^{-1}\circ\Abb^{\rd}_{\nbigv,\mp})\circ\rho_z(v) 
=
 \\
 (2\pi\sqrt{-1})^{-2}
 \rho_z\circ
 s_{\vecF_{!,\mp}(\nbigv),\pm}
 \Bigl(
 w^{-1}
 \Ftilde^{\gbigf}_{\ast,\pm}
 \bigl(
 s_{\nbigv[!0],\mp}
 (-\del_z\otimes\Ftilde_{!,\mp}(v))
 \bigr)
 \Bigr).
\end{multline}
Note that
$w^{-1}
\Ftilde^{\gbigf}_{\ast,\pm}
\bigl(
s_{\nbigv[!0],\mp}
(-\del_z\otimes\Ftilde_{!,\mp}(v))
\bigr)
=\mp
\Ftilde^{\gbigf}_{\ast,\pm}
\bigl(
s_{\nbigv[!0],\mp}
(\Ftilde_{!,\mp}(v))$.
We obtain
\begin{multline}
  s_{\vecF_{!,\mp}(\nbigv),\pm}
 \Bigl(
 w^{-1}
 \Ftilde^{\gbigf}_{\ast,\pm}
 \bigl(
 s_{\nbigv[!0],\mp}
 (-\del_z\otimes\Ftilde_{!,\mp}(v))
 \bigr)
 \Bigr)
 =
\\ \mp
 s_{\vecF_{!,\mp}(\nbigv),\pm}
 \bigl(
 \Ftilde^{\gbigf}_{\ast,\pm}
 \circ
 s_{\nbigv[!0],\mp}
 (\Ftilde_{!,\mp}(v))
 \bigr)
=
 \\
 s_{\vecF_{!,\mp}(\nbigv),\pm}
 \Bigl(
 -\int_{\Gamma_{\ast,\pm}}
 \exp(-N\log\eta)
 e^{\mp\eta}d\eta
 \cdot
 s_{\nbigv(!0),\mp}
 \Bigl(
 \int_{\Gamma_{!\mp}}
 \exp(N\log\zeta)e^{\pm\zeta}\frac{d\zeta}{\zeta}
 \Bigr)(v)
 \Bigr).
\end{multline}
We obtain (\ref{eq;25.3.11.50}) and (\ref{eq;25.3.11.51})
by using Lemma \ref{lem;25.3.11.60} below.
Similarly, we have
\begin{multline}
 (c^{-1}\circ\Abb^{\rd}_{\vecF_{\ast,\mp}\nbigv,\pm})\circ
 (c^{-1}\circ\Abb^{\mg}_{\nbigv,\mp})\circ\rho_z(v)
 =
\\
 (2\pi\sqrt{-1})^{-2}
 \rho_z\Bigl(
 s_{\vecF_{\ast,\mp}\nbigv,\pm}
 \Ftilde^{\gbigf}_{!,\pm}
 \bigl(
 -\del_w\otimes s_{\nbigv,\mp}(z^{-1}\Ftilde_{\ast,\mp}(v))
 \bigr)
 \Bigr)
 =\\
 \pm
 (2\pi\sqrt{-1})^{-2}
 \rho_z\Bigl(
 s_{\vecF_{\ast,\mp}\nbigv,\pm}
 \Ftilde^{\gbigf}_{!,\pm}
 \bigl(
 s_{\nbigv,\mp}(\Ftilde_{\ast,\mp}(v))
 \bigr)
 \Bigr).
\end{multline}
We have
\[
 \pm\Ftilde^{\gbigf}_{!,\pm}
 \circ\Ftilde_{\ast,\mp}
=-\int_{\Gamma_{!,\pm}}
 \exp(-N\log\eta)e^{\mp\eta}\frac{d\eta}{\eta}
 \int_{\Gamma_{\ast,\mp}}
 \exp((\id+N)\log\zeta)e^{\pm\zeta}\frac{d\zeta}{\zeta}.
\]
Then, we obtain
the equalities (\ref{eq;25.3.11.52})
and (\ref{eq;25.3.11.53})
by using Lemma \ref{lem;25.3.11.60} below.
\hfill\qed

\begin{lem}
\label{lem;25.3.11.60}
We obtain the following equalities from 
Lemma {\rm\ref{lem;25.3.11.40}}
and Corollary {\rm\ref{cor;25.3.11.41}}.
\begin{equation}
 \int_{\Gamma_{\ast,+}}
  \exp
  \bigl((\id-N)\log\eta\bigr)e^{-\eta}\frac{d\eta}{\eta}
  \int_{\Gamma_{!,-}}
  \exp\bigl(N\log\zeta\bigr)e^{\zeta}
  \frac{d\zeta}{\zeta}
=2\pi\sqrt{-1}\id,
 \end{equation}
\begin{equation}
 \int_{\Gamma_{\ast,-}}
  \exp
  \bigl((\id-N)\log\eta\bigr)e^{-\eta}\frac{d\eta}{\eta}
  \int_{\Gamma_{!,+}}
  \exp\bigl(N\log\zeta\bigr)e^{\zeta}
  \frac{d\zeta}{\zeta}
=-2\pi\sqrt{-1}e^{-2\pi\sqrt{-1}N},
\end{equation}
\begin{equation}
 \int_{\Gamma_{!,+}}
  \exp
  \bigl(-N\log\eta\bigr)e^{-\eta}\frac{d\eta}{\eta}
  \int_{\Gamma_{\ast,-}}
  \exp\bigl((\id+N)\log\zeta\bigr)e^{\zeta}
  \frac{d\zeta}{\zeta}
=-2\pi\sqrt{-1}e^{2\pi\sqrt{-1}N},
\end{equation}
\begin{equation}
 \int_{\Gamma_{!,-}}
  \exp
  \bigl(-N\log\eta\bigr)e^{-\eta}\frac{d\eta}{\eta}
  \int_{\Gamma_{\ast,+}}
  \exp\bigl((\id+N)\log\zeta\bigr)e^{\zeta}
  \frac{d\zeta}{\zeta}
=2\pi\sqrt{-1}\id.
\end{equation}
\hfill\qed
\end{lem}

\subsection{Proof of Proposition \ref{prop;25.3.10.110}}
\label{subsection;25.3.16.2}

\subsubsection{Complexes and representatives}

We shall use the notation in
\S\ref{subsection;25.3.12.1} and \S\ref{subsection;25.3.12.2}.
There exist the following natural isomorphisms
\begin{equation}
 R^1\pi_{w\ast}
 \bigl(
 \nbigctilde^{\bullet}_{\pm,0}(\nbigv)
 \bigr)
 (\ast 0)
 \simeq
\Fourier_{\pm}(\nbigv[!0])(\ast 0),
\end{equation}
\begin{equation}
 R^1\pi_{w\ast}
 \bigl(
 \nbigctilde^{\bullet}_{\pm,1}(\nbigv)
 \bigr)
 (\ast 0)
 \simeq
\Fourier_{\pm}(\nbigv)(\ast 0).
\end{equation}

There exist the isomorphisms
\begin{equation}
\label{eq;25.3.10.51}
 \Fourier_{\pm}(\nbigv(!0))(\ast 0)_{|w}
 \simeq
 H^1\bigl(
 \proj^1,
 \nbigc^{\bullet}_{\pm,0}(\nbigv)_w
\bigr),
\end{equation}
\begin{equation}
\label{eq;25.3.10.50}
 \Fourier_{\pm}(\nbigv)(\ast 0)_{|w}
 \simeq
 H^1\bigl(
 \proj^1,
 \nbigc^{\bullet}_{\pm,1}(\nbigv)_w
\bigr).
\end{equation}
\begin{lem}
For $v\in A$,
$s_{\nbigv,\pm}(z^{-1}v)_{|w}$
are represented by
$v\otimes dz/z$ of
$\nbigc^1_{\pm,1}(\nbigv)_w$
under the isomorphisms {\rm(\ref{eq;25.3.10.50})}.
\hfill\qed
\end{lem}

Let $C^{\bullet}_{\pm,0,C^{\infty}}(\nbigv)_u$
denote the Dolbeault resolution
of $C^{\bullet}_{\pm,0}(\nbigv)_u$.
Let $\chi_0:\proj^1\to [0,1]$
be a $C^{\infty}$-function such that
$\chi_0=1$ around $0$
and $\chi_0=0$ on $\{|z|\geq 1\}$.
Let $\chi_{\infty}:\proj^1\to[0,1]$
be a $C^{\infty}$-function such that
$\chi_{\infty}=1$ around $\infty$
and $\chi_{\infty}=0$ on $\{|z|\leq 10\}$.
For $v\in A$,
we define
$B_{\pm}(v)\in C^{1}_{\pm,0,C^{\infty}}(\nbigv)_w$
by
\[
B_{\pm}(v)
=(-\del_z\otimes v)\,dz
+(\nabla\pm w\,dz)(\chi_0v)
\mp(\nabla\pm w\,dz)(\chi_{\infty}w^{-1}z^{-1}N(v)).
\]
\begin{lem}
$s_{\nbigv[!0],\pm}(-\del_z\otimes v)_{|w}$
are represented by 
 $B_{\pm}(v)$
under the isomorphisms {\rm(\ref{eq;25.3.10.51})}.
\hfill\qed
\end{lem}

\subsubsection{Duality}

Let $A^{\lor}$ denote the dual space,
which is equipped with the dual endomorphism $N^{\lor}$.
Let $\langle\cdot,\cdot\rangle$ denote the natural pairing of
$A$ and $A^{\lor}$.
We obtain the dual bundle
$\nbigv^{\lor}=
A^{\lor}\otimes
\nbigo_{\proj^1}(\ast\{0,\infty\})$
with the induced connection $d+N^{\lor}\frac{dz}{z}$.

There exist the natural pairings
\[
\Tot\Bigl(
\nbigctilde^{\bullet}_{\pm,0}(\nbigv)
\otimes
\nbigctilde^{\bullet}_{\mp,1}(\nbigv^{\lor})
\Bigr)
\lrarr
\pi_z^{\ast}\Omega^{\bullet}_{\proj^1}.
\]
They induce the following perfect pairings
\[
\langle\cdot,\cdot\rangle_{\pm,0}:
 R^1\pi_{w\ast}\Bigl(
\nbigctilde^{\bullet}_{\pm,0}(\nbigv)
 \Bigr)(\ast 0)
 \otimes
 R^1\pi_{w\ast}\Bigl(
 \nbigctilde^{\bullet}_{\mp,1}(\nbigv^{\lor})
 \Bigr)(\ast 0)
 \lrarr
 \nbigo_{\cnum_w}(\ast 0).
\]
Similarly, we obtain the perfect pairings
\[
\langle\cdot,\cdot\rangle_{\pm,1}:
 R^1\pi_{w\ast}\Bigl(
\nbigctilde^{\bullet}_{\pm,1}(\nbigv)
 \Bigr)(\ast 0)
 \otimes
 R^1\pi_{w\ast}\Bigl(
 \nbigctilde^{\bullet}_{\mp,0}(\nbigv^{\lor})
 \Bigr)(\ast 0)
 \lrarr
 \nbigo_{\cnum_w}(\ast 0).
\]

\begin{lem}
$\bigl\langle
s_{\nbigv[!0],\pm}(-\del_z\otimes v),
s_{\nbigv^{\lor},\mp}(z^{-1}v^{\lor})
\bigr\rangle_{\pm,0}
=-2\pi\sqrt{-1}\langle
v,v^{\lor}
\rangle$. 
\end{lem}
\pf
We obtain
\begin{multline}
 \bigl\langle s_{\nbigv[!0],\pm}(-\del_z\otimes v),
 s_{\nbigv,\mp}(z^{-1}v^{\lor})
 \bigr\rangle_{\pm,0}
 =\int_{\proj^1}
 \bigl\langle
 B_{\pm}(v),
 v^{\lor}(dz/z)
 \bigr\rangle
 =\\
 \int\delbar\chi_0\langle v,v^{\lor}\rangle
 \frac{dz}{z}
 \mp\int
 \delbar(\chi_{\infty})
 w^{-1}
 \langle N(v),v^{\lor}\rangle
 \frac{dz}{z^2}
=-2\pi\sqrt{-1}
 \langle v,v^{\lor}\rangle.
\end{multline}
Thus, we are done.
\hfill\qed

\vspace{.1in}
Similarly, we obtain the following.
\begin{lem}
$\bigl\langle
 s_{\nbigv,\pm}(z^{-1}v),
 s_{\nbigv^{\lor}[!0],\mp}(-\del_z\otimes v^{\lor})
\bigr\rangle_{\pm,1}=2\pi\sqrt{-1}\langle v,v^{\lor}\rangle$.
\hfill\qed
\end{lem}

\subsubsection{Proof of Proposition \ref{prop;25.3.10.110}}

We have
\begin{multline}
\label{eq;25.3.10.111}
\bigl\langle
 c^{-1}\circ\Abb^{\rd}_{\nbigv,\pm}(\rho_z(v)),
 s_{\nbigv^{\lor},\pm}(z^{-1}v^{\lor})
 \bigr\rangle_{\pm,0}
 =\int_{\Gamma_{!,\pm,\theta^u}}
 \bigl\langle
 \exp(N\log z)v,v^{\lor}e^{\mp wz}\frac{dz}{z}
 \bigr\rangle
\\
 = \int_{\Gamma_{!,\pm}}
 \bigl\langle
 \exp(N\log (\zeta w^{-1}))v,v^{\lor}
 \bigr\rangle
  e^{\mp\zeta}\frac{d\zeta}{\zeta}.
\end{multline}
We rewrite (\ref{eq;25.3.10.111}) as 
\begin{multline}
 \bigl\langle
 \exp(-N\log w)
 \Ftilde_{!,\pm}(v),v^{\lor}
 \bigr\rangle
 =\\
 \frac{-1}{2\pi\sqrt{-1}}
 \Bigl\langle
 \rho_w\circ
 s_{\nbigv(!0),\pm}
 (-\del_z\otimes \Ftilde_{!,\pm}(v)),
 s_{\nbigv^{\lor},\mp}(z^{-1}v^{\lor})
 \Bigr\rangle_{\pm,0}.
\end{multline}
Hence, we obtain (\ref{eq;25.3.10.130}).
Similarly, we have
\begin{multline}
\label{eq;25.3.10.120}
\Bigl\langle
c^{-1}\circ\Abb^{\mg}_{\nbigv,\pm}\bigl(
\rho_z(v)
\bigr),
s_{\nbigv^{\lor}(!0),\mp} (-\del_z\otimes v^{\lor})
\Bigr\rangle_{\pm,1}
=\\
 \int_{\Gamma_{\ast,\pm,\theta^u}}
 \bigl\langle
 \exp(N\log z)v,
 -\del_z\otimes v^{\lor}
 \bigr\rangle
 e^{\mp wz}\,dz
 =
 \\
 \int_{\Gamma_{\ast,\pm,\theta}}
 \bigl\langle
 \exp(N\log z)v,
 v^{\lor}
 \bigr\rangle
 (\mp w) e^{\mp wz}dz
 =\mp
 \int_{\Gamma_{\ast,\pm}}
 \bigl\langle
 \exp(N\log(\zeta w^{-1}))v,v^{\lor}
 \bigr\rangle
 e^{\mp\zeta}d\zeta.
\end{multline}
We rewrite (\ref{eq;25.3.10.120}) as
\begin{multline}
 \mp\Bigl\langle
 \exp(-N\log w)
 \int_{\Gamma_{\ast,\pm}}
 \exp(N\log \zeta)v e^{\mp\zeta}\,d\zeta,v^{\lor}
 \Bigr\rangle
 =\\
 \frac{1}{2\pi\sqrt{-1}}
 \Bigl\langle
 s_{\nbigv,\pm}\bigl(
 z^{-1}\Ftilde_{\ast,\pm}(v)
 \bigr),
 s_{\nbigv^{\lor}(!0),\mp}(-\del_z\otimes v^{\lor})
\Bigr\rangle_{1,\pm}.
\end{multline}
Thus, we obtain (\ref{eq;25.3.10.131}).
\hfill\qed

\section[Stokes structure of Fourier transform of
holonomic $\nbigd$-modules]{Stokes structure of Fourier transform of
holonomic $\nbigd$-modules at $\infty$}

\label{subsection;18.6.2.120}

Let $D\subset\cnum$ be a finite subset.
We set $\Dbar=D\cup\{\infty\}$.
Let $\Hol(\proj^1,D,\infty)$
denote the category of 
holonomic $\nbigd_{\proj^1}$-modules $\nbigm$
such that
$\nbigm(\ast D)$ is a meromorphic flat bundle
on $(\proj^1,\Dbar)$.
\index{category $\Hol(\proj^1,D,\infty)$}

Let $\nbigm\in \Hol(\proj^1,D,\infty)$.
We obtain the $\nbigd$-module $\nbigm^{\gbigf}:=\Fourier_+(\nbigm)$
on $\proj^1_w$.
\index{$\nbigd$-module $\nbigm^{\gbigf}$}
(See \S\ref{subsection;25.2.12.20} for the Fourier transform
$\Fourier_+(\nbigm)$.)
There exists a neighbourhood $U_{\infty}$ of $\infty$
such that
$\nbigm^{\gbigf}_{|U_{\infty}}$
is a meromorphic flat bundle
on $(U_{\infty},\infty)$.
Let $u=w^{-1}$ be the coordinate of $\proj^1_w$
around $\infty$.
We set
$\nbigi_D=\{\alpha u^{-1}\,|\,\alpha\in D\}$
and
$\nbigitilde^{\circ}:=
\nbigi_{\infty}(\nbigm^{\gbigf})
\cup
\nbigi_D$.
We shall study the $2\pi\seisuu$-equivariant local system
with Stokes structure
$(\gbigl^{\gbigf}(\nbigm),\vecnbigf)
\in \Loc^{\St}(\nbigitilde^{\circ})$
corresponding to $\nbigm^{\gbigf}_{|U_{\infty}}$.

\subsection{The formal structure of the Fourier transform at infinity}

\label{subsection;18.6.24.1}

The formal structure of $\nbigm^{\gbigf}$ at $\infty$
was studied in \cite{Sabbah-stationary}.

We set
\[
 \nbigp(u^{-1}):=
 \bigcup_{e\in\seisuu_{>0}} u^{-1/e}\cnum[u^{-1/e}].
\]
For any non-zero element 
$f=\sum_{j=1}^Nf_ju^{-j/e}$ of $\nbigp(u^{-1})$,
we set
$\ord_{u^{-1}}(f)=\max\{j/e\,|\,f_j\neq 0\}$.
Let 
$\nbigp_{a}(u^{-1}):=
 \{f\in\nbigp(u^{-1})\,|\,f\neq 0,\,\,\ord_{u^{-1}}(f)=a\}$.
We set 
$\nbigp_{\leq a}(u^{-1}):=
 \{0\}\cup\bigcup_{b\leq a}\nbigp_b(u^{-1})$
and 
$\nbigp_{>a}(u^{-1}):=
 \bigcup_{b> a}\nbigp_b(u^{-1})$.
The following lemma is easy and well known.

\begin{lem}
\label{lem;24.4.9.2}
There exists the decomposition
\begin{equation}
\label{eq;25.2.12.30}
 \nbigm^{\gbigf}_{|\widehat{\infty}}
=\bigoplus_{\alpha \in D}
 \nbigm^{\gbigf}_{\widehat{\infty},\alpha}
\oplus
 \gbigv(\nbigm),
\end{equation}
such that the sets of ramified irregular values of
$\nbigm^{\gbigf}_{\widehat{\infty},\alpha}$
are contained in
$\bigl\{
\alpha u^{-1}+f\,\big|\,
 f\in\nbigp_{<1}(u^{-1})
 \bigr\}$,
and the set of ramified irregular values
of $\gbigv(\nbigm)$ is contained in
 $\nbigp_{>1}(u^{-1})$.
\end{lem}
\pf
We explain an outline of the proof for the convenience of the readers.
For the $\nbigd_{\proj_z^1}$-module
$\nbigg_{\alpha}=\nbigd_{\proj^1}\big/\nbigd_{\proj^1}(z-\alpha)$,
there exists a natural isomorphism
$\nbigg_{\alpha}^{\gbigf}
\simeq
\nbige(\alpha w):=
\bigl(\nbigo_{\proj_w^1}(\ast\infty),d+\alpha\,dw\bigr)$.
If the support of $\nbigm$ is contained in $D$,
there exists an isomorphism
$\nbigm\simeq \bigoplus_{\alpha\in D}\nbigg_{\alpha}^{\oplus m(\alpha)}$
for some non-negative integers $m(\alpha)$.
Hence,
we obtain
$\nbigm^{\gbigf}\simeq
 \bigoplus \nbige(\alpha w)^{m(\alpha)}$.
In general, the kernel and the cokernel of 
$\nbigm\to\nbigm(\ast D)$ are contained in $D$.
Hence, we obtain the claim of the lemma
by using the results for meromorphic flat bundles
(see \S\ref{section;18.6.3.12}).
\hfill\qed

\vspace{.1in}
We obtain the following lemma similarly.
\begin{lem}
\label{lem;25.2.12.112}
Let $D_1\subsetneq D$.
\begin{itemize}
 \item For any $\alpha\in D\setminus D_1$,
the induced morphisms
$\nbigm(!D_1)^{\gbigf}_{\inftyhat,\alpha}
 \to
 \nbigm^{\gbigf}_{\inftyhat,\alpha}
 \to
 \nbigm(\ast D_1)^{\gbigf}_{\inftyhat,\alpha}$
are isomorphisms. 
 \item For any $\alpha\in D_1$,
       the induced morphisms
       $\nbigm(!D)^{\gbigf}_{\inftyhat,\alpha}
       \to
       \nbigm(!D_1)^{\gbigf}_{\inftyhat,\alpha}$
       and
        $\nbigm(\ast D_1)^{\gbigf}_{\inftyhat,\alpha}
       \to
       \nbigm(\ast D)^{\gbigf}_{\inftyhat,\alpha}$
       are isomorphisms.
\hfill\qed
\end{itemize}
\end{lem}

There exists the decomposition
\begin{equation}
\label{eq;20.9.8.10}
  \nbigm^{\gbigf}_{\widehat{\infty},\alpha}
 \otimes
 \bigl(
 \cnum(\!(u)\!),
 d-d(\alpha u^{-1})
 \bigr)
=\bigl(
 \nbigm^{\gbigf}_{\widehat{\infty},\alpha}
 \bigr)_1
\oplus
 \bigl(
 \nbigm^{\gbigf}_{\widehat{\infty},\alpha}
 \bigr)_2
\end{equation}
into the regular part and the irregular part.
The following lemma is also easy and well known,
which can be proved 
by the argument in the proof of Lemma \ref{lem;24.4.9.2}.
\begin{lem}
\label{lem;25.2.12.100}
Note that the natural morphisms
$\nbigm(!D)\lrarr\nbigm\lrarr\nbigm(\ast D)$
induce the following isomorphisms.
\begin{itemize}
\item
$\gbigv(\nbigm(!D))\simeq
 \gbigv(\nbigm)\simeq
 \gbigv(\nbigm(\ast D))$.
\item
$\bigl(
 \nbigm(!D)^{\gbigf}_{\widehat{\infty},\alpha}
 \bigr)_2
\simeq
\bigl(
 \nbigm^{\gbigf}_{\widehat{\infty},\alpha}
 \bigr)_2
\simeq
\bigl(
 \nbigm(\ast D)^{\gbigf}_{\widehat{\infty},\alpha}
 \bigr)_2$.
\item
$\Gr^V_{\gamma}\Bigl(
\bigl(
 \nbigm(!D)^{\gbigf}_{\widehat{\infty},\alpha}
 \bigr)_1
\Bigr)
 \simeq
 \Gr^V_{\gamma}\Bigl(
 \bigl(
 \nbigm^{\gbigf}_{\widehat{\infty},\alpha}
 \bigr)_1
\Bigr)
 \simeq
 \Gr^V_{\gamma}\Bigl(
 \bigl(
 \nbigm(\ast D)^{\gbigf}_{\widehat{\infty},\alpha}
 \bigr)_1
 \Bigr)$
for any $\gamma\not\in \seisuu$.
\hfill\qed
\end{itemize}
\end{lem}

\begin{cor}
$\nbigitilde^{\circ}
 =\nbigi_{\infty}\bigl(
 \nbigm(\ast D)^{\gbigf}
 \bigr)
 \cup \nbigi_D$.
 \hfill\qed 
\end{cor}

For each $\alpha\in D$,
there exists 
$\nbigm_{\alpha}\in\Hol(\proj^1,\alpha,\infty)$
such that 
$\nbigm_{\alpha|U_{\alpha}}\simeq \nbigm_{|U_{\alpha}}$,
and that $\nbigm_{\alpha}$ is regular singular at $\infty$.
Such $\nbigm_{\alpha}$ is unique up to isomorphisms.
The following lemma is also standard.
We shall explain a proof in \S\ref{subsection;25.2.12.110}
for the convenience of readers.
\begin{lem}
\label{lem;18.6.1.21}
For each $\alpha\in D$,
and for any $a\in\cnum$
with $-1<_{\cnum}a\leq_{\cnum}0$
(see {\rm\S\ref{subsection;24.4.11.1}}
for a total order $\leq_{\cnum}$),
there exists a natural isomorphism
\[
 \Gr^V_{a}
 (\nbigm_{\alpha})
 \simeq
  \Gr^V_{a-1}\Bigl(
 \bigl(
 \nbigm^{\gbigf}_{\inftyhat,\alpha}
 \bigr)_1
 \Bigr).
\] 
The induced operator
$-\del_uu$  on 
$\Gr^V_{a}\Bigl(
 \bigl(
 \nbigm^{\gbigf}_{\inftyhat,\alpha}
 \bigr)_1
 \Bigr)$
equals 
 $-\del_z(z-\alpha)$
 under the isomorphism.
\end{lem}

\subsection{Comparison of the graded pieces of the Stokes structure}

\subsubsection{The isomorphisms in the general parts}

We obtain the following proposition
from Lemma \ref{lem;25.2.12.100}
about the formal structure of $\nbigm^{\gbigf}$ at $\infty$.

\begin{prop}
\label{prop;25.2.12.111}
For any $\gminia\in \nbigitilde^{\circ}\setminus\nbigi_D$,
the induced morphisms 
\[
 \Gr^{\vecnbigf}_{\gminia}
 \gbigl^{\gbigf}(\nbigm(!D))
 \lrarr
  \Gr^{\vecnbigf}_{\gminia}
 \gbigl^{\gbigf}(\nbigm)
 \lrarr
 \Gr^{\vecnbigf}_{\gminia}
 \gbigl^{\gbigf}(\nbigm(\ast D))
\]
are isomorphisms. 
\hfill\qed
\end{prop}

We obtain the following proposition from 
Lemma \ref{lem;25.2.12.112}.
\begin{prop}
Let $D_1\subset D$.
\begin{itemize}
 \item For any $\alpha\in D\setminus D_1$,
       the morphisms
\[
 \Gr^{\vecnbigf}_{\alpha u^{-1}}
 \gbigl^{\gbigf}(\nbigm(!D_1))
 \lrarr
  \Gr^{\vecnbigf}_{\alpha u^{-1}}
 \gbigl^{\gbigf}(\nbigm)
 \lrarr
 \Gr^{\vecnbigf}_{\alpha u^{-1}}
 \gbigl^{\gbigf}(\nbigm(\ast D_1))
\]
are isomorphisms. 
\item For any $\alpha\in D_1$,
       the morphisms
 $\Gr^{\vecnbigf}_{\alpha u^{-1}}
 \gbigl^{\gbigf}(\nbigm(!D))
\to
 \Gr^{\vecnbigf}_{\alpha u^{-1}}
 \gbigl^{\gbigf}(\nbigm(!D_1))$
and
$\Gr^{\vecnbigf}_{\alpha u^{-1}}
 \gbigl^{\gbigf}(\nbigm(\ast D_1))
\to
 \Gr^{\vecnbigf}_{\alpha u^{-1}}
 \gbigl^{\gbigf}(\nbigm(\ast D))$
are isomorphisms.
\hfill\qed
\end{itemize}
\end{prop}

Let $\alpha\in D$.
We have the generalized eigen decomposition
with respect to the monodromy automorphism
\[
\Gr^{\vecnbigf}_{\alpha u^{-1}}
\gbigl^{\gbigf}(\nbigm)
=\bigoplus_{b\in\cnum^{\ast}}
\Gr^{\vecnbigf}_{\alpha u^{-1}}
\gbigl^{\gbigf}(\nbigm)_b.
\]
We also obtain the following proposition
from Lemma \ref{lem;25.2.12.100}.
\begin{prop}
\label{prop;25.2.12.130}
The natural morphisms
\[
 \Gr^{\vecnbigf}_{\alpha u^{-1}}
 \gbigl^{\gbigf}(\nbigm(!D))_b
 \lrarr
  \Gr^{\vecnbigf}_{\alpha u^{-1}}
 \gbigl^{\gbigf}(\nbigm)_b
 \lrarr
 \Gr^{\vecnbigf}_{\alpha u^{-1}}
 \gbigl^{\gbigf}(\nbigm(\ast D))_b
\] 
are isomorphisms.
\hfill\qed
\end{prop}

\subsubsection{The graded pieces corresponding to $\alpha u^{-1}\in\nbigi_D$}
\label{subsection;25.2.12.120}

For any $\alpha\in D$,
let $U_{\alpha}$ be a small neighbourhood of $\alpha$.
There exists
$\nbigm_{\alpha}\in \Hol(\proj^1,\alpha,\infty)$
such that
$\nbigm_{\alpha}\simeq \nbigm_{|U_{\alpha}}$
and that $\nbigm_{\alpha}$ is regular singular at $\infty$.
We set $\nbigv_{\alpha}=\nbigm_{\alpha}(\ast\alpha)$.
We obtain the regular meromorphic flat bundle
$\Gr^{\vecnbigf}_0(\nbigv_{\alpha})$ on $(U_{\alpha},\alpha)$
as the graduation of $\nbigv_{\alpha}$ with respect to
the Stokes structure.
It naturally extends to the regular meromorphic flat bundle
on $(\proj^1,\{\alpha,\infty\})$,
which is also denoted by
$\Gr^{\vecnbigf}_0(\nbigv_{\alpha})$.
We have
$\Gr^V_a(\Gr^{\vecnbigf}_0(\nbigv_{\alpha}))
=\Gr^V_a(\nbigv_{\alpha})$.
Let $\Gr^{\vecnbigf}_0(\nbigm_{\alpha})'$
be the regular holonomic $\nbigd$-module
on $U_{\alpha}$
corresponding to
$\Gr^{\vecnbigf}_0(\nbigv_{\alpha})$
and the morphisms
\[
\Gr^V_{-1}(\Gr^{\vecnbigf}_0\nbigv_{\alpha})
=\Gr^V_{-1}(\nbigv_{\alpha})
\to
\Gr^V_0(\nbigm_{\alpha})\to
\Gr^V_{-1}(\nbigv_{\alpha})
=\Gr^V_{-1}(\Gr^{\vecnbigf}_0\nbigv_{\alpha}).
\]
We obtain the $\nbigd_{\proj^1}$-module
$\Gr^{\vecnbigf}_0(\nbigm_{\alpha})$
such that
$\Gr^{\vecnbigf}_0(\nbigm_{\alpha})_{|U_{\alpha}}
=\Gr^{\vecnbigf}_0(\nbigm_{\alpha})'$,
and that $\Gr^{\vecnbigf}_0(\nbigm_{\alpha})$
is regular singular at $\infty$.
We have $\Gr^{\vecnbigf}_0(\nbigm_{\alpha})(\ast \alpha)
=\Gr^{\vecnbigf}_0(\nbigv_{\alpha})$.

\begin{prop}
\label{prop;24.4.14.30}
For any $\nbigm\in\Hol(\proj^1,D,\infty)$,
there exist the natural isomorphisms
\begin{equation}
\label{eq;24.4.14.20}
 \Gr^{\vecnbigf}_{\alpha u^{-1}}\bigl(
 \gbigl^{\gbigf}(\nbigm),\vecnbigf
 \bigr)
 \simeq
 \bigl(
 \gbigl^{\gbigf}(\Gr^{\vecnbigf}_0\nbigm_{\alpha}),\vecnbigf
 \bigr)
 \quad(\alpha\in D).
\end{equation}
The morphisms are functorial in the sense that
for any morphism $\nbigm_1\to\nbigm_2$ in $\Hol(\proj^1,D,\infty)$,
the following diagrams are commutative:
\begin{equation}
\begin{CD}
 \Gr^{\vecnbigf}_{\alpha u^{-1}}\bigl(
 \gbigl^{\gbigf}(\nbigm_1),\vecnbigf
 \bigr)
 @>{\simeq}>>
 \bigl(
 \gbigl^{\gbigf}(\Gr^{\vecnbigf}_0\nbigm_{1,\alpha}),\vecnbigf
 \bigr)
 \\
 @VVV @VVV \\
  \Gr^{\vecnbigf}_{\alpha u^{-1}}\bigl(
 \gbigl^{\gbigf}(\nbigm_2),\vecnbigf
 \bigr)
 @>{\simeq}>>
 \bigl(
 \gbigl^{\gbigf}(\Gr^{\vecnbigf}_0\nbigm_{2,\alpha}),\vecnbigf
 \bigr).
\end{CD}
\end{equation}
\end{prop}
\pf
For any $D_1\subset D$,
we consider the full subcategory
$\Hol(\proj^1,D,\infty)_{D_1}$
of objects $\nbigm$ in $\Hol(\proj^1,D,\infty)$
satisfying the following condition.
\begin{itemize}
 \item $\nbigm_{\alpha}=\nbigm_{\alpha}(!\alpha)$
       or
       $\nbigm_{\alpha}=\nbigm_{\alpha}(\ast\alpha)$
for each $\alpha\in D\setminus D_1$.
\end{itemize}
We have already obtained the functorial isomorphisms
(\ref{eq;24.4.14.20})
for objects in $\Hol(\proj^1,D,\infty)_{\emptyset}$.
Let $D_1\subset D$.
Let $\alpha\in D_1$
and $D_2:=D_1\setminus\{\alpha\}$.
Suppose that we have already obtained the functorial
isomorphisms
(\ref{eq;24.4.14.20})
for objects in $\Hol(\proj^1,D,\infty)_{D_2}$.

Let $\alpha\in\cnum$.
We set
$\gbigi^{a,b}_{\alpha}
=\nbigo_{\proj^1}(\ast\{\alpha,\infty\})\otimes A^{a,b}$
which is equipped with the connection
$\nabla_{\alpha}=d+N_A\cdot\frac{dz}{z-\alpha}$.
We may apply the construction in \S\ref{subsection;24.4.14.1}
to objects $\nbigm$ in $\Hol(\proj^1,D,\infty)$
by using $\gbigi^{a,b}_{\alpha}$.
The obtained functors are denoted by
$\Pi^{a,b}_{\alpha,\star}(\nbigm)$ $(\star=!,\ast)$,
$\Pi^{a,b}_{\alpha,\ast!}(\nbigm)$,
$\Xi_{\alpha}(\nbigm)$,
$\psi^{(a)}_{\alpha}(\nbigm)$
and $\phi_{\alpha}(\nbigm)$.
We have
$\psi^{(a)}_{\alpha}(\nbigm)
=\psi^{(a)}_{\alpha}(\nbigm_{\alpha})$
and
$\phi_{\alpha}(\nbigm)
=\phi_{\alpha}(\nbigm_{\alpha})$.

Because
$\Pi^{a,b}_{\alpha,\ast !}(\nbigm)$
is the cokernel of
$\Pi^{b,N}_{\alpha,!}(\nbigm)
\to
\Pi^{a,N}_{\alpha,\ast}(\nbigm)$,
\[
\Gr^{\vecnbigf}_0\bigl(
 \gbigl^{\gbigf}(\Pi^{a,b}_{\alpha,\ast !}(\nbigm)),\vecnbigf
 \bigr)
\]
is the cokernel of 
$\Gr^{\vecnbigf}_0\bigl(
 \gbigl^{\gbigf}(\Pi^{b,N}_{\alpha,!}(\nbigm)),\vecnbigf
 \bigr)
 \to
\Gr^{\vecnbigf}_0\bigl(
 \gbigl^{\gbigf}(\Pi^{a,N}_{\alpha,\ast}(\nbigm)),\vecnbigf
 \bigr)$.
Similarly,
\[
\bigl(
 \gbigl^{\gbigf}(\Pi^{a,b}_{\alpha,\ast !}(\Gr^{\vecnbigf}_0\nbigm_{\alpha})),\vecnbigf
 \bigr) 
\]
is the cokernel of 
$\bigl(
 \gbigl^{\gbigf}(\Pi^{b,N}_{\alpha,!}(\Gr^{\vecnbigf}_0\nbigm_{\alpha})),
 \vecnbigf
 \bigr) 
 \to
 \bigl(
 \gbigl^{\gbigf}(\Pi^{a,N}_{\alpha,\ast}(\Gr^{\vecnbigf}_0\nbigm_{\alpha})),
 \vecnbigf
 \bigr)$.
Hence, we obtain the morphisms
(\ref{eq;24.4.14.20})
for $\Pi^{a,b}_{\alpha,\ast !}(\nbigm)$.
In particular,
we obtain the morphisms
(\ref{eq;24.4.14.20})
for $\Xi_{\alpha}(\nbigm)$
and $\psi^{(a)}_{\alpha}(\nbigm)$.
Note that
the isomorphisms for
$\psi^{(a)}_{\alpha}(\nbigm)$
equal the identity by Lemma \ref{lem;24.4.14.21}.

We can reconstruct $\nbigm$
as the cohomology of
\[
 \psi^{(1)}_{\alpha}(\nbigm)
 \to
 \Xi_{\alpha}(\nbigm)
 \oplus
 \phi_{\alpha}(\nbigm)
 \to
 \psi^{(0)}_{\alpha}(\nbigm).
\]
Hence,
$\Gr^{\vecnbigf}_{\alpha u^{-1}}\gbigl^{\gbigf}(\nbigm)$
is reconstructed as the cohomology of
\[
\gbigl^{\gbigf}
 (\psi^{(1)}_{\alpha}(\nbigm))
 \to
  \Gr^{\vecnbigf}_{\alpha u^{-1}}\bigl(
 \gbigl^{\gbigf}\bigl(
 \Xi_{\alpha}(\nbigm)
 \bigr)
 \bigr)
 \oplus
 \gbigl^{\gbigf}\bigl(
 \phi_{\alpha}(\nbigm)
 \bigr)
 \to
 \gbigl^{\gbigf}\bigl(
 \psi^{(0)}_{\alpha}(\nbigm)
 \bigr).
\]
Similarly,
$\gbigl^{\gbigf}\Gr^{\vecnbigf}_0(\nbigm_{\alpha})$
is reconstructed as the cohomology of
{\small
\[
\gbigl^{\gbigf}(
\psi^{(1)}_{\alpha}(\Gr^{\vecnbigf}_0(\nbigm_{\alpha})))
 \to
\gbigl^{\gbigf}(\Xi_{\alpha}(\Gr^{\vecnbigf}_0(\nbigm_{\alpha})))
 \oplus
\gbigl^{\gbigf}(\phi_{\alpha}(\Gr^{\vecnbigf}_0(\nbigm_{\alpha})))
\to
\gbigl^{\gbigf}(
 \psi^{(0)}_{\alpha}(\Gr^{\vecnbigf}_0(\nbigm))).
\]}
Hence, we obtain the isomorphism for $\nbigm$.
\hfill\qed

\subsubsection{Description in terms of the $V$-filtrations}

Let $\rho_{\alpha}:\proj^1\to\proj^1$
be defined by $\rho_{\alpha}(z)=z+\alpha$.
We set
$\nbigm_{\alpha}^0:=
\rho_{\alpha}^{\ast}\Gr^{\vecnbigf}_0(\nbigm_{\alpha})
\in\Hol(\proj^1,0,\infty)$.
We also set $\nbigv_{\alpha}^0=\nbigm_{\alpha}^0(\ast 0)$.
We use the notation in \S\ref{subsection;25.3.12.42}.
By Proposition \ref{prop;25.3.11.200},
we obtain the following proposition.
\begin{prop}
There exist the following commutative diagram:
\begin{equation}
\begin{CD}
 \psitilde(\nbigm_{\alpha}^0)
 @>{\cantilde_{\nbigm_{\alpha}^0}\circ \Phi_{!,\pm}}>>
 \phitilde(\nbigm_{\alpha}^0)
 @>{(\Phi_{\ast,\pm})^{-1}\circ\var_{\nbigm_{\alpha}^0}}>>
 \psitilde(\nbigm_{\alpha}^0)
 \\
 @V{\simeq}V{\Abb^{\rd}_{\nbigv,\pm}\circ\rhotilde_z}V
 @V{\simeq}V{\Psi_{\nbigm^0_{\alpha},\pm}}V
 @V{\simeq}V{\Abb^{\mg}_{\nbigv,\pm}\circ\rhotilde_z}V \\
 H^0(\real,\gbigl^{\gbigf}_{!}(\nbigv_{\alpha}^0))
 @>>>
 H^0(\real,\gbigl^{\gbigf}(\nbigm_{\alpha}^0))
 @>>>
 H^0(\real,\gbigl^{\gbigf}_{\ast}(\nbigv_{\alpha}^0)).
\end{CD}
\end{equation}
Here, the lower horizontal arrows are the natural morphisms.
The monodromy automorphisms of
$H^0(\real,\gbigl^{\gbigf}_!(\nbigv_{\alpha}^0))$,
$H^0(\real,\gbigl^{\gbigf}(\nbigm_{\alpha}^0))$
and
$H^0(\real,\gbigl^{\gbigf}(\nbigv_{\alpha}^0))$
are equal to
$M_{\psitilde(\nbigm^0_{\alpha})}$,
$M_{\phitilde(\nbigm^0_{\alpha})}$,
and $M_{\psitilde(\nbigm^0_{\alpha})}$,
respectively.
\hfill\qed
\end{prop}

\subsection{Description of $(\gbigl^{\gbigf}(\nbigm),\vecnbigf)$
as an extension}

Recall that $\Csf(D)$
denotes the category of maps $D\lrarr \{!,\circ,\ast\}$,
where 
\[
 \Hom_{\Csf(D)}(\varrho_1,\varrho_2):=
 \prod_{\alpha\in D}
 \Hom_{\Csf_1}\bigl(
 \varrho_1(\alpha),\varrho_2(\alpha)
 \bigr).
\]

Let $\nbigm\in\Hol(\proj^1,D,\infty)$.
For any $\varrho\in\Csf(D)$,
let $\nbigm(\varrho)$ be the $\nbigd_{\proj^1}$-module
determined by
the conditions
$\nbigm(\varrho)\otimes\nbigo(\ast D)=\nbigm(\ast D)$
and
\[
 \nbigm(\varrho)_{|U_{\alpha}}
\simeq
\left\{
 \begin{array}{ll}
 \nbigm_{|U_{\alpha}}
 & (\varrho(\alpha)=\circ), \\
\nbigm(\ast\alpha)
 & (\varrho(\alpha)=\ast), \\
\nbigm(!\alpha)
 & (\varrho(\alpha)=!).
 \end{array}
\right.
\]
We set
$\nbige(\varrho):=
 \bigl(
 \gbigl^{\gbigf}(\nbigm(\varrho)),\vecnbigf
 \bigr)$
in $\Loc^{\St}(\nbigitilde^{\circ})$.
It induces a functor $\nbige$
from $\Csf(D)$
to $\Loc^{\St}(\nbigitilde^{\circ})$.
By Proposition \ref{prop;25.2.12.111},
$\Gr^{\vecnbigf}_{\gminia}(\nbige(\varrho))$
are independent of $\varrho$ if $\gminia\not\in\nbigi_D$,
and that
$\Gr^{\vecnbigf}_{\alpha u^{-1}}(\nbige(\varrho))$
($\alpha\in D$)
depend only on $\varrho(\alpha)$.

Let
$\iota^{\ast}(\nbige):
\Dsf(D)\lrarr \Loc^{\St}(\nbigitilde^{\circ})$
denote the naturally defined functor.
We set $\nbigv=\nbigm(\ast D)$.
We have
$\iota^{\ast}(\nbige)
=(\gbigl^{\gbigf}_{\varrho}(\nbigv),\vecnbigf)$.
We obtain the following theorem
from Proposition \ref{prop;24.4.14.30}.
(See \S\ref{section;21.4.30.1} for the notion of extensions
of local systems with Stokes structure.)
\begin{thm}
\label{thm;24.4.14.40}
The functor $\nbige$ is
the extension of the base tuple
$\gbigl^{\gbigf}_{\varrho}(\nbigv)$
$(\varrho\in\Dsf(D))$
by the morphisms of $2\pi\seisuu$-equivariant local systems
\begin{equation}
\label{eq;25.2.12.140}
 \gbigl^{\gbigf}\bigl(
 \Gr^{\vecnbigf}_0(\nbigm_{\alpha})(!\alpha)
 \bigr)
 \lrarr
 \gbigl^{\gbigf}\bigl(
 \Gr^{\vecnbigf}_0(\nbigm_{\alpha})
 \bigr)
 \lrarr
 \gbigl^{\gbigf}\bigl(
 \Gr^{\vecnbigf}_0(\nbigm_{\alpha})(\ast\alpha)
 \bigr).
\end{equation}
Here, $\Gr^{\vecnbigf}_0(\nbigm_{\alpha})$
are the regular holonomic $\nbigd$-modules
as in {\rm\S\ref{subsection;25.2.12.120}}.
\hfill\qed
\end{thm}

\subsection{Reductions}

\subsubsection{Reductions at $0$}

Let $\nbigm\in\Hol(\proj^1,0,\infty)$
which is regular singular at $\infty$.
We set $\nbigv=\nbigm(\ast 0)$
and $\omega=-\ord\nbigi(\nbigv)$.
We set
$V=\nbigs_{\omega}(\nbigv)$
and $\omega^{\circ}=(1+\omega)^{-1}\omega$.
We obtain the following corollary
from Theorem \ref{thm;24.3.24.1},
Theorem \ref{thm;24.3.25.50},
and Theorem \ref{thm;24.4.14.40}.
\begin{cor}
\label{cor;24.4.15.10}
There exists the functorial isomorphism
\[
 \nbigt_{\omega^{\circ}}(\gbigl^{\gbigf}(\nbigm),\vecnbigf)
 \simeq
 (\gbigl^{\gbigf}(\nbigt_{\omega}\nbigm),\vecnbigf).
\]
The $2\pi\seisuu$-equivariant local system with Stokes structure
$\nbigs_{\omega^{\circ}}\bigl(
\gbigl^{\gbigf}(\nbigm),\vecnbigf
 \bigr)$
is obtained as the extension of
the base tuple
$(\gbigl^{\gbigf}_!(V),\vecnbigf)
 \to
 (\gbigl^{\gbigf}_{\ast}(V),\vecnbigf)$
by the morphisms of the $2\pi\seisuu$-equivariant local systems
\begin{equation}
 \gbigl^{\gbigf}_!(\nbigt_{\omega}(V))
\lrarr
 \gbigl^{\gbigf}(\nbigt_{\omega}(\nbigm))
\lrarr
 \gbigl^{\gbigf}_{\ast}(\nbigt_{\omega}(V)).
\end{equation}
\hfill\qed
\end{cor}

\subsubsection{Reductions at finite place}

Let $\nbigm\in\Hol(\proj^1,D,\infty)$
which is regular singular at $\infty$.
Let $V$ denote the regular singular meromorphic flat bundle
on $(\proj^1,D\cup\{\infty\})$
associated with the local system corresponding to $\nbigm(\ast D)$.
We set $\vecnbigf^{(1)}=\pi_{1\ast}(\vecnbigf)$
on $\gbigl^{\gbigf}(\nbigm)$.
We obtain the following corollary
from Proposition \ref{prop;24.3.30.1},
Proposition \ref{prop;24.4.14.50},
and Theorem \ref{thm;24.4.14.40}.
\begin{cor}
\label{cor;24.4.15.13}
There exist the functorial isomorphisms
\[
 \Gr^{\vecnbigf^{(1)}}_{\alpha u^{-1}}
 (\gbigl^{\gbigf}(\nbigm),\vecnbigf)
 \simeq
 (\gbigl^{\gbigf}(\nbigm_{\alpha}),\vecnbigf).
\]
Here, $\nbigm_{\alpha}\in\Hol(\proj^1,\alpha,\infty)$
are the $\nbigd$-modules as in {\rm\S\ref{subsection;18.6.24.1}}.
The $2\pi\seisuu$-equivariant local system with Stokes structure
$(\gbigl^{\gbigf}(\nbigm),\vecnbigf^{(1)})$
is the extension of 
the base tuple
$(\gbigl^{\gbigf}_{\varrho}(V),\vecnbigf)$
$(\varrho\in \Dsf(D))$
by the $2\pi\seisuu$-equivariant local systems
\[
 \gbigl^{\gbigf}_{!}(V_{\alpha})
\lrarr
 \gbigl^{\gbigf}(\nbigm_{\alpha})
\lrarr
 \gbigl^{\gbigf}_{\ast}(V_{\alpha}).
\] 
Here, $V_{\alpha}$ denote the regular singular meromorphic 
flat bundles on $(\proj^1,\{\alpha,\infty\})$
obtained as the extension of the restriction of
$V$ to a neighbourhood of $\alpha$. 
\hfill\qed
\end{cor}

\subsubsection{Reductions at infinity}

Let $\nbigm\in\Hol(\proj^1,D,\infty)$.
We obtain the following proposition
from Proposition \ref{prop;24.3.17.121}
and Theorem \ref{thm;24.4.14.40}.
\begin{prop}
\label{prop;24.4.15.40}
There exists the functorial isomorphism
\[
 \bigl(
 \gbigl^{\gbigf}(\nbigstilde^{\infty}_1\nbigm),\vecnbigf
 \bigr)
 \simeq
  \bigl(
 \gbigl^{\gbigf}(\nbigm),\vecnbigf
 \bigr).
\]
\hfill\qed
\end{prop}

We set
$\omega=
\min\bigl\{
 -\ord(\gminia)\,|\,\gminia\in\nbigi_{\infty}(\nbigm)
 \bigr\}$.
Suppose $\omega>1$
and put $\omega^{\circ}=(\omega-1)^{-1}\omega$.
We obtain the meromorphic flat bundle
$V_{\infty}=\nbigttilde_{\omega}(\nbigm)$
on $(\proj^1,\{0,\infty\})$.
Let $V^{\reg}_{\infty}=\nbigstilde^{\infty}_{\omega}(V_{\infty})$.
We obtain the following corollary from
Theorem \ref{thm;24.3.29.40},
Theorem \ref{thm;24.3.29.41},
and Theorem \ref{thm;24.4.14.40}.
\begin{cor}
\label{cor;24.4.15.51}
There exists the functorial isomorphism
\[
 \nbigt_{\omega^{\circ}}\bigl(
 \gbigl^{\gbigf}(\nbigm),\vecnbigf
 \bigr)
\simeq
 \bigl(
 \gbigl^{\gbigf}(\nbigstilde^{\infty}_{\omega}\nbigm),
 \vecnbigf
 \bigr).
\]
The $2\pi\seisuu$-equivariant local system
with Stokes structure
$\nbigs_{\omega^{\circ}}\bigl(
\gbigl^{\gbigf}(\nbigm),\vecnbigf
 \bigr)$
is obtained as the extension of
$(\gbigl^{\gbigf}_{!}(V_{\infty}),\vecnbigf)
 \to
(\gbigl^{\gbigf}_{\ast}(V_{\infty}),\vecnbigf)$
by the following natural morphisms of
$2\pi\seisuu$-equivariant local systems:
\begin{equation}
 \gbigl^{\gbigf}_!(V^{\reg}_{\infty})
 \to
 \gbigl^{\gbigf}(\nbigstilde^{\infty}_{\omega}(\nbigm))
 \to
 \gbigl^{\gbigf}_{\ast}(V^{\reg}_{\infty}).
\end{equation}
\hfill\qed
\end{cor}

\subsection{Appendix: Proof of Lemma \ref{lem;18.6.1.21}}
\label{subsection;25.2.12.110}
We explain only an outline of the proof
just for the convenience of the readers.
For each $\alpha\in D$,
there exists the decomposition
\[
 \nbigm_{|\alphahat}
=\nbigmhat_{\alpha}^{\reg}
\oplus
 \nbigmhat_{\alpha}^{\irr},
\]
where $\nbigmhat_{\alpha}^{\reg}$ is regular singular,
and $\nbigmhat_{\alpha}^{\irr}$ is isomorphic
to the formal completion of
a meromorphic flat bundle
whose set of ramified irregular values
does not contain $0$.
There exists a good lattice pair
for $\nbigmhat_{\alpha}^{\irr}$
in the sense of \cite{Bloch-Esnault1},
i.e.,
sub-lattices
$\nbigchat_{\alpha}^{0,\irr}
\subset
\nbigchat_{\alpha}^{1,\irr}
\subset\nbigmhat_{\alpha}^{\irr}$
such that
(i) $\del_z\bigl(
 \nbigchat_{\alpha}^{0,\irr}
 \bigr)
\subset
 \nbigchat_{\alpha}^{1,\irr}$,
(ii)
$\nbigchat_{\alpha}^{0,\irr}
\stackrel{\del_z}{\lrarr}
\nbigchat_{\alpha}^{1,\irr}$
is naturally quasi isomorphic to
$\nbigmhat_{\alpha}^{\irr}
\stackrel{\del_z}{\lrarr}
 \nbigmhat_{\alpha}^{\irr}$.
We set
\[
 \nbigchat^0_{\alpha}:=
 V_{-1}(\nbigmhat^{\reg}_{\alpha})
\oplus 
 \nbigchat^{0,\irr}_{\alpha},
\quad\quad
 \nbigchat^1_{\alpha}:=
 V_{0}(\nbigmhat^{\reg}_{\alpha})
 \oplus 
 \nbigchat^{1,\irr}_{\alpha}.
\]

There exists the decomposition
$\nbigm_{|\inftyhat}
=\nbigg_1\oplus\nbigg_2$,
where the set of ramified irregular values of
$\nbigg_1$ is contained in $\nbigp_{\leq 1}(z)$,
and the set of ramified irregular values of $\nbigg_2$
is contained in $\nbigp_{>1}(z)$.
Let $\nbigg_1^{\DM}\subset\nbigg_1$
denote the Deligne-Malgrange lattice.
For any $\ell\in\seisuu$,
$z^{\ell}\nbigg_1^{\DM}
\stackrel{\del_z+w}{\lrarr}
z^{\ell}\nbigg_1^{\DM}$
is naturally quasi isomorphic to
$\nbigg_1
\stackrel{\del_z+w}{\lrarr}
 \nbigg_1$
if $|w|$ is sufficiently large.
There exist lattices 
$\nbigchat_{\infty}^{0}(\nbigg_2)
\subset
\nbigchat_{\infty}^{1}(\nbigg_2)
\subset\nbigg_2$
such that 
(i) 
$\del_z\nbigchat_{\infty}^{0}(\nbigg_2)
\subset
\nbigchat_{\infty}^{1}(\nbigg_2)$,
(ii) the complexes
$\nbigchat_{\infty}^0(\nbigg_2)
\stackrel{\del_z}{\lrarr}
\nbigchat_{\infty}^1(\nbigg_2)$
and
$\nbigg_2
\stackrel{\del_z}{\lrarr}
\nbigg_2$
are naturally quasi-isomorphic.
Note that
$z^{\ell}\nbigchat_{\infty}^0(\nbigg_2)
\stackrel{\del_z+w}{\lrarr}
z^{\ell}\nbigchat_{\infty}^1(\nbigg_2)$
and
$\nbigg_2
\stackrel{\del_z+w}{\lrarr}
\nbigg_2$
are naturally quasi-isomorphic
for any $w$.
We set
\[
 \nbigp_{\ell}\nbigchat_{\infty}^0:=
 z^{\ell}
 \Bigl(
 \nbigg_1^{\DM}
\oplus
 \nbigchat^{0}_{\infty}(\nbigg_2)
\Bigr),
\quad\quad
  \nbigp_{\ell}\nbigchat_{\infty}^1:=
  z^{\ell}
  \Bigl(
 \nbigg_1^{\DM}
\oplus
 \nbigchat^{1}_{\infty}(\nbigg_2)
 \Bigr).
\]

There exist the $\nbigo_{\proj^1}$-coherent submodules
$\nbigp_{\ell}\nbigc^0\subset
\nbigp_{\ell}\nbigc^1\subset\nbigm$
such that
$\nbigp_{\ell}\nbigc^i_{|\alphahat}
\simeq
 \nbigchat^{i}_{\alpha}$
for any $\alpha$,
and 
$\nbigp_{\ell}\nbigc^i_{|\inftyhat}
\simeq
 \nbigp_{\ell}\nbigchat^{i}_{\infty}$.
The complex
$\nbigp_{\ell}\nbigc^0
 \stackrel{\del_z+w}{\lrarr}
\nbigp_{\ell}\nbigc^1$
is naturally quasi-isomorphic to
$\nbigm\stackrel{\del_z+w}{\lrarr}\nbigm$
if $|w|$ is sufficiently large.
If $\ell$ is sufficiently large,
we may assume that
$H^1(\proj^1,\nbigp_{\ell}\nbigc^i)=0$.

We take a small neighbourhood 
$U$ of $\infty$ in $\proj^1_w$
with the coordinate $u=w^{-1}$.
Let $p_z:\proj^1_z\times U\lrarr\proj^1_z$
and $p_w:\proj^1_z\times U\lrarr U$
denote the projections.

We obtain the following complexes
$\nbigp_{\ell}\nbigctilde^{\bullet}(\nbigm)$
on $\proj^1_z\times\proj^1_w$:
\[
\begin{CD}
 p_z^{\ast}\nbigp_{\ell}\nbigc^0\otimes
 p_w^{\ast}\nbigo_{U}(-\{\infty\})
 @>{\del_z+u^{-1}}>>
 p_z^{\ast}\nbigp_{\ell}\nbigc^1.
\end{CD}
\]
If $\ell$ is large,
we obtain $R^1p_{w\ast}\nbigp_{\ell}\nbigctilde^{i}=0$.
Hence,
$\nbige_{\ell}:=
 R^1p_{w\ast}
 \nbigp_{\ell}\nbigctilde^{\bullet}$
on $U$
is obtained as the cokernel of
the morphism of coherent $\nbigo_U$-modules
$p_{w\ast}(\nbigp_{\ell}\nbigctilde^0)
\lrarr
p_{w\ast}(\nbigp_{\ell}\nbigctilde^1)$. 
There exists a natural isomorphism
$\nbige_{\ell}(\ast \infty)
\simeq
 \Fourier_+(\nbigm)_{|U}$.
As in \cite{Bloch-Esnault1},
we obtain
\[
 \nbige_{\ell|\inftyhat}
\simeq
 \bigoplus_{\alpha\in D}
 \Cok\Bigl(
 u\nbigchat^{0}_{\alpha}[\![u]\!]
\stackrel{\del_z+u^{-1}}\lrarr
 \nbigchat^{1}_{\alpha}[\![u]\!]
 \Bigr)
\oplus
 \Cok\Bigl(
 u\nbigp_{\ell}\nbigchat^{0}_{\infty}[\![u]\!]
\stackrel{\del_z+u^{-1}}\lrarr
 \nbigp_{\ell}\nbigchat^1_{\infty}[\![u]\!]
 \Bigr).
\]
There exists the decomposition
\begin{multline}
 \Cok\Bigl(
 u\nbigchat^{0}_{\alpha}[\![u]\!]
\stackrel{\del_z+u^{-1}}\lrarr
 \nbigchat^{1}_{\alpha}[\![u]\!]
 \Bigr)
\simeq \\
 \Cok\Bigl(
 uV_{-1}(\nbigmhat^{\reg}_{\alpha})[\![u]\!]
\stackrel{\del_z+u^{-1}}{\lrarr}
 V_0(\nbigmhat^{\reg}_{\alpha})[\![u]\!]
 \Bigr)
\oplus
 \Cok\Bigl(
 u\nbigchat^{0,\irr}_{\alpha}[\![u]\!]
\stackrel{\del_z+u^{-1}}\lrarr
  \nbigchat^{1,\irr}_{\alpha}[\![u]\!]
 \Bigr).
\end{multline}
It induces the decomposition (\ref{eq;20.9.8.10}),
i.e.,
\begin{equation}
\label{eq;20.9.8.11}
  \Cok\Bigl(
 uV_{-1}(\nbigmhat^{\reg}_{\alpha})[\![u]\!]
\stackrel{\del_z+u^{-1}}{\lrarr}
 V_0(\nbigmhat^{\reg}_{\alpha})[\![u]\!]
 \Bigr)
 (\ast\infty)
\simeq
 (\nbigm^{\gbigf}_{\inftyhat,\alpha})_1
 \otimes
 \bigl(\cnum(\!(u)\!),d+d(\alpha u^{-1})\bigr),
\end{equation}
\begin{equation}
  \Cok\Bigl(
 u\nbigchat^{0,\irr}_{\alpha}[\![u]\!]
\stackrel{\del_z+u^{-1}}\lrarr
  \nbigchat^{1,\irr}_{\alpha}[\![u]\!]
  \Bigr)(\ast\infty)
 \simeq
   (\nbigm^{\gbigf}_{\inftyhat,\alpha})_2
 \otimes
 \bigl(\cnum(\!(u)\!),d+d(\alpha u^{-1})\bigr).
\end{equation}

We take a vector subspace
$H_{\alpha}\subset V_0\nbigmhat_{\alpha}^{\reg}$
such that
$H_{\alpha}\oplus
 V_{-1}(\nbigmhat^{\reg}_{\alpha})
=V_{0}(\nbigmhat^{\reg}_{\alpha})$
as a $\cnum$-vector space
such that
$\del_z(z-\alpha)H_{\alpha}\subset H_{\alpha}$.
It is easy to see that
$V_0(\nbigmhat^{\reg}_{\alpha})[\![u]\!]
=H_{\alpha}[\![u]\!]
\oplus\Image(\del_z+u^{-1})$.
Hence, we obtain the following 
$\cnum[\![u]\!]$-isomorphism
\[
H_{\alpha}[\![u]\!]
\simeq
  \Cok\Bigl(
 uV_{-1}(\nbigmhat^{\reg}_{\alpha})[\![u]\!]
\stackrel{\del_z+u^{-1}}{\lrarr}
 V_0(\nbigmhat^{\reg}_{\alpha})[\![u]\!]
 \Bigr)
=:\nbigl'_{\alpha}(\nbigm).
\]
For any $f\in H_{\alpha}\subset\nbigl'_{\alpha}(\nbigm)$,
we obtain
\begin{equation}
\label{eq;20.9.8.12}
 u(\del_u
 +\alpha u^{-2})f=
 \del_z(z-\alpha)f
 \in \nbigl_{\alpha}(\nbigm)'
\end{equation}
holds
in the $\cnum(\!(u)\!)$-module (\ref{eq;20.9.8.11}).
We can easily check that
$V_{-1}(\nbigmhat_{\alpha}^{\reg})
\subset
\Image(\del_z+u^{-1})
+uV_0(\nbigmhat_{\alpha}^{\reg})[\![u]\!]$.

Let $\nbigl_{\alpha}(\nbigm)$
denote the $\cnum[\![u]\!]$-lattice of
$\bigl(
 \nbigm^{\gbigf}_{\inftyhat,\alpha}
 \bigr)_1$
 induced by
 $\nbigl'_{\alpha}(\nbigm)$
 and the isomorphism (\ref{eq;20.9.8.11}).
By (\ref{eq;20.9.8.12}),
we obtain
$u\del_u\nbigl_{\alpha}
\subset\nbigl_{\alpha}$.
Moreover, the endomorphism of
\[
\nbigl_{\alpha}(\nbigm)/u\nbigl_{\alpha}(\nbigm)
\simeq
H_{\alpha}\simeq
 V_0(\nbigmhat^{\reg}_{\alpha})
  \big/V_{-1}(\nbigmhat^{\reg}_{\alpha})
\]
induced by $-u\del_u$
is identified with
the endomorphism
induced by $-\del_z(z-\alpha)$.
Hence, we obtain
$\nbigl_{\alpha}(\nbigm)
=V_{0}\Bigl(
 \bigl(
 \nbigm^{\gbigf}_{\inftyhat,\alpha}
 \bigr)_1
 \Bigr)$,
 which implies the claim of the lemma.
 \hfill\qed

\vspace{.1in}
 
We obtain the following commutative diagrams
for $\alpha\in D$
and $-1<_{\cnum}a\leq_{\cnum}0$:
\[
 \begin{CD}
 \Gr^V_{a}(\nbigm(!D)_{\alpha})
@>>>
 \Gr^V_{a}(\nbigm_{\alpha})
@>>> 
 \Gr^V_a(\nbigm(\ast D)_{\alpha})\\
@V{\simeq}VV @V{\simeq}VV @V{\simeq}VV \\
  \Gr^{V}_{a-1}\bigl(
  \bigl(
 \nbigm(!D)^{\gbigf}_{\inftyhat,\alpha}
 \bigr)_1
 \bigr)
@>>>
 \Gr^{V}_{a-1}\bigl(
  \bigl(
 \nbigm^{\gbigf}_{\inftyhat,\alpha}
 \bigr)_1
 \bigr)
@>>>
 \Gr^{V}_{a-1}\bigl(
  \bigl(
 \nbigm(\ast D)^{\gbigf}_{\inftyhat,\alpha}
 \bigr)_1
 \bigr).
 \end{CD}
\]
If $a\neq 0$,
the horizontal arrows are also isomorphisms.
If $a=0$,
the horizontal morphisms are identified with
$\Gr^V_{-1}(\nbigm_{\alpha})
\lrarr
 \Gr^V_{0}(\nbigm_{\alpha})
\lrarr
 \Gr^V_{-1}(\nbigm_{\alpha})$.
 
\section{Local systems with Stokes structure at finite place}

Let $\nbigm\in\Hol(\proj^1,D,\infty)$.
We obtain the following collection 
$\LS^{\fin}(\nbigm)$
of the data associated with $\nbigm$:
\index{data $\LS^{\fin}(\nbigm)$}
\begin{itemize}
 \item The local system $\nbigl(\nbigm)$ on $\cnum\setminus D$.
 \item $2\pi\seisuu$-equivariant local systems with Stokes structure
       $(L_{\alpha}(\nbigm),\vecnbigf)$ $(\alpha\in D)$.
 \item $\psi_{z-\alpha}(\nbigm)$
       and $\phi_{z-\alpha}(\nbigm)$
       and the morphisms
\[
       \psi_{z-\alpha}(\nbigm)
       \to
       \phi_{z-\alpha}(\nbigm)
       \to
       \psi_{z-\alpha}(\nbigm)
\]
for $\alpha\in D$.
\end{itemize}

Let $\nbigm^{\gbigf}=\Fourier_+(\nbigm)$
be the holonomic $\nbigd$-module
on $\proj^1_w$ as in \S\ref{subsection;18.6.2.120}.
In \S\ref{subsection;18.6.2.120},
we have explained how to compute
$(\gbigl^{\gbigf}(\nbigm),\vecnbigf)
=(L_{\infty}(\nbigm^{\gbigf}),\vecnbigf)$
from
$\LS^{\fin}(\nbigm)$
and $\nbigstilde_1(L_{\infty}(\nbigm),\vecnbigf)$.
Let us complement how to obtain
the rest of $\LS^{\fin}(\nbigm^{\gbigf})$ from
$\LS^{\fin}(\nbigm)$
and 
$(L_{\infty}(\nbigm),\vecnbigf)$.

\subsection{Fourier transform and constructible complexes}
\label{subsection;25.3.5.1}

Let $D\subset \cnum$ be any finite subset.
Let $\nbign\in\Hol(\proj^1_z,D,\infty)$.
We set
$\nbign^{\gbigf_{\pm}}=\Fourier_{\pm}(\nbign)$
which are $\nbigd$-modules on $\proj^1_w$.
Let us study the perverse sheaves
$\DR_{\cnum}\nbign^{\gbigf_{\pm}}=
\nbign^{\gbigf_{\pm}}_{|\cnum}\otimes\Omega^{\bullet}_{\cnum}[1]$
on $\cnum$.

\subsubsection{}

We set $X^{(0)}=\proj^1_z\times\cnum_w$.
We set $H^{(0)}_D=D\times\cnum_w$,
$H^{(0)}_{\infty}=\{\infty\}\times \cnum_w$
and $H^{(0)}=H^{(0)}_D\cup H^{(0)}_{\infty}$.
Let $p_1:X^{(0)}\to \proj^1_z$ and $p_2:X^{(0)}\to\cnum_w$
denote the projections.

Let $U_{\infty}\subset\proj^1_z$ be a neighbourhood of $\infty$.
We set $\nbigu^{(0)}_{\infty}=U_{\infty}\times\cnum_w$.

We set
$\nbige(\pm zw)
=(\nbigo_{X^{(0)}}(\ast H^{(0)}_{\infty}),d\pm d(zw))$
on $(X^{(0)},H^{(0)}_{\infty})$.
We obtain the $\nbigd_{X^{(0)}}(\ast H^{(0)}_{\infty})$-modules
$\nbign^{(0)}_{\pm}=p_1^{\ast}(\nbign)\otimes \nbige(\pm zw)$.
We have
$p^0_{2+}\bigl(
\nbign_{\pm}^{(0)}
\bigr)=\nbign^{\gbigf_{\pm}}_{|\cnum}$.

\subsubsection{}

Note that
$\nbign^{(0)}_{\pm|\nbigu^{(0)}_{\infty}}$
are meromorphic flat bundles on
$(\nbigu^{(0)}_{\infty},H_{\infty})$.
There exists a projective morphism of complex manifolds
$\rho:X^{(1)}\to X^{(0)}$ 
such that
(i) $H^{(1)}_{\infty}=\rho^{-1}(H^{(0)}_{\infty})$ is
a simple normal crossing hypersurface,
(ii) $\rho$ induces an isomorphism
$X^{(1)}\setminus H^{(1)}_{\infty}
\simeq
X^{(0)}\setminus H^{(0)}_{\infty}$,
(iii) $\nbign^{(1)}_{\pm|\nbigu^{(1)}_{\infty}}$
are good meromorphic flat bundles on
$(\nbigu^{(1)}_{\infty},H^{(1)}_{\infty})$,
where
we set
$\nbign_{\pm}^{(1)}=\rho^{\ast}(\nbign_{\pm}^{(0)})$
and
$\nbigu^{(1)}_{\infty}=\rho^{-1}(\nbigu_{\infty}^{(0)})$.
We set
$H^{(1)}_D=\rho^{-1}(H^{(0)}_D)$
and 
$H^{(1)}=H^{(1)}_{\infty}\cup H^{(1)}_D$.

Let $\varpi:\Xtilde^{(1)}(H^{(1)})\to X^{(1)}$
denote the oriented real blow up along $H^{(1)}$.
Let $\nbiga^{\mg}_{\Xtilde^{(1)}(H^{(1)})}$
denote the sheaf of holomorphic functions with moderate growth
on $\Xtilde^{(1)}(H^{(1)})$.
We obtain the following cohomologically constructible complex
on $\Xtilde^{(1)}(H^{(1)})$:
\[
 \DR^{\leq 0}_{\Xtilde^{(1)}(H^{(1)})} (\nbign_{\pm}^{(1)})
 :=\nbiga^{\mg}_{\Xtilde^{(1)}(H^{(1)})}
 \otimes
 _{\varpi^{-1}(\nbigo_X)}
 \varpi^{-1}
 \Bigl(
 \Omega^{\bullet}_{X^{(1)}}\otimes_{\nbigo_{X^{(1)}}}
 \nbign^{(1)}_{\pm}
 \Bigr)[2].
\]
We set $\ptilde_2=p_2\circ\rho\circ\varpi$.
There exists a quasi-isomorphism
\[
 R\ptilde_{2\ast}\Bigl(
 \DR^{\leq 0}_{\Xtilde^{(1)}(H^{(1)})} (\nbign_{\pm}^{(1)})
 \Bigr)
 \simeq
 \DR_{\cnum_w}(\nbign^{\gbigf_{\pm}}).
\]

\subsubsection{}

We set
$\nbigutilde_{\infty}^{(1)}=\varpi^{-1}(\nbigu_{\infty}^{(1)})$.
We obtain the local systems
$\nbigl(\nbign^{(1)}_{\pm})$ on $\nbigutilde_{\infty}^{(1)}$
associated with $\nbign_{\pm}^{(1)}$.
We obtain the constructible subsheaf
$\nbigl(\nbign^{(1)}_{\pm})^{\leq 0}
\subset\nbigl(\nbign_{\pm}^{(1)})$
associated with the Stokes structure of
$\nbign_{\pm|\nbigutilde^{(1)}_{\infty}}$.
There exists the following natural quasi-isomorphism
\[
 \nbigl(\nbign_{\pm}^{(1)})^{\leq 0}[2]
 \lrarr
 \DR^{\leq 0}_{\Xtilde^{(1)}(H^{(1)})} (\nbign_{\pm}^{(1)})
 _{|\nbigutilde^{(1)}_{\infty}}.
\]

\subsection{Canonical morphisms}
\label{subsection;25.3.12.60}

We continue to use the notation in \S\ref{subsection;25.3.5.1}.
We obtain the $2\pi\seisuu$-equivariant local system with Stokes structure
$(L_{\infty}(\nbign),\vecnbigf)$
corresponding to $\nbign$ at $\infty$.
We obtain the meromorphic flat bundle 
$\nbigttilde^{\infty}_1(\nbign)$
on $(\proj^1,\{0,\infty\})$
such that
(i) it is regular singular at $0$,
(ii)
$(L_{\infty}(\nbigttilde^{\infty}_1(\nbign)),\vecnbigf)
=\nbigttilde_1(L_{\infty}(\nbign),\vecnbigf)$.

\begin{prop}
\label{prop;25.3.5.30}
There exist the following natural morphisms
of perverse sheaves:
\begin{equation}
\label{eq;25.3.4.1}
 \DR_{\cnum}\bigl(
 \nbigttilde_1(\nbign)(!0)^{\gbigf_{\pm}}
 \bigr)
 \stackrel{c_1}{\lrarr}
 \DR_{\cnum}(\nbign^{\gbigf_{\pm}})
 \stackrel{c_2}{\lrarr}
 \DR_{\cnum}\bigl(
 \nbigttilde_1(\nbign)(\ast 0)^{\gbigf_{\pm}}
 \bigr).
 \end{equation}
The kernel and the cokernel of
$c_i$ are cohomologically locally constant sheaves.
\end{prop}

We shall prove the proposition in
\S\ref{subsection;25.3.6.1}--\S\ref{subsection;25.3.6.2}.
To simplify the description,
we explain the case $\nbign^{\gbigf_+}$,
and we omit to denote $+$.

\subsubsection{}
\label{subsection;25.3.6.1}

We set $H^{(1)}_0=\rho^{-1}(\{0\}\times\cnum_w)$.
We obtain the good meromorphic flat bundle
\[
\nbigttilde_1^{\infty}(\nbign)^{(1)}
=\rho^{\ast}\Bigl(
 p_1^{\ast}(\nbigttilde_1^{\infty}(\nbign))
 \otimes\nbige(zw)
\Bigr)
\]
on $(X^{(1)},H^{(1)}_0\cup H^{(1)}_{\infty})$.
We have
\[
  \DR_{\cnum}\bigl(
 \nbigttilde_1(\nbign)^{\gbigf}(\star 0)
 \bigr)
 \simeq
 R\ptilde_{2\ast}
 \Bigl(
 \DR^{\leq 0}_{\Xtilde^{(1)}(H^{(1)})}
 \bigl(
 \nbigttilde_1^{\infty}(\nbign)^{(1)}(\star H^{(1)}_0)
 \bigr)
 \Bigr).
\]

\subsubsection{}

We obtain the local system
$\nbigl(\nbigttilde^{\infty}_1(\nbign)^{(1)})$
on $\nbigutilde^{(1)}_{\infty}$.
It is equipped with the Stokes structure.
We obtain the constructible subsheaf
$\nbigl(\nbigttilde^{\infty}_1(\nbign)^{(1)})^{\leq 0}
\subset
\nbigl(\nbigttilde^{\infty}_1(\nbign)^{(1)})$.
There exists a natural quasi-isomorphism
\[
 \nbigl(\nbigttilde_1^{\infty}(\nbign)^{(1)})^{\leq 0}[2]
 \lrarr
 \DR^{\leq 0}_{\Xtilde^{(1)}(H^{(1)})}
 \bigl(\nbigttilde^{\infty}_1(\nbign)^{(1)}\bigr)
 _{|\nbigutilde^{(1)}_{\infty}}.
\]

Let $\iota:\nbigutilde^{(1)}_{\infty}\lrarr
\Xtilde^{(1)}(H^{(1)})$
denote the inclusion.
There exists the following morphism:
\begin{equation}
\label{eq;25.3.5.2}
 \iota_!\Bigl(
 \nbigl(\nbigttilde_1^{\infty}(\nbign)^{(1)})^{\leq 0}[2]
 \Bigr)
 \lrarr
 \DR^{\leq 0}_{\Xtilde^{(1)}(H^{(1)})}
 \bigl(\nbigttilde^{\infty}_1(\nbign)^{(1)}(!H_0)\bigr).
\end{equation}
There also exist the following morphisms:
\begin{multline}
\label{eq;25.3.5.4}
\DR^{\leq 0}_{\Xtilde^{(1)}(H^{(1)})}
  (\nbigttilde^{\infty}_1(\nbign)^{(1)}(\ast H_{\infty}))
  \lrarr
R\iota_{\ast}\Bigl(
 \DR^{\leq 0}_{\Xtilde^{(1)}(H^{(1)})}
  (\nbigttilde^{\infty}_1(\nbign)^{(1)})_{|\nbigutilde^{(1)}_{\infty}}
 \Bigr)
 \\
\llarr
R\iota_{\ast}\Bigl(
  \nbigl(\nbigttilde_1^{\infty}(\nbign)^{(1)})^{\leq 0}
\Bigr)
=\iota_{\ast}\Bigl(
  \nbigl(\nbigttilde_1^{\infty}(\nbign)^{(1)})^{\leq 0}
\Bigr).
\end{multline}
The following lemma is standard and easy to see.
\begin{lem}
The morphism {\rm(\ref{eq;25.3.5.2})}
induces an isomorphism
\[
 R\ptilde_{2\ast}\circ
 \iota_!\Bigl(
 \nbigl(\nbigttilde_1^{\infty}(\nbign)^{(1)})^{\leq 0}[2]
 \Bigr)
\simeq
 \DR_{\cnum}\bigl(
 \nbigttilde^{\infty}_1(\nbign)(!0)^{\gbigf}
 \bigr).
\]
The morphisms {\rm(\ref{eq;25.3.5.4})}
induce an isomorphism
\[
  R\ptilde_{2\ast}\circ
 \iota_{\ast}\Bigl(
 \nbigl(\nbigttilde_1^{\infty}(\nbign)^{(1)})^{\leq 0}[2]
 \Bigr)
\simeq
 \DR_{\cnum}\bigl(
 \nbigttilde^{\infty}_1(\nbign)(\ast 0)^{\gbigf}
 \bigr).
\] 
\hfill\qed
\end{lem}

\subsubsection{}

Let $\varpi_1:\Utilde_{\infty}\to U_{\infty}$
denote the real oriented blow up along $\infty$.
We also set
$U_{\infty}^{\ast}=U_{\infty}\setminus\{\infty\}
=\Utilde_{\infty}\setminus\varpi_1^{-1}(\infty)$.
Let $q_{\Utilde}:\Utilde_{\infty}\to \varpi^{-1}_1(\infty)$
and $q_{U^{\ast}}:U_{\infty}^{\ast}\to \varpi^{-1}_1(\infty)$
denote the projections
with respect to the polar decomposition induced by
the complex coordinate $z^{-1}$.

We obtain the holonomic $\nbigd$-module
$\nbigstilde^{\infty}_1(\nbign)\in\Hol(\proj^1,D,\infty)$.
We set
$\nbigv=\nbigstilde^{\infty}_1(\nbign)_{|U_{\infty}}$,
which is a meromorphic flat bundle on $(U_{\infty},\infty)$.
Let $\nbigl(\nbigv)$ denote the local system on
$\Utilde_{\infty}$ corresponding to $\nbigv$.
It is equipped with the Stokes structure.

We obtain the local system
$L_{1,S^1}=\nbigl(\nbigv)_{|\varpi_1^{-1}(\infty)}$
on $\varpi^{-1}(\infty)$.
We obtain the constructible subsheaves
$L_{1,S^1}^{<0}
\subset
L_{1,S^1}^{\leq 0}
\subset
L_{1,S^1}$.
We identify
$q_{\Utilde}^{-1}(L_{1,S^1})=\nbigl(\nbigv)$.
We obtain the constructible subsheaves
\[
q_{\Utilde}^{-1}(L_{1,S^1}^{<0})
\subset
q_{\Utilde}^{-1}(L_{1,S^1}^{\leq 0})
\subset
\nbigl(\nbigv),
\quad\quad
q_{U^{\ast}}^{-1}(L_{1,S^1}^{<0})
\subset
q_{U^{\ast}}^{-1}(L_{1,S^1}^{\leq 0})
\subset
\nbigl(\nbigv)_{|U^{\ast}_{\infty}}.
\]

We set $\nbigv_0=\nbigttilde^{\infty}_1(\nbign)_{|U_{\infty}}$.
Let $\nbigl(\nbigv_0)$
denote the associated local system on
$\Utilde_{\infty}$.
It is equipped with the Stokes structure.
There exists the natural isomorphism
\[
 \nbigl(\nbigv_0)_{|\Utilde_{\infty}}
 \simeq
 q_{\Utilde}^{-1}(L_{1,S^1}^{\leq 0})
 \big/
 q_{\Utilde}^{-1}(L_{1,S^1}^{< 0}).
\]

\subsubsection{}

Let $f_{\Utilde}:\nbigutilde^{(1)}_{\infty}\to \Utilde_{\infty}$
denote the induced map,
which satisfies
$\varpi_1\circ f_{\Utilde}=p_1\circ\rho\circ\varpi$.
We set
$\nbigv^{(1)}=
(p_1\circ\rho)^{\ast}(\nbigv)$
which is a good meromorphic flat bundle on
$(\nbigu^{(1)}_{\infty},H^{(1)}_{\infty})$.
Let $\nbigl(\nbigv^{(1)})$ denote the associated local system on
$\nbigutilde^{(1)}_{\infty}$
corresponding to $\nbigv^{(1)}$.

We set
$\nbigutilde^{(1)\ast}_{\infty}=
\nbigu^{(1)}_{\infty}\setminus H^{(1)}_{\infty}
=\nbigutilde^{(1)}_{\infty}\setminus
\varpi^{-1}(H^{(1)}_{\infty})$.
The multiplication of $\exp(-wz)$ induces isomorphisms
$\nbigl(\nbigv^{(1)})_{|\nbigutilde^{(1)\ast}_{\infty}}
\simeq\nbigl(\nbign^{(1)})_{|\nbigutilde^{(1)\ast}_{\infty}}$.
It induces
$\nbigl(\nbigv^{(1)})
\simeq\nbigl(\nbign^{(1)})$ on $\nbigutilde^{(1)}_{\infty}$.
Similarly, we obtain
$f_{\Utilde}^{-1}(\nbigl(\nbigv_0))
\simeq
 \nbigl(\nbigttilde^{\infty}_1(\nbign)^{(1)})$.

\subsubsection{}

There exist the constructible subsheaves
\[
 f_{\Utilde}^{-1}(q_{\Utilde}^{-1}(L_{1,S^1}^{<0}))
 \subset
 f_{\Utilde}^{-1}(q_{\Utilde}^{-1}(L_{1,S^1}^{\leq 0}))
 \subset
 \nbigl(\nbign^{(1)}).
\]
Note that
\[
 f_{\Utilde}^{-1}(q_{\Utilde}^{-1}(L_{1,S^1}^{\leq 0}))\Big/
 f_{\Utilde}^{-1}(q_{\Utilde}^{-1}(L_{1,S^1}^{<0}))
 \simeq
 \nbigl(\nbigttilde^{\infty}_1(\nbign)^{(1)})
\]
There exist the following exact sequence:
\begin{multline}
 0\lrarr
 f_{\Utilde}^{-1}(q_{\Utilde}^{-1}(L_{1,S^1}^{<0}))
 _{|\varpi^{-1}(H^{(1)}_{\infty})}
 \lrarr
 \nbigl(\nbign^{(1)})^{\leq 0}_{|\varpi^{-1}(H^{(1)}_{\infty})}
 \\
 \lrarr
 \nbigl(\nbigttilde^{\infty}_1(\nbign)^{(1)})^{\leq 0}
 _{|\varpi^{-1}(H^{(1)}_{\infty})}
 \lrarr 0.
\end{multline}

\subsubsection{}

Let $N_{!}\subset \nbigl(\nbign^{(1)})^{\leq 0}$
denote the constructible subsheaf determined by
the following exact sequence:
\[
 0\lrarr f_{\Utilde}^{-1}(q_{\Utilde}^{-1}(L_{1,S^1}^{<0}))
 \lrarr
 N_{!}
 \lrarr
 \nbigl(\nbigttilde^{\infty}_1(\nbign)^{(1)})^{\leq 0}
 \lrarr 0.
\]

\begin{lem}
\label{lem;25.3.5.11}
$R(p_2\circ\rho\circ\varpi)_{\ast}\Bigl(
 \iota_{!}f_{\Utilde}^{-1}\bigl(q_{\Utilde}^{-1}(L_{1,S^1}^{<0})\bigr)
 \Bigr)=0$.
\end{lem}
\pf
Let $j:\nbigutilde^{(1)\ast}_{\infty}\lrarr \nbigutilde^{(1)}_{\infty}$
denote the inclusion.
Let $f_{U^{\ast}}:
\nbigu^{(1)\ast}_{\infty}\to U_{\infty}^{\ast}$
denote the projection.
Note that
\[
f_{\Utilde}^{-1}(q_{\Utilde}^{-1}(L_{1,S^1}^{<0}))
=j_{\ast}\bigl(
f_{U^{\ast}}^{-1}(q_{U^{\ast}}^{-1}(L_{1,S^1}^{<0}))
\bigr)
=Rj_{\ast}\bigl(
f_{U^{\ast}}^{-1}(q_{U^{\ast}}^{-1}(L_{1,S^1}^{<0}))
\bigr).
\]
Let $\varpi^{(0)}:\Xtilde^{(0)}(H^{(0)})\to X^{(0)}$
denote the oriented real blow up along $H^{(0)}$.
We have the induced map $\rhotilde:\Xtilde^{(1)}(H^{(1)})
\to \Xtilde^{(0)}(H^{(0)})$
which satisfies
$\varpi^{(0)}\circ\rhotilde
=\rho\circ\varpi$.
Let $j^{(0)}:
\nbigu_{\infty}^{(1)\ast}
=U_{\infty}^{\ast}\times\cnum_w
\to
U_{\infty}\times\cnum_w=
\nbigu^{(0)}_{\infty}$
denote the inclusion.
We obtain
\[
 R\rhotilde_{\ast}
 f_{\Utilde}^{-1}(q_{\Utilde}^{-1}(L_{1,S^1}^{<0}))
 =
 Rj^{(0)}_{\ast}\bigl(
 f_{U^{\ast}}^{-1}(q_{U^{\ast}}^{-1}(L_{1,S^1}^{<0}))
 \bigr)
 =j^{(0)}_{\ast}\bigl(
 f_{U^{\ast}}^{-1}(q_{U^{\ast}}^{-1}(L_{1,S^1}^{<0}))
 \bigr).
\]
Let $\iota^{(0)}:\nbigu^{(0)}_{\infty}\to X^{(0)}$
denote the inclusion.
It is easy to see that
\[
 R(p_2\circ\varpi^{(0)})_{\ast}
 \iota^{(0)}_!
j^{(0)}_{\ast}
\bigl(
 f_{U^{\ast}}^{-1}(q_1^{-1}(L_{1,S^1}^{<0}))
 \bigr)=0.
\]
Then, we obtain the claim of Lemma \ref{lem;25.3.5.11}.
\hfill\qed

\vspace{.1in}
By Lemma \ref{lem;25.3.5.11},
the monomorphism 
$\iota_{!}(N_{!})
\to
\iota_!
 (\nbigl^{(1)})^{\leq 0}$
induces an isomorphism
\[
 R(p_2\circ\rho\circ\varpi)_{\ast}
 \Bigl(
 \iota_{!}(N_{!})[2]
 \Bigr)
 \simeq
 \DR_{\cnum}
 (\nbigttilde^{\infty}_1(\nbign)(!0)^{\gbigf}).
\]

There exist the natural morphisms:
\begin{equation}
\label{eq;25.3.5.12}
 \iota_!N_{!}[2]
 \lrarr
 \iota_!\nbigl(\nbign^{(1)})^{\leq 0}[2]
 \lrarr
 \DR^{\leq 0}_{\Xtilde^{(1)}(H^{(1)})}
 (\nbign^{(1)}).
\end{equation}
Thus, we obtain the following morphisms
\[
 \DR_{\cnum}(\nbigttilde^{\infty}_1(\nbign)(!0)^{\gbigf})
 \simeq
 R(p_2\circ\rho\circ\varpi)_{\ast}
 \bigl(
 \iota_!N_{!}[2]
 \bigr)
 \lrarr
 \DR_{\cnum}(\nbign^{\gbigf}).
\]
Thus, we obtain the first morphism $c_1$ in (\ref{eq;25.3.4.1}).

\subsubsection{}
\label{subsection;25.3.5.20}

We put $W_0=\proj^1\setminus U_{\infty}$
and $\nbigw^{(1)}_0=(p_2\circ\rho)^{-1}(W)=W_0\times \cnum_w$.
We naturally regard $\nbigw^{(1)}_0$
as a subset of $\Xtilde^{(1)}(H^{(1)})$.
The inclusion
$\nbigw^{(1)}_0\subset\Xtilde^{(1)}(H^{(1)})$
is denoted by $\iota_W$.
Let $f_W:\nbigw^{(1)}_0\to W_0$ denote the projection.
There exists a natural quasi-isomorphism
\[
 \iota_W^{-1}
  \DR^{\leq 0}_{\Xtilde^{(1)}(H^{(1)})}
  (\nbign^{(1)})_{|\nbigw^{(1)}_0}[-2]
  \simeq
  f_W^{-1}\DR(\nbign)_{|W_0}[-1].
\]
Hence,
$R^j(p_2\circ\rho)_{\ast}
R\iota_{W\ast}
\iota_W^{-1}
\DR^{\leq 0}_{\Xtilde^{(1)}(H^{(1)})}
(\nbign^{(1)})$
are constant sheaves on $\cnum$.

Because
\[
 \nbigl(\nbign^{(1)})^{\leq 0}
 \big/
 N_{!}
 \simeq
 j_!\Bigl(
 f_{U^{\ast}}^{-1}\bigl(
 q_{U^{\ast}}^{-1}(L_{1,S^1}/L_{1,S^1}^{\leq 0})
 \bigr)
 \Bigr),
\]
$R^j(p_2\circ\rho\circ\varpi)_{\ast}
 \iota_!\Bigl(
 \nbigl(\nbign^{(1)})^{\leq 0}
 \big/
 N_{!}
 \Bigr)$
are constant sheaves.
Hence, the kernel and the cokernel of $c_1$
are cohomologically locally constant.

\subsubsection{}

There exists the following exact sequence
of constructible sheaves on
$\nbigutilde^{(1)}_{\infty}$:
\begin{multline}
 0\lrarr
 \nbigl(\nbigttilde^{\infty}_1(\nbign)^{(1)})^{\leq 0}
 \lrarr
 \nbigl(\nbign^{(1)})^{\leq 0}\big/
  f_{\Utilde}^{-1}(q_{\Utilde}^{-1}(L_{1,S^1}^{<0}))
\\
 \lrarr
 j_!\Bigl(
 f_{U^{\ast}}\bigl(q_{U^{\ast}}^{-1}(L_{1,S^1}/L_{1,S^1}^{\leq 0})\bigr)
 \Bigr)
 \lrarr 0.
\end{multline}

\begin{lem}
 $R(p_2\circ\rho\circ\varpi)_{\ast}
 R\iota_{\ast}
 \Bigl(
j_!\Bigl(
 f_{U^{\ast}}\bigl(q_{U^{\ast}}^{-1}(L_{1,S^1}/L_{1,S^1}^{\leq 0})\bigr)
 \Bigr)
 \Bigr)=0$.
\end{lem}
\pf
We have
\[
R\rhotilde_{\ast}
R\iota_{\ast}
\Bigl(
j_!\Bigl(
 f_{U^{\ast}}\bigl(q_{U^{\ast}}^{-1}(L_{1,S^1}/L_{1,S^1}^{\leq 0})\bigr)
 \Bigr)
 =R\iota^{(0)}_{\ast}
 j^{(0)}_!
 \Bigl(
 f_{U^{\ast}}\bigl(q_{U^{\ast}}^{-1}(L_{1,S^1}/L_{1,S^1}^{\leq 0})\bigr)
 \Bigr).
\]
It is easy to see that
$R(p_2\circ\varpi)_{\ast}
R\iota^{(0)}_{\ast}
j^{(0)}_!
 \Bigl(
 f_{U^{\ast}}\bigl(q_{U^{\ast}}^{-1}(L_{1,S^1}/L_{1,S^1}^{\leq 0})\bigr)
 \Bigr)=0$.
\hfill\qed

\vspace{.1in}
There exist the following natural morphisms:
\begin{multline}
\DR^{\leq 0}_{\Xtilde^{(1)}(H^{(1)})}
\bigl(\nbign^{(1)}\bigr) 
\lrarr
 R\iota_{\ast}\Bigl(
 \DR^{\leq 0}_{\Xtilde^{(1)}(H^{(1)})}
\bigl(\nbign^{(1)}\bigr) _{|\nbigutilde^{(1)}_{\infty}}
 \Bigr)
 \\
 \llarr
 \iota_{\ast}
 \nbigl(\nbign^{(1)})^{\leq 0}
 \lrarr
 \iota_{\ast}\Bigl(
 \nbigl(\nbign^{(1)})^{\leq 0}
 \big/f_U^{-1}(q_0^{-1}(L_{1,S^1}^{< 0}))
 \Bigr)
 \llarr
 \iota_{\ast}\bigl(
 \nbigl(\nbigttilde^{\infty}_1(\nbign)^{(1)})
 \bigr).
\end{multline}
They induce the following morphism
\[
 \DR_{\cnum}(\nbign^{\gbigf})
 \lrarr
 \DR_{\cnum}(\nbigttilde^{\infty}_1(\nbign)(\ast 0)^{\gbigf}).
\]
Thus, we obtain the second morphism in (\ref{eq;25.3.4.1}).

\subsubsection{}
\label{subsection;25.3.6.2}

Let $W_0^{\circ}$ denote the interior part of $W_0$.
We set
$\nbigw^{(1)\circ}_0=(p_2\circ\rho)^{-1}(W_0)=W_0\times\cnum$.
Let $f_{W^{\circ}}:\nbigw^{(1)\circ}_0\to W_0$
denote the projection.
We naturally regard $\nbigw^{(1)\circ}_0$
as a subset of $\Xtilde^{(1)}(H^{(1)})$.
Let $\iota_{W^{\circ}}:\nbigw^{(1)\circ}_0\to \Xtilde^{(1)}(H^{(1)})$
denote the inclusion.
As in the case of \S\ref{subsection;25.3.5.20},
$R^j(p_2\circ\rho)_{\ast}
R\iota_{W^{\circ}!}
\iota_{W^{\circ}}^{-1}
\DR^{\leq 0}_{\Xtilde^{(1)}(H^{(1)})}
(\nbign^{(1)})$
are constant sheaves on $\cnum$.
Note that
\[
 R\rhotilde_{\ast}
 R\iota_{\ast}
 \bigl(
 f_{\Utilde}^{-1}(q_{\Utilde}^{-1}(L_1^{<0}))
 \bigr)
=R\iota^{(0)}_{\ast}
 Rj^{(0)}_{\ast}
 \bigl(
 f_{U^{\ast}}^{-1}\bigl(
 q_{U^{\ast}}^{-1}(L_1^{<0})
 \bigr)
 \bigr)
=\iota^{(0)}_{\ast}
 j^{(0)}_{\ast}
 \bigl(
 f_{U^{\ast}}^{-1}\bigl(
 q_{U^{\ast}}^{-1}(L_1^{<0})
 \bigr)
 \bigr).
\]
Hence, 
$R^j(p_2\circ\rho)_{\ast}
R\iota_{\ast}
\bigl(
f_{\Utilde}^{-1}(q_{\Utilde}^{-1}(L_1^{<0}))
\bigr)$
are locally constant sheaves.
We obtain that
the kernel and the cokernel of $c_2$
are cohomologically locally constant.
Thus, the proof of Proposition \ref{prop;25.3.5.30}
is completed.
\hfill\qed

\subsection{The case $\nbigstilde_1(\nbign)=\nbign$}

Let us consider the case $\nbigstilde_1(\nbign)=\nbign$,
i.e.,
$\nbigttilde_1(\nbigi_{\infty}(\nbign))=\{0\}$.
It implies that
$\nbigttilde^{\infty}_1(\nbign)(\star 0)$ are regular singular.
It is easy to see 
\[
 \nbigttilde^{\infty}_1(\nbign)(!0)^{\gbigf_{\pm}}
 \simeq
 \bigl(
 \nbigttilde^{\infty}_1(\nbign)(!0)^{\gbigf_{\pm}}
 \bigr)(\ast 0),
\quad\quad
 \nbigttilde^{\infty}_1(\nbign)(\ast 0)^{\gbigf_{\pm}}
 \simeq
 \bigl(
 \nbigttilde^{\infty}_1(\nbign)(\ast 0)^{\gbigf_{\pm}}
 \bigr)(!0).
\]
The $\nbigd$-modules
$\nbign^{\gbigf_{\pm}}(\ast\{0\})$
and
$\nbigttilde^{\infty}_1(\nbign)(\star 0)^{\gbigf_{\pm}}(\ast\{0\})$
are meromorphic flat bundles on $(\proj^1,\{0,\infty\})$,
and regular singular at $0$.
As explained in \S\ref{subsection;25.3.12.60},
there exist the natural morphisms
\[
 \DR_{\cnum}\nbigttilde^{\infty}_1(\nbign)(!0)^{\gbigf_{\pm}}
 \lrarr
 \DR_{\cnum}\nbign^{\gbigf_{\pm}}
 \lrarr
 \DR_{\cnum}\nbigttilde^{\infty}_1(\nbign)(\ast 0)^{\gbigf_{\pm}}.
\]
The kernel and the cokernel of the morphisms are constant sheaves.

We obtain the induced morphisms of
$2\pi\seisuu$-equivariant local systems
\begin{equation}
\label{eq;25.3.17.1}
L_{\infty}(\nbigttilde^{\infty}_1(\nbign)(! 0)^{\gbigf_{\pm}})
\stackrel{a}{\lrarr}
L_{\infty}(\nbign^{\gbigf_{\pm}})
\stackrel{b}{\lrarr}
L_{\infty}(\nbigttilde^{\infty}_1(\nbign)(\ast 0)^{\gbigf_{\pm}}).
\end{equation}
By using the maps $\Abb^{\rd}_{\nbigttilde^{\infty}_1(\nbign),\pm}$
and $\Abb^{\mg}_{\nbigttilde^{\infty}_1(\nbign),\pm}$,
we obtain the isomorphisms:
\begin{equation}
\label{eq;25.3.5.50}
 L_0(\nbigttilde^{\infty}_1(\nbign)(\star 0))\simeq
 L_{\infty}(\nbigttilde^{\infty}_1(\nbign)(\star 0)^{\gbigf_{\pm}}).
\end{equation}
We rewrite (\ref{eq;25.3.17.1}) as 
\begin{equation}
\label{eq;25.3.15.20}
 L_0(\nbigttilde^{\infty}_1(\nbign)(!0))
\stackrel{a}{\lrarr}
L_{\infty}(\nbign^{\gbigf_{\pm}})
\stackrel{b}{\lrarr}
 L_0(\nbigttilde^{\infty}_1(\nbign)(\ast 0)).
\end{equation}
\begin{rem}
In the case of $+$,
these are the connecting morphisms
in {\rm\S\ref{subsection;18.6.2.120}}.
\hfill\qed
\end{rem}
 Note that
$L_0(\nbigttilde^{\infty}_1(\nbign)(!0))=
L_0(\nbigttilde^{\infty}_1(\nbign)(\ast 0))$.
Let $M_{L_{\infty}(\nbign^{\gbigf_{\pm}})}$
and $M_{L_0(\nbigttilde_1^{\infty}(\nbign))}$
denote the monodromy automorphisms of
$L_{\infty}(\nbign^{\gbigf_{\pm}})$
and $L_0(\nbigttilde^{\infty}_1(\nbign))$, respectively.
We obtain the following from Proposition \ref{prop;25.3.5.51}.

\begin{prop}
$b\circ a=\id-M_{L_0(\nbigttilde^{\infty}_1(\nbign))}^{-1}$
and 
$a\circ b=\id-M_{L_{\infty}(\nbign^{\gbigf_{\pm}})}^{-1}$
under the isomorphisms {\rm(\ref{eq;25.3.5.50})}.
\hfill\qed
\end{prop}

\begin{rem}
We can recover
$\psitilde(\nbign^{\gbigf_{\pm}})\to
 \phitilde(\nbign^{\gbigf_{\pm}})\to
\psitilde(\nbign^{\gbigf_{\pm}})$ 
from {\rm(\ref{eq;25.3.15.20})}
as explained in 
{\rm\S\ref{subsection;25.3.13.12}}.
\hfill\qed
\end{rem}

\subsection{Some isomorphisms}
\label{subsection;25.3.4.3}

Let $\check{D}=\bigl\{
\alpha\in\cnum\,\big|\,
-\alpha u^{-1}\in
\pi_1(\nbigttilde_1(\nbigi_{\infty}(\nbign)))
\bigr\}$.
Let $R>0$ such that
$R>\!>|\alpha|$ for any $\alpha\in \check{D}$.
We set $Y(R)=\{|w|>R\}$.
The restrictions
$\nbign^{\gbigf_{\pm}}_{|Y(R)}$
and
$\nbigttilde^{\infty}_1(\nbign)(\star 0)^{\gbigf_{\pm}}_{|Y(R)}$
are flat bundles.

\begin{lem}
\label{lem;25.3.14.1}
On $Y(R)$,
there exist the following natural commutative diagrams of
the local systems:
{\small
\[
\begin{CD}
\DR_{\cnum}
\nbigttilde^{\infty}_1(\nbign)(! 0)^{\gbigf_{\pm}}_{|Y(R)}
@>>>
\DR_{\cnum}(\nbign^{\gbigf_{\pm}})_{|Y(R)}
@>>>
\DR_{\cnum}
 \nbigttilde^{\infty}_1(\nbign)(\ast 0)^{\gbigf_{\pm}}
 _{|Y(R)}
 \\
 @V{\simeq}VV @V{\simeq}VV @V{\simeq}VV \\
\DR_{\cnum}
 \nbigstilde^{\infty}_1
 \nbigttilde^{\infty}_1
 (\nbign)(! 0)^{\gbigf_{\pm}}_{|Y(R)}
@>>>
 \DR_{\cnum}
 \nbigstilde^{\infty}_1(\nbign)^{\gbigf_{\pm}}_{|Y(R)}
@>>>
\DR_{\cnum}
 \nbigstilde^{\infty}_1
 \nbigttilde^{\infty}_1
 (\nbign)(\ast 0)^{\gbigf_{\pm}}_{|Y(R)}.
 \end{CD}
\]}
\end{lem}
\pf
Because
$\nbigl(\nbigttilde^{\infty}_1(\nbign_{\pm})^{(1)})^{\leq 0}
=\nbigl(\nbigttilde^{\infty}_1
\nbigstilde_1^{\infty}(\nbign_{\pm})^{(1)})^{\leq 0}$
on $(\rho\circ\varpi)^{-1}(U_{\infty}\times Y(R))$,
we obtain the claim of the lemma.
\hfill\qed

\vspace{.1in}

Let $\gbigf^c_{\pm}$ denote
the transforms for constructible complexes
explained in \S\ref{subsection;25.3.14.30} below.
Because the kernels and the cokernels
of the horizontal arrows are constant sheaves,
we obtain the isomorphisms
of the local systems on $\cnum^{\ast}$:
\begin{equation}
\label{eq;25.3.15.40}
 \gbigf^c_{\mp}\bigl(
 \DR_{\cnum}
 \nbigttilde^{\infty}_1(\nbign)(! 0)^{\gbigf_{\pm}}
 \bigr)_{|\cnum^{\ast}}
\simeq
 \gbigf^c_{\mp}\bigl(
 \DR_{\cnum}(\nbign^{\gbigf_{\pm}})
 \bigr)_{|\cnum^{\ast}}
\simeq
 \gbigf^c_{\mp}\bigl(
 \DR_{\cnum}
 \nbigttilde^{\infty}_1(\nbign)(\ast 0)^{\gbigf_{\pm}}
 \bigr)_{|\cnum^{\ast}}.
\end{equation}
By Lemma \ref{lem;25.3.14.15},
we obtain the isomorphisms
\begin{equation}
 \label{eq;25.3.15.41}
 \gbigf^c_{\mp}\bigl(
 \DR_{\cnum}
 \nbigttilde^{\infty}_1
 (\nbign)(! 0)^{\gbigf_{\pm}}
 \bigr)
 \simeq
 \gbigf^c_{\mp}\bigl(
 \iota_{R\ast}
 \iota_R^{-1}
 \DR_{\cnum}
 \nbigttilde^{\infty}_1
 (\nbign)(! 0)^{\gbigf_{\pm}}
 \bigr),
\end{equation}
\begin{equation}
 \label{eq;25.3.15.42}
 \gbigf^c_{\mp}\bigl(
 \DR_{\cnum}
 \nbigttilde^{\infty}_1
 (\nbign)(\ast 0)^{\gbigf_{\pm}}
 \bigr)
 \simeq
 \gbigf^c_{\mp}\bigl(
 \iota_{R!}
 \iota_R^{-1}
 \DR_{\cnum}
 \nbigttilde^{\infty}_1
 (\nbign)(\ast 0)^{\gbigf_{\pm}}
 \bigr).
\end{equation}
Hence, we obtain isomorphisms of local systems on $\cnum^{\ast}$:
\begin{equation}
  \gbigf^c_{\mp}\bigl(
 \iota_{R\ast}
 \iota_R^{-1}
 \DR_{\cnum}
 \nbigttilde^{\infty}_1
 (\nbign)(! 0)^{\gbigf_{\pm}}
 \bigr)_{|\cnum^{\ast}}
 \simeq
  \gbigf^c_{\mp}\bigl(
 \iota_{R!}
 \iota_R^{-1}
 \DR_{\cnum}
 \nbigttilde^{\infty}_1
 (\nbign)(\ast 0)^{\gbigf_{\pm}}
 \bigr)_{|\cnum^{\ast}}.
\end{equation}
Similarly,
we obtain the following isomorphisms of local systems
on $\cnum^{\ast}$:
\begin{multline}
 \gbigf^c_{\mp}\bigl(
 \iota_{R\ast}
 \iota_R^{-1}
 \DR_{\cnum}
 \nbigstilde^{\infty}_1
 \nbigttilde^{\infty}_1
 (\nbign)(! 0)^{\gbigf_{\pm}}
 \bigr)_{|\cnum^{\ast}}
\simeq
 \gbigf^c_{\mp}\bigl(
 \DR_{\cnum}
 \nbigstilde^{\infty}_1
 \nbigttilde^{\infty}_1
 (\nbign)(! 0)^{\gbigf_{\pm}}
 \bigr)
 _{|\cnum^{\ast}}
 \\
\simeq
 \gbigf^c_{\mp}\bigl(
 \DR_{\cnum}
 \nbigstilde^{\infty}_1(\nbign)^{\gbigf_{\pm}}
 \bigr)
 _{|\cnum^{\ast}}
\simeq
 \gbigf^c_{\mp}\bigl(
 \DR_{\cnum}
 \nbigstilde^{\infty}_1
 \nbigttilde^{\infty}_1
 (\nbign)(\ast 0)^{\gbigf_{\pm}}
 \bigr)
 _{|\cnum^{\ast}}
\\
 \simeq
  \gbigf^c_{\mp}\bigl(
 \iota_{R!}
 \iota_R^{-1}
 \DR_{\cnum}
 \nbigstilde^{\infty}_1
 \nbigttilde^{\infty}_1
 (\nbign)(\ast 0)^{\gbigf_{\pm}}
 \bigr)_{|\cnum^{\ast}}.
\end{multline}
By Proposition \ref{prop;25.3.14.100} below,
we obtain the following proposition.
\begin{prop}
\label{prop;25.3.15.50}
The following diagram {\rm(\ref{eq;25.3.14.101})}
is commutative.
{\small
\begin{equation}
\label{eq;25.3.14.101}
 \begin{CD}
 \gbigf^c_{\mp}\bigl(
 \iota_{R\ast}
 \iota_R^{-1}
 \DR_{\cnum}
 \nbigttilde^{\infty}_1
 (\nbign)(! 0)^{\gbigf_{\pm}}
 \bigr)_{|\cnum^{\ast}}
 @>{\simeq}>>
  \gbigf^c_{\mp}\bigl(
 \iota_{R!}
 \iota_R^{-1}
 \DR_{\cnum}
 \nbigttilde^{\infty}_1
 (\nbign)(\ast 0)^{\gbigf_{\pm}}
  \bigr)_{|\cnum^{\ast}}
  \\
@V{\simeq}VV @V{\simeq}VV \\
 \gbigf^c_{\mp}\bigl(
 \iota_{R\ast}
 \iota_R^{-1}
 \DR_{\cnum}
 \nbigstilde^{\infty}_1
 \nbigttilde^{\infty}_1
 (\nbign)(! 0)^{\gbigf_{\pm}}
  \bigr)_{|\cnum^{\ast}}
  @>{\simeq}>>
  \gbigf^c_{\mp}\bigl(
 \iota_{R!}
 \iota_R^{-1}
 \DR_{\cnum}
 \nbigstilde^{\infty}_1
 \nbigttilde^{\infty}_1
 (\nbign)(\ast 0)^{\gbigf_{\pm}}
 \bigr)_{|\cnum^{\ast}}.
 \end{CD}
 \end{equation}}
\hfill\qed
\end{prop}
For $\star=!,\ast$, we set $\star'=\ast,!$.
There exist the following isomorphisms:
\begin{multline}
\label{eq;25.3.15.30}
 \gbigf^c_{\mp}\bigl(
 \iota_{R\star'}
 \iota_R^{-1}
 \DR_{\cnum}
 \nbigstilde^{\infty}_1
 \nbigttilde^{\infty}_1
 (\nbign)(\star 0)^{\gbigf_{\pm}}
  \bigr)
  \simeq
  \gbigf^c_{\mp}\Bigl(
 \DR_{\cnum}
 \Bigl(
 \nbigstilde^{\infty}_1
 \nbigttilde^{\infty}_1
 (\nbign)(\star 0)^{\gbigf_{\pm}}(\star' 0)
 \Bigr)
 \Bigr)
 \\
\simeq
 \gbigf^c_{\mp}\Bigl(
 \DR_{\cnum}
 \Bigl(
 \nbigstilde^{\infty}_1
 \nbigttilde^{\infty}_1
 (\nbign)(\star 0)^{\gbigf_{\pm}}
 \Bigr)
 \Bigr)
 \simeq
 \DR_{\cnum}
 \Bigl(
 \Bigl(
 \nbigstilde^{\infty}_1
 \nbigttilde^{\infty}_1
 (\nbign)(\star 0)^{\gbigf_{\pm}}
 \Bigr)^{\gbigf_{\mp}}
 \Bigr)
\\
 \simeq
 \DR_{\cnum}\Bigl(
 \nbigstilde^{\infty}_1
 \nbigttilde^{\infty}_1
 (\nbign)(\star 0)
 \Bigr).
\end{multline}

In summary,
we obtain the following isomorphisms
\begin{equation}
\label{eq;25.3.15.51}
L_{\infty}\bigl(
\nbigs^{\infty}_{\infty}
(\nbign^{\gbigf_{\pm}})^{\gbigf_{\mp}}
\bigr)
\simeq
L_{\infty}\Bigl(
\bigl(
\nbigstilde^{\infty}_1
\nbigttilde^{\infty}_1
\nbign(\star 0)\bigr)^{\gbigf_{\pm}}
(\star'0)^{\gbigf_{\mp}}
\Bigr)
\simeq
L_{\infty}\bigl(
\nbigstilde^{\infty}_1\nbigttilde^{\infty}_1\nbign
\bigr).
\end{equation}
The compositions of the morphisms
are independent of the choice of
$(\star,\star')=(!,\ast),(\ast,!)$
by Proposition \ref{prop;25.3.15.50}.

\subsection{Some commutative diagrams}
\label{subsection;25.3.18.101}

We obtain $\nbigs_{\infty}^{\infty}\nbigm
\in \Hol(\proj^1,D,\infty)$,
i.e.,
$\nbigs_{\infty}^{\infty}\nbigm$
is characterized by the conditions that
(i) $\nbigs^{\infty}_{\infty}(\nbigm)_{|\cnum}=\nbigm_{|\cnum}$,
(ii) $\nbigs^{\infty}_{\infty}(\nbigm)$ is regular at $\infty$.
Let $V_{\infty}$ be a regular singular meromorphic flat bundle
on $(\proj^1,\{0,\infty\})$
with an isomorphism
$V_{\infty}\simeq \nbigs^{\infty}_{\infty}(\nbigm)$
on a neighbourhood of $\infty$.
There exist the following natural morphisms
\[
\DR_{\cnum}\bigl(
 V_{\infty}(!0)^{\gbigf_+}
\bigr)
\lrarr
\DR_{\cnum}\bigl(
\nbigs^{\infty}_{\infty}(\nbigm)^{\gbigf_+}\bigr)
\lrarr
\DR_{\cnum}\bigl(
 V_{\infty}(\ast 0)^{\gbigf_+}
\bigr).
\]
There exist the isomorphisms
induced by $\Abb^{\rd}_+$ and $\Abb^{\mg}_+$:
\[
 c^{-1}L_{\infty}(\nbigm)=
 c^{-1}L_{\infty}(V_{\infty}(\star 0))\simeq
 L_{\infty}\bigl(V_{\infty}(\star 0)^{\gbigf_+}\bigr).
\]
We obtain the following morphisms
\begin{equation}
\label{eq;25.3.15.2}
 c^{-1}L_{\infty}(\nbigm)
 \stackrel{a_1}{\lrarr}
 L_{\infty}
 \bigl(
 \nbigs^{\infty}_{\infty}(\nbigm)^{\gbigf_+}
 \bigr)
 \stackrel{a_2}{\lrarr}
 c^{-1}L_{\infty}(\nbigm).
\end{equation}

\subsubsection{}

We set
$\nbign=\nbigm^{\gbigf_+}$.
We obtain the following morphisms:
\begin{equation}
\label{eq;25.3.15.10}
\DR_{\cnum}\bigl(
\nbigttilde^{\infty}_1\nbigstilde^{\infty}_1(\nbign(!0))^{\gbigf_-}
\bigr)
\lrarr
\DR_{\cnum}\bigl(
 (\nbigstilde^{\infty}_1\nbign)^{\gbigf_-}
 \bigr)
 \lrarr
  \DR_{\cnum}\bigl(
 \nbigttilde^{\infty}_1\nbigstilde^{\infty}_1(\nbign(\ast 0))^{\gbigf_-}
 \bigr).
\end{equation}
By using $\Abb^{\rd}_-$ and $\Abb^{\mg}_-$,
we obtain the following isomorphisms:
\[
 L_{\infty}\bigl(
 \nbigttilde^{\infty}_1
 \nbigstilde^{\infty}_1
 (\nbign)
 \bigr)
= L_{\infty}\bigl(
 \nbigttilde^{\infty}_1
 \nbigstilde^{\infty}_1
 (\nbign)(\star 0)
 \bigr)
 \simeq
 c^{-1}L_{\infty}\bigl(
\nbigttilde^{\infty}_1
 \nbigstilde^{\infty}_1
 (\nbign)(\star 0)^{\gbigf_-} 
 \bigr).
\]
We obtain the following morphisms
\[
 L_{\infty}\bigl(
 \nbigttilde^{\infty}_1
 \nbigstilde^{\infty}_1
 (\nbign)
 \bigr)
\stackrel{b_1}{\lrarr}
c^{-1}L_{\infty}\bigl(
  (\nbigstilde^{\infty}_1\nbign)^{\gbigf_-}
\bigr)
\stackrel{b_2}{\lrarr}
  L_{\infty}\bigl(
 \nbigttilde^{\infty}_1
 \nbigstilde^{\infty}_1
 (\nbign)
 \bigr).
\]

\subsubsection{}
\label{subsection;25.3.18.30}

Because $\nbigm=\nbign^{\gbigf_-}$,
by Lemma \ref{lem;25.3.14.1},
there exists $R>0$ such that
\[
 \DR_{\cnum}\nbigm_{|Y(R)}
 \simeq
 \DR_{\cnum}(\nbign^{\gbigf_-})_{|Y(R)}
 \simeq
 \DR_{\cnum}(\nbigstilde^{\infty}_1(\nbign)^{\gbigf_-})_{|Y(R)}.
\]
In particular, we obtain the following isomorphism:
\[
 d_1:
L_{\infty}\bigl(
 \nbigstilde^{\infty}_1(\nbign)^{\gbigf_-}
 \bigr)
 \stackrel{\simeq}{\lrarr}
 L_{\infty}(\nbigm).
\]
From (\ref{eq;25.3.15.51}),
we obtain the following isomorphism:
\[
 d_2:
L_{\infty}(\nbigstilde^{\infty}_1\nbigttilde^{\infty}_1\nbign)
\stackrel{\simeq}{\lrarr}
 L_{\infty}(\nbigs^{\infty}_{\infty}(\nbigm)^{\gbigf_+}).
\]
Let $M$ denote the monodromy automorphism of
$L_{\infty}(\nbigttilde^{\infty}_1\nbigstilde^{\infty}_1\nbign)$.

\begin{prop}
\label{prop;25.3.15.100}
The following diagrams are commutative:
{\small
\begin{equation}
\begin{CD}
 c^{-1}
 L_{\infty}\bigl(
 (\nbigstilde^{\infty}_1\nbign)^{\gbigf_-}
 \bigr)
 @>{(2\pi\sqrt{-1})^{-1}M\cdot b_2}>>
 L_{\infty}
 (\nbigttilde^{\infty}_1\nbigstilde^{\infty}_1(\nbign))
 @>{-(2\pi\sqrt{-1})b_1}>>
 c^{-1}
 L_{\infty}\bigl(
 (\nbigstilde^{\infty}_1\nbign)^{\gbigf_-}
 \bigr)
 \\
 @V{d_1}VV @V{d_2}VV
 @V{d_1}VV
 \\
 c^{-1}L_{\infty}(\nbigm)
 @>{a_1}>>
 L_{\infty}\bigl(
 \nbigs^{\infty}_{\infty}(\nbigm)^{\gbigf_+}
 \bigr)
 @>{a_2}>>
 c^{-1}L_{\infty}(\nbigm).
\end{CD}
\end{equation}}
\end{prop}
\pf
Let $\varpi:\projtilde^1\to\proj^1$ denote
the oriented real blow up along $\infty$.
We use the polar coordinate $(r,e^{\sqrt{-1}\theta})$
around $\varpi^{-1}(\infty)$
induced by the polar decomposition
$z^{-1}=|z^{-1}|e^{\sqrt{-1}\theta}$.
Let $\varphi_{\infty}:\real\to\varpi^{-1}(\infty)$
be defined by
$\varphi_{\theta}(\theta)=(0,e^{\sqrt{-1}\theta})$.
Let $j_{\infty}:\cnum\to \proj^1$ denote the inclusion.
Let $\nbigp$ be a perverse sheaf on $\cnum$
such that $\nbigp_{|\cnum\setminus D'}=\nbigl_{\cnum\setminus D'}[1]$
for a finite subset $D'\subset\cnum$
and a local system $\nbigl_{\cnum\setminus D'}$
on $\cnum\setminus D'$.
For such $\nbigp$,
we obtain a local system
$j_{\infty\ast}\nbigl_{\cnum\setminus D'}$ on $\projtilde^1\setminus D'$.
We set
$\nbigp_{\infty}=
H^0(\real,\varphi_{\infty}j_{\infty\ast}\nbigl_{\cnum\setminus D'})$.
It is equipped with the induced $2\pi\seisuu$-action.

There exists the following commutative diagram
of vector spaces with a $2\pi\seisuu$-action:
{\tiny
\begin{equation}
\begin{CD}
 \gbigf^c_+
 \bigl(
 \bigl(\nbigttilde^{\infty}_1\nbign(!0)\bigr)^{\gbigf_-}
 \bigr)_{\infty}
 @>{\simeq}>>
 \gbigf^c_+
 \bigl(
 \nbign^{\gbigf_-}
 \bigr)_{\infty}
 @>{\simeq}>>
 \gbigf^c_+
 (\DR\nbigm)_{\infty}\\
@V{\simeq}VV @VVV @VVV \\
 \gbigf^c_+\bigl(
 \iota_{R\ast}\iota_R^{-1}
 \DR\bigl(
 (\nbigttilde^{\infty}_1\nbign(!0))^{\gbigf_-}
 \bigr)
 \bigr)_{\infty}
 @>>>
 \gbigf^c_+\bigl(
 \iota_{R\ast}\iota_R^{-1}
 \DR\nbign^{\gbigf_-}
 \bigr)_{\infty}
 @>{\simeq}>>
 \gbigf^c_+\bigl(
 \iota_{R\ast}\iota_R^{-1}
 \DR\nbigm
 \bigr)_{\infty}\\
@V{\simeq}VV @V{\simeq}VV @V{=}VV \\ 
\gbigf^c_+\bigl(
 \iota_{R\ast}\iota_R^{-1}
\DR(\nbigttilde^{\infty}_1\nbigstilde^{\infty}_1\nbign)(!0)^{\gbigf_-}
 \bigr)_{\infty}
 @>>>
 \gbigf^c_+
 \bigl(
 \iota_{R\ast}\iota_R^{-1}
 \DR\bigl(
  (\nbigstilde^1_{\infty}\nbign)^{\gbigf_-}
 \bigr)
 \bigr)_{\infty}
 @>{\simeq}>>
  \gbigf^c_+\bigl(
 \iota_{R\ast}\iota_R^{-1}
 \DR\nbigm
 \bigr)_{\infty}\\
 @A{\simeq}A{\Abb^{\mg}_+\circ c^{-1}}A
 @A{\simeq}A{\Abb^{\mg}_+\circ c^{-1}}A
 @A{\simeq}A{\Abb^{\mg}_+\circ c^{-1}}A\\
 \DR
 \bigl(
 \bigl(\nbigttilde^{\infty}_1
 \nbigstilde^{\infty}_1(\nbign)(!0)
 \bigr)^{\gbigf_-}
 \bigr)_{\infty}
 @>>>
 \DR\bigl(
 (\nbigstilde_1^{\infty}\nbign)^{\gbigf_-}
 \bigr)_{\infty}
 @>{\simeq}>>
 \DR\nbigm_{\infty}\\
 @A{\simeq}A{\Abb^{\rd}_-\circ c^{-1}}A \\
 \DR(\nbigstilde^{\infty}_1\nbigttilde^{\infty}_1\nbign(!0))_{\infty}
\end{CD}
\end{equation}
}
Note that
\[
\gbigf^c_{+}\DR(\nbigm)_{\infty}
=H^0(\real,L_{\infty}(\nbigs^{\infty}_{\infty}\nbigm)^{\gbigf_+}),
\quad
 \DR\nbigm_{\infty}
 =H^0(\real,L_{\infty}(\nbigm)),
\]
We also have
 \[
DR(\nbigstilde^{\infty}_1\nbigttilde^{\infty}_1\nbign(\star 0))_{\infty}
=H^0(\real,L_{\infty}(\nbigttilde^{\infty}_1\nbigstilde^{\infty}_1\nbign)),
\]
\[
  \DR\bigl(
 (\nbigstilde_1^{\infty}\nbign)^{\gbigf_-}
 \bigr)_{\infty}
 =H^0(\real,L_{\infty}(
 (\nbigstilde_1^{\infty}\nbign)^{\gbigf_-}
 )).
\]
For $(\star,\star')=(!,\ast),(\ast,!)$,
we recall (\ref{eq;25.3.15.51}):
\begin{multline}
\gbigf^c_+\bigl(
 \iota_{R\star'}\iota_R^{-1}
\DR(\nbigttilde^{\infty}_1\nbigstilde^{\infty}_1\nbign)(\star 0)^{\gbigf_-}
 \bigr)_{\infty}
 =
 \DR\Bigl(
 \bigl(\nbigttilde^{\infty}_1\nbigstilde^{\infty}_1\nbign)(\star 0)^{\gbigf_-}
 \bigr)(\star' 0)^{\gbigf_+}
 \Bigr)_{\infty}
\\
= \DR\bigl(
 (\nbigttilde^{\infty}_1\nbigstilde^{\infty}_1\nbign)(\star 0)
 \bigr)_{\infty}
=H^0(\real,L_{\infty}(\nbigttilde^{\infty}_1\nbigstilde^{\infty}_1\nbign)).
\end{multline}
By a diagram chasing and Corollary \ref{cor;25.3.15.1},
we obtain
\[
 d_1\circ b_1
=a_2\circ d_2\circ
 \bigl(
 \Abb^{\mg}_+\circ c^{-1}
 \bigr)
 \circ
 \bigl(
 \Abb^{\rd}_-\circ c^{-1}
 \bigr)
=(-2\pi\sqrt{-1})^{-1}a_2\circ d_2.
\]
Similarly, there exists the following commutative diagram:
{\tiny
\begin{equation}
\begin{CD}
 @. @.
\DR(\nbigstilde^{\infty}_1\nbigttilde^{\infty}_1\nbign(\ast 0))_{\infty}
\\
 @. @.
 @V{\Abb^{\mg}_-\circ c^{-1}}V{\simeq}V
\\
\DR(\nbigm)_{\infty}
@>{\simeq}>>
\DR
\bigl(
\bigl(
 \nbigstilde^{\infty}_1\nbign
 \bigr)^{\gbigf_-}
\bigr)_{\infty}
@>>>
\DR\bigl(
(\nbigstilde^{\infty}_1\nbigttilde^{\infty}_1\nbign(\ast 0))^{\gbigf_-}
\bigr)_{\infty}
\\
@V{\Abb^{\rd}_+\circ c^{-1}}V{\simeq}V
@V{\Abb^{\rd}_+\circ c^{-1}}V{\simeq}V
@V{\Abb^{\rd}_+\circ c^{-1}}V{\simeq}V
\\
\gbigf^c_+\bigl(
\iota_{R!}\iota_R^{-1}
\DR\nbigm
\bigr)_{\infty}
@>{\simeq}>>
 \gbigf^c_+
\iota_{R!}\iota_R^{-1}
 \DR
 \bigl(
 \bigl(
 \nbigstilde^{\infty}_1\nbign
 \bigr)^{\gbigf_-}
\bigr)_{\infty}
@>>>
 \gbigf^c_+
\iota_{R!}\iota_R^{-1}
 \DR\bigl(
(\nbigstilde^{\infty}_1\nbigttilde^{\infty}_1\nbign(\ast 0))^{\gbigf_-}
 \bigr)_{\infty}
 \\
 @V{\simeq}VV @V{\simeq}VV @V{\simeq}VV \\
\gbigf^c_+\bigl(
\iota_{R!}\iota_R^{-1}
\DR\nbigm
\bigr)_{\infty}
 @>{\simeq}>>
 \gbigf^c_+\bigl(
 \iota_{R!}\iota_R^{-1}
 \DR\nbign^{\gbigf_-}
 \bigr)_{\infty}
 @>>>
 \gbigf^c_+
\iota_{R!}\iota_R^{-1}
 \DR\bigl(
(\nbigttilde^{\infty}_1\nbign(\ast 0))^{\gbigf_-}
 \bigr)_{\infty}\\
 @VVV @VVV @V{\simeq}VV \\
\gbigf^c_+\bigl(
\DR\nbigm
 \bigr)_{\infty}
 @>{\simeq}>>
 \gbigf^c_+\bigl(
 \DR\nbign^{\gbigf_-}
 \bigr)_{\infty}
 @>{\simeq}>>
  \gbigf^c_+
 \DR\bigl(
(\nbigttilde^{\infty}_1\nbign(\ast 0))^{\gbigf_-}
 \bigr)_{\infty}
\end{CD}
\end{equation}
}
By a diagram chasing and Corollary \ref{cor;25.3.15.1},
we obtain
\[
a_1\circ d_1
=
d_2\circ
\bigl(
\Abb^{\rd}_+\circ c^{-1}
\bigr)
\circ
\bigl(
\Abb^{\mg}_-\circ c^{-1}
\bigr)
\circ b_2
=d_2\circ\bigl(
(2\pi\sqrt{-1})^{-1}M
\bigr)\circ b_2
\]
Thus, we obtain Proposition \ref{prop;25.3.15.100}.
\hfill\qed

\begin{rem}
By Proposition {\rm\ref{prop;25.3.15.100}},
with the results in 
{\rm\S\ref{subsection;25.3.17.2}}
and {\rm\S\ref{subsection;25.3.17.3}},
we can recover
$\LS^{\fin}(\nbigm)$
from 
the $2\pi\seisuu$-equivariant local system
with Stokes structure
 $\nbigttilde_1\bigl(
 L_{\infty}(\nbign),\vecnbigf
\bigr)$,
morphisms
 {\small
 \[
\begin{CD}
 c^{-1}
 L_{\infty}\bigl(
 (\nbigstilde^{\infty}_1\nbign)^{\gbigf_-}
 \bigr)
 @>{(2\pi\sqrt{-1})^{-1}M\cdot b_2}>>
 L_{\infty}
 (\nbigttilde^{\infty}_1\nbigstilde^{\infty}_1(\nbign))
 @>{-(2\pi\sqrt{-1})b_1}>>
 c^{-1}
 L_{\infty}\bigl(
 (\nbigstilde^{\infty}_1\nbign)^{\gbigf_-}
 \bigr)
\end{CD}.
\]}
and an isomorphism
$L_{\infty}(\nbigm)
\simeq
L_{\infty}\bigl(
 (\nbigstilde^{\infty}_1\nbign)^{\gbigf_-}
 \bigr)$.
\hfill\qed 
\end{rem}

\subsection{Computation of $\LS^{\fin}(\Fourier_+\nbigm)$}
\label{subsection;25.3.18.100}

Let $\nbigm\in\Hol(\proj^1,D,\infty)$.
We obtain
$\pi_1(\nbigttilde_1(\nbigi_{\infty}(\nbigm)))
\subset z^{-1}\cnum$.
We set
$D'=\{\alpha\in\cnum\,|\,-\alpha z\in
\pi_1(\nbigttilde_1(\nbigi_{\infty}(\nbigm)))
\}$.
We have
$\nbigm^{\gbigf_+}\in\Hol(\proj^1,D',\infty)$
and
$\nbigm=(\nbigm^{\gbigf_+})^{\gbigf_-}$.
Let $h:\proj^1\to \proj^1$ be defined by $h(z)=-z$.
We have
$(\nbigm^{\gbigf_+})^{\gbigf_+}=h^{\ast}\nbigm$.

There exists an isomorphism
of $2\pi\seisuu$-equivariant local systems
\[
 L_{\infty}(\nbigm^{\gbigf_+})
 \simeq
 L_{\infty}\bigl(
 \nbigstilde_1^{\infty}(h^{\ast}\nbigm)^{\gbigf_-}
 \bigr).
\]
Moreover,
the tuple
$\LS^{\fin}(\nbigm^{\gbigf_+})$
is computed from 
the $2\pi\seisuu$-equivariant local system
$\nbigttilde_1\bigl(
 L_{\infty}(h^{\ast}\nbigm),\vecnbigf
 \bigr)$
and
{\footnotesize
\begin{equation}
\label{eq;25.3.15.110}
\begin{CD}
 c^{-1}L_{\infty}
 \bigl(
 (\nbigstilde^{\infty}_1h^{\ast}\nbigm)^{\gbigf_-}
 \bigr)
 @>{(2\pi\sqrt{-1})^{-1}M\cdot b_2}>>
 L_{\infty}\bigl(
 \nbigttilde^{\infty}_1
 (h^{\ast}\nbigm)
 \bigr)
 @>{-(2\pi\sqrt{-1})b_1}>>
c^{-1}L_{\infty}
 \bigl(
 (\nbigstilde^{\infty}_1h^{\ast}\nbigm)^{\gbigf_-}
\bigr).
 \end{CD}
\end{equation}}
We can compute
$(L_{\infty}(h^{\ast}\nbigm),\vecnbigf)$
and (\ref{eq;25.3.15.110})
from
$(L_{\infty}(\nbigm),\vecnbigf)$
and $\LS^{\fin}(\nbigm)$.
Hence, we can also compute
$\LS^{\fin}(\nbigm^{\gbigf})$
from
$(L_{\infty}(\nbigm),\vecnbigf)$
and $\LS^{\fin}(\nbigm)$.

\subsection{Appendix: Fourier transforms for constructible sheaves}
\label{subsection;25.3.14.30}

We recall Fourier transforms for
cohomologically constructible complexes.
Let $\Kbb$ be any field
though we are interested in only the case $\Kbb=\cnum$.
On a topological space $Y$,
let $\Kbb_Y$ denote the sheaf of
$\Kbb$-valued locally constant functions.
Let $D^b(\Kbb_Y)$ denote the derived category of
cohomologically bounded $\Kbb_Y$-complexes.

Let $\varpi_1:\projtilde^1_z\to \proj^1_z$
denote the oriented real blow up along $\infty$.
Let $j_{\infty}:\cnum\to \projtilde^1$ denote the inclusion.
Let $A^{\bullet}$ be a complex of $\Kbb_{\cnum}$-modules
such that
(i) the cohomology sheaves $\nbigh^j(A^{\bullet})$ are
$\real$-constructible,
(ii) $\nbigh^{j}(A^{\bullet})=\nbigh^{-j}(A^{\bullet})=0$
for any sufficiently large $j$,
(iii) there exists a finite subset $D$
such that $\nbigh^j(A^{\bullet})_{|\cnum\setminus D}$
are locally constant sheaves.

Let $\rho_Z:Z\to \proj^1_z\times\cnum_w$
denote the complex blow up at the point $(\infty,0)$.
We set $H_Z=\rho_Z^{-1}(\{\infty\}\times\cnum_w)$.
Let $\varpi_Z:\Ztilde\to Z$ denote the oriented real blow up
of $Z$ along $H_Z$.
Let $W_{\pm}\subset \varpi_Z^{-1}(H_Z)$
denote the open subset determined by the following condition:
\begin{itemize}
 \item $P\in W$ if and only if
       $(\varpi_Z\circ\rho_Z)^{-1}\exp(\mp zw)$
       is bounded around $P$.
\end{itemize}
We obtain the open subsets
$\Ztilde^{\circ}_{\pm}=(\cnum_z\times\cnum_w)\cup W_{\pm}
\subset
\Ztilde$.
Let $j_{\Ztilde^{\circ}_{\pm},1}:
\Ztilde^{\circ}_{\pm}\to \Ztilde$
denote the inclusions.
Let 
$j_{\Ztilde^{\circ}_{\pm},2}:
\cnum_z\times\cnum_w
\to 
\Ztilde^{\circ}_{\pm}$
denote the inclusion.

Let $p_1:\cnum_z\times\cnum_w\to \cnum_z$
denote the projection.
Let $\ptilde_2:\Ztilde\to \cnum_w$
denote the composition of
$\rho_{Z}\circ\varpi_Z$
and the projection $\proj^1_z\times \cnum_w\to \cnum_w$.
We obtain
\[
 \gbigf^{c}_{\pm}(A^{\bullet})
 :=
 R\ptilde_{2\ast}\Bigl(
 Rj_{\Ztilde^{\circ}_{\pm},1!}
 Rj_{\Ztilde^{\circ}_{\pm},2\ast}
 p_1^{-1}(A^{\bullet})
 \Bigr)
 \in
 D^b(\Kbb_{\cnum_w}).
\]
It is cohomologically $\real$-constructible.

\subsubsection{Holonomic $\nbigd$-modules}

In the following, we set $\Kbb=\cnum$.

\begin{lem}
Let $M\in \Hol(\proj^1,D,\infty)$
such that $M$ is regular at $\infty$.
Then,  
$\gbigf^c_{\pm}\DR_{\cnum_z}(M)=
\DR_{\cnum_w}\Fourier_{\pm}(M)$.
\hfill\qed
\end{lem}

For any $R\geq 0$,
let $Y(R)=\{|z|>R\}$.
Let $\iota_R:Y(R)\to \cnum$ denote the inclusion.

\begin{lem}
For any regular meromorphic flat bundle $V$
on $(\proj^1,\{0,\infty\})$,
$\gbigf^c_{\pm}
\iota_{R\star}\DR_{\cnum}(V)_{|Y(R)}$
are naturally isomorphic to
$\DR_{\cnum}\Fourier_{\pm}(V(\star 0))$. 
\hfil\qed
\end{lem}

\subsubsection{Some meromorphic flat bundles}
\label{subsection;25.3.14.13}

Let $V$ be a meromorphic flat bundle on $(\proj^1,\{0,\infty\})$
such that
(i) $V$ is regular singular at $0$,
(ii) $\nbigttilde_1(\nbigi_{\infty}(V))=\nbigi_{\infty}(V)$.
We set $V(\star 0)^{\gbigf_{\pm}}=\Fourier_{\pm}(V(\star 0))$,
which is regular singular at $\infty$.

As in Lemma \ref{lem;25.3.14.1},
there exists $R>0$ such that
$V(\star 0)^{\gbigf_{\pm}}_{|Y(R)}$
are flat bundles on $Y(R)$.
Let $\nbigl_{R,\pm,\star}$ denote the local systems
on $Y(R)$
obtained as the sheaves of flat sections of
$V(\star 0)^{\gbigf_{\pm}}_{|Y(R)}$.

There exist the following natural morphisms:
\[
\begin{CD}
 \iota_{R!}\nbigl_{R,\pm,\star}[1]
@>{a_{\pm,\star,1}}>>
 \DR_{\cnum}
 V(\star 0)^{\gbigf_{\pm}}
@>{a_{\pm,\star,2}}>>
 \iota_{R\ast}\nbigl_{R,\pm,\star}[1].
\end{CD}
\]
We obtain the following lemma
as a special case of Proposition 
\ref{prop;25.3.5.30}.
\begin{lem}
\label{lem;25.3.14.10}
The kernel and the cokernel 
of the induced morphisms
\[
\begin{CD}
 \gbigf^c_{\mp}\bigl(
 \iota_{R!}\nbigl_{R,\pm,\star}[1]
 \bigr)
@>{\gbigf^{c}_{\mp}(a_{\pm,\star,1})}>>
 \gbigf^{c}_{\mp}\DR_{\cnum}
 V(\star 0)^{\gbigf_{\pm}}
@>{\gbigf^c_{\mp}(a_{\pm,\star,2})}>>
 \gbigf^c_{\mp}\bigl(
 \iota_{R\ast}\nbigl_{R,\pm,\star}[1]
 \bigr)
\end{CD}
\]
are constant sheaves.
\hfill\qed
\end{lem}

\begin{lem}
\label{lem;25.3.14.15}
The following induced morphisms
\begin{equation}
\label{eq;25.3.14.11}
\gbigf^{c}_{\mp}(a_{\pm,\ast,1}):
 \gbigf^c_{\mp}
 \bigl(
 \iota_{R!}\nbigl_{R,\pm,\ast}[1]
 \bigr)
 \lrarr
 \gbigf^c_{\mp}
 \bigl(
 \DR_{\cnum}
 V(\ast 0)^{\gbigf_{\pm}}
 \bigr),
\end{equation}
\begin{equation}
\label{eq;25.3.14.12}
 \gbigf^c_{\mp}(a_{\pm,!,2}):
 \gbigf^c_{\mp}
 \bigl(
 \DR_{\cnum}
 V(!0)^{\gbigf_{\pm}}
 \bigr)
 \lrarr
 \gbigf^c_{\mp}
 \bigl(
 \iota_{R\ast}\nbigl_{R,\pm,!}[1]
 \bigr)
\end{equation}
are isomorphisms. 
\end{lem}
\pf
Because $V(\star 0)^{\gbigf_{\pm}}$
are regular singular at $\infty$,
we have
\[
 \gbigf^{c}_{\mp}
 \DR_{\cnum}(V(\star 0)^{\gbigf_{\pm}})
 \simeq
 \DR_{\cnum}\bigl(
 \Fourier_{\mp}(V(\star 0)^{\gbigf}_{\pm})
 \bigr)
 \simeq
 \DR_{\cnum}(V(\star 0)).
\]
There exist regular singular meromorphic flat bundles
$V_{1,\pm,\star}$ on $(\proj^1,\{0,\infty\})$
such that $\nbigl_{R,\pm,\star}$
are the sheaves of flat sections of
$V_{1,\pm,\star}$ on $Y(R)$.
We have
\[
 \gbigf^c_{\mp}
 \bigl(
 \iota_{R,\star'}\nbigl_{R,\pm,\star}[1]
 \bigr)
 \simeq
 \DR_{\cnum}
 \Fourier_{\mp}
 (V_{1,\pm,\star}(\star'0)).
\]
If $(\star,\star')=(\ast,!)$ or
$(\star,\star')=(!,\ast)$,
then we have
$\Fourier_{\mp}
(V_{1,\pm,\star}(\star'0))(\star 0)
=\Fourier_{\mp}
(V_{1,\pm,\star}(\star'0))$.
By Lemma \ref{lem;25.3.14.10},
the morphisms
(\ref{eq;25.3.14.11})
and
(\ref{eq;25.3.14.12})
induces isomorphisms of the vanishing cycle sheaves.
Then, the claim of the lemma follows.
\hfill\qed

\vspace{.1in}
By Lemma \ref{lem;25.3.14.15},
we obtain the following isomorphisms
of the local systems on $\cnum^{\ast}$:
\begin{multline}
 \gbigf^c_{\mp}
 \bigl(
 \iota_{R\ast}\nbigl_{R,\pm,!}[1]
 \bigr)_{|\cnum^{\ast}}
 \stackrel{\simeq}{\llarr}
 \gbigf^c_{\mp}
 \bigl(
 \DR_{\cnum}
 V(!0)^{\gbigf_{\pm}}
 \bigr)_{|\cnum^{\ast}}
 \stackrel{\simeq}{\lrarr}
\\
 \gbigf^c_{\mp}
 \bigl(
 \DR_{\cnum}
 V(\ast 0)^{\gbigf_{\pm}}
 \bigr)_{|\cnum^{\ast}}
 \stackrel{\simeq}{\llarr}
 \gbigf^c_{\mp}
 \bigl(
 \iota_{R!}\nbigl_{R,\pm,\ast}[1]
 \bigr)_{|\cnum^{\ast}}.
\end{multline}
Thus, we obtain the isomorphism
\begin{equation}
 c_{V,\pm}:
  \gbigf^c_{\mp}
 \bigl(
 \iota_{R\ast}\nbigl_{R,\pm,!}[1]
 \bigr)_{|\cnum^{\ast}}
 \simeq
  \gbigf^c_{\mp}
 \bigl(
 \iota_{R!}\nbigl_{R,\pm,\ast}[1]
 \bigr)_{|\cnum^{\ast}}.
\end{equation}

\subsubsection{The regularization}

We continue to use the notation in
\S\ref{subsection;25.3.14.13}.
We have the regular singular meromorphic flat bundle
$V^{\reg}=\nbigstilde^{\infty}_1(V)$
on $(\proj^1,\{0,\infty\})$
with the equality $V^{\reg}=V_{|\cnum}$.

By Lemma \ref{lem;25.3.14.1},
there exist the natural isomorphisms
\[
\DR_{\cnum}(V^{\reg}(\star 0)^{\gbigf_{\pm}})
\simeq
\DR_{\cnum}(V(\star 0)^{\gbigf_{\pm}}).
\]
We obtain the following natural morphisms:
\[
\begin{CD}
 \iota_{R!}\nbigl_{R,\star}[1]
@>{b_{\pm,\star,1}}>>
\DR_{\cnum}\bigl(
 V^{\reg}(\star 0)^{\gbigf_{\pm}}
 \bigr)
@>{b_{\pm,\star,2}}>>
 \iota_{R\ast}\nbigl_{R,\star}[1].
\end{CD}
\]
By Lemma \ref{lem;25.3.14.10},
The kernel and the cokernel 
of the induced morphisms
\[
\begin{CD}
 \gbigf^c_{\mp}\bigl(
 \iota_{R!}\nbigl_{R,\pm,\star}[1]
 \bigr)
@>{\gbigf^{c}_{\mp}(b_{\pm,\star,1})}>>
 \gbigf^{c}_{\mp}\DR_{\cnum}
 V^{\reg}(\star 0)^{\gbigf_{\pm}}
@>{\gbigf^c_{\mp}(b_{\pm,\star,2})}>>
 \gbigf^c_{\mp}\bigl(
 \iota_{R\ast}\nbigl_{R,\pm,\star}[1]
 \bigr)
\end{CD}
\]
are constant sheaves.
By Lemma \ref{lem;25.3.14.15},
The following induced morphisms
\begin{equation}
 \gbigf^{c}_{\mp}(b_{\pm,\ast,1}):
 \gbigf^c_{\mp}
 \bigl(
 \iota_{R!}\nbigl_{R,\pm,\ast}[1]
 \bigr)
 \lrarr
 \gbigf^c_{\mp}
 \bigl(
 \DR_{\cnum}
 V^{\reg}(\ast 0)^{\gbigf_{\pm}}
 \bigr),
\end{equation}
\begin{equation}
 \gbigf^c_{\mp}(b_{\pm,!,2}):
 \gbigf^c_{\mp}
 \bigl(
 \DR_{\cnum}
 V^{\reg}(!0)^{\gbigf_{\pm}}
 \bigr)
 \lrarr
 \gbigf^c_{\mp}
 \bigl(
 \iota_{R\ast}\nbigl_{R,\pm,!}[1]
 \bigr)
\end{equation}
are isomorphisms. 
As in the case of $V$,
we obtain the following isomorphisms
of the local systems on $\cnum^{\ast}$:
\begin{equation}
 c_{V^{\reg},\pm}:
 \gbigf^c_{\mp}
 \bigl(
 \iota_{R\ast}\nbigl_{R,\pm,!}[1]
 \bigr)_{|\cnum^{\ast}}
\simeq
 \gbigf^c_{\mp}
 \bigl(
 \iota_{R!}\nbigl_{R,\pm,\ast}[1]
 \bigr)_{|\cnum^{\ast}}.
\end{equation}
\begin{prop}
\label{prop;25.3.14.100}
$c_{V,\pm}=c_{V^{\reg},\pm}$.
\end{prop}
\pf
Let $f_{\pm,1}$ denote the following morphisms induced by
$a_{\pm,\ast,1}$ and $-b_{\pm,\ast,1}$:
\[
\iota_{R!}\nbigl_{R,\pm,\ast}[1]
\lrarr
\DR_{\cnum}
V(\ast 0)^{\gbigf_{\pm}}
\oplus
\DR_{\cnum}
V^{\reg}(\ast 0)^{\gbigf_{\pm}}.
\]
Let $f_{\pm,2}$ denote the following morphism induced by
$a_{\pm,!,2}$ and $-b_{\pm,!,2}$:
\[
\DR_{\cnum}
V(! 0)^{\gbigf_{\pm}}
\oplus
\DR_{\cnum}
V^{\reg}(! 0)^{\gbigf_{\pm}}
\lrarr
\iota_{R\ast}\nbigl_{R,\pm,!}[1].
\]
We set
$\nbigg_{1,\pm}=C(f_{1,\pm})$ and
$\nbigg_{2,\pm}=C(f_{2,\pm})[-1]$,
where $C(f_{i,\pm})$ denote the mapping cones of $f_{i,\pm}$.

There exist the natural commutative diagrams:
\[
\begin{CD}
\DR_{\cnum}
V(!0)^{\gbigf_{\pm}}
@<<<
\nbigg_{2,\pm}
@>>>
\DR_{\cnum}
V^{\reg}(!0)^{\gbigf_{\pm}}
 \\
 @VVV @VVV @VVV
 \\
 \DR_{\cnum}V(\ast 0)^{\gbigf_{\pm}}
@>>>
\nbigg_{1,\pm}
@<<<
\DR_{\cnum}
V^{\reg}(\ast 0)^{\gbigf_{\pm}}.
\end{CD}
\]
They induce the following commutative diagrams of
isomorphisms:
\[
\begin{CD}
 \gbigf^c_{\mp}
 \DR_{\cnum}
V(!0)^{\gbigf_{\pm}}
@<{\simeq}<<
 \gbigf^c_{\mp}\nbigg_{2,\pm}
@>{\simeq}>>
 \gbigf^c_{\mp}\DR_{\cnum}
V^{\reg}(!0)^{\gbigf_{\pm}}
 \\
 @V{\simeq}VV @V{\simeq}VV @V{\simeq}VV
 \\
\gbigf^c_{\mp}\DR_{\cnum}V(\ast 0)^{\gbigf_{\pm}}
@>{\simeq}>>
 \gbigf^c_{\mp}\nbigg_{1,\pm}
@<{\simeq}<<
 \gbigf^c_{\mp}\DR_{\cnum}
V^{\reg}(\ast 0)^{\gbigf_{\pm}}.
\end{CD}
\]

There exist the following natural morphisms:
\begin{equation}
\label{eq;25.3.14.20}
 \nbigg_{2,\pm}
 \lrarr
 \DR_{\cnum}
 V(!0)^{\gbigf_{\pm}}
 \lrarr
 \iota_{R\ast}\nbigl_{R,\pm,!}[1],
\end{equation}
\begin{equation}
\label{eq;25.3.14.21}
  \nbigg_{2,\pm}
 \lrarr
 \DR_{\cnum}
 V^{\reg}(!0)^{\gbigf_{\pm}}
 \lrarr
 \iota_{R\ast}\nbigl_{R,\pm,!}[1].
\end{equation}
By the construction of $\nbigg_{2,\pm}$,
we can check that the morphisms
$\gbigf^c_{\mp}
\bigl(
\nbigg_{2,\pm}
\bigr)
\to
\gbigf^c_{\mp}\bigl(
\iota_{R\ast}\nbigl_{R,\pm,!}[1]
\bigr)$
induced by (\ref{eq;25.3.14.20})
and (\ref{eq;25.3.14.21})
are equal.
Similarly, the natural morphisms
\begin{equation}
 \iota_{R!}\nbigl_{R,\pm,\ast}[1]
 \lrarr
 \DR_{\cnum}
 V(\ast 0)^{\gbigf_{\pm}}
 \lrarr
  \nbigg_{1,\pm},
\end{equation}
\begin{equation}
 \iota_{R!}\nbigl_{R,\pm,\ast}[1]
  \lrarr
 \DR_{\cnum}
 V^{\reg}(\ast 0)^{\gbigf_{\pm}}
  \lrarr
  \nbigg_{1,\pm}
\end{equation}
induce the same morphisms
$
\gbigf^c_{\mp}\bigl(
\iota_{R!}\nbigl_{R,\pm,\ast}[1]
\bigr)
\to
\gbigf^c_{\mp}
\bigl(
\nbigg_{1,\pm}
\bigr)$.
Then, we obtain $c_{V,\pm}=c_{V^{\reg},\pm}$.
\hfill\qed

\backmatter

\printindex

\end{document}